\documentclass[12pt,a4paper,twoside,openright]{book}

\usepackage[T1]{fontenc}
\usepackage{amsmath,amssymb,amsthm,mathtools}
\usepackage{mathrsfs}
\usepackage{geometry}
\usepackage{booktabs,longtable,array,tabularx,multirow}
\usepackage{enumitem}
\usepackage{xcolor}
\usepackage[expansion=false]{microtype}
\usepackage{graphicx}
\usepackage{float}
\usepackage{tikz}
\usepackage{tikz-cd}
\usepackage{pgfplots}
\usepackage[hidelinks]{hyperref}
\pgfplotsset{compat=1.18}
\usepgfplotslibrary{groupplots}
\usetikzlibrary{arrows.meta,calc,decorations.pathreplacing,positioning}

\definecolor{curveblue}{RGB}{31,103,178}
\definecolor{curveorange}{RGB}{230,126,34}
\definecolor{curvegreen}{RGB}{32,140,92}
\definecolor{curvered}{RGB}{192,57,43}
\definecolor{softgray}{RGB}{110,118,125}

\hypersetup{
  colorlinks=true,
  linkcolor=blue!45!black,
  citecolor=green!35!black,
  urlcolor=blue!60!black,
  pdftitle={A New Characteristic-Uniform Model for Elliptic Curves --
  Theory, Arithmetic, and Applications},
  pdfauthor={Hongfeng Wu}
}
\setlist{nosep}
\allowdisplaybreaks[3]
\makeatletter
\renewcommand*\l@section{\@dottedtocline{1}{1.5em}{3.5em}}
\renewcommand*\l@subsection{\@dottedtocline{2}{5.0em}{4.5em}}
\renewcommand{\cleardoublepage}{%
  \clearpage
  \if@twoside
    \ifodd\c@page
    \else
      \hbox{}%
      \thispagestyle{empty}%
      \newpage
      \if@twocolumn
        \hbox{}%
        \newpage
      \fi
    \fi
  \fi}
\makeatother

\pgfplotsset{
  CdCurve/.style={
    blue!70!black,
    very thick,
    smooth,
    no marks
  },
  CdAsymptote/.style={
    gray!60,
    densely dashed,
    thin
  }
}

\newtheorem{theorem}{Theorem}[chapter]
\newtheorem{proposition}[theorem]{Proposition}
\newtheorem{lemma}[theorem]{Lemma}
\newtheorem{corollary}[theorem]{Corollary}
\theoremstyle{definition}
\newtheorem{definition}[theorem]{Definition}
\newtheorem{example}[theorem]{Example}
\theoremstyle{remark}
\newtheorem{remark}[theorem]{Remark}

\newcommand{\F}{\mathbb F}
\newcommand{\PP}{\mathbb P}
\newcommand{\C}{\mathcal C}
\newcommand{\T}{\mathcal T}
\newcommand{\R}{\mathcal R}
\newcommand{\Q}{\mathcal Q}
\newcommand{\JQ}{\mathcal J}
\newcommand{\wpmap}{\wp}
\newcommand{\Jac}{\operatorname{Jac}}
\newcommand{\Pic}{\operatorname{Pic}}
\newcommand{\Cd}{\mathcal C_d}
\newcommand{\Erho}{\mathcal E_\rho}
\newcommand{\MAB}{\mathcal M_{A,\beta}}
\newcommand{\Tr}{\operatorname{Tr}}
\newcommand{\ord}{\operatorname{ord}}
\newcommand{\Div}{\operatorname{div}}
\newcommand{\charac}{\operatorname{char}}
\newcommand{\M}{\mathbf M}
\newcommand{\Sqr}{\mathbf S}
\newcommand{\Dpar}{\mathbf D}
\newcommand{\Cmul}{\mathbf C}

\newcommand{\Dconst}[1]{\mathbf D_{#1}}
\newcommand{\Dnine}{\mathbf D_{9}}
\newcommand{\Inv}{\mathbf I}
\newcommand{\mBase}{\mathbf m_t}
\newcommand{\mCurve}{\mathbf m_c}
\newcommand{\Cube}{\mathbf C_3}
\newcommand{\Full}{\mathsf F}
\newcommand{\Ffin}{\mathsf F_{\mathrm{fin}}}
\newcommand{\Ktwo}{\mathsf K_2}
\newcommand{\Kfour}{\mathsf K_{4,\pm}}
\newcommand{\hexalpha}[1]{\langle\mathtt{0x#1}\rangle_{\alpha}}
\newcolumntype{L}[1]{>{\raggedright\arraybackslash}p{#1}}
\newcommand{\partoverview}[1]{%
  \thispagestyle{plain}
  \vspace*{0.08\textheight}
  \begin{center}
  {\Large\bfseries Part overview}
  \end{center}
  \vspace{1.5em}
  \begin{quote}\large #1\end{quote}
  \clearpage}

\newcommand{\MonographTitle}{A New Characteristic-Uniform Model for Elliptic Curves}
\newcommand{\MonographSubtitle}{Theory, Arithmetic, and Applications}
\newcommand{\AuthorName}{Hongfeng Wu}
\newcommand{\AuthorEmail}{whfmath@gmail.com}
\newcommand{\AuthorAffiliation}{College of Science, North China University
of Technology, Beijing, China}

\begin{document}

\frontmatter
\hypersetup{pageanchor=false}
\begin{titlepage}
\thispagestyle{empty}
\centering
\vspace*{0.13\textheight}

{\LARGE\bfseries \MonographTitle\par}
\vspace{1.5em}
{\Large \MonographSubtitle\par}

\vfill

{\large \AuthorName\par}
\vspace{0.8em}
{\normalsize E-mail address: \href{mailto:\AuthorEmail}{\texttt{\AuthorEmail}}\par}
\vspace{0.8em}
{\normalsize \AuthorAffiliation\par}

\vfill

{\normalsize \today\par}
\end{titlepage}
\cleardoublepage

\thispagestyle{empty}
\null
\vspace*{0.28\textheight}
\begin{center}
  {\Large\itshape \par}
\end{center}
\vfill
\null
\cleardoublepage

\hypersetup{pageanchor=true}
\setcounter{page}{1}

\chapter*{Preface}
\addcontentsline{toc}{chapter}{Preface}
Explicit elliptic-curve arithmetic is often described as a choice among
normal forms: Weierstrass equations provide universality, Montgomery curves
provide an efficient quotient ladder, and Edwards curves provide symmetric
and often complete addition laws.  The family is useful, but it leaves a
more structural question.  Can a single low-degree model preserve its own
coordinates, marked geometry, quotient, and input/output interface while
still exposing the different simplifications required by odd characteristic,
characteristic three, and characteristic two?

This monograph provides such a model
\[
             \mathcal C_d:\qquad (u^2+u)(v^2+v)=d.
\]
We call this the \(\mathcal C_d\)-curve family and each smooth member a
\(\mathcal C_d\) curve (or a \emph{Cd curve} in plain text).  The interaction
between Artin--Schreier geometry and Edwards arithmetic naturally leads to
this equation.  Its factor \(z^2+z\) becomes a centered even quadratic in
odd characteristic and remains the Artin--Schreier operator in characteristic
two.  More importantly, the smooth \((2,2)\)-completion has the same marked
four-point boundary in every characteristic, inversion fixes the same native
function \(u\), and the boundary selects the same projective Kummer
coordinate \((u+1:u)\).

The equation also carries a literal square symmetry.  Its eight natural
affine transformations form the dihedral group \(D_8\) of order eight.  If
the characteristic is not two, the shift
\((r,s)=(2u+1,2v+1)\) identifies them with all signed permutations of two
coordinates; over \(\mathbb R\) these are exactly the rigid symmetries of a
square centered at \((-1/2,-1/2)\) in the original affine plane.  In every
characteristic the same symmetry has the uniform algebraic invariant
\[
        J=(u^2+u)+(v^2+v),
        \qquad k(\mathcal C_d)^{D_8}=k(J).
\]
This geometric and invariant-theoretic structure is developed intrinsically
before any passage to an Edwards or Weierstrass model.

The isomorphisms with classical models are made explicit throughout.
Edwards, Montgomery, \(Z/4\mathbb Z\)-normal, and Weierstrass equations are used as
proof devices, optimization bridges, and standards for comparison.  Our
convention is that they do not replace the object under study: every formula
is translated back to a full point of
\(\mathcal C_d\), a point on its own smooth completion, or a value on its
declared native quotient.  This distinction between an abstract elliptic
curve and an arithmetic model is central throughout the book.

The reciprocal presentation
\[
              (u+1)(v+1)=du^2v^2
\]
is treated differently from an auxiliary classical model.  It is obtained
from \(\mathcal C_d\) by coordinate inversion and factor exchange on
\(\PP^1\times\PP^1\), so it is an internal affine chart of the same marked
biquadratic model.  In that chart the native Kummer ratio becomes simply
\(v+1\).  In this chart, \(v+1\) gives the corresponding Montgomery or binary-Weierstrass quotient coordinate and allows the known difference to be normalized directly when the external point is already stored in reciprocal coordinates.

The exposition follows one recurring chain:
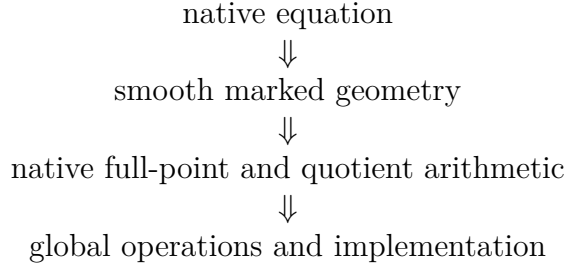
\[
\begin{gathered}
\text{native equation}\\[-1mm]
\Downarrow\\[-1mm]
\text{smooth marked geometry}\\[-1mm]
\Downarrow\\[-1mm]
\text{native full-point and quotient arithmetic}\\[-1mm]
\Downarrow\\[-1mm]
\text{global operations and implementation}
\end{gathered}
\]
Geometry is therefore established before formulas are optimized, and the
domain of a formula is established before its operation count is compared.
This order is especially important for completeness: no finite affine chart
can contain \(P+(-P)=O\), whereas the smooth native completion does support a
complete addition-law atlas and, on the appropriate finite-field subfamily, a
single complete formula.

The monograph develops native full-point, differential, and Kummer arithmetic for \(\mathcal C_d\), 
together with complete addition atlases, scalar multiplication, recovery, halving, tripling, and \(2P+Q\). 
It also develops division polynomials, isogenies, trace moments, CM endomorphisms, and Tate and Weil pairings. 
The later extension chapters treat \(\mathcal T_{a,d}\), reciprocal \(\mathcal C\)-curves, 
and the symmetric QRT envelope.
For isogenies, the native odd-degree \(u\)-Kummer quotient is integrated
into the general theory and repeated in the cryptographic chapter.  In the
ordinary binary branch, the one-sided and full product extensions support a
twist-stable separable two-isogeny
\(W_{A,d}\to W_{A,\sqrt d}\); its native \(C_d\) and
\(\mathcal T_{a,d}\) realizations require only one multiplication or one
squaring and are compared with binary Edwards using equal endpoints.
For the QRT envelope, the adjacent-state analysis includes an exact descent
criterion, a ground-field group-law statement with origin alignment, a
one-formula complete Montgomery-conjugate state update, ordinary-binary
completeness, full-point recovery, and coordinatewise isogeny functoriality.
In particular, the Binary pointed state model theorem gives every ordinary
pointed elliptic curve over a perfect binary field its ground-field state
equation, oriented shift, complete doubling update, and Artin--Schreier
twist data in one statement.

Every operation count is attached to stated coordinates, inputs, outputs,
and completeness hypotheses.  Symbolic derivations establish rational-map
identities, while finite-field examples test affine points, boundary points,
exceptional charts, and ordinary as well as supersingular behavior.  These
counts are algebraic cost models rather than claims about a particular
processor or security parameter.  The principal contribution is a characteristic-uniform arithmetic 
system in which local formulas, quotient arithmetic, and global operations are realized on the same marked biquadratic model.

\chapter*{Synopsis}
\addcontentsline{toc}{chapter}{Synopsis}
An explicit elliptic-curve model is determined not only by its abstract
isomorphism class, but also by its embedding, marked identity, boundary,
coordinate functions, exceptional divisors, and arithmetic circuits.  This
monograph develops all of these structures for the characteristic-uniform
family
\[
   \Cd:\qquad (u^2+u)(v^2+v)=d.
\]
The monograph proves that the complete group law, Kummer quotient, scalar multiplication, division theory, isogenies, 
and pairings of \(\mathcal C_d\) can all be expressed within the native model or its explicitly declared quotients.
 
\medskip
\noindent\textbf{Geometry and the native quotient.}
The natural completion in \(\PP^1\times\PP^1\) is a smooth genus-one curve
exactly when \(d(1-16d)\ne0\).  With \(O=(0,\infty)\) as identity, one has in
every characteristic
\[
  -(u,v)=(u,-v-1),\qquad
  \kappa(u,v)=(u+1:u).
\]
Moreover, the four boundary points form a rational cyclic subgroup of order four, and
\(k(u)\) is the fixed field of inversion.  The two factor involutions and
coordinate exchange generate the full eight-element square group
\(D_8\).  When \(2\ne0\), its true affine center is
\((-1/2,-1/2)\), and the centered coordinates
\((r,s)=(2u+1,2v+1)\) turn the action into the signed permutation group
\(W(B_2)\).  Algebraically, in every characteristic,
\[
 k(\mathcal C_d)^{D_8}
   =k\bigl((u^2+u)+(v^2+v)\bigr).
\]
The four boundary points are the vertices of the corresponding square
permutation representation.  Thus \(\kappa\) is not an imported Montgomery
coordinate: it is the boundary-normalized native Kummer map, distinct from
the further quotient by the full square group.

\medskip
\noindent\textbf{Odd-characteristic arithmetic and completeness.}
In odd characteristic, put \(\rho=1-16d\).  Centered reciprocal coordinates
\[
 (\xi,\eta)=\bigl((2v+1)^{-1},(2u+1)^{-1}\bigr)
\]
give the Edwards equation
\(\xi^2+\eta^2=1+\rho\xi^2\eta^2\), while
\[
 U=\frac{u+1}{u},\qquad
 V=2(2v+1)U
\]
give the Montgomery equation
\[
 \frac1{16d}V^2
   =U^3+\left(\frac1{4d}-2\right)U^2+U.
\]
Addition and doubling are first derived in \(u,v\); tripling and closed
\(2P+Q\) are likewise stated as maps \(\mathcal C_d\to\mathcal C_d\).
The explicit Edwards and Montgomery coordinate correspondences are then
used to factor and optimize native shifted-projective full-point formulas,
Kummer \(XZ\) formulas, multiplication--squaring tradeoffs, halving,
coordinate recovery, and complete scalar multiplication. 
The ordinary two-to-one quotient \(U\) is kept distinct from the
degree-eight post-isogeny Kummer coordinate
\[
 \omega=\frac{\rho}{(2u+1)^2(2v+1)^2}.
\]
The latter can support a complete differential law on the complete Edwards
subfamily, but it retains the eight-point ambiguity
\(\{\pm P\}+\langle T_4\rangle\), where
\(T_4=(1,0)\) is the Edwards image of the marked native point \(R\).

Completeness is proved intrinsically on \(\mathcal C_d\).  A finite affine
\((u,v)\)-formula cannot represent every sum because
\(P+(-P)=O=(0,\infty)\).  On the smooth native completion, the centered
Segre coordinates
\[
(X:Y:T:Z)=(RV_0:SU_0:U_0V_0:RS)
\]
carry the native law
\[
                  P_1+P_2=(EF:GH:EH:FG).
\]
It costs \(9\M+\Dpar\), its mixed form costs \(8\M+\Dpar\), and its
specialized double costs \(4\M+4\Sqr\).  Over an odd finite field this
single formula is \(\F_q\)-complete exactly when \(\rho\) is a nonsquare.
For all other smooth parameters, and also in characteristic two, finitely
many native formulas form a base-point-free complete addition-law atlas.
The native differential formulas form complete atlases after the
exceptional difference graphs are included.

On the native Kummer input \((X:Z)=(u+1:u)\), the identity \(X-Z=1\)
reduces the first double from \(2\M+2\Sqr+\Dpar\) to
\(\M+\Sqr+\Dpar\).  This is an initialization saving, not a recurring
per-bit count.  A dedicated odd-characteristic full-point double costs
\(3\M+4\Sqr+\Dpar\), and native shifted-projective mixed addition costs
\(8\M+\Sqr+\Dpar\).  The operation tables distinguish these specialized
circuits from the complete Segre law and from formulas on auxiliary models.

In characteristic different from two and three, two CM subfamilies support
additional low-cost endomorphisms directly in the
\(\mathcal C_d\)-coordinates.  At \(d=1/8\), the \(j=1728\) map
\[
 (u,v)\longmapsto
 \left(-\frac{u}{2u+1},
       \frac{-1-i(2v+1)}2\right)
\]
satisfies \(\phi^2=[-1]\) and acts on the native Kummer line as
\((X:Z)\mapsto(-X:Z)\).  When
\(A=(4d)^{-1}-2\) satisfies \(A^2=3\), the \(j=0\) map induces
\[
       (X:Z)\longmapsto(\zeta X+cZ:Z),
       \qquad c=A(\zeta-1)/3,
\]
and satisfies \(\psi^2+\psi+[1]=[0]\).  These maps give the eigenvalue
equations used for two-dimensional GLV scalar decomposition.

\medskip
\noindent\textbf{The reciprocal chart.}
The ambient involution
\[
 (u_0,v_0)\longmapsto
 \left(\frac1{v_0},\frac1{u_0}\right)
\]
turns \(\mathcal C_d\) into
\((u+1)(v+1)=du^2v^2\) without changing its smooth completion, parameter,
marked torsion, or Kummer quotient.  In the reciprocal coordinates the
quotient is \((v+1:1)\).  Thus an already reciprocal input gives the
odd-characteristic affine-difference ladder cost
\(5\M+4\Sqr+\Dpar\) without preprocessing inversion and gives the two
normalized binary ladder costs directly.  A strict interface comparison
also records the converse: the original native representative retains the
better odd-characteristic first double, while converting a raw native point
to the reciprocal quotient costs the same inversion as ordinary Kummer
normalization.  The reciprocal equation therefore serves as a quotient-friendly internal
chart of the same marked \(\mathcal C_d\) model.  It supplies a distinct
stored-input interface in which the known difference can be normalized
directly, while preserving the same smooth completion, parameter, marked
torsion, and native Kummer quotient.

\medskip
\noindent\textbf{Characteristic three.}
Characteristic three receives a separate treatment.  Native Kummer
differential addition has the Montgomery-shaped circuit
\(4\M+2\Sqr\), or
\(3\M+2\Sqr\) for an affine known difference.  The native Kummer tripling
map has
cost \(2\M+2\Sqr+2\Dpar+2\Cube\) and factors as relative Frobenius followed
by Verschiebung.  The first factor is purely inseparable; the second is
separable on an ordinary curve and purely inseparable on a supersingular
curve.  Neither factor is to be confused with a V\'elu quotient by a
chosen reduced cyclic subgroup.  This tripling formula is not obtained by
reducing a classical Hessian normal form: in characteristic three the
diagonal Hessian cubic degenerates because
\(X^3+Y^3+Z^3=(X+Y+Z)^3\).

\medskip
\noindent\textbf{Characteristic two.}
In characteristic two the same equation is an Artin--Schreier double cover.
The native affine double is
\[
 u_{2P}=\frac{(u^2+u)^2}{(u^2+u)^2+d},\qquad
 v_{2P}=u^2+v^2,
\]
and a direct native addition law is obtained before introducing any
auxiliary model.  The linear substitutions
\(x=u/(u+1)\), \(y=v/(v+1)\) then give the
\(Z/4\mathbb Z\)-normal form
\[
        (1+x)^2(1+y)^2=d^{-1}xy,
\]
and
\[
       X=d\frac{u+1}{u},\qquad Y=Xv
\]
give \(Y^2+XY=X^3+d^2X\).  This yields native binary Kummer ladder steps
costing
\[
 4\M+5\Sqr+\mBase+\mCurve
 \quad\text{or}\quad
 4\M+4\Sqr+\mBase+2\mCurve,
\]
after the fixed difference has been normalized once.  A raw native
projective difference instead gives
\(6\M+5\Sqr+\mCurve\) per step and requires no input inversion.
These recurring costs exclude the explicitly stated one-time
normalization.  The same theory also gives point recovery, tripling,
complete differential charts,
model-preserving odd-degree isogenies, and native Miller functions.  The
binary complete Segre chart is
\[
       AD=BC,\qquad D(A+B+C+D)=dA^2,
\]
so completeness in characteristic two is again a property of the smooth
\(\mathcal C_d\)-model, not of an auxiliary equation.

\medskip
\noindent\textbf{Global constructions.}
The global theory determines geometric and \(\F_q\)-isomorphism classes,
exact first and second Frobenius-trace moments with their parameter weights,
and division-fiber mean values.
Classical Weierstrass and Edwards division functions are pulled back to the
native variables through order seven.  Independently, paired native Kummer
recursions construct degree-\(n^2\) division pairs directly from
\(\mathcal C_d\) doubling and differential addition; the same addition-chain
circuit performs scalar multiplication, torsion testing, and division
evaluation.  For integers \(n\) prime to the characteristic, the
characteristic-uniform trace identity
\[
 \frac1{n^2}\sum_{[n]P=Q}\frac{u(P)+1}{u(P)}
   =\frac{u(Q)+1}{u(Q)}
\]
is proved from generalized division polynomials.  Odd characteristic also
has reciprocal-coordinate means, with characteristic three included whenever
\(3\nmid n\).

Model-preserving V\'elu-type isogenies are constructed in the odd and binary
branches, together with kernel-polynomial evaluation, Galois descent,
relative Frobenius and Verschiebung, composition, and dual-isogeny formulas.
Sparse native Miller numerators and vertical denominators depending only on
the native Kummer coordinate lead to Tate and Weil pairing wrappers written
in \(u,v\).  A precise subfield criterion permits denominator elimination,
and product pairings share the final exponentiation.  Closed \(2P+Q\)
circuits are useful in selected affine and fixed-point settings, and the two
CM subfamilies support two-dimensional scalar decomposition.  Canonical point
encodings use a square-and-sign criterion in odd characteristic and a
trace criterion supplemented by a linear functional in
characteristic two; four reserved tags represent the native boundary, and a
masked decode-and-validate schedule keeps exceptional inputs out of
secret-dependent control flow.

\medskip
\noindent\textbf{The 25519 specialization.}
A concrete 25519 specialization tests the theory against an established
cryptographic parameter.  Over \(p=2^{255}-19\), put
\[
 d_{\rm C}=(16\cdot121666)^{-1}.
\]
Then
\[
 \mathcal C_{25519}:\quad
 (u^2+u)(v^2+v)=d_{\rm C}
\]
is isomorphic to Edwards25519 by
\[
 x=-\frac{\imath}{2v+1},\qquad
 y=\frac1{2u+1},\qquad \imath^2=-1.
\]
Its native Kummer ratio
\[
 U=\frac{u+1}{u}=\frac{1+y}{1-y}
\]
is exactly the Curve25519/X25519 Montgomery coordinate.  A diagonal
\(\imath\)-twist of the native Segre embedding gives a single complete
Cd25519 law at \(8\M+\Dpar\), mixed addition at
\(7\M+\Dpar\), and doubling at \(4\M+4\Sqr\).  On the quotient line, the
native first double costs \(\M+\Sqr+\Dpar\); a fixed difference \(U=9\)
reduces the base-point ladder step to
\(4\M+4\Sqr+\Dpar+\Dnine\).  The analysis also separates main-curve
Cd25519 points from the curve-or-twist input semantics required by X25519.

\medskip
\noindent\textbf{Comparison, examples, and checks.}
The model-by-model comparison separates full-point, Kummer, mixed-addition,
tripling, and binary costs, and records the precise coordinate and output
type in every row.  It therefore uses the general formulas and the Cd25519
specialization already established, rather than anticipating either one.
Explicit examples over
\(\F_9\), \(\F_{27}\), \(\F_{97}\), \(\F_{101}\),
\(\F_{101^2}\), \(\F_{1009}\), and \(\F_{2^8}\), together with symbolic
and finite-field checks, accompany the proofs.  The Cd25519 chapter also
checks the standard \(p=2^{255}-19\) base point and an RFC~7748 X25519
vector.  These examples verify not only
identities on generic charts but also boundary points,
exceptional-difference charts, completeness claims, scalar chains, the
ordinary/supersingular degree-three distinction, and the stated arithmetic
output conventions.

\medskip
\noindent\textbf{The symmetric three-parameter extension.}
The final extension Part begins with
\[
       \mathcal C_{a,b,d}:\quad
       (u^2+u+a)(v^2+v+b)=d.
\]
In odd characteristic, with
\(\lambda=1-4a\), \(\mu=1-4b\),
\(\delta=16d\), and \(\rho=\lambda\mu-\delta\), centered reciprocal
coordinates give
\[
             \lambda x^2+\mu y^2=1+\rho x^2y^2,
\]
and the base-field Jacobian is
\[
       Y^2=X\bigl(X^2+2(2\lambda\mu-\delta)X+\delta^2\bigr).
\]
The subfamily \(b=0\) is a twisted Edwards curve over the ground field;
\(a=b\) has the additional coordinate-exchange symmetry; and \(d=ab\)
has four explicit rational affine points and a base-field Huff--Weierstrass
map.  In characteristic two the Jacobian is
\[
       Y^2+XY=X^3+(a+b)X^2+d^2X,
\]
so \(a+b\) is the Artin--Schreier twist parameter while \(d\) determines
the geometric \(j\)-invariant.  This chapter identifies the rigidity and the rationally marked boundary of
the specialization \(a=b=0\).

\medskip
\noindent\textbf{Twisted, reciprocal, and QRT extensions.}
The remaining chapters of the same Part develop the one-sided, reciprocal,
and QRT forms.  The one-sided subfamily
\[
 \mathcal T_{a,d}:\qquad (u^2+u+a)(v^2+v)=d
\]
represents the Jacobian classes of the full product family; over finite
fields it represents the curves themselves after an origin is chosen.  In
odd characteristic it is exactly a twisted-Edwards/Montgomery normal form,
with \(a=1/2\) producing the native coefficient \(A=-1\); in characteristic
two it gives
\(Y^2+XY=X^3+aX^2+d^2X\) and represents every ordinary elliptic curve over
a perfect field.

The reciprocal family
\[
 \mathcal R_{\tau,\sigma,\kappa}:\qquad
 (x^2-\tau)(y^2-\sigma)=\kappa xy
\]
connects reciprocal quotient involutions, Jacobi quartics, symmetric QRT
curves, and split or nonsplit four-torsion normal forms.  In characteristic
two it is M\"obius-equivalent to
\(\mathcal C_{\sqrt{\tau\sigma}/\kappa}\).  Its complete projective
interface is obtained from explicit Weierstrass/Jacobi atlases and Kummer
differential arithmetic.

Finally, the symmetric envelope
\[
 \mathcal Q_{\alpha,\beta,\gamma}:\qquad
 x^2y^2+\alpha(x^2+y^2)+\beta xy+\gamma=0
\]
is treated as an Euler--Chasles symmetric biquadratic and as a QRT invariant
fiber.  Its Vieta involutions generate the McMillan map
\(T(x,y)=(y,-x-\beta y/(y^2+\alpha))\), which is independent of
\(\gamma\) and acts as a fixed elliptic translation on every smooth fiber.
After a pointed elliptic curve \((E,D)\) is fixed, the same geometry becomes
the adjacent Kummer state
\((\kappa(P),\kappa(P+D))\); state doubling and differential addition then
give a logarithmic binary ladder rather than linear repeated QRT iteration.
The same state calculus yields an elliptic Lucas system with nonlinear
addition and exact fast-index doubling, QRT state-division polynomials, and
proved bridges to Ward elliptic divisibility sequences, sigma functions, and
elliptic nets.  The exact-even odd-characteristic state update is
PGL$_2$-conjugate to compiled Montgomery \texttt{xDBLADD} and costs
\(4\M+4\Sqr+2\Dpar\); each split odd branch and each ordinary-binary
branch is represented by one complete projective formula.  The exact
field-of-definition obstruction, three-chart recovery atlas, fixed-state
MSM reachability theorem, and coordinatewise isogeny law are stated separately so
that these cryptographic conclusions are not conflated with the existence
of the QRT recurrence alone.  Over a perfect binary field, the Binary
pointed state model theorem applies to every ordinary pointed elliptic
curve and records the state equation, the original marked translation, the
complete update, and the residual Artin--Schreier twist simultaneously.

\medskip
\noindent\textbf{Keywords:}
\(\mathcal C_d\) curves; elliptic curves; elliptic-curve cryptography; characteristic-uniform models;
dihedral symmetry; invariant theory; Artin--Schreier covers; Edwards curves; Montgomery
arithmetic; reciprocal chart; characteristic two; characteristic three; Kummer lines;
differential addition; division polynomials; mean values; isogenies;
V\'elu formulas; endomorphisms; complex multiplication; GLV decomposition;
Tate pairing; Weil pairing; symmetric biquadratic curves;
Euler--Chasles correspondence; Jacobi quartics; reciprocal \(C\)-curves;
QRT maps; McMillan iteration; adjacent Kummer states; Cd25519; Curve25519;
Edwards25519; X25519.

\tableofcontents

\mainmatter

\part{Motivation, Context, and Native Geometry}
\partoverview{This part explains why the factorized model is worth studying,
fixes the distinction between an abstract curve and an arithmetic model, and
constructs the geometric objects on which every later algorithm depends.
The progression is deliberate: research problem and literature, cost and
output conventions, the smooth marked \((2,2)\)-completion, its square
symmetry and invariant algebra, the native Kummer coordinate, and only then
the odd-characteristic dictionaries.}

\chapter{Introduction}
\label{ch:introduction}
\section{The arithmetic-model problem and the central claim}

The subject of this monograph is the one-parameter biquadratic family
\begin{equation}\label{eq:model}
 \Cd:\qquad (u^2+u)(v^2+v)=d.
\end{equation}
The starting problem is broader than finding another equation for an
elliptic curve.  Fast formulas normally arise because a model makes one
piece of the group structure visible: rational two-torsion gives a
Montgomery quotient, a symmetric embedding gives Edwards addition, and
special level structures give Hessian or theta coordinates.  These
advantages are usually studied in separate normal forms and, in low
characteristic, sometimes in entirely separate theories.  What is missing is
a model in which the equation, the marked torsion, the native quotient, and
the implementation interface survive a change of characteristic.

Equation~\eqref{eq:model} is a natural candidate.  It is symmetric in \(u\)
and \(v\), is invariant under
\[
 u\longmapsto-u-1,\qquad v\longmapsto-v-1,
\]
and makes sense without alteration in every characteristic.  Together with
coordinate exchange, these involutions produce eight transformations.  If
\(2\ne0\), the shift \((r,s)=(2u+1,2v+1)\) turns them into
\[
 (r,s),(s,r),(-r,s),(s,-r),(r,-s),(-s,r),(-r,-s),(-s,-r),
\]
the complete rigid-motion group of a square.  Thus the geometric center in
the original coordinates is \((-1/2,-1/2)\), not the origin.  The center
belongs only to the degenerate member \(d=1/16\).  In characteristic two
the Euclidean center disappears, but the same eight algebraic maps and the
same invariant-field quotient remain.  The factor
\(z^2+z\) is the key.  Away from characteristic two it becomes an even
quadratic after the linear shift \(z\mapsto z+1/2\); in characteristic two
it becomes the additive Artin--Schreier polynomial.  Thus the same factorization
predicts Edwards-type bilinear arithmetic in one branch and
Artin--Schreier/\(Z/4\mathbb Z\) arithmetic in the other.  In characteristic
three, its interaction with Frobenius cubing produces a third specialization.
The interaction between Artin--Schreier geometry and Edwards arithmetic
naturally leads to the characteristic-uniform biquadratic family
\(\mathcal C_d\): its fundamental factor \(z^2+z\) becomes a centered even
quadratic in odd characteristic and retains its intrinsic
Artin--Schreier form in characteristic two.  This structural synthesis
motivates the designation \(\mathcal C_d\), used throughout this monograph
for the family and its smooth members.

The central claim is therefore:
\begin{quote}
\emph{The family \(\mathcal C_d\) is a characteristic-uniform explicit
elliptic-curve model whose native geometry supports full-point arithmetic,
Kummer arithmetic, division theory, isogenies, and pairings without making an
auxiliary isomorphic model the declared input or output interface.}
\end{quote}
This statement must be proved at four levels: smooth geometry, intrinsic
quotients, explicit formulas with specified domains and costs, and comparisons
that separate structural advantages from operation-count advantages.

It is essential to distinguish four notions:
\begin{enumerate}[label=\textup{(\roman*)}]
 \item an abstract elliptic curve over \(k\);
 \item a projective model and its marked identity and torsion data;
 \item a coordinate chart on that model;
 \item an arithmetic circuit acting on the chosen coordinates.
\end{enumerate}
Every elliptic curve with a rational point admits a Weierstrass equation.
Edwards, Montgomery, Hessian, and Jacobi quartic models nevertheless retain
their own embeddings, coordinates, marked data, and arithmetic interfaces.
The same principle applies to \eqref{eq:model}.
Throughout this monograph, \(\mathcal C_d\) is therefore treated as a
distinct explicit elliptic-curve model: it
has its own defining equation, smooth completion, marked identity, boundary
divisor, native coordinates, and addition-law spaces.  An isomorphism to an
Edwards or Weierstrass equation is a computational bridge between models;
it is not an identification of their equations, embeddings, coordinate
functions, or arithmetic interfaces.

\paragraph{Model and coordinate dictionaries.}
Throughout this monograph, a \emph{dictionary} between two curve models
means an explicit interface consisting of parameter correspondences,
forward and inverse coordinate maps at the declared information level,
the images of the marked identity, boundary, and torsion points, and the
induced translation of full-point or quotient arithmetic.  The domain of
each coordinate formula, its quotient degree, its recovery data, and its
extension to the smooth completions are regarded as part of the dictionary.
The surrounding \emph{framework} organizes the geometry, models, and
algorithms, whereas a dictionary is one of the concrete bridges inside
that framework.

A \emph{Montgomery dictionary} records the parameter specialization to a
Montgomery equation and identifies the corresponding Kummer coordinate,
marked points, ladder constants, differential operations, and recovery
maps.  A \emph{Möbius dictionary} is the special case in which a
\(\mathbb P^1\)-coordinate is changed by a projective linear
transformation
\[
   (X:Z)\longmapsto(aX+bZ:cX+dZ),
   \qquad ad-bc\ne0;
\]
on a model in \(\mathbb P^1\times\mathbb P^1\), such a transformation may
be applied in either factor.  The induced arithmetic maps are transported
by the corresponding conjugation.

\begin{definition}[\(\mathcal C_d\) curves]\label{def:Cd-curve-name}
We call the family
\[
             \mathcal C_d:\qquad (u^2+u)(v^2+v)=d
\]
the \(\mathcal C_d\)-curve family and its members \(\mathcal C_d\)
curves.  A parameter \(d\) satisfying the smoothness condition
\(d(1-16d)\ne0\) determines a smooth \(\mathcal C_d\) curve.
In plain-text titles and searchable metadata, the notation is written
\emph{Cd curve}.
\end{definition}

Writing \(\Phi(z)=z^2+z\), the equation takes the compact form
\[
             \mathcal C_d:\qquad \Phi(u)\Phi(v)=d.
\]

\begin{definition}[Arithmetic closure]
An explicit elliptic-curve model is \emph{arithmetically closed} if its
principal arithmetic objects and operations---complete group law, Kummer
quotient, scalar multiplication, isogenies, and pairings---can all be
expressed and executed intrinsically within the same model, without requiring
an auxiliary curve model as the computational interface.
\end{definition}

The results proved in this monograph establish that the
\(\mathcal C_d\)-curve model is arithmetically closed in this sense.

Figures 1.1--1.4 illustrate the intrinsic symmetry and real loci of the curve
\[
 \mathcal C_d:\qquad (u^2+u)(v^2+v)=d.
\]

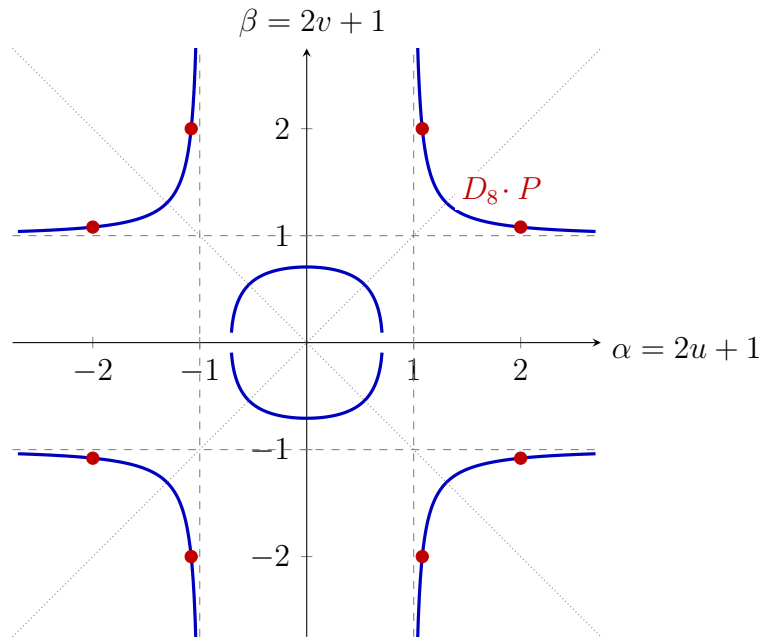
\begin{figure}[htbp]
\centering
\begin{tikzpicture}
\begin{axis}[
  width=0.78\textwidth,
  height=0.62\textwidth,
  axis equal image,
  xmin=-2.75,xmax=2.75,
  ymin=-2.75,ymax=2.75,
  axis lines=middle,
  xlabel={$\alpha=2u+1$},
  ylabel={$\beta=2v+1$},
  xlabel style={
    at={(axis description cs:1,0.49)},
    anchor=west
  },
  ylabel style={
    at={(axis description cs:0.51,1)},
    anchor=south
  },
  xtick={-2,-1,0,1,2},
  ytick={-2,-1,0,1,2},
  samples=180,
  clip=true
]
  \foreach \xstart/\xend in
    {-2.7/-1.015,-0.704/0.704,1.015/2.7}
  {
    \addplot[blue!75!black,very thick,domain=\xstart:\xend]
      {sqrt((x^2-0.5)/(x^2-1))};
    \addplot[blue!75!black,very thick,domain=\xstart:\xend]
      {-sqrt((x^2-0.5)/(x^2-1))};
  }

  \addplot[gray,dashed]
    coordinates {(-1,-2.75) (-1,2.75)};
  \addplot[gray,dashed]
    coordinates {(1,-2.75) (1,2.75)};
  \addplot[gray,dashed]
    coordinates {(-2.75,-1) (2.75,-1)};
  \addplot[gray,dashed]
    coordinates {(-2.75,1) (2.75,1)};

  \addplot[gray,densely dotted]
    coordinates {(-2.75,-2.75) (2.75,2.75)};
  \addplot[gray,densely dotted]
    coordinates {(-2.75,2.75) (2.75,-2.75)};

  \addplot[
    only marks,
    mark=*,
    mark size=2.3pt,
    red!75!black
  ]
  coordinates {
    (2,1.080123)
    (2,-1.080123)
    (-2,1.080123)
    (-2,-1.080123)
    (1.080123,2)
    (-1.080123,2)
    (1.080123,-2)
    (-1.080123,-2)
  };

  \node[
    fill=white,
    inner sep=2pt,
    text=red!75!black
  ]
  at (axis cs:1.82,1.42)
  {$D_8\!\cdot P$};
\end{axis}
\end{tikzpicture}
\caption{The real locus of \(\mathcal C_{1/32}\) in the centered coordinates
\((\alpha,\beta)\). The coordinate axes and the two diagonals are
reflection axes for the \(D_8\)-action. The marked points constitute
a generic eight-element orbit, while the dashed lines
\(\alpha=\pm1\) and \(\beta=\pm1\) are asymptotes.}
\label{fig:cd-d8-symmetry}
\end{figure}

\begin{figure}[htbp]
\centering
\begin{tikzpicture}
\begin{groupplot}[
  group style={
    group size=3 by 1,
    horizontal sep=0.65cm
  },
  width=0.31\textwidth,
  height=0.31\textwidth,
  axis equal image,
  xmin=-2.55,xmax=2.55,
  ymin=-2.55,ymax=2.55,
  axis lines=middle,
  xtick={-1,0,1},
  ytick={-1,0,1},
  samples=150,
  clip=true,
  title style={font=\small},
  tick label style={font=\scriptsize}
]

\nextgroupplot[
  title={$d=1/32$ \((16d=1/2)\)}
]
  \foreach \xstart/\xend in
    {-2.5/-1.018,-0.704/0.704,1.018/2.5}
  {
    \addplot[blue!75!black,thick,domain=\xstart:\xend]
      {sqrt((x^2-0.5)/(x^2-1))};
    \addplot[blue!75!black,thick,domain=\xstart:\xend]
      {-sqrt((x^2-0.5)/(x^2-1))};
  }
  \node[
    font=\scriptsize,
    fill=white,
    inner sep=1pt
  ]
  at (axis cs:0,0)
  {central oval};

\nextgroupplot[
  title={$d=1/8$ \((16d=2)\)}
]
  \foreach \xstart/\xend in
    {-2.5/-1.018,1.018/2.5}
  {
    \addplot[orange!85!black,thick,domain=\xstart:\xend]
      {sqrt((x^2+1)/(x^2-1))};
    \addplot[orange!85!black,thick,domain=\xstart:\xend]
      {-sqrt((x^2+1)/(x^2-1))};
  }
  \node[
    font=\scriptsize,
    fill=white,
    inner sep=1pt
  ]
  at (axis cs:0,0)
  {no compact oval};

\nextgroupplot[
  title={$d=-1/16$ \((16d=-1)\)}
]
  \foreach \xstart/\xend in
    {-2.5/-1.418,-0.982/0.982,1.418/2.5}
  {
    \addplot[green!60!black,thick,domain=\xstart:\xend]
      {sqrt((x^2-2)/(x^2-1))};
    \addplot[green!60!black,thick,domain=\xstart:\xend]
      {-sqrt((x^2-2)/(x^2-1))};
  }
  \node[
    font=\scriptsize,
    fill=white,
    inner sep=1pt
  ]
  at (axis cs:0,0)
  {alternating regions};

\end{groupplot}
\end{tikzpicture}
\caption{Representative nonsingular real loci of \(\mathcal C_d\).
From left to right: \(0<d<1/16\), \(d>1/16\), and \(d<0\).
For \(0<d<1/16\), a compact central oval coexists with four
unbounded branches. For \(d>1/16\), the central oval disappears.
For \(d<0\), the branches occupy alternating horizontal and
vertical regions.}
\label{fig:cd-real-loci}
\end{figure}
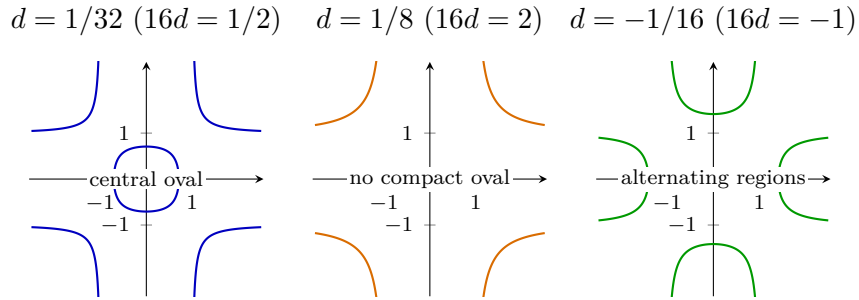

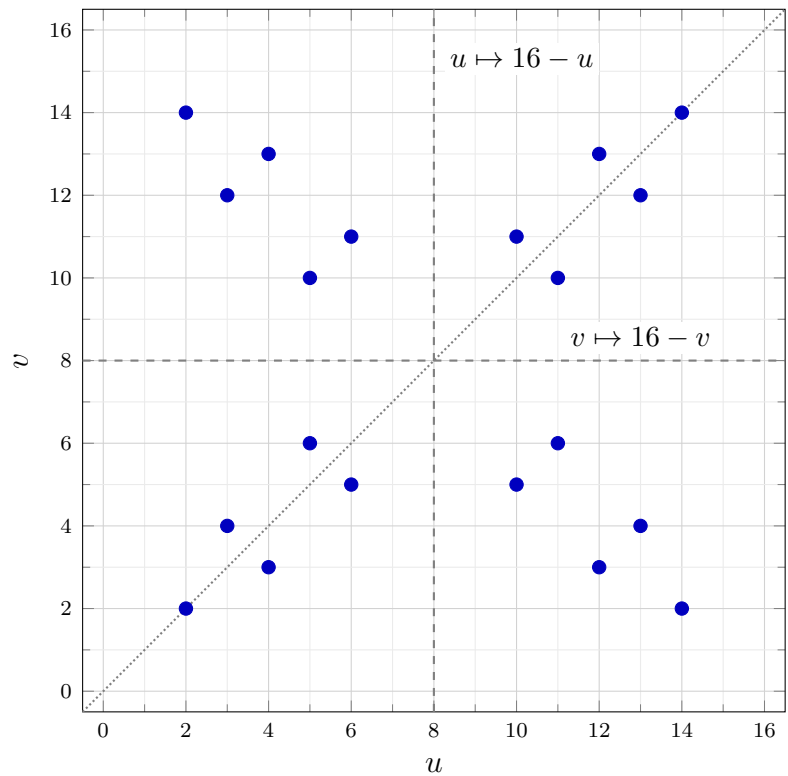
\begin{figure}[htbp]
\centering
\begin{tikzpicture}
\begin{axis}[
  width=0.72\textwidth,
  height=0.72\textwidth,
  axis equal image,
  xmin=-0.5,xmax=16.5,
  ymin=-0.5,ymax=16.5,
  xlabel={$u$},
  ylabel={$v$},
  xtick={0,2,...,16},
  ytick={0,2,...,16},
  minor x tick num=1,
  minor y tick num=1,
  grid=both,
  major grid style={draw=gray!35},
  minor grid style={draw=gray!16},
  tick label style={font=\scriptsize}
]
  \addplot[gray,dashed,thick]
    coordinates {(8,-0.5) (8,16.5)};
  \addplot[gray,dashed,thick]
    coordinates {(-0.5,8) (16.5,8)};
  \addplot[gray,densely dotted,thick]
    coordinates {(-0.5,-0.5) (16.5,16.5)};

  \addplot[
    only marks,
    mark=*,
    mark size=2.5pt,
    blue!75!black
  ]
  coordinates {
    (2,2)
    (2,14)
    (3,4)
    (3,12)
    (4,3)
    (4,13)
    (5,6)
    (5,10)
    (6,5)
    (6,11)
    (10,5)
    (10,11)
    (11,6)
    (11,10)
    (12,3)
    (12,13)
    (13,4)
    (13,12)
    (14,2)
    (14,14)
  };

  \node[
    anchor=west,
    fill=white,
    inner sep=2pt,
    font=\small
  ]
  at (axis cs:8.25,15.3)
  {$u\mapsto16-u$};

  \node[
    anchor=south,
    fill=white,
    inner sep=2pt,
    font=\small
  ]
  at (axis cs:13.0,8.15)
  {$v\mapsto16-v$};
\end{axis}
\end{tikzpicture}
\caption{The complete affine point set of
\(\mathcal C_2:(u^2+u)(v^2+v)=2\) over \(\mathbb F_{17}\).
The twenty blue points are the exact
finite-field solutions. The dashed lines and the diagonal display
the finite-field realization of the same \(D_8\)-symmetry.}
\label{fig:cd-f17}
\end{figure}

\begin{figure}[htbp]
\centering
\begin{tikzpicture}[>=Latex,font=\small]
  \draw[very thick,blue!75!black]
    (-5,1.25)--(1.0,1.25);
  \node[left] at (-5,1.25) {sheet \(+\)};

  \draw[very thick,orange!85!black]
    (-5,-1.25)--(1.0,-1.25);
  \node[left] at (-5,-1.25) {sheet \(-\)};

  \foreach \x/\branchlabel in {
    -4.3/$-1$,
    -3.2/$-\sqrt{1-16d}$,
    -0.8/$\sqrt{1-16d}$,
     0.3/$1$
  }
  {
    \fill[red!75!black]
      (\x,1.25) circle (2.1pt);
    \fill[red!75!black]
      (\x,-1.25) circle (2.1pt);
    \draw[gray,densely dashed]
      (\x,-1.25)--(\x,1.25);
    \node[above=3pt]
      at (\x,1.25)
      {\branchlabel};
  }

  \draw[
    line width=3pt,
    green!55!black
  ]
    (-4.3,1.25)--(-3.2,1.25);
  \draw[
    line width=3pt,
    green!55!black
  ]
    (-0.8,1.25)--(0.3,1.25);
  \draw[
    line width=3pt,
    green!55!black
  ]
    (-4.3,-1.25)--(-3.2,-1.25);
  \draw[
    line width=3pt,
    green!55!black
  ]
    (-0.8,-1.25)--(0.3,-1.25);

  \node[
    align=center,
    text=gray!80!black
  ]
  at (-2.0,0)
  {crosswise gluing\\along two cuts};

  \draw[->,very thick]
    (1.3,0)--(2.15,0)
    node[midway,above] {glue};

  \begin{scope}[xshift=5.2cm]
    \shade[
      ball color=blue!35,
      opacity=.85
    ]
      (0,0) ellipse (1.65cm and 1.08cm);

    \shade[ball color=white]
      (0,0.06) ellipse (0.72cm and 0.38cm);

    \draw[
      blue!70!black,
      very thick
    ]
      (0,0) ellipse (1.65cm and 1.08cm);

    \draw[
      blue!70!black,
      thick
    ]
      (0,0.06) ellipse (0.72cm and 0.38cm);

    \draw[
      red!75!black,
      very thick,
      ->
    ]
      (-1.38,0.18)
      .. controls (-0.6,0.72) and (0.6,0.72) ..
      (1.38,0.18);

    \draw[
      green!55!black,
      very thick,
      ->
    ]
      (0.05,-0.40)
      .. controls (0.86,-0.62) and (1.02,0.35) ..
      (0.37,0.71);

    \node[
      below=8pt,
      align=center
    ]
    at (0,-1.08)
    {\(\overline{\mathcal C}_d(\mathbb C)\):
     a genus-one Riemann surface\\
     (schematic topological model)};
  \end{scope}
\end{tikzpicture}
\caption{The complex projective curve \(\overline{\mathcal C}_d\) as a
double cover of \(\mathbb P^1\) branched at
\(\alpha=\pm1\) and \(\alpha=\pm\sqrt{1-16d}\).
The branch-point configuration is algebraically exact, whereas the
torus is a schematic representation of the resulting genus-one
Riemann surface.}
\label{fig:cd-complex-cover}
\end{figure}
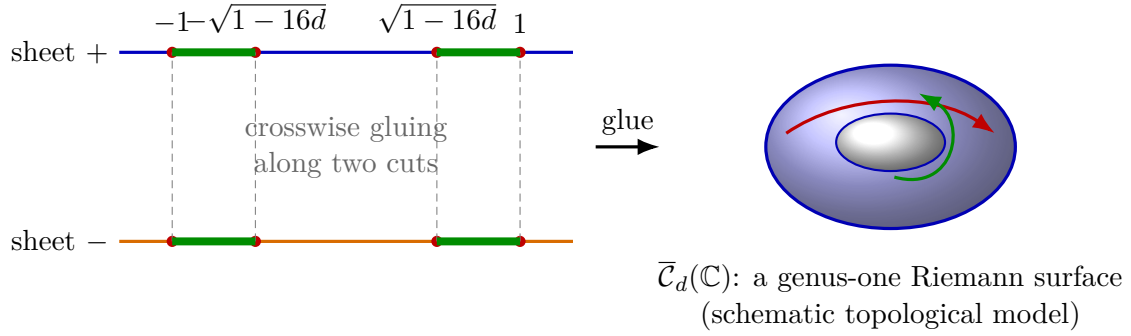

\section[The factorized equation]{Why this particular factorized equation is worth studying}

The family \(\mathcal C_d\) is supported by six mutually reinforcing
geometric, arithmetic, and computational features.
\begin{enumerate}[label=\textup{(\roman*)}]
 \item The equation is a smooth \((2,2)\)-model in
       \(\PP^1\times\PP^1\), rather than an affine quartic whose singular
       plane closure must be mistaken for the elliptic curve.
 \item Its four boundary points are rational and form a marked cyclic
       subgroup.  The same marked divisor controls completeness questions,
       quotient maps, isogeny kernels, and exceptional cases.
 \item The two factor involutions and coordinate exchange extend across the
       smooth completion and generate an intrinsic \(D_8\)-action.  In odd
       characteristic this is literally the square group centered at
       \((-1/2,-1/2)\); in every characteristic its full invariant field is
       generated by \((u^2+u)+(v^2+v)\).  Hence inversion and translation by
       the marked two- and four-torsion points remain visible as sparse
       coordinate symmetries.
 \item Inversion fixes the rational function \(u\).  Consequently the
       Kummer line is not imported from Montgomery coordinates: it is
       generated by a native coordinate and normalized by the boundary to
       \((u+1:u)\).
 \item The odd and binary theories are not unrelated accidents.  The same
       two factors \(u^2+u\) and \(v^2+v\) explain both the centered
       projective variables \(R=2U+Z\), \(S=2V+Z\) and the binary
       Artin--Schreier simplifications.
 \item The parameter \(d\) simultaneously controls smoothness, isomorphism
       classes, completeness, ladder constants, and binary invariants.  A
       single parameter-selection problem can therefore be joined to a
       single arithmetic interface.
\end{enumerate}

These observations create a coherent research program.  First determine
the marked geometry and its intrinsic quotient.  Next derive the group law
on \(\mathcal C_d\) itself.  Only then use auxiliary models to optimize
dependency graphs and pull the results back.  Finally test whether the same
native functions continue to organize division polynomials, isogenies,
Miller functions, and finite-field statistics.  The subsequent parts follow
this order.

\section[Why Cd is a distinct and useful arithmetic model]
{Why \texorpdfstring{\(\mathcal C_d\)}{Cd} is a distinct and useful arithmetic model}

The preceding section motivates the equation from its factorization and
marked geometry.  This section has the narrower task of stating the
arithmetic consequences proved later; it is a summary of benefits rather
than a second motivation section.

Within the arithmetic framework and the corresponding cost models developed
in this monograph, \(\mathcal C_d\) curves attain the optimized Edwards
field-operation counts on the loci where the two dependency graphs are
linearly identified, while retaining additional native advantages.  These
include a characteristic-uniform defining
equation, a natural \(D_8\)-symmetry, four rational boundary points, native
Kummer coordinates, preservation of the same equation form in characteristic
\(2\), and a more natural unifying framework encompassing QRT, reciprocal,
twisted, and symmetric biquadratic extensions.

These advantages arise from the interaction between the factorized equation,
the marked boundary, the intrinsic quotient, and the characteristic-uniform
geometry.  The principal features established throughout the monograph are
summarized below.

\begin{enumerate}[label=\textup{(\alph*)}]
 \item \textbf{One equation in every characteristic.}
       The factorized equation
       \((u^2+u)(v^2+v)=d\) requires no change of model in
       characteristics two or three.  Only the interpretation of the two
       quadratic factors changes.
 \item \textbf{A smooth symmetric \((2,2)\)-completion.}
       The two coordinates occur symmetrically, while the chosen identity
       and inverse are visible from the boundary:
       \(O=(0,\infty)\) and
       \(-(u,v)=(u,-v-1)\).  The eight natural transformations realize all
       square symmetries; the odd-characteristic center is
       \((-1/2,-1/2)\), and the characteristic-uniform algebraic quotient is
       generated by \((u^2+u)+(v^2+v)\).
 \item \textbf{Marked rational four-torsion.}
       The four boundary points form a distinguished cyclic subgroup.
       This marking simultaneously explains the Edwards-type arithmetic in
       odd characteristic and the \(Z/4\mathbb Z\)-normal geometry in
       characteristic two.
 \item \textbf{A genuinely native Kummer coordinate.}
       The fixed field of inversion is \(k(u)\), and the boundary
       normalization selects
       \(\kappa_d(P)=(u+1:u)\).  Thus scalar multiplication, differential
       addition, division polynomials, isogenies, and pairing denominators
       can share one coordinate native to \(\mathcal C_d\).
 \item \textbf{Complete internal addition laws.}
       The model carries its own internal arithmetic rather than merely
       importing formulas from auxiliary isomorphic models.  On its smooth
       \((2,2)\)-completion, \(\mathcal C_d\) admits a native complete
       addition law: a single formula is \(\F_q\)-complete when
       \(\rho=1-16d\) is a nonsquare, while a finite base-point-free complete
       addition-law atlas exists in general.  Its intrinsic Kummer quotient
       also supports direct differential-addition formulas.  In
       characteristic two, the Artin--Schreier structure of the defining
       equation yields native complete full-point and differential-addition
       atlases without requiring arithmetic on an auxiliary Edwards,
       Montgomery, or Weierstrass model.
 \item \textbf{Characteristic-specific simplification without changing the
       declared model.}
       In odd characteristic the linear shifts \(2u+1,2v+1\) expose
       Edwards-type bilinear laws; in characteristic two the same factors
       are Artin--Schreier polynomials; in characteristic three cubing is
       Frobenius.  The optimized formulas differ, but every input and output
       remains a \(\mathcal C_d\)-point or its declared native quotient.
 \item \textbf{Composable arithmetic.}
       Full addition, Kummer ladders, halving, tripling, closed \(2P+Q\),
       division polynomials, model-preserving isogenies, and Miller
       functions are expressed through a common \((u,v,d)\)-interface.  This avoids
       treating conversion to another model as the mathematical endpoint.
 \item \textbf{Transparent parameter control.}
       The single parameter \(d\) controls smoothness, the Edwards
       completeness parameter \(\rho=1-16d\), the Montgomery ladder
       constant, the binary \(j\)-invariant, and the marked-torsion
       subfamily.  Parameter selection can therefore be analyzed together
       with the arithmetic rather than after choosing unrelated models.
\end{enumerate}

\section{Architecture and organizing themes}
\label{sec:scope-organizing-themes}

The theory developed in this monograph is organized around six
interconnected themes.  Together they describe the geometric,
arithmetic, algorithmic, and cryptographic architecture of the
\(\mathcal C_d\) model.

\paragraph{Intrinsic geometry and symmetry.}
The first theme is the geometry of the smooth \((2,2)\)-completion of
\[
             \mathcal C_d:\qquad (u^2+u)(v^2+v)=d.
\]
This includes the exact smoothness criterion, the marked identity and
boundary, rational four-torsion, the inverse map, and the eight natural
coordinate transformations.  In odd characteristic these transformations
become the signed symmetries of a square after centering, while in every
characteristic their complete invariant field can be described explicitly.
The reciprocal presentation
\[
             (u+1)(v+1)=du^2v^2
\]
is treated as a complementary affine chart of the same completed curve,
with its marked points and automorphisms obtained by explicit conjugation.

\paragraph{Native full-point and quotient arithmetic.}
The second theme is arithmetic whose declared input and output remain on
\(\mathcal C_d\), its smooth completion, or its intrinsic Kummer quotient.
The monograph constructs native addition, doubling, mixed addition,
differential addition, halving, tripling, closed \(2P+Q\), full-point
recovery, and scalar-multiplication procedures.  It distinguishes formulas
derived directly from the biquadratic equation from circuits obtained by
transport through Edwards, Montgomery, Jacobi, or Weierstrass coordinates.
The native Kummer coordinate and its normalization are characterized
intrinsically, and the exceptional loci of the differential formulas are
completed by explicit projective charts.

\paragraph{Characteristic-dependent realizations of one model.}
The third theme is characteristic uniformity without imposing one arithmetic
circuit on every field.  In odd characteristic, centered coordinates expose
Edwards and Montgomery structures and give optimized full-point and Kummer
arithmetic.  The resulting \(\mathcal C_d\) implementations attain the
corresponding optimized Edwards or Montgomery dependency graphs on the
linearly identified loci while retaining the native equation, boundary,
four-torsion, symmetry, and quotient interface.  In characteristic two, the
same defining equation retains its Artin--Schreier structure and leads to
native binary full-point, Kummer, recovery, and complete differential
atlases.  Characteristic three is treated within the odd branch, with
Frobenius, Verschiebung, scalar tripling, and separable three-isogenies kept
mathematically distinct.

\paragraph{Global arithmetic constructions.}
The fourth theme concerns structures extending beyond scalar multiplication.
The monograph determines geometric and finite-field isomorphism classes,
parameter multiplicities, point-count moments, division polynomials,
division-point mean values, model-preserving isogenies, and native Miller
functions.  It also develops Tate and Weil pairing wrappers whose running
points and declared outputs remain in the \(\mathcal C_d\) interface.
Special \(j=0\) and \(j=1728\) loci yield explicit low-cost complex
multiplication maps, together with their Kummer actions and the precise
field conditions under which GLV decompositions are available.

\paragraph{Concrete arithmetic and the Cd25519 specialization.}
The fifth theme is the conversion of the structural theory into exact
algorithms and parameter choices.  Operation counts are attached to
specified coordinate schedules and declared outputs, with conversion,
normalization, validation, and recovery costs separated from recurring loop
costs.  For
\[
 d_{\rm C}=(16\cdot121666)^{-1}
 \qquad\text{over}\qquad
 \mathbb F_{2^{255}-19},
\]
the resulting Cd25519 model has native Kummer coordinate equal to the X25519
\(U\)-coordinate.  Its full-point arithmetic attains the optimized
\(A=-1\) Edwards dependency graph after an invertible native recoding, while
its original biquadratic presentation retains additional marked geometry,
symmetry, and native-input optimizations.

\paragraph{Extensions, reciprocal families, and QRT state arithmetic.}
The sixth theme places \(\mathcal C_d\) inside a broader biquadratic
framework.  The symmetric family
\[
 \mathcal C_{a,b,d}:
 \qquad (u^2+u+a)(v^2+v+b)=d
\]
separates the characteristic-uniform product geometry from the additional
symmetries and rational torsion of the specialization
\(\mathcal C_d=\mathcal C_{0,0,d}\).  The one-sided family
\(\mathcal T_{a,d}\), reciprocal \(C\)-curves, and the symmetric QRT family
\[
 \mathcal Q_{\alpha,\beta,\gamma}:
 \qquad
 x^2y^2+\alpha(x^2+y^2)+\beta xy+\gamma=0
\]
are then developed through explicit Weierstrass, Jacobi, Montgomery, and
binary models.  In the QRT setting, the intrinsic McMillan displacement is
separated from the choice of group origin, and adjacent Kummer states turn
a fixed translation into logarithmic scalar multiplication through state
doubling.  Elliptic Lucas sequences, elliptic divisibility sequences,
division sections, elliptic nets, and coordinatewise isogenies are thereby
assembled into one compatible state-arithmetic framework.

Here a \emph{QRT state} (or \emph{genus-one state}) is a pair of adjacent
Kummer values constrained to lie on a smooth genus-one curve; after an
origin is chosen, state evolution becomes elliptic translation.  This
terminology distinguishes the constrained two-coordinate state from an
arbitrary pair of quotient-line values.

These six themes determine the logical organization of the monograph.  The
next section records the principal results in the order in which their
mathematical prerequisites are established.

\section{Main contributions}

The principal results are listed below in dependency order.  The list is
intended as a theorem-level synopsis of the monograph; precise hypotheses,
field-of-definition conditions, exceptional loci, and proofs are given in
the corresponding chapters.

\begin{enumerate}[label=\textup{(\arabic*)}]
 \item The completion of \(\Cd\) is smooth precisely when
       \(d(1-16d)\ne0\), has genus one, and possesses a characteristic-free
       cyclic four-torsion subgroup on the boundary.  Three elementary
       symmetries generate a faithful intrinsic \(D_8\)-action on that
       boundary and on the smooth completion.  In odd characteristic the
       action is the signed permutation group about
       \((-1/2,-1/2)\); its four boundary points form the square-vertex
       representation.  In every characteristic the full fixed field is
       generated by \((u^2+u)+(v^2+v)\).
 \item The fixed function field under inversion is \(k(u)\).  Hence every
       degree-one Kummer coordinate is a M\"obius transform of \(u\).  The
       conditions \(O\mapsto\infty\), \(T\mapsto0\) leave the unavoidable
       target scale \(c\in k^\times\); fixing \(c=1\) selects
       \((u+1:u)\).
 \item The reciprocal presentation
       \((u+1)(v+1)=du^2v^2\) is proved to be an ambient reciprocal chart of
       the same smooth \((2,2)\)-model.  Its marked points, inverse, and all
       eight \(D_8\)-transformations are obtained by explicit conjugation.
       Its quotient coordinate \((v+1:1)\) is proved equal to the native
       Kummer point \((u_0+1:u_0)\).  Odd Montgomery, reduced binary
       Weierstrass, general, mixed, and fully affine differential formulas
       are derived, and their endpoint and recurring costs are compared
       without counting the chart as a second model.
 \item In odd characteristic, explicit Edwards, inverted Edwards,
       Montgomery, and extended Edwards dictionaries are given and proved.
       In characteristic two, explicit \(Z/4\mathbb Z\), binary Weierstrass,
       and Kummer dictionaries are given and proved.
 \item Before those dictionaries are used for optimization, explicit
       addition and doubling formulas are proved directly in the native
       \(u,v\) coordinates in odd characteristic, characteristic three,
       and characteristic two.
 \item Affine incompleteness is separated from completeness on the smooth
       model.  A native Segre law is proved directly from
       \(\mathcal C_d\), its exact nonsquare-\(\rho\) completeness
       criterion is established, and a finite complete native atlas is
       proved to exist in every characteristic.
 \item Standard and multiplication--squaring-tradeoff differential
       formulas are stated as operations on the native quotient
       \((u+1:u)\).  Their domains, quotient degrees, recovery maps, and
       operation counts are stated separately.
 \item The native-input identity \(X-Z=1\) gives a first Kummer double at
       cost \(\M+\Sqr+\Dpar\), compared with
       \(2\M+2\Sqr+\Dpar\) for a general projective Kummer input.
       Dedicated native full-point doubling and mixed-addition schedules
       are stated separately from the strongly unified law.
 \item Odd-characteristic halving is factored into two quadratic equations.
       A square-class criterion and a linear recovery formula for the
       ordinate are obtained.
 \item Characteristic three admits both a native full-point tripling formula
       and a Frobenius-specialized Kummer tripling map.  The latter is
       identified as the Kummer shadow of
       \([3]=V_{3,d}\circ F_{3,d}\); the ordinary and supersingular cases,
       as well as a separable V\'elu quotient by a reduced kernel, are
       distinguished.  Classical Hessian formulas do not survive by naive
       reduction.
 \item Exact geometric and rational isomorphism counts, parameter-weighted
       first and second Frobenius-trace moments, division-polynomial
       recurrences, and fiber averages are established in both
       characteristic branches.
 \item Odd-order separable kernels give model-preserving isogenies.  In odd
       characteristic the formulas are multiplicative Edwards products; in
       characteristic two they reduce to half-kernel sums with
       \(d'=d+\Sigma_K^{1/2}+\Sigma_K^{1/4}\), where
       \(\Sigma_K\) is the half-kernel abscissa sum.
 \item The odd-degree quotient specializes on the native \(u\)-Kummer line
       to two monic products and cost \(2m\M+2\Sqr\) for
       \(\ell=2m+1\); the same target parameter is one eighth-power kernel
       product.  The theorem is placed in the general isogeny chapter and
       rederived in the self-contained cryptographic chapter.
 \item In characteristic two, the reduced two-torsion quotient
       \(W_{A,d}\to W_{A,\sqrt d}\) preserves \(A\).  It costs one
       multiplication on native \(C_d\), and one multiplication or one
       square on \(\mathcal T_{a,d}\); it also gives Jacobian closure for
       \(\mathcal C_{a,b,d}\).  Exact binary Edwards endpoint costs,
       target re-embedding equations, a reference benchmark, and exhaustive
       finite-field checks accompany the proof.
 \item Native tripling, native closed \(2P+Q\), pulled-back Miller functions,
       reduced Tate pairings, and Weil pairings are derived, proved, assigned
       explicit costs, and verified over explicit finite fields.
 \item The \(j=1728\) member \(d=1/8\) and the \(j=0\) locus
       \(A^2=3\) admit explicit fourth- and third-order endomorphisms in
       \(u,v\).  Their Kummer actions, exact costs, eigenvalue equations,
       and finite-field GLV examples are established.
 \item For
       \(d_{\rm C}=(16\cdot121666)^{-1}\) over
       \(\F_{2^{255}-19}\), the exact Cd25519 dictionaries are derived.
       The native quotient equals the X25519 \(U\)-line; an
       \(\imath\)-twisted native Segre recoding gives complete addition at
       \(8\M+\Dpar\) and mixed addition at \(7\M+\Dpar\); and
       leading-bit, affine-difference, and \(U=9\) fixed-base ladder
       schedules are proved and assigned explicit costs.  Encoding, twist,
       cofactor, and
       constant-time boundaries are stated explicitly.
 \item The symmetric family
       \(\mathcal C_{a,b,d}:(u^2+u+a)(v^2+v+b)=d\) is placed in the same
       geometric framework.  Its smoothness condition, characteristic-free
       factor symmetries, odd-characteristic diagonal Edwards equation,
       odd and binary Weierstrass Jacobians, and the specializations
       \(b=0\), \(a=b\), and \(d=ab\) are organized so that
       \(\mathcal C_d=\mathcal C_{0,0,d}\) is recovered without changing
       notation or origin conventions.
 \item The one-sided family \(\mathcal T_{a,d}\) is proved to represent the
       Jacobians of the full product family, with actual curve isomorphisms
       over finite fields.  Its odd-characteristic twisted-Edwards and
       Montgomery dictionaries, the distinguished value \(a=1/2\), and its
       characteristic-two universal ordinary Weierstrass form are made
       explicit.
 \item Reciprocal \(C\)-curves
       \((x^2-\tau)(y^2-\sigma)=\kappa xy\) are developed through their
       smooth completion, quotient involutions, four-torsion action, Jacobi
       and Weierstrass models, complete full-point atlases, Kummer
       differential arithmetic, and characteristic-two identification with
       \(\mathcal C_d\).
 \item The symmetric QRT envelope is given exact odd- and
       characteristic-two smoothness criteria, Weierstrass and Jacobi
       models, finite-field model and abstract-isomorphism classifications,
       explicit two-, three-, and odd-prime-degree isogeny transport, and
       complete arithmetic interfaces with all field-of-definition and
       exceptional-fiber hypotheses stated.
 \item The Vieta involutions and coordinate exchange are proved to produce
       the McMillan translation on every smooth QRT fiber.  The intrinsic
       displacement is separated from the origin-dependent point
       \(D=T(O)\); pointed curves admit adjacent-Kummer-state models, explicit
       state doubling, differential recovery, and binary or signed-digit
       ladder recurrences.  Exact evenization is characterized by a split
       rational two-torsion Kummer involution and needs an extension of
       degree at most six.  On the split locus, a complete one-core state
       update has the proved bound \(4\M+4\Sqr+2\Dpar\); the ordinary-binary
       counterpart is universal for ordinary pointed curves over perfect
       binary fields, carries the full Artin--Schreier twist data, and costs
       \(4\M+5\Sqr+3\Dpar\).  Recovery, origin alignment, generic-MSM
       extensions, and isogeny functoriality are included with their
       precise hypotheses.
\end{enumerate}

\section{Guide to the monograph}

The logical dependencies are more important than the chronological order in
which individual formulas were discovered:
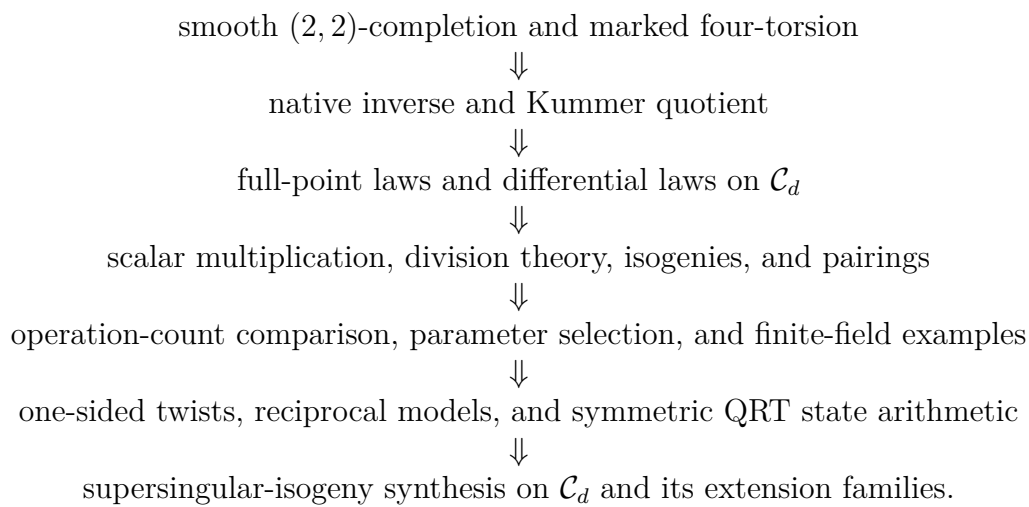
\[
\begin{array}{c}
\text{smooth \((2,2)\)-completion and marked four-torsion}\\
\Downarrow\\
\text{native inverse and Kummer quotient}\\
\Downarrow\\
\text{full-point laws and differential laws on \(\mathcal C_d\)}\\
\Downarrow\\
\text{scalar multiplication, division theory, isogenies, and pairings}\\
\Downarrow\\
\text{operation-count comparison, parameter selection, and finite-field examples}\\
\Downarrow\\
\text{one-sided twists, reciprocal models, and symmetric QRT state arithmetic}\\
\Downarrow\\
\text{supersingular-isogeny synthesis on \(\mathcal C_d\) and its extension families}.
\end{array}
\]

The dependency blocks used throughout the book are summarized here.  A
right-hand block uses only objects constructed in blocks to its left or in
earlier rows.
\begin{center}
\small
\setlength{\tabcolsep}{4pt}
\renewcommand{\arraystretch}{1.18}
\begin{tabularx}{\textwidth}{@{}
  >{\raggedright\arraybackslash}X
  >{\raggedright\arraybackslash}X
  >{\raggedright\arraybackslash}X@{}}
\toprule
\textbf{Foundation} & \textbf{Arithmetic layer} & \textbf{Global layer}\\
\midrule
smooth completion, boundary, inverse
 & native affine and Segre laws
 & complete scalar multiplication\\
native Kummer field and normalization
 & \(x\mathrm{DBL}\), \(x\mathrm{ADD}\), recovery
 & division theory and mean values\\
odd and binary dictionaries
 & tripling, halving, CM maps
 & isogenies and pairings\\
cost and output conventions
 & parameter-specific circuits
 & comparison and implementation\\
\bottomrule
\end{tabularx}
\end{center}

Part I explains the motivation, literature, common projective geometry,
square symmetry and invariant algebra, and coordinate dictionaries.  It
establishes the objects used everywhere else.
Part II develops the odd-characteristic branch, always beginning with
native \(u,v\)-formulas before optimized projective or Kummer formulas.  It
then studies the reciprocal chart as an internal affine normalization of
the same Kummer line, and isolates the two low-cost CM endomorphism loci
before treating characteristic three in its final chapter.  Part III repeats
the same sequence in characteristic two: native group law,
Artin--Schreier geometry, binary Kummer line, recovery, and full-point
arithmetic.  Part IV asks whether the native interface survives beyond
scalar multiplication, treating moduli, moments, division polynomials,
isogenies, and pairings.  Part V turns the theory into a decision framework
through parameter selection, the Cd25519 specialization, model-by-model
comparisons, and finite-field examples.  Part VI begins with the symmetric
three-parameter product family and then develops four connected chapters:
the one-sided family in odd and binary characteristic, reciprocal
\(C\)-curves, and the symmetric QRT envelope with its McMillan and
adjacent-Kummer-state arithmetic.  The Cd25519 chapter
deliberately precedes the comparison chapter, so every parameter-specific
cost cited in the comparison has already been proved.  Part VII then uses
the completed \(\mathcal C_d\), full-product, one-sided, reciprocal, and QRT theories
to develop a self-contained supersingular-isogeny interface and an
equal-output comparison with Edwards arithmetic.

Readers interested primarily in implementation may read Chapters~
\ref{ch:geometry}, \ref{ch:odd-full}, \ref{ch:odd-kummer},
\ref{ch:reciprocal-chart},
\ref{ch:binary-geometry}, \ref{ch:binary-kummer}, and
\ref{ch:Cd25519}, then Chapter~\ref{ch:comparison}.  The fixed-translation
and adjacent-state viewpoint is developed in
Chapters~\ref{ch:reciprocal-C-curves} and
\ref{ch:symmetric-QRT-envelope}; concrete native templates remain in
Appendix~\ref{app:algorithms}.  Readers
interested in arithmetic geometry may instead continue from
Chapter~\ref{ch:geometry} to Chapters~\ref{ch:moduli}--\ref{ch:pairings}.
For isogeny cryptography, Chapter~\ref{ch:Cd-isogeny-cryptography} may be
read after Chapters~\ref{ch:isogenies}, \ref{ch:generalized-Cabd},
\ref{ch:one-sided-twist}, \ref{ch:reciprocal-C-curves}, and
\ref{ch:symmetric-QRT-envelope}; its central formulas are nevertheless
restated so that the cryptographic comparison does not require repeated
backward lookup.  In all routes the auxiliary coordinate dictionaries are
tools, not a change in the object being studied.

\chapter{Background and State of the Art}
\label{ch:background}
\section{Weierstrass equations and arithmetic models}

A Weierstrass equation is the universal algebraic interface for an elliptic
curve with a chosen origin.  For explicit computation, however, the degrees
of the embedding and the action of rational torsion strongly influence the
shape and completeness of the group law.  Montgomery's model
\[
  BV^2=U^3+AU^2+U,\qquad B(A^2-4)\ne0,
\]
made the quotient \(E/\{\pm1\}\) computationally central: the \(U\)-coordinate
supports differential addition and a regular ladder without an ordinate
\cite{Montgomery1987,CostelloSmith2018}.  Edwards introduced a symmetric
quartic normal form, and Bernstein--Lange converted its symmetry into fast,
unified addition formulas \cite{Edwards2007,BernsteinLange2007}.
Twisted, inverted, projective, and extended Edwards coordinates subsequently
optimized different full-point tasks
\cite{BernsteinEtAl2008,BernsteinLangeInverted,HisilEtAl2008}.

The terminology ``new model'' therefore concerns more than an abstract
isomorphism class.  A model may be valuable because its embedding linearizes
a torsion action, because an addition law is complete on rational points, or
because its Kummer quotient has a sparse pseudo-addition law.  The present
family is assessed by exactly these criteria.

\section{Characteristic two}

Naively reducing odd-characteristic formulas modulo two usually destroys the
separable double cover underlying the ordinary Montgomery form.  Binary
Edwards curves supplied complete addition formulas for ordinary binary
elliptic curves \cite{BinaryEdwards2008}.  Kohel's
\(Z/4\mathbb Z\)- and \(\mu_4\)-normal forms instead organize the arithmetic
around a rational four-torsion point and its theta group
\cite{Kohel2012,KohelTwistedMu4}.  Gaudry--Lubicz developed the associated
Kummer pseudo-addition law in characteristic two
\cite{GaudryLubicz2009}.  A central result of this monograph is the linear equivalence of
\eqref{eq:model} with the split \(Z/4\mathbb Z\)-normal form over the
base field.

\section{Characteristic three}

In characteristic three, the same model equation
\[
             \mathcal C_d:\quad (u^2+u)(v^2+v)=d
\]
is smooth precisely when \(d\ne0,1\), and \(2\) is still invertible.
Consequently
the centered Edwards dictionary and the Montgomery quotient dictionary of
Chapter~\ref{ch:odd-dictionary} remain valid.  In particular, ordinary
full-point addition and doubling are obtained by specializing the
odd-characteristic formulas, while Montgomery \(x\)DBL and \(x\)ADD retain
their usual projective form.  This observation is important: differential
addition does not fail merely because the field has characteristic three.

What changes is the multiplication-by-three map.  For any elliptic curve in
characteristic three,
\[
             [3]=V\circ F,
\]
where \(F:E\to E^{(3)}\) is the relative Frobenius and
\(V:E^{(3)}\to E\) is the Verschiebung.  The Frobenius part is purely
inseparable of degree \(3\); on an ordinary curve the Verschiebung is
separable of degree \(3\), whereas on a supersingular curve it is also
purely inseparable.  A point-sum V\'elu construction requires separate data:
it needs a specified reduced cyclic kernel and does not sum
over the full group scheme \(E[3]\).  Thus cubing can be treated as a
Frobenius operation rather than a general multiplication, but only after
the separable and inseparable factors have been identified on the chosen
model.

Farashahi, Wu, and Zhao developed low-cost doubling, tripling, and addition
for characteristic-three elliptic curves by exploiting exactly this
Frobenius structure \cite{FarashahiWuZhao2013}.  Their work provides the
proper comparison point for the ternary branch of this monograph.  It does
not, however, justify reducing a diagonal Hessian model unchanged, since
\[
 X^3+Y^3+Z^3=(X+Y+Z)^3.
\]
The diagonal cubic therefore loses the geometry on which the usual Hessian
addition law depends.

The approach taken here is model-preserving.  Chapter~\ref{ch:char3} first
writes addition and doubling on \(\mathcal C_d\), then retains the
two-to-one quotient
\(\kappa_d=(u+1:u)\) for differential addition, and finally derives a
Frobenius-aware quotient tripling formula from the multiplication-by-three
map.  The result is compared separately with full-point tripling: the former
returns only \(\kappa_d([3]P)\), whereas the latter returns
\((u_{3P},v_{3P})\).  The chapter then states the
Frobenius--Verschiebung factorization on the smooth native completion and
identifies exactly when a separable three-point kernel can enter V\'elu's
construction.  The characteristic-three analysis therefore keeps the curve model, quotient
map, separability, and output representation separate.

\section{Efficient endomorphisms and GLV decomposition}

For elliptic curves in characteristic different from \(2,3\), the
origin-preserving geometric automorphism group is generically
\(\{\pm1\}\), and it enlarges at \(j=1728\) and \(j=0\)
\cite{Silverman2009}.  When the corresponding roots of unity lie in the
base field, these automorphisms can become efficient endomorphisms for
scalar multiplication.  Gallant, Lambert, and Vanstone showed how an
endomorphism eigenvalue on a large prime-order subgroup produces a
two-dimensional scalar decomposition
\cite{GallantLambertVanstone2001}.

For the present model, complex multiplication on the abstract curve is only
one part of the problem.  The endomorphism must also have a sparse expression
in the declared \(u,v\) coordinates, descend to the native Kummer line, and
satisfy the required field-of-definition and subgroup conditions.
Chapter~\ref{ch:native-endomorphisms} establishes the corresponding
arithmetic results for the
\(j=1728\) member \(d=1/8\) and the \(j=0\) locus
\(((4d)^{-1}-2)^2=3\).

\section[Arithmetic extensions]{Differential addition, division polynomials, isogenies, and pairings}

Differential-addition formulas for twisted Edwards and Jacobi quartic
models trade multiplications for squarings and parameter multiplications
\cite{WuSong2022}.  For Jacobi quartics, the work of
Wu and Song was received in August 2021, finalized in October 2021, and
published in 2022.  It gives mixed differential-addition/doubling costs
\(5\M+4\Sqr+\Dpar\) and, under an additional square-class condition,
\(3\M+6\Sqr+3\Dpar\).  The twisted-Edwards tradeoff also includes a
\(3\M+7\Sqr+\Dpar\) schedule; the relevant identities are derived
directly in Chapter~\ref{ch:tradeoffs}.  This monograph identifies the
exact native quotient maps on which the transported circuits act on
\(\Cd\), rather than relying on a coordinate-level analogy.

Division polynomials on Weierstrass, Edwards, and Jacobi quartic models
encode torsion and multiplication maps
\cite{Silverman2009,HittMcGuireMoloney2008,MoodyJacobiDivision2011}.
They also give coordinate averages over division fibers
\cite{MoodyMeanEdwards2011,FengWuMean2014}.  V\'elu's formulas compute
normalized isogenies from a finite separable kernel
\cite{Velu1971}; Moody--Shumow derived lower-complexity analogues on
Edwards and Huff curves \cite{MoodyShumow2016}.

For pairings, Miller's divisor recursion remains the fundamental
computational mechanism \cite{Miller2004}.  Edwards conic constructions
provide an established benchmark for pairing arithmetic
\cite{AreneLangeNaehrigRitzenthaler2011}.  The native
\(\mathcal C_d\) framework developed here evaluates sparse Miller
numerators directly in the \(u,v\)-coordinates, while its vertical
denominators depend only on the native Kummer coordinate.  The same
characteristic-uniform wrapper supports the odd-characteristic,
characteristic-three, and characteristic-two Artin--Schreier branches.
The resulting end-to-end cost ledger records the embedding degree, twist,
extension-field arithmetic, loop length, denominator elimination, final
exponentiation, and the native point operations used inside the Miller loop.

Complete addition laws on projective embeddings are naturally formulated
through their base loci and the associated line bundles
\cite{AreneKohelRitzenthaler2012}.  The base-locus analysis below
distinguishes a single finite-field-complete tuple from a geometrically
complete finite atlas on the smooth \((2,2)\)-model.

\section{The 25519 specialization as a test case}

Curve25519 and its Edwards form provide a stringent test of whether the
model-internal viewpoint has practical content.  The standard Montgomery
coefficient is \(A=486662\), the quotient base point is \(U=9\), and the
field is \(\F_{2^{255}-19}\)
\cite{BernsteinCurve25519,RFC7748}.  Solving the native parameter relation
\[
                 A=\frac1{4d}-2
\]
selects
\[
                 d=(16\cdot121666)^{-1}.
\]
The question is then more precise than asking whether the abstract curve is
isomorphic to Edwards25519.  One must determine whether the isomorphism
returns full points to \(u,v\), whether the quotient coordinate agrees with
the X25519 wire coordinate, whether the complete law can be accelerated in
coordinates intrinsic to the native completion, and which X25519 protocol
rules survive unchanged.

Chapter~\ref{ch:Cd25519} develops these constructions and establishes the corresponding
arithmetic results.  It derives the two-way
Cd25519--Edwards25519 dictionary, proves equality of the native Kummer ratio
and the Curve25519 \(U\)-coordinate, constructs an
\(\imath\)-twisted native complete law, specializes \(x\)DBL and \(x\)ADD,
and analyzes leading-bit and fixed-base acceleration.  It also records two
protocol facts: generic X25519 inputs range over the main curve and its
quadratic twist, while Cd25519 uses its dedicated ladder and complete-law
accelerations independently of the two CM-special parameter loci developed
earlier.

\section{Contributions and relation to established models}

This monograph develops \(\mathcal C_d\) as an independent,
characteristic-uniform biquadratic arithmetic model with a single smooth
projective completion and a native Kummer coordinate valid in every
characteristic.  Edwards, Montgomery, Weierstrass, and related models are
connected to \(\mathcal C_d\) through explicit dictionaries and serve as
proof frameworks and optimization bridges; every transported construction
is returned explicitly to the native \(u,v\)-coordinates or to a declared
native quotient.  Transported formulas retain attribution to their original
sources.

The contributions of this monograph include the characteristic-uniform
\(\mathcal C_d\) framework, native full-point and quotient arithmetic,
characteristic-specific circuits where direct formulas are required, and a
structural synthesis connecting marked geometry, complete addition laws,
Kummer arithmetic, division theory, isogenies, endomorphisms, pairings, and
adjacent-state QRT dynamics.  For every arithmetic circuit, the domain of
validity, quotient degree, retained output information, completeness
properties, and operation count are stated explicitly.

\chapter{Conventions and Computational Cost}
\label{ch:conventions}
\section{Fields, points, and operation counts}

Unless stated otherwise, \(k\) is a field, \(\bar k\) an algebraic closure,
and \(\F_q\) a finite field.  We use \(\M,\Sqr,\Dpar,\Inv\) for one
base-field multiplication, squaring, multiplication by a fixed curve
constant, and inversion.  When the multiplier matters,
\(\Dconst{c}\) denotes multiplication by the displayed constant \(c\), and
\(\Dnine\) abbreviates \(\Dconst{9}\).  Additions, subtractions, negations,
and multiplications by small integers realized by additions are normally
omitted from headline costs.  The one explicit exception is \(\Dnine\),
which is retained in the Cd25519 fixed-difference ledger so that the
specialization \(U=9\) remains visible and can be assigned a
platform-specific cost.  Uncharged small-integer operations are still
retained in dependency graphs when they affect data flow.  In binary
Kummer arithmetic, \(\mBase\) denotes
multiplication by the normalized base-point coordinate and \(\mCurve\)
multiplication by \(d^{-1}\).  In characteristic three, \(\Cube\) denotes
Frobenius cubing.  It may be a coordinate permutation in a normal basis, but
it is not assigned zero cost unless explicitly stated.

Unless an extension field is explicitly named, these symbols count
operations in the declared base field.  In the pairing chapter, base-field
point arithmetic and extension-field accumulator arithmetic are therefore
reported separately.  A base-field count alone is never presented as the
cost of an entire pairing, encoding, or protocol.

For platform-dependent comparisons we write
\[
  \gamma_{\rm S}=\frac{\Sqr}{\M},\qquad
  \gamma_{\rm D}=\frac{\Dpar}{\M}.
\]
Hardware evaluation supplements these symbolic counts with reduction
latency, instruction scheduling, memory traffic, and register pressure.

\section{Coordinate states, preprocessing, and charged work}

A projective tuple is always interpreted up to a common nonzero scalar.  A
headline cost assumes that every input is already in the working coordinates
named in the statement and that the output remains in those coordinates.
Thus an affine-to-projective lift, a diagonal recoding, a Kummer
normalization, full-point recovery, and final affine conversion are charged
only when the statement says that they are included.  This convention is
particularly important when a low-cost loop uses a different representation
from its external \((u,v)\)-interface.

The term \emph{mixed addition} means that the running point is projective
and the second point is affine or precomputed with its scale coordinate equal
to one.  A table entry quoted as a mixed cost does not include table
construction, constant-time table selection, or conversion of an arbitrary
second input into the normalized form.  Similarly, batch inversion,
fixed-base precomputation, and parameter setup are recorded separately from
the recurring loop cost.  When a scalar loop has \(N\) iterations, the
accounting principle is
\[
 C_{\rm total}
   =C_{\rm input}+C_{\rm pre}+N C_{\rm loop}+C_{\rm recovery}+C_{\rm output}.
\]
Any vanishing term must be justified by the stated input or output format.
Conditional swaps, memory traffic, and side-channel countermeasures are not
field operations, but they remain part of an implementation and are
discussed separately in Chapter~\ref{ch:implementation}.

Single capital letters introduced as intermediate quantities inside a
displayed formula, dependency graph, or algorithm are local to that block.
Likewise, the letters in a newly declared projective tuple are local
coordinate names.  Such declarations may shadow a global model parameter;
outside their declared block, \(A=(4d)^{-1}-2\) retains its meaning as the
odd-characteristic Montgomery coefficient.  This convention prevents local
temporaries such as \(A,B,C,D\) from creating new global notation.  In
particular, \(T\) and \(R\) in a declared coordinate tuple are coordinates;
outside that tuple the same letters denote the marked boundary points fixed
in Chapter~\ref{ch:geometry}.

\section{Output categories}

The following output types are not interchangeable:
\[
\begin{array}{c|l}
\Full & \text{a full point on the smooth completed curve},\\
\Ffin& \text{a full point only while it remains on the finite native chart},\\
\Ktwo & \{P,-P\}\text{ on the ordinary Kummer line},\\
\Kfour& \{\pm P\}+\langle T_4\rangle\text{ on the Kummer line of the
order-four quotient}.
\end{array}
\]
Here \(T_4\) denotes the marked point of order four; under the
odd-characteristic Edwards dictionary it is \((1,0)\), the image of the
native boundary point \(R\).
A circuit producing \(\Ffin\), \(\Ktwo\), or \(\Kfour\)
cannot replace a
full-point operation unless an explicit recovery step is available and its
cost is included.

Here a \(\Full\)-output is a point of the smooth completed curve; it need not
lie in the finite affine \((u,v)\)-chart.  A \(\Ffin\)-output retains both
coordinates only on that finite chart; a one-scale plane-quartic tuple at
infinity records merely the image of a branch under the singular plane
model.  By contrast, \(\Ktwo\) forgets one sign bit, while the higher-degree
quotients forget additional torsion data.
Equality of projective quotient coordinates therefore does not certify
equality of full points.

\section{Domains, completeness, and model transport}

The following terminology is used consistently throughout the monograph.
\begin{enumerate}[label=\textup{(\arabic*)}]
 \item A rational formula is an identity on the dense open set on which its
       denominators, or its homogeneous output tuple, are defined.
 \item A formula is \emph{unified} when the same dependency graph handles
       addition and doubling.  Unified does not mean complete.
 \item A homogeneous addition law is \emph{\(k\)-complete} when it has no
       base pair in \(E(k)\times E(k)\), and is \emph{geometrically
       complete} when the same holds over \(\bar k\).
 \item A finite collection of laws is a \emph{complete atlas} when their
       common base locus is empty.  Individual charts in such an atlas may
       have exceptional pairs.
\end{enumerate}
In particular, a finite affine \((u,v)\)-formula cannot represent
\(P+(-P)=O=(0,\infty)\).  Completeness claims in this book therefore refer
to a homogeneous law on the smooth completion, to an explicitly stated
finite-field domain, or to a complete atlas; Chapter~\ref{ch:native-completeness}
proves the corresponding assertions for \(\mathcal C_d\).

The terms \emph{model dictionary}, \emph{coordinate dictionary}, and
\emph{quotient dictionary} are used in the sense fixed in the Introduction;
the declared information level and all required recovery data form part of
the dictionary.

The word \emph{native} means that the declared input and output are a point
of \(\mathcal C_d\), a point of its own smooth completion, or one of its
declared intrinsic quotients.  An invertible linear recoding of its
projective coordinates is an internal working representation.  A formula
proved on an Edwards, Montgomery, \(Z/4\mathbb Z\), or Weierstrass equation
is called \emph{transported} until the inverse dictionary has returned its
output to the native interface.  This convention prevents an inexpensive
operation on an auxiliary model from being counted as a full
\(\mathcal C_d\)-operation without its required conversion or recovery.

\section{Dependency and local-restatement conventions}

The first occurrence of a symbol in a proof is either a definition in that
proof or a backward reference to an earlier definition.  Forward references
in the Preface, Background, chapter guides, and concluding transitions are
navigational only; no proof of a theorem relies on a result established later.
This convention is useful in a monograph containing several parallel
characteristic branches.

At the start of a chapter or a characteristic-specific section, the full
equation of \(\mathcal C_d\), the characteristic hypothesis, and the local
parameter abbreviations are often restated.  Such a restatement fixes the
local domain and does not redefine the global model.  By contrast, a new
coordinate state---for example, centered Segre, ordinary Kummer, or binary
Kummer coordinates---is introduced once together with its recovery map or
output interpretation; later algorithmic chapters refer back to that
construction instead of repeating its correctness proof.

Proofs that use an auxiliary equation follow the same three-step order:
declare the dictionary, prove or cite the identity in the auxiliary
coordinates, and return the result to the declared native output.  Cost
analysis follows only after those three steps.  These ordering rules are
used to distinguish a logical prerequisite from a comparison or historical
remark.

\section{Generalized Weierstrass equations}

We use
\begin{equation}\label{eq:general-W}
 E:\quad Y^2+a_1XY+a_3Y=X^3+a_2X^2+a_4X+a_6.
\end{equation}
For \(R=(x_1,y_1)\), \(S=(x_2,y_2)\), the chord or tangent slope is
\begin{equation}\label{eq:general-slope}
\lambda=
\begin{cases}
\dfrac{y_2-y_1}{x_2-x_1},&R\ne S,\\[3mm]
\dfrac{3x_1^2+2a_2x_1+a_4-a_1y_1}
      {2y_1+a_1x_1+a_3},&R=S.
\end{cases}
\end{equation}
Whenever the indicated denominator is nonzero, put
\[
       \nu=y_1-\lambda x_1.
\]
The inverse and the sum are then
\begin{equation}\label{eq:general-Weierstrass-addition}
 -(x,y)=(x,-y-a_1x-a_3),
\end{equation}
and
\begin{equation}\label{eq:general-Weierstrass-sum}
 \begin{aligned}
 x_3&=\lambda^2+a_1\lambda-a_2-x_1-x_2,\\
 y_3&=-(\lambda+a_1)x_3-\nu-a_3,\\
 R+S&=(x_3,y_3).
 \end{aligned}
\end{equation}
Indeed, the line \(Y=\lambda X+\nu\) has a third intersection with
\eqref{eq:general-W}, and \eqref{eq:general-Weierstrass-addition} reflects
that intersection across the generalized Weierstrass inverse.  Vertical and
boundary cases are interpreted on the smooth projective curve.  These
formulas, rather than short-Weierstrass formulas with \(2\) in a denominator,
are used in low characteristic and in the divisor calculations of the
pairing chapter.

\chapter{The Smooth Biquadratic Model}
\label{ch:geometry}
\section{Projective completion}

The affine equation
\[
             \mathcal C_d:\qquad (u^2+u)(v^2+v)=d
\]
has the following natural bihomogeneous completion in
\(\PP^1_u\times\PP^1_v\):
\begin{equation}\label{eq:homogeneous-thesis}
 U_1(U_1+U_0)V_1(V_1+V_0)=dU_0^2V_0^2,
\end{equation}
where \(u=U_1/U_0\) and \(v=V_1/V_0\).  The four rational boundary points are
\begin{equation}\label{eq:boundary-thesis}
 O=(0,\infty),\quad T=(-1,\infty),\quad
 R=(\infty,0),\quad -R=(\infty,-1).
\end{equation}
They remain distinct as projective points in characteristic two.

\begin{theorem}[Smoothness and genus]\label{thm:smooth-thesis}
The completion \eqref{eq:homogeneous-thesis} is smooth if and only if
\begin{equation}\label{eq:smooth-condition}
          d(1-16d)\ne0.
\end{equation}
When this holds, it is a genus-one curve with rational point \(O\).
Thus the condition is \(d\ne0,1/16\) in odd characteristic and simply
\(d\ne0\) in characteristic two.
\end{theorem}

\begin{proof}
Set \(F=(u^2+u)(v^2+v)-d\).  Its affine partial derivatives are
\[
 F_u=(2u+1)(v^2+v),\qquad
 F_v=(2v+1)(u^2+u).
\]
If \(d\ne0\), neither factor \(u^2+u\) nor \(v^2+v\) vanishes on the
curve.  In characteristic two the partial derivatives become
\(v^2+v\) and \(u^2+u\), so there is no affine singularity.  In odd
characteristic both derivatives can vanish only at
\((u,v)=(-1/2,-1/2)\), which lies on the curve exactly when \(d=1/16\).

It remains to check the boundary.  Near \(u=\infty\), use \(t=1/u\).
The local equation is
\[
       (1+t)(v^2+v)=dt^2.
\]
At \(t=0\), the two points have \(v=0,-1\), and the derivative with
respect to \(v\) is, respectively, \(1\) and \(-1\); both values are
nonzero, including in characteristic two where \(-1=1\).

Near \(v=\infty\), put \(s=1/v\).  After multiplication by \(s^2\), the
local equation becomes
\[
       (u^2+u)(1+s)=ds^2.
\]
At \(s=0\), the two boundary points have \(u=0,-1\), and the derivative
with respect to \(u\) is again \(1\) or \(-1\), hence nonzero in every
characteristic.  Finally, at the projective point
\(((U_0:U_1),(V_0:V_1))=((0:1),(0:1))\), corresponding to
\((u,v)=(\infty,\infty)\), the bihomogeneous equation has left-hand side
\(1\) and right-hand side \(0\).  Hence this point does not lie on the
completion.  Thus no boundary
singularity occurs.  For \(d=0\) the equation is reducible.  This proves
\eqref{eq:smooth-condition}.  A smooth curve of bidegree \((2,2)\) in
\(\PP^1\times\PP^1\) has genus \((2-1)(2-1)=1\), and \(O\) is rational.
\end{proof}

Let \(F\) be a field of arbitrary characteristic.  For \(D\in F\), write
\[
  \C_D:\qquad (x^2+x)(y^2+y)=D.
\]

\begin{proposition}[Explicit isomorphism with a \(\mathcal C_D\) curve]
\label{prop:general-biquadratic-to-Cd}
Let \(a,b,c,d,e\in F\) with
\[
  abcde\ne 0,
\]
and consider the affine biquadratic curve
\[
  \mathcal B_{a,b,c,d,e}:\qquad
  (au^2+bu)(cv^2+dv)=e.
\]
Set
\[
  D=\frac{ace}{b^2d^2}.
\]
Then \(\mathcal B_{a,b,c,d,e}\) is \(F\)-isomorphic to
\[
  \C_D:\qquad
  (x^2+x)(y^2+y)=\frac{ace}{b^2d^2}.
\]
More precisely, the isomorphism and its inverse are
\begin{align*}
  \Phi:\mathcal B_{a,b,c,d,e}&\longrightarrow \C_D,
  & (u,v)&\longmapsto
  \left(\frac{a}{b}u,\frac{c}{d}v\right),\\
  \Phi^{-1}:\C_D&\longrightarrow \mathcal B_{a,b,c,d,e},
  & (x,y)&\longmapsto
  \left(\frac{b}{a}x,\frac{d}{c}y\right).
\end{align*}
These maps are defined over \(F\) and are valid in every characteristic.
\end{proposition}

\begin{proof}
Put
\[
  x=\frac{a}{b}u,
  \qquad
  y=\frac{c}{d}v;
  \qquad\text{equivalently}\qquad
  u=\frac{b}{a}x,
  \quad
  v=\frac{d}{c}y.
\]
Direct substitution gives
\begin{align*}
  au^2+bu
  &=a\left(\frac{b}{a}x\right)^2
    +b\left(\frac{b}{a}x\right)
    =\frac{b^2}{a}(x^2+x),\\
  cv^2+dv
  &=c\left(\frac{d}{c}y\right)^2
    +d\left(\frac{d}{c}y\right)
    =\frac{d^2}{c}(y^2+y).
\end{align*}
Consequently,
\[
  (au^2+bu)(cv^2+dv)
  =\frac{b^2d^2}{ac}(x^2+x)(y^2+y).
\]
Thus the equation \((au^2+bu)(cv^2+dv)=e\) is equivalent to
\[
  (x^2+x)(y^2+y)=\frac{ace}{b^2d^2}=D.
\]
Since \(a,b,c,d\) are nonzero, both coordinate substitutions are invertible
over \(F\).  Hence \(\Phi\) is an \(F\)-isomorphism with the displayed inverse.
The computation uses only multiplication and division by nonzero elements;
in particular, it does not divide by \(2\) or extract square roots.  It
therefore remains valid in characteristics \(2\) and \(3\), as well as in every
other characteristic.
\end{proof}

\begin{remark}[Extension to the natural \((2,2)\)-completions]
Write \(u=U_1/U_0\), \(v=V_1/V_0\), \(x=X_1/X_0\), and \(y=Y_1/Y_0\).
The preceding affine isomorphism is induced by the automorphism of
\(\PP^1\times\PP^1\)
\[
  \bigl((U_0:U_1),(V_0:V_1)\bigr)
  \longmapsto
  \bigl((bU_0:aU_1),(dV_0:cV_1)\bigr).
\]
Its inverse is
\[
  \bigl((X_0:X_1),(Y_0:Y_1)\bigr)
  \longmapsto
  \bigl((aX_0:bX_1),(cY_0:dY_1)\bigr).
\]
Therefore the isomorphism extends from the affine equations to their natural
bidegree-\((2,2)\) projective completions and, whenever necessary, to their
smooth projective models.
\end{remark}

\section{The group law and rational four-torsion}

\begin{proposition}\label{prop:inverse-torsion}
Take \(O=(0,\infty)\) as identity.  Then
\[
        -(u,v)=(u,-v-1).
\]
Moreover \(T\) has order two, \(2R=T\), and
\(\langle R\rangle\simeq\mathbb Z/4\mathbb Z\) is rational over \(k\).
\end{proposition}

\begin{proof}
Write temporarily
\[
        R'=(\infty,-1),
\]
so that the notation \(-R\) is not used before it has been justified by
the group law.  We first compute four principal divisors directly on the
smooth \((2,2)\)-completion.

Near \(R=(\infty,0)\), put \(t=1/u\).  The local equation is
\begin{equation}\label{eq:local-at-R-for-torsion}
             (1+t)(v^2+v)=dt^2.
\end{equation}
Its derivative with respect to \(v\) at \((t,v)=(0,0)\) is \(1\).
Consequently \(t\) is a uniformizer and
\(\operatorname{ord}_R(v)=2\).  At \(R'=(\infty,-1)\), write
\(w=v+1\).  Then \(v^2+v=w^2-w\), and
\eqref{eq:local-at-R-for-torsion} becomes
\[
       (1+t)(w^2-w)=dt^2.
\]
The derivative with respect to \(w\) at \((t,w)=(0,0)\) is \(-1\),
which is nonzero in every characteristic.  Hence \(t\) is also a
uniformizer at \(R'\), and comparison of the lowest-order terms gives
\(w= -dt^2+O(t^3)\), so
\(\operatorname{ord}_{R'}(v+1)=2\).  There are no other zeros of \(v\)
or \(v+1\), since a finite zero of either function would make the
left-hand side of the affine equation vanish while \(d\ne0\).

Near \(O=(0,\infty)\), put \(z=1/v\).  After multiplication by \(z^2\),
the local equation becomes
\begin{equation}\label{eq:local-at-O-for-torsion}
              (u^2+u)(1+z)=dz^2.
\end{equation}
The derivative with respect to \(u\) at \((u,z)=(0,0)\) is \(1\), so
\(z\) is a uniformizer, \(u\) has order two, and \(v=z^{-1}\) has a
simple pole.  Replacing \(u\) by \(-1+w\) gives the identical conclusion
at \(T=(-1,\infty)\): \(u+1\) has order two and \(v\) has a simple
pole.  Finally, \(u=1/t\) has a simple pole at each of \(R,R'\).
We have therefore proved
\begin{align}
 \Div(u)&=2(O)-(R)-(R'),&
 \Div(u+1)&=2(T)-(R)-(R'),\label{eq:u-divisors-four-torsion}\\
 \Div(v)&=2(R)-(O)-(T),&
 \Div(v+1)&=2(R')-(O)-(T).
 \label{eq:v-divisors-four-torsion}
\end{align}
These local calculations remain valid in characteristic two: every linear
coefficient used above is \(1\) or \(-1=1\), and hence is nonzero.

Identify the smooth curve with its Jacobian by
\(P\mapsto[(P)-(O)]\).  A principal divisor maps to the identity.
The first relation in \eqref{eq:u-divisors-four-torsion} gives
\[
                       R+R'=O.
\]
Thus \(R'\) is indeed the group inverse of \(R\).  The second relation
then gives \(2T=O\), while the first relation in
\eqref{eq:v-divisors-four-torsion} gives
\[
                       2R=T.
\]
Since the four boundary points are distinct, \(T\ne O\); hence \(T\)
has exact order two and \(R\) has exact order four.

It remains to identify inversion on an arbitrary point.  The map
\[
              \iota(u,v)=(u,-v-1)
\]
preserves the equation, fixes \(O\), and is the nontrivial deck
transformation of the projection to the \(u\)-line.  For a generic
\(a\in\overline{k}\), the fiber of \(u=a\) consists of
\(P\) and \(\iota(P)\), and
\[
        \Div(u-a)=(P)+(\iota(P))-(R)-(R').
\]
The already proved equality \(R+R'=O\) therefore implies
\(P+\iota(P)=O\).  Hence \(\iota=[-1]\) on a dense open subset.  Both
maps are morphisms of the smooth projective curve, so they agree
everywhere.  This proves
\( -(u,v)=(u,-v-1)\) and all the asserted torsion relations without
using an auxiliary curve model.
\end{proof}

\begin{figure}[H]
\centering
\begin{tikzpicture}[scale=1.0,every node/.style={font=\small},>=Latex]
  \coordinate (O) at (0,1.6);
  \coordinate (R) at (1.8,0);
  \coordinate (T) at (0,-1.6);
  \coordinate (Rm) at (-1.8,0);
  \draw[rounded corners=8pt,gray!40] (-2.4,-2.1) rectangle (2.4,2.1);
  \draw[thick] (O)--(R)--(T)--(Rm)--cycle;
  \foreach \P/\Lab in {O/O,R/R,T/T,Rm/-R}
    \fill (\P) circle (2pt);
  \node[above] at (O) {$O=(0,\infty)$};
  \node[right] at (R) {$R=(\infty,0)$};
  \node[below] at (T) {$T=(-1,\infty)$};
  \node[left] at (Rm) {$-R=(\infty,-1)$};
  \draw[->,thick,blue!70!black] (1.05,1.0) arc[start angle=55,end angle=-35,radius=1.45];
  \draw[->,thick,blue!70!black] (0.8,-0.95) arc[start angle=-35,end angle=-125,radius=1.45];
  \draw[->,thick,blue!70!black] (-0.8,-0.95) arc[start angle=-125,end angle=-215,radius=1.45];
  \draw[->,thick,blue!70!black] (-1.05,1.0) arc[start angle=145,end angle=55,radius=1.45];
  \node[blue!70!black] at (0,0.25) {$+R$};
  \node at (0,2.35) {$2R=T,\qquad 4R=O$};
  \node at (0,-2.35) {$\langle R\rangle=\{O,R,T,-R\}\cong \mathbb Z/4\mathbb Z$};
\end{tikzpicture}
\caption{The four rational boundary points form a cyclic subgroup of order four.  The native inverse fixes the first coordinate and sends $R$ to $-R=(\infty,-1)$, while addition by $R$ cycles through the boundary in the displayed order.}
\label{fig:Cd-boundary-four-torsion}
\end{figure}
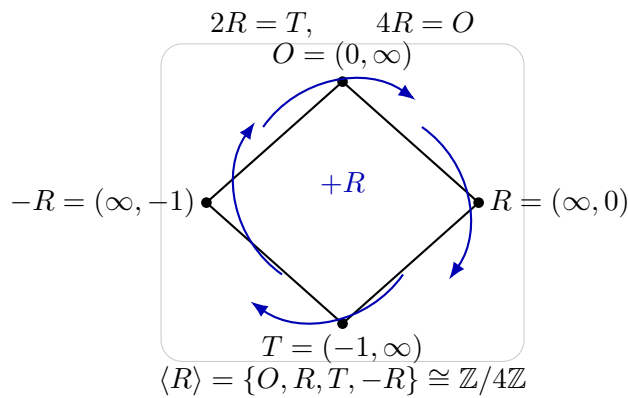

\section[Geometric addition on the (2,2)-completion]
{Geometric addition on the \texorpdfstring{\((2,2)\)}{(2,2)}-completion}
\label{sec:geometric-addition-Cd}

The chord-and-tangent construction for a plane cubic uses a line, whose
intersection number with the cubic is three.  The native completion of
\(\mathcal C_d\) is instead a \((2,2)\)-curve in
\(\PP^1\times\PP^1\).  A curve of bidegree \((1,1)\) meets it in four
points, counted with multiplicity.  The marked point
\(T=(-1,\infty)\) supplies the fourth fixed incidence needed to turn this
four-point intersection into an addition construction.

Let \(F_u\) and \(F_v\) denote the divisor classes of a vertical and a
horizontal fiber, respectively.  The principal divisors already computed in
\eqref{eq:u-divisors-four-torsion}--\eqref{eq:v-divisors-four-torsion}
give
\begin{equation}\label{eq:fiber-classes-geometric-addition}
 F_u\sim(R)+(-R)\sim2(O),\qquad
 F_v\sim(O)+(T).
\end{equation}
Consequently the divisor cut out on \(\overline{\mathcal C}_d\) by a
\((1,1)\)-curve has the fixed class
\begin{equation}\label{eq:one-one-class-Cd}
 H:=F_u+F_v\sim3(O)+(T).
\end{equation}

\begin{theorem}[The \((1,1)\)-intersection construction]
\label{thm:geometric-addition-one-one}
Let \(P,Q\in\overline{\mathcal C}_d(\bar k)\) be in general position.
There is a unique \((1,1)\)-curve \(\Gamma_{P,Q,T}\) through
\(P,Q,T\).  If its fourth intersection with
\(\overline{\mathcal C}_d\) is \(W\), with all intersections counted with
multiplicity, then
\begin{equation}\label{eq:geometric-addition-Cd}
             \boxed{\qquad P+Q=-W.\qquad}
\end{equation}
For \(P=Q\), the two incidence conditions at \(P\) are replaced by a
tangency condition, and the same construction gives \(2P=-W\).
\end{theorem}

\begin{proof}
The projective space of bihomogeneous forms of bidegree \((1,1)\) has
dimension three.  Three independent point conditions therefore determine a
unique member.  B\'ezout's theorem on
\(\PP^1\times\PP^1\) gives
\[
 (1,1)\mathbin{\cdot}(2,2)=1\cdot2+1\cdot2=4,
\]
so its intersection divisor with the native completion is
\[
       (P)+(Q)+(T)+(W)\sim H\sim3(O)+(T).
\]
Passing to \(\operatorname{Pic}^0\) and using \(O\) as the group identity
gives
\[
 (P-O)+(Q-O)+(T-O)+(W-O)=(T-O).
\]
Canceling \((T-O)\) yields \(P+Q+W=O\), which is
\eqref{eq:geometric-addition-Cd}.  When \(P=Q\), the intersection divisor
contains \(2(P)\), and the same divisor-class calculation proves the tangent
statement.
\end{proof}

In homogeneous coordinates, a general \((1,1)\)-curve through
\(T=((1:-1),(0:1))\) can be written
\begin{equation}\label{eq:one-one-through-T}
 aU_1V_1+bU_1V_0+aU_0V_1+eU_0V_0=0.
\end{equation}
On the finite chart this is the bi-affine equation
\[
             a(uv+v)+bu+e=0.
\]
The two remaining linear conditions obtained by substituting \(P\) and
\(Q\) determine \((a:b:e)\) generically.  Substitution of
\eqref{eq:one-one-through-T} into the native biquadratic equation then
determines the fourth intersection \(W\), after which the native inverse
\[
             -W=(u_W,-v_W-1)
\]
produces the sum.  This gives a geometric explanation for the denominator
and exceptional-divisor phenomena in the affine formulas of
Chapter~\ref{ch:odd-full}: a chosen \((1,1)\)-member may cease to be unique
or may meet a boundary point with higher multiplicity, although the group-law
morphism remains defined.

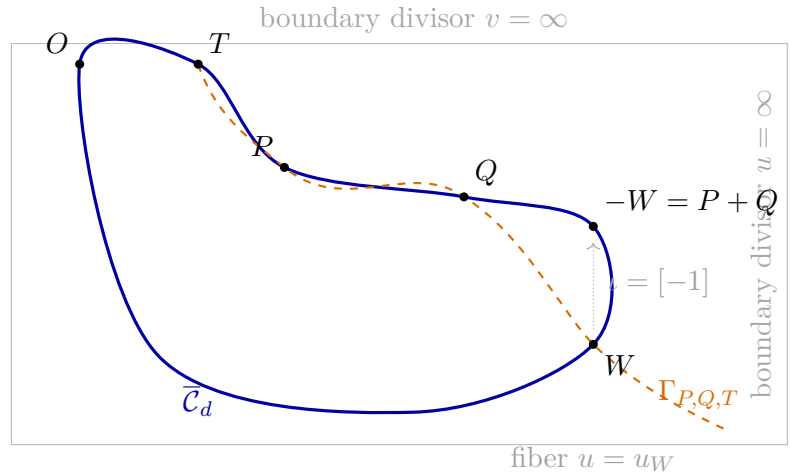
\begin{figure}[H]
\centering
\begin{tikzpicture}[x=0.95cm,y=0.78cm,
                    every node/.style={font=\small}]
  \draw[gray!55] (0.4,0.3) rectangle (11.2,7.1);
  \node[gray!70,above] at (6.0,7.1) {boundary divisor \(v=\infty\)};
  \node[gray!70,rotate=90,above] at (11.2,3.7)
       {boundary divisor \(u=\infty\)};

  \draw[blue!65!black,very thick]
    plot[smooth cycle,tension=0.58] coordinates
    {(1.35,6.75) (3.0,6.75) (4.2,5.0) (6.7,4.5)
     (8.5,4.0) (8.5,2.0) (6.0,0.85) (2.45,1.8)};
  \node[blue!65!black] at (3.0,1.05)
       {\(\overline{\mathcal C}_d\)};

  \draw[orange!85!black,thick,dashed]
    (3.0,6.75)
      .. controls (3.25,5.85) and (3.75,5.45) .. (4.2,5.0)
      .. controls (5.05,4.15) and (5.95,5.15) .. (6.7,4.5)
      .. controls (7.55,3.85) and (8.05,2.55) .. (8.5,2.0)
      .. controls (9.05,1.35) and (9.65,0.95) .. (10.35,0.55);
  \node[orange!85!black,right] at (9.25,1.15)
       {\(\Gamma_{P,Q,T}\)};

  \draw[gray!75,densely dotted,->]
       (8.5,2.25)--(8.5,3.72);
  \node[gray!75,right] at (8.55,3.05)
       {\(\iota=[-1]\)};
  \node[gray!75,below] at (8.5,0.45) {fiber \(u=u_W\)};

  \foreach \x/\y in {1.35/6.75,3.0/6.75,4.2/5.0,6.7/4.5,8.5/2.0,8.5/4.0}
       \fill[black] (\x,\y) circle (1.7pt);
  \node[above left]  at (1.35,6.75) {\(O\)};
  \node[above right] at (3.0,6.75)  {\(T\)};
  \node[above left]  at (4.2,5.0)   {\(P\)};
  \node[above right] at (6.7,4.5)   {\(Q\)};
  \node[below right] at (8.5,2.0)   {\(W\)};
  \node[above right] at (8.5,4.0)   {\(-W=P+Q\)};
\end{tikzpicture}
\caption{Schematic \((1,1)\)-intersection construction of the group law.
The dashed curve meets \(\overline{\mathcal C}_d\) at
\(P,Q,T,W\); applying the native inverse to the fourth intersection gives
\(P+Q=-W\).  The drawing represents incidence in
\(\PP^1\times\PP^1\), not a metrically accurate real affine locus.}
\label{fig:geometric-addition-Cd}
\end{figure}

Figure~\ref{fig:geometric-addition-Cd} is the native analogue of the
chord-and-tangent picture for a plane cubic.  The replacement of a line by a
\((1,1)\)-curve and the required passage through \(T\) are dictated by
the embedding class \eqref{eq:one-one-class-Cd}; they are not artifacts of
an Edwards or Weierstrass transformation.  For exceptional triples one uses
another member of the addition-law atlas from
Chapter~\ref{ch:native-completeness}.  Thus the picture explains the generic
law, while the atlas supplies global completeness.

\section{Intrinsic dihedral symmetry}

The factorized equation has more structure than the single involution used
for inversion.  Define
\begin{equation}\label{eq:dihedral-generators-Cd}
\begin{aligned}
 \iota_u(u,v)&=(-u-1,v),\\
 \iota_v(u,v)&=(u,-v-1),\\
 \sigma(u,v)&=(v,u).
\end{aligned}
\end{equation}

\begin{theorem}[The intrinsic dihedral subgroup]
\label{thm:Cd-dihedral-symmetry}
For every characteristic and every smooth parameter \(d\), the maps in
\eqref{eq:dihedral-generators-Cd} extend to automorphisms of the smooth
\((2,2)\)-completion and generate a subgroup
\[
          \langle\iota_u,\iota_v,\sigma\rangle\simeq D_8
\]
of order eight.  Their action on the marked boundary is
\[
\begin{array}{c|cccc}
 &O&T&R&-R\\ \hline
\iota_u&T&O&R&-R\\
\iota_v&O&T&-R&R\\
\sigma&R&-R&O&T
\end{array}
\]
and is faithful.

If \(\charac k\ne2\), let \(\tau_Q(P)=P+Q\).  Relative to the group law
with identity \(O\), these same automorphisms satisfy
\begin{equation}\label{eq:dihedral-group-interpretation}
 \iota_v=[-1],\qquad
 \iota_u=\tau_T\circ[-1],\qquad
 \sigma=\tau_R\circ[-1],
\end{equation}
and consequently
\[
        \iota_u\iota_v=\tau_T,\qquad
        \sigma\iota_v=\tau_R.
\]
\end{theorem}

\begin{proof}
The polynomial \(z^2+z\) is invariant under \(z\mapsto-z-1\), and the
equation is symmetric in its two factors.  Thus all three maps preserve the
affine equation.  In homogeneous coordinates, \(\iota_u\) is induced by
\[
 (U_0:U_1)\longmapsto(U_0:-U_1-U_0),
\]
\(\iota_v\) is induced by
\[
 (V_0:V_1)\longmapsto(V_0:-V_1-V_0),
\]
and
\(\sigma\) exchanges the two \(\PP^1\)-factors.  They therefore extend
across the boundary and restrict to automorphisms of the smooth
completion.

Direct composition gives
\[
 \iota_u^2=\iota_v^2=\sigma^2=1,\qquad
 \iota_u\iota_v=\iota_v\iota_u,\qquad
 \sigma\iota_u\sigma=\iota_v.
\]
Hence the generated group is a quotient of
\((C_2\times C_2)\rtimes C_2\), where the last involution exchanges the
two factors.  Evaluating the maps at the four boundary points gives the
displayed table.  In permutation notation the three rows are
\[
 (O\,T),\qquad (R\,{-R}),\qquad (O\,R)(T\,{-R}).
\]
They generate eight distinct permutations, so no additional relation is
present and the group is \(D_8\).  Since the induced boundary
permutations of these eight group elements are pairwise distinct, the
kernel of the action on the boundary is trivial.  Hence the action is
faithful.

It remains to identify the group-theoretic meaning in odd characteristic.
Use the centered reciprocal coordinates
\[
       (\xi,\eta)=\bigl((2v+1)^{-1},(2u+1)^{-1}\bigr).
\]
Then
\[
 \iota_v(\xi,\eta)=(-\xi,\eta),\qquad
 \iota_u(\xi,\eta)=(\xi,-\eta),\qquad
 \sigma(\xi,\eta)=(\eta,\xi),
\]
while
\[
 O=(0,1),\qquad T=(0,-1),\qquad R=(1,0).
\]
The Edwards addition law gives
\[
       P+T=(-\xi,-\eta),\qquad P+R=(\eta,-\xi),
\]
and inversion is \((\xi,\eta)\mapsto(-\xi,\eta)\).  Applying these three
identities first to \(P\) and then to \(-P\) yields
\eqref{eq:dihedral-group-interpretation} and the two translation
identities.
\end{proof}

\begin{remark}
The theorem identifies an intrinsic symmetry subgroup of the marked
factorized model, not necessarily the full automorphism group of the
abstract elliptic curve.  At the CM parameters \(j=0\) and \(j=1728\) the
abstract automorphism group can be larger.  The \(D_8\)-action is present
for every smooth \(d\) and explains why inversion, translation by the
marked two-torsion point, translation by a marked four-torsion point, and
exchange of the two factors remain visible in the same coordinates.
\end{remark}

\begin{example}[A full \(D_8\)-orbit over \(\F_{101}\)]
On \(\mathcal C_1\), the point \(P=(6,42)\) has the eight-point orbit
\[
\begin{aligned}
\mathcal O_{D_8}(P)=\{&
(6,42),(94,42),(6,58),(94,58),\\
&(42,6),(42,94),(58,6),(58,94)\}.
\end{aligned}
\]
Every listed pair satisfies
\((u^2+u)(v^2+v)=1\pmod {101}\).  In particular,
\[
 \iota_u\iota_v(P)=(94,58)=P+T,\qquad
 \sigma\iota_v(P)=(58,6)=P+R.
\]
Thus the coordinate involutions exhibit both marked translations without
leaving the \(\mathcal C_1\)-coordinates.
\end{example}

\section{The native Kummer coordinate}

\begin{theorem}[Classification of degree-one Kummer coordinates]
\label{thm:kummer-class-thesis}
For every smooth member,
\[
       k(\Cd)^{[-1]}=k(u).
\]
Consequently every degree-one generator of the Kummer function field is
\[
       \frac{au+b}{cu+e},\qquad ae-bc\ne0.
\]
Every normalization sending \(O\) to infinity and \(T\) to zero has the
form
\[
       t_c=c\frac{u+1}{u},\qquad c\in k^\times.
\]
After fixing the target scale by choosing \(c=1\), its
homogeneous representative is
\begin{equation}\label{eq:native-kummer-thesis}
       \kappa(P)=(X:Z)=(u+1:u).
\end{equation}
\end{theorem}

\begin{proof}
Let \(K=k(\Cd)\) be the function field of the smooth completion.  The
first projection
\[
       \pi_u:\overline{\mathcal C}_d\longrightarrow\PP^1,
       \qquad (u,v)\longmapsto u,
\]
has degree two.  Indeed, a generic vertical fiber of the equation of
bidegree \((2,2)\) consists of the two roots of
\begin{equation}\label{eq:quadratic-over-ku}
        Z^2+Z-\frac{d}{u^2+u}=0.
\end{equation}
Equivalently, \([K:k(u)]=2\).  The substitution
\(v\mapsto-v-1\) interchanges the two roots and is nontrivial, so
\(K/k(u)\) is a quadratic Galois extension once separability has been
checked.

In characteristic two, the derivative of
\eqref{eq:quadratic-over-ku} with respect to \(Z\) is \(1\).  In odd
characteristic it is \(2Z+1\), which is not the zero element of
\(K\); otherwise \(v=-1/2\) would be constant and the relation
\(v^2+v=d/(u^2+u)\) would force the nonconstant function \(u^2+u\) to
be constant.  Thus the extension is separable in every characteristic.
By Proposition~\ref{prop:inverse-torsion}, its nontrivial automorphism is
exactly inversion.  The fixed-field theorem for a degree-two Galois
extension now gives
\[
       K^{[-1]}=K^{\langle v\mapsto-v-1\rangle}=k(u).
\]

The quotient \(\overline{\mathcal C}_d/\{\pm1\}\) is therefore the
projective line with rational coordinate \(u\).  If \(t\) is another
degree-one generator of the same rational function field, then the maps
\(u,t:\PP^1\to\PP^1\) both have degree one.  Hence their transition map
is an automorphism of \(\PP^1\), and consequently
\[
                  t=\frac{au+b}{cu+e},
                  \qquad ae-bc\ne0.
\]
Conversely every such fractional linear function generates \(k(u)\), so
this lists all degree-one Kummer coordinates.

For the marked normalization, Proposition~\ref{prop:inverse-torsion} and
\eqref{eq:u-divisors-four-torsion} give
\[
                   u(O)=0,\qquad u(T)=-1.
\]
A degree-one function on the quotient having a pole at the image of \(O\)
and a zero at the image of \(T\) must therefore be
\[
                   c\frac{u+1}{u},\qquad c\in k^\times.
\]
Choosing \(c=1\) yields the homogeneous coordinate
\[
                       (X:Z)=(u+1:u).
\]
Multiplying both entries by the same nonzero scalar does not change a
projective representative.  Changing \(c\) changes the automorphism of the target \(\PP^1\).  Thus
the normalization \eqref{eq:native-kummer-thesis} is selected by the
additional choice \(c=1\), whereas common projective scaling leaves the
represented point unchanged.
\end{proof}

\begin{remark}
The word ``native'' is geometric here.  It means that the coordinate is
determined by the fixed field and the marked boundary divisor.  Montgomery
and \(Z/4\mathbb Z\) coordinates later appear as linear bases of this same
projective line.
\end{remark}

\section{Square geometry and invariant algebra}
\label{sec:Cd-square-symmetry}

Theorem~\ref{thm:Cd-dihedral-symmetry} constructs the intrinsic dihedral
subgroup from three elementary maps and identifies its group-law meaning.
We now make two further features explicit.  First, away from characteristic
two this action is literally the full rigid-motion group of a square whose
center, in the original affine coordinates, is
\((-1/2,-1/2)\).  Second, in every characteristic the same action has a
simple algebraic quotient generated by a characteristic-uniform invariant.
The square interpretation therefore has an algebraic counterpart in the invariant theory of the \(D_8\)-action.

Write
\[
             \Phi(z)=z^2+z,
             \qquad
             \mathcal C_d:\ \Phi(u)\Phi(v)=d,
\]
and let
\[
       G=\langle\iota_u,\iota_v,\sigma\rangle\simeq D_8,
\]
where \(D_8\) denotes the dihedral group of order eight.  Products of maps
below are composed from right to left.

\begin{proposition}[The eight natural transformations and the true center]
\label{prop:Cd-eight-square-transformations}
The elements of \(G\) are exactly the following eight affine
transformations:
\begin{equation}\label{eq:Cd-eight-square-transformations}
\begin{aligned}
1:\quad &(u,v)\longmapsto (u,v),\\
\sigma:\quad &(u,v)\longmapsto (v,u),\\
\iota_u:\quad &(u,v)\longmapsto (-1-u,v),\\
\sigma\iota_u:\quad &(u,v)\longmapsto (v,-1-u),\\
\iota_v:\quad &(u,v)\longmapsto (u,-1-v),\\
\sigma\iota_v:\quad &(u,v)\longmapsto (-1-v,u),\\
\iota_u\iota_v:\quad &(u,v)\longmapsto (-1-u,-1-v),\\
\sigma\iota_u\iota_v:\quad &(u,v)\longmapsto (-1-v,-1-u).
\end{aligned}
\end{equation}
They preserve \(\mathcal C_d\) and extend to its smooth
\((2,2)\)-completion in every characteristic.

Suppose that \(\charac k\ne2\), and put
\begin{equation}\label{eq:Cd-square-centered-coordinates}
                 r=2u+1,\qquad s=2v+1.
\end{equation}
Then the eight maps become the eight signed permutations of \((r,s)\), as
shown in Table~\ref{tab:Cd-square-transformations}.  Over \(\mathbb R\)
they are precisely the four rotations and four reflections of a square
centered at \((r,s)=(0,0)\).  Consequently their common geometric center in
the original \((u,v)\)-plane is
\begin{equation}\label{eq:Cd-true-geometric-center}
                         c=\left(-\frac12,-\frac12\right).
\end{equation}
The point \(c\) lies on the affine equation only when \(d=1/16\), which is
exactly the singular odd-characteristic parameter.  Hence the center is not
a point of any smooth member of the family.
\end{proposition}

\begin{table}[H]
\centering
\small
\setlength{\tabcolsep}{3.5pt}
\renewcommand{\arraystretch}{1.12}
\caption{The eight natural transformations as square symmetries}
\label{tab:Cd-square-transformations}
\begin{tabularx}{\textwidth}{@{}L{1.55cm}L{3.35cm}L{2.55cm}X@{}}
\toprule
element & action on \((u,v)\) & action on \((r,s)\) & Euclidean meaning\\
\midrule
\(1\) & \((u,v)\) & \((r,s)\) & identity\\
\(\sigma\) & \((v,u)\) & \((s,r)\) & reflection in \(r=s\), i.e.\ \(u=v\)\\
\(\iota_u\) & \((-1-u,v)\) & \((-r,s)\) & reflection in \(r=0\), i.e.\ \(u=-1/2\)\\
\(\sigma\iota_u\) & \((v,-1-u)\) & \((s,-r)\) & clockwise quarter-turn\\
\(\iota_v\) & \((u,-1-v)\) & \((r,-s)\) & reflection in \(s=0\), i.e.\ \(v=-1/2\)\\
\(\sigma\iota_v\) & \((-1-v,u)\) & \((-s,r)\) & counterclockwise quarter-turn\\
\(\iota_u\iota_v\) & \((-1-u,-1-v)\) & \((-r,-s)\) & half-turn\\
\(\sigma\iota_u\iota_v\) & \((-1-v,-1-u)\) & \((-s,-r)\) & reflection in \(r+s=0\), i.e.\ \(u+v=-1\)\\
\bottomrule
\end{tabularx}
\end{table}

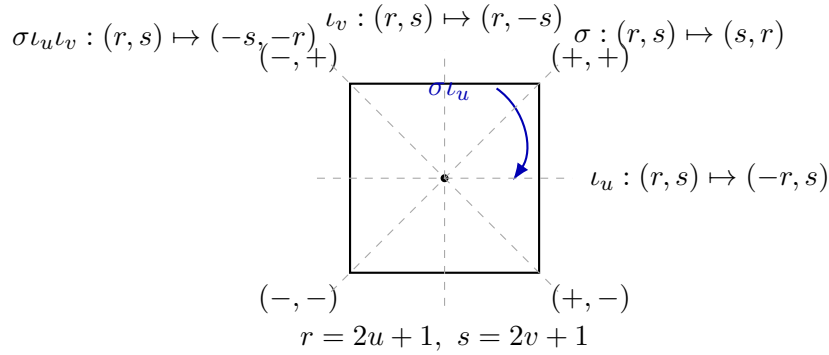
\begin{figure}[H]
\centering
\begin{tikzpicture}[scale=1.25,every node/.style={font=\small},>=Latex]
  \draw[thick] (-1,-1)--(1,-1)--(1,1)--(-1,1)--cycle;
  \fill (0,0) circle (1.2pt);
  \node[below left] at (-1,-1) {$(-,-)$};
  \node[below right] at (1,-1) {$(+,-)$};
  \node[above right] at (1,1) {$(+,+)$};
  \node[above left] at (-1,1) {$(-,+)$};
  \draw[->,thick,blue!70!black] (0.55,0.95) arc[start angle=55,end angle=-35,radius=0.7];
  \node[blue!70!black] at (0.05,0.92) {$\sigma\iota_u$};
  \draw[dashed,gray!70] (-1.35,0)--(1.35,0);
  \draw[dashed,gray!70] (0,-1.35)--(0,1.35);
  \draw[dashed,gray!70] (-1.2,-1.2)--(1.2,1.2);
  \draw[dashed,gray!70] (-1.2,1.2)--(1.2,-1.2);
  \node[right] at (1.42,0) {$\iota_u:(r,s)\mapsto(-r,s)$};
  \node[above] at (0,1.42) {$\iota_v:(r,s)\mapsto(r,-s)$};
  \node[above right] at (1.25,1.25) {$\sigma:(r,s)\mapsto(s,r)$};
  \node[above left] at (-1.23,1.22) {$\sigma\iota_u\iota_v:(r,s)\mapsto(-s,-r)$};
  \node at (0,-1.7) {$r=2u+1,\ s=2v+1$};
\end{tikzpicture}
\caption{In centered coordinates $(r,s)=(2u+1,2v+1)$, the eight natural automorphisms of $\mathcal C_d$ become the rigid symmetries of a square.  This is the geometric meaning of the intrinsic dihedral subgroup $D_8$.}
\label{fig:Cd-D8-square}
\end{figure}

\begin{proof}
The elementary identity
\begin{equation}\label{eq:Phi-reflection-invariance}
                 \Phi(-1-z)=(-1-z)^2+(-1-z)=z^2+z=\Phi(z)
\end{equation}
holds over every field.  It proves invariance under \(\iota_u\) and
\(\iota_v\), while symmetry of the product proves invariance under
\(\sigma\).  The extension to the completion was established in
Theorem~\ref{thm:Cd-dihedral-symmetry}.  Taking all products of the three
generators gives the eight maps in
\eqref{eq:Cd-eight-square-transformations}; their distinctness follows
either from their signed-permutation actions when \(2\ne0\), or from their
faithful boundary action in arbitrary characteristic.

Under \eqref{eq:Cd-square-centered-coordinates}, the maps
\(u\mapsto-1-u\), \(v\mapsto-1-v\), and \((u,v)\mapsto(v,u)\) become
\[
       (r,s)\mapsto(-r,s),\qquad
       (r,s)\mapsto(r,-s),\qquad
       (r,s)\mapsto(s,r),
\]
respectively.  Their products are all signed permutation matrices of rank
two.  These matrices form the Weyl group \(W(B_2)\), equivalently the full
orthogonal symmetry group of a square.  More concretely, if
\[
       \varrho=\sigma\iota_u,\qquad \varepsilon=\sigma,
\]
then
\[
       \varrho(r,s)=(s,-r),\qquad
       \varepsilon(r,s)=(s,r),
\]
and
\[
       \varrho^4=\varepsilon^2=1,
       \qquad
       \varepsilon\varrho\varepsilon=\varrho^{-1}.
\]
This is the standard presentation of the symmetry group of a square.

A point fixed by both coordinate reflections must satisfy
\(2u+1=2v+1=0\).  When \(2\) is invertible this has the unique solution
\eqref{eq:Cd-true-geometric-center}.  Substitution gives
\[
       \Phi(-1/2)^2=(-1/4)^2=1/16,
\]
so the center lies on \(\mathcal C_d\) precisely for \(d=1/16\).
Theorem~\ref{thm:smooth-thesis} identifies this as the singular parameter.
In characteristic two, a common affine fixed point would have to satisfy
\(u+1=u\), which would imply \(1=0\).  Thus the eight algebraic
automorphisms persist, but the affine-center description does not; the
coordinate reflections become Artin--Schreier translations.
\end{proof}

\subsection{Centered factorization and a scaled square family}
\label{subsec:Cd-centered-factorization}

Translation to the true center found in Proposition~\ref{prop:Cd-eight-square-transformations}
gives the following factorized equation.  The next result separates
the translation to that center from an optional uniform dilation.  In this
form the four linear factors are the four signed distances from the two
pairs of parallel lines bounding the reference square.

\begin{proposition}[Centered and scaled square forms]
\label{prop:Cd-centered-scaled-square-forms}
Assume that \(\charac k\ne2\).  Under the translation
\begin{equation}\label{eq:Cd-half-centered-translation}
             x=u+\frac12,\qquad y=v+\frac12,
\end{equation}
the curve \(\mathcal C_d\) has the equivalent equations
\begin{align}
 \left(x^2-\frac14\right)
 \left(y^2-\frac14\right)&=d,
 \label{eq:Cd-half-centered-quadratic-form}\\
 \left(x-\frac12\right)\left(x+\frac12\right)
 \left(y-\frac12\right)\left(y+\frac12\right)&=d.
 \label{eq:Cd-half-centered-linear-form}
\end{align}
Its eight natural symmetries become the signed permutations
\begin{equation}\label{eq:Cd-centered-eight-signed-permutations}
\begin{aligned}
1:\quad &(x,y)\longmapsto (x,y),\\
\sigma:\quad &(x,y)\longmapsto (y,x),\\
\iota_u:\quad &(x,y)\longmapsto (-x,y),\\
\sigma\iota_u:\quad &(x,y)\longmapsto (y,-x),\\
\iota_v:\quad &(x,y)\longmapsto (x,-y),\\
\sigma\iota_v:\quad &(x,y)\longmapsto (-y,x),\\
\iota_u\iota_v:\quad &(x,y)\longmapsto (-x,-y),\\
\sigma\iota_u\iota_v:\quad &(x,y)\longmapsto (-y,-x).
\end{aligned}
\end{equation}

More generally, let \(a\in k^\times\) and consider
\begin{equation}\label{eq:scaled-square-family}
 \mathcal Q_{a,d}:\qquad (x^2-a^2)(y^2-a^2)=d.
\end{equation}
Then \(\mathcal Q_{a,d}\) is \(k\)-isomorphic to
\begin{equation}\label{eq:scaled-square-associated-Cd}
 \mathcal C_D:\qquad (U^2+U)(V^2+V)=D,
 \qquad D=\frac{d}{16a^4},
\end{equation}
by
\begin{equation}\label{eq:scaled-square-to-Cd-map}
 U=\frac{x}{2a}-\frac12,\qquad
 V=\frac{y}{2a}-\frac12,
\end{equation}
with inverse
\begin{equation}\label{eq:Cd-to-scaled-square-map}
 x=a(2U+1),\qquad y=a(2V+1).
\end{equation}
In particular, \(\mathcal Q_{a,d}\) is smooth precisely when
\begin{equation}\label{eq:scaled-square-smoothness}
                         d(a^4-d)\ne0.
\end{equation}
For \(a=1/2\), one has \(D=d\), and
\eqref{eq:scaled-square-family} is exactly the centered presentation
\eqref{eq:Cd-half-centered-quadratic-form} of \(\mathcal C_d\).
\end{proposition}

\begin{proof}
From \eqref{eq:Cd-half-centered-translation} one has
\[
 u=x-\frac12,\qquad v=y-\frac12.
\]
Completing the square, without introducing any rational map or exceptional
affine locus, gives
\[
 u^2+u=\left(x-\frac12\right)^2+
        \left(x-\frac12\right)=x^2-\frac14,
\]
and
\[
 v^2+v=\left(y-\frac12\right)^2+
        \left(y-\frac12\right)=y^2-\frac14.
\]
Substitution into the defining equation
of \(\mathcal C_d\) proves
\eqref{eq:Cd-half-centered-quadratic-form}; factoring each difference of
squares proves \eqref{eq:Cd-half-centered-linear-form}.

The involution \(u\mapsto-1-u\) sends
\[
 x=u+\frac12\longmapsto -1-u+\frac12=-x,
\]
The involution \(v\mapsto-1-v\) sends
\[
 y=v+\frac12\longmapsto -1-v+\frac12=-y.
\]
Coordinate
exchange sends \((x,y)\) to \((y,x)\).  Composing these three operations
from right to left gives exactly the eight maps in
\eqref{eq:Cd-centered-eight-signed-permutations}.  Each map preserves the
unordered pair \(\{x^2,y^2\}\), so it preserves
\((x^2-1/4)(y^2-1/4)\).

For the scaled family, substitute the inverse formulas
\eqref{eq:Cd-to-scaled-square-map}.  A direct calculation gives
\[
 x^2-a^2
 =a^2\bigl((2U+1)^2-1\bigr)
 =4a^2(U^2+U),
\]
and likewise \(y^2-a^2=4a^2(V^2+V)\).  Hence
\[
 (x^2-a^2)(y^2-a^2)
 =16a^4(U^2+U)(V^2+V).
\]
Because \(2a\ne0\), formulas
\eqref{eq:scaled-square-to-Cd-map} and
\eqref{eq:Cd-to-scaled-square-map} are mutually inverse affine linear
changes of coordinates.  Dividing the last identity by \(16a^4\) proves
\eqref{eq:scaled-square-associated-Cd} and the asserted isomorphism of the
affine curves; the same linear changes homogenize to an isomorphism of
their smooth \((2,2)\)-completions.

By Theorem~\ref{thm:smooth-thesis}, \(\mathcal C_D\) is smooth exactly when
\(D(1-16D)\ne0\).  With \(D=d/(16a^4)\) and \(a\ne0\), this condition is
equivalent to
\[
 \frac{d}{16a^4}\left(1-\frac d{a^4}\right)\ne0,
\]
which is precisely \eqref{eq:scaled-square-smoothness}.  Finally, setting
\(a=1/2\) gives \(16a^4=1\), so \(D=d\) and the displayed specialization
follows.
\end{proof}

\begin{remark}[Scale-free parameter]
If \(X=x/a\) and \(Y=y/a\), then
\begin{equation}\label{eq:scaled-square-normalized-form}
                 (X^2-1)(Y^2-1)=\frac d{a^4}.
\end{equation}
Thus the intrinsic dimensionless parameter of the scaled square family is
\(d/a^4\).  The form \eqref{eq:scaled-square-family} is geometrically
useful because the four distinguished lines are visibly
\(x=\pm a\) and \(y=\pm a\), while the original \(\mathcal C_d\) equation
remains the characteristic-uniform model: the centered and scaled formulas
in this subsection require \(2a\in k^\times\).
\end{remark}

The square is also visible on the marked boundary.  The quarter-turn
\(\varrho=\sigma\iota_u\) acts by the four-cycle
\begin{equation}\label{eq:Cd-boundary-square-cycle}
                 O\longmapsto -R\longmapsto T
                 \longmapsto R\longmapsto O.
\end{equation}
Hence the four boundary points are the four vertices of the intrinsic
square, in the precise sense of permutation representations.

\begin{figure}[H]
\centering
\begin{tikzpicture}[scale=1.0,>=Latex,every node/.style={font=\small}]
  \coordinate (O)  at (0,2.05);
  \coordinate (mR) at (2.05,0);
  \coordinate (T)  at (0,-2.05);
  \coordinate (R)  at (-2.05,0);

  \draw[curveblue,very thick,->] (O)--(mR);
  \draw[curveblue,very thick,->] (mR)--(T);
  \draw[curveblue,very thick,->] (T)--(R);
  \draw[curveblue,very thick,->] (R)--(O);
  \draw[softgray,densely dashed] (-2.45,2.45)--(2.45,-2.45);

  \fill[curvered] (O) circle (2.2pt) node[above=4pt] {\(O\)};
  \fill[curvered] (mR) circle (2.2pt) node[right=5pt] {\(-R\)};
  \fill[curvered] (T) circle (2.2pt) node[below=4pt] {\(T\)};
  \fill[curvered] (R) circle (2.2pt) node[left=5pt] {\(R\)};

  \node[curveblue,fill=white,inner sep=2pt] at (1.72,1.72)
       {\(\varrho=\sigma\iota_u\)};
  \node[softgray,fill=white,inner sep=2pt,rotate=-45] at (-1.55,1.55)
       {reflection axis for \(\sigma\)};
\end{tikzpicture}
\caption{The exact permutation action of the quarter-turn
\(\varrho=\sigma\iota_u\) on the four marked boundary points.  The square
is a representation of the boundary action, not an affine placement of
these points in \(\PP^1\times\PP^1\).}
\label{fig:Cd-boundary-square}
\end{figure}
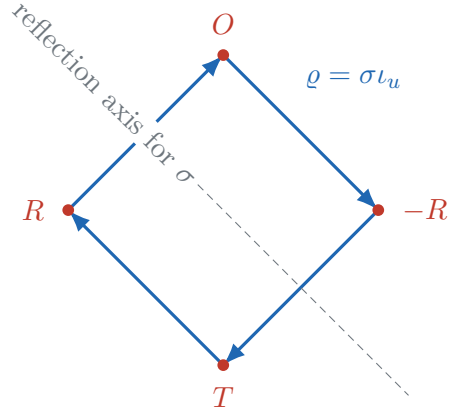

\begin{theorem}[The invariant field of the square action]
\label{thm:Cd-square-invariant-field}
Let \(K=k(\overline{\mathcal C}_d)\) be the function field of a smooth
member, and let \(G\simeq D_8\) be the group above.  Define
\begin{equation}\label{eq:Cd-full-square-invariant}
       J=\Phi(u)+\Phi(v)
        =u^2+u+v^2+v.
\end{equation}
Then
\begin{equation}\label{eq:Cd-square-fixed-field}
                         K^G=k(J).
\end{equation}
Equivalently, the quotient map by all eight natural symmetries is the
degree-eight map
\[
       \pi_G:\overline{\mathcal C}_d\longrightarrow\PP^1_J,
       \qquad P\longmapsto J(P).
\]

If \(\charac k\ne2\), put \(\rho=1-16d\) and use the centered coordinates
\eqref{eq:Cd-square-centered-coordinates}.  Then
\begin{equation}\label{eq:Cd-square-invariants-centered}
  \mathcal C_d:\quad r^2s^2-r^2-s^2+\rho=0,
  \qquad
  J=\frac{r^2+s^2-2}{4}.
\end{equation}
Moreover,
\begin{equation}\label{eq:B2-polynomial-invariants}
       k[r,s]^G=k[p,q],\qquad
       p=r^2+s^2,\qquad q=r^2s^2,
\end{equation}
and restriction to \(\mathcal C_d\) imposes the single relation
\(q-p+\rho=0\).  Thus the affine invariant coordinate ring is a polynomial
ring in \(p\), consistently with \(K^G=k(J)\).
\end{theorem}

\begin{proof}
Set
\[
                  A=\Phi(u),\qquad B=\Phi(v).
\]
The curve relation is \(AB=d\), and hence \(B=d/A\) in \(K\).  Let
\[
                  H=\langle\iota_u,\iota_v\rangle
                    \simeq C_2\times C_2.
\]
The four elements of \(H\) fix \(A\) and \(B\).  Conversely, over
\(k(A)\) the two native coordinates satisfy
\[
             u^2+u=A,
             \qquad
             v^2+v=\frac d A.
\]
Therefore \([K:k(A)]\le4\).  The four distinct automorphisms in \(H\)
fix \(k(A)\), while Artin's fixed-field theorem gives
\([K:K^H]=|H|=4\).  It follows that
\begin{equation}\label{eq:Cd-first-square-quotient}
                         K^H=k(A).
\end{equation}
This argument is separable in every characteristic: in characteristic two
the two quadratics are Artin--Schreier equations with derivative \(1\), and
in odd characteristic their derivatives are nonzero rational functions.

The remaining involution induced by \(\sigma\) exchanges \(A\) and \(B\),
so on \(k(A)\) it acts as
\[
                         A\longmapsto\frac d A.
\]
Its invariant \(A+d/A\) is exactly \(J\).  Since \(A\) satisfies
\begin{equation}\label{eq:Cd-square-quotient-quadratic}
                         X^2-JX+d=0
\end{equation}
over \(k(J)\), the extension \(k(A)/k(J)\) has degree at most two.  The
involution is nontrivial, so its fixed field has index two and must be
\(k(J)\).  Combining this with
\eqref{eq:Cd-first-square-quotient} proves
\eqref{eq:Cd-square-fixed-field}.  It also exhibits the quotient as the
two-stage tower
\begin{equation}\label{eq:Cd-two-stage-square-quotient}
 \overline{\mathcal C}_d
 \longrightarrow \{AB=d\}
 \longrightarrow \PP^1_J,
 \qquad
 \deg=4\ \text{followed by}\ \deg=2.
\end{equation}

Assume now that \(2\) is invertible.  Since
\[
       \Phi(u)=\frac{r^2-1}{4},
       \qquad
       \Phi(v)=\frac{s^2-1}{4},
\]
the native equation and \(J\) become
\eqref{eq:Cd-square-invariants-centered}.  A polynomial fixed by the two
independent sign changes is a polynomial in \(r^2,s^2\), and invariance
under exchange makes it a symmetric polynomial in those two variables.
The fundamental theorem of symmetric polynomials therefore gives
\eqref{eq:B2-polynomial-invariants}.  On the curve,
\(q-p+\rho=0\).  Since \(|G|=8\) is invertible in \(k\), averaging over
\(G\) shows that taking invariants commutes with this invariant quotient;
hence
\[
 \bigl(k[r,s]/(q-p+\rho)\bigr)^G
       \simeq k[p,q]/(q-p+\rho)
       \simeq k[p].
\]
Finally \(p=4J+2\), so this polynomial quotient agrees with the
function-field description.
\end{proof}

\begin{corollary}[Generic eight-point orbits]
\label{cor:Cd-generic-square-orbits}
Over an algebraic closure, every unramified fibre of \(\pi_G\) is one
eight-element \(G\)-orbit.  An orbit has fewer than eight points exactly
when its stabilizer is nontrivial; geometrically these exceptional points
occur on the reflection loci or on their projective boundary
specializations.
\end{corollary}

\begin{proof}
Artin's theorem and Theorem~\ref{thm:Cd-square-invariant-field} give
\([K:k(J)]=|G|=8\) with automorphism group \(G\).  Away from the ramification
locus, a fibre therefore consists of eight distinct points permuted simply
transitively by \(G\).  The orbit--stabilizer formula gives the final
statement.
\end{proof}

\begin{remark}[Geometric and algebraic meanings of the same symmetry]
In odd characteristic, \(p=r^2+s^2\) is the algebraic analogue of squared
distance from the true center, and the equation expresses the second basic
square invariant \(q=r^2s^2\) as \(q=p-\rho\).  In characteristic two
there is no affine center and no nontrivial sign change, but
\(A=\Phi(u)\) and \(B=\Phi(v)\) remain the Artin--Schreier invariants of the
two translations.  The same formula \(J=A+B\) therefore describes the full
dihedral quotient in every characteristic.  This is the algebraic reason
that the visible square symmetry survives even when its Euclidean picture
does not.
\end{remark}

The equation
\[
\mathcal C_d:\qquad (u^2+u)(v^2+v)=d
\]
is the pullback of the hyperbola
\[
XY=d
\]
under the two quadratic maps
\[
X=u^2+u,\qquad Y=v^2+v.
\]
The two deck transformations
\[
u\longmapsto -1-u,\qquad
v\longmapsto -1-v,
\]
together with interchange of the two factors, generate the intrinsic
\(D_8\)-action described above.  In odd characteristic, after the translation
\[
r=2u+1,\qquad s=2v+1,
\]
this action becomes the group of signed permutations of \((r,s)\), hence the
full symmetry group of a square centered at
\((-1/2,-1/2)\) in the original affine coordinates.  In characteristic two,
the same automorphisms are expressed by the two Artin--Schreier translations
and factor exchange.

The \(D_8\)-action is reflected in several structures used later in the
monograph.  It permutes the marked boundary, contains the native inverse,
determines the quotient invariant
\[
J=(u^2+u)+(v^2+v),
\qquad
k(\mathcal C_d)^{D_8}=k(J),
\]
and is compatible with the native Kummer quotient and the arithmetic
formulas derived from the two quadratic factors.  Thus the geometric square
symmetry and the invariant-theoretic description give two realizations of
the same automorphism structure on \(\mathcal C_d\).

Proposition~\ref{prop:general-biquadratic-to-Cd} gives the exact coverage of
the corresponding factorized normal form.  More generally, let \(e\in F^\times\) and
\[
       Q_1(x)Q_2(y)=e
\]
with \(Q_1,Q_2\in F[T]\) separable quadratic polynomials.  If each
\(Q_i\) has a chosen simple root in \(F\), translate those roots to zero.
Then
\[
       Q_1(x_0+u)=a u^2+b u,
       \qquad
       Q_2(y_0+v)=c v^2+d v
\]
with \(abcd\ne0\), and Proposition~\ref{prop:general-biquadratic-to-Cd}
identifies the resulting curve with an explicitly determined
\(\mathcal C_D\) over \(F\).  If one of the two quadratic factors has no
\(F\)-rational root, such a factor-preserving affine reduction need not
exist over \(F\); it becomes available after a splitting extension.  Thus
\(\mathcal C_d\) is a characteristic-uniform normal form for the split,
rationally rooted factorized subclass, rather than for every equation
\(Q_1(x)Q_2(y)=e\) over an arbitrary ground field.

\subsection{Visible eight-point orbits}
\label{subsec:Cd-visible-eight-point-orbits}

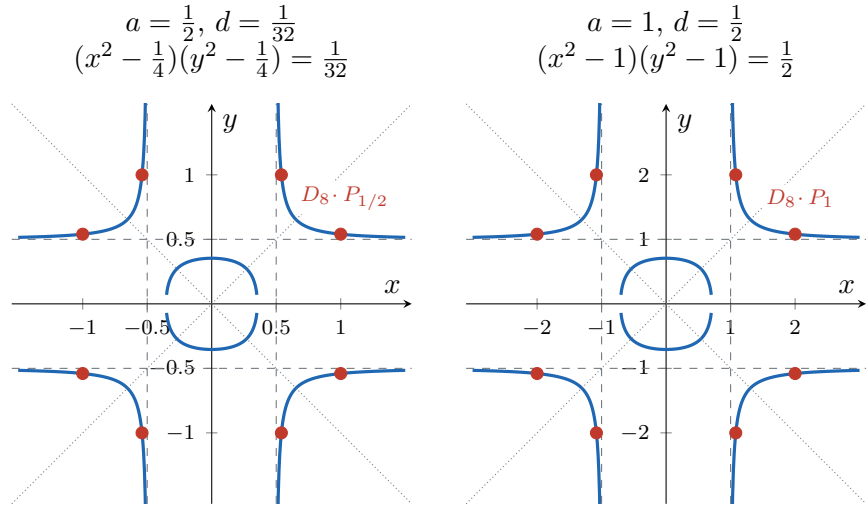
\begin{figure}[H]
\centering
\begin{tikzpicture}
\begin{groupplot}[
  group style={group size=2 by 1,horizontal sep=0.72cm},
  width=0.455\textwidth,
  height=0.455\textwidth,
  axis equal image,
  axis lines=middle,
  xlabel={$x$},
  ylabel={$y$},
  samples=180,
  clip=true,
  tick label style={font=\scriptsize},
  label style={font=\small},
  title style={font=\small,align=center}
]

\nextgroupplot[
  xmin=-1.55,xmax=1.55,
  ymin=-1.55,ymax=1.55,
  xtick={-1,-0.5,0,0.5,1},
  ytick={-1,-0.5,0,0.5,1},
  title={$a=\tfrac12$, $d=\tfrac1{32}$\\
         $(x^2-\tfrac14)(y^2-\tfrac14)=\tfrac1{32}$}
]
  \foreach \xa/\xb in {-1.5/-0.505,-0.35/0.35,0.505/1.5}{
    \addplot[curveblue,very thick,domain=\xa:\xb]
      {sqrt(0.25+0.03125/(x^2-0.25))};
    \addplot[curveblue,very thick,domain=\xa:\xb]
      {-sqrt(0.25+0.03125/(x^2-0.25))};
  }
  \addplot[softgray,dashed] coordinates {(-0.5,-1.55) (-0.5,1.55)};
  \addplot[softgray,dashed] coordinates {(0.5,-1.55) (0.5,1.55)};
  \addplot[softgray,dashed] coordinates {(-1.55,-0.5) (1.55,-0.5)};
  \addplot[softgray,dashed] coordinates {(-1.55,0.5) (1.55,0.5)};
  \addplot[softgray,densely dotted] coordinates {(-1.55,-1.55) (1.55,1.55)};
  \addplot[softgray,densely dotted] coordinates {(-1.55,1.55) (1.55,-1.55)};
  \addplot[only marks,mark=*,mark size=2.25pt,curvered]
    coordinates {
      (1,0.540061725) (1,-0.540061725)
      (-1,0.540061725) (-1,-0.540061725)
      (0.540061725,1) (-0.540061725,1)
      (0.540061725,-1) (-0.540061725,-1)
    };
  \node[fill=white,inner sep=1.5pt,text=curvered,font=\scriptsize]
    at (axis cs:1.03,0.82) {$D_8\!\cdot P_{1/2}$};

\nextgroupplot[
  xmin=-3.1,xmax=3.1,
  ymin=-3.1,ymax=3.1,
  xtick={-2,-1,0,1,2},
  ytick={-2,-1,0,1,2},
  title={$a=1$, $d=\tfrac12$\\
         $(x^2-1)(y^2-1)=\tfrac12$}
]
  \foreach \xa/\xb in {-3/-1.01,-0.70/0.70,1.01/3}{
    \addplot[curveblue,very thick,domain=\xa:\xb]
      {sqrt(1+0.5/(x^2-1))};
    \addplot[curveblue,very thick,domain=\xa:\xb]
      {-sqrt(1+0.5/(x^2-1))};
  }
  \addplot[softgray,dashed] coordinates {(-1,-3.1) (-1,3.1)};
  \addplot[softgray,dashed] coordinates {(1,-3.1) (1,3.1)};
  \addplot[softgray,dashed] coordinates {(-3.1,-1) (3.1,-1)};
  \addplot[softgray,dashed] coordinates {(-3.1,1) (3.1,1)};
  \addplot[softgray,densely dotted] coordinates {(-3.1,-3.1) (3.1,3.1)};
  \addplot[softgray,densely dotted] coordinates {(-3.1,3.1) (3.1,-3.1)};
  \addplot[only marks,mark=*,mark size=2.25pt,curvered]
    coordinates {
      (2,1.080123450) (2,-1.080123450)
      (-2,1.080123450) (-2,-1.080123450)
      (1.080123450,2) (-1.080123450,2)
      (1.080123450,-2) (-1.080123450,-2)
    };
  \node[fill=white,inner sep=1.5pt,text=curvered,font=\scriptsize]
    at (axis cs:2.08,1.64) {$D_8\!\cdot P_1$};

\end{groupplot}
\end{tikzpicture}
\caption{Two centered real members of the scaled square family and their
exact eight-point orbits from
\eqref{eq:scaled-square-exact-eight-point-orbit}.  The dashed lines are
\(x=\pm a\) and \(y=\pm a\), and the dotted diagonals are the other two
reflection axes.  The right-hand panel is the dilation by a factor of two
of the left-hand panel.  Decimal values are used only to place the red
markers; their exact coordinates are given by
\eqref{eq:scaled-square-exact-eight-point-orbit}.}
\label{fig:Cd-centered-eight-point-orbits}
\end{figure}

The scaled family also supplies an exact one-parameter construction of a
generic real orbit.  For \(a>0\), take \(d=a^4/2\) and
\begin{equation}\label{eq:scaled-square-distinguished-real-point}
             P_a=\left(2a,a\sqrt{\frac76}\right).
\end{equation}
Indeed,
\[
 \bigl((2a)^2-a^2\bigr)
 \left(\left(a\sqrt{\frac76}\right)^2-a^2\right)
 =3a^2\cdot\frac{a^2}{6}=\frac{a^4}{2}.
\]
Because \(2a\ne a\sqrt{7/6}\ne0\), none of the signed-coordinate
reflections or coordinate exchanges fixes \(P_a\).  Its orbit therefore
has all eight points
\begin{equation}\label{eq:scaled-square-exact-eight-point-orbit}
 D_8\!\cdot P_a
 =\left\{
   \left(\pm2a,\,\pm a\sqrt{\frac76}\right),
   \left(\pm a\sqrt{\frac76},\,\pm2a\right)
  \right\},
\end{equation}
where the two signs in each ordered pair are independent.  Figure~
\ref{fig:Cd-centered-eight-point-orbits} plots this exact orbit for
\(a=1/2\) and for its dilation \(a=1\).

\chapter{Odd-Characteristic Coordinate Dictionaries}
\label{ch:odd-dictionary}
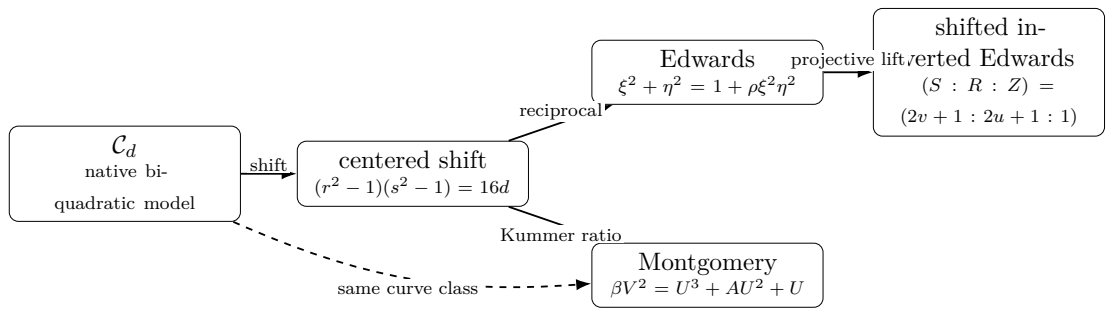
\begin{figure}[H]
\centering
\resizebox{0.96\textwidth}{!}{%
\begin{tikzpicture}[x=1cm,y=1cm,>=Latex, box/.style={draw,rounded corners,align=center,inner sep=4pt,font=\small,text width=3.25cm}]
  \node[box] (Cd) at (0,0) {$\mathcal C_d$\\[-1mm]\scriptsize native biquadratic model};
  \node[box] (shift) at (4.4,0) {centered shift\\[-1mm]\scriptsize $(r^2-1)(s^2-1)=16d$};
  \node[box] (Ed) at (8.9,1.55) {Edwards\\[-1mm]\scriptsize $\xi^2+\eta^2=1+\rho\xi^2\eta^2$};
  \node[box] (InvEd) at (13.2,1.55) {shifted inverted Edwards\\[-1mm]\scriptsize $(S:R:Z)=(2v+1:2u+1:1)$};
  \node[box] (Mont) at (8.9,-1.55) {Montgomery\\[-1mm]\scriptsize $\beta V^2=U^3+AU^2+U$};
  \draw[->,thick] (Cd) -- node[above,font=\scriptsize,fill=white,inner sep=1pt]{shift} (shift);
  \draw[->,thick] (shift) -- node[above,font=\scriptsize,fill=white,inner sep=1pt]{reciprocal} (Ed);
  \draw[->,thick] (Ed) -- node[above,font=\scriptsize,fill=white,inner sep=1pt]{projective lift} (InvEd);
  \draw[->,thick] (shift) -- node[below,font=\scriptsize,fill=white,inner sep=1pt]{Kummer ratio} (Mont);
  \draw[->,dashed,thick] (Cd) to[bend right=14] node[below,font=\scriptsize,fill=white,inner sep=1pt]{same curve class} (Mont);
\end{tikzpicture}%
}
\caption{The odd-characteristic dictionaries used throughout the monograph.  The arrows denote explicit coordinate interfaces on dense charts.  The dashed arrow records equivalence of abstract elliptic curves, not equality of arithmetic models.}
\label{fig:Cd-odd-dictionary-map}
\end{figure}

Throughout this chapter, \(\charac k\ne2\) and \(d\ne0,1/16\).  We begin
from the native affine equation
\begin{equation}\label{eq:odd-dictionary-native-start}
 \boxed{\qquad
 \mathcal C_d:\quad (u^2+u)(v^2+v)=d.
 \qquad}
\end{equation}
The restrictions on \(d\) are exactly the smoothness restrictions from
Chapter~\ref{ch:geometry}.  Since \(2\in k^\times\), the polynomial shift
\begin{equation}\label{eq:shift-dictionary}
       r=2u+1,\qquad s=2v+1,\qquad D=16d,\qquad \rho=1-D
\end{equation}
is invertible, with \(u=(r-1)/2\) and \(v=(s-1)/2\).  Moreover
\[
 u^2+u=\frac{r^2-1}{4},\qquad
 v^2+v=\frac{s^2-1}{4}.
\]
Consequently \eqref{eq:odd-dictionary-native-start} becomes
\begin{equation}\label{eq:shifted-biquadratic}
  (r^2-1)(s^2-1)=D,\qquad
  r^2+s^2=r^2s^2+\rho.
\end{equation}
The first equality is obtained by multiplying the native equation by
\(16\); expanding it and using \(\rho=1-D\) gives the second.  Thus the
shift itself discards no affine point.  Denominators enter only in the
reciprocal Edwards and affine Montgomery charts introduced below.  Notice
also that smoothness gives \(D\rho\ne0\).  Figure~\ref{fig:Cd-odd-dictionary-map}
summarizes the interfaces that will be used repeatedly in the odd-characteristic
arithmetic chapters.

\section{Edwards and inverted Edwards}

\begin{theorem}[Edwards identification]\label{thm:edwards-id-thesis}
The rational map
\begin{equation}\label{eq:edwards-map-thesis}
 (\xi,\eta)=\left(\frac1s,\frac1r\right)
 =\left(\frac1{2v+1},\frac1{2u+1}\right)
\end{equation}
extends uniquely to an origin-preserving isomorphism from the smooth
completion of \(\Cd\) to the smooth completion of the Edwards affine model
\begin{equation}\label{eq:edwards-thesis}
 \Erho:\qquad \xi^2+\eta^2=1+\rho\xi^2\eta^2.
\end{equation}
\end{theorem}

\begin{proof}
Divide the second equation in \eqref{eq:shifted-biquadratic} by
\(r^2s^2\):
\[
  s^{-2}+r^{-2}=1+\rho r^{-2}s^{-2}.
\]
This is \eqref{eq:edwards-thesis}.  On \(\xi\eta\ne0\), the inverse is
\[
 u=\frac{\eta^{-1}-1}{2},\qquad
 v=\frac{\xi^{-1}-1}{2}.
\]
The curves are therefore birational.  A birational map between smooth
projective curves extends uniquely to an isomorphism.  At the native identity
\(O=(0,\infty)\), one has \(r=1\), while \(s\) has a pole; hence
\((\xi,\eta)=(0,1)\), the Edwards identity.  The extended isomorphism is
therefore origin-preserving.
\end{proof}

\begin{proposition}[The shifted inverted-Edwards chart]
\label{prop:shifted-inverted-Edwards-dictionary}
On the finite native chart, set
\begin{equation}\label{eq:shifted-inverted-Edwards-lift}
        (S:R:Z)=(s:r:1)=(2v+1:2u+1:1).
\end{equation}
Then
\begin{equation}\label{eq:shifted-inverted-Edwards-equation}
       (S^2+R^2)Z^2=S^2R^2+\rho Z^4.
\end{equation}
On the chart \(Z\ne0\), the native coordinates are recovered by
\begin{equation}\label{eq:shifted-inverted-Edwards-recovery}
       u=\frac{R-Z}{2Z},\qquad
       v=\frac{S-Z}{2Z}.
\end{equation}
The lift \eqref{eq:shifted-inverted-Edwards-lift} is division-free, but the
single \(Z\ne0\) chart is not the whole smooth completion.
\end{proposition}

\begin{proof}
In inverted Edwards notation,
\((\xi,\eta)=(Z/S,Z/R)\).  Multiplying
\(\xi^2+\eta^2=1+\rho\xi^2\eta^2\) by \(S^2R^2\) gives
\eqref{eq:shifted-inverted-Edwards-equation}.  Alternatively, substituting
\((S:R:Z)=(s:r:1)\) reduces that equation directly to the second identity
in \eqref{eq:shifted-biquadratic}.  Finally \(r=R/Z\) and \(s=S/Z\), so
the inverse shift gives \eqref{eq:shifted-inverted-Edwards-recovery}.
Points with \(Z=0\) require the boundary charts of the smooth completion,
which is why division-free entry does not imply completeness of this one
computational chart.
\end{proof}

\begin{corollary}[Model-theoretic distinction from the Edwards presentation]
\label{cor:model-not-class}
In odd characteristic, \(\Cd\) is a distinct biquadratic elliptic-curve model
that admits a shifted inverted-Edwards dictionary.  The dictionary does not
identify the two models: the \((2,2)\)-embedding, boundary marking,
characteristic-free equation, native functions \(u,v\), and native Kummer
map of \(\mathcal C_d\) differ from the usual Edwards presentation.
\end{corollary}

\begin{proof}
Inverted Edwards coordinates \((S:R:Z)\) represent
\((\xi,\eta)=(Z/S,Z/R)\).  Taking \((S:R:Z)=(s:r:1)\) yields
the chart of Proposition~\ref{prop:shifted-inverted-Edwards-dictionary} and
recovers \eqref{eq:edwards-map-thesis}.  The coordinate shift and scaling by \(2\)
do not alter the principal multiplication graph, but the ambient projective
embeddings and boundary divisors are different.  Moreover
\eqref{eq:edwards-map-thesis} is unavailable in characteristic two while
\eqref{eq:model} remains smooth.  Thus equivalence of abstract curves does
not identify the two arithmetic models.  This is the same distinction by
which Edwards, Montgomery, Hessian, Jacobi quartic, and Weierstrass equations
are recognized as different elliptic-curve models even when explicit
isomorphisms exist between their smooth completions.
\end{proof}

\section{Montgomery form}

\begin{theorem}[Explicit Montgomery isomorphism]\label{thm:mont-thesis}
Set
\begin{equation}\label{eq:mont-parameters-thesis}
 \beta=\frac1{16d},\qquad A=\frac1{4d}-2.
\end{equation}
Then \(\Cd\) is isomorphic over \(k\), with the isomorphism preserving the
origins, to
\begin{equation}\label{eq:mont-thesis}
 \MAB:\qquad \beta V^2=U^3+AU^2+U
\end{equation}
through
\begin{equation}\label{eq:mont-map-thesis}
 U=\frac{r+1}{r-1}=\frac{u+1}{u},\qquad
 V=2sU.
\end{equation}
The inverse on the dense affine chart is
\begin{equation}\label{eq:mont-inverse-thesis}
 u=\frac1{U-1},\qquad v=\frac{V}{4U}-\frac12.
\end{equation}
\end{theorem}

\begin{proof}
Solve \(r=(U+1)/(U-1)\) and \(s=V/(2U)\).  Then
\[
 r^2-1=\frac{4U}{(U-1)^2},\qquad
 s^2-1=\frac{V^2-4U^2}{4U^2}.
\]
Substitution into \((r^2-1)(s^2-1)=16d\) gives
\[
 V^2=16d\,U^3+(4-32d)U^2+16d\,U.
\]
Division by \(16d\) proves \eqref{eq:mont-thesis}.  Also
\[
 A^2-4=\frac{1-16d}{16d^2}=\frac{\rho}{16d^2}\ne0,
\]
so the Montgomery curve is nonsingular.  The inverse formulas are obtained
by solving \eqref{eq:mont-map-thesis}; smoothness then extends the birational
map to the projective completions.  At \(O=(0,\infty)\), the function
\(U=(u+1)/u\) has a pole, so the image is the unique Montgomery point at
infinity.  Thus the extended isomorphism preserves the chosen identities.
\end{proof}

\begin{remark}[Full-point information and the native quotient]
The first Montgomery coordinate is not an arbitrary imported abscissa:
\[
       U=\frac{u+1}{u}
\]
is invariant under the native inverse
\((u,v)\mapsto(u,-v-1)\) and is the affine ratio of the native Kummer map
\((u+1:u)\) from Chapter~\ref{ch:geometry}.  The coordinate \(V=2sU\)
restores the sign information on a dense chart.  Thus \((U,V)\) is a
full-point dictionary, whereas retaining only \(U\) gives a degree-two
quotient output of type \(\Ktwo\).
\end{remark}

\section{Dictionary domains and recovery}

For later arithmetic it is useful to separate an isomorphism of smooth
completions from the affine formula by which it is evaluated.  The relevant
dense-chart data are
\begin{table}[H]
\centering
\small
\setlength{\tabcolsep}{3.5pt}
\begin{tabular}{L{2.3cm}L{3.2cm}L{4.4cm}L{2.8cm}}
\toprule
presentation & forward coordinates & native recovery & affine restriction\\
\midrule
centered shift
 & \((r,s)=(2u+1,2v+1)\)
 & \(\displaystyle
    \begin{gathered}u=(r-1)/2,\\v=(s-1)/2\end{gathered}\)
 & none\\
Edwards
 & \((\xi,\eta)=(s^{-1},r^{-1})\)
 & \(\displaystyle
    \begin{gathered}u=(\eta^{-1}-1)/2,\\v=(\xi^{-1}-1)/2\end{gathered}\)
 & \(rs\ne0\)\\
inverted Edwards
 & \((S:R:Z)=(s:r:1)\)
 & \(\displaystyle
    \begin{gathered}u=(R-Z)/(2Z),\\v=(S-Z)/(2Z)\end{gathered}\)
 & \(Z\ne0\) for recovery\\
Montgomery
 & \((U,V)=((u+1)/u,\,2sU)\)
 & \(\displaystyle
    \begin{gathered}u=1/(U-1),\\v=V/(4U)-1/2\end{gathered}\)
 & no additional restriction on finite \(\mathcal C_d\)-points\\
\bottomrule
\end{tabular}
\caption{Odd-characteristic coordinate dictionaries on their dense affine charts}
\label{tab:odd-dictionary-domains}
\end{table}
The restrictions in the last column describe only the displayed affine
inverse formulas.  They are not singularity conditions on the curve: each
rational dictionary extends across the missing points after passing to the
appropriate smooth completion.  In particular, an operation count obtained
in one row must include recovery to \((u,v)\) if a native full-point output
is required, whereas a Kummer algorithm deliberately retains only the
projective ratio \((u+1:u)\).
For the Montgomery row, the absence of an additional restriction follows from
\(d\ne0\): equation~\eqref{eq:odd-dictionary-native-start} forces
\(u(u+1)\ne0\), and hence \(U\ne0,1\).

The endpoint costs also explain the working coordinates used later.  The
centered shift and the inverted lift require only additions.  Computing both
Edwards reciprocals separately costs \(2\Inv\), or \(\Inv+3\M\) by one
simultaneous inversion.  A full Montgomery image can be formed as
\(U=1+u^{-1}\), \(V=2sU\), at cost \(\Inv+\M\), whereas the quotient input
\((u+1:u)\) costs no multiplication or inversion.  These are endpoint costs,
not recurring costs of a projective scalar-multiplication loop.

\section{Invariants and completeness}

\begin{proposition}\label{prop:j-thesis}
The \(j\)-invariant is
\begin{equation}\label{eq:j-thesis}
 j(\Cd)=16\frac{(\rho^2+14\rho+1)^3}{\rho(1-\rho)^4}.
\end{equation}
Over a finite field of odd characteristic, the Edwards affine addition law
on \(\Erho\) is complete when \(\rho\) is a nonsquare.
\end{proposition}

\begin{proof}
Substitute \(A=(4d)^{-1}-2\) into the Montgomery invariant
\[
 j=256\frac{(A^2-3)^3}{A^2-4}
\]
and use \(d=(1-\rho)/16\).  Then
\[
 A=\frac{2(1+\rho)}{1-\rho},\qquad
 A^2-3=\frac{\rho^2+14\rho+1}{(1-\rho)^2},\qquad
 A^2-4=\frac{16\rho}{(1-\rho)^2}.
\]
Substitution gives
\[
 256\frac{(A^2-3)^3}{A^2-4}
 =16\frac{(\rho^2+14\rho+1)^3}
          {\rho(1-\rho)^4},
\]
which is \eqref{eq:j-thesis}.  The completeness criterion for
\(a\xi^2+\eta^2=1+d_E\xi^2\eta^2\) over a finite field is that \(a\) be a
square and \(d_E\) a nonsquare.  Here \(a=1\) and \(d_E=\rho\).
\end{proof}

\begin{remark}
Completeness of the Edwards affine law does not imply completeness of the
inverted coordinate chart: the latter omits points with \(\xi\eta=0\),
including the identity and boundary torsion.
\end{remark}

\part{Native Arithmetic in Odd Characteristic}
\partoverview{This part starts and ends on
\(\mathcal C_d:(u^2+u)(v^2+v)=d\).  It first records addition and doubling in
the original \(u,v\) coordinates, then introduces native shifted-projective
coordinates as an internal linear recoding.  A separate completeness chapter
then constructs the native Segre law and the odd and binary differential
atlases before the native Kummer line, its reciprocal affine chart,
arithmetic tradeoffs, halving, recovery, tripling, closed \(2P+Q\), and
characteristic-three Frobenius
specialization are developed.  The \(j=1728\) and \(j=0\) endomorphisms and
their Kummer-compatible GLV actions receive their own chapter.}

\chapter[Native Odd-Characteristic Arithmetic]
{Native Full-Point Arithmetic on \(\mathcal C_d\) in Odd Characteristic}
\label{ch:odd-full}
\section[Native affine group-law formulas]
{Native affine group-law formulas on \(\mathcal C_d\)}

Throughout this chapter,
\[
 \mathcal C_d:\qquad (u^2+u)(v^2+v)=d,\qquad
 \charac k\ne2,\qquad d(1-16d)\ne0.
\]
The group law is recorded first as formulas, independently of any
implementation schedule.  With
\[
 O=(0,\infty),\qquad -(u,v)=(u,-v-1),\qquad
 r_i=2u_i+1,\quad s_i=2v_i+1,\quad \rho=1-16d,
\]
the sum \(P_3=P_1+P_2=(u_3,v_3)\) is
\begin{equation}\label{eq:native-affine-add-thesis}
\boxed{
\begin{aligned}
u_3&=
\frac{r_1r_2s_1s_2-\rho-s_1s_2+r_1r_2}
     {2(s_1s_2-r_1r_2)},\\
v_3&=
\frac{r_1r_2s_1s_2+\rho-s_1r_2-r_1s_2}
     {2(s_1r_2+r_1s_2)}.
\end{aligned}}
\end{equation}
Expanding every shifted symbol gives the same formula entirely in the
original coordinates:
\begin{equation*}\tag{\thechapter.1-uv}\label{eq:native-affine-add-uv}
\boxed{
\begin{aligned}
u_3={}&
\frac{
(2u_1+1)(2u_2+1)(2v_1+1)(2v_2+1)-(1-16d)}
{2\bigl((2v_1+1)(2v_2+1)-(2u_1+1)(2u_2+1)\bigr)}
\\[-1mm]
&\quad+
\frac{(2u_1+1)(2u_2+1)-(2v_1+1)(2v_2+1)}
{2\bigl((2v_1+1)(2v_2+1)-(2u_1+1)(2u_2+1)\bigr)},\\
v_3={}&
\frac{
(2u_1+1)(2u_2+1)(2v_1+1)(2v_2+1)+(1-16d)}
{2\bigl((2v_1+1)(2u_2+1)+(2u_1+1)(2v_2+1)\bigr)}
\\[-1mm]
&\quad-
\frac{(2v_1+1)(2u_2+1)+(2u_1+1)(2v_2+1)}
{2\bigl((2v_1+1)(2u_2+1)+(2u_1+1)(2v_2+1)\bigr)}.
\end{aligned}}
\end{equation*}
The numerators in formula~\eqref{eq:native-affine-add-uv} are split only for
readability; each coordinate is a single rational function after combining the two
terms.  In particular, no coordinate belonging to an isomorphic curve
occurs in this version.
Equivalently,
\begin{equation}\label{eq:native-shift-add-thesis}
\boxed{\qquad
 r_3=\frac{r_1r_2s_1s_2-\rho}{s_1s_2-r_1r_2},
 \qquad
 s_3=\frac{r_1r_2s_1s_2+\rho}{s_1r_2+r_1s_2}.
\qquad}
\end{equation}
For \(P_1=P_2=P=(u,v)\), with \(r=2u+1\) and \(s=2v+1\),
\begin{equation}\label{eq:native-affine-dbl-thesis}
\boxed{
\begin{aligned}
 r_{2P}&=\frac{r^2s^2-\rho}{s^2-r^2},&
 s_{2P}&=\frac{r^2s^2+\rho}{2rs},\\
 u_{2P}&=\frac{r^2s^2-\rho-s^2+r^2}
                 {2(s^2-r^2)},&
 v_{2P}&=\frac{r^2s^2+\rho-2rs}{4rs}.
\end{aligned}}
\end{equation}
Without the abbreviations \(r,s\), this becomes
\begin{equation*}\tag{\thechapter.3-uv}\label{eq:native-affine-double-uv}
\boxed{
\begin{aligned}
u_{2P}
&=\frac{(2u+1)^2(2v+1)^2-(1-16d)
        -(2v+1)^2+(2u+1)^2}
       {2\bigl((2v+1)^2-(2u+1)^2\bigr)},\\
v_{2P}
&=\frac{(2u+1)^2(2v+1)^2+(1-16d)
        -2(2u+1)(2v+1)}
       {4(2u+1)(2v+1)}.
\end{aligned}}
\end{equation*}
These are equations in the native coordinates \(u,v\); no Edwards,
Montgomery, or Weierstrass coordinate is part of their statement.

\medskip
\noindent\textbf{Ordinary native projective form.}
For implementations that begin with the usual one-scale homogenization,
write
\[
 \mathbf P=(U:V:Z),\qquad u=U/Z,\quad v=V/Z,
\]
so that the native plane chart is
\begin{equation}\label{eq:ordinary-native-projective-curve}
 U(U+Z)V(V+Z)=dZ^4.
\end{equation}
Let
\(\nu_{\rm pl}:\overline{\mathcal C}_d\to C_d^{\rm pl}\)
denote the normalization map to this plane quartic.  It is an isomorphism
over \(Z\ne0\), but its two singular points at infinity each identify two
branches of the smooth \((2,2)\)-completion.  Thus a nonzero tuple with
\(Z=0\) records \(\nu_{\rm pl}(P)\), not the branch data needed to recover
the full point \(P\).
For \(\mathbf P_i=(U_i:V_i:Z_i)\), define the centered linear forms
\[
 R_i=2U_i+Z_i,\qquad S_i=2V_i+Z_i
\]
and
\[
\begin{aligned}
A&=Z_1Z_2,&B&=\rho A^2,&C&=S_1S_2,&D&=R_1R_2,\\
E&=CD,&H&=C-D,&I&=S_1R_2+R_1S_2.
\end{aligned}
\]
Then, whenever the displayed output tuple is nonzero, an ordinary projective
representative of \(\nu_{\rm pl}(P_1+P_2)\) is
\begin{equation}\label{eq:ordinary-native-projective-add}
\boxed{
\begin{aligned}
U_3&=(E-B)I-AHI,\\
V_3&=(E+B)H-AHI,\\
Z_3&=2AHI.
\end{aligned}}
\end{equation}
For doubling, put
\[
\begin{aligned}
R&=2U+Z,&S&=2V+Z,\\
A_0&=S^2,&B_0&=R^2,&C_0&=A_0+B_0,\\
D_0&=A_0-B_0,&E_0&=(S+R)^2-C_0,
\end{aligned}
\]
and
\[
 \widehat S_2=C_0D_0,\qquad
 \widehat R_2=E_0(C_0-2\rho Z^2),\qquad
 \widehat Z_2=D_0E_0.
\]
Whenever it is nonzero, the ordinary projective image of the double is
\begin{equation}\label{eq:ordinary-native-projective-double}
\boxed{\quad
(U_2:V_2:Z_2)=
(\widehat R_2-\widehat Z_2:
 \widehat S_2-\widehat Z_2:
 2\widehat Z_2).
\quad}
\end{equation}
The direct addition costs \(9\M+\Sqr+\Dpar\), its mixed
projective--affine specialization costs \(8\M+\Sqr+\Dpar\), and
\eqref{eq:ordinary-native-projective-double} costs
\(3\M+4\Sqr+\Dpar\).  These counts exclude additions and multiplication
by the small integer \(2\).  Equation
\eqref{eq:ordinary-native-projective-curve} is an arithmetic chart; the
smooth completion remains the native \((2,2)\)-curve in
\(\PP^1\times\PP^1\).  Consequently these tuples have output type
\(\Ffin\) when their final scale is nonzero; a boundary result must be
retained in, or lifted to, the native Segre completion before it is called
a \(\Full\)-output.

\section{Derivation and exceptional affine divisors}

\begin{theorem}[Native affine addition]\label{thm:native-affine-add-thesis}
For \(P_i=(u_i,v_i)\in\mathcal C_d(k)\), put
\[
\begin{aligned}
 C&=s_1s_2,&D&=r_1r_2,&E&=CD,\\
 H&=C-D,&I&=s_1r_2+r_1s_2.
\end{aligned}
\]
If \(HI\ne0\), formulas
\eqref{eq:native-affine-add-thesis} and
\eqref{eq:native-shift-add-thesis} give the group-theoretic sum
\(P_1+P_2\) on the smooth completion of \(\mathcal C_d\).
With two separate inversions the literal affine schedule costs
\(6\M+2\Inv\).  With Montgomery's simultaneous inversion of \(H\) and
\(I\), it costs \(9\M+\Inv\).
\end{theorem}

\begin{proof}
By Theorem~\ref{thm:edwards-id-thesis}, the map
\[
 \xi_i=s_i^{-1},\qquad \eta_i=r_i^{-1}
\]
is an isomorphism from \(\mathcal C_d\) to
\(\xi^2+\eta^2=1+\rho\xi^2\eta^2\).  The two numerators in the Edwards
addition law simplify as
\[
\begin{aligned}
\xi_1\eta_2+\eta_1\xi_2
 &=\frac{s_1r_2+r_1s_2}{r_1r_2s_1s_2}
   =\frac I E,\\
\eta_1\eta_2-\xi_1\xi_2
 &=\frac{s_1s_2-r_1r_2}{r_1r_2s_1s_2}
   =\frac H E,
\end{aligned}
\]
while
\[
 \rho\xi_1\xi_2\eta_1\eta_2=\frac{\rho}{E}.
\]
Consequently
\[
 \xi_3=\frac{\xi_1\eta_2+\eta_1\xi_2}
              {1+\rho\xi_1\xi_2\eta_1\eta_2},\qquad
 \eta_3=\frac{\eta_1\eta_2-\xi_1\xi_2}
              {1-\rho\xi_1\xi_2\eta_1\eta_2}
\]
become
\[
       \xi_3=\frac I{E+\rho}=\frac1{s_3},\qquad
       \eta_3=\frac H{E-\rho}=\frac1{r_3}.
\]
Thus \eqref{eq:native-shift-add-thesis} is the pullback of the group law.
Curve preservation alone would not establish this identification.  Finally
\(u_3=(r_3-1)/2\) and
\(v_3=(s_3-1)/2\), giving \eqref{eq:native-affine-add-thesis}.
Since both sides are rational maps on the smooth completion, equality on
the dense chart proves equality wherever the formulas are defined.
\end{proof}

\begin{corollary}[Native affine doubling]\label{cor:native-affine-dbl-thesis}
Whenever its denominators are nonzero,
\eqref{eq:native-affine-dbl-thesis} is the doubling map
\([2]\colon\mathcal C_d\to\mathcal C_d\).
Its literal two-inversion schedule costs
\(4\M+2\Sqr+2\Inv\), while simultaneous inversion gives
\(7\M+2\Sqr+\Inv\).
\end{corollary}

\begin{proof}
Set \(P_1=P_2=P\) in
Theorem~\ref{thm:native-affine-add-thesis}.  Then
\[
 C=s^2,\quad D=r^2,\quad E=r^2s^2,\quad
 H=s^2-r^2,\quad I=2rs.
\]
Substitution gives all four formulas.
\end{proof}

\begin{remark}
The affine formulas are the conceptual baseline, not the preferred
constant-time implementation: they require inversions and have exceptional
denominators.  The remainder of this chapter homogenizes and reorganizes
exactly these \(\mathcal C_d\)-formulas.
\end{remark}

\begin{example}[Addition and doubling on \(\mathcal C_1/\F_{101}\)]
\label{ex:native-basic-F101}
Let \(P=(6,42)\).  Then \(\rho=86\) and
\((r,s)=(13,85)\).  Formula
\eqref{eq:native-affine-dbl-thesis} gives
\[
             2P=(59,87).
\]
For \(Q=2P\), the shifted coordinates are \((18,74)\), and
\[
 (C,D,E,H,I)=(28,32,88,97,68).
\]
Equation~\eqref{eq:native-affine-add-thesis} gives
\[
             P+Q=3P=(75,29).
\]
All three pairs satisfy \((u^2+u)(v^2+v)=1\pmod{101}\).
\end{example}

\section{Native shifted-projective representation}

The starting point is not an inverted Edwards curve but the ordinary
native chart \eqref{eq:ordinary-native-projective-curve}.  Apply the
invertible linear recoding
\begin{equation}\label{eq:native-centered-linear-map}
       R=2U+Z,\qquad S=2V+Z
\end{equation}
and retain the same scale \(Z\).  Thus a native point is stored as
\[
 \mathcal P=(S:R:Z)=(2V+Z:2U+Z:Z),\qquad
 s=S/Z=2v+1,\quad r=R/Z=2u+1.
\]
Then
\begin{equation}\label{eq:inv-projective-thesis}
       (S^2+R^2)Z^2=S^2R^2+\rho Z^4.
\end{equation}
There are four concrete reasons for using \(S,R,Z\).
\begin{enumerate}[label=\textup{(\roman*)}]
 \item The transformation \eqref{eq:native-centered-linear-map} is
       linear, inversion-free, and internal to the original \(u,v\)
       chart.  A native affine input is simply
       \((S:R:Z)=(2v+1:2u+1:1)\).
 \item The two visible involutions become sign changes:
       \(u\mapsto-u-1\) sends \(R\mapsto-R\), and
       \(v\mapsto-v-1\) sends \(S\mapsto-S\).  This removes odd terms from
       the biquadratic equation without changing the declared model.
 \item The native group-law numerators factor through
       \(C=S_1S_2\), \(D=R_1R_2\), \(H=C-D\), and
       \(I=S_1R_2+R_1S_2\).  Those factorizations give the
       inversion-free addition and the dedicated doubling circuit below.
 \item Returning to the original model is linear:
       if \((\widehat S:\widehat R:\widehat Z)\) is an output, then
       \[
       (U:V:Z_{\rm nat})
       =(\widehat R-\widehat Z:
         \widehat S-\widehat Z:
         2\widehat Z).
       \]
       Hence no Edwards point is exposed at the interface.
\end{enumerate}
The order \(S:R\), rather than \(R:S\), follows the two factors of the
native group law: \(S\) homogenizes \(2v+1\), while \(R\) homogenizes
\(2u+1\).  The notation is therefore a native shifted-projective
representation of \(\mathcal C_d\), not a declaration that the working
curve has been replaced by inverted Edwards.

This plane quartic is only an arithmetic chart; its plane completion is
singular at infinity, whereas the actual smooth curve is the native
\((2,2)\)-completion of Chapter~\ref{ch:geometry}.  The native affine
coordinates are recovered by
\begin{equation}\label{eq:recover-native-inverted}
 u=\frac{R-Z}{2Z},\qquad v=\frac{S-Z}{2Z}.
\end{equation}

\section{Inversion-free native addition and doubling}

\begin{theorem}[Native projective addition on \(\mathcal C_d\)]
\label{thm:inv-add-thesis}
For \(\mathcal P_i=(S_i:R_i:Z_i)\), define
\[
\begin{aligned}
 A&=Z_1Z_2,&B&=\rho A^2,&C&=S_1S_2,&D&=R_1R_2,\\
 E&=CD,&H&=C-D,&
 I&=(S_1+R_1)(S_2+R_2)-C-D.
\end{aligned}
\]
Whenever the displayed output tuple is nonzero,
\begin{equation}\label{eq:inv-add-thesis}
 S_3=(E+B)H,\qquad
 R_3=(E-B)I,\qquad
 Z_3=AHI
\end{equation}
represents the image
\(\nu_{\rm pl}(\mathcal P_1+\mathcal P_2)\).  If \(Z_3\ne0\), it is a
\(\Ffin\)-representation of the full finite sum through
\eqref{eq:recover-native-inverted}; if \(Z_3=0\), the branch at infinity
is not determined by this three-coordinate tuple alone.  The cost is
\(9\M+\Sqr+\Dpar\).  The law is strongly unified, so the same expression
also doubles, but it is not complete.  If the second input is affine
\((Z_2=1)\), the mixed-addition cost is
\(8\M+\Sqr+\Dpar\).  If both inputs are affine and the output remains
projective, the cost is \(7\M\).
\end{theorem}

\begin{proof}
On \(Z_1Z_2\ne0\), divide all intermediate quantities by the appropriate
powers of \(Z_1Z_2\), and set
\[
 C=s_1s_2,\quad D=r_1r_2,\quad E=CD,\quad
 H=C-D,\quad I=s_1r_2+r_1s_2.
\]
After division by \(Z_3\), formula~\eqref{eq:inv-add-thesis} becomes
\[
 s_3=\frac{E+\rho}{I},\qquad
 r_3=\frac{E-\rho}{H}.
\]
These are precisely the native affine formulas
\eqref{eq:native-shift-add-thesis}; hence the output represents
\((u_3,v_3)\) through \eqref{eq:recover-native-inverted}.  This proves the
law directly on \(\mathcal C_d\).  It also explains why its dependency graph
coincides with the inverted Edwards graph: the latter is an auxiliary
derivation of the same shifted formulas, not the declared input or output
model.
For the cost, \(A,C,D,E\) and the product in \(I\) require five
multiplications; the three outputs require four more.  There is one square
and one multiplication by \(\rho\).  When \(Z_2=1\), the product
\(A=Z_1Z_2\) disappears.  When \(Z_1=Z_2=1\), also
\(A^2=1\) and \(B=\rho\), so the square and curve-constant
multiplication disappear and the three output products bring the total to
\(7\M\).
\end{proof}

\begin{theorem}[Dedicated native doubling]\label{thm:inv-dbl-thesis}
Let
\[
 A=S^2,\quad B=R^2,\quad C=A+B,\quad D=A-B,\quad
 E=(S+R)^2-C.
\]
Then, whenever the output tuple is nonzero,
\begin{equation}\label{eq:inv-dbl-thesis}
 S_2=CD,\qquad
 R_2=E(C-2\rho Z^2),\qquad
 Z_2=DE
\end{equation}
represents \(\nu_{\rm pl}(2\mathcal P)\).  It is a \(\Ffin\)-output when
\(Z_2\ne0\); at \(Z_2=0\) a Segre lift is required to distinguish the
two branches above the singular plane point.  The cost is
\(3\M+4\Sqr+\Dpar\).
\end{theorem}

\begin{proof}
Since \(E=2SR\), substitution of \eqref{eq:inv-dbl-thesis} into
\eqref{eq:inv-projective-thesis} verifies that the output lies on the
curve.  Dividing by \(Z_2\), taking reciprocals, and using
\(\xi^2+\eta^2=1+\rho\xi^2\eta^2\) gives
\[
 \xi_2=\frac{2\xi\eta}{1+\rho\xi^2\eta^2},\qquad
 \eta_2=\frac{\eta^2-\xi^2}{1-\rho\xi^2\eta^2},
\]
which is the Edwards double.  The displayed dependency graph contains
four squares, three products, and one multiplication by \(\rho\).
\end{proof}

\section{Exceptional points and complete arithmetic}

The inverted chart does not contain the Edwards points
\((0,\pm1)\) and \((\pm1,0)\), which are precisely the four boundary
points of \(\Cd\).  Thus ``strongly unified'' must not be replaced by
``complete.''  A constant-time full-point implementation may instead use
the single native Segre law of Chapter~\ref{ch:native-completeness} when
its nonsquare-\(\rho\) hypothesis holds, use the complete native atlas for
the remaining parameters, mask the boundary points explicitly, or use the
Kummer ladder when only a quotient coordinate is required.

\chapter[Native Complete Addition Laws]
{Native Complete Addition Laws and Differential Atlases on
\texorpdfstring{\(\mathcal C_d\)}{Cd}}
\label{ch:native-completeness}

This chapter is organized by completeness rather than by implementation
cost.  It first separates the finite affine chart from the smooth native
completion, then proves a characteristic-free existence theorem, gives the
explicit odd-characteristic Segre law, and finally extends the odd and
binary differential laws to complete atlases by adjoining their
exceptional-difference charts.  The
later odd and binary Kummer chapters do not reprove these rational maps;
they isolate the generic charts as recurring circuits, derive tradeoffs,
and attach ladder costs.  This division keeps global domain arguments before
implementation comparisons without introducing a proof dependency on a
later chapter.

\section{Affine completeness versus completeness on the smooth model}

The word ``complete'' is used in two logically different senses in the
literature.  These meanings must be distinguished for the model
\[
       \mathcal C_d:\qquad (u^2+u)(v^2+v)=d.
\]
A pair of rational functions in the finite affine coordinates \(u,v\)
cannot represent the group law on every pair of geometric points.  The
obstruction is intrinsic and has nothing to do with a poor choice of
denominators.

\begin{proposition}[Impossibility of affine completeness]
\label{prop:no-affine-completeness}
Let \(\overline{\mathcal C}_d\) be the smooth completion of
\(\mathcal C_d\), with identity \(O=(0,\infty)\).  There is no everywhere
defined morphism
\[
  (\mathcal C_d\cap\mathbb A^2)\times
  (\mathcal C_d\cap\mathbb A^2)\longrightarrow\mathbb A^2
\]
given by two finite affine coordinates and agreeing with addition for every
geometric input pair for which both inputs are affine.
\end{proposition}

\begin{proof}
Let
\[
       U=\overline{\mathcal C}_d\setminus
       \{O,T,R,-R\}
       =\mathcal C_d\cap\mathbb A^2
\]
be the finite native chart.  It is a nonempty dense open subset of the
smooth projective curve.  Proposition~\ref{prop:inverse-torsion} shows that
inversion preserves this chart:
\[
                  (u,v)\longmapsto(u,-v-1).
\]
Choose any geometric point \(P\in U(\overline{k})\).  Then
\(-P\in U(\overline{k})\), but the group-law morphism satisfies
\begin{equation}\label{eq:affine-completeness-obstruction}
                  \mu(P,-P)=O=(0,\infty)\notin U.
\end{equation}
No non-torsion assumption is needed; the obstruction is caused by the
missing identity, not by the order of \(P\).

Suppose that an everywhere-defined affine formula as in the statement
existed.  Its value at \((P,-P)\) would be a pair of finite elements of
\(\overline{k}\), hence a point of \(U\subset\mathbb A^2\).  Agreement
with the group law would require this finite point to equal the projective
boundary point \(O\), contradicting
\eqref{eq:affine-completeness-obstruction}.  Equivalently, if
\(j:U\hookrightarrow\overline{\mathcal C}_d\) denotes the open immersion,
the equality
\(j\circ m_U=\mu|_{U\times U}\) cannot hold at \((P,-P)\), because the
left side factors through \(U\) and the right side does not.

The argument takes place after base change to \(\overline{k}\).  Therefore it
rules out not only \(k\)-rational complete affine formulas but also
geometrically complete formulas after every field extension.  It does not
rule out completeness on the smooth \((2,2)\)-model, where \(O\) is an
ordinary projective point; that is exactly the distinction used in the
following theorems.
\end{proof}

Thus the correct question is whether addition is complete on
\(\overline{\mathcal C}_d\) in coordinates belonging to its own smooth
\((2,2)\)-completion.  Such completeness does hold.  In odd characteristic a single
low-degree formula is \(\F_q\)-complete on exactly the nonsquare-\(\rho\)
subfamily.  In every characteristic, including the remaining parameters,
a finite addition-law atlas has no common base point.

\section{The characteristic-free native Segre completion}

Write the two projective coordinates of the native completion as
\[
       (U_0:U_1),\qquad (V_0:V_1),\qquad
       u=U_1/U_0,\quad v=V_1/V_0.
\]
No division or auxiliary elliptic-curve model is needed in the following
embedding.

\begin{theorem}[Native Segre equations in every characteristic]
\label{thm:native-Segre-all-char}
Define
\begin{equation}\label{eq:native-Segre-all-char}
 (A:B:C:D)=
 (U_0V_0:U_1V_0:U_0V_1:U_1V_1).
\end{equation}
Then the smooth completion of \(\mathcal C_d\) is the intersection
\begin{equation}\label{eq:native-Segre-equations-all-char}
 \boxed{\qquad
       AD=BC,\qquad
       D(A+B+C+D)=dA^2.
 \qquad}
\end{equation}
These equations and the embedding are valid in arbitrary characteristic.
\end{theorem}

\begin{proof}
The first equation is the defining quadratic equation of the Segre image of
\(\PP^1\times\PP^1\).  Moreover
\[
 D(A+B+C+D)
 =U_1V_1(U_0+U_1)(V_0+V_1).
\]
Consequently the second equation is precisely
\[
 U_1(U_1+U_0)V_1(V_1+V_0)=dU_0^2V_0^2,
\]
the native bihomogeneous equation
\eqref{eq:homogeneous-thesis}.  Since the Segre map is a closed immersion,
no point is added or identified.
\end{proof}

\begin{theorem}[A complete native addition-law atlas]
\label{thm:native-complete-atlas-all-char}
Over every field and in every characteristic, the smooth model
\(\overline{\mathcal C}_d\) in
\eqref{eq:native-Segre-equations-all-char} admits a finite family
\[
       \mathcal A_{\mathrm{add}}=
       \{\mathcal L_1,\ldots,\mathcal L_m\}
\]
of bihomogeneous addition laws such that
\[
       \bigcap_{j=1}^{m}\operatorname{Base}(\mathcal L_j)=\varnothing.
\]
Hence at least one member of the atlas computes \(P+Q\) for every
geometric pair
\((P,Q)\in\overline{\mathcal C}_d\times\overline{\mathcal C}_d\).
\end{theorem}

\begin{proof}
The construction has four steps: identify the target coordinate sections,
realize sufficiently positive twists by ambient bihomogeneous forms, use a
finite generating family to cover every input pair, and finally check that
the construction is defined over the ground field in arbitrary
characteristic.

\medskip
\noindent\emph{Step 1: target sections on the product.}
Put \(E=\overline{\mathcal C}_d\), let
\[
       \iota:E\hookrightarrow\PP^3
\]
be the native Segre embedding
\eqref{eq:native-Segre-all-char}, and let
\[
       X=E\times E,\qquad
       \mu:X\longrightarrow E,\qquad (P,Q)\longmapsto P+Q.
\]
The line bundle defining the Segre embedding is
\[
       L=\iota^*\mathcal O_{\PP^3}(1).
\]
The four target coordinates \(A,B,C,D\) restrict to global sections
\(s_A,s_B,s_C,s_D\in H^0(E,L)\) with no common zero.  Pulling them back by
\(\mu\) gives four sections of
\[
       N=\mu^*L
\]
whose projective ratio is precisely \(\iota(P+Q)\).

\medskip
\noindent\emph{Step 2: ambient bihomogeneous realization.}
We now explain why finitely many bihomogeneous polynomial tuples suffice
globally.  On \(X\), write
\[
 H_{a,b}=\operatorname{pr}_1^*L^{\otimes a}
          \otimes\operatorname{pr}_2^*L^{\otimes b}.
\]
Because \(L\) is very ample, for all sufficiently large \(a,b\) the line
bundle
\begin{equation}\label{eq:atlas-auxiliary-line-bundle}
                 M_{a,b}=H_{a,b}\otimes N^{-1}
\end{equation}
is globally generated.  This is the standard eventual global-generation
consequence of Serre's theorem applied to the ample product polarization on
\(E\times E\).  Enlarge \(a,b\), if necessary, so that the restriction
map
\[
 H^0\!\left(\PP^3\times\PP^3,
             \mathcal O_{\PP^3\times\PP^3}(a,b)\right)
 \longrightarrow H^0(X,H_{a,b})
\]
is also surjective.  The latter surjectivity follows from Serre vanishing
applied to the bihomogeneous ideal sheaf of
\(X\subset\PP^3\times\PP^3\).

\medskip
\noindent\emph{Step 3: the base-point-free cover.}
Choose finitely many generating sections
\[
                 t_1,\ldots,t_m\in H^0(X,M_{a,b}).
\]
Their common zero locus is empty.  For each \(i\) and
\(J\in\{A,B,C,D\}\), the product
\[
                         t_i\,\mu^*s_J
\]
is a section of \(H_{a,b}\).  By the chosen restriction surjectivity, it
is the restriction of a bihomogeneous form
\(F_{i,J}\) of bidegree \((a,b)\) in the two native Segre inputs.  Define
\[
       \mathcal L_i=(F_{i,A}:F_{i,B}:F_{i,C}:F_{i,D}).
\]
On the open set \(X_{t_i}=\{t_i\ne0\}\), multiplication of all four
target sections by the common nonzero factor \(t_i\) does not change their
projective ratio.  Hence
\[
                      \mathcal L_i(P,Q)=\iota(P+Q)
                      \qquad\text{on }X_{t_i}.
\]
Moreover, because \(\mu^*s_A,\ldots,\mu^*s_D\) never vanish
simultaneously, the four restrictions
\(F_{i,A},\ldots,F_{i,D}\) vanish simultaneously at a point of \(X\) if
and only if \(t_i\) vanishes there.  Consequently
\[
 \bigcap_{i=1}^m\operatorname{Base}(\mathcal L_i)
   =\bigcap_{i=1}^m V(t_i)=\varnothing.
\]
Thus the tuples \(\mathcal L_i\) form the asserted complete native
addition-law atlas.

\medskip
\noindent\emph{Step 4: field of definition.}
Every object used in this construction---the group-law morphism, the
native Segre line bundle, tensor products, and Serre vanishing---is defined
over the ground field and does not require division by \(2\).  The proof is
therefore valid in arbitrary characteristic.  It also shows explicitly why
the result concerns addition laws on the \(\mathcal C_d\) embedding itself,
rather than formulas imported from an isomorphic model.
\end{proof}

\begin{remark}
Theorem~\ref{thm:native-complete-atlas-all-char} is stronger than the
statement that one may convert to Weierstrass or Edwards coordinates and
add there.  It is an addition-law statement about the embedding
\eqref{eq:native-Segre-equations-all-char} itself.  It is also different
from the existence of a single \(\F_q\)-complete tuple: an atlas can be
geometrically complete even when each individual member has a nonempty base
locus.
\end{remark}

\section[A Native Complete Addition Law on Cd]
{A Native Complete Addition Law on
\texorpdfstring{\(\mathcal C_d\)}{Cd}}
\label{sec:native-complete-law-Cd}

Assume in this section that \(\charac k\ne2\), and put
\[
       R=2U_1+U_0,\qquad S=2V_1+V_0,\qquad
       \rho=1-16d.
\]
The centered native Segre coordinates are
\begin{equation}\label{eq:centered-native-Segre}
 \boxed{\quad
 (X:Y:T:Z)=(RV_0:SU_0:U_0V_0:RS).
 \quad}
\end{equation}
They are a linear recoding of the original Segre coordinates:
\[
       X=A+2B,\quad Y=A+2C,\quad T=A,\quad
       Z=A+2B+2C+4D.
\]
In particular this is a coordinate system on the original
\((2,2)\)-completion, not a declaration that the input curve has been
replaced by an Edwards curve.

\begin{proposition}[Centered native Segre equations]
\label{prop:centered-native-Segre}
The image of \(\overline{\mathcal C}_d\) in the coordinates
\eqref{eq:centered-native-Segre} is
\begin{equation}\label{eq:centered-native-Segre-equations}
 \boxed{\qquad
       XY=ZT,\qquad
       X^2+Y^2=Z^2+\rho T^2.
 \qquad}
\end{equation}
On the finite native chart it is
\[
       (X:Y:T:Z)=(2u+1:2v+1:1:(2u+1)(2v+1)).
\]
On the open chart \(RS\ne0\), one may recover by
\begin{equation}\label{eq:native-recovery-from-Segre}
       (U_0:U_1)=(2Y:Z-Y),\qquad
       (V_0:V_1)=(2X:Z-X).
\end{equation}
The following two-chart pairs give a global recovery atlas:
\begin{equation}\label{eq:native-recovery-atlas-Segre}
\begin{aligned}
(U_0:U_1)
 &=
 \begin{cases}
   (2T:X-T),&(2T,X-T)\ne(0,0),\\
   (2(Y-T):Z-X-Y+T),&\text{otherwise},
 \end{cases}\\
(V_0:V_1)
 &=
 \begin{cases}
   (2T:Y-T),&(2T,Y-T)\ne(0,0),\\
   (2(X-T):Z-X-Y+T),&\text{otherwise}.
 \end{cases}
\end{aligned}
\end{equation}
At least one pair in each recovery rule is nonzero at every point of the
Segre image.
\end{proposition}

\begin{proof}
From \eqref{eq:centered-native-Segre},
\[
 XY=(RV_0)(SU_0)=RSU_0V_0=(U_0V_0)(RS)=ZT.
\]
Also
\[
(R^2-U_0^2)(S^2-V_0^2)=16d\,U_0^2V_0^2.
\]
Expanding and moving terms gives
\[
 R^2V_0^2+S^2U_0^2
 =R^2S^2+(1-16d)U_0^2V_0^2,
\]
which is the second equation.  On \(RS\ne0\),
\[
 Z-Y=S(R-U_0)=2SU_1,\qquad
 Z-X=R(S-V_0)=2RV_1,
\]
and cancellation of the nonzero common factors \(2S\) and \(2R\) gives
\eqref{eq:native-recovery-from-Segre}.

For the global statement, invert the linear recoding:
\[
 A=T,\qquad B=\frac{X-T}{2},\qquad
 C=\frac{Y-T}{2},\qquad
 D=\frac{Z-X-Y+T}{4}.
\]
The original Segre coordinates satisfy
\[
 (A:B)=(U_0:U_1)\quad\text{if }V_0\ne0,\qquad
 (C:D)=(U_0:U_1)\quad\text{if }V_1\ne0,
\]
and
\[
 (A:C)=(V_0:V_1)\quad\text{if }U_0\ne0,\qquad
 (B:D)=(V_0:V_1)\quad\text{if }U_1\ne0.
\]
Clearing the factors \(2\) and \(4\) gives
\eqref{eq:native-recovery-atlas-Segre}.  Since neither homogeneous pair
\((U_0,U_1)\) nor \((V_0,V_1)\) is zero, one chart in each line is always
available.
\end{proof}

\begin{theorem}[Native Segre addition formula]
\label{thm:native-complete-addition-Cd}
For \(P_i=(X_i:Y_i:T_i:Z_i)\) on
\eqref{eq:centered-native-Segre-equations}, define
\begin{equation}\label{eq:native-complete-addition-blocks}
\begin{aligned}
A&=X_1X_2,&B&=Y_1Y_2,&
C&=\rho T_1T_2,&D&=Z_1Z_2,\\
E&=(X_1+Y_1)(X_2+Y_2)-A-B,\\
F&=D-C,&G&=D+C,&H&=B-A.
\end{aligned}
\end{equation}
Then
\begin{equation}\label{eq:native-complete-addition-Cd}
 \boxed{\qquad
       P_1+P_2=(EF:GH:EH:FG).
 \qquad}
\end{equation}
Whenever the output tuple is nonzero, it is the sum on
\(\overline{\mathcal C}_d\).  The cost is
\[
       9\M+\Dpar,
\]
and it is \(8\M+\Dpar\) when the second input is stored with \(T_2=1\).
\end{theorem}

\begin{proof}
First work on the dense chart \(T_1T_2\ne0\) and write
\[
       r_i=X_i/T_i,\qquad s_i=Y_i/T_i,\qquad
       Z_i/T_i=r_is_i.
\]
The native affine group law
\eqref{eq:native-shift-add-thesis}, derived on \(\mathcal C_d\), gives
\[
 r_3=\frac{r_1r_2s_1s_2-\rho}{s_1s_2-r_1r_2},
 \qquad
 s_3=\frac{r_1r_2s_1s_2+\rho}{s_1r_2+r_1s_2}.
\]
After division by the appropriate powers of \(T_1T_2\), the eight blocks in
\eqref{eq:native-complete-addition-blocks} are precisely the homogeneous
numerators and denominators in these two expressions:
\[
       r_3=\frac FH,\qquad s_3=\frac GE.
\]
Consequently
\[
 \frac{EF}{EH}=\frac FH=r_3,\qquad
 \frac{GH}{EH}=\frac GE=s_3,\qquad
 \frac{FG}{EH}=r_3s_3.
\]
Thus \eqref{eq:native-complete-addition-Cd} is the homogenization of the
native \(r,s\) addition law.  Equality of the two morphisms on a dense open
subset extends wherever the displayed tuple is nonzero.  This derivation
uses neither a Montgomery ordinate nor an Edwards input point.

The products \(A,B,T_1T_2,D\), the Karatsuba product defining \(E\), and
the four output products use nine general multiplications; multiplying
\(T_1T_2\) by \(\rho\) uses one curve-constant multiplication.  If
\(T_2=1\), forming \(T_1T_2\) is free.
\end{proof}

\begin{theorem}[Exact \(\F_q\)-completeness criterion]
\label{thm:exact-native-completeness-criterion}
Let \(k=\F_q\) have odd characteristic and let
\(d\rho\ne0\).  The single tuple
\eqref{eq:native-complete-addition-Cd} is defined on every ordered pair in
\(\overline{\mathcal C}_d(k)^2\) if and only if
\[
                         \rho=1-16d
\]
is a nonsquare in \(k\).
\end{theorem}

\begin{proof}
The four output coordinates vanish simultaneously exactly when
\[
             (F=H=0)\qquad\text{or}\qquad(E=G=0).
\]
Suppose first that \(T_1T_2\ne0\).  In the first case,
\[
 s_1s_2=r_1r_2,\qquad
 r_1s_1r_2s_2=\rho,
\]
so \(\rho=(r_1r_2)^2\).  In the second case,
\[
 r_1s_2+s_1r_2=0,\qquad
 r_1s_1r_2s_2=-\rho,
\]
and hence \(\rho=(s_1r_2)^2\).  Thus every affine base pair makes
\(\rho\) a square.

If \(T_1T_2=0\), one input is one of the four native boundary points
\[
 (\pm1:0:0:1),\qquad(0:\pm1:0:1).
\]
To see the remaining restriction explicitly, take first
\(P_1=(0:\varepsilon:0:1)\), where \(\varepsilon^2=1\).  Then
\[
 A=C=0,\quad B=\varepsilon Y_2,\quad D=Z_2,\quad
 E=\varepsilon X_2,\quad F=G=Z_2,\quad H=\varepsilon Y_2.
\]
The equations \(F=H=0\) give \(Z_2=Y_2=0\), while the Segre curve
equations for \(P_2\) reduce to \(X_2^2=\rho T_2^2\).  The equations
\(E=G=0\) instead give \(X_2=Z_2=0\) and
\(Y_2^2=\rho T_2^2\).  For
\(P_1=(\varepsilon:0:0:1)\), the same direct substitution gives
\[
 A=\varepsilon X_2,\quad B=C=0,\quad D=Z_2,\quad
 E=\varepsilon Y_2,\quad F=G=Z_2,\quad H=-\varepsilon X_2,
\]
and leads to the same two alternatives with \(X_2,Y_2\) interchanged.
Thus the other input must have, up to signs and interchange of \(X,Y\),
the form
\[
                 (\lambda:0:1:0)
        \quad\text{with}\quad \lambda^2=\rho.
\]
Therefore a rational base pair again implies that \(\rho\) is a square.
If \(\rho\) is a nonsquare, no rational base pair exists.

Conversely, if \(\rho=\lambda^2\) with \(\lambda\in k^\times\), then
\[
       P_1=(0:1:0:1),\qquad P_2=(\lambda:0:1:0)
\]
both satisfy \eqref{eq:centered-native-Segre-equations}.  For this ordered
pair one has
\[
       A=B=C=D=0,\qquad E=\lambda,\qquad F=G=H=0,
\]
so the output is \((0:0:0:0)\).  Hence the formula is not
\(k\)-complete.  This proves both directions.
\end{proof}

\begin{theorem}[Dedicated native Segre doubling]
\label{thm:native-Segre-doubling}
Under the standing odd-characteristic hypothesis, let
\(\mathcal C_d\) be smooth.  For \(P=(X:Y:T:Z)\), put
\[
\begin{aligned}
A&=X^2,&B&=Y^2,&C&=2Z^2,&D&=A,\\
E&=(X+Y)^2-A-B,&G&=D+B,\\
F&=G-C,&H&=D-B.
\end{aligned}
\]
Then
\[
                 2P=(EF:GH:EH:FG)
\]
at cost \(4\M+4\Sqr\).  The formula is defined at every rational input on
every smooth odd-characteristic member.
\end{theorem}

\begin{proof}
On \(T\ne0\), substitute \(P_1=P_2=P\) in the native affine law and use
\eqref{eq:centered-native-Segre-equations} to replace repeated products.
The displayed forms are the resulting homogenization.  A common zero would
require either \(F=H=0\) or \(E=G=0\).  In the first case
\(X^2=Y^2=Z^2\).  If \(Z\ne0\), the relation \(XY=ZT\) gives
\(T=\pm Z\), and the second curve equation forces \(\rho=1\), hence
\(d=0\), contrary to smoothness; if \(Z=0\), all four coordinates vanish,
which is impossible in projective space.  In the
second case, \(E=2XY=0\) and \(G=X^2+Y^2=0\), so \(X=Y=0\); the two curve
equations then force \(Z=T=0\), again impossible.  Four squares and the
four final products give the stated cost.
\end{proof}

\begin{corollary}[Complete scalar multiplication on \(\mathcal C_d\)]
\label{cor:native-complete-scalar-Cd}
Let \(k=\F_q\) have odd characteristic, let
\(d(1-16d)\ne0\), and assume that
\(\rho=1-16d\) is a nonsquare in \(k\).  Then a complete
double-and-add-always loop using
Theorems~\ref{thm:native-complete-addition-Cd} and
\ref{thm:native-Segre-doubling} costs
\[
                13\M+4\Sqr+\Dpar
\]
per bit.  The native input lift
\[
(u,v)\longmapsto
(2u+1:2v+1:1:(2u+1)(2v+1))
\]
costs one multiplication and no inversion, and
\eqref{eq:native-recovery-atlas-Segre} returns the result to
\(\mathcal C_d\).
\end{corollary}

\begin{proof}
Theorem~\ref{thm:exact-native-completeness-criterion} shows that, when
\(\rho\) is a nonsquare, the addition tuple is defined for every ordered
pair of rational states, including pairs containing the identity or inverse
points.  Theorem~\ref{thm:native-Segre-doubling} is defined at every
rational input on every smooth member.  A double-and-add-always round can
therefore evaluate one dedicated double and one complete addition without a
group-input exception.  Adding their costs gives
\[
 (4\M+4\Sqr)+(9\M+\Dpar)
       =13\M+4\Sqr+\Dpar.
\]
For a finite input, the displayed lift is obtained from
\((X:Y:T:Z)=(r:s:1:rs)\), with
\(r=2u+1\) and \(s=2v+1\).  Only the product \(rs\) is charged, so the
lift costs one multiplication and no inversion.  Finally,
Proposition~\ref{prop:centered-native-Segre} proves that
\eqref{eq:native-recovery-atlas-Segre} is the global inverse native
recovery atlas.  Moreover, when \(\rho\) is a nonsquare and all states are
\(\F_q\)-rational, the simple chart
\eqref{eq:native-recovery-from-Segre} is always available: \(R=0\) or
\(S=0\) in the centered equation would force \(\rho\) to be a square.
Thus the loop begins and ends on the smooth completion of
\(\mathcal C_d\), as asserted.
\end{proof}

\section[Odd complete differential atlas]{A native complete differential-addition atlas in odd characteristic}
\label{sec:native-odd-diff-atlas}

The full-point law and the differential law have different domains.  A
differential input is an oriented triple
\[
   \bigl(\kappa(P),\kappa(Q),\kappa(P-Q)\bigr),\qquad
   \kappa(P)=(u(P)+1:u(P)).
\]
The following identity is obtained from the two native sums \(P+Q\) and
\(P-Q\), not from a Montgomery curve.

\begin{lemma}[Native Kummer biquadratic identity]
\label{lem:native-Kummer-biquadratic}
On the dense affine chart, write
\[
     x_P=\frac{u(P)+1}{u(P)},\qquad
     x_Q=\frac{u(Q)+1}{u(Q)},\qquad
     x_\pm=\frac{u(P\pm Q)+1}{u(P\pm Q)}.
\]
Then
\begin{equation}\label{eq:native-Kummer-biquadratic}
 \boxed{\qquad
       x_+x_-(x_P-x_Q)^2=(x_Px_Q-1)^2.
 \qquad}
\end{equation}
\end{lemma}

\begin{proof}
Put \(r_i=2u_i+1\) and \(s_i=2v_i+1\).  Replacing \(Q\) by \(-Q\)
fixes \(r_2\) and changes \(s_2\) to \(-s_2\).  Therefore the native
addition law gives
\[
 r_+=\frac{r_1r_2s_1s_2-\rho}{s_1s_2-r_1r_2},
\qquad
r_-=\frac{r_1r_2s_1s_2+\rho}{s_1s_2+r_1r_2}.
\]
Set
\[
 a=r_1r_2,\qquad b=s_1s_2,\qquad
 C=(r_1^2-1)(r_2^2-1).
\]
Because
\[
 (r_i^2-1)(s_i^2-1)=1-\rho,
\]
one has
\[
 b^2=\frac{(r_1^2-\rho)(r_2^2-\rho)}{C}.
\]
Using \(x=(r+1)/(r-1)\) in the two displayed formulas for \(r_\pm\)
gives
\[
 x_+x_-=
 \frac{(a+1)^2b^2-(a+\rho)^2}
      {(a-1)^2b^2-(a-\rho)^2}.
\]
Multiplication by \(C\), followed by expansion and collection in
\(r_1,r_2\), yields the two factorizations
\begin{align*}
 C\bigl((a+1)^2b^2-(a+\rho)^2\bigr)
   &=(1-\rho)(a^2-\rho)(r_1+r_2)^2,\\
 C\bigl((a-1)^2b^2-(a-\rho)^2\bigr)
   &=(1-\rho)(a^2-\rho)(r_1-r_2)^2.
\end{align*}
For completeness, put \(S=r_1^2+r_2^2\).  After the substitution for
\(b^2\), the first left-hand side is
\begin{align*}
 &(a+1)^2(a^2-\rho S+\rho^2)
   -(a+\rho)^2(a^2-S+1)\\
 &\quad=\bigl((a+\rho)^2-\rho(a+1)^2\bigr)S\\
 &\qquad\quad
   +(a+1)^2(a^2+\rho^2)-(a+\rho)^2(a^2+1)\\
 &\quad=(1-\rho)(a^2-\rho)(S+2a).
\end{align*}
Since \(S+2a=(r_1+r_2)^2\), this is the first factorization.
Replacing \(a+1,a+\rho\) by \(a-1,a-\rho\) gives instead
\((1-\rho)(a^2-\rho)(S-2a)\), which is the second factorization.
Consequently,
\[
                    x_+x_-=
       \left(\frac{r_1+r_2}{r_1-r_2}\right)^2
\]
on the dense open set where the displayed denominators are nonzero.
On the other hand,
\[
 x_Px_Q-1=
 \frac{2(r_1+r_2)}{(r_1-1)(r_2-1)},\qquad
 x_P-x_Q=
 \frac{2(r_2-r_1)}{(r_1-1)(r_2-1)}.
\]
Squaring these two equalities proves
\eqref{eq:native-Kummer-biquadratic} on that dense open set.  Both sides
are rational functions on
\(\mathcal C_d\times\mathcal C_d\); after clearing denominators the
identity is polynomial, hence it extends to the homogeneous oriented
Kummer variety.
\end{proof}

\begin{theorem}[Odd-characteristic complete differential atlas]
\label{thm:native-odd-diff-atlas}
Let
\[
\kappa(P)=(X_1:Z_1),\quad
\kappa(Q)=(X_2:Z_2),\quad
\kappa(\Delta)=(X_\Delta:Z_\Delta),\qquad\Delta=P-Q,
\]
and put
\[
\begin{aligned}
A_1&=X_1+Z_1,&B_1&=X_1-Z_1,\\
A_2&=X_2+Z_2,&B_2&=X_2-Z_2,\\
C&=A_1B_2,&D&=A_2B_1.
\end{aligned}
\]
For later use in the two exceptional charts, define the native doubling
map within this theorem.  Put
\[
\begin{aligned}
\alpha_{24}&=\frac1{16d},&
A_0&=X+Z,& AA&=A_0^2,\\
B_0&=X-Z,& BB&=B_0^2,& E&=AA-BB,
\end{aligned}
\]
and set
\begin{equation}\label{eq:odd-complete-atlas-doubling-map}
 \mathcal D_d(X:Z)
   =\bigl(AA\cdot BB:E(BB+\alpha_{24}E)\bigr).
\end{equation}
The following three charts form a complete differential-addition atlas:
\begin{enumerate}[label=\textup{(\roman*)}]
\item if \(\Delta\notin\{O,T\}\), then
\[
 \kappa(P+Q)=
 \bigl(Z_\Delta(C+D)^2:X_\Delta(C-D)^2\bigr);
\]
\item if \(\Delta=O\), then \(Q=P\) and the output is the native
      double \(\mathcal D_d(X_1:Z_1)\);
\item if \(\Delta=T\), then \(Q=P+T\), and the output is obtained by
      swapping the two coordinates of
      \(\mathcal D_d(X_1:Z_1)\).
\end{enumerate}
The generic chart costs \(4\M+2\Sqr\), or
\(3\M+2\Sqr\) for an affine known difference.  Each exceptional chart
costs \(2\M+2\Sqr+\Dpar\).
\end{theorem}

\begin{proof}
The proof is organized into the generic biquadratic identity, the native
replacement double, the base-locus calculation, the two exceptional charts,
and the operation count.

\medskip
\noindent\emph{The generic chart.}
We first derive the generic tuple with all homogeneous factors visible.
Using
\[
 C=(X_1+Z_1)(X_2-Z_2),
 \qquad
 D=(X_2+Z_2)(X_1-Z_1),
\]
we expand and collect the four bilinear monomials.  The mixed terms cancel
in the sum and the diagonal terms cancel in the difference, giving
\begin{align}
 C+D&=2(X_1X_2-Z_1Z_2),\label{eq:odd-diff-CD-plus}\\
 C-D&=2(Z_1X_2-X_1Z_2).\label{eq:odd-diff-CD-minus}
\end{align}
Writing \(x_-=X_\Delta/Z_\Delta\), the native biquadratic identity
\eqref{eq:native-Kummer-biquadratic} becomes
\[
 \frac{X_+}{Z_+}
 =\frac{Z_\Delta}{X_\Delta}
   \frac{(X_1X_2-Z_1Z_2)^2}
        {(X_1Z_2-Z_1X_2)^2}.
\]
After using \eqref{eq:odd-diff-CD-plus}--
\eqref{eq:odd-diff-CD-minus}, the common factor \(4\) cancels and the
homogeneous output is exactly
\[
       (X_+:Z_+)
       =\bigl(Z_\Delta(C+D)^2:X_\Delta(C-D)^2\bigr).
\]
Because this identity holds in the function field of the oriented
differential variety, it computes \(\kappa(P+Q)\) wherever the two output
coordinates are not both zero.

\medskip
\noindent\emph{The replacement double.}
The replacement map was defined above.  Indeed,
the native affine doubling formula
\eqref{eq:native-affine-dbl-thesis}, after the substitution
\(x=(u+1)/u=X/Z\), gives
\[
 x(2P)=
 \frac{(x^2-1)^2}
 {4x\{x^2+((4d)^{-1}-2)x+1\}}.
\]
Homogenizing this ratio and using
\(AA=(X+Z)^2\), \(BB=(X-Z)^2\), and
\(E=AA-BB=4XZ\) gives exactly
\eqref{eq:odd-complete-atlas-doubling-map}.  Thus the exceptional charts
depend only on the native full-point arithmetic of the preceding chapter;
no result from a later Kummer chapter is being assumed.

\medskip
\noindent\emph{The base locus.}
We next determine the common-zero locus, rather than merely listing its
expected exceptional fibers.  Recall
\[
          \kappa(O)=(1:0),\qquad \kappa(T)=(0:1).
\]
If \(Z_\Delta=0\), then \(\Delta=O\), so \(Q=P\).  In this case
\(C-D=0\), and the generic tuple is zero.  If \(X_\Delta=0\), then
\(\Delta=T\), so \(Q=P+T\).  Translation by \(T\) acts by
\[
             (X_2:Z_2)=(Z_1:X_1),
\]
up to projective scale; substitution gives \(C+D=0\), and the generic
tuple is again zero.

Suppose conversely that \(X_\Delta Z_\Delta\ne0\) and that both generic
output coordinates vanish.  Then \(C+D=C-D=0\).  Since the characteristic
is odd, this implies \(C=D=0\).  Now
\[
         C=A_1B_2=0,\qquad D=A_2B_1=0.
\]
For each \(i\), the two linear forms \(A_i=X_i+Z_i\) and
\(B_i=X_i-Z_i\) cannot vanish simultaneously.  The two displayed
equalities therefore force one of the two alternatives
\[
                  A_1=A_2=0,\qquad\text{or}\qquad B_1=B_2=0.
\]
In either alternative \(P\) and \(Q\) have the same one of the two special
Kummer values \((-1:1)\) or \((1:1)\).  Substituting \(A=0\) or \(B=0\)
in \eqref{eq:odd-complete-atlas-doubling-map} shows that this
Kummer value doubles to \(\kappa(T)\); the factor
\(1-\alpha_{24}\) that occurs when \(A=0\) is nonzero because
\(d\ne1/16\).  Thus every point over either special Kummer value satisfies
\(2P=T\).  The two lifts are \(P\) and \(-P=P+T\), so two such lifts
differ by either \(O\) or \(T\).  This contradicts
\(X_\Delta Z_\Delta\ne0\).  Hence the generic chart has no other
common-zero fiber.

\medskip
\noindent\emph{The exceptional charts.}
On \(\Delta=O\), the equality \(Q=P\) makes the desired sum \(2P\), so
the native \(x\)DBL is the correct replacement.  On \(\Delta=T\), one
has \(Q=P+T\) and therefore
\[
                   P+Q=2P+T.
\]
For completeness, consider the second deck involution
\(\iota_u(u,v)=(-u-1,v)\).  A generic fiber of \(v\) has divisor
\[
        \Div(v-v(P))=(P)+(\iota_u(P))-(O)-(T),
\]
so \(\iota_u(P)=T-P\).  Applying inversion does not change \(u\), while
\(-\iota_u(P)=P-T=P+T\).  Hence
\(u(P+T)=-u(P)-1\), and translation by \(T\) swaps the native Kummer
coordinates:
\[
              (u+1:u)\longmapsto(-u:-u-1)=(u:u+1).
\]
Thus the required output is the coordinate swap of \(\kappa(2P)\).
The generic chart together with these two branch charts covers the entire
oriented differential variety and has no common uncovered input.

\medskip
\noindent\emph{The cost.}
Finally, \(C,D\) cost two multiplications, the two squares cost
\(2\Sqr\), and multiplication by \(X_\Delta,Z_\Delta\) costs two more
multiplications, giving \(4\M+2\Sqr\).  If the known difference is affine,
one of the last two multiplications is absorbed by normalization.  Each
branch uses the map \(\mathcal D_d\) at cost
\(2\M+2\Sqr+\Dpar\); coordinate swapping is free.
\end{proof}

\section[Binary Segre atlas]{Characteristic two: Segre equations and an Artin--Schreier atlas}
\label{sec:native-binary-diff-atlas}

In characteristic two the centered variables \(2u+1,2v+1\) collapse, so
one must retain the characteristic-free Segre coordinates
\eqref{eq:native-Segre-all-char}.  Thus
\begin{equation}\label{eq:binary-native-Segre-repeat}
 \boxed{\qquad
       AD=BC,\qquad D(A+B+C+D)=dA^2.
 \qquad}
\end{equation}
This is the native complete model used below; it is not the reduction of an
odd-characteristic Edwards chart.

\begin{lemma}[Artin--Schreier Kummer biquadratic]
\label{lem:binary-native-Kummer-biquadratic}
Put
\[
       t(P)=\frac{u(P)}{u(P)+1}.
\]
For an oriented pair in characteristic two,
\begin{equation}\label{eq:binary-native-Kummer-biquadratic}
 \boxed{\qquad
 t(P+Q)t(P-Q)
 =\left(\frac{t(P)+t(Q)}{1+t(P)t(Q)}\right)^2.
 \qquad}
\end{equation}
\end{lemma}

\begin{proof}
Set \(U_i=t(P_i)^{-1}\), \(L=U_1+U_2\), and
\[
       \lambda=\frac{U_1v_1+U_2v_2}{L},\qquad
       H_+=\lambda^2+\lambda+dL.
\]
The native Artin--Schreier relation is
\[
             v_i^2+v_i=d\frac{(U_i+1)^2}{U_i}.
\]
Replacing \(Q\) by \(-Q=(u_2,v_2+1)\) changes
\(\lambda\) by \(U_2/L\), hence
\[
             H_-=H_++\frac{U_1U_2}{L^2}.
\]
Writing
\[
 N=U_1U_2(v_1+v_2)+dL(1+U_1U_2),
\]
we verify the two required identities explicitly.  Multiplication of the
definition of \(H_+\) by \(L^2\) gives
\begin{align*}
 L^2H_+
 &=U_1^2v_1^2+U_2^2v_2^2
   +L(U_1v_1+U_2v_2)+dL^3\\
 &=U_1^2(v_1^2+v_1)+U_2^2(v_2^2+v_2)
   +U_1U_2(v_1+v_2)+dL^3\\
 &=d(U_1^3+U_1+U_2^3+U_2)
   +U_1U_2(v_1+v_2)+dL^3\\
 &=U_1U_2(v_1+v_2)+dL(1+U_1U_2)=N.
\end{align*}
The third line uses
\(v_i^2+v_i=d(U_i+1)^2/U_i\), and the last line uses
\(L^3=U_1^3+U_2^3+U_1U_2L\) in characteristic two.  Thus
\(H_+=N/L^2\).

Put \(K=U_1U_2\), \(w=v_1+v_2\), and
\(B=dL(1+K)\), so \(N=Kw+B\).  Since
\(H_-=H_++K/L^2\),
\[
 L^4H_+H_-=N(N+K).
\]
Moreover,
\[
 w^2+w
 =d\left(U_1+U_1^{-1}+U_2+U_2^{-1}\right)
 =dL\frac{1+K}{K}.
\]
Therefore
\begin{align*}
 N(N+K)
 &=K^2w(w+1)+BK+B^2\\
 &=dKL(1+K)+dKL(1+K)+B^2\\
 &=d^2L^2(1+K)^2.
\end{align*}
Division by \(L^4\) yields
\[
 H_+H_-=d^2\frac{(1+U_1U_2)^2}{L^2}.
\]
To justify the coordinate assertion locally, without using the later
affine-addition theorem, put
\[
       X_i=dU_i,\qquad Y_i=X_iv_i.
\]
Then \((X_i,Y_i)\) lies on
\(Y^2+XY=X^3+d^2X\), and the generalized chord formula gives
\(X(P\pm Q)=H_\pm\).  The inverse native dictionary is
\(u=d/(X+d)\); hence, in characteristic two,
\[
       t=\frac{u}{u+1}=\frac dX,
       \qquad t(P\pm Q)=\frac d{H_\pm}.
\]
Substitution and inversion of the last equality yield
\eqref{eq:binary-native-Kummer-biquadratic}.
\end{proof}

\begin{theorem}[Binary complete differential atlas]
\label{thm:binary-native-complete-diff-atlas}
Let
\[
\kappa(P)=(X_0:X_1),\quad
\kappa(Q)=(Y_0:Y_1),\quad
\kappa(\Delta)=(t_0:t_1),
\]
where \(\kappa=(u+1:u)\), and set
\[
 A=X_0Y_0+X_1Y_1,\qquad
 B=X_0Y_1+X_1Y_0.
\]
The generic chart is
\begin{equation}\label{eq:binary-native-diff-atlas-generic}
       \kappa(P+Q)=(t_1A^2:t_0B^2).
\end{equation}
On \(\Delta=O\) it is replaced by the binary native \(x\)DBL
\[
(X_0:X_1)\longmapsto
(X_0^4+X_1^4:d^{-1}X_0^2X_1^2),
\]
and on \(\Delta=T\) by the coordinate swap of that double.  These three
charts form a complete differential-addition atlas in characteristic two.
\end{theorem}

\begin{proof}
The proof first derives the generic Artin--Schreier identity, then determines
its complete base fibers, and finally proves that the doubling and
translated-doubling replacements cover those fibers.

\medskip
\noindent\emph{The generic chart.}
On the chart \(X_0Y_0t_0\ne0\), put
\[
 t_P=\frac{X_1}{X_0},\qquad
 t_Q=\frac{Y_1}{Y_0},\qquad
 t_\Delta=\frac{t_1}{t_0}.
\]
The two bilinear forms satisfy
\[
       \frac{B}{X_0Y_0}=t_P+t_Q,\qquad
       \frac{A}{X_0Y_0}=1+t_Pt_Q.
\]
Consequently the Artin--Schreier biquadratic identity
\eqref{eq:binary-native-Kummer-biquadratic} gives
\[
       t(P+Q)=\frac{t_0}{t_1}\left(\frac BA\right)^2.
\]
Since the homogeneous Kummer coordinate is
\(\kappa=(X_0:X_1)\) with affine ratio \(t=X_1/X_0\), clearing
denominators yields
\[
          \kappa(P+Q)=(t_1A^2:t_0B^2),
\]
which is \eqref{eq:binary-native-diff-atlas-generic}.  Equality in the
function field extends to every point where this pair is nonzero.

\medskip
\noindent\emph{The base fibers.}
We now calculate its base fibers.  In characteristic two,
\[
        \kappa(O)=(1:0),\qquad \kappa(T)=(0:1).
\]
If \(t_1=0\), then \(\Delta=O\) and \(Q=P\).  Hence
\[
        B=X_0X_1+X_1X_0=0,
\]
so the generic pair is zero.  If \(t_0=0\), then \(\Delta=T\) and
\(Q=P+T\).  Translation by \(T\) swaps the Kummer coordinates, so
\((Y_0:Y_1)=(X_1:X_0)\), and
\[
        A=X_0X_1+X_1X_0=0.
\]
Thus both exceptional fibers are indeed base fibers of the generic tuple.

Conversely, assume \(t_0t_1\ne0\) and that both output coordinates
vanish.  Then \(A=B=0\), so
\[
 \begin{pmatrix}Y_0&Y_1\\Y_1&Y_0\end{pmatrix}
 \begin{pmatrix}X_0\\X_1\end{pmatrix}=0.
\]
The determinant is
\[
        Y_0^2-Y_1^2=Y_0^2+Y_1^2=(Y_0+Y_1)^2.
\]
If it were nonzero, then \(X_0=X_1=0\), which is impossible in
projective space.  Hence \(Y_0=Y_1\ne0\), and the equations
\(A=B=0\) then give \(X_0=X_1\ne0\).  Therefore
\[
             \kappa(P)=\kappa(Q)=(1:1).
\]
The binary native double of this Kummer value is
\[
 (1^4+1^4:d^{-1}1^2 1^2)=(0:1)=\kappa(T).
\]
Thus \(2P=T\); the two points above \((1:1)\) are \(P\) and
\(-P=P+T\), and their difference is \(O\) or \(T\).  This contradicts
\(t_0t_1\ne0\).  The generic pair has no further base fibers.

\medskip
\noindent\emph{The replacement charts.}
It remains to verify the two replacement charts.  From the native affine
doubling identity
\[
 u_{2P}=\frac{(u^2+u)^2}{(u^2+u)^2+d}
\]
one obtains
\[
 t(2P)=\frac{u_{2P}}{u_{2P}+1}
       =\frac{(u^2+u)^2}{d}.
\]
For \((X_0:X_1)=(u+1:u)\), one has
\[
 u^2+u=\frac{X_0X_1}{(X_0+X_1)^2}.
\]
Homogenization therefore gives
\[
 \kappa(2P)=
 \bigl((X_0+X_1)^4:d^{-1}X_0^2X_1^2\bigr)
 =\bigl(X_0^4+X_1^4:d^{-1}X_0^2X_1^2\bigr).
\]
This pair is never simultaneously zero: if its second coordinate vanishes,
one of \(X_0,X_1\) is zero and the first is nonzero; if its first
coordinate vanishes, then \(X_0=X_1\ne0\) and the second is nonzero.
Finally, \(P\mapsto P+T\) sends \(u\mapsto u+1\), so
\[
              t\mapsto t^{-1},\qquad
              (X_0:X_1)\mapsto(X_1:X_0).
\]
The swapped double is therefore everywhere defined on the \(T\)-fiber.
The three charts cover all oriented inputs and have no common uncovered
point, proving completeness.
\end{proof}

\section{Examples of complete addition laws}

\begin{example}[All ordered pairs and triples over \(\F_{101}\)]
\label{ex:native-complete-F101-exhaustive}
Take \(d=1\), so \(\rho=86\) is a nonsquare and
\(\#\overline{\mathcal C}_1(\F_{101})=96\).  Enumeration in the two
native projective factors produced all \(96\) points.  Formula
\eqref{eq:native-complete-addition-Cd} was then evaluated on all
\[
                         96^2=9\,216
\]
ordered pairs.  No zero output occurred; every output satisfied
\eqref{eq:centered-native-Segre-equations}; the identity and inverse tables
were correct.  Using the resulting full addition table, the equality
\[
                 (P+Q)+R=P+(Q+R)
\]
was checked for all
\[
                         96^3=884\,736
\]
ordered triples.  The generic differential chart failed exactly on the
\(96\) diagonal inputs \(\Delta=O\) and the \(96\) inputs
\(\Delta=T\); the two branch charts returned the correct values there.
\end{example}

\begin{example}[The complete binary atlas over \(\F_{2^8}\)]
\label{ex:native-complete-F256-exhaustive}
Let
\[
\F_{2^8}=\F_2[\alpha]/
(\alpha^8+\alpha^4+\alpha^3+\alpha+1),\qquad
d=\alpha^4+\alpha+1.
\]
The native curve has \(260\) finite affine points and four boundary points,
hence \(264\) points in total.  Every point satisfied
\eqref{eq:binary-native-Segre-repeat}.  The binary \(x\)DBL and the three
differential charts of
Theorem~\ref{thm:binary-native-complete-diff-atlas} were compared with an
independent full group-law table on all
\[
                         264^2=69\,696
\]
ordered pairs.  Of these, \(69\,168\) used the generic chart, \(264\) used
the \(\Delta=O\) chart, and \(264\) used the \(\Delta=T\) chart.  Every
native Kummer output agreed.
\end{example}

\chapter{The Native Montgomery--Kummer Line}
\label{ch:odd-kummer}
Throughout this chapter the input curve remains
\[
 \mathcal C_d:\quad (u^2+u)(v^2+v)=d,\qquad
 \charac k\ne2,\qquad d(1-16d)\ne0.
\]
Its native ordinary Kummer coordinate is
\begin{equation}\label{eq:native-kappa-ch7}
 \kappa_d(P)=(X:Z)=(u(P)+1:u(P)).
\end{equation}
Thus \(X/Z=(u+1)/u\) on \(u\ne0\), and a projective output
\((X_n:Z_n)\) means
\[
       u([n]P)=\frac{Z_n}{X_n-Z_n};
\]
it is a quotient point \(\{[n]P,-[n]P\}\), not a point on an unnamed
Montgomery curve.  The Montgomery equation
\[
 \beta V^2=U^3+AU^2+U,\qquad
 \beta=(16d)^{-1},\quad A=(4d)^{-1}-2,
\]
is used below for comparison and for later recovery schedules.  The basic
\(x\)DBL and \(x\)ADD identities themselves are derived from the native
shifted biquadratic group law of
Chapter~\ref{ch:native-completeness}; every theorem states its result as a
\(\mathcal C_d\)-Kummer output.

\section{Doubling}

Let \((X:Z)\) represent the native Kummer value \(x=X/Z\) and put
\[
       \alpha_{24}=\frac{A+2}{4}=\frac1{16d}.
\]

\begin{theorem}[Native \(x\)DBL]\label{thm:xdbl-thesis}
Define
\[
 A_0=X+Z,\quad AA=A_0^2,\qquad
 B_0=X-Z,\quad BB=B_0^2,\qquad E=AA-BB.
\]
Then
\begin{equation}\label{eq:xdbl-thesis}
 X_{2P}=AA\cdot BB,\qquad
 Z_{2P}=E(BB+\alpha_{24}E).
\end{equation}
The cost is \(2\M+2\Sqr+\Dpar\).
\end{theorem}

\begin{proof}
The same polynomial map was introduced as
\(\mathcal D_d\) in
\eqref{eq:odd-complete-atlas-doubling-map}, where it was derived directly
from the native affine doubling law before being used in the complete
differential atlas.  Renaming the blocks in that map gives
\eqref{eq:xdbl-thesis}.  The present theorem isolates the map as an
implementation circuit: \(AA,BB\) use two squarings, the product
\(AA\cdot BB\) and the final product defining \(Z_{2P}\) use two general
multiplications, and \(\alpha_{24}E\) uses one curve-constant
multiplication.  This proves the stated cost without repeating the earlier
rational-map derivation.
\end{proof}

\begin{corollary}[First-double saving]\label{cor:first-dbl-thesis}
For a native input \((X:Z)=(u+1:u)\), one has \(X-Z=1\), and
\[
 AA=(2u+1)^2,\quad E=AA-1,\quad
 (X_2:Z_2)=\bigl(AA:E(1+\alpha_{24}E)\bigr).
\]
The initialization costs exactly
\[
             \boxed{\M+\Sqr+\Dpar}.
\]
Here \(AA\) costs one squaring, \(\alpha_{24}E\) costs one
curve-constant multiplication, and the final product for \(Z_2\) costs
one general multiplication; \(X_2=AA\) is already available.
\end{corollary}

\begin{proof}
For the native representative \((X:Z)=(u+1:u)\), one has
\(X-Z=1\), and hence \(BB=(X-Z)^2=1\).  Substitution in
Theorem~\ref{thm:xdbl-thesis} gives
\[
 X_2=AA\cdot BB=AA,\qquad
 Z_2=E(BB+\alpha_{24}E)=E(1+\alpha_{24}E).
\]
Moreover \(X+Z=2u+1\), so \(AA=(2u+1)^2\) and
\(E=AA-BB=AA-1\).  Thus the dependency graph contains one squaring,
one multiplication by \(\alpha_{24}\), and one general multiplication.
No multiplication is needed for \(X_2\), proving the formula and count.
\end{proof}

\begin{remark}[Native-input initialization saving]
This saving belongs specifically to the original
\(\mathcal C_d\)-input interface: the two homogeneous coordinates
\((u+1:u)\) differ by \(1\) before any inversion or model conversion.
It is a one-time first-double optimization.  A general ladder state no
longer satisfies \(X-Z=1\), so subsequent doubles cost
\(2\M+2\Sqr+\Dpar\).  If \(\alpha_{24}=1/(16d)\) is selected as a
machine-level small constant, \(\Dpar\) may be cheaper in practice, but it
is retained in every symbolic count.
\end{remark}

\section{Bihomogeneous differential addition}

\begin{theorem}[Native \(x\)ADD]\label{thm:xadd-thesis}
Given
\[
\kappa(P)=(X_1:Z_1),\quad
\kappa(Q)=(X_2:Z_2),\quad
\kappa(P-Q)=(X_\Delta:Z_\Delta),
\]
put
\[
\begin{aligned}
A_1&=X_1+Z_1,&B_1&=X_1-Z_1,\\
A_2&=X_2+Z_2,&B_2&=X_2-Z_2,\\
C&=A_1B_2,&D&=A_2B_1.
\end{aligned}
\]
Then, on the oriented differential-addition variety,
\begin{equation}\label{eq:xadd-thesis}
 X_{P+Q}=Z_\Delta(C+D)^2,\qquad
 Z_{P+Q}=X_\Delta(C-D)^2.
\end{equation}
The general cost is \(4\M+2\Sqr\), and it is
\(3\M+2\Sqr\) when the known difference is affine.
\end{theorem}

\begin{proof}
Formula~\eqref{eq:xadd-thesis} is the generic chart of
Theorem~\ref{thm:native-odd-diff-atlas}.  Its derivation there starts from
Lemma~\ref{lem:native-Kummer-biquadratic}, and that theorem also proves
that the only omitted fibers are \(\Delta=O\) and \(\Delta=T\).
Accordingly, this implementation theorem does not repeat the base-locus
argument.  It only extracts the dependency graph: \(C,D\) cost two
multiplications, the two displayed squares cost \(2\Sqr\), and
multiplication by \(X_\Delta,Z_\Delta\) costs two further
multiplications.  When the known difference is affine, one of the last
two products is absorbed by its fixed scale.  This gives the two stated
costs.
\end{proof}

\begin{remark}
A differential addition law is not an everywhere-defined ternary operation
on arbitrary affine Kummer values.  Its correct domain is the oriented
subvariety consisting of triples that arise from \((P,Q,P-Q)\), and its
correct expression is bihomogeneous on \((\PP^1)^3\).
\end{remark}

\section{The standard ladder}

\begin{figure}[H]
\centering
\begin{tikzpicture}[x=1cm,y=1cm,>=Latex,every node/.style={font=\small}]
  \node at (0,1.65) {Native Kummer ladder on $\kappa(P)=(u+1:u)$};
  \node[draw,rounded corners,inner sep=5pt] (R0) at (-3.6,0.35) {$R_0=\kappa([n]P)$};
  \node[draw,rounded corners,inner sep=5pt] (R1) at (3.6,0.35) {$R_1=\kappa([n+1]P)$};
  \node[draw,rounded corners,inner sep=5pt] (B0) at (-3.6,-2.05) {$\kappa([2n]P)$};
  \node[draw,rounded corners,inner sep=5pt] (B1) at (3.6,-2.05) {$\kappa([2n+1]P)$};
  \node[draw,rounded corners,inner sep=4pt,align=center] (D) at (0,-0.95) {$\kappa(P)$ fixed\\known difference};
  \draw[->,thick] (R0) -- node[above,font=\scriptsize,fill=white,inner sep=1pt] {$x\mathrm{ADD}$} (R1);
  \draw[->,thick] (R0) -- node[left,font=\scriptsize,fill=white,inner sep=1pt] {$x\mathrm{DBL}$} (B0);
  \draw[->,thick] (R1) -- node[right,font=\scriptsize,fill=white,inner sep=1pt] {$x\mathrm{ADD}$} (B1);
  \draw[dashed] (-4.8,-2.8) -- (4.8,-2.8);
  \node at (0,-3.25) {a constant-time swap selects the ordered pair for the next bit};
\end{tikzpicture}
\caption{The standard native Kummer ladder maintains the adjacent pair $(\kappa([n]P),\kappa([n+1]P))$.  This is the same adjacency pattern that later reappears in the QRT state model of Chapter~\ref{ch:symmetric-QRT-envelope}.}
\label{fig:Cd-native-Kummer-ladder}
\end{figure}
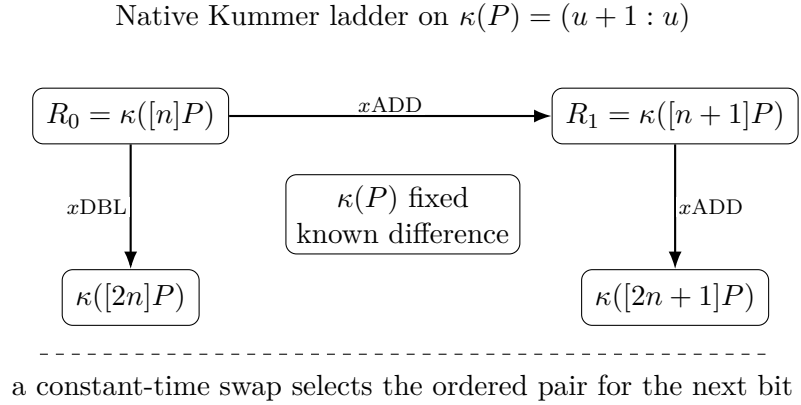

Combining Theorems~\ref{thm:xdbl-thesis} and~\ref{thm:xadd-thesis} gives
\[
\begin{array}{c|c}
\text{known difference}&x\mathrm{DBL}+x\mathrm{ADD}\\ \hline
\text{projective}&6\M+4\Sqr+\Dpar\\
\text{affine}&5\M+4\Sqr+\Dpar.
\end{array}
\]
For a secret scalar \(m=(m_{\ell-1}\cdots m_0)_2\), initialize
\[
 R_0=(1:0),\qquad R_1=(u_P+1:u_P).
\]
This is the inversion-free native initialization and uses the projective
difference row, hence \(6\M+4\Sqr+\Dpar\) per bit.  Alternatively, when
\(u_P\ne0\), one may pay one preprocessing inversion, form
\[
             U_P=1+u_P^{-1},\qquad R_1=(U_P:1),
\]
and then use the affine-difference row at
\(5\M+4\Sqr+\Dpar\) per bit.  The preprocessing inversion is not included
in that recurring cost.
At every bit, use a constant-time conditional swap and compute one
\(x\)DBL and one \(x\)ADD while maintaining
\[
 R_0=\kappa([n]P),\qquad R_1=\kappa([n+1]P).
\]
The fixed difference is always \(P\).  Secret-dependent branches and memory
addresses are forbidden.  Figure~\ref{fig:Cd-native-Kummer-ladder} shows the state pattern preserved by the native ladder.

\section{Linear transport on a Kummer line}

\begin{proposition}[Linear transport principle]\label{prop:linear-transport-thesis}
Let \(\kappa:E\to K\simeq\PP^1\) be a Kummer quotient and let
\(\lambda\in\operatorname{PGL}_2(k)\).  If \((F_0:F_1)\) is a homogeneous
doubling or oriented differential-addition formula for
\(\lambda\circ\kappa\), then
\[
 \lambda^{-1}\circ(F_0:F_1)\circ(\lambda,\ldots,\lambda)
\]
is the corresponding formula for \(\kappa\).  The transport introduces no
inversion and preserves degree and global domain.
\end{proposition}

\begin{proof}
Choose a matrix representative
\[
       L=\begin{pmatrix}a&b\\ c&d\end{pmatrix}\in
       \operatorname{GL}_2(k)
\]
for \(\lambda\).  Thus
\[
       \lambda(X:Z)=(aX+bZ:cX+dZ),
\]
and \(\lambda^{-1}\) is represented by the inverse matrix \(L^{-1}\),
up to a nonzero scalar.  Both maps are everywhere-defined automorphisms of
\(\PP^1\); in homogeneous coordinates they use only additions and
multiplications by field constants, not inversion of an input-dependent
quantity.

Let \(r=1\) for doubling and \(r=3\) for oriented differential addition,
and let \(V\subseteq(\PP^1)^r\) denote the appropriate input variety.  The
coordinate change sends it isomorphically to
\[
       V'=\lambda^{\times r}(V).
\]
By hypothesis, the homogeneous tuple \(F=(F_0:F_1)\) satisfies
\[
       F\bigl(\lambda\kappa(P_1),\ldots,
               \lambda\kappa(P_r)\bigr)
       =\lambda\kappa(R)
\]
on the dense domain on which it is nonzero, where \(R=2P_1\) in the
doubling case and \(R=P_1+P_2\) in the differential-addition case.
Applying \(\lambda^{-1}\) gives
\[
 \lambda^{-1}F\bigl(\lambda\kappa(P_1),\ldots,
                     \lambda\kappa(P_r)\bigr)=\kappa(R),
\]
which is exactly the transported formula in the statement.  Equality on a
dense open subset proves equality as rational maps, and hence at every point
where the transported tuple is defined.

It remains to verify the structural assertions.  Substitution of independent
linear forms for each input variable preserves the degree in every input
block.  Post-composition by \(L^{-1}\) takes invertible linear combinations
of the two output forms and therefore also preserves their multidegree.  If
the transported output were \((0,0)\) at an input \(w\in V\), invertibility
of \(L^{-1}\) would imply
\[
                   F(\lambda^{\times r}w)=(0,0).
\]
Conversely, if \(w\in\operatorname{Base}(F)\), then
\(F(\lambda^{\times r}(\lambda^{-\times r}w))=F(w)=(0,0)\), so
\(\lambda^{-\times r}w\in\operatorname{Base}(F^\lambda)\).  Hence
\[
 \operatorname{Base}(F^{\lambda})
   =(\lambda^{\times r})^{-1}\bigl(\operatorname{Base}(F)\bigr).
\]
In particular, an empty base locus, a prescribed exceptional fiber, or a
complete collection of charts is carried bijectively to the corresponding
domain for \(\kappa\).  No affine denominator is introduced at any stage.
The exact number of additions or constant multiplications may change with
the chosen matrix \(L\), but degree, inversion-freeness, and global domain
are preserved.
\end{proof}

\begin{corollary}
All degree-one Kummer coordinates for \(\Cd\) form one
\(\operatorname{PGL}_2(k)\)-orbit.  A different ``bilinear-then-square''
formula can select a cheaper basis for constants or squarings, but does not
create a new Kummer isomorphism class.
\end{corollary}

\begin{proof}
Theorem~\ref{thm:kummer-class-thesis} identifies the Kummer function field
with \(k(u)\) and shows that every degree-one generator has the form
\((au+b)/(cu+e)\), with \(ae-bc\ne0\).  Such substitutions are precisely
the action of \(\operatorname{PGL}_2(k)\) on the quotient line.  Applying
Proposition~\ref{prop:linear-transport-thesis} to any two of these bases
transports their doubling and differential-addition laws by an invertible
linear change.  Hence a basis can change a dependency graph and its
constant multipliers, but it represents the same quotient morphism up to a
projective-line automorphism.
\end{proof}

\chapter[The Reciprocal Chart]
{The Reciprocal Chart and Its Exposed Montgomery--Kummer Structure}
\label{ch:reciprocal-chart}

This chapter studies the equation
\begin{equation}\label{eq:reciprocal-chart-equation}
             \mathcal R_d:\qquad
             (u+1)(v+1)=du^2v^2.
\end{equation}
It is important to state at the outset how \(\mathcal R_d\) is being used.
It is not introduced as an unrelated elliptic-curve family, nor as a second
competitor to \(\mathcal C_d\).  It is the reciprocal affine chart of the
same marked \((2,2)\)-model.  The chart deserves a separate short chapter
because the native Kummer ratio becomes the affine linear function \(v+1\),
so the Montgomery quotient and its normalization are visible directly in
the defining coordinates.  This gives a useful input interface and several
low-cost differential schedules, while preserving the distinction between
quotient arithmetic and full-point arithmetic.

To avoid ambiguity, \((u_0,v_0)\) denotes the original native coordinates
on \(\mathcal C_d\) throughout this chapter, whereas \((u,v)\) denotes the
reciprocal coordinates on \(\mathcal R_d\).  The order of the two reciprocal
coordinates is chosen deliberately: it makes \(v+1\), rather than \(u+1\),
the marked Kummer coordinate.

\section{The reciprocal transformation on
\texorpdfstring{\(\PP^1\times\PP^1\)}{P1 x P1}}
\label{sec:reciprocal-global-map}

\begin{proposition}[The reciprocal chart as an ambient involution]
\label{prop:reciprocal-global-isomorphism}
Let \(d(1-16d)\ne0\).  On the dense affine locus define
\begin{equation}\label{eq:Cd-to-reciprocal-map}
 \Psi:\mathcal C_d\dashrightarrow\mathcal R_d,
 \qquad
 (u_0,v_0)\longmapsto
 (u,v)=\left(\frac1{v_0},\frac1{u_0}\right).
\end{equation}
Then \(\Psi\) extends to an isomorphism of the natural smooth
\((2,2)\)-completions.  More precisely, if
\[
 u_0=\frac{U_1}{U_0},\quad v_0=\frac{V_1}{V_0},
 \qquad
 u=\frac{R_1}{R_0},\quad v=\frac{S_1}{S_0},
\]
then the target completion is
\begin{equation}\label{eq:reciprocal-homogeneous-completion}
 R_0(R_1+R_0)S_0(S_1+S_0)=dR_1^2S_1^2,
\end{equation}
and the extended map is the ambient automorphism
\begin{equation}\label{eq:reciprocal-ambient-involution}
 \bigl((U_0:U_1),(V_0:V_1)\bigr)
 \longmapsto
 \bigl((V_1:V_0),(U_1:U_0)\bigr).
\end{equation}
It is an involution after the source and target factors are identified.
Consequently \(\mathcal R_d\) is smooth under exactly the same condition
\(d(1-16d)\ne0\), and the parameter \(d\) is unchanged.
\end{proposition}

\begin{proof}
On \(u_0v_0\ne0\), put \(u=v_0^{-1}\) and \(v=u_0^{-1}\).  Then
\[
 u_0^2+u_0=\frac{v+1}{v^2},
 \qquad
 v_0^2+v_0=\frac{u+1}{u^2}.
\]
Thus
\[
 (u_0^2+u_0)(v_0^2+v_0)=d
 \quad\Longleftrightarrow\quad
 (u+1)(v+1)=du^2v^2,
\]
which proves the affine identity.  Notice that every finite native point
has \(u_0v_0(u_0+1)(v_0+1)\ne0\), because \(d\ne0\).  Hence the displayed
reciprocals introduce no exceptional point on the finite native chart.

The map in \eqref{eq:reciprocal-ambient-involution} is the product of
coordinate inversion on each projective line and exchange of the two
factors.  It is therefore an everywhere-defined automorphism of
\(\PP^1\times\PP^1\).  Substituting
\[
 (R_0:R_1)=(V_1:V_0),
 \qquad
 (S_0:S_1)=(U_1:U_0)
\]
into \eqref{eq:homogeneous-thesis} gives exactly
\eqref{eq:reciprocal-homogeneous-completion}.  Applying the same operation
twice restores both projective factors.  The map therefore restricts to an
isomorphism of the two completions.  Smoothness and the value of the
parameter are preserved by this isomorphism, proving the final assertion.
\end{proof}

\begin{remark}[The reciprocal presentation as a chart]
Equation~\eqref{eq:reciprocal-chart-equation} is obtained by an ambient
\(\operatorname{PGL}_2\times\operatorname{PGL}_2\) transformation together
with factor exchange; it preserves the bidegree, the parameter, and the
marked smooth completion.  It is therefore a quotient-friendly internal
chart of the Cd model, with its own stored-input interface and normalization
advantages.  The reciprocal presentation exposes an affine Kummer
normalization while retaining the same marked \((2,2)\)-curve.
\end{remark}

\section{Boundary points, the identity, and marked torsion}
\label{sec:reciprocal-marked-points}

Transport the group law through \(\Psi\), so that \(\Psi\) is an
origin-preserving group isomorphism.  The marked points then have the
following coordinates.

\begin{proposition}[Marked-point dictionary and inversion]
\label{prop:reciprocal-marked-points}
Under \(\Psi\),
\begin{equation}\label{eq:reciprocal-marked-point-dictionary}
\begin{array}{c|c|c}
\text{point}&\mathcal C_d\text{ coordinates}&
             \mathcal R_d\text{ coordinates}\\ \hline
O&(0,\infty)&(0,\infty)\\
T&(-1,\infty)&(0,-1)\\
R&(\infty,0)&(\infty,0)\\
-R&(\infty,-1)&(-1,0)
\end{array}
\end{equation}
Thus the reciprocal affine chart contains \(T\) and \(-R\), while its two
points at infinity are \(O\) and \(R\).  The identity is still
\(O=(0,\infty)\), the point \(T=(0,-1)\) has order two, and
\(2R=T\).  The inverse map on the reciprocal chart is
\begin{equation}\label{eq:reciprocal-group-inverse}
             -(u,v)=\left(-\frac{u}{u+1},v\right),
\end{equation}
interpreted projectively when \(u=-1\).
\end{proposition}

\begin{proof}
The four point correspondences follow immediately from
\((u,v)=(v_0^{-1},u_0^{-1})\), with \(0^{-1}=\infty\) and
\(\infty^{-1}=0\) on \(\PP^1\).  They also show directly that the only
points outside the finite reciprocal chart are \(O\) and \(R\).  Orders and
the relation \(2R=T\) are preserved because the group law was transported
through \(\Psi\).

On \(\mathcal C_d\), inversion is
\((u_0,v_0)\mapsto(u_0,-1-v_0)\).  Hence the first reciprocal coordinate
changes by
\[
 \frac1{v_0}\longmapsto\frac1{-1-v_0}
 =-\frac{u}{u+1},
\]
while \(v=u_0^{-1}\) is fixed.  This proves
\eqref{eq:reciprocal-group-inverse}.  The fractional transformation is the
projective automorphism
\[
       (R_0:R_1)\longmapsto(R_0+R_1:-R_1),
\]
so the apparent denominator at \(u=-1\) is only an affine-chart boundary.
\end{proof}

\section{The conjugated eight-element symmetry}
\label{sec:reciprocal-dihedral-symmetry}

Put
\begin{equation}\label{eq:reciprocal-tau-involution}
                 \tau(z)=-\frac{z}{z+1}.
\end{equation}
This is a projective involution in every characteristic; in characteristic
two the minus sign disappears but the same formula still satisfies
\(\tau^2=1\).

\begin{proposition}[The eight reciprocal transformations]
\label{prop:reciprocal-eight-transformations}
The conjugate of the intrinsic \(D_8\)-action of
Theorem~\ref{thm:Cd-dihedral-symmetry} consists of the following eight
transformations, displayed one per line:
\begin{equation}\label{eq:reciprocal-eight-transformations}
\begin{aligned}
1:\quad &(u,v)\longmapsto(u,v),\\
\sigma:\quad &(u,v)\longmapsto(v,u),\\
\widehat\iota_u:\quad &(u,v)\longmapsto(\tau(u),v),\\
\sigma\widehat\iota_u:\quad &(u,v)\longmapsto(v,\tau(u)),\\
\widehat\iota_v:\quad &(u,v)\longmapsto(u,\tau(v)),\\
\sigma\widehat\iota_v:\quad &(u,v)\longmapsto(\tau(v),u),\\
\widehat\iota_u\widehat\iota_v:\quad
 &(u,v)\longmapsto(\tau(u),\tau(v)),\\
\sigma\widehat\iota_u\widehat\iota_v:\quad
 &(u,v)\longmapsto(\tau(v),\tau(u)).
\end{aligned}
\end{equation}
They extend to the smooth completion, preserve
\eqref{eq:reciprocal-chart-equation}, and form a faithful group isomorphic
to \(D_8\).  Their correspondence with the native generators is
\begin{equation}\label{eq:reciprocal-generator-conjugacy}
 \widehat\iota_u=\Psi\iota_v\Psi^{-1},\qquad
 \widehat\iota_v=\Psi\iota_u\Psi^{-1},\qquad
 \sigma=\Psi\sigma\Psi^{-1}.
\end{equation}
In odd characteristic,
\[
 \widehat\iota_u=[-1],\qquad
 \widehat\iota_v=\tau_T\circ[-1],\qquad
 \sigma=\tau_R\circ[-1].
\]
\end{proposition}

\begin{proof}
Conjugating the native transformation
\(v_0\mapsto-1-v_0\) changes \(u=v_0^{-1}\) to \(\tau(u)\) and fixes
\(v\); conjugating \(u_0\mapsto-1-u_0\) changes \(v\) to \(\tau(v)\)
and fixes \(u\).  Native coordinate exchange remains coordinate exchange
because \(\Psi\) already inverts and exchanges both factors symmetrically.
This proves \eqref{eq:reciprocal-generator-conjugacy}.  Composing the three
generators from right to left gives the eight lines in
\eqref{eq:reciprocal-eight-transformations}.

The maps extend projectively because \(\tau\in\operatorname{PGL}_2(k)\),
and they preserve the curve because they are conjugates of automorphisms of
\(\overline{\mathcal C}_d\).  Conjugation preserves all relations and
faithfulness, so the generated group is again \(D_8\).  The final
group-law identities are the conjugates of
\eqref{eq:dihedral-group-interpretation}.
\end{proof}

\section{The exposed Kummer coordinate \texorpdfstring{\(v+1\)}{v+1}}
\label{sec:reciprocal-kummer-coordinate}

\begin{theorem}[The reciprocal Kummer quotient]
\label{thm:reciprocal-kummer-coordinate}
Let \(K=k(\overline{\mathcal R}_d)\).  Then
\begin{equation}\label{eq:reciprocal-kummer-fixed-field}
                         K^{[-1]}=k(v).
\end{equation}
The boundary normalization sending \(O\) to infinity and \(T\) to zero is
\begin{equation}\label{eq:reciprocal-kummer-map}
 \kappa_{\mathcal R}(P)=(X:Z)=(v(P)+1:1).
\end{equation}
In homogeneous reciprocal coordinates it is
\begin{equation}\label{eq:reciprocal-kummer-homogeneous}
                         (X:Z)=(S_1+S_0:S_0).
\end{equation}
Under the reciprocal dictionary this coordinate is exactly the native
Kummer point:
\begin{equation}\label{eq:reciprocal-native-kummer-equality}
 (v+1:1)
 =\left(\frac1{u_0}+1:1\right)
 =(u_0+1:u_0)=\kappa_d(P).
\end{equation}
\end{theorem}

\begin{proof}
By \eqref{eq:reciprocal-group-inverse}, inversion fixes \(v\).  The
projection to the \(v\)-line has degree two because it is conjugate through
\(\Psi\) to the native projection to the \(u_0\)-line.  Its nontrivial deck
transformation is inversion, so the fixed-field theorem gives
\eqref{eq:reciprocal-kummer-fixed-field}.

Proposition~\ref{prop:reciprocal-marked-points} gives
\(v(O)=\infty\) and \(v(T)=-1\).  Thus \(v+1\) has the required pole and
zero on the quotient line, and its homogeneous representative is
\((S_1+S_0:S_0)\).  Finally \(v=u_0^{-1}\), and clearing the common
denominator \(u_0\) gives
\eqref{eq:reciprocal-native-kummer-equality}.
\end{proof}

\begin{remark}[Quotient information in reciprocal coordinates]
A finite reciprocal point supplies the affine Kummer representative
\((v+1:1)\) using one addition and no inversion.  The second full coordinate
\(u\) is still essential for the sign of the elliptic-curve point.  Hence
the reciprocal chart makes quotient normalization free; it does not turn a
Kummer output into a full point.
\end{remark}

\section{The odd-characteristic Montgomery dictionary}
\label{sec:reciprocal-odd-Montgomery}

Assume in this section that \(\charac k\ne2\).  Retain
\[
 \beta=\frac1{16d},\qquad
 A=\frac1{4d}-2.
\]

\begin{theorem}[Montgomery form exposed by the reciprocal chart]
\label{thm:reciprocal-Montgomery-dictionary}
On the dense affine locus \(u\ne0\), define
\begin{equation}\label{eq:reciprocal-to-Montgomery-map}
 U=v+1,\qquad
 V=2\frac{u+2}{u}(v+1)=2\left(1+\frac2u\right)U.
\end{equation}
Then
\begin{equation}\label{eq:reciprocal-Montgomery-equation}
             \beta V^2=U^3+AU^2+U.
\end{equation}
The inverse on the dense locus \(V-2U\ne0\) is
\begin{equation}\label{eq:Montgomery-to-reciprocal-map}
             v=U-1,\qquad
             u=\frac{4U}{V-2U}.
\end{equation}
The birational dictionary extends to an origin-preserving isomorphism of
smooth completions, and its first coordinate \(U=v+1\) is exactly
\(\kappa_{\mathcal R}\).
\end{theorem}

\begin{proof}
Put \(U=v+1\), so \(v=U-1\).  The reciprocal equation becomes
\begin{equation}\label{eq:reciprocal-relation-for-Montgomery-proof}
                 (u+1)U=du^2(U-1)^2.
\end{equation}
Using the definition of \(V\), multiply the desired identity
\eqref{eq:reciprocal-Montgomery-equation} by \(4du^2\).  Its left side
becomes \((u+2)^2U^2\), while its right side becomes
\begin{align*}
 4du^2\left(U(U-1)^2+\frac{U^2}{4d}\right)
 &=4du^2U(U-1)^2+u^2U^2\\
 &=\bigl(4(u+1)+u^2\bigr)U^2\\
 &=(u+2)^2U^2.
\end{align*}
The middle equality uses
\eqref{eq:reciprocal-relation-for-Montgomery-proof}.  Since \(du\ne0\)
on the stated affine locus, this proves
\eqref{eq:reciprocal-Montgomery-equation} without dividing by \(u+1\).

Solving \(U=v+1\) gives \(v=U-1\), while
\(V=2U+4U/u\) gives \(u=4U/(V-2U)\).  Thus the two maps are inverse on a
dense open set.  Smoothness extends them uniquely to the projective
completions.  Since \(v\) has a pole at
\(O=(0,\infty)\), the point \(O\) maps to the Montgomery point at infinity.
Theorem~\ref{thm:reciprocal-kummer-coordinate} identifies the reciprocal
Kummer point on this affine chart with \((v+1:1)\).  Since the present
coordinate is \(U=v+1\), it represents that same Kummer point; equality on
the smooth completion follows from equality on the dense affine chart.
\end{proof}

\begin{remark}[Endpoint costs]
From a reciprocal affine input, the quotient coordinate \(U=v+1\) is free
in the operation-count convention.  Forming the full coordinate \(V\)
costs \(\Inv+\M\): invert \(u\), form \(1+2u^{-1}\), and multiply by
\(U\); multiplication by the small integer \(2\) is uncharged.  Conversely,
recovering \((u,v)\) from a full Montgomery point costs \(\Inv+\M\) on the
dense chart.  These endpoint costs are absent when only the Kummer output
\(U\) is required.
\end{remark}

\section{The characteristic-two reduced Weierstrass dictionary}
\label{sec:reciprocal-binary-Weierstrass}

The phrase ``short Weierstrass form'' must be used carefully in
characteristic two.  The appropriate reduced binary equation retains the
\(XY\) term.

\begin{theorem}[Reduced binary Weierstrass form]
\label{thm:reciprocal-binary-Weierstrass}
Assume \(\charac k=2\) and \(d\ne0\).  On \(u\ne0\), put
\begin{equation}\label{eq:reciprocal-to-binary-W-map}
                 X=d(v+1),\qquad Y=\frac{X}{u}.
\end{equation}
Then
\begin{equation}\label{eq:reciprocal-binary-W-equation}
                 Y^2+XY=X^3+d^2X.
\end{equation}
The inverse on the dense locus \(dY\ne0\) is
\begin{equation}\label{eq:binary-W-to-reciprocal-map}
                 u=\frac{X}{Y},\qquad
                 v=\frac{X}{d}+1.
\end{equation}
The discriminant is \(d^4\), the \(j\)-invariant is \(d^{-4}\), and the
Weierstrass abscissa satisfies
\begin{equation}\label{eq:binary-W-exposed-kummer}
                         \frac Xd=v+1.
\end{equation}
\end{theorem}

\begin{proof}
Because \(Y=X/u\),
\[
 Y^2+XY=X^2\left(\frac1{u^2}+\frac1u\right)
          =X^2\frac{u+1}{u^2}.
\]
The reciprocal equation gives
\[
                  \frac{u+1}{u^2}
                  =\frac{dv^2}{v+1}.
\]
Since \(X=d(v+1)\), it follows that
\[
 Y^2+XY=d^2Xv^2.
\]
In characteristic two,
\[
 X^3+d^2X
 =d^2X\bigl((v+1)^2+1\bigr)
 =d^2Xv^2,
\]
which proves the equation.  Solving the definitions gives the inverse
formulas.  For the generalized Weierstrass coefficients one has
\[
 a_1=1,\qquad a_2=a_3=a_6=0,\qquad a_4=d^2.
\]
Consequently
\[
 b_2=1,\qquad b_4=b_6=0,\qquad b_8=-a_4^2=d^4
\]
because minus and plus agree in characteristic two.  The standard invariant
formulas therefore give
\[
 c_4=b_2^2-24b_4=1,
 \qquad
 \Delta=-b_2^2b_8=d^4,
 \qquad
 j=\frac{c_4^3}{\Delta}=d^{-4}.
\]
Finally, the definition \(X=d(v+1)\) gives \(X/d=v+1\), which is
\eqref{eq:binary-W-exposed-kummer}.
\end{proof}

\begin{remark}[Endpoint costs]
The quotient coordinate \(v+1\) again needs only an addition.  A full
forward map uses one multiplication by the curve constant \(d\), followed
by \(\Inv+\M\) to form \(Y=X/u\).  The affine inverse uses multiplication
by \(d^{-1}\) together with \(\Inv+\M\) for \(X/Y\).  As in odd
characteristic, the full-point dictionary has an inversion even though the
quotient dictionary does not.
\end{remark}

\section{General and mixed differential addition}
\label{sec:reciprocal-differential-addition}

The reciprocal chart supplies an affine normalization of the native quotient
law.  Its operation counts therefore record the stored-input and
normalization advantages of this interface on the same native Kummer line.

\begin{proposition}[Odd-characteristic reciprocal \(x\)ADD]
\label{prop:reciprocal-odd-xadd}
Assume \(\charac k\ne2\).  Let
\[
 K(P)=(X_1:Z_1),\quad K(Q)=(X_2:Z_2),\quad
 K(P-Q)=(X_\Delta:Z_\Delta)
\]
be reciprocal Kummer points.  Define
\[
\begin{aligned}
A_1&=X_1+Z_1,&B_1&=X_1-Z_1,\\
A_2&=X_2+Z_2,&B_2&=X_2-Z_2,\\
C&=A_1B_2,&D&=A_2B_1.
\end{aligned}
\]
Then the generic differential chart is
\begin{equation}\label{eq:reciprocal-odd-general-xadd}
 K(P+Q)=\bigl(Z_\Delta(C+D)^2:
                    X_\Delta(C-D)^2\bigr).
\end{equation}
Its general projective cost is \(4\M+2\Sqr\).  If the known difference is
the reciprocal affine input
\[
 K(P-Q)=(v_\Delta+1:1),
\]
the cost is \(3\M+2\Sqr\).

If all three quotient inputs arise from finite reciprocal coordinates
\(v_1,v_2,v_\Delta\), the same identity simplifies to
\begin{equation}\label{eq:reciprocal-odd-affine-xadd}
 K(P+Q)=
 \left(
   (v_1v_2+v_1+v_2)^2:
   (v_\Delta+1)(v_2-v_1)^2
 \right),
\end{equation}
at cost \(2\M+2\Sqr\).  Formula
\eqref{eq:reciprocal-odd-affine-xadd} is a fully affine-input differential
formula, not a recurring projective ladder step.
\end{proposition}

\begin{proof}
Equation~\eqref{eq:reciprocal-odd-general-xadd} is
\eqref{eq:xadd-thesis} written in the identical projective coordinate
identified by \eqref{eq:reciprocal-native-kummer-equality}.  The two products
\(C,D\), the two squares, and the two multiplications by the difference
coordinates give \(4\M+2\Sqr\).  Setting \(Z_\Delta=1\) removes one of
the last multiplications and gives \(3\M+2\Sqr\).

For three affine reciprocal quotient inputs, put
\(X_i=v_i+1\), \(Z_i=1\).  Then
\[
 C=(v_1+2)v_2,\qquad D=(v_2+2)v_1,
\]
so
\[
 C+D=2(v_1v_2+v_1+v_2),\qquad
 C-D=2(v_2-v_1).
\]
The common square factor \(4\) cancels projectively, proving
\eqref{eq:reciprocal-odd-affine-xadd}.  One multiplication forms
\(v_1v_2\), one forms the second output coordinate after its square, and
the two displayed squares give the stated cost.
\end{proof}

\begin{proposition}[Characteristic-two reciprocal \(x\)ADD]
\label{prop:reciprocal-binary-xadd}
Assume \(\charac k=2\).  Write
\[
 K(P)=(X_0:X_1),\qquad K(Q)=(Y_0:Y_1),\qquad
 K(P-Q)=(T_0:T_1)
\]
and put
\[
 A=X_0Y_0+X_1Y_1,\qquad
 B=X_0Y_1+X_1Y_0.
\]
Then
\begin{equation}\label{eq:reciprocal-binary-general-xadd}
                 K(P+Q)=(T_1A^2:T_0B^2)
\end{equation}
costs \(5\M+2\Sqr\) for a general projective difference.  For a finite
reciprocal difference outside the exceptional fibre \(\Delta=T\), let
\[
                 U_\Delta=v_\Delta+1,\qquad
                 K(P-Q)=(U_\Delta:1).
\]
Then the two equivalent mixed circuits are
\begin{align}
 K(P+Q)&=(A^2:U_\Delta B^2),
 &&3\M+2\Sqr+\mBase,
 \label{eq:reciprocal-binary-mixed-xadd-a}\\
 &=\bigl(U_\Delta B^2+d^{-1}C:B^2\bigr),
 &&3\M+\Sqr+\mBase+\mCurve,
 \label{eq:reciprocal-binary-mixed-xadd-b}
\end{align}
where
\[
                         C=(X_0Y_1)(X_1Y_0).
\]
Combining these circuits with binary \(x\)DBL gives the two recurring
reciprocal-input ladder costs
\begin{equation}\label{eq:reciprocal-binary-ladder-costs}
 \boxed{4\M+5\Sqr+\mBase+\mCurve},
 \qquad
 \boxed{4\M+4\Sqr+\mBase+2\mCurve}.
\end{equation}
\end{proposition}

\begin{proof}
The generic tuple is the binary differential chart
\eqref{eq:binary-native-diff-atlas-generic}, proved earlier in
Theorem~\ref{thm:binary-native-complete-diff-atlas}.  Karatsuba forms
\(A\) and \(B\) with the three products
\[
 X_0Y_0,\qquad X_1Y_1,
 \qquad (X_0+X_1)(Y_0+Y_1).
\]
Two squarings and the two final multiplications by \(T_1,T_0\) give
\(5\M+2\Sqr\).  If \((T_0:T_1)=(U_\Delta:1)\), direct substitution in
\eqref{eq:reciprocal-binary-general-xadd} gives
\((A^2:U_\Delta B^2)\), proving
\eqref{eq:reciprocal-binary-mixed-xadd-a} and its cost.

For the second circuit we derive the required oriented relation from the
Artin--Schreier calculation already used in
Lemma~\ref{lem:binary-native-Kummer-biquadratic}.  Put
\[
 x=\frac{X_1}{X_0},\qquad y=\frac{Y_1}{Y_0},
 \qquad z=\frac{T_1}{T_0}
\]
on the dense chart where the denominators are nonzero.  In the notation of
that lemma, \(U_1=x^{-1}\), \(U_2=y^{-1}\),
\(L=U_1+U_2\), and the Weierstrass abscissae \(H_+=X(P+Q)\),
\(H_-=X(P-Q)\) satisfy
\[
 H_-=H_++\frac{U_1U_2}{L^2},
 \qquad
 H_+H_-=d^2\frac{(1+U_1U_2)^2}{L^2}.
\]
Because the native Kummer ratio is \(t=d/X\), one has \(H_-=d/z\).
Eliminating \(H_+\) from the two displayed identities and substituting
\[
 \frac{U_1U_2}{L^2}=\frac{xy}{(x+y)^2},
 \qquad
 \frac{(1+U_1U_2)^2}{L^2}
   =\frac{(1+xy)^2}{(x+y)^2}
\]
gives, after multiplication by \(z^2(x+y)^2/d^2\),
\[
 (x+y)^2+z^2(1+xy)^2+d^{-1}zxy=0.
\]
Homogenizing this identity yields
\begin{equation}\label{eq:reciprocal-binary-oriented-local}
 T_0^2B^2+T_1^2A^2
 +d^{-1}T_0T_1X_0X_1Y_0Y_1=0.
\end{equation}
The oriented differential variety is the closure of the image of the
irreducible surface \(E\times E\) under
\((P,Q)\mapsto(\kappa(P),\kappa(Q),\kappa(P-Q))\); hence it is irreducible.
Since both sides are bihomogeneous, equality on the dense chart therefore
proves the identity everywhere on that variety.

Specializing \eqref{eq:reciprocal-binary-oriented-local} to
\((T_0:T_1)=(U_\Delta:1)\) gives
\[
 U_\Delta^2B^2+A^2
 +d^{-1}U_\Delta X_0X_1Y_0Y_1=0.
\]
With \(C=(X_0Y_1)(X_1Y_0)=X_0X_1Y_0Y_1\), this identity implies
\[
 A^2=U_\Delta\bigl(U_\Delta B^2+d^{-1}C\bigr).
\]
Replacing the first coordinate of \((A^2:U_\Delta B^2)\) by the right-hand
side and cancelling the common nonzero factor \(U_\Delta\) gives
\eqref{eq:reciprocal-binary-mixed-xadd-b}.  The three products
\(X_0Y_1\), \(X_1Y_0\), and their product \(C\) cost \(3\M\); the
remaining charged operations are exactly those displayed.

Finally, the doubling chart in
Theorem~\ref{thm:binary-native-complete-diff-atlas} is
\[
 (X_0:X_1)\longmapsto
 (X_0^4+X_1^4:d^{-1}X_0^2X_1^2).
\]
It is evaluated by one product, three squarings, and one multiplication by
\(d^{-1}\), hence costs \(\M+3\Sqr+\mCurve\).  Adding this cost to the
two mixed circuits proves \eqref{eq:reciprocal-binary-ladder-costs}.
\end{proof}

\begin{remark}[Differential completeness]
The generic formulas in this section retain the exceptional difference
fibres of the native formulas.  In odd characteristic the
\(\Delta=O,T\) charts of
Theorem~\ref{thm:native-odd-diff-atlas}, and in characteristic two the
charts of Theorem~\ref{thm:binary-native-complete-diff-atlas}, transport
through \(\Psi\) without changing their base loci.  The result is a complete
differential-addition atlas on the reciprocal presentation: the transported
generic chart together with the two exceptional charts
\(\Delta=O,T\) has empty common base locus.
\end{remark}

\section{A strict comparison with the native
\texorpdfstring{\(\mathcal C_d\)}{Cd} interface}
\label{sec:reciprocal-strict-cost-comparison}

The following comparison charges the coordinates actually supplied by the
external interface.  It therefore distinguishes a native input
\((u_0,v_0)\), a reciprocal input \((u,v)\), and a normalized Kummer input.

\paragraph{Odd characteristic.}
For a reciprocal input put \(W=v+1\).  Besides the general doubling circuit,
one may exploit \((X:Z)=(W:1)\) to write
\begin{equation}\label{eq:reciprocal-odd-first-double}
 K(2P)=
 \left((W^2-1)^2:
 4W\bigl((W-1)^2+4\alpha_{24}W\bigr)\right).
\end{equation}
This has two useful schedules:
\begin{equation}\label{eq:reciprocal-odd-first-double-costs}
          \M+2\Sqr+\Dpar
          \qquad\text{or}\qquad
          2\M+\Sqr+\Dpar.
\end{equation}
The first computes \(W^2\) and then \((W^2-1)^2\); the second puts
\(W=v+1\), reuses \(v^2\), and evaluates
\[
 \bigl((v+2)^2v^2:
       4(v+1)(v^2+4\alpha_{24}(v+1))\bigr).
\]

\begin{table}[H]
\centering
\footnotesize
\renewcommand{\arraystretch}{1.18}
\begin{tabularx}{\textwidth}{@{}L{2.65cm}X X@{}}
\toprule
operation or state & native \(\mathcal C_d\) input
                   & reciprocal-chart input\\
\midrule
Kummer lift
 & \((u_0+1:u_0)\), additions only; difference is projective
 & \((v+1:1)\), additions only; difference is already affine\\
first \(x\)DBL
 & \(\M+\Sqr+\Dpar\)
 & \(\M+2\Sqr+\Dpar\) or
   \(2\M+\Sqr+\Dpar\)\\
\(x\)ADD, fixed difference
 & \(4\M+2\Sqr\) without normalization
 & \(3\M+2\Sqr\) directly\\
\(x\)DBLADD, fixed difference
 & \(6\M+4\Sqr+\Dpar\) without normalization
 & \(5\M+4\Sqr+\Dpar\) directly\\
fully affine \(x\)ADD
 & requires first forming affine quotient values
 & \(2\M+2\Sqr\) by
   \eqref{eq:reciprocal-odd-affine-xadd}\\
normalization from a raw native input
 & \(\Inv\) gives \(U=1+u_0^{-1}\)
 & no cost only when the point is already stored reciprocally\\
\bottomrule
\end{tabularx}
\caption{Exact odd-characteristic interface comparison on the ordinary
Kummer line}
\label{tab:reciprocal-odd-cost-comparison}
\end{table}

Thus the native representative has the strictly cheaper specialized first
double: it saves one squaring relative to the first reciprocal schedule and
one multiplication relative to the second.  An input already stored in the
reciprocal chart supplies an affine fixed difference for free and saves one
multiplication in every ordinary ladder step.  After a native input has been
normalized once, the two interfaces use the same recurring Kummer circuit;
the native interface continues to retain its first-double advantage and its
direct interpretation on \(\mathcal C_d\).

\paragraph{Characteristic two.}
For a reciprocal input, the binary Kummer representative is
\((X_0:X_1)=(v+1:1)\).  The doubling chart in
Theorem~\ref{thm:binary-native-complete-diff-atlas} becomes
\begin{equation}\label{eq:reciprocal-binary-first-double}
 K(2P)=\bigl(v^4:d^{-1}(v^2+1)\bigr),
\end{equation}
which costs \(2\Sqr+\mCurve\): compute \(v^2\), square once more, and
reuse \(v^2+1=(v+1)^2\).  For a native input
\((X_0:X_1)=(u_0+1:u_0)\), characteristic two gives
\(X_0+X_1=1\), and the same doubling chart specializes to
\begin{equation}\label{eq:native-binary-first-double-local}
 \kappa(2P)=\bigl(1:d^{-1}(u_0^2+u_0)^2\bigr).
\end{equation}
One squaring forms \(u_0^2\), a second squares \(u_0^2+u_0\), and the
remaining charged operation is multiplication by \(d^{-1}\).  Thus the
native and reciprocal first doubles both cost \(2\Sqr+\mCurve\).

\begin{table}[H]
\centering
\footnotesize
\renewcommand{\arraystretch}{1.18}
\begin{tabularx}{\textwidth}{@{}L{2.65cm}X X@{}}
\toprule
operation or state & native \(\mathcal C_d\) input
                   & reciprocal-chart input\\
\midrule
Kummer lift
 & \((u_0+1:u_0)\), projective difference
 & \((v+1:1)\), affine difference\\
first \(x\)DBL
 & \(2\Sqr+\mCurve\)
 & \(2\Sqr+\mCurve\)\\
\(x\)ADD, projective difference
 & \(5\M+2\Sqr\)
 & the same general formula\\
\(x\)ADD, affine difference
 & requires one normalization inversion from raw native input
 & \(3\M+2\Sqr+\mBase\), or
   \(3\M+\Sqr+\mBase+\mCurve\) directly\\
\(x\)DBLADD, affine difference
 & same recurring count after one normalization inversion
 & \(4\M+5\Sqr+\mBase+\mCurve\), or
   \(4\M+4\Sqr+\mBase+2\mCurve\) directly\\
\bottomrule
\end{tabularx}
\caption{Exact characteristic-two interface comparison on the native
Kummer line}
\label{tab:reciprocal-binary-cost-comparison}
\end{table}

\begin{proposition}[Interface-level conclusion]
\label{prop:reciprocal-interface-conclusion}
The reciprocal chart has the following interface-dependent arithmetic properties.
\begin{enumerate}[label=\textup{(\roman*)}]
 \item If the external point is already represented on \(\mathcal R_d\),
       the fixed-difference Kummer coordinate is affine without inversion.
       The odd ladder therefore uses
       \(5\M+4\Sqr+\Dpar\) per bit, and the binary ladder uses either cost
       in \eqref{eq:reciprocal-binary-ladder-costs}.
 \item If the external point is given in native \(\mathcal C_d\)
       coordinates, forming only the reciprocal Kummer value costs one
       inversion and is exactly the already-known normalization
       \(U=1+u_0^{-1}\).  Forming both reciprocal affine coordinates costs
       two inversions, or \(\Inv+3\M\) by simultaneous inversion.
  \item The native input supplies the lower-cost specialized first double in
       odd characteristic, while the native and reciprocal inputs attain
       the same specialized first-double cost in characteristic two.

 \item The reciprocal chart is optimized for points already stored in
       quotient-friendly coordinates and for fully affine differential
       additions with an affine known difference.  The native chart retains
       its specialized initial doubling and its direct marked-boundary
       interface.  After the one-time Kummer normalization of a native
       input, the two charts use the same recurring differential-addition
       core, so an implementation can select the endpoint representation
       best suited to its surrounding protocol.
\end{enumerate}
\end{proposition}

\begin{proof}
Part~(i) is the mixed-difference count of
Propositions~\ref{prop:reciprocal-odd-xadd} and
\ref{prop:reciprocal-binary-xadd} combined with the corresponding doubling
costs.  For part~(ii), \(v=u_0^{-1}\), so the quotient value \(v+1\) needs
one inversion.  To form both \(u=v_0^{-1}\) and \(v=u_0^{-1}\), separate
inversion costs \(2\Inv\).  The usual two-input simultaneous inversion
forms \(u_0v_0\), inverts it, and multiplies twice, for
\(\Inv+3\M\).  For part~(iii), the odd reciprocal costs are
\eqref{eq:reciprocal-odd-first-double-costs}, whereas the native odd cost is
\(\M+\Sqr+\Dpar\) by Corollary~\ref{cor:first-dbl-thesis}.  In
characteristic two, equations
\eqref{eq:reciprocal-binary-first-double} and
\eqref{eq:native-binary-first-double-local} both cost
\(2\Sqr+\mCurve\).  This proves the comparison without appealing to a
later implementation chapter.  Part~(iv) is the resulting combination of
the endpoint costs and the recurring costs already proved in this section.
\end{proof}

\begin{example}[The reciprocal chart over \(\F_{101}\)]
\label{ex:reciprocal-F101}
On \(\mathcal C_1\), take the point \(P_0=(6,42)\) from
the full-orbit example in Chapter~\ref{ch:geometry}.  Since
\[
                    42^{-1}=89,\qquad 6^{-1}=17
                    \quad\text{in }\F_{101},
\]
its reciprocal image is
\[
                         P=(u,v)=(89,17).
\]
Indeed,
\[
 (u+1)(v+1)=90\cdot18\equiv4
 \equiv89^2\cdot17^2=du^2v^2\pmod{101}.
\]
The exposed Kummer coordinate is
\[
                 U=v+1=18
                 =\frac{6+1}{6},
\]
which is exactly the native Kummer ratio of \(P_0\).  Here
\(\beta=16^{-1}=19\) and \(A=4^{-1}-2=74\).  Formula
\eqref{eq:reciprocal-to-Montgomery-map} gives \(V=30\), and
\[
 19\cdot30^2\equiv31
 \equiv18^3+74\cdot18^2+18\pmod{101}.
\]
Thus the reciprocal point, the native quotient, and the Montgomery
abscissa agree exactly in a concrete finite-field computation.
\end{example}

The reciprocal presentation strengthens the internal arithmetic picture of
\(\mathcal C_d\).  The same smooth marked curve supports an inversion-free
native projective input with a specialized first-double saving and a
reciprocal affine input with a freely normalized fixed difference.  Thus
\[
 \mathcal R_d:(u+1)(v+1)=du^2v^2
\]
is a distinguished computational chart of
\(\mathcal C_d:(u^2+u)(v^2+v)=d\) in which the Montgomery--Kummer structure
is explicit.  The reciprocal chart is optimized for affine fixed differences
and stored quotient inputs, while the native chart retains the specialized
first double, the marked boundary, and the characteristic-uniform
input/output interface.  After normalization, both use the same recurring
differential-addition core.

\chapter[Arithmetic Tradeoffs and Torsion Quotients]
{Multiplication--Squaring Tradeoffs and Torsion Quotients}
\label{ch:tradeoffs}
Throughout this chapter,
\[
 \charac k\ne2,\qquad d(1-16d)\ne0,\qquad
 \rho=1-16d,\qquad \alpha_{24}=\frac1{16d}.
\]
Chapter~\ref{ch:odd-kummer} fixed the ordinary native Kummer map and its
baseline doubling and differential-addition circuits.  The present chapter
asks two logically separate optimization questions.  The first section
changes the multiplication--squaring balance without changing the quotient;
the second passes through the marked four-torsion quotient and therefore
changes the information retained by the output.  The final thresholds keep
these two kinds of tradeoff from being compared as though they had identical
semantics.

\section{Tradeoff on the ordinary Kummer line}

Assume that the known difference is affine.  Retain \(A_i,B_i,C,D\) from
Theorem~\ref{thm:xadd-thesis} and compute
\[
 X_{1+2}=(C+D)^2,\qquad
 Z_{1+2}=U_\Delta(C-D)^2.
\]

\begin{theorem}\label{thm:trade-seven-thesis}
Put \(\lambda=1-\alpha_{24}\), and define
\[
 E=A_1^2-B_1^2,\qquad
 F=A_1^4+B_1^4-E^2.
\]
Then
\begin{equation}\label{eq:trade-seven-thesis}
 X_{2P_1}=F,\qquad
 Z_{2P_1}=2(A_1^4-\lambda E^2)-F
\end{equation}
together with the preceding \(x\)ADD costs
\[
             3\M+7\Sqr+\Dpar.
\]
\end{theorem}

\begin{proof}
Write \(A=A_1\) and \(B=B_1\).  Expanding the square gives
\[
 F=A^4+B^4-(A^2-B^2)^2=2A^2B^2.
\]
The baseline Kummer doubling formula
\eqref{eq:xdbl-thesis} has first coordinate \(A^2B^2\).  Hence the
displayed first coordinate is exactly twice the baseline coordinate.
For the second coordinate, put \(E=A^2-B^2\) and use
\(\lambda=1-\alpha_{24}\).  Then
\begin{align*}
 2(A^4-\lambda E^2)-F
 &=2A^4-2(1-\alpha_{24})E^2-2A^2B^2\\
 &=2\bigl(A^4-A^2B^2-E^2+\alpha_{24}E^2\bigr)\\
 &=2\bigl(EA^2-E(A^2-B^2)+\alpha_{24}E^2\bigr)\\
 &=2E(B^2+\alpha_{24}E).
\end{align*}
This is twice the second baseline coordinate.  Thus both projective
coordinates have been multiplied by the same nonzero scalar \(2\), so
\eqref{eq:trade-seven-thesis} represents the same doubled Kummer point.

The differential-addition part uses \(3\M+2\Sqr\).  For the double,
the squares \(A^2,B^2,E^2,A^4,B^4\) cost \(5\Sqr\), and the product
\(\alpha_{24}E^2\) costs one \(\Dpar\); all remaining operations are
additions or multiplications by the small integer \(2\).  The combined
cost is therefore \(3\M+7\Sqr+\Dpar\).
\end{proof}

\begin{theorem}\label{thm:trade-six-thesis}
Suppose \(\rho=\gamma^2\) in \(k\).  Define
\[
\begin{aligned}
H_+&=(\gamma^{-1}A_1^2+B_1^2)^2,&
H_-&=(\gamma^{-1}A_1^2-B_1^2)^2,\\
G&=H_++H_-,&K&=H_+-H_-,\\
S&=\gamma K,&T&=\gamma^{-1}K.
\end{aligned}
\]
Then
\[
       X_{2P_1}=T-S,\qquad Z_{2P_1}=2G-S-T
\]
gives a mixed ladder step of cost
\(3\M+6\Sqr+3\Dpar\).
\end{theorem}

\begin{proof}
Write \(A=A_1\) and \(B=B_1\).  Since \(\rho=\gamma^2\), direct
expansion of the two squares gives
\begin{align*}
 G&=2\bigl(\gamma^{-2}A^4+B^4\bigr),\\
 K&=4\gamma^{-1}A^2B^2,\\
 S&=4A^2B^2,\\
 T&=4\gamma^{-2}A^2B^2.
\end{align*}
Consequently
\begin{equation}
 T-S=\frac{4(1-\rho)}{\rho}A^2B^2.
 \label{eq:trade-six-X-expanded}
\end{equation}
For the second coordinate,
\begin{equation}
\begin{aligned}
 2G-S-T
 &=4\gamma^{-2}A^4+4B^4
   -4A^2B^2-4\gamma^{-2}A^2B^2\\
 &=\frac4{\rho}
   \bigl(A^4-(1+\rho)A^2B^2+\rho B^4\bigr)\\
 &=\frac4{\rho}(A^2-B^2)(A^2-\rho B^2).
\end{aligned}
\label{eq:trade-six-Z-expanded}
\end{equation}
On the other hand, since
\(\alpha_{24}=1/(1-\rho)\), the baseline double is
\[
 (X_{2P_1}:Z_{2P_1})
 =\left(
 A^2B^2:
 \frac{(A^2-B^2)(A^2-\rho B^2)}{1-\rho}
 \right).
\]
Multiplying this pair by \(4(1-\rho)/\rho\) gives exactly
\eqref{eq:trade-six-X-expanded} and
\eqref{eq:trade-six-Z-expanded}.  The new pair therefore represents the
same doubled point.

The two input squares \(A^2,B^2\), followed by the two squares defining
\(H_+\) and \(H_-\), give \(4\Sqr\); together with the two output squares
from differential addition, this is \(6\Sqr\).  Differential addition
uses \(3\M\).  The fixed products by \(\gamma^{-1}\), \(\gamma\), and
\(\gamma^{-1}\) in the definitions of \(H_\pm,S,T\) give
\(3\Dpar\).  Hence the complete mixed step costs
\(3\M+6\Sqr+3\Dpar\), as asserted.
\end{proof}

\begin{remark}
The complete Edwards range requires \(\rho\) to be a nonsquare, whereas
Theorem~\ref{thm:trade-six-thesis} requires it to be a square.  The two
claims cannot be combined on the ordinary Kummer line.
\end{remark}

\section{The Kummer line after the four-torsion quotient}

Define
\begin{equation}\label{eq:omega-thesis}
 \omega(P)=\rho\xi^2\eta^2
 =\frac{\rho}{(2u+1)^2(2v+1)^2}.
\end{equation}
For the Edwards four-torsion point \(T_4=(1,0)\), translation by
\(\langle T_4\rangle\) and inversion preserve \(\omega\).  In
coordinates, the resulting eight-point orbit is
\[
 \{(\pm\xi,\pm\eta),(\pm\eta,\pm\xi)\}.
\]
Generically these eight points are distinct.  Indeed, if
\(a=\xi^2\), \(b=\eta^2\), and \(w=\omega(P)\), then
\begin{equation}\label{eq:omega-fiber-quadratic}
       a+b=1+w,\qquad ab=\frac{w}{\rho}.
\end{equation}
Thus \(a,b\) are the two roots of
\[
       T^2-(1+w)T+\frac{w}{\rho}=0.
\]
For a generic \(w\), the two root orderings and the four independent sign
choices give eight points.  Consequently
\[
             \deg(\omega)=8,
\]
and \(\omega\) determines
\[
       \{\pm P\}+\langle T_4\rangle,
\]
not the four-point coset \(P+\langle T_4\rangle\) and not the ordinary
Kummer pair \(\{P,-P\}\).  Equivalently, \(\omega\) is the ordinary
Kummer coordinate on the elliptic quotient
\(\mathcal C_d/\langle T_4\rangle\): the order-four isogeny has degree
four and the subsequent sign quotient has degree two.
Native initialization is inversion-free:
\[
       (W:Z)=\bigl(\rho:(rs)^2\bigr).
\]

\begin{theorem}[Fast steps on the post-isogeny Kummer line]\label{thm:omega-fast-thesis}
Let \(\omega_i=W_i/Z_i\), let the difference be affine, and set
\[
\begin{aligned}
A_1&=W_1+Z_1,&B_1&=W_1-Z_1,\\
A_2&=W_2+Z_2,&B_2&=W_2-Z_2,\\
C&=A_1B_2,&D&=A_2B_1,\\
E&=A_1^2-B_1^2,&F&=A_1^4+B_1^4-E^2.
\end{aligned}
\]
Then
\[
\begin{aligned}
W_{2P_1}&=2(A_1^4-\rho^{-1}E^2)-F,&Z_{2P_1}&=F,\\
W_{P_1+P_2}&=(C-D)^2,&
Z_{P_1+P_2}&=\omega_\Delta(C+D)^2
\end{aligned}
\]
costs \(3\M+7\Sqr+\Dpar\).

If \(d=c^2\), let \(\tau=4c\), so \(\tau^2=1-\rho\), and put
\[
\begin{aligned}
H_+&=(\tau A_1^2+B_1^2)^2,&
H_-&=(\tau A_1^2-B_1^2)^2,\\
G&=H_++H_-,&K&=H_+-H_-,\\
S&=\tau^{-1}K,&T&=\tau K.
\end{aligned}
\]
Replacing the double by
\[
       W_{2P_1}=2G-S-T,\qquad Z_{2P_1}=T-S
\]
gives \(3\M+6\Sqr+3\Dpar\).
\end{theorem}

\begin{proof}
On the Edwards model
\[
 x^2+y^2=1+\rho x^2y^2,
 \qquad w=\rho x^2y^2,
\]
the doubling formulas are
\[
 x(2P)=\frac{2xy}{1+w},
 \qquad
 y(2P)=\frac{y^2-x^2}{1-w}.
\]
Moreover,
\[
 (y^2-x^2)^2=(x^2+y^2)^2-4x^2y^2
 =(1+w)^2-\frac{4w}{\rho}.
\]
It follows that
\begin{equation}
 w(2P)=
 \frac{4w\bigl((1+w)^2-4w/\rho\bigr)}{(1-w^2)^2}.
 \label{eq:omega-affine-double-derived}
\end{equation}

Now put \(w=W/Z\), \(A=W+Z\), \(B=W-Z\), and
\(E=A^2-B^2\).  The identity
\[
 F=A^4+B^4-E^2=2A^2B^2=2(W^2-Z^2)^2
\]
shows that the displayed square-heavy denominator represents
\(2(1-w^2)^2\).  A direct expansion gives
\begin{align*}
 2(A^4-\rho^{-1}E^2)-F
 &=8WZ\left((W+Z)^2-\frac{4WZ}{\rho}\right).
\end{align*}
The ratio of these two homogeneous coordinates is therefore exactly
\eqref{eq:omega-affine-double-derived}.

For the six-square formula, use \(\tau^2=1-\rho\).  Expanding \(H_+\)
and \(H_-\) gives
\[
 G=2(\tau^2A^4+B^4),
 \qquad K=4\tau A^2B^2,
\]
and hence
\[
 S=4A^2B^2,\qquad T=4\tau^2A^2B^2.
\]
Thus \(T-S=-4\rho A^2B^2=-2\rho F\).  For the other coordinate,
\begin{align*}
 2G-S-T
 &=4(1-\rho)A^4+4B^4-4(2-\rho)A^2B^2\\
 &=-2\rho\bigl(2(A^4-\rho^{-1}E^2)-F\bigr).
\end{align*}
Both six-square coordinates are therefore obtained from the preceding
pair by the same nonzero scalar \(-2\rho\).  This proves the two doubling
identities without changing the represented Kummer point.

It remains to verify differential addition.  On the Edwards chart write
\(P_i=(x_i,y_i)\), set
\[
 w_i=\rho x_i^2y_i^2,\qquad
 w_\pm=w(P_1\pm P_2),
\]
and put
\[
 a=x_1^2,\quad b=y_1^2,\quad c=x_2^2,\quad e=y_2^2.
\]
Multiplying the Edwards addition and subtraction formulas gives
\[
 x_+x_-=\frac{ae-bc}{1-w_1w_2},\qquad
 y_+y_-=\frac{be-ac}{1-w_1w_2}.
\]
The source equations imply
\[
 a+b=1+w_1,\qquad c+e=1+w_2,\qquad
 \rho ab=w_1,\qquad \rho ce=w_2.
\]
Consequently,
\begin{align*}
 \rho(ae-bc)(be-ac)
 &=\rho ab(e^2+c^2)-\rho ce(a^2+b^2)\\
 &=w_1\left((1+w_2)^2-\frac{2w_2}{\rho}\right)
   -w_2\left((1+w_1)^2-\frac{2w_1}{\rho}\right)\\
 &=(w_1-w_2)(1-w_1w_2).
\end{align*}
It follows that
\[
            w_+w_-=
 \left(\frac{w_1-w_2}{1-w_1w_2}\right)^2.
\]
If \(w_i=W_i/Z_i\), the definitions of \(C,D\) give
\[
 C-D=2Z_1Z_2(w_2-w_1),\qquad
 C+D=2Z_1Z_2(w_1w_2-1).
\]
Since the fixed difference satisfies \(w_-=\omega_\Delta\), clearing
denominators yields
\[
 (W_{P_1+P_2}:Z_{P_1+P_2})
   =\bigl((C-D)^2:\omega_\Delta(C+D)^2\bigr).
\]
The derivation is valid on the dense open set on which the affine
denominators and the displayed output pair are nonzero.  Clearing
denominators gives a homogeneous identity, so it extends to the entire
nondegenerate differential-addition chart.

For the first schedule, the differential addition uses
\(3\M+2\Sqr\).  The double uses the five squares
\[
 A^2,\quad B^2,\quad A^4,\quad B^4,\quad E^2
\]
and one fixed multiplication by \(\rho^{-1}\), so the combined cost is
\(3\M+7\Sqr+\Dpar\).  In the second schedule the double uses the two
input squares \(A^2,B^2\), the two squares defining \(H_+,H_-\), and
the three fixed products by \(\tau,\tau^{-1},\tau\); together with the
differential addition this gives \(3\M+6\Sqr+3\Dpar\).
Finally, \(\tau^2=16d\), so the second schedule is defined over \(k\)
exactly when \(d\) is a square in \(k\).
\end{proof}

\begin{theorem}[Complete differential law on the post-isogeny Kummer line]
\label{thm:omega-complete-thesis}
Let \(k=\F_q\), let \(\rho\) be a nonsquare, and put
\(e_\omega=4/\rho\).
For oriented triples arising from \(P_1,P_2\in E_\rho(k)\), define
\[
\begin{aligned}
A_1&=W_1+Z_1,&B_1&=W_1-Z_1,&E_1&=A_1^2-B_1^2,\\
A_2&=W_1W_2,&B_2&=Z_1Z_2,&
E_2&=(W_1+Z_1)(W_2+Z_2),\\
C&=A_2+B_2,&D&=A_2-B_2.
\end{aligned}
\]
Then
\[
\begin{aligned}
W_{2P_1}&=A_1^4-B_1^4+(1-e_\omega/2)E_1^2,\\
Z_{2P_1}&=A_1^4+B_1^4-E_1^2,\\
W_{P_1+P_2}&=C(2E_2-e_\omega C)
 -(\omega_\Delta-e_\omega+2)D^2,\\
Z_{P_1+P_2}&=D^2
\end{aligned}
\]
is everywhere defined on all such rational oriented inputs and costs
\(5\M+6\Sqr+2\Dpar\).
\end{theorem}

\begin{proof}
The doubling pair is the square-heavy pair of
Theorem~\ref{thm:omega-fast-thesis}, because
\(1-e_\omega/2=1-2/\rho\).  Hence it represents
\(w(2P_1)\).  We derive the differential pair directly.

Work first on the affine Edwards chart and write
\(P_i=(x_i,y_i)\).  Put
\[
 a=x_1x_2,\qquad b=y_1y_2,\qquad
 c=x_1y_2,\qquad d=y_1x_2,\qquad h=\rho ab.
\]
The Edwards addition and subtraction laws give
\begin{align*}
 w(P_1+P_2)
 &=\frac{\rho(c+d)^2(b-a)^2}{(1-h^2)^2},\\
 w(P_1-P_2)
 &=\frac{\rho(c-d)^2(b+a)^2}{(1-h^2)^2}.
\end{align*}
Adding these two fractions and expanding their numerators gives
\begin{equation}
\begin{aligned}
 &(c+d)^2(b-a)^2+(c-d)^2(b+a)^2\\
 &\qquad=2(c^2+d^2)(a^2+b^2)-8a^2b^2.
\end{aligned}
\label{eq:omega-differential-expansion-one}
\end{equation}
Let \(w_i=\rho x_i^2y_i^2\), \(p=w_1w_2=h^2\), and
\(s=w_1+w_2\).  Since
\(x_i^2+y_i^2=1+w_i\), another direct expansion gives
\begin{equation}
 (c^2+d^2)(a^2+b^2)
 =\frac{(1+p)s+4p}{\rho}-\frac{4p}{\rho^2}.
 \label{eq:omega-differential-expansion-two}
\end{equation}
Substitution of
\eqref{eq:omega-differential-expansion-two} into
\eqref{eq:omega-differential-expansion-one}, followed by multiplication
by \(\rho/(1-p)^2\), yields
\begin{equation}
 w(P_1+P_2)+w(P_1-P_2)
 =\frac{2(1+p)s+8p-16p/\rho}{(p-1)^2}.
 \label{eq:omega-affine-differential-sum}
\end{equation}

In the affine specialization \(Z_1=Z_2=1\), the theorem's
intermediates satisfy
\[
 C=p+1,\qquad D=p-1,\qquad E_2=p+s+1.
\]
Using \(e_\omega=4/\rho\), expand the left-hand side:
\begin{align*}
 &C(2E_2-e_\omega C)+(e_\omega-2)D^2\\
 &\quad
 =2(p+1)(p+s+1)-e_\omega(p+1)^2
   +(e_\omega-2)(p-1)^2\\
 &\quad
 =2(1+p)s+8p-4e_\omega p\\
 &\quad
 =2(1+p)s+8p-\frac{16p}{\rho}.
\end{align*}
Thus
\[
 C(2E_2-e_\omega C)+(e_\omega-2)D^2
 =2(1+p)s+8p-16p/\rho.
\]
Equation~\eqref{eq:omega-affine-differential-sum} therefore becomes
\[
 w(P_1+P_2)
 =\frac{C(2E_2-e_\omega C)
 -(w(P_1-P_2)-e_\omega+2)D^2}{D^2},
\]
which is exactly the displayed projective differential pair.  Clearing
the two input denominators gives a polynomial identity on a dense open
subset of the variety parametrized by \((P_1,P_2)\).  The identity
therefore holds identically and extends to every point at which the
projective output pair is nonzero.

It remains to exclude a rational base point.  The doubling denominator is
\(Z_{2P_1}=2A_1^2B_1^2\).  If \(A_1=0\), substitution into its numerator
gives \(-2B_1^4/\rho\ne0\).  If \(B_1=0\), the numerator is
\(2(1-1/\rho)A_1^4\ne0\), because a nonsquare \(\rho\) is neither zero
nor one.  Thus the doubling pair is base-point-free.

Because \(\rho\) is a nonsquare, the \(k\)-rational points of the complete
Edwards model lie in its affine chart.  Indeed, at \(x=\infty\) the
bihomogeneous Edwards equation reduces to
\(Y_0^2=\rho Y_1^2\), and at \(y=\infty\) it reduces to
\(X_0^2=\rho X_1^2\); neither equation has a projective \(k\)-solution.
We may therefore put \(Z_1=Z_2=1\) when
checking a rational base point of the differential pair.  A common zero
would give
\[
 w_1w_2=1,\qquad w_1+w_2=e_\omega-2.
 \label{eq:omega-exceptional-pair}
\]
For a rational lift \(P_i=(x_i,y_i)\), the two values
\(x_i^2,y_i^2\) are roots of
\[
 T^2-(1+w_i)T+\frac{w_i}{\rho}=0.
\]
The two values \(w_1,w_2\) are the roots of
\[
 W^2-(e_\omega-2)W+1=0.
\]
Consequently the discriminant of the preceding quadratic in \(T\) is
\[
 (1+w_i)^2-\frac{4w_i}{\rho}
 =w_i^2+(2-e_\omega)w_i+1=0.
\]
Hence \(x_i^2=y_i^2=a_i=(1+w_i)/2\), and substitution in the Edwards
equation gives \(\rho a_i^2-2a_i+1=0\).  Moreover,
\[
 a_1a_2
 =\frac{(1+w_1)(1+w_2)}4
 =\frac{2+w_1+w_2}{4}
 =\frac1{\rho}.
\]
Both \(a_1\) and \(a_2\) are squares in \(k\), whereas \(1/\rho\) is a
nonsquare, a contradiction.  The differential pair is therefore defined
on every rational oriented input.

Finally, the three products \(W_1W_2\), \(Z_1Z_2\), and
\((W_1+Z_1)(W_2+Z_2)\), the product by \(C\), and the product
\(\omega_\Delta D^2\) give \(5\M\).  The four powers used by the
double, the square \(E_1^2\), and \(D^2\) give \(6\Sqr\); the occurrences
of \(e_\omega/2\) and \(e_\omega\) give \(2\Dpar\).  This proves the
claimed cost and completes the proof.
\end{proof}

\section{Cost thresholds}

With \(\Sqr=\gamma_{\rm S}\M\),
\(\Dpar=\gamma_{\rm D}\M\), the standard, seven-square,
six-square, and complete-quotient costs are
\[
\begin{aligned}
C_0&=5+4\gamma_{\rm S}+\gamma_{\rm D},&
C_1&=3+7\gamma_{\rm S}+\gamma_{\rm D},\\
C_2&=3+6\gamma_{\rm S}+3\gamma_{\rm D},&
C_{\rm comp}&=5+6\gamma_{\rm S}+2\gamma_{\rm D}.
\end{aligned}
\]
Therefore
\[
 C_1<C_0\iff\gamma_{\rm S}<\frac23,\qquad
 C_2<C_0\iff\gamma_{\rm S}+\gamma_{\rm D}<1.
\]
The complete quotient law is symbolically more expensive for positive
costs; its purpose is to remove exceptional inputs.  Any end-to-end use of
\(\omega\) must also account for the eight-point ambiguity
\(\{\pm P\}+\langle T_4\rangle\).

\chapter{Halving and Complete Scalar Multiplication}
\label{ch:halving}
Throughout this chapter,
\[
 \mathcal C_d:\quad (u^2+u)(v^2+v)=d,\qquad
 \charac k\ne2,\qquad d(1-16d)\ne0,
\]
and
\[
             \beta=\frac1{16d},\qquad A=\frac1{4d}-2.
\]
The divisions by \(2\) and the separability of the doubling map used below
are the reasons for keeping the odd-characteristic hypothesis explicit.

\section{Halving through two quadratic equations}

\begin{theorem}[Native halving: the abscissa stage]
\label{thm:halving-thesis}
Let \(Q=(u_Q,v_Q)\) be a finite point of
\[
       \mathcal C_d:\quad (u^2+u)(v^2+v)=d
\]
distinct from the identity, and define the native pullback quantities
\[
 q=\frac{u_Q+1}{u_Q},\qquad
 V_Q=2(2v_Q+1)q.
\]
For every choice of roots
\[
\begin{aligned}
R^2&=q^2+Aq+1,\\
T&=2q+2R,\qquad S^2=T^2-4,\\
U&=\frac{T+S}{2},
\end{aligned}
\]
the value \(U\) is the Kummer coordinate \((u_P+1)/u_P\) of a solution
of \(2P=Q\).  Allowing the two signs of \(R\) and \(S\) gives all four
geometric preimages, with multiplicity.
\end{theorem}

\begin{proof}
The duplication equation is
\[
 q=\frac{(U^2-1)^2}{4U(U^2+AU+1)}.
\]
For \(U\ne0\), divide by \(U^2\) and set \(T=U+U^{-1}\).  This gives
\[
 T^2-4=4q(T+A),
\]
or
\[
 T^2-4qT-4(Aq+1)=0.
\]
The quadratic formula gives \(T=2q\pm2\sqrt{q^2+Aq+1}\).  The equation
\(U+U^{-1}=T\) is \(U^2-TU+1=0\), whence
\(U=(T\pm\sqrt{T^2-4})/2\).  For every value obtained in this way,
\(U\ne0\) because the constant term of \(U^2-TU+1\) is one, and
\[
  U+U^{-1}=T.
\]
Hence
\[
  (U^2-1)^2=U^2(T^2-4),
  \qquad
  U(U^2+AU+1)=U^2(T+A).
\]
Moreover \(T+A\ne0\).  Indeed, if \(T=-A\), then
\(T^2-4=4q(T+A)\) gives \(A^2-4=0\), contrary to smoothness.  Therefore
\[
 \frac{(U^2-1)^2}{4U(U^2+AU+1)}
 =\frac{T^2-4}{4(T+A)}=q,
\]
which verifies the duplication equation directly.  In particular no
resulting value has \(f(U)=U(U^2+AU+1)=0\).  Thus, over \(\bar k\), there are
two points above each candidate \(U\).  Their doubles have Kummer
coordinate \(q\) and are negatives of one another.  They are therefore
\(Q\) and \(-Q\), with the two values coinciding when \(Q\) is
two-torsion.  Hence every candidate is the Kummer coordinate of a half
of \(Q\).  Finally, the
duplication map on the Kummer line has degree four, so these roots, counted
with multiplicity, form the complete list.
\end{proof}

\begin{proposition}[Geometric meaning of the four branches]
\label{prop:halving-branch-structure}
Over \(\bar k\), if \(P_0\) is one half of \(Q\), then
\[
             [2]^{-1}(Q)=P_0+\mathcal C_d[2].
\]
If \(Q\notin\mathcal C_d[2]\), the four branches in
Theorem~\ref{thm:halving-thesis} give four distinct Kummer coordinates.
If \(Q\in\mathcal C_d[2]\setminus\{O\}\), they give two Kummer
coordinates, each with multiplicity two.
\end{proposition}

\begin{proof}
Because \(\charac k\ne2\), the morphism \([2]\) is separable of degree
four and its kernel is \(\mathcal C_d[2]\).  Every nonempty geometric fiber
is therefore the displayed torsor.  Two points in the fiber have the same
ordinary Kummer coordinate exactly when they are negatives.  If both \(P\)
and \(-P\) double to \(Q\), then \(Q=-Q\).  Hence no identification occurs
when \(Q\notin\mathcal C_d[2]\), while a nonzero two-torsion target pairs
the four halves into two inverse pairs.  This also explains the
``with multiplicity'' clause in Theorem~\ref{thm:halving-thesis}.
\end{proof}

\begin{corollary}[Rationality criterion for a Kummer half]
\label{cor:halving-criterion-thesis}
For a finite point \(Q=(u_Q,v_Q)\in\mathcal C_d(k)\), put
\(q=(u_Q+1)/u_Q\).  At least one Kummer preimage \(U\) belongs to \(k\)
if and only if \(q^2+Aq+1\) is a square in \(k\) and, for at least one
root \(R\), \((2q+2R)^2-4\) is a square in \(k\).
\end{corollary}

\begin{proof}
Suppose first that a Kummer preimage \(U\in k\) exists.  Since \(U\ne0\),
\[
             T=U+U^{-1}\in k.
\]
The first quadratic in the proof of Theorem~\ref{thm:halving-thesis} is
\[
             T^2-4qT-4(Aq+1)=0,
\]
whose discriminant is
\[
             16(q^2+Aq+1).
\]
Because \(4\in k^\times\), the existence of \(T\in k\) implies that
\(q^2+Aq+1\) is a square in \(k\).  Writing
\(T=2q+2R\) for one of its two square roots, the equation
\(U+U^{-1}=T\) is
\[
             U^2-TU+1=0,
\]
with discriminant \(T^2-4=(2q+2R)^2-4\).  Since \(U\in k\), this second
discriminant is also a square in \(k\).

Conversely, suppose that \(q^2+Aq+1=R^2\) for some \(R\in k\) and that,
for one choice of \(R\), the element
\[
             S^2=(2q+2R)^2-4
\]
has a square root \(S\in k\).  Put
\[
             T=2q+2R,
             \qquad
             U=\frac{T+S}{2}.
\]
Then \(U\in k\) and \(U^2-TU+1=0\).  Its constant term is one, so
\(U\ne0\) and division by \(U\) gives \(U+U^{-1}=T\).  The identity
\(R^2=q^2+Aq+1\) implies
\[
\begin{aligned}
 T^2-4qT-4(Aq+1)
 &=4\bigl((q+R)^2-2q(q+R)-Aq-1\bigr)\\
 &=4\bigl(R^2-q^2-Aq-1\bigr)=0.
\end{aligned}
\]
Thus \(T^2-4=4q(T+A)\).  If \(T+A=0\), then this equality gives
\(T^2=4\), hence \(A^2=4\), contrary to the smoothness condition
\(A^2-4\ne0\).  Consequently
\[
\begin{aligned}
 \frac{(U^2-1)^2}{4U(U^2+AU+1)}
 &=\frac{U^2(T^2-4)}{4U^2(T+A)}\\
 &=q.
\end{aligned}
\]
This is exactly the Kummer duplication equation of
Theorem~\ref{thm:halving-thesis}; hence \(U\) is the Kummer coordinate of
a half of \(Q\).
\end{proof}

\begin{proposition}[Inversion-free full native recovery of a half]
\label{prop:halving-full-thesis}
Put
\[
 f(U)=U^3+AU^2+U,\qquad g(U)=3U^2+2AU+1.
\]
For a candidate \(U\) and \(V_Q\ne0\), define
\begin{equation}\label{eq:halving-projective-blocks}
       N_U=g(U)(U-q)-2f(U),\qquad
       \Delta_Q=2\beta V_Q.
\end{equation}
Then
\begin{equation}\label{eq:halving-V-thesis}
                  V=\frac{N_U}{\Delta_Q}
\end{equation}
is the unique auxiliary ordinate for the chosen Kummer root.  The full
half can be returned without an inversion in the native projective pairs
\begin{equation}\label{eq:halving-native-projective-recovery}
 \boxed{\qquad
 (U_0:U_1)=(U-1:1),\qquad
 (V_0:V_1)=(4U\Delta_Q:N_U-2U\Delta_Q).
 \qquad}
\end{equation}
If affine output is required, this becomes
\[
        u=\frac1{U-1},\qquad v=\frac{V}{4U}-\frac12.
\]
\end{proposition}

\begin{proof}
The tangent slope is \(\lambda=g(U)/(2\beta V)\), and the ordinate
addition formula is \(V_Q=\lambda(U-q)-V\).  Multiply by \(2\beta V\)
and use \(\beta V^2=f(U)\); solving the resulting linear equation in \(V\)
formally gives \eqref{eq:halving-V-thesis}.  We now verify the reverse
implication, which is needed because \(U\) was obtained only from the
abscissa-duplication equation.

For a candidate root \(U\), that equation is equivalent to
\begin{equation}\label{eq:halving-g-square-identity}
       g(U)^2=4f(U)(q+A+2U).
\end{equation}
Indeed, the tangent formula on
\(\beta V^2=f(U)\) gives
\[
 q=\beta\left(\frac{g(U)}{2\beta V}\right)^2-A-2U
   =\frac{g(U)^2}{4f(U)}-A-2U,
\]
and clearing the denominator gives
\eqref{eq:halving-g-square-identity}.  Conversely, assume
\eqref{eq:halving-g-square-identity}.  The argument below shows that
\(f(U)\ne0\); division by \(4f(U)\) then gives
\[
  q+A+2U=\frac{g(U)^2}{4f(U)},
\]
which is precisely the Montgomery abscissa-doubling equation.  Thus the
square identity and the abscissa equation are equivalent for every candidate
root.  No candidate satisfying
\eqref{eq:halving-g-square-identity} has \(f(U)=0\): smoothness gives
\(\gcd(f,g)=1\), whereas \(f(U)=0\) in that identity would also force
\(g(U)=0\).

The cubic Taylor identity
\[
 f(q)=f(U)+g(U)(q-U)
        +(q-U)^2(q+2U+A)
\]
together with \eqref{eq:halving-g-square-identity} now yields
\begin{equation}\label{eq:halving-norm-identity}
 \bigl(g(U)(U-q)-2f(U)\bigr)^2=4f(U)f(q).
\end{equation}
Since \(Q\) lies on the curve,
\(f(q)=\beta V_Q^2\).  Therefore the value
\[
       V=\frac{N_U}{2\beta V_Q}
\]
satisfies
\[
       \beta V^2
       =\frac{N_U^2}{4\beta V_Q^2}
       =\frac{4f(U)f(q)}{4f(q)}
       =f(U).
\]
It is consequently an ordinate of a point on the auxiliary Montgomery
curve.  Moreover
\[
       2\beta VV_Q=N_U
       =g(U)(U-q)-2\beta V^2,
\]
which rearranges to
\[
       V_Q=\frac{g(U)}{2\beta V}(U-q)-V.
\]
Thus its tangent double has both coordinates \((q,V_Q)\).  Because
\(V_Q\ne0\), this last equation is linear with nonzero coefficient
\(2\beta V_Q\), so the ordinate is unique for the chosen \(U\).

The native formulas are
\eqref{eq:mont-inverse-thesis}.  Before dividing by \(\Delta_Q\), they give
\[
 v=\frac{N_U}{4U\Delta_Q}-\frac12
   =\frac{N_U-2U\Delta_Q}{4U\Delta_Q},
\]
which proves the projective pairs
\eqref{eq:halving-native-projective-recovery}.
\end{proof}

\begin{corollary}[Full rational-halving criterion]
\label{cor:full-rational-halving-Cd}
Suppose \(Q=(u_Q,v_Q)\in\mathcal C_d(k)\) is finite, put
\(q=(u_Q+1)/u_Q\) and \(V_Q=2(2v_Q+1)q\), and assume \(V_Q\ne0\).  Then
\(Q\in[2]\mathcal C_d(k)\) if and only if the two-stage square conditions
of Corollary~\ref{cor:halving-criterion-thesis} hold for at least one
branch.  Every successful branch gives a full \(k\)-rational half through
\eqref{eq:halving-native-projective-recovery}.
\end{corollary}

\begin{proof}
A \(k\)-rational half has a \(k\)-rational Kummer coordinate, so necessity
follows from Corollary~\ref{cor:halving-criterion-thesis}.  Conversely, a
successful square-root branch gives \(U\in k\).  Since
\(q,V_Q,\beta,A\in k\) and \(V_Q\ne0\), both entries in each projective pair
of \eqref{eq:halving-native-projective-recovery} lie in \(k\), and neither
pair is zero.  The tangent calculation in
Proposition~\ref{prop:halving-full-thesis} then verifies that the resulting
point doubles to \(Q\).
\end{proof}

\begin{remark}[The two-torsion target]
The condition \(V_Q=0\) means that the finite Montgomery image of \(Q\) is
a nonzero two-torsion point.  Formula~\eqref{eq:halving-V-thesis} cannot
select an ordinate in that case because both \(P\) and \(-P\) double to the
same \(Q\).  For each Kummer candidate \(U\in k\), a full rational lift
exists exactly when
\[
                   \frac{f(U)}{\beta}
\]
is a square in \(k\).  Choosing either square root \(V\), converting by
\eqref{eq:mont-inverse-thesis}, and verifying \([2]P=Q\) handles this
exceptional fiber.  Proposition~\ref{prop:halving-branch-structure} explains
why the two signs have the same double.
\end{remark}

\subsection{Cost and appropriate use of halving}

Starting from a native affine point, forming
\(q=1+u_Q^{-1}\) costs one inversion; if a normalized Kummer input \(q\)
is already available, this cost disappears.  The first discriminant uses
\(\Sqr+\Dconst{A}\).  Enumerating all four Kummer candidates then requires
one square-root extraction for \(R\), at most two further squarings, and at
most two further square-root extractions for the two values of \(T\).
Each retained affine value \(U=(T+S)/2\) also uses one
\(\Dconst{1/2}\).
Square-root extraction is field- and representation-dependent and is not
hidden inside the \(\M,\Sqr,\Dpar\) count.

When full recovery is requested, forming \(V_Q\) costs one multiplication.
For each successful \(U\), the blocks
\eqref{eq:halving-projective-blocks} and the inversion-free native output
\eqref{eq:halving-native-projective-recovery} can be formed in
\[
                  3\M+\Sqr+\Dconst{A},
\]
in addition to one shared \(\Dconst{\beta}\).  An explicit schedule first
forms \(U^2\) and \(AU\), then reuses them in
\[
 f(U)=U(U^2+AU+1),\qquad g(U)=3U^2+2AU+1.
\]
The three general products are the product defining \(f(U)\), the product
in \(N_U\), and the shared product \(U\Delta_Q\) used in both entries of
the second native projective pair.  Additions and small
integer multiples are omitted as in Chapter~\ref{ch:conventions}.  Affine
normalization is a separate endpoint cost and can be batched across several
candidates.  Because square roots and inversions dominate on many platforms,
halving is principally useful for divisibility tests, decompression,
precomputation, and protocols that already require a half; it is not a
generic replacement for a regular secret-scalar ladder.

\section{Complete scalar multiplication from the native Segre law}

Chapter~\ref{ch:native-completeness} proves, directly on the smooth
\((2,2)\)-completion of \(\mathcal C_d\), the centered native Segre
equations
\[
       XY=ZT,\qquad X^2+Y^2=Z^2+\rho T^2
\]
and the native Segre addition tuple
\[
       (X_3:Y_3:T_3:Z_3)=(EF:GH:EH:FG).
\]
This is Theorem~\ref{thm:native-complete-addition-Cd}; over a finite base
field its exact completeness criterion is the nonsquare condition stated in
Theorem~\ref{thm:exact-native-completeness-criterion}.  It is not an
instruction to convert the scalar-multiplication state to an Edwards curve.
The native
input lift and output recovery are
\[
\begin{aligned}
(u,v)&\longmapsto
(2u+1:2v+1:1:(2u+1)(2v+1)),\\
(U_0:U_1)&=(2Y:Z-Y),\qquad
(V_0:V_1)=(2X:Z-X).
\end{aligned}
\]
The finite input lift costs one multiplication, namely the product
\((2u+1)(2v+1)\), and no inversion.  The displayed output recovery already
returns two native projective pairs and uses only additions.  If affine
\((u,v)\) output is required, simultaneously normalizing both pairs costs
one inversion and five multiplications by the standard product-inversion
trick.

When \(k=\F_q\) and \(\rho\) is a nonsquare,
double-and-add-always costs
\[
       13\M+4\Sqr+\Dpar
\]
per bit.  More precisely, a uniform \(\ell\)-round schedule has the loop
bound
\begin{equation}\label{eq:complete-double-add-total-cost}
       C_{\rm daa}(\ell)
       =\ell(13\M+4\Sqr+\Dpar),
\end{equation}
before its one-time input and output costs.  For a width-\(w\) table that
includes the identity, the conservative single-law bound is
\begin{equation}\label{eq:complete-window-total-cost}
 C_{\rm win}^{\rm gen}(\ell,w)
 \le \ell(4\M+4\Sqr)
    +\left\lceil\frac{\ell}{w}\right\rceil(9\M+\Dpar),
\end{equation}
or asymptotically
\[
       4\M+4\Sqr+\frac{9\M+\Dpar}{w}
\]
per bit, excluding amortized precomputation and table selection.  The reason
for using the general \(9\M+\Dpar\) addition count here is that the identity
has \(T=0\) in the centered native Segre embedding and cannot be normalized
to the mixed condition \(T_2=1\).

For nonzero table digits, precomputed points with \(T_2=1\) use the mixed
cost \(8\M+\Dpar\).  If \(N_{\ne0}\) windows are nonzero, a mixed schedule
therefore has arithmetic part
\begin{equation}\label{eq:complete-window-mixed-qualified}
 \ell(4\M+4\Sqr)+N_{\ne0}(8\M+\Dpar)+C_{0},
\end{equation}
where \(C_0\) is the cost of a separately specified constant-time zero-digit
treatment.  It is legitimate to set \(C_0=0\) only for a public or
variable-time recoding that skips zero windows.  Completeness removes
exceptional group inputs from the chosen addition graph; it does not by
itself make table selection or zero-digit handling constant time.

The exact necessity and sufficiency of the nonsquare condition are proved
in Theorem~\ref{thm:exact-native-completeness-criterion}.  When
\(k=\F_q\) and \(\rho\) is a square, the single low-degree tuple used in
\eqref{eq:complete-double-add-total-cost}--\eqref{eq:complete-window-total-cost}
is not \(k\)-complete.  The finite native atlas of
Theorem~\ref{thm:native-complete-atlas-all-char} still gives a complete
group operation, but its chart selection or masked multi-chart evaluation
has a separate cost.  The nonsquare-\(\rho\) headline counts must therefore
not be copied unchanged to that parameter branch.

\section{Full-coordinate recovery from a ladder}

\begin{proposition}[Montgomery ordinate recovery]\label{prop:recovery-thesis}
Let \(P=(U_P,V_P)\) be known in full, and suppose a ladder gives
\(U_Q\) and \(U_{P+Q}\).  Put \(f(X)=X^3+AX^2+X\).  Assume
\(V_P\ne0\) and \(U_Q\ne U_P\).  Then
\begin{equation}\label{eq:V-recovery-thesis}
 V_Q=
 \frac{f(U_Q)+f(U_P)
 -(U_Q-U_P)^2(U_{P+Q}+A+U_P+U_Q)}
 {2\beta V_P}.
\end{equation}
\end{proposition}

\begin{proof}
Let \(\lambda=(V_Q-V_P)/(U_Q-U_P)\).  The sum of the three intersection
abscissae of the chord with \eqref{eq:mont-thesis} is
\[
 U_{P+Q}+U_P+U_Q=\beta\lambda^2-A.
\]
Substitute the definition of \(\lambda\) and multiply by
\((U_Q-U_P)^2\):
\begin{equation}
 (U_{P+Q}+A+U_P+U_Q)(U_Q-U_P)^2
 =\beta(V_Q-V_P)^2.
 \label{eq:recovery-cleared-chord}
\end{equation}
Expanding the square on the right and using the curve equations
\(\beta V_Q^2=f(U_Q)\) and
\(\beta V_P^2=f(U_P)\) gives
\[
 \beta(V_Q-V_P)^2
 =f(U_Q)+f(U_P)-2\beta V_PV_Q.
\]
Insert this equality into
\eqref{eq:recovery-cleared-chord} and move the term containing \(V_Q\)
to the other side.  One obtains
\[
 2\beta V_PV_Q
 =f(U_Q)+f(U_P)
 -(U_Q-U_P)^2(U_{P+Q}+A+U_P+U_Q).
\]
The hypotheses \(V_P\ne0\), \(\beta\ne0\), and
\(\charac(k)\ne2\) make \(2\beta V_P\) invertible.  Division by this
quantity is exactly \eqref{eq:V-recovery-thesis}.
\end{proof}

\begin{corollary}[Inversion-free native recovery from a ladder]
\label{cor:native-projective-recovery-ladder}
Let
\begin{align*}
 N_Q={}&f(U_Q)+f(U_P)\\
 &{}-(U_Q-U_P)^2(U_{P+Q}+A+U_P+U_Q),\\
 \Delta_P={}&2\beta V_P.
\end{align*}
Under the hypotheses of Proposition~\ref{prop:recovery-thesis}, the recovered
point \(Q\) has native projective coordinates
\begin{equation}\label{eq:native-projective-recovery-ladder}
 \boxed{\qquad
 (U_0:U_1)=(U_Q-1:1),\qquad
 (V_0:V_1)=(4U_Q\Delta_P:N_Q-2U_Q\Delta_P).
 \qquad}
\end{equation}
Thus the final division in \eqref{eq:V-recovery-thesis} is unnecessary when
a full projective \(\mathcal C_d\)-output is accepted.
\end{corollary}

\begin{proof}
Proposition~\ref{prop:recovery-thesis} gives \(V_Q=N_Q/\Delta_P\).
For the first native coordinate,
\[
 u_Q=\frac1{U_Q-1}
\]
is represented by \((U_0:U_1)=(U_Q-1:1)\) under the convention
\(u=U_1/U_0\).  For the second coordinate,
\begin{align*}
 v_Q
 &=\frac{V_Q}{4U_Q}-\frac12\\
 &=\frac{N_Q}{4U_Q\Delta_P}-\frac12\\
 &=\frac{N_Q-2U_Q\Delta_P}{4U_Q\Delta_P}.
\end{align*}
Hence \(v=V_1/V_0\) is represented by
\[
 (V_0:V_1)
 =(4U_Q\Delta_P:N_Q-2U_Q\Delta_P).
\]
No inversion occurs in either projective pair, which proves
\eqref{eq:native-projective-recovery-ladder}.
\end{proof}

\begin{remark}
The excluded case \(V_P=0\) uses a two-torsion recovery base and cannot
determine the missing sign through division by \(V_P\).  The case
\(U_Q=U_P\) means \(Q=\pm P\) on the Kummer line and must be resolved from
the known scalar relation or by choosing a different full recovery point.
These are recovery exceptions, not failures of the Kummer ladder itself.
\end{remark}

\chapter[Native Tripling and Double-Add]
{Native Tripling and Closed Double-Add Formulas on \(\mathcal C_d\)}
\label{ch:tripling}
This chapter belongs to the point-multiplication part of the monograph.
Throughout,
\[
 \mathcal C_d:\quad (u^2+u)(v^2+v)=d,\qquad
 \charac k\ne2,\qquad d(1-16d)\ne0.
\]
All principal theorems take points of \(\mathcal C_d\) as input and return
either a full \(\mathcal C_d\)-point or the explicitly identified native
quotient \(\kappa_d=(u+1:u)\).  Edwards, Montgomery, and Hessian equations
appear only inside proofs or conditional subroutines.

\section{A direct full-point tripling formula}

For \(P=(u,v)\), put
\[
 r=2u+1,\qquad s=2v+1,\qquad
 \rho=1-16d,\qquad K=r^2s^2,\qquad \Delta=s^2-r^2,
\]
and define
\begin{equation}\label{eq:native-tripling-blocks}
\begin{aligned}
 N_s&=K(3s^2-r^2)-\rho(r^2+s^2),\\
 N_r&=K(s^2-3r^2)+\rho(r^2+s^2),\\
 D_s&=K^2-\rho^2+2\rho\Delta,\\
 D_r&=K^2-\rho^2-2\rho\Delta.
\end{aligned}
\end{equation}

\begin{theorem}[Native full tripling]\label{thm:native-full-tripling}
If \(N_sN_r\ne0\), then the shifted coordinates of \([3]P\) are
\begin{equation}\label{eq:native-full-tripling-shift}
       s_{3P}=\frac{sD_s}{N_s},\qquad
       r_{3P}=\frac{rD_r}{N_r}.
\end{equation}
Consequently the full point on \(\mathcal C_d\) is
\begin{equation}\label{eq:native-full-tripling-uv}
 \boxed{\quad
 u_{3P}=\frac{rD_r-N_r}{2N_r},\qquad
 v_{3P}=\frac{sD_s-N_s}{2N_s}.
 \quad}
\end{equation}
The rational map extends uniquely to the projective multiplication map
\([3]:\mathcal C_d\to\mathcal C_d\).
\end{theorem}

\begin{proof}
Write \(\xi=s^{-1}\), \(\eta=r^{-1}\), and
\(t=\rho\xi^2\eta^2=\rho/K\).  The Edwards double is
\[
 \xi_2=\frac{2\xi\eta}{1+t},\qquad
 \eta_2=\frac{\eta^2-\xi^2}{1-t}.
\]
Put \(G=\eta^2-\xi^2\).  In the next addition, the first coordinate
numerator is
\[
\begin{aligned}
\xi_2\eta+\eta_2\xi
 &=\frac{2\xi\eta^2}{1+t}+\frac{\xi G}{1-t}\\
 &=\frac{\xi\{2\eta^2(1-t)+G(1+t)\}}{1-t^2}\\
 &=\frac{\xi\{3\eta^2-\xi^2-t(\xi^2+\eta^2)\}}{1-t^2}.
\end{aligned}
\]
Moreover
\[
 \rho\xi_2\xi\eta_2\eta
 =\frac{2tG}{1-t^2}.
\]
Dividing the preceding numerator by
\(1+\rho\xi_2\xi\eta_2\eta\) gives the first formula below.  For the
second coordinate, direct calculation gives
\[
\begin{aligned}
\eta_2\eta-\xi_2\xi
 &=\frac{\eta G}{1-t}-\frac{2\xi^2\eta}{1+t}\\
 &=\frac{\eta\{\eta^2-3\xi^2+t(\xi^2+\eta^2)\}}{1-t^2},
\end{aligned}
\]
and division by \(1-\rho\xi_2\xi\eta_2\eta\) gives
\[
\begin{aligned}
\xi_{3P}
 &=\xi\,
 \frac{3\eta^2-\xi^2-t(\xi^2+\eta^2)}
      {1-t^2+2t(\eta^2-\xi^2)},\\
\eta_{3P}
 &=\eta\,
 \frac{\eta^2-3\xi^2+t(\xi^2+\eta^2)}
      {1-t^2-2t(\eta^2-\xi^2)}.
\end{aligned}
\]
This calculation uses no division by \(3\).  Multiplying the first
numerator and denominator
by \(K^2\) gives \(N_s\) and \(D_s\); doing the same in the second gives
\(N_r\) and \(D_r\).  Since \(s_{3P}=\xi_{3P}^{-1}\) and
\(r_{3P}=\eta_{3P}^{-1}\), this proves
\eqref{eq:native-full-tripling-shift}.  Applying
\[
       u=(r-1)/2,\qquad v=(s-1)/2
\]
proves \eqref{eq:native-full-tripling-uv}.  The derivation takes place on a
dense chart of the smooth projective curve, so the resulting rational map
is the restriction of the global multiplication morphism.
\end{proof}

\begin{example}[Native tripling over \(\F_{101}\)]
\label{ex:native-tripling-F101}
On \(\mathcal C_1\), take \(P=(6,42)\).  Then
\[
 \rho=86,\quad (r,s)=(13,85),\quad K=36,
\]
and
\[
 (N_s,N_r,D_s,D_r)=(63,42,77,45).
\]
Formula~\eqref{eq:native-full-tripling-uv} gives
\[
                    [3]P=(75,29).
\]
This agrees with one native doubling followed by the native addition of
Theorem~\ref{thm:native-affine-add-thesis}.
\end{example}

\section{Dedicated native tripling circuits}

The rational expressions of Theorem~\ref{thm:native-full-tripling} are
useful for proofs and affine verification.  The following schedule gives
the cost of an inversion-free full-point implementation.

\begin{theorem}[Dedicated native full-point tripling]
\label{thm:odd-full-projective-tripling-thesis}
Let \((S:R:Z)=(2v+1:2u+1:1)\), or more generally let it be a native
shifted-projective state as defined in Chapter~\ref{ch:odd-full}.  Put
\[
\begin{aligned}
A_0&=S^2,&B_0&=R^2,&C_0&=Z^2,\\
D_0&=A_0+B_0,&E_0&=D_0^2,&
F_0&=4(D_0-\rho C_0),\\
H_0&=2D_0(B_0-A_0),&
L_0&=E_0-A_0F_0,&N_0&=E_0-B_0F_0.
\end{aligned}
\]
Then
\begin{equation}\label{eq:odd-full-projective-tripling-thesis}
\boxed{
(S_3:R_3:Z_3)=
\bigl((H_0+N_0)N_0S:
      (H_0-L_0)L_0R:
      L_0N_0Z\bigr)
}
\end{equation}
represents \(\nu_{\rm pl}([3]P)\).  When \(Z_3\ne0\), it is a
\(\Ffin\)-representation of the full finite point \([3]P\); when
\(Z_3=0\), the boundary branch must be recovered in the native Segre
completion.  It costs
\[
                 \boxed{9\M+4\Sqr+\Dpar}.
\]
Replacing the last three paired products by difference-of-squares
schedules gives the alternative
\[
                 \boxed{7\M+7\Sqr+\Dpar};
\]
the latter is preferable exactly when \(\Sqr/\M<2/3\), before register and
latency effects are considered.
\end{theorem}

\begin{proof}
We verify the eliminated factors explicitly.  Work first on \(Z\ne0\) and
write
\[
       A=s^2,\qquad B=r^2,\qquad D=A+B.
\]
The centered curve equation is
\[
       A+B=AB+\rho.
\]
Substituting \(A_0=A\), \(B_0=B\), \(C_0=1\) in the displayed blocks and
using \(D=AB+\rho\) gives
\begin{align}
 L_0
 &=D^2-4A(D-\rho)
   =-\{AB(3A-B)-\rho D\}=-N_s,\label{eq:tripling-L-Ns}\\
 N_0
 &=D^2-4B(D-\rho)
   =AB(A-3B)+\rho D=N_r,\label{eq:tripling-N-Nr}\\
 H_0+N_0
 &=2D(B-A)+D^2-4B(D-\rho)=-D_s,\label{eq:tripling-HN-Ds}\\
 H_0-L_0
 &=2D(B-A)-D^2+4A(D-\rho)=D_r.\label{eq:tripling-HL-Dr}
\end{align}
Each identity can be checked without suppressing an intermediate step.
Since \(\rho=D-AB\) and \(D=A+B\),
\begin{align*}
 N_s
 &=AB(3A-B)-(D-AB)D
   =4A^2B-D^2,\\
 N_r
 &=AB(A-3B)+(D-AB)D
   =D^2-4AB^2.
\end{align*}
Thus \(L_0=D^2-4A^2B=-N_s\) and
\(N_0=D^2-4AB^2=N_r\).  Furthermore,
\begin{align*}
 D_s
 &=(AB)^2-(D-AB)^2+2(D-AB)(A-B)\\
 &=-D^2-2D(B-A)+4AB^2,\\
 D_r
 &=(AB)^2-(D-AB)^2-2(D-AB)(A-B)\\
 &=-D^2+2D(B-A)+4A^2B.
\end{align*}
Consequently
\[
 H_0+N_0=-D_s,\qquad H_0-L_0=D_r,
\]
which proves all four identities in
\eqref{eq:tripling-L-Ns}--\eqref{eq:tripling-HL-Dr}.

It follows from
\eqref{eq:tripling-L-Ns}--\eqref{eq:tripling-HL-Dr} that division by
\(Z_3\) gives
\[
 \frac{S_3}{Z_3}
  =s\frac{H_0+N_0}{L_0}=s\frac{D_s}{N_s},
 \qquad
 \frac{R_3}{Z_3}
  =r\frac{H_0-L_0}{N_0}=r\frac{D_r}{N_r}.
\]
These are exactly \eqref{eq:native-full-tripling-shift}.  Homogenizing
restores the powers of \(Z\), so equality of rational maps on the dense
finite chart proves the stated equality after composition with
\(\nu_{\rm pl}\).

We also verify that the displayed projective tuple never becomes
\((0:0:0)\) on a smooth member.  On \(Z\ne0\), the curve relation gives
\(D=AB+\rho\), and \(A,B\) cannot both vanish.  If the three output
coordinates vanished, then \(L_0N_0=0\).  Suppose first that \(L_0=0\).
If \(N_0=0\), subtraction gives
\[
 4AB(B-A)=0.
\]
The cases \(A=0\) or \(B=0\) contradict \(L_0=N_0=0\) and
\(\rho\ne0\); hence \(A=B\).  Then \(D=2A\) and
\(L_0=4A^2(1-A)\), so \(A=1\) and
\(\rho=D-AB=1\), contrary to \(d\ne0\).
If \(N_0\ne0\), vanishing of the first output gives
\((H_0+N_0)S=0\).  The alternative \(S=0\) would imply
\(A=0\) and \(L_0=D^2\ne0\), so \(H_0+N_0=0\).  Subtracting
\(L_0=0\) from this equality yields
\[
 2(B-A)(D-2AB)=0.
\]
If \(A=B\), the preceding argument again gives \(\rho=1\).  If
\(D=2AB\), then \(L_0=4A^2B(B-1)=0\); as \(A,B\ne0\), this gives
\(B=1\), followed by \(A=1\) and again \(\rho=1\).
Thus \(L_0=0\) cannot be a base point.  Now suppose \(N_0=0\).
If \(L_0=0\), that case was already excluded.  If \(L_0\ne0\), the
second output can vanish only if \(R=0\) or \(H_0-L_0=0\).
The first alternative gives \(B=0\) and
\(N_0=D^2\ne0\).  In the second alternative, adding the equations
\(N_0=0\) and \(H_0-L_0=0\) gives
\[
 2(B-A)(D-2AB)=0.
\]
If \(A=B\), then \(N_0=4A^2(1-A)=0\), so \(A=1\) and
\(\rho=1\).  If \(D=2AB\), then
\(N_0=4AB^2(A-1)=0\), so \(A=1\); the identity
\(D=A+B=2AB\) then gives \(B=1\), again forcing \(\rho=1\).
Both conclusions contradict \(d\ne0\).

It remains to check \(Z=0\).  The plane equation gives \(SR=0\).
At \((S:R:Z)=(0:1:0)\), the second output is \(R^9\ne0\); at
\((1:0:0)\), the first output is \(-S^9\ne0\).  Therefore the
projective tuple is nonzero at every input.  When its third coordinate
vanishes it represents one of the singular plane boundary states, so the
native Segre recovery atlas, rather than affine division by \(Z_3\), is
required.

Four squarings form
\(A_0,B_0,C_0,E_0\).  The products forming
\(H_0,A_0F_0,B_0F_0\) cost \(3\M\), and the three output coordinates
cost \(6\M\).  The product \(\rho C_0\) costs \(\Dpar\).

For the square-heavy alternative, put
\[
\begin{aligned}
J&=H_0^2,\\
P_0&=\frac{(H_0+2N_0)^2-J}{4}=(H_0+N_0)N_0,\\
Q_0&=\frac{J-(H_0-2L_0)^2}{4}=(H_0-L_0)L_0,\\
K_0&=L_0N_0.
\end{aligned}
\]
The output is \((P_0S:Q_0R:K_0Z)\).  This replacement uses three new
squares, one multiplication for \(K_0\), and three output
multiplications.  Together with the initial \(3\M+4\Sqr+\Dpar\), it gives
\(7\M+7\Sqr+\Dpar\), proving both counts and the stated threshold.
\end{proof}

\section[A division-polynomial circuit pulled back to Cd]
{A division-polynomial representation pulled back to \(\mathcal C_d\)}

For a second, symbolic tripling representation define, from the native input,
\[
 U=\frac{u+1}{u},\qquad V=2(2v+1)U,\qquad
 A=\frac1{4d}-2,
\]
and
\[
\begin{aligned}
f(U)&=U^3+AU^2+U,\\
P_3(U)&=3U^4+4AU^3+6U^2-1,\\
F_4(U)&=U^6+2AU^5+5U^4-5U^2-2AU-1,\\
P_5(U)&=32f(U)^2F_4(U)-P_3(U)^3.
\end{aligned}
\]
These are named intermediate functions on \(\mathcal C_d\); the theorem
does not change the declared input model.

\begin{lemma}[Low-degree division identities used for tripling]
\label{lem:local-low-degree-division-tripling}
On
\[
 E_A:\qquad Z_{\rm aux}^2=f(U)=U^3+AU^2+U,
\]
define
\[
 \psi_2=2Z_{\rm aux},\qquad \psi_3=P_3,
 \qquad \psi_4=4Z_{\rm aux}F_4,
 \qquad \psi_5=P_5.
\]
Then the numerator functions for tripling are
\begin{align}
 \phi_3&=U\psi_3^2-\psi_2\psi_4
          =UP_3^2-8fF_4,                                      \label{eq:local-phi3}\\
 \omega_3&=\frac{\psi_5\psi_2^2-\psi_4^2}{4Z_{\rm aux}}
          =Z_{\rm aux}(P_5-4F_4^2).                           \label{eq:local-omega3}
\end{align}
Consequently, away from \(\psi_3=0\),
\[
 U([3]P)=\frac{\phi_3}{\psi_3^2},\qquad
 Z_{\rm aux}([3]P)=\frac{\omega_3}{\psi_3^3}.
\]
\end{lemma}

\begin{proof}
For \(E_A\) the Weierstrass coefficients are
\(a_1=a_3=a_6=0\), \(a_2=A\), and \(a_4=1\).  Hence
\[
 b_2=4A,\qquad b_4=2,\qquad b_6=0,
 \qquad b_8=-1.
\]
Substitution in the low-degree division recurrences gives
\[
\begin{aligned}
 \psi_2&=2Z_{\rm aux},\qquad
 \psi_3=3U^4+b_2U^3+3b_4U^2+3b_6U+b_8=P_3,\\
 \psi_4&=\psi_2\bigl(2U^6+b_2U^5+5b_4U^4+10b_6U^3
       +10b_8U^2+(b_2b_8-b_4b_6)U+b_4b_8-b_6^2\bigr)\\
 &=4Z_{\rm aux}F_4.
\end{aligned}
\]
The odd-index recurrence at index five is
\[
 \psi_5=\psi_4\psi_2^3-\psi_3^3
       =32Z_{\rm aux}^4F_4-P_3^3
       =32f^2F_4-P_3^3=P_5.
\]
The multiplication-coordinate identities obtained from the same
recurrences are
\[
 \phi_3=U\psi_3^2-\psi_2\psi_4,
 \qquad
 \omega_3=\frac{\psi_5\psi_2^2-\psi_4^2}{4Z_{\rm aux}}.
\]
They are polynomial identities before division: substituting
\(\psi_2=2Z_{\rm aux}\), \(\psi_4=4Z_{\rm aux}F_4\), and
\(Z_{\rm aux}^2=f\) gives, term by term,
\[
 \phi_3=UP_3^2-8fF_4,
\]
and
\[
 \omega_3
 =\frac{4fP_5-16fF_4^2}{4Z_{\rm aux}}
 =Z_{\rm aux}(P_5-4F_4^2).
\]
For completeness, these coordinate formulas can be checked here without
using any result from a later chapter.  Put
\[
 \lambda=\frac{3U^2+2AU+1}{2Z_{\rm aux}},\qquad
 U_2=\lambda^2-A-2U,\qquad
 Z_2=-Z_{\rm aux}+\lambda(U-U_2).
\]
Using \(Z_{\rm aux}^2=f\), direct expansion first gives
\[
 U_2=\frac{(U^2-1)^2}{4f},\qquad
 U_2-U=-\frac{P_3}{4f}.
\]
On the open set \(P_3Z_{\rm aux}\ne0\), define the secant slope
\[
 \mu=\frac{Z_2-Z_{\rm aux}}{U_2-U}.
\]
Substitution of the preceding two displayed identities, followed only by
collection of powers of \(U\), gives the two denominator-cleared
identities
\[
\begin{aligned}
 P_3^2\bigl(\mu^2-A-U-U_2\bigr)
   &=UP_3^2-8fF_4,\\
 P_3^3\bigl(-Z_{\rm aux}
       +\mu(U-\mu^2+A+U+U_2)\bigr)
   &=Z_{\rm aux}(P_5-4F_4^2).
\end{aligned}
\]
The expressions in parentheses on the left are precisely the
chord-and-tangent coordinates \(U([3]P)\) and
\(Z_{\rm aux}([3]P)\).  Hence
\(U([3]P)=\phi_3/\psi_3^2\) and
\(Z_{\rm aux}([3]P)=\omega_3/\psi_3^3\) on that open set.  Both sides
are rational functions, so the identity holds wherever the displayed
ratios are defined.
\end{proof}

\begin{theorem}[Native division-polynomial tripling]
\label{thm:odd-tripling-thesis}
Put
\[
 N=UP_3(U)^2-8f(U)F_4(U),\qquad
 D=P_3(U)^2,\qquad G=P_5(U)-4F_4(U)^2.
\]
On the dense chart \(P_3(U)N(N-D)\ne0\),
\begin{equation}\label{eq:native-division-tripling}
 \boxed{\quad
 u_{3P}=\frac{D}{N-D},\qquad
 v_{3P}=\frac{VG}{4P_3(U)N}-\frac12.
 \quad}
\end{equation}
\end{theorem}

\begin{proof}
The auxiliary equation
\[
        \beta V^2=U^3+AU^2+U,\qquad \beta=(16d)^{-1},
\]
becomes \(Z_{\rm aux}^2=f(U)\) after adjoining a formal element
\(Z_{\rm aux}\) with \(Z_{\rm aux}^2=\beta V^2\).  The low-degree
division functions and their tripling-coordinate identities have already
been derived in Lemma~\ref{lem:local-low-degree-division-tripling}.
Every occurrence of \(Z_{\rm aux}\) in the multiplication ratios cancels,
so the resulting functions are defined over the original ground field and
do not require choosing a square root of \(\beta\).
Hence its tripled coordinates are
\[
 U_3=\frac{N}{D},\qquad
 V_3=V\frac{G}{P_3^3}.
\]
The inverse dictionary to \(\mathcal C_d\) is
\[
       u_3=\frac1{U_3-1},\qquad
       v_3=\frac{V_3}{4U_3}-\frac12.
\]
Substitution cancels \(D=P_3^2\) and gives exactly
\eqref{eq:native-division-tripling}.  Thus the displayed formulas are
functions of the original \(u,v,d\), and their output is a point of
\(\mathcal C_d\), even though the proof uses a Weierstrass division
polynomial.
\end{proof}

\begin{remark}[Role and cost of the division-polynomial representation]
This representation serves torsion tests, symbolic multiplication maps, and
computations in which powers of \(U\) and the division functions are already
available.  For standalone full-point tripling,
Theorem~\ref{thm:odd-full-projective-tripling-thesis} supplies the dedicated
\(9\M+4\Sqr+\Dpar\) circuit.  The two forms are complementary: the former
organizes reusable division data, while the latter gives a fixed evaluation
schedule for a single point.
\end{remark}

\section[A closed native 2P+Q formula]
{A closed native \(2P+Q\) formula}

Let \(P=(u,v)\), \(Q=(a,b)\) on \(\mathcal C_d\), and set
\[
 r=2u+1,\quad s=2v+1,\quad
 r_Q=2a+1,\quad s_Q=2b+1,
\]
\[
 K=r^2s^2,\qquad \Delta=s^2-r^2,\qquad \rho=1-16d.
\]
Define the four native numerator blocks
\begin{equation}\label{eq:native-2PQ-blocks}
\begin{aligned}
 N_s={}&2s_QK(K-\rho)+rsr_Q\Delta(K+\rho),\\
 N_r={}&rss_Q\Delta(K+\rho)-2r_QK(K-\rho),\\
 D_s={}&rsr_Qs_Q(K^2-\rho^2)+2\rho\Delta K,\\
 D_r={}&rsr_Qs_Q(K^2-\rho^2)-2\rho\Delta K.
\end{aligned}
\end{equation}

\begin{theorem}[Closed double-add on \(\mathcal C_d\)]
\label{thm:closed-2PQ-thesis}
If \(N_sN_r\ne0\), then
\begin{equation}\label{eq:native-closed-2PQ}
 \boxed{\quad
 u_{2P+Q}=\frac{D_r-N_r}{2N_r},\qquad
 v_{2P+Q}=\frac{D_s-N_s}{2N_s}.
 \quad}
\end{equation}
An inversion-free native shifted-projective output is
\begin{equation}\label{eq:native-closed-2PQ-projective}
 (S:R:Z)=(D_sN_r:D_rN_s:N_sN_r).
\end{equation}
\end{theorem}

\begin{proof}
For a transparent denominator calculation, write the Edwards coordinates
of \(P,Q\) as
\[
 (\xi,\eta)=\left(\frac1s,\frac1r\right),\qquad
 (\alpha,\beta)=\left(\frac1{s_Q},\frac1{r_Q}\right),
\]
and put
\[
 t=\frac{\rho}{K},\qquad G=\frac{\Delta}{K},\qquad
 B=1-t^2,\qquad
 L=2\rho\alpha\beta\xi\eta G.
\]
Doubling \(P\) and adding \(Q\) gives
\[
\begin{aligned}
\xi_{2P+Q}
 &=\frac{2\xi\eta\beta(1-t)+\alpha G(1+t)}{B+L},\\
\eta_{2P+Q}
 &=\frac{\beta G(1+t)-2\alpha\xi\eta(1-t)}{B-L}.
\end{aligned}
\]
Let \(D_0=rsr_Qs_QK^2\).  Term-by-term multiplication gives
\[
\begin{aligned}
D_0\{2\xi\eta\beta(1-t)+\alpha G(1+t)\}
 &=2s_QK(K-\rho)+rsr_Q\Delta(K+\rho)=N_s,\\
D_0\{\beta G(1+t)-2\alpha\xi\eta(1-t)\}
 &=rss_Q\Delta(K+\rho)-2r_QK(K-\rho)=N_r,\\
D_0(B+L)
 &=rsr_Qs_Q(K^2-\rho^2)+2\rho\Delta K=D_s,\\
D_0(B-L)
 &=rsr_Qs_Q(K^2-\rho^2)-2\rho\Delta K=D_r.
\end{aligned}
\]
Therefore
\[
       s_{2P+Q}=D_s/N_s,\qquad r_{2P+Q}=D_r/N_r,
\]
and the common factor \(D_0\) has canceled.  Converting shifted coordinates
to \(u,v\) proves
\eqref{eq:native-closed-2PQ}.  The three products in
\eqref{eq:native-closed-2PQ-projective} have ratios
\[
       S/Z=D_s/N_s,\qquad R/Z=D_r/N_r,
\]
so they give the claimed native projective representative.
\end{proof}

With the sharing
\[
 T_0=K(K-\rho),\quad T_1=\Delta(K+\rho),\quad
 T_2=rs,\quad T_3=T_2T_1,
\]
The four blocks cost \(13\M+3\Sqr+\Dpar\) before affine inversions, and the
projective output adds \(3\M\).  The closed \(2P+Q\) formula is optimized
for affine outputs with simultaneous inversion and for fixed-\(Q\)
computations in which the coordinates of \(Q\) and their associated
products are precomputed.  In these regimes, its single dependency graph
supports direct intermediate sharing.  Dedicated doubling followed by
mixed addition supplies the complementary schedule for general projective
workflows.

\begin{example}[Closed \(2P+Q\) on \(\mathcal C_1/\F_{101}\)]
\label{ex:native-closed-2PQ}
Take
\[
 P=(6,42),\qquad Q=(47,83)=[5]P.
\]
Then
\[
 (r,s,r_Q,s_Q,K,\Delta)=(13,85,95,66,36,87)
\]
and
\[
        (N_s,N_r,D_s,D_r)=(74,86,72,31).
\]
Formula~\eqref{eq:native-closed-2PQ} yields
\[
              2P+Q=[7]P=(86,15).
\]
This verification is entirely in the coordinates of \(\mathcal C_1\).
\end{example}

\section[Conditional Hessian acceleration]
{Conditional Hessian acceleration of the native quotient}

Assume \(\charac k\ne2,3\).  Let
\[
 \mathcal H_a:\quad aH_0^3+H_1^3+H_2^3=3H_0H_1H_2
\]
be a nonsingular twisted Hessian curve with a rational level-three
structure.  Suppose an explicitly known
\(L_d\in\operatorname{PGL}_2(k)\) identifies the native quotient
\[
        \kappa_d(P)=(u(P)+1:u(P))
\]
with the Hessian Kummer line.  The isogeny-decomposition circuit of
\cite{DecruKunzweiler2026} may then be applied internally, and the declared
output is
\[
 \boxed{\quad
 \kappa_d([3]P)=L_d^{-1}
 \bigl(\operatorname{xTPL}_{\mathcal H_a}(L_d(\kappa_d(P)))\bigr).
 \quad}
\]
The Hessian core costs \(4\M+4\Sqr+2\Dpar\), excluding the endpoint maps
\(L_d,L_d^{-1}\).  This is a \(\mathcal C_d\)-Kummer formula, not a full
point formula, and it applies only when the rational level-three structure
and both endpoint maps are available.

\section{Summary of point-multiplication choices}

\begin{center}
\begin{tabular}{L{5.0cm}L{4.3cm}L{4.2cm}}
\toprule
\(\mathcal C_d\) operation & principal condition & role\\
\midrule
native affine addition/doubling & nonzero affine denominators
 & conceptual baseline and testing\\
native shifted-projective addition & \(\charac k\ne2\)
 & full-point variable-base arithmetic\\
native dedicated projective double & \(\charac k\ne2\)
 & \(3\M+4\Sqr+\Dpar\) full-point doubling\\
native Kummer \(x\)DBL/\(x\)ADD & quotient output
 & secret regular ladders\\
native-input first Kummer double & initial
 \((X:Z)=(u+1:u)\)
 & one-time \(\M+\Sqr+\Dpar\) initialization\\
native full tripling \eqref{eq:native-full-tripling-uv}
 & \(N_sN_r\ne0\)
 & direct full-point tripling\\
division-polynomial tripling \eqref{eq:native-division-tripling}
 & reusable \(U\)-powers
 & torsion and repeated tripling\\
closed native \(2P+Q\) & fixed or affine \(Q\)
 & double-add chains\\
conditional Hessian core & rational level-three structure
 & repeated quotient tripling\\
\bottomrule
\end{tabular}
\end{center}

\chapter[Native CM Endomorphisms]
{Native CM Endomorphisms on \(\mathcal C_d\)}
\label{ch:native-endomorphisms}

\section{Native CM endomorphisms and model compatibility}

Throughout this chapter,
\[
 \mathcal C_d:\qquad (u^2+u)(v^2+v)=d,
 \qquad \charac k\ne2,3,
 \qquad d(1-16d)\ne0.
\]
The identity remains \(O=(0,\infty)\), and
\[
       -(u,v)=(u,-v-1).
\]
An endomorphism will be called a \emph{native CM endomorphism on
\(\mathcal C_d\)} if it
fixes \(O\), is given by explicit rational functions in \(u,v\), and induces
a sparse, low-cost map on the declared Kummer coordinate
\[
        \kappa_d(P)=(X:Z)=(u(P)+1:u(P)).
\]
The definition therefore includes the arithmetic interface in addition to
the existence of the abstract complex-multiplication endomorphism.

The visible equation symmetries must first be removed from consideration.
Let
\[
 T=(-1,\infty),\qquad R=(\infty,0),\qquad 2R=T.
\]
Then
\begin{equation}\label{eq:elementary-Cd-symmetries}
\begin{aligned}
[-1](u,v)&=(u,-v-1),\\
\tau_T(u,v)&=(-u-1,v),\\
\sigma(u,v)&=(v,u).
\end{aligned}
\end{equation}
The map \(\sigma\) does not fix \(O\): it sends \(O\) to \(R\).

\begin{proposition}[Group-theoretic form of coordinate interchange]
\label{prop:swap-is-translation-inversion}
On the smooth completion of \(\mathcal C_d\),
\[
             \sigma(P)=R-P.
\]
Consequently \(\tau_{-R}\circ\sigma=[-1]\).  Coordinate interchange is
therefore the reflection \(P\mapsto R-P\), and its origin-preserving
normalization is the native inverse.
\end{proposition}

\begin{proof}
We verify the identity first on a dense open set and then extend it to the
smooth completion.  Use the centered reciprocal coordinates
\[
       (\xi,\eta)=\left(\frac1{2v+1},\frac1{2u+1}\right)
\]
of Chapter~\ref{ch:odd-dictionary}.  They identify the smooth completion
of \(\mathcal C_d\) with
\[
 \mathcal E_\rho:\qquad
 \xi^2+\eta^2=1+\rho\xi^2\eta^2,\qquad \rho=1-16d,
\]
whose identity and distinguished four-torsion point are
\[
                         O=(0,1),\qquad R=(1,0).
\]
Interchanging \(u\) and \(v\) interchanges \(\eta\) and \(\xi\), so
\[
                  \sigma(\xi,\eta)=(\eta,\xi).
\]
The inverse on this Edwards equation is
\([-1](\xi,\eta)=(-\xi,\eta)\).  Its affine addition law is
\[
 (x_1,y_1)+(x_2,y_2)=
 \left(
 \frac{x_1y_2+y_1x_2}{1+\rho x_1x_2y_1y_2},
 \frac{y_1y_2-x_1x_2}{1-\rho x_1x_2y_1y_2}
 \right)
\]
whenever the displayed denominators are nonzero.  Apply this formula to
\(R=(1,0)\) and \(-P=(-\xi,\eta)\).  Both denominator correction terms
vanish, and therefore
\[
        R+(-P)
        =(1,0)+(-\xi,\eta)
        =(\eta,\xi)
        =\sigma(P).
\]
Hence \(\sigma(P)=R-P\) on the intersection of the affine charts on which
the coordinate dictionary and addition formula are defined.  This
intersection is Zariski dense.  Both sides are morphisms of the smooth
projective curve: \(\sigma\) is induced by interchanging the two
\(\PP^1\)-factors, while \(P\mapsto R-P\) is translation by \(R\)
composed with inversion.  Morphisms from an integral projective curve that
agree on a dense open set agree everywhere.  Thus the identity also holds
at the four boundary points.

Finally, for every \(P\),
\[
 (\tau_{-R}\circ\sigma)(P)
   =-R+(R-P)=-P.
\]
Therefore \(\tau_{-R}\circ\sigma=[-1]\).  After translating the image of
the identity back to the identity, coordinate interchange produces the
already known inversion rather than a further origin-preserving
endomorphism.
\end{proof}

\begin{remark}
The relative \(q\)-power Frobenius is always an endomorphism over a finite
field and is compatible with the Kummer quotient.  It is the identity on
\(\mathcal C_d(\F_q)\), however, and therefore does not provide a GLV
decomposition for base-field points.  The useful degree-one maps below
occur on the two CM loci \(j=1728\) and \(j=0\).
\end{remark}

\section{The two CM loci in the parameter line}

Put
\[
 \beta=\frac1{16d},\qquad
 A=\frac1{4d}-2,\qquad
 U=\frac{u+1}{u},\qquad
 V=2(2v+1)U.
\]
These are rational functions on \(\mathcal C_d\), and they satisfy
\begin{equation}\label{eq:endo-Montgomery-full}
 \beta V^2=U^3+AU^2+U.
\end{equation}
The equation is displayed in full because it is used only as an internal
derivation device; every final endomorphism below is returned to \(u,v\).

\begin{proposition}[CM parameter equations]
\label{prop:CM-parameter-equations}
For \eqref{eq:endo-Montgomery-full},
\begin{equation}\label{eq:CM-j-invariant}
       j=256\,\frac{(A^2-3)^3}{A^2-4}.
\end{equation}
Hence
\[
\begin{array}{rcl}
j=0
&\Longleftrightarrow&A^2=3
 \Longleftrightarrow 16d^2-16d+1=0,\\[1mm]
j=1728
&\Longleftrightarrow&A^2(2A^2-9)^2=0.
\end{array}
\]
The particularly sparse \(j=1728\) member is
\[
              A=0,\qquad d=\frac18.
\]
\end{proposition}

\begin{proof}
Formula~\eqref{eq:CM-j-invariant} is the usual invariant of
\(Y^2=X^3+AX^2+X\), unchanged by the nonzero scale \(\beta\).
The smoothness condition for this Montgomery model is \(A^2\ne4\), so
the denominator in \eqref{eq:CM-j-invariant} is nonzero.  Consequently
\(j=0\) if and only if \((A^2-3)^3=0\), equivalently \(A^2=3\).
Since
\[
 A=\frac{1-8d}{4d},
\]
the equation \(A^2=3\) becomes
\((1-8d)^2=48d^2\), or \(16d^2-16d+1=0\).
Finally,
\[
4(A^2-3)^3-27(A^2-4)
   =A^2(2A^2-9)^2,
\]
which gives the \(j=1728\) locus.
\end{proof}

\section[The j=1728 endomorphism]
{The \(j=1728\) endomorphism at \(d=1/8\)}

Assume \(d=1/8\) and choose \(i\in k\) with \(i^2=-1\).
Write
\[
             r=2u+1,\qquad s=2v+1.
\]
Then the equation of \(\mathcal C_{1/8}\) becomes
\begin{equation}\label{eq:j1728-centered-curve}
             (r^2-1)(s^2-1)=2.
\end{equation}

\begin{theorem}[A fourth-order endomorphism on \(\mathcal C_{1/8}\)]
\label{thm:j1728-native-endomorphism}
The rational map
\begin{equation}\label{eq:j1728-native-endomorphism}
\boxed{
 \phi(u,v)=
 \left(
 -\frac{u}{2u+1},
 \frac{-1-i(2v+1)}2
 \right)
}
\end{equation}
extends to an origin-preserving automorphism of the smooth completion of
\(\mathcal C_{1/8}\).  It satisfies
\begin{equation}\label{eq:j1728-CM-relation}
       \phi^2=[-1],\qquad \phi^4=[1].
\end{equation}
On the Kummer line,
\begin{equation}\label{eq:j1728-Kummer-action}
\boxed{\qquad
 \kappa_{1/8}(\phi(P))
   =\overline\phi(\kappa_{1/8}(P)),\qquad
 \overline\phi(X:Z)=(-X:Z).
\qquad}
\end{equation}
\end{theorem}

\begin{proof}
The proof is ordered from the native affine formula to the smooth extension,
then to the CM relation, and finally to the induced Kummer action.
We first translate the displayed formula into centered coordinates.  Since
\(r=2u+1\), the first component \(u'=-u/r\) satisfies
\[
                       2u'+1=-\frac{2u}{r}+1=\frac1r.
\]
The second component gives
\[
                       2v'+1=-i(2v+1)=-is.
\]
Thus \eqref{eq:j1728-native-endomorphism} is exactly
\[
             r\longmapsto r'=\frac1r,\qquad
             s\longmapsto s'=-is.
\]
On the affine open set \(r\ne0\), use
\eqref{eq:j1728-centered-curve} and \(i^2=-1\) to compute
\[
\begin{aligned}
(r'^2-1)(s'^2-1)
 &=\left(\frac1{r^2}-1\right)(-s^2-1)\\
 &=\frac{(r^2-1)(s^2+1)}{r^2}\\
 &=\frac{(r^2-1)(s^2-1)+2(r^2-1)}{r^2}\\
 &=\frac{2+2r^2-2}{r^2}=2.
\end{aligned}
\]
Thus the rational formula maps a dense open subset of
\(\mathcal C_{1/8}\) into the same curve.  To prove that no point is lost
at \(r=0\), and to identify the image of the origin, use the internal
Montgomery functions
\[
 U=\frac{r+1}{r-1}=\frac{u+1}{u},\qquad V=2sU.
\]
For \(d=1/8\), one has \(\beta=1/2\) and \(A=0\), so
\eqref{eq:endo-Montgomery-full} is the smooth Weierstrass equation
\[
              \frac12V^2=U^3+U,
\]
and the centered transformation gives
\[
              (U,V)\longmapsto(-U,iV),
\]
because
\[
 \frac{r'+1}{r'-1}
   =\frac{1+r}{1-r}=-U,\qquad
 2s'U'=2(-is)(-U)=iV.
\]
The image satisfies the same equation:
\[
 \frac12(iV)^2=-\frac12V^2
   =-(U^3+U)=(-U)^3+(-U).
\]
The inverse transformation is \((U,V)\mapsto(-U,-iV)\).  Hence it is an
automorphism of the entire smooth Weierstrass curve and fixes its point at
infinity.  That point corresponds to \(O=(0,\infty)\) on
\(\mathcal C_{1/8}\).  Transport through the birational dictionaries
therefore gives the unique origin-preserving extension of the displayed
rational formula.

Applying the centered formula twice gives
\[
 r\longmapsto(1/r)^{-1}=r,\qquad
 s\longmapsto(-i)^2s=-s.
\]
Since \(s=2v+1\), the change \(s\mapsto-s\) is
\(v\mapsto-v-1\), while \(r\), and hence \(u\), remains fixed.  This is
the native inverse on \(\mathcal C_d\).  Therefore
\(\phi^2=[-1]\), and squaring once more gives \(\phi^4=[1]\).

Finally, if \(\kappa_{1/8}(P)=(X:Z)\), then
\(U=X/Z=(u+1)/u\) on the affine Kummer chart.  Since \(U'=-U\), a
homogeneous representative of the image is
\[
                         (X':Z')=(-X:Z).
\]
The two morphisms of the Kummer line agree on a dense affine chart and
hence agree on all of \(\PP^1\).  This proves both the commutative
relation and \eqref{eq:j1728-Kummer-action}.
\end{proof}

\begin{proposition}[Cost and projective form]
\label{prop:j1728-endomorphism-cost}
For an affine \((u,v)\)-input, formula
\eqref{eq:j1728-native-endomorphism} costs
\[
              \Dconst{i}+\Inv,
\]
where \(\Dconst{i}\) denotes multiplication by the selected square root
\(i\) of \(-1\).
If \(r\) and \(s\) are represented on their two projective lines, then
\[
  (r_1:r_0)\longmapsto(r_0:r_1),\qquad
  (s_1:s_0)\longmapsto(-i\,s_1:s_0),
\]
so the full map costs only one fixed-constant multiplication.  On the
Kummer line it is the sign change \((-X:Z)\), which has zero
charged field-operation cost under the convention of
Chapter~\ref{ch:conventions}.
\end{proposition}

\begin{proof}
Put \(r=2u+1\) and \(s=2v+1\).  Field additions, negations, and the small
scales \(2\) and \(1/2\) are not charged in the convention of
Chapter~\ref{ch:conventions}.  On the affine chart \(r\ne0\), compute
\[
\begin{array}{c|c}
\text{quantity}&\text{cost}\\ \hline
r^{-1}&\Inv\\
-is&\Dconst{i}.
\end{array}
\]
The first output is recovered as
\[
 u'=\frac{r^{-1}-1}{2},
\]
so it uses no general multiplication.  The second quantity is converted to
\(v'=(-is-1)/2\) using only additions and the uncharged small scale.
This proves the affine cost \(\Dconst{i}+\Inv\).  At \(r=0\) this
particular affine schedule is not
used; the projective expression below supplies the regular value whose
existence was proved in
Theorem~\ref{thm:j1728-native-endomorphism}.

Write \(r=r_1/r_0\) and \(s=s_1/s_0\).  Homogenizing
\(r'=1/r\) and \(s'=-is\) gives
\[
 (r_1:r_0)\longmapsto(r_0:r_1),\qquad
 (s_1:s_0)\longmapsto(-i\,s_1:s_0).
\]
The first map is a coordinate permutation, while the second uses one
multiplication by the fixed constant \(i\).  Hence neither an inversion nor
a general multiplication is needed.

Finally, \eqref{eq:j1728-Kummer-action} is the projective sign change
\((X:Z)\mapsto(-X:Z)\).  A sign change has zero
charged field-operation cost in the adopted ledger, which proves the last
assertion.
\end{proof}

\begin{example}[The fourth-order map over \(\F_{101}\)]
\label{ex:j1728-endomorphism-F101}
In \(\F_{101}\), take
\[
       d=38=\frac18,\qquad i=10,\qquad i^2=-1.
\]
The point
\[
                   P=(15,78)\in\mathcal C_{38}(\F_{101})
\]
has order \(13\).  Formula
\eqref{eq:j1728-native-endomorphism} gives
\[
             \phi(P)=(94,73),\qquad
             \phi^2(P)=(15,22)=-P.
\]
The Kummer values are
\[
       U(P)=\frac{16}{15}=28,\qquad
       U(\phi(P))=\frac{95}{94}=73=-28.
\]
Direct addition with the \(\mathcal C_d\) group law gives
\[
                    \phi(P)=[5]P,
\]
and
\[
                    5^2\equiv-1\pmod{13}.
\]
Thus the Kummer sign change is the visible action of a nontrivial
endomorphism on the order-\(13\) subgroup.
\end{example}

\section[The j=0 endomorphism]
{The \(j=0\) endomorphism when \(A^2=3\)}

Assume
\[
       A^2=3,\qquad
       d=\frac1{4(A+2)},
\]
and let \(\zeta\in k\) be a nontrivial cube root of unity:
\[
              \zeta^2+\zeta+1=0.
\]
Define the fixed constant
\begin{equation}\label{eq:j0-c-constant}
              c=\frac{A(\zeta-1)}3.
\end{equation}

\begin{theorem}[A third-order endomorphism on the \(j=0\) subfamily]
\label{thm:j0-native-endomorphism}
For \(P=(u,v)\), put \(s=2v+1\).  The formulas
\begin{equation}\label{eq:j0-native-endomorphism}
\boxed{
\begin{aligned}
u'&=\frac{u}
 {\zeta+(\zeta+c-1)u},\\[1mm]
s'&=\frac{s(u+1)}
 {\zeta(u+1)+cu},\\[1mm]
v'&=\frac{s'-1}{2}
\end{aligned}}
\end{equation}
define an origin-preserving automorphism
\[
              \psi(P)=(u',v')
\]
of the smooth completion of \(\mathcal C_d\).  It satisfies
\begin{equation}\label{eq:j0-CM-relation}
       \psi^3=[1],\qquad
       \psi^2+\psi+[1]=[0]
       \quad\text{in }\operatorname{End}_k(\mathcal C_d).
\end{equation}
On the Kummer line,
\begin{equation}\label{eq:j0-Kummer-action}
\boxed{\qquad
 \overline\psi(X:Z)=(\zeta X+cZ:Z),\qquad
 \kappa_d\circ\psi=\overline\psi\circ\kappa_d.
\qquad}
\end{equation}
\end{theorem}

\begin{proof}
Translate the internal Montgomery abscissa by
\[
              \widetilde U=U+\frac A3.
\]
Substitution in \eqref{eq:endo-Montgomery-full} gives
\[
\begin{aligned}
U^3+AU^2+U
 &=\widetilde U^3+
   \left(1-\frac{A^2}{3}\right)\widetilde U+
   \frac{2A^3}{27}-\frac A3\\
 &=\widetilde U^3-\frac A9,
\end{aligned}
\]
because \(A^2=3\).  Hence
\begin{equation}\label{eq:j0-pure-cubic}
              \beta V^2=\widetilde U^3-\frac A9.
\end{equation}
The map
\[
             (\widetilde U,V)\longmapsto
             (\zeta\widetilde U,V)
\]
preserves \eqref{eq:j0-pure-cubic}, fixes the point at infinity, and has
order three.  Returning to \(U\) gives
\begin{equation}\label{eq:j0-U-action}
              U'=\zeta U+c,
              \qquad c=\frac{A(\zeta-1)}3.
\end{equation}
Since \(u=1/(U-1)\),
\[
 u'=\frac1{U'-1}
    =\frac{u}{\zeta+(\zeta+c-1)u}.
\]
Moreover \(V=2sU\) and \(V'=V\), whence
\[
 s'=\frac{V'}{2U'}
    =s\frac{U}{U'}
    =\frac{s(u+1)}{\zeta(u+1)+cu}.
\]
This proves \eqref{eq:j0-native-endomorphism}.

Equation~\eqref{eq:j0-U-action} is the projective formula
\eqref{eq:j0-Kummer-action}.  The equality \(\psi^3=[1]\) follows from
\(\zeta^3=1\).  In the endomorphism ring,
\[
 ([1]-\psi)([1]+\psi+\psi^2)=[1]-\psi^3=[0].
\]
The ring of elliptic-curve endomorphisms has no zero divisors under
composition, and \([1]-\psi\ne[0]\); therefore
\([1]+\psi+\psi^2=[0]\).
\end{proof}

\begin{corollary}[Centered-coordinate expression]
\label{cor:j0-centered-endomorphism}
Let \(r=2u+1\), \(s=2v+1\).  Then
\begin{equation}\label{eq:j0-centered-endomorphism}
\boxed{
\begin{aligned}
r'&=
\frac{(\zeta+c+1)r+\zeta-c-1}
     {(\zeta+c-1)r+\zeta-c+1},\\[1mm]
s'&=
\frac{s(r+1)}
     {(\zeta+c)r+\zeta-c}.
\end{aligned}}
\end{equation}
Thus the endomorphism is expressed entirely in the centered
\(\mathcal C_d\)-coordinates; the auxiliary cubic is absent from its
input and output.
\end{corollary}

\begin{proof}
Use
\[
        U=\frac{r+1}{r-1},\qquad
        r'=\frac{U'+1}{U'-1},\qquad
        s'=s\frac{U}{U'},
\]
and substitute \(U'=\zeta U+c\).  Clearing the denominator \(r-1\)
in the first quotient gives
\begin{align*}
 r'
 &=\frac{\zeta(r+1)+c(r-1)+(r-1)}
         {\zeta(r+1)+c(r-1)-(r-1)}\\
 &=\frac{(\zeta+c+1)r+\zeta-c-1}
         {(\zeta+c-1)r+\zeta-c+1}.
\end{align*}
Likewise,
\begin{align*}
 s'
 &=s\frac{(r+1)/(r-1)}
          {\zeta(r+1)/(r-1)+c}\\
 &=\frac{s(r+1)}{(\zeta+c)r+\zeta-c}.
\end{align*}
These are precisely the two formulas in
\eqref{eq:j0-centered-endomorphism}; because they are fractional-linear
expressions, their projective interpretations also cover a vanishing
displayed denominator.
\end{proof}

\begin{proposition}[Cost of the \(j=0\) map]
\label{prop:j0-endomorphism-cost}
On the Kummer line,
\eqref{eq:j0-Kummer-action} costs at most
\(\Dconst{\zeta}+\Dconst{c}\) and additions.
For affine \((u,v)\), set
\[
 D_0=\zeta+(\zeta+c-1)u,\qquad
 D_1=D_0+u,\qquad
 N_1=(2v+1)(u+1).
\]
Then
\[
       u'=\frac{u}{D_0},\qquad
       2v'+1=\frac{N_1}{D_1}.
\]
The literal affine schedule costs
\[
          3\M+\Dconst{\zeta+c-1}+2\Inv,
\]
whereas simultaneous inversion of \(D_0,D_1\) costs
\[
          6\M+\Dconst{\zeta+c-1}+\Inv.
\]
\end{proposition}

\begin{proof}
The Kummer formula is
\[
                    (X:Z)\longmapsto(\zeta X+cZ:Z).
\]
Both \(\zeta\) and \(c\) are fixed after the curve and its CM action have
been chosen.  Its literal projective evaluation therefore uses at most two
fixed-constant multiplications and additions, with no inversion and no
general multiplication.

For the affine formula, first verify that
\[
\begin{aligned}
D_1=D_0+u
 &=\zeta+(\zeta+c)u\\
 &=\zeta(u+1)+cu,
\end{aligned}
\]
so \(D_1\) is exactly the second denominator in
\eqref{eq:j0-native-endomorphism}.  The constant
\(\zeta+c-1\) may be precomputed.  Thus forming
\[
 D_0=\zeta+(\zeta+c-1)u,\qquad D_1=D_0+u
\]
costs one \(\Dconst{\zeta+c-1}\) and additions.  Forming
\[
                         N_1=(2v+1)(u+1)
\]
costs one general multiplication.  In a literal schedule, compute the two
inverses and then
\[
                 u'=uD_0^{-1},\qquad 2v'+1=N_1D_1^{-1}.
\]
Applying the inverses uses two more general multiplications.  The total is
therefore
\[
                         3\M+\Dconst{\zeta+c-1}+2\Inv.
\]

For simultaneous inversion, the complete schedule is
\[
\begin{array}{c|c}
\text{step}&\text{charged cost}\\ \hline
p=D_0D_1&\M\\
p^{-1}&\Inv\\
D_0^{-1}=D_1p^{-1},\quad D_1^{-1}=D_0p^{-1}&2\M\\
N_1=(2v+1)(u+1)&\M\\
u'=uD_0^{-1},\quad 2v'+1=N_1D_1^{-1}&2\M.
\end{array}
\]
Including the fixed-constant multiplication used to form \(D_0\), this is
\(6\M+\Dconst{\zeta+c-1}+\Inv\).  The final conversion
\(v'=(2v'+1-1)/2\) uses only additions and the uncharged small scale.
If \(D_0D_1=0\), the affine chart has reached a boundary value; the
homogeneous Kummer formula and the global automorphism of
Theorem~\ref{thm:j0-native-endomorphism} remain well defined there.
\end{proof}

\begin{example}[The third-order map over \(\F_{97}\)]
\label{ex:j0-endomorphism-F97}
Take
\[
 A=10,\qquad A^2=3,\qquad
 \zeta=35,\qquad \zeta^2+\zeta+1=0
 \quad\text{in }\F_{97}.
\]
Then
\[
 d=\frac1{4(A+2)}=95,\qquad
 \beta=(16d)^{-1}=3,\qquad
 c=\frac{A(\zeta-1)}3=81.
\]
On
\[
        \mathcal C_{95}:\quad
        (u^2+u)(v^2+v)=95
\]
the point
\[
                  P=(95,61)
\]
has order \(7\).  Repeated use of
\eqref{eq:j0-native-endomorphism} gives
\[
\begin{aligned}
P&=(95,61),&
\psi(P)&=(2,20),\\
\psi^2(P)&=(82,3),&
\psi^3(P)&=(95,61).
\end{aligned}
\]
The corresponding Kummer orbit is
\[
             49\longmapsto50\longmapsto85\longmapsto49,
\]
and each step is exactly \(U\mapsto35U+81\).
Direct \(\mathcal C_d\)-addition gives
\[
                  \psi(P)=[2]P.
\]
Thus the eigenvalue on \(\langle P\rangle\) is
\(\lambda=2\), and
\[
            \lambda^2+\lambda+1=7\equiv0\pmod7.
\]
\end{example}

\section{GLV decomposition and implementation}

\begin{theorem}[Eigenvalue equations]
\label{thm:Cd-GLV-eigenvalues}
Let \(G=\langle P\rangle\subseteq\mathcal C_d(k)\) be a cyclic subgroup
of prime order \(n\), stable under the indicated endomorphism.
\begin{enumerate}[label=\textup{(\roman*)}]
 \item On the \(j=1728\) subfamily, if
       \(\phi(P)=[\lambda]P\), then
       \[
                    \lambda^2+1\equiv0\pmod n.
       \]
 \item On the \(j=0\) subfamily, if
       \(\psi(P)=[\lambda]P\), then
       \[
                    \lambda^2+\lambda+1\equiv0\pmod n.
       \]
\end{enumerate}
Consequently a scalar \(m\) may be decomposed as
\[
       m\equiv m_1+m_2\lambda\pmod n,
       \qquad |m_1|,|m_2|=O(\sqrt n),
\]
and evaluated as
\[
       [m]P=[m_1]P+[m_2]\vartheta(P),
       \qquad \vartheta=\phi\ \text{or}\ \psi.
\]
\end{theorem}

\begin{proof}
The proof first obtains the eigenvalue congruences, then constructs a short
lattice basis from the relevant CM order, and finally reduces the scalar
coset and verifies the joint-multiplication identity.
Because \(G\) is cyclic of prime order and is stable under the indicated
endomorphism, the restriction of that endomorphism to \(G\) is
multiplication by a unique
\(\lambda\in\mathbb Z/n\mathbb Z\).  On the \(j=1728\) locus,
\[
                         [\lambda^2]P=\phi^2(P)=-P.
\]
Thus \([\lambda^2+1]P=O\).  Since \(P\) has exact order \(n\), this is
equivalent to
\[
                         \lambda^2+1\equiv0\pmod n.
\]
On the \(j=0\) locus, applying
\(\psi^2+\psi+[1]=[0]\) to \(P\) gives
\[
                         [\lambda^2+\lambda+1]P=O,
\]
and hence
\[
                         \lambda^2+\lambda+1\equiv0\pmod n.
\]
This proves the two eigenvalue equations.

We next justify the balanced decomposition.  Define
\[
 L_\lambda=
 \{(a,b)\in\mathbb Z^2:a+b\lambda\equiv0\pmod n\}.
\]
The homomorphism
\[
 \pi:\mathbb Z^2\longrightarrow\mathbb Z/n\mathbb Z,\qquad
 (a,b)\longmapsto a+b\lambda
\]
is surjective because \(\pi(1,0)=1\), and its kernel is \(L_\lambda\).
The first isomorphism theorem therefore gives
\[
             [\mathbb Z^2:L_\lambda]=n,\qquad
             \det L_\lambda=n.
\]

The CM equation supplies more information than the determinant alone.  In
the fourth-order case identify \((a,b)\) with
\(a+bi\in\mathbb Z[i]\).  Evaluation at \(i=\lambda\) defines a
surjective ring map
\[
 \mathbb Z[i]\longrightarrow\mathbb F_n,\qquad
 a+bi\longmapsto a+b\lambda.
\]
Its kernel is the ideal corresponding to \(L_\lambda\), and the quotient
has \(n\) elements.  The Gaussian integers are Euclidean, so this ideal is
principal, say \((\alpha)\) with \(\alpha=a+bi\).  Its ideal norm is \(n\);
therefore
\[
                         a^2+b^2=n.
\]
As a \(\mathbb Z\)-lattice it has the orthogonal basis
\[
                         (a,b),\qquad(-b,a),
\]
whose two vectors both have length \(\sqrt n\).  The exceptional prime
\(n=2\) is bounded and can be handled directly.

In the third-order case use the Eisenstein ring
\(\mathbb Z[\varpi]\), with \(\varpi^2+\varpi+1=0\), and evaluate
\(\varpi\) at \(\lambda\).  The kernel again has index and ideal norm
\(n\).  Since \(\mathbb Z[\varpi]\) is Euclidean, it has a generator
\(\alpha=a+b\varpi\) satisfying
\[
                         a^2-ab+b^2=n.
\]
Multiplication by \(\varpi\) shows that the corresponding lattice has
basis
\[
                         (a,b),\qquad(-b,a-b).
\]
Under the usual complex embedding these vectors have length \(\sqrt n\)
and angle \(120^\circ\).  The ramified prime \(n=3\) is another bounded
case and may be treated directly.

For a given scalar \(m\), all pairs satisfying
\[
                         m_1+m_2\lambda\equiv m\pmod n
\]
form the coset \((m,0)+L_\lambda\).  Express \((m,0)\) in one of the
preceding real bases and round both coefficients to nearest integers.
Subtracting the resulting lattice vector leaves a representative
\((m_1,m_2)\) in a fundamental parallelogram whose diameter is at most the
sum of the two basis lengths.  In the Eisenstein case, the identity
\[
 |x+y\varpi|^2=x^2-xy+y^2
    \ge\frac12(x^2+y^2)
\]
shows that the complex norm and the ordinary coordinate norm are uniformly
equivalent; the Gaussian case is already orthonormal.  Hence, for an
absolute constant \(C\),
\[
                         |m_1|,|m_2|\le C\sqrt n.
\]
This is the asserted \(O(\sqrt n)\) bound.  In implementations, Gauss
reduction and nearest-plane reduction produce the same type of
decomposition directly from \(n\) and \(\lambda\), as in the standard GLV
procedure \cite{GallantLambertVanstone2001}.

Finally, because \(\vartheta(P)=[\lambda]P\),
\[
\begin{aligned}
[m_1]P+[m_2]\vartheta(P)
  &=[m_1]P+[m_2\lambda]P\\
  &=[m_1+m_2\lambda]P=[m]P.
\end{aligned}
\]
This proves the decomposition and the joint-multiplication identity.
\end{proof}

\begin{table}[H]
\centering
\small
\setlength{\tabcolsep}{3pt}
\caption{Low-cost endomorphisms available on \(\mathcal C_d\)}
\label{tab:Cd-endomorphism-comparison}
\begin{tabular}{L{2.2cm}L{2.7cm}L{3.3cm}L{2.4cm}L{2.3cm}}
\toprule
locus & field condition & Kummer action & affine cost & CM relation\\
\midrule
\(d=1/8\), \(j=1728\)
& \(i^2=-1\in k\)
& \((X:Z)\mapsto(-X:Z)\)
& \(\Dconst{i}+\Inv\)
& \(\phi^2+1=0\)\\
\(A^2=3\), \(j=0\)
& \(\zeta^2+\zeta+1=0\) in \(k\)
& \((X:Z)\mapsto(\zeta X+cZ:Z)\)
& \(3\M+\Dconst{\zeta+c-1}+2\Inv\), or
  \(6\M+\Dconst{\zeta+c-1}+\Inv\)
& \(\psi^2+\psi+1=0\)\\
\bottomrule
\end{tabular}
\end{table}

For \(k=\F_q\), the first map is rational when \(-1\) is a square in
\(\F_q\); the second requires both a solution of \(A^2=3\) and a
nontrivial cube root of unity.  If these constants lie only in an
extension field, the formulas remain geometric endomorphisms but do not
give a base-field GLV map.  Completeness of the scalar-multiplication
routine still depends on the addition formulas used for the final joint
multiplication.  The endomorphism itself does not turn an incomplete
full-point chart into a complete one.  In particular, on the \(j=1728\)
locus one has \(\rho=-1\); whenever \(i\in k\), this is the square
\(i^2\).  The exact criterion of
Theorem~\ref{thm:exact-native-completeness-criterion} therefore excludes
the single nonsquare-\(\rho\) Segre tuple in precisely that base-field CM
setting.  A constant-time implementation must use the native complete
atlas or another full-point law whose completeness hypotheses have been
verified.

\chapter{Characteristic Three}
\label{ch:char3}
Characteristic three belongs to the odd-characteristic branch for the basic
group law, but Frobenius changes the organization of tripling.  We therefore
start by specializing native addition and differential addition, then derive
the Frobenius-sensitive tripling maps, and only afterward distinguish scalar
tripling from separable degree-three quotients.  This order prevents an
inexpensive cubing stage from being mistaken for a complete multiplication
map or for a V\'elu isogeny.

\section{Basic native addition and doubling}

Let
\[
 \mathcal C_d:\quad (u^2+u)(v^2+v)=d
\]
over a field of characteristic three, with \(d\ne0,1\).  Here
\[
 \rho=1-d,\qquad r=2u+1=1-u,\qquad s=2v+1=1-v.
\]
For \(P_i=(u_i,v_i)\), define
\[
 C=s_1s_2,\quad D=r_1r_2,\quad E=CD,\quad
 H=C-D,\quad I=s_1r_2+r_1s_2.
\]
The basic \(\mathcal C_d\)-addition law specializes to
\begin{equation}\label{eq:char3-native-add}
 \boxed{\quad
 u_3=1-\frac{E-\rho}{H},\qquad
 v_3=1-\frac{E+\rho}{I}.
 \quad}
\end{equation}
For doubling, put \(K=r^2s^2\); then
\begin{equation}\label{eq:char3-native-double}
 \boxed{\quad
 u_{2P}=1-\frac{K-\rho}{s^2-r^2},\qquad
 v_{2P}=1+\frac{K+\rho}{rs}.
 \quad}
\end{equation}
The second plus sign occurs because \(2=-1\).  These are formulas on
\(\mathcal C_d\), not reductions of an equation whose definition required
division by three.
Using simultaneous inversion, the affine addition costs
\(9\M+\Inv\), and the affine double costs
\(7\M+2\Sqr+\Inv\).  The corresponding two-inversion counts are
\(6\M+2\Inv\) and \(4\M+2\Sqr+2\Inv\), respectively.

\begin{proof}
Equations~\eqref{eq:char3-native-add} and
\eqref{eq:char3-native-double} are obtained from
\eqref{eq:native-shift-add-thesis} and
\eqref{eq:native-affine-dbl-thesis} by using
\[
 16=1,\quad 2=-1,\quad u=1-r,\quad v=1-s.
\]
Their derivation uses only division by \(2\), which is valid in
characteristic three.  Hence no inseparable \([3]\)-operation is hidden in
the basic group law.
\end{proof}

\begin{example}[Native arithmetic over \(\F_9\)]
\label{ex:char3-native-F9}
Let
\[
 \F_9=\F_3[\alpha]/(\alpha^2+1),\qquad d=2,
\]
so \(\rho=2\).  On \(\mathcal C_2\), take
\[
             P=(\alpha,2\alpha).
\]
The native doubling formula gives
\[
             2P=(1+2\alpha,1),
\]
and applying \eqref{eq:char3-native-add} to \(P\) and \(2P\) gives
\[
             3P=(\alpha,2+\alpha).
\]
Direct substitution verifies all three points on
\((u^2+u)(v^2+v)=2\).
\end{example}

\section{Differential addition in characteristic three}

The model \(\Cd\) is smooth in characteristic three for
\(d\ne0,1\), because \(1/16=1\) in \(\F_3\).  All the odd-characteristic
coordinate dictionaries remain valid.  Nevertheless, arithmetic involving
\([3]\), Hessian equations, and three-isogenies changes qualitatively:
the inseparable degree of multiplication by three is nontrivial, and cubing
is Frobenius.

The diagonal Hessian identity
\[
 aX^3+Y^3+Z^3=3XYZ
\]
cannot be used as an ordinary nonsingular Hessian equation after reduction
modulo three.  Its right side vanishes and
\(Y^3+Z^3=(Y+Z)^3\).  Thus a characteristic-not-three Hessian
three-isogeny decomposition may be transported to \(\Cd\) only before
reduction and only when its level-three structure is defined over the base
field.  It is not a universal characteristic-three formula.

\begin{theorem}[Characteristic-three native \(x\)ADD]\label{thm:char3-xadd}
Let \(k\) have characteristic three and \(d\ne0,1\).  In the native
Kummer coordinate \((X:Z)=(u+1:u)\),
Theorem~\ref{thm:xadd-thesis} remains valid:
\[
\begin{aligned}
C&=(X_1+Z_1)(X_2-Z_2),\\
D&=(X_2+Z_2)(X_1-Z_1),\\
(X_{P+Q}:Z_{P+Q})
 &=\bigl(Z_\Delta(C+D)^2:X_\Delta(C-D)^2\bigr).
\end{aligned}
\]
Its cost is \(4\M+2\Sqr\), or \(3\M+2\Sqr\) for an affine known
difference.  Together with \(x\)DBL, an affine-difference ladder step costs
\(5\M+4\Sqr+\Dpar\).
\end{theorem}

\begin{proof}
The derivation of Theorem~\ref{thm:xadd-thesis} uses only
\(\charac k\ne2\), the nonsingularity \(A^2-4\ne0\), and homogeneous
polynomial identities.  In characteristic three,
\[
 A^2-4=A^2-1=\frac{\rho}{d^2}
\]
up to a nonzero square, and \(\rho=1-d\ne0\).  Hence the Montgomery model
is smooth and every step of the derivation remains valid.  The cost graph
contains no division by three.  Therefore both the formula and its count
survive unchanged.
\end{proof}

\begin{corollary}[First native double in characteristic three]
For \((X:Z)=(u+1:u)\), the first double still costs
\(\M+\Sqr+\Dpar\), because \(X-Z=1\).
\end{corollary}

\begin{proof}
The identity \(X-Z=(u+1)-u=1\) is unchanged in characteristic three, so
\(BB=(X-Z)^2=1\) in the native \(x\)DBL formula.  Therefore
\[
 (X_2:Z_2)=
 \bigl(AA:E(1+\alpha_{24}E)\bigr),
 \qquad AA=(2u+1)^2,\quad E=AA-1.
\]
Here \(\alpha_{24}=(16d)^{-1}=d^{-1}\) is defined because \(d\ne0\).
One squaring, one parameter multiplication, and one final product give the
stated cost.
\end{proof}

\section{Frobenius-specialized tripling}

On the auxiliary curve
\[
        Z^2=f(U)=U^3+AU^2+U,
\]
let
\[
        P_3(U)=3U^4+4AU^3+6U^2-1.
\]
In characteristic three this reduces to
\begin{equation}\label{eq:P3-char3}
        P_3(U)=AU^3-1.
\end{equation}

\begin{theorem}[Characteristic-three Kummer tripling]
\label{thm:char3-tripling}
For a projective Kummer input \((X:Z)\), put
\[
       T=X^3,\qquad W=Z^3.
\]
Then
\begin{equation}\label{eq:char3-tripling}
 \kappa([3]P)=
 \bigl(T(T-AW)^2:W(AT-W)^2\bigr).
\end{equation}
The cost is
\[
       2\M+2\Sqr+2\Dpar+2\Cube.
\]
\end{theorem}

\begin{proof}
Apply Lemma~\ref{lem:local-low-degree-division-tripling} before reducing
its integral identities modulo three.  In characteristic three,
\[
 P_3=AU^3-1,\qquad
 F_4=U^6-AU^5+2U^4+U^2+AU+2,
\]
and \(-8=1\).  Therefore
\[
 \phi_3=UP_3^2+fF_4.
\]
The two summands expand as
\[
\begin{aligned}
 UP_3^2&=A^2U^7+AU^4+U,\\
 fF_4&=U^9-A^2U^7+AU^6-AU^4+A^2U^3-U.
\end{aligned}
\]
Adding them cancels the terms of degrees seven, four, and one and gives
\[
 \phi_3=U^9+AU^6+A^2U^3.
\]
Since \(-2=1\) in characteristic three,
\[
 U^9+AU^6+A^2U^3
   =U^3(U^3-A)^2,
\]
whereas
\[
 \psi_3^2=(AU^3-1)^2.
\]
Therefore
\[
 U([3]P)=
 \frac{U^3(U^3-A)^2}{(AU^3-1)^2}.
\]
Homogenizing with \(U=X/Z\) gives
\eqref{eq:char3-tripling}.  Two Frobenius cubing operations form \(T,W\); two
constant products form \(AW,AT\); two squares and two final products form
the output.
\end{proof}

\section{Full-point tripling and double-base chains}

The Frobenius circuit \eqref{eq:char3-tripling} returns only
\(\kappa_d([3]P)\).  A full native output follows by specializing
Theorem~\ref{thm:native-full-tripling}.  Namely, for
\[
 r=1-u,\quad s=1-v,\quad K=r^2s^2,\quad
 \Delta=s^2-r^2,\quad \rho=1-d,
\]
put
\begin{equation}\label{eq:char3-full-tripling-blocks}
\begin{aligned}
 N_s&=-Kr^2-\rho(r^2+s^2),&
 D_s&=K^2-\rho^2-\rho\Delta,\\
 N_r&=Ks^2+\rho(r^2+s^2),&
 D_r&=K^2-\rho^2+\rho\Delta.
\end{aligned}
\end{equation}
Then, when \(N_sN_r\ne0\),
\begin{equation}\label{eq:char3-full-tripling-native}
 \boxed{\quad
 [3](u,v)=
 \left(1-\frac{rD_r}{N_r},
       1-\frac{sD_s}{N_s}\right).
 \quad}
\end{equation}
Thus characteristic three has both a very small native Kummer tripling
and an explicit full-point formula on \(\mathcal C_d\).  The former cannot
by itself supply a Miller line or choose between the two points in a Kummer
fiber.

\begin{proof}
Reduce \eqref{eq:native-tripling-blocks} modulo three.  Since
\(3=0\) and \(2=-1\), the four blocks become
\eqref{eq:char3-full-tripling-blocks}.  Moreover \(u=1-r\) and
\(v=1-s\), so \eqref{eq:native-full-tripling-shift} becomes
\eqref{eq:char3-full-tripling-native}.
\end{proof}

For secret variable-base scalars, the bitwise native Kummer ladder in
characteristic three provides
a regular constant-pattern \(x\)DBLADD schedule.  For public or fixed
scalars, a sparse \(\{2,3\}\)-double-base chain can exploit the
characteristic-three Kummer tripling formula together with ordinary
doubling and differential addition.  The especially small cost of
\eqref{eq:char3-tripling} makes this mixed \(x\)DBL/\(x\)TPL strategy
particularly effective when Frobenius cubing is inexpensive.  The two
scalar engines are therefore assigned to the following complementary
tasks:
\[
\begin{array}{c|c}
\text{task}&\text{native arithmetic engine}\\ \hline
\text{secret variable-base scalar}
  &\text{bitwise native Kummer ladder}\\
\text{public/fixed sparse \(2,3\)-chain}
  &\text{mixed \(x\)DBL/\(x\)TPL chain}
\end{array}
\]

\section{Degree-three maps and separability}
\label{sec:char3-degree-three-maps}

The inexpensive cubing operations in Theorem~\ref{thm:char3-tripling} do not mean
that multiplication by three has become a degree-three separable map.  In
characteristic three one must distinguish three objects:
\[
 \begin{array}{c|c|c}
 \text{map}&\text{degree}&\text{kernel type}\\
 \hline
 [3]:E\longrightarrow E&9&\text{the full group scheme }E[3],\\
 F_{3,d}:E\longrightarrow E^{(3)}&3&\text{local; purely inseparable},\\
 \varphi_K:E\longrightarrow E/K&3&
     \text{\'etale when }K\simeq\mathbb Z/3\mathbb Z.
 \end{array}
\]
Here \(E=\overline{\mathcal C}_d\), \(E^{(3)}\) is its Frobenius twist,
and the last row exists only after a separable cyclic kernel has been
specified.  The purpose of this section is to identify which of these maps
is being evaluated by each formula in this chapter.

\begin{proposition}[The characteristic-three factorization of tripling]
\label{prop:char3-FV-distinction}
Let \(k\) be a perfect field of characteristic three and let
\(d\ne0,1\).  The native coordinate map
\begin{equation}\label{eq:char3-native-relative-Frobenius}
 F_{3,d}:\overline{\mathcal C}_d\longrightarrow
          \overline{\mathcal C}_{d^3},
 \qquad (u,v)\longmapsto(u^3,v^3)
\end{equation}
is the relative Frobenius and is purely inseparable of degree three.  If
\(V_{3,d}:\overline{\mathcal C}_{d^3}\to
\overline{\mathcal C}_d\) denotes its dual Verschiebung, then
\begin{equation}\label{eq:char3-FV-factorization}
             [3]=V_{3,d}\circ F_{3,d}.
\end{equation}
Moreover:
\begin{enumerate}[label=\textup{(\roman*)}]
 \item if \(E\) is ordinary, \(V_{3,d}\) is separable of degree three;
       over \(\bar k\), \(E[3]\) has one local factor and one reduced
       cyclic factor of order three;
 \item if \(E\) is supersingular, \(V_{3,d}\) is also purely inseparable,
       \(E[3](\bar k)=\{O\}\), and \([3]\) has inseparable degree nine;
 \item a V\'elu sum over point coordinates defines a degree-three quotient
       only for a reduced cyclic subgroup
       \(K=\{O,Q,-Q\}\).  It neither sums over the scheme-theoretic
       \(E[3]\) nor by itself represents \([3]\).
\end{enumerate}
If \(d\in\F_3\), source and target in
\eqref{eq:char3-native-relative-Frobenius} have the same equation, but
\(F_{3,d}\) is still not \([3]\): their degrees are three and nine,
respectively.
\end{proposition}

\begin{proof}
In characteristic three the binomial identity gives
\[
 \bigl((u^3)^2+u^3\bigr)\bigl((v^3)^2+v^3\bigr)
 =\bigl((u^2+u)(v^2+v)\bigr)^3=d^3.
\]
Thus \eqref{eq:char3-native-relative-Frobenius} has the stated target.
On the bihomogeneous completion it is obtained by cubing each coordinate
in both \(\PP^1\)-factors, so it extends across the four boundary points
and sends the marked identity to the marked identity.  The induced
function-field inclusion is
\[
 k\bigl(u^3,v^3\bigr)\subset k(u,v).
\]
For a smooth curve over a perfect field this extension is purely
inseparable of degree three.  Hence \(F_{3,d}\) is an isogeny of degree
three.  Its dual isogeny exists and, by the defining property of the dual,
satisfies
\[
 V_{3,d}\circ F_{3,d}=[\deg F_{3,d}]=[3],
\]
which proves \eqref{eq:char3-FV-factorization}.

For an ordinary elliptic curve in characteristic three, the connected--
\'etale sequence of the three-torsion has geometric form
\[
       0\longrightarrow\mu_3\longrightarrow E[3]
       \longrightarrow\mathbb Z/3\mathbb Z\longrightarrow0.
\]
The kernel of the relative Frobenius is the connected factor; the kernel
of the dual Verschiebung is reduced.  Therefore \(F_{3,d}\) is
inseparable and \(V_{3,d}\) is separable, each of degree three.  In the
supersingular case the three-torsion has no reduced nonidentity geometric
point.  Both degree-three factors are then purely inseparable, so their
composition has inseparable degree nine.

Finally, the point-sum form of V\'elu's construction requires a finite
separable kernel: its summands are indexed by the distinct nonidentity
geometric points of that kernel.  The reduced subgroup
\(K=\{O,Q,-Q\}\) meets this requirement, whereas a local kernel has no
such three-point indexing set and the full group scheme \(E[3]\) has
scheme-theoretic length nine.  This proves all three distinctions.  When
\(d^3=d\), the Frobenius twist is represented by the same
\(\mathcal C_d\)-equation, but the degree calculation is unchanged and
precludes \(F_{3,d}=[3]\).
\end{proof}

Theorem~\ref{thm:char3-tripling} makes the factorization visible on the
native Kummer line.  The coordinatewise cubing map
\[
                  (X:Z)\longmapsto(T:W)=(X^3:Z^3)
\]
is the Kummer shadow of \(F_{3,d}\); the remaining degree-three tuple
\[
        (T:W)\longmapsto
        \bigl(T(T-AW)^2:W(AT-W)^2\bigr)
\]
is the Kummer shadow of the Verschiebung after the stated parameter
identification.  Consequently cheap cubing reduces the cost of the first
factor, but does not remove the second factor.  For isogeny construction,
one must instead exhibit a reduced kernel and use the separable formulas of
Chapter~\ref{ch:isogenies}.  Example~\ref{ex:F27-Frob-isogeny} gives an
explicit ordinary/supersingular comparison over \(\F_{27}\).

\part{Native Arithmetic in Characteristic Two}
\partoverview{The defining equation is unchanged, but its two quadratic
factors are now Artin--Schreier polynomials.  This part therefore restarts
from the native group law instead of reducing odd-characteristic formulas.
It develops the \(Z/4\mathbb Z\) geometry, the native binary Kummer ladder,
full-coordinate recovery, tripling, and a precise comparison with binary
Edwards and related normal forms.}

\chapter[Artin--Schreier and Z/4Z Geometry]
{Artin--Schreier and \texorpdfstring{\(Z/4\mathbb Z\)}{Z/4Z} Geometry}
\label{ch:binary-geometry}
Throughout Part III, \(\charac k=2\) and \(d\ne0\).  The affine equation is
\begin{equation}\label{eq:AS-thesis}
       v^2+v=\frac{d}{u^2+u}.
\end{equation}
It is an Artin--Schreier double cover with poles above \(u=0,1\).  This
description explains both the genus and the linear nature of point
decompression.

\section{Basic native addition and doubling}

The identity is again the boundary point \(O=(0,\infty)\), and
\[
              -(u,v)=(u,v+1).
\]
Thus \(u\) generates the Kummer function field; the globally normalized
projective coordinate remains \(\kappa=(u+1:u)\).  The following formulas
give the full group operation before any use of the \(Z/4\mathbb Z\)-normal
form.

\begin{theorem}[Native binary affine addition]
\label{thm:binary-native-affine-add}
Let \(P_i=(u_i,v_i)\in\mathcal C_d(k)\), with
\[
        U_i=\frac{u_i+1}{u_i}.
\]
Assume \(u_1u_2(U_1+U_2)\ne0\), and put
\begin{equation}\label{eq:binary-native-add-blocks}
\lambda=\frac{U_1v_1+U_2v_2}{U_1+U_2},\qquad
\nu=\frac{dU_1U_2(v_1+v_2)}{U_1+U_2},
\end{equation}
\[
              H=\lambda^2+\lambda+d(U_1+U_2).
\]
If \(H(H+d)\ne0\), then
\begin{equation}\label{eq:binary-native-affine-add}
 \boxed{\quad
 u_3=\frac d{H+d},\qquad
 v_3=\lambda+1+\frac{\nu}{H}.
 \quad}
\end{equation}
If \(U_1,U_2\) are already available, a literal schedule using one
inversion for \(U_1+U_2\) and a simultaneous inversion of \(H,H+d\)
costs
\[
             10\M+\Sqr+3\Dpar+2\Inv.
\]
Starting from raw \(u_1,u_2\), batch formation of
\(U_i=(u_i+1)/u_i=1+u_i^{-1}\) adds \(3\M+\Inv\).
\end{theorem}

\begin{proof}
Introduce, only for the proof,
\[
 X_i=dU_i,\qquad Y_i=X_iv_i.
\]
Then \(P_i\) lies on the explicitly stated binary Weierstrass equation
\[
             W_d^+:\qquad Y^2+XY=X^3+d^2X.
\]
For \(X_1\ne X_2\), its chord parameters are
\[
 \lambda_W=\frac{Y_1+Y_2}{X_1+X_2}=\lambda,\qquad
 \nu_W=\frac{X_1Y_2+X_2Y_1}{X_1+X_2}=\nu,
\]
and the sum satisfies
\[
 X_3=\lambda^2+\lambda+X_1+X_2=H,\qquad
 Y_3=(\lambda+1)X_3+\nu.
\]
The inverse dictionary is
\[
             u_3=\frac d{X_3+d},\qquad v_3=\frac{Y_3}{X_3}.
\]
Substitution yields \eqref{eq:binary-native-affine-add}.  Hence both the
input and the stated output are native \(\mathcal C_d\)-coordinates.
For the cost, the two numerators in
\eqref{eq:binary-native-add-blocks}, their common division, and the
formation of \(H\) use \(6\M+\Sqr+2\Dpar+\Inv\).
Batch-inverting \(H\) and \(H+d\), multiplying the two outputs, and
forming \(d/(H+d)\) add \(4\M+\Dpar+\Inv\).
The standard two-input batch inversion first forms
\(p=u_1u_2\), then computes
\(u_1^{-1}=u_2p^{-1}\) and \(u_2^{-1}=u_1p^{-1}\).
It therefore forms \(U_i=1+u_i^{-1}\) from raw native abscissae in
\(3\M+\Inv\).
\end{proof}

\begin{theorem}[Native binary affine doubling]
\label{thm:binary-native-affine-double}
Let \(F=u^2+u\).  Whenever the denominators are nonzero,
\begin{equation}\label{eq:binary-native-affine-double}
 \boxed{\quad
 u_{2P}=\frac{F^2}{F^2+d},\qquad
 v_{2P}=u^2+v^2.
 \quad}
\end{equation}
Reusing the square \(u^2\), the literal affine schedule costs
\[
                 \M+3\Sqr+\Inv.
\]
\end{theorem}

\begin{proof}
On \(W_d^+:Y^2+XY=X^3+d^2X\), the tangent slope at
\((X,Y)=(d(u+1)/u,Xv)\) is
\[
 \lambda=\frac{Y+X^2+d^2}{X}
         =v+d\frac{U^2+1}{U}
         =v+\frac d{u^2+u}=v^2.
\]
Therefore
\[
 X(2P)=\lambda^2+\lambda=(v^2+v)^2=\frac{d^2}{F^2},
\]
and the inverse dictionary gives \(u_{2P}=F^2/(F^2+d)\).
The tangent intercept is \(X^2+d^2\), so
\[
 \frac{Y(2P)}{X(2P)}
 =\lambda+1+\frac{X^2+d^2}{X(2P)}
 =v^2+1+(U^2+1)F^2.
\]
Since \((U^2+1)F^2=(u+1)^2\), the last expression is
\(u^2+v^2\), proving the second formula.
The squares \(u^2,F^2,v^2\), one inversion of \(F^2+d\), and one
final multiplication give the stated cost.
\end{proof}

\begin{example}[Native binary addition over \(\F_{2^8}\)]
\label{ex:binary-native-basic}
Use the polynomial basis
\[
\F_{2^8}=\F_2[\alpha]/(\alpha^8+\alpha^4+\alpha^3+\alpha+1)
\]
and fix the following polynomial-basis encoding.  If an eight-bit integer
\(h=\sum_{i=0}^7h_i2^i\), with \(h_i\in\{0,1\}\), is written as the
hexadecimal byte \(\mathtt{0xHH}\), then
\[
 \langle\mathtt{0xHH}\rangle_\alpha
 :=\sum_{i=0}^7h_i\alpha^i\in\F_{2^8}.
\]
For example,
\[
 \hexalpha{46}=\alpha^6+\alpha^2+\alpha,
 \qquad
 \hexalpha{5C}=\alpha^6+\alpha^4+\alpha^3+\alpha^2.
\]
Thus the coordinates below are field elements, not decimal integers or
formal two-character strings.  Take
\[
 d=\hexalpha{13}=\alpha^4+\alpha+1,
 \qquad
 P=(u,v)=(\hexalpha{46},\hexalpha{5C}).
\]
The curve equation may already be checked at the input from
\[
 u^2+u=\hexalpha{F9},\qquad v^2+v=\hexalpha{BC},
 \qquad \hexalpha{F9}\hexalpha{BC}=\hexalpha{13}.
\]
Formula~\eqref{eq:binary-native-affine-double} gives
\[
             2P=(\hexalpha{E2},\hexalpha{5F}).
\]
For the addition \(P+2P\), the intermediate values are
\[
 U_1=\hexalpha{F4},\quad U_2=\hexalpha{D7},\quad
 \lambda=\hexalpha{BE},\quad \nu=\hexalpha{50},\quad
 H=\hexalpha{64}.
\]
Equation~\eqref{eq:binary-native-affine-add} then gives
\[
             3P=(\hexalpha{A9},\hexalpha{9E}).
\]
\end{example}

\section[Linear equivalence with the Z/4Z-normal form]
{Linear equivalence with the
\texorpdfstring{\(Z/4\mathbb Z\)}{Z/4Z}-normal form}

\begin{theorem}\label{thm:z4-thesis}
Put \(d_{\rm K}=d^{-1}\) and make the projective linear substitutions
\[
(\widetilde U_0:\widetilde U_1)=(U_0+U_1:U_1),\qquad
(\widetilde V_0:\widetilde V_1)=(V_0+V_1:V_1).
\]
On the affine charts, these substitutions are
\[
             x=\frac{u}{u+1},\qquad y=\frac{v}{v+1}.
\]
Equation~\eqref{eq:homogeneous-thesis} becomes
\begin{equation}\label{eq:z4-p1-thesis}
(\widetilde U_0+\widetilde U_1)^2
(\widetilde V_0+\widetilde V_1)^2
=d_{\rm K}\widetilde U_0\widetilde U_1
  \widetilde V_0\widetilde V_1,
\end{equation}
or \((1+x)^2(1+y)^2=d_{\rm K}xy\).

Under the skew Segre embedding
\[
(X_0:X_1:X_2:X_3)=
(\widetilde U_0\widetilde V_0:
 \widetilde U_1\widetilde V_0:
 \widetilde U_1\widetilde V_1:
 \widetilde U_0\widetilde V_1),
\]
the equations are
\begin{equation}\label{eq:z4-p3-thesis}
X_0X_2=X_1X_3,\qquad
(X_0+X_1+X_2+X_3)^2
=d_{\rm K}X_0X_2=d_{\rm K}X_1X_3.
\end{equation}
\end{theorem}

\begin{proof}
Since \(u=x/(1+x)\) in characteristic two,
\[
       u^2+u=\frac{x}{(1+x)^2},
\qquad v^2+v=\frac{y}{(1+y)^2}.
\]
Substitution into \eqref{eq:model} gives
\eqref{eq:z4-p1-thesis}.  The skew Segre image satisfies
\(X_0X_2=X_1X_3\), and
\[
 X_0+X_1+X_2+X_3
 =(\widetilde U_0+\widetilde U_1)
  (\widetilde V_0+\widetilde V_1).
\]
Squaring and applying \eqref{eq:z4-p1-thesis} proves
\eqref{eq:z4-p3-thesis}.
\end{proof}

\begin{remark}
If the original Segre coordinates are ordered as
\[
(S_0:S_1:S_2:S_3)
=(U_0V_0:U_1V_0:U_1V_1:U_0V_1),
\]
the ambient linear map is
\[
(X_0:X_1:X_2:X_3)
=(S_0+S_1+S_2+S_3:S_1+S_2:S_2:S_3+S_2).
\]
Thus the addition-law space and Kummer projections of the normal form
transport linearly, not through a costly general birational map.
\end{remark}

\section{Binary Weierstrass form}

\begin{theorem}\label{thm:binary-W-thesis}
On \(u\ne0\), define
\begin{equation}\label{eq:binary-W-map-thesis}
        X=d\frac{u+1}{u},\qquad Y=Xv.
\end{equation}
Then
\begin{equation}\label{eq:binary-W-thesis}
        Y^2+XY=X^3+d^2X.
\end{equation}
The inverse is
\[
        u=\frac{d}{X+d},\qquad v=\frac{Y}{X}.
\]
The discriminant is \(d^4\) and
\begin{equation}\label{eq:binary-j-thesis}
        j(\Cd)=d^{-4}.
\end{equation}
\end{theorem}

\begin{proof}
By \eqref{eq:AS-thesis},
\[
Y^2+XY=X^2(v^2+v)=X^2\frac{d}{u^2+u}.
\]
Since \(u^2+u=u(u+1)\), substitution of
\(X=d(u+1)/u\) gives
\begin{align*}
 X^2\frac d{u^2+u}
 &=\frac{d^3(u+1)}{u^3},\\
 X^3+d^2X
 &=\frac{d^3(u+1)^3}{u^3}
   +\frac{d^3(u+1)}u\\
  =\frac{d^3(u+1)}{u^3},
\end{align*}
where the last equality uses
\((u+1)^2+u^2=1\) in characteristic two.  This proves the Weierstrass
equation.  From \(Xu=d(u+1)\) one obtains
\((X+d)u=d\), hence \(u=d/(X+d)\); the definition \(Y=Xv\) gives
\(v=Y/X\).  These formulas are inverse on their common dense affine
open and therefore extend uniquely to the smooth projective completions.
For
\(a_1=1,a_4=d^2\), the generalized Weierstrass invariants are
\(c_4=1\) and \(\Delta=d^4\), so \(j=c_4^3/\Delta=d^{-4}\).
\end{proof}

\begin{corollary}\label{cor:binary-coverage-thesis}
Every \(\Cd\) in characteristic two is ordinary and has a distinguished
rational point of order four.  Conversely, an ordinary elliptic curve with
a marked rational four-torsion point admits a \(Z/4\mathbb Z\)-normal
embedding of the form \eqref{eq:z4-p3-thesis}.  The family therefore covers
the rational-four-torsion subfamily, not every ordinary binary elliptic
curve.
\end{corollary}

\begin{proof}
Theorem~\ref{thm:binary-W-thesis} gives
\(j(\mathcal C_d)=d^{-4}\ne0\).  An elliptic curve in characteristic two
is ordinary exactly when its \(j\)-invariant is nonzero, so every smooth
\(\mathcal C_d\) is ordinary.  Proposition~\ref{prop:inverse-torsion}
exhibits \(R=(\infty,0)\) with \(2R=T\ne O\) and \(4R=O\), proving that
the distinguished point has exact order four.

Conversely, the normal-form theorem for an ordinary elliptic curve with a
marked rational point of exact order four supplies a
\(Z/4\mathbb Z\)-normal embedding of the form
\eqref{eq:z4-p3-thesis} \cite{Kohel2012}.  The inverse ambient linear map
recorded after Theorem~\ref{thm:z4-thesis} returns this embedding to the
native Segre presentation of a member \(\mathcal C_d\).  Thus the family
covers the locus with the stated marking.  An ordinary binary elliptic
curve need not possess a rational point of order four, so ordinarity alone
is insufficient.
\end{proof}

\section{Decompression}

Given \(u\notin\{0,1\}\), put
\[
       a=\frac{d}{u^2+u}.
\]
The equation \(v^2+v=a\) has a solution in \(\F_{2^m}\) if and only if
\[
       \Tr_{\F_{2^m}/\F_2}(a)=0.
\]
Its two solutions are \(v\) and \(v+1\), exactly the two points in one
Kummer fiber.  For odd \(m\), a half-trace computes a root; for general
\(m\), the Artin--Schreier map is a fixed linear transformation that can be
solved by a precomputed linear circuit.  Thus binary decompression replaces
a generic square root by a linearized equation.

\chapter{The Native Binary Kummer Ladder}
\label{ch:binary-kummer}
Throughout this chapter,
\[
       \mathcal C_d:\quad (u^2+u)(v^2+v)=d,\qquad
       \charac k=2,\quad d\ne0.
\]
Use its characteristic-free native coordinate
\[
\kappa(P)=(X_0:X_1)=(u+1:u),\qquad
t(P)=X_1/X_0=u/(u+1).
\]
In characteristic two,
\[
                         -(u,v)=(u,v+1),
\]
so \(u\), and hence \(\kappa\), is constant on the pair
\(\{P,-P\}\).  Since the projection
\(\overline{\mathcal C}_d\to\PP^1_u\) has degree two, \(\kappa\) is the
native Kummer quotient rather than merely a convenient rational function.
The reciprocal affine coordinate
\[
                         U=\frac{X_0}{X_1}=\frac{u+1}{u}
\]
is used when a formula has lower degree in \(U\); the ladder itself uses
\(t=U^{-1}\).  The boundary normalization is
\[
                         \kappa(O)=(1:0),\qquad
                         \kappa(T)=(0:1).
\]
This coordinate is also the first projection of the
\(Z/4\mathbb Z\)-normal form, but every output below is interpreted on
\(\mathcal C_d\):
\[
        u([n]P)=\frac{X_1^{(n)}}{X_0^{(n)}+X_1^{(n)}}.
\]
Thus the computational state never ceases to represent the original
\(\mathcal C_d\)-point modulo sign.

\section{Doubling}

\begin{theorem}[Binary \(x\)DBL]\label{thm:binary-xdbl-thesis}
Let \(d_{\rm K}=d^{-1}\).  Then
\begin{equation}\label{eq:binary-xdbl-thesis}
(X_0:X_1)\longmapsto
\bigl(X_0^4+X_1^4:d_{\rm K}X_0^2X_1^2\bigr)
\end{equation}
represents \(\kappa(2P)\).  The circuit
\[
A=X_0+X_1,\quad B=X_0X_1,\quad
X'_0=A^4,\quad X'_1=d_{\rm K}B^2
\]
costs \(\M+3\Sqr+\mCurve\).
\end{theorem}

\begin{proof}
The same map was derived from the native affine double and proved
base-point-free in
Theorem~\ref{thm:binary-native-complete-diff-atlas}.  The present theorem
isolates that exceptional-chart map as the recurring binary
\(x\)DBL circuit.  To evaluate it, form
\(A=X_0+X_1\), \(B=X_0X_1\), square \(A\) twice to obtain \(A^4\),
square \(B\) once, and multiply the latter by \(d_{\rm K}\).  This uses
\(\M+3\Sqr+\mCurve\).  Thus the correctness proof remains in the complete
atlas chapter, while this chapter records only the implementation
dependency graph.
\end{proof}

\section{Differential addition}

\begin{theorem}[Binary \(x\)ADD]\label{thm:binary-xadd-thesis}
Let
\[
\kappa(P)=(X_0:X_1),\quad
\kappa(Q)=(Y_0:Y_1),\quad
\kappa(P-Q)=(t_0:t_1),
\]
and put
\[
A=X_0Y_0+X_1Y_1,\qquad
B=X_0Y_1+X_1Y_0.
\]
Then
\begin{equation}\label{eq:binary-xadd-thesis}
       \kappa(P+Q)=\bigl(t_1A^2:t_0B^2\bigr).
\end{equation}
The oriented inputs satisfy
\begin{equation}\label{eq:binary-oriented-thesis}
t_0^2B^2+t_1^2A^2
+d_{\rm K}t_0t_1X_0X_1Y_0Y_1=0.
\end{equation}
When \((t_0:t_1)=(1:t)\) with \(t\ne0\), equivalent circuits are
\begin{align}
\kappa(P+Q)&=(tA^2:B^2),
 &&3\M+2\Sqr+\mBase,\label{eq:binary-xadd-a-thesis}\\
&=(A^2:tA^2+d_{\rm K}(X_0Y_0)(X_1Y_1)),
 &&3\M+\Sqr+\mBase+\mCurve.
\label{eq:binary-xadd-b-thesis}
\end{align}
\end{theorem}

\begin{proof}
The generic formula~\eqref{eq:binary-xadd-thesis}, its two exceptional
difference fibers, and their replacement charts were established in
Theorem~\ref{thm:binary-native-complete-diff-atlas} from the
Artin--Schreier biquadratic identity.  It remains here to derive the
oriented relation and the two optimized affine-difference circuits.
Substituting the generic output into the native
Segre equation~\eqref{eq:binary-native-Segre-repeat}, or equivalently
clearing the same Artin--Schreier identity before solving for the output,
gives \eqref{eq:binary-oriented-thesis}.

In the affine-difference chart, \eqref{eq:binary-oriented-thesis} becomes
\[
 B^2=t^2A^2+d_{\rm K}t(X_0Y_0)(X_1Y_1).
\]
Substituting this relation into the first representative and removing a
common factor \(t\) yields the second.  The restriction \(t\ne0\) is
necessary for this cancellation.  At \(t=0\), the difference is \(O\),
so \(P=Q\) and \(B=0\); the first affine-difference pair degenerates to
\((0:0)\).  The second pair, however, specializes to
\[
 \bigl((X_0+X_1)^4:d_{\rm K}X_0^2X_1^2\bigr),
\]
which is exactly the binary \(x\)DBL map
\eqref{eq:binary-xdbl-thesis}.  Thus the second circuit extends regularly
to the diagonal and supplies its replacement chart; no equivalence at
\(t=0\) is asserted.  Karatsuba computes
\(X_0Y_0\), \(X_1Y_1\), and
\((X_0+X_1)(Y_0+Y_1)\) with three multiplications.  The remaining
operations give the displayed counts.
\end{proof}

\begin{corollary}[Binary ladder costs]\label{cor:binary-ladder-cost}
Combining binary \(x\)DBL with the two \(x\)ADD circuits gives
\[
\boxed{4\M+5\Sqr+\mBase+\mCurve},
\qquad
\boxed{4\M+4\Sqr+\mBase+2\mCurve}.
\]
For a native input, \(X_0+X_1=1\), so
\[
        \kappa(2P)=\bigl(1:d_{\rm K}(u^2+u)^2\bigr)
\]
costs \(2\Sqr+\mCurve\) and no general multiplication.
\end{corollary}

\begin{proof}
Theorem~\ref{thm:binary-xdbl-thesis} costs
\(\M+3\Sqr+\mCurve\).  Adding the two alternative \(x\)ADD costs of
Theorem~\ref{thm:binary-xadd-thesis} gives, respectively,
\[
\begin{aligned}
 &4\M+5\Sqr+\mBase+\mCurve,\\
 &4\M+4\Sqr+\mBase+2\mCurve.
\end{aligned}
\]
For a native input \((X_0:X_1)=(u+1:u)\), characteristic two gives
\(X_0+X_1=1\).  Hence the first coordinate in
\eqref{eq:binary-xdbl-thesis} is \(1\), while the second is
\[
 d^{-1}X_0^2X_1^2=d^{-1}(u^2+u)^2.
\]
One squaring forms \(u^2\), a second squares \(u^2+u\), and the remaining
charged operation is multiplication by \(d^{-1}\).  This proves the
specialized first-double count.
\end{proof}

\section{Constant-time ladder and recovery}

There are two mathematically distinct initialization choices.  The
inversion-free native choice is
\[
       R_0=(1:0),\qquad R_1=(u_P+1:u_P).
\]
With this projective known difference, the general formula
\eqref{eq:binary-xadd-thesis} costs
\(5\M+2\Sqr\) for \(x\)ADD and hence
\[
       6\M+5\Sqr+\mCurve
\]
per \(x\)DBLADD step.  Alternatively, one may first normalize the fixed
difference to
\[
       R_1=(1:t_P),\qquad t_P=\frac{u_P}{u_P+1}.
\]
Starting from a raw native input, this normalization costs one inversion
and additions, because in characteristic two
\(t_P=1+(u_P+1)^{-1}\).  The two affine-difference circuits of
Corollary~\ref{cor:binary-ladder-cost} then give the lower recurring costs
displayed there.  This preprocessing cost is not part of those per-step
figures.

At every secret bit, perform a constant-time swap, one binary \(x\)ADD,
one binary \(x\)DBL, and the inverse swap.  The invariant is
\[
 R_0=\kappa([n]P),\qquad R_1=\kappa([n+1]P).
\]
The difference between the two registers is always the fixed input
\(\kappa(P)\); hence the differential-addition input is public and does
not depend on the processed scalar prefix.  The formulas are homogeneous,
so neither option requires normalization or inversion inside the ladder.
The specialized first double of
Corollary~\ref{cor:binary-ladder-cost} is evaluated in the unnormalized
native representative \((u_P+1:u_P)\); if the affine-difference option is
selected, its separate normalization is charged as just stated.

For full-point recovery only, convert the fixed native input
\(P=(u_P,v_P)\) to the temporary binary Weierstrass coordinates
\[
 x_P=d\frac{u_P+1}{u_P},\qquad y_P=x_Pv_P.
\]
The ladder abscissae are converted by the same
\(x=d(u+1)/u\) dictionary.

\begin{proposition}[Binary full-point recovery]\label{prop:binary-recovery-thesis}
Let \(P=(x_P,y_P)\) lie on \eqref{eq:binary-W-thesis} with \(x_P\ne0\).
If the ladder gives \(x_Q=X([m]P)\) and
\(x_R=X([m+1]P)\), with \(x_Q\) and \(x_R\) finite, then
\begin{equation}\label{eq:binary-recovery-thesis}
\lambda=
\frac{(x_P+x_Q)x_R+y_P+x_Px_Q+d^2}{x_P},
\qquad
y_Q=y_P+\lambda(x_P+x_Q).
\end{equation}
If in addition \(x_Q(x_Q+d)\ne0\), the finite native point is
\[
        u_Q=\frac{d}{x_Q+d},\qquad v_Q=\frac{y_Q}{x_Q}.
\]
When \(x_Q=0\) or \(x_Q=d\), the recovered Weierstrass point is instead
returned through the reserved native boundary encoding; the displayed
affine inverse dictionary is not evaluated.
\end{proposition}

\begin{proof}
Put \(s=x_P+x_Q\).  First assume \(s\ne0\), and let
\[
                         \lambda=\frac{y_P+y_Q}{s}.
\]
Because \(R=P+Q\), the binary addition law gives
\[
                         x_R=\lambda^2+\lambda+s
\]
and \(y_Q=y_P+\lambda s\).  Substitution in the equation of \(Q\) gives
\[
 (y_P+\lambda s)^2+x_Q(y_P+\lambda s)
   =x_Q^3+d^2x_Q.
\]
In characteristic two the cross term in the square vanishes.  Use
\[
                         y_P^2+x_Py_P=x_P^3+d^2x_P
\]
to eliminate \(y_P^2\), bring all terms to the left, and divide by \(s\).
Since
\[
 x_P^3+x_Q^3=s(x_P^2+x_Px_Q+x_Q^2),
\]
the result is
\[
 \lambda^2s+\lambda x_Q+
 x_P^2+x_Px_Q+x_Q^2+d^2+y_P=0.
\]
Replace \(\lambda^2\) by \(x_R+\lambda+s\).  The coefficient of
\(\lambda\) becomes \(s+x_Q=x_P\), while
\[
 s^2+x_P^2+x_Px_Q+x_Q^2=x_Px_Q.
\]
Consequently
\[
 \lambda x_P=(x_P+x_Q)x_R+y_P+x_Px_Q+d^2.
\]
The assumption \(x_P\ne0\) now gives the first formula in
\eqref{eq:binary-recovery-thesis}; the second is
\(y_Q=y_P+\lambda(x_P+x_Q)\).

When \(s=0\) and \(R\) is finite, the applicable case is \(Q=P\).  The
same formula reduces to
\[
                 \lambda=\frac{y_P+x_P^2+d^2}{x_P},
\]
which is the tangent slope, and it returns \(y_Q=y_P\).  The alternative
\(Q=-P\) has \(R=O\) and is a boundary output, not a finite recovery
input.  Thus the formula covers every finite ladder-recovery case allowed
by the statement.

Finally, under the additional hypotheses \(x_Q(x_Q+d)\ne0\), the inverse
dictionary in
Theorem~\ref{thm:binary-W-thesis} is
\[
        u_Q=\frac{d}{x_Q+d},\qquad v_Q=\frac{y_Q}{x_Q},
\]
which returns the recovered point to \(\mathcal C_d\).
\end{proof}

\paragraph{Boundary handling.}
The projective ladder remains meaningful when one register is
\(\kappa(O)\) or \(\kappa(T)\).  The affine recovery formula is invoked
only after masks have verified that \(x_P,x_Q\), and the required
inverse-dictionary denominators are nonzero.  A protocol that admits a
boundary output returns its reserved encoding directly; it does not feed
an infinite Weierstrass abscissa to
\eqref{eq:binary-recovery-thesis}.

\chapter[Binary Full-Point and Tripling]
{Characteristic-Two Full-Point and Tripling Arithmetic on
\texorpdfstring{\(\mathcal C_d\)}{Cd}}
\label{ch:binary-comparison}
The preceding binary Kummer chapter produced quotient arithmetic and
recovery.  This chapter first records how a linearly equivalent
\(Z/4\mathbb Z\)-chart can be used while retaining full
\(\mathcal C_d\)-inputs and outputs, then returns to the native Kummer line
for tripling, and finally restores the missing full coordinate.  Its closing
section states precisely which geometric objects are shared across
characteristics and which circuits are characteristic-specific.

\section[Full-point arithmetic through the Z/4Z chart]
{Full-point arithmetic on \(\mathcal C_d\) through its
\texorpdfstring{\(Z/4\mathbb Z\)}{Z/4Z} chart}

Throughout this chapter,
\[
 \mathcal C_d:\qquad (u^2+u)(v^2+v)=d,\qquad
 \charac k=2,\qquad d\ne0.
\]
This first section records how a full-point addition law is evaluated while
\(\mathcal C_d\) remains the declared input and output model.
Chapter~\ref{ch:comparison} subsequently gives the general comparison of
the characteristic-specific implementations.

Let \(\mathcal L\) be the ambient projective linear map displayed after
Theorem~\ref{thm:z4-thesis}.  It sends the native Segre completion of
\(\mathcal C_d\) to the \(Z/4\mathbb Z\)-normal embedding using only
coordinate additions.  If \(\mathcal A_{Z/4}\) is a homogeneous
full-point addition tuple on that normal form, then
\[
 \mathcal A_{\mathcal C_d}(P,Q)
   =\mathcal L^{-1}
      \bigl(\mathcal A_{Z/4}(\mathcal L(P),\mathcal L(Q))\bigr)
\]
is a full-point addition tuple on the original smooth completion of
\(\mathcal C_d\).  Because \(\mathcal L\) and \(\mathcal L^{-1}\) are
linear and multiplication-free in characteristic two, this pullback has
the same charged \(\M,\Sqr,\mCurve\) count.  Thus the
\(Z/4\mathbb Z\) row below is also an available
\(\mathcal C_d\)-full-point implementation after linear recoding.

For a long sequence of full additions, one may similarly convert once to a
split or twisted \(\mu_4\)-normal chart and remain there until the final
linear return to \(\mathcal C_d\).  This is different from the ladder of
Chapter~\ref{ch:binary-kummer}: the latter stores only a Kummer point,
whereas every row of the following table stores a full point.

The operation counts in the table are a comparison of established
homogeneous circuits, not an assertion preceding their formulas.  The
binary-Edwards row is the complete law of
\cite{BinaryEdwards2008}; the \(Z/4\mathbb Z\)-normal and split
\(\mu_4\)-normal rows are the circuits of \cite{Kohel2012}; and the
twisted row is the corresponding circuit of \cite{KohelTwistedMu4}.
Theorem~\ref{thm:z4-thesis} and its displayed ambient inverse prove that
transporting the \(Z/4\mathbb Z\) row to and from \(\mathcal C_d\) uses
only additions, so its charged multiplication and squaring count is
unchanged.  Thus every cost below has both a specified full-point endpoint
and an identified formula source.  All rows use the same ledger:
\(\M\) is a general field multiplication, \(\Sqr\) is a field squaring,
and \(\mCurve\) is multiplication by a fixed curve coefficient; field
additions and multiplication by \(0,1\), or \(-1\) are uncharged.  The
listed costs are per group operation after the indicated chart has been
selected.  A one-time change of chart is therefore not charged to each
row; for the displayed \(Z/4\mathbb Z\) pullback it is, in any event,
linear and uses no charged operation.

\begin{table}[htbp]
\centering
\caption{Representative characteristic-two full-point costs}
\label{tab:binary-full-thesis}
\begin{tabular}{L{4.5cm}L{2.2cm}L{3.1cm}L{3.1cm}}
\toprule
model & output & addition & doubling\\
\midrule
binary Edwards, complete & \(\Full\)
 & \(16\M+\Sqr+4\mCurve\)
 & \(2\M+5\Sqr+2\mCurve\)\\
\(\mathcal C_d\) via its \(Z/4\mathbb Z\) chart & \(\Full\)
 & \(12\M\) & \(7\M+2\Sqr\)\\
\(\mathcal C_d\) via a split \(\mu_4\) chart & \(\Full\)
 & \(7\M+2\Sqr+2\mCurve\)
 & \(2\M+5\Sqr+2\mCurve\)\\
\(\mathcal C_d\) via a twisted \(\mu_4\) chart & \(\Full\)
 & \(9\M+2\Sqr\)
 & \(2\M+5\Sqr+2\mCurve\)\\
\bottomrule
\end{tabular}
\end{table}

The table separates two implementation choices that were blurred in the
previous formulation. A protocol requiring arbitrary full additions can use a pulled-back
normal-form law and return the result explicitly to
\(\mathcal C_d\).  Quotient scalar multiplication uses the native Kummer
ladder followed by one full-point recovery when a terminal full point is
required.  The explicit \(Z/4\mathbb Z\)-normal, twisted
\(\mu_4\)-normal, and binary Edwards dictionaries extend the accompanying
full-point arithmetic framework to the remaining ordinary twist classes,
while the marked rational-four-torsion locus retains the native
\(\mathcal C_d\) input and output interface.

\section{Native binary Kummer tripling}

We now return explicitly to the native curve
\[
 \mathcal C_d:\qquad (u^2+u)(v^2+v)=d,
 \qquad \charac k=2,\qquad d\ne0.
\]
Its inverse is \(-(u,v)=(u,v+1)\), so the quotient by sign is represented
by
\begin{equation}\label{eq:binary-tripling-Kummer-coordinate}
 \boxed{\qquad
 \kappa_d(P)=(X_0:Z_0)=(u(P)+1:u(P)),\qquad
 U(P)=\frac{X_0}{Z_0}=\frac{u(P)+1}{u(P)}.
 \qquad}
\end{equation}
The symbols \((X_0:Z_0)\) are used here to emphasize the affine coordinate
\(U=X_0/Z_0\).  They represent the same projective Kummer point denoted
\((X_0:X_1)\) in Chapter~\ref{ch:binary-kummer}, where the reciprocal
ladder coordinate was \(t=X_1/X_0=U^{-1}\).  The boundary values are
\[
                         \kappa_d(O)=(1:0),\qquad
                         \kappa_d(T)=(0:1).
\]
Accordingly, the theorem below computes
\(\kappa_d([3]P)\) directly; it does not output a point on an auxiliary
Weierstrass or Edwards model.

\begin{lemma}[Low-degree generalized division identities in characteristic two]
\label{lem:binary-low-degree-division-tripling}
On
\[
 W_d^+:\qquad Y^2+XY=X^3+d^2X,\qquad d\ne0,
\]
the generalized division functions needed for the abscissa of \([3]P\)
are
\[
 \psi_2=X,\qquad
 \psi_3=X^4+X^3+d^4,\qquad
 \psi_4=X^6+d^4X^2.
\]
Thus
\begin{equation}\label{eq:binary-local-phi3}
 \phi_3=X\psi_3^2-\psi_2\psi_4
       =X^9+d^4X^3+d^8X,
\end{equation}
and, wherever \(\psi_3\ne0\),
\[
                         X([3]P)=\frac{\phi_3}{\psi_3^2}.
\]
\end{lemma}

\begin{proof}
The Weierstrass coefficients are
\[
 a_1=1,\qquad a_2=a_3=a_6=0,\qquad a_4=d^2.
\]
Consequently, after reduction in characteristic two,
\[
 b_2=1,\qquad b_4=b_6=0,\qquad b_8=d^4.
\]
The low-degree generalized division recurrences now give
\[
 \psi_2=2Y+a_1X+a_3=X
\]
and
\[
 \psi_3=3X^4+b_2X^3+3b_4X^2+3b_6X+b_8
       =X^4+X^3+d^4.
\]
For the fourth function the unreduced recurrence is
\[
 \psi_4=\psi_2\bigl(2X^6+b_2X^5+5b_4X^4+10b_6X^3
 +10b_8X^2+(b_2b_8-b_4b_6)X+b_4b_8-b_6^2\bigr).
\]
Reducing its coefficients modulo two leaves
\[
 \psi_4=X(X^5+d^4X)=X^6+d^4X^2.
\]
Finally, subtraction equals addition in characteristic two, and hence
\[
\begin{aligned}
 \phi_3
 &=X(X^4+X^3+d^4)^2+X(X^6+d^4X^2)\\
 &=X(X^8+X^6+d^8)+X^7+d^4X^3\\
 &=X^9+d^4X^3+d^8X.
\end{aligned}
\]
The generalized division-coordinate identity
\(X([3]P)=\phi_3/\psi_3^2\) then gives the last assertion.
\end{proof}

\begin{theorem}[Binary Kummer tripling]\label{thm:binary-tripling-thesis}
On the native affine Kummer coordinate,
\begin{equation}\label{eq:binary-tripling-affine-thesis}
 U([3]P)=
 \frac{d^2U^9+U^3+d^2U}{(dU^4+U^3+d)^2}.
\end{equation}
For a projective input \(U=X_0/Z_0\), define
\[
\begin{aligned}
X_2&=X_0^2,&Z_2&=Z_0^2,&X_3&=X_2X_0,\\
X_4&=X_2^2,&Z_4&=Z_2^2,&T&=X_4+Z_4,\\
H&=dT+X_3Z_0,&Z_6&=Z_2Z_4,\\
N&=d^2X_0T^2+X_3Z_6.
\end{aligned}
\]
Then
\begin{equation}\label{eq:binary-tripling-projective-thesis}
       \kappa([3]P)=(N:Z_0H^2)
\end{equation}
at cost \(6\M+6\Sqr+2\Dpar\).
\end{theorem}

\begin{proof}
Lemma~\ref{lem:binary-low-degree-division-tripling} gives the required
generalized division functions and \(\phi_3\) on
\eqref{eq:binary-W-thesis}.
Since \(X=dU\),
\[
U([3]P)=\frac{\phi_3(dU)}{d\psi_3(dU)^2},
\]
and
\[
\begin{aligned}
\phi_3(dU)
 &=d^7\{d^2U^9+U^3+d^2U\},\\
d\psi_3(dU)^2
 &=d^7\{dU^4+U^3+d\}^2.
\end{aligned}
\]
Cancellation of the nonzero common factor \(d^7\) proves
\eqref{eq:binary-tripling-affine-thesis}.  Homogenizing its numerator and
denominator to the same degree gives
\[
d^2X_0(X_0^8+Z_0^8)+X_0^3Z_0^6,
\quad
Z_0\{d(X_0^4+Z_0^4)+X_0^3Z_0\}^2.
\]
Because \(T=X_0^4+Z_0^4\) and
\(T^2=X_0^8+Z_0^8\), the first form is
\[
 N=d^2X_0T^2+X_3Z_6,
\]
and the second is \(Z_0H^2\).  The circuit uses six squares for
\(X_2,Z_2,X_4,Z_4,T^2,H^2\); six general multiplications for
\(X_3,X_3Z_0,Z_6,X_0T^2,X_3Z_6,Z_0H^2\); and two fixed-parameter
multiplications by \(d\) and \(d^2\).  This proves the homogeneous formula
and the cost \(6\M+6\Sqr+2\Dpar\).
\end{proof}

\section[Full-point tripling returned to Cd]
{Full-point tripling returned to \texorpdfstring{\(\mathcal C_d\)}{Cd}}

Let \(P=(u,v)\in\mathcal C_d(k)\) and put
\[
                         U=\frac{u+1}{u}.
\]
When the ordinate of \([3]P\) is required, introduce only as temporary
variables the coordinates
\[
 x_1=dU,\qquad y_1=x_1v
\]
on the locally stated binary Weierstrass equation
\[
                         y^2+xy=x^3+d^2x.
\]
The tangent parameters at \((x_1,y_1)\) are
\[
 \lambda=v+dU+\frac dU,\qquad
 \nu_2=x_1^2+d^2=d^2(U^2+1).
\]
Therefore the double is
\[
 x_2=\lambda^2+\lambda,\qquad
 y_2=(\lambda+1)x_2+\nu_2.
\]
For the chord through \(2P\) and \(P\), set
\[
 \mu=\frac{y_2+y_1}{x_2+x_1},\qquad
 \nu=\frac{x_2y_1+x_1y_2}{x_2+x_1}.
\]
The generalized binary chord law gives
\[
 x_3=\mu^2+\mu+x_2+x_1,\qquad
 y_3=(\mu+1)x_3+\nu.
\]
The output is then returned explicitly to the original model:
\begin{equation}\label{eq:binary-full-tripling-return-Cd}
 \boxed{\qquad
       u([3]P)=\frac{d}{x_3+d},\qquad
       v([3]P)=\frac{y_3}{x_3}.
 \qquad}
\end{equation}
Thus the Weierstrass equation is used to factor the derivation, but neither
the declared input nor the declared output has changed.

With \(U\) already available, a literal affine schedule costs
\[
          11\M+3\Sqr+4\Dpar+3\Inv.
\]
The charged inversions are the inversion of \(U\), the chord denominator
\(x_2+x_1\), and one batched inversion used for \(x_3\) and \(x_3+d\).
Starting from raw \(u\), forming
\(U=(u+1)/u=1+u^{-1}\) adds \(\Inv\) and additions.
Vanishing tangent, chord, or return denominators are handled by the
projective formulas or by constant-time boundary masks.  The full affine
schedule supplies the full-coordinate recovery companion to the native
Kummer tripling map
\eqref{eq:binary-tripling-projective-thesis}.  A quotient-only
scalar-multiplication chain remains entirely on the native Kummer line and
therefore avoids the ordinate and all affine inversions.  When a full point
is required at the endpoint, the native Kummer computation is followed by
one full-point recovery, giving a unified quotient-to-full-point arithmetic
pipeline.

\section[Uniform geometry and arithmetic]{Uniform geometry and characteristic-specific arithmetic}

The purpose of this section is to state exactly how the characteristic-two
theory fits the odd-characteristic chapters.  The uniformity of
\(\mathcal C_d\) is a uniformity of the curve model, its marked geometry,
and its arithmetic interface.  It does not mean that one polynomial
circuit can be reduced coefficientwise from odd characteristic to
characteristic two.

\begin{center}
\begin{tabular}{L{3.1cm}L{4.7cm}L{4.7cm}}
\toprule
layer & characteristic-uniform statement & characteristic-specific realization\\
\midrule
curve model
& \((u^2+u)(v^2+v)=d\) and its smooth
  \((2,2)\)-completion
& \(d(1-16d)\ne0\) in odd characteristic;
  \(d\ne0\) in characteristic two\\
marked group data
& \(O=(0,\infty)\), rational four-torsion boundary, and
  \(-(u,v)=(u,-v-1)\)
& in characteristic two the inverse becomes
  \((u,v+1)\)\\
Kummer interface
& \(\kappa_d(P)=(u(P)+1:u(P))\)
& Montgomery-type biquadratics in odd characteristic;
  Artin--Schreier biquadratics in characteristic two\\
full-point charts
& inputs and final outputs are points of \(\mathcal C_d\)
& centered reciprocal/Segre charts when \(2\ne0\);
  \(Z/4\mathbb Z\) and binary Weierstrass derivation charts when
  \(2=0\)\\
scalar multiplication
& quotient ladder followed, when necessary, by one full-point recovery
& the actual \(x\)DBL, \(x\)ADD, and recovery circuits are derived
  separately in the two characteristics\\
\bottomrule
\end{tabular}
\end{center}

This distinction has three practical consequences.  First, the same
external representation \((u+1:u)\) and the same boundary normalization
can be used by characteristic-independent protocol code.  Second, the
internal arithmetic backend must select the formulas proved for the actual
characteristic; setting \(2=0\) in an odd-characteristic denominator is
not a derivation of the binary law.  Third, a Kummer output represents the
pair \(\{P,-P\}\), so a protocol requiring a signed full point must invoke
the recovery procedure of Chapter~\ref{ch:binary-kummer}.

Finally, the binary family considered here is exactly the ordinary
rational-four-torsion subfamily described in
Corollary~\ref{cor:binary-coverage-thesis}.  An ordinary binary curve
without a rational point of order four requires a different marked model.
Likewise, operation counts in \(\M,\Sqr,\Dpar\) compare algebraic circuits;
they become implementation timings only after field representation,
reduction, memory traffic, and processor architecture have been fixed.

\part{Moduli, Division Theory, Isogenies, and Pairings}
\partoverview{The purpose of this part is to test whether the native
\(\mathcal C_d\) interface remains useful beyond scalar multiplication.
The chapters move from parameter spaces and point-count moments to division
polynomials and division-fiber means, then to model-preserving isogenies and
native Miller functions.  Auxiliary Weierstrass or Edwards equations are
restated locally and every final formula is returned to \(u,v,d\).}

\chapter{Parameters, Moduli, and Isomorphism Classes}
\label{ch:moduli}
\section{Definitions and the distinction between models and curves}

Let
\[
\mathcal D_q=
\begin{cases}
\F_q\setminus\{0,1/16\},&q\ \text{odd},\\
\F_q^\times,&q\ \text{even}.
\end{cases}
\]
Write
\[
\begin{aligned}
J_{\C}(q)
 &=\#\{\Cd\otimes\overline{\F}_q:d\in\mathcal D_q\}/
            \cong_{\overline{\F}_q},\\
I_{\C}(q)
 &=\#\{\Cd:d\in\mathcal D_q\}/\cong_{\F_q}.
\end{aligned}
\]
These counts forget the \((2,2)\)-embedding and the selected generator of
the rational four-torsion subgroup.  Arithmetic equivalence of marked models
is finer than elliptic-curve isomorphism.

\section{Odd characteristic}

The parameter change \(d\mapsto\rho=1-16d\) is a bijection
\(\mathcal D_q\to\F_q\setminus\{0,1\}\), and
\[
       j_d=\frac{(16d^2-16d+1)^3}{d^4(1-16d)}.
\]
Thus two parameters are geometrically isomorphic if and only if
\begin{equation}\label{eq:geometric-j-thesis}
(\rho^2+14\rho+1)^3\rho'(1-\rho')^4
=({\rho'}^2+14\rho'+1)^3\rho(1-\rho)^4.
\end{equation}
The map \(\rho\mapsto j\) has degree six; the shorter fibers at
\(j=0,1728\) account for the correction terms in the exact counts.

For a rational criterion, put \(c=4d\) and use
\begin{equation}\label{eq:E-c-thesis}
 E_c:\quad y^2=x^3+(1-2c)x^2+c^2x.
\end{equation}

\begin{theorem}[Rational isomorphism criterion]\label{thm:Fq-iso-thesis}
For \(d,d'\in\mathcal D_q\), put
\[
a=1-2c,\quad b=c^2,\qquad
a'=1-2c',\quad b'={c'}^2.
\]
Then \(\Cd\cong_{\F_q}\C_{d'}\) if and only if there are
\(\lambda\in\F_q^\times\), \(r\in\F_q\) such that
\begin{equation}\label{eq:Fq-iso-system-thesis}
\begin{aligned}
r(r^2+ar+b)&=0,\\
3r+a&=\lambda^2a',\\
3r^2+2ar+b&=\lambda^4b'.
\end{aligned}
\end{equation}
Over \(\overline{\F}_q\), the same system is equivalent to
\eqref{eq:geometric-j-thesis}.
\end{theorem}

\begin{proof}
We first make explicit why the Weierstrass equations in
\eqref{eq:E-c-thesis} test isomorphism of the original \(\mathcal C_d\)
models.  On the dense affine open set on which \(u\ne0\), put
\[
 U=\frac{u+1}{u},\qquad V=2(2v+1)U,
 \qquad x=cU,\qquad y=\frac{cV}{2}.
\]
The equation of \(\mathcal C_d\), with \(c=4d\), then gives
\[
 \frac{1}{4c}V^2=U^3+(c^{-1}-2)U^2+U
\]
and hence
\[
 y^2=x^3+(1-2c)x^2+c^2x.
\]
The inverse rational functions are obtained successively from
\(U=x/c\), \(u=(U-1)^{-1}\), and
\(2v+1=V/(2U)=y/x\).  Since both equations have smooth projective
completions, this birational map extends uniquely to an isomorphism of
elliptic curves taking the declared identity to the point at infinity.

We now determine all isomorphisms between \(E_{c'}\) and \(E_c\).  A
general admissible change of Weierstrass variables has the form
\[
 x=\lambda^2x'+r,\qquad
 y=\lambda^3y'+s\lambda^2x'+t,
 \qquad \lambda\ne0.
\]
Both equations have \(a_1=a_3=0\).  The transformation formula for
\(a_1\), or simply the coefficient of \(x'y'\) after substitution, gives
\(2s=0\); the coefficient of \(y'\) then gives \(2t=0\).  Because the
characteristic is odd, \(s=t=0\).  Expanding the right-hand side of
\(E_c\) gives
\[
\begin{aligned}
 &(\lambda^2x'+r)^3+a(\lambda^2x'+r)^2
       +b(\lambda^2x'+r)\\
 &\quad=\lambda^6{x'}^3
   +\lambda^4(3r+a){x'}^2
   +\lambda^2(3r^2+2ar+b)x'
   +r^3+ar^2+br.
\end{aligned}
\]
On the other hand,
\(\lambda^6{y'}^2=\lambda^6({x'}^3+a'{x'}^2+b'x')\).
Equality of the constant, quadratic, and linear coefficients is therefore
equivalent, respectively, to
\[
 r(r^2+ar+b)=0,\qquad
 3r+a=\lambda^2a',\qquad
 3r^2+2ar+b=\lambda^4b'.
\]
This proves necessity.  If the three equations hold, the displayed
substitution has inverse
\(x'=\lambda^{-2}(x-r)\), \(y'=\lambda^{-3}y\).  In the displayed
coefficient expansion, the constant term vanishes by
\(r(r^2+ar+b)=0\), while the \(x'^2\)- and \(x'\)-coefficients become
\(\lambda^6a'\) and \(\lambda^6b'\), respectively.  Hence the
substitution carries \(E_{c'}\) onto \(E_c\), proving sufficiency.

Finally extend scalars to \(\overline{\F}_q\).  The curves remain smooth,
and two elliptic curves over an algebraically closed field are isomorphic
if and only if their \(j\)-invariants agree.  Formula
\eqref{eq:j-thesis} turns this equality into
\eqref{eq:geometric-j-thesis}.  Applying the preceding coefficient
argument over \(\overline{\F}_q\) proves the last assertion, including the
exceptional values \(j=0\) and \(j=1728\).
\end{proof}

\begin{theorem}[Exact odd-characteristic class counts]
\label{thm:odd-class-counts-thesis}
For an odd prime power \(q\),
\begin{equation}\label{eq:J-odd-thesis}
J_{\C}(q)=
\begin{cases}
\left\lfloor\dfrac{5q+7}{12}\right\rfloor,&q\equiv1\pmod4,\\[2mm]
\left\lfloor\dfrac{3q-1}{8}\right\rfloor,&q\equiv3\pmod4,
\end{cases}
\end{equation}
and
\begin{equation}\label{eq:I-odd-thesis}
I_{\C}(q)=
\begin{cases}
(2q+1)/3,&q\equiv1\pmod{12},\\
(2q-1)/3,&q\equiv5\pmod{12},\\
2q/3,&q\equiv9\pmod{12},\\
(3q-5)/4,&q\equiv3\pmod4.
\end{cases}
\end{equation}
\end{theorem}

\begin{proof}
The proof has four stages: transport to the Edwards parameter \(\rho\);
count geometric \(j\)-fibers; count rational classes with nonsquare
\(\rho\); and apply Burnside's lemma to the square-\(\rho\) stratum.
Separating the last two strata is essential because their rational
two-torsion behaves differently.
Throughout the proof, \(q\) is an odd prime power and \(\chi\) is the
quadratic character of \(\F_q\), extended by \(\chi(0)=0\).  Fixed points
and roots are counted as distinct elements of \(\F_q\), not with algebraic
multiplicity.  In particular, a double root in characteristic three
contributes one element to a fixed-point set.

\smallskip
\noindent\emph{Stage 1: parameter transport.}
Put
\[
 \xi=(2v+1)^{-1},\qquad \eta=(2u+1)^{-1}.
\]
Multiplying the equation of \(\mathcal C_d\) by
\(16\xi^2\eta^2\) gives
\[
             \xi^2+\eta^2=1+\rho\xi^2\eta^2,
             \qquad \rho=1-16d.
\]
The inverse formulas are
\(u=(\eta^{-1}-1)/2\) and \(v=(\xi^{-1}-1)/2\); extension across the
finitely many omitted points identifies the smooth completions.  Thus the
parameter set is exactly \(\rho\in\F_q\setminus\{0,1\}\), with neither
duplication nor loss of a boundary point.

\smallskip
\noindent\emph{Stage 2: geometric classes.}
Substitution in
\eqref{eq:j-thesis} gives the rational function
\[
 j(\rho)=16\frac{(\rho^2+14\rho+1)^3}
                    {\rho(1-\rho)^4}.
\]
Its generic geometric fiber has six points.  The finite-field fiber
enumeration is most conveniently carried out after introducing the
associated Legendre parameter \(z\).  On
\(\F_q\setminus\{0,1\}\), the three transpositions may be represented by
\[
 z\longmapsto 1-z,\qquad z\longmapsto z^{-1},
 \qquad z\longmapsto \frac{z}{z-1}.
\]
Their fixed-point equations are, respectively,
\[
 2z=1,\qquad z^2=1,\qquad z(z-2)=0.
\]
After the excluded values \(0,1\) are removed, the fixed points are
\(1/2,-1,2\).  The two nonidentity three-cycles are represented by
\(z\mapsto(1-z)^{-1}\) and \(z\mapsto(z-1)/z\); both have fixed-point
equation
\[
                         z^2-z+1=0.
\]
These fixed-point incidences account for the shorter fibers above
\(j=1728\) and \(j=0\), respectively.  In characteristic three some of
the displayed values coalesce, but Burnside's lemma counts the fixed set
of each group element separately, so no incidence is lost.  Tracking
which of these auxiliary
parameters descend to an Edwards parameter in \(\F_q\) uses
\[
 \chi(-1)=(-1)^{(q-1)/2},\qquad
 \chi(2)=(-1)^{(q^2-1)/8},\qquad
 \#\{z:z^2-z+1=0\}=1+\chi(-3),
\]
where \(\chi(0)=0\).  In characteristic three the second polynomial has
one double root, so the same expression records the coalescence rather than
silently counting two roots.  We now invoke the exact finite-field
anharmonic-orbit enumeration proved in
\cite{FarashahiMoodyWu2012}; its hypotheses are exactly
\(\rho\in\F_q\setminus\{0,1\}\) and \(q\) odd, and its five cases are
\[
\begin{array}{c|c}
\text{condition on }q&\#j(\F_q\setminus\{0,1\})\\ \hline
q\equiv1\pmod{12}&(5q+7)/12\\
q\equiv5\pmod{12}&(5q-1)/12\\
q\equiv9\pmod{12}&(5q+3)/12\\
q\equiv3\pmod 8&(3q-1)/8\\
q\equiv7\pmod 8&(3q-5)/8.
\end{array}
\]
For \(q\equiv1\pmod4\), the first three rows are exactly
\(\lfloor(5q+7)/12\rfloor\); for \(q\equiv3\pmod4\), the last two are
\(\lfloor(3q-1)/8\rfloor\).  This proves
\eqref{eq:J-odd-thesis}.  Notice that the row \(q\equiv9\pmod{12}\)
is the characteristic-three case with even extension degree.

\smallskip
\noindent\emph{Stage 3: the nonsquare rational stratum.}
Let
\(c=(1-\rho)/4\).  In the Weierstrass equation \(E_c\),
\[
 a=\frac{1+\rho}{2},\qquad
 b=\frac{(1-\rho)^2}{16},\qquad a^2-4b=\rho.
\]
If \(\chi(\rho)=-1\), then \((0,0)\) is the unique nonzero rational
two-torsion point.  Every rational isomorphism must fix it, so the
translation parameter in Theorem~\ref{thm:Fq-iso-thesis} is \(r=0\).
The two remaining equations imply either \(\rho'=\rho\), or
\(\rho'=\rho^{-1}\) with \(\lambda^2=\rho\).  The latter alternative is
impossible for a nonsquare.  (When \(\rho=-1\), the equation \(a=0\)
immediately gives \(\rho'=-1\), so no division by \(1+\rho\) is being
used.)  Hence distinct nonsquare parameters are not \(\F_q\)-isomorphic,
and they contribute
\[
                         I_{\rm ns}(q)=\frac{q-1}{2}.
\]

\smallskip
\noindent\emph{Stage 4: the square rational stratum.}
Write \(\rho=\theta^2\) and
factor
\[
 x^3+(1-2c)x^2+c^2x
   =x(x+s^2)(x+t^2),\qquad
 s=\frac{1+\theta}{2},\quad t=\frac{1-\theta}{2}.
\]
After scaling the three rational roots, this is a Legendre equation.  When
\(-1\) is a square its parameter belongs to
\[
 Q=\{z\in\F_q\setminus\{0,1\}:\chi(z)=1\}.
\]
Here \(E_{L,z}\) denotes
\[
                   E_{L,z}:\quad Y^2=X(X-1)(X-z).
\]
More precisely, one may take \(z=(t/s)^2\) when \(-1\) is a square.  When
\(-1\) is a nonsquare, one obtains instead the Legendre parameter
\(1-z\), where \(z=(t/s)^2\in Q\).  In the latter case
\(E_{L,1-z}\) is the nonsquare quadratic twist of \(E_{L,z}\); twisting
every member gives a bijection on rational isomorphism classes, so it is
still enough to count the corresponding parameters in \(Q\).  Replacing
\(\theta\) by \(-\theta\) replaces \(z\) by \(z^{-1}\), and therefore
does not duplicate an orbit.  The six possible Legendre
parameters in one geometric class are
\[
 z,\quad z^{-1},\quad 1-z,\quad(1-z)^{-1},\quad
 \frac{z}{z-1},\quad\frac{z-1}{z}.
\]

We now perform the Burnside count on \(Q\), using the rational
Legendre-transition criterion from the same cited orbit lemma.  If
\(q\equiv3\pmod4\), only the reciprocal pair remains rationally
isomorphic inside this square-parameter family.  The involution
\(z\mapsto z^{-1}\) has no fixed point in \(Q\), because its only
admissible fixed point is \(-1\), a nonsquare.  Since
\(\#Q=(q-3)/2\), the number of classes is \((q-3)/4\).

Suppose now that \(q\equiv1\pmod4\), and set
\[
 Q_+=\{z\in Q:\chi(1-z)=1\},\qquad
 Q_-=\{z\in Q:\chi(1-z)=-1\}.
\]
The standard quadratic-character sum
\(\sum_z\chi(z(1-z))=-1\) gives
\[
                 \#Q_+=\frac{q-5}{4},\qquad
                 \#Q_-=\frac{q-1}{4}.
\]
Here is the calculation, including the boundary terms.  Since
\(q\equiv1\pmod4\), \(\chi(-1)=1\).  The substitution \(x=2z-1\),
together with
\[
 \sum_{x\in\F_q}\chi(x^2-1)=-1,
\]
gives \(\sum_z\chi(z(1-z))=-1\).  To verify the displayed auxiliary
sum, count pairs satisfying \(y^2=x^2-1\): the map
\((x,y)\mapsto(x-y,x+y)\) is a bijection onto pairs
\((r,s)\in(\F_q^\times)^2\) with \(rs=1\), so there are \(q-1\) pairs;
on the other hand their number is
\(q+\sum_x\chi(x^2-1)\).  On the set \(z\ne0,1\), one also has
\[
 \sum\chi(z)=-1,\qquad \sum\chi(1-z)=-1.
\]
Consequently
\[
\begin{aligned}
4\#Q_+
 &=\sum_{z\ne0,1}(1+\chi(z))(1+\chi(1-z))=q-5,\\
4\#Q_-
 &=\sum_{z\ne0,1}(1+\chi(z))(1-\chi(1-z))=q-1.
\end{aligned}
\]
The full anharmonic group \(S_3\) acts on \(Q_+\), whereas only the
reciprocal involution acts on \(Q_-\).  Let
\(\epsilon_2=1\) if \(2\) is a square and \(0\) otherwise, and let
\[
 \nu_3=
 \begin{cases}
 2,&q\equiv1\pmod{12},\\
 0,&q\equiv5\pmod{12},\\
 1,&q\equiv9\pmod{12}.
 \end{cases}
\]
Here \(\nu_3\) is the number, with distinct field elements rather than
algebraic multiplicity, of roots of \(z^2-z+1\) in \(Q_+\).  The three
transposition fixed points are \(1/2,-1,2\), as calculated in Stage~2.
For \(q\equiv1\pmod4\), each belongs to \(Q_+\) exactly when \(2\) is a
square; hence the three transpositions have total
\(3\epsilon_2\) fixed points on \(Q_+\).  Each root of
\(z^2-z+1\) is fixed by both three-cycles, so those cycles contribute
\(2\nu_3\).  If the characteristic is not three, such roots occur for
\(q\equiv1\pmod{12}\), and they are squares because their order divides
six while \(12\mid(q-1)\); moreover
\(z(1-z)=1\), so they lie in \(Q_+\).  For
\(q\equiv5\pmod{12}\) there is no root in \(\F_q\).  In characteristic
three with \(q\equiv9\pmod{12}\), there is the single root
\(z=1/2=-1\), which lies in \(Q_+\) because \(-1\) is a square.  This
proves the three stated values of \(\nu_3\).

On \(Q_-\), inversion can fix only \(-1\).  Because
\(\chi(-1)=1\), this element lies in \(Q_-\) precisely when
\(\chi(1-(-1))=\chi(2)=-1\); its fixed-point contribution is therefore
\(1-\epsilon_2\).  Burnside's lemma now
gives
\begin{equation}\label{eq:square-orbit-count-expanded}
 I_{\rm sq}(q)=
 \frac{\#Q_++3\epsilon_2+2\nu_3}{6}
 +\frac{\#Q_-+1-\epsilon_2}{2}.
\end{equation}
Substitution yields
\[
I_{\rm sq}(q)=
\begin{cases}
(q+5)/6,&q\equiv1\pmod{12},\\
(q+1)/6,&q\equiv5\pmod{12},\\
(q+3)/6,&q\equiv9\pmod{12},\\
(q-3)/4,&q\equiv3\pmod4.
\end{cases}
\]
The characteristic-three specialization is included: for even extension
degree the polynomial \(z^2-z+1\) has one double root as recorded by
\(\nu_3=1\), and for odd extension degree one is in the
\(q\equiv3\pmod4\) case.  Finally
\(I_{\mathcal C}(q)=I_{\rm ns}(q)+I_{\rm sq}(q)\), and adding
\((q-1)/2\) in the four cases gives exactly
\eqref{eq:I-odd-thesis}.
\end{proof}

If \(\chi\) is the quadratic character, the complete Edwards subfamily
\(\chi(\rho)=-1\) has
\[
        I_{\rm comp}(q)=\frac{q-1}{2},
\]
and distinct nonsquare parameters are never \(\F_q\)-isomorphic.  The
square-parameter subfamily has the four-case value of
\(I_{\rm sq}(q)\) computed immediately above from
\eqref{eq:square-orbit-count-expanded}.
Their sum is \eqref{eq:I-odd-thesis}.

\section{Characteristic two}

\begin{theorem}[Radicial parameter map]\label{thm:binary-classes-thesis}
For \(q=2^m\) and \(d,d'\in\F_q^\times\),
\[
\Cd\cong_{\overline{\F}_q}\C_{d'}
\iff d=d'
\iff \Cd\cong_{\F_q}\C_{d'}.
\]
Consequently
\[
           J_{\C}(q)=I_{\C}(q)=q-1.
\]
\end{theorem}

\begin{proof}
Assume first that
\(\mathcal C_d\cong_{\overline{\F}_q}\mathcal C_{d'}\).  Isomorphic
elliptic curves have equal \(j\)-invariants, and
\eqref{eq:binary-j-thesis} therefore gives
\[
                         d^{-4}={d'}^{-4}.
\]
Both parameters are nonzero, so multiplication by \(d^4{d'}^4\) gives
\(d^4={d'}^4\).  In characteristic two,
\[
                         (d-d')^4=d^4-{d'}^4=0.
\]
A field has no nonzero nilpotents, whence \(d=d'\).  This proves that even
geometric isomorphism forces equality of parameters.  Equality of
parameters, in turn, gives the identity isomorphism over \(\F_q\), and an
\(\F_q\)-isomorphism is certainly a geometric one.  Thus the three
conditions in the statement are equivalent.

There are exactly \(q-1\) admissible parameters.  Equivalently, the map
\(d\mapsto d^{-4}\) is a bijection on \(\F_q^\times\): inversion is a
bijection and \(x\mapsto x^4\) is the second iterate of Frobenius, hence an
automorphism of the finite field.  Consequently the \(q-1\) parameters give
\(q-1\) distinct geometric classes and, a fortiori, \(q-1\) distinct
rational classes.
\end{proof}

\begin{remark}
In odd characteristic, the parameter-to-moduli map is separable of generic
degree six.  In characteristic two, \(d\mapsto j=d^{-4}\) is radicial and
injective.  Thus a binary parameter is already an isomorphism-class
identifier for this family.
\end{remark}

\section{Average parameter multiplicity}

The exact average numbers of parameters per geometric and rational class are
\[
\overline m_{\rm geom}(q)=\frac{\#\mathcal D_q}{J_{\C}(q)},\qquad
\overline m_{\F_q}(q)=\frac{\#\mathcal D_q}{I_{\C}(q)}.
\]
Here \(\#\mathcal D_q=q-2\) for odd \(q\), and both averages are \(1\)
for even \(q\).  These orbit averages must not be confused with point-count
moments or averages over division fibers.

\chapter{Exact First and Second Moments of the Family}
\label{ch:first-moment}
\chaptermark{Exact Moments of the Family}
Chapter~\ref{ch:moduli} counted the geometric and rational isomorphism
fibers of the parameter map.  Here the parameters themselves remain the
sample space: we first express each point count directly from the factorized
equation, then compute the first moment, then the second moment, and only at
the end translate the formulas into class-weighted language.  Division-fiber
averages use a different sample space and are postponed to
Chapter~\ref{ch:mean-values}.

Let
\[
N_d=\#\Cd(\F_q),\qquad a_d=q+1-N_d,
\]
where \(\Cd\) denotes the smooth projective completion.

\section{Native point-count character sums}

The factorized equation gives a point-count expression before any change of
model.  This expression is also useful for checking individual parameters.

\begin{proposition}[Exact native character-sum formulas]
\label{prop:native-point-count-character-sum}
If \(q\) is odd, \(d\ne0,1/16\), and \(\chi\) is the quadratic character
extended by \(\chi(0)=0\), then
\begin{equation}\label{eq:native-point-count-odd}
 N_d=q+2+
 \sum_{u\in\F_q\setminus\{0,-1\}}
 \chi\!\left(1+\frac{4d}{u^2+u}\right).
\end{equation}
If \(q=2^m\), \(d\ne0\), and
\(\psi(z)=(-1)^{\operatorname{Tr}_{\F_q/\F_2}(z)}\), then
\begin{equation}\label{eq:native-point-count-binary}
 N_d=q+2+
 \sum_{u\in\F_q\setminus\{0,1\}}
 \psi\!\left(\frac{d}{u^2+u}\right).
\end{equation}
\end{proposition}

\begin{proof}
Fix \(u\) with \(u^2+u\ne0\), and put
\(c=d/(u^2+u)\).  In odd characteristic, the equation
\(v^2+v=c\) has discriminant \(1+4c\), and therefore has
\(1+\chi(1+4c)\) solutions.  There are \(q-2\) admissible values of
\(u\), and the smooth completion contributes the four rational boundary
points.  Summing gives \(4+(q-2)\) plus the character sum in
\eqref{eq:native-point-count-odd}.

In characteristic two, the Artin--Schreier equation
\(v^2+v=c\) has two solutions when \(\operatorname{Tr}(c)=0\) and no
solution when \(\operatorname{Tr}(c)=1\).  Its number of solutions is
therefore \(1+\psi(c)\).  Again there are \(q-2\) admissible values of
\(u\) and four boundary points, proving
\eqref{eq:native-point-count-binary}.
\end{proof}

\begin{lemma}[The quadratic correction]\label{lem:first-moment-character-sum}
For every odd \(q\),
\begin{equation}\label{eq:quadratic-correction-first-moment}
              \sum_{r\in\F_q}\chi(r^2-1)=-1.
\end{equation}
\end{lemma}

\begin{proof}
Count the affine solutions of \(y^2=r^2-1\) in two ways.  For fixed
\(r\), the number of \(y\) is \(1+\chi(r^2-1)\), so the total is
\(q+\sum_r\chi(r^2-1)\).  On the other hand,
\[
             (r-y)(r+y)=1.
\]
Because \(2\) is invertible, every \(t\in\F_q^\times\) gives exactly one
solution
\(r=(t+t^{-1})/2\), \(y=(t^{-1}-t)/2\), and every solution arises this
way.  The total is \(q-1\), which proves
\eqref{eq:quadratic-correction-first-moment}.
\end{proof}

\section{The parameter-weighted first moment}

\begin{theorem}[Parameter-weighted first moment]\label{thm:first-moment-thesis}
If \(q\) is odd and \(\chi\) is the quadratic character, then
\begin{align}
\sum_{d\in\mathcal D_q}N_d
 &=q^2-q-1+\chi(-1),\label{eq:sum-N-odd-thesis}\\
\frac1{q-2}\sum_{d\in\mathcal D_q}N_d
 &=q+1+\frac{1+\chi(-1)}{q-2},\\
\sum_{d\in\mathcal D_q}a_d&=-1-\chi(-1).
\end{align}
If \(q=2^m\), then
\begin{align}
\sum_{d\in\F_q^\times}N_d&=q^2,\label{eq:sum-N-even-thesis}\\
\frac1{q-1}\sum_{d\in\F_q^\times}N_d
 &=q+1+\frac1{q-1},\\
\sum_{d\in\F_q^\times}a_d&=-1.
\end{align}
\end{theorem}

\begin{proof}
Let \(f(z)=z^2+z\).  Every smooth completion contributes four rational
boundary points.

Suppose first that \(q\) is odd.  There are \(q-2\) values of \(u\) with
\(f(u)\ne0\), so \((q-2)^2\) affine pairs have nonzero product.  We must
remove the pairs belonging to the singular parameter \(d=1/16\).  With
\(r=2u+1\), \(s=2v+1\), these pairs satisfy
\[
       (r^2-1)(s^2-1)=1.
\]
For \(r\ne0,\pm1\), the number of \(s\) is
\(1+\chi(r^2-1)\), and \(r=0\) contributes one pair.  Hence
\[
\begin{aligned}
T_q
 &=1+\sum_{r\ne0,\pm1}\{1+\chi(r^2-1)\}\\
 &=1+(q-3)+\{-1-\chi(-1)\}\\
 &=q-3-\chi(-1),
\end{aligned}
\]
where Lemma~\ref{lem:first-moment-character-sum} was used and the terms at
\(r=\pm1\) vanish.  The affine total over smooth parameters is
\((q-2)^2-T_q\).  There are \(q-2\) smooth parameters, so adding their
four boundary points gives
\[
\begin{aligned}
 \sum_{d\in\mathcal D_q}N_d
  &=(q-2)^2-(q-3-\chi(-1))+4(q-2)\\
  &=q^2-q-1+\chi(-1),
\end{aligned}
\]
which is
\eqref{eq:sum-N-odd-thesis}.  Division by \(q-2\) and
\(a_d=q+1-N_d\) give the remaining odd formulas.

If \(q\) is even, \(f(z)=0\) still has exactly the two roots \(0,1\), and
every nonzero product is a smooth parameter.  Thus the affine total is
\((q-2)^2\), while the boundary total is \(4(q-1)\).  Their sum is
\(q^2\), proving \eqref{eq:sum-N-even-thesis} and its consequences.
Equivalently, summing \eqref{eq:native-point-count-binary} over
\(d\in\F_q^\times\) uses
\(\sum_{d\ne0}\psi(d/c)=-1\) for every \(c\ne0\), and gives
\[
 (q-1)(q+2)-(q-2)=q^2.
\]
\end{proof}

\begin{remark}
This proof uses the factorized equation directly.  It avoids a transformation
to Weierstrass form and exposes the small correction to the Hasse center as
a boundary-and-singular-fiber contribution.
\end{remark}

\section{Exact second moments of Frobenius traces}

The same native character sums determine the quadratic fluctuation about
the Hasse center.  This gives information that is invisible in the first
moment and, in particular, provides a sensitive check on the singular-fiber
correction in odd characteristic.

\begin{theorem}[Exact parameter-weighted second moment]
\label{thm:second-moment-thesis}
Let \(a_d=q+1-N_d\).  Then
\begin{equation}\label{eq:second-moment-thesis}
 \sum_{d\in\mathcal D_q}a_d^2
 =
 \begin{cases}
   q^2-2q-3=(q-3)(q+1),&q\text{ odd},\\
   q^2-q-1,&q\text{ even}.
 \end{cases}
\end{equation}
In the even case \(\mathcal D_q=\F_q^\times\), as in
Chapter~\ref{ch:moduli}.
\end{theorem}

\begin{proof}
The proof has two characteristic branches.  In the odd branch we first
sum over every parameter and only then remove the two singular parameters;
in the binary branch every nonzero parameter is smooth, so additive
orthogonality gives the answer directly.

\smallskip
\noindent\emph{Odd characteristic: the all-parameter sum.}
Put
\[
 D=\F_q\setminus\{0,-1\},\qquad h_u=u^2+u,
 \qquad
 S_d=\sum_{u\in D}\chi\!\left(1+\frac{4d}{h_u}\right).
\]
Proposition~\ref{prop:native-point-count-character-sum} gives
\(a_d=-1-S_d\) for every smooth \(d\).  For the moment, allow \(d\) to
range over all of \(\F_q\), and write \(c_u=4/h_u\).  Expanding the square
and interchanging summations gives
\begin{equation}\label{eq:all-d-S-square-odd}
 \sum_{d\in\F_q}S_d^2
 =\sum_{u,w\in D}\sum_{d\in\F_q}
       \chi\bigl((1+c_ud)(1+c_wd)\bigr).
\end{equation}
If \(c_u=c_w\), the inner sum is \(q-1\).  If \(c_u\ne c_w\), the
quadratic polynomial has two distinct roots and its character sum is
\(-\chi(c_uc_w)=-\chi(h_uh_w)\).  To verify the identity, note that
\[
 \sum_{d\in\F_q}\chi((d-a)(d-b))=-1\qquad(a\ne b):
\]
indeed, the affine change
\(r=(2d-a-b)/(b-a)\) changes the summand into
\(\chi(r^2-1)\), and
Lemma~\ref{lem:first-moment-character-sum} applies.  Multiplication by
\(\chi(c_uc_w)\) proves the inner-sum assertion.  Next,
\[
 h_u=h_w
 \quad\Longleftrightarrow\quad
 (u-w)(u+w+1)=0.
\]
Thus \(w=u\) or \(w=-u-1\).  The involution \(u\mapsto-u-1\) has the
single fixed point \(-1/2\) in \(D\), so the number of ordered equal pairs
is
\[
                 E=2|D|-1=2q-5.
\]
Moreover
\[
 H:=\sum_{u\in D}\chi(h_u)=-1.
\]
Indeed, completing the square changes this sum into
\(\sum_{r\in\F_q}\chi(r^2-1)\), with \(r=2u+1\); the omitted roots
contribute zero, and
Lemma~\ref{lem:first-moment-character-sum} applies.  Hence the sum of
\(\chi(h_uh_w)\) over the unequal ordered pairs is \(H^2-E=1-E\).
Equation~\eqref{eq:all-d-S-square-odd} therefore becomes
\begin{equation}\label{eq:all-d-S-square-odd-value}
       \sum_{d\in\F_q}S_d^2
       =E(q-1)-(1-E)=Eq-1=2q^2-5q-1.
\end{equation}

\smallskip
\noindent\emph{Odd characteristic: removal of singular parameters.}
Two parameters in the preceding sum are not smooth.  At \(d=0\), one has
\(S_0=q-2\).  At \(d_*=1/16\),
\[
 1+\frac{4d_*}{h_u}
   =\frac{(2u+1)^2}{4h_u}.
\]
Every term except the one at \(u=-1/2\) is therefore \(\chi(h_u)\),
whereas that exceptional term is zero.  Since
\(\chi(h_{-1/2})=\chi(-1)\),
\[
                   S_{d_*}=-1-\chi(-1).
\]
Subtracting these two contributions from
\eqref{eq:all-d-S-square-odd-value} yields
\begin{equation}\label{eq:smooth-S-square-odd}
 \sum_{d\in\mathcal D_q}S_d^2
       =q^2-q-7-2\chi(-1).
\end{equation}
The first-moment formula and \(a_d=-1-S_d\) also give
\[
 \sum_{d\in\mathcal D_q}S_d=3-q+\chi(-1).
\]
Using \(|\mathcal D_q|=q-2\), we finally obtain
\[
 \sum_{d\in\mathcal D_q}a_d^2
 =(q-2)+2\sum_dS_d+\sum_dS_d^2
 =q^2-2q-3.
\]

\smallskip
\noindent\emph{Characteristic two.}
Now let \(q\) be even and put
\[
 D=\F_q\setminus\{0,1\},\qquad h_u=u^2+u,
 \qquad S_d=\sum_{u\in D}\psi(d/h_u).
\]
Again \(a_d=-1-S_d\).  Additive-character orthogonality shows that
\[
 \sum_{d\in\F_q}S_d^2=qE,
\]
where \(E\) counts ordered pairs \((u,w)\in D^2\) with \(h_u=h_w\).
Indeed, the inner sum over \(d\) is
\[
 \sum_{d\in\F_q}\psi\!\left(
 d\left(\frac1{h_u}+\frac1{h_w}\right)\right),
\]
which is \(q\) when \(h_u=h_w\) and zero otherwise.
In characteristic two,
\[
 h_u=h_w\quad\Longleftrightarrow\quad w=u\ \text{or}\ w=u+1,
\]
and the two choices are distinct; hence \(E=2(q-2)\).  Removing \(d=0\),
where \(S_0=q-2\), gives
\[
 \sum_{d\ne0}S_d^2=2q(q-2)-(q-2)^2=q^2-4.
\]
Also, additive-character orthogonality gives
\(\sum_{d\ne0}S_d=-(q-2)=2-q\).  Since there are \(q-1\) smooth
parameters,
\[
 \sum_{d\ne0}a_d^2
 =(q-1)+2(2-q)+(q^2-4)=q^2-q-1,
\]
which completes the proof.
\end{proof}

\begin{corollary}[Average square and extension-field counts]
\label{cor:second-moment-and-zeta}
The parameter averages of \(a_d^2\) are
\[
 \frac1{|\mathcal D_q|}\sum_{d\in\mathcal D_q}a_d^2
 =\begin{cases}
 q-\dfrac{3}{q-2},&q\text{ odd},\\[4pt]
 q-\dfrac1{q-1},&q\text{ even}.
 \end{cases}
\]
For each smooth parameter, set \(A_{d,0}=2\), \(A_{d,1}=a_d\), and
\begin{equation}\label{eq:Frobenius-trace-recurrence-Cd}
 A_{d,m}=a_dA_{d,m-1}-qA_{d,m-2}\qquad(m\ge2).
\end{equation}
Then
\[
 \#\mathcal C_d(\F_{q^m})=q^m+1-A_{d,m},\qquad
 Z(\mathcal C_d/\F_q,T)
  =\frac{1-a_dT+qT^2}{(1-T)(1-qT)}.
\]
\end{corollary}

\begin{proof}
The two averages follow by dividing
\eqref{eq:second-moment-thesis} by \(q-2\) and \(q-1\), respectively:
\[
 \frac{q^2-2q-3}{q-2}
   =q-\frac3{q-2},\qquad
 \frac{q^2-q-1}{q-1}
   =q-\frac1{q-1}.
\]
If \(\pi_d,\bar\pi_d\) are the roots of
\(X^2-a_dX+q\), then \(A_{d,m}=\pi_d^m+\bar\pi_d^m\); this proves the
recurrence, the extension-field point count, and the zeta-function formula.
\end{proof}

\section{Weights, class averages, and small checks}

Let \(\mathcal I_q\) be the set of \(\F_q\)-isomorphism classes represented
by the family, and let \(m(C)\) be the number of smooth parameters that
represent \(C\).  Because rationally isomorphic curves have the same point
count, Theorem~\ref{thm:first-moment-thesis} is equivalently the exact
weighted class identity
\begin{equation}\label{eq:first-moment-class-weight}
 \sum_{C\in\mathcal I_q}m(C)N(C)
 =\begin{cases}
 q^2-q-1+\chi(-1),&q\text{ odd},\\
 q^2,&q\text{ even}.
 \end{cases}
\end{equation}
Thus the theorem is parameter-weighted, not the unweighted average over
isomorphism classes.  The multiplicities computed in
Chapter~\ref{ch:moduli} are indispensable when the two averages are
compared.

\begin{example}[The field \(\F_5\)]
Here the singular value is \(d=1\), so the smooth parameters are
\(2,3,4\).  Directly from the factor table for \(z^2+z\),
\[
       N_2=8,\qquad N_3=4,\qquad N_4=8.
\]
Their sum is \(20\), agreeing with
\(5^2-5-1+\chi(-1)=20\); moreover
\(a_2+a_3+a_4=-2=-1-\chi(-1)\), and
\[
 a_2^2+a_3^2+a_4^2=12=5^2-2\cdot5-3.
\]
\end{example}

\begin{example}[The field \(\F_4\)]
Write \(\F_4=\{0,1,\alpha,\alpha+1\}\) with
\(\alpha^2+\alpha+1=0\).  The two nonzero inputs of
\(z^2+z\) both have value \(1\).  Hence
\[
       N_1=8,\qquad N_\alpha=N_{\alpha+1}=4,
\]
and \(8+4+4=16=q^2\), as predicted.
The corresponding traces are \(-3,1,1\), whose squared sum is
\(11=4^2-4-1\).
\end{example}

\chapter{Division Polynomials}
\label{ch:division}
This chapter presents two complementary constructions and keeps their roles
separate.  Sections~\ref{sec:division-generalized-recursion}--
\ref{sec:division-edwards-recurrence} begin with established division
functions on generalized Weierstrass or Edwards models and pull them back to
the declared \(\mathcal C_d\)-coordinates; this supplies normalized formulas
and convenient low-degree expansions.  Section~\ref{sec:native-scalar-division-polynomials}
then starts instead from the native doubling and differential-addition maps
and constructs the same division fibers by an addition-chain recursion.  The
second construction is therefore not a repetition of the first: it explains
how division evaluation, torsion testing, and scalar multiplication share one
native circuit.

\section{The generalized Weierstrass recursion}
\label{sec:division-generalized-recursion}

For \eqref{eq:general-W}, define
\[
\begin{aligned}
b_2&=a_1^2+4a_2,&b_4&=a_1a_3+2a_4,\\
b_6&=a_3^2+4a_6,&
b_8&=a_1^2a_6+4a_2a_6-a_1a_3a_4
      +a_2a_3^2-a_4^2.
\end{aligned}
\]
The division polynomials are
\begin{align}
\psi_0&=0,\quad \psi_1=1,\quad
\psi_2=2Y+a_1X+a_3,\label{eq:psi12-thesis}\\
\psi_3&=3X^4+b_2X^3+3b_4X^2+3b_6X+b_8,\label{eq:psi3-thesis}\\
\psi_4&=\psi_2\bigl(2X^6+b_2X^5+5b_4X^4+10b_6X^3
 +10b_8X^2\notag\\
&\hspace{24mm}+(b_2b_8-b_4b_6)X+b_4b_8-b_6^2\bigr),\\
\psi_{2m+1}
 &=\psi_{m+2}\psi_m^3-\psi_{m-1}\psi_{m+1}^3,\label{eq:psi-odd-thesis}\\
\psi_{2m}
 &=\frac{\psi_m}{\psi_2}
(\psi_{m+2}\psi_{m-1}^2-\psi_{m-2}\psi_{m+1}^2).
\label{eq:psi-even-thesis}
\end{align}
Set
\[
       \phi_n=X\psi_n^2-\psi_{n-1}\psi_{n+1}.
\]
If \((n,\charac k)=1\), then
\[
X([n]P)=\frac{\phi_n(P)}{\psi_n(P)^2},\qquad
[n]P=O\iff\psi_n(P)=0\quad(P\ne O).
\]
The generalized value \(\psi_2=2Y+a_1X+a_3\) is indispensable in
characteristic two.

\section{Odd-characteristic pullback}

Scale the Montgomery model by
\[
       X=\beta U,\qquad Y=\beta^2V
\]
to obtain
\begin{equation}\label{eq:odd-W-division-thesis}
 W_d^-:\quad Y^2=X^3+\beta AX^2+\beta^2X.
\end{equation}
The native division function is
\[
\Psi_{n,d}^{(-)}(u,v)=
\psi_n\left(
\beta\frac{u+1}{u},
2\beta^2(2v+1)\frac{u+1}{u}\right).
\]
Clearing chart denominators and common curve factors gives a genuine native
division polynomial whose zero divisor is the closure of the nonzero
\(n\)-torsion in that chart.

For compact formulas put
\[
\begin{aligned}
f(U)&=U^3+AU^2+U,\\
P_3(U)&=3U^4+4AU^3+6U^2-1,\\
F_4(U)&=U^6+2AU^5+5U^4-5U^2-2AU-1.
\end{aligned}
\]
Then
\[
\beta^{-4}\psi_3(\beta U)=P_3,\qquad
\frac{\psi_4(\beta U,\beta^2V)}{2\beta^8V}=2F_4.
\]
The fifth polynomial is
\begin{align*}
P_5(U)={}&5U^{12}+20AU^{11}+(16A^2+62)U^{10}
+80AU^9-105U^8-360AU^7\\
&-(240A^2+300)U^6-(64A^3+368A)U^5\\
&-(160A^2+125)U^4-140AU^3-50U^2+1,
\end{align*}
and the shorter check is
\begin{equation}\label{eq:P5-factor-thesis}
        P_5=32f^2F_4-P_3^3.
\end{equation}
On \(Z^2=f(U)\), the normalized polynomials through seven are
\begin{align}
\vartheta_2&=2Z,&
\vartheta_3&=P_3,&
\vartheta_4&=4ZF_4,&
\vartheta_5&=P_5,\notag\\
\vartheta_6&=2ZP_3(P_5-4F_4^2),\label{eq:theta6-thesis}\\
\vartheta_7&=P_5P_3^3-128f^2F_4^3.\label{eq:theta7-thesis}
\end{align}
These factorizations follow immediately from
\eqref{eq:psi-odd-thesis}--\eqref{eq:psi-even-thesis} and are preferable
to full expansion for torsion and small-kernel computations.

\section{Characteristic-two pullback}

Here the auxiliary equation is restated in full:
\[
       W_d^+:\qquad Y^2+XY=X^3+d^2X.
\]
For this curve,
\[
a_1=1,\quad a_2=a_3=a_6=0,\quad a_4=d^2,
\]
so
\begin{equation}\label{eq:binary-small-psi-thesis}
\psi_2=X,\qquad
\psi_3=X^4+X^3+d^4,\qquad
\psi_4=X^6+d^4X^2.
\end{equation}
For odd \(n\), let \(m_n=(n^2-1)/2\) and define
\[
        D_{n,d}^{(+)}(U)=d^{-m_n}\psi_n(dU).
\]
Then
\begin{align}
D_{3,d}^{(+)}(U)&=U^4+d^{-1}U^3+1,\label{eq:D3-binary-thesis}\\
D_{5,d}^{(+)}(U)&=U^{12}+d^{-1}U^{11}+d^{-2}U^{10}
+U^8+d^{-2}U^6\notag\\
&\quad+d^{-3}U^5+U^4+d^{-1}U^3+1,\label{eq:D5-binary-thesis}\\
D_{7,d}^{(+)}(U)&=U^{24}+d^{-2}U^{22}+d^{-3}U^{21}
+(1+d^{-4})U^{16}\notag\\
&\quad+d^{-5}U^{15}+d^{-6}U^{14}
+(1+d^{-4})U^8+d^{-5}U^7\notag\\
&\quad+d^{-2}U^6+d^{-3}U^5+1.
\label{eq:D7-binary-thesis}
\end{align}
For example, \(\psi_5=\psi_4\psi_2^3+\psi_3^3\) in characteristic two,
and substitution \(X=dU\) yields \eqref{eq:D5-binary-thesis}.
The sparsity through degrees \(4,12,24\) is a practical advantage of the
native coordinate.

\section{An Edwards recurrence in native variables}
\label{sec:division-edwards-recurrence}

Assume in this section that \(\charac k\ne2,3\).  On the Edwards curve
\[
 E_\rho:\qquad \xi^2+\eta^2=1+\rho\xi^2\eta^2,
 \qquad \rho\ne0,1,
\]
one may use the division rational functions of
\cite[Theorem~5.1 and Corollary~5.2]{HittMcGuireMoloney2008}.  We record
the normalization because it is needed to compare their zero divisors with
the native pullback.  Under the isomorphism from \(\mathcal C_d\),
\[
\xi=(2v+1)^{-1},\qquad\eta=(2u+1)^{-1}.
\]
With \(e_0=0,e_1=1\), define
\begin{align*}
e_2&=\frac{(1-\rho)(\eta+1)}{2\xi(1-\eta)},\\
e_3&=\frac{(1-\rho)^3(-\rho\eta^4-2\rho\eta^3+2\eta+1)}
           {(2(1-\eta))^4},\\
e_4&=\frac{2(1-\rho)^6\eta(1+\eta)(1-\rho\eta^4)}
           {\xi(2(1-\eta))^7},
\end{align*}
For \(m\ge2\) in the odd recurrence and \(m\ge3\) in the even recurrence,
define
\begin{align*}
e_{2m+1}&=e_{m+2}e_m^3-e_{m-1}e_{m+1}^3,\\
e_{2m}&=\frac{e_m}{e_2}
 \bigl(e_{m+2}e_{m-1}^2-e_{m-2}e_{m+1}^2\bigr).
\end{align*}
These are identities in the function field: the displayed quotient is the
regular division function supplied by the recurrence, rather than a
pointwise instruction to divide by a possibly zero value of \(e_2\).  Put
\[
\Phi_n=\frac{1+\eta}{1-\eta}e_n^2
 -\frac4{1-\rho}e_{n-1}e_{n+1},\qquad
\Omega_n=\frac{2e_{2n}}{(1-\rho)e_n}.
\]
Then
\[
[n](\xi,\eta)=
\left(\frac{\Phi_ne_n}{\Omega_n},
\frac{\Phi_n-e_n^2}{\Phi_n+e_n^2}\right).
\]
Indeed, the cited Edwards division functions have parameters
\(a=1\) and \(d=\rho\).  Their first four terms are therefore exactly the
four displayed functions above, and their multiplication formula is the
last display.  Composing that formula with
\(u=(\eta^{-1}-1)/2\), \(v=(\xi^{-1}-1)/2\) gives the native
multiplication map on the common dense chart.  Both sides are rational maps
between smooth projective curves, so equality on that chart proves equality
globally.  Finally, these \(e_n\) are the pullbacks, up to the displayed
nonzero normalization factors, of the generalized Weierstrass division
sections.  Pullback preserves their zero divisors, including multiplicity.
Consequently denominator clearing gives a second native division system
whose nonzero \(n\)-torsion divisor is the same divisor as in the
Weierstrass pullback.  Characteristic three is treated separately in
Chapter~\ref{ch:char3}; no short-Weierstrass division recurrence is invoked
there.

\section[Native division polynomials]{Division polynomials from native scalar multiplication}
\label{sec:native-scalar-division-polynomials}

The preceding constructions begin with a classical division-polynomial
recursion and pull it back to \(\mathcal C_d\).  There is a second,
model-internal construction: iterate the native doubling and differential
addition laws first, and read the division polynomial from the coordinate
that vanishes at the identity.  This construction is useful for three
reasons.  It proves that division theory can be developed from
\(\mathcal C_d\)-arithmetic alone, it produces an addition-chain circuit
without expanding large polynomials, and it gives an independent check on
the Weierstrass and Edwards pullbacks above.

\subsection{The odd-characteristic native multiplication pair}

Assume \(\charac k\ne2\), put
\[
       \alpha_{24}=\frac1{16d},
\]
and write a native Kummer point as
\[
       \mathbf P=(X:Z)=\kappa(P)=(u(P)+1:u(P)).
\]
For a homogeneous pair \(\mathbf P=(X:Z)\), define the native doubling
operator
\begin{equation}\label{eq:native-division-doubling-operator}
\begin{aligned}
 E_{\mathbf P}&=(X+Z)^2-(X-Z)^2,\\
 \mathcal D_d(\mathbf P)
   &=
 \bigl((X+Z)^2(X-Z)^2:
       E_{\mathbf P}\bigl((X-Z)^2+\alpha_{24}E_{\mathbf P}\bigr)\bigr).
\end{aligned}
\end{equation}
For three oriented Kummer points
\[
 \mathbf P=(X_1:Z_1),\qquad
 \mathbf Q=(X_2:Z_2),\qquad
 \boldsymbol\Delta=(X_\Delta:Z_\Delta),
\]
put
\begin{equation}\label{eq:native-division-addition-operator}
\begin{aligned}
 C&=(X_1+Z_1)(X_2-Z_2),\\
 D&=(X_2+Z_2)(X_1-Z_1),\\
 \mathcal A(\mathbf P,\mathbf Q;\boldsymbol\Delta)
   &=
 \bigl(Z_\Delta(C+D)^2:X_\Delta(C-D)^2\bigr).
\end{aligned}
\end{equation}
Equations~\eqref{eq:native-division-doubling-operator} and
\eqref{eq:native-division-addition-operator} are precisely the native
\(x\)DBL and \(x\)ADD derived in
Theorems~\ref{thm:xdbl-thesis} and~\ref{thm:xadd-thesis}.

If \((F:G)\) is a pair of homogeneous polynomials, write
\(\operatorname{prim}(F:G)\) for the pair obtained by dividing by its
greatest common homogeneous factor and then multiplying both entries by one
nonzero scalar to fix a chosen leading coefficient.  Define recursively
\begin{equation}\label{eq:native-Kummer-multiplication-recursion}
\begin{aligned}
 \mathbf K_0&=(1:0),&
 \mathbf K_1&=(X:Z),\\
 \mathbf K_{2m}
   &=\operatorname{prim}\bigl(\mathcal D_d(\mathbf K_m)\bigr),
   &&m\ge1,\\
 \mathbf K_{2m+1}
   &=\operatorname{prim}\bigl(
       \mathcal A(\mathbf K_{m+1},\mathbf K_m;\mathbf K_1)\bigr),
   &&m\ge1.
\end{aligned}
\end{equation}
The recursion is well founded because \(m,m+1<2m+1\).  Write
\[
                 \mathbf K_n=(\mathsf X_n:\mathsf Z_n).
\]

\begin{theorem}[Native scalar-multiplication division pair]
\label{thm:native-scalar-division-pair}
For every \(n\ge0\), the primitive pair
\((\mathsf X_n:\mathsf Z_n)\) in
\eqref{eq:native-Kummer-multiplication-recursion} represents the
multiplication-by-\(n\) map on the native Kummer line:
\begin{equation}\label{eq:native-Kn-represents-nP}
 \boxed{\qquad
   \kappa([n]P)
   =\bigl(\mathsf X_n(\kappa(P)):
           \mathsf Z_n(\kappa(P))\bigr).
 \qquad}
\end{equation}
The two entries are homogeneous of common degree \(n^2\).  In particular,
\[
       \Delta^{\mathrm{nat}}_{n,d}(u)
       =\mathsf Z_n(u+1,u)
\]
is a native Kummer division polynomial: for every finite affine point
\(P=(u,v)\) at which the curve is smooth,
\begin{equation}\label{eq:native-division-zero-criterion}
       [n]P=O
       \quad\Longleftrightarrow\quad
       \Delta^{\mathrm{nat}}_{n,d}(u(P))=0.
\end{equation}
Boundary zeros are interpreted on the smooth \((2,2)\)-completion.
\end{theorem}

\begin{proof}
Let \(\pi=\kappa:\overline{\mathcal C}_d\to\PP^1\) be the native Kummer
quotient.  The native formulas prove the two commutative identities
\[
 \mathcal D_d(\pi(P))=\pi(2P),\qquad
 \mathcal A(\pi(P),\pi(Q);\pi(P-Q))=\pi(P+Q)
\]
as identities of rational maps.  They were derived from the
\(\mathcal C_d\) group law, so no auxiliary curve is used here.

We prove \eqref{eq:native-Kn-represents-nP} by induction on \(n\).
For \(n=0\), the initial pair represents \(\pi(O)\), and for \(n=1\) it
represents \(\pi(P)\) by definition.  If the assertion is known for \(m\),
the first
identity gives
\[
       \mathcal D_d(\mathbf K_m(\pi(P)))=\pi([2m]P).
\]
If it is known for \(m\) and \(m+1\), then
\[
       [m+1]P-[m]P=P,
\]
so the second identity, with known difference \(\mathbf K_1=\pi(P)\),
gives
\[
 \mathcal A(\mathbf K_{m+1}(\pi(P)),
            \mathbf K_m(\pi(P));\mathbf K_1(\pi(P)))
       =\pi([2m+1]P).
\]
Removing a common polynomial factor does not alter a projective rational
map.  This proves the induction on the dense generic input.  The primitive
pair has no common projective zero and therefore extends uniquely to the
Kummer morphism induced by \([n]\).

The descended morphism \([n]_{\mathrm K}:\PP^1\to\PP^1\) satisfies
\[
                \pi\circ[n]=[n]_{\mathrm K}\circ\pi.
\]
Since \(\deg(\pi)=2\) and \(\deg([n])=n^2\), multiplicativity of degrees
gives
\[
       2n^2=\deg(\pi\circ[n])
            =\deg([n]_{\mathrm K})\deg(\pi)
            =2\deg([n]_{\mathrm K}).
\]
Thus \(\deg([n]_{\mathrm K})=n^2\).  A primitive homogeneous pair
representing a morphism of \(\PP^1\) has common degree equal to the degree
of that morphism, proving the degree assertion.

Finally,
\[
       \kappa(O)=(1:0).
\]
Therefore the second coordinate of
\(\kappa([n]P)\) vanishes exactly when \([n]P=O\).  Substituting the
native input \((X:Z)=(u+1:u)\) into \(\mathsf Z_n\) proves
\eqref{eq:native-division-zero-criterion}.  When
\((n,\charac k)=1\), this also has the divisor-theoretic interpretation
\begin{equation}\label{eq:native-division-zero-divisor}
 \operatorname{div}_0\!
 \bigl(\mathsf Z_n(\kappa(P))\bigr)
       =[n]^*(2(O))
       =2\sum_{Q\in\overline{\mathcal C}_d[n]}(Q),
\end{equation}
because \(O\) is a branch point of the double cover \(\kappa\).
Thus the recurrence has exactly the nonzero \(n\)-torsion zeros, together
with the projective identity and the multiplicities imposed by the Kummer
quotient.
\end{proof}

\begin{example}[The first native multiplication pairs]
\label{ex:first-native-multiplication-pairs}
For \(n=2\), the recurrence gives
\[
\begin{aligned}
 \mathsf X_2&=(X+Z)^2(X-Z)^2,\\
 \mathsf Z_2&=E\bigl((X-Z)^2+\alpha_{24}E\bigr),
 \qquad E=(X+Z)^2-(X-Z)^2.
\end{aligned}
\]
At the native input \((X:Z)=(u+1:u)\), one has
\[
       X-Z=1,\qquad E=4u(u+1),
\]
and hence, up to the nonzero scalar \(d^{-1}\),
\begin{equation}\label{eq:native-Delta2-explicit}
 \boxed{\qquad
       \Delta^{\mathrm{nat}}_{2,d}(u)
       \doteq u(u+1)(u^2+u+4d).
 \qquad}
\end{equation}
The roots \(u=0,-1\) are the marked boundary points \(O,T\); the remaining
quadratic gives the other geometric two-torsion Kummer values.

For \(n=3\), direct execution of
\[
        \mathbf K_3
        =\operatorname{prim}
          \mathcal A(\mathbf K_2,\mathbf K_1;\mathbf K_1)
\]
gives
\begin{equation}\label{eq:native-K3-denominator}
       \mathsf Z_3(X,Z)
       \doteq Z\bigl(P_3^{\,h}(X,Z)\bigr)^2,
\end{equation}
where
\[
 P_3^{\,h}(X,Z)
   =3X^4+4AX^3Z+6X^2Z^2-Z^4,
 \qquad A=\frac1{4d}-2.
\]
Thus
\[
 \Delta^{\mathrm{nat}}_{3,d}(u)
   \doteq
 u\bigl(P_3^{\,h}(u+1,u)\bigr)^2.
\]
The factor \(u\) records the projective identity, while the quartic
\(P_3^{\,h}(u+1,u)\) gives the nonzero three-torsion Kummer values.
The next two pairs are obtained without expansion from
\[
 \mathbf K_4=\operatorname{prim}\mathcal D_d(\mathbf K_2),
 \qquad
 \mathbf K_5=\operatorname{prim}
 \mathcal A(\mathbf K_3,\mathbf K_2;\mathbf K_1);
\]
their primitive degrees are \(16\) and \(25\), respectively.  These
factorized circuits are generally preferable to fully expanded
polynomials.
\end{example}

\subsection{The characteristic-two native recursion}

In characteristic two, retain the same Kummer coordinate
\[
             \kappa(P)=(X_0:X_1)=(u+1:u)
\]
and define
\begin{equation}\label{eq:binary-native-division-double-operator}
 \mathcal D_d^+(X_0:X_1)
   =\bigl(X_0^4+X_1^4:d^{-1}X_0^2X_1^2\bigr).
\end{equation}
For
\(\mathbf P=(X_0:X_1)\), \(\mathbf Q=(Y_0:Y_1)\), and
\(\boldsymbol\Delta=(t_0:t_1)\), put
\begin{equation}\label{eq:binary-native-division-add-operator}
\begin{aligned}
 A&=X_0Y_0+X_1Y_1,\qquad
 B=X_0Y_1+X_1Y_0,\\
 \mathcal A^+(\mathbf P,\mathbf Q;\boldsymbol\Delta)
   &=\bigl(t_1A^2:t_0B^2\bigr).
\end{aligned}
\end{equation}
These are the two native Artin--Schreier Kummer operations of
Theorem~\ref{thm:binary-native-complete-diff-atlas}.  Define
\(\mathbf K_n^+=(\mathsf X_n^+:\mathsf Z_n^+)\) by the same recursion
\eqref{eq:native-Kummer-multiplication-recursion}, replacing
\(\mathcal D_d,\mathcal A\) by
\(\mathcal D_d^+,\mathcal A^+\).

\begin{corollary}[Binary native division pair]
\label{cor:binary-native-division-pair}
The pair \(\mathbf K_n^+\) represents
\(\kappa([n]P)\) and, after primitive cancellation, has degree \(n^2\).
For every \(n\ge1\),
\[
       \Delta_{n,d}^{\mathrm{nat},+}(u)
       =\mathsf Z_n^+(u+1,u)
\]
vanishes exactly at the native \(n\)-division fiber, with the
scheme-theoretic multiplicities supplied by the Kummer quotient.  In
particular,
\[
\begin{aligned}
 \Delta_{2,d}^{\mathrm{nat},+}(u)
   &\doteq d^{-1}u^2(u+1)^2,\\
 \Delta_{3,d}^{\mathrm{nat},+}(u)
   &\doteq
   u\bigl((u+1)^4+d^{-1}(u+1)^3u+u^4\bigr)^2.
\end{aligned}
\]
\end{corollary}

\begin{proof}
The characteristic-two maps satisfy
\[
 \mathcal D_d^+(\kappa(P))=\kappa(2P),\qquad
 \mathcal A^+(\kappa(P),\kappa(Q);\kappa(P-Q))
   =\kappa(P+Q)
\]
as rational maps.  The assertion is true for \(n=0,1\) by the initial
conditions.  If \(\mathbf K_m^+\) represents \(\kappa([m]P)\), the
first identity shows that \(\mathbf K_{2m}^+\) represents
\(\kappa([2m]P)\).  If \(\mathbf K_m^+\) and
\(\mathbf K_{m+1}^+\) have the asserted meanings, then
\[
                 [m+1]P-[m]P=P,
\]
and the second identity gives
\[
 \mathcal A^+(\mathbf K_{m+1}^+,\mathbf K_m^+;\mathbf K_1^+)
       =\kappa([2m+1]P).
\]
This proves the recursion by induction.  Removing a common homogeneous
factor does not change the represented rational map; the resulting
primitive pair extends uniquely to the descended morphism
\([n]_{\mathrm K}^+:\PP^1\to\PP^1\).

The commutative relation
\[
           \kappa\circ[n]=[n]_{\mathrm K}^+\circ\kappa
\]
and \(\deg\kappa=2\) give
\[
 2n^2=\deg(\kappa\circ[n])
     =2\deg([n]_{\mathrm K}^+),
\]
so \(\deg([n]_{\mathrm K}^+)=n^2\).  Moreover,
\(\kappa^*(X_1=0)=2(O)\).  Therefore the zero divisor of the second
coordinate is
\[
 ([n]_{\mathrm K}^+\circ\kappa)^*(X_1=0)
      =[n]^*(2(O)).
\]
This divisor identity includes the scheme-theoretic multiplicities when
\([n]\) has an inseparable part.

It remains to check the two displayed low-degree polynomials.
Substitution of \((X_0:X_1)=(u+1:u)\) in
\eqref{eq:binary-native-division-double-operator} gives the first displayed
polynomial, because
\[
             (u+1)^4+u^4=1.
\]
For \(n=3\), primitive evaluation of
\(\mathcal A^+(\mathbf K_2^+,\mathbf K_1^+;\mathbf K_1^+)\) gives
\[
 \mathsf Z_3^+(X_0,X_1)
 \doteq X_1
 \bigl(X_0^4+d^{-1}X_0^3X_1+X_1^4\bigr)^2,
\]
which yields the stated native polynomial.  This agrees with
\eqref{eq:D3-binary-thesis}, but here it has been obtained from
\(\mathcal C_d\)-multiplication rather than pulled back from
\(W_d^+\).
\end{proof}

\subsection{Full-point certification and addition-chain evaluation}

The Kummer denominator is the most compact single division polynomial, but
the full native complete laws give a second certification.  In odd
characteristic the identity in centered native Segre coordinates is
\[
                    O=(0:1:0:1).
\]
If repeated native addition and doubling produce
\[
        [n]P=(X_n:Y_n:T_n:Z_n),
\]
then
\begin{equation}\label{eq:full-native-division-ideal-odd}
        [n]P=O
        \quad\Longleftrightarrow\quad
        X_n=T_n=Y_n-Z_n=0.
\end{equation}
In characteristic two, using the characteristic-free Segre coordinates
\((A:B:C:D)\), one has
\[
                    O=(0:0:1:0),
\]
and therefore
\begin{equation}\label{eq:full-native-division-ideal-binary}
        [n]P=O
        \quad\Longleftrightarrow\quad
        A_n=B_n=D_n=0.
\end{equation}
Equations~\eqref{eq:full-native-division-ideal-odd} and
\eqref{eq:full-native-division-ideal-binary} are native division ideals on
the smooth completion.  Eliminating the redundant full-point coordinates
recovers the Kummer polynomial
\(\Delta^{\mathrm{nat}}_{n,d}\).

For evaluation, it is unnecessary to expand any degree-\(n^2\) polynomial.
A binary addition chain evaluates the paired recursion
\((\mathbf K_m,\mathbf K_{m+1})\) with one native \(x\)DBL and one native
\(x\)ADD per processed bit.  In odd characteristic the projective
difference schedule costs
\[
                   6\M+4\Sqr+\Dpar
\]
per bit, and the first native double has the special cost
\(\M+\Sqr+\Dpar\) from
Corollary~\ref{cor:first-dbl-thesis}.  The output second coordinate is the
value of the native division polynomial up to a nonzero projective scale.
Thus the same circuit serves scalar multiplication, torsion testing, and
division-polynomial evaluation while keeping every input and output on the
\(\mathcal C_d\) interface.

\chapter{Mean Values on Division Fibers}
\label{ch:mean-values}
This chapter changes the sample space used in Chapter~\ref{ch:first-moment}.
There the averaging variable was the parameter \(d\); here \(d\) is fixed
and the averaging variable is a point in one multiplication fiber.  The
first section proves a characteristic-uniform trace identity for
multiplication maps of degree prime to the characteristic.  The second uses the odd-characteristic Edwards dictionary for
two reciprocal centered coordinates.  Low-degree checks are postponed until
both trace identities are available.

Let \(n\) be prime to \(p=\charac\F_q\), let \(Q\) be generic, and put
\[
\Lambda_{n,Q}=\{P:[n]P=Q\},\qquad
\langle h\rangle_{n,Q}
=\frac1{n^2}\sum_{P\in\Lambda_{n,Q}}h(P).
\]
If a fiber meets a pole, the identities below are interpreted as
function-field trace identities and extended from a dense open set.

\section{The prime-to-characteristic native Kummer trace}

\begin{theorem}[Prime-to-characteristic native mean]\label{thm:Kummer-mean-thesis}
Assume \(n\ge1\) and \(p\nmid n\), as in the standing convention of this
chapter.  For
\[
            U(P)=\frac{u(P)+1}{u(P)},
\]
one has
\begin{equation}\label{eq:Kummer-mean-thesis}
       \langle U\rangle_{n,Q}=U(Q),
\qquad
       \sum_{[n]P=Q}U(P)=n^2U(Q).
\end{equation}
\end{theorem}

\begin{proof}
If \(n=1\), the fiber consists of the single point \(Q\), and both sides
of \eqref{eq:Kummer-mean-thesis} equal \(U(Q)\).  Hence assume \(n\ge2\).
The argument has two stages.  We first prove a universal trace identity for
the Weierstrass abscissa by comparing the two leading coefficients of the
division polynomial.  We then apply the characteristic-specific dictionary
and return the result to the native Kummer coordinate \(U=(u+1)/u\).

\medskip
\noindent\emph{Stage 1: the universal Weierstrass trace.}
We first prove the coefficient statement needed below on the general
Weierstrass equation \eqref{eq:general-W}.
Let \(b_2=a_1^2+4a_2\), and use the generalized division polynomials of
Chapter~\ref{ch:division}.  Their two leading terms are
\begin{align}
 \psi_n(T)^2
   &=n^2T^{n^2-1}
     +\frac{n^2(n^2-1)}{12}b_2T^{n^2-2}
     +O(T^{n^2-3}),\label{eq:psi-square-leading-mean}\\
 \psi_{n-1}(T)\psi_{n+1}(T)
   &=(n^2-1)T^{n^2}
     +\frac{n^2(n^2-1)}{12}b_2T^{n^2-1}
     +O(T^{n^2-2}).\label{eq:psi-neighbours-leading-mean}
\end{align}
We now justify these coefficients in the universal coefficient ring.
Work first over its rationalization.  The admissible change
\[
 x'=T+\frac{b_2}{12},
 \qquad
 y'=Y+\frac{a_1T+a_3}{2}
\]
puts the generalized Weierstrass equation into short form
\[
 y'^2=x'^3-\frac{c_4}{48}x'-\frac{c_6}{864}.
\]
The change has scale factor one, so the division functions \(\psi_m\) are
unchanged as rational functions.  On the short equation, assign weights
\(2,4,6\) to \(x'\), \(-c_4/48\), and \(-c_6/864\), respectively.
The polynomial \(\psi_n^2\) is weighted homogeneous of weight
\(2(n^2-1)\) and has leading term \(n^2x'^{n^2-1}\).  A term
\(x'^{n^2-2}\) would require a coefficient of weight two, but the
coefficient ring generated by the short-Weierstrass parameters of weights
four and six contains no nonzero element of weight two.  Hence the
coefficient of \(x'^{n^2-2}\) is zero.  The product
\(\psi_{n-1}\psi_{n+1}\) has weight \(2n^2\) and leading term
\((n^2-1)x'^{n^2}\).  A coefficient of \(x'^{n^2-1}\) would again have
weight two, which is absent from the same coefficient ring; that
coefficient is therefore zero as well.

Substituting \(x'=T+b_2/12\) now gives
\[
 \psi_n(T)^2
 =n^2T^{n^2-1}
  +\frac{n^2(n^2-1)}{12}b_2T^{n^2-2}
  +O(T^{n^2-3}),
\]
and
\[
 \psi_{n-1}(T)\psi_{n+1}(T)
 =(n^2-1)T^{n^2}
  +\frac{n^2(n^2-1)}{12}b_2T^{n^2-1}
  +O(T^{n^2-2}).
\]
The integer \(n^2(n^2-1)\) is divisible by \(12\): divisibility by
three follows from one of \(n-1,n,n+1\), and divisibility by four follows
from \(n^2\) when \(n\) is even and from \(n^2-1\) when \(n\) is odd.
Both sides are coefficients of integral universal division polynomials, so
the identities obtained after tensoring with \(\mathbb Q\) already hold in
the universal integral coefficient ring.  They may therefore be reduced to
any characteristic.  Consequently
\begin{equation}\label{eq:phi-leading-mean-proof}
 \phi_n(T)=T\psi_n(T)^2-\psi_{n-1}(T)\psi_{n+1}(T)
          =T^{n^2}+0\,T^{n^2-1}+O(T^{n^2-2}).
\end{equation}

Fix first a point \(Q\) that is not two-torsion and for which the fiber
avoids the poles of \(X\).  Since \(p\nmid n\), multiplication by \(n\) is
finite \'etale of degree \(n^2\).  The polynomial
\begin{equation}\label{eq:mean-fiber-polynomial}
 F_{n,Q}(T)=\phi_n(T)-X(Q)\psi_n(T)^2
\end{equation}
is monic of degree \(n^2\).  Its equation is exactly
\(X([n]P)=X(Q)\), hence \([n]P=Q\) or \([n]P=-Q\).  The two sets are
interchanged by \(P\mapsto-P\), and \(X(P)=X(-P)\).  Moreover, when
\(Q\ne-Q\), two points of \([n]^{-1}(Q)\) cannot have the same
\(X\)-coordinate: equality would make them negatives and would force
\(Q=-Q\).  It follows that
\begin{equation}\label{eq:mean-fiber-factorization}
       F_{n,Q}(T)=\prod_{[n]P=Q}\bigl(T-X(P)\bigr).
\end{equation}
By \eqref{eq:psi-square-leading-mean} and
\eqref{eq:phi-leading-mean-proof}, the coefficient of \(T^{n^2-1}\) in
\eqref{eq:mean-fiber-polynomial} is \(-n^2X(Q)\).  Vi\`ete's formula
applied to \eqref{eq:mean-fiber-factorization} therefore gives
\begin{equation}\label{eq:weierstrass-X-trace-mean}
             \sum_{[n]P=Q}X(P)=n^2X(Q).
\end{equation}
Both sides are rational functions of the generic point \(Q\).  Hence
\eqref{eq:weierstrass-X-trace-mean} is an identity in the function field
of \(E\), and it extends to the omitted fibers in the trace sense specified
before the theorem.

\medskip
\noindent\emph{Stage 2: return to the native quotient.}
We now return to \(\mathcal C_d\).  In odd characteristic use the full
equation
\[
 W_d^-:\quad Y^2=X^3+\beta AX^2+\beta^2X,
 \qquad X=\beta U,
\]
where \(\beta=(16d)^{-1}\).  In characteristic two use
\[
 W_d^+:\quad Y^2+XY=X^3+d^2X,
 \qquad X=dU.
\]
The smoothness assumptions make \(\beta\) and \(d\) nonzero.  Applying
\eqref{eq:weierstrass-X-trace-mean} and dividing by the corresponding
fixed scale gives
\[
             \sum_{[n]P=Q}U(P)=n^2U(Q).
\]
Division by \(n^2\) proves \eqref{eq:Kummer-mean-thesis}.  In
characteristic two the use of
\(\psi_2=2Y+a_1X+a_3=X\), rather than the short-Weierstrass expression
\(2Y\), is precisely what makes the preceding universal argument valid.
\end{proof}

\section[Reciprocal trace identities]{Reciprocal centered-coordinate traces in odd characteristic}

\begin{theorem}[Odd-characteristic reciprocal means]
\label{thm:reciprocal-means-thesis}
Let \(p\ne2\) and \(p\nmid n\).  For odd \(n\), set
\(\epsilon_n=(-1)^{(n-1)/2}\); then
\begin{align}
\frac1{n^2}\sum_{[n]P=Q}\frac1{2v(P)+1}
 &=\frac1{n(2v(Q)+1)},\label{eq:mean-v-thesis}\\
\frac1{n^2}\sum_{[n]P=Q}\frac1{2u(P)+1}
 &=\frac{\epsilon_n}{n(2u(Q)+1)}.
\label{eq:mean-u-thesis}
\end{align}
For even \(n\), both averages are zero.
\end{theorem}

\begin{proof}
Put
\[
 \xi=(2v+1)^{-1},\qquad \eta=(2u+1)^{-1}.
\]
Then \(\mathcal C_d\) is written as the Edwards curve
\begin{equation}\label{eq:Edwards-mean-proof}
 E_\rho:\quad \xi^2+\eta^2=1+\rho\xi^2\eta^2,
 \qquad \rho=1-16d,
\end{equation}
with identity \((0,1)\).  Assume first that \(n=2k+1\) is odd.
The Edwards mean-value identity of
\cite{MoodyMeanEdwards2011}, in the normalization
\eqref{eq:Edwards-mean-proof}, is the function-field trace formula
\begin{equation}\label{eq:eta-trace-mean-proof}
       \sum_{[n]P=Q}\eta(P)=(-1)^k n\eta(Q).
\end{equation}
More precisely, multiplication by \(n\) is finite \'etale of degree
\(n^2\), and the left side is
\(\operatorname{Tr}_{[n]}(\eta)\).  The cited identity is an equality on
the universal smooth Edwards family over
\[
  \mathbb Z[1/2,n^{-1},\rho,\rho^{-1},(1-\rho)^{-1}].
\]
Function-field trace commutes with specialization under this finite
\'etale base change.  It therefore applies in every odd characteristic
prime to \(n\), including characteristic three when \(3\nmid n\).
This invocation supplies both the coefficient and multiplicity statements,
including the leading coefficient used in the displayed trace identity.

To obtain the other coordinate without a second elimination, let
\(R=(1,0)\), a point of order four on \(E_\rho\).  The Edwards law gives
\[
       P+R=(\eta(P),-\xi(P)),\qquad
       P-R=(-\eta(P),\xi(P)).
\]
Translation by \(R\) bijects the fiber over \(Q\) with the fiber over
\(Q+[n]R\).  If \(n\equiv1\pmod4\), then \([n]R=R\) and
\(\eta(Q+R)=-\xi(Q)\); if \(n\equiv3\pmod4\), then
\([n]R=-R\) and \(\eta(Q-R)=\xi(Q)\).  Applying
\eqref{eq:eta-trace-mean-proof} to the translated fiber in the two cases
gives, in both cases,
\begin{equation}\label{eq:xi-trace-mean-proof}
                  \sum_{[n]P=Q}\xi(P)=n\xi(Q).
\end{equation}

Now let \(n\) be even.  The nontrivial point
\(T=(0,-1)\in E_\rho[2]\) acts by
\[
                    P+T=(-\xi(P),-\eta(P)).
\]
Because \(T\in E_\rho[n]\), this fixed-point-free translation preserves
every \([n]\)-fiber.  Pairing \(P\) with \(P+T\) proves directly that
both coordinate sums are zero.

The trace identity and the translation argument are identities of rational
functions on the smooth Edwards family.  As in
Theorem~\ref{thm:Kummer-mean-thesis}, equality on the generic fiber extends
to pole fibers as a function-field trace identity.  Finally substitute
\(\xi=(2v+1)^{-1}\), \(\eta=(2u+1)^{-1}\), and divide the odd sums by
\(n^2\).  This proves \eqref{eq:mean-v-thesis}--
\eqref{eq:mean-u-thesis} and the even case.
\end{proof}

\section{Low-degree checks}

For \(\xi_Q=(2v(Q)+1)^{-1}\), \(\eta_Q=(2u(Q)+1)^{-1}\),
\[
\begin{array}{c|c|c|c}
n&\langle U\rangle&
\langle(2v+1)^{-1}\rangle&
\langle(2u+1)^{-1}\rangle\\ \hline
2&U(Q)&0&0\\
3&U(Q)&\xi_Q/3&-\eta_Q/3\\
5&U(Q)&\xi_Q/5&\eta_Q/5\\
7&U(Q)&\xi_Q/7&-\eta_Q/7\\
9&U(Q)&\xi_Q/9&\eta_Q/9
\end{array}
\]
whenever the index is prime to the characteristic.

For \(U(Q)=q_0\), the odd-characteristic third-division fiber polynomial is
\begin{equation}\label{eq:odd-3fiber-thesis}
\mathfrak F_{3,q_0}^{(-)}(T)
=(T-q_0)P_3(T)^2-8f(T)F_4(T).
\end{equation}
Its \(T^8\)-coefficient is \(-9q_0\).  In characteristic two,
\begin{equation}\label{eq:binary-3fiber-thesis}
\mathfrak F_{3,q_0}^{(+)}(T)
=T^9+q_0T^8+d^{-2}q_0T^6+d^{-2}T^3+T+q_0.
\end{equation}
Its root sum is \(q_0=9q_0\), confirming the same mean identity.

\begin{example}
Over \(\F_{101}\), take \(d=1\), so \(A=74\), and let \(U(Q)=13\).
Then \eqref{eq:odd-3fiber-thesis} is
\[
\begin{aligned}
T^9+85T^8+29T^7+20T^6+11T^5
+31T^4+4T^3+13T^2+9T+88.
\end{aligned}
\]
The nine roots sum to \(-85=16=9\cdot13\pmod{101}\).
\end{example}

\chapter{Model-Preserving Isogenies}
\label{ch:isogenies}
The goal of this chapter is an isogeny
\[
 \varphi:\mathcal C_d\longrightarrow\mathcal C_{d'},\qquad
 (u,v)\longmapsto(u',v'),
\]
with both source and target written as
\[
 (u^2+u)(v^2+v)=d,\qquad
 ({u'}^2+u')({v'}^2+v')=d'.
\]
The generalized Weierstrass formulas in the first section are the proof
engine for an arbitrary separable kernel.  The later theorems explicitly
return \(d'\), \(u'\), and \(v'\); an isogeny whose codomain is left only
in Weierstrass or Edwards form is not called model-preserving here.

The dependency order is as follows.  Chapter~\ref{ch:division} supplies
kernel and torsion data; the odd and binary dictionaries supply temporary
derivation coordinates.  We first state the common V\'elu interface, then
derive the odd and binary model-preserving formulas separately.  Only after
both characteristic branches are available do we discuss descent,
characteristic-three inseparability, composition, duality, and examples.
Thus the structural results near the end of the chapter depend on the
explicit maps already constructed near its beginning.

\section{The generalized V\'elu interface}
\label{sec:generalized-Velu-interface}

Let
\[
E:\quad y^2+a_1xy+a_3y=x^3+a_2x^2+a_4x+a_6
\]
and let an odd-order separable kernel be
\[
K=\{O,\pm Q_1,\ldots,\pm Q_m\},\qquad \ell=2m+1,
\]
with \(Q_i=(x_i,y_i)\).  Put
\[
\begin{aligned}
f_{x,i}&=3x_i^2+2a_2x_i+a_4-a_1y_i,\\
f_{y,i}&=-2y_i-a_1x_i-a_3,\\
t_i&=2f_{x,i}-a_1f_{y,i},\qquad
w_i=f_{y,i}^2+x_it_i,\\
t_K&=\sum_i t_i,\qquad w_K=\sum_iw_i.
\end{aligned}
\]
The normalized V\'elu map is characterized by
\begin{align}
\widehat x(P)&=x(P)+
\sum_{Q\in K\setminus\{O\}}\{x(P+Q)-x(Q)\},\label{eq:velu-xsum-thesis}\\
\widehat y(P)&=y(P)+
\sum_{Q\in K\setminus\{O\}}\{y(P+Q)-y(Q)\},
\label{eq:velu-ysum-thesis}
\end{align}
and the codomain is
\begin{equation}\label{eq:velu-codomain-thesis}
\widehat y^2+a_1\widehat x\widehat y+a_3\widehat y
=\widehat x^3+a_2\widehat x^2
+(a_4-5t_K)\widehat x
+a_6-(a_1^2+4a_2)t_K-7w_K.
\end{equation}
This formulation remains valid in low characteristic.  A short-Weierstrass
formula cannot be reduced modulo two in its place.

\section{Odd-characteristic product formulas}

Assume \(\charac k\ne2\), let \(d(1-16d)\ne0\), and let
\[
 K=\{O,\pm Q_1,\ldots,\pm Q_m\}
\]
be a finite separable subgroup of odd order \(\ell=2m+1\).  Put
\[
r_P=2u(P)+1,\qquad s_P=2v(P)+1.
\]
For \(Q_i=(u_i,v_i)\), write
\[
r_i=2u_i+1,\qquad s_i=2v_i+1,\qquad R_K=\prod_{i=1}^m r_i,
\]
and define
\begin{align}
\Delta_i(P)&=(r_Pr_is_Ps_i)^2-\rho^2,\label{eq:velu-delta-thesis}\\
A_i(P)&=s_i^2s_P^2-r_i^2r_P^2,\label{eq:velu-A-thesis}\\
B_i(P)&=s_i^2r_P^2-r_i^2s_P^2.\label{eq:velu-B-thesis}
\end{align}

Every nonidentity kernel point occurring here satisfies \(r_is_i\ne0\).
Indeed, the four boundary points have order dividing four, so a nonidentity
odd-order kernel point is affine.  If \(r_i=0\), then its native Kummer
value \(x=(r_i+1)/(r_i-1)\) equals \(-1\).  At the representative
\((X:Z)=(-1:1)\), native \(x\)DBL has first coordinate zero and second
coordinate \(16(\alpha_{24}-1)\).  Since
\(\alpha_{24}=1/(1-\rho)\) and \(\rho\ne0\), this second coordinate is
nonzero.  Thus \(\kappa(2Q_i)=\kappa(T)\), whence \(2Q_i=T\) and
\(Q_i\) has order four, a contradiction.  If \(s_i=0\), then
\(H=\sigma(Q_i)=R-Q_i\) has \(r(H)=0\), so \(2H=T\).  Since
\(2R=T\), it follows that \(2Q_i=O\), again contradicting the nonzero
odd order of \(Q_i\).  Therefore the reciprocal Edwards coordinates used
below are defined, \(R_K\ne0\), and every displayed inverse involving
\(r_i\) or \(s_i\) is legitimate.

\begin{theorem}[Odd-degree product isogeny]\label{thm:odd-velu-thesis}
Define
\begin{equation}\label{eq:odd-velu-parameter-thesis}
        \rho'=\frac{\rho^\ell}{R_K^8},\qquad
        d'=\frac{1-\rho'}{16}.
\end{equation}
On a dense open set, the model-preserving isogeny with kernel \(K\) is
\begin{align}
r_{\varphi(P)}
 &=\frac{r_P}{R_K^2}\prod_{i=1}^m
       \frac{\Delta_i(P)}{A_i(P)},\label{eq:velu-r-thesis}\\
s_{\varphi(P)}
 &=\frac{s_P}{R_K^2}\prod_{i=1}^m
       \frac{\Delta_i(P)}{B_i(P)}.\label{eq:velu-s-thesis}
\end{align}
The native output is
\[
       u'=\frac{r_{\varphi(P)}-1}{2},\qquad
       v'=\frac{s_{\varphi(P)}-1}{2}.
\]
Projective extension maps the kernel to the identity.
The extension is a separable isogeny of degree \(\ell\), its kernel is
exactly \(K\), and \(\rho'\ne0,1\); hence \(d'(1-16d')\ne0\).
\end{theorem}

\begin{proof}
We separate the proof into the quotient construction, identification of the
kernel and degree, determination of the target parameter, and return to the
native coordinates.

\smallskip
\noindent\emph{Step 1: the quotient map in centered reciprocal coordinates.}
In Edwards coordinates,
\[
(\xi,\eta)=(s_P^{-1},r_P^{-1}),\qquad
(\alpha_i,\beta_i)=(s_i^{-1},r_i^{-1}),
\]
and \(B_K=\prod_i\beta_i=R_K^{-1}\).
The Moody--Shumow product isogeny is
\[
\begin{aligned}
\xi'&=\frac{\xi}{B_K^2}
\prod_i\frac{\beta_i^2\xi^2-\alpha_i^2\eta^2}
 {1-\rho^2\alpha_i^2\beta_i^2\xi^2\eta^2},\\
\eta'&=\frac{\eta}{B_K^2}
\prod_i\frac{\beta_i^2\eta^2-\alpha_i^2\xi^2}
 {1-\rho^2\alpha_i^2\beta_i^2\xi^2\eta^2}
\end{aligned}
\]
\cite{MoodyShumow2016}.  Before simplification, these two products are
obtained by multiplying the translates of the Edwards coordinate functions
over \(K\).  Translation by any \(Q\in K\) merely permutes the factors.
Thus the map is \(K\)-invariant.  It sends
\(O=(0,1)\) to \((0,1)\), and its first coordinate has a nonzero linear
term at \(O\); consequently it is nonconstant.  A rational map from a
smooth projective curve to a projective curve extends across the finitely
many apparent poles.  The resulting morphism preserves the identity, and
is therefore a homomorphism of elliptic curves.

\smallskip
\noindent\emph{Step 2: kernel and degree.}
The unsimplified translation product shows that \(\xi'(P)=0\) precisely
when \(\xi(P+Q)=0\) for some \(Q\in K\).  The two points of
\(\mathcal E_\rho\) with first coordinate zero are
\(O=(0,1)\) and \(T_2=(0,-1)\).  If \(P+Q=O\), then
\(P=-Q\in K\).  If \(P+Q=T_2\), translation invariance and the displayed
product give
\[
       \Psi_K(P)=\Psi_K(T_2-Q)=\Psi_K(T_2)=T_2\ne O;
\]
indeed, at \(T_2=(0,-1)\) the product formula gives
\[
 \xi'(T_2)=0,\qquad
 \eta'(T_2)=\frac{-1}{B_K^2}\prod_i\beta_i^2=-1,
\]
because \(B_K=\prod_i\beta_i\).  Thus \(\Psi_K(T_2)=T_2\), and here
oddness of \(\#K\) is essential: it ensures that the nonzero two-torsion
point \(T_2\) does not belong to \(K\).  Hence
\(\ker\Psi_K=K\).  Since \(K\) consists of \(\ell\) distinct geometric
points, it is \'etale, the quotient map is separable, and
\[
                 \deg\Psi_K=\#\ker\Psi_K=\ell.
\]

\smallskip
\noindent\emph{Step 3: the target parameter.}
Substitution of the Edwards addition law in the translation products shows
that the image has an ordinary Edwards equation
\[
        {\xi'}^2+{\eta'}^2=1+\rho'{\xi'}^2{\eta'}^2
\]
for one constant \(\rho'\).  We determine that constant without expanding
the full products.  Work temporarily over a field containing
\(\lambda\) and \(\imath\) with \(\lambda^4=\rho\) and
\(\imath^2=-1\), and take
\[
              P_*=(\lambda^{-1},\imath\lambda^{-1}).
\]
Then \(P_*\in\mathcal E_\rho\).  For each kernel point, the identity
\(\alpha_i^2+\beta_i^2=1+\rho\alpha_i^2\beta_i^2\)
reduces the corresponding factor, and the product evaluates to
\[
 \Psi_K(P_*)=
 \left(\frac1{B_K^2\lambda^\ell},
       \frac{(-1)^m\imath}{B_K^2\lambda^\ell}\right).
\]
The sum of the squares of these two coordinates is zero.  Substitution in
the image equation therefore gives
\[
 0=1-\frac{\rho'}{B_K^8\lambda^{4\ell}},
 \qquad\text{so}\qquad
 \rho'=B_K^8\rho^\ell.
\]
The calculation took place after a scalar extension, but the resulting
identity has coefficients in the original field.  Therefore it descends,
and the displayed value is the target parameter \(\rho'\).

\smallskip
\noindent\emph{Step 4: return to \(\mathcal C_d\).}
Since \(B_K=R_K^{-1}\), the preceding equality is
\[
             \rho'=\frac{\rho^\ell}{R_K^8}.
\]
In every factor of the two Edwards products, multiply numerator and
denominator by \(r_P^2s_P^2r_i^2s_i^2\).  For the first coordinate this
changes the numerator and denominator into \(B_i(P)\) and
\(\Delta_i(P)\), respectively; for the second it gives
\(A_i(P)\) and \(\Delta_i(P)\).  Taking reciprocals
\(r'=1/\eta'\), \(s'=1/\xi'\) now yields
\eqref{eq:velu-r-thesis}--\eqref{eq:velu-s-thesis}.  Finally
\(d'=(1-\rho')/16\) gives the asserted native target.

The image of a separable quotient of a smooth elliptic curve is again a
smooth elliptic curve.  Hence its ordinary Edwards parameter cannot be
\(0\) or \(1\); equivalently \(\rho'\ne0,1\), and
\(d'(1-16d')\ne0\).  Apparent zeros in the affine numerators and
denominators are resolved by the projective extension constructed in
Step~1, which sends every point of \(K\) to \(O\).

Replacing \(Q_i\) by \(-Q_i\) changes \(s_i\) to \(-s_i\) and fixes
\(r_i\).  Every displayed factor uses \(s_i^2\), while \(R_K\) is
unchanged.  Thus the formula is independent of the chosen half-kernel
representatives, as required.
\end{proof}

\begin{remark}[Normalization]
The product formulas preserve the ordinary Edwards subfamily and hence the
shape \(\C_{d'}\), but they need not be V\'elu-normalized.  If
\(A_K=\prod_i\alpha_i\), \(B_K=\prod_i\beta_i\), the normalized
Moody--Shumow codomain is twisted Edwards with
\[
\widehat a=A_K^4/B_K^4,\qquad
\widehat d_E=A_K^4B_K^4\rho^\ell.
\]
A final scaling returns to ordinary Edwards parameter
\(\widehat d_E/\widehat a=B_K^8\rho^\ell\).  Kernel normalization and
choice of a canonical family representative are distinct operations.
\end{remark}

\section[Native u-Kummer specialization]
{The native \texorpdfstring{\(u\)}{u}-Kummer specialization}
\label{sec:native-u-Kummer-isogeny-general}

The full centered-coordinate products above prove the global quotient and
return both native coordinates.  When only the quotient by sign is needed,
their specialization to the native \(u\)-Kummer line is substantially smaller.
This specialization belongs in the general isogeny theory because it is not a
protocol-dependent optimization: it is the rational map induced on
\(\mathcal C_d/\{\pm1\}\).  Chapter~\ref{ch:Cd-isogeny-cryptography}
deliberately derives the formula again, with its implementation schedule, so
that the cryptographic chapter remains independently readable.

\begin{theorem}[Native odd-degree \(u\)-Kummer quotient]
\label{thm:native-u-Kummer-isogeny-general}
Retain the hypotheses and notation of
Theorem~\ref{thm:odd-velu-thesis}.  For one representative
\(Q_i=(u_i,v_i)\) of each pair \(\{Q_i,-Q_i\}\), put
\begin{equation}
 \mathcal A_K(u)=\prod_{i=1}^m(u+u_i+1),\qquad
 \mathcal B_K(u)=\prod_{i=1}^m(u-u_i),
\label{eq:native-u-Kummer-products-general}
\end{equation}
and
\[
                 N=(u+1)\mathcal A_K(u)^2,\qquad
                 D=u\mathcal B_K(u)^2.
\]
Then the quotient isogeny to \(\mathcal C_{d'}\) induces
\begin{equation}
 \boxed{\qquad
       u'=\frac{D}{N-D},\qquad
       (u'_0:u'_1)=(D:N-D).
 \qquad}
\label{eq:native-u-Kummer-isogeny-general}
\end{equation}
Equivalently, for projective native Kummer pairs
\[
 (X:Z)=(u+1:u),\qquad (X_i:Z_i)=(u_i+1:u_i),
\]
define
\begin{equation}
 P_K=\prod_{i=1}^m(XX_i-ZZ_i),\qquad
 Q_K=\prod_{i=1}^m(XZ_i-ZX_i).
\label{eq:native-u-Kummer-projective-products-general}
\end{equation}
Then an inversion-free target Kummer pair is
\begin{equation}
       \boxed{\qquad
       (X':Z')=(XP_K^2:ZQ_K^2).
       \qquad}
\label{eq:native-u-Kummer-projective-general}
\end{equation}
All three formulas depend only on the unsigned kernel classes.
\end{theorem}

\begin{proof}
The exact Montgomery dictionary of this monograph is
\[
                    x=\frac{u+1}{u},\qquad
                    x_i=\frac{u_i+1}{u_i}.
\]
For a normalized cyclic quotient of odd order \(2m+1\), its Kummer map is
\begin{equation}
 f(x)=x\prod_{i=1}^m
       \left(\frac{xx_i-1}{x-x_i}\right)^2
 \cite{CostelloHisil2017}.
\label{eq:Montgomery-Kummer-product-general}
\end{equation}
The two linear factors simplify without any root extraction:
\begin{align*}
 xx_i-1
 &=\frac{(u+1)(u_i+1)-uu_i}{uu_i}
   =\frac{u+u_i+1}{uu_i},\\
 x-x_i
 &=\frac{(u+1)u_i-(u_i+1)u}{uu_i}
   =-\frac{u-u_i}{uu_i}.
\end{align*}
The sign and the common denominator disappear after squaring.  Hence
\[
 f(x)=\frac{u+1}{u}
       \left(\frac{\mathcal A_K(u)}{\mathcal B_K(u)}\right)^2
      =\frac ND.
\]
The target dictionary is \(u'=1/(x'-1)\).  Substitution of
\(x'=N/D\) gives \(u'=D/(N-D)\), proving
\eqref{eq:native-u-Kummer-isogeny-general}.

For the projective form, homogenize each factor of
\eqref{eq:Montgomery-Kummer-product-general}.  The numerator factor is
\(XX_i-ZZ_i\), the denominator factor is \(XZ_i-ZX_i\), and the leading
factor \(x=X/Z\) gives exactly
\eqref{eq:native-u-Kummer-projective-general}.  Negation fixes \(u_i\),
so replacing any \(Q_i\) by \(-Q_i\) changes none of the displayed data.
\end{proof}

\begin{proposition}[Exact native Kummer evaluation costs]
\label{prop:native-u-Kummer-isogeny-cost-general}
Assume \(m\ge1\), and ignore additions, subtractions, and multiplication
by \(2\).
\begin{enumerate}[label=\textup{(\roman*)}]
 \item Affine \(u,u_i\) to the projective \(u'\)-pair in
       \eqref{eq:native-u-Kummer-isogeny-general} costs
       \begin{equation}
                         \boxed{2m\M+2\Sqr}.
       \label{eq:native-u-Kummer-affine-cost-general}
       \end{equation}
 \item Normalizing that output costs one additional \(\M+\Inv\).
 \item General projective source and kernel Kummer pairs to
       \eqref{eq:native-u-Kummer-projective-general} cost
       \begin{equation}
                         \boxed{4m\M+2\Sqr}.
       \label{eq:native-u-Kummer-projective-cost-general}
       \end{equation}
\end{enumerate}
\end{proposition}

\begin{proof}
For \(m=0\), the kernel is \(\{O\}\), both products are empty, and the
quotient map is the identity; after suppressing operations on the unit
constant, its evaluation cost is zero.  The remaining count therefore
assumes \(m\ge1\).

For affine kernel values, the two products of \(m\) monic linear factors
cost \(2(m-1)\M\); the two squares and the two multiplications by \(u+1\)
and \(u\) give \(2\Sqr+2\M\).  This proves part~(i), and a single
projective normalization proves part~(ii).

For a projective factor pair, compute
\[
 a_i=(X-Z)(X_i+Z_i),\qquad
 b_i=(X+Z)(X_i-Z_i).
\]
Then
\[
 a_i+b_i=2(XX_i-ZZ_i),\qquad
 a_i-b_i=2(XZ_i-ZX_i).
\]
The projectively common powers of \(2\) may be omitted.  The \(m\) pairs
cost \(2m\M\), the two product accumulators cost \(2(m-1)\M\), and final
assembly costs \(2\M+2\Sqr\).  Their sum is \(4m\M+2\Sqr\).
\end{proof}

\begin{corollary}[Target parameter from the same native kernel]
\label{cor:native-u-Kummer-target-general}
With \(\rho=1-16d\), the target of
Theorem~\ref{thm:native-u-Kummer-isogeny-general} is the smooth member
\(\mathcal C_{d'}\) with
\begin{equation}
 \boxed{\qquad
 \rho'=\frac{\rho^{2m+1}}
  {\left(\prod_{i=1}^m(2u_i+1)\right)^8},\qquad
 d'=\frac{1-\rho'}{16}.
 \qquad}
\label{eq:native-u-Kummer-target-general}
\end{equation}
For projective pairs \((X_i:Z_i)=(u_i+1:u_i)\), the same identity is
\[
 \rho'=\rho^{2m+1}
   \left(\prod_i\frac{X_i-Z_i}{X_i+Z_i}\right)^8.
\]
\end{corollary}

\begin{proof}
Theorem~\ref{thm:odd-velu-thesis} gives the first expression with
\(R_K=\prod_i(2u_i+1)\).  Moreover
\[
 \frac{X_i-Z_i}{X_i+Z_i}=\frac1{2u_i+1},
\]
which gives the projective form.  The quotient is the same smooth quotient
already constructed in Theorem~\ref{thm:odd-velu-thesis}, so its parameter
satisfies \(\rho'\ne0,1\) and no additional closure argument is required.
\end{proof}

\section{Kernel polynomials and projective evaluation}

Define
\[
\mathcal K_1(T,Z)=\prod_i(\beta_i^2T-\alpha_i^2Z),\qquad
\mathcal K_0(T)=\prod_i(1-\rho^2\alpha_i^2\beta_i^2T).
\]
Then
\[
\xi'=\frac{\xi\mathcal K_1(\xi^2,\eta^2)}
 {B_K^2\mathcal K_0(\xi^2\eta^2)},\qquad
\eta'=\frac{\eta\mathcal K_1(\eta^2,\xi^2)}
 {B_K^2\mathcal K_0(\xi^2\eta^2)}.
\]
This directly connects the isogeny to the division polynomials of
Chapter~\ref{ch:division}.

For \(r=R/Z\), \(s=S/Z\), set
\[
\begin{aligned}
\Delta_i^h&=(RSr_is_i)^2-\rho^2Z^4,\\
A_i^h&=s_i^2S^2-r_i^2R^2,\qquad
B_i^h=s_i^2R^2-r_i^2S^2,
\end{aligned}
\]
and let \(\Delta^h,A^h,B^h\) be their products.  An inversion-free
representative is
\[
(R':S':Z')=
\bigl(R\Delta^hB^h:S\Delta^hA^h:
R_K^2Z^{2m+1}A^hB^h\bigr).
\]
All coordinates have degree \(6m+1\).

\subsection{Evaluation costs and kernel preprocessing}

The products should be separated into kernel preprocessing and per-point
evaluation.  Let \(\mathbf C_K\) denote multiplication by a precomputed
kernel-dependent constant; additions are ignored.  A direct schedule gives
the following conservative upper bounds.
\begin{table}[H]
\centering
\caption{Odd-degree \(\mathcal C_d\)-isogeny evaluation, \(\ell=2m+1\)}
\label{tab:Cd-isogeny-evaluation-cost}
\begin{tabular}{L{3.6cm}L{3.0cm}L{5.8cm}}
\toprule
input and output & per-point upper bound & assumptions\\
\midrule
affine \((r,s)\mapsto(r',s')\)
& \((3m+5)\M+2\Sqr+(5m+2)\mathbf C_K+\Inv\)
& the two denominator products use one simultaneous inversion\\
projective \((R:S:Z)\mapsto(R':S':Z')\)
& \((4m+4)\M+5\Sqr+\Dpar+(5m+1)\mathbf C_K\)
& sequential products and sequential evaluation of \(Z^{2m+1}\)\\
\bottomrule
\end{tabular}
\end{table}

For the affine row, compute \(r^2,s^2,r^2s^2\) once, form the three
families \(\Delta_i,A_i,B_i\), and accumulate their products.  The two
final denominator products are inverted with Montgomery's
simultaneous-inversion trick, using one inversion and three
multiplications.  For the projective row, compute
\(R^2,S^2,(RS)^2,Z^2,Z^4\) once and form the reusable term
\(\rho^2Z^4\) with one fixed-parameter multiplication; the three product
accumulators cost \(3(m-1)\M\), and coordinate assembly gives the stated
bound.  These transparent generic bounds provide uniform evaluation
schedules for every odd degree \(\ell=2m+1\); the low-degree cases may be
specialized further by expanding their kernel products.

When a kernel is reused, the coefficients of
\(\mathcal K_0\) and \(\mathcal K_1\), the values
\(r_i^2,s_i^2,(r_is_i)^2\), and \(R_K^{\pm2}\) are precomputed once.
Horner evaluation of the kernel polynomials is linear in \(m\), while
product trees reduce kernel construction to quasi-linear polynomial
arithmetic for large \(\ell\).  This representation also avoids storing
both \(Q_i\) and \(-Q_i\).

\section{Characteristic-two half-kernel formulas}

\subsection{The distinguished separable two-isogeny}
\label{subsec:Cd-char2-two-isogeny-general}

The odd-order formulas below enumerate pairs \(\{\pm Q_i\}\).  In
characteristic two the reduced two-torsion instead gives a separate
degree-two quotient.  It is particularly small on the native \(u\)-Kummer
line and therefore must be recorded before the odd half-kernel construction.

\begin{theorem}[Native binary two-isogeny]
\label{thm:Cd-char2-two-isogeny-general}
Let \(k\) be a perfect field of characteristic two, let \(d\in k^\times\),
and let \(e\in k^\times\) be the unique element with \(e^2=d\).  On
\[
 W_d^+:\qquad y^2+xy=x^3+d^2x
\]
put \(T=(0,0)\).  The rational functions
\begin{equation}
 \begin{aligned}
 X&=x+\frac{d^2}{x},\\
 Y&=y+\frac{d^2(y+x)}{x^2}+d
 \end{aligned}
\label{eq:Cd-char2-two-isogeny-W-general}
\end{equation}
extend to a separable isogeny
\[
        \Phi_{2,d}:W_d^+\longrightarrow W_e^+
\]
of degree two with kernel \(\{O,T\}\).

Under the native dictionary
\[
 x=d\frac{u+1}{u},\qquad y=xv,
\]
the same isogeny is
\begin{equation}
 \boxed{\qquad
 \Phi_{2,d}^{(u)}(u,v)=
 \left(\frac{u(u+1)}{u(u+1)+e},\ u+v\right)
 \in\mathcal C_e .
 \qquad}
\label{eq:Cd-char2-two-isogeny-native-general}
\end{equation}
The first coordinate is understood projectively at its boundary.  If
\[
 z=\frac{u+1}{u}=\frac{K_0}{K_1},
\]
the induced Kummer map is
\begin{equation}
 \boxed{\qquad
 (K'_0:K'_1)=
 \bigl(e(K_0+K_1)^2:K_0K_1\bigr).
 \qquad}
\label{eq:Cd-char2-two-isogeny-Kummer-general}
\end{equation}
For the raw native input \((K_0:K_1)=(u+1:u)\), this specializes to
\begin{equation}
             (K'_0:K'_1)=(e:u(u+1)).
\label{eq:Cd-char2-two-isogeny-native-Kummer-general}
\end{equation}
\end{theorem}

\begin{proof}
\noindent\emph{Step 1: translation by the reduced two-torsion point.}
Let \(P=(x,y)\) with \(x\ne0\).  The line through \(P\) and \(T\) has
slope \(y/x\).  The generalized binary addition formula gives
\begin{equation}
 x(P+T)=\frac{d^2}{x},\qquad
 y(P+T)=\frac{d^2(y+x)}{x^2}.
\label{eq:Cd-char2-translation-by-T-general}
\end{equation}
Indeed, division of the curve equation by \(x^2\) gives
\[
 \left(\frac yx\right)^2+\frac yx=x+\frac{d^2}{x},
\]
and substitution in the chord formula yields the two displayed values.
Thus the unshifted functions
\[
 X=x+x(P+T),\qquad Y_0=y+y(P+T)
\]
are invariant under translation by \(T\).

\smallskip
\noindent\emph{Step 2: identify the quotient equation.}
Put \(c=d^2\).  Direct expansion, with every equality in characteristic
two, gives
\begin{align*}
 Y_0^2+XY_0
 &=(y^2+xy)\left(1+\frac{c^2}{x^4}\right)+c\\
 &=(x^3+cx)\left(1+\frac{c^2}{x^4}\right)+c\\
 &=X^3+c.
\end{align*}
After the shift \(Y=Y_0+d\),
\[
 Y^2+XY=X^3+c+d^2+dX=X^3+dX.
\]
Since \(e^2=d\), this is precisely \(W_e^+\).

\smallskip
\noindent\emph{Step 3: degree, separability, and kernel.}
The relation
\[
                    x^2+Xx+d^2=0
\]
shows that \(k(W_d^+)\) has degree at most two over \(k(X,Y_0)\).
Conversely, translation by \(T\) is nontrivial and fixes \(X,Y_0\), so
the degree is exactly two.  The derivative of the displayed polynomial
with respect to \(x\) is the nonzero function \(X\); hence the extension
is separable.  The two points \(O,T\) form the fiber over the target
identity, and the degree count leaves no additional kernel point.

\smallskip
\noindent\emph{Step 4: return to native coordinates.}
Write \(z=(u+1)/u\), so \(x=dz\), \(y=dzv\), and
\[
 X=d(z+z^{-1})=\frac{d}{u(u+1)}.
\]
The inverse target dictionary is \(u'=e/(X+e)\), \(v'=Y/X\).
Consequently
\[
 u'=\frac{u(u+1)}{u(u+1)+e}.
\]
Moreover
\[
 \frac YX=
 \frac{zv+(v+1)z^{-1}+1}{z+z^{-1}}=u+v;
\]
multiplication of numerator and denominator by \(u(u+1)\) verifies the
last equality term by term.  This proves
\eqref{eq:Cd-char2-two-isogeny-native-general}.

Finally \(z'=X/e=e(z+z^{-1})\).  Homogenizing \(z=K_0/K_1\) gives
\[
 z'=\frac{e(K_0^2+K_1^2)}{K_0K_1}
    =\frac{e(K_0+K_1)^2}{K_0K_1},
\]
which proves \eqref{eq:Cd-char2-two-isogeny-Kummer-general}.  For the
native pair \(K_0=u+1,K_1=u\), its sum is \(1\), proving the final
specialization.
\end{proof}

\begin{corollary}[Dual and exact costs]
\label{cor:Cd-char2-two-isogeny-cost-general}
The dual is relative Frobenius,
\begin{equation}
 \widehat\Phi_{2,d}:W_e^+\longrightarrow W_d^+,\qquad
             (X,Y)\longmapsto(X^2,Y^2),
\label{eq:Cd-char2-two-isogeny-dual-general}
\end{equation}
or natively \((u,v)\mapsto(u^2,v^2)\).  Both compositions are
multiplication by two.  Ignoring additions:
\begin{enumerate}[label=\textup{(\roman*)}]
 \item a general projective Kummer input in
       \eqref{eq:Cd-char2-two-isogeny-Kummer-general} costs
       \[
                         \M+\Sqr+\Dpar;
       \]
 \item a raw native affine input costs exactly
       \[
                         \boxed{\M}
       \]
       for the projective Kummer output
       \eqref{eq:Cd-char2-two-isogeny-native-Kummer-general};
 \item the full native map to the two projective target coordinates
       \[
        (u'_0:u'_1)=(u(u+1):u(u+1)+e),\qquad v'=u+v
       \]
       also costs exactly \(\boxed{\M}\);
 \item affine normalization of \(u'\) adds \(\M+\Inv\), while the dual
       costs two squarings.
\end{enumerate}
\end{corollary}

\begin{proof}
Squaring the equation of \(W_e^+\) shows immediately that
\((X^2,Y^2)\) lies on \(W_d^+\).  To verify the first composition
directly, put \(t=y/x\) and \(q=d^2/x\).  Then
\[
 t^2+t=x+q,\qquad
 \lambda=\frac{y+x^2+d^2}{x}=t^2
\]
for the tangent at \(P\).  Hence
\[
 x([2]P)=\lambda^2+\lambda=x^2+q^2
        =\left(x+\frac{d^2}{x}\right)^2.
\]
The tangent intercept is \(\nu=x^2+d^2\), so
\[
 y([2]P)=(t^2+1)(x^2+q^2)+x^2+d^2
         =t^6+t^4+q^2+d^2.
\]
On the other hand, the ordinate in
\eqref{eq:Cd-char2-two-isogeny-W-general} is
\[
 Y=t^3+t^2+q+d,
\]
whose square is the same expression.  Thus
\(\widehat\Phi_{2,d}\circ\Phi_{2,d}=[2]\); the identical calculation on
the target gives the reverse composition.  The native Frobenius formula
follows by squaring the two rational dictionaries.

For the costs, the general Kummer circuit forms \(K_0K_1\), squares
\(K_0+K_1\), and multiplies the square by \(e\).  On a native input the
sum is \(1\), so only \(u(u+1)\) remains.  The same product supplies the
full native first coordinate, and the second coordinate uses additions
only.  Projective normalization and the two coordinate squarings give
the remaining rows.
\end{proof}

On \(W_d^+:y^2+xy=x^3+d^2x\), write
\[
Q_i=(a_i,b_i),\qquad -Q_i=(a_i,b_i+a_i),
\]
and put
\[
        \Sigma_K=\sum_{i=1}^m a_i.
\]
Assume that \(\Sigma_K\) has a fourth root in the coefficient field, and
denote it by \(\tau_K\).  Such a root is necessarily unique in
characteristic two, and its existence is automatic over every perfect
field, in particular over every finite field.

\begin{theorem}[Binary model-preserving isogeny]
\label{thm:binary-velu-thesis}
For \(P=(x,y)\notin K\), define
\[
h_i=\frac{a_i}{x+a_i},\qquad
\lambda_i=\frac{y+b_i}{x+a_i},
\]
and
\begin{align}
\widehat x&=x+\sum_i\frac{a_ix}{(x+a_i)^2},
\label{eq:binary-velu-x-thesis}\\
\widehat y&=y+\sum_i
h_i(x+1+\lambda_i^2+h_i\lambda_i+h_i^2).
\label{eq:binary-velu-y-thesis}
\end{align}
Let
\begin{equation}\label{eq:binary-dprime-thesis}
d'=d+\tau_K^2+\tau_K,\qquad
x'=\widehat x,\qquad y'=\widehat y+\tau_K^2.
\end{equation}
Then \(d'\ne0\), \((x',y')\in W_{d'}^+\), and the displayed rational
functions extend to a separable isogeny of degree \(\ell=2m+1\) with
kernel exactly \(K\).  Moreover
\[
       u'=\frac{d'}{x'+d'},\qquad v'=\frac{y'}{x'}
\]
returns to \(\C_{d'}\).
\end{theorem}

\begin{proof}
\noindent\emph{Step 1: specialize the generalized V\'elu data.}
For
\[
 a_1=1,\qquad a_2=a_3=a_6=0,\qquad a_4=d^2,
\]
the definitions in Section~\ref{sec:generalized-Velu-interface} give, in
characteristic two,
\[
f_{x,i}=a_i^2+d^2+b_i,\qquad
f_{y,i}=a_i,\qquad t_i=a_i,\qquad w_i=0.
\]
Indeed, \(3=1\), \(2=0\), and
\[
 t_i=2f_{x,i}-f_{y,i}=a_i,\qquad
 w_i=f_{y,i}^2+a_it_i=a_i^2+a_i^2=0.
\]
Hence \(t_K=\Sigma_K\) and \(w_K=0\).
Thus \eqref{eq:velu-codomain-thesis} gives
\begin{equation}\label{eq:binary-raw-codomain-thesis}
\widehat y^2+\widehat x\widehat y
=\widehat x^3+(d^2+\Sigma_K)\widehat x+\Sigma_K.
\end{equation}

\smallskip
\noindent\emph{Step 2: simplify one kernel pair.}
Fix \(i\), abbreviate
\[
 a=a_i,\quad b=b_i,\quad
 h=\frac{a}{x+a},\quad \lambda=\frac{y+b}{x+a},
\]
and write \(Q=(a,b)\), \(-Q=(a,b+a)\).  The two line slopes are
\(\lambda\) and \(\lambda+h\).  The generalized Weierstrass addition law
on \(W_d^+\) gives
\[
\begin{aligned}
 x_+&=x(P+Q)=\lambda^2+\lambda+x+a,\\
 x_-&=x(P-Q)=(\lambda+h)^2+(\lambda+h)+x+a.
\end{aligned}
\]
Therefore the paired contribution to the V\'elu abscissa is
\[
 \{x_+-a\}+\{x_--a\}
 =h^2+h
 =\frac{a^2+a(x+a)}{(x+a)^2}
 =\frac{ax}{(x+a)^2},
\]
which proves the \(i\)th summand of
\eqref{eq:binary-velu-x-thesis}.

For the ordinates, put
\[
 \nu_+=y+\lambda x,\qquad \nu_-=\nu_++hx.
\]
The two third-intersection reflections are
\[
 y_+=(\lambda+1)x_++\nu_+,\qquad
 y_-=(\lambda+h+1)x_-+\nu_-.
\]
Since the ordinates of \(Q\) and \(-Q\) are \(b\) and \(b+a\), their
paired V\'elu contribution is \(y_++y_-+a\).  Substitution of \(x_+\)
and \(x_-\), followed by collection in characteristic two, gives
\[
\begin{aligned}
 y_++y_-+a
 &=h\lambda^2+h^2\lambda+ha+h^3+h+a\\
 &=h(x+1+\lambda^2+h\lambda+h^2).
\end{aligned}
\]
For the last equality we used \(h(x+a)=a\).  This is precisely the
\(i\)th summand of \eqref{eq:binary-velu-y-thesis}.

Replacing \(Q_i\) by \(-Q_i\) replaces
\(\lambda_i\) by \(\lambda_i+h_i\).  Directly,
\[
 (\lambda_i+h_i)^2+h_i(\lambda_i+h_i)+h_i^2
 =\lambda_i^2+h_i\lambda_i+h_i^2,
\]
so the simplified formula is independent of the chosen representative of
the pair.

\smallskip
\noindent\emph{Step 3: normalize the target.}
By definition \(\tau_K^4=\Sigma_K\).  With
 \(x'=\widehat x\) and \(y'=\widehat y+\tau_K^2\), the left side of the
target equation is
\[
 {y'}^2+x'y'
 =\widehat y^2+\widehat x\widehat y
   +\tau_K^4+\tau_K^2\widehat x.
\]
Insert \eqref{eq:binary-raw-codomain-thesis}.  Its constant term cancels
because \(\Sigma_K=\tau_K^4\), and the coefficient of \(x'\) becomes
\[
d^2+\tau_K^4+\tau_K^2
=(d+\tau_K^2+\tau_K)^2=d'^2.
\]
Thus
\[
                  {y'}^2+x'y'={x'}^3+d'^2x',
\]
which is \(W_{d'}^+\).

\smallskip
\noindent\emph{Step 4: kernel, degree, and return to the native model.}
The formulas \((\widehat x,\widehat y)\) were obtained by pairwise
simplification of the normalized V\'elu sums
\eqref{eq:velu-xsum-thesis}--\eqref{eq:velu-ysum-thesis}; hence they define
the quotient morphism by \(K\).  In particular all apparent poles at
\(x=a_i\) resolve on the smooth projective curve, every point of \(K\)
maps to the identity, and no point outside \(K\) does.  Because \(K\) is
separable and has \(2m+1=\ell\) geometric points, the quotient is separable,
\[
                    \ker\varphi=K,\qquad \deg\varphi=\ell.
\]
Translation of the ordinate by \(\tau_K^2\) is an isomorphism fixing the
identity, so it changes neither the kernel nor the degree.

The quotient of a smooth elliptic curve is smooth.  The curve \(W_{d'}^+\)
would be singular for \(d'=0\); therefore \(d'\ne0\).  Finally the binary
dictionary
\[
  u'=\frac{d'}{x'+d'},\qquad v'=\frac{y'}{x'}
\]
is the inverse of \eqref{eq:binary-W-map-thesis} on a dense open set and
extends to the smooth completions.  It sends the codomain to
\(\mathcal C_{d'}\), completing the proof.
\end{proof}

\begin{corollary}[Binary three-isogeny]\label{cor:binary-3iso-thesis}
For \(K=\langle Q\rangle\), \(Q=(a,b)\) of order three, let
\(\tau^4=a\).  Then
\[
d'=d+\tau^2+\tau,\quad
h=\frac a{x+a},\quad
\lambda=\frac{y+b}{x+a},
\]
\[
x'=x+\frac{ax}{(x+a)^2},\qquad
y'=y+h(x+1+\lambda^2+h\lambda+h^2)+\tau^2.
\]
If the input and kernel generator are stored on \(\mathcal C_d\) as
\[
 P=(u,v),\qquad Q=(u_Q,v_Q),\qquad
 U=\frac{u+1}{u},\quad U_Q=\frac{u_Q+1}{u_Q},
\]
then the same three-isogeny can be evaluated without first materializing a
separate curve point by
\begin{equation}\label{eq:binary-native-3iso-wrapper}
\boxed{
\begin{aligned}
 a&=dU_Q,& \tau^4&=dU_Q,& d'&=d+\tau^2+\tau,\\
 h&=\frac{U_Q}{U+U_Q},&
 \lambda&=\frac{Uv+U_Qv_Q}{U+U_Q},\\
 x'&=dU+\frac{UU_Q}{(U+U_Q)^2},\\
 y'&=dUv+h(dU+1+\lambda^2+h\lambda+h^2)+\tau^2,\\
 u'&=\frac{d'}{x'+d'},&v'&=\frac{y'}{x'}.
\end{aligned}}
\end{equation}
The projective extension of \eqref{eq:binary-native-3iso-wrapper} has
degree three, kernel \(\langle Q\rangle\), and codomain
\(\mathcal C_{d'}\).
\end{corollary}

\begin{proof}
For a three-element kernel one has \(m=1\),
\(\Sigma_K=a\), and \(\tau_K=\tau\).  Substitution in
Theorem~\ref{thm:binary-velu-thesis} immediately gives the first displayed
formulas and proves that the map has degree three, kernel
\(\{O,Q,-Q\}\), and target \(W_{d'}^+\).

The formula does not depend on which of \(Q\) and \(-Q\) is chosen.  Indeed,
negation on \(W_d^+\) replaces \(b\) by \(b+a\), hence replaces
\(\lambda\) by \(\lambda+h\).  In characteristic two,
\[
 (\lambda+h)^2+h(\lambda+h)+h^2
   =\lambda^2+h\lambda+h^2,
\]
so the ordinate formula is unchanged.

It remains to verify the native wrapper.  The binary dictionary is
\[
        x=dU,\qquad y=dUv,\qquad
        a=dU_Q,\qquad b=dU_Qv_Q.
\]
Therefore
\[
 \frac{a}{x+a}=\frac{U_Q}{U+U_Q},\qquad
 \frac{y+b}{x+a}=\frac{Uv+U_Qv_Q}{U+U_Q},
\]
and
\[
 x+\frac{ax}{(x+a)^2}
   =dU+\frac{UU_Q}{(U+U_Q)^2}.
\]
These substitutions give every line of
\eqref{eq:binary-native-3iso-wrapper}.  Finally, the inverse binary
dictionary
\(u'=d'/(x'+d')\), \(v'=y'/x'\) proves directly that the output is a
point of \(\mathcal C_{d'}\).  At vanishing affine denominators the
degree-three morphism supplied by Theorem~\ref{thm:binary-velu-thesis}
provides the unique projective extension.
\end{proof}

\section{Field of definition and descent}

\begin{theorem}[Descent from a Galois-stable kernel]
\label{thm:isogeny-kernel-descent}
Let \(k=\F_q\), and let \(K\) be a finite separable odd-order subgroup of
\(\mathcal C_d(\overline{k})\) stable under
\(\operatorname{Gal}(\overline{k}/k)\).  Then the model-preserving
isogenies of Theorems~\ref{thm:odd-velu-thesis} and
\ref{thm:binary-velu-thesis} are defined over \(k\), even when the
individual half-kernel points are not \(k\)-rational.  Their target
parameters \(d'\) belong to \(k\).
\end{theorem}

\begin{proof}
Write \(\Gamma=\operatorname{Gal}(\overline{k}/k)\).  We prove not only
that the target parameter is fixed by \(\Gamma\), but also that the two
coordinate functions of the quotient lie in the native function field
\(k(\mathcal C_d)\).

Suppose first that \(\charac k\ne2\).  Since \(K\) has odd order, no
nonidentity point of \(K\) is equal to its negative.  Thus
\[
 K\setminus\{O\}=\coprod_{i=1}^m\{Q_i,-Q_i\}
\]
is an intrinsic partition into unordered pairs.  Galois stability of
\(K\) means that, for every \(\gamma\in\Gamma\), there are a permutation
\(\pi_\gamma\) of \(\{1,\ldots,m\}\) and signs
\(\epsilon_{\gamma,i}\in\{\pm1\}\) such that
\[
       \gamma(Q_i)=\epsilon_{\gamma,i}Q_{\pi_\gamma(i)}.
\]
In the centered coordinates \(r=2u+1\), \(s=2v+1\), the inverse map is
\((r,s)\mapsto(r,-s)\).  Consequently
\[
 \gamma(r_i)=r_{\pi_\gamma(i)},\qquad
 \gamma(s_i)=\epsilon_{\gamma,i}s_{\pi_\gamma(i)}.
\]
It follows term by term that
\[
\begin{aligned}
 \gamma(R_K)&=R_K,&
 \gamma\!\left(\prod_i\Delta_i\right)&=\prod_i\Delta_i,\\
 \gamma\!\left(\prod_i A_i\right)&=\prod_i A_i,&
 \gamma\!\left(\prod_i B_i\right)&=\prod_i B_i.
\end{aligned}
\]
Here the possible signs disappear because \(A_i,B_i,\Delta_i\) contain
\(r_i,s_i\) only through squares.  By definition, every coefficient of
\(\mathcal K_0\) and \(\mathcal K_1\) is a sum of products of these
quantities over the full index set.  The permutation
\(\pi_\gamma\) merely reorders those summands and factors, so every such
coefficient is fixed by \(\Gamma\).  Since \(r,s,\rho\) are themselves
defined over \(k\), the
right sides of \eqref{eq:velu-r-thesis} and
\eqref{eq:velu-s-thesis} are fixed elements of
\(\overline{k}(\mathcal C_d)\).  The standard fixed-field identity
\[
       \overline{k}(\mathcal C_d)^\Gamma=k(\mathcal C_d)
\]
therefore places both coordinate functions in \(k(\mathcal C_d)\).
Moreover
\[
       \gamma(\rho')
       =\gamma\!\left(\frac{\rho^\ell}{R_K^8}\right)
       =\frac{\rho^\ell}{R_K^8}=\rho',
\]
so \(\rho'\), and hence \(d'=(1-\rho')/16\), lies in \(k\).

Now suppose that \(\charac k=2\).  On
\(W_d^+:y^2+xy=x^3+d^2x\), negation sends
\((a_i,b_i)\) to \((a_i,b_i+a_i)\); in particular it fixes the
abscissa \(a_i\).  Hence \(\Gamma\) permutes the multiset
\(\{a_1,\ldots,a_m\}\), and
\[
                    \Sigma_K=\sum_{i=1}^m a_i\in k.
\]
The paired summands in \eqref{eq:binary-velu-x-thesis} and
\eqref{eq:binary-velu-y-thesis} are independent of the choice between
\(Q_i\) and \(-Q_i\): for the ordinate this is the identity
\[
 (\lambda_i+h_i)^2+h_i(\lambda_i+h_i)+h_i^2
       =\lambda_i^2+h_i\lambda_i+h_i^2.
\]
Galois therefore only permutes these paired summands.  Their sums
\(\widehat x,\widehat y\) are fixed elements of
\(\overline{k}(W_d^+)\), and hence belong to \(k(W_d^+)\).

Because \(k=\F_q\) is finite of characteristic two, the map
\(z\mapsto z^4\) is a field automorphism: it is the second iterate of
Frobenius, or equivalently it is injective and therefore bijective on the
finite set \(k\).  Thus the equation
\(\tau_K^4=\Sigma_K\) has a unique solution \(\tau_K\in k\).  It follows
that
\[
       d'=d+\tau_K^2+\tau_K\in k,
       \qquad x'=\widehat x,\qquad
       y'=\widehat y+\tau_K^2
\]
are defined over \(k\).  The inverse binary dictionary
\[
       u'=\frac{d'}{x'+d'},\qquad v'=\frac{y'}{x'}
\]
also has coefficients in \(k\).

In either characteristic the resulting rational map is defined over
\(k\) on a dense open subset.  Over \(\overline{k}\), the preceding
product or V\'elu construction extends it to the quotient morphism by
\(K\).  Since a rational map from a smooth projective curve to a proper
curve has a unique extension across the missing points, this extension is
fixed by \(\Gamma\) and therefore descends to \(k\).  Its kernel, computed
after base change to \(\overline{k}\), is still exactly \(K\).  This proves
the asserted descent of both the isogeny and its target parameter.
\end{proof}

\section{Characteristic-three separability boundary}

In characteristic three, Theorem~\ref{thm:odd-velu-thesis} applies to a
finite separable odd-order kernel.  A rational cyclic subgroup scheme
isomorphic to \(\mathbb Z/3\mathbb Z\) is admissible.  The full \(E[3]\)
of an ordinary curve has a non-\'etale component and cannot be represented by
a list of nine geometric kernel points.  An inseparable isogeny must instead
be handled through Frobenius and Verschiebung.  This is why the
characteristic-three tripling formula factors through Frobenius but the
ordinary V\'elu point sum does not describe the whole \([3]\)-map.

\begin{proposition}[Relative Frobenius in native coordinates]
\label{prop:Cd-relative-Frobenius}
Let \(k\) have characteristic \(p>0\).  The map
\begin{equation}\label{eq:Cd-relative-Frobenius}
 F_{p,d}:\mathcal C_d\longrightarrow\mathcal C_{d^p},
 \qquad (u,v)\longmapsto(u^p,v^p)
\end{equation}
extends to the smooth completions and is purely inseparable of degree
\(p\).  Its dual Verschiebung
\(V_{p,d}:\mathcal C_{d^p}\to\mathcal C_d\) satisfies
\[
       V_{p,d}\circ F_{p,d}=[p]_{\mathcal C_d},\qquad
       F_{p,d}\circ V_{p,d}=[p]_{\mathcal C_{d^p}}.
\]
If \(d\in\F_p\), then \(F_{p,d}\) is a native endomorphism of
\(\mathcal C_d\).
\end{proposition}

\begin{proof}
We verify successively the target equation, extension to the completion,
degree, group-homomorphism property, and the two duality identities.

In characteristic \(p\), the binomial identity
\((a+b)^p=a^p+b^p\) gives
\[
\begin{aligned}
 (u^{2p}+u^p)(v^{2p}+v^p)
  &=\bigl((u^2+u)(v^2+v)\bigr)^p\\
  &=d^p.
\end{aligned}
\]
Thus \((u^p,v^p)\) satisfies the equation of \(\mathcal C_{d^p}\).
On the projective completion
\eqref{eq:homogeneous-thesis}, the same map is
\[
 ((U_0:U_1),(V_0:V_1))
 \longmapsto
 ((U_0^p:U_1^p),(V_0^p:V_1^p)).
\]
Neither pair of output coordinates vanishes simultaneously, so this is a
morphism on all of \(\PP^1\times\PP^1\).  Raising the homogeneous equation
to the \(p\)th power shows that it restricts to the stated morphism of
smooth completions.  It also sends
\(O=((1:0),(0:1))\) to the identity on the target.

Let \(L=k(\mathcal C_d)=k(u,v)\) and let \(L_0\) be the image of
\(k(\mathcal C_{d^p})\) under pullback.  Then
\[
                    L_0=k(u^p,v^p).
\]
Every element of \(L\) has its \(p\)th power in \(L_0\), so
\(L/L_0\) is purely inseparable of exponent one.  Because
\(\mathcal C_d\) is a smooth geometrically integral curve, its function
field is separably generated of transcendence degree one over \(k\).  The
one-element \(p\)-basis theorem therefore gives
\[
             [L:kL^p]=p.
\]
Here \(kL^p=k(u^p,v^p)=L_0\), and hence
\([\!L:L_0\!]=p\).  Equivalently, the displayed morphism is finite,
radicial, and purely inseparable of degree \(p\).

The relative Frobenius of a group scheme is compatible with multiplication
and inversion.  In the present curve language one can see the same fact
from the rigidity lemma: a morphism between elliptic curves that sends the
identity to the identity and is the relative Frobenius is a group
homomorphism.  Thus \(F_{p,d}\) is an isogeny.  Every isogeny
\(\psi:E\to E'\) of degree \(n\) has a unique dual isogeny
\(\widehat\psi:E'\to E\) characterized by
\[
       \widehat\psi\circ\psi=[n]_E,
       \qquad
       \psi\circ\widehat\psi=[n]_{E'}.
\]
Applying this theorem to \(F_{p,d}\), whose degree is \(p\), defines
\(V_{p,d}=\widehat{F}_{p,d}\) and yields exactly the two displayed
composition identities.

Finally, if \(d\in\F_p\), then \(d^p=d\).  The source and target are
therefore the same marked equation, and the coordinate polynomials
\(u^p,v^p\) are already defined over \(\F_p\).  Hence \(F_{p,d}\) is an
endomorphism written directly in the native \(\mathcal C_d\) coordinates.
\end{proof}

\section{Composition and dual isogenies}

\begin{proposition}[Composition law for the target parameter]
\label{prop:Cd-isogeny-composition-parameter}
In odd characteristic, suppose a degree-\(\ell_1\) product isogeny with
half-kernel product \(R_1\) sends \(\mathcal C_d\) to
\(\mathcal C_{d_1}\), and a degree-\(\ell_2\) product isogeny with
half-kernel product \(R_2\) sends \(\mathcal C_{d_1}\) to
\(\mathcal C_{d_2}\).  If \(\rho_i=1-16d_i\), then
\begin{equation}\label{eq:Cd-isogeny-composition-parameter}
 \rho_2
 =\frac{\rho_1^{\ell_2}}{R_2^8}
 =\frac{\rho^{\ell_1\ell_2}}
        {R_1^{8\ell_2}R_2^8}.
\end{equation}
The composite is a model-preserving isogeny of degree
\(\ell_1\ell_2\), and its kernel fits into
\[
 0\longrightarrow K_1\longrightarrow
 \ker(\varphi_2\circ\varphi_1)
 \longrightarrow K_2\longrightarrow0.
\]
\end{proposition}

\begin{proof}
Write \(\varphi_i\) for the two product isogenies and
\(K_i=\ker\varphi_i\).  Applying
\eqref{eq:odd-velu-parameter-thesis} to \(\varphi_1\) gives
\[
                    \rho_1=\frac{\rho^{\ell_1}}{R_1^8}.
\]
The same formula on the second source curve \(\mathcal C_{d_1}\) gives
\[
                    \rho_2=\frac{\rho_1^{\ell_2}}{R_2^8}.
\]
Substituting the first equality into the second, and keeping track of the
power applied to the first denominator, yields
\[
 \rho_2
 =\frac{(\rho^{\ell_1}/R_1^8)^{\ell_2}}{R_2^8}
 =\frac{\rho^{\ell_1\ell_2}}
        {R_1^{8\ell_2}R_2^8},
\]
which proves \eqref{eq:Cd-isogeny-composition-parameter}.  In particular
the source and output of the composite are again written in the
\(\mathcal C_d\) family.

We next prove the kernel statement rather than merely count degrees.  Over
an algebraic closure define
\[
 \theta:\ker(\varphi_2\circ\varphi_1)\longrightarrow K_2,
 \qquad P\longmapsto\varphi_1(P).
\]
This is well defined because
\(\varphi_2(\varphi_1(P))=O\).  Its kernel consists of those \(P\) for
which \(\varphi_1(P)=O\), and is therefore exactly \(K_1\).
It is also surjective.  Indeed, if \(Q\in K_2\), surjectivity of the
isogeny \(\varphi_1\) supplies a point \(P\) with
\(\varphi_1(P)=Q\); then
\(\varphi_2(\varphi_1(P))=\varphi_2(Q)=O\), so
\(P\in\ker(\varphi_2\circ\varphi_1)\).  We have proved the exact sequence
\[
 0\longrightarrow K_1\longrightarrow
 \ker(\varphi_2\circ\varphi_1)
 \xrightarrow{\,\theta\,}K_2\longrightarrow0.
\]
Because the product isogenies are separable, these are finite \'etale
kernel group schemes; the geometric-point argument is therefore equivalent
to exactness of the group schemes.

Finally, degrees of finite morphisms multiply.  Equivalently, the exact
sequence gives
\[
 \#\ker(\varphi_2\circ\varphi_1)
       =\#K_1\,\#K_2=\ell_1\ell_2.
\]
The composite is separable, so its degree equals this kernel cardinality.
Hence
\(\deg(\varphi_2\circ\varphi_1)=\ell_1\ell_2\), as asserted.
\end{proof}

\begin{proposition}[Dual isogeny and return to the family]
\label{prop:Cd-dual-isogeny}
Assume \(\charac k\ne2\).  Let \(\ell\) be odd and invertible in \(k\),
and let
\(\varphi:\mathcal C_d\to\mathcal C_{d'}\) be a separable
model-preserving \(\ell\)-isogeny.  There is a unique dual isogeny
\[
 \widehat\varphi:\mathcal C_{d'}\longrightarrow\mathcal C_d,
 \qquad
 \widehat\varphi\circ\varphi=[\ell],\quad
 \varphi\circ\widehat\varphi=[\ell].
\]
Its kernel is
\begin{equation}\label{eq:Cd-dual-kernel}
       \ker\widehat\varphi=\varphi(\mathcal C_d[\ell]),
\end{equation}
and has order \(\ell\).  Hence it can be evaluated by the same odd-degree
product formulas, followed, if necessary, by the family isomorphism that
selects the original representative \(\mathcal C_d\).  If that
normalization makes the product target parameter exactly \(\rho=1-16d\),
then its half-kernel product \(R_{\widehat K}\) satisfies
\begin{equation}\label{eq:Cd-dual-parameter-check}
       R_{\widehat K}^{\,8}=\frac{{\rho'}^{\ell}}{\rho}.
\end{equation}
\end{proposition}

\begin{proof}
Put \(E=\mathcal C_d\), \(E'=\mathcal C_{d'}\), and
\(K=\ker\varphi\).  The dual-isogeny theorem gives a unique isogeny
\(\widehat\varphi:E'\to E\) for which
\[
 \widehat\varphi\circ\varphi=[\ell]_E,
 \qquad
 \varphi\circ\widehat\varphi=[\ell]_{E'}.
\]
Taking degrees in the first identity gives
\[
 \deg\widehat\varphi\cdot\deg\varphi
       =\deg[\ell]=\ell^2.
\]
Since \(\deg\varphi=\ell\), one obtains
\(\deg\widehat\varphi=\ell\).  The assumption that \(\ell\) is
invertible in \(k\) implies
that \([\ell]\), \(\varphi\), and \(\widehat\varphi\) are separable;
in particular their kernels are finite \'etale group schemes.

We now identify the dual kernel.  For every
\(P\in E[\ell]\),
\[
       \widehat\varphi(\varphi(P))=[\ell]P=O,
\]
and hence
\[
                  \varphi(E[\ell])\subseteq\ker\widehat\varphi.
\]
Because \(\ell\) is invertible in \(k\), the group
\(E[\ell](\overline{k})\) has
\(\ell^2\) points.  The subgroup \(K\) has order \(\ell\), and every
element of \(K\) is killed by \(\ell\); thus \(K\subseteq E[\ell]\).
The kernel of the restricted homomorphism
\[
       \varphi|_{E[\ell]}:E[\ell]\longrightarrow E'[\ell]
\]
is therefore exactly \(K\).  The first isomorphism theorem gives
\[
       \#\varphi(E[\ell])
       =\frac{\#E[\ell]}{\#K}
       =\frac{\ell^2}{\ell}=\ell.
\]
On the other hand, separability of \(\widehat\varphi\) gives
\(\#\ker\widehat\varphi=\deg\widehat\varphi=\ell\).  The inclusion of
two finite subgroups of the same order is equality, proving
\eqref{eq:Cd-dual-kernel}.

Since \(\ell\) is odd, the kernel
\(\widehat K=\ker\widehat\varphi\) is eligible for the product formulas
of Theorem~\ref{thm:odd-velu-thesis}.  Applying them with source parameter
\(\rho'=1-16d'\) gives a quotient
\[
 \psi:E'\longrightarrow\mathcal C_{d''},
 \qquad
 \rho''=1-16d''=\frac{{\rho'}^\ell}
                         {R_{\widehat K}^{\,8}}.
\]
The maps \(\psi\) and \(\widehat\varphi\) have the same kernel.  By the
uniqueness of the quotient by a finite subgroup, there is a unique
identity-preserving isomorphism
\(\iota:\mathcal C_{d''}\to\mathcal C_d\) such that
\[
                     \widehat\varphi=\iota\circ\psi.
\]
This is the family isomorphism mentioned in the statement.  If the
normalization is chosen so that the product codomain itself is the original
representative \(\mathcal C_d\), then \(\rho''=\rho=1-16d\).  Substitution
in the preceding parameter formula yields
\[
       \rho=\frac{{\rho'}^\ell}{R_{\widehat K}^{\,8}},
       \qquad\text{hence}\qquad
       R_{\widehat K}^{\,8}=\frac{{\rho'}^\ell}{\rho},
\]
which is \eqref{eq:Cd-dual-parameter-check}.
\end{proof}

\section{Examples}

\begin{example}[Odd three-isogeny over \(\F_{101}\)]
Let \(d=1\), so \(\rho=86\).  A point
\(P=(u,v)=(6,42)\) maps to Edwards \((82,70)\) and has order \(96\).
Let \(Q=[32]P\), whose Edwards coordinates are
\((\alpha,\beta)=(61,56)\).  Thus \(Q\) has order \(3\), and
\[
r_1=\beta^{-1}=92,\qquad s_1=\alpha^{-1}=53.
\]
Formula~\eqref{eq:odd-velu-parameter-thesis} gives
\[
\rho'=\frac{86^3}{92^8}=10,\qquad d'=31.
\]
Evaluating the product isogeny at \(P\) gives
\[
             \varphi(P)=(u',v')=(47,37),
\]
and direct substitution verifies
\((47^2+47)(37^2+37)=31\) in \(\F_{101}\).
\end{example}

\begin{example}[Binary three-isogeny over \(\F_{2^8}\)]
\label{ex:binary-three-isogeny-F256}
Let
\[
\F_{2^8}=\F_2[\alpha]/(\alpha^8+\alpha^4+\alpha^3+\alpha+1).
\]
For this example, encode a hexadecimal byte \(\mathtt{0xHH}\) by
\[
 \langle\mathtt{0xHH}\rangle_\alpha
 =\sum_{i=0}^7h_i\alpha^i
 \quad\text{when}\quad
 \mathtt{0xHH}=\sum_{i=0}^7h_i2^i,\qquad h_i\in\{0,1\}.
\]
For instance,
\[
\begin{aligned}
 \hexalpha{DE}&=\alpha^7+\alpha^6+\alpha^4+\alpha^3+\alpha^2+\alpha,\\
 \hexalpha{97}&=\alpha^7+\alpha^4+\alpha^2+\alpha+1,\\
 \hexalpha{05}&=\alpha^2+1.
\end{aligned}
\]
In particular, every two-digit expression below denotes the corresponding
polynomial-basis field element.  Set
\[
             d=\hexalpha{13}=\alpha^4+\alpha+1.
\]
On the explicitly stated auxiliary curve
\[
       W_d^+:\qquad Y^2+XY=X^3+d^2X,
\]
the binary dictionary
\[
 X=d\frac{u+1}{u},\qquad Y=Xv
\]
sends the native point
\[
 P_{\mathcal C}=(u,v)=(\hexalpha{46},\hexalpha{5C})
\]
to
\[
 P=(X,Y)=(\hexalpha{DE},\hexalpha{97}).
\]
Repeated use of the group law gives
\[
\begin{aligned}
 [24]P&=(\hexalpha{49},\hexalpha{CC}),&
 [88]P&=(\hexalpha{05},\hexalpha{FD}),\\
 [132]P&=(\hexalpha{00},\hexalpha{00}),&
 [264]P&=O.
\end{aligned}
\]
Since \(264=2^3\cdot3\cdot11\) and the three points with multipliers
\(264/2\), \(264/3\), and \(264/11\) are nonzero, it follows that
\(\ord(P)=264\).  Let
\[
 Q=[88]P=(\hexalpha{05},\hexalpha{FD}),
\]
where the scalar \(88\) is written in decimal.  Hence \(Q\ne O\) and
\([3]Q=[264]P=O\), so \(Q\) has order three.  The two displayed points
are readily checked on \(W_d^+\):
\[
\begin{aligned}
 Y(P)^2+X(P)Y(P)
  &=X(P)^3+d^2X(P)=\hexalpha{B9},\\
 Y(Q)^2+X(Q)Y(Q)
  &=X(Q)^3+d^2X(Q)=\hexalpha{33}.
\end{aligned}
\]
The unique \(\tau\in\F_{2^8}\) satisfying
\(\tau^4=X(Q)=\hexalpha{05}\) is
\(\tau=\hexalpha{FB}\).  Therefore
\[
d'=d+\tau^2+\tau
   =\hexalpha{13}+\hexalpha{FB}^{\,2}+\hexalpha{FB}
   =\hexalpha{EB}.
\]
The three-isogeny of Corollary~\ref{cor:binary-3iso-thesis} gives
\[
 (x',y')=(\hexalpha{D4},\hexalpha{55})\in W_{d'}^+.
\]
Indeed,
\[
 y'^2+x'y'=x'^3+{d'}^2x'=\hexalpha{76}.
\]
Applying the inverse binary dictionary gives
\[
 u'=\frac{d'}{x'+d'}=\hexalpha{C0},\qquad
 v'=\frac{y'}{x'}=\hexalpha{A5}.
\]
Consequently
\[
 (u'^2+u')(v'^2+v')=\hexalpha{EB}=d',
\]
which verifies the return to \(\mathcal C_{d'}\) without any ambiguity in
the coordinate representation.
\end{example}

\begin{remark}
The model-preserving statements above apply to odd-order separable kernels.
For an even-order kernel, its intersection with the marked four-torsion
subgroup enters the construction; an inseparable kernel is represented by
its finite group scheme.  These data determine the appropriate affine or
projective isogeny formulation in each characteristic.
\end{remark}

\chapter[Tate and Weil Pairings]
{Tate and Weil Pairings on \(\mathcal C_d\)}
\label{ch:pairings}
Throughout this chapter, assume
\[
                         d(1-16d)\ne0,
\]
so that the completed curve \(\mathcal C_d\) is smooth.  In
characteristic two this condition reduces to \(d\ne0\).
The chapter follows the dependency order required by a Miller proof.  It
first determines divisors of native coordinate functions, then constructs
odd- and binary-characteristic line factors, proves the accumulator
recursion, and only afterward packages the normalized Tate and Weil
pairings.  Sparse evaluation, denominator elimination, simultaneous
pairings, and the numerical example are consequences of that divisor
foundation rather than independent formula lists.

\section[Coordinate divisors and Miller functions]
{Coordinate divisors and Miller functions on \(\mathcal C_d\)}

All functions in this chapter belong to the function field of the smooth
completion of
\[
           \mathcal C_d:\quad (u^2+u)(v^2+v)=d.
\]
Its identity and inverse are
\[
        O=(0,\infty),\qquad -(u,v)=(u,-v-1).
\]
To avoid collisions with the running points used in Miller steps, denote
the remaining boundary points in this chapter by
\[
 B_2=(-1,\infty),\qquad
 B_4^+=(\infty,0),\qquad B_4^-=(\infty,-1).
\]
In the global notation of Chapter~\ref{ch:geometry}, these are
\[
                    B_2=T,\qquad B_4^+=R,\qquad B_4^-=-R.
\]
Thus \(B_2\) has order two and \(B_4^-=-B_4^+\).

\begin{proposition}[Divisors of the basic coordinates]
\label{prop:pairing-coordinate-divisors}
On the smooth completion of \(\mathcal C_d\),
\begin{equation}\label{eq:basic-coordinate-divisors}
\begin{aligned}
\Div(u)&=2(O)-(B_4^+)-(B_4^-),\\
\Div(u+1)&=2(B_2)-(B_4^+)-(B_4^-),\\
\Div(v)&=2(B_4^+)-(O)-(B_2),\\
\Div(v+1)&=2(B_4^-)-(O)-(B_2).
\end{aligned}
\end{equation}
Consequently,
\begin{equation}\label{eq:kummer-coordinate-divisor}
 \Div\!\left(\frac{u+1}{u}\right)=2(B_2)-2(O).
\end{equation}
If \(A=(u_A,v_A)\) is affine, then
\begin{equation}\label{eq:vertical-u-divisors}
\begin{aligned}
\Div(u-u_A)
  &=(A)+(-A)-(B_4^+)-(B_4^-),\\
\Div\!\left(\frac{u-u_A}{u}\right)
  &=(A)+(-A)-2(O).
\end{aligned}
\end{equation}
\end{proposition}

\begin{proof}
The proof has three stages: compute the four boundary orders in local
parameters, assemble the global coordinate divisors, and then treat a
general fiber of the degree-two function \(u\).

\medskip
\noindent\emph{Boundary orders.}
We compute the orders locally; this also verifies the multiplicities that
are invisible in the set-theoretic fibers.  Recall the bihomogeneous
equation
\[
 U_1(U_1+U_0)V_1(V_1+V_0)=dU_0^2V_0^2,
 \qquad u=\frac{U_1}{U_0},\quad v=\frac{V_1}{V_0}.
\]
At \(O=(0,\infty)\), use the chart \(U_0=V_1=1\) and put
\(w=V_0/V_1=1/v\).  The local equation is
\begin{equation}\label{eq:pairing-local-O}
             u(u+1)(1+w)=dw^2.
\end{equation}
The partial derivative with respect to \(u\) is \(1\) at \((u,w)=(0,0)\).
Hence \(w\) is a uniformizer at \(O\), and
\eqref{eq:pairing-local-O} gives
\[
                \ord_O(u)=2,\qquad \ord_O(v)=-1.
\]
More precisely, the leading term is \(u=dw^2+\cdots\), so no cancellation
can increase the order.

At \(B_2=(-1,\infty)\), in the same chart put \(x=u+1\).  The equation is
\[
             (-1+x)x(1+w)=dw^2.
\]
Its derivative with respect to \(x\) at \((0,0)\) is \(-1\), which is
nonzero in every characteristic.  Thus \(w\) is again a uniformizer and
\[
             \ord_{B_2}(u+1)=2,\qquad \ord_{B_2}(v)=-1.
\]

At \(B_4^+=(\infty,0)\), use \(U_1=V_0=1\), set
\(z=U_0/U_1=1/u\), and retain \(v=V_1/V_0\).  The local equation becomes
\[
                 (1+z)v(v+1)=dz^2.
\]
Its derivative with respect to \(v\) at \((z,v)=(0,0)\) is \(1\).
Therefore \(z\) is a uniformizer and
\[
             \ord_{B_4^+}(u)=-1,\qquad \ord_{B_4^+}(v)=2.
\]
Finally, at \(B_4^-=(\infty,-1)\), put \(y=v+1\) in the same chart.
Then
\[
                 (1+z)(-1+y)y=dz^2,
\]
so \(z\) is a uniformizer and
\[
             \ord_{B_4^-}(u)=-1,\qquad
             \ord_{B_4^-}(v+1)=2.
\]

\medskip
\noindent\emph{Global coordinate divisors.}
There are no further zeros or poles.  Indeed, setting \(u=0\) in the
homogeneous equation forces \(V_0=0\), because \(d\ne0\); hence the
zero fiber of \(u\) is supported only at \(O\).  Setting
\(U_1+U_0=0\), equivalently \(u=-1\), forces \(V_0=0\) and hence gives
only \(B_2\); setting \(U_0=0\) gives precisely \(B_4^+\) and
\(B_4^-\).
The four local order computations consequently give
\[
\begin{aligned}
\Div(u)&=2(O)-(B_4^+)-(B_4^-),\\
\Div(u+1)&=2(B_2)-(B_4^+)-(B_4^-).
\end{aligned}
\]
The coordinate interchange \(\sigma(u,v)=(v,u)\) is an automorphism
of the smooth completion.  It sends
\[
  O\longleftrightarrow B_4^+,
  \qquad
  B_2\longleftrightarrow B_4^-,
\]
and satisfies \(\sigma^*u=v\) and \(\sigma^*(u+1)=v+1\).
Pulling back the two divisors already computed therefore gives
\[
\begin{aligned}
  \Div(v)&=2(B_4^+)-(O)-(B_2),\\
  \Div(v+1)&=2(B_4^-)-(O)-(B_2),
\end{aligned}
\]
which are the remaining two formulas.
Subtracting \(\Div(u)\) from \(\Div(u+1)\) proves
\eqref{eq:kummer-coordinate-divisor}.

\medskip
\noindent\emph{A general \(u\)-fiber.}
Let now \(A=(u_A,v_A)\) be affine.  Since \(d\ne0\), one has
\(u_A\ne0,-1\).  The fiber above \(u_A\) is cut out by
\[
       v^2+v=\frac{d}{u_A(u_A+1)}.
\]
Its two points, counted with scheme-theoretic multiplicity, are
\[
       A=(u_A,v_A),\qquad
       -A=(u_A,-v_A-1).
\]
If the two points coincide, the notation \((A)+(-A)\) means \(2(A)\);
the degree-two projection still supplies the correct multiplicity.
The poles of \(u-u_A\) are the same simple poles as those of \(u\).
Therefore
\[
 \Div(u-u_A)=(A)+(-A)-(B_4^+)-(B_4^-).
\]
Subtracting the already established divisor of \(u\) gives
\[
 \Div\!\left(\frac{u-u_A}{u}\right)
 =(A)+(-A)-2(O),
\]
which completes the proof.
\end{proof}

\begin{corollary}[Divisors of the pairing coordinates]
\label{cor:pairing-special-coordinate-divisors}
Assume first that \(\charac k\ne2\), put
\(\rho=1-16d\), and work over a field containing \(\sqrt\rho\).  Define
\[
\begin{aligned}
 P_u^\pm&=\left(-\frac12,\frac{-1\pm\sqrt\rho}{2}\right),\\
 P_v^\pm&=\left(\frac{-1\pm\sqrt\rho}{2},-\frac12\right).
\end{aligned}
\]
Then
\begin{equation}\label{eq:centered-coordinate-divisors}
\begin{aligned}
\Div(2u+1)&=(P_u^+)+(P_u^-)-(B_4^+)-(B_4^-),\\
\Div(2v+1)&=(P_v^+)+(P_v^-)-(O)-(B_2),\\
\Div\!\left((2u+1)^{-1}\right)
  &=(B_4^+)+(B_4^-)-(P_u^+)-(P_u^-),\\
\Div\!\left((2v+1)^{-1}\right)
  &=(O)+(B_2)-(P_v^+)-(P_v^-).
\end{aligned}
\end{equation}
Put \(\beta=(16d)^{-1}\).  For the odd-characteristic functions
\[
 x_d^-=\beta\frac{u+1}{u},\qquad
 y_d^-=2\beta^2(2v+1)\frac{u+1}{u},
\]
which will be used for line evaluation below,
\begin{equation}\label{eq:odd-pairing-coordinate-divisors}
\Div(x_d^-)=2(B_2)-2(O),\qquad
\Div(y_d^-)=(P_v^+)+(P_v^-)+(B_2)-3(O).
\end{equation}

In characteristic two, put
\[
 x_d^+=d\frac{u+1}{u},\qquad
 y_d^+=dv\frac{u+1}{u}.
\]
Then
\begin{equation}\label{eq:binary-pairing-coordinate-divisors}
\begin{aligned}
\Div(x_d^+)&=2(B_2)-2(O),\\
\Div(y_d^+)&=2(B_4^+)+(B_2)-3(O),\\
\Div(x_d^++y_d^+)&=2(B_4^-)+(B_2)-3(O).
\end{aligned}
\end{equation}
\end{corollary}

\begin{proof}
The odd-characteristic centered coordinates, the odd pairing coordinates,
and the binary pairing coordinates are treated in that order.

\medskip
\noindent\emph{Centered coordinates in odd characteristic.}
Assume first that \(2\ne0\).  At a zero of \(2u+1\) one has
\(u=-1/2\), and therefore
\[
             u^2+u=-\frac14.
\]
Substitution in the equation of \(\mathcal C_d\) gives
\[
             v^2+v+4d=0.
\]
This quadratic has discriminant
\(\rho=1-16d\), so its two roots over \(k(\sqrt\rho)\) are
\[
             v=\frac{-1\pm\sqrt\rho}{2}.
\]
The smoothness assumption gives \(\rho\ne0\); hence the two roots are
distinct, and both zeros of \(2u+1\) are simple.  They are precisely
\(P_u^+\) and \(P_u^-\).  Multiplication by the nonzero constant \(2\)
does not change a divisor, and
Proposition~\ref{prop:pairing-coordinate-divisors} shows that \(u\), hence also
\(2u+1\), has one simple pole at each of \(B_4^+\) and \(B_4^-\).
Thus
\[
 \Div(2u+1)=(P_u^+)+(P_u^-)-(B_4^+)-(B_4^-).
\]

At a zero of \(2v+1\), substitution of \(v=-1/2\) gives
\(u^2+u+4d=0\), whose two distinct roots are the points
\(P_v^+\) and \(P_v^-\).  Proposition~\ref{prop:pairing-coordinate-divisors}
gives one simple pole of \(v\), and hence of \(2v+1\), at each of
\(O\) and \(B_2\).  Therefore
\[
 \Div(2v+1)=(P_v^+)+(P_v^-)-(O)-(B_2).
\]
Negating these two divisors proves the reciprocal identities in
\eqref{eq:centered-coordinate-divisors}.

\medskip
\noindent\emph{Odd-characteristic pairing coordinates.}
Constants have trivial divisor.
Proposition~\ref{prop:pairing-coordinate-divisors} gives this divisor as
\eqref{eq:kummer-coordinate-divisor}.
Consequently
\[
\begin{aligned}
\Div(x_d^-)
 &=2(B_2)-2(O),\\
\Div(y_d^-)
 &=\Div(2v+1)+2(B_2)-2(O)\\
 &=(P_v^+)+(P_v^-)+(B_2)-3(O),
\end{aligned}
\]
which proves \eqref{eq:odd-pairing-coordinate-divisors}.

\medskip
\noindent\emph{Binary pairing coordinates.}
Suppose now that \(\charac k=2\).  Again constants \(d\ne0\) have
trivial divisor.  Using the four identities of
\eqref{eq:basic-coordinate-divisors} and
\eqref{eq:kummer-coordinate-divisor}, we obtain
\[
\begin{aligned}
\Div(x_d^+)
 &=2(B_2)-2(O),\\
\Div(y_d^+)
 &=\Div(v)+\Div\!\left(\frac{u+1}{u}\right)\\
 &=\{2(B_4^+)-(O)-(B_2)\}
   +\{2(B_2)-2(O)\}\\
 &=2(B_4^+)+(B_2)-3(O).
\end{aligned}
\]
Since subtraction and addition coincide in characteristic two,
\[
 x_d^++y_d^+
 =d(v+1)\frac{u+1}{u}.
\]
Replacing \(\Div(v)\) by
\(\Div(v+1)=2(B_4^-)-(O)-(B_2)\) in the preceding calculation gives
\[
 \Div(x_d^++y_d^+)=2(B_4^-)+(B_2)-3(O).
\]
This proves every identity in
\eqref{eq:binary-pairing-coordinate-divisors}.
\end{proof}

Because the curve is fixed throughout this chapter, no model superscript is
needed on a Miller function.  For \(P\in\mathcal C_d\), choose \(f_{n,P}\)
with
\[
 \Div(f_{n,P})
 =n(P)-([n]P)-(n-1)(O).
\]
For each admissible pair of running points \(R,S\in\mathcal C_d\), write
\(g_{R,S}\) for a rational function with
\begin{equation}\label{eq:native-Miller-factor-divisor}
 \Div(g_{R,S})
 =(R)+(S)-(R+S)-(O).
\end{equation}

\begin{proposition}[Miller divisor]\label{prop:Miller-divisor-thesis}
If \eqref{eq:native-Miller-factor-divisor} holds, then, up to a nonzero
normalization constant,
\begin{equation}\label{eq:Miller-recursion-thesis}
 f_{m+n,P}=f_{m,P}f_{n,P}g_{[m]P,[n]P}.
\end{equation}
\end{proposition}

\begin{proof}
Set
\[
       H=f_{m,P}f_{n,P}g_{[m]P,[n]P}.
\]
The divisor map is additive under multiplication, so the hypotheses give
\[
\begin{aligned}
 \Div(H)
 &=m(P)-([m]P)-(m-1)(O)\\
 &\quad+n(P)-([n]P)-(n-1)(O)\\
 &\quad+([m]P)+([n]P)-([m+n]P)-(O).
\end{aligned}
\]
The terms \(([m]P)\) and \(([n]P)\) cancel.  The coefficient of \(O\)
is
\[
       -(m-1)-(n-1)-1=-(m+n-1).
\]
Therefore
\[
 \Div(H)
 =(m+n)(P)-([m+n]P)-(m+n-1)(O)
 =\Div(f_{m+n,P}).
\]
This calculation takes place in the free abelian group of divisors, so it
remains valid when some of the displayed points coincide or when
\(m=n\).

It follows that \(H/f_{m+n,P}\) has divisor zero.  A rational function
with zero divisor on a smooth projective geometrically integral curve is
both regular and invertible everywhere, and hence is a nonzero constant.
Thus
\[
        f_{m+n,P}=c\,f_{m,P}f_{n,P}g_{[m]P,[n]P}
\]
for some \(c\in\overline{k}^{\times}\).  Choosing the prescribed Miller
normalization absorbs \(c\), proving
\eqref{eq:Miller-recursion-thesis} up to the stated constant.
\end{proof}

\section[Odd-characteristic line factors in uv]
{Odd-characteristic line factors in \(u,v\)}

In this section
\[
 \mathcal C_d:\quad (u^2+u)(v^2+v)=d,\qquad
 \charac k\ne2,\qquad d(1-16d)\ne0.
\]
Put
\[
 \beta=\frac1{16d},\qquad A=\frac1{4d}-2,\qquad s(P)=2v_P+1,
\]
and define two rational functions on \(\mathcal C_d\):
\begin{equation}\label{eq:pairing-odd-native-dictionary}
\boxed{\quad
 x_d^-(u,v)=\beta\frac{u+1}{u},\qquad
 y_d^-(u,v)=2\beta^2(2v+1)\frac{u+1}{u}.
\quad}
\end{equation}
These are rational functions of \(u,v\); they do not replace the running
\(\mathcal C_d\)-point by a point on another declared model.  They satisfy
the auxiliary relation
\begin{equation}\label{eq:pairing-Wminus}
 W_d^-:\qquad
 Y^2=X^3+\beta AX^2+\beta^2X.
\end{equation}

Indeed, if \(U=(u+1)/u\) and \(V=2(2v+1)U\), the defining equation of
\(\mathcal C_d\) gives
\[
       \beta V^2=U^3+AU^2+U.
\]
Multiplication by \(\beta^3\) proves directly that
\((x_d^-,y_d^-)\) satisfies \eqref{eq:pairing-Wminus}.

Let
\[
 R=(u_R,v_R),\quad S=(u_S,v_S),\quad
 T=R+S=(u_T,v_T),\quad Q=(u_Q,v_Q)
\]
be points of \(\mathcal C_d\).  The secant or tangent coefficient is
\begin{equation}\label{eq:odd-native-slope}
\lambda_{R,S}^-=
\begin{cases}
\displaystyle
2\beta\,
\frac{s(S)(u_S+1)u_R-s(R)(u_R+1)u_S}{u_R-u_S},
&S\ne R,-R,\\[4mm]
\displaystyle
\frac{3(u_R+1)^2+2Au_R(u_R+1)+u_R^2}
     {4s(R)u_R(u_R+1)},&S=R.
\end{cases}
\end{equation}

When \(S=-R\), the secant is vertical; define its evaluation factor before
it is used in the domain statement:
\begin{equation}\label{eq:odd-native-vertical}
 g_{R,-R}^{-}(Q)
 =x_d^-(Q)-x_d^-(R)
 =\beta\frac{u_R-u_Q}{u_Qu_R}.
\end{equation}

\begin{proposition}[Derivation and domain of the odd slope]
\label{prop:odd-native-slope-derivation}
Equation~\eqref{eq:odd-native-slope} is the secant or tangent coefficient
of the functions \((x_d^-,y_d^-)\) on \(\mathcal C_d\).  The first branch
applies when \(R,S\) are distinct and not negatives.  If \(S=-R\), the
line is vertical and \eqref{eq:odd-native-vertical} is used.  If an affine
tangent denominator vanishes at a boundary or torsion point, its homogeneous
tangent specialization gives the projective value.
\end{proposition}

\begin{proof}
The argument proceeds through the ordinary secant, the vertical fiber, the
tangent specialization, and finally the projective boundary interpretation.
Put
\[
             U_P=\frac{u_P+1}{u_P}.
\]
For every affine point of \(\mathcal C_d\), the equality \(d\ne0\)
implies \(u_P\ne0,-1\), so the following denominators are legitimate on
the stated affine domain.

For distinct points \(R,S\) with different cubic abscissas, the secant
coefficient on \eqref{eq:pairing-Wminus} is
\[
 \lambda=\frac{y_d^-(S)-y_d^-(R)}
               {x_d^-(S)-x_d^-(R)}.
\]
The numerator and denominator are, respectively,
\[
\begin{aligned}
 y_d^-(S)-y_d^-(R)
   &=2\beta^2\{s(S)U_S-s(R)U_R\},\\
 x_d^-(S)-x_d^-(R)
   &=\beta(U_S-U_R).
\end{aligned}
\]
Moreover,
\[
 U_S-U_R
 =\left(1+\frac1{u_S}\right)
  -\left(1+\frac1{u_R}\right)
 =\frac{u_R-u_S}{u_Ru_S}.
\]
Substitution and multiplication of numerator and denominator by
\(u_Ru_S\) give
\[
 \lambda
 =2\beta\,
 \frac{s(S)(u_S+1)u_R-s(R)(u_R+1)u_S}{u_R-u_S},
\]
which is the first branch of \eqref{eq:odd-native-slope}.

The equality \(x_d^-(S)=x_d^-(R)\) is equivalent to
\(U_S=U_R\), and hence to \(u_S=u_R\).  The degree-two \(u\)-fiber
consists of \(R\) and
\(-R=(u_R,-v_R-1)\), counted with multiplicity.  Therefore, for distinct
points, the secant denominator vanishes exactly when \(S=-R\).
In centered coordinates negation sends \(s(R)=2v_R+1\) to \(-s(R)\);
thus \(y_d^-(-R)=-y_d^-(R)\) while the abscissa is unchanged.  The
corresponding line is the vertical line \(X=x_d^-(R)\).

For the tangent, implicit differentiation of
\[
       Y^2=X^3+\beta AX^2+\beta^2X
\]
at a point with \(2Y_R\ne0\) gives
\[
 2Y_R\lambda
   =3X_R^2+2\beta AX_R+\beta^2.
\]
Now substitute
\[
 X_R=\beta\frac{u_R+1}{u_R},\qquad
 Y_R=2\beta^2s(R)\frac{u_R+1}{u_R}.
\]
After factoring \(\beta^2/u_R^2\) from the numerator and multiplying
numerator and denominator by \(u_R^2\), one obtains
\[
\begin{aligned}
\lambda
 &=\frac{\beta^2
  \{3(u_R+1)^2+2Au_R(u_R+1)+u_R^2\}/u_R^2}
 {4\beta^2s(R)(u_R+1)/u_R}\\
 &=\frac{3(u_R+1)^2+2Au_R(u_R+1)+u_R^2}
        {4s(R)u_R(u_R+1)}.
\end{aligned}
\]
In the first displayed fraction the factor \(\beta^2\) is common to the
whole bracketed numerator; cancelling it gives the second line.  This is
the tangent branch of \eqref{eq:odd-native-slope}.

For a finite affine point, \(u_R(u_R+1)\ne0\).  Hence the only possible
vanishing factor in the tangent denominator is \(s(R)\), equivalently
\(Y_R=0\).  Such a tangent is vertical and is represented by the vertical
Miller factor rather than by a finite affine slope.  At boundary points
the formulas are interpreted after homogenizing the tangent line.  Since
the affine identities hold on a dense open subset, their homogeneous
specializations are the unique projective extensions.
\end{proof}

Define directly on \(\mathcal C_d\)
\[
\begin{aligned}
\ell_{R,S}^{-}(Q)
 &=y_d^-(Q)-y_d^-(R)
   -\lambda_{R,S}^-\bigl(x_d^-(Q)-x_d^-(R)\bigr),\\
\nu_T^{-}(Q)&=x_d^-(Q)-x_d^-(T).
\end{aligned}
\]
The corresponding Miller factor is
\begin{equation}\label{eq:odd-native-Miller-thesis}
\boxed{
g_{R,S}^{-}(Q)
=\frac{\ell_{R,S}^{-}(Q)}{\nu_T^{-}(Q)}
=\frac{u_T\,\mathcal L_{R,S}^-(Q)}
       {u_R(u_T-u_Q)},
}
\end{equation}
\[
\mathcal L_{R,S}^-(Q)
=2\beta\!\left[
 u_Rs(Q)(u_Q+1)-u_Qs(R)(u_R+1)
 \right]
-\lambda_{R,S}^-(u_R-u_Q).
\]

\begin{proposition}[Odd-characteristic Miller factor]
\label{prop:odd-native-Miller-divisor}
Equations~\eqref{eq:odd-native-Miller-thesis} and
\eqref{eq:odd-native-vertical} satisfy
\eqref{eq:native-Miller-factor-divisor}.
More explicitly, for \(T=R+S\),
\[
\begin{aligned}
\Div(\ell_{R,S}^{-})
  &=(R)+(S)+(-T)-3(O),\\
\Div(\nu_T^{-})
  &=(T)+(-T)-2(O).
\end{aligned}
\]
\end{proposition}

\begin{proof}
The proof first derives the native rational expression, then identifies the
smooth cubic used to read the intersection divisor, and finally treats the
vertical specialization.

\medskip
\noindent\emph{The rational expression.}
We first verify the displayed rational formula and then compute its divisor.
For any affine \(Q\) and \(R\),
\[
\begin{aligned}
 x_d^-(Q)-x_d^-(R)
 &=\beta\left(\frac{u_Q+1}{u_Q}
             -\frac{u_R+1}{u_R}\right)
  =\beta\frac{u_R-u_Q}{u_Qu_R},\\
 x_d^-(Q)-x_d^-(T)
 &=\beta\frac{u_T-u_Q}{u_Qu_T}.
\end{aligned}
\]
Likewise,
\[
\begin{aligned}
y_d^-(Q)-y_d^-(R)
 &=2\beta^2\left\{
     s(Q)\frac{u_Q+1}{u_Q}
    -s(R)\frac{u_R+1}{u_R}\right\}\\
 &=\frac{2\beta^2}{u_Qu_R}
   \{u_Rs(Q)(u_Q+1)-u_Qs(R)(u_R+1)\}.
\end{aligned}
\]
It follows that
\[
\begin{aligned}
\ell_{R,S}^-(Q)
 &=y_d^-(Q)-y_d^-(R)
   -\lambda_{R,S}^-\{x_d^-(Q)-x_d^-(R)\}\\
 &=\frac{\beta}{u_Qu_R}
 \Bigl(
 2\beta\{u_Rs(Q)(u_Q+1)-u_Qs(R)(u_R+1)\}\\
 &\hspace{38mm}
 -\lambda_{R,S}^-(u_R-u_Q)\Bigr)\\
 &=\frac{\beta}{u_Qu_R}\,
       \mathcal L_{R,S}^-(Q).
\end{aligned}
\]
On the other hand,
\[
       \nu_T^-(Q)
       =\frac{\beta(u_T-u_Q)}{u_Qu_T}.
\]
Dividing these identities cancels \(\beta\) and \(u_Q\) and gives
\[
 \frac{\ell_{R,S}^-(Q)}{\nu_T^-(Q)}
 =\frac{u_T\mathcal L_{R,S}^-(Q)}
        {u_R(u_T-u_Q)},
\]
which proves \eqref{eq:odd-native-Miller-thesis} as an identity in
\(k(\mathcal C_d)\).

\medskip
\noindent\emph{The divisor calculation.}
We now justify the divisor calculation, including the identification of
the two smooth models.  The rational functions
\[
 X=x_d^-=\beta\frac{u+1}{u},\qquad
 Y=y_d^-=2\beta^2(2v+1)\frac{u+1}{u}
\]
satisfy \eqref{eq:pairing-Wminus}.  On the dense open subset where
\(X(X-\beta)\ne0\), they have the inverse
\[
        u=\frac{\beta}{X-\beta},\qquad
        2v+1=\frac{Y}{2\beta X}.
\]
Thus the map is birational.  Both source and the smooth projective
completion of \eqref{eq:pairing-Wminus} are smooth projective curves.
For the target this follows explicitly from
\[
       \beta\ne0,\qquad
       A^2-4=\frac{1-16d}{16d^2}=\frac{\rho}{16d^2}\ne0:
\]
the cubic polynomial
\(X(X^2+\beta AX+\beta^2)\) has three distinct roots.  Hence the
birational map extends uniquely to an isomorphism.  By
\eqref{eq:odd-pairing-coordinate-divisors}, \(X\) has a pole of order
two and \(Y\) a pole of order three at \(O\); hence \(O\) maps to the
point at infinity of the cubic.  The isomorphism is therefore
identity-preserving.

For a nonvertical secant or tangent, the affine line
\[
       Y-Y_R-\lambda_{R,S}^-(X-X_R)=0
\]
meets the cubic at \(R\), \(S\), and a third point, with tangent
multiplicities when \(R=S\).  By the chord-and-tangent definition of the
group law, the third point is \(-T\), where \(T=R+S\).  The restricted
line function has a pole of order three at the point at infinity, because
its \(Y\)-term has order \(-3\).  B\'ezout's theorem, with intersection
multiplicities, therefore gives
\[
        (R)+(S)+(-T)-3(O),
\]
as the divisor of \(\ell_{R,S}^-\).

The function \(X-X_T\) has the two zeros \(T\) and \(-T\), counted with
multiplicity, because negation on the cubic fixes \(X\).  Since \(X\) has
a double pole at \(O\),
\[
        \Div(\nu_T^-)=(T)+(-T)-2(O).
\]
Subtracting the two divisors yields
\[
\begin{aligned}
\Div(g_{R,S}^-)
 &=\{(R)+(S)+(-T)-3(O)\}\\
 &\quad-\{(T)+(-T)-2(O)\}\\
 &=(R)+(S)-(T)-(O),
\end{aligned}
\]
which is \eqref{eq:native-Miller-factor-divisor}.

\medskip
\noindent\emph{The vertical specialization.}
If \(S=-R\), then \(T=O\) and the secant is vertical.  The correct factor
is simply
\[
 X(Q)-X(R)=\beta\frac{u_R-u_Q}{u_Qu_R},
\]
whose divisor is
\((R)+(-R)-2(O)\), exactly the specialization
\((R)+(S)-(R+S)-(O)\).  All arguments use intersection multiplicities,
so they also cover tangencies and torsion specializations.  Since \(2\)
is invertible in characteristic three, the proof applies there without
change.
\end{proof}

\section[Characteristic-two line factors in uv]
{Characteristic-two line factors in \(u,v\)}

Here
\[
 \mathcal C_d:\quad (u^2+u)(v^2+v)=d,\qquad
 \charac k=2,\qquad d\ne0.
\]
Define the following rational functions on \(\mathcal C_d\):
\begin{equation}\label{eq:pairing-binary-native-dictionary}
\boxed{\quad
 x_d^+(u,v)=d\frac{u+1}{u},\qquad
 y_d^+(u,v)=d\,v\frac{u+1}{u}.
\quad}
\end{equation}
They satisfy the full auxiliary equation
\begin{equation}\label{eq:pairing-Wplus}
 W_d^+:\qquad Y^2+XY=X^3+d^2X.
\end{equation}
To verify this directly in \(k(\mathcal C_d)\), put
\(U=(u+1)/u\).  In characteristic two,
\[
 U+1=\frac1u,\qquad U^3+U=U(U+1)^2=\frac{u+1}{u^3}.
\]
The equation of \(\mathcal C_d\) gives
\[
\begin{aligned}
(y_d^+)^2+x_d^+y_d^+
 &=d^2U^2(v^2+v)
  =d^3\frac{U^2}{u(u+1)}\\
 &=d^3(U^3+U)
  =(x_d^+)^3+d^2x_d^+.
\end{aligned}
\]

Let
\[
 R=(u_R,v_R),\quad S=(u_S,v_S),\quad
 T=R+S=(u_T,v_T),\quad Q=(u_Q,v_Q).
\]
The secant or tangent coefficient in \(u,v\) is
\begin{equation}\label{eq:binary-native-slope}
\lambda_{R,S}^+=
\begin{cases}
\displaystyle
\frac{v_S(u_S+1)u_R-v_R(u_R+1)u_S}{u_R-u_S},
&S\ne R,-R,\\[4mm]
\displaystyle
v_R+d\frac{u_R+1}{u_R}
   +d\frac{u_R}{u_R+1},&S=R.
\end{cases}
\end{equation}

When \(S=-R\), define the vertical factor at once:
\begin{equation}\label{eq:binary-native-vertical}
 g_{R,-R}^{+}(Q)
 =x_d^+(Q)-x_d^+(R)
 =d\frac{u_R-u_Q}{u_Qu_R}.
\end{equation}
All displayed minus signs may of course be evaluated as plus signs in
characteristic two.

\begin{proposition}[Derivation and domain of the binary slope]
\label{prop:binary-native-slope-derivation}
Equation~\eqref{eq:binary-native-slope} is the secant or tangent coefficient
of \((x_d^+,y_d^+)\).  Its vertical branch is exactly
\eqref{eq:binary-native-vertical}; boundary and torsion specializations are
understood on the smooth projective model.
\end{proposition}

\begin{proof}
The proof follows the same four-stage order as its odd-characteristic
counterpart: secant, vertical fiber, tangent, and projective boundary.
For an affine point of \(\mathcal C_d\), the relation \(d\ne0\) forces
\(u\ne0,1\); here \(1=-1\) in characteristic two.  Thus
\[
       U_P=\frac{u_P+1}{u_P},\qquad
       X_P=dU_P,\qquad Y_P=dU_Pv_P
\]
are well defined.

Let first \(R\) and \(S\) be distinct with \(X_R\ne X_S\).  The secant
coefficient is
\[
\begin{aligned}
\lambda
 &=\frac{Y_S-Y_R}{X_S-X_R}\\
 &=\frac{d\{v_SU_S-v_RU_R\}}
         {d(U_S-U_R)}\\
 &=\frac{v_S(u_S+1)u_R-v_R(u_R+1)u_S}
         {u_R-u_S},
\end{aligned}
\]
where in the last step we used
\[
          U_S-U_R=\frac{u_R-u_S}{u_Ru_S}.
\]
This is the first branch of \eqref{eq:binary-native-slope}; its displayed
minus signs may be read as plus signs.

Because \(U=1+1/u\), equality \(X_R=X_S\) is equivalent to
\(u_R=u_S\).  The two points of that degree-two fiber are
\[
       R=(u_R,v_R),\qquad
       -R=(u_R,v_R+1),
\]
since negation in characteristic two sends \(v\) to \(v+1\).
Consequently a secant through distinct points has zero denominator exactly
when \(S=-R\), in which case the line is vertical.

For the tangent, recall that a generalized Weierstrass equation
\[
Y^2+a_1XY+a_3Y=X^3+a_2X^2+a_4X+a_6
\]
has tangent coefficient
\[
 \lambda
 =\frac{3X_R^2+2a_2X_R+a_4-a_1Y_R}
        {2Y_R+a_1X_R+a_3}.
\]
For \(W_d^+:Y^2+XY=X^3+d^2X\), one has
\(a_1=1\), \(a_2=a_3=a_6=0\), and \(a_4=d^2\).  Reducing the preceding
formula in characteristic two gives
\[
 \lambda
 =\frac{X_R^2+d^2+Y_R}{X_R}
 =X_R+\frac{d^2}{X_R}+\frac{Y_R}{X_R}.
\]
Substitution of the native expressions gives
\[
\begin{aligned}
\lambda
 &=d\frac{u_R+1}{u_R}
   +d\frac{u_R}{u_R+1}+v_R,
\end{aligned}
\]
which is the second branch of
\eqref{eq:binary-native-slope}.

For a finite affine \(\mathcal C_d\)-point, \(X_R\ne0\), so the tangent
coefficient is finite.  The denominator \(X_R\) can vanish only at a
projective boundary specialization; there the homogeneous tangent is
vertical and the factor \eqref{eq:binary-native-vertical} supplies the
corresponding Miller function.  The affine identities hold on a dense
open subset and hence extend uniquely to the smooth projective model.
\end{proof}

Define directly on \(\mathcal C_d\)
\[
\begin{aligned}
\ell_{R,S}^{+}(Q)
 &=y_d^+(Q)-y_d^+(R)
   -\lambda_{R,S}^+\bigl(x_d^+(Q)-x_d^+(R)\bigr),\\
\nu_T^{+}(Q)&=x_d^+(Q)-x_d^+(T).
\end{aligned}
\]
The factor used in the Miller recurrence is
\begin{equation}\label{eq:binary-native-Miller-thesis}
\boxed{
g_{R,S}^{+}(Q)
=\frac{\ell_{R,S}^{+}(Q)}{\nu_T^{+}(Q)}
=\frac{u_T\,\mathcal L_{R,S}^+(Q)}
       {u_R(u_T-u_Q)},
}
\end{equation}
\[
\mathcal L_{R,S}^+(Q)
=u_Rv_Q(u_Q+1)-u_Qv_R(u_R+1)
-\lambda_{R,S}^+(u_R-u_Q).
\]

\begin{proposition}[Characteristic-two Miller factor]
\label{prop:binary-native-Miller-divisor}
Equations~\eqref{eq:binary-native-Miller-thesis} and
\eqref{eq:binary-native-vertical} satisfy
\eqref{eq:native-Miller-factor-divisor}.
More explicitly, for \(T=R+S\),
\[
\begin{aligned}
\Div(\ell_{R,S}^{+})
  &=(R)+(S)+(-T)-3(O),\\
\Div(\nu_T^{+})
  &=(T)+(-T)-2(O).
\end{aligned}
\]
\end{proposition}

\begin{proof}
The proof first derives the native rational expression, then identifies the
binary cubic and its line divisor, and finally treats the vertical
specialization.

\medskip
\noindent\emph{The rational expression.}
We first derive the boxed expression directly in the function field.
For affine \(Q,R,T\), the identity
\[
       \frac{u+1}{u}=1+\frac1u
\]
gives
\[
\begin{aligned}
x_d^+(Q)-x_d^+(R)
 &=d\frac{u_R-u_Q}{u_Qu_R},\\
x_d^+(Q)-x_d^+(T)
 &=d\frac{u_T-u_Q}{u_Qu_T}.
\end{aligned}
\]
Moreover,
\[
\begin{aligned}
y_d^+(Q)-y_d^+(R)
 &=d\left\{
      v_Q\frac{u_Q+1}{u_Q}
     -v_R\frac{u_R+1}{u_R}\right\}\\
 &=\frac{d}{u_Qu_R}
   \{u_Rv_Q(u_Q+1)-u_Qv_R(u_R+1)\}.
\end{aligned}
\]
It follows that
\[
\begin{aligned}
\ell_{R,S}^+(Q)
 &=y_d^+(Q)-y_d^+(R)
   -\lambda_{R,S}^+\{x_d^+(Q)-x_d^+(R)\}\\
 &=\frac{d}{u_Qu_R}
 \Bigl(
 u_Rv_Q(u_Q+1)-u_Qv_R(u_R+1)
 -\lambda_{R,S}^+(u_R-u_Q)\Bigr)\\
 &=\frac{d}{u_Qu_R}\mathcal L_{R,S}^+(Q),
\end{aligned}
\]
whereas
\[
       \nu_T^+(Q)=\frac{d(u_T-u_Q)}{u_Qu_T}.
\]
Their quotient is therefore
\[
 \frac{\ell_{R,S}^+(Q)}{\nu_T^+(Q)}
 =\frac{u_T\mathcal L_{R,S}^+(Q)}
        {u_R(u_T-u_Q)},
\]
which is \eqref{eq:binary-native-Miller-thesis}.

\medskip
\noindent\emph{The divisor calculation.}
We next verify the divisor statement.  The functions
\[
 X=x_d^+=d\frac{u+1}{u},\qquad
 Y=y_d^+=dv\frac{u+1}{u}
\]
satisfy \(W_d^+:Y^2+XY=X^3+d^2X\).  Conversely, on the dense open set
\(X(X+d)\ne0\),
\[
              u=\frac{d}{X+d},\qquad v=\frac{Y}{X}.
\]
Thus the map is birational and extends to an isomorphism of the smooth
projective completions.  For completeness, the target is nonsingular:
in the affine chart the simultaneous equations
\[
 X=0,\qquad Y+X^2+d^2=0,\qquad
 Y^2+XY+X^3+d^2X=0
\]
would force \(X=Y=0\) and then \(d=0\), a contradiction; at the
projective point at infinity the derivative with respect to the
homogenizing coordinate equals \(1\).  From
\eqref{eq:binary-pairing-coordinate-divisors}, \(X\) and \(Y\) have
poles of orders two and three, respectively, at \(O\).  Hence \(O\)
maps to the point at infinity and the isomorphism preserves the identity.

On \(W_d^+\), negation is
\[
              (X,Y)\longmapsto(X,X+Y).
\]
A nonvertical secant or tangent line through \(R,S\) meets the cubic a
third time at \(-T\), where \(T=R+S\).  The restricted line function has
a triple pole at \(O\); B\'ezout's theorem, including tangent
multiplicities, gives
\[
       (R)+(S)+(-T)-3(O),
\]
for \(\Div(\ell_{R,S}^+)\).  Since \(X-X_T\) vanishes at \(T\) and
\(-T\), counted with multiplicity, and has a double pole at \(O\),
\[
       (T)+(-T)-2(O).
\]
is the divisor of \(\nu_T^+\).  Subtraction gives
\[
\begin{aligned}
\Div(g_{R,S}^+)
 &=\{(R)+(S)+(-T)-3(O)\}\\
 &\quad-\{(T)+(-T)-2(O)\}\\
 &=(R)+(S)-(T)-(O),
\end{aligned}
\]
which proves \eqref{eq:native-Miller-factor-divisor}.

\medskip
\noindent\emph{The vertical specialization.}
If \(S=-R\), then \(T=O\) and the line is vertical.  Its function
\[
       X(Q)-X(R)=d\frac{u_R-u_Q}{u_Qu_R}
\]
has divisor \((R)+(-R)-2(O)\), the required specialization of the
Miller divisor.  Since all divisor computations use intersection
multiplicity, the conclusion also covers tangents and torsion points.
\end{proof}

\section[Miller accumulation on Cd]
{Miller accumulation on \(\mathcal C_d\)}

Choose the sign \(\epsilon=-\) in odd characteristic and \(\epsilon=+\)
in characteristic two.  Fix \(P\in\mathcal C_d\) and suppress it from the
running accumulator:
\[
       R_1=P,\qquad F_1(Q)=1.
\]
For any addition chain and positive integers \(a,b\), close the state by
\begin{equation}\label{eq:native-Miller-closure}
\boxed{
\begin{aligned}
R_{a+b}&=R_a+R_b
       &&\text{by the group law on \(\mathcal C_d\)},\\
F_{a+b}(Q)
 &=F_a(Q)F_b(Q)g_{R_a,R_b}^{\epsilon}(Q).
\end{aligned}}
\end{equation}
The factor in the second line is exactly
\eqref{eq:odd-native-Miller-thesis} or
\eqref{eq:binary-native-Miller-thesis}.  Thus every running point is stored
as a \(\mathcal C_d\)-point and every evaluation is a rational expression
in its \(u,v\) coordinates.

The boxed affine factors are to be used only when the displayed affine
coordinates and their denominators are defined.  If a running point is a
boundary point, the line and vertical functions must instead be evaluated
in the corresponding projective Weierstrass model (or by an explicitly
homogenized native formula).  In particular, an expression containing
\(u_T/u_R\) is not to be evaluated when \(u_R=0\) or \(u_R=\infty\).

No extension-field inversion is required at each step.  Write
\[
g_{R,S}^{\epsilon}(Q)
=\frac{\ell_{R,S}^{\epsilon}(Q)}
       {\nu_{R+S}^{\epsilon}(Q)},\qquad
F_a(Q)=\frac{N_a(Q)}{D_a(Q)}.
\]
Then the same closure is
\begin{equation}\label{eq:native-Miller-ND}
\begin{aligned}
N_{a+b}&=N_aN_b\,\ell_{R_a,R_b}^{\epsilon}(Q),&
D_{a+b}&=D_aD_b\,\nu_{R_{a+b}}^{\epsilon}(Q),\\
N_{2a}&=N_a^2\,\ell_{R_a,R_a}^{\epsilon}(Q),&
D_{2a}&=D_a^2\,\nu_{R_{2a}}^{\epsilon}(Q).
\end{aligned}
\end{equation}
Vertical steps use the factors
\eqref{eq:odd-native-vertical} or
\eqref{eq:binary-native-vertical} with denominator \(1\).
Only the final quotient \(N_n/D_n\) must be formed, unless a valid
denominator-elimination argument removes \(D_n\).

\begin{theorem}[A native denominator-elimination criterion]
\label{thm:Cd-denominator-elimination}
Let \(k_{\rm emb}\) be even, let \(n\mid q^{k_{\rm emb}/2}+1\), and
suppose all running
Miller points lie in \(\mathcal C_d(\F_q)\).  Choose the evaluation point
\(Q\in\mathcal C_d(\F_{q^{k_{\rm emb}}})\) so that
\[
                   u_Q\in\F_{q^{k_{\rm emb}/2}}.
\]
Then every vertical denominator appearing in
\eqref{eq:native-Miller-ND} belongs to
\(\F_{q^{k_{\rm emb}/2}}^\times\), and it is killed by the reduced Tate final
exponent.  Consequently the denominator accumulator may be omitted for the
reduced Tate pairing.  The same conclusion holds for shifted evaluations
provided the shifted arguments retain this subfield condition on their
\(u\)-coordinates.
\end{theorem}

\begin{proof}
Put \(h=k_{\rm emb}/2\) and
\[
             K_0=\F_{q^h}\subset K=\F_{q^{k_{\rm emb}}}.
\]
Because every running point is \(\F_q\)-rational, every finite affine
running abscissa \(u_T\) lies in
\(\F_q\subset K_0\).  The evaluation hypothesis gives \(u_Q\in K_0\).
For a nonvertical Miller step the denominator is
\[
\nu_T^\epsilon(Q)
=x_d^\epsilon(Q)-x_d^\epsilon(T)
=c_\epsilon\frac{u_T-u_Q}{u_Tu_Q},
\]
where
\[
        c_-=\beta\in\F_q^\times
        \quad\text{and}\quad
        c_+=d\in\F_q^\times.
\]
Thus \(\nu_T^\epsilon(Q)\in K_0\).  Admissibility of the evaluation
excludes its being zero and excludes \(u_Tu_Q=0\), so in fact
\(\nu_T^\epsilon(Q)\in K_0^\times\).

If a step has \(T=O\), it is a vertical step of the special form
\eqref{eq:odd-native-vertical} or
\eqref{eq:binary-native-vertical}; by the convention in
\eqref{eq:native-Miller-ND}, that factor is placed wholly in the numerator
and contributes denominator \(1\).  If
\(T=B_4^+\) or \(B_4^-\), then \(u_T=\infty\), and the homogeneous limit
of the displayed formula is
\[
                  \nu_T^\epsilon(Q)=\frac{c_\epsilon}{u_Q}\in K_0.
\]
The remaining boundary point \(B_2\) has finite nonzero \(u\)-coordinate,
so the preceding affine argument already applies.  Boundary
specializations therefore introduce no additional denominator outside
\(K_0\).

Starting with \(D_1=1\), the recurrences
\[
       D_{a+b}=D_aD_b\,\nu_{R_{a+b}}^\epsilon(Q),
       \qquad
       D_{2a}=D_a^2\,\nu_{R_{2a}}^\epsilon(Q)
\]
use only multiplication and squaring in \(K_0\).  Induction along the
addition chain consequently gives
\[
                       D_n(Q)\in K_0^\times.
\]

Now \(n\mid q^h+1\), so
\[
\begin{aligned}
\frac{q^{k_{\rm emb}}-1}{n}
 &=\frac{q^{2h}-1}{n}\\
 &=(q^h-1)\frac{q^h+1}{n}.
\end{aligned}
\]
Every \(z\in K_0^\times\) satisfies \(z^{q^h-1}=1\).  Hence
\[
       z^{(q^{k_{\rm emb}}-1)/n}=1
       \qquad\text{for every }z\in K_0^\times,
\]
and in particular the final exponentiation sends \(D_n(Q)\) to \(1\).
Therefore
\[
\left(\frac{N_n(Q)}{D_n(Q)}\right)^{
                  (q^{k_{\rm emb}}-1)/n}
=N_n(Q)^{(q^{k_{\rm emb}}-1)/n},
\]
which proves that the denominator accumulator may be omitted.

For each shifted argument \(Q_j\) whose \(u\)-coordinate belongs to
\(K_0\), the displayed formula for
\(\nu_T^\epsilon(Q_j)\), including its homogeneous boundary limits,
places every denominator factor in \(K_0^\times\).  Products of these
factors therefore remain in \(K_0^\times\) and are killed by the same
final exponent.  Thus the shifted conclusion holds exactly under the
coordinate hypothesis stated in the theorem.  The decisive
model-specific fact is that
\(\nu_T^\epsilon\) depends on the evaluation point only through its
quotient coordinate \(u\), equivalently through the normalized native
Kummer coordinate \(U=(u+1)/u\).
\end{proof}

\subsection{Sparse evaluation and accumulator costs}

Fix \(Q\) during a Miller loop and precompute
\[
 A_Q^-=(2v_Q+1)(u_Q+1),\qquad
 A_Q^+=v_Q(u_Q+1),\qquad B_Q=u_Q.
\]
The two line numerators can be rearranged as
\begin{align}
 \mathcal L_{R,S}^-(Q)
  &=(2\beta u_R)A_Q^-
   +\{\lambda_{R,S}^--2\beta s(R)(u_R+1)\}B_Q
   -\lambda_{R,S}^-u_R,\label{eq:odd-line-sparse-eval}\\
 \mathcal L_{R,S}^+(Q)
  &=u_RA_Q^+
   +\{\lambda_{R,S}^+-v_R(u_R+1)\}B_Q
   -\lambda_{R,S}^+u_R.
\label{eq:binary-line-sparse-eval}
\end{align}
Thus, after the one-time computation of \(A_Q^\epsilon\), each line
numerator uses only two multiplications of an extension-field element by a
base-field scalar.  Multiplying by the scale \(u_T/u_R\) uses one more if
that scale is not absorbed into the accumulator.  This count excludes the
base-field group operation and slope computation, which can share
intermediates with point addition or doubling.

Let \(\mathbf M_K,\mathbf S_K\) denote multiplication and squaring in the
extension field.  With separate numerator and denominator accumulators, a
doubling update in \eqref{eq:native-Miller-ND} costs
\[
                 2\mathbf S_K+2\mathbf M_K
\]
in accumulator arithmetic, and a general \(a+b\) update costs
\(4\mathbf M_K\), in addition to the sparse line and vertical evaluations.
Under Theorem~\ref{thm:Cd-denominator-elimination}, these become
\(\mathbf S_K+\mathbf M_K\) and \(2\mathbf M_K\), respectively; an
addition by a state whose accumulator is \(1\) saves one further
\(\mathbf M_K\).  These figures state exactly what is counted and avoid
mixing base-field point arithmetic with extension-field accumulation.

\begin{definition}\label{def:native-Miller-wrapper}
The raw Miller evaluation on \(\mathcal C_d\) is
\begin{equation}\label{eq:native-Miller-wrapper}
 \mathsf M_d^{\rm raw}(n;P,Q)
 :=F_n(Q)=\frac{N_n(Q)}{D_n(Q)},
\end{equation}
where \(R_a,N_a,D_a\) are obtained from
\eqref{eq:native-Miller-closure}--\eqref{eq:native-Miller-ND}.
If \([n]P=O\) and \(P\ne O\), fix a local parameter \(t\) at \(O\) and
write
\[
       F_n=c_{n,P}^{(t)}t^{-n}(1+O(t)),
       \qquad c_{n,P}^{(t)}\ne0.
\]
The locally normalized direct evaluation is
\begin{equation}\label{eq:native-Miller-local-normalization}
 \mathsf M_{d,t}(n;P,Q)
   :=\bigl(c_{n,P}^{(t)}\bigr)^{-1}
      \mathsf M_d^{\rm raw}(n;P,Q).
\end{equation}
Ratios of two evaluations of the same Miller function may use either
normalization, since the scalar cancels.  A single direct evaluation must
use \(\mathsf M_{d,t}\), not an unspecified raw scalar multiple.
\end{definition}

\begin{proposition}\label{prop:native-Miller-wrapper-divisor}
For every positive integer \(n\),
\[
\Div(F_n)
=n(P)-([n]P)-(n-1)(O).
\]
\end{proposition}

\begin{proof}
We prove simultaneously that
\[
                 R_a=[a]P
\]
and that
\begin{equation}\label{eq:native-Miller-wrapper-induction}
 \Div(F_a)=a(P)-([a]P)-(a-1)(O)
\end{equation}
for every index \(a\) constructed by the chosen addition chain.

For \(a=1\), the initialization gives \(R_1=P=[1]P\) and \(F_1=1\).
Therefore
\[
       \Div(F_1)=0=(P)-([1]P)-0(O),
\]
so both assertions hold at the base of the chain.

Suppose they hold for already constructed indices \(a\) and \(b\).
The first line of \eqref{eq:native-Miller-closure} uses the group law on
\(\mathcal C_d\), and hence
\[
       R_{a+b}=R_a+R_b=[a]P+[b]P=[a+b]P.
\]
For the accumulator,
Proposition~\ref{prop:odd-native-Miller-divisor} in odd characteristic and
Proposition~\ref{prop:binary-native-Miller-divisor} in characteristic two give
\[
\Div(g_{R_a,R_b}^\epsilon)
=([a]P)+([b]P)-([a+b]P)-(O).
\]
Using the induction hypotheses and additivity of divisors,
\[
\begin{aligned}
\Div(F_{a+b})
 &=\Div(F_a)+\Div(F_b)
   +\Div(g_{R_a,R_b}^\epsilon)\\
 &=a(P)-([a]P)-(a-1)(O)\\
 &\quad+b(P)-([b]P)-(b-1)(O)\\
 &\quad+([a]P)+([b]P)-([a+b]P)-(O)\\
 &=(a+b)(P)-([a+b]P)-(a+b-1)(O).
\end{aligned}
\]
This is \eqref{eq:native-Miller-wrapper-induction} for \(a+b\).
The doubling update is the special case \(b=a\); the separate
numerator--denominator recurrences in \eqref{eq:native-Miller-ND} form
the same quotient and therefore have the same divisor.

Induction over the finite sequence of chain operations reaches \(n\) and
proves the displayed assertion.  Different addition chains may change
\(F_n\) by a nonzero scalar, but
Proposition~\ref{prop:Miller-divisor-thesis} shows that they all have the same divisor;
degree-zero pairing evaluation removes this normalization ambiguity.
\end{proof}

\section[Pairing wrappers on Cd]
{Pairing wrappers on \(\mathcal C_d\)}

\begin{theorem}[Reduced Tate pairing on \(\mathcal C_d\)]
\label{thm:Tate-thesis}
Let \(K=\F_{q^{k_{\rm emb}}}\), let \(n>1\), and let
\[
P\in\mathcal C_d(\F_q)[n]\setminus\{O\},\qquad
Q\in\mathcal C_d(K),\qquad
\gcd(n,q)=1,\qquad n\mid q^{k_{\rm emb}}-1.
\]
If the evaluation avoids zeros and poles, then
\begin{equation}\label{eq:Tate-thesis}
\boxed{\qquad
t_{n,d}^{\mathcal C_d}(P,[Q])
=\mathsf M_{d,t}(n;P,Q)^{(q^{k_{\rm emb}}-1)/n}
\in\mu_n.
\qquad}
\end{equation}
Here \(t\) is the declared local parameter at \(O\) used in
\eqref{eq:native-Miller-local-normalization}.
For an arbitrary admissible shift \(T\in\mathcal C_d(K)\), the
normalization-free wrapper is
\begin{equation}\label{eq:Tate-shifted-Cd}
\boxed{\qquad
t_{n,d}^{\mathcal C_d}(P,[Q])
=\left(
\frac{\mathsf M_d^{\rm raw}(n;P,Q+T)}
     {\mathsf M_d^{\rm raw}(n;P,T)}
\right)^{(q^{k_{\rm emb}}-1)/n}.
\qquad}
\end{equation}
All sums in \eqref{eq:Tate-shifted-Cd} are sums on
\(\mathcal C_d:(u^2+u)(v^2+v)=d\).
For the omitted trivial input \(P=O\), both pairing wrappers are defined
to have value \(1\).
\end{theorem}

\begin{proof}
The proof has four parts: choose a disjoint degree-zero representative,
prove independence of that representative, prove invariance modulo
\(nE(K)\), and finally compare shifted evaluation with the locally
normalized direct evaluation.
Write \(E=\mathcal C_d\) and
\[
                  D_P=(P)-(O).
\]
Because \(P\in E[n]\), one has \([n]P=O\).
Proposition~\ref{prop:native-Miller-wrapper-divisor} therefore gives
\[
             \Div(F_n)=n(P)-n(O)=nD_P.
\]

Choose \(T\) so that all supports used below are disjoint from the zeros
and poles of \(F_n\), and set
\[
                  D_{Q,T}=(Q+T)-(T).
\]
Under the standard identification
\[
 E(K)\longrightarrow\operatorname{Pic}^0(E),\qquad
 R\longmapsto[(R)-(O)],
\]
the divisor \(D_{Q,T}\) represents the point
\((Q+T)-T=Q\).  Hence it represents the same divisor class as
\((Q)-(O)\), but its support can be moved away from \(P\) and \(O\).
By the definition of evaluation of a function on a divisor,
\[
 F_n(D_{Q,T})=\frac{F_n(Q+T)}{F_n(T)},
\]
which is the unexponentiated quotient in
\eqref{eq:Tate-shifted-Cd}.  If \(F_n\) is replaced by \(cF_n\), then
\[
       (cF_n)(D_{Q,T})
       =c^{\deg D_{Q,T}}F_n(D_{Q,T})
       =F_n(D_{Q,T}),
\]
because \(\deg D_{Q,T}=0\).  Thus the shifted quotient has no
normalization ambiguity.

We next prove well-definedness in
\(K^\times/(K^\times)^n\).  First replace \(D_{Q,T}\) by a linearly
equivalent degree-zero divisor
\[
                 D_{Q,T}'=D_{Q,T}+\Div(h)
\]
whose support is also admissible.  Weil reciprocity gives
\[
 \frac{F_n(D_{Q,T}')}{F_n(D_{Q,T})}
 =F_n(\Div(h))
 =h(\Div(F_n))
 =h(nD_P)
 =h(D_P)^n.
\]
Thus changing a divisor representative changes the evaluation only by an
\(n\)th power.

Now replace \(Q\) by a point \(Q'=Q+[n]R\).  Choose an admissible
degree-zero divisor \(D_R\) representing \(R\).  In
\(\operatorname{Pic}^0(E)\),
\[
          [D_{Q'}]-[D_Q]=n[D_R].
\]
Consequently there is a rational function \(h\) such that, after moving
supports if necessary,
\[
              D_{Q'}-D_Q=nD_R+\Div(h).
\]
Evaluating \(F_n\) and applying Weil reciprocity once more gives
\[
\begin{aligned}
\frac{F_n(D_{Q'})}{F_n(D_Q)}
 &=F_n(D_R)^nF_n(\Div(h))\\
 &=F_n(D_R)^n h(\Div(F_n))\\
 &=F_n(D_R)^n h(D_P)^n\\
 &=\{F_n(D_R)h(D_P)\}^n.
\end{aligned}
\]
Hence the class of \(F_n(D_Q)\) in
\(K^\times/(K^\times)^n\) depends only on
\([Q]\in E(K)/nE(K)\).

Put
\[
             h_{\rm red}=\frac{q^{k_{\rm emb}}-1}{n}.
\]
If an evaluation is multiplied by an \(n\)th power \(b^n\), its
\(h_{\rm red}\)-th power is unchanged because
\[
          (b^n)^{h_{\rm red}}=b^{q^{k_{\rm emb}}-1}=1.
\]
Furthermore, for \(a\in K^\times\),
\[
          (a^{h_{\rm red}})^n=a^{q^{k_{\rm emb}}-1}=1.
\]
Thus final exponentiation gives a well-defined element of \(\mu_n\).
This proves the shifted formula
\eqref{eq:Tate-shifted-Cd}.

It remains to relate the shifted divisor evaluation to the direct wrapper.
For a local parameter \(t\) at \(O\), write
\[
                 F_n=c_{n,P}^{(t)}t^{-n}(1+O(t)).
\]
Define the regularized evaluation on the intersecting divisor
\((Q)-(O)\) by
\[
 F_n[(Q)-(O)]_t
 :=\frac{F_n(Q)}{c_{n,P}^{(t)}}
 =\mathsf M_{d,t}(n;P,Q).
\]
Put \(D_0=(Q)-(O)\), and choose \(h\) with
\(D_{Q,T}=D_0+\Div(h)\).  The tame-symbol form of Weil reciprocity gives,
with leading coefficients used at the common support,
\[
 \frac{F_n(D_{Q,T})}{F_n[D_0]_t}
       =(-1)^n h(D_P)^n.
\]
Here \(h(D_P)\) is interpreted by the same leading-coefficient convention
at any common support.
The sign is itself an \(n\)th power: it is \(1\) when \(n\) is even,
and is \((-1)^n\) when \(n\) is odd.  Consequently the two evaluations
determine the same class in
\(K^\times/(K^\times)^n\).  Moreover, replacing \(t\) by
\(t'=at+O(t^2)\) multiplies the regularized value by \(a^{-n}\), again an
\(n\)th power.  Final exponentiation therefore gives
\eqref{eq:Tate-thesis}, independently of the chosen local parameter.  The
shifted degree-zero quotient cancels every nonzero scalar already before
the final exponentiation.
\end{proof}

\begin{theorem}[Weil pairing on \(\mathcal C_d\)]
\label{thm:Weil-thesis}
Let \(\gcd(n,q)=1\), and let
\(P,Q\in\mathcal C_d(\overline{\F}_q)[n]\) be independent.  Choose \(T\)
so that all four evaluations below avoid zeros and poles.  Then
\begin{equation}\label{eq:Weil-shifted-Cd}
\boxed{
e_{n,d}^{\mathcal C_d}(P,Q)
=
\frac{\mathsf M_d^{\rm raw}(n;P,Q+T)}
     {\mathsf M_d^{\rm raw}(n;P,T)}
\frac{\mathsf M_d^{\rm raw}(n;Q,-T)}
     {\mathsf M_d^{\rm raw}(n;Q,P-T)}.
}
\end{equation}
With the common local-parameter normalization at \(O\), and when direct
evaluation at \(P,Q\) is admissible, this reduces to
\begin{equation}\label{eq:Weil-thesis}
\boxed{\qquad
e_{n,d}^{\mathcal C_d}(P,Q)
=(-1)^n
\frac{\mathsf M_{d,t}(n;P,Q)}
     {\mathsf M_{d,t}(n;Q,P)}.
\qquad}
\end{equation}
\end{theorem}

\begin{proof}
We first obtain a normalization-free formula from disjoint divisors and then
take the moving divisor to the identity to recover the sign in the direct
locally normalized formula.
Put \(E=\mathcal C_d\), and choose Miller functions
\[
\begin{aligned}
\Div(F_P)&=n(P)-n(O),&
F_P&=\mathsf M_d^{\rm raw}(n;P,\mathord{\cdot}),\\
\Div(F_Q)&=n(Q)-n(O),&
F_Q&=\mathsf M_d^{\rm raw}(n;Q,\mathord{\cdot}).
\end{aligned}
\]
We first derive the shifted expression from disjoint divisors, so that no
evaluation at a common zero or pole is hidden.

Let
\[
        D_P=(P)-(O),\qquad
        D_{Q,T}=(Q+T)-(T).
\]
The second divisor represents the class of \(Q\), and the choice of
\(T\) makes the two supports disjoint and all evaluations admissible.
Define the translated function
\[
               \widetilde F_Q(X)=F_Q(X-T).
\]
Translation is an automorphism of the curve, and therefore
\[
\begin{aligned}
\Div(\widetilde F_Q)
 &=n(Q+T)-n(T)\\
 &=nD_{Q,T}.
\end{aligned}
\]
The divisor definition of the Weil pairing for the two disjoint
representatives is
\begin{equation}\label{eq:Weil-divisor-definition-expanded}
 e_n(P,Q)
 =\frac{F_P(D_{Q,T})}
        {\widetilde F_Q(D_P)}.
\end{equation}
The numerator is
\[
       F_P(D_{Q,T})
       =\frac{F_P(Q+T)}{F_P(T)},
\]
while the denominator is
\[
       \widetilde F_Q(D_P)
       =\frac{\widetilde F_Q(P)}{\widetilde F_Q(O)}
       =\frac{F_Q(P-T)}{F_Q(-T)}.
\]
Substituting these two evaluations into
\eqref{eq:Weil-divisor-definition-expanded} gives
\[
 e_n(P,Q)
 =\frac{F_P(Q+T)}{F_P(T)}
  \frac{F_Q(-T)}{F_Q(P-T)},
\]
which is exactly \eqref{eq:Weil-shifted-Cd}.  Multiplying either Miller
function by a nonzero scalar does not change its evaluation on a
degree-zero divisor.  Thus both normalization constants cancel
independently.

We now derive the sign in the direct formula.  Choose one local parameter
\(t\) at \(O\) and normalize both Miller functions by
\[
            F_P=t^{-n}(1+O(t)),\qquad
            F_Q=t^{-n}(1+O(t)).
\]
The inversion map on an elliptic curve has differential \(-1\) at the
identity, so
\begin{equation}\label{eq:pairing-local-inversion-parameter}
                    t(-T)=-t(T)+O(t(T)^2)
                    \qquad (T\longrightarrow O).
\end{equation}
The shifted expression is independent of the admissible moving point
\(T\), because it is the divisor definition of \(e_n(P,Q)\).  We may
therefore compute its value from its Laurent expansion as \(T\) tends to
\(O\).  Admissibility of the direct evaluations gives
\[
\begin{aligned}
F_P(Q+T)&=F_P(Q)+O(t(T)),\\
F_Q(P-T)&=F_Q(P)+O(t(T)).
\end{aligned}
\]
The two factors near the common pole satisfy
\[
\begin{aligned}
\frac{F_Q(-T)}{F_P(T)}
 &=
 \frac{t(-T)^{-n}(1+O(t(T)))}
      {t(T)^{-n}(1+O(t(T)))}\\
 &=
 \left(\frac{t(T)}{-t(T)+O(t(T)^2)}\right)^n
   (1+O(t(T)))\\
 &=(-1)^n+O(t(T)).
\end{aligned}
\]
Taking the constant term of the shifted identity consequently yields
\[
              e_n(P,Q)
              =(-1)^n\frac{F_P(Q)}{F_Q(P)},
\]
which is \eqref{eq:Weil-thesis}.  In characteristic two the sign is,
of course, \(1\).

Finally, \eqref{eq:Weil-divisor-definition-expanded} is the standard
divisor construction of the Weil pairing.  Multiplication of functions
and addition of degree-zero divisor classes give bilinearity, the
commutator construction gives alternation, and prime-to-characteristic
Cartier duality makes the pairing nondegenerate on \(E[n]\).  Since the
native wrappers have exactly the required divisors by
Proposition~\ref{prop:native-Miller-wrapper-divisor}, the displayed formulas compute
the Weil pairing itself and not an undetermined scalar multiple.
\end{proof}

\section{Simultaneous and product pairings}

Let \((P_j,Q_j)\), \(1\le j\le h\), be admissible Tate inputs of the same
order \(n\) over \(K=\F_{q^{k_{\rm emb}}}\).  Multiplicativity gives the product
wrapper
\begin{equation}\label{eq:Cd-product-Tate-pairing}
 \prod_{j=1}^h t_{n,d}^{\mathcal C_d}(P_j,[Q_j])
 =\left(\prod_{j=1}^h\mathsf M_{d,t}(n;P_j,Q_j)\right)^{
       (q^{k_{\rm emb}}-1)/n}.
\end{equation}
Thus only one final exponentiation is required.  If all \(P_j\) are equal,
one native addition chain, one sequence of running points, and one sequence
of slopes can be shared; the sparse expressions
\eqref{eq:odd-line-sparse-eval}--\eqref{eq:binary-line-sparse-eval} are
then evaluated at the different \(Q_j\).  If the \(P_j\) differ, the line
states remain separate but the final exponentiation is still shared.

The same principle applies to products of Weil pairings by multiplying the
two shifted Miller quotients for every pair before taking the final product.
By Theorem~\ref{thm:Tate-thesis}, each locally normalized factor is the
corresponding Tate pairing.  Multiplying those identities and distributing
the common exponent over their product gives
\eqref{eq:Cd-product-Tate-pairing}; each direct factor uses the local
normalization \eqref{eq:native-Miller-local-normalization}.  Alternatively,
each factor may be replaced by its shifted raw quotient, in which case the
scalar cancels on a degree-zero divisor.

\section[Tripling and 2P+Q inside Miller loops]
{Tripling and \texorpdfstring{\(2P+Q\)}{2P+Q} inside Miller loops}

From \eqref{eq:Miller-recursion-thesis},
\begin{equation}\label{eq:Miller-tripling-thesis}
F_{3m}(Q)
=F_m(Q)^3
g_{[m]P,[m]P}^{\epsilon}(Q)
g_{[2m]P,[m]P}^{\epsilon}(Q).
\end{equation}
Applying \eqref{eq:Miller-recursion-thesis} first to \(2m\) and then
to \(2m+n\) gives
\begin{equation}\label{eq:Miller-doubleadd-thesis}
F_{2m+n}(Q)
=F_m(Q)^2F_n(Q)
g_{[m]P,[m]P}^{\epsilon}(Q)
g_{[2m]P,[n]P}^{\epsilon}(Q).
\end{equation}
Thus the Kummer tripling core supplies the updated running point, while the
Miller accumulator records the two line factors with their distinct
divisors.  A closed \(2P+Q\) coordinate formula can share base-field
intermediates with these evaluations, and the two divisor functions complete
the corresponding Miller step.

\section{Comparison with Edwards pairing formulas}

The models and the evaluated functions must be named before costs are
compared.  The relevant interfaces are as follows.
\begin{table}[H]
\centering
\caption{Pairing interfaces: declared model, running point, and function}
\label{tab:pairing-interface-comparison}
\begin{tabular}{L{2.8cm}L{4.2cm}L{2.8cm}L{3.2cm}}
\toprule
declared model & full equation & running point & Miller object\\
\midrule
\(\mathcal C_d\), odd
& \((u^2+u)(v^2+v)=d\)
& \((u_R,v_R)\) on \(\mathcal C_d\)
& \(\mathcal L^-_{R,S}(Q)/(u_T-u_Q)\),
  \eqref{eq:odd-native-Miller-thesis}\\
\(\mathcal C_d\), binary
& \((u^2+u)(v^2+v)=d\)
& \((u_R,v_R)\) on \(\mathcal C_d\)
& \(\mathcal L^+_{R,S}(Q)/(u_T-u_Q)\),
  \eqref{eq:binary-native-Miller-thesis}\\
Edwards
& \(\xi^2+\eta^2=1+\rho\xi^2\eta^2\)
& Edwards extended coordinates
& conic divided by vertical factors\\
Jacobian Weierstrass
& \(Y^2=X^3+a_4X+a_6\)
& Jacobian coordinates
& tangent/secant line divided by a vertical line\\
\bottomrule
\end{tabular}
\end{table}

For a general twisted Edwards equation
\(a\xi^2+\eta^2=1+d_E\xi^2\eta^2\), direct conic formulas in a common
even-embedding-degree setting have base-field coefficient costs
\[
6\M+5\Sqr+\Dconst{a}
\quad\text{for doubling},\qquad
12\M+\Dconst{a}
\quad\text{for mixed addition},
\]
in addition to extension-field accumulator updates
\cite{AreneLangeNaehrigRitzenthaler2011}.  For the row displayed above,
\(a=1\), so \(\Dconst{a}\) is free and the specialized costs are
\(6\M+5\Sqr\) and \(12\M\).  Representative optimized
Jacobian Weierstrass costs are
\[
\M+11\Sqr+\Dconst{a_4},
\qquad
6\M+6\Sqr,
\]
again before the extension-field terms are specialized.

The direct \(\mathcal C_d\) formulas have four concrete advantages.
\begin{enumerate}[label=\textup{(\roman*)}]
 \item The running points, additions, doublings, and final pairing arguments
       all remain in \(u,v\); an Edwards or Weierstrass point is never the
       declared loop state.
 \item The line numerators are sparse.  In odd characteristic
       \(\mathcal L^-_{R,S}(Q)\) is linear in
       \(s(Q)(u_Q+1)\) and \(u_Q\); in characteristic two
       \(\mathcal L^+_{R,S}(Q)\) is linear in
       \(v_Q(u_Q+1)\) and \(u_Q\).
 \item Every vertical denominator depends only on the Kummer
       coordinate \(u_Q\).  Consequently each Miller update places the
       sparse line numerator in the numerator accumulator and the vertical
       factor in the denominator accumulator, which gives
       \eqref{eq:native-Miller-ND} and permits the usual denominator
       elimination whenever the twist and final exponent satisfy its
       hypotheses.
 \item The same raw wrapper
       \(\mathsf M_d^{\rm raw}(n;P,Q)\), followed by either local
       normalization or a shifted degree-zero quotient, works in characteristic three
       through the odd branch and in characteristic two through the
       Artin--Schreier branch.
\end{enumerate}

Taken together, these properties give the native \(\mathcal C_d\) Miller
framework a unified end-to-end structure: sparse line numerators are
evaluated directly in the native \(u,v\)-coordinates, vertical denominators
depend only on the native Kummer coordinate, and the numerator--denominator
organization supports denominator elimination and shared final
exponentiation under the stated twist and subfield hypotheses.  The same
characteristic-uniform point interface supplies the Miller updates in odd
characteristic, in characteristic three, and in the characteristic-two
Artin--Schreier branch.  Its complete implementation cost is determined
explicitly by the embedding degree, twist, extension-field basis,
denominator-elimination rule, addition chain, and final exponentiation.

\section{A pairing example}

\begin{example}[A third-order pairing over \(\F_{101^2}\)]
Let
\[
\F_{101^2}=\F_{101}[\iota]/(\iota^2-2),\qquad d=1.
\]
Work throughout on
\[
\mathcal C_1:\quad (u^2+u)(v^2+v)=1,
\qquad
\#\mathcal C_1(\F_{101})=96.
\]
The point
\[
             P=(96,26)
\]
has order \(3\), and
\[
       2P=-P=(96,74).
\]
Take the Tate argument
\[
       R=(8+95\iota,\,90+67\iota)\in\mathcal C_1(\F_{101^2}).
\]
Use the standard local parameter \(t=-X/Y\) at the point at infinity of
\(W_1^-\).  For the two-step order-three chain used below, the first
line-over-vertical factor has leading term \(-t^{-1}\) and the final
vertical factor has leading term \(t^{-2}\).  Hence the raw accumulator
has \(c_{3,P}^{(t)}=-1\), and its locally normalized value is the negative
of the displayed raw value.  The Tate exponent \(3400\) is even, and the
same factor \(-1\) occurs in both order-three Weil accumulators, so neither
final value below is changed.
Here
\[
 \beta=19,\qquad A=74,\qquad s(P)=53,
\qquad \lambda_{P,P}^-=49.
\]
Substitution of these \(\mathcal C_1\)-coordinates into
\eqref{eq:odd-native-Miller-thesis} and
\eqref{eq:odd-native-vertical} gives
\[
g_{P,P}^{-}(R)=50+87\iota,\qquad
g_{2P,P}^{-}(R)
=19\,\frac{96-(8+95\iota)}{96(8+95\iota)}
=5+11\iota.
\]
Consequently the Miller evaluation is
\[
\mathsf M_1^{\rm raw}(3;P,R)
=(50+87\iota)(5+11\iota)
=43+76\iota.
\]
Since the embedding degree is \(2\) and
\((101^2-1)/3=3400\), the direct
\(\mathcal C_1\)-Tate wrapper gives
\[
t_{3,1}^{\mathcal C_1}(P,[R])
=(43+76\iota)^{3400}=50+54\iota,
\qquad (50+54\iota)^3=1.
\]

For the Weil pairing, take another independent order-three point
\[
       Q=(57,50+92\iota)\in\mathcal C_1[3].
\]
The same \(u,v\)-formulas give
\[
\mathsf M_1^{\rm raw}(3;P,Q)=57+7\iota,\qquad
\mathsf M_1^{\rm raw}(3;Q,P)=30+6\iota.
\]
Formula~\eqref{eq:Weil-thesis} therefore yields
\[
e_{3,1}^{\mathcal C_1}(P,Q)
=-\frac{57+7\iota}{30+6\iota}
=50+47\iota,\qquad (50+47\iota)^3=1.
\]
Every point addition, tangent coefficient, vertical factor, and evaluation
in this calculation is expressed in the coordinates of
\(\mathcal C_1\).  Interchanging \(P,Q\) gives the inverse, checking
alternation.
\end{example}

\part[Implementation and comparison]
{Implementation, Specialization, Comparison, and Examples}
\partoverview{This part converts the preceding formulas into a
decision framework.  It separates full-point, Kummer, and torsion-quotient
outputs and first records parameter-selection, encoding, and constant-time
requirements.  The Cd25519 case study is then developed before its costs are
used in the model-by-model comparison.  Representative finite-field examples
follow the comparison and close the implementation-oriented part.}

\chapter[Parameters, Encoding, and Constant Time]
{Parameter Selection, Encoding, and Constant-Time Use}
\label{ch:implementation}
This chapter turns the algebraic formulas into admissible interfaces before
any parameter-specific benchmark is discussed.  It first states smoothness,
completeness, and constant-cost criteria for parameter selection in the two
characteristic branches; it then defines the affine and boundary encodings;
finally it separates mathematical decoding from subgroup policy and
constant-time failure handling.  The Cd25519 specialization in the next
chapter uses these conventions unchanged.

\section{Odd-characteristic parameter selection}

To choose a cheap Montgomery constant, start with
\(c=\alpha_{24}\ne0,1\).  Then
\begin{equation}\label{eq:reverse-parameter-thesis}
       d=\frac1{16c},\qquad
       \rho=1-\frac1c,\qquad
       A=4c-2.
\end{equation}

\begin{proposition}[Admissibility in the \(c\)-parameter]
\label{prop:odd-implementation-parameter}
Let \(q\) be odd.  The substitution
\[
 c=\alpha_{24}=\frac1{16d}
\]
is a bijection between smooth parameters
\(d\in\F_q\setminus\{0,1/16\}\) and
\(c\in\F_q\setminus\{0,1\}\).  Under this bijection,
\[
 A^2-4=16c(c-1),\qquad \rho=\frac{c-1}{c}.
\]
Consequently the native Segre addition tuple of
Theorem~\ref{thm:native-complete-addition-Cd} is a single
\(\F_q\)-complete law exactly when \((c-1)/c\) is a nonsquare.
\end{proposition}

\begin{proof}
Because the characteristic is odd, \(16\in\F_q^\times\).  Hence
\[
       \Phi(d)=\frac{1}{16d},\qquad
       \Psi(c)=\frac{1}{16c}
\]
are mutually inverse on \(\F_q^\times\).  Under this inversion,
\(d=1/16\) is sent to \(c=1\), and conversely.  Therefore \(\Phi\)
restricts to the claimed bijection
\[
 \F_q\setminus\{0,1/16\}\longleftrightarrow
 \F_q\setminus\{0,1\}.
\]
This also verifies the smoothness condition directly: using
\(d=(16c)^{-1}\),
\[
 d(1-16d)=\frac{1}{16c}\left(1-\frac1c\right)
          =\frac{c-1}{16c^2},
\]
which is nonzero exactly when \(c(c-1)\ne0\).

The two parameter identities follow by calculation:
\[
 \rho=1-16d=1-\frac1c=\frac{c-1}{c}
\]
and, since \(A=4c-2\),
\[
 A^2-4=(4c-2)^2-4=16c^2-16c=16c(c-1).
\]
The tuple referred to in the statement is the addition tuple
\eqref{eq:native-complete-addition-Cd}, expressed in the centered native
Segre coordinates of \(\overline{\mathcal C}_d\).  By
Theorem~\ref{thm:exact-native-completeness-criterion}, it has no
\(\F_q\)-rational base pair if and only if \(\rho\) is a nonsquare.
Substituting \(\rho=(c-1)/c\) gives exactly the final criterion.  Thus the
assertion concerns a law on the original smooth \((2,2)\)-completion and
does not require changing the computational model.
\end{proof}

If complete Edwards addition is required, \(\rho\) must be a nonsquare.
If the six-square post-four-isogeny Kummer coordinate is required
simultaneously, one may
instead choose
\[
       d=h^2,\qquad \rho=1-16h^2
\]
with \(\rho\) nonsquare and \(\tau=4h\) cheap.  These objectives need not
make \(\alpha_{24}\), \(\tau\), and \(\tau^{-1}\) cheap at the same time.

If a GLV endomorphism is required, parameter selection is restricted to one
of two loci:
\[
\begin{array}{c|c|c}
\text{locus}&\text{parameter choice}&\text{base-field constant}\\ \hline
j=1728&d=1/8&i^2=-1,\\
j=0&A^2=3,\quad d=1/(4(A+2))
&\zeta^2+\zeta+1=0.
\end{array}
\]
For cryptographic deployment, the chosen large prime-order subgroup is
verified to be stable under the endomorphism, the corresponding eigenvalue
\(\lambda\) is computed, and both the main curve and its twist are selected
with the required subgroup structure.

An implementation parameter search must also verify:
\begin{enumerate}[label=\textup{(\arabic*)}]
 \item a sufficiently large prime-order subgroup and a controlled cofactor;
 \item twist security and the absence of dangerous small factors;
 \item exclusion of anomalous and low-embedding-degree curves when the
       protocol requires it;
 \item a field modulus suitable for constant-time reduction;
 \item the exact encoding and validation rules for boundary values.
\end{enumerate}
A small curve constant is an implementation filter, not a security
criterion.

\section{Characteristic-two parameter selection}

The Kummer constant is \(d_{\rm K}=d^{-1}\), and
\[
       d=d_{\rm K}^{-1},\qquad j=d_{\rm K}^4.
\]
One may therefore choose a sparse or otherwise cheap \(d_{\rm K}\), then recover
\(d\).  Because every member has rational four-torsion, its group order is
divisible by four; cofactor handling must be part of the protocol.  If the
target ordinary curve has no rational four-torsion, a twisted
\(\mu_4\)-normal form or general binary Edwards model is required instead.

For \(q=2^m\), every \(d_{\rm K}\in\F_q^\times\) gives
\(d=d_{\rm K}^{-1}\ne0\) and hence a
smooth member of the family.  Distinct implementation goals should still be
tested separately: sparsity of \(d_{\rm K}\) controls multiplication by the Kummer
constant, the factorization of \(N_d=\#\mathcal C_d(\F_q)\) controls subgroup
selection, and the order of the quadratic twist controls twist security.
None of these latter two properties follows from the Hamming weight of
\(d_{\rm K}\).

\section{Point encoding}

\begin{theorem}[Canonical affine encoding and decoding]
\label{thm:Cd-canonical-encoding}
Fix a canonical representation of \(\F_q\).

\begin{enumerate}[label=\textup{(\roman*)}]
\item If \(q\) is odd and \(u\notin\{0,-1\}\), put
\[
       w_u=1+\frac{4d}{u^2+u}.
\]
There is an affine point above \(u\) if and only if \(w_u\) is a square.
For a square root \(z^2=w_u\), the two possible ordinates are
\[
                    v_\pm=\frac{-1\pm z}{2}.
\]
A sign function satisfying
\(\operatorname{sgn}(-z)=1-\operatorname{sgn}(z)\) for \(z\ne0\)
therefore gives an injective encoding \((u,\operatorname{sgn}(2v+1))\);
when \(z=0\), only the sign bit \(0\) is canonical.

\item If \(q=2^m\) and \(u\notin\{0,1\}\), put
\[
                    a_u=\frac{d}{u^2+u}.
\]
There is an affine point above \(u\) if and only if
\(\operatorname{Tr}_{\F_q/\F_2}(a_u)=0\).  Its two ordinates are
\(v\) and \(v+1\).  Choose an \(\F_2\)-linear functional
\(L:\F_q\to\F_2\) with \(L(1)=1\).  Then
\((u,L(v))\) is an injective encoding.
\end{enumerate}

In both cases the four boundary points are encoded by four reserved tags
\[
 O=(0,\infty),\quad B_2=(-1,\infty),\quad
 B_4^+=(\infty,0),\quad B_4^-=(\infty,-1),
\]
and no affine byte string is permitted to alias a reserved tag.
\end{theorem}

\begin{proof}
We prove the two characteristic cases separately.

Assume first that \(q\) is odd.  Since \(d\ne0\), the affine equation has
no point with \(u^2+u=0\); thus every affine point has
\(u\notin\{0,-1\}\), and division by \(u^2+u\) is legitimate.  The equation
becomes
\[
                    v^2+v=\frac{d}{u^2+u}.
\]
Multiplication by \(4\) and addition of \(1\) give the equivalent equation
\begin{equation}\label{eq:odd-decompression-thesis}
       (2v+1)^2=1+\frac{4d}{u^2+u}=w_u.
\end{equation}
Consequently a point above \(u\) exists if and only if \(w_u\) is a
square.  If \(z^2=w_u\), solving the linear equations
\(2v+1=z\) and \(2v+1=-z\) gives
\[
                  v_+=\frac{-1+z}{2},\qquad
                  v_-=\frac{-1-z}{2}.
\]
For \(z\ne0\) these are distinct, and they are the only roots of the
quadratic.  The group inverse on \(\mathcal C_d\) is
\((u,v)\mapsto(u,-v-1)\), so it interchanges \(v_+\) and \(v_-\) and sends
\(2v+1\) to its negative.  The defining property of
\(\operatorname{sgn}\) therefore assigns opposite bits to the two roots.
If \(z=0\), the two displayed roots coincide at \(v=-1/2\); admitting both
bits would create two strings for the same point, so the convention that
only bit \(0\) is accepted is both necessary and sufficient for
canonicity.  It follows that the pair
\((u,\operatorname{sgn}(2v+1))\) determines exactly one affine point.

Now let \(q=2^m\).  Here \(u^2+u=0\) exactly for \(u=0,1\), and again
these values cannot occur on the affine curve because \(d\ne0\).  For every
other \(u\), the fiber equation is
\[
                         v^2+v=a_u.
\]
Consider the \(\F_2\)-linear Artin--Schreier map
\(\wp(z)=z^2+z\).  Its kernel is \(\{0,1\}\), so rank--nullity gives
\(\#\operatorname{im}(\wp)=q/2\).  Moreover
\[
 \operatorname{Tr}(\wp(z))
 =\operatorname{Tr}(z^2)+\operatorname{Tr}(z)
 =\operatorname{Tr}(z)^2+\operatorname{Tr}(z)=0.
\]
The trace-zero elements also form an \(\F_2\)-hyperplane of size \(q/2\).
The inclusion just proved is therefore an equality:
\[
 \operatorname{im}(\wp)=\ker\bigl(\operatorname{Tr}_{\F_q/\F_2}\bigr).
\]
This proves the trace criterion.  Whenever it holds, one solution \(v\)
produces the full fiber \(v+\ker(\wp)=\{v,v+1\}\).

An \(\F_2\)-linear functional with \(L(1)=1\) always exists: extend
\(1\) to an \(\F_2\)-basis of \(\F_q\), assign value \(1\) to the first
basis vector and arbitrary fixed values to the others, and extend linearly.
Then
\[
                         L(v+1)=L(v)+1,
\]
so the two roots receive opposite bits.  Hence \((u,L(v))\) also determines
exactly one affine point.

Finally, the smooth completion lies in \(\PP^1\times\PP^1\).  Its four
boundary points each have one infinite coordinate and are the four
distinct projective points listed in the statement (in characteristic two,
\(-1=1\), but the two choices in each factor remain distinct points of
\(\PP^1\)).  None is represented by a finite affine pair.  Four disjoint
reserved tags therefore encode them injectively, and forbidding those tags
in the affine format prevents every boundary/affine alias.  Together with
the preceding fiber arguments, this proves the theorem.
\end{proof}

\begin{remark}[Deterministic Artin--Schreier solvers]
When \(m\) is odd, the half-trace
\[
 \operatorname{HT}(a)=
 \sum_{i=0}^{(m-1)/2}a^{\,2^{2i}}
\]
satisfies
\(\operatorname{HT}(a)^2+\operatorname{HT}(a)=a\) for trace-zero \(a\).
For even \(m\), one may use a fixed linearized-polynomial solver or a
precomputed linear map in the chosen basis.  The choice is part of the
field implementation and does not change the encoding theorem.
\end{remark}

\begin{proposition}[Canonical decoder invariant]
\label{prop:Cd-decoder-invariant}
Suppose a decoder accepts only canonical field elements, enforces the
square or trace test of Theorem~\ref{thm:Cd-canonical-encoding}, checks the
root-selection bit, separates the four boundary tags, verifies the curve
equation, and performs the protocol's subgroup test.  Then every accepted
string determines exactly one permitted \(\mathcal C_d\)-point, and
re-encoding that point returns the original string.
\end{proposition}

\begin{proof}
Let \(s\) be an accepted string.  There are two disjoint cases.

If \(s\) is one of the four reserved tags, the decoder returns the unique
boundary point assigned to that tag.  The tags are pairwise distinct and the
affine parser is forbidden to accept them, so no other decoding path can
return the same point.  Re-encoding uses the same table and returns \(s\).

Otherwise \(s\) consists of a canonical field representation of \(u\) and
a root-selection bit.  Canonical parsing supplies one field element, rather
than several byte representatives of the same element.  In odd
characteristic, the square test rejects precisely the values for which
\eqref{eq:odd-decompression-thesis} has no solution.  If \(w_u\ne0\), the
two square roots are \(z\) and \(-z\), and the sign rule selects exactly one
of them; the equation \(v=(-1+z)/2\) then determines exactly one ordinate.
If \(w_u=0\), there is one ordinate and the decoder accepts only the
canonical bit \(0\).  Thus no odd-characteristic affine point has two
accepted strings.

In characteristic two, the trace test is equivalent to solvability by
Theorem~\ref{thm:Cd-canonical-encoding}.  Starting from either solution of
\(v^2+v=a_u\), the complete solution set is \(\{v,v+1\}\), and
\(L(v+1)=L(v)+1\).  Exactly one of these two elements has the requested
bit, so the ordinate is again unique.  This also shows that the selected
root is independent of which preliminary Artin--Schreier solution a fixed
solver happens to return.

The explicit curve-equation check can only reject the uniquely reconstructed
candidate; it cannot create a second candidate.  Likewise, the protocol's
subgroup test restricts the accepted point set but does not identify two
points.  Therefore every accepted string decodes to exactly one permitted
point.  In each case the stored field element is already canonical and the
root bit is, by construction, the canonical bit of the recovered ordinate,
so encoding the point returns the original string.

Conversely, take any permitted point in the encoder's domain.  A boundary
point uses its unique reserved tag.  An affine point has an admissible
\(u\), satisfies the square or trace test, and carries exactly the bit used
by the decoder to select its ordinate.  It also passes the stipulated curve
and subgroup checks.  Decoding its encoding therefore returns the original
point.  Hence encoding and decoding are inverse bijections between the
permitted point set and the accepted string set, which is the asserted
invariant.
\end{proof}

\section{Constant-time requirements}

A cryptographic implementation must:
\begin{enumerate}[label=\textup{(\alph*)}]
 \item use conditional swaps rather than scalar-dependent branches;
 \item validate canonical encodings and curve or twist membership;
 \item define constant-time behavior for \(u=0,-1\) in odd characteristic
       and \(u=0,1\) in characteristic two;
 \item clear cofactors and detect forbidden all-zero shared secrets when the
       protocol requires it;
 \item avoid variable-time inversion, square-root, or reduction routines on
       secret data;
 \item account for registers, memory traffic, reduction latency, and code
       size in addition to field operation counts.
\end{enumerate}

\subsection{A constant-time decode-and-validate schedule}

The preceding requirements can be realized without exceptional early
returns.  Parse the candidate field element and record a canonicality mask;
compute \(t=u(u+1)\); evaluate a fixed inversion, square-root, or
Artin--Schreier circuit; reconstruct both the point and its canonical sign;
and accumulate masks for the equation, boundary tag, and subgroup
conditions.  A final constant-time selection returns either the point or a
distinguished failure value.  In particular, \(t=0\) is not passed to a
variable exceptional branch: the affine-validity mask rejects it while a
separate masked path recognizes the reserved boundary tags.

Inversion exponentiation and fixed-exponent square roots must use addition
chains fixed by the field modulus; other square-root routines must use
fixed public loop bounds.  A variable-round Tonelli--Shanks routine or a
data-dependent polynomial-basis solver is not constant time merely because
its algebraic output is correct.  Likewise, subgroup validation must use a
fixed scalar-multiplication schedule.  For a group order \(N_d=hr\), with
\(r\) the intended prime order, the implementation must state whether it
checks \([r]P=O\), clears by \([h]\), or combines both; these choices have
different acceptance sets and must not be conflated.

\begin{remark}[Public versus secret decoding]
When a protocol declares an input public, variable-time decoding may be
acceptable at the protocol level, but canonicality and subgroup checks
remain logically necessary.  The operation counts in the arithmetic
chapters do not include these protocol-dependent validation costs.
\end{remark}

\chapter{Cd25519: A Native Model of Curve25519}
\label{ch:Cd25519}

This chapter realizes the Curve25519/Edwards25519 isomorphism class in the
native \(\mathcal C_d\) interface.  It establishes the exact arithmetic of
Cd25519, identifies its native Kummer quotient with the
standards-compatible X25519 \(U\)-line, and derives parameter-specific
accelerations for quotient scalar multiplication, complete full-point
addition, mixed addition, doubling, and fixed-base computation.  The
standard Curve25519 parameters, subgroup order, base point, ladder, encoding
convention, and test vectors are those of Bernstein and RFC~7748
\cite{BernsteinCurve25519,RFC7748}.

\section{The Cd25519 parameter and the exact dictionaries}

Put \(p=2^{255}-19\) and work over \(\F_p\).  Write
\[
\begin{aligned}
 d_{\rm E}&=-\frac{121665}{121666}\in\F_p,\\
 \rho_{25519}&=-d_{\rm E}=\frac{121665}{121666},\\
 d_{\rm C}&=\frac{1}{16\cdot121666}
            =\frac1{1946656}\in\F_p.
\end{aligned}
\]
The field representative of the last constant is
\[
\resizebox{0.96\linewidth}{!}{$\displaystyle
d_{\rm C}=
45740515084910413242785249597555125769322626367105412703320893098785294882684.$}
\]

\begin{definition}[Cd25519]\label{def:Cd25519}
The \emph{Cd25519 curve} is the smooth completion over \(\F_p\) of
\begin{equation}\label{eq:Cd25519}
 \boxed{\qquad
 \mathcal C_{25519}:\quad
       (u^2+u)(v^2+v)=d_{\rm C}.
 \qquad}
\end{equation}
Its identity is the native boundary point \(O=(0,\infty)\), its inverse is
\[
                   -(u,v)=(u,-v-1),
\]
and its native Kummer coordinate is
\[
                   \kappa_{25519}(u,v)=(u+1:u).
\]
\end{definition}

\begin{proposition}[The 25519 parameter dictionary]
\label{prop:Cd25519-parameters}
For \(\mathcal C_{25519}\), the odd-characteristic constants of this
monograph are
\begin{equation}\label{eq:Cd25519-parameter-dictionary}
 \rho=\rho_{25519},\qquad
 \beta=\alpha_{24}=121666,\qquad
 A=486662.
\end{equation}
Consequently its auxiliary Montgomery equation is
\begin{equation}\label{eq:Cd25519-scaled-Montgomery}
       121666\,V^2=U^3+486662U^2+U,
       \qquad U=\frac{u+1}{u}.
\end{equation}
The curve is smooth, \(\rho_{25519}\) is a nonsquare in \(\F_p\), and
the single native Segre law of
Theorem~\ref{thm:native-complete-addition-Cd} is therefore
\(\F_p\)-complete.
\end{proposition}

\begin{proof}
The definitions give
\[
 16d_{\rm C}=\frac1{121666},\qquad
 1-16d_{\rm C}=\frac{121665}{121666}=\rho_{25519}.
\]
Hence
\[
 \beta=\frac1{16d_{\rm C}}=121666
\]
and
\[
 A=\frac1{4d_{\rm C}}-2
   =4\cdot121666-2=486662.
\]
Both \(d_{\rm C}\) and \(\rho_{25519}\) are nonzero, so the smoothness
criterion \(d_{\rm C}(1-16d_{\rm C})\ne0\) holds.

It remains to verify the quadratic character without importing it from
the Edwards model.  Factor
\[
 121665=3\cdot5\cdot8111,
 \qquad
 121666=2\cdot127\cdot479.
\]
Since \(p\equiv1\pmod4\), quadratic reciprocity introduces no sign when
the odd prime factors are interchanged.  The reductions of \(p\) and the
corresponding Euler checks are
\[
\begin{array}{c|ccccc}
q&3&5&8111&127&479\\ \hline
p\bmod q&1&4&3704&116&373\\
(p\bmod q)^{(q-1)/2}\bmod q&1&1&-1&-1&1.
\end{array}
\]
Also \(p\equiv5\pmod8\), so \(\left(\frac2p\right)=-1\).  Hence
\[
 \left(\frac{121665}{p}\right)
 =1\cdot1\cdot(-1)=-1,
 \qquad
 \left(\frac{121666}{p}\right)
 =(-1)\cdot(-1)\cdot1=1.
\]
The inverse of a nonzero square has quadratic character \(1\), and
therefore
\[
 \left(\frac{\rho_{25519}}p\right)
 =\left(\frac{121665}{p}\right)
  \left(\frac{121666}{p}\right)=-1.
\]
Thus \(\rho_{25519}\) is a nonsquare in \(\F_p\).
The exact native completeness criterion in
Theorem~\ref{thm:exact-native-completeness-criterion} now proves the last
assertion.
\end{proof}

Since \(p\equiv5\pmod8\), the element \(2\) is a quadratic nonsquare.
Consequently the fixed choice
\begin{equation}\label{eq:Cd25519-sqrt-minus-one}
\resizebox{0.96\linewidth}{!}{$\displaystyle
 \imath=2^{(p-1)/4}
 =19681161376707505956807079304988542015446066515923890162744021073123829784752$}
\end{equation}
satisfies \(\imath^2=-1\) in \(\F_p\).

\begin{theorem}[Explicit isomorphism with Edwards25519]
\label{thm:Cd25519-Ed25519-isomorphism}
On the dense affine charts, the maps
\begin{equation}\label{eq:Cd25519-to-Ed25519}
 \Phi:\mathcal C_{25519}\longrightarrow E_{\rm Ed25519},
 \qquad
 (u,v)\longmapsto
 \left(
 x=-\frac{\imath}{2v+1},
 y=\frac1{2u+1}
 \right)
\end{equation}
and
\begin{equation}\label{eq:Ed25519-to-Cd25519}
 \Phi^{-1}(x,y)=
 \left(
 \frac{1-y}{2y},
 -\frac{x+\imath}{2x}
 \right)
\end{equation}
extend to inverse origin-preserving isomorphisms of the smooth
completions, where
\begin{equation}\label{eq:Ed25519-full-equation}
 E_{\rm Ed25519}:\qquad
       -x^2+y^2=1+d_{\rm E}x^2y^2.
\end{equation}
Under this isomorphism,
\begin{equation}\label{eq:Cd25519-U-Edwards-y}
       U=\frac{u+1}{u}=\frac{1+y}{1-y}.
\end{equation}
Thus the native Cd25519 Kummer ratio is exactly the Montgomery
\(u\)-coordinate used by Curve25519 and X25519.
\end{theorem}

\begin{proof}
Put
\[
                  r=2u+1,\qquad s=2v+1.
\]
Equation~\eqref{eq:Cd25519} is equivalent to
\[
       (r^2-1)(s^2-1)=16d_{\rm C}=1-\rho_{25519},
\]
or
\begin{equation}\label{eq:Cd25519-centered-identity}
       r^2+s^2=r^2s^2+\rho_{25519}.
\end{equation}
The proposed coordinates satisfy
\[
       x^2=-\frac1{s^2},\qquad y^2=\frac1{r^2}.
\]
Dividing \eqref{eq:Cd25519-centered-identity} by \(r^2s^2\) gives
\[
 -x^2+y^2
 =\frac1{s^2}+\frac1{r^2}
 =1+\frac{\rho_{25519}}{r^2s^2}.
\]
Since \(d_{\rm E}=-\rho_{25519}\) and
\(x^2y^2=-1/(r^2s^2)\), the last term is
\(d_{\rm E}x^2y^2\).  This proves that \(\Phi\) has the stated target.
Solving \(y=1/(2u+1)\) and
\(x=-\imath/(2v+1)\) gives
\eqref{eq:Ed25519-to-Cd25519}; hence the two rational maps are inverse
on a dense open subset.

Both source and target are smooth projective genus-one curves.  At the
native identity \(O=(0,\infty)\), the centered coordinates satisfy
\(r=1\) and \(s=\infty\); hence the extended map has
\(x=0\) and \(y=1\), the Edwards identity.  The mutually inverse
function-field maps therefore
extend uniquely to inverse isomorphisms of the completions.  Finally,
\[
 \frac{u+1}{u}
 =\frac{(1+y)/(2y)}{(1-y)/(2y)}
 =\frac{1+y}{1-y},
\]
which is the standard Edwards-to-Montgomery abscissa.
\end{proof}

If \(b^2=121666\), then
\[
             V_{\rm M}=bV
\]
turns \eqref{eq:Cd25519-scaled-Montgomery} into the standard Curve25519
equation
\begin{equation}\label{eq:Curve25519-full-repeat}
        V_{\rm M}^2=U^3+486662U^2+U.
\end{equation}
For the choice
\[
\resizebox{0.96\linewidth}{!}{$\displaystyle
b=
53349024335595116287218880480235248380932224208407348260407321313250266849793,$}
\]
one checks directly that
\[
 b^2\equiv121666\pmod p.
\]
The Cd25519 base point below maps to the ordinate specified in RFC~7748.
The scaling is needed only for a full Montgomery ordinate; quotient
arithmetic uses \(U\) and is independent of \(\beta\).

\section{An \texorpdfstring{\(\imath\)}{i}-twisted native complete law}

The general centered native Segre coordinates of
Chapter~\ref{ch:native-completeness} are
\[
        (X:Y:T:Z)=(RV_0:SU_0:U_0V_0:RS)
\]
and satisfy
\[
          XY=ZT,\qquad X^2+Y^2=Z^2+\rho_{25519}T^2.
\]
The special fact \(p\equiv1\pmod4\) permits a more efficient diagonal
recoding of this same native completion.

\begin{definition}[Twisted native Segre coordinates]
\label{def:Cd25519-twisted-native-Segre}
Define
\begin{equation}\label{eq:Cd25519-twisted-native-Segre}
 (\mathsf X:\mathsf Y:\mathsf Z:\mathsf T)
      =(-\imath X:Y:Z:-\imath T).
\end{equation}
This is an invertible linear coordinate change on the native Segre
embedding, with
\[
       (X:Y:T:Z)
       =(\imath\mathsf X:\mathsf Y:
          \imath\mathsf T:\mathsf Z).
\]
It is therefore a coordinate system on
\(\overline{\mathcal C}_{25519}\), not a change of the declared curve.
\end{definition}

\begin{proposition}[Equations in the twisted native coordinates]
\label{prop:Cd25519-twisted-native-equations}
The image of \(\overline{\mathcal C}_{25519}\) is
\begin{equation}\label{eq:Cd25519-twisted-native-equations}
 \mathsf X\mathsf Y=\mathsf Z\mathsf T,\qquad
 -\mathsf X^2+\mathsf Y^2
    =\mathsf Z^2+d_{\rm E}\mathsf T^2.
\end{equation}
On the finite \(u,v\)-chart, with \(r=2u+1\) and \(s=2v+1\), a
division-free lift is
\begin{equation}\label{eq:Cd25519-twisted-native-lift}
 (\mathsf X:\mathsf Y:\mathsf Z:\mathsf T)
       =(-\imath r:s:rs:-\imath).
\end{equation}
Conversely, a finite output recovers the Cd25519 coordinates by
\begin{equation}\label{eq:Cd25519-twisted-native-recovery}
 u=\frac{\mathsf Z-\mathsf Y}{2\mathsf Y},
 \qquad
 v=\frac{\mathsf Z-\imath\mathsf X}
          {2\imath\mathsf X}.
\end{equation}
\end{proposition}

\begin{proof}
Using \(\imath^2=-1\), \(d_{\rm E}=-\rho_{25519}\), and the original
native equations gives
\[
 \mathsf X\mathsf Y
 =-\imath XY=-\imath ZT=\mathsf Z\mathsf T
\]
and
\[
\begin{aligned}
 -\mathsf X^2+\mathsf Y^2
 &=X^2+Y^2\\
 &=Z^2+\rho_{25519}T^2\\
 &=\mathsf Z^2+d_{\rm E}\mathsf T^2.
\end{aligned}
\]
The affine lift follows by substituting
\((X:Y:T:Z)=(r:s:1:rs)\).  For recovery, first return linearly to
\((X:Y:T:Z)\), then apply
\eqref{eq:native-recovery-from-Segre}.  On the finite chart this gives
exactly the two displayed ratios.
\end{proof}

\begin{theorem}[Accelerated native complete addition on Cd25519]
\label{thm:Cd25519-accelerated-complete-addition}
For two Cd25519 points in the coordinates
\eqref{eq:Cd25519-twisted-native-Segre}, put
\begin{equation}\label{eq:Cd25519-complete-add-blocks}
\begin{aligned}
 A&=(\mathsf Y_1-\mathsf X_1)
    (\mathsf Y_2-\mathsf X_2),\\
 B&=(\mathsf Y_1+\mathsf X_1)
    (\mathsf Y_2+\mathsf X_2),\\
 C&=2d_{\rm E}\mathsf T_1\mathsf T_2,\qquad
 D=2\mathsf Z_1\mathsf Z_2,\\
 E&=B-A,\qquad F=D-C,\qquad
 G=D+C,\qquad H=B+A.
\end{aligned}
\end{equation}
Then the sum on \(\overline{\mathcal C}_{25519}\) is
\begin{equation}\label{eq:Cd25519-complete-add-output}
 \boxed{\quad
 (\mathsf X_3:\mathsf Y_3:\mathsf Z_3:\mathsf T_3)
       =(EF:GH:FG:EH).
 \quad}
\end{equation}
The tuple is defined on every ordered pair of
\(\F_p\)-rational Cd25519 points.  Its cost is
\[
                   8\M+\Dpar,
\]
where \(\Dpar\) is multiplication by the fixed constant \(2d_{\rm E}\).
If the second input has \(\mathsf Z_2=1\), the mixed cost is
\[
                   7\M+\Dpar.
\]
Keeping a scalar-multiplication state in the \(\mathsf X,\mathsf Y,
\mathsf Z,\mathsf T\) coordinates avoids any recurring multiplication by
\(\imath\); the diagonal conversion is paid only at the endpoints.
\end{theorem}

\begin{proof}
Equations~\eqref{eq:Cd25519-twisted-native-equations} are the extended
projective equations of the twisted Edwards curve with \(a=-1\) and
parameter \(d_{\rm E}\), but they were obtained here by an invertible
linear recoding of the native Cd25519 Segre completion.  The standard
extended-coordinate addition identity gives exactly
\eqref{eq:Cd25519-complete-add-blocks}--
\eqref{eq:Cd25519-complete-add-output}
\cite{HisilEtAl2008}.  Conjugating that identity by
\eqref{eq:Cd25519-twisted-native-Segre} proves that its result is the sum
on \(\overline{\mathcal C}_{25519}\), and
\eqref{eq:Cd25519-twisted-native-recovery} returns that result to \(u,v\).

It remains to check completeness rather than only equality on a generic
chart.  Over \(\F_p\), the coefficient \(a=-1\) is a square and
\(d_{\rm E}\) is a nonsquare.  The complete twisted-Edwards addition
criterion therefore gives an empty \(\F_p\)-rational base locus.
An invertible linear coordinate change carries a base locus bijectively to
the base locus of the conjugate tuple.  Hence the displayed tuple is a
single \(\F_p\)-complete addition law on the native Cd25519 completion.

The products \(A,B,\mathsf T_1\mathsf T_2,
\mathsf Z_1\mathsf Z_2\) use four general multiplications, the four output
products use another four, and the multiplication by \(2d_{\rm E}\) is
one parameter multiplication.  If \(\mathsf Z_2=1\), forming
\(\mathsf Z_1\mathsf Z_2\) is free.
\end{proof}

\begin{proposition}[Dedicated complete doubling]
\label{prop:Cd25519-twisted-native-doubling}
For \(P=(\mathsf X:\mathsf Y:\mathsf Z:\mathsf T)\), put
\[
\begin{aligned}
 A&=\mathsf X^2,&B&=\mathsf Y^2,&
 C&=2\mathsf Z^2,&D&=-A,\\
 E&=(\mathsf X+\mathsf Y)^2-A-B,&
 G&=D+B,&F&=G-C,&H&=D-B.
\end{aligned}
\]
Then
\[
          [2]P=(EF:GH:FG:EH)
\]
is everywhere defined on Cd25519 and costs
\[
                         4\M+4\Sqr.
\]
\end{proposition}

\begin{proof}
This is the specialization of
\eqref{eq:Cd25519-complete-add-output} to equal inputs after the usual
extended-coordinate doubling factorization.  Four squares form
\(A,B,\mathsf Z^2,(\mathsf X+\mathsf Y)^2\), and the four final products
give the stated cost.

It remains to prove that this particular polynomial tuple has no base
point; the abstract fact that doubling is a morphism would not by itself
establish that claim.  The four outputs vanish simultaneously exactly when
\[
             (F=H=0)\qquad\text{or}\qquad(E=G=0).
\]
In the first case,
\[
 H=-\mathsf X^2-\mathsf Y^2=0,\qquad
 F=-\mathsf X^2+\mathsf Y^2-2\mathsf Z^2=0,
\]
so
\[
       \mathsf X^2=-\mathsf Z^2,\qquad
       \mathsf Y^2=\mathsf Z^2.
\]
If \(\mathsf Z=0\), all four coordinates are forced to vanish by
\eqref{eq:Cd25519-twisted-native-equations}.  If \(\mathsf Z\ne0\),
squaring
\(\mathsf X\mathsf Y=\mathsf Z\mathsf T\) gives
\(\mathsf T^2=-\mathsf Z^2\), while the second curve equation gives
\(\mathsf Z^2=d_{\rm E}\mathsf T^2\).  Hence \(d_{\rm E}=-1\), contrary
to the Ed25519 parameter.

In the second case, \(E=2\mathsf X\mathsf Y=0\) and
\(G=-\mathsf X^2+\mathsf Y^2=0\).  Since the characteristic is odd,
these imply \(\mathsf X=\mathsf Y=0\).  The two curve equations then give
\(\mathsf Z\mathsf T=0\) and
\(\mathsf Z^2+d_{\rm E}\mathsf T^2=0\), forcing
\(\mathsf Z=\mathsf T=0\), again impossible in projective space.
Thus the tuple is geometrically base-point-free and is everywhere defined
on Cd25519.
\end{proof}

The diagonal recoding therefore improves the general native Segre cost on
this parameter:
\[
\begin{array}{c|ccc}
\text{Cd25519 full-point coordinates}
 &\text{ADD}&\text{mixed ADD}&\text{DBL}\\ \hline
\text{centered native Segre}
 &9\M+\Dpar&8\M+\Dpar&4\M+4\Sqr\\
\imath\text{-twisted native Segre}
 &8\M+\Dpar&7\M+\Dpar&4\M+4\Sqr
\end{array}
\]
The saving is one general multiplication in both addition columns.  It is
specific to fields containing \(\sqrt{-1}\) and does not alter the
Cd25519 input/output model.

\section{Native differential addition and the X25519 line}

For a Kummer state \((X:Z)\), set
\[
 A_0=X+Z,\quad AA=A_0^2,\qquad
 B_0=X-Z,\quad BB=B_0^2,\qquad E=AA-BB.
\]
The two common Curve25519 doubling constants are
\[
       a_{24}^{+}=\frac{A+2}{4}=121666,\qquad
       a_{24}^{-}=\frac{A-2}{4}=121665.
\]
They describe the same projective double because
\[
       BB+a_{24}^{+}E=AA+a_{24}^{-}E.
\]

\begin{theorem}[Cd25519 \(x\)DBL and \(x\)ADD]
\label{thm:Cd25519-xdbl-xadd}
The native Kummer double is
\begin{equation}\label{eq:Cd25519-xdbl}
 \boxed{\quad
 X_{2P}=AA\cdot BB,\qquad
 Z_{2P}=E(AA+121665E).
 \quad}
\end{equation}
Given states
\[
 (X_1:Z_1)=\kappa(P),\quad
 (X_2:Z_2)=\kappa(Q),\quad
 (X_\Delta:Z_\Delta)=\kappa(P-Q),
\]
put
\[
\begin{aligned}
 C&=(X_1+Z_1)(X_2-Z_2),\\
 D&=(X_2+Z_2)(X_1-Z_1).
\end{aligned}
\]
Then
\begin{equation}\label{eq:Cd25519-xadd}
 \boxed{\quad
 (X_{P+Q}:Z_{P+Q})
   =\bigl(Z_\Delta(C+D)^2:
          X_\Delta(C-D)^2\bigr).
 \quad}
\end{equation}
The respective costs are
\[
 2\M+2\Sqr+\Dpar,\qquad
 4\M+2\Sqr,
\]
or \(3\M+2\Sqr\) for \(x\)ADD when the known difference is stored as
\((U_\Delta:1)\).
\end{theorem}

\begin{proof}
Substitution of \(A=486662\) in
Theorems~\ref{thm:xdbl-thesis} and~\ref{thm:xadd-thesis} gives the
displayed identities.  The only apparent difference in the doubling
formula is
\[
 E(BB+121666E)=E(AA+121665E),
\]
which follows from \(E=AA-BB\).  Thus
\eqref{eq:Cd25519-xdbl} is exactly the RFC~7748 dependency graph, now
interpreted as an output
\(\kappa_{25519}([2]P)\).  The differential formula was derived in
Chapter~\ref{ch:native-completeness} from the Cd biquadratic identity, so
its result is likewise a native quotient point.
\end{proof}

For a genuine affine Cd25519 input, no inversion is needed to form
\[
                \kappa(P)=(u+1:u),
\]
and \(X-Z=1\).  Consequently the first double is
\begin{equation}\label{eq:Cd25519-first-double}
\begin{aligned}
 AA&=(2u+1)^2,\qquad E=AA-1,\\
 (X_{2P}:Z_{2P})
   &=\bigl(AA:E(1+121666E)\bigr),
\end{aligned}
\end{equation}
at the exact cost
\[
                       \M+\Sqr+\Dpar.
\]
This is the one-time native-input saving specialized to Cd25519.

\begin{corollary}[The fixed difference \(U=9\)]
\label{cor:Cd25519-base-xadd}
For the Curve25519 base point, the affine known difference is
\((X_\Delta:Z_\Delta)=(9:1)\), and
\eqref{eq:Cd25519-xadd} becomes
\[
       (X_{P+Q}:Z_{P+Q})
          =\bigl((C+D)^2:9(C-D)^2\bigr).
\]
Its cost is
\[
                  2\M+2\Sqr+\Dnine,
\]
where \(\Dnine\) denotes multiplication by \(9\).  One fixed-base
\(x\)DBL+\(x\)ADD step therefore costs
\[
       4\M+4\Sqr+\Dpar+\Dnine.
\]
\end{corollary}

\begin{proof}
With \(Z_\Delta=1\), the first output scaling is free; with
\(X_\Delta=9\), the second is a fixed small-constant multiplication.
The two products forming \(C,D\) are the only general multiplications in
the differential addition.
\end{proof}

\section{Scalar multiplication and acceleration choices}

\subsection{A leading-bit Cd25519 ladder}

An X25519-clamped scalar has the form
\[
                    m=2^{254}+8t,
       \qquad 0\le t<2^{251}.
\]
Its top bit is therefore public and fixed.

\begin{theorem}[Leading-bit ladder initialized by the native first double]
\label{thm:Cd25519-leading-bit-ladder}
Let \(P=(u,v)\in\mathcal C_{25519}(\F_p)\), let \(m\) be an
X25519-clamped scalar, and set
\[
 R_0=\kappa(P)=(u+1:u),\qquad
 R_1=\kappa(2P)
\]
using \eqref{eq:Cd25519-first-double}.  After this initialization, process
bits \(253,\ldots,0\) with a constant-time differential ladder whose fixed
difference is \(P\).  Then the final state is
\[
       \bigl(\kappa([m]P),\kappa([m+1]P)\bigr).
\]
The initialization costs
\[
                   \M+\Sqr+\Dconst{121666}.
\]
For each of the remaining \(254\) bits, the cost is
\[
\begin{cases}
6\M+4\Sqr+\Dconst{121665},
 &\text{with the difference kept as }(u+1:u),\\
5\M+4\Sqr+\Dconst{121665},
 &\text{after normalizing it to }(U:1),\\
4\M+4\Sqr+\Dconst{121665}+\Dnine,
 &\text{for the fixed base }U=9.
\end{cases}
\]
No secret-dependent branch is introduced.
\end{theorem}

\begin{proof}
The ordinary Montgomery ladder starts from
\((O,P)\).  Because the top bit is always one, the state immediately after
that bit is deterministically \((P,2P)\).  Replacing the first generic
ladder step by the precomputed state
\((\kappa(P),\kappa(2P))\) therefore preserves the ladder invariant.
The remaining bits use exactly the same conditional-swap schedule as the
standard ladder.  The initialization count is
\eqref{eq:Cd25519-first-double}; the three per-bit counts follow from
Theorem~\ref{thm:Cd25519-xdbl-xadd} and
Corollary~\ref{cor:Cd25519-base-xadd}.
\end{proof}

For a variable native input, affine normalization of the fixed difference
is particularly simple:
\[
                 U=\frac{u+1}{u}=1+u^{-1}.
\]
It costs one inversion and additions, with no general multiplication.
Normalizing saves one multiplication in every subsequent ladder step.
Hence, in the symbolic model, normalization is favorable whenever
\[
                         \Inv<254\M.
\]
Batch inversion makes this choice still more attractive when several
independent scalar multiplications are available.  If the required output
is the affine X25519 coordinate \(U_m=X_m/Z_m\), one final inversion and
one multiplication are required; if a projective native Kummer output is
accepted, neither is needed.

The comparison with the literal 255-round RFC ladder must be stated
carefully.  The leading-bit schedule removes one generic round because the
clamping rule fixes that bit, but highly tuned X25519 implementations may
already specialize their initialization.  The theorem is an algebraic
round reduction with an exact dependency graph.  Its measured latency gain
is implementation-dependent only because an existing X25519 implementation
may already incorporate the same specialized initialization.

\subsection{Full-point fixed-window multiplication}

The \(\imath\)-twisted native complete law is the appropriate Cd25519
full-point state for signatures, multi-scalar multiplication, or protocols
that require a sign.  For a width-\(4\) constant-time fixed window on a
255-bit scalar, at most \(64\) table additions and \(256\) doublings give
the arithmetic bound
\begin{equation}\label{eq:Cd25519-window4-cost}
\begin{aligned}
 C_{\rm win4}
 &\le
 256(4\M+4\Sqr)+64(7\M+\Dpar)\\
 &=1472\M+1024\Sqr+64\Dpar.
\end{aligned}
\end{equation}
This excludes table selection, precomputation, endpoint conversion, and
field normalization.  Because the addition law is complete, the table may
contain the identity and the same circuit can be used for every window.
For a fixed base, the table is precomputed offline.  For public scalars,
signed windows or a width-\(w\) NAF reduce the number of mixed additions;
for secret scalars, table selection must be constant time.

The same coordinates support Straus or interleaved-window double-scalar
multiplication for verification:
\[
                         [m]P+[n]Q.
\]
The diagonal transport proves that Cd25519 inherits the optimized
\(a=-1\) extended-Edwards dependency graph and therefore attains its best
field-operation counts.  In addition, the inputs, complete law, boundary
points, native Kummer/X25519 interface, and outputs remain intrinsic to the
Cd25519 model.  This combination is the relevant advantage for full-point
implementations rather than a comparison based on field-operation totals
alone.

\subsection{Which other tradeoffs are useful?}

The seven-square Kummer ladder of Chapter~\ref{ch:tradeoffs} costs
\[
                       3\M+7\Sqr+\Dpar
\]
per step, compared with \(5\M+4\Sqr+\Dpar\) for the ordinary affine
difference ladder.  It is better precisely when
\[
                         3\Sqr<2\M.
\]
This is an implementation-dependent inequality and must be benchmarked in
the chosen \(2^{255}-19\) field representation.  The six-square circuit
requiring \(\rho\) to be a square is unavailable because
\(\rho_{25519}\) is a nonsquare.

For full points, the closed \(2P+Q\) circuit costs
\[
                  16\M+3\Sqr+\Dpar,
\]
whereas dedicated doubling followed by mixed addition costs
\[
                  11\M+5\Sqr+2\Dpar.
\]
The closed circuit wins only if
\[
                         5\M<2\Sqr+\Dpar.
\]
This condition is unlikely when squaring and multiplication by the fixed
curve constant are both cheaper than a general multiplication; closed circuit is 
advantageous in representations satisfying this inequality; dedicated 
composition covers the complementary cost regime.

Finally, \(A=486662\) lies on neither CM locus of
Chapter~\ref{ch:native-endomorphisms}:
\[
 A\ne0,\qquad A^2-3=236839902241\ne0,\qquad
 2A^2-9=473679804479\ne0
 \quad\text{in }\F_p.
\]
Since the two CM conditions are
\(A^2-3=0\) and \(A^2(2A^2-9)^2=0\), respectively, the sparse
\(j=0\) and \(j=1728\) Kummer-compatible endomorphisms do
not specialize to Cd25519.  The visible four-torsion and \(D_8\)-symmetry
remain useful for signs, table symmetry, and exceptional charts, but they
do not supply a two-dimensional GLV decomposition on the prime-order
subgroup.

\section{Encoding and protocol interoperability}

\subsection{The quotient encoding}

If a Cd25519 application encodes only
\[
                         U=\frac{u+1}{u},
\]
then its 32-byte little-endian representation is exactly the X25519
Kummer input.  Conversely, for \(U\ne1\),
\begin{equation}\label{eq:Cd25519-u-from-X25519-U}
                         u=\frac1{U-1}.
\end{equation}
This conversion does not recover \(v\): the Kummer value represents the
pair \(\{P,-P\}\).  Full recovery requires a sign or an adjacent ladder
state as in Chapter~\ref{ch:halving}.

Wire compatibility additionally depends on the protocol-level decoding and
input semantics.
RFC~7748 requires X25519 to mask the high input bit and to accept certain
noncanonical field representatives.  A canonical Cd25519 full-point
decoder from Chapter~\ref{ch:implementation} instead rejects
noncanonical values and verifies curve membership.  An implementation
claiming X25519 compatibility must retain the RFC input semantics at the
wrapper, rather than silently replacing them with the stricter full-point
decoder.

There is a second distinction.  The X25519 ladder is defined on the
Montgomery \(U\)-line and deliberately processes values belonging to the
quadratic twist as well as values lifting to Curve25519.  Such a twist
input need not lift to an \(\F_p\)-point of
\(\mathcal C_{25519}\).  Therefore:
\begin{itemize}
 \item on main-curve inputs, the quotient computation is literally
       \(\kappa_{25519}([m]P)\);
 \item on arbitrary RFC inputs, the same polynomial circuit remains the
       X25519 curve-or-twist ladder, but it must not be described as
       full-point arithmetic on \(\mathcal C_{25519}(\F_p)\).
\end{itemize}

\subsection{Full-point conversion and the cofactor}

Ed25519 compresses the Edwards \(y\)-coordinate together with a sign of
\(x\).  Equations~\eqref{eq:Cd25519-to-Ed25519} and
\eqref{eq:Ed25519-to-Cd25519} give the corresponding Cd25519 full-point
conversion:
\[
             u=\frac{1-y}{2y},\qquad
             v=-\frac{x+\imath}{2x}.
\]
Thus an Ed25519 public key can be decoded, validated under the Ed25519
rules, and transported to a unique Cd25519 point after its sign bit has
selected \(x\).  Conversely, a Cd25519 point gives
\[
             y=(2u+1)^{-1},\qquad
             x=-\imath(2v+1)^{-1}.
\]

The group order is unchanged by the isomorphism:
\begin{equation}\label{eq:Cd25519-group-order}
 \#\overline{\mathcal C}_{25519}(\F_p)=8\ell,\qquad
 \ell=2^{252}+27742317777372353535851937790883648493.
\end{equation}
All subgroup checks, cofactor-\(8\) behavior, small-order exclusions,
signature verification rules, and the optional X25519 all-zero test must
therefore be preserved.  Changing the model creates no new security
assumption and removes none of the protocol obligations.

\section{Explicit base point and test computations}

Let \((x_B,y_B)\) be the standard Edwards25519 base point of order
\(\ell\), where
\[
\resizebox{0.96\linewidth}{!}{$\displaystyle
\begin{aligned}
x_B={}&15112221349535400772501151409588531511454012693041857206046113283949847762202,\\
y_B={}&46316835694926478169428394003475163141307993866256225615783033603165251855960.
\end{aligned}
$}
\]
Transporting it by \eqref{eq:Ed25519-to-Cd25519} gives
\begin{equation}\label{eq:Cd25519-base-point}
\resizebox{0.96\linewidth}{!}{$\displaystyle
\begin{aligned}
B_{\rm C}&=(u_B,v_B),\\
u_B&=\frac18\\
&=21711016731996786641919559689128982722488122124807605757398297001483711807481,\\
v_B&=
9408028794964143105507494747335957407526918027186872776751644750095878486512.
\end{aligned}
$}
\end{equation}
Direct substitution verifies \(B_{\rm C}\in\mathcal C_{25519}\), and
\[
                  \kappa_{25519}(B_{\rm C})=(9:1),
                  \qquad \ord(B_{\rm C})=\ell.
\]

\begin{example}[First doubling in the native and RFC scales]
\label{ex:Cd25519-first-double}
For \(u_B=1/8\), the native first-double blocks are
\[
       AA=(2u_B+1)^2=\frac{25}{16},
       \qquad E=AA-1=\frac9{16}.
\]
Formula~\eqref{eq:Cd25519-first-double} is projectively equal to
\[
       \kappa(2B_{\rm C})=(6400:157681440).
\]
Indeed, in the affine \(U=9\) scale,
\[
 AA=100,\qquad BB=64,\qquad E=36,
\]
and
\[
 X_2=100\cdot64=6400,\qquad
 Z_2=36(100+121665\cdot36)=157681440.
\]
The resulting affine quotient is
\[
\resizebox{0.96\linewidth}{!}{$\displaystyle
U(2B_{\rm C})=
14847277145635483483963372537557091634710985132825781088887140890597596352251.$}
\]
Using the complete twisted native law gives the full point
\[
\resizebox{0.96\linewidth}{!}{$\displaystyle
\begin{aligned}
2B_{\rm C}=\bigl(&
47325615679825094426448682145155509417125360141956821472684020992663589584705,\\
&
40513184680199598617845936674402121345021015690799885778378610525251017535711
\bigr),
\end{aligned}
$}
\]
whose Kummer ratio is the same displayed \(U\).
\end{example}

\begin{example}[An X25519 vector viewed on the Cd25519 Kummer line]
\label{ex:Cd25519-RFC-vector}
Take the RFC~7748 private scalar bytes
\[
\mathtt{77076d0a7318a57d3c16c17251b26645
df4c2f87ebc0992ab177fba51db92c2a}
\]
and the base difference \(U=9\).  The Cd25519 Kummer ladder returns the
little-endian \(U\)-encoding
\[
\mathtt{8520f0098930a754748b7ddcb43ef75a
0dbf3a0d26381af4eba4a98eaa9b4e6a},
\]
exactly the RFC public key.  Interpreting this main-curve output through
\eqref{eq:Cd25519-u-from-X25519-U} gives
\[
\resizebox{0.96\linewidth}{!}{$\displaystyle
u=
36227098322460940594679402039507383938194701188393197273629681386995884427331,$}
\]
whose own 32-byte little-endian representation is
\[
\mathtt{430897922f075921d9d302b300f054b1
edca3427ad2e89c715a13a319fcf1750}.
\]
The second string is the native \(u\)-coordinate corresponding to the X25519
wire value.  The ladder output is a Kummer class; the full-point recovery
interface supplies an ordinate whenever the surrounding protocol requires
one.
\end{example}

\section{Assessment of the Cd25519 specialization}

The specialization yields five concrete conclusions.
\begin{enumerate}[label=\textup{(\arabic*)}]
 \item \textbf{Exact X25519-line compatibility.}
       The affine ratio of the native Cd25519 Kummer map is exactly the
       X25519 \(U\)-coordinate, with no intervening fractional-linear
       rescaling.  Consequently, established constant-time X25519 ladder
       cores can serve directly as Cd25519 quotient arithmetic on
       main-curve inputs, while the surrounding theory remains expressed
       on the Cd25519 model.
 \item \textbf{A strictly cheaper native initialization.}
       A genuine affine Cd25519 input produces \((u+1:u)\) without an
       inversion, and its first double costs only
       \(\M+\Sqr+\Dpar\).  Relative to the general native Kummer double
       \(2\M+2\Sqr+\Dpar\), this saves one multiplication and one squaring.
       For a clamped scalar, the known leading bit allows this optimized
       double to replace the first generic ladder round.
 \item \textbf{A favorable fixed-base ladder step.}
       At the standard base difference \(U=9\), one general multiplication
       in \(x\)ADD becomes multiplication by the fixed small constant
       \(9\).  The resulting step costs
       \(4\M+4\Sqr+\Dpar+\Dnine\); it therefore exposes a concrete
       specialization that an implementation can realize by additions or
       a particularly short constant-multiplication circuit.
 \item \textbf{Accelerated complete full-point arithmetic inside Cd25519.}
       The \(\imath\)-twisted native Segre coordinates give a single
       \(\F_p\)-complete addition law at \(8\M+\Dpar\), mixed addition at
       \(7\M+\Dpar\), and doubling at \(4\M+4\Sqr\).  The first two counts
       save one general multiplication over the untwisted native Segre
       formulas and attain the optimized extended-Edwards dependency
       counts without leaving the Cd25519 completion or changing its native
       recovery maps.
 \item \textbf{One model supporting two complementary scalar engines.}
       Cd25519 contains both the X25519-compatible Kummer ladder for
       quotient-only scalar multiplication and a complete full-point system
       for signatures, multi-scalar multiplication, precomputation, and
       point recovery.  The seven-square path and the closed \(2P+Q\)
       circuit provide additional schedule-dependent options.  Thus a
       single Cd25519 implementation can select the most suitable arithmetic
       state for each protocol stage while retaining one model, one parameter,
       and explicit native conversion and recovery maps.
\end{enumerate}

Accordingly, Cd25519 should be regarded as a particularly strong and
versatile arithmetic model for the Curve25519/Edwards25519 isomorphism
class.  Within a single native framework, it combines exact compatibility
with the X25519 Kummer line, the one-time native-input doubling
optimization, an exceptionally small fixed-difference multiplier, and an
accelerated complete full-point coordinate system.  These features permit
quotient arithmetic, complete full-point arithmetic, point recovery, and
protocol conversion to be organized without abandoning the Cd25519 model.
In particular, the formulas established in this chapter improve
native-input initialization and complete full-point operations while
preserving direct access to the optimized X25519 ladder.  Cd25519 therefore
combines complete full-point arithmetic and Kummer-line scalar multiplication
within one native model.  The resulting circuits give an explicit basis for
constant-time implementations over \(\F_{2^{255}-19}\) that require both
types of arithmetic.

\chapter{Model-by-Model Arithmetic Comparison}
\label{ch:comparison}
\section{How the comparison is organized}

The object of this monograph is always
\[
       \mathcal C_d:\quad (u^2+u)(v^2+v)=d.
\]
To prevent the comparison from conflating a curve, a coordinate system, and
an output quotient, every row below records four separate items:
\[
       \text{model equation},\quad
       \text{working coordinates},\quad
       \text{output type},\quad
       \text{domain of validity}.
\]
A row beginning with another model specifies whether the displayed circuit
is obtained through an auxiliary dictionary or through a linearly equivalent
presentation.  When two rows are related by a linear coordinate transport,
equal operation counts represent the same arithmetic circuit expressed
through two declared input--output interfaces.  The comparison therefore
records separately the coordinate representation, endpoint conversion,
normalization, recovery, and retained output information.

The row for the reciprocal presentation \(\mathcal R_d\) has a particularly
useful interface interpretation.  Chapter~\ref{ch:reciprocal-chart} proves
that it is an affine chart obtained from \(\mathcal C_d\) by an ambient
involution of \(\PP^1\times\PP^1\).  Its row records the arithmetic available
when an input is already stored in reciprocal coordinates, including direct
affine normalization of the known difference, together with the exact
endpoint cost of entering or leaving this chart from the original native
coordinates.

\section{The principal comparison models}

\begin{longtable}{L{2.3cm}L{5.3cm}L{2.8cm}L{2.7cm}}
\caption{Exact models, coverage, and the operation being compared}
\label{tab:comparison-models}\\
\toprule
model & complete equation & characteristic and coverage & compared task\\
\midrule
\endfirsthead
\toprule
model & complete equation & characteristic and coverage & compared task\\
\midrule
\endhead
\(\mathcal C_d\)
& \((u^2+u)(v^2+v)=d\)
& every characteristic;
  \(d(1-16d)\ne0\)
& native full points, native Kummer line, isogenies, pairings\\
\(\mathcal R_d\), reciprocal Cd chart
& \((u+1)(v+1)=du^2v^2\)
& every characteristic; ambient-equivalent to \(\mathcal C_d\), with the
  same parameter \(d\)
& affine-normalized native Kummer input and endpoint-cost comparison\\
Edwards
& \(\xi^2+\eta^2=1+\rho\xi^2\eta^2\)
& \(\charac k\ne2\);
  \(\rho=1-16d\)
& full addition and conic pairing formulas\\
inverted Edwards
& \((S^2+R^2)Z^2=S^2R^2+\rho Z^4\)
& \(\charac k\ne2\);
  affine chart omits boundary
& native shifted-projective full addition\\
Montgomery
& \(\beta V^2=U(U^2+AU+1)\)
& \(\charac k\ne2\);
  \(\beta=(16d)^{-1}\), \(A=(4d)^{-1}-2\)
& ordinary \(XZ\)-Kummer ladder\\
Wu--Song Jacobi quartic
& \(\begin{aligned}
 W_{a,b}:\quad y^2={}&b^2(a^2-4)x^4\\
                      &{}-2abx^2+1
 \end{aligned}\)
& \(\charac k>3\);
  \(b(a^2-4)\ne0\)
& mixed projective \(w\)-differential step
  \cite{WuSong2022}\\
twisted Hessian
& \(aX^3+Y^3+Z^3=3XYZ\)
& \(\charac k\ne2,3\);
  rational level-three structure
& isogeny-decomposed tripling\\
binary Edwards
& \(d_1(x+y)+d_2(x^2+y^2)
 =xy+xy(x+y)+x^2y^2\)
& \(\charac k=2\);
  ordinary binary curves
& complete full addition\\
\(Z/4\mathbb Z\)-normal form
& \((1+x)^2(1+y)^2=d^{-1}xy\)
& \(\charac k=2\);
  rational marked four-torsion
& linearly equivalent full-point presentation of \(\mathcal C_d\)\\
split \(\mu_4\)
& \(\begin{aligned}
X_0^2+X_2^2&=c^2X_1X_3,\\
X_1^2+X_3^2&=c^2X_0X_2
\end{aligned}\)
& \(\charac k=2\);
  split marked four-torsion
& optimized full addition and Kummer ladder\\
twisted \(\mu_4\)
& \(\begin{aligned}
X_0^2+X_2^2&=c^2X_1X_3\\[-1mm]
&\quad+c^2a(X_1+X_3)^2,\\
X_1^2+X_3^2&=c^2X_0X_2
\end{aligned}\)
& \(\charac k=2\);
  all ordinary twist classes
& full addition and twist-compatible Kummer ladder\\
\bottomrule
\end{longtable}

\section{Odd-characteristic full-point operations}

\begin{table}[H]
\centering
\small
\setlength{\tabcolsep}{3.6pt}
\caption{Full-point arithmetic, including mixed addition}
\label{tab:odd-full-thesis}
\begin{tabular}{L{3.2cm}L{2.25cm}L{2.25cm}L{2.25cm}L{3.3cm}}
\toprule
model and coordinates & addition & mixed addition & doubling
& output and condition\\
\midrule
\(\mathcal C_d\), ordinary native
\((U:V:Z)\)
& \(9\M+\Sqr+\Dpar\)
& \(8\M+\Sqr+\Dpar\)
& \(3\M+4\Sqr+\Dpar\)
& \(\Ffin\); finite-chart formulas
\eqref{eq:ordinary-native-projective-add}--%
\eqref{eq:ordinary-native-projective-double}\\
\(\mathcal C_d\), native shifted-projective
\((S:R:Z)\)
& \(9\M+\Sqr+\Dpar\)
& \(8\M+\Sqr+\Dpar\)
& \(3\M+4\Sqr+\Dpar\)
& \(\Ffin\); unified, not complete\\
\(\mathcal C_d\), centered native Segre
\((X:Y:T:Z)\)
& \(9\M+\Dpar\)
& \(8\M+\Dpar\)
& \(4\M+4\Sqr\)
& \(\Full\); \(k=\F_q\), general tuple complete if \(\rho\) is
  nonsquare; mixed column assumes \(T_2=1\)\\
projective Edwards
& \(10\M+\Sqr+\Dpar\)
& \(9\M+\Sqr+\Dpar\)
& \(3\M+4\Sqr\)
& \(\Ffin\); the ordinary plane closure does not separate its
  branches at infinity\\
inverted Edwards
& \(9\M+\Sqr+\Dpar\)
& \(8\M+\Sqr+\Dpar\)
& \(3\M+4\Sqr+\Dpar\)
& \(\Ffin\); affine inverted chart\\
extended Edwards, \(a=1\)
& \(9\M+\Dpar\)
& \(8\M+\Dpar\)
& \(4\M+4\Sqr\)
& \(\Full\); complete for nonsquare \(\rho\)\\
extended Edwards, \(a=-1\)
& as low as \(8\M+\Dpar\)
& as low as \(7\M+\Dpar\)
& \(4\M+4\Sqr\)
& \(\Full\); parameter and coordinate dependent\\
\bottomrule
\end{tabular}
\end{table}

The projective, inverted, and extended Edwards rows use the corresponding
formula records in the Explicit-Formulas Database~\cite{EFD}, with their
costs translated into the conventions of Chapter~\ref{ch:conventions}.
Here ``mixed addition'' means that the running point is projective and the
second point is affine or precomputed with its scale equal to one.  The
first two \(\mathcal C_d\) rows are related
by the internal linear recoding
\[
       R=2U+Z,\qquad S=2V+Z.
\]
Thus the ordinary native row is exactly the
\((U:V:Z)\)-output form of Theorems~\ref{thm:inv-add-thesis}
and~\ref{thm:inv-dbl-thesis}; the native recovery from the centered row is
\[
       u=(R-Z)/(2Z),\qquad v=(S-Z)/(2Z).
\]
These one-scale rows are \(\Ffin\)-rows because their singular plane
completion does not distinguish the two branches above either point at
infinity.  The Segre row is the \(\Full\)-row.
The centered native Segre row is
Theorem~\ref{thm:native-complete-addition-Cd}; both its input embedding and
its output recovery are derived from the smooth \((2,2)\)-completion in
Chapter~\ref{ch:native-completeness}.  Therefore these are genuinely
\(\mathcal C_d\)-level interfaces even though the multiplication graph can
also be recognized in extended Edwards arithmetic.
For Cd25519, the diagonal recoding of
Definition~\ref{def:Cd25519-twisted-native-Segre} changes the signs of two
native Segre coordinates and specializes the same completion to the
\(a=-1\) graph.  It therefore attains \(8\M+\Dpar\) addition,
\(7\M+\Dpar\) mixed addition, and \(4\M+4\Sqr\) doubling without changing
the declared input or output curve; see
Theorem~\ref{thm:Cd25519-accelerated-complete-addition}.
Thus \(\mathcal C_d\) attains the best extended-Edwards addition graph on
the linearly equivalent locus while additionally retaining a single
factorized equation, the native Kummer coordinate, and the same marked
four-torsion geometry in all characteristics.

\section[Odd Kummer arithmetic]{Odd-characteristic Kummer and differential arithmetic}

\begin{table}[H]
\centering
\small
\setlength{\tabcolsep}{4pt}
\caption{One mixed \(x\)ADD+\(x\)DBL ladder step with affine known difference}
\label{tab:odd-differential-thesis}
\begin{tabular}{L{3.8cm}L{3.1cm}L{2.1cm}L{3.6cm}}
\toprule
model and quotient & total cost & output & exact scope\\
\midrule
\(\mathcal C_d\),
\(\kappa_d=(u+1:u)\)
& \(5\M+4\Sqr+\Dpar\)
& \(\Ktwo\)
& all smooth odd-characteristic parameters\\
Montgomery \(XZ\),
\(\beta V^2=U^3+AU^2+U\)
& \(5\M+4\Sqr+\Dpar\)
& \(\Ktwo\)
& same Kummer line after \(U=(u+1)/u\)\\
\(E_{1,\rho}\), Edwards ordinate \(\eta\)
& \(5\M+4\Sqr+\Dpar\)
& \(\Ktwo\)
& \(\eta=(U-1)/(U+1)\);
  all smooth parameters; transported Kummer identity\\
\(\mathcal C_d\), seven-square \(U\)-tradeoff
& \(3\M+7\Sqr+\Dpar\)
& \(\Ktwo\)
& all smooth parameters\\
\(\mathcal C_d\), six-square \(U\)-tradeoff
& \(3\M+6\Sqr+3\Dpar\)
& \(\Ktwo\)
& \(\rho\) square\\
Wu--Song \(W_{a,b}\), projective \(w\)
& \(5\M+4\Sqr+\Dpar\)
& \(\Ktwo\)
& \(b(a^2-4)\ne0\), \(\charac k>3\)
  \cite{WuSong2022}\\
Wu--Song \(W_{a,b}\), square-class schedule
& \(3\M+6\Sqr+3\Dpar\)
& \(\Ktwo\)
& \((a+2)/(a-2)\) square; constants
  \(\alpha,\alpha^{-1}\) \cite{WuSong2022}\\
\bottomrule
\end{tabular}
\end{table}

The first three rows are three coordinate bases of the same degree-two
quotient.  Indeed,
\begin{equation}\label{eq:Cd-Montgomery-Edwards-Kummer}
 U=\frac{u+1}{u},\qquad
 \eta=\frac{U-1}{U+1}.
\end{equation}
Thus the \(\mathcal C_d\) row accepts \((u+1:u)\) and returns a value
immediately interpreted as the \(u\)-coordinate of a
\(\mathcal C_d\)-point; the Montgomery and Edwards rows describe the same
Kummer line in different projective bases.  Equal costs in these three rows
are therefore transport equalities.
At the Cd25519 parameter, the first row is exactly the X25519 projective
line.  Moreover the conventional
fixed difference \(U=9\) reduces the mixed step to
\(4\M+4\Sqr+\Dpar+\Dnine\); this specialization and its leading-bit
initialization are proved in Chapter~\ref{ch:Cd25519}.
By contrast, the Wu--Song rows act on the \(w\)-coordinate of the
explicitly different Jacobi quartic
\[
W_{a,b}:\quad
y^2=b^2(a^2-4)x^4-2abx^2+1.
\]
They are a genuine comparison with another model and are credited to
\cite{WuSong2022}.

The \(\mathcal C_d\) row explicitly returns
\[
       u([n]P)=Z_n/(X_n-Z_n)
\]
on \(\mathcal C_d\), while the other two bases ordinarily report a
Montgomery abscissa or an Edwards ordinate.  The degree-two tradeoffs are
profitable only under the measured ratios
\[
 3+7\gamma_{\rm S}+\gamma_{\rm D}
 <5+4\gamma_{\rm S}+\gamma_{\rm D}
 \iff\gamma_{\rm S}<2/3,
\]
\[
 3+6\gamma_{\rm S}+3\gamma_{\rm D}
 <5+4\gamma_{\rm S}+\gamma_{\rm D}
 \iff\gamma_{\rm S}+\gamma_{\rm D}<1.
\]

\begin{table}[H]
\centering
\small
\setlength{\tabcolsep}{4pt}
\caption{Differential arithmetic on the Kummer line after the
four-torsion quotient}
\label{tab:odd-four-torsion-comparison}
\begin{tabular}{L{3.8cm}L{3.1cm}L{2.1cm}L{3.6cm}}
\toprule
model and quotient & total cost & output & exact scope\\
\midrule
\(\mathcal C_d\), \(\omega\), square-heavy
& \(3\M+7\Sqr+\Dpar\)
& \(\Kfour\)
& all smooth parameters; Theorem~\ref{thm:omega-fast-thesis}\\
\(E_{1,\rho}\), \(w=\rho\xi^2\eta^2\), square-heavy
& \(3\M+7\Sqr+\Dpar\)
& \(\Kfour\)
& same rational function after pullback;
  Theorem~\ref{thm:omega-fast-thesis}\\
\(\mathcal C_d\), \(\omega\), six-square
& \(3\M+6\Sqr+3\Dpar\)
& \(\Kfour\)
& \(d\) square, equivalently \(1-\rho\) square\\
\(E_{1,\rho}\), \(w=\rho\xi^2\eta^2\), six-square
& \(3\M+6\Sqr+3\Dpar\)
& \(\Kfour\)
& \(1-\rho\) square; Theorem~\ref{thm:omega-fast-thesis}\\
\(\mathcal C_d\), \(\omega\), complete
& \(5\M+6\Sqr+2\Dpar\)
& \(\Kfour\)
& \(k=\F_q\), \(\rho\) nonsquare;
  Theorem~\ref{thm:omega-complete-thesis}\\
\(E_{1,\rho}\), \(w=\rho\xi^2\eta^2\), complete
& \(5\M+6\Sqr+2\Dpar\)
& \(\Kfour\)
& complete Edwards subfamily, \(\rho\) nonsquare;
  Theorem~\ref{thm:omega-complete-thesis}\\
\bottomrule
\end{tabular}
\end{table}

The two models in each pair of
Table~\ref{tab:odd-four-torsion-comparison} use exactly the same rational
function:
\begin{equation}\label{eq:Cd-Edwards-four-quotient}
 w=\rho\xi^2\eta^2
   =\frac{\rho}{(2u+1)^2(2v+1)^2}
   =\omega.
\end{equation}
The table nevertheless lists both interfaces so that an implementation
claim can state whether its inputs and outputs are Edwards coordinates or
\(\mathcal C_d\)-coordinates.  None of these rows returns the ordinary
Kummer pair \(\{P,-P\}\); it returns the larger orbit
\(\{\pm P\}+\langle T_4\rangle\).  Moreover, the six-square condition for the
degree-two \(\mathcal C_d\) schedule is that \(\rho\) be a square, whereas
the condition for the four-torsion schedule is that
\(1-\rho=16d\) be a square.  They are different subfamilies and cannot be
merged into one completeness claim.

\section{CM endomorphism comparison}

\begin{table}[H]
\centering
\small
\caption{Low-cost CM maps in their declared quotient coordinates}
\label{tab:CM-model-comparison}
\begin{tabular}{L{3.5cm}L{3.8cm}L{2.3cm}L{3.7cm}}
\toprule
model and CM locus & quotient action & cost & interpretation\\
\midrule
\(\mathcal C_{1/8}\), \(j=1728\)
& \((X:Z)\mapsto(-X:Z)\)
& sign only
& native \(\kappa=(u+1:u)\) output\\
Montgomery \(A=0\)
& \(U\mapsto-U\)
& sign only
& same Kummer line in affine basis\\
short Weierstrass \(y^2=x^3+ax\)
& \(x\mapsto-x\)
& sign only
& standard \(j=1728\) CM presentation\\
\(\mathcal C_d\), \(A^2=3\), \(j=0\)
& \((X:Z)\mapsto(\zeta X+cZ:Z)\)
& \(\Dconst{\zeta}+\Dconst{c}\)
& native Kummer basis\\
centered Montgomery abscissa
\(\widetilde U=U+A/3\)
& \(\widetilde U\mapsto\zeta\widetilde U\)
& \(\Dconst{\zeta}\)
& translated basis of the same Kummer line\\
short Weierstrass \(y^2=x^3+b\)
& \(x\mapsto\zeta x\)
& \(\Dconst{\zeta}\)
& standard \(j=0\) CM presentation\\
\bottomrule
\end{tabular}
\end{table}

The \(j=1728\) Kummer action is equally sparse in all three presentations.
At \(j=0\), translation of the Kummer basis removes the constant term and
saves one fixed-constant multiplication, while the native basis preserves
direct access to the declared input and output interface.  In particular,
the same native \(\mathcal C_d\) Kummer coordinate supports the scalar
ladder, full-point recovery, division polynomials, model-preserving
isogenies, and an explicit sparse CM action.  Thus the CM endomorphism is
integrated into the complete native arithmetic framework rather than
requiring a separate quotient coordinate.

\section{Tripling and characteristic three}

\begin{table}[H]
\centering
\caption{Tripling operations relevant to \(\mathcal C_d\)}
\label{tab:tripling-thesis}
\begin{tabular}{L{4.5cm}L{3.4cm}L{2.5cm}L{2.8cm}}
\toprule
method & condition & output & principal cost\\
\midrule
\(\mathcal C_d\) native full tripling
\eqref{eq:odd-full-projective-tripling-thesis}
& \(\charac k\ne2\), native shifted-projective chart
& \(\Ffin\)
& \(9\M+4\Sqr+\Dpar\), or
  \(7\M+7\Sqr+\Dpar\)\\
\(\mathcal C_d\) division-polynomial tripling
\eqref{eq:native-division-tripling}
& reusable \(U\)-powers
& \(\Ffin\)
& symbolic representation; use a dedicated row for a standalone point\\
extended Edwards \(E_{1,\rho}\), DBL then mixed ADD
& \(\charac k\ne2\), affine copy of \(P\) retained
& \(\Full\)
& \(12\M+4\Sqr+\Dpar\); unspecialized chain baseline\\
\(\mathcal C_d\) Frobenius \(x\)TPL
\eqref{eq:char3-tripling}
& \(\charac k=3\)
& \(\Ktwo\)
& \(2\M+2\Sqr+2\Dpar+2\Cube\)\\
\(\mathcal C_d\) binary native \(x\)TPL
& \(\charac k=2\)
& \(\Ktwo\)
& \(6\M+6\Sqr+2\Dpar\)\\
twisted Hessian isogeny core
& rational level-three structure,
  \(\charac k\ne2,3\)
& transported \(\Ktwo\)
& \(4\M+4\Sqr+2\Dpar\), endpoints excluded\\
\bottomrule
\end{tabular}
\end{table}

For ordinary odd characteristic, the dedicated
\(\mathcal C_d\) tripling schedule saves \(3\M\) over the displayed
extended-Edwards double-then-mixed-add baseline at the same
\(4\Sqr+\Dpar\) level.  This is a finite full-coordinate comparison on the
stated native shifted-projective chart; it
must not be combined with the lower-cost Hessian or characteristic-three
rows, which return only quotient data.

In characteristic three the relevant distinction is between
\eqref{eq:char3-full-tripling-native}, which returns \((u_{3P},v_{3P})\),
and \eqref{eq:char3-tripling}, which returns only
\((u_{3P}+1:u_{3P})\).  The diagonal Hessian equation is singular after
naive reduction modulo three, so it is not listed as a characteristic-three
competitor.  For secret scalars, the regular binary ladder remains the
baseline; the Frobenius tripling is most useful for public or fixed
\(\{2,3\}\)-chains.

\begin{table}[H]
\centering
\small
\caption{Full-point realization of \(2P+Q\) in odd characteristic}
\label{tab:double-add-comparison}
\begin{tabular}{L{4.6cm}L{3.2cm}L{2.0cm}L{3.3cm}}
\toprule
model and method & cost & output & interpretation\\
\midrule
\(\mathcal C_d\), closed projective formula
\eqref{eq:native-closed-2PQ}
& \(16\M+3\Sqr+\Dpar\)
& \(\Ffin\)
& one closed dependency graph in \(u,v\)\\
\(\mathcal C_d\), dedicated DBL then mixed ADD
& \(11\M+5\Sqr+2\Dpar\)
& \(\Ffin\)
& native shifted-projective or ordinary projective coordinates\\
extended Edwards \(E_{1,\rho}\), DBL then mixed ADD
& \(12\M+4\Sqr+\Dpar\)
& \(\Full\)
& \(a=1\), affine \(Q\) retained\\
\bottomrule
\end{tabular}
\end{table}

The closed \(\mathcal C_d\) dependency graph realizes the exact tradeoff
\[
        5\M \quad\longleftrightarrow\quad 2\Sqr+\Dpar
\]
relative to the dedicated doubling followed by mixed addition: it uses five
additional general multiplications while saving two squarings and one
parameter multiplication.  Consequently, the closed circuit is advantageous
whenever
\[
                         5\M<2\Sqr+\Dpar.
\]
It is especially well suited to field representations and fixed-point
schedules that reward this exchange, while the dedicated two-step
\(\mathcal C_d\) schedule supplies the complementary arithmetic regime.

\section{Characteristic-two comparison}

The two-isogeny rows below distinguish a Kummer output obtained from a
full native point from a map whose declared input is already Kummer-only.
The \(\mathcal T_{a,d}\) entries preview
Corollary~\ref{cor:T-char2-two-isogeny-native}; their formulas and costs
are proved there before they are used subsequently.

\begin{table}[H]
\centering
\caption{Characteristic-two full and quotient arithmetic}
\label{tab:binary-differential-thesis}
\begin{tabular}{L{4.8cm}L{4.3cm}L{4.0cm}}
\toprule
model & operation and cost & coverage/output\\
\midrule
\(\mathcal C_d\), native affine
& \eqref{eq:binary-native-affine-add} and
  \eqref{eq:binary-native-affine-double}
& \(\Ffin\); basic formulas with inversions\\
\(\mathcal C_d\), native Kummer, circuit I
& \(x\)DBL+\(x\)ADD:
  \(4\M+5\Sqr+\mBase+\mCurve\)
& \(\Ktwo\); rational four-torsion family\\
\(\mathcal C_d\), native Kummer, circuit II
& \(4\M+4\Sqr+\mBase+2\mCurve\)
& \(\Ktwo\); same family\\
\(Z/4\mathbb Z\)-normal form
& full addition \(12\M\), doubling \(7\M+2\Sqr\)
& \(\Full\); linearly equivalent to \(\mathcal C_d\)\\
split \(\mu_4\)-normal form
& full addition \(7\M+2\Sqr+2\mCurve\)
& \(\Full\); marked split four-torsion\\
twisted \(\mu_4\) Kummer
& \(4\M+4\Sqr+\mBase+2\mCurve\)
& \(\Ktwo\); all ordinary twist classes\\
binary Edwards
& complete addition \(16\M+\Sqr+4\mCurve\)
& \(\Full\); broader ordinary-curve coverage\\
\(\mathcal C_d\), native separable two-isogeny, Kummer
& raw affine \(u\) to projective target Kummer: \(\M\)
& \(\Ktwo\); forms \(u(u+1)\)\\
\(\mathcal C_d\), native separable two-isogeny, full
& affine \(u,v\) to two projective target coordinates: \(\M\)
& \(\Full\); \(d\mapsto\sqrt d\)\\
\(\mathcal T_{a,d}\), image \(u\)-Kummer from a full native point
& \(u'=u+v+1\): additions only
& \(\Ktwo\); the declared source contains both \(u\) and \(v\)\\
\(\mathcal T_{a,d}\), native separable two-isogeny
& full projective output: \(\min(\M,\Sqr)\)
& \(\Full\); \(a\) fixed, \(d\mapsto\sqrt d\)\\
binary Edwards through its standard \(W\)-dictionary
& one-way Kummer two-isogeny:
  \(2\M+\Sqr+3\mCurve\)
& projective \(W\)-Kummer output; target binary-Edwards
  re-embedding is separate\\
\bottomrule
\end{tabular}
\end{table}

The last five rows concern isogeny evaluation rather than addition or
doubling.  Their derivation appears first in
Theorem~\ref{thm:Cd-char2-two-isogeny-general} and
Corollary~\ref{cor:T-char2-two-isogeny-native}, and is repeated with the
binary Edwards endpoint dictionary and benchmark in
Section~\ref{sec:Cd-crypto-binary-two-isogeny}.  The comparison is
equal-output: a binary Edwards implementation that persists in the same
\(W\)-Kummer coordinate uses the general
\(\M+\Sqr+\mCurve\) core.  A raw \(C_d\) input exposes its one-product
specialization.  On \(\mathcal T_{a,d}\), the additions-only \(u'\) row
uses a full native source; a Kummer-only source does not contain the
missing \(v\)-coordinate and must not be assigned that cost.  Producing
the second target coordinate costs \(\min(\M,\Sqr)\).

If \(\mBase=\M\), \(\Sqr=\gamma_{2,{\rm S}}\M\), and
\(\mCurve=\gamma_{2,{\rm C}}\M\), the two native ladder steps have
normalized costs
\[
 C_{2,1}=5+5\gamma_{2,{\rm S}}+\gamma_{2,{\rm C}},\qquad
 C_{2,2}=5+4\gamma_{2,{\rm S}}+2\gamma_{2,{\rm C}}.
\]
The second is better exactly when
\(\gamma_{2,{\rm C}}<\gamma_{2,{\rm S}}\).  The
\(\mathcal C_d\) family gives a native characteristic-two model for the
ordinary binary locus carrying marked rational four-torsion.  The explicit
\(Z/4\mathbb Z\)-normal, twisted \(\mu_4\)-normal, and binary Edwards
dictionaries extend the accompanying arithmetic framework to the remaining
ordinary twist classes.  These compatible interfaces therefore organize both
the marked-four-torsion locus and the broader ordinary binary family within
one explicit system of coordinate transformations and arithmetic schedules.

\section[Consolidated cost ledger for Cd]
{Consolidated cost ledger for \texorpdfstring{\(\mathcal C_d\)}{Cd}}
\label{sec:Cd-cost-ledger}

This section collects every point-arithmetic circuit for which the monograph
makes an efficiency claim.  Additions, subtractions, sign changes, and
multiplication by the small integers \(2,3,4\) are omitted.  The symbols
\(\M,\Sqr,\Dpar,\Inv,\mBase,\mCurve,\Cube\) have the meanings fixed in
Chapter~\ref{ch:conventions}.  Symbolic division-polynomial identities,
kernel-size-dependent isogenies, and complete pairing loops are not assigned
a misleading constant cost: their totals depend on the polynomial
evaluation schedule, kernel size, embedding degree, addition chain, twist,
and extension-field basis.  Their local circuits are analyzed in their own
chapters.

\begin{longtable}{L{4.2cm}L{3.2cm}L{1.9cm}L{3.4cm}}
\caption{Native \(\mathcal C_d\) costs in odd characteristic}
\label{tab:Cd-odd-cost-ledger}\\
\toprule
operation and coordinates & exact cost & output & scope\\
\midrule
\endfirsthead
\toprule
operation and coordinates & exact cost & output & scope\\
\midrule
\endhead
affine ADD, two inversions
& \(6\M+2\Inv\) & \(\Ffin\)
& equation~\eqref{eq:native-affine-add-thesis} and its displayed
  expansion~\eqref{eq:native-affine-add-uv}\\
affine ADD, simultaneous inversion
& \(9\M+\Inv\) & \(\Ffin\)
& same rational map\\
affine DBL, two inversions
& \(4\M+2\Sqr+2\Inv\) & \(\Ffin\)
& equation~\eqref{eq:native-affine-dbl-thesis} and its displayed
  expansion~\eqref{eq:native-affine-double-uv}\\
affine DBL, simultaneous inversion
& \(7\M+2\Sqr+\Inv\) & \(\Ffin\)
& same rational map\\
ordinary or native shifted-projective ADD
& \(9\M+\Sqr+\Dpar\) & \(\Ffin\)
& unified, not complete\\
ordinary or native shifted-projective mixed ADD
& \(8\M+\Sqr+\Dpar\) & \(\Ffin\)
& second point affine\\
affine--affine ADD to projective output
& \(7\M\) & \(\Ffin\)
& both input scales equal to one\\
dedicated ordinary or native shifted-projective DBL
& \(3\M+4\Sqr+\Dpar\) & \(\Ffin\)
& Theorem~\ref{thm:inv-dbl-thesis}\\
native-to-extended input lift
& \(\M\) & \(\Full\)
& forms \(Z=(2u+1)(2v+1)\), no inversion\\
complete native Segre ADD
& \(9\M+\Dpar\) & \(\Full\)
& \(k=\F_q\), \(\rho\) nonsquare\\
native Segre mixed ADD, restricted complete tuple
& \(8\M+\Dpar\) & \(\Full\)
& \(k=\F_q\), \(\rho\) nonsquare; second point in the \(T_2=1\) chart\\
dedicated native Segre DBL
& \(4\M+4\Sqr\) & \(\Full\)
& Segre output retained\\
complete double-and-add-always
& \(13\M+4\Sqr+\Dpar\) & \(\Full\)
& \(k=\F_q\), \(\rho\) nonsquare; per scalar bit\\
native \(x\)DBL, general Kummer input
& \(2\M+2\Sqr+\Dpar\) & \(\Ktwo\)
& Theorem~\ref{thm:xdbl-thesis}\\
\textbf{native-input first \(x\)DBL}
& \(\boldsymbol{\M+\Sqr+\Dpar}\) & \(\Ktwo\)
& \(X-Z=1\); one-time initialization\\
native \(x\)ADD, projective difference
& \(4\M+2\Sqr\) & \(\Ktwo\)
& oriented differential input\\
native \(x\)ADD, affine difference
& \(3\M+2\Sqr\) & \(\Ktwo\)
& fixed difference\\
standard \(x\)DBLADD, affine difference
& \(5\M+4\Sqr+\Dpar\) & \(\Ktwo\)
& one ladder bit\\
seven-square \(x\)DBLADD
& \(3\M+7\Sqr+\Dpar\) & \(\Ktwo\)
& all smooth odd parameters\\
six-square \(x\)DBLADD
& \(3\M+6\Sqr+3\Dpar\) & \(\Ktwo\)
& \(\rho\) square\\
fast post-isogeny Kummer step
& \(3\M+7\Sqr+\Dpar\) & \(\Kfour\)
& eight-point orbit, not an ordinary Kummer pair\\
complete post-isogeny differential step
& \(5\M+6\Sqr+2\Dpar\) & \(\Kfour\)
& \(k=\F_q\), \(\rho\) nonsquare\\
native full TPL, multiplication schedule
& \(9\M+4\Sqr+\Dpar\) & \(\Ffin\)
& native shifted-projective output\\
native full TPL, square-heavy schedule
& \(7\M+7\Sqr+\Dpar\) & \(\Ffin\)
& advantageous if \(\Sqr/\M<2/3\)\\
closed \(2P+Q\), projective output
& \(16\M+3\Sqr+\Dpar\) & \(\Ffin\)
& \(13\M+3\Sqr+\Dpar\) blocks plus \(3\M\) output\\
closed \(2P+Q\), affine output
& \(18\M+3\Sqr+\Dpar+\Inv\) & \(\Ffin\)
& simultaneous inversion of \(N_s,N_r\)\\
native halving, one successful branch
& \(4\M+3\Sqr+2\Dconst{A}+\Dconst{\beta}
     +\Dconst{1/2}+\Inv+2\sqrt{\phantom{x}}\) & \(\Full\)
& raw affine input to projective native output\\
\(j=1728\) affine endomorphism \(\phi\)
& \(\Dconst{i}+\Inv\) & \(\Ffin\)
& \(d=1/8\), \(i^2=-1\);
  Theorem~\ref{thm:j1728-native-endomorphism}\\
\(j=1728\) Kummer endomorphism
& sign change only & \(\Ktwo\)
& \((X:Z)\mapsto(-X:Z)\)\\
\(j=0\) affine endomorphism \(\psi\)
& \(3\M+\Dconst{\zeta+c-1}+2\Inv\), or
  \(6\M+\Dconst{\zeta+c-1}+\Inv\) & \(\Ffin\)
& \(A^2=3\), \(\zeta^2+\zeta+1=0\)\\
\(j=0\) Kummer endomorphism
& \(\Dconst{\zeta}+\Dconst{c}\) & \(\Ktwo\)
& \((X:Z)\mapsto(\zeta X+cZ:Z)\)\\
characteristic-three first \(x\)DBL
& \(\M+\Sqr+\Dpar\) & \(\Ktwo\)
& same native initialization\\
characteristic-three Frobenius \(x\)TPL
& \(2\M+2\Sqr+2\Dpar+2\Cube\) & \(\Ktwo\)
& \(d\ne0,1\)\\
\bottomrule
\end{longtable}

For the halving row, the count is for one prescribed successful choice of
the two square roots.  Forming \(q\) costs \(\Inv\), forming \(V_Q\) costs
\(\M\), and the two quadratic stages cost
\(2\Sqr+\Dconst{A}+\Dconst{1/2}
  +2\sqrt{\phantom{x}}\).  The projective full-recovery
blocks then cost \(3\M+\Sqr+\Dconst{A}\), together with the shared
\(\Dconst{\beta}\).  Their sum is the entry in the table.  The identity
\(U^{-1}=T-U\) avoids another inversion.  If affine native coordinates are
required, simultaneous normalization of the two projective pairs adds
\(4\M+\Inv\).  Enumerating both signs of the first root and all retained
second-root branches has a variable cost and is therefore not folded into
this single-branch ledger entry.

\begin{longtable}{L{4.2cm}L{3.3cm}L{2.1cm}L{3.2cm}}
\caption{Native \(\mathcal C_d\) costs in characteristic two}
\label{tab:Cd-binary-cost-ledger}\\
\toprule
operation and coordinates & exact cost & output & scope\\
\midrule
\endfirsthead
\toprule
operation and coordinates & exact cost & output & scope\\
\midrule
\endhead
native affine ADD, \(U_i\) available
& \(10\M+\Sqr+3\Dpar+2\Inv\) & \(\Ffin\)
& simultaneous inversion of \(H,H+d\)\\
formation of both \(U_i\) from raw \(u_i\)
& \(3\M+\Inv\) & prep.
& one batch inversion\\
native affine DBL
& \(\M+3\Sqr+\Inv\) & \(\Ffin\)
& equation~\eqref{eq:binary-native-affine-double}\\
binary native \(x\)DBL
& \(\M+3\Sqr+\mCurve\) & \(\Ktwo\)
& general projective input\\
\textbf{binary native-input first \(x\)DBL}
& \(\boldsymbol{2\Sqr+\mCurve}\) & \(\Ktwo\)
& \(X_0+X_1=1\), no general multiplication\\
binary known-difference normalization
& \(\Inv\) & prep.
& \((u+1:u)\mapsto(1:u/(u+1))\); paid once\\
binary native \(x\)ADD, projective difference
& \(5\M+2\Sqr\) & \(\Ktwo\)
& inversion-free raw native difference\\
binary native \(x\)ADD, circuit I
& \(3\M+2\Sqr+\mBase\) & \(\Ktwo\)
& affine known difference\\
binary native \(x\)ADD, circuit II
& \(3\M+\Sqr+\mBase+\mCurve\) & \(\Ktwo\)
& affine known difference\\
binary native \(x\)DBLADD, circuit I
& \(4\M+5\Sqr+\mBase+\mCurve\) & \(\Ktwo\)
& one ladder bit\\
binary native \(x\)DBLADD, circuit II
& \(4\M+4\Sqr+\mBase+2\mCurve\) & \(\Ktwo\)
& one ladder bit\\
binary native \(x\)DBLADD, projective difference
& \(6\M+5\Sqr+\mCurve\) & \(\Ktwo\)
& inversion-free raw native fixed difference\\
binary native \(x\)TPL
& \(6\M+6\Sqr+2\Dpar\) & \(\Ktwo\)
& Theorem~\ref{thm:binary-tripling-thesis}\\
binary full affine TPL, \(U\) available
& \(11\M+3\Sqr+4\Dpar+3\Inv\) & \(\Ffin\)
& tangent--chord schedule; not the preferred quotient circuit\\
linearly transported \(Z/4\mathbb Z\) ADD
& \(12\M\) & \(\Full\)
& full state kept in the transported coordinates\\
linearly transported \(Z/4\mathbb Z\) DBL
& \(7\M+2\Sqr\) & \(\Full\)
& same qualification\\
\bottomrule
\end{longtable}

\section[Distinct advantages of Cd]
{Distinct advantages of \texorpdfstring{\(\mathcal C_d\)}{Cd}}
\label{sec:Cd-distinct-advantages}

The advantages of \(\mathcal C_d\) are not confined to a sentence in the
conclusion.  They arise from the equation, the marked completion, and the
way these structures interact with arithmetic.  The following table states
each advantage together with its exact mechanism and its legitimate scope.

\begin{longtable}{L{3.1cm}L{4.0cm}L{3.8cm}L{1.9cm}}
\caption{Structural and arithmetic advantages of the model
\(\mathcal C_d:(u^2+u)(v^2+v)=d\)}
\label{tab:Cd-advantages}\\
\toprule
advantage & intrinsic mechanism on \(\mathcal C_d\)
& concrete consequence & boundary\\
\midrule
\endfirsthead
\toprule
advantage & intrinsic mechanism on \(\mathcal C_d\)
& concrete consequence & boundary\\
\midrule
\endhead
one equation in every characteristic
& the two factors \(z^2+z\) remain meaningful without changing the
  defining equation
& odd characteristic, characteristic three, and characteristic two are
  studied in one model and one notation
& optimized circuits still depend on the characteristic\\
new explicit arithmetic model
& its own \((2,2)\)-completion, identity \(O=(0,\infty)\), boundary
  divisor, native coordinates, and addition-law spaces
& isomorphisms with Edwards, Montgomery, and Weierstrass models serve as
  explicit proof dictionaries and optimization bridges while preserving
  the native \(\mathcal C_d\) arithmetic interface
& its arithmetic identity and operation counts are established through
  explicit native formulas, coordinates, and interfaces\\
native \(u\)-Kummer coordinate
& \(k(\mathcal C_d)^{[-1]}=k(u)\) and the boundary normalization gives
  \(\kappa_d=(u+1:u)\)
& input needs one addition and no multiplication or inversion
& quotient output determines \(\{P,-P\}\), not a full point\\
first-double saving
& the original input satisfies \(X-Z=1\)
& the first odd-characteristic \(x\)DBL costs
  \(\M+\Sqr+\Dpar\); the binary first double costs
  \(2\Sqr+\mCurve\)
& the identity is normally lost after the first double\\
native shifted-projective recoding
& \(R=2U+Z,\ S=2V+Z\) is a linear change inside the native chart
& inversion-free ADD, mixed ADD, and dedicated DBL with an additions-only
  return to \((U:V:Z)\)
& its singular plane closure does not distinguish the boundary branches,
  and the formulas are not complete\\
rational four-torsion built into the boundary
& the four boundary points form a marked cyclic subgroup
& Kummer normalization and halving are base-field operations, and the
  marked subgroups of orders two and four give base-field quotient kernels
& the family represents the marked-four-torsion locus\\
two useful quotient layers
& \(\kappa_d\) has degree two, while
  \(\omega=\rho/((2u+1)^2(2v+1)^2)\) has degree eight
& ordinary ladders and complete-parameter quotient laws can coexist
& a \(\Kfour\) output must not be reported as \(\Ktwo\)\\
cheap complete-coordinate entry
& \((u,v)\mapsto(2u+1:2v+1:1:(2u+1)(2v+1))\)
& complete native Segre arithmetic begins with one multiplication and no
  inversion when \(\rho\) is nonsquare
& persistent four-coordinate storage is required for the stated loop costs\\
Artin--Schreier binary branch
& \(v^2+v=d/(u^2+u)\)
& decompression is a linearized equation; the first double has no general
  multiplication
& only the marked-four-torsion ordinary subfamily is covered\\
Frobenius ternary branch
& in characteristic three, cubing is Frobenius and
  \(\psi_3\) collapses
& native \(x\)TPL costs
  \(2\M+2\Sqr+2\Dpar+2\Cube\)
& useful mainly for public or fixed \(2,3\)-chains\\
CM endomorphism subfamilies
& \(d=1/8\) gives \(\phi^2=[-1]\), while
  \(A^2=3\) gives \(\psi^2+\psi+[1]=[0]\)
& the Kummer actions are respectively
  \((-X:Z)\) and \((\zeta X+cZ:Z)\)
& requires the CM constants in the base field and a stable subgroup\\
one reusable native interface
& \(u,v,d\) and \(\kappa_d\) are shared by group laws, division
  polynomials, isogenies, and Miller factors
& conversions are localized inside proofs; outputs return to
  \(\mathcal C_d\)
& task-specific internal representations are selected while every
  declared endpoint remains in the native \(\mathcal C_d\) interface\\
one tunable recurring constant
& \(\alpha_{24}=1/(16d)\) in odd characteristic and
  \(d_{\rm K}=d^{-1}\) in
  characteristic two
& \(d\) may be selected so the recurring constant multiplication is cheap
& security, twist, completeness, and group-order constraints come first\\
\bottomrule
\end{longtable}

These properties place \(\mathcal C_d\) among the explicit arithmetic
models that support full-point and quotient computations within one
low-degree, characteristic-uniform equation.  The same native interface
contains specialized doubling and tripling, halving, isogenies, and pairings
while retaining the original \(u,v\) coordinates.

\section{Operation-specific comparison and model selection}

For pairings, the compared method is the raw Miller evaluation
\(\mathsf M_d^{\rm raw}(n;P,Q)\) of
\eqref{eq:native-Miller-wrapper}.  Its odd and binary factors are the
explicit \(u,v\)-formulas
\eqref{eq:odd-native-Miller-thesis} and
\eqref{eq:binary-native-Miller-thesis}; no auxiliary curve is the declared
Miller-loop state.

\begin{longtable}{L{4.8cm}L{4.7cm}L{3.6cm}}
\caption{Recommended role of \(\mathcal C_d\) by task}
\label{tab:task-recommendations-thesis}\\
\toprule
task & \(\mathcal C_d\) method & final interpretation\\
\midrule
\endfirsthead
\toprule
task & \(\mathcal C_d\) method & final interpretation\\
\midrule
\endhead
secret odd-characteristic scalar
& native Kummer ladder, then recovery if required
& \(\Ktwo\), or a full point of the native completion after recovery\\
complete odd-characteristic full addition
& centered native Segre coordinates
& over \(k=\F_q\), \(\Full\) when \(\rho\) is nonsquare; affine
  \((u,v)\) only off the boundary\\
CM-enabled odd-characteristic scalar
& Theorem~\ref{thm:Cd-GLV-eigenvalues} with joint multiplication through
  a complete native atlas or a separately verified complete law
& \(\Full\); only on the \(j=1728\) or \(j=0\) subfamilies\\
public characteristic-three \(2,3\)-chain
& native \(x\)DBL plus Frobenius \(x\)TPL
& native Kummer quotient\\
full characteristic-three tripling
& \eqref{eq:char3-full-tripling-native}
& \((u_{3P},v_{3P})\)\\
secret characteristic-two scalar
& native binary Kummer ladder plus optional recovery
& \(\Ktwo\), or a full point of the native completion after recovery\\
fixed-point double-add chain
& native closed \(2P+Q\)
& finite full \(\mathcal C_d\)-point; use the native Segre completion
  for a boundary output\\
small isogeny kernels and torsion
& native division polynomials in
  \(\kappa_d=(u+1:u)\)
& parameter and point returned to \(\mathcal C_d\)\\
Tate or Weil pairing
& Miller evaluation \eqref{eq:native-Miller-wrapper} and wrappers
  \eqref{eq:Tate-shifted-Cd}, \eqref{eq:Weil-shifted-Cd}
& pairing computed from original \((u,v)\)-points\\
\bottomrule
\end{longtable}

The comparison identifies the following advantages of \(\mathcal C_d\):
\begin{enumerate}[label=\textup{(\arabic*)}]
 \item one factorized model equation, one identity, one inverse, and one
       native Kummer coordinate in every characteristic;
 \item direct full-point addition and doubling formulas in \(u,v\), followed
       by inversion-free native realizations when needed;
 \item the simultaneous availability of full-point arithmetic, two-to-one
       Kummer ladders, a post-four-isogeny Kummer line, halving, tripling, and closed
       \(2P+Q\) in a single coordinate language;
 \item exceptional simplifications from Artin--Schreier arithmetic in
       characteristic two and Frobenius cubing in characteristic three;
 \item native division polynomials, model-preserving isogenies, and the
       direct Miller evaluation of Chapter~\ref{ch:pairings}, all sharing
       \(\kappa_d=(u+1:u)\);
 \item explicit recovery to \((u,v)\), so an isomorphic Edwards,
       Montgomery, Hessian, Jacobi, or Weierstrass model is a proof or
       optimization device rather than the endpoint of the theory.
\end{enumerate}
Individual competitors may have a lower count for a selected operation on
their preferred subfamily.  The distinctive contribution of
\(\mathcal C_d\) is that these arithmetic layers coexist on one new
biquadratic model across characteristics two, three, and greater than
three.

\chapter{Finite-Field Examples}
\label{ch:finite-field-examples}
\section[A small Montgomery constant over F101]
{A small Montgomery constant over \texorpdfstring{\(\F_{101}\)}{F101}}

\begin{example}\label{ex:F101-small}
Take \(d=60\).  Then
\[
16d=51,\qquad \rho=51,\qquad (16d)^{-1}=2.
\]
Since \(51^{50}=-1\pmod{101}\), \(\rho\) is a nonsquare and the Edwards
law is complete.  The three models are
\[
(u^2+u)(v^2+v)=60,
\]
\[
\xi^2+\eta^2=1+51\xi^2\eta^2,
\]
\[
2V^2=U^3+6U^2+U,\qquad \alpha_{24}=2.
\]
There are \(100\) rational points on the smooth completion, \(96\) affine
and four at the boundary.

For \(P=(4,17)\), one has \((U,V)=(77,37)\) and
\(\ord(P)=20\).  Repeated addition gives
\[
\begin{array}{c|c|c}
n&(U([n]P),V([n]P))&(u([n]P),v([n]P))\\ \hline
1&(77,37)&(4,17)\\
2&(93,72)&(56,73)\\
3&(19,52)&(73,56)\\
4&(7,25)&(17,4)\\
5&(1,2)&(\infty,0)\\
10&(0,0)&(-1,\infty)\\
20&O&O.
\end{array}
\]
This table explicitly displays the rational four- and two-torsion boundary.
\end{example}

\section{Complete arithmetic and halving}

\begin{example}\label{ex:F101-complete-thesis}
Take \(d=1\) over \(\F_{101}\).  Then
\[
\rho=86,\qquad \chi(\rho)=-1,\qquad
\tau=4,\qquad\alpha_{24}=19,
\]
and \(\#\Cd(\F_{101})=96\).  Let \(P=(6,42)\).  Its Edwards image is
\((82,70)\), its order is \(96\), and
\[
\omega(P)=8,\qquad (U_P,V_P)=(18,30).
\]
The native Kummer input is \((X:Z)=(7:6)\), so \(X-Z=1\).
Corollary~\ref{cor:first-dbl-thesis} gives
\[
 AA=13^2=68,\qquad E=67,\qquad
 (X_2:Z_2)=\bigl(68:67(1+19\cdot67)\bigr)=(68:13),
\]
and therefore \(X_2/Z_2=13\).  This numerical calculation uses exactly
one squaring, one multiplication by \(\alpha_{24}=19\), and one general
multiplication, illustrating the \(1\M+1\Sqr+\Dpar\) first-double count.
The ordinary native projective formulas
\eqref{eq:ordinary-native-projective-add} and
\eqref{eq:ordinary-native-projective-double}, applied to
\((U:V:Z)=(6:42:1)\), return \((59:87:1)\) up to common scale.
The native double is
\[
           2P=(59,87)\in\mathcal C_1(\F_{101}),
\]
whose auxiliary Montgomery coordinates are \(Q_M=(13,5)\).  With \(A=74\),
\[
\Delta_Q=13^2+74\cdot13+1=21.
\]
Choose \(R=83\) with \(R^2=21\).  Then
\[
T=2\cdot13+2\cdot83=91,\qquad T^2-4=96.
\]
Taking \(S=46\), one obtains \(U=(T+S)/2=18\), and
\eqref{eq:halving-V-thesis} returns \(V=30\).  Thus the halving formulas
recover the original native point \(P=(6,42)\) exactly.
\end{example}

\section{A binary example}

\begin{example}\label{ex:F256-thesis}
Let
\[
\F_{2^8}=\F_2[\alpha]/(\alpha^8+\alpha^4+\alpha^3+\alpha+1)
\]
and use the same explicit polynomial-basis encoding locally in this
example: if \(\mathtt{0xHH}=\sum_{i=0}^7h_i2^i\), with
\(h_i\in\{0,1\}\), then
\[
 \hexalpha{HH}:=\sum_{i=0}^7h_i\alpha^i.
\]
For example,
\[
 \hexalpha{13}=\alpha^4+\alpha+1,\qquad
 \hexalpha{46}=\alpha^6+\alpha^2+\alpha,\qquad
 \hexalpha{5C}=\alpha^6+\alpha^4+\alpha^3+\alpha^2.
\]
Take
\[
d=\hexalpha{13},\qquad d_{\rm K}=d^{-1}=\hexalpha{4B}.
\]
The curve has \(264\) points.  The native point
\[
P=(u,v)=(\hexalpha{46},\hexalpha{5C})
\]
has order \(264\) and
\[
\kappa(P)=(\hexalpha{47}:\hexalpha{46}).
\]
Since \(u^2+u=\hexalpha{F9}\), the first-double formula gives
\[
\kappa(2P)=(\hexalpha{01}:\hexalpha{EA}).
\]
The relevant auxiliary equation in this example is
\[
       W_d^+:\qquad Y^2+XY=X^3+d^2X.
\]
The Weierstrass image is
\[
(X,Y)=(\hexalpha{DE},\hexalpha{97}),
\]
and doubling on \(W_d^+\) gives
\[
       (X(2P),Y(2P))=(\hexalpha{BD},\hexalpha{30}).
\]
Converting back, the native \(\mathcal C_d\)-point is
\[
                 2P=(\hexalpha{E2},\hexalpha{5F}),
\]
in agreement with
Theorem~\ref{thm:binary-native-affine-double}.  Its Kummer pair
\((\hexalpha{BD}:\hexalpha{13})\) is projectively equal to
\((\hexalpha{01}:\hexalpha{EA})\).  Moreover, in native coordinates,
with the scalars \(132,264,3\) written in decimal,
\[
[132]P=(0,0),\qquad [264]P=O,\qquad
[3]P=(\hexalpha{A9},\hexalpha{9E}).
\]
\end{example}

\section[A complete Segre computation over F1009]
{A complete Segre computation over
\texorpdfstring{\(\F_{1009}\)}{F1009}}

\begin{example}[Full scalar multiplication across the native boundary]
\label{ex:F1009-complete-Segre}
Take \(d=2\).  Then
\[
        \rho=1-16d=978,
        \qquad 978^{504}=-1\pmod {1009},
\]
so the single native Segre tuple of
Theorem~\ref{thm:native-complete-addition-Cd} is
\(\F_{1009}\)-complete.  Direct enumeration gives
\[
        \#\overline{\mathcal C}_2(\F_{1009})=1064.
\]
Choose the two affine native points
\[
             P=(2,108),\qquad Q=(9,178).
\]
Their centered native Segre representatives are
\[
 \widetilde P=(5:217:1:76),\qquad
 \widetilde Q=(19:357:1:729).
\]
For this pair, the blocks
\eqref{eq:native-complete-addition-blocks} are
\[
\begin{array}{c|rrrrrrrr}
 &A&B&C&D&E&F&G&H\\ \hline
 \widetilde P,\widetilde Q
 &95&785&978&918&863&949&887&690.
\end{array}
\]
The native complete law therefore returns
\[
       \widetilde P+\widetilde Q
       =(688:576:160:257).
\]
Applying the native recovery pairs
\((U_0:U_1)=(2Y:Z-Y)\) and
\((V_0:V_1)=(2X:Z-X)\) gives
\[
              P+Q=(153,304),
\]
and direct substitution yields
\[
       (153^2+153)(304^2+304)=2\pmod {1009}.
\]

Now use only the same complete addition tuple and its dedicated double in a
binary chain.  For \(37=(100101)_2\), the result, normalized by its first
Segre coordinate, is
\[
       [37]\widetilde P=(1:16:667:472),
       \qquad [37]P=(771,740).
\]
Again
\((771^2+771)(740^2+740)=2\pmod {1009}\).  Repeated native arithmetic
also gives
\[
\begin{array}{c|c}
n&[n]\widetilde P\text{ (first-coordinate normalization)}\\ \hline
28 &(1:836:110:613)\\
76 &(1:519:735:672)\\
133&(1:0:0:1)\\
266&(0:1:0:-1)\\
532&(0:1:0:1)=O.
\end{array}
\]
Since \(532=2^2\cdot7\cdot19\), the nonidentity values at
\(532/2=266\), \(532/7=76\), and \(532/19=28\) prove that
\(\ord(P)=532\).  The multiples \(133P\) and \(266P\) are boundary
points rather than affine pairs; nevertheless the same polynomial tuple
evaluates without an exceptional branch.  This is the concrete distinction
between completeness on the smooth \((2,2)\)-model and formulas confined
to its affine chart.  A double-and-add-always realization has the bound
\(13\M+4\Sqr+\Dpar\) per processed bit, with one additional multiplication
for the initial native Segre lift.
\end{example}

\section[Frobenius and three-torsion over F27]
{Frobenius, Verschiebung, and three-torsion over
\texorpdfstring{\(\F_{27}\)}{F27}}

\begin{example}[Ordinary and supersingular degree-three behavior]
\label{ex:F27-Frob-isogeny}
Let
\[
       \F_{27}=\F_3[\alpha]/(\alpha^3-\alpha-1),
       \qquad \alpha^3=\alpha+1.
\]
We first take the ordinary parameter \(d=\alpha\).  Enumeration of the
native fibers gives
\[
       \#\overline{\mathcal C}_{\alpha}(\F_{27})=24,
       \qquad a_{\alpha}=27+1-24=4.
\]
Since \(3\nmid a_{\alpha}\), the curve is ordinary.  The point
\[
                       Q=(2\alpha,\alpha)
\]
lies on \(\mathcal C_{\alpha}\), because
\[
 (Q_u^2+Q_u)(Q_v^2+Q_v)
  =(\alpha^2+2\alpha)(\alpha^2+\alpha)=\alpha.
\]
The native addition formula gives
\[
       2Q=-Q=(2\alpha,2\alpha+2),\qquad 3Q=O.
\]
Thus \(K=\{O,Q,-Q\}\) is a reduced rational cyclic kernel and is an
admissible input to the characteristic-three separable V\'elu construction.

The relative Frobenius has a different effect.  It sends
\[
 F_{3,\alpha}(Q)=Q^{(3)}
   =(2\alpha+2,\alpha+1)
   \in\mathcal C_{\alpha+1}(\F_{27}),
\]
because \(\alpha^3=\alpha+1\).  This point is nonidentity and has order
three.  On the other hand \([3]Q=O\); hence the two maps cannot be
identified.  The factorization
\(V_{3,\alpha}\circ F_{3,\alpha}=[3]\) gives
\[
          V_{3,\alpha}(2\alpha+2,\alpha+1)=O.
\]
Since the curve is ordinary, \(V_{3,\alpha}\) is separable of degree
three, and its kernel is exactly
\(\langle(2\alpha+2,\alpha+1)\rangle\).  This calculation exhibits both
a valid separable three-isogeny kernel and the distinct inseparable
Frobenius factor.

For contrast, take \(d=2\in\F_3\).  Here
\[
       \#\overline{\mathcal C}_2(\F_{27})=28,
       \qquad a_2=0,
\]
so the curve is supersingular.  Because \(\rho=1-d=2\) is a nonsquare in
\(\F_{27}\), its native Segre law is a single complete formula.  The point
\[
        P=(\alpha^2+2,\alpha^2+\alpha+1)
\]
has order \(28\); native multiplication gives
\[
\begin{aligned}
 [3]P&=(\alpha^2+\alpha,\,2\alpha^2+\alpha+1),\\
 [7]P&=(1:0:0:1),\\
 [14]P&=(0:1:0:2),\\
 [28]P&=(0:1:0:1)=O,
\end{aligned}
\]
where the last three values are centered native Segre coordinates.  In the
same native affine coordinates,
\[
          F_{3,2}(P)=P^{(3)}
          =(\alpha^2+2\alpha,\alpha^2)
          \ne[3]P.
\]
Although \(d^3=d\), so \(F_{3,2}\) is an endomorphism of the same written
equation, it remains a degree-three purely inseparable map.  Indeed
\(F_{3,2}^3\) fixes every \(\F_{27}\)-rational coordinate, whereas
\([3]^3P=[27]P=-P\ne P\).  In this supersingular case the Verschiebung is
also purely inseparable and there is no reduced geometric subgroup of
order three to which a three-point V\'elu sum could be applied.  The two
halves of the example therefore realize all alternatives in
Proposition~\ref{prop:char3-FV-distinction}.
\end{example}

\section{Class-number and moment checks}

Direct enumeration gives
\[
\begin{array}{c|rrrrrr}
q&3&5&7&11&13&17\\ \hline
J_{\C}(q)&1&2&2&4&6&7\\
I_{\C}(q)&1&3&4&7&9&11\\
\sum_dN_d&4&20&40&108&156&272\\
\sum_da_d^2&0&12&32&96&140&252
\end{array}
\]
in odd characteristic, and
\[
\begin{array}{c|rrrr}
q&2&4&8&16\\ \hline
J_{\C}(q)=I_{\C}(q)&1&3&7&15\\
\sum_{d\ne0}N_d&4&16&64&256\\
\sum_{d\ne0}a_d^2&1&11&55&239
\end{array}
\]
in characteristic two.  These data agree with
Chapters~\ref{ch:moduli} and~\ref{ch:first-moment}.

\part[Twisted, reciprocal, and QRT extensions]
{Three-Parameter, Twisted, Reciprocal, and Symmetric QRT Extensions}
\label{part:twisted-reciprocal-qrt}
\partoverview{This part begins with the full symmetric product family
\(\mathcal C_{a,b,d}:(u^2+u+a)(v^2+v+b)=d\), establishing its exact
smoothness conditions, characteristic-uniform symmetries, Edwards and
Weierstrass dictionaries, rational-point strata, and the boundary of the
results proved for the unrestricted family.  The next two chapters develop
the one-sided twisted family
\(\mathcal T_{a,d}:(u^2+u+a)(v^2+v)=d\) in odd and binary
characteristic.  The fourth chapter studies reciprocal \(C\)-curves,
including their projective completion, four-torsion symmetry, explicit
identification with twisted \(C_d\) models, full-point and Kummer arithmetic,
and characteristic-two collapse to \(\mathcal C_d\).  The fifth chapter
places all of these forms in a symmetric QRT envelope and develops its
Euler--Chasles, Jacobi, Edwards, McMillan, finite-field, isogeny, and
adjacent-Kummer-state structures, including elliptic Lucas recurrences,
state-division polynomials, and the EDS--QRT bridge.  Every completeness and cost statement
retains the coordinate, field, and exceptional-input hypotheses under which
it was proved.  These extension-model results supply the structural
prerequisites for the unified supersingular-isogeny treatment in Part~VII.}

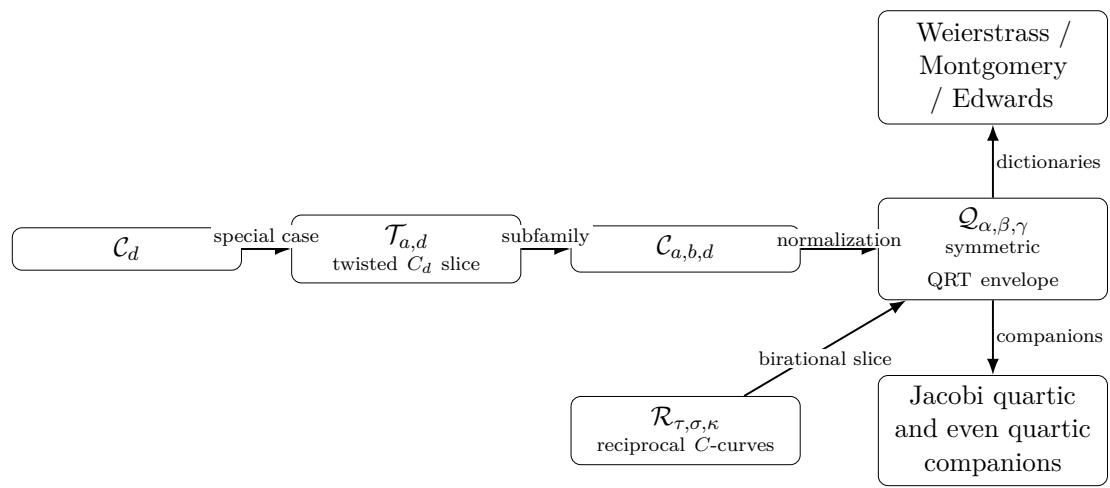
\begin{figure}[H]
\centering
\resizebox{0.96\textwidth}{!}{%
\begin{tikzpicture}[x=1cm,y=1cm,>=Latex, box/.style={draw,rounded corners,align=center,inner sep=4pt,font=\small,text width=3.0cm}]
  \node[box] (Cd) at (0,0) {$\mathcal C_d$};
  \node[box] (Tad) at (4.0,0) {$\mathcal T_{a,d}$\\[-1mm]\scriptsize twisted $C_d$ slice};
  \node[box] (Cabd) at (8.0,0) {$\mathcal C_{a,b,d}$};
  \node[box] (Q) at (12.4,0) {$\mathcal Q_{\alpha,\beta,\gamma}$\\[-1mm]\scriptsize symmetric QRT envelope};
  \node[box] (R) at (8.0,-2.6) {$\mathcal R_{\tau,\sigma,\kappa}$\\[-1mm]\scriptsize reciprocal $C$-curves};
  \node[box] (J) at (12.4,-2.6) {Jacobi quartic\\and even quartic companions};
  \node[box] (W) at (12.4,2.6) {Weierstrass /\\Montgomery / Edwards};
  \draw[->,thick] (Cd) -- node[above,font=\scriptsize,fill=white,inner sep=1pt]{special case} (Tad);
  \draw[->,thick] (Tad) -- node[above,font=\scriptsize,fill=white,inner sep=1pt]{subfamily} (Cabd);
  \draw[->,thick] (Cabd) -- node[above,font=\scriptsize,fill=white,inner sep=1pt]{normalization} (Q);
  \draw[->,thick] (R) -- node[below,font=\scriptsize,fill=white,inner sep=1pt]{birational slice} (Q);
  \draw[->,thick] (Q) -- node[right,font=\scriptsize,fill=white,inner sep=1pt]{companions} (J);
  \draw[->,thick] (Q) -- node[right,font=\scriptsize,fill=white,inner sep=1pt]{dictionaries} (W);
\end{tikzpicture}%
}
\caption{Part~\ref{part:twisted-reciprocal-qrt} organizes the later families around a common biquadratic envelope.  The arrows indicate either inclusion as a distinguished subfamily or an explicit birational/interface dictionary, depending on the label.}
\label{fig:partVI-family-map}
\end{figure}

Figure~\ref{fig:partVI-family-map} gives a roadmap for the families developed in this Part and indicates how the later chapters fit into a single biquadratic framework.

\chapter[A Symmetric Three-Parameter Generalization]
{A Symmetric Three-Parameter Generalization of \texorpdfstring{\(\mathcal C_d\)}{Cd}}
\label{ch:generalized-Cabd}

The one-parameter family studied in the preceding chapters sits naturally
inside the symmetric three-parameter family
\begin{equation}\label{eq:Cabd-definition}
 \mathcal C_{a,b,d}:\qquad
 (u^2+u+a)(v^2+v+b)=d.
\end{equation}
This chapter establishes the geometry of the full product family and its
dictionaries with Edwards and Weierstrass equations.  It does not assign
general operation counts to every \(\mathcal C_{a,b,d}\), because the marked
rational boundary and hence the appropriate arithmetic interface can change
with \(a\) and \(b\).  The following four chapters of this Part develop
complete arithmetic for the one-sided subfamily, reciprocal \(C\)-curves,
and a symmetric QRT envelope; those results do not silently extend to every
torsor in the full three-parameter family.  The first point requiring care is that a smooth
member of \eqref{eq:Cabd-definition} is always a genus-one curve, but it is an
elliptic curve over the ground field only when a ground-field point is chosen
as the identity.

Throughout the chapter put
\[
                 f_a(T)=T^2+T+a.
\]
The original family is \(\mathcal C_d=\mathcal C_{0,0,d}\).  We use
\(\overline{\mathcal C}_{a,b,d}\) for the smooth completion in
\(\PP^1\times\PP^1\).

\section{Completion, smoothness, and the genus-one qualification}

In bihomogeneous coordinates
\((U_0:U_1)\times(V_0:V_1)\), the natural completion is
\begin{equation}\label{eq:Cabd-homogeneous}
\begin{aligned}
 &(U_1^2+U_1U_0+aU_0^2)
  (V_1^2+V_1V_0+bV_0^2)\\
 &\hspace{42mm}=dU_0^2V_0^2.
\end{aligned}
\end{equation}

\begin{proposition}[Exact smoothness conditions]
\label{prop:Cabd-smoothness}
Let \(k\) be a field.
\begin{enumerate}[label=\textup{(\roman*)}]
 \item If \(\charac k\ne2\), put
 \begin{equation}\label{eq:Cabd-four-parameters}
   \lambda=1-4a,\qquad \mu=1-4b,\qquad
   \delta=16d,\qquad \rho=\lambda\mu-\delta.
 \end{equation}
 Then \eqref{eq:Cabd-homogeneous} is smooth if and only if
 \begin{equation}\label{eq:Cabd-smooth-odd}
                  \lambda\mu\delta\rho\ne0.
 \end{equation}
 \item If \(\charac k=2\), it is smooth if and only if \(d\ne0\).
\end{enumerate}
Every smooth member has genus one.  It is an elliptic curve over \(k\) if
and only if it has a \(k\)-rational point, one of which is selected as the
identity.
\end{proposition}

\begin{proof}
Suppose first that \(\charac k\ne2\).  On the affine chart set
\[
             r=2u+1,\qquad s=2v+1.
\]
Since
\[
 f_a(u)=\frac{r^2-\lambda}{4},\qquad
 f_b(v)=\frac{s^2-\mu}{4},
\]
equation~\eqref{eq:Cabd-definition} becomes
\begin{equation}\label{eq:Cabd-centered}
              (r^2-\lambda)(s^2-\mu)=\delta.
\end{equation}
If \(\delta=0\), this equation is reducible.  Assume \(\delta\ne0\).
At an affine singular point, neither factor in
\eqref{eq:Cabd-centered} is zero, while the two partial derivatives are
\[
       2r(s^2-\mu),\qquad 2s(r^2-\lambda).
\]
They vanish simultaneously exactly when \(r=s=0\); the point lies on the
curve exactly when \(\delta=\lambda\mu\), or equivalently \(\rho=0\).

It remains to inspect the boundary.  Near \(r=\infty\), put \(t=r^{-1}\).
Equation~\eqref{eq:Cabd-centered} becomes
\begin{equation}\label{eq:Cabd-boundary-local}
          (1-\lambda t^2)(s^2-\mu)=\delta t^2.
\end{equation}
At \(t=0\) one has \(s^2=\mu\).  If \(\mu\ne0\), then
\(s\ne0\) and the derivative with respect to \(s\) is \(2s\ne0\), so
both boundary points are smooth.  If \(\mu=0\), then \((t,s)=(0,0)\)
is singular.  Interchanging the two factors proves that the boundary over
\(s=\infty\) is smooth exactly when \(\lambda\ne0\).  This proves
\eqref{eq:Cabd-smooth-odd}.

Now let \(\charac k=2\).  Write
\(F=f_a(u)f_b(v)-d\).  Then
\[
                      F_u=f_b(v),\qquad F_v=f_a(u).
\]
When \(d\ne0\), the equation \(F=0\) implies that both factors are
nonzero, so there is no affine singularity.  Near \(u=\infty\), with
\(t=u^{-1}\), the equation is
\[
       (1+t+at^2)f_b(v)=dt^2.
\]
At \(t=0\), the derivative with respect to \(v\) is
\(2v+1=1\), hence the two geometric boundary points are smooth.  The
other boundary is identical.  If \(d=0\), the equation is the reducible
union \(f_a(u)f_b(v)=0\).  This proves the binary assertion.

Finally, a smooth divisor of bidegree \((2,2)\) in
\(\PP^1\times\PP^1\) has genus
\((2-1)(2-1)=1\) by adjunction.  A genus-one curve has a group law over
\(k\) only after a \(k\)-rational base point has been specified; in the
absence of such a point its Jacobian is an elliptic curve but the original
curve can be a nontrivial torsor under that Jacobian.
\end{proof}

Two useful sources of rational points should be recorded immediately.
If \(b=0\), the two points over \(u=\infty\) with
\(v=0,-1\) are rational in every characteristic.  If \(d=ab\ne0\), then
\begin{equation}\label{eq:Cabd-four-affine-points}
       (0,0),\quad(0,-1),\quad(-1,0),\quad(-1,-1)
\end{equation}
are rational points; in characteristic two, \(-1=1\), and these are still
four distinct ordered pairs.  Thus both subfamilies are elliptic rather than
merely genus-one after one of these points is selected.

\section{Characteristic-uniform symmetry}

\begin{proposition}[Factor involutions and parameter exchange]
\label{prop:Cabd-symmetry}
For every characteristic, the maps
\begin{align}
 \iota_u(u,v)&=(-u-1,v),&
 \iota_v(u,v)&=(u,-v-1)\label{eq:Cabd-factor-involutions}
\end{align}
extend to automorphisms of the smooth completion and generate a Klein
four subgroup.  Coordinate exchange gives an isomorphism
\begin{equation}\label{eq:Cabd-swap}
 \tau:\mathcal C_{a,b,d}\longrightarrow\mathcal C_{b,a,d},
             \qquad (u,v)\longmapsto(v,u).
\end{equation}
When \(a=b\), the map \(\tau\) is an automorphism and
\[
       \langle\iota_u,\iota_v,\tau\rangle
       \simeq (C_2\times C_2)\rtimes C_2\simeq D_8.
\]
\end{proposition}

\begin{proof}
The identity
\[
       (-T-1)^2+(-T-1)+c=T^2+T+c
\]
holds in every characteristic.  Hence each map in
\eqref{eq:Cabd-factor-involutions} preserves the equation; its homogeneous
form is a projective linear transformation on the corresponding
\(\PP^1\), so it extends across the boundary.  The two maps are commuting
involutions.  Equation~\eqref{eq:Cabd-swap} follows by interchanging the
two factors.  If \(a=b\), then
\(\tau\iota_u\tau^{-1}=\iota_v\), which gives the stated semidirect
product and the standard presentation of \(D_8\).
\end{proof}

In odd characteristic the involutions become the independent sign changes
\((r,s)\mapsto(-r,s)\) and \((r,s)\mapsto(r,-s)\) in
\eqref{eq:Cabd-centered}.  In characteristic two they become the two
Artin--Schreier translations \(u\mapsto u+1\) and \(v\mapsto v+1\).
Thus the same abstract symmetry is visible as reflection geometry in odd
characteristic and as additive Galois geometry in characteristic two.

\section{Odd characteristic: the diagonal Edwards dictionary}

Assume throughout this section that \(\charac k\ne2\) and that
\eqref{eq:Cabd-smooth-odd} holds.  Retain the four parameters in
\eqref{eq:Cabd-four-parameters}.

\begin{theorem}[Centered reciprocal equation]
\label{thm:Cabd-diagonal-Edwards}
On the dense open set \((2u+1)(2v+1)\ne0\), define
\begin{equation}\label{eq:Cabd-reciprocal-map}
         x=\frac1{2u+1},\qquad y=\frac1{2v+1}.
\end{equation}
Then \(\mathcal C_{a,b,d}\) is birational to the diagonal Edwards
biquadratic
\begin{equation}\label{eq:Cabd-diagonal-Edwards}
       \boxed{\qquad
       \lambda x^2+\mu y^2
       =1+\rho x^2y^2.
       \qquad}
\end{equation}
The inverse map is
\begin{equation}\label{eq:Cabd-reciprocal-inverse}
       u=\frac{x^{-1}-1}{2},\qquad
       v=\frac{y^{-1}-1}{2}.
\end{equation}
The birational maps extend uniquely to the smooth projective completions.

If \(\mu=\epsilon^2\) for some \(\epsilon\in k^\times\), put
\(Y=\epsilon y\).  Then the equation becomes the standard twisted Edwards
curve
\begin{equation}\label{eq:Cabd-twisted-Edwards}
        \lambda x^2+Y^2
        =1+\frac{\rho}{\mu}x^2Y^2.
\end{equation}
If instead \(\lambda\) is a square, the same conclusion follows after
interchanging the two coordinates and the two parameters.
\end{theorem}

\begin{proof}
Starting from \eqref{eq:Cabd-centered}, substitute
\(r=x^{-1}\) and \(s=y^{-1}\), then multiply by \(x^2y^2\):
\[
 (1-\lambda x^2)(1-\mu y^2)=\delta x^2y^2.
\]
Expanding and using \(\rho=\lambda\mu-\delta\) gives
\eqref{eq:Cabd-diagonal-Edwards}.  Solving the two reciprocal definitions
gives \eqref{eq:Cabd-reciprocal-inverse}.  Since both completed curves are
smooth projective genus-one curves, a birational map between them extends
uniquely to an isomorphism.  Finally \(Y^2=\mu y^2\), so substitution in
\eqref{eq:Cabd-diagonal-Edwards} gives
\eqref{eq:Cabd-twisted-Edwards}.  Its two twisted Edwards parameters are
\[
                  a_E=\lambda,\qquad d_E=\rho/\mu.
\]
Their nonvanishing and inequality follow from
\(\lambda\mu\rho\delta\ne0\), because
\(a_E-d_E=\delta/\mu\ne0\).
\end{proof}

The square-root condition in Theorem~\ref{thm:Cabd-diagonal-Edwards} is a
field-of-definition statement, not a geometric obstruction.  Over
\(k(\sqrt\mu)\) every smooth member has the displayed twisted Edwards
equation.  Over \(k\), the original genus-one curve may have no rational
point at all, in which case it cannot be \(k\)-isomorphic to any elliptic
curve equation with a declared rational identity.  Even when it has a
rational point, the absence of a square root of \(\mu\) means only that this
particular boundary-normalized diagonal scaling is not defined over \(k\).

\begin{corollary}[The subfamily \(b=0\)]
\label{cor:Cabd-b-zero-TE}
If \(b=0\), then \(\mu=1\), and the map
\eqref{eq:Cabd-reciprocal-map} is defined over \(k\) and gives exactly
\begin{equation}\label{eq:Cabd-b-zero-TE}
 \boxed{\quad
 (1-4a)x^2+y^2
 =1+(1-4a-16d)x^2y^2.
 \quad}
\end{equation}
For \(a=b=0\), this reduces to
\begin{equation}\label{eq:Cabd-Cd-Edwards}
             x^2+y^2=1+(1-16d)x^2y^2,
\end{equation}
the centered reciprocal Edwards equation of \(\mathcal C_d\).
\end{corollary}

\begin{proof}
Put \(\mu=1\) in
\eqref{eq:Cabd-diagonal-Edwards}; then
\(\rho=\lambda-\delta=1-4a-16d\).  Putting \(a=0\) gives the second
formula.  No further rescaling or extension of the ground field is used.
\end{proof}

\begin{corollary}[The symmetric subfamily \(a=b\)]
\label{cor:Cabd-a-equals-b}
Let \(a=b\), put \(\lambda=1-4a\), and assume
\(\lambda d(\lambda^2-16d)\ne0\).  Then
\begin{equation}\label{eq:Cabd-symmetric-diagonal}
        \lambda(x^2+y^2)
        =1+(\lambda^2-16d)x^2y^2.
\end{equation}
If \(\lambda=\theta^2\in k^{\times2}\), the scaling
\(X=\theta x\), \(Y=\theta y\) gives the ordinary Edwards equation
\begin{equation}\label{eq:Cabd-symmetric-Edwards}
       X^2+Y^2
       =1+\frac{\lambda^2-16d}{\lambda^2}X^2Y^2.
\end{equation}
If \(\lambda\) is a nonsquare, equation
\eqref{eq:Cabd-symmetric-diagonal} is the corresponding symmetric quadratic
twist; the coordinate-exchange automorphism remains defined over \(k\).
\end{corollary}

\begin{proof}
When \(a=b\), one has \(\mu=\lambda=1-4a\) and
\(\rho=\lambda^2-16d\).  Substitution in
\eqref{eq:Cabd-diagonal-Edwards} gives
\eqref{eq:Cabd-symmetric-diagonal}.  If
\(\lambda=\theta^2\), set \(X=\theta x\) and \(Y=\theta y\).  Then
\(x^2+y^2=(X^2+Y^2)/\lambda\) and
\(x^2y^2=X^2Y^2/\lambda^2\); substitution gives
\eqref{eq:Cabd-symmetric-Edwards}.  If \(\lambda\) is a nonsquare, this
simultaneous scaling is defined only after adjoining \(\sqrt\lambda\), so
the base-field equation is its symmetric quadratic twist.
Proposition~\ref{prop:Cabd-symmetry} independently shows that coordinate exchange is
defined over \(k\), regardless of this square class.
\end{proof}

\section{Odd-characteristic Weierstrass and Jacobian models}

The distinction between the curve and its Jacobian is essential for the
general parameters.  The next theorem gives a Weierstrass equation over
\(k\) even when \(\mathcal C_{a,b,d}(k)\) is empty.

\begin{theorem}[The Jacobian in odd characteristic]
\label{thm:Cabd-odd-Jacobian}
Put
\begin{equation}\label{eq:Cabd-c-parameter}
                  c=\lambda\mu+\rho=2\lambda\mu-\delta.
\end{equation}
Then the Jacobian of \(\mathcal C_{a,b,d}\) is
\begin{equation}\label{eq:Cabd-odd-Jacobian}
 \boxed{\qquad
 E_{a,b,d}^{J}:\quad
 Y^2=X\bigl(X^2+2cX+\delta^2\bigr).
 \qquad}
\end{equation}
Its discriminant is
\begin{equation}\label{eq:Cabd-Jacobian-discriminant}
              256\delta^4\lambda\mu\rho,
\end{equation}
and
\begin{equation}\label{eq:Cabd-Jacobian-j}
 c_4=16(4c^2-3\delta^2),\qquad
 j=16\frac{(4c^2-3\delta^2)^3}
                  {\delta^4\lambda\mu\rho}.
\end{equation}
Thus its smoothness is equivalent to
\eqref{eq:Cabd-smooth-odd}.

Let \(K=k(\epsilon)\), where \(\epsilon^2=\mu\).  Over \(K\), select as
origin the boundary point characterized in centered coordinates by
\(r=\infty\) and \(s=\epsilon\), equivalently by \(x=0\) and
\(y_E=1\).  The isomorphism of this pointed curve with its Jacobian can
then be written as
\begin{equation}\label{eq:Cabd-to-Jacobian-map}
\begin{aligned}
 x&=(2u+1)^{-1},& y&=(2v+1)^{-1},\\
 y_E&=\epsilon y,&
 U&=\frac{1+y_E}{1-y_E},\\
 V&=\frac{U}{x},&
 X&=\delta U,\qquad Y=2\epsilon\delta V.
\end{aligned}
\end{equation}
On the open set where the displayed denominators are nonzero, the inverse is
\begin{equation}\label{eq:Cabd-Jacobian-inverse}
\begin{aligned}
 U&=X/\delta,&V&=Y/(2\epsilon\delta),\\
 x&=U/V=2\epsilon X/Y,&
 y_E&=\frac{U-1}{U+1}=\frac{X-\delta}{X+\delta},\\
 u&=\frac{x^{-1}-1}{2},&
 v&=\frac{\epsilon/y_E-1}{2}.
\end{aligned}
\end{equation}
If \(\mu\) is a square in \(k\), these are \(k\)-rational maps.  If not,
\eqref{eq:Cabd-odd-Jacobian} is still the \(k\)-rational Jacobian, but the
original curve need not be a trivial torsor under it.
\end{theorem}

\begin{proof}
The proof proceeds by eliminating one centered coordinate, obtaining a
quartic genus-one equation, passing to a cubic, and then normalizing the
coefficients to the displayed Jacobian.  The inverse substitutions are
recorded at the end so that the construction is a dictionary rather than
only an invariant calculation.
First eliminate the second centered coordinate from
\eqref{eq:Cabd-centered}.  Put
\begin{equation}\label{eq:Cabd-even-quartic-z}
                       z=(r^2-\lambda)s.
\end{equation}
Since
\[
 (r^2-\lambda)s^2
 =\mu(r^2-\lambda)+\delta
 =\mu r^2-\rho,
\]
we obtain the even quartic
\begin{equation}\label{eq:Cabd-even-quartic}
 z^2=(r^2-\lambda)(\mu r^2-\rho)
    =\mu r^4-cr^2+\lambda\rho.
\end{equation}
For an even quartic
\[
             z^2=Ar^4+Br^2+C,
\]
the binary-quartic Jacobian construction gives
\begin{equation}\label{eq:even-quartic-Jacobian-template}
       Y^2=X\{X^2-2BX+(B^2-4AC)\}.
\end{equation}
For completeness, the coefficients in this construction can be checked by
passing to a field where the leading and constant coefficients are squares,
scaling to a Jacobi quartic \(w^2=t^4+2qt^2+1\), and using
\[
  2V^2=U\bigl(U^2-2qU+q^2-1\bigr).
\]
Undoing the two scalings gives
\eqref{eq:even-quartic-Jacobian-template}; because all its coefficients are
polynomials in \(A,B,C\), the equation descends to the original field.  This
is the canonical binary-quartic Jacobian, rather than an arbitrarily chosen
quadratic twist.

In \eqref{eq:Cabd-even-quartic}, one has
\[
 A=\mu,\qquad B=-c,\qquad C=\lambda\rho,\qquad
 B^2-4AC=(\lambda\mu-\rho)^2=\delta^2.
\]
Substitution in \eqref{eq:even-quartic-Jacobian-template} proves
\eqref{eq:Cabd-odd-Jacobian}.

We now verify the explicit map.
Theorem~\ref{thm:Cabd-diagonal-Edwards} gives over \(K\)
\[
       \lambda x^2+y_E^2
       =1+\frac{\rho}{\mu}x^2y_E^2.
\]
The standard Edwards-to-Montgomery substitution
\[
       U=\frac{1+y_E}{1-y_E},\qquad V=\frac Ux
\]
gives
\begin{equation}\label{eq:Cabd-Montgomery}
 \frac{4\mu}{\delta}V^2
 =U^3+\frac{2c}{\delta}U^2+U.
\end{equation}
Indeed, substituting
\(y_E=(U-1)/(U+1)\) and \(x=U/V\), clearing
\((U+1)^2V^2\), and using
\(\lambda-\rho/\mu=\delta/\mu\) gives exactly
\eqref{eq:Cabd-Montgomery}.  Finally set
\(X=\delta U\) and \(Y=2\epsilon\delta V\).  Multiplying
\eqref{eq:Cabd-Montgomery} by \(\delta^3\) yields
\eqref{eq:Cabd-odd-Jacobian}.  On the dense open set where the denominators
are nonzero, the inverse formulas follow successively from
\[
 U=X/\delta,
 \qquad
 V=Y/(2\epsilon\delta),
 \qquad
 x=U/V,
 \qquad
 y_E=(U-1)/(U+1).
\]
Since \(x=(2u+1)^{-1}\) and
\(y_E=\epsilon(2v+1)^{-1}\), solving these two linear-fractional
relations gives
\[
 u=\frac{x^{-1}-1}{2},
 \qquad
 v=\frac{\epsilon/y_E-1}{2}.
\]
These are exactly the formulas in \eqref{eq:Cabd-Jacobian-inverse}.

For the discriminant, equation~\eqref{eq:Cabd-odd-Jacobian} has
\(a_2=2c\), \(a_4=\delta^2\), and \(a_1=a_3=a_6=0\).  Hence
\[
 \Delta_E=16a_4^2(a_2^2-4a_4)
 =64\delta^4(c^2-\delta^2).
\]
Since
\[
 c^2-\delta^2
 = (2\lambda\mu-\delta)^2-\delta^2
 =4\lambda\mu(\lambda\mu-\delta)
 =4\lambda\mu\rho,
\]
formula~\eqref{eq:Cabd-Jacobian-discriminant} follows.
The same generalized Weierstrass formulas give
\(c_4=16(4c^2-3\delta^2)\); division of \(c_4^3\) by the discriminant
gives \eqref{eq:Cabd-Jacobian-j}.
\end{proof}

\begin{corollary}[Direct \(b=0\) Weierstrass dictionary]
\label{cor:Cabd-b-zero-Weierstrass}
Assume \(b=0\), put
\[
 \lambda=1-4a,\qquad \delta=16d,\qquad
 c=2\lambda-\delta,
\]
and let \((x,y)\) be the twisted Edwards coordinates in
\eqref{eq:Cabd-reciprocal-map}.  Then
\begin{equation}\label{eq:Cabd-b-zero-to-W}
 U=\frac{1+y}{1-y},\qquad V=\frac Ux,\qquad
 X=\delta U,\qquad Y=2\delta V
\end{equation}
gives
\begin{equation}\label{eq:Cabd-b-zero-W}
       Y^2=X\bigl(X^2+2(2\lambda-\delta)X+\delta^2\bigr).
\end{equation}
The inverse is obtained by
\[
 U=X/\delta,\qquad V=Y/(2\delta),\qquad
 x=U/V,\qquad y=(U-1)/(U+1),
\]
followed by \eqref{eq:Cabd-reciprocal-inverse}.  For \(a=0\), this is the
direct Weierstrass dictionary of \(\mathcal C_d\), expressed in the scaling
used in the present chapter.
\end{corollary}

\begin{proof}
This is Theorem~\ref{thm:Cabd-odd-Jacobian} with
\(\mu=\epsilon=1\).  All substitutions and their inverses are therefore
defined over \(k\), and no square-root extension is present.
\end{proof}

\begin{remark}[Curve, Jacobian, and twisted Edwards form]
The three equations in this section have different logical roles.
Equation~\eqref{eq:Cabd-diagonal-Edwards} is birational to the original
genus-one curve over \(k\).  Equation~\eqref{eq:Cabd-twisted-Edwards} is a
standard twisted Edwards equation over a field containing \(\sqrt\mu\).
Equation~\eqref{eq:Cabd-odd-Jacobian} is the Jacobian over \(k\).  If the
original curve has a \(k\)-point and that point is chosen as identity, it is
\(k\)-isomorphic to the Jacobian, but the isomorphism need not be the
boundary-normalized map \eqref{eq:Cabd-to-Jacobian-map} when
\(\mu\notin k^{\times2}\).
\end{remark}

\section[Short forms outside characteristic two]
{Short forms in characteristic greater than three and in characteristic three}

\begin{proposition}[Short Weierstrass form for \(\charac k>3\)]
\label{prop:Cabd-short-greater-three}
Assume \(\charac k>3\).  In
\eqref{eq:Cabd-odd-Jacobian} set
\begin{equation}\label{eq:Cabd-short-change}
          \mathsf X=X+\frac{2c}{3},\qquad \mathsf Y=Y.
\end{equation}
Then the Jacobian has the short form
\begin{equation}\label{eq:Cabd-short-greater-three}
          \mathsf Y^2=\mathsf X^3+A_4\mathsf X+A_6,
\end{equation}
where
\begin{equation}\label{eq:Cabd-short-coefficients}
 A_4=\delta^2-\frac{4c^2}{3},\qquad
 A_6=\frac{16c^3}{27}-\frac{2c\delta^2}{3}.
\end{equation}
The inverse change is \(X=\mathsf X-2c/3\), \(Y=\mathsf Y\).
\end{proposition}

\begin{proof}
Write \eqref{eq:Cabd-odd-Jacobian} as
\[
                 Y^2=X^3+2cX^2+\delta^2X.
\]
Substitute \(X=\mathsf X-2c/3\).  The coefficient of
\(\mathsf X^2\) is
\(-3(2c/3)+2c=0\).  Expanding the remaining two coefficients gives
\[
 \delta^2-\frac{(2c)^2}{3},
 \qquad
 \frac{2(2c)^3}{27}-\frac{(2c)\delta^2}{3},
\]
which are exactly \eqref{eq:Cabd-short-coefficients}.
\end{proof}

\begin{proposition}[The characteristic-three reduced form]
\label{prop:Cabd-char-three-reduced}
Suppose \(\charac k=3\).  Then
\begin{equation}\label{eq:Cabd-char-three-parameters}
 \lambda=1-a,\qquad \mu=1-b,\qquad
 \delta=d,\qquad \rho=(1-a)(1-b)-d,
\end{equation}
and smoothness is
\[
             d(1-a)(1-b)\bigl((1-a)(1-b)-d\bigr)\ne0.
\]
Put \(a_2=2c\) and \(a_4=\delta^2\).  If \(a_2=0\), the Jacobian is
already the short equation
\begin{equation}\label{eq:Cabd-char-three-short-special}
                    Y^2=X^3+a_4X.
\end{equation}
If \(a_2\ne0\), put
\begin{equation}\label{eq:Cabd-char-three-shift}
       r_0=\frac{a_4}{a_2},\qquad X=Z+r_0.
\end{equation}
Then
\begin{equation}\label{eq:Cabd-char-three-reduced}
 Y^2=Z^3+a_2Z^2+B_6,\qquad
 B_6=r_0^3-\frac{a_4^2}{a_2}.
\end{equation}
In general the nonzero \(Z^2\)-coefficient cannot be removed by a
translation in characteristic three.
\end{proposition}

\begin{proof}
The parameter reductions follow from \(4=16=1\) in characteristic three.
The Jacobian remains
\(Y^2=X^3+a_2X^2+a_4X\).  If \(a_2=0\), this is
\eqref{eq:Cabd-char-three-short-special}.  Otherwise substitute
\(X=Z+r_0\).  Since \((Z+r_0)^3=Z^3+r_0^3\) and
\((Z+r_0)^2=Z^2-r_0Z+r_0^2\), the coefficient of \(Z\) is
\(-a_2r_0+a_4=0\).  The constant term is
\[
 r_0^3+a_2r_0^2+a_4r_0
 =r_0^3+2\frac{a_4^2}{a_2}
 =r_0^3-\frac{a_4^2}{a_2},
\]
because \(2=-1\).  A translation cannot change the coefficient of
\(Z^2\), since the missing binomial coefficient is \(3=0\).
\end{proof}

For example, the \(b=0\) equation in characteristic three is
\begin{equation}\label{eq:Cabd-b-zero-char-three}
       (1-a)x^2+y^2
       =1+(1-a-d)x^2y^2,
\end{equation}
with no formal division by a vanished integer.  This is the direct
characteristic-three specialization of
\eqref{eq:Cabd-b-zero-TE}, whereas the reduced Weierstrass form must be
handled by Proposition~\ref{prop:Cabd-char-three-reduced}, not by the
characteristic-greater-than-three shift.

\section[Characteristic two: binary Weierstrass form]{Characteristic two: Artin--Schreier descent and binary Weierstrass form}

Assume throughout this section that \(\charac k=2\) and \(d\ne0\).

\begin{theorem}[Binary Jacobian and its twist parameter]
\label{thm:Cabd-binary-Jacobian}
The Jacobian of \(\mathcal C_{a,b,d}\) has the sparse ordinary binary
Weierstrass equation
\begin{equation}\label{eq:Cabd-binary-Jacobian}
 \boxed{\qquad
 E_{a,b,d}^{J}:\quad
 Y^2+XY=X^3+(a+b)X^2+d^2X.
 \qquad}
\end{equation}
It has
\begin{equation}\label{eq:Cabd-binary-invariants}
             \Delta_E=d^4,\qquad c_4=1,\qquad j=d^{-4}.
\end{equation}
Over a field containing elements \(r,s\) with
\begin{equation}\label{eq:Cabd-AS-roots}
                   r^2+r=a,\qquad s^2+s=b,
\end{equation}
the translation
\begin{equation}\label{eq:Cabd-AS-translation}
                   U=u+r,\qquad V=v+s
\end{equation}
identifies \(\mathcal C_{a,b,d}\) with \(\mathcal C_d\).
Thus \(a+b\) records the Artin--Schreier twist of the Jacobian, while
\(d\) determines its geometric \(j\)-invariant.
\end{theorem}

\begin{proof}
In characteristic two,
\[
 U^2+U=u^2+u+r^2+r=u^2+u+a,
\]
and
\[
 V^2+V=v^2+v+s^2+s=v^2+v+b.
\]
This proves
\eqref{eq:Cabd-AS-translation}.  Over the same splitting field the
Jacobian is therefore the binary curve
\[
               W_d^+:\quad Y_0^2+X Y_0=X^3+d^2X
\]
from Theorem~\ref{thm:binary-W-thesis}.

Put \(t=r+s\).  Then
\[
                       t^2+t=a+b.
\]
The change \(Y_0=Y+tX\) transforms
\eqref{eq:Cabd-binary-Jacobian} into \(W_d^+\), because
\[
 (Y+tX)^2+X(Y+tX)
 =Y^2+XY+(t^2+t)X^2.
\]
Consequently \eqref{eq:Cabd-binary-Jacobian} is precisely the
Artin--Schreier twist that descends the split Jacobian.  The separate changes
\(r\mapsto r+1\) and \(s\mapsto s+1\) alter the genus-one trivialization by
the two factor involutions; on \(\operatorname{Pic}^0\) their common linear
part is the binary inverse \((X,Y)\mapsto(X,Y+X)\).  Hence the sum
\(r+s\), and therefore the class of \(a+b\), is the descent parameter.

For \(a_1=1\), \(a_2=a+b\), \(a_4=d^2\), and
\(a_3=a_6=0\), the generalized Weierstrass invariant formulas give
\(c_4=1\) and \(\Delta_E=d^4\), proving
\eqref{eq:Cabd-binary-invariants}.  Notice that replacing \(a+b\) by
\(a+b+q^2+q\) gives an isomorphic curve through \(Y\mapsto Y+qX\), as
required for an Artin--Schreier twist parameter.
\end{proof}

Over a finite field \(k=\F_{2^m}\), equations
\eqref{eq:Cabd-AS-roots} can be solved in \(k\) exactly when
\(\Tr(a)=\Tr(b)=0\).  If this happens, the original curve itself is
\(k\)-isomorphic to \(\mathcal C_d\) through
\eqref{eq:Cabd-AS-translation}.  If not, it can still have a rational
point, but the isomorphism must be constructed from that point rather than
from two separate Artin--Schreier roots.

\begin{theorem}[A direct binary map when \(b=0\)]
\label{thm:Cabd-binary-b-zero}
Let \(b=0\), put \(F=u^2+u+a\), and define
\begin{equation}\label{eq:Cabd-binary-b-zero-map}
                X=Fv+d,\qquad Y=uX.
\end{equation}
Then
\begin{equation}\label{eq:Cabd-binary-b-zero-W}
       Y^2+XY=X^3+aX^2+d^2X.
\end{equation}
On the dense open set \(X F\ne0\), the inverse is
\begin{equation}\label{eq:Cabd-binary-b-zero-inverse}
             u=Y/X,\qquad v=(X+d)/(u^2+u+a).
\end{equation}
Thus the binary \(b=0\) curve is already an elliptic curve over the base
field and is isomorphic there to the Jacobian in
\eqref{eq:Cabd-binary-Jacobian}.
\end{theorem}

\begin{proof}
Put \(T=Fv\), so \(X=T+d\).  The curve equation
\(F(v^2+v)=d\) gives
\[
       T^2+FT=F^2(v^2+v)=dF,
\]
hence \(T^2=F(T+d)=FX\).  Since
\(T^2=(X+d)^2=X^2+d^2\), one obtains
\[
                         FX=X^2+d^2.
\]
Multiplying \(F=u^2+u+a\) by \(X^2\), and using \(Y=uX\), gives
\[
 Y^2+XY+aX^2=X^2F=X^3+d^2X,
\]
which is \eqref{eq:Cabd-binary-b-zero-W}.  The inverse formulas follow
directly from \(Y=uX\) and \(Fv=X+d\).  At the boundary, the point
\((u,v)=(\infty,1)\) maps to the point at infinity of the Weierstrass
model and therefore supplies the identity under which the birational map
extends to an elliptic-curve isomorphism.  The other boundary point
\((\infty,0)\) maps to \((X,Y)=(0,0)\).
\end{proof}

Equation~\eqref{eq:Cabd-binary-Jacobian} is the appropriate binary short
form for this family.  The term \(XY\) must be retained: deleting it would
destroy the separable degree-two map to the \(X\)-line.  Thus the phrase
``short Weierstrass form'' has a characteristic-dependent meaning; the
equation \(Y^2=X^3+AX+B\) used in characteristic greater than three is not
a valid universal target in characteristic two.

\section[The rational-point subfamily d=ab]
{The rational-point subfamily \texorpdfstring{\(d=ab\)}{d=ab}}

Assume \(d=ab\ne0\).  The point \((0,0)\) is rational, so the genus-one
curve is an elliptic curve after choosing this point as identity.  This
subfamily admits especially explicit base-field maps which do not require
square roots of \(1-4a\) or \(1-4b\).

\subsection{Odd characteristic}

Assume \(\charac k\ne2\).  Smoothness becomes
\begin{equation}\label{eq:Cabd-dab-smooth-odd}
 ab(1-4a)(1-4b)(1-4a-4b)\ne0.
\end{equation}
Indeed,
\(\rho=(1-4a)(1-4b)-16ab=1-4a-4b\).

\begin{theorem}[Huff-type cubic and explicit Weierstrass map]
\label{thm:Cabd-dab-Huff}
Put
\begin{equation}\label{eq:Cabd-dab-rs}
       r=\frac{u}{u+1},\qquad s=\frac{v}{v+1},
       \qquad c_0=1-2a-2b.
\end{equation}
Then \(\mathcal C_{a,b,ab}\) is birational to the cubic
\begin{equation}\label{eq:Cabd-dab-cubic}
 b r(1+s^2)+a s(1+r^2)+c_0rs=0.
\end{equation}
Set
\begin{equation}\label{eq:Cabd-dab-Huff-parameters}
 x=-br,\qquad y=as,\qquad
 A=-a^{-2},\quad B=-b^{-2},\quad C=\frac{c_0}{ab}.
\end{equation}
Then the cubic is the generalized Huff equation
\begin{equation}\label{eq:Cabd-generalized-Huff}
            y-x=xy(Bx-Ay+C).
\end{equation}
Define
\begin{equation}\label{eq:Cabd-Huff-to-W}
 U=\frac1{xy},\qquad
 V=2U(U+B)x+CU.
\end{equation}
The resulting Weierstrass equation is
\begin{equation}\label{eq:Cabd-dab-Weierstrass}
 \boxed{\quad
 V^2=C^2U^2+4U(U+A)(U+B).
 \quad}
\end{equation}
Equivalently, with \(W=V/2\),
\begin{equation}\label{eq:Cabd-dab-Weierstrass-monic}
 W^2=U^3+\left(A+B+\frac{C^2}{4}\right)U^2+ABU.
\end{equation}
The inverse rational map is
\begin{equation}\label{eq:Cabd-Huff-inverse}
\begin{aligned}
 x&=\frac{V-CU}{2U(U+B)},&
 y&=\frac{V+CU}{2U(U+A)},\\
 u&=-\frac{x}{b+x},&
 v&=\frac{y}{a-y}.
\end{aligned}
\end{equation}
\end{theorem}

\begin{proof}
Since \(u=r/(1-r)\), one has
\[
 u^2+u+a
 =\frac{ar^2+(1-2a)r+a}{(1-r)^2},
\]
while direct substitution of \(v=s/(1-s)\) gives
\[
 v^2+v+b
 =\frac{bs^2+(1-2b)s+b}{(1-s)^2}.
\]
Substitute these two
expressions in
\((u^2+u+a)(v^2+v+b)=ab\), multiply by
\((1-r)^2(1-s)^2\), and expand.  The quartic, pure quadratic, and constant
terms on the two sides cancel; the remaining terms are
\[
       br+as+c_0rs+ar^2s+brs^2=0,
\]
which is \eqref{eq:Cabd-dab-cubic}.

Now substitute \(x=-br\), \(y=as\).  The five terms become
\[
 -x+y-\frac{xy^2}{a^2}+\frac{x^2y}{b^2}
       -\frac{c_0xy}{ab}=0.
\]
Moving the nonlinear terms to the right and using the definitions of
\(A,B,C\) gives \eqref{eq:Cabd-generalized-Huff}.

From \eqref{eq:Cabd-generalized-Huff},
\[
 U=\frac{Bx-Ay+C}{y-x}=\frac1{xy}.
\]
Rearranging the numerator identity gives
\[
              (A+U)y=(B+U)x+C.
\]
Using \(xy=U^{-1}\), eliminate \(y\) to obtain
\begin{equation}\label{eq:Cabd-Huff-quadratic}
       U(U+B)x^2+CUx-(U+A)=0.
\end{equation}
The discriminant of this quadratic is
\[
             C^2U^2+4U(U+A)(U+B).
\]
The expression for \(V\) in \eqref{eq:Cabd-Huff-to-W} is exactly
\(2\) times the leading coefficient times \(x\), plus the linear
coefficient; squaring it and using
\eqref{eq:Cabd-Huff-quadratic} proves
\eqref{eq:Cabd-dab-Weierstrass}.  Solving the same quadratic gives the first
formula in \eqref{eq:Cabd-Huff-inverse}.  Moreover,
\[
 (V-CU)(V+CU)=4U(U+A)(U+B),
\]
so the displayed expression for \(y\) satisfies \(xy=U^{-1}\).
Finally, solving \(r=u/(u+1)=-x/b\) and
\(s=v/(v+1)=y/a\) gives the last two inverse formulas.  Every step is
reversible on a dense open set, and smoothness extends the birational map to
an isomorphism of the completed elliptic curves.
\end{proof}

\subsection{Characteristic two}

Now let \(\charac k=2\).  The same substitutions
\(r=u/(u+1)\), \(s=v/(v+1)\) give a different, genuinely binary
cubic.

\begin{theorem}[Explicit binary Weierstrass map for \(d=ab\)]
\label{thm:Cabd-dab-binary}
Put
\begin{equation}\label{eq:Cabd-dab-binary-scaled}
  x=br,\qquad y=as,\qquad
  A=a^{-2},\quad B=b^{-2},\quad C=(ab)^{-1}.
\end{equation}
Then
\begin{equation}\label{eq:Cabd-dab-binary-Huff}
          x+y=xy(Bx+Ay+C).
\end{equation}
Define
\begin{equation}\label{eq:Cabd-dab-binary-UV}
       U=\frac1{xy},\qquad V=U(U+B)x.
\end{equation}
Then
\begin{equation}\label{eq:Cabd-dab-binary-W-raw}
      V^2+CUV=U(U+A)(U+B).
\end{equation}
After the scaling
\begin{equation}\label{eq:Cabd-dab-binary-normalization}
                  X=U/C^2=d^2U,\qquad Y=V/C^3=d^3V,
\end{equation}
one obtains
\begin{equation}\label{eq:Cabd-dab-binary-W}
 \boxed{\quad
 Y^2+XY=X^3+(a+b)^2X^2+d^2X.
 \quad}
\end{equation}
Finally, replacing \(Y\) by \(Y+(a+b)X\) changes the \(X^2\)-coefficient
from \((a+b)^2\) to \(a+b\), and hence gives the Jacobian form
\eqref{eq:Cabd-binary-Jacobian}.
On the dense open set on which the displayed denominators are nonzero,
an inverse map is obtained as follows.  If \(Y_0\) denotes the ordinate
in \eqref{eq:Cabd-dab-binary-W} before the final replacement and
\(Y_J=Y_0+(a+b)X\), put
\[
 U=C^2X,\qquad V=C^3Y_0=C^3\bigl(Y_J+(a+b)X\bigr),
\]
and then
\begin{equation}\label{eq:Cabd-dab-binary-inverse}
 x=\frac{V}{U(U+B)},\qquad y=\frac1{Ux},\qquad
 r=\frac{x}{b},\quad s=\frac{y}{a},\qquad
 u=\frac{r}{1+r},\quad v=\frac{s}{1+s}.
\end{equation}
\end{theorem}

\begin{proof}
In characteristic two,
\[
 u^2+u+a=\frac{ar^2+r+a}{(1+r)^2},
\qquad
 v^2+v+b=\frac{bs^2+s+b}{(1+s)^2}.
\]
After multiplication by
\((1+r)^2(1+s)^2\), the equation with \(d=ab\) is
\[
 (ar^2+r+a)(bs^2+s+b)=ab(1+r^2)(1+s^2).
\]
The terms \(abr^2s^2,abr^2,abs^2,ab\) cancel, leaving
\[
              ar^2s+brs^2+rs+br+as=0.
\]
Substitution of \eqref{eq:Cabd-dab-binary-scaled} gives
\eqref{eq:Cabd-dab-binary-Huff}.

The latter equation implies
\[
       Bx+Ay+C=U(x+y),\qquad xy=U^{-1}.
\]
Eliminating \(y\) gives
\[
              U(U+B)x^2+CUx+(U+A)=0.
\]
Multiplying this equality by \(U(U+B)\) and using
\(V=U(U+B)x\) yields
\[
       V^2+CUV=U(U+A)(U+B),
\]
which proves \eqref{eq:Cabd-dab-binary-W-raw}.  Substitute
\(U=C^2X\), \(V=C^3Y\) and divide by \(C^6\).  Since
\[
 \frac{A+B}{C^2}=a^2+b^2=(a+b)^2,
 \qquad
 \frac{AB}{C^4}=a^2b^2=d^2,
\]
equation~\eqref{eq:Cabd-dab-binary-W} follows.  If
\(Y'=Y+(a+b)X\), then
\[
 {Y'}^2+XY'=Y^2+XY+\bigl((a+b)^2+(a+b)\bigr)X^2.
\]
Adding this coefficient to the existing \((a+b)^2X^2\) leaves
\((a+b)X^2\), which proves the final assertion.  Conversely,
\eqref{eq:Cabd-dab-binary-inverse} first reverses the translation and
scaling, then reverses \(V=U(U+B)x\), \(U=(xy)^{-1}\), and finally the
definitions of \(x,y,r,s\).  Thus it is inverse to the forward map on a
dense open set.  Smoothness then extends the resulting birational map to
an isomorphism of the completed elliptic curves.
\end{proof}

\section{Specializations and model properties}

The principal cases are summarized in
Table~\ref{tab:Cabd-specializations}.  A statement that a standard Edwards
equation is obtained over \(k\) includes the indicated square condition;
without it, the diagonal equation remains defined over \(k\) and becomes
standard Edwards after the stated quadratic extension.

\begin{longtable}{L{0.16\textwidth}L{0.20\textwidth}L{0.27\textwidth}L{0.25\textwidth}}
\caption{Distinguished subfamilies of \(\mathcal C_{a,b,d}\)}
\label{tab:Cabd-specializations}\\
\toprule
parameters & rational-point structure & odd-characteristic model & additional structure \\
\midrule
\endfirsthead
\toprule
parameters & rational-point structure & odd-characteristic model & additional structure \\
\midrule
\endhead
General \(a,b,d\)
& May be a nontrivial genus-one torsor
& \(\lambda x^2+\mu y^2=1+\rho x^2y^2\); standard twisted Edwards over
  \(k(\sqrt\mu)\)
& Klein four factor symmetry; Jacobian given by
  \eqref{eq:Cabd-odd-Jacobian} or \eqref{eq:Cabd-binary-Jacobian} \\

\(b=0\)
& Two rational boundary points
& \((1-4a)x^2+y^2=1+(1-4a-16d)x^2y^2\)
& A genuine twisted Edwards curve over \(k\); direct binary map
  \eqref{eq:Cabd-binary-b-zero-map} \\

\(a=b=0\)
& The marked boundary of \(\mathcal C_d\)
& \(x^2+y^2=1+(1-16d)x^2y^2\)
& The original \(\mathcal C_d\) arithmetic model and all preceding chapters \\

\(a=b\)
& Not automatic for general \(d\)
& \(\lambda(x^2+y^2)=1+(\lambda^2-16d)x^2y^2\)
& Full coordinate-exchange symmetry; \(D_8\) acts over \(k\) \\

\(d=ab\ne0\)
& Four explicit affine points in \eqref{eq:Cabd-four-affine-points}
& Generalized Huff cubic and the explicit Weierstrass maps
  \eqref{eq:Cabd-dab-Weierstrass}--\eqref{eq:Cabd-Huff-inverse}
& Rational base point without extracting \(\sqrt\lambda\) or
  \(\sqrt\mu\) \\

\(a=b,\ d=a^2\)
& Four explicit affine points
& Symmetric diagonal equation plus the \(d=ab\) Huff map
& Simultaneously rationally pointed and coordinate symmetric \\
\bottomrule
\end{longtable}

The characteristic-dependent Weierstrass targets are collected in
Table~\ref{tab:Cabd-characteristic-forms}.

\begin{table}[H]
\centering
\caption{Weierstrass forms attached to \(\mathcal C_{a,b,d}\)}
\label{tab:Cabd-characteristic-forms}
\begin{tabular}{L{0.17\textwidth}L{0.36\textwidth}L{0.35\textwidth}}
\toprule
characteristic & base-field Jacobian & reduced or short form \\
\midrule
\(>3\)
& \(Y^2=X^3+2cX^2+\delta^2X\)
& \(\mathsf Y^2=\mathsf X^3+A_4\mathsf X+A_6\), with
  \eqref{eq:Cabd-short-coefficients} \\
\(3\)
& The same equation, with \(4=16=1\)
& \(Y^2=Z^3+a_2Z^2+B_6\) if \(a_2\ne0\); ordinary short form only
  on the locus \(a_2=0\) \\
\(2\)
& \(Y^2+XY=X^3+(a+b)X^2+d^2X\)
& Already a sparse ordinary binary form; the \(XY\)-term must remain \\
\bottomrule
\end{tabular}
\end{table}

\section{Why the generalization is structurally natural}

The equation~\eqref{eq:Cabd-definition} is not obtained by appending two
arbitrary constants to a convenient formula.  Each factor is the same
quadratic cover of \(\PP^1\), and the equation couples the two covers only
through their product.  The parameters \(a\) and \(b\) move the branch
divisors of the two covers, while \(d\) controls their multiplicative
coupling.  The exchange
\[
       (a,u)\longleftrightarrow(b,v)
\]
is therefore built into the parameter space rather than imposed after a
normal-form transformation.

In odd characteristic, centering turns the two factors into the diagonal
quadratics \(r^2-\lambda\) and \(s^2-\mu\); reciprocal coordinates then
turn their product into the diagonal Edwards equation
\eqref{eq:Cabd-diagonal-Edwards}.  In characteristic two, exactly the same
factors remain separable Artin--Schreier polynomials, and \(a+b\) becomes
the twist parameter of the binary Jacobian.  Thus reflection, reciprocal
Edwards geometry, and Artin--Schreier translation are three realizations of
one factorized symmetry.

The special loci have transparent intrinsic meanings.  The condition
\(b=0\) makes one branch divisor rational and produces a twisted Edwards
identity without extending the field.  The condition \(a=b\) identifies the
two covers and upgrades parameter exchange to an automorphism.  The condition
\(d=ab\) forces the four rational points
\eqref{eq:Cabd-four-affine-points} and gives a base-field Huff--Weierstrass
map.  Finally \(a=b=0\) selects the member with the rational marked boundary
and characteristic-uniform arithmetic developed in the rest of this book.
The generalization therefore explains why \(\mathcal C_d\) is especially
suited to explicit arithmetic while placing it inside a larger symmetric
genus-one geometry.

\section{Research program for the full product family}
\label{sec:Cabd-further-directions}

The three-parameter family
\[
 \mathcal C_{a,b,d}:\qquad
 (u^2+u+a)(v^2+v+b)=d
\]
contains a rich collection of arithmetic interfaces.  This chapter
establishes the structural layer:
smoothness, symmetry, the distinction between a genus-one torsor and its
Jacobian, and the explicit Edwards and Weierstrass dictionaries.  The next
Part develops substantial arithmetic theories for
\(\mathcal T_{a,d}=\mathcal C_{a,0,d}\), reciprocal \(C\)-curves, and the
symmetric QRT envelope.  The following items form a coherent research
program for extending these constructions across the full
\(\mathcal C_{a,b,d}\) parameter space.

\begin{enumerate}[label=\textup{(\arabic*)},leftmargin=2.5em]

 \item Construct model-internal affine, projective, mixed, and complete
 addition laws on the smooth completion of
 \(\mathcal C_{a,b,d}\).  In particular, determine precisely which
 parameter subfamilies admit a single low-degree \(k\)-complete addition
 law and which require a finite addition-law atlas.

 \item Determine the most efficient doubling, tripling, differential
 addition, closed \(2P+Q\), halving, and scalar-multiplication formulas
 in the original \((u,v)\)-coordinates.  The resulting costs should be
 compared with those of twisted Edwards, Montgomery, Huff, Jacobi
 quartic, and other explicitly identified models under matching input,
 output, normalization, and completeness assumptions.

 \item Develop the intrinsic Kummer theory of the generalized family.
 The two factor involutions suggest several natural quotient functions,
 but the optimal quotient can depend on the rationality of the boundary
 divisors and on whether the original genus-one curve is a trivial torsor
 under its Jacobian.

 \item Classify the \(k\)-isomorphism classes, geometric isomorphism
 classes, twists, and torsor classes represented by
 \(\mathcal C_{a,b,d}\).  The loci
 \[
        b=0,\qquad a=b,\qquad d=ab
 \]
 should be viewed as distinguished strata of the parameter space rather
 than isolated substitutions.  Their intersections, automorphism jumps,
 CM loci, torsion structures, and fields of definition deserve a
 systematic moduli-theoretic treatment.

 \item Construct division polynomials, model-preserving isogenies,
 explicit endomorphisms, Miller functions, Tate pairings, and Weil
 pairings directly in the generalized coordinates.  It is especially
 interesting to determine which of these constructions descend to a
 genus-one torsor before a rational base point is chosen.

 \item Establish exact first and higher moments of point counts over
 finite fields, average parameter multiplicities, average torsion
 distributions, and statistics of the distinguished subfamilies.  The
 additional parameters \(a\) and \(b\) may make possible averaging
 arguments that are unavailable in the one-parameter family.

 \item Investigate characteristic-specific phenomena beyond the initial
 dictionaries of this chapter.  In characteristic two, the interaction
 between the two Artin--Schreier classes \(a\) and \(b\) should be studied
 at the level of complete arithmetic, isogenies, and descent.  In
 characteristic three, one should examine whether special parameter
 strata yield new Frobenius--Verschiebung factorizations or especially
 sparse tripling formulas.

 \item Study cryptographic parameter generation, encoding, decoding,
 twist security, subgroup structure, side-channel resistance, and
 constant-time implementations for selected subfamilies.  Such an
 investigation must distinguish the arithmetic advantages of a model
 from the security properties of a particular curve and protocol.

\end{enumerate}

\subsection{Elliptic-curve factorization}

A particularly promising direction is the elliptic curve method for
integer factorization.  Both \(\mathcal C_d\) and
\(\mathcal C_{a,b,d}\) should be investigated as arithmetic models for
ECM over \(\mathbb Z/N\mathbb Z\).

The one-parameter family \(\mathcal C_d\) already provides several
features relevant to an ECM implementation:

\begin{itemize}[leftmargin=2.2em]

 \item a factorized equation with inexpensive parameter handling;

 \item a native Kummer coordinate supporting differential addition and
 constant-time ladder structures;

 \item a native-input first doubling requiring only
 \(\M+\Sqr+\Dpar\) under the hypotheses established earlier;

 \item specialized doubling, mixed-addition, tripling, and closed
 \(2P+Q\) formulas that may be combined with ECM addition chains;

 \item rational marked torsion and explicit transformations to Edwards
 and Montgomery models, which permit a direct comparison with existing
 Edwards-ECM and Montgomery-ECM implementations.

\end{itemize}

The generalized family introduces two further parameters with which one may
attempt to control torsion, small constants, discriminants, twists, and the
probability that the group order modulo an unknown prime divisor is suitably
smooth.  In particular, the subfamilies \(b=0\), \(a=b\), and \(d=ab\)
offer different compromises among rational torsion, symmetry, parameter
freedom, and arithmetic cost.  This additional flexibility may lead to
useful analogues of the parameterizations employed in Montgomery and
Edwards ECM.

These structural features make \(\mathcal C_d\) and
\(\mathcal C_{a,b,d}\) strong foundations for a dedicated ECM theory and
implementation program.  Their factorized equations, native Kummer
coordinates, inexpensive initial doubling, specialized differential
arithmetic, rational marked torsion, and flexible parameter spaces can be
developed into an end-to-end factorization framework through the following
five closely connected components:

\begin{enumerate}[label=\textup{(\roman*)},leftmargin=2.5em]

 \item construct torsion-aware parameterizations that generate large
 admissible families with efficient parameter selection;

 \item derive exact Stage~1 and Stage~2 operation ledgers, including curve
 generation, coordinate conversion, normalization, recovery, and
 failure-detection costs;

 \item determine the distribution of smooth group orders modulo the unknown
 prime divisors and identify the parameter strata that maximize the relevant
 torsion and smoothness probabilities;

 \item develop optimized constant-pattern and variable-time implementations
 over \(\mathbb Z/N\mathbb Z\), using the native Kummer ladder,
 specialized first doubling, mixed arithmetic, batch inversion, and
 product-tree techniques;

 \item produce reproducible implementation comparisons with Montgomery-ECM
 and Edwards-ECM software at identical bounds, integer sizes, parameter
 families, and output requirements.

\end{enumerate}

Together, these components turn the native arithmetic and parameter
flexibility of \(\mathcal C_d\) and \(\mathcal C_{a,b,d}\) into a concrete
ECM platform.  In particular, their marked torsion, characteristic-uniform
quotient structure, inexpensive native initialization, and explicit
model dictionaries provide a substantial foundation for new Stage~1 and
Stage~2 algorithms and for reproducible high-performance implementations.

\subsection{Structural economy of the factorized family}

The structural content of these curves is not exhausted by operation counts.  Their
equations exhibit an unusual concentration of compatible structures.
The same quadratic expression
\[
                         T^2+T+c
\]
appears in both coordinates, and the two copies are coupled by a single
multiplicative parameter.  From this simple construction arise the
involutions
\[
 (u,v)\longmapsto(-u-1,v),\qquad
 (u,v)\longmapsto(u,-v-1),
\]
together with the parameter-exchange symmetry
\[
              (a,u)\longleftrightarrow(b,v).
\]
When \(a=b\), parameter exchange becomes an automorphism and completes the
visible \(D_8\)-symmetry.

In odd characteristic, centering transforms the equation into
\[
 (r^2-\lambda)(s^2-\mu)=\delta,
\]
and reciprocal coordinates reveal the diagonal Edwards equation
\[
             \lambda x^2+\mu y^2
             =1+\rho x^2y^2.
\]
In characteristic two, the same factors do not degenerate into an awkward
reduction of an odd-characteristic formula.  They become genuine
Artin--Schreier operators, and the combination \(a+b\) records the twist of
the binary Jacobian.  Thus reflection symmetry in odd characteristic and
additive Artin--Schreier symmetry in characteristic two are two
manifestations of the same equation.

The special parameter conditions also have intrinsic geometric meanings:
\begin{center}
\begin{tabular}{L{2.0cm}L{10.5cm}}
\toprule
condition & geometric effect\\
\midrule
\(b=0\)
   & one boundary divisor becomes rational and a ground-field twisted
     Edwards chart appears;\\
\(a=b\)
   & the two quadratic covers coincide and coordinate exchange becomes an
     automorphism;\\
\(d=ab\)
   & four explicit affine rational points appear;\\
\(a=b=0\)
   & the rational marked boundary and the native arithmetic of
     \(\mathcal C_d\) are recovered.\\
\bottomrule
\end{tabular}
\end{center}
These cancellations identify visible arithmetic strata in a symmetric
parameter space.

For these reasons, \(\mathcal C_d\) and
\(\mathcal C_{a,b,d}\) provide a compact family in which symmetry, marked
torsion, reciprocal Edwards geometry, Artin--Schreier descent, Kummer
quotients, and explicit arithmetic can be studied in the original
coordinates.  The economy of this presentation makes the family a compact
laboratory in which these structures can be developed together.  The
structural theory and explicit arithmetic of the distinguished strata supply
the foundations for extending complete formulas, isogenies, and
characteristic-specific circuits throughout the three-parameter family.

\chapter{The One-Sided Twisted Family}
\label{ch:one-sided-twist}

\section{Coverage of the full product family}
Let $k$ be a field.  Retain the polynomial
\[
  f_c(t)=t^2+t+c
\]
and the three-parameter family $\mathcal C_{a,b,d}$ from
Chapter~\ref{ch:generalized-Cabd}.  Its one-sided twisted subfamily is
\begin{equation}
  \T_{a,d}=\C_{a,0,d}:\qquad
  (u^2+u+a)(v^2+v)=d.
  \label{eq:Tad}
\end{equation}
Within this monograph, \(\T_{a,d}\) is also called the
\emph{one-sided twisted \(C_d\) curve}, or simply a \emph{twisted
\(C_d\) curve}; the original family is the specialization
\(C_d=\T_{0,d}\).
All genus statements refer to the smooth completion in
$\PP^1\times\PP^1$.  A smooth curve of bidegree $(2,2)$ has genus one.
Whenever a $k$-rational point has been selected, the curve is regarded as
an elliptic curve; without such a point it is a genus-one torsor under its
Jacobian.

The phrase that the family \eqref{eq:Tad} \emph{covers} the full family
\eqref{eq:Cabd-definition} will be used in the following precise sense.

\begin{enumerate}[label=(\roman*)]
  \item After a splitting extension, every smooth member of
  $\C_{a,b,d}$ is isomorphic, by a separate fractional-linear change in
  the two $\PP^1$ coordinates, to a member of $\T_{a,d}$.

  \item Over the ground field, the Jacobian of every smooth member of
  $\C_{a,b,d}$ is represented by a member of $\T_{a,d}$; consequently,
  the same is true for the curve itself whenever its torsor class is
  trivial.

  \item Over a finite field every smooth genus-one curve has a rational
  point.  Thus, in the finite-field setting, the preceding Jacobian
  statement becomes an actual ground-field isomorphism statement.
\end{enumerate}

This distinction is essential.  A parameter such as $b$ may record the
splitting behavior of a marked degree-two boundary divisor or a nontrivial
genus-one torsor, even when it does not produce a new Jacobian isomorphism
class.

\section{Geometry and basic symmetries}
The model $\T_{a,d}$ retains the two coordinate involutions
\begin{equation}
  \iota_u(u,v)=(-1-u,v),
  \qquad
  \iota_v(u,v)=(u,-1-v).
  \label{eq:T-involutions}
\end{equation}
In characteristic $2$ these become $u\mapsto u+1$ and $v\mapsto v+1$.
The second quadratic factor is split over every ground field:
$v^2+v=v(v+1)$.  Consequently, one of the two boundary degree-two
divisors is always rational and split.  The first factor carries the
remaining quadratic or Artin--Schreier twist information.  Thus the two
parameters of $\T_{a,d}$ have conceptually different roles: one controls
the geometric moduli, while the other records a ground-field twist.
On every smooth completion, the product $\iota_u\iota_v$ is fixed-point
free.  After a geometric origin is chosen, it is translation on the
genus-one curve by a nonzero two-torsion displacement class
\(Q\in\Pic^0[2]\).  Because the automorphism is defined over \(k\), this
class is \(k\)-rational; the induced pullback action on \(\Pic^0\), and
hence on the Jacobian, is the identity.  Thus the one-sided twist preserves the basic
reflection--reflection--translation mechanism of $C_d$, even though the
coordinate interchange $u\leftrightarrow v$ is generally lost when
$a\ne0$.

\section{Direct reduction of the three-parameter family}

\begin{proposition}[Coordinate-preserving reduction]
\label{prop:direct-reduction}
Let $\C_{a,b,d}$ be smooth.

\begin{enumerate}[label=(\alph*)]
  \item Suppose $\charac(k)\ne2$.  Set
  \[
    r=2u+1,\qquad s=2v+1,
    \qquad
    \lambda=1-4a,\quad \mu=1-4b,\quad \delta=16d.
  \]
  Then
  \begin{equation}
    (r^2-\lambda)(s^2-\mu)=\delta.
    \label{eq:centered-Cab}
  \end{equation}
  If $\mu$ is a square in $k$, smoothness forces $\mu\ne0$, so write
  $\mu=m^2$ with $m\in k^\times$.  Then $s=mS$ gives
  \begin{equation}
    (r^2-\lambda)(S^2-1)=\frac{\delta}{\mu},
    \label{eq:direct-odd-reduction}
  \end{equation}
  which is a member of the centered one-sided family.

  \item Suppose $\charac(k)=2$.  If $b=\beta^2+\beta$ for some
  $\beta\in k$, then $V=v+\beta$ gives
  \[
    v^2+v+b=V^2+V,
  \]
  and hence $\C_{a,b,d}\cong\T_{a,d}$ over $k$.

  \item Over a separable closure, the required square root in (a), or
  the required Artin--Schreier root in (b), always exists.  Thus every
  smooth $\C_{a,b,d}$ is geometrically a member of the one-sided family.
\end{enumerate}
\end{proposition}

\begin{proof}
Part (a) follows from
\[
  u^2+u+a=\frac{r^2-\lambda}{4},
  \qquad
  v^2+v+b=\frac{s^2-\mu}{4}.
\]
Substituting $s=mS$ in \eqref{eq:centered-Cab} and dividing by $\mu$
gives \eqref{eq:direct-odd-reduction}.  In characteristic $2$ one has
$(v+\beta)^2+(v+\beta)=v^2+v+\beta^2+\beta$, proving (b).  For the final assertion, first note that smoothness in odd
characteristic forces $\mu\ne0$: if $\mu=0$, the bihomogeneous form of
\eqref{eq:centered-Cab} is
\[
  (R^2-\lambda Z^2)S^2=16dZ^2W^2,
\]
and the point $((R:Z),(S:W))=((1:0),(0:1))$ is singular.  Hence
$X^2-\mu$ is separable.  In characteristic $2$, the polynomial
$X^2+X+b$ has derivative $1$ and is therefore separable for every
$b$.  The necessary roots consequently exist over a separable closure,
which proves the geometric reduction.
\end{proof}

\begin{lemma}[Exact smoothness of the centered product in odd characteristic]
\label{lem:Cab-smoothness-odd}
Assume $\charac(k)\ne2$ and write the centered completion as
\begin{equation}
  (R^2-\lambda Z^2)(S^2-\mu W^2)=\delta Z^2W^2
  \quad\text{in }\PP^1_{(R:Z)}\times\PP^1_{(S:W)}.
  \label{eq:Cab-centered-bihomogeneous}
\end{equation}
It is a smooth curve of genus one if and only if
\begin{equation}
  \lambda\mu\delta(\lambda\mu-\delta)\ne0.
  \label{eq:Cab-exact-smoothness}
\end{equation}
\end{lemma}

\begin{proof}
Let $F$ be the left-hand side of
\eqref{eq:Cab-centered-bihomogeneous} minus its right-hand side.  On the
affine chart $Z=W=1$, the equation is
$(r^2-\lambda)(s^2-\mu)=\delta$ and
\[
  F_r=2r(s^2-\mu),
  \qquad
  F_s=2s(r^2-\lambda).
\]
If $\delta\ne0$, neither factor $r^2-\lambda$ nor $s^2-\mu$ vanishes on
the curve.  Hence an affine singular point must have $r=s=0$, and this
point lies on the curve exactly when $\lambda\mu=\delta$.

On the boundary $Z=0$, one has $R\ne0$ and
$S^2=\mu W^2$.  The case $W=0$ is impossible.  If $\mu\ne0$, then
$S\ne0$ and
\[
  F_S=2S(R^2-\lambda Z^2)=2SR^2\ne0.
\]
If $\mu=0$, the point $((1:0),(0:1))$ lies on the curve and all four
first derivatives of $F$ vanish there.  The boundary $W=0$ is symmetric:
it is smooth when $\lambda\ne0$, whereas $\lambda=0$ produces the
singular point $((0:1),(1:0))$.  The corner $Z=W=0$ is not on the curve,
because there $F=R^2S^2\ne0$.  Finally, if $\delta=0$, the equation is
the union of the two degree-two fibre divisors and is singular at their
geometric intersections.  These cases are exhaustive and prove the
criterion.  Under the criterion, a smooth curve of bidegree $(2,2)$ has
genus $(2-1)(2-1)=1$.
\end{proof}

The obstruction in Proposition~\ref{prop:direct-reduction} is a genuine
ground-field obstruction.  If $\mu$ is nonsquare in odd characteristic,
or if $b\notin\wpmap(k)$ in characteristic $2$, where
\[
  \wpmap(k)=\{z^2+z:z\in k\},
\]
then the second marked double cover does not split over $k$.  The next
results show that this extra boundary datum does not require a larger
family of Jacobians.

\section[Odd-characteristic Edwards and Montgomery forms]
{Odd characteristic: exact relation with twisted Edwards and Montgomery forms}

Assume throughout this section that $\charac(k)\ne2$.  For the
one-sided model, set
\begin{equation}
  r=2u+1,\qquad s=2v+1,
  \qquad
  \lambda=1-4a,\qquad \delta=16d.
  \label{eq:T-centered-parameters}
\end{equation}
Then
\begin{equation}
  \T_{a,d}:\qquad (r^2-\lambda)(s^2-1)=\delta.
  \label{eq:T-centered}
\end{equation}
On the open set $rs\ne0$, put
\begin{equation}
  x=\frac1r,\qquad y=\frac1s.
  \label{eq:T-to-TE-map}
\end{equation}
Multiplication of \eqref{eq:T-centered} by $x^2y^2$ gives
\begin{equation}
  \lambda x^2+y^2
  =1+(\lambda-\delta)x^2y^2.
  \label{eq:T-to-TE}
\end{equation}
Thus $\T_{a,d}$ is birational over $k$ to the twisted Edwards curve
\begin{equation}
  E_{A,D}:\qquad Ax^2+y^2=1+Dx^2y^2
  \label{eq:twisted-Edwards}
\end{equation}
with
\begin{equation}
  A=\lambda=1-4a,
  \qquad
  D=\lambda-\delta=1-4a-16d.
  \label{eq:T-TE-parameters}
\end{equation}
Conversely, every twisted Edwards equation \eqref{eq:twisted-Edwards}
with $AD(A-D)\ne0$ is obtained from
\begin{equation}
  a=\frac{1-A}{4},
  \qquad
  d=\frac{A-D}{16}.
  \label{eq:TE-to-T-parameters}
\end{equation}
The exceptional points of the affine reciprocal map are precisely the
points that are restored on the smooth $(2,2)$ completion.

The standard twisted-Edwards--Montgomery dictionary
\cite{BernsteinEtAl2008} sends \eqref{eq:twisted-Edwards} to
\begin{equation}
  B_M Y^2=X^3+A_MX^2+X,
  \label{eq:Montgomery}
\end{equation}
where
\begin{equation}
  A_M=\frac{2(A+D)}{A-D},
  \qquad
  B_M=\frac{4}{A-D}.
  \label{eq:TE-Montgomery-parameters}
\end{equation}
For $\T_{a,d}$ this becomes
\begin{equation}
  A_M=\frac{4\lambda}{\delta}-2,
  \qquad
  B_M=\frac4\delta.
  \label{eq:T-Montgomery-parameters}
\end{equation}
Consequently, the one-sided family is exactly the $C$-curve realization
of the class of Montgomery curves admitting a twisted Edwards model,
rather than merely the
untwisted Edwards subfamily represented by $C_d$.  Indeed, $a=0$ gives
$A=1$ in \eqref{eq:T-TE-parameters}, whereas arbitrary $a$ permits an
arbitrary nonzero twisted Edwards coefficient $A$.

The smoothness condition is
\begin{equation}
  \lambda\delta(\lambda-\delta)\ne0.
  \label{eq:T-smooth-odd}
\end{equation}
The $j$-invariant can be read from the twisted Edwards formula:
\begin{equation}
  j(\T_{a,d})
  =16\frac{\bigl(\lambda^2+14\lambda(\lambda-\delta)
       +(\lambda-\delta)^2\bigr)^3}
       {\lambda(\lambda-\delta)\delta^4}
  =16\frac{(16\lambda^2-16\lambda\delta+\delta^2)^3}
       {\lambda(\lambda-\delta)\delta^4}.
  \label{eq:T-j-odd}
\end{equation}

\section{The Jacobian of the full family in odd characteristic}

We first recall the even-quartic Jacobian formula.

\begin{lemma}[Jacobian of an even quartic]
\label{lem:even-quartic-jacobian}
Let $\charac(k)\ne2$, and let
\begin{equation}
  \mathcal X:\qquad z^2=Ax^4+Bx^2+C_0,
  \label{eq:even-quartic-general}
\end{equation}
where $AC_0(B^2-4AC_0)\ne0$.  Then
\begin{equation}
  \Jac(\mathcal X):\qquad
  V^2=U\bigl(U^2-2BU+B^2-4AC_0\bigr).
  \label{eq:even-quartic-jacobian}
\end{equation}
\end{lemma}

\begin{proof}
For the binary quartic
$AX^4+BX^2Z^2+C_0Z^4$, put
\[
  I=B^2+12AC_0,
  \qquad
  J=72ABC_0-2B^3.
\]
The binary-quartic Jacobian theorem states that the Jacobian of this
smooth quartic is
\[
  V^2=U\bigl(U^2-2BU+(B^2-4AC_0)\bigr);
\]
see, for example, the invariant-theoretic construction in
\cite{CremonaEtAlDescent}.  The theorem applies because
$AC_0(B^2-4AC_0)\ne0$.  For completeness, we verify the normalization of
this model by calculating its invariants and discriminant.  Write
\[
  D=B^2-4AC_0.
\]
For
\[
  V^2=U^3-2BU^2+DU,
\]
the generalized Weierstrass coefficients are
$a_1=a_3=a_6=0$, $a_2=-2B$, and $a_4=D$.  Therefore
\[
  b_2=-8B,
  \qquad b_4=2D,
  \qquad b_6=0,
  \qquad b_8=-D^2.
\]
Substitution in the Weierstrass invariant formulas gives
\begin{align*}
  c_4
  &=b_2^2-24b_4
    =64B^2-48(B^2-4AC_0)
    =16(B^2+12AC_0)=16I,\\
  c_6
  &=-b_2^3+36b_2b_4-216b_6\\
  &=512B^3-576B(B^2-4AC_0)
    =32(72ABC_0-2B^3)=32J,
\end{align*}
and
\begin{align*}
  \Delta
  &=-b_2^2b_8-8b_4^3-27b_6^2+9b_2b_4b_6\\
  &=64B^2D^2-64D^3
    =64D^2(B^2-D)
    =256AC_0(B^2-4AC_0)^2.
\end{align*}
The hypothesis $AC_0(B^2-4AC_0)\ne0$ makes this discriminant nonzero.
Thus the cited binary-quartic Jacobian theorem applies without a
singular or repeated-root exception, and the calculation verifies the
normalization used in \eqref{eq:even-quartic-jacobian}.
\end{proof}

\begin{theorem}[One-sided representative of every odd-characteristic Jacobian]
\label{thm:odd-Jacobian-coverage}
Assume $\charac(k)\ne2$, and write
\[
  \lambda=1-4a,
  \qquad
  \mu=1-4b,
  \qquad
  \delta=16d.
\]
Suppose
\begin{equation}
  \lambda\mu\delta(\lambda\mu-\delta)\ne0,
  \label{eq:Cab-smooth-odd}
\end{equation}
which is exactly the smoothness condition by
Lemma~\ref{lem:Cab-smoothness-odd}.  Then its Jacobian is birational over
$k$ to the one-sided curve $\T_{a_0,d_0}$ determined by
\begin{equation}
  1-4a_0=\frac{4\lambda\mu}{\delta^2},
  \qquad
  16d_0=\frac4\delta.
  \label{eq:odd-Jac-T-centered-parameters}
\end{equation}
Equivalently,
\begin{equation}
  a_0=\frac14-\frac{(1-4a)(1-4b)}{256d^2},
  \qquad
  d_0=\frac{1}{64d}.
  \label{eq:odd-Jac-T-original-parameters}
\end{equation}
If $\C_{a,b,d}(k)\ne\varnothing$, then $\C_{a,b,d}$ itself is
$k$-isomorphic, after choosing a rational origin, to this member of the
one-sided family.
\end{theorem}

\begin{proof}
The centered equation is
\[
  (r^2-\lambda)(s^2-\mu)=\delta.
\]
Define
\[
  z=(r^2-\lambda)s.
\]
Because the curve equation and $\delta\ne0$ imply
$r^2-\lambda\ne0$ in its function field, one recovers
$s=z/(r^2-\lambda)$.  Thus this substitution is birational.  Moreover,
\begin{align*}
  z^2
  &=(r^2-\lambda)\bigl(\mu r^2+\delta-\lambda\mu\bigr)\\
  &=\mu r^4+(\delta-2\lambda\mu)r^2
    +\lambda(\lambda\mu-\delta).
\end{align*}
For this even quartic, set
\[
  A=\mu,
  \qquad
  B=\delta-2\lambda\mu,
  \qquad
  C=\lambda(\lambda\mu-\delta).
\]
The discriminant term is a square because
\begin{align*}
  B^2-4AC
  &=(\delta-2\lambda\mu)^2
    -4\mu\lambda(\lambda\mu-\delta)\\
  &=\delta^2-4\delta\lambda\mu+4\lambda^2\mu^2
    -4\lambda^2\mu^2+4\delta\lambda\mu\\
  &=\delta^2.
\end{align*}
Thus
\begin{equation}
  B^2-4AC=\delta^2.
  \label{eq:perfect-square-quartic-disc}
\end{equation}
Lemma~\ref{lem:even-quartic-jacobian} therefore yields
\begin{equation}
  V^2
  =U\Bigl(U^2+(4\lambda\mu-2\delta)U+\delta^2\Bigr).
  \label{eq:Cab-Jac-Montgomery-pre}
\end{equation}
Put $U=\delta X$ and $V=\delta^2Y$.  Equation
\eqref{eq:Cab-Jac-Montgomery-pre} becomes
\begin{equation}
  \delta Y^2
  =X^3+\left(\frac{4\lambda\mu}{\delta}-2\right)X^2+X.
  \label{eq:Cab-Jac-Montgomery}
\end{equation}
This is a Montgomery equation with
\[
  A_M=\frac{4\lambda\mu}{\delta}-2,
  \qquad B_M=\delta.
\]
By the inverse of \eqref{eq:TE-Montgomery-parameters}, the associated
twisted Edwards parameters are
\begin{equation}
  A=\frac{A_M+2}{B_M}
    =\frac{4\lambda\mu}{\delta^2},
  \qquad
  D=\frac{A_M-2}{B_M}
    =\frac{4(\lambda\mu-\delta)}{\delta^2}.
  \label{eq:Cab-Jac-TE}
\end{equation}
A one-sided curve with centered parameters
\[
  \lambda_0=\frac{4\lambda\mu}{\delta^2},
  \qquad
  \delta_0=\frac4\delta
\]
has twisted Edwards parameters
\[
  A=\lambda_0,
  \qquad
  D=\lambda_0-\delta_0
    =\frac{4(\lambda\mu-\delta)}{\delta^2}.
\]
Thus its twisted Edwards model is exactly
\eqref{eq:Cab-Jac-TE}, proving
\eqref{eq:odd-Jac-T-centered-parameters}.  Formula
\eqref{eq:odd-Jac-T-original-parameters} follows by substituting
$\delta=16d$.  Finally, a genus-one curve with a $k$-rational point is
isomorphic to its Jacobian once that point is chosen as the origin.
\end{proof}

\begin{remark}
The identity \eqref{eq:perfect-square-quartic-disc} explains why the
three-parameter product family has no larger Jacobian class than the
one-sided family.  The discriminant term that would normally be an
independent square class collapses to the square $\delta^2$, forcing the
Jacobian into Montgomery form over the ground field.
\end{remark}

\section{The distinguished value \texorpdfstring{$a=\tfrac12$}{a=1/2}}

Suppose $\charac(k)\ne2$ and take
\begin{equation}
  a=\frac12.
  \label{eq:a-half}
\end{equation}
Then $\lambda=1-4a=-1$, so \eqref{eq:T-to-TE} becomes
\begin{equation}
  -x^2+y^2=1+D x^2y^2,
  \qquad
  D=-1-16d.
  \label{eq:a-minus-one-TE}
\end{equation}
Thus the coefficient $A=-1$ occurs natively, without adjoining a square
root of $-1$ in order to rescale an $A=1$ Edwards equation.  The
smoothness condition is
\[
  d\ne0,
  \qquad
  -1-16d\ne0.
\]

The importance of $A=-1$ is visible in extended twisted Edwards
coordinates.  Write a projective point as $(X:Y:Z:T)$ with
$x=X/Z$, $y=Y/Z$, and $T=XY/Z$.  For two input points, the standard
$A=-1$ addition dependency graph is
\begin{align*}
  A_0&=(Y_1-X_1)(Y_2-X_2),
  &B_0&=(Y_1+X_1)(Y_2+X_2),\\
  C_0&=2D T_1T_2,
  &D_0&=2Z_1Z_2,\\
  E_0&=B_0-A_0,
  &F_0&=D_0-C_0,\\
  G_0&=D_0+C_0,
  &H_0&=B_0+A_0,\\
  X_3&=E_0F_0,
  &Y_3&=G_0H_0,\\
  T_3&=E_0H_0,
  &Z_3&=F_0G_0.
\end{align*}
If multiplication by the fixed curve parameter $D$ is denoted by
$\Dpar$, this dependency graph costs
\[
  8\M+1\Dpar
\]
for general addition and
\[
  7\M+1\Dpar
\]
for mixed addition with the second input affine.  When multiplication by
$D$ is absorbed into the fixed-constant model, these are often quoted as
$8\M$ and $7\M$, respectively
\cite{HisilEtAl2008}.

Doubling can be organized as
\begin{align*}
  A_0&=X_1^2,
  &B_0&=Y_1^2,
  &C_0&=2Z_1^2,
  &D_0&=-A_0,\\
  E_0&=(X_1+Y_1)^2-A_0-B_0,
  &G_0&=D_0+B_0,\\
  F_0&=G_0-C_0,
  &H_0&=D_0-B_0,\\
  X_3&=E_0F_0,
  &Y_3&=G_0H_0,
  &T_3&=E_0H_0,
  &Z_3&=F_0G_0,
\end{align*}
with cost $4\M+4\Sqr$.

The specialized dependency graph and the completeness of a single
addition law are logically distinct.  A standard sufficient condition
for completeness of the affine twisted Edwards law is that $A$ be a
square and $D$ a nonsquare.  Hence for $A=-1$ one obtains this criterion
when $-1$ is a square and $D=-1-16d$ is a nonsquare.  If $-1$ is a
nonsquare, the model remains valid and the specialized formulas remain
useful, but completeness must be supplied by a suitable addition-law
atlas or by restrictions appropriate to the chosen subgroup.

\chapter[One-Sided Twists in Characteristic Two]{The One-Sided Twisted Family in Characteristic Two}
\label{ch:one-sided-twist-char2}

\section{Native Weierstrass form and parameter separation}
Throughout this chapter assume that $\charac(k)=2$.  Put
\[
  F=u^2+u+a.
\]
Then $\T_{a,d}$ is
\[
  F(v^2+v)=d.
\]
The characteristic-two geometry is intrinsically Artin--Schreier rather
than a degeneration of the odd-characteristic reciprocal substitution.
Nevertheless, the one-sided model admits an especially simple
Weierstrass dictionary.

\begin{theorem}[Native Weierstrass model in characteristic two]
\label{thm:T-char2-Weierstrass}
Let $k$ have characteristic $2$ and let $d\ne0$.  The curve
$\T_{a,d}$ is birational to
\begin{equation}
  E_{a,d}:\qquad
  Y^2+XY=X^3+aX^2+d^2X.
  \label{eq:T-char2-Weierstrass}
\end{equation}
A birational map is
\begin{equation}
  X=Fv+d,
  \qquad
  Y=uX.
  \label{eq:T-char2-forward}
\end{equation}
On the open set where the expressions are defined, the inverse is
\begin{equation}
  u=\frac{Y}{X},
  \qquad
  v=\frac{X}{X+d}.
  \label{eq:T-char2-inverse}
\end{equation}
Moreover,
\begin{equation}
  \Delta(E_{a,d})=d^4,
  \qquad
  j(E_{a,d})=d^{-4}.
  \label{eq:T-char2-invariants}
\end{equation}
\end{theorem}

\begin{proof}
Let $T=Fv$.  Multiplying $F(v^2+v)=d$ by $F$ gives
\[
  T^2+FT=dF.
\]
Since $X=T+d$, one obtains
\begin{equation}
  T^2=F(T+d)=FX.
  \label{eq:T-square-FX}
\end{equation}
On the other hand, $T=X+d$, so
\[
  X^2+d^2=FX.
\]
Multiplying by $X$ yields
\begin{equation}
  FX^2=X^3+d^2X.
  \label{eq:FX2}
\end{equation}
Because $Y=uX$ and $u^2+u=F+a$ in characteristic $2$,
\begin{align*}
  Y^2+XY
  &=(u^2+u)X^2\\
  &=(F+a)X^2\\
  &=X^3+aX^2+d^2X,
\end{align*}
which proves \eqref{eq:T-char2-Weierstrass}.

Conversely, on \eqref{eq:T-char2-Weierstrass}, set $u=Y/X$.  Then
\begin{align*}
  u^2+u+a
  &=\frac{Y^2+XY+aX^2}{X^2}\\
  &=\frac{X^3+d^2X}{X^2}
    =\frac{X^2+d^2}{X}
    =\frac{(X+d)^2}{X}.
\end{align*}
For $v=X/(X+d)$,
\[
  v^2+v
  =\frac{X^2}{(X+d)^2}+\frac{X}{X+d}
  =\frac{dX}{(X+d)^2}.
\]
Their product is $d$, proving the inverse formula.

For \eqref{eq:T-char2-Weierstrass}, the Weierstrass coefficients are
$a_1=1$, $a_2=a$, $a_3=0$, $a_4=d^2$, and $a_6=0$.  Hence
$c_4=1$ and $\Delta=d^4$, which gives $j=d^{-4}$.
\end{proof}

\section{A twist-stable separable two-isogeny}
\label{sec:T-char2-two-isogeny}

The sparse Weierstrass equation does more than classify the binary twist:
it makes the reduced two-torsion quotient independent of the twist
coefficient.  This yields a degree-two operation whose online circuit is
shared by the one-sided and full product extensions.

\begin{theorem}[Twist-stable binary two-isogeny]
\label{thm:T-char2-two-isogeny}
Let \(k\) be perfect of characteristic two, let \(d\ne0\), and let
\(e^2=d\).  For any \(A\in k\), put
\[
 W_{A,d}:\qquad y^2+xy=x^3+Ax^2+d^2x .
\]
Then
\begin{equation}
 \begin{aligned}
 X&=x+\frac{d^2}{x},\\
 Y&=y+\frac{d^2(y+x)}{x^2}+d
 \end{aligned}
\label{eq:T-char2-two-isogeny-W}
\end{equation}
extends to a separable isogeny
\begin{equation}
        \Phi_{A,d}:W_{A,d}\longrightarrow W_{A,e}
\label{eq:T-char2-two-isogeny-arrow}
\end{equation}
of degree two and kernel \(\{O,(0,0)\}\).  In particular, \(A\) is
preserved exactly and \(d\) is updated by the inverse Frobenius
\[
                            d\longmapsto e=\sqrt d.
\]
Its dual is
\begin{equation}
 \widehat\Phi_{A,d}:W_{A,e}\longrightarrow W_{A,d},\qquad
 (X,Y)\longmapsto\bigl(X^2,Y^2+AX^2\bigr),
\label{eq:T-char2-two-isogeny-dual}
\end{equation}
and the two compositions are the corresponding multiplication-by-two maps.
\end{theorem}

\begin{proof}
\noindent\emph{Step 1: the kernel translation is independent of \(A\).}
For \(T=(0,0)\) and \(P=(x,y)\), the chord slope is \(y/x\).  Since
\[
 \left(\frac yx\right)^2+\frac yx=x+A+\frac{d^2}{x},
\]
the generalized binary addition law gives
\[
 x(P+T)=\frac{d^2}{x},\qquad
 y(P+T)=\frac{d^2(y+x)}{x^2}.
\]
The term \(A\) cancels from the abscissa because it occurs once in the
curve relation and once in the chord formula.  Consequently
\(X=x+x(P+T)\) and
\(Y_0=y+y(P+T)\) are invariant under the order-two translation.

\smallskip
\noindent\emph{Step 2: compute the target equation.}
Set \(c=d^2\).  Expansion gives
\begin{align*}
 Y_0^2+XY_0
 &=(y^2+xy)\left(1+\frac{c^2}{x^4}\right)+c\\
 &=(x^3+Ax^2+cx)\left(1+\frac{c^2}{x^4}\right)+c\\
 &=X^3+AX^2+c.
\end{align*}
The last equality follows from
\[
 X^3=x^3+cx+\frac{c^2}{x}+\frac{c^3}{x^3},
 \qquad
 AX^2=Ax^2+\frac{Ac^2}{x^2}.
\]
Putting \(Y=Y_0+d\) cancels the constant \(c=d^2\) and gives
\[
                  Y^2+XY=X^3+AX^2+dX
                            =X^3+AX^2+e^2X.
\]
This is \(W_{A,e}\).

\smallskip
\noindent\emph{Step 3: degree and separability.}
The equation \(x^2+Xx+d^2=0\) has derivative \(X\), so it is generically
separable.  The nontrivial translation by \(T\) fixes \(X,Y_0\), while
the displayed quadratic recovers \(x\) and then
\[
                    y=\frac{xY_0+d^2}{X}.
\]
Thus the induced function-field extension has degree exactly two.  The
source abscissa has order \(-2\) at \(O\) and order \(2\) at
\(T=(0,0)\): near \(T\), the nonzero linear term \(d^2x\) in the curve
equation gives \(x=d^{-2}y^2+\) higher-order terms.  Consequently
\(X=x+d^2/x\) has a pole at both \(O\) and \(T\), so both points map to
the target identity.  Since the map has degree two, these two points,
with their fiber multiplicities, exhaust that fiber.  Hence the kernel
is exactly \(\{O,T\}\).

\smallskip
\noindent\emph{Step 4: verify the dual and both compositions.}
If \((X,Y)\in W_{A,e}\), then
\begin{align*}
 &(Y^2+AX^2)^2+X^2(Y^2+AX^2)\\
 &\quad=(Y^2+XY)^2+(A^2+A)X^4\\
 &\quad=X^6+AX^4+d^2X^2,
\end{align*}
which proves \eqref{eq:T-char2-two-isogeny-dual}.

For a direct composition check on \(W_{A,d}\), set
\(t=y/x\) and \(q=d^2/x\).  Then
\[
 x+q=t^2+t+A,\qquad
 \lambda=\frac{y+x^2+d^2}{x}=t^2+A.
\]
The tangent formula yields
\[
 x([2]P)=\lambda^2+\lambda+A
         =t^4+t^2+A^2
         =\left(x+\frac{d^2}{x}\right)^2.
\]
The intercept is again \(\nu=x^2+d^2\).  The ordinate from
\eqref{eq:T-char2-two-isogeny-W} is
\[
                  Y=t^3+t^2+At+q+d.
\]
Using \(x^2+q^2=t^4+t^2+A^2\), direct expansion gives
\[
 y([2]P)
  =(t^2+A+1)x([2]P)+x^2+d^2
  =Y^2+A\,x([2]P).
\]
This is precisely the ordinate produced by
\(\widehat\Phi_{A,d}\circ\Phi_{A,d}\), so this composition is \([2]\)
on \(W_{A,d}\).  The nonconstant isogeny \(\Phi_{A,d}\) is surjective.
Therefore
\[
 (\Phi_{A,d}\circ\widehat\Phi_{A,d})\circ\Phi_{A,d}
 =\Phi_{A,d}\circ[2]
 =[2]\circ\Phi_{A,d}.
\]
Surjectivity permits cancellation on the right, proving
\(\Phi_{A,d}\circ\widehat\Phi_{A,d}=[2]\) on \(W_{A,e}\).
\end{proof}

\begin{corollary}[Native one-sided formula and two optimal schedules]
\label{cor:T-char2-two-isogeny-native}
On
\[
 \mathcal T_{a,d}:\qquad
 F(v^2+v)=d,\qquad F=u^2+u+a,
\]
the isogeny of Theorem~\ref{thm:T-char2-two-isogeny} is
\begin{equation}
 \boxed{\qquad
 (u,v)\longmapsto
 \left(u+v+1,\ \frac{F}{F+e}\right)
 =
 \left(u+v+1,\ \frac{e}{e+v(v+1)}\right)
 \in\mathcal T_{a,e}.
 \qquad}
\label{eq:T-char2-two-isogeny-native}
\end{equation}
Its dual in native coordinates is
\begin{equation}
             (u',v')\longmapsto({u'}^2+a,{v'}^2).
\label{eq:T-char2-two-isogeny-native-dual}
\end{equation}
For raw affine input and projective output in the two
\(\mathbb P^1\)-factors, additions being free, there are two exact
schedules:
\begin{align}
 F&=u^2+u+a,&
 (v'_0:v'_1)&=(F:F+e)
 &&\text{cost }\boxed{\Sqr},\label{eq:T-char2-two-isogeny-S-schedule}\\
 h&=v(v+1),&
 (v'_0:v'_1)&=(e:e+h)
 &&\text{cost }\boxed{\M}.
\label{eq:T-char2-two-isogeny-M-schedule}
\end{align}
In both rows \(u'=u+v+1\).  Thus the implementation may choose the
cheaper of one multiplication and one squaring for its field basis.
\end{corollary}

\begin{proof}
The dictionary of Theorem~\ref{thm:T-char2-Weierstrass} gives
\[
 x=Fv+d=\frac{dv}{v+1},\qquad y=ux.
\]
Put \(z=x/d=v/(v+1)\).  Then
\[
 X=d(z+z^{-1})=\frac{d}{v(v+1)}=F,
\]
and the target inverse dictionary gives \(v'=X/(X+e)=F/(F+e)\).
Furthermore
\[
 \frac YX
 =\frac{zu+(u+1)z^{-1}+1}{z+z^{-1}}
 =u+v+1,
\]
as is checked after multiplying by \(v(v+1)\).  Since
\(Fv(v+1)=d=e^2\), cross multiplication gives
\[
                  \frac F{F+e}=\frac e{e+v(v+1)}.
\]
This proves the two schedules and their costs.

For the dual, square the target Weierstrass coordinates and apply the
\(A=a\) shear in \eqref{eq:T-char2-two-isogeny-dual}.  The inverse
dictionary then gives \(u''={u'}^2+a\) and \(v''={v'}^2\).
\end{proof}

\begin{corollary}[Full-product Jacobian closure under the two-isogeny]
\label{cor:Cabd-char2-two-isogeny-Jacobian}
For the smooth full product family,
\[
 \operatorname{Jac}(\mathcal C_{a,b,d})
   =W_{a+b,d}
 \quad\longrightarrow\quad
 W_{a+b,e}
   =\operatorname{Jac}(\mathcal C_{a,b,e}).
\label{eq:Cabd-char2-two-isogeny-Jacobian}
\]
Thus both displayed extension models retain their written parameter shape:
\(a\), respectively \(a+b\), is unchanged and only
\(d\mapsto\sqrt d\) is performed.  If both
\(\mathcal C_{a,b,d}(k)\) and \(\mathcal C_{a,b,e}(k)\) are nonempty,
choose a rational origin on each curve.  Composition with the two
pointed curve--Jacobian isomorphisms then gives the corresponding
degree-two map of pointed genus-one curves.  Over a finite field the two
nonemptiness hypotheses hold automatically.
\end{corollary}

\begin{proof}
Theorem~\ref{thm:Cabd-binary-Jacobian} identifies the source with
\(W_{a+b,d}\).  Apply Theorem~\ref{thm:T-char2-two-isogeny} with
\(A=a+b\).  Its target is \(W_{a+b,e}\), which the same Jacobian theorem
identifies with \(\operatorname{Jac}(\mathcal C_{a,b,e})\).  A selected
rational origin gives an Abel--Jacobi isomorphism for each of the source
and target curves.  Composing the source isomorphism, the displayed
Jacobian isogeny, and the inverse target isomorphism proves the pointed
statement.  The finite-field assertion follows from the existence of a
rational point on every genus-one curve over a finite field.
\end{proof}

Theorem~\ref{thm:T-char2-Weierstrass} separates the two arithmetic
roles of the parameters with unusual clarity:
\begin{equation}
  \begin{aligned}
    d&\quad\hbox{determines the geometric invariant }j=d^{-4},\\
    a&\quad\hbox{determines the Artin--Schreier twist class.}
  \end{aligned}
  \label{eq:char2-parameter-separation}
\end{equation}
Indeed, the change $Y\mapsto Y+tX$ replaces $a$ by
\[
  a+t^2+t.
\]
Thus the twist parameter naturally belongs to the quotient
$k/\wpmap(k)$.  Over $\F_{2^m}$ this quotient has two elements,
distinguished by the absolute trace.

\section{Universality for ordinary binary elliptic curves}
\begin{theorem}[Universality for ordinary elliptic curves]
\label{thm:T-char2-universal}
Let $k$ be a perfect field of characteristic $2$.  Every ordinary
elliptic curve over $k$ is birational over $k$ to some $\T_{a,d}$.
In particular, this holds for every ordinary elliptic curve over a finite
binary field.
\end{theorem}

\begin{proof}
Start with a generalized Weierstrass equation
\[
 y_0^2+a_1x_0y_0+a_3y_0
 =x_0^3+a_2x_0^2+a_4x_0+a_6.
\]
Ordinarity in characteristic two is equivalent to \(a_1\ne0\).
The scaling \(x_0=a_1^2x_1\), \(y_0=a_1^3y_1\), followed by division
by \(a_1^6\), normalizes the coefficient of \(x_1y_1\) to one.  Write
\[
 c=\frac{a_3}{a_1^3},\qquad
 A_2=\frac{a_2}{a_1^2},\qquad
 A_4=\frac{a_4}{a_1^4},\qquad
 A_6=\frac{a_6}{a_1^6}.
\]
The normalized equation is
\[
 y_1^2+x_1y_1+cy_1
 =x_1^3+A_2x_1^2+A_4x_1+A_6.
\]
Set \(x_1=x_2+c\).  The two occurrences of \(cy_1\) cancel, and
\[
 y_1^2+x_2y_1
 =x_2^3+(A_2+c)x_2^2+(A_4+c^2)x_2+B_0,
\]
where \(B_0=c^3+A_2c^2+A_4c+A_6\).  Now put
\(y_1=y_2+t\) with \(t=A_4+c^2\).  Since
\[
 (y_2+t)^2+x_2(y_2+t)
 =y_2^2+x_2y_2+t^2+tx_2,
\]
the linear term cancels and the curve takes the form
\begin{equation}
  y_2^2+x_2y_2=x_2^3+Ax_2^2+B,
  \qquad B\ne0.
  \label{eq:ordinary-char2-normal}
\end{equation}
Here \(A=A_2+c\) and \(B=B_0+t^2\).  The condition \(B\ne0\) follows
from smoothness: if \(B=0\), then \((x_2,y_2)=(0,0)\) is singular,
whereas for this normal form the Weierstrass discriminant equals \(B\).

Because $k$ is perfect, the Frobenius is bijective, so there is a unique
$d\in k^\times$ satisfying $d^4=B$.  Put $Y=y_2+d^2$.  Since
\[
  (Y+d^2)^2+x_2(Y+d^2)
  =Y^2+x_2Y+d^4+d^2x_2,
\]
equation \eqref{eq:ordinary-char2-normal} becomes
\[
  Y^2+x_2Y=x_2^3+Ax_2^2+d^2x_2.
\]
This is \eqref{eq:T-char2-Weierstrass} with $a=A$, and
Theorem~\ref{thm:T-char2-Weierstrass} completes the proof.
\end{proof}

The preceding theorem is consistent with the complete coverage of
ordinary binary curves by binary Edwards and twisted $\mu_4$ models
\cite{BinaryEdwards2008,Kohel2012,KohelTwistedMu4}.  Here the coverage is
expressed directly in the $C$-curve product geometry.  In particular,
the birational dictionary \eqref{eq:T-char2-forward}--\eqref{eq:T-char2-inverse}
places every $\T_{a,d}$ in the same isomorphism class as models for which
highly optimized characteristic-two arithmetic is already available.
For orientation, binary Edwards differential addition together with
doubling can be organized in $5\M+4\Sqr+2\Dpar$ when the known
difference is affine, while the twisted $\mu_4$ framework gives general
addition in $9\M+2\Sqr$, doubling in
$2\M+5\Sqr+2\Dpar$, and a Kummer ladder with point recovery whose
per-bit cost is $4\M+4\Sqr+1\mBase+2\mCurve$
\cite{BinaryEdwards2008,KohelTwistedMu4}.  Through the explicit Weierstrass
and normal-form dictionaries, these established circuits become available
to the \(\mathcal T_{a,d}\) isomorphism class.  The intrinsic advantage of
\eqref{eq:Tad} is that this
coverage coexists with the product equation, its two Artin--Schreier
reflections, and the clean parameter separation
\eqref{eq:char2-parameter-separation}.

\section{The full product family and finite-field coverage}
We can now identify the Jacobian of the full three-parameter family in
characteristic $2$.

\begin{theorem}[Characteristic-two collapse of the two twist parameters]
\label{thm:char2-general-Jacobian}
Let $k$ have characteristic $2$ and let $d\ne0$.  Then
\begin{equation}
  \Jac(\C_{a,b,d})\cong\T_{a+b,d}.
  \label{eq:char2-Jac-collapse}
\end{equation}
More precisely, the equality of $a+b$ is an equality in the
Artin--Schreier quotient $k/\wpmap(k)$.  If
$\C_{a,b,d}(k)\ne\varnothing$, then the curve itself is $k$-isomorphic,
after a choice of origin, to $\T_{a+b,d}$.
\end{theorem}

\begin{proof}
Let $k_s$ be a separable closure.  Choose
$\alpha,\beta\in k_s$ such that
\[
  \alpha^2+\alpha=a,
  \qquad
  \beta^2+\beta=b.
\]
Then
\[
  U=u+\alpha,
  \qquad
  V=v+\beta
\]
identifies $\C_{a,b,d}$ over $k_s$ with the split curve
\[
  C_d:\qquad (U^2+U)(V^2+V)=d.
\]
This split completion is smooth when $d\ne0$.  On the affine chart its
partial derivatives are $V^2+V$ and $U^2+U$, and neither can vanish on
the curve because their product is $d$.  On the boundary, for example
at $Z=0$ in
\[
  (U^2+UZ)(V^2+VW)=dZ^2W^2,
\]
one has $U\ne0$, $W\ne0$, and
$\partial F/\partial V=U^2W\ne0$; the other boundary is symmetric and
the corner $Z=W=0$ is absent.  Smoothness descends from $k_s$ to $k$,
so the Jacobian in the theorem is well defined.
Let
\[
  \sigma_U(U,V)=(U+1,V),
  \qquad
  \sigma_V(U,V)=(U,V+1).
\]
For $g\in\operatorname{Gal}(k_s/k)$, put
\[
  \epsilon_a(g)=g(\alpha)+\alpha\in\F_2,
  \qquad
  \epsilon_b(g)=g(\beta)+\beta\in\F_2.
\]
The descent cocycle of $\C_{a,b,d}$ relative to $C_d$ is
\begin{equation}
  g\longmapsto
  \sigma_U^{\epsilon_a(g)}\sigma_V^{\epsilon_b(g)}.
  \label{eq:char2-descent-cocycle}
\end{equation}

Each of $\sigma_U$ and $\sigma_V$ is a nontrivial involution with fixed
points on the smooth completion.  On a genus-one curve, such an
involution has the form $P\mapsto Q-P$ after a geometric origin is
chosen.  Its induced action on $\Pic^0$ is therefore $[-1]$; the
translation by $Q$ acts trivially on $\Pic^0$.  Consequently, the
cocycle induced by \eqref{eq:char2-descent-cocycle} on the Jacobian is
\[
  g\longmapsto[-1]^{\epsilon_a(g)+\epsilon_b(g)}.
\]
But $\alpha+\beta$ satisfies
\[
  (\alpha+\beta)^2+(\alpha+\beta)=a+b,
\]
so $\epsilon_a+\epsilon_b$ is precisely the Artin--Schreier character
of $a+b$.  The curve $\T_{a+b,d}$ is the corresponding twist of $C_d$.
This proves \eqref{eq:char2-Jac-collapse}.  The final assertion again
uses the fact that a genus-one curve with a rational point is isomorphic
to its Jacobian after choosing that point as the origin.
\end{proof}

\begin{corollary}[Finite-field coverage]
\label{cor:finite-field-coverage}
Let $k$ be finite.

\begin{enumerate}[label=(\alph*)]
  \item In odd characteristic, every smooth $\C_{a,b,d}$ is
  $k$-isomorphic to the member $\T_{a_0,d_0}$ of
  Theorem~\ref{thm:odd-Jacobian-coverage}.

  \item In characteristic $2$, every smooth $\C_{a,b,d}$ is
  $k$-isomorphic to $\T_{a+b,d}$.
\end{enumerate}
\end{corollary}

\begin{proof}
Every smooth genus-one curve over a finite field has a rational point;
this follows, for example, from Lang's theorem.  Apply
Theorems~\ref{thm:odd-Jacobian-coverage} and
\ref{thm:char2-general-Jacobian}.
\end{proof}

Thus the one-sided family retains all elliptic isomorphism classes that
occur in the full product family over finite fields, while using only
one explicit twist parameter.  In odd characteristic it is exactly a
$C$-curve incarnation of twisted Edwards/Montgomery geometry; in
characteristic $2$ it is a characteristic-native normal form for all
ordinary elliptic curves, with $d$ and $a$ cleanly separating the
$j$-invariant and the Artin--Schreier twist.

\chapter{Reciprocal \texorpdfstring{$C$}{C}-Curves}
\label{ch:reciprocal-C-curves}

The notation \(\mathcal R_{\tau,\sigma,\kappa}\) in this chapter denotes a
three-parameter reciprocal \(C\)-curve.  It must not be confused with the
one-parameter reciprocal chart
\(\mathcal R_d:(u+1)(v+1)=du^2v^2\) of
Chapter~\ref{ch:reciprocal-chart}, which is an internal chart of
\(\mathcal C_d\).  The precise relation between the two reciprocal
constructions is established below.

\section{Definition, compactification, and quotient involutions}

Let $k$ be a field and let $\tau,\sigma,\kappa\in k$.  The reciprocal
$C$-curve is the biquadratic curve
\begin{equation}
  \R_{\tau,\sigma,\kappa}:\qquad
  (x^2-\tau)(y^2-\sigma)=\kappa xy.
  \label{eq:reciprocal-general}
\end{equation}
Its natural completion in $\PP^1\times\PP^1$ is the zero locus of
\begin{equation}
\begin{split}
  F(X,Z;Y,W)
  ={}&X^2Y^2-\sigma X^2W^2-\tau Z^2Y^2\\
     &{}-\kappa XZYW+\tau\sigma Z^2W^2.
  \label{eq:R-bihomogeneous}
\end{split}
\end{equation}
The affine equation has a group law only after a rational point has been
selected.  When no such point has been selected, or when no rational
point exists, $\R_{\tau,\sigma,\kappa}$ is to be regarded as a genus-one
torsor under the Jacobian described below.

On the dense torus $xy\ne0$, division by $xy$ rewrites
\eqref{eq:reciprocal-general} as
\begin{equation}
  \left(x-\frac{\tau}{x}\right)
  \left(y-\frac{\sigma}{y}\right)=\kappa.
  \label{eq:reciprocal-quotient-product}
\end{equation}
Thus the equation is the product of two degree-two quotient functions,
in exact analogy with the product
$(u^2+u)(v^2+v)$ on $C_d$.

Assume first that $\charac(k)\ne2$.  The two deck transformations are
\begin{equation}
  \iota_x(x,y)=\left(-\frac{\tau}{x},y\right),
  \qquad
  \iota_y(x,y)=\left(x,-\frac{\sigma}{y}\right).
  \label{eq:reciprocal-involutions}
\end{equation}
They extend to the projective linear transformations
\begin{equation}
  (X:Z)\longmapsto(-\tau Z:X),
  \qquad
  (Y:W)\longmapsto(-\sigma W:Y).
  \label{eq:reciprocal-PGL2}
\end{equation}
The matrices representing these transformations square to scalar
matrices, so both transformations have order two in
$\operatorname{PGL}_2(k)$.  Their fixed divisors are defined by
$x^2=-\tau$ and $y^2=-\sigma$.  Consequently, the individual
involutions can be nonsplit over $k$, even though the involutions
themselves are defined over $k$.

\section{Exact smoothness criterion in odd characteristic}

Regard \eqref{eq:reciprocal-general} as a quadratic equation in $y$ and
put
\begin{equation}
  W_R=2(x^2-\tau)y-\kappa x.
  \label{eq:R-W-definition}
\end{equation}
Multiplying the quadratic equation
\[
  (x^2-\tau)y^2-\kappa xy-\sigma(x^2-\tau)=0
\]
by $4(x^2-\tau)$ and completing the square gives
\begin{equation}
  W_R^2
  =4\sigma x^4+(\kappa^2-8\tau\sigma)x^2
   +4\sigma\tau^2.
  \label{eq:R-even-quartic}
\end{equation}
The inverse on the open set $x^2\ne\tau$ is
\begin{equation}
  y=\frac{W_R+\kappa x}{2(x^2-\tau)}.
  \label{eq:R-even-quartic-inverse}
\end{equation}

\begin{theorem}[Smoothness of the reciprocal family]
\label{thm:R-smoothness-odd}
Suppose $\charac(k)\ne2$.  The completion of
$\R_{\tau,\sigma,\kappa}$ in $\PP^1\times\PP^1$ is a smooth genus-one
curve if and only if
\begin{equation}
  \tau\sigma\kappa(\kappa^2-16\tau\sigma)\ne0.
  \label{eq:R-smoothness}
\end{equation}
\end{theorem}

\begin{proof}
The right-hand side of \eqref{eq:R-even-quartic} is an even quartic
$Ax^4+Bx^2+C$ with
\[
  A=4\sigma,
  \qquad B=\kappa^2-8\tau\sigma,
  \qquad C=4\sigma\tau^2.
\]
Its constant and leading coefficients are nonzero exactly when
$\tau\sigma\ne0$.  Moreover,
\begin{align*}
  B^2-4AC
  &=(\kappa^2-8\tau\sigma)^2
    -4(4\sigma)(4\sigma\tau^2)\\
  &=\kappa^4-16\kappa^2\tau\sigma
    +64\tau^2\sigma^2-64\tau^2\sigma^2\\
  &=\kappa^2(\kappa^2-16\tau\sigma).
\end{align*}
Hence the quartic has four distinct geometric roots exactly under
\eqref{eq:R-smoothness}.  In particular, it is not a square in
$\overline{k}(x)$, so both the quartic double cover and the reciprocal
biquadratic are geometrically integral.  In that case its smooth
projective completion
is a double cover of $\PP^1_x$ branched at four distinct points and has
genus one.  Equations \eqref{eq:R-W-definition} and
\eqref{eq:R-even-quartic-inverse} identify its function field with that
of the reciprocal biquadratic.  The projective biquadratic has
arithmetic genus $(2-1)(2-1)=1$.  Its normalization already has genus
one, so the sum of its local $\delta$-invariants is
$p_a-g=1-1=0$.  Every local $\delta$-invariant is nonnegative and
vanishes exactly at a nonsingular point; consequently, the
biquadratic completion is itself smooth and is isomorphic to the
quartic completion.

It remains to account for the excluded parameters.  If $\kappa=0$, the
bihomogeneous equation factors geometrically into the union of the two
vertical and two horizontal branch fibres.  If $\tau=0$ or $\sigma=0$,
the affine point $(0,0)$ is singular.  Finally, if
$\kappa^2=16\tau\sigma$, the even quartic has a repeated root and its
normalization has genus zero.  Thus every excluded parameter produces a
singular biquadratic completion, proving necessity as well as
sufficiency.
\end{proof}

The commuting product
\begin{equation}
  \epsilon_R=\iota_x\iota_y:
  (x,y)\longmapsto
  \left(-\frac{\tau}{x},-\frac{\sigma}{y}\right)
  \label{eq:R-product-involution}
\end{equation}
is fixed-point free under the hypotheses of
Theorem~\ref{thm:R-smoothness-odd}.  Indeed, a fixed point would satisfy
$x^2=-\tau$ and $y^2=-\sigma$.  Substitution in
\eqref{eq:reciprocal-general} gives
$16\tau\sigma=\kappa^2$.  After a geometric origin is chosen, this
fixed-point-free involution is translation on the genus-one curve by a
nonzero two-torsion displacement class in \(\Pic^0\).  The class is
\(k\)-rational because \(\epsilon_R\) is defined over \(k\), whereas the
induced action of the translation on \(\Pic^0\) itself is the identity.

\section[The order-four action]{An order-four action and the four-torsion interpretation}

Suppose $\tau/\sigma$ is a square in $k$ and choose $c\in k^\times$
with
\begin{equation}
  c^2=\frac{\tau}{\sigma}.
  \label{eq:R-c-scaling}
\end{equation}
The scaled interchange
\[
  S_c(x,y)=\left(cy,\frac{x}{c}\right)
\]
preserves \eqref{eq:reciprocal-general}.  Indeed,
\[
  (c^2y^2-\tau)\left(\frac{x^2}{c^2}-\sigma\right)
  =c^2(y^2-\sigma)\frac{x^2-\tau}{c^2}
  =(x^2-\tau)(y^2-\sigma),
\]
where $c^2\sigma=\tau$, while $(cy)(x/c)=xy$.  Define
\begin{equation}
  \rho_R=S_c\circ\iota_x,
  \qquad
  \rho_R(x,y)=\left(cy,-\frac{\tau}{cx}\right).
  \label{eq:R-order-four-map}
\end{equation}
Applying this map twice gives
\begin{align*}
  \rho_R^2(x,y)
  &=\rho_R\left(cy,-\frac{\tau}{cx}\right)\\
  &=\left(-\frac{\tau}{x},
          -\frac{\tau}{c^2y}\right)\\
  &=\left(-\frac{\tau}{x},-\frac{\sigma}{y}\right)
   =\epsilon_R(x,y),
\end{align*}
where the equality $\tau/c^2=\sigma$ uses
\eqref{eq:R-c-scaling}.  Consequently,
\begin{equation}
  \rho_R^2=\epsilon_R,
  \qquad \rho_R^4=1.
  \label{eq:R-rho-square}
\end{equation}
A fixed point of $\rho_R$ would satisfy $x=cy$ and
$x^2=-\tau$; substitution again forces
$\kappa^2=16\tau\sigma$.  Hence $\rho_R$ is fixed-point free on every
smooth member.  After a rational origin is selected, $\rho_R$ is
translation by a point of exact order four.

For the normalized split member
\begin{equation}
  \R_\kappa:\qquad
  (x^2-1)(y^2-1)=\kappa xy,
  \qquad \kappa(\kappa^2-16)\ne0,
  \label{eq:R-normalized}
\end{equation}
one may select $(1,0)$ as origin and take
\[
  \rho_R(x,y)=\left(y,-\frac1x\right).
\]
On the smooth projective completion this gives the four-cycle
\begin{equation}
  (1,0)\longmapsto(0,-1)\longmapsto(-1,\infty)
  \longmapsto(\infty,1)\longmapsto(1,0).
  \label{eq:R-four-cycle}
\end{equation}
Thus $(0,-1)$ has order four for this choice of origin.  A different
origin will translate this marked four-torsion configuration; this
point must therefore not be confused with the identity selected below
by the Jacobi-quartic arithmetic model.

If $\tau/\sigma$ is nonsquare, the order-four generator becomes visible
after adjoining $c$, whereas its square $\epsilon_R$ is already defined
over $k$.  This is the precise descent-theoretic sense in which the
model is a non-split four-torsion normal form.  The two reflections and
the scaled interchange generate a dihedral group of order eight after
the splitting extension.

\section{Jacobi quartic, Weierstrass Jacobian, and invariant}

The even quartic \eqref{eq:R-even-quartic} and
Lemma~\ref{lem:even-quartic-jacobian} give the Jacobian over the ground
field without any square-root hypothesis.  Put
\begin{equation}
  q_R=\kappa^2-8\tau\sigma,
  \qquad
  \Omega_R=\kappa^2(\kappa^2-16\tau\sigma).
  \label{eq:R-q-Omega}
\end{equation}
Then
\begin{equation}
  E_R:\qquad
  V^2=U\bigl(U^2-2q_RU+\Omega_R\bigr)
  \label{eq:R-Jacobian}
\end{equation}
is $\Jac(\R_{\tau,\sigma,\kappa})$.  In particular, $E_R$ has the
rational two-torsion point $(0,0)$.  If
$\R_{\tau,\sigma,\kappa}(k)$ is nonempty and an origin is chosen, the
genus-one curve itself becomes $k$-isomorphic to $E_R$.

Over an extension containing square roots of $\tau$ and $\sigma$, write
\begin{equation}
  x=\sqrt{\tau}\,X,
  \qquad y=\sqrt{\sigma}\,Y,
  \qquad
  \kappa_0=\frac{\kappa}{\sqrt{\tau\sigma}}.
  \label{eq:R-normalization-scaling}
\end{equation}
Then the curve becomes
$(X^2-1)(Y^2-1)=\kappa_0XY$.  With
\[
  Z=(X^2-1)Y-\frac{\kappa_0}{2}X
\]
one obtains the Jacobi quartic
\begin{equation}
  \JQ_{\kappa_0}:\qquad
  Z^2=X^4+2a_JX^2+1,
  \qquad
  a_J=\frac{\kappa_0^2}{8}-1.
  \label{eq:R-Jacobi-quartic}
\end{equation}
The inverse on $X^2\ne1$ is
\begin{equation}
  Y=\frac{Z+(\kappa_0/2)X}{X^2-1}.
  \label{eq:R-Jacobi-inverse}
\end{equation}
The identity $(0,1)$ of the Jacobi model corresponds to the reciprocal
point $(0,-1)$.  Hence the Jacobi arithmetic below uses $(0,-1)$ as
identity, rather than the origin $(1,0)$ used in
\eqref{eq:R-four-cycle}.

Set
\begin{equation}
  t=\frac{\kappa^2}{\tau\sigma}=\kappa_0^2.
  \label{eq:R-t-parameter}
\end{equation}
For the even quartic, the binary-quartic invariant
$I=B^2+12AC$ equals
\[
  (\tau\sigma)^2(t^2-16t+256),
\]
up to the common scaling that cancels in $j$, while
$B^2-4AC=(\tau\sigma)^2t(t-16)$.  Substitution in the Jacobian
invariant therefore yields
\begin{equation}
  j(\R_{\tau,\sigma,\kappa})
  =\frac{(t^2-16t+256)^3}{t^2(t-16)^2}.
  \label{eq:R-j-invariant}
\end{equation}
The poles $t=0$ and $t=16$ coincide with the singular values in
Theorem~\ref{thm:R-smoothness-odd}.

\section{The reciprocal model as a twist of \texorpdfstring{$C_d$}{Cd}}
\label{subsec:R-explicit-twisted-Cd}

Choose square roots
\begin{equation}
  \alpha_0^2=-\tau,
  \qquad
  \beta_0^2=-\sigma,
  \qquad
  \varrho=\alpha_0\beta_0,
  \qquad
  \varrho^2=\tau\sigma.
  \label{eq:R-C-square-roots}
\end{equation}
Apply the Cayley transformations
\begin{equation}
  x=\alpha_0\frac{1+r}{1-r},
  \qquad
  y=\beta_0\frac{1+s}{1-s},
  \label{eq:R-C-Cayley}
\end{equation}
whose inverse formulas are
\begin{equation}
  r=\frac{x-\alpha_0}{x+\alpha_0},
  \qquad
  s=\frac{y-\beta_0}{y+\beta_0}.
  \label{eq:R-C-Cayley-inverse}
\end{equation}
The reciprocal involutions become $r\mapsto-r$ and $s\mapsto-s$.
Indeed,
\[
  \alpha_0\frac{1-r}{1+r}
  =\frac{\alpha_0^2}{x}=-\frac{\tau}{x},
\]
and the second coordinate is identical.

The factors in \eqref{eq:reciprocal-general} transform as
\[
  x^2-\tau
  =\frac{2\alpha_0^2(1+r^2)}{(1-r)^2},
  \qquad
  y^2-\sigma
  =\frac{2\beta_0^2(1+s^2)}{(1-s)^2},
\]
whereas
\[
  xy=\varrho\frac{(1+r)(1+s)}{(1-r)(1-s)}.
\]
After multiplication by $(1-r)^2(1-s)^2/\varrho$, the transformed
equation is
\begin{equation}
  4\varrho(1+r^2)(1+s^2)
  =\kappa(1-r^2)(1-s^2).
  \label{eq:R-C-after-Cayley}
\end{equation}
Put
\begin{equation}
  \eta=\frac{\kappa}{4\varrho},
  \qquad
  A_C=\frac{1+\eta}{1-\eta}
      =\frac{4\varrho+\kappa}{4\varrho-\kappa}.
  \label{eq:R-C-eta}
\end{equation}
Then \eqref{eq:R-C-after-Cayley} is equivalent to
\begin{equation}
  r^2s^2+A_C(r^2+s^2)+1=0.
  \label{eq:R-C-intermediate-QRT}
\end{equation}
The smoothness condition
$\kappa^2\ne16\tau\sigma=16\varrho^2$ implies
$4\varrho\pm\kappa\ne0$; hence $A_C$ is defined, nonzero, and different
from $\pm1$.

\begin{theorem}[Explicit identification with a twisted \texorpdfstring{$C_d$}{Cd} curve]
\label{thm:R-explicit-twisted-Cd}
Assume $\charac(k)\ne2$ and
\(\tau\sigma\kappa(\kappa^2-16\tau\sigma)\ne0\).
Let $L/k$ contain the elements $\alpha_0,\beta_0$ of
\eqref{eq:R-C-square-roots} and an element $c_C\in L^\times$ satisfying
\begin{equation}
  c_C^2=-A_C
  =-\frac{4\varrho+\kappa}{4\varrho-\kappa}.
  \label{eq:R-twisted-Cd-c}
\end{equation}
Define
\begin{align}
  a_{\mathrm{tw}}
  &=\frac{1+A_C}{4}
    =\frac{2\varrho}{4\varrho-\kappa},
  \label{eq:R-twisted-Cd-a}\\
  d_{\mathrm{tw}}
  &=\frac{A_C^{-1}-A_C}{16}
    =\frac{\kappa\varrho}{\kappa^2-16\tau\sigma}.
  \label{eq:R-twisted-Cd-d}
\end{align}
Then the smooth projective completions are $L$-isomorphic:
\begin{equation}
  \boxed{
  \bigl(\R_{\tau,\sigma,\kappa}\bigr)_L
  \simeq
  \bigl(\T_{a_{\mathrm{tw}},d_{\mathrm{tw}}}\bigr)_L:
  (u^2+u+a_{\mathrm{tw}})(v^2+v)=d_{\mathrm{tw}}.}
  \label{eq:R-to-explicit-twisted-Cd}
\end{equation}
On the dense affine chart on which the displayed denominators are nonzero,
an isomorphism is
\begin{equation}
\begin{split}
  u&=\frac12\left(\frac{x-\alpha_0}{x+\alpha_0}-1\right)
     =-\frac{\alpha_0}{x+\alpha_0},\\
  v&=\frac12\left(
       \frac{y-\beta_0}{c_C(y+\beta_0)}-1
     \right).
\end{split}
\label{eq:R-to-twisted-Cd-map}
\end{equation}
The inverse map is
\begin{equation}
\begin{split}
  x&=-\alpha_0\frac{u+1}{u},\\
  y&=\beta_0
      \frac{1+c_C(2v+1)}{1-c_C(2v+1)}.
\end{split}
\label{eq:R-to-twisted-Cd-inverse}
\end{equation}
Both coordinate changes are projective linear transformations on the two
\(\PP^1\) factors, so the formulas extend across all omitted affine
points.
\end{theorem}

\begin{proof}
Let
\[
  R=2u+1=r,
  \qquad
  S=2v+1=\frac{s}{c_C}.
\]
Substituting $s=c_CS$ and $c_C^2=-A_C$ into
\eqref{eq:R-C-intermediate-QRT}, then dividing by $-A_C$, gives
\begin{equation}
  R^2S^2-R^2+A_CS^2-A_C^{-1}=0.
  \label{eq:R-twisted-Cd-centered-expanded}
\end{equation}
Set
\[
  \lambda_{\mathrm{tw}}=1-4a_{\mathrm{tw}}=-A_C=c_C^2,
  \qquad
  \delta_{\mathrm{tw}}=16d_{\mathrm{tw}}
  =A_C^{-1}-A_C.
\]
Then
\begin{align*}
  &(R^2-\lambda_{\mathrm{tw}})(S^2-1)
     -\delta_{\mathrm{tw}}\\
  &\qquad=R^2S^2-R^2+A_CS^2-A_C^{-1},
\end{align*}
so \eqref{eq:R-twisted-Cd-centered-expanded} is exactly the centered
form of \(\T_{a_{\mathrm{tw}},d_{\mathrm{tw}}}\).
The rational formulas \eqref{eq:R-to-twisted-Cd-map} are obtained by
combining \eqref{eq:R-C-Cayley-inverse} with
$R=2u+1$ and $S=2v+1$; solving these equations gives
\eqref{eq:R-to-twisted-Cd-inverse}.  The parameter identities follow from
\[
  \frac{1+A_C}{4}
  =\frac14\left(1+\frac{4\varrho+\kappa}
                         {4\varrho-\kappa}\right)
  =\frac{2\varrho}{4\varrho-\kappa}
\]
and
\begin{align*}
  \frac{A_C^{-1}-A_C}{16}
  &=\frac1{16}\left(
      \frac{4\varrho-\kappa}{4\varrho+\kappa}
      -\frac{4\varrho+\kappa}{4\varrho-\kappa}
    \right)\\
  &=\frac{\kappa\varrho}{\kappa^2-16\varrho^2}
   =\frac{\kappa\varrho}{\kappa^2-16\tau\sigma}.
\end{align*}
Finally,
$\lambda_{\mathrm{tw}}\delta_{\mathrm{tw}}
(\lambda_{\mathrm{tw}}-\delta_{\mathrm{tw}})\ne0$ because
$A_C\ne0,\pm1$ and
$\lambda_{\mathrm{tw}}-\delta_{\mathrm{tw}}=-A_C^{-1}$.
Thus the target twisted \(C_d\) member is smooth, and the birational
maps extend uniquely to an isomorphism of the smooth projective
completions.
\end{proof}

\begin{remark}[Splitting data]
\label{rem:R-twisted-Cd-splitting}
The isomorphism in Theorem~\ref{thm:R-explicit-twisted-Cd} is defined over
$k$ precisely when the chosen splitting data
$\alpha_0,\beta_0,c_C$ lie in $k$.  Changing the signs of the square
roots changes the written parameters and the displayed coordinate map,
but not the resulting $L$-isomorphism class, because every such target is
isomorphic to the same reciprocal curve over $L$.  Over the ground field,
the failure of these square roots to descend is exactly the reciprocal
split/nonsplit twisting data described in the preceding sections.
\end{remark}

A second, symmetric scaling identifies the same reciprocal curve with an
untwisted member of the original family.  Put
\[
  r=c_CR,
  \qquad
  s=c_CS.
\]
Division of \eqref{eq:R-C-intermediate-QRT} by $A_C^2$ gives
\begin{equation}
  (R^2-1)(S^2-1)=\Delta_C,
  \qquad
  \Delta_C=1-\frac1{A_C^2}
          =\frac{4\eta}{(1+\eta)^2}.
  \label{eq:R-to-C-centered}
\end{equation}
Thus, with
\begin{equation}
  U=\frac{R-1}{2},
  \qquad
  V=\frac{S-1}{2},
  \label{eq:R-to-Cd-affine-final}
\end{equation}
the target is
\begin{equation}
  C_{d_C}:\qquad (U^2+U)(V^2+V)=d_C,
  \qquad
  d_C=\frac{\Delta_C}{16}
      =\frac{\kappa\varrho}{(4\varrho+\kappa)^2}.
  \label{eq:R-to-C-d-parameter}
\end{equation}
Equations \eqref{eq:R-C-Cayley-inverse},
\eqref{eq:R-twisted-Cd-c}, and \eqref{eq:R-to-Cd-affine-final} therefore
give a fully explicit coordinate map to
$C_{d_C}=\T_{0,d_C}$ as well.  The two targets
$\T_{a_{\mathrm{tw}},d_{\mathrm{tw}}}$ and
$\T_{0,d_C}$ are different written models of the same curve over $L$.

Thus the reciprocal family is a reciprocal twist of $C_d$ in a precise
sense: after the splitting extensions and the one-sided diagonal
normalization by $c_C$, it is the explicitly identified twisted curve
$\T_{a_{\mathrm{tw}},d_{\mathrm{tw}}}$; applying the same scaling
symmetrically to both centered coordinates yields $C_{d_C}$ itself.  Its
ground-field descent
retains reciprocal rather than additive reflections.  The equality of the
$j$-invariants follows by substituting $t=16\eta^2$ and
$\Delta_C=4\eta/(1+\eta)^2$ into
\begin{equation}
  j(C_{d_C})
  =16\frac{(\Delta_C^2-16\Delta_C+16)^3}
  {(1-\Delta_C)\Delta_C^4}.
  \label{eq:C-j-centered}
\end{equation}

\section{Characteristic two: identification with \texorpdfstring{$C_d$}{Cd}}
\label{subsubsec:R-char2}

The characteristic-two member is not an independent nonsplit
reciprocal twist.  The semisimple split/nonsplit language used above
must be replaced by Artin--Schreier geometry.

\begin{theorem}[Reciprocal curves in characteristic two]
\label{thm:R-char2-Cd}
Let $k$ be a perfect field of characteristic $2$, and suppose
$\tau\sigma\kappa\ne0$.  Let
\[
  \alpha_2^2=\tau,
  \qquad \beta_2^2=\sigma,
  \qquad
  d_R=\frac{\alpha_2\beta_2}{\kappa}
      =\frac{\sqrt{\tau\sigma}}{\kappa}.
\]
Then
\begin{equation}
  \R_{\tau,\sigma,\kappa}
  \cong C_{d_R}:
  (u^2+u)(v^2+v)=d_R.
  \label{eq:R-char2-Cd-isomorphism}
\end{equation}
An isomorphism of the two $\PP^1\times\PP^1$ models is given by
\begin{equation}
  u=\frac{\alpha_2}{x+\alpha_2},
  \qquad
  v=\frac{\beta_2}{y+\beta_2},
  \label{eq:R-char2-Mobius}
\end{equation}
with inverse
\begin{equation}
  x=\alpha_2\frac{u+1}{u},
  \qquad
  y=\beta_2\frac{v+1}{v}.
  \label{eq:R-char2-Mobius-inverse}
\end{equation}
The corresponding Weierstrass model is
\begin{equation}
  E_R^{(2)}:\qquad
  Y^2+XY=X^3+\frac{\tau\sigma}{\kappa^2}X.
  \label{eq:R-char2-Weierstrass}
\end{equation}
Its invariants are
\begin{equation}
  \Delta=\frac{\tau^2\sigma^2}{\kappa^4},
  \qquad
  j=\frac{\kappa^4}{\tau^2\sigma^2}.
  \label{eq:R-char2-invariants}
\end{equation}
Consequently, the biquadratic completion is smooth if and only if
$\tau\sigma\kappa\ne0$.
\end{theorem}

\begin{proof}
In characteristic $2$ the defining equation is
\[
  (x^2+\tau)(y^2+\sigma)=\kappa xy.
\]
Using \eqref{eq:R-char2-Mobius-inverse}, one obtains
\[
  x^2+\tau=\frac{\tau}{u^2},
  \qquad
  y^2+\sigma=\frac{\sigma}{v^2},
\]
and
\[
  xy=\alpha_2\beta_2
  \frac{(u+1)(v+1)}{uv}.
\]
Substitution and multiplication by $u^2v^2$ give
\[
  (\alpha_2\beta_2)^2
  =\kappa\alpha_2\beta_2u(u+1)v(v+1).
\]
All three parameters are nonzero, so division by
$\kappa\alpha_2\beta_2$ yields
\[
  (u^2+u)(v^2+v)=\frac{\alpha_2\beta_2}{\kappa}=d_R.
\]
This proves the isomorphism of the smooth completions because
\eqref{eq:R-char2-Mobius} is a pair of projective linear coordinate
changes.

Apply Theorem~\ref{thm:T-char2-Weierstrass} with $a=0$ and $d=d_R$.
The result is
$Y^2+XY=X^3+d_R^2X$, and
$d_R^2=\tau\sigma/\kappa^2$.  Formula
\eqref{eq:T-char2-invariants} gives
$\Delta=d_R^4$ and $j=d_R^{-4}$, which are exactly
\eqref{eq:R-char2-invariants}.  Thus the completion is smooth whenever
$\tau\sigma\kappa\ne0$.

For the converse, the vanishing of any parameter produces a singular
completion.  If $\tau=0$ or $\sigma=0$, the affine point $(0,0)$ lies on
the curve and both affine partial derivatives
$F_x=\kappa y$ and $F_y=\kappa x$ vanish there.  If $\kappa=0$, then,
over the perfect field $k$,
\[
  (x^2+\tau)(y^2+\sigma)
  =(x+\alpha_2)^2(y+\beta_2)^2,
\]
so the projective curve is nonreduced and cannot be smooth.  This proves
that smoothness is equivalent to $\tau\sigma\kappa\ne0$.
\end{proof}

Under \eqref{eq:R-char2-Mobius}, the reciprocal involutions become
\begin{equation}
  \iota_x:(u,v)\longmapsto(u+1,v),
  \qquad
  \iota_y:(u,v)\longmapsto(u,v+1).
  \label{eq:R-char2-AS-involutions}
\end{equation}
For example,
\[
  \frac{\alpha_2}{\tau/x+\alpha_2}
  =\frac{x}{x+\alpha_2}
  =\frac{\alpha_2}{x+\alpha_2}+1.
\]
The fixed equation $x^2=\tau$ has one geometric root of multiplicity
two, rather than two distinct fixed points.  Thus the reciprocal
involutions are wild Artin--Schreier involutions, not semisimple
elements of nonsplit one-dimensional tori.

The QRT root exchanges remain
\[
  (x,y)\longmapsto\left(\frac{\tau}{x},y\right),
  \qquad
  (x,y)\longmapsto\left(x,\frac{\sigma}{y}\right).
\]
Their product has no fixed point on the smooth curve: a fixed point
would have $x^2=\tau$ and $y^2=\sigma$, making the left side of the
curve equation zero and the right side
$\kappa\alpha_2\beta_2\ne0$.

Put $c_2=\alpha_2/\beta_2$, so that
$c_2^2=\tau/\sigma$.  The scaled interchange
$(x,y)\mapsto(c_2y,x/c_2)$ preserves the characteristic-two reciprocal
equation by the same factor calculation used in odd characteristic,
and $x\mapsto\tau/x$ is a root exchange.  Their composition is
therefore an automorphism.  Explicitly,
\begin{equation}
  \rho_R^{(2)}(x,y)
  =\left(c_2y,\frac{\tau}{c_2x}\right),
  \qquad
  (\rho_R^{(2)})^2(x,y)
  =\left(\frac{\tau}{x},\frac{\sigma}{y}\right).
  \label{eq:R-char2-rho}
\end{equation}
With $(\alpha_2,0)$ as origin, this automorphism gives the rational
four-cycle
\begin{equation}
  (\alpha_2,0)\longmapsto(0,\beta_2)
  \longmapsto(\alpha_2,\infty)
  \longmapsto(\infty,\beta_2)
  \longmapsto(\alpha_2,0).
  \label{eq:R-char2-four-cycle}
\end{equation}
The normalized equation
\begin{equation}
  (X+1)^2(Y+1)^2=\lambda XY,
  \qquad
  \lambda=\frac{\kappa}{\sqrt{\tau\sigma}},
  \label{eq:R-char2-Z4-normal}
\end{equation}
is therefore a biquadratic realization of the characteristic-two
$\mathbb Z/4\mathbb Z$-normal form.  On
\eqref{eq:R-char2-Weierstrass}, writing $d=d_R$, the points
\begin{equation}
  T_2=(0,0),
  \qquad
  T_4=(d,0),
  \qquad
  -T_4=(d,d)
  \label{eq:R-char2-torsion-points}
\end{equation}
have orders $2,4,4$, respectively.  Indeed, after the shift
$\widehat Y=Y+d^2$, the tangent slope at $(d,d^2)$ is zero, so its
double is $(0,d^2)$ in shifted coordinates, namely $(0,0)$ in the
original coordinates.

If $k$ is not perfect, the same reduction holds over the purely
inseparable extension $k(\sqrt\tau,\sqrt\sigma)$.  This is a
purely inseparable failure of splitting and must not be identified with
the quadratic separable nonsplit-torus phenomenon of odd
characteristic.

\section[Arithmetic conventions and the Weierstrass atlas]
{Arithmetic conventions and a complete odd-characteristic Weierstrass atlas}
\label{subsec:odd-Weierstrass-toolkit}

The arithmetic formulas in the remainder of this section use the
following cost notation.  A general field multiplication, a field
squaring, a multiplication by a fixed curve parameter, and a field
inversion are denoted by
\[
  \M,\qquad \Sqr,\qquad \Dpar,\qquad \Inv,
\]
respectively.  Additions, subtractions, negations, and multiplication by
small integers are not included.  The stated costs are straight-line
upper bounds; a specialized implementation may merge operations or
replace a fixed-parameter multiplication by additions when the
parameter has a special representation.  Equality tests, conditional
branches, conditional swaps, and coordinate copies are not included in
the field-operation counts; they must nevertheless be implemented in a
complete constant-time routine.

Throughout this section, assume $\charac(k)\ne2$.
We first record a Weierstrass toolkit that applies both to reciprocal
curves and to the symmetric QRT family.  Let
\begin{equation}
  E_{A,B}:\qquad v^2=u^3+Au^2+Bu,
  \qquad B(A^2-4B)\ne0,
  \label{eq:EAB}
\end{equation}
and let $O$ denote its point at infinity.  The point $(0,0)$ is a
nonzero point of order two.

\begin{proposition}[Affine group law on $E_{A,B}$]
\label{prop:EAB-affine-law}
Let $P=(u_1,v_1)$ and $Q=(u_2,v_2)$ be affine points of
\eqref{eq:EAB}.

\begin{enumerate}[label=(\alph*)]
  \item If $u_1\ne u_2$, put
  \begin{equation}
    \lambda=\frac{v_2-v_1}{u_2-u_1}.
    \label{eq:EAB-affine-lambda-add}
  \end{equation}
  Then $P+Q=(u_3,v_3)$, where
  \begin{equation}
  \begin{split}
    u_3&=\lambda^2-A-u_1-u_2,\\
    v_3&=\lambda(u_1-u_3)-v_1.
  \end{split}
  \label{eq:EAB-affine-add}
  \end{equation}

  \item If $P=Q$ and $v_1\ne0$, put
  \begin{equation}
    \lambda=\frac{3u_1^2+2Au_1+B}{2v_1}.
    \label{eq:EAB-affine-lambda-double}
  \end{equation}
  Then $2P=(u_3,v_3)$, where
  \begin{equation}
  \begin{split}
    u_3&=\lambda^2-A-2u_1,\\
    v_3&=\lambda(u_1-u_3)-v_1.
  \end{split}
  \label{eq:EAB-affine-double}
  \end{equation}

  \item If $u_1=u_2$ and $v_2=-v_1$, then $P+Q=O$.  This branch
  includes the case \(P=Q\) and \(v_1=0\).
\end{enumerate}
Together with $P+O=P$, these cases are mutually exclusive and exhaust
all pairs of points.
\end{proposition}

\begin{proof}
For $u_1\ne u_2$, the line through $P$ and $Q$ is
$v=\lambda(u-u_1)+v_1$.  Substitution in \eqref{eq:EAB} gives a cubic
polynomial in $u$ whose three roots are $u_1,u_2,u_3$.  The coefficient
of $u^2$ is $A-\lambda^2$, so Vieta's formula gives
$u_1+u_2+u_3=\lambda^2-A$.  Reflection in the $u$-axis changes the
third intersection point into $P+Q$, which gives the formula for $v_3$.
For doubling, the derivative of
$v^2-u^3-Au^2-Bu$ at $P$ gives the tangent slope
\eqref{eq:EAB-affine-lambda-double}; the same Vieta argument gives
\eqref{eq:EAB-affine-double}.  If two affine points have the same
$u$-coordinate, the curve equation implies $v_2=\pm v_1$.  The two
possibilities are therefore precisely the inverse and tangent cases.
Finally, a point with $v=0$ is equal to its inverse and hence has order
two.  This proves both the formulas and the exhaustion assertion.
\end{proof}

With one inversion used for each slope, the affine addition schedule in
part~(a) costs
\begin{equation}
  2\M+1\Sqr+1\Inv,
  \label{eq:EAB-affine-add-cost}
\end{equation}
and affine doubling in part~(b) costs
\begin{equation}
  2\M+2\Sqr+1\Dpar+1\Inv.
  \label{eq:EAB-affine-double-cost}
\end{equation}
For addition, one multiplication forms the slope after inversion and
one forms the output ordinate.  For doubling, the fixed multiplication
is $Au_1$; the two squares are $u_1^2$ and $\lambda^2$.

For inversion-free arithmetic, represent $u=X/Z^2$ and $v=Y/Z^3$ by
Jacobian coordinates $(X:Y:Z)$, with $O=(1:1:0)$.  For two inputs with
$H\ne0$, where $H$ is defined below, set
\begin{align}
  Z_1^{(2)}&=Z_1^2,& Z_2^{(2)}&=Z_2^2,\notag\\
  U_1&=X_1Z_2^{(2)},& U_2&=X_2Z_1^{(2)},\notag\\
  S_1&=Y_1Z_2^{(2)}Z_2,&
  S_2&=Y_2Z_1^{(2)}Z_1,\notag\\
  H&=U_2-U_1,& R&=S_2-S_1,\notag\\
  I&=H^2,& J&=HI,& V_0&=U_1I,\notag\\
  Z_3&=HZ_1Z_2,\notag\\
  X_3&=R^2-AZ_3^2-J-2V_0,\notag\\
  Y_3&=R(V_0-X_3)-S_1J.
  \label{eq:EAB-Jacobian-add}
\end{align}
We now verify the dehomogenization.  Put
\[
  u_i=\frac{X_i}{Z_i^2},
  \qquad
  v_i=\frac{Y_i}{Z_i^3},
  \qquad
  h=u_2-u_1,
  \qquad
  r=v_2-v_1.
\]
For finite inputs,
\[
  H=Z_1^2Z_2^2h,
  \qquad
  R=Z_1^3Z_2^3r,
  \qquad
  Z_3=Z_1^3Z_2^3h.
\]
Moreover,
\[
  J=Z_1^6Z_2^6h^3,
  \qquad
  V_0=Z_1^6Z_2^6u_1h^2.
\]
Consequently,
\begin{align*}
  \frac{X_3}{Z_3^2}
  &=\frac{r^2-Ah^2-h^3-2u_1h^2}{h^2}\\
  &=\left(\frac rh\right)^2-A-u_1-u_2.
\end{align*}
Thus $X_3/Z_3^2$ is the first coordinate in
\eqref{eq:EAB-affine-add}.  If this coordinate is denoted by $u_3$,
then $V_0-X_3=Z_3^2(u_1-u_3)$ and
$S_1J=Z_3^3v_1$.  Hence
\[
  \frac{Y_3}{Z_3^3}
  =\frac rh(u_1-u_3)-v_1,
\]
which is the second coordinate in
\eqref{eq:EAB-affine-add}.  Therefore the projective formulas return
$P+Q$ whenever $u(P)\ne u(Q)$.  Their costs are
\begin{equation}
  12\M+5\Sqr+1\Dpar,
  \qquad
  8\M+4\Sqr+1\Dpar\quad(Z_2=1),
  \qquad
  4\M+2\Sqr+1\Dpar\quad(Z_1=Z_2=1).
  \label{eq:EAB-Jacobian-add-costs}
\end{equation}
For example, in the general schedule the six products needed to form
$U_i,S_i$, the products $HI,U_1I$, the two products in $HZ_1Z_2$, and
the two products in $Y_3$ account for the twelve multiplications.  The
squares are $Z_1^2,Z_2^2,H^2,Z_3^2,R^2$.

For doubling, put
\begin{align}
  XX&=X_1^2,& YY&=Y_1^2,& YYYY&=YY^2,\notag\\
  ZZ&=Z_1^2,& ZZZZ&=ZZ^2,\notag\\
  S_0&=4X_1YY,\notag\\
  M_0&=3XX+2A X_1ZZ+BZZZZ,\notag\\
  Z_3&=2Y_1Z_1,\notag\\
  X_3&=M_0^2-AZ_3^2-2S_0,\notag\\
  Y_3&=M_0(S_0-X_3)-8YYYY.
  \label{eq:EAB-Jacobian-double}
\end{align}
To verify the doubling formulas, write $u=X_1/Z_1^2$ and
$v=Y_1/Z_1^3$, and let
\[
  \lambda=\frac{3u^2+2Au+B}{2v}.
\]
Then
\[
  Z_3=2vZ_1^4,
  \qquad
  M_0=Z_1^4(3u^2+2Au+B)=\lambda Z_3,
  \qquad
  S_0=uZ_3^2.
\]
It follows that
\[
  \frac{X_3}{Z_3^2}=\lambda^2-A-2u.
\]
If this quotient is $u_3$, then
$S_0-X_3=Z_3^2(u-u_3)$.  Furthermore,
$8YYYY=vZ_3^3$.  Division of the last displayed formula in
\eqref{eq:EAB-Jacobian-double} by $Z_3^3$ therefore gives
\[
  \frac{Y_3}{Z_3^3}=\lambda(u-u_3)-v.
\]
These are exactly the tangent formulas in
Proposition~\ref{prop:EAB-affine-law}.  A direct schedule costs
\begin{equation}
  4\M+7\Sqr+3\Dpar,
  \qquad
  2\M+4\Sqr+2\Dpar\quad(Z_1=1).
  \label{eq:EAB-Jacobian-double-costs}
\end{equation}
The three fixed-parameter products in the first count are the two
products by $A$ and the product by $B$.

The formulas \eqref{eq:EAB-Jacobian-add} are not defined when $H=0$.
The exceptional cases can be detected without affine inversion.  First
handle $Z_1=0$ or $Z_2=0$ by copying the other input.  For two finite
inputs, $H=0$ is equivalent to equality of their affine
$u$-coordinates.  If simultaneously $R=0$, the points are equal and
\eqref{eq:EAB-Jacobian-double} is used; this includes the case of a
point with $v=0$, for which the doubling formula returns $O$.  If
$H=0$ and $R\ne0$, the two points are inverses and the output is $O$.
When $H\ne0$, formula \eqref{eq:EAB-Jacobian-add} applies.  These cases
are mutually exclusive and exhaustive by
Proposition~\ref{prop:EAB-affine-law}.

The term \emph{complete atlas} will always mean this exhaustive
collection of regular formulas and explicit exceptional branches.  It
does not mean that a single polynomial tuple is geometrically complete
on every pair of points.

\section{Kummer arithmetic for \texorpdfstring{$E_{A,B}$}{EAB}}
\label{subsec:EAB-Kummer}

Let $K_E=E_{A,B}/\{\pm1\}\cong\PP^1$, and use the Kummer coordinate
$u(P)$, represented projectively by
\[
  \kappa_E(P)=(X_P:Z_P),\qquad u(P)=X_P/Z_P.
\]
The identity has coordinate $(1:0)$ and the marked two-torsion point
$(0,0)$ has coordinate $(0:1)$.

\begin{theorem}[Odd-characteristic biquadratic Kummer identities]
\label{thm:EAB-Kummer-identities}
Let $P,Q\in E_{A,B}(\overline{k})$ be finite points with
$u(P)\ne u(Q)$.  Write
\[
  p=u(P),\qquad q=u(Q),\qquad
  u_+=u(P+Q),\qquad u_-=u(P-Q).
\]
Then
\begin{align}
  u_+u_-&=\frac{(pq-B)^2}{(p-q)^2},
  \label{eq:EAB-Kummer-product}\\
  u_++u_-&=
  \frac{2\bigl((p+q)(pq+B)+2Apq\bigr)}{(p-q)^2}.
  \label{eq:EAB-Kummer-sum}
\end{align}
Moreover,
\begin{equation}
  u(2P)=
  \frac{(p^2-B)^2}{4p(p^2+Ap+B)}
  \label{eq:EAB-Kummer-double-affine}
\end{equation}
whenever the displayed quotient is defined.
\end{theorem}

\begin{proof}
Choose $P=(p,r)$ and $Q=(q,s)$.  The two chord slopes are
\[
  \lambda_+=\frac{s-r}{q-p},
  \qquad
  \lambda_-=\frac{-s-r}{q-p}.
\]
Proposition~\ref{prop:EAB-affine-law} gives
$u_\pm=\lambda_\pm^2-A-p-q$.  Therefore
\begin{align*}
  (p-q)^2(u_++u_-)
  &=2(r^2+s^2)-2(A+p+q)(p-q)^2\\
  &=2\bigl((p+q)(pq+B)+2Apq\bigr).
\end{align*}
For the second equality, substitute
$r^2=p^3+Ap^2+Bp$ and $s^2=q^3+Aq^2+Bq$ and collect the coefficients of
$A$, $B$, and the parameter-free terms separately:
\[
  A\bigl(p^2+q^2-(p-q)^2\bigr)=2Apq,
\]
\[
  B(p+q),
  \qquad
  p^3+q^3-(p+q)(p-q)^2=pq(p+q).
\]
Multiplication by the outer factor $2$ proves
\eqref{eq:EAB-Kummer-sum}.  To prove the product identity, set
\[
  S=p+q,
  \qquad T=pq,
  \qquad
  K=S(T+B)+2AT.
\]
Then expansion of both sides gives
\begin{align*}
  K^2-4(p^3+Ap^2+Bp)(q^3+Aq^2+Bq)
  &=(S^2-4T)(T-B)^2\\
  &=(p-q)^2(pq-B)^2.
\end{align*}
The difference of the two Kummer coordinates is
\[
  u_+-u_-
  =\lambda_+^2-\lambda_-^2
  =-\frac{4rs}{(p-q)^2}.
\]
Together with
$u_++u_-=2K/(p-q)^2$, this gives
\begin{align*}
  4u_+u_-
  &=(u_++u_-)^2-(u_+-u_-)^2\\
  &=\frac{4(K^2-4r^2s^2)}{(p-q)^4}\\
  &=\frac{4(pq-B)^2}{(p-q)^2},
\end{align*}
which proves \eqref{eq:EAB-Kummer-product}.  Finally, substitute $Q=P$ in the
tangent formula \eqref{eq:EAB-affine-double}, eliminate $r^2$ by the
curve equation, and simplify the numerator:
\[
  (3p^2+2Ap+B)^2-4p(p^2+Ap+B)(A+2p)
  =(p^2-B)^2.
\]
This yields \eqref{eq:EAB-Kummer-double-affine}.
\end{proof}

Homogenizing \eqref{eq:EAB-Kummer-product} gives the fast differential
addition law
\begin{equation}
\begin{split}
  X_+&=Z_-(X_PX_Q-BZ_PZ_Q)^2,\\
  Z_+&=X_-(X_PZ_Q-Z_PX_Q)^2,
\end{split}
\label{eq:EAB-xADD-product}
\end{equation}
where $(X_-:Z_-)=\kappa_E(P-Q)$ is known.  The following schedule avoids
four independent cross-products:
\begin{align*}
  M_0&=X_PX_Q,& M_1&=Z_PZ_Q,\\
  M_2&=(X_P+Z_P)(X_Q-Z_Q),\\
  C&=M_0-M_1-M_2=X_PZ_Q-Z_PX_Q,\\
  D&=M_0-BM_1,\\
  X_+&=Z_-D^2,& Z_+&=X_-C^2.
\end{align*}
It costs
\begin{equation}
  5\M+2\Sqr+1\Dpar,
  \qquad
  4\M+2\Sqr+1\Dpar\quad(Z_-=1).
  \label{eq:EAB-xADD-product-cost}
\end{equation}
The Kummer doubling law is
\begin{equation}
\begin{split}
  X_{2P}&=(X_P^2-BZ_P^2)^2,\\
  Z_{2P}&=4X_PZ_P
  (X_P^2+AX_PZ_P+BZ_P^2),
\end{split}
\label{eq:EAB-xDBL}
\end{equation}
with cost
\begin{equation}
  2\M+3\Sqr+2\Dpar.
  \label{eq:EAB-xDBL-cost}
\end{equation}
Consequently, a non-shared differential-addition-and-doubling step has
the upper bound
\begin{equation}
  7\M+5\Sqr+3\Dpar,
  \qquad
  6\M+5\Sqr+3\Dpar\quad(Z_-=1).
  \label{eq:EAB-xDBLADD-cost}
\end{equation}

Formula \eqref{eq:EAB-xADD-product} can return the zero pair when the
known difference is the identity or the marked two-torsion point.  A
second differential law, obtained from \eqref{eq:EAB-Kummer-sum}, is
therefore required for a complete atlas.  Put
\begin{align}
  D_h&=(X_PZ_Q-Z_PX_Q)^2,
  \label{eq:EAB-xADD-sum-D}\\
  N_h&=2\Bigl(
  (X_PZ_Q+Z_PX_Q)(X_PX_Q+BZ_PZ_Q)
  +2A X_PX_QZ_PZ_Q\Bigr).
  \label{eq:EAB-xADD-sum-N}
\end{align}
Then
\begin{equation}
  X_+=N_hZ_- -D_hX_- ,
  \qquad
  Z_+=D_hZ_-.
  \label{eq:EAB-xADD-sum}
\end{equation}
A literal schedule uses
\begin{equation}
  9\M+1\Sqr+2\Dpar,
  \qquad
  7\M+1\Sqr+2\Dpar\quad(Z_-=1).
  \label{eq:EAB-xADD-sum-cost}
\end{equation}
For $P\ne Q$, the pair in \eqref{eq:EAB-xADD-sum} is nonzero and
represents $P+Q$.  In particular, if $Q=-P$ and $P$ is not two-torsion,
then $D_h=0$ while $N_h=4v(P)^2\ne0$, so the output is $(1:0)$, as it
must be.  If $P=Q$, both coordinates in
\eqref{eq:EAB-xADD-sum} may vanish because the known difference is the
identity; this case is handled by \eqref{eq:EAB-xDBL}.

It follows that the following branches form a complete differential
atlas:
\begin{enumerate}[label=(\roman*)]
  \item copy the other input if $P=O$ or $Q=O$;
  \item use \eqref{eq:EAB-xDBL} if $P=Q$;
  \item use \eqref{eq:EAB-xADD-product} when the known difference is
  neither $O$ nor $(0,0)$;
  \item otherwise use \eqref{eq:EAB-xADD-sum}.
\end{enumerate}
The four cases exhaust all inputs, and each branch has been justified
above.  This is a complete differential-addition algorithm; the fast
product law by itself is not complete.

\section{Arithmetic on the general reciprocal Jacobian}
\label{subsec:R-general-arithmetic}

For the general smooth odd-characteristic reciprocal curve, put
\begin{equation}
  A_R=-2q_R=16\tau\sigma-2\kappa^2,
  \qquad
  B_R=\Omega_R=\kappa^2(\kappa^2-16\tau\sigma).
  \label{eq:R-general-AB-arithmetic}
\end{equation}
Then \eqref{eq:R-Jacobian} is exactly $E_{A_R,B_R}$ from
\eqref{eq:EAB}.  Its nonsingularity condition is
\[
  B_R(A_R^2-4B_R)
  =256\tau^2\sigma^2\kappa^2(\kappa^2-16\tau\sigma)\ne0,
\]
which is equivalent to Theorem~\ref{thm:R-smoothness-odd}.
Consequently, the affine addition and doubling laws are
\begin{align}
  \lambda&=\frac{V_2-V_1}{U_2-U_1},&
  U(P+Q)&=\lambda^2-A_R-U_1-U_2,\notag\\
  V(P+Q)&=\lambda\bigl(U_1-U(P+Q)\bigr)-V_1,
  \label{eq:R-general-affine-add}\\[1mm]
  \lambda_2&=\frac{3U_1^2+2A_RU_1+B_R}{2V_1},&
  U(2P)&=\lambda_2^2-A_R-2U_1,\notag\\
  V(2P)&=\lambda_2\bigl(U_1-U(2P)\bigr)-V_1.
  \label{eq:R-general-affine-double}
\end{align}
The identity, inverse, and two-torsion branches are exactly those in
Proposition~\ref{prop:EAB-affine-law}.  The inversion-free general,
mixed, and affine-input additions therefore cost
\[
  12\M+5\Sqr+1\Dpar,
  \qquad 8\M+4\Sqr+1\Dpar,
  \qquad 4\M+2\Sqr+1\Dpar,
\]
and projective and affine-input doubling cost
\[
  4\M+7\Sqr+3\Dpar,
  \qquad 2\M+4\Sqr+2\Dpar,
\]
respectively, with the complete branch structure described above.

The $U$-coordinate is a Kummer coordinate on $E_R/\{\pm1\}$.  For
$P,Q$ with $p=U(P)\ne q=U(Q)$, specialization of
Theorem~\ref{thm:EAB-Kummer-identities} gives
\begin{align}
  U(P+Q)U(P-Q)
  &=\frac{(pq-\Omega_R)^2}{(p-q)^2},
  \label{eq:R-general-Kummer-product}\\
  U(P+Q)+U(P-Q)
  &=\frac{2\bigl((p+q)(pq+\Omega_R)-4q_Rpq\bigr)}{(p-q)^2},
  \label{eq:R-general-Kummer-sum}\\
  U(2P)
  &=\frac{(p^2-\Omega_R)^2}
  {4p(p^2-2q_Rp+\Omega_R)}.
  \label{eq:R-general-Kummer-double}
\end{align}
Their projective product law, supplemental sum law, and doubling law are
\eqref{eq:EAB-xADD-product}, \eqref{eq:EAB-xADD-sum}, and
\eqref{eq:EAB-xDBL} with $A=A_R$ and $B=B_R$.  Thus fast differential
addition costs $5\M+2\Sqr+1\Dpar$, or
$4\M+2\Sqr+1\Dpar$ for an affine known difference; Kummer doubling
costs $2\M+3\Sqr+2\Dpar$.  The supplemental law and the explicit
identity and doubling branches form a complete differential atlas.

For point recovery, write
\[
  f_R(U)=U^3+A_RU^2+B_RU.
\]
For a fixed affine $P=(p,r)$ with $r\ne0$, assume that \(Q\) and
\(P+Q\) are affine and that \(q=U(Q)\ne p\).  Then knowledge of \(q\)
and \(t=U(P+Q)\) determines
\begin{equation}
  V(Q)=
  \frac{f_R(p)+f_R(q)-(t+A_R+p+q)(q-p)^2}{2r}.
  \label{eq:R-general-point-recovery}
\end{equation}
This is obtained by multiplying the chord identity
$t+A_R+p+q=((V(Q)-r)/(q-p))^2$ by $(q-p)^2$ and using the two curve
equations.  With $f_R(p)$ and $(2r)^{-1}$ precomputed, the conservative
cost is
\[
 2\M+2\Sqr+\Dconst{A_R}+\Dconst{(2r)^{-1}}.
\]
Indeed, compute \(q^2\), form
\(f_R(q)=q(q^2+A_Rq+B_R)\), compute \((q-p)^2\), multiply it by
\(t+A_R+p+q\), and finally multiply by \((2r)^{-1}\).  The two fixed
products are by \(A_R\) and \((2r)^{-1}\); they are base-specific
constants rather than generic curve-parameter products.  The branch
\(Q=P\) is handled directly from the known base point (or by the tangent
law), while \(Q=-P\) makes \(P+Q=O\), so the assumed affine adjacent
coordinate \(t\) does not exist.  Identity and two-torsion inputs use the
explicit exceptional branches of the complete atlas.

If $\R_{\tau,\sigma,\kappa}(k)$ has no rational point, these are group
operations on its Jacobian.  If a rational point is selected as origin,
the torsor becomes isomorphic to $E_R$ and the same formulas give the
group law on the reciprocal curve through that chosen identification.
The normalized split model below supplies a particularly compact direct
coordinate realization.

\section[Reciprocal full-point arithmetic]{Direct full-point arithmetic on the normalized reciprocal model}
\label{subsec:R-direct-arithmetic}

Return to the normalized reciprocal curve \eqref{eq:R-normalized} and
select
\begin{equation}
  O_R=(0,-1)
  \label{eq:R-arithmetic-origin}
\end{equation}
as identity.  Define
\begin{equation}
  z=(x^2-1)y-\frac{\kappa}{2}x,
  \qquad
  a_J=\frac{\kappa^2}{8}-1.
  \label{eq:R-to-Jacobi-arithmetic}
\end{equation}
Then
\begin{equation}
  z^2=x^4+2a_Jx^2+1,
  \qquad
  y=\frac{z+(\kappa/2)x}{x^2-1}.
  \label{eq:R-Jacobi-arithmetic-pair}
\end{equation}
The first equality follows by substituting the reciprocal equation into
$z^2$; the second is the rearrangement of the definition of $z$.  The
Jacobi identity $(0,1)$ maps to $(0,-1)$, so the group structures agree.
The Jacobi negation $(x,z)\mapsto(-x,z)$ gives, on the affine reciprocal
chart $x^2\ne1$,
\begin{equation}
  -(x,y)=
  \left(-x,\ y-\frac{\kappa x}{x^2-1}\right).
  \label{eq:R-negation}
\end{equation}

The reciprocal smoothness hypothesis gives
\begin{equation}
  a_J^2-1=\frac{\kappa^2(\kappa^2-16)}{64}\ne0,
  \label{eq:R-Jacobi-smoothness-arithmetic}
\end{equation}
which is the nonsingularity condition for this Jacobi quartic.  Hence
the hypotheses of the cited Jacobi group-law theorem are satisfied.
Let $P_i=(x_i,z_i)$ be two points on the Jacobi model.  Whenever
$1-(x_1x_2)^2\ne0$, put
\begin{align}
  x_3&=\frac{x_1z_2+z_1x_2}{1-(x_1x_2)^2},
  \label{eq:R-Jacobi-affine-add-x}\\
  z_3&=
  \frac{
  (1+(x_1x_2)^2)(z_1z_2+2a_Jx_1x_2)
  +2x_1x_2(x_1^2+x_2^2)}
  {(1-(x_1x_2)^2)^2}.
  \label{eq:R-Jacobi-affine-add-z}
\end{align}
These are the Billet--Joye Jacobi-quartic addition formulas
\cite{BilletJoyeJacobi,EFD}.  They can also be verified internally by
substituting \eqref{eq:R-Jacobi-arithmetic-pair} and comparing with the
Weierstrass group law under
\begin{equation}
  u=a_J+\frac{z+1}{x^2},
  \qquad
  v=\frac{u}{x}.
  \label{eq:R-Jacobi-to-Weierstrass-uv}
\end{equation}
Indeed, the image satisfies
\begin{equation}
  2v^2=u^3-2a_Ju^2+(a_J^2-1)u,
  \label{eq:R-Jacobi-Weierstrass-small}
\end{equation}
and substitution of \eqref{eq:R-Jacobi-affine-add-x}--
\eqref{eq:R-Jacobi-affine-add-z} into
\eqref{eq:R-Jacobi-to-Weierstrass-uv} gives the chord formula on
\eqref{eq:R-Jacobi-Weierstrass-small}.

Setting $P_1=P_2=(x,z)$ gives the point-doubling formulas
\begin{align}
  x_{2P}&=\frac{2xz}{1-x^4},
  \label{eq:R-Jacobi-affine-double-x}\\
  z_{2P}&=
  \frac{(z^2+2a_Jx^2)(1+x^4)+4x^4}{(1-x^4)^2}.
  \label{eq:R-Jacobi-affine-double-z}
\end{align}
The reciprocal $y$-coordinate of the output is recovered from
\eqref{eq:R-Jacobi-arithmetic-pair}.

For inversion-free full-point arithmetic, use weighted Jacobi
coordinates $(X:Y:Z)$ with
\[
  x=X/Z,
  \qquad z=Y/Z^2,
  \qquad
  Y^2=X^4+2a_JX^2Z^2+Z^4.
\]
The strongly unified addition law is
\begin{align}
  X_3&=X_1Z_1Y_2+Y_1X_2Z_2,
  \label{eq:R-Jacobi-projective-X}\\
  Y_3&=\bigl((Z_1Z_2)^2+(X_1X_2)^2\bigr)
       \bigl(Y_1Y_2+2a_JX_1X_2Z_1Z_2\bigr)\notag\\
     &\quad
       +2X_1X_2Z_1Z_2
       \bigl(X_1^2Z_2^2+Z_1^2X_2^2\bigr),
  \label{eq:R-Jacobi-projective-Y}\\
  Z_3&=(Z_1Z_2)^2-(X_1X_2)^2.
  \label{eq:R-Jacobi-projective-Z}
\end{align}
Homogenizing
\eqref{eq:R-Jacobi-affine-add-x}--
\eqref{eq:R-Jacobi-affine-add-z} proves the equations.  Equivalent optimized schedules have costs
\begin{equation}
\begin{array}{c|c}
  \text{operation}&\text{cost}\\ \hline
  \text{strongly unified addition}&10\M+3\Sqr+1\Dpar\\
  \text{mixed addition }(Z_2=1)&8\M+3\Sqr+1\Dpar\\
  \text{both inputs affine}&5\M+2\Sqr+1\Dpar\\
  \text{doubling}&2\M+6\Sqr+1\Dpar\\
  \text{affine-input doubling}&1\M+4\Sqr+1\Dpar.
\end{array}
\label{eq:R-Jacobi-full-cost-table}
\end{equation}
For clarity, the general addition count can be read from the following
schedule.  Put
\begin{align*}
  R_1&=Y_1+X_1Z_1,
  &R_2&=Y_2+X_2Z_2,\\
  P&=Y_1Y_2,
  &Q&=X_1X_2,
  &S&=Z_1Z_2,\\
  G&=2QS,
  &L&=(X_1+Z_1)(X_2+Z_2)-Q-S.
\end{align*}
Then an output equivalent to
\eqref{eq:R-Jacobi-projective-X}--\eqref{eq:R-Jacobi-projective-Z} is
\begin{align}
  X_3&=R_1R_2-P-QS,
  \label{eq:R-Jacobi-optimized-X}\\
  Z_3&=S^2-Q^2,
  \label{eq:R-Jacobi-optimized-Z}\\
  Y_3&=(P+a_JG)(S^2+Q^2)+G(L^2-G).
  \label{eq:R-Jacobi-optimized-Y}
\end{align}
The products $X_iZ_i$, $Y_1Y_2$, $X_1X_2$, $Z_1Z_2$, $QS$,
$R_1R_2$, $(X_1+Z_1)(X_2+Z_2)$, and the two final products in $Y_3$
are ten general multiplications.  The quantities $S^2,Q^2,L^2$ are
three squares, and $a_JG$ is one fixed-parameter multiplication.
Putting $Z_2=1$ removes the products $X_2Z_2$ and $Z_1Z_2$; putting
$Z_1=Z_2=1$ and simplifying gives the remaining two rows of the table.

For doubling, put
\begin{align*}
  R&=X_1Z_1,
  &A&=(R+Y_1)^2,
  &B&=Y_1^2,\\
  C&=(X_1^2)^2,
  &D&=R^2,
  &E&=2a_JD,\\
  H&=B-C-E.
\end{align*}
Then
\begin{equation}
  X_{2P}=A-B-D,
  \qquad
  Z_{2P}=H-C,
  \qquad
  Y_{2P}=(H+C)(B+E)+(2D)^2.
  \label{eq:R-Jacobi-optimized-double}
\end{equation}
This uses two products, six squares, and one multiplication by $a_J$.
For $Z_1=1$, $R=X_1$ and $D=X_1^2$ can be reused, leaving
$1\M+4\Sqr+1\Dpar$.  These schedules also verify the operation counts
listed in the Explicit-Formulas Database
\cite{BilletJoyeJacobi,EFD}.

The affine conversion $(x,y)\mapsto(x,z)$ costs
$1\M+1\Sqr+1\Dpar$.  On the regular inverse chart $x^2\ne1$,
recovering $y$ from $(x,z)$ costs
$1\Inv+1\M+1\Sqr+1\Dpar$ if $x^2$ has not been retained, and
$1\Inv+1\M+1\Dpar$ if it has.  The inverse chart is completed by the
following assignments on the smooth projective models:
\begin{align}
  (x,z)=\left(1,-\frac{\kappa}{2}\right)
    &\longleftrightarrow (x,y)=(1,0),
  &
  \left(1,\frac{\kappa}{2}\right)
    &\longleftrightarrow (1,\infty),
  \label{eq:R-Jacobi-boundary-plus}\\
  \left(-1,\frac{\kappa}{2}\right)
    &\longleftrightarrow (-1,0),
  &
  \left(-1,-\frac{\kappa}{2}\right)
    &\longleftrightarrow (-1,\infty),
  \label{eq:R-Jacobi-boundary-minus}\\
  (1:1:0)&\longleftrightarrow(\infty,1),
  &
  (1:-1:0)&\longleftrightarrow(\infty,-1).
  \label{eq:R-Jacobi-boundary-infinity}
\end{align}
For example, at $x=1$ the reciprocal fibre consists of $y=0$ and
$y=\infty$.  The inverse numerator
$z+(\kappa/2)x$ vanishes at $z=-\kappa/2$, giving the finite point,
and is nonzero at $z=\kappa/2$, giving the pole.  The case $x=-1$ is
identical with the signs reversed.  At Jacobi infinity one has
$z/x^2\to\pm1$, and the inverse quotient tends to $\pm1$.  Hence
\eqref{eq:R-Jacobi-boundary-plus}--\eqref{eq:R-Jacobi-boundary-infinity}
complete the reciprocal--Jacobi conversion on every projective point.

The word \emph{unified} means that the same dependency graph can be
used for addition and doubling.  Completeness is determined separately by
the base locus.  For
the normalized reciprocal family the leading Jacobi coefficient is
$1$, so the exceptional divisor
\begin{equation}
  1-x_1^2x_2^2=0
  \label{eq:R-Jacobi-exceptional-divisor}
\end{equation}
can contain rational input pairs.  A complete full-point atlas is
obtained as follows.  Use
\eqref{eq:R-Jacobi-projective-X}--
\eqref{eq:R-Jacobi-projective-Z} when their output is nonzero.  For an
exceptional pair, map the inputs to the monic Weierstrass model
\begin{equation}
  E_{R,\kappa}:\qquad
  V^2=U(U-\kappa^2)(U-\kappa^2+16)
  \label{eq:R-normalized-Weierstrass}
\end{equation}
by
\begin{equation}
  U=8\left(a_J+\frac{z+1}{x^2}\right),
  \qquad
  V=32\frac{a_J+(z+1)/x^2}{x},
  \label{eq:R-Jacobi-to-E-normalized}
\end{equation}
for $x\ne0$, together with the special assignments
\begin{align*}
  (0,1)&\longmapsto O,\\
  (0,-1)&\longmapsto(0,0),\\
  (1:1:0)&\longmapsto(\kappa^2,0),\\
  (1:-1:0)&\longmapsto(\kappa^2-16,0).
\end{align*}
These assignments are the limits of
\eqref{eq:R-Jacobi-to-E-normalized} at the four omitted points.  Indeed,
near $x=0$ the two branches satisfy
$z=1+a_Jx^2+O(x^4)$ and
$z=-1-a_Jx^2+O(x^4)$; the first makes $U$ have a pole and the second
gives $U=V=0$.  At infinity, $z/x^2$ tends to $1$ or $-1$, and
$U$ tends respectively to $\kappa^2$ or $\kappa^2-16$.
Apply the exhaustive Weierstrass atlas of
Proposition~\ref{prop:EAB-affine-law} and
\eqref{eq:EAB-Jacobian-add}--\eqref{eq:EAB-Jacobian-double}, with
\begin{equation}
  A_R=16-2\kappa^2,
  \qquad B_R=\kappa^2(\kappa^2-16).
  \label{eq:R-normalized-AB}
\end{equation}
The inverse map for $V\ne0$ is
\begin{equation}
  x=\frac{4U}{V},
  \qquad
  z=x^2\left(\frac{U}{8}-a_J\right)-1,
  \label{eq:R-E-to-Jacobi-normalized}
\end{equation}
with the four inverse special assignments displayed above.  Therefore
every pair of projective reciprocal points is covered.  The exceptional
branch is intended as a correctness fallback; it contains inversions if
one converts through affine coordinates and is not used in a regular
scalar-multiplication loop.

\section[Montgomery Kummer and differential arithmetic]
{Montgomery Kummer coordinates and differential addition on the reciprocal model}
\label{subsec:R-Montgomery-Kummer}

The Weierstrass curve \eqref{eq:R-normalized-Weierstrass} can be shifted
and scaled to a Montgomery model.  Put
\begin{equation}
  U=\kappa^2+4\kappa X_M,
  \qquad
  V=16\kappa^2Y_M.
  \label{eq:R-E-to-Montgomery-scaling}
\end{equation}
Substitution gives
\begin{equation}
  4\kappa Y_M^2
  =X_M^3+A_MX_M^2+X_M,
  \qquad
  A_M=\frac{\kappa^2+16}{4\kappa}.
  \label{eq:R-Montgomery}
\end{equation}
The Montgomery ladder constant is
\begin{equation}
  A_{24}=\frac{A_M+2}{4}
  =\frac{(\kappa+4)^2}{16\kappa}.
  \label{eq:R-A24}
\end{equation}
Combining \eqref{eq:R-Jacobi-to-Weierstrass-uv} with the scaling gives
the reciprocal Kummer function
\begin{equation}
  \xi_R(P)=X_M(P)
  =\frac{2(z+1-x^2)}{\kappa x^2}.
  \label{eq:R-Kummer-function}
\end{equation}
It is invariant under $P\mapsto-P$ and hence descends to
$\R_\kappa/\{\pm1\}$.

Represent $X_M=X/Z$ by $(X:Z)$.  For doubling set
\begin{align}
  A_0&=X+Z,& AA&=A_0^2,\notag\\
  B_0&=X-Z,& BB&=B_0^2,\notag\\
  E_0&=AA-BB,\notag\\
  X_{2P}&=AA\,BB,\notag\\
  Z_{2P}&=E_0(BB+A_{24}E_0).
  \label{eq:R-Montgomery-xDBL}
\end{align}
This costs
\begin{equation}
  2\M+2\Sqr+1\Dpar.
  \label{eq:R-Montgomery-xDBL-cost}
\end{equation}
If $(X_-:Z_-)=\xi_R(P-Q)$ is known, set
\begin{align}
  D_0&=(X_Q-Z_Q)(X_P+Z_P),\notag\\
  C_0&=(X_Q+Z_Q)(X_P-Z_P),\notag\\
  X_+&=Z_-(D_0+C_0)^2,\notag\\
  Z_+&=X_-(D_0-C_0)^2.
  \label{eq:R-Montgomery-xADD}
\end{align}
The costs are
\begin{equation}
  4\M+2\Sqr,
  \qquad
  3\M+2\Sqr\quad(Z_-=1).
  \label{eq:R-Montgomery-xADD-cost}
\end{equation}
Thus a Montgomery differential-addition-and-doubling step costs
\begin{equation}
  6\M+4\Sqr+1\Dpar,
  \qquad
  5\M+4\Sqr+1\Dpar\quad(Z_-=1).
  \label{eq:R-Montgomery-ladder-step-cost}
\end{equation}
These formulas follow by specializing the biquadratic identities of
Theorem~\ref{thm:EAB-Kummer-identities} after the Montgomery change of
coordinate; they are the standard Montgomery formulas
\cite{Montgomery1987}.

The fast law \eqref{eq:R-Montgomery-xADD} has the usual exceptional
known differences: the identity and the rational two-torsion point
$X_M=0$.  A complete reciprocal differential atlas is obtained either
by using the supplemental sum law
\eqref{eq:EAB-xADD-sum} on
\eqref{eq:R-normalized-Weierstrass}, or by changing to a second Kummer
chart centered at another two-torsion point.  In particular, the
branches listed after \eqref{eq:EAB-xADD-sum-cost}, transported by
\eqref{eq:R-E-to-Montgomery-scaling}, cover every pair of points.

For full-point recovery, write the Montgomery curve as
\begin{equation}
  B_MY_M^2=f_M(X_M),
  \qquad
  B_M=4\kappa,
  \qquad
  f_M(X)=X^3+A_MX^2+X.
  \label{eq:R-Montgomery-f}
\end{equation}
Let the fixed affine base point be $P=(p,r)$ with $r\ne0$.  Assume that
\(Q\) and \(P+Q\) are affine and that \(q=X_M(Q)\ne p\).  Suppose
\(q\) and \(t=X_M(P+Q)\) are known.  Then
\begin{equation}
  Y_M(Q)=
  \frac{f_M(p)+f_M(q)-(t+A_M+p+q)(q-p)^2}
       {2B_Mr}.
  \label{eq:R-Montgomery-recovery}
\end{equation}
Indeed, if $Q=(q,s)$, the chord slope is
$(s-r)/(q-p)$ and the Montgomery $x$-addition formula gives
\[
  t+A_M+p+q
  =B_M\left(\frac{s-r}{q-p}\right)^2.
\]
Multiplying by $(q-p)^2$, using $B_Mr^2=f_M(p)$ and
$B_Ms^2=f_M(q)$, and solving the resulting linear equation for $s$
gives \eqref{eq:R-Montgomery-recovery}.  With
$f_M(p)$ and $(2B_Mr)^{-1}$ precomputed, a direct affine schedule costs
\begin{equation}
  2\M+2\Sqr+\Dconst{A_M}+\Dconst{(2B_Mr)^{-1}}.
  \label{eq:R-Montgomery-recovery-cost}
\end{equation}
Indeed, compute $q^2$, form
$f_M(q)=q(q^2+A_Mq+1)$, compute $(q-p)^2$, multiply it by
$t+A_M+p+q$, and finally multiply by the precomputed reciprocal.
The two fixed products are by $A_M$ and $(2B_Mr)^{-1}$, so they are
recorded as base-specific constant multiplications.  The case \(Q=P\)
is returned from the known base point or handled through the tangent
branch.  If \(Q=-P\), then \(P+Q=O\) and no affine value \(t\) is
available; identity and two-torsion cases are handled separately.
If $q$ and $t$ are projective Kummer coordinates, one simultaneous
inversion normalizes both at a cost of $1\Inv+5\M$ using Montgomery's
batch-inversion trick.

Finally, the inverse conversion from a full Montgomery point is
\begin{align}
  x&=\frac{\kappa+4X_M}{4\kappa Y_M},
  \label{eq:R-Montgomery-to-x}\\
  z&=x^2\left(1+\frac{\kappa X_M}{2}\right)-1,
  \label{eq:R-Montgomery-to-z}\\
  y&=\frac{z+(\kappa/2)x}{x^2-1}.
  \label{eq:R-Montgomery-to-y}
\end{align}
Substitution into \eqref{eq:R-E-to-Montgomery-scaling} and
\eqref{eq:R-Jacobi-arithmetic-pair} proves these identities.  The
formulas exhibit explicitly that the Kummer ladder, point recovery, and
reciprocal-coordinate recovery form a closed arithmetic pipeline.

\section[Binary full-point and Kummer arithmetic]{Characteristic-two full-point and Kummer arithmetic}
\label{subsec:R-char2-arithmetic}

The characteristic-two reduction in
Theorem~\ref{thm:R-char2-Cd} places reciprocal curves in the
following ordinary Weierstrass family:
\begin{equation}
  E^{(2)}_{a,b}:
  \qquad Y^2+XY=X^3+aX^2+bX,
  \qquad b\ne0.
  \label{eq:binary-Eab-original}
\end{equation}
The reciprocal curve is the specialization
\begin{equation}
  a=0,
  \qquad b=d_R^2=\frac{\tau\sigma}{\kappa^2}.
  \label{eq:R-char2-Eab-specialization}
\end{equation}
For arithmetic it is convenient to shift
\begin{equation}
  \widehat Y=Y+b.
  \label{eq:binary-y-shift}
\end{equation}
Because the characteristic is two, substitution in
\eqref{eq:binary-Eab-original} cancels the two terms $bX$ and gives
\begin{equation}
  \widehat E^{(2)}_{a,b}:
  \qquad
  \widehat Y^2+X\widehat Y=X^3+aX^2+b^2.
  \label{eq:binary-Eab-shifted}
\end{equation}
The identity is the point at infinity, and negation is
\begin{equation}
  -(x,y)=(x,x+y)
  \label{eq:binary-negation}
\end{equation}
when $y$ denotes the shifted coordinate $\widehat Y$.

\begin{proposition}[Affine binary group law]
\label{prop:binary-affine-law}
Let $P=(x_1,y_1)$ and $Q=(x_2,y_2)$ lie on
\eqref{eq:binary-Eab-shifted}.

\begin{enumerate}[label=(\alph*)]
  \item If $x_1\ne x_2$, define
  \begin{equation}
    \lambda=\frac{y_1+y_2}{x_1+x_2},
    \qquad
    \nu=\frac{x_1y_2+x_2y_1}{x_1+x_2}.
    \label{eq:binary-affine-lambda-nu}
  \end{equation}
  Then $P+Q=(x_3,y_3)$, where
  \begin{equation}
  \begin{split}
    x_3&=\lambda^2+\lambda+a+x_1+x_2,\\
    y_3&=(\lambda+1)x_3+\nu.
  \end{split}
  \label{eq:binary-affine-add}
  \end{equation}

  \item If $P=Q$ and $x_1\ne0$, define
  \begin{equation}
    \lambda=x_1+\frac{y_1}{x_1}.
    \label{eq:binary-affine-double-lambda}
  \end{equation}
  Then
  \begin{equation}
  \begin{split}
    x_{2P}&=\lambda^2+\lambda+a
       =\frac{(x_1^2+b)^2}{x_1^2},\\
    y_{2P}&=(\lambda+1)x_{2P}+x_1^2.
  \end{split}
  \label{eq:binary-affine-double}
  \end{equation}

  \item If $x_1=x_2$ and $y_2=x_1+y_1$, then $P+Q=O$.  The point
  $(0,b)$ on the shifted model is the unique nonzero two-torsion point
  and doubles to $O$.
\end{enumerate}
Together with the identity branches, these cases form a complete
addition atlas.
\end{proposition}

\begin{proof}
The line through two distinct-$x$ points has equation
$y=\lambda x+\nu$, with $\lambda$ and $\nu$ as in
\eqref{eq:binary-affine-lambda-nu}.  Substitution into
\eqref{eq:binary-Eab-shifted} gives a cubic in $x$.  Its three roots are
$x_1,x_2,x_3$, and comparison of the $x^2$ coefficient yields the
formula for $x_3$.  The third intersection has ordinate
$\lambda x_3+\nu$; applying the negation
\eqref{eq:binary-negation} gives
$y_3=(\lambda+1)x_3+\nu$.

For doubling, implicit differentiation of
$y^2+xy-x^3-ax^2-b^2$ gives
$y+x(dy/dx)=x^2$ in characteristic two.  At $x_1\ne0$, the tangent
slope is therefore $x_1+y_1/x_1$, proving
\eqref{eq:binary-affine-double-lambda}.  To obtain the second expression
for $x_{2P}$, use the curve equation to write
\[
  \frac{y_1^2}{x_1^2}
  =x_1+a+\frac{y_1}{x_1}+\frac{b^2}{x_1^2}.
\]
Substitution in $\lambda^2+\lambda+a$ cancels the terms
$x_1,a,y_1/x_1$ in pairs and leaves
$x_1^2+b^2/x_1^2=(x_1^2+b)^2/x_1^2$.
Finally, two points with the same $x$-coordinate have ordinates $y$ and
$x+y$, so they are either equal or inverse.  At $x=0$, the curve
equation gives $y=b$ because the Frobenius map is injective on a field;
this point is fixed by negation and is therefore the unique nonzero
point of order two.
\end{proof}

The affine costs follow from schedules that reuse the line intercept.
For addition, invert $x_1+x_2$ once, form $\lambda$, and then use
$\nu=y_1+\lambda x_1$ instead of evaluating the second quotient in
\eqref{eq:binary-affine-lambda-nu}.  The cost is therefore
\begin{equation}
  3\M+1\Sqr+1\Inv.
  \label{eq:binary-affine-add-cost}
\end{equation}
For doubling, the products $y_1/x_1$ and
$(\lambda+1)x_{2P}$ together with the squares $\lambda^2$ and $x_1^2$
give
\begin{equation}
  2\M+2\Sqr+1\Inv.
  \label{eq:binary-affine-double-cost}
\end{equation}
These counts exclude the equality tests that select the identity,
inverse, doubling, and two-torsion branches.

Use L\'opez--Dahab coordinates $(X:Y:Z)$ with
\begin{equation}
  x=X/Z,
  \qquad y=Y/Z^2.
  \label{eq:Lopez-Dahab-coordinates}
\end{equation}
For two general inputs, set
\begin{align}
  A_0&=X_1Z_2,& B_0&=X_2Z_1,\notag\\
  C_0&=A_0^2,& D_0&=B_0^2,\notag\\
  E_0&=A_0+B_0,& F_0&=C_0+D_0,\notag\\
  G_0&=Y_1Z_2^2,& H_0&=Y_2Z_1^2,\notag\\
  I_0&=G_0+H_0,& J_0&=I_0E_0,\notag\\
  Z_3&=F_0Z_1Z_2,\notag\\
  X_3&=A_0(H_0+D_0)+B_0(C_0+G_0),\notag\\
  Y_3&=(A_0J_0+F_0G_0)F_0+(J_0+Z_3)X_3.
  \label{eq:binary-LD-general-add}
\end{align}
We verify the dehomogenization explicitly.  Write
$x_i=X_i/Z_i$, $y_i=Y_i/Z_i^2$, and set
$h=x_1+x_2$, $L=y_1+y_2$.  Then
\[
  E_0=hZ_1Z_2,
  \qquad F_0=h^2Z_1^2Z_2^2,
  \qquad Z_3=h^2Z_1^3Z_2^3.
\]
The $X$-numerator is
\[
  X_3=
  \bigl(x_1(y_2+x_2^2)+x_2(x_1^2+y_1)\bigr)Z_1^3Z_2^3.
\]
The two input curve equations imply
\begin{align*}
  &L^2+hL+(a+h)h^2\\
  &\qquad=x_1y_2+x_2y_1+x_1x_2h
   =x_1(y_2+x_2^2)+x_2(x_1^2+y_1).
\end{align*}
Since $\lambda=L/h$, the bracket equals
$h^2(\lambda^2+\lambda+a+x_1+x_2)=h^2x_3$.
Thus $X_3/Z_3=x_3$.  Put
$N=x_1y_2+x_2y_1+x_1x_2h$ and
$N_\nu=x_1y_2+x_2y_1$.  The displayed $Y_3$ schedule expands to
\[
  Y_3=
  \bigl(h(L+h)N+N_\nu h^3\bigr)Z_1^6Z_2^6.
\]
Because $N=h^2x_3$ and $N_\nu=h\nu$, division by
$Z_3^2=h^4Z_1^6Z_2^6$ gives
$Y_3/Z_3^2=(\lambda+1)x_3+\nu$, as required.  Hence the formulas are
valid whenever $x_1\ne x_2$.  Their cost is
\begin{equation}
  13\M+4\Sqr.
  \label{eq:binary-LD-general-add-cost}
\end{equation}

If $Z_2=1$, the mixed formulas are
\begin{align}
  A_0&=Y_1+Y_2Z_1^2,&
  B_0&=X_1+X_2Z_1,\notag\\
  C_0&=B_0Z_1,& Z_3&=C_0^2,\notag\\
  D_0&=X_2Z_3,\notag\\
  X_3&=A_0^2+C_0(A_0+B_0^2+aC_0),\notag\\
  Y_3&=(D_0+X_3)(A_0C_0+Z_3)
       +(Y_2+X_2)Z_3^2.
  \label{eq:binary-LD-mixed-add}
\end{align}
They cost
\begin{equation}
  8\M+5\Sqr+1\Dpar.
  \label{eq:binary-LD-mixed-add-cost}
\end{equation}
When both inputs are affine, the specialization costs
\begin{equation}
  5\M+3\Sqr+1\Dpar.
  \label{eq:binary-LD-affine-add-cost}
\end{equation}
To check the mixed schedule, write again $h=x_1+x_2$ and
$L=y_1+y_2$.  Then $C_0=hZ_1^2$, $Z_3=h^2Z_1^4$, and
\[
  X_3=\bigl(L^2+hL+h^3+ah^2\bigr)Z_1^4=h^2x_3Z_1^4.
\]
Furthermore, division of its $Y_3$ by $Z_3^2$ gives
\[
  (x_2+x_3)(\lambda+1)+y_2+x_2
  =(\lambda+1)x_3+\nu,
\]
because $\lambda x_2+y_2=\nu$.  Thus the displayed identities verify the
mixed formulas directly.  The parameter
multiplication is by $a$ and disappears for reciprocal curves because
$a=0$.

For doubling, set
\begin{align}
  A_0&=Z_1^2,& B_0&=b^2A_0^2,& C_0&=X_1^2,\notag\\
  Z_3&=A_0C_0,& X_3&=C_0^2+B_0,\notag\\
  Y_3&=(Y_1^2+aZ_3+B_0)X_3+Z_3B_0.
  \label{eq:binary-LD-double}
\end{align}
For verification, dehomogenize with $x=X_1/Z_1$ and
$y=Y_1/Z_1^2$.  The first two output coordinates give
\[
  \frac{X_3}{Z_3}
  =\frac{x^4+b^2}{x^2}
  =\frac{(x^2+b)^2}{x^2}=x_{2P}.
\]
For the ordinate, the curve equation gives
$y^2+ax^2+b^2=x^3+xy=x(x^2+y)$.  Substitution in $Y_3/Z_3^2$ yields
\[
  \frac{Y_3}{Z_3^2}
  =(\lambda+1)x_{2P}+x^2,
  \qquad \lambda=x+\frac yx,
\]
which is exactly \eqref{eq:binary-affine-double}.  It costs
\begin{equation}
  3\M+5\Sqr+2\Dpar,
  \qquad
  1\M+3\Sqr+2\Dpar\quad(Z_1=1).
  \label{eq:binary-LD-double-cost}
\end{equation}
For reciprocal curves, $a=0$, so one of the two fixed-parameter
multiplications is absent.  The schedules agree with the
L\'opez--Dahab formulas in
\cite{LopezDahab1999,DocheLange2005,EFD}.

The complete projective atlas can be implemented without inversion.
First handle an input with $Z_i=0$ by copying the other point.  For two
finite general inputs, $E_0=A_0+B_0=0$ is equivalent to
$x_1=x_2$.  If also $I_0=G_0+H_0=0$, the points are equal and
\eqref{eq:binary-LD-double} is used; otherwise they are inverse and the
output is $O$.  When $E_0\ne0$, use
\eqref{eq:binary-LD-general-add}.  In the mixed schedule,
$B_0=X_1+X_2Z_1=(x_1+x_2)Z_1$, so $B_0=0$ is exactly the condition
$x_1=x_2$.  On this locus,
$A_0=Y_1+Y_2Z_1^2=(y_1+y_2)Z_1^2$.  Thus $A_0=0$ means that the two
points are equal and the doubling formula is used; if $A_0\ne0$, the
two points are the two distinct points above the same abscissa and are
therefore inverses, so the output is $O$.  These branches are exactly
the identity, equality, and inverse cases in
Proposition~\ref{prop:binary-affine-law}; hence the generic addition
formula alone is not complete, but the stated atlas is complete.

The binary Kummer line again uses the $x$-coordinate.  The identities
are particularly simple.

\begin{theorem}[Binary Kummer identities]
\label{thm:binary-Kummer-identities}
Let $P,Q$ be finite points of \eqref{eq:binary-Eab-shifted} with
$p=x(P)\ne q=x(Q)$.  Then
\begin{align}
  x(P+Q)x(P-Q)&=\frac{(pq+b)^2}{(p+q)^2},
  \label{eq:binary-Kummer-product}\\
  x(P+Q)+x(P-Q)&=\frac{pq}{(p+q)^2},
  \label{eq:binary-Kummer-sum}\\
  x(2P)&=\frac{(p^2+b)^2}{p^2}.
  \label{eq:binary-Kummer-double-affine}
\end{align}
\end{theorem}

\begin{proof}
Write $P=(p,r)$ and $Q=(q,s)$, and put $h=p+q\ne0$.  The slope
for $P+Q$ is $\lambda=(r+s)/h$.  Since $-Q=(q,q+s)$, the slope for
$P-Q$ is
\[
  \lambda_-=\frac{r+s+q}{h}=\lambda+\frac qh.
\]
Let $x_+=x(P+Q)$ and $x_-=x(P-Q)$.  By
\eqref{eq:binary-affine-add},
\[
  x_+=\lambda^2+\lambda+a+h,
  \qquad
  x_-=\lambda_-^2+\lambda_-+a+h.
\]
Their sum is
\begin{align*}
  x_++x_-
  &=\left(\frac qh\right)^2+\frac qh
    =\frac{q^2+q(p+q)}{h^2}
    =\frac{pq}{h^2},
\end{align*}
which is \eqref{eq:binary-Kummer-sum}.

For the product, define $N=h^2x_+$.  Using
\[
  r^2+pr=p^3+ap^2+b^2,
  \qquad
  s^2+qs=q^3+aq^2+b^2,
\]
one obtains, with every equality taken in characteristic two,
\begin{align*}
  N
  &=(r+s)^2+h(r+s)+(a+h)h^2\\
  &=pq(p+q)+ps+qr.
\end{align*}
The sum identity already proved gives $h^2x_-=N+pq$.  A fully expanded
cancellation is
\begin{align*}
  N(N+pq)+(pq+b)^2h^2
  ={}&q^2\bigl(r^2+pr+p^3+ap^2+b^2\bigr)\\
     &{}+p^2\bigl(s^2+qs+q^3+aq^2+b^2\bigr)=0.
\end{align*}
Hence $h^4x_+x_-=(pq+b)^2h^2$, and division by $h^4$ proves
\eqref{eq:binary-Kummer-product}.  The doubling identity is the first
coordinate in \eqref{eq:binary-affine-double}.
\end{proof}

For projective Kummer coordinates $(X:Z)$, the fast product law is
\begin{equation}
\begin{split}
  X_+&=Z_-(X_PX_Q+bZ_PZ_Q)^2,\\
  Z_+&=X_-(X_PZ_Q+Z_PX_Q)^2.
\end{split}
\label{eq:binary-xADD-product}
\end{equation}
Using the Karatsuba cross-product schedule, its costs are
\begin{equation}
  5\M+2\Sqr+1\Dpar,
  \qquad
  4\M+2\Sqr+1\Dpar\quad(Z_-=1).
  \label{eq:binary-xADD-product-cost}
\end{equation}
The doubling law is
\begin{equation}
  X_{2P}=(X_P^2+bZ_P^2)^2,
  \qquad
  Z_{2P}=X_P^2Z_P^2,
  \label{eq:binary-xDBL}
\end{equation}
with cost
\begin{equation}
  1\M+3\Sqr+1\Dpar.
  \label{eq:binary-xDBL-cost}
\end{equation}
Thus a binary differential-addition-and-doubling step has the upper
bounds
\begin{equation}
  6\M+5\Sqr+2\Dpar,
  \qquad
  5\M+5\Sqr+2\Dpar\quad(Z_-=1).
  \label{eq:binary-xDBLADD-cost}
\end{equation}

The supplemental sum law is obtained by homogenizing
\eqref{eq:binary-Kummer-sum}.  Put
\begin{equation}
  D_h=(X_PZ_Q+Z_PX_Q)^2,
  \qquad
  N_h=X_PX_QZ_PZ_Q.
  \label{eq:binary-xADD-sum-DN}
\end{equation}
Then
\begin{equation}
  X_+=N_hZ_-+D_hX_-,
  \qquad
  Z_+=D_hZ_-.
  \label{eq:binary-xADD-sum}
\end{equation}
Computing the two cross-products first and using their product for
$N_h$ gives the costs
\begin{equation}
  6\M+1\Sqr,
  \qquad
  4\M+1\Sqr\quad(Z_-=1).
  \label{eq:binary-xADD-sum-cost}
\end{equation}
If $Q=-P$ and $P$ is not two-torsion, then $D_h=0$ and
$N_h=x(P)^2\ne0$, so \eqref{eq:binary-xADD-sum} returns the identity.
If $P=Q$, use \eqref{eq:binary-xDBL}.  Together with the identity
branches and the fast law when its known difference is regular, these
form a complete binary differential atlas.

A full point can be recovered from two adjacent Kummer coordinates on
the regular recovery chart.  Let the fixed affine point on
\eqref{eq:binary-Eab-shifted} be \(P=(p,r)\) with \(p\ne0\), and put
\begin{equation}
  f_2(x)=x^3+ax^2+b^2.
  \label{eq:binary-f2}
\end{equation}
Assume that \(Q\) and \(P+Q\) are affine and that
\(p+q\ne0\), where \(q=x(Q)\) and \(t=x(P+Q)\).  Then
\begin{equation}
  y(Q)=
  \frac{(t+a+p+q)(p+q)^2+qr+f_2(p)+f_2(q)}{p}.
  \label{eq:binary-point-recovery}
\end{equation}
To prove the formula, multiply
$t+a+p+q=\lambda^2+\lambda$ by $(p+q)^2$, substitute
$\lambda=(r+y(Q))/(p+q)$, and replace $r^2$ and $y(Q)^2$ by the two
curve equations.  The resulting identity is
\[
  p\,y(Q)
  =(t+a+p+q)(p+q)^2+qr+f_2(p)+f_2(q),
\]
which is \eqref{eq:binary-point-recovery}.  With $p^{-1}$, $f_2(p)$,
$p^2$, and $b^2$ precomputed, a direct schedule is
\begin{equation}
  2\M+1\Sqr+\Dconst{a}+\Dconst{r}+\Dconst{p^{-1}}.
  \label{eq:binary-point-recovery-cost}
\end{equation}
Compute $q^2$, then $f_2(q)=q q^2+a q^2+b^2$, use
$(p+q)^2=p^2+q^2$, multiply this square by $t+a+p+q$, and finally
multiply the numerator by $p^{-1}$.  Thus the two general products are
$q q^2$ and the product by $(p+q)^2$; the three base-specific fixed
products are by $a$, $r$, and $p^{-1}$.  For the reciprocal
specialization \(a=0\), the first of these products disappears.  The
branch \(Q=P\), for which \(p+q=0\), is returned from the known base
point or handled by doubling.  The branch \(Q=-P\) has \(P+Q=O\) and no
affine adjacent coordinate; identity and two-torsion inputs are handled
separately.  Joint normalization of two projective Kummer coordinates
adds $1\Inv+5\M$.

For the characteristic-two reciprocal curve, insert
$a=0$ and $b=\tau\sigma/\kappa^2$ in all formulas above.  An arithmetic
output is first obtained on the shifted model as
$(X_W,\widehat Y_W)$.  Undo \eqref{eq:binary-y-shift} by setting
$Y_W=\widehat Y_W+b$, and then use the inverse of
Theorem~\ref{thm:R-char2-Cd}:
\begin{equation}
  u=\frac{Y_W}{X_W},
  \qquad
  v=\frac{X_W}{X_W+d_R},
  \qquad
  x_R=\alpha_2\frac{u+1}{u},
  \qquad
  y_R=\beta_2\frac{v+1}{v}.
  \label{eq:R-char2-arithmetic-inverse}
\end{equation}
The first step is essential: the birational map
\eqref{eq:T-char2-inverse} uses the unshifted ordinate $Y_W$, not the
L\'opez--Dahab ordinate $\widehat Y_W$.  The next two equalities are the
inverse $C_{d_R}$ map, and the last two are the inverse M\"obius
substitutions.

The rational expressions in \eqref{eq:R-char2-arithmetic-inverse} omit
four projective points.  They are supplied by the following conversion
atlas on the unshifted model \eqref{eq:R-char2-Weierstrass}:
\begin{equation}
\begin{array}{c|c}
  \text{Weierstrass point}&\text{reciprocal point}\\ \hline
  O_E&(\alpha_2,0)\\
  (0,0)&(\alpha_2,\infty)\\
  (d_R,0)&(\infty,\beta_2)\\
  (d_R,d_R)&(0,\beta_2).
\end{array}
\label{eq:R-char2-complete-boundary-atlas}
\end{equation}
To verify the table, first work on $C_{d_R}$.  The inverse functions
$u=Y_W/X_W$ and $v=X_W/(X_W+d_R)$ give
\begin{equation}
\begin{array}{c|c}
  O_E&(\infty,1)\\
  (0,0)&(\infty,0)\\
  (d_R,0)&(0,\infty)\\
  (d_R,d_R)&(1,\infty).
\end{array}
\label{eq:Cd-char2-Weierstrass-boundary-atlas}
\end{equation}
For $(X_W,Y_W)=(d_R,0)$ the inverse functions give
$u=0$ and $v=\infty$, whereas for $(X_W,Y_W)=(d_R,d_R)$ they give
$u=1$ and $v=\infty$; these are the two middle boundary assignments in
\eqref{eq:Cd-char2-Weierstrass-boundary-atlas}.  At
$O_E$, the Weierstrass functions have orders
$\operatorname{ord}_{O_E}(X_W)=-2$ and $\operatorname{ord}_{O_E}(Y_W)=-3$, so $u$ has a pole and
$v\to1$.  At $(0,0)$, the partial derivative with respect to $X_W$ is
$d_R^2\ne0$; hence $Y_W$ is a local parameter,
$X_W=Y_W^2/d_R^2+O(Y_W^3)$, so again $u$ has a pole while $v\to0$.
Applying the projective M\"obius maps
\[
  x_R=\alpha_2\frac{u+1}{u},
  \qquad
  y_R=\beta_2\frac{v+1}{v}
\]
to \eqref{eq:Cd-char2-Weierstrass-boundary-atlas} gives
\eqref{eq:R-char2-complete-boundary-atlas}.  Consequently, the dense
formulas together with this finite table constitute a complete
projective interface.  They supply full-point, Kummer, differential,
and recovery arithmetic directly for every smooth reciprocal curve
over a perfect field of characteristic two.

\chapter{The Symmetric Biquadratic QRT Envelope}
\label{ch:symmetric-QRT-envelope}

\begin{figure}[H]
\centering
\resizebox{0.94\textwidth}{!}{%
\begin{tikzpicture}[x=1cm,y=1cm,>=Latex,
  box/.style={draw,rounded corners,align=center,inner sep=4pt,font=\small,text width=3.05cm}]
  \node[box] (Cd)   at (0,  2.4) {$\mathcal C_d$};
  \node[box] (Cabd) at (0,  0.0) {$\mathcal C_{a,b,d}$};
  \node[box] (R)    at (0, -2.4) {$\mathcal R_{\tau,\sigma,\kappa}$};

  \node[box] (Q)    at (5.3, 0.0) {$\mathcal Q_{\alpha,\beta,\gamma}$};

  \node[box] (J)    at (11.0,  3.0) {Jacobi quartic};
  \node[box] (W)    at (11.0,  1.0) {rational $2$-torsion\\Weierstrass model};
  \node[box] (M)    at (11.0, -1.0) {Montgomery};
  \node[box] (E)    at (11.0, -3.0) {Edwards /\\twisted Edwards};

  \draw[->,thick]
    (Cd.east) to[bend left=10]
    node[above,font=\scriptsize,fill=white,inner sep=1pt]{slice}
    (Q.north west);
  \draw[->,thick]
    (Cabd.east) --
    node[above,font=\scriptsize,fill=white,inner sep=1pt]{envelope}
    (Q.west);
  \draw[->,thick]
    (R.east) to[bend right=10]
    node[below,font=\scriptsize,fill=white,inner sep=1pt]{reciprocal slice}
    (Q.south west);
  \draw[->,thick]
    (Q.east) to[bend left=12]
    node[above,font=\scriptsize,fill=white,inner sep=1pt]{full-point companion}
    (J.west);
  \draw[->,thick]
    (Q.east) --
    node[above,font=\scriptsize,fill=white,inner sep=1pt]{Jacobian}
    (W.west);
  \draw[->,thick]
    (W.south) --
    node[right,font=\scriptsize,fill=white,inner sep=1pt]{if $\Omega_Q$ splits}
    (M.north);
  \draw[->,thick]
    (M.south) --
    node[right,font=\scriptsize,fill=white,inner sep=1pt]{standard dictionary}
    (E.north);
  \draw[->,dashed,thick]
    (J.south east) to[bend left=14]
    node[right,font=\scriptsize,fill=white,inner sep=1pt]{two-isogeny / dense opens}
    (E.north east);
\end{tikzpicture}%
}
\caption{A three-column overview of the principal interfaces surrounding the symmetric QRT family.  The left column records distinguished factorized biquadratic families, the middle column isolates the symmetric envelope $\mathcal Q_{\alpha,\beta,\gamma}$, and the right column collects its principal companion and dictionary models.}
\label{fig:QRT-model-relations-total}
\end{figure}
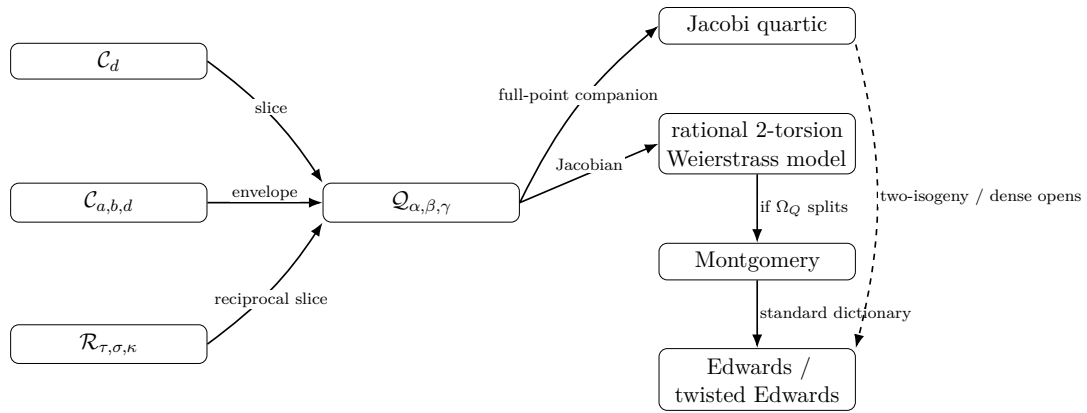

This chapter studies the symmetric biquadratic normal form that contains
both the centered \(\mathcal C_d\) slice and the reciprocal \(C\)-curve
slice.  In addition to its elliptic-curve models and arithmetic, the
chapter records the intrinsic QRT displacement before an origin is chosen
and the adjacent-Kummer-state interpretation after a pointed elliptic curve
has been fixed.  Figure~\ref{fig:QRT-model-relations-total} summarizes the
main interfaces among the families and companion models used in this chapter.

\section{Definition and root-exchange involutions}

Let $k$ be a field and let $\alpha,\beta,\gamma\in k$.  The symmetric
QRT biquadratic is
\begin{equation}
  \Q_{\alpha,\beta,\gamma}:
  \qquad
  x^2y^2+\alpha(x^2+y^2)+\beta xy+\gamma=0.
  \label{eq:symmetric-QRT}
\end{equation}
In some symmetric-biquadratic literature the same family is denoted
$C_{\alpha,\beta,\gamma}$; the symbol $\Q$ is used here to avoid a
conflict with the $C$-curve notation used for the product families.
Its completion in $\PP^1\times\PP^1$ is
\begin{equation}
\begin{split}
  F_Q(X,Z;Y,W)={}&X^2Y^2
  +\alpha(X^2W^2+Z^2Y^2)\\
  &{}+\beta XZYW+\gamma Z^2W^2=0.
  \label{eq:QRT-bihomogeneous}
\end{split}
\end{equation}
The equation is invariant under interchange of the two factors.  In
odd characteristic it is also invariant under simultaneous negation
$(x,y)\mapsto(-x,-y)$; in characteristic two this latter map is the
identity.

For fixed $x$, the equation is quadratic in $y$:
\begin{equation}
  (x^2+\alpha)y^2+\beta xy+(\alpha x^2+\gamma)=0.
  \label{eq:QRT-quadratic-y}
\end{equation}
If $y$ and $y'$ are the two roots and neither denominator vanishes,
Vieta's formula gives
\[
  yy'=\frac{\alpha x^2+\gamma}{x^2+\alpha}.
\]
The vertical and horizontal root exchanges are therefore
\begin{align}
  \jmath_y(x,y)
  &=\left(x,
  \frac{\alpha x^2+\gamma}{(x^2+\alpha)y}\right),
  \label{eq:QRT-vertical-involution}\\
  \jmath_x(x,y)
  &=\left(
  \frac{\alpha y^2+\gamma}{(y^2+\alpha)x},y\right).
  \label{eq:QRT-horizontal-involution}
\end{align}
On the regular affine locus the same two involutions have the
complementary Vieta-sum representatives
\begin{align}
  \jmath_y(x,y)
  &=\left(x,-y-\frac{\beta x}{x^2+\alpha}\right),
  \label{eq:QRT-vertical-involution-sum}\\
  \jmath_x(x,y)
  &=\left(-x-\frac{\beta y}{y^2+\alpha},y\right).
  \label{eq:QRT-horizontal-involution-sum}
\end{align}
Indeed, the sum of the two roots of
\eqref{eq:QRT-quadratic-y} is
$-\beta x/(x^2+\alpha)$, while their product is
$(\alpha x^2+\gamma)/(x^2+\alpha)$.  Hence
\eqref{eq:QRT-vertical-involution} and
\eqref{eq:QRT-vertical-involution-sum} agree on the curve wherever both
are defined.  Viewing the equation as a quadratic in \(x\) gives root sum
\(-\beta y/(y^2+\alpha)\) and root product
\((\alpha y^2+\gamma)/(y^2+\alpha)\); these are precisely
\eqref{eq:QRT-horizontal-involution} and
\eqref{eq:QRT-horizontal-involution-sum}.  The sum forms involve only one multiplication by the curve parameter
$\beta$ and do not contain $\gamma$.  The
product forms are needed as complementary representatives on fibres
where a sum denominator vanishes.

The rational expressions extend uniquely to automorphisms of every
smooth projective completion: on a fibre on which a displayed quotient
is indeterminate, the automorphism is the exchange of the two points in
the corresponding degree-two fibre.  Each exchange has order two, and
interchange of the two $\PP^1$ factors conjugates $\jmath_x$ to
$\jmath_y$.

\begin{definition}[QRT map of a biquadratic pencil]
\label{def:intrinsic-QRT-map}
Let
\[
  \mathcal P=\{F_0+\lambda F_1=0\}_{\lambda\in\PP^1}
\]
be a pencil of curves of bidegree $(2,2)$ in
$\PP^1\times\PP^1$, with no common curve component.  On a regular point
of a regular fibre, the \emph{vertical root exchange} fixes the first
projective coordinate and exchanges the two intersection points with the
corresponding vertical ruling; the \emph{horizontal root exchange} is
defined analogously.  A \emph{QRT map} is the birational self-map obtained
by composing these two exchanges in one of the two possible orders.
Equivalently, it is the fibre-preserving birational map generated by the
two degree-two deck involutions of the biquadratic pencil.
\end{definition}

This definition is intrinsic: it does not depend on an affine chart or on
a choice of elliptic-curve origin.  By construction, a QRT map preserves
the pencil parameter.  On a smooth genus-one fibre, after an origin is
chosen, the two root exchanges become reflections and their composition
becomes translation by the difference of the two degree-two projection
classes.  The corresponding displacement class in $\Pic^0$ is derived
explicitly after the boundary and projection line bundles have been introduced.

For the pencil containing \eqref{eq:symmetric-QRT}, the standard full-step
map is
\begin{equation}
  \Phi_Q=\jmath_x\circ\jmath_y.
  \label{eq:QRT-map}
\end{equation}
Because the pencil is symmetric, the coordinate exchange
$S(x,y)=(y,x)$ also produces the natural McMillan half-step
\begin{equation}
  \mathcal M_Q=S\circ\jmath_x,
  \qquad
  \mathcal M_Q(x,y)
  =\left(y,-x-\frac{\beta y}{y^2+\alpha}\right).
  \label{eq:QRT-McMillan-half-step-preview}
\end{equation}
After the elliptic translation model has been constructed, direct
composition of the two Vieta exchanges proves
$\mathcal M_Q^2=\Phi_Q^{-1}$.  Thus the terms \emph{QRT full step} and
\emph{McMillan half-step} refer to two related but distinct maps.  Both
are fibrewise elliptic translations; the half-step is the map that shifts
an adjacent-state pair by one place.  This is the geometric source of QRT
integrability \cite{QRT1988,QRT1989,DuistermaatQRT}.

A smooth member of \eqref{eq:symmetric-QRT} need not possess
a $k$-rational point over an arbitrary field.  In that case it is a
genus-one torsor and the group operations below apply to its Jacobian.
Whenever an explicit rational point is available, the curve itself is
identified with its Jacobian after that point is chosen as origin.

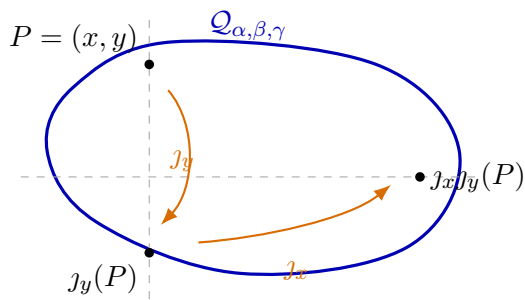
\begin{figure}[H]
\centering
\begin{tikzpicture}[x=1.0cm,y=0.85cm,every node/.style={font=\small},>=Latex]
  \draw[blue!70!black,very thick]
    plot[smooth cycle,tension=0.62] coordinates
    {(0.9,3.6) (2.2,4.3) (4.9,3.95) (6.0,2.6)
     (5.3,1.1) (3.1,0.7) (1.2,1.5) (0.55,2.6)};
  \node[blue!70!black] at (3.2,4.55) {$\mathcal Q_{\alpha,\beta,\gamma}$};
  \draw[dashed,gray!70] (1.9,0.3)--(1.9,4.7);
  \draw[dashed,gray!70] (0.2,2.2)--(6.3,2.2);
  \fill (1.9,3.95) circle (1.8pt) node[above left] {$P=(x,y)$};
  \fill (1.9,1.02) circle (1.8pt) node[below left] {$\jmath_y(P)$};
  \fill (5.48,2.2) circle (1.8pt) node[right] {$\jmath_x\jmath_y(P)$};
  \draw[->,thick,orange!85!black] (2.15,3.55) .. controls (2.55,3.0) and (2.5,2.0) .. (2.05,1.45);
  \node[orange!85!black,left] at (2.65,2.45) {$\jmath_y$};
  \draw[->,thick,orange!85!black] (2.55,1.18) .. controls (3.4,1.25) and (4.55,1.55) .. (5.1,2.08);
  \node[orange!85!black,below] at (3.85,1.05) {$\jmath_x$};
  \node at (3.35,-0.15) {$\mathcal M_Q=S\circ\jmath_x$ and $\Phi_Q=\jmath_x\circ\jmath_y$ arise from the two degree-two rulings.};
\end{tikzpicture}
\caption{The symmetric biquadratic carries two degree-two projections.  Each Vieta involution exchanges the two points in one ruling, and their compositions produce the QRT full step and the McMillan half-step.}
\label{fig:QRT-Vieta-geometry}
\end{figure}

\section[Euler--Chasles and the symmetric normal form]{The even symmetric biquadratic normal form and the Euler viewpoint}
\label{subsec:QRT-natural-normal-form}

The model \eqref{eq:symmetric-QRT} is not an arbitrary three-parameter
specialization of a general biquadratic.  It is the natural even normal
form singled out by the two most visible symmetries of a symmetric
biquadratic curve.

\begin{definition}[Even symmetric biquadratic]
\label{def:even-symmetric-biquadratic}
Assume first that $\charac(k)\ne2$.  A polynomial
$F(x,y)\in k[x,y]$ of bidegree at most $(2,2)$ is called an
\emph{even symmetric biquadratic} if
\[
  F(y,x)=F(x,y),
  \qquad
  F(-x,-y)=F(x,y).
\]
Two such equations are considered equivalent as affine equations if
one is obtained from the other by multiplication by a nonzero scalar.
In characteristic two, where simultaneous negation is trivial, the
same phrase will mean the specialization of this monomial normal form,
not a characterization by the second identity.
\end{definition}

\begin{proposition}[Invariant-theoretic normal form]
\label{prop:QRT-invariant-normal-form}
Let $\charac(k)\ne2$, and let $F$ be an even symmetric biquadratic.
Then there exist $A,B,C,D\in k$ such that
\begin{equation}
  F(x,y)=A x^2y^2+B(x^2+y^2)+Cxy+D.
  \label{eq:QRT-invariant-monomials}
\end{equation}
If $A\ne0$, division by $A$ identifies $F=0$ with a unique member
$\Q_{\alpha,\beta,\gamma}$, where
\[
  \alpha=B/A,\qquad \beta=C/A,\qquad \gamma=D/A.
\]
\end{proposition}

\begin{proof}
Write a general polynomial of bidegree at most $(2,2)$ as
$\sum_{0\le i,j\le2}a_{ij}x^iy^j$.  The identity
$F(-x,-y)=F(x,y)$ forces $a_{ij}=0$ whenever $i+j$ is odd, because
$2\ne0$ in $k$.  The surviving monomials are
\[
  1,\ x^2,\ xy,\ y^2,\ x^2y^2.
\]
The identity $F(y,x)=F(x,y)$ then forces the coefficients of $x^2$ and
$y^2$ to agree.  This gives exactly
\eqref{eq:QRT-invariant-monomials}.  If $A\ne0$, normalizing the
coefficient of $x^2y^2$ to one gives the stated parameters, and their
uniqueness follows by comparison of coefficients.
\end{proof}

The classical symmetric biquadratic is more general than
\eqref{eq:symmetric-QRT}.  Its affine equation can be written as
\begin{equation}
\begin{split}
  B(x,y)={}&a_0x^2y^2+a_1xy(x+y)+a_2(x^2+y^2)\\
  &{}+a_3xy+a_4(x+y)+a_5=0.
  \label{eq:general-symmetric-biquadratic}
\end{split}
\end{equation}
The six displayed monomials form the complete invariant space for the
coordinate interchange $(x,y)\mapsto(y,x)$ among polynomials of
bidegree at most $(2,2)$.  The terms with coefficients $a_1$ and $a_4$
are removed by the additional common involution used in
Definition~\ref{def:even-symmetric-biquadratic}; thus
\eqref{eq:symmetric-QRT} is the even symmetric normal form rather than
the most general symmetric equation.

The historical reason for the importance of
\eqref{eq:general-symmetric-biquadratic} is the Euler--Chasles
correspondence.  The following formulation is over $\mathbb C$, where
elliptic-function uniformization is available.

\begin{theorem}[Euler--Chasles parametrization]
\label{thm:Euler-Chasles-symmetric-biquadratic}
Let $C\subset\PP^1\times\PP^1$ be a smooth irreducible symmetric
biquadratic over $\mathbb C$, and assume that both projections have
degree two.  Then there are a complex torus $E=\mathbb C/\Lambda$, an
even elliptic function $f$ of degree two, and a fixed class
$\delta\in E$ such that, after fractional-linear changes of the two
projective coordinates, the normalization map can be written as
\begin{equation}
  t\longmapsto \bigl(f(t),f(t+\delta)\bigr).
  \label{eq:Euler-Chasles-parametrization}
\end{equation}
Replacing $\delta$ by $-\delta$ reverses the orientation of the
correspondence.  Degenerate biquadratics are excluded because their
normalizations need not be elliptic curves.
\end{theorem}

\begin{proof}
Let $E$ be the normalization of $C$, and denote the two degree-two
coordinate functions by $x$ and $y$.  The deck involution $\jmath_x$ of
$x:E\to\PP^1$ has a fixed point.  Choose that point as the origin of
$E$; after analytic uniformization, $\jmath_x$ is $t\mapsto-t$ and
$f(t)=x(t)$ is even.  Let $S$ be the involution induced by interchanging
the two coordinates.  The diagonal has intersection number four with
the $(2,2)$ curve and is not a component of a smooth irreducible
member, so $S$ has geometric fixed points and is also a reflection.
The composition $T=S\circ\jmath_x$ is therefore the product of two
reflections of the elliptic curve, and hence is a translation
$t\mapsto t+\delta$ for a fixed $\delta\in E$.  Since $S=T\circ\jmath_x$, one has
\[
  y(t)=x(S(t))=f(-t+\delta)=f(t-\delta).
\]
Replacing $\delta$ by $-\delta$ gives
\eqref{eq:Euler-Chasles-parametrization}.  Conversely, eliminating $t$
from two degree-two elliptic functions $f(t)$ and $f(t+\delta)$ gives a
relation of bidegree at most $(2,2)$, and symmetry follows from
$t\mapsto-t-\delta$.  This is the geometric content of the classical
Euler--Chasles theorem; algebraic and analytic treatments are given in
\cite{BuchstaberVeselov1996,IatrouRoberts2001,IatrouRoberts2002}.
\end{proof}

A standard normalization used by Baxter and in the later literature
reduces a nondegenerate symmetric biquadratic, over a splitting field
and outside the exceptional normal-form strata, to
\begin{equation}
  X^2Y^2+c_1(X^2+Y^2)+2d_1XY+1=0.
  \label{eq:Baxter-symmetric-biquadratic}
\end{equation}
This is exactly the shape of \eqref{eq:symmetric-QRT} with constant
term one.  For the present family no general M\"obius reduction is
needed to obtain this normalization.  If $\gamma\ne0$ and
$\rho^4=\gamma$ in the chosen field, the common scaling
$x=\rho X$, $y=\rho Y$ changes \eqref{eq:symmetric-QRT} into
\eqref{eq:Baxter-symmetric-biquadratic} with
\begin{equation}
  c_1=\frac{\alpha}{\rho^2},
  \qquad
  d_1=\frac{\beta}{2\rho^2}.
  \label{eq:QRT-to-Baxter-parameters}
\end{equation}
The requirement that $\rho$ and $2^{-1}$ exist is part of the
field-of-definition statement.  Over a ground field in which the
necessary fourth root does not exist, the unscaled equation
\eqref{eq:symmetric-QRT} is the corresponding twisted Baxter form.
The historical normal form and its elliptic parametrizations are
discussed in \cite{Baxter1982,IatrouRoberts2001,IatrouRoberts2002}.

For comparison with the preceding theorem, the Weierstrass addition
theorem gives a completely explicit symmetric biquadratic.  With the
convention
\[
  \wp'(u)^2=4\wp(u)^3-g_2\wp(u)-g_3,
\]
put
\[
  X=\wp(u),\qquad Y=\wp(u+\eta),\qquad Z=\wp(\eta).
\]
The Frobenius--Stickelberger form of the addition theorem gives
\begin{equation}
  \left(XY+YZ+ZX+\frac{g_2}{4}\right)^2
  =(X+Y+Z)(4XYZ-g_3).
  \label{eq:Euler-symmetric-biquadratic}
\end{equation}
Here is a direct verification of the normalization in
\eqref{eq:Euler-symmetric-biquadratic}.  The three points
\[
  (X,\wp'(u)),\qquad
  (Z,\wp'(\eta)),\qquad
  (Y,-\wp'(u+\eta))
\]
are collinear because their parameters sum to zero.  Write their common
line as $v=mt+n$.  The three abscissas $X,Y,Z$ are therefore the roots
of
\[
  4t^3-g_2t-g_3-(mt+n)^2=0.
\]
If
\[
  S_1=X+Y+Z,\qquad S_2=XY+YZ+ZX,\qquad S_3=XYZ,
\]
comparison of coefficients gives
\[
  S_1=\frac{m^2}{4},\qquad
  S_2=-\frac{2mn+g_2}{4},\qquad
  S_3=\frac{n^2+g_3}{4}.
\]
Consequently
\[
  \left(S_2+\frac{g_2}{4}\right)^2
  =\frac{m^2n^2}{4}
  =S_1(4S_3-g_3),
\]
which is precisely \eqref{eq:Euler-symmetric-biquadratic}.  For fixed
$Z$, this is a symmetric polynomial of degree two in each of $X$ and
$Y$.  Thus a translation $u\mapsto u+\eta$ becomes a root-exchange
dynamics on a symmetric biquadratic.  The even form
\eqref{eq:symmetric-QRT} is obtained after a common involution has been
normalized and, when the required fixed divisors split, after a common
fractional-linear change of coordinates.  Over a nonsplitting ground
field one obtains a twisted symmetric form rather than a literal
member of \eqref{eq:symmetric-QRT}.

\section{Why a smooth member is a genus-one curve}
\label{subsec:QRT-why-elliptic}

The genus-one property is not an accidental consequence of a quartic
substitution.  It is built into the ambient surface, the two degree-two
projections, and the anticanonical geometry that underlies every QRT
pencil.

\begin{proposition}[Anticanonical realization and canonical differential]
\label{prop:QRT-anticanonical}
Let $C\subset\PP^1\times\PP^1$ be a smooth curve of bidegree $(2,2)$.
Then
\begin{equation}
  C\sim-K_{\PP^1\times\PP^1},
  \qquad
  K_C\simeq\mathcal O_C.
  \label{eq:QRT-anticanonical}
\end{equation}
Consequently, $C$ has genus one.  If $F(x,y)=0$ is an affine equation
for $C$, the local expressions
\begin{equation}
  \omega_C=\frac{dx}{F_y}
  =-\frac{dy}{F_x}
  \label{eq:QRT-residue-differential}
\end{equation}
on the open sets $F_y\ne0$ and $F_x\ne0$ patch to a nonzero regular
differential on the projective curve.  Each degree-two root exchange
pulls $\omega_C$ back to $-\omega_C$, whereas every QRT translation
preserves $\omega_C$.
\end{proposition}

\begin{proof}
Let $H_1$ and $H_2$ be the two ruling classes.  One has
\[
  K_{\PP^1\times\PP^1}=-2H_1-2H_2,
  \qquad
  [C]=2H_1+2H_2,
\]
so $C$ is anticanonical.  Adjunction gives
\[
  K_C=(K_{\PP^1\times\PP^1}+C)|_C\simeq\mathcal O_C.
\]
Equivalently,
\[
  2p_a(C)-2
  =(2H_1+2H_2)\cdot0=0,
\]
and smoothness identifies the arithmetic and geometric genera, so
$g(C)=1$.

The Poincar\'e residue of a local meromorphic two-form with simple pole
along the anticanonical divisor gives a section of $K_C$.  In affine
coordinates this residue is represented by
$dx/F_y=-dy/F_x$, because differentiation of $F=0$ gives
$F_xdx+F_ydy=0$.  Smoothness guarantees that $F_x$ and $F_y$ do not
vanish simultaneously, so the two expressions cover the affine curve;
the anticanonical construction shows that they extend across the
boundary.  The resulting section is nonzero and, since $K_C$ has degree
zero, it has no zeros.

After choosing a geometric origin, a degree-two deck involution has the
form $P\mapsto A-P$ for a suitable point $A$.  It is a translation
followed by inversion, so it acts by $-1$ on the invariant differential.
A composition of two such reflections is a translation and therefore
acts trivially on the differential.  The statement remains valid in
characteristic two, where $-\omega_C=\omega_C$.
\end{proof}

The two rulings contain more information than the genus alone.  Their
restrictions
\begin{equation}
  \mathcal L_x=\mathcal O_C(1,0),
  \qquad
  \mathcal L_y=\mathcal O_C(0,1)
  \label{eq:QRT-two-degree-two-line-bundles}
\end{equation}
are degree-two line bundles.  They encode the intrinsic displacement of
the QRT map before any group origin has been selected.

\begin{proposition}[The projection classes determine the QRT displacement]
\label{prop:QRT-line-bundle-displacement}
Let $C$ be a smooth $(2,2)$ curve, and choose a geometric point $O$ as
origin.  Let $S_x,S_y\in C$ be the points corresponding to
$\mathcal L_x,\mathcal L_y\in\Pic^2(C)$ under the identification
\begin{equation}
  \mathcal O_C(P+Q)\longmapsto P+Q
  \label{eq:QRT-Pic2-identification}
\end{equation}
with the elliptic-curve group law.  Then the deck involutions of the two
projections satisfy
\begin{equation}
  \jmath_y(P)=S_x-P,
  \qquad
  \jmath_x(P)=S_y-P,
  \label{eq:QRT-reflection-centres}
\end{equation}
and hence
\begin{equation}
  \Phi_Q(P)=\jmath_x\jmath_y(P)
  =P+(S_y-S_x).
  \label{eq:QRT-Picard-translation}
\end{equation}
The intrinsic displacement class is therefore
\begin{equation}
  \delta_{\Phi_Q}
  =[\mathcal L_y\otimes\mathcal L_x^{-1}]
  \in\Pic^0(C),
  \label{eq:QRT-line-bundle-displacement}
\end{equation}
which is independent of the chosen origin.
\end{proposition}

\begin{proof}
A fibre of the $x$-projection is a divisor
$P+\jmath_y(P)$ in the complete linear system $|\mathcal L_x|$.
All divisors in this linear system have the same sum in the elliptic-curve
group, namely $S_x$.  Hence
$P+\jmath_y(P)=S_x$, which gives the first formula in
\eqref{eq:QRT-reflection-centres}.  A fibre of the \(y\)-projection is
the divisor \(P+\jmath_x(P)\) in \(|\mathcal L_y|\); its group sum is
\(S_y\), so \(P+\jmath_x(P)=S_y\) and
\(\jmath_x(P)=S_y-P\).  Their composition is
\[
  \jmath_x(\jmath_y(P))
  =S_y-(S_x-P)=P+(S_y-S_x),
\]
proving \eqref{eq:QRT-Picard-translation}.  Under the standard
identification of $\Pic^0(C)$ with the elliptic curve, the tensor quotient
$\mathcal L_y\otimes\mathcal L_x^{-1}$ corresponds exactly to
$S_y-S_x$.  The line bundle quotient is defined before an origin is
chosen, so the resulting displacement class is intrinsic.
\end{proof}

For the symmetric half-step of
\eqref{eq:QRT-McMillan-half-step-preview}, suppose a chosen origin makes
$\mathcal M_Q$ translation by $D$.  Since
$\mathcal M_Q^2=\Phi_Q^{-1}$, Proposition~\ref{prop:QRT-line-bundle-displacement}
gives
\begin{equation}
  [\mathcal L_y\otimes\mathcal L_x^{-1}]=[-2D].
  \label{eq:QRT-half-displacement-line-bundle}
\end{equation}
Thus the $(2,2)$ embedding intrinsically records twice the McMillan
translation.  Choosing the oriented half-step amounts to choosing a
solution of $[2]D=-\delta_{\Phi_Q}$; any two such solutions differ by a
two-torsion point, while reversing the orientation replaces $D$ by
$-D$.  This explains why the pair consisting of the curve and its QRT
dynamics contains more information than the abstract elliptic curve
alone.

The second genus calculation comes from either projection.  Equation
\eqref{eq:QRT-quadratic-y} makes the curve a degree-two cover of the
$x$-line.  In odd characteristic its discriminant is the even quartic
\[
  -4\alpha x^4+q_Qx^2-4\alpha\gamma.
\]
When this quartic has four distinct roots, Riemann--Hurwitz gives
\[
  2g-2=2(-2)+4,
\]
again yielding $g=1$.  The four branch points are the ramification divisor
of the Kummer quotient; their cross-ratio determines the elliptic modulus
and is the reason that the Jacobi-quartic and Legendre-type descriptions
appear naturally rather than through an unrelated change of variables.

There is also a surface-theoretic explanation of QRT integrability.  Two
members of a biquadratic pencil have intersection number
\begin{equation}
  (2H_1+2H_2)^2=8.
  \label{eq:QRT-eight-base-points}
\end{equation}
Thus a pencil has a base scheme of length eight, counted with
multiplicity, when its generators have no common component.  Blowing up
that base scheme, including any required infinitely near points, resolves
the pencil into a genus-one fibration on a rational surface.  The two root
exchanges and their QRT composition lift to fibre-preserving
birational automorphisms of this surface.  On every smooth fibre they are
the reflections and translations described above.  In the present
symmetric pencil, the visible involutions organize the eight base
conditions into symmetric boundary orbits; this is the global geometric
content behind the elementary Vieta formulas.

Finally, the bidegree formula
\[
  g=(m-1)(n-1)
\]
shows why $(2,2)$ is the first symmetric bidegree in
$\PP^1\times\PP^1$ that has positive genus, and why its genus is exactly
one.  Bidegrees $(1,n)$ are rational, while increasing either degree
beyond two generally leaves the elliptic category.  A smooth member is
therefore a genus-one curve for structural, not accidental, reasons.  It
becomes an elliptic curve only after a rational point is selected.  Over a
finite field such a point always exists, so every smooth
$\Q_{\alpha,\beta,\gamma}/\F_q$ becomes an elliptic curve after choosing
an origin; over a general field it may remain a nontrivial torsor under
its Jacobian; an explicit Weierstrass equation for that Jacobian is derived
after the boundary analysis.

\section{The boundary and the visible automorphisms}
\label{subsec:QRT-boundary-symmetry}

Assume first that $\charac(k)\ne2$ and that the curve is smooth.  Over a
field containing a square root $r$ of $-\alpha$, the four boundary
points are
\begin{equation}
  (\infty,r),\quad(\infty,-r),\quad
  (r,\infty),\quad(-r,\infty).
  \label{eq:QRT-four-boundary-points}
\end{equation}
Indeed, setting $Z=0$ in \eqref{eq:QRT-bihomogeneous} gives
$Y^2+\alpha W^2=0$, and setting $W=0$ gives the symmetric equation.
The corner $Z=W=0$ is not on the curve because its first term is then
$X^2Y^2\ne0$.  If $-\alpha$ is not a square in $k$, each pair in
\eqref{eq:QRT-four-boundary-points} is a rational degree-two divisor
rather than two rational points.

The manifest automorphisms in odd characteristic are
\begin{equation}
  S(x,y)=(y,x),
  \qquad
  N(x,y)=(-x,-y),
  \label{eq:QRT-S-N}
\end{equation}
together with the two root exchanges $\jmath_x$ and $\jmath_y$.
They satisfy
\begin{equation}
  S^2=N^2=\jmath_x^2=\jmath_y^2=1,
  \qquad
  S\jmath_xS=\jmath_y,
  \qquad
  N\jmath_x=\jmath_xN,
  \qquad
  N\jmath_y=\jmath_yN.
  \label{eq:QRT-basic-relations}
\end{equation}
After an origin is chosen, each root exchange is a reflection
$R\mapsto T-R$ for a suitable $T$, and a product of two such exchanges
is a translation.  This gives a geometric, rather than merely formal,
explanation of the explicit symmetry group.

In characteristic two, if $k$ is perfect and $r^2=\alpha$, the boundary
consists of the two geometric points
\begin{equation}
  (\infty,r),\qquad(r,\infty).
  \label{eq:QRT-char2-boundary-two-points}
\end{equation}
The equation on a boundary fibre has a double root, but the total curve
is smooth when $\alpha\beta\gamma\ne0$: for example, at
$(\infty,r)$ the derivative in the transverse $Z$-direction is
$\beta XYW\ne0$.  Thus inseparability of the boundary fibre must not be
confused with singularity of the curve.

\section[Odd-characteristic QRT geometry]
{Odd characteristic: product-family slices, exact smoothness, and Weierstrass model}
\label{subsec:QRT-odd-geometry}

Assume in this section that $\charac(k)\ne2$.  Define
\begin{equation}
  W_Q=2(x^2+\alpha)y+\beta x.
  \label{eq:QRT-W}
\end{equation}
Multiplying \eqref{eq:QRT-quadratic-y} by $4(x^2+\alpha)$ and completing
the square gives
\begin{equation}
  W_Q^2=-4\alpha x^4+q_Qx^2-4\alpha\gamma,
  \label{eq:QRT-even-quartic}
\end{equation}
where
\begin{equation}
  q_Q=\beta^2-4\alpha^2-4\gamma,
  \qquad
  \Omega_Q=q_Q^2-64\alpha^2\gamma.
  \label{eq:QRT-q-Omega}
\end{equation}
On the open set $x^2+\alpha\ne0$, the inverse is
\begin{equation}
  y=\frac{W_Q-\beta x}{2(x^2+\alpha)}.
  \label{eq:QRT-even-quartic-inverse}
\end{equation}

\subsection{The product family as the \texorpdfstring{$\beta=0$}{beta=0} QRT locus}

The relation with the three-parameter product family
$\C_{a,b,d}:(u^2+u+a)(v^2+v+b)=d$ is exact once the
allowed coordinate changes are stated.  Put
\begin{equation}
  \lambda=1-4a,
  \qquad
  \mu=1-4b,
  \qquad
  \delta=16d.
  \label{eq:Cabd-to-QRT-parameters}
\end{equation}
In centered coordinates $r=2u+1$ and $s=2v+1$, the product family is
\begin{equation}
  r^2s^2-\mu r^2-\lambda s^2+\lambda\mu-\delta=0.
  \label{eq:Cabd-centered-as-diagonal-biquadratic}
\end{equation}
This is an even diagonal biquadratic, but its two pure-square coefficients
need not be equal.  The following proposition identifies precisely when it
is a member of the symmetric normal form over the ground field.

\begin{proposition}[Exact product-family locus in the symmetric QRT family]
\label{prop:Cabd-as-QRT-beta-zero}
Assume $\charac(k)\ne2$ and that
$\C_{a,b,d}$ is smooth, so
$\lambda\mu\delta(\lambda\mu-\delta)\ne0$.

\begin{enumerate}[label=\textup{(\roman*)}]
  \item In the unchanged centered coordinates $(x,y)=(r,s)$,
  \eqref{eq:Cabd-centered-as-diagonal-biquadratic} is a member of
  $\Q_{\alpha,\beta,\gamma}$ if and only if
  \(\lambda=\mu\), equivalently \(a=b\).  In that case
  \begin{equation}
    \boxed{
    \C_{a,a,d}\simeq
    \Q_{4a-1,\,0,\,(1-4a)^2-16d}}
    \label{eq:Caa-as-QRT}
  \end{equation}
  through $x=2u+1$ and $y=2v+1$.

  \item More generally, a centered diagonal equivalence over $k$ to a
  member of the $\beta=0$ QRT slice exists if and only if
  \begin{equation}
    \frac{\lambda}{\mu}\in(k^\times)^2.
    \label{eq:Cabd-QRT-square-class-condition}
  \end{equation}
  If $c^2=\lambda/\mu$, then
  \begin{equation}
    x=\frac{2u+1}{c},
    \qquad
    y=2v+1
    \label{eq:Cabd-to-QRT-map}
  \end{equation}
  gives
  \begin{equation}
    \boxed{
    \C_{a,b,d}\simeq
    \Q_{-\mu,\,0,\,\mu(\lambda\mu-\delta)/\lambda}.}
    \label{eq:Cabd-to-QRT-target}
  \end{equation}

  \item Over a separable closure, condition
  \eqref{eq:Cabd-QRT-square-class-condition} always holds.  Hence every
  smooth odd-characteristic product curve is geometrically a member of
  the $\beta=0$ QRT locus.  The obstruction over $k$ is exactly the
  square class of $\lambda/\mu$ for centered diagonal equivalence.
\end{enumerate}
\end{proposition}

\begin{proof}
Equation \eqref{eq:Cabd-centered-as-diagonal-biquadratic} already has
coefficient one on $r^2s^2$ and no $rs$ term.  It has equal coefficients
on $r^2$ and $s^2$ exactly when $-\mu=-\lambda$, proving part~(i) and
\eqref{eq:Caa-as-QRT}.

For a general centered diagonal change $r=px$, $s=qy$, division by
$p^2q^2$ gives
\[
  x^2y^2-\frac{\mu}{q^2}x^2-\frac{\lambda}{p^2}y^2
  +\frac{\lambda\mu-\delta}{p^2q^2}=0.
\]
The pure-square coefficients are equal if and only if
\[
  \frac{\mu}{q^2}=\frac{\lambda}{p^2},
  \qquad\text{equivalently}\qquad
  \left(\frac pq\right)^2=\frac{\lambda}{\mu}.
\]
This proves the necessity in part~(ii).  Conversely, if
$c^2=\lambda/\mu$, take $p=c$ and $q=1$.  The resulting equation is
\[
  x^2y^2-\mu(x^2+y^2)
  +\frac{\mu(\lambda\mu-\delta)}{\lambda}=0,
\]
which is exactly \eqref{eq:Cabd-to-QRT-target}.  Part~(iii) follows
because $\lambda/\mu$ has a square root over a separable closure.
\end{proof}

\begin{remark}[Important specializations and the role of \texorpdfstring{$\beta$}{beta}]
\label{rem:Cabd-QRT-specializations}
The original family satisfies
\begin{equation}
  C_d\simeq\Q_{-1,0,1-16d}
  \label{eq:Cd-as-QRT}
\end{equation}
through $x=2u+1$, $y=2v+1$.  For the one-sided twisted curve
$\T_{a,d}=\C_{a,0,d}$, if
$c^2=1-4a$, then
\begin{equation}
  \T_{a,d}\simeq
  \Q_{-1,0,1-16d/(1-4a)}.
  \label{eq:Tad-as-QRT-beta-zero}
\end{equation}
When $1-4a$ is nonsquare, the same statement holds only after the
corresponding quadratic extension; over $k$ the one-sided curve retains
that square-class twist.

The smoothness criterion of the target in
\eqref{eq:Cabd-to-QRT-target} agrees exactly with the product-family
criterion.  Indeed, for
\[
  \alpha_Q=-\mu,
  \qquad
  \gamma_Q=\frac{\mu(\lambda\mu-\delta)}{\lambda},
\]
one has
\[
  \alpha_Q^2-\gamma_Q=\frac{\mu\delta}{\lambda},
\]
so
$\alpha_Q\gamma_Q(\alpha_Q^2-\gamma_Q)\ne0$ if and only if
$\lambda\mu\delta(\lambda\mu-\delta)\ne0$.

Thus the centered product family occupies the $\beta=0$ stratum of the
symmetric QRT geometry, up to the explicit square-class obstruction.
The parameter $\beta\ne0$ is not produced by factor-preserving centered
diagonal changes of $\C_{a,b,d}$; it is the additional mixed coupling
that makes the McMillan update
$y_{n+1}+y_{n-1}=-\beta y_n/(y_n^2+\alpha)$ nontrivial.  This statement
concerns the written biquadratic model and its marked projections, not
merely the abstract isomorphism class of its Jacobian.
\end{remark}

\begin{theorem}[Exact odd-characteristic smoothness criterion]
\label{thm:QRT-smoothness-odd}
Suppose $\charac(k)\ne2$.  The completion of
$\Q_{\alpha,\beta,\gamma}$ is a smooth genus-one curve if and only if
\begin{equation}
  \alpha\gamma\Omega_Q
  =\alpha\gamma
  \bigl((\beta+2\alpha)^2-4\gamma\bigr)
  \bigl((\beta-2\alpha)^2-4\gamma\bigr)
  \ne0.
  \label{eq:QRT-smoothness-odd}
\end{equation}
\end{theorem}

\begin{proof}
The right-hand side of \eqref{eq:QRT-even-quartic} has coefficients
\[
  A=-4\alpha,
  \qquad B=q_Q,
  \qquad C=-4\alpha\gamma.
\]
Thus $AC\ne0$ is equivalent to $\alpha\gamma\ne0$, and
$B^2-4AC=\Omega_Q$.  Under
\eqref{eq:QRT-smoothness-odd}, the quartic has degree four, nonzero
constant term, and four distinct geometric roots.  It is therefore not
a square in $\overline{k}(x)$, so the quartic double cover and the QRT
biquadratic are geometrically integral.  Its smooth
projective completion is therefore a genus-one double cover of
$\PP^1_x$.  The maps \eqref{eq:QRT-W} and
\eqref{eq:QRT-even-quartic-inverse} identify the function fields.  Since
a projective $(2,2)$ curve has arithmetic genus one, its normalization
has the same genus as its arithmetic genus; hence the total
$\delta$-invariant is zero and the QRT completion is already smooth.
The two smooth projective curves are consequently isomorphic.

For necessity, if $\gamma=0$, the affine point $(0,0)$ lies on the curve
and both affine partial derivatives vanish there.  If $\alpha=0$, let $F$ denote the bihomogeneous polynomial in
\eqref{eq:QRT-bihomogeneous}.  The point
$((X:Z),(Y:W))=((1:0),(0:1))$ lies on $F=0$.  At that point,
\[
  F_X=2XY^2+\beta ZYW=0,\quad
  F_Z=\beta XYW+2\gamma ZW^2=0,
\]
\[
  F_Y=2X^2Y+\beta XZW=0,\quad
  F_W=\beta XZY+2\gamma Z^2W=0.
\]
Thus the projective QRT curve is singular.  If $\Omega_Q=0$ with $\alpha\gamma\ne0$, the quartic has
a repeated root, so its double cover is singular and its normalization
has genus zero.  These three alternatives exhaust the failure of
\eqref{eq:QRT-smoothness-odd}.
\end{proof}

\begin{corollary}[The centered $\beta=0$ slice]
\label{cor:QRT-smoothness-beta-zero}
If $\charac(k)\ne2$ and $\beta=0$, then
\begin{equation}
  \Omega_Q=16(\alpha^2-\gamma)^2.
  \label{eq:QRT-Omega-beta-zero}
\end{equation}
Consequently, the particularly compact condition
\begin{equation}
  \alpha\gamma(\alpha^2-\gamma)\ne0
  \label{eq:QRT-smoothness-beta-zero}
\end{equation}
is necessary and sufficient for smoothness on this slice.  It is not
the smoothness criterion for the general three-parameter family when
$\beta\ne0$.
\end{corollary}

\begin{proof}
For $\beta=0$, one has
$q_Q=-4(\alpha^2+\gamma)$.  Therefore
\[
  \Omega_Q
  =16(\alpha^2+\gamma)^2-64\alpha^2\gamma
  =16(\alpha^2-\gamma)^2.
\]
Substitution in Theorem~\ref{thm:QRT-smoothness-odd} proves the claim.
\end{proof}

Lemma~\ref{lem:even-quartic-jacobian} now gives the ground-field
Jacobian
\begin{equation}
  E_Q:
  \qquad
  V^2=U\bigl(U^2-2q_QU+\Omega_Q\bigr).
  \label{eq:QRT-Jacobian}
\end{equation}
Thus the universal Weierstrass parameters of
\eqref{eq:EAB} are
\begin{equation}
  A_Q=-2q_Q,
  \qquad B_Q=\Omega_Q.
  \label{eq:QRT-AB}
\end{equation}
The discriminant and $j$-invariant are
\begin{align}
  \Delta(E_Q)&=4096\alpha^2\gamma\Omega_Q^2,
  \label{eq:QRT-Delta-odd}\\
  j(E_Q)&=
  \frac{(q_Q^2+192\alpha^2\gamma)^3}
       {\alpha^2\gamma\,(q_Q^2-64\alpha^2\gamma)^2}.
  \label{eq:QRT-j}
\end{align}
To verify these formulas, use
$c_4=16(A_Q^2-3B_Q)$ and
$\Delta=16B_Q^2(A_Q^2-4B_Q)$ for
\eqref{eq:EAB}, and substitute
$A_Q=-2q_Q$, $B_Q=q_Q^2-64\alpha^2\gamma$.

Without any square-root hypothesis, full-point arithmetic on the
Jacobian \eqref{eq:QRT-Jacobian} is already explicit.  For affine points
$P_i=(U_i,V_i)$ with $U_1\ne U_2$, put
\begin{equation}
  \lambda_Q=\frac{V_2-V_1}{U_2-U_1}.
  \label{eq:QRT-affine-lambda}
\end{equation}
Then
\begin{equation}
\begin{split}
  U(P_1+P_2)&=\lambda_Q^2+2q_Q-U_1-U_2,\\
  V(P_1+P_2)&=\lambda_Q
  \bigl(U_1-U(P_1+P_2)\bigr)-V_1.
\end{split}
\label{eq:QRT-affine-full-add}
\end{equation}
For $P_1=P_2$ with $V_1\ne0$, put
\begin{equation}
  \lambda_{Q,2}
  =\frac{3U_1^2-4q_QU_1+\Omega_Q}{2V_1}.
  \label{eq:QRT-affine-double-slope}
\end{equation}
Then
\begin{equation}
\begin{split}
  U(2P_1)&=\lambda_{Q,2}^2+2q_Q-2U_1,\\
  V(2P_1)&=\lambda_{Q,2}
  \bigl(U_1-U(2P_1)\bigr)-V_1.
\end{split}
\label{eq:QRT-affine-full-double}
\end{equation}
The identity, inverse, and two-torsion cases are those of
Proposition~\ref{prop:EAB-affine-law}.  Substitution of
$A=-2q_Q$ and $B=\Omega_Q$ in
\eqref{eq:EAB-Jacobian-add}--\eqref{eq:EAB-Jacobian-double} gives a
complete inversion-free atlas.  General, mixed, and affine-input
addition cost, respectively,
\[
  12\M+5\Sqr+1\Dpar,
  \qquad 8\M+4\Sqr+1\Dpar,
  \qquad 4\M+2\Sqr+1\Dpar,
\]
and projective and affine-input doubling cost
\[
  4\M+7\Sqr+3\Dpar,
  \qquad 2\M+4\Sqr+2\Dpar.
\]
These formulas operate over $k$ on the Jacobian for every smooth QRT
member, whether or not the QRT torsor has a $k$-rational point.

An explicit birational map is available when the constant term of the
even quartic is a square.  Suppose
\begin{equation}
  c^2=-4\alpha\gamma
  \label{eq:QRT-c-square}
\end{equation}
has a solution $c\in k^\times$.  Then
\begin{align}
  U&=q_Q+\frac{2c(W_Q+c)}{x^2},
  \label{eq:QRT-explicit-U}\\
  V&=\frac{2cU}{x}
  \label{eq:QRT-explicit-V}
\end{align}
map the even quartic to \eqref{eq:QRT-Jacobian}.  Conversely, for
$UV\ne0$,
\begin{align}
  x&=\frac{2cU}{V},
  \label{eq:QRT-explicit-inverse-x}\\
  W_Q&=\frac{(U-q_Q)x^2}{2c}-c,
  \label{eq:QRT-explicit-inverse-W}\\
  y&=\frac{W_Q-\beta x}{2(x^2+\alpha)}.
  \label{eq:QRT-explicit-inverse-y}
\end{align}
For completeness, the algebraic verification is as follows.  Substitute
\eqref{eq:QRT-explicit-U}--\eqref{eq:QRT-explicit-V} into the left-hand
side of \eqref{eq:QRT-Jacobian}; after multiplication by $x^6$, the
difference factors as
\begin{equation}
  4c^2
  \bigl(q_Qx^2+2c^2+2cW_Q\bigr)
  \bigl(-4\alpha x^4+q_Qx^2+c^2-W_Q^2\bigr).
  \label{eq:QRT-map-factorization}
\end{equation}
The second factor vanishes by \eqref{eq:QRT-even-quartic}; hence the
whole difference is zero.  Substituting the inverse formulas
recovers $x,W_Q,y$, so the maps are mutually inverse on dense opens.
They therefore extend to an isomorphism of the smooth projective
curves.  In particular, condition \eqref{eq:QRT-c-square} supplies a
$k$-rational point and turns the QRT curve itself into an elliptic
curve over $k$.

\section[Jacobi QRT arithmetic]{Extended Jacobi coordinates and full-point QRT arithmetic}
\label{subsec:QRT-extended-Jacobi}

Continue to assume \eqref{eq:QRT-c-square}.  Put
\begin{equation}
  z=\frac{W_Q}{c},
  \qquad
  \epsilon=\frac1\gamma,
  \qquad
  \delta=-\frac{q_Q}{8\alpha\gamma}.
  \label{eq:QRT-Jacobi-parameters}
\end{equation}
Dividing \eqref{eq:QRT-even-quartic} by
$c^2=-4\alpha\gamma$ gives the extended Jacobi quartic
\begin{equation}
  \JQ_{\epsilon,\delta}:
  \qquad
  z^2=\epsilon x^4+2\delta x^2+1.
  \label{eq:QRT-extended-Jacobi}
\end{equation}
Moreover,
\begin{equation}
  \epsilon(\delta^2-\epsilon)
  =\frac{\Omega_Q}{64\alpha^2\gamma^3}\ne0,
  \label{eq:QRT-Jacobi-smoothness}
\end{equation}
so the Jacobi model is smooth.  The identity $(0,1)$ corresponds to
\begin{equation}
  O_Q=\left(0,\frac{c}{2\alpha}\right)
  \label{eq:QRT-origin-odd}
\end{equation}
on \eqref{eq:symmetric-QRT}, while $(0,-1)$ is the marked two-torsion
point.  The conversion formulas are
\begin{equation}
  z=\frac{2(x^2+\alpha)y+\beta x}{c},
  \qquad
  y=\frac{cz-\beta x}{2(x^2+\alpha)}.
  \label{eq:QRT-Jacobi-conversion}
\end{equation}

\begin{proposition}[Complete QRT--Jacobi interface]
\label{prop:QRT-Jacobi-original-atlas}
Under the assumptions of this section, the two rational formulas
in \eqref{eq:QRT-Jacobi-conversion}, supplemented by the cases below,
define mutually inverse maps between the smooth projective QRT curve
and its extended Jacobi model.

Over $\overline{k}$, choose $x_0$ with $x_0^2=-\alpha$.  If
$\beta\ne0$, then
\begin{align}
  \left(x_0,\frac{\beta x_0}{c}\right)
  &\longleftrightarrow
  \left(x_0,\frac{\alpha^2-\gamma}{\beta x_0}\right),
  \notag\\
  \left(x_0,-\frac{\beta x_0}{c}\right)
  &\longleftrightarrow (x_0,\infty).
  \label{eq:QRT-Jacobi-x0-beta}
\end{align}
If $\beta=0$, then smoothness implies $\gamma\ne\alpha^2$, and the
single ramified point is
\begin{equation}
  (x_0,0)\longleftrightarrow(x_0,\infty).
  \label{eq:QRT-Jacobi-x0-beta0}
\end{equation}
Finally, if $\epsilon=e^2$ in an extension field, the two points at
Jacobi infinity satisfy
\begin{equation}
  (0:\pm e:1:0)
  \longleftrightarrow
  \left(\infty,\pm\frac{ce}{2}\right).
  \label{eq:QRT-Jacobi-original-infinity}
\end{equation}
On the regular affine locus, the forward and inverse interface costs
are, respectively,
\begin{equation}
  1\M+1\Sqr+2\Dpar,
  \qquad
  1\Inv+1\M+1\Sqr+2\Dpar.
  \label{eq:QRT-Jacobi-interface-costs}
\end{equation}
Here multiplication by $c^{-1}$, $c$, and $\beta$ is classified as
fixed-parameter multiplication.
\end{proposition}

\begin{proof}
Put $A_x=x^2+\alpha$.  The fibre equation is
\[
  A_xy^2+\beta xy+\alpha x^2+\gamma=0.
\]
At $x=x_0$, one has $A_x=0$.  If $\beta\ne0$, its finite point is
therefore
$y=(\alpha^2-\gamma)/(\beta x_0)$.  At that point
$W_Q=2A_xy+\beta x=\beta x_0$, which gives the first line of
\eqref{eq:QRT-Jacobi-x0-beta}.  The second root tends to infinity.  By
Vieta's formula its leading term satisfies
$A_xy\to-\beta x_0$, and hence
$W_Q=2A_xy+\beta x\to-\beta x_0$.  This gives the second line.

If $\beta=0$, then
\[
  \Omega_Q=16(\alpha^2-\gamma)^2.
\]
The smoothness condition $\Omega_Q\ne0$ therefore gives
$\gamma\ne\alpha^2$.  The homogeneous fibre at $x=x_0$ is
$(\gamma-\alpha^2)W^2=0$, so it consists of the single double point at
$y=\infty$.  The quartic coordinate satisfies $W_Q=0$ there, proving
\eqref{eq:QRT-Jacobi-x0-beta0}.

At a Jacobi point at infinity,
$z/x^2\to\pm e$.  Dividing the inverse formula in
\eqref{eq:QRT-Jacobi-conversion} by $x^2$ gives
\[
  y=\frac{c(z/x^2)-\beta/x}{2(1+\alpha/x^2)}
  \longrightarrow\pm\frac{ce}{2},
\]
which proves \eqref{eq:QRT-Jacobi-original-infinity}.  Conversely,
$(ce/2)^2=-\alpha$, because $c^2=-4\alpha\gamma$ and
$e^2=1/\gamma$; these are exactly the two points on the QRT boundary
$x=\infty$.

For the costs, the forward map squares $x$, multiplies
$(x^2+\alpha)$ by $y$, and performs fixed multiplications by $\beta$
and $c^{-1}$.  The inverse map uses the same square, one inversion of
$x^2+\alpha$, one general multiplication, and fixed multiplications by
$c$ and $\beta$.  The listed exceptional fibres are precisely the
zeros of $x^2+\alpha$ together with the points at infinity, so the
atlas is exhaustive.
\end{proof}

For two affine Jacobi points $(x_i,z_i)$, define
\begin{align}
  x_3&=\frac{x_1z_2+x_2z_1}
  {1-\epsilon x_1^2x_2^2},
  \label{eq:QRT-Jacobi-affine-add-x}\\
  z_3&=
  \frac{(z_1z_2+2\delta x_1x_2)
        (1+\epsilon x_1^2x_2^2)
       +2\epsilon x_1x_2(x_1^2+x_2^2)}
       {(1-\epsilon x_1^2x_2^2)^2}.
  \label{eq:QRT-Jacobi-affine-add-z}
\end{align}
These formulas are valid whenever the denominator is nonzero.  Setting
the two inputs equal gives
\begin{align}
  x_{2P}&=\frac{2xz}{1-\epsilon x^4},
  \label{eq:QRT-Jacobi-affine-double-x}\\
  z_{2P}&=
  \frac{(z^2+2\delta x^2)(1+\epsilon x^4)
        +4\epsilon x^4}
       {(1-\epsilon x^4)^2}.
  \label{eq:QRT-Jacobi-affine-double-z}
\end{align}
Equations
\eqref{eq:QRT-Jacobi-affine-add-x}--
\eqref{eq:QRT-Jacobi-affine-double-z} are the extended Jacobi group law
\cite{HisilCarterDawson,WuSong2022}.  Their correctness can
be checked without appealing to the citation: substitute
\eqref{eq:QRT-Jacobi-parameters} in the Weierstrass map
\eqref{eq:QRT-explicit-U}--\eqref{eq:QRT-explicit-V}; both sides then
reduce to the affine chord-and-tangent formulas of
Proposition~\ref{prop:EAB-affine-law}.

For inversion-free arithmetic, use extended coordinates
$(X:Y:T:Z)$ satisfying
\begin{equation}
  X^2=TZ,
  \qquad
  Y^2=\epsilon T^2+2\delta X^2+Z^2,
  \label{eq:QRT-extended-projective-equations}
\end{equation}
where an affine point is represented by $(x:z:x^2:1)$.  Put
\begin{align}
  A_0&=X_1Y_2+X_2Y_1,\notag\\
  B_0&=Z_1Z_2,\notag\\
  C_0&=T_1T_2,\notag\\
  D_-&=B_0-\epsilon C_0,
  \qquad
  D_+=B_0+\epsilon C_0,\notag\\
  P_0&=X_1X_2,\notag\\
  E_0&=Y_1Y_2+2\delta P_0,\notag\\
  X_3&=A_0D_-,
  \label{eq:QRT-projective-X}\\
  Y_3&=E_0D_+
  +2\epsilon P_0(T_1Z_2+T_2Z_1),
  \label{eq:QRT-projective-Y}\\
  T_3&=A_0^2,
  \qquad
  Z_3=D_-^2.
  \label{eq:QRT-projective-TZ}
\end{align}
To dehomogenize, write
$x_i=X_i/Z_i$, $z_i=Y_i/Z_i$, and $t_i=T_i/Z_i=x_i^2$.  Then
\[
  \frac{A_0}{Z_1Z_2}=x_1z_2+x_2z_1,
  \qquad
  \frac{D_-}{Z_1Z_2}=1-\epsilon x_1^2x_2^2,
\]
so $X_3/Z_3$ is exactly
\eqref{eq:QRT-Jacobi-affine-add-x}.  Likewise,
\[
  \frac{E_0}{Z_1Z_2}=z_1z_2+2\delta x_1x_2,
  \quad
  \frac{D_+}{Z_1Z_2}=1+\epsilon x_1^2x_2^2,
\]
and
\[
  \frac{T_1Z_2+T_2Z_1}{Z_1Z_2}=x_1^2+x_2^2.
\]
Division of $Y_3$ by $Z_3$ therefore gives
\eqref{eq:QRT-Jacobi-affine-add-z}.  Finally,
$X_3^2=A_0^2D_-^2=T_3Z_3$.  Since the affine output satisfies
\eqref{eq:QRT-extended-Jacobi}, homogenization by $Z_3^2$ proves the
second relation in \eqref{eq:QRT-extended-projective-equations}.
A straight schedule gives
\begin{equation}
\begin{array}{c|c}
  \text{operation}&\text{cost}\\ \hline
  \text{projective addition}&11\M+2\Sqr+3\Dpar\\
  \text{mixed addition }(T_2=x_2^2,Z_2=1)&9\M+2\Sqr+3\Dpar\\
  \text{doubling}&5\M+6\Sqr+3\Dpar.
\end{array}
\label{eq:QRT-projective-cost-table}
\end{equation}
The three parameter products are by $\epsilon$, $\delta$, and
$\epsilon$ in the final cross term.  Caching a fixed affine input can
reduce the mixed count further, but
\eqref{eq:QRT-projective-cost-table} is a uniform upper bound.

The QRT point-addition formula in the original coordinates is now
fully explicit.  Convert each input by the regular formula
\eqref{eq:QRT-Jacobi-conversion} or, on an exceptional fibre, by
Proposition~\ref{prop:QRT-Jacobi-original-atlas}; apply either the
affine formulas \eqref{eq:QRT-Jacobi-affine-add-x}--
\eqref{eq:QRT-Jacobi-affine-add-z} or the projective formulas
\eqref{eq:QRT-projective-X}--\eqref{eq:QRT-projective-TZ}; and convert
the output back by the same complete interface.  The same procedure
with equal inputs gives point doubling.  The affine interface uses an
inversion in the reverse direction; for repeated arithmetic one stores
points in Jacobi or Weierstrass coordinates and converts only at the
external interface.

The Weierstrass conversion also has a complete boundary atlas.  The
two affine Jacobi points above $x=0$ extend as
\begin{equation}
  (0,1)\longmapsto O,
  \qquad
  (0,-1)\longmapsto(0,0).
  \label{eq:QRT-Jacobi-special-zero}
\end{equation}
If $\epsilon=e^2$ with $e\in k^\times$, the two rational points at
infinity in the extended coordinates are
\[
  \mathcal P_\pm^\infty=(0:\pm e:1:0).
\]
They map to the remaining two-torsion points
\begin{equation}
  \mathcal P_\pm^\infty
  \longmapsto
  \bigl(q_Q\pm2c^2e,0\bigr).
  \label{eq:QRT-Jacobi-special-infinity}
\end{equation}
Indeed, near $x=0$ the two branches of
\eqref{eq:QRT-extended-Jacobi} are
$z=1+\delta x^2+O(x^4)$ and
$z=-1-\delta x^2+O(x^4)$.  Substitution in
\eqref{eq:QRT-explicit-U} sends the first branch to the point at
infinity and the second to $U=0$.  At
$\mathcal P_\pm^\infty$ one has $z/x^2\to\pm e$, and hence
$U\to q_Q\pm2c^2e$.  Finally,
\[
  (2c^2e)^2=64\alpha^2\gamma=q_Q^2-\Omega_Q,
\]
so these two limits are exactly the roots of
$U^2-2q_QU+\Omega_Q$.  Equations
\eqref{eq:QRT-Jacobi-special-zero} and
\eqref{eq:QRT-Jacobi-special-infinity}, together with the dense-open
maps, remove every ambiguity at the conversion boundary.

There is one important completeness criterion.  Let $k=\F_q$ have odd
characteristic.  If $\epsilon$ is a nonsquare, then
$1-\epsilon x_1^2x_2^2$ cannot vanish for $k$-rational affine inputs,
because such a vanishing would express $\epsilon$ as the square
$(x_1x_2)^{-2}$.  Moreover, the two points at infinity of
\eqref{eq:QRT-extended-Jacobi} are not $k$-rational when $\epsilon$ is a
nonsquare.  Hence
\eqref{eq:QRT-projective-X}--\eqref{eq:QRT-projective-TZ} is a single
$k$-complete addition law in this case
\cite{HisilCarterDawson}.  Since $\epsilon=1/\gamma$, this criterion is
equivalent to $\gamma$ being a nonsquare.  When $\epsilon$ is a square,
a complete algorithmic atlas is obtained by using the Jacobi law on its
regular domain and, for exceptional pairs, the conversion atlas
\eqref{eq:QRT-Jacobi-special-zero}--
\eqref{eq:QRT-Jacobi-special-infinity} followed by the exhaustive
Weierstrass atlas on \eqref{eq:QRT-Jacobian}.  This distinction
prevents the unified Jacobi formulas from being incorrectly described
as geometrically complete in every parameter regime.

The preceding affine and projective formulas already determine the group
law on the QRT model once an origin has been selected.  The following
proposition records this conclusion without repeating the derivation of the
QRT--Jacobi interface.

\begin{proposition}[Transported full-point group law on the pointed QRT curve]
\label{prop:QRT-full-point-group-law-summary}
\label{prop:QRT-transported-affine-group-law}
Let \(k\) be a field with \(\operatorname{char}(k)\ne2\), and let
\[
  \mathcal Q_{\alpha,\beta,\gamma}:
  \qquad
  x^2y^2+\alpha(x^2+y^2)+\beta xy+\gamma=0
\]
be smooth.  Assume that \(c\in k^\times\) satisfies
\[
  c^2=-4\alpha\gamma.
\]
Define
\[
  q_Q=\beta^2-4\alpha^2-4\gamma,
  \qquad
  \epsilon=\frac1\gamma,
  \qquad
  \delta=-\frac{q_Q}{8\alpha\gamma},
\]
and let
\[
  \mathcal J_{\epsilon,\delta}:
  \qquad
  z^2=\epsilon x^4+2\delta x^2+1
\]
be the associated Jacobi quartic, with identity
\[
  O_J=(0,1).
\]

The affine QRT--Jacobi correspondence
\[
  \Phi_c:
  \mathcal Q_{\alpha,\beta,\gamma}
  \dashrightarrow
  \mathcal J_{\epsilon,\delta},
  \qquad
  (x,y)\longmapsto
  \left(
    x,\,
    \frac{2(x^2+\alpha)y+\beta x}{c}
  \right)
\]
extends to an isomorphism of the smooth projective completions.  Its affine
inverse on the chart \(x^2+\alpha\ne0\) is
\[
  \Psi_c(x,z)
  =
  \left(
    x,\,
    \frac{cz-\beta x}{2(x^2+\alpha)}
  \right).
\]

In particular,
\[
  O_Q:=\Psi_c(O_J)
  =
  \left(0,\frac{c}{2\alpha}\right)
\]
is a \(k\)-rational point of
\(\mathcal Q_{\alpha,\beta,\gamma}\).  Equip the smooth projective
completion of \(\mathcal Q_{\alpha,\beta,\gamma}\) with \(O_Q\) as its
identity.

Then the full-point addition law developed above is precisely the
transport of the Jacobi group law through \(\Phi_c\).  Explicitly, for
\(P_1,P_2\in\mathcal Q_{\alpha,\beta,\gamma}(k)\),
\begin{equation}
  P_1\oplus_Q P_2
  =
  \Psi_c\!\left(
    \Phi_c(P_1)\oplus_J\Phi_c(P_2)
  \right),
  \label{eq:QRT-transported-full-point-law}
\end{equation}
where \(\oplus_J\) is the group law on
\(\mathcal J_{\epsilon,\delta}\) with identity \(O_J\).

More explicitly, if
\[
  P_i=(x_i,y_i),
  \qquad
  z_i=
  \frac{2(x_i^2+\alpha)y_i+\beta x_i}{c},
  \qquad i=1,2,
\]
and the Jacobi addition formulas give
\[
  (x_3,z_3)
  =
  (x_1,z_1)\oplus_J(x_2,z_2),
\]
then, whenever \(x_3^2+\alpha\ne0\),
\begin{equation}
  P_1\oplus_Q P_2
  =
  \left(
    x_3,\,
    \frac{cz_3-\beta x_3}{2(x_3^2+\alpha)}
  \right).
  \label{eq:QRT-affine-full-point-sum}
\end{equation}
When \(x_3^2+\alpha=0\), the affine expression in
\eqref{eq:QRT-affine-full-point-sum} is replaced by the corresponding
boundary chart of the complete QRT--Jacobi interface.  Thus the
projective addition atlas established above covers every pair of input
points.
\end{proposition}

\begin{proof}
Proposition~\ref{prop:QRT-Jacobi-original-atlas} gives mutually inverse
rational maps
\[
  \Phi_c:
  \mathcal Q_{\alpha,\beta,\gamma}
  \dashrightarrow
  \mathcal J_{\epsilon,\delta}
\]
and
\[
  \Psi_c:
  \mathcal J_{\epsilon,\delta}
  \dashrightarrow
  \mathcal Q_{\alpha,\beta,\gamma}.
\]
Because both projective completions are smooth projective curves, these
birational maps extend uniquely to mutually inverse projective
isomorphisms.

Substituting \(O_J=(0,1)\) into the inverse affine formula gives
\[
  \Psi_c(O_J)
  =
  \left(
    0,\,
    \frac{c}{2\alpha}
  \right)
  =
  O_Q.
\]
Hence \(\Phi_c\) is an isomorphism of pointed curves
\[
  \Phi_c:
  \bigl(\mathcal Q_{\alpha,\beta,\gamma},O_Q\bigr)
  \xrightarrow{\;\sim\;}
  \bigl(\mathcal J_{\epsilon,\delta},O_J\bigr).
\]

Transporting the Jacobi group operation through this pointed
isomorphism defines
\[
  P_1\oplus_QP_2
  =
  \Psi_c\!\left(
    \Phi_c(P_1)\oplus_J\Phi_c(P_2)
  \right).
\]
Since \(\Phi_c\) and \(\Psi_c\) are mutually inverse, this transported
operation inherits associativity, commutativity, inverses, and the
identity \(O_Q\) from the Jacobi group law.  It is therefore the
elliptic-curve group law on the pointed QRT curve.

If the Jacobi sum is represented by an affine point \((x_3,z_3)\) with
\(x_3^2+\alpha\ne0\), applying the affine inverse \(\Psi_c\) gives
\eqref{eq:QRT-affine-full-point-sum}.  If \(x_3^2+\alpha=0\), only that
particular affine inverse chart is undefined.  The projective
isomorphism remains defined, and the complementary boundary chart of
the complete interface gives the same group sum.
\end{proof}

\begin{remark}[Dependence on the chosen origin]
\label{rem:QRT-full-point-law-origin}
The preceding proposition concerns the pointed curve
\[
  \left(
    \mathcal Q_{\alpha,\beta,\gamma},
    O_Q
  \right),
  \qquad
  O_Q=\left(0,\frac{c}{2\alpha}\right).
\]
The unpointed genus-one curve
\(\mathcal Q_{\alpha,\beta,\gamma}\) does not by itself determine a
group identity.  Replacing \(c\) by \(-c\) replaces \(O_Q\) by
\[
  \left(0,-\frac{c}{2\alpha}\right)
\]
and therefore gives the group law corresponding to the other natural
choice of origin above \(x=0\).

If \(c\notin k\), the same construction is defined over
\(k(c)\).  Over the original field \(k\), the unpointed QRT curve and
its QRT dynamics remain defined, but the displayed pointed Jacobi
identification need not descend to \(k\).
\end{remark}

\section[Jacobi, Montgomery, and Edwards links]
{Direct links with Jacobi quartics, Montgomery curves, and Edwards curves}
\label{subsec:QRT-Jacobi-Edwards-links}

The QRT model is related to a Jacobi quartic by the explicit birational
construction below.  The coordinate
\[
  W_Q=2(x^2+\alpha)y+\beta x
\]
removes the linear term of the quadratic in $y$ and gives
\[
  W_Q^2=-4\alpha x^4+q_Qx^2-4\alpha\gamma.
\]
When $c^2=-4\alpha\gamma$, division by $c^2$ produces
\[
  z^2=\epsilon x^4+2\delta x^2+1,
  \qquad
  \epsilon=\gamma^{-1},
  \qquad
  \delta=-\frac{q_Q}{8\alpha\gamma},
\]
which is exactly the extended Jacobi quartic
\eqref{eq:QRT-extended-Jacobi}.  Thus the QRT $x$-projection is the
standard degree-two Jacobi projection and the four branch points are
the roots of the even quartic.  The inequality
\[
  \epsilon(\delta^2-\epsilon)\ne0
\]
is equivalent to the QRT smoothness condition.

There are two different links with Edwards curves, and they should not
be conflated.

\subsection{A degree-two Jacobi--Edwards isogeny.}
Suppose $\epsilon=h^2$ with $h\in k^\times$, and put
\begin{equation}
  d_E=\frac{h-\delta}{2}.
  \label{eq:QRT-Jacobi-Edwards-d}
\end{equation}
Then
\begin{equation}
  \phi_{J,E}(x,z)
  =\left(
     \frac{2x}{1+hx^2},
     \frac{1-hx^2}{z}
   \right)
  \label{eq:QRT-Jacobi-Edwards-map}
\end{equation}
lands on the twisted Edwards curve
\begin{equation}
  \operatorname{TE}_{h,d_E}:
  \qquad
  hX^2+Y^2=1+d_EX^2Y^2.
  \label{eq:QRT-Jacobi-Edwards-target}
\end{equation}
To verify the equation, put $D=1+hx^2$.  After multiplying the
difference between the two sides of
\eqref{eq:QRT-Jacobi-Edwards-target} by $D^2z^2$, substitution of
$z^2=h^2x^4+2\delta x^2+1$ makes all coefficients cancel.  The map is
invariant under
\begin{equation}
  \iota_h(x,z)
  =\left(\frac1{hx},-\frac{z}{hx^2}\right),
  \label{eq:QRT-Jacobi-Edwards-deck}
\end{equation}
which is a nontrivial involution of the Jacobi quartic.  Conversely,
$X=2x/(1+hx^2)$ gives the quadratic equation
\[
  hXx^2-2x+X=0
\]
for $x$ over $k(X,Y)$, and its two roots are interchanged by
\eqref{eq:QRT-Jacobi-Edwards-deck}.  Hence the map has degree two.  It
is separable because $2hXx-2$ is not the zero rational function.  On a
smooth member the involution has no geometric fixed point: a fixed
point would satisfy $hx^2=1$ and $z=0$, which would force
$\delta=-h$ and hence $\delta^2=\epsilon$, contrary to smoothness.
Thus the completed map is an unramified separable double cover between
genus-one curves.  With compatible origins it is a degree-two isogeny,
not a birational equivalence in general.  The distinction matters when
transferring point counts, completeness statements, or operation
costs.  This isogeny is the direct form of the standard connection
between Jacobi quartics and twisted Edwards curves
\cite{WuSong2022}.

\subsection{A birational Montgomery--Edwards chain.}
Suppose in addition that
\begin{equation}
  s^2=\Omega_Q
  \label{eq:QRT-Omega-square-for-Montgomery}
\end{equation}
has a solution $s\in k^\times$.  Put
\[
  U=sX_M,
  \qquad
  V=sY_M.
\]
Then \eqref{eq:QRT-Jacobian} becomes the Montgomery equation
\begin{equation}
  B_MY_M^2=X_M^3+A_MX_M^2+X_M,
  \qquad
  A_M=-\frac{2q_Q}{s},
  \qquad
  B_M=\frac1s.
  \label{eq:QRT-Montgomery-model}
\end{equation}
The standard Montgomery--Edwards transformation gives the twisted
Edwards parameters
\begin{equation}
  a_E=\frac{A_M+2}{B_M}=2(s-q_Q),
  \qquad
  d_E'=\frac{A_M-2}{B_M}=-2(q_Q+s),
  \label{eq:QRT-Montgomery-Edwards-parameters}
\end{equation}
and the birational coordinates
\begin{equation}
  X_E=\frac{U}{V},
  \qquad
  Y_E=\frac{U-s}{U+s}.
  \label{eq:QRT-Montgomery-Edwards-map}
\end{equation}
This chain is an isomorphism on dense opens, unlike
\eqref{eq:QRT-Jacobi-Edwards-map}.  Its additional square condition is
also different: it requires the rational two-torsion cubic to split in
the manner needed for a Montgomery parameter.

Consequently, the relations among the four models may be summarized as
\[
  \Q_{\alpha,\beta,\gamma}
  \dashrightarrow \JQ_{\epsilon,\delta}
  \xrightarrow{\;2:1\;} \operatorname{TE}_{h,d_E},
\]
and, when \eqref{eq:QRT-Omega-square-for-Montgomery} holds,
\[
  \Q_{\alpha,\beta,\gamma}
  \dashrightarrow E_Q
  \dashrightarrow M_{A_M,B_M}
  \dashrightarrow \operatorname{TE}_{a_E,d_E'}.
\]
The first arrow in each line is an isomorphism only after the stated
rational-point or square hypotheses have been checked; the middle
Jacobi--Edwards arrow in the first line is a genuine two-isogeny.

\section{The McMillan map as an elliptic translation}
\label{subsec:QRT-McMillan-translation}

The same symmetric biquadratic is also invariant under a second-order birational recurrence, which gives the QRT presentation its dynamical interpretation.

Define
\begin{equation}
  \mathcal M_Q=S\circ\jmath_x.
  \label{eq:QRT-McMillan-definition}
\end{equation}
On the regular affine locus, Vieta's formula for the two $x$-roots in a
fixed horizontal fibre gives
\begin{equation}
  \mathcal M_Q(x,y)
  =\left(y,-x-\frac{\beta y}{y^2+\alpha}\right).
  \label{eq:QRT-McMillan-sum-form}
\end{equation}
The equivalent product form is
\begin{equation}
  \mathcal M_Q(x,y)
  =\left(y,
  \frac{\alpha y^2+\gamma}{(y^2+\alpha)x}\right).
  \label{eq:QRT-McMillan-product-form}
\end{equation}
Each formula has a base locus, but the two expressions are equal on the
curve wherever both are defined and extend to the same automorphism of
the smooth projective completion.

If $(x_n,x_{n+1})$ denotes the $n$th point of an orbit, then
\eqref{eq:QRT-McMillan-sum-form} is the McMillan recurrence
\cite{McMillan1971}
\begin{equation}
  x_{n+2}+x_n
  =-\frac{\beta x_{n+1}}{x_{n+1}^2+\alpha}.
  \label{eq:QRT-McMillan-recurrence}
\end{equation}
The biquadratic is a first integral:
\begin{equation}
  x_n^2x_{n+1}^2
  +\alpha(x_n^2+x_{n+1}^2)
  +\beta x_nx_{n+1}+\gamma=0
  \label{eq:QRT-McMillan-invariant}
\end{equation}
for every $n$ for which the orbit is defined.  This follows either by
substitution into \eqref{eq:QRT-McMillan-recurrence}, or more
conceptually because $\mathcal M_Q$ is the composition of a root
exchange with the coordinate interchange.  The inverse map is
\begin{equation}
  \mathcal M_Q^{-1}(x,y)
  =\left(-y-\frac{\beta x}{x^2+\alpha},x\right).
  \label{eq:QRT-McMillan-inverse}
\end{equation}

The parameter $\gamma$ selects the level of the invariant pencil.  Define
\begin{equation}
  H_{\alpha,\beta}(x,y)
  =x^2y^2+\alpha(x^2+y^2)+\beta xy.
  \label{eq:QRT-McMillan-Hamiltonian}
\end{equation}
The map $\mathcal M_Q$ depends only on $\alpha$ and $\beta$, while
$\Q_{\alpha,\beta,\gamma}$ is the level
$H_{\alpha,\beta}=-\gamma$.  This independence must not be read as
saying that the elliptic translating point is independent of
$\gamma$.  The Jacobian, the selected fibre, and the intrinsic
displacement class $\delta_{\mathcal M_Q}$ all vary with the level;
only the rational expression defining the shift is common to the whole
pencil.

\begin{proposition}[Invariant QRT pencil and symplectic form]
\label{prop:QRT-McMillan-pencil-symplectic}
On its regular affine locus, the McMillan map satisfies
\begin{equation}
  H_{\alpha,\beta}\circ\mathcal M_Q
  =H_{\alpha,\beta}
  \label{eq:QRT-McMillan-H-invariance}
\end{equation}
and
\begin{equation}
  \mathcal M_Q^*(dx\wedge dy)=dx\wedge dy.
  \label{eq:QRT-McMillan-symplectic}
\end{equation}
Thus the same birational map preserves every member of the pencil
$\{\Q_{\alpha,\beta,\gamma}\}_{\gamma}$ on its regular locus.

If $\charac(k)\ne2$ and $\alpha\ne0$, the finite singular members of
this pencil occur exactly at
\begin{equation}
  \gamma=0,
  \qquad
  \gamma=\frac{(\beta+2\alpha)^2}{4},
  \qquad
  \gamma=\frac{(\beta-2\alpha)^2}{4}.
  \label{eq:QRT-McMillan-singular-levels}
\end{equation}
The member at $\gamma=\infty$ is $Z^2W^2=0$.  The four boundary points
in \eqref{eq:QRT-four-boundary-points} are the common base points of
the pencil, each counted with intersection multiplicity two.
\end{proposition}

\begin{proof}
Fix $(x,y)$ and let $z$ be the second coordinate of
$\mathcal M_Q(x,y)$.  For the level
$h=H_{\alpha,\beta}(x,y)$, the two possible first coordinates in the
horizontal fibre through $y$ are the roots of
\[
  (y^2+\alpha)T^2+\beta yT+\alpha y^2-h=0.
\]
Their sum is $-\beta y/(y^2+\alpha)$.  Since one root is $x$, the other
is precisely
$z=-x-\beta y/(y^2+\alpha)$.  Hence
$H_{\alpha,\beta}(y,z)=h$, proving
\eqref{eq:QRT-McMillan-H-invariance} without assuming in advance that
$(x,y)$ lies on a particular member.

Write $\mathcal M_Q(x,y)=(y,z)$.  Its Jacobian matrix has the form
\[
  \begin{pmatrix}
    0&1\\
    -1&\partial z/\partial y
  \end{pmatrix},
\]
whose determinant is one.  This proves
\eqref{eq:QRT-McMillan-symplectic}.

The finite singular levels follow by factoring the exact smoothness
discriminant as in \eqref{eq:QRT-smoothness-odd}.  Two different
members of the pencil differ by a nonzero multiple of $Z^2W^2$;
hence every common point lies on $Z=0$ or $W=0$.  The equations on
those two boundary divisors give the four points
\eqref{eq:QRT-four-boundary-points}.  Because the pencil parameter
occurs through $Z^2W^2$, each is a double base point in the
intersection-theoretic count.  Finally, the parameter value at infinity
is represented by $Z^2W^2=0$ itself.
\end{proof}

Moreover, since $S\jmath_xS=\jmath_y$,
\begin{equation}
  \mathcal M_Q^2=\jmath_y\jmath_x
  =(\jmath_x\jmath_y)^{-1}=\Phi_Q^{-1}.
  \label{eq:QRT-McMillan-square}
\end{equation}
Thus the one-step McMillan map is a natural square root, in the group of
curve automorphisms, of the inverse of the standard two-involution QRT
map.

The next theorem identifies its translation point explicitly.  The
formula is useful both for dynamics and for scalar multiplication.
Figure~\ref{fig:QRT-McMillan-translation} records the conceptual meaning:
a McMillan step is the adjacent-state shadow of the elliptic translation
\(P\mapsto P+D\).

\begin{figure}[H]
\centering
\begin{tikzcd}[column sep=huge,row sep=large]
P \arrow[r,"+D"] \arrow[d,"\mathcal S_D"']
  & P+D \arrow[d,"\mathcal S_D"] \\
\mathcal S_D(P)
  \arrow[r,"T_D"]
  & \mathcal S_D(P+D)
\end{tikzcd}
\[
\mathcal S_D(P)=\bigl(\kappa(P),\kappa(P+D)\bigr),
\qquad
\mathcal M_Q(x,y)=\left(y,-x-\frac{\beta y}{y^2+\alpha}\right).
\]
\caption{The McMillan map is the state-level realization of elliptic translation by the marked displacement $D$.  In symmetric biquadratic coordinates, the map $T_D$ is represented by the McMillan update $\mathcal M_Q$.}
\label{fig:QRT-McMillan-translation}
\end{figure}

\begin{theorem}[Explicit McMillan translation point]
\label{thm:QRT-McMillan-point}
Assume $\charac(k)\ne2$ and the smoothness condition
\eqref{eq:QRT-smoothness-odd}.  Work first over a field extension
$L/k$ containing an element $c$ with $c^2=-4\alpha\gamma$, and choose
\begin{equation}
  O_Q=\left(0,\frac{c}{2\alpha}\right)
  \label{eq:QRT-McMillan-origin}
\end{equation}
as the origin over $L$.  When $c\in k$, this is a $k$-rational
elliptic origin; when $c\notin k$, the calculation below still
identifies the induced translation on the ground-field Jacobian.  Put
\begin{equation}
  D_+=(\beta+2\alpha)^2-4\gamma,
  \qquad
  D_-=(\beta-2\alpha)^2-4\gamma.
  \label{eq:QRT-Dpm}
\end{equation}
Then
\begin{equation}
  \Omega_Q=D_+D_-.
  \label{eq:QRT-Omega-factorized}
\end{equation}
The automorphism $\mathcal M_Q$ is translation by the point $P_M$ whose
coordinates in the three models are
\begin{align}
  P_M^{Q}
  &=\left(
     \frac{c}{2\alpha},
     -\frac{\beta c}{2(\alpha^2-\gamma)}
     \right),
  \label{eq:QRT-McMillan-point-Q}\\
  P_M^{J}
  &=\left(\frac{c}{2\alpha},-\frac{\beta}{2\alpha}\right),
  \label{eq:QRT-McMillan-point-J}\\
  P_M^{E}
  &=(D_-,4\alpha D_-)
  \label{eq:QRT-McMillan-point-E}
\end{align}
whenever the displayed affine coordinates are finite.  The Weierstrass point in
\eqref{eq:QRT-McMillan-point-E} is $k$-rational and independent of the
choice of the square root $c$.  Hence the induced translation on the
Jacobian is defined over the ground field even when the displayed QRT
origin exists only over $L$.

Its double is
\begin{equation}
  2P_M^E
  =\bigl(\beta^2,-4\beta(\alpha^2-\gamma)\bigr).
  \label{eq:QRT-McMillan-double}
\end{equation}
In particular, under the smoothness hypothesis, $P_M$ has exact order
four if and only if
\begin{equation}
  \beta(\alpha^2-\gamma)=0.
  \label{eq:QRT-McMillan-order-four}
\end{equation}
\end{theorem}

\begin{proof}
First expand the product in \eqref{eq:QRT-Dpm}:
\begin{align*}
  D_+D_-
  &=(\beta^2+4\alpha^2-4\gamma)^2-16\alpha^2\beta^2\\
  &=(\beta^2-4\alpha^2-4\gamma)^2-64\alpha^2\gamma
   =\Omega_Q.
\end{align*}
Thus both $D_+$ and $D_-$ are nonzero on a smooth member.

By definition $P_M^Q=\mathcal M_Q(O_Q)$.  Substitution of
$O_Q$ in \eqref{eq:QRT-McMillan-sum-form} gives
\[
  x(P_M^Q)=\frac{c}{2\alpha}.
\]
Since
\[
  \left(\frac{c}{2\alpha}\right)^2+\alpha
  =-\frac{\gamma}{\alpha}+\alpha
  =\frac{\alpha^2-\gamma}{\alpha},
\]
the second coordinate is the one in
\eqref{eq:QRT-McMillan-point-Q}; when $\alpha^2=\gamma$ the image is the
corresponding boundary point and the projective interpretation is used.

The Jacobi ordinate is $z=W_Q/c$.  Substituting
\eqref{eq:QRT-McMillan-point-Q} into
$W_Q=2(x^2+\alpha)y+\beta x$ gives
$z=-\beta/(2\alpha)$, proving
\eqref{eq:QRT-McMillan-point-J}.  Applying
\eqref{eq:QRT-explicit-U}--\eqref{eq:QRT-explicit-V} and using
$c^2=-4\alpha\gamma$ gives
\[
  U(P_M)=D_-,\qquad V(P_M)=4\alpha D_-,
\]
which proves \eqref{eq:QRT-McMillan-point-E}.  Since the latter
coordinates belong to $k$, the Jacobian translation descends even if
the chosen origin is defined only after a quadratic extension.

On $E_Q$, the tangent slope at $P_M^E$ is
\[
  \lambda
  =\frac{3D_-^2-4q_QD_-+\Omega_Q}{8\alpha D_-}
  =4\alpha-\beta.
\]
The last equality follows after substituting
$\Omega_Q=D_+D_-$ and expanding $D_\pm$.  The affine doubling law then
gives
\[
  U(2P_M)=\lambda^2+2q_Q-2D_-=\beta^2
\]
and
\[
  V(2P_M)=\lambda(D_--\beta^2)-4\alpha D_-
  =-4\beta(\alpha^2-\gamma),
\]
which is \eqref{eq:QRT-McMillan-double}.

A finite point on an odd-characteristic Weierstrass curve is
$2$-torsion exactly when its ordinate is zero.  Hence $P_M$ has order
four exactly when the second coordinate in
\eqref{eq:QRT-McMillan-double} vanishes.  If $\beta=0$, the double is
$(0,0)$; if $\alpha^2=\gamma$, smoothness forces $\beta\ne0$, and the
double is another nonzero two-torsion point.  The two alternatives are
therefore sufficient.  Conversely, order four forces the ordinate of
$2P_M$ to vanish, proving necessity.
\end{proof}

\begin{corollary}[Iteration is addition]
\label{cor:QRT-iteration-addition}
Let $\psi_Q$ denote any elliptic identification of the smooth QRT curve
with its Jacobian for which $O_Q$ maps to the identity.  Then
\begin{equation}
  \psi_Q\bigl(\mathcal M_Q^n(R)\bigr)
  =\psi_Q(R)+[n]P_M
  \label{eq:QRT-iteration-is-addition}
\end{equation}
for every integer $n$.  The McMillan map has finite order $m$ if and
only if $P_M$ is an $m$-torsion point; equivalently,
$[m]P_M=O$ and $[r]P_M\ne O$ for $0<r<m$.
\end{corollary}

\begin{proof}
The automorphism $\mathcal M_Q=S\circ\jmath_x$ sends the origin to
$P_M$.  Over $\overline{k}$, the root exchange $\jmath_x$ has fixed
points at the ramification points of the degree-two $y$-projection and
therefore induces $[-1]$ on $\Pic^0$.  The interchange $S$ is also a
nonidentity involution with fixed points: its fixed locus is obtained by
intersecting the smooth curve with the diagonal $x=y$.  A nonidentity
involution of a genus-one curve with a geometric fixed point is a
reflection and likewise induces $[-1]$ on $\Pic^0$.  Consequently
$S\circ\jmath_x$ induces the identity on $\Pic^0$.  An automorphism of
a genus-one curve acting trivially on $\Pic^0$ is a translation, and
its value at the origin shows that it is translation by $P_M$.
Iterating gives \eqref{eq:QRT-iteration-is-addition}.  A translation has
order $m$ precisely when its translating point has exact order $m$.
\end{proof}

The corollary is the precise content of the statement that QRT or
McMillan iteration is elliptic-curve addition.  It is stronger than the
observation that the invariant level set has genus one: it identifies
the actual increment on the group.
\section[Intrinsic displacement and origin]{The intrinsic displacement, the choice of origin, and the point \texorpdfstring{$D$}{D}}
\label{subsec:QRT-origin-and-D}

The equality ``the McMillan map is addition by a point'' requires a
choice of elliptic origin.  Before that choice, the intrinsic object is
the pair $(C,T)$ consisting of a smooth genus-one curve and a
translation-type automorphism.  There is nevertheless a canonical
translation class on the Jacobian.

\begin{theorem}[Origin-independent displacement class]
\label{thm:QRT-origin-independent-displacement}
Let $C/k$ be a smooth genus-one curve, and let $T:C\to C$ be a
$k$-automorphism inducing the identity on $\Pic^0(C)$.  Then
\begin{equation}
  \delta_T=[T(P)-P]\in\Pic^0(C)(k)
  \label{eq:QRT-intrinsic-displacement}
\end{equation}
is independent of the geometric point $P$.  If $O\in C(k)$ is chosen
as an origin and the corresponding group law is denoted by
$\oplus_O$, then
\begin{equation}
  T(P)=P\oplus_O D_O,
  \qquad
  D_O=T(O),
  \label{eq:QRT-translation-with-origin}
\end{equation}
and the identification
$\iota_O:C\to\Pic^0(C)$, $P\mapsto[P-O]$, sends $D_O$ to
$\delta_T$.

If $A\in C(k)$ is selected as a new origin, then
\begin{align}
  P\oplus_A R&=P\oplus_O R\ominus_O A,
  \label{eq:QRT-origin-change-law}\\
  D_A&=T(A)=A\oplus_O D_O.
  \label{eq:QRT-origin-change-D}
\end{align}
The map
\begin{equation}
  \varphi_A:(C,\oplus_A)\longrightarrow(C,\oplus_O),
  \qquad
  P\longmapsto P\ominus_O A,
  \label{eq:QRT-origin-change-isomorphism}
\end{equation}
is a group isomorphism and satisfies
$\varphi_A(D_A)=D_O$.  In particular, the order of the translation
increment is independent of the selected origin.
\end{theorem}

\begin{proof}
We separate the argument into four steps.

\emph{Step 1: independence of the geometric point.}
Work first over a separable closure $k_s$ of $k$.  Because $T$ is an
automorphism, the push-forward action of $T$ on divisor classes is the
inverse of its pullback action.  The hypothesis that $T$ induces the
identity on $\Pic^0(C)$ therefore implies that, for all
$P,R\in C(k_s)$,
\begin{equation}
  [T(P)-T(R)]=[P-R]
  \label{eq:QRT-displacement-difference-invariance}
\end{equation}
in $\Pic^0(C_{k_s})$.  Subtracting the classes attached to two possible
base points gives
\begin{align*}
  [T(P)-P]-[T(R)-R]
  &=[T(P)-T(R)]-[P-R]\\
  &=0
\end{align*}
by \eqref{eq:QRT-displacement-difference-invariance}.  Hence the class
$[T(P)-P]$ does not depend on $P\in C(k_s)$.

\emph{Step 2: descent to $k$.}
Let $\sigma\in\operatorname{Gal}(k_s/k)$.  Since $T$ is defined over
$k$, one has $\sigma\circ T=T\circ\sigma$.  Consequently,
\[
  \sigma\bigl([T(P)-P]\bigr)
  =[T(\sigma P)-\sigma P].
\]
The class on the right is the same point-independent class obtained in
Step~1.  Thus this class is fixed by the Galois group and defines
$\delta_T\in\Pic^0(C)(k)$.

\emph{Step 3: interpretation after choosing an origin.}
Fix $O\in C(k)$.  The Abel--Jacobi map
\[
  \iota_O:C\longrightarrow\Pic^0(C),
  \qquad P\longmapsto[P-O],
\]
is a $k$-isomorphism sending $O$ to the identity class.  For every
geometric point $P$,
\begin{align*}
  \iota_O(T(P))
  &=[T(P)-O]\\
  &=[T(P)-T(O)]+[T(O)-O]\\
  &=[P-O]+\delta_T
       &&\text{by \eqref{eq:QRT-displacement-difference-invariance}}\\
  &=\iota_O(P)+\delta_T.
\end{align*}
Because $\iota_O$ is a $k$-isomorphism and $\delta_T$ is $k$-rational,
there is a unique point $D_O\in C(k)$ satisfying
$\iota_O(D_O)=\delta_T$.  Transporting the addition law of
$\Pic^0(C)$ through $\iota_O$ turns the preceding equality into
\[
  T(P)=P\oplus_OD_O.
\]
Putting $P=O$ gives $D_O=T(O)$, and then
$\iota_O(D_O)=[T(O)-O]=\delta_T$.

\emph{Step 4: change of origin.}
Let $A\in C(k)$ and define
\[
  \varphi_A(P)=P\ominus_OA.
\]
This map is bijective, with inverse
$Q\mapsto Q\oplus_OA$.  Transporting the old group law through
$\varphi_A$ gives, by definition,
\begin{align*}
  P\oplus_AR
  &=\varphi_A^{-1}
      \bigl(\varphi_A(P)\oplus_O\varphi_A(R)\bigr)\\
  &=\bigl((P\ominus_OA)\oplus_O(R\ominus_OA)\bigr)
       \oplus_OA\\
  &=P\oplus_OR\ominus_OA.
\end{align*}
In particular $A$ is the identity for $\oplus_A$.  Since
$T(A)=A\oplus_OD_O$, the point translating $T$ in the new group law is
\[
  D_A=T(A)=A\oplus_OD_O.
\]
Indeed, for every $P$,
\begin{align*}
  P\oplus_AD_A
  &=P\oplus_OT(A)\ominus_OA\\
  &=P\oplus_O(A\oplus_OD_O)\ominus_OA\\
  &=P\oplus_OD_O\\
  &=T(P).
\end{align*}
Finally,
\[
  \varphi_A(D_A)
  =T(A)\ominus_OA
  =(A\oplus_OD_O)\ominus_OA
  =D_O.
\]
The defining transport identity also shows directly that
$\varphi_A(P\oplus_AR)=\varphi_A(P)\oplus_O\varphi_A(R)$; hence
$\varphi_A$ is a group isomorphism.  Group isomorphisms preserve exact
orders, so the order of the translating element is independent of the
chosen origin.
\end{proof}

The theorem gives a precise interpretation of the symbol $D$.  The
class $\delta_T\in\Jac(C)(k)$ is intrinsic to $(C,T)$, whereas the point
$D_O\in C(k)$ depends on the chosen identification of the genus-one
curve with its Jacobian.  Thus a point of $C$ is not an absolute
translation increment until an origin has been fixed.

\begin{corollary}[Making a prescribed curve point the translation point]
\label{cor:QRT-prescribed-D-by-origin}
Under the hypotheses of
Theorem~\ref{thm:QRT-origin-independent-displacement}, let
$Q\in C(k)$ be arbitrary.  There is a unique origin for which $T$ is
translation by the curve point $Q$, namely
\begin{equation}
  O_Q=T^{-1}(Q).
  \label{eq:QRT-prescribed-origin}
\end{equation}
With the group law having identity $O_Q$,
\begin{equation}
  T(P)=P\oplus_{O_Q}Q
  \label{eq:QRT-prescribed-translation}
\end{equation}
for every $P\in C$.
\end{corollary}

\begin{proof}
For an origin $A$, Theorem~\ref{thm:QRT-origin-independent-displacement}
shows that the translating point is $T(A)$.  It equals $Q$ if and only
if $A=T^{-1}(Q)$.  Since $T$ is an automorphism, this point exists and
is unique.
\end{proof}

For the affine McMillan map, if $Q=(u,v)$ and $u^2+\alpha\ne0$, the
origin in \eqref{eq:QRT-prescribed-origin} is
\begin{equation}
  O_Q=
  \left(-v-\frac{\beta u}{u^2+\alpha},u\right).
  \label{eq:QRT-prescribed-origin-affine}
\end{equation}
This formula is obtained by solving
$\mathcal M_Q(x,y)=(u,v)$: the first coordinate gives $y=u$, and the
second then gives the displayed value of $x$.  The condition
$u^2+\alpha\ne0$ is only the domain condition for this affine
representative, not a condition for the existence of $O_Q$.

More precisely, write
$Q=((U:Z),(V:W))\in\PP^1\times\PP^1$.  Two complementary homogeneous
representatives for the first coordinate of
$\mathcal M_Q^{-1}(Q)$ are
\begin{align}
  (X_s:Z_s)
  &={}
  \bigl(-V(U^2+\alpha Z^2)-\beta UZW:
        W(U^2+\alpha Z^2)\bigr),
  \label{eq:QRT-prescribed-origin-projective-sum}\\
  (X_p:Z_p)
  &={}
  \bigl((\alpha U^2+\gamma Z^2)W:
        (U^2+\alpha Z^2)V\bigr).
  \label{eq:QRT-prescribed-origin-projective-product}
\end{align}
The second coordinate is $(U:Z)$.  Equations
\eqref{eq:QRT-prescribed-origin-projective-sum} and
\eqref{eq:QRT-prescribed-origin-projective-product} are the homogeneous
sum and product forms of the same Vieta root.  Where one representative
is indeterminate, the other or the boundary atlas of the smooth
projective automorphism is used.  Thus a zero denominator in
\eqref{eq:QRT-prescribed-origin-affine} does not destroy the geometric
construction.

The effect on the Euler--Chasles parameter is equally explicit.  If
\[
  P(t)=\bigl(f(t),f(t+\delta)\bigr),
  \qquad T(P(t))=P(t+\delta),
\]
and $Q=P(q)$, define
\[
  \widetilde P(s)=P(s+q-\delta).
\]
Then
\[
  \widetilde P(0)=P(q-\delta)=T^{-1}(Q),
  \qquad
  \widetilde P(\delta)=P(q)=Q.
\]
Hence reselecting the origin is exactly a translation of the elliptic
parameter.

\begin{remark}[Changing the origin versus encoding a pointed curve]
\label{rem:QRT-reorigining-versus-pointing}
Corollary~\ref{cor:QRT-prescribed-D-by-origin} does not change the
intrinsic displacement $\delta_T$.  Under the group isomorphism
\eqref{eq:QRT-origin-change-isomorphism}, the newly named point $Q$ is
sent to the original translation element.  Its order, the order of
$T$, and the isomorphism class of the pointed Jacobian element are
unchanged.  This operation is therefore different from the following
construction: start with a fixed elliptic curve $(E,O)$ and a specified
point $D\in E(k)$, and choose coordinates in which the natural QRT
shift represents the original translation $P\mapsto P+D$.  The latter
is a genuine coordinate model of the pointed object $(E,D)$ and is the
one relevant when a cryptographic base point or fixed difference is
prescribed in advance.
\end{remark}

The intrinsic origin statement is easiest to interpret after separating a
mere change of origin from a genuine pointed-state construction.  The
following discussion explains exactly what a cryptographer may and may not
infer from the freedom to reselect the origin.  The operation count for the
McMillan step is treated in the next section.

\subsection{Effect of reorigining on the pointed state}
\label{subsec:QRT-reorigining-preserves-changes}

The observation \(O_Q=T^{-1}(Q)\) is elementary, but it prevents several
misinterpretations that are important in cryptographic arithmetic.  There are
two different meanings of the phrase ``preassign the point \(D\)''.

\paragraph{Changing the origin on a fixed pair \((C,T)\).}
For a fixed genus-one curve with its QRT automorphism, every
\(Q\in C(k)\) can be made the translating curve point by choosing
\(O_Q=T^{-1}(Q)\).  Under the resulting group law,
\[
  T(P)=P\oplus_{O_Q}Q.
\]
This is only a change of the torsor origin.  The isomorphism
\(\varphi_{O_Q}(P)=P\ominus_OO_Q\) sends the newly named increment \(Q\)
to the original increment \(D_O\).  Hence the Jacobian class
\(\delta_T\), the exact order of the translation, and the conjugacy class of
\(T\) do not change.

\paragraph{Encoding a point on an already pointed elliptic curve.}
A stronger construction starts with an elliptic curve \((E,O)\) and a point
\(D\in E(k)\) in the original group law, and then chooses coordinates in
which the natural QRT shift is the genuine map \(P\mapsto P+D\).  The
adjacent-state construction developed later realizes this by using
\((\kappa(P),\kappa(P+D))\).  In that setting the fixed difference is encoded
by the pointed coordinate model, rather than produced by renaming the origin.
Thus the two operations have different invariants:
\begin{enumerate}[label=(\alph*)]
  \item reorigining changes the representative of \(\delta_T\) on the torsor
  but preserves \(\delta_T\) and \(\operatorname{ord}(T)\);
  \item translation-adapted modelling changes the coordinates of the fixed
  object \((E,O,D)\) but preserves that pointed elliptic curve.
\end{enumerate}

\paragraph{Iteration, homomorphisms, and scalar multiplication.}
Fix an origin \(A\) and put \(D_A=T(A)\).  Since
\(T(P)=P\oplus_AD_A\), induction gives
\begin{equation}
  T^n(P)=P\oplus_A[n]_AD_A
  \qquad(n\in\mathbb Z).
  \label{eq:QRT-reorigin-iterate}
\end{equation}
For \(n\geq0\), the induction step is
\[
  T^{n+1}(P)
  =T^n(P)\oplus_AD_A
  =P\oplus_A[n+1]_AD_A;
\]
for \(m\geq0\), induction applied to
\(T^{-1}(P)=P\oplus_A[-1]_AD_A\) gives
\[
 T^{-m}(P)=P\oplus_A[-m]_AD_A.
\]
Together the two inductions prove the formula for every integer \(n\).
Thus a cheap McMillan step gives a cheap fixed
translation, but \(n\) successive shifts still require \(O(n)\) steps.

The translation \(T\) is a group homomorphism for \((C,\oplus_A)\) only in
the trivial case.  Indeed,
\begin{align*}
  T(P\oplus_AR)&=P\oplus_AR\oplus_AD_A,\\
  T(P)\oplus_AT(R)&=P\oplus_AR\oplus_A[2]_AD_A.
\end{align*}
Equality for all \(P,R\) is equivalent to
\(D_A=[2]_AD_A\), hence to \(D_A=A\), the identity of the new group.
Consequently a nontrivial QRT translation is not a GLV/GLS-type endomorphism.
The logarithmic algorithm requires the separate multiplication map \([2]\),
or its state conjugate \(\Delta_D\), rather than repeated use of \(T\) alone.

\paragraph{Inverse and doubling after reorigining.}
We derive both formulas from
\eqref{eq:QRT-origin-change-law}.  The inverse \(R=[-1]_A(P)\) is determined
by \(P\oplus_AR=A\), which is equivalent to
\(P\oplus_OR\ominus_OA=A\), and therefore
\(R=2A\ominus_OP\).  For doubling, the definition of \(\oplus_A\) gives
\[
 P\oplus_AP=P\oplus_OP\ominus_OA=[2]_O(P)\ominus_OA.
\]
Hence
\begin{equation}
  [-1]_A(P)=2A\ominus_OP,
  \qquad
  [2]_A(P)=[2]_O(P)\ominus_OA.
  \label{eq:QRT-reorigin-inverse-double}
\end{equation}
These identities show exactly why state doubling is not attached to the
unpointed pair \((C,T)\): multiplication by two becomes definite only after an
origin is selected.

\paragraph{Cryptographic reading.}
If a protocol has selected a point \(Q\), the choice
\(O_Q=T^{-1}(Q)\) makes the low-degree McMillan update equal to ``add \(Q\)''
in the reoriginated group.  This can simplify the interpretation of a fixed
shift, but it neither changes the intrinsic translation order nor supplies a
fast scalar multiplication by itself.  The genuinely pointed construction
\(P\mapsto(\kappa(P),\kappa(P+D))\), together with its state-doubling map,
is the stronger statement relevant to a fixed cryptographic base point.
The adjacent-state construction plays a role analogous to a pointed normal form: a Tate
normal form simplifies the coordinates of the marked point, whereas the
symmetric biquadratic state model simplifies the action of translation by the
marked point.  The ground-field classification of these state forms is given
in Sections~\ref{subsec:QRT-finite-field-classification-char2} and
\ref{subsec:QRT-abstract-isomorphism-odd}.

\section{Projective McMillan iteration and its cost}
\label{subsec:QRT-McMillan-projective}

Write a point as
\[
  ((X:Z),(Y:W))\in\PP^1\times\PP^1.
\]
The product-root form of the McMillan step has the inversion-free
projective realization
\begin{align}
  A&=Y^2,& B&=W^2,\notag\\
  H&=A+\alpha B,& K&=\alpha A+\gamma B,
  \label{eq:QRT-McMillan-projective-precomp}\\
  \mathcal M_Q((X:Z),(Y:W))
  &=((Y:W),(KZ:HX)).
  \label{eq:QRT-McMillan-projective-product}
\end{align}
It costs
\begin{equation}
  2\M+2\Sqr+3\Dpar.
  \label{eq:QRT-McMillan-product-cost}
\end{equation}
Here the three fixed multiplications are by $\alpha$ twice and by
$\gamma$ once.  Formula
\eqref{eq:QRT-McMillan-projective-product} is exceptionally small, but
it is not a single complete law: it becomes indeterminate on fibres for
which the product form has a zero numerator and denominator.

A complementary projective form comes from the sum of the two roots.
Put
\begin{equation}
  C=YW,
  \qquad
  N=-XH-\beta ZC,
  \qquad
  D=ZH.
  \label{eq:QRT-McMillan-projective-sum-precomp}
\end{equation}
Then
\begin{equation}
  \mathcal M_Q((X:Z),(Y:W))
  =((Y:W),(N:D)),
  \label{eq:QRT-McMillan-projective-sum}
\end{equation}
with cost
\begin{equation}
  4\M+2\Sqr+2\Dpar.
  \label{eq:QRT-McMillan-sum-cost}
\end{equation}
The two fixed multiplications are by $\alpha$ and $\beta$.  On their
common domain, equality of
\eqref{eq:QRT-McMillan-projective-product} and
\eqref{eq:QRT-McMillan-projective-sum} follows from the defining
biquadratic.  Together with the boundary tables already established for
the QRT--Jacobi and QRT--Weierstrass interfaces, these two formulas form
a complete iteration atlas.  Calling the first formula alone complete
would be incorrect.

Sequential computation of $\mathcal M_Q^n(R)$ by
\eqref{eq:QRT-McMillan-projective-product} takes $O(n)$ field
operations.  This fixed-step recurrence therefore makes no logarithmic
complexity claim.  A later Kummer construction uses the same fixed
difference $P_M$ to obtain logarithmic scalar multiplication.
Because
\begin{equation}
  U(P_M)=D_-,
  \label{eq:QRT-McMillan-known-difference}
\end{equation}
the affine-difference ladder step costs
\begin{equation}
  6\M+5\Sqr+3\Dpar
  \label{eq:QRT-McMillan-ladder-cost-odd}
\end{equation}
with the formulas used in this chapter.  Thus the direct McMillan step
is attractive for short fixed-step iteration, while the Kummer ladder
is the appropriate method for computing $[n]P_M$ with $O(\log n)$
steps.  For an arbitrary starting point $R$, first compute $[n]P_M$ and
then add it to $R$.

The use of QRT dynamics for ECM is structurally natural: scalar
multiplication on the genus-one fibre is implemented by a birational
recurrence in $\PP^1\times\PP^1$.  Projective and parallel Lyness/QRT
variants of ECM have been developed in the literature
\cite{HoneLynessECM}, providing an established
implementation benchmark for the complete QRT-state arithmetic constructed
here.  An end-to-end evaluation combines state doubling, differential
addition, quotient or full-point recovery, exceptional-fibre handling,
curve generation, and Stage~1 and Stage~2 integration.  This common
accounting framework permits reproducible comparisons with Montgomery-ECM
at identical bounds and integer sizes while making the state constraint,
compiled fixed difference, and parallel structure of the QRT interface
fully visible.

\section{Pointed elliptic curves and adjacent Kummer states}
\label{subsec:QRT-adjacent-Kummer-states}

The preceding origin discussion separates two constructions.  The
first starts from the intrinsic pair $(C,T)$ and changes the origin.
The second starts from a pointed elliptic curve $(E,D)$ with its
original identity fixed and constructs a biquadratic state model for
the genuine translation $P\mapsto P+D$.  The latter construction is a
Kummer-state version of the Euler--Chasles correspondence.

The structural simplicity of the symmetric biquadratic model is explained by
this adjacent-state construction.  Its two coordinates are not independent
affine coordinates on an elliptic curve: they are two consecutive values of
the same degree-two Kummer function,
\[
  x=\kappa(P),\qquad y=\kappa(P+D),
\]
separated by a fixed elliptic displacement $D$.  Since both coordinate
functions have degree two, their joint image is naturally a curve of
bidegree $(2,2)$ in $\PP^1\times\PP^1$.  The symmetry under coordinate
exchange reflects the identity
\[
  \mathcal S_D(-P-D)
  =
  \bigl(\kappa(P+D),\kappa(P)\bigr),
\]
while translation by $D$ becomes the shift
\[
  \bigl(\kappa(P),\kappa(P+D)\bigr)
  \longmapsto
  \bigl(\kappa(P+D),\kappa(P+2D)\bigr).
\]
Thus the symmetry, the Vieta involutions, and the low-degree QRT/McMillan
translation all arise from the same degree-two quotient together with a
fixed displacement on the elliptic curve.  Whenever the two Kummer
coordinates admit a common fractional-linear normalization over the ground
field, this state curve is represented by the even symmetric form
\[
  x^2y^2+\alpha(x^2+y^2)+\beta xy+\gamma=0.
\]
The following theorem makes this construction precise.

\begin{theorem}[Adjacent Kummer state curve]
\label{thm:QRT-adjacent-Kummer-state}
Let $(E,O)$ be an elliptic curve over a field $k$, let
$\kappa:E\to\PP^1$ be a separable degree-two quotient satisfying
\begin{equation}
  \kappa(P)=\kappa(-P)
  \label{eq:QRT-Kummer-even-function}
\end{equation}
with generic fibres $\{P,-P\}$, and let
$D\in E(k)$ satisfy $2D\ne O$.  Define
\begin{equation}
  \mathcal S_D:E\longrightarrow\PP^1\times\PP^1,
  \qquad
  \mathcal S_D(P)=\bigl(\kappa(P),\kappa(P+D)\bigr).
  \label{eq:QRT-adjacent-state-map}
\end{equation}
Then the following statements hold.

\begin{enumerate}[label=(\roman*)]
  \item The map $\mathcal S_D$ is an isomorphism from $E$ onto a smooth
  curve $\mathcal B_D\subset\PP^1\times\PP^1$ of bidegree $(2,2)$.

  \item The image is symmetric under
  $\sigma(x,y)=(y,x)$, and
  \begin{equation}
    \sigma\bigl(\mathcal S_D(P)\bigr)
    =\mathcal S_D(-P-D).
    \label{eq:QRT-adjacent-state-swap}
  \end{equation}

  \item Translation by $D$ is represented by the state shift
  \begin{equation}
    T_D\bigl(\mathcal S_D(P)\bigr)
    =\mathcal S_D(P+D)
    =\bigl(\kappa(P+D),\kappa(P+2D)\bigr).
    \label{eq:QRT-adjacent-state-shift}
  \end{equation}
  Whenever a common fractional-linear normalization of the two Kummer
  coordinates to \eqref{eq:symmetric-QRT} exists over the field under
  consideration, $T_D$ becomes one of the two orientations
  $\mathcal M_Q^{\pm1}$ of the McMillan map.  Over a nonsplitting
  ground field the state curve remains a twisted symmetric
  biquadratic.
\end{enumerate}
\end{theorem}

\begin{proof}
Suppose two generic points $P,R\in E$ have the same image under
$\mathcal S_D$.  Equality of the first coordinates and the generic
fibre description of $\kappa$ give $R=P$ or $R=-P$.  In the second
case, equality of the second coordinates gives
\[
  \kappa(D-P)=\kappa(D+P),
\]
so generically either $D-P=D+P$, which would force $2P=O$, or
$D-P=-D-P$, which would force $2D=O$.  The first condition holds only
on a finite set of $P$, and the second is excluded.  Hence
$\mathcal S_D$ is generically injective and therefore birational onto
its image.

The first projection of the image, composed with $\mathcal S_D$, is
$\kappa$ and has degree two; the second is
$\kappa\circ\tau_D$, where $\tau_D(P)=P+D$, and also has degree two.
Because $\mathcal S_D$ is birational, the image has bidegree $(2,2)$.
An integral $(2,2)$ curve in $\PP^1\times\PP^1$ has arithmetic genus
one.  Its normalization is $E$, which also has genus one, so the total
$\delta$-invariant is zero.  Thus the image is smooth, and the
birational morphism from the smooth projective curve $E$ is an
isomorphism.

Using the evenness of $\kappa$,
\[
  \mathcal S_D(-P-D)
  =\bigl(\kappa(P+D),\kappa(P)\bigr),
\]
which proves \eqref{eq:QRT-adjacent-state-swap}.  Moreover,
\[
 \mathcal S_D(P+D)
 =\bigl(\kappa(P+D),\kappa(P+2D)\bigr)
 =T_D\bigl(\mathcal S_D(P)\bigr),
\]
which proves \eqref{eq:QRT-adjacent-state-shift}.  On any symmetric
biquadratic normal form, the second
coordinate of the shifted state is the second root in the horizontal
fibre through the first state.  Hence the shift is the coordinate swap
followed by the corresponding Vieta exchange, which is precisely a
McMillan orientation.
\end{proof}

The hypothesis $2D\ne O$ is necessary for the stated birational model.
If $2D=O$, then
$\mathcal S_D(P)=\mathcal S_D(-P)$ generically, so the state map factors
through the Kummer line and the image is no longer a smooth
biquadratic model of $E$.

The theorem makes the Euler--Chasles parameterization transparent.  If
$E(\mathbb C)=\mathbb C/\Lambda$, if $\kappa(P(t))=f(t)$ with $f$ even,
and if $D$ corresponds to $\delta$, then
\begin{equation}
  \mathcal S_D(P(t))
  =\bigl(f(t),f(t+\delta)\bigr),
  \qquad
  T_D:t\longmapsto t+\delta.
  \label{eq:QRT-adjacent-state-Euler}
\end{equation}
The symmetric biquadratic relation records the unordered displacement
$\{D,-D\}$; choosing $T_D$ rather than $T_D^{-1}$ selects its
orientation.

\begin{proposition}[Ground-field state group law and origin alignment]
\label{prop:QRT-pointed-state-ground-field-law}
Under the hypotheses of
Theorem~\ref{thm:QRT-adjacent-Kummer-state}, the formula
\begin{equation}
 R_1+_D R_2
 =\mathcal S_D\bigl(
   \mathcal S_D^{-1}(R_1)+\mathcal S_D^{-1}(R_2)
  \bigr)
 \label{eq:QRT-pointed-state-ground-field-law}
\end{equation}
defines a group law over \(k\) on \(\mathcal B_D\), with identity
\begin{equation}
 O_S=\mathcal S_D(O).
 \label{eq:QRT-pointed-state-identity}
\end{equation}
It is a morphism on the smooth projective state curve and admits a finite
projective formula atlas.  An unpointed smooth QRT equation, by contrast,
is only a genus-one torsor unless a rational origin is supplied.

Suppose additionally that \(\charac(k)\ne2\), that the state has been put
in the exact form \(\Q_{\alpha,\beta,\gamma}\), and that
\(c_J^2=-4\alpha\gamma\) in the working field.  The Jacobi law of
Section~\ref{subsec:QRT-extended-Jacobi} has identity
\[
 O_J^Q=\left(0,\frac{c_J}{2\alpha}\right).
\]
If \(\oplus_J\) denotes that transported law on the QRT curve, then
\begin{equation}
 R_1\oplus_JR_2
 =R_1+_DR_2-_DO_J^Q.
 \label{eq:QRT-Jacobi-origin-alignment-law}
\end{equation}
Thus a fixed origin-alignment translation is required to turn a Jacobi-core
addition into the originally prescribed pointed-state addition unless
\(O_J^Q=O_S\).
\end{proposition}

\begin{proof}
The map \(\mathcal S_D\) is a \(k\)-isomorphism by
Theorem~\ref{thm:QRT-adjacent-Kummer-state}.  Conjugating the addition
morphism \(E\times E\to E\) by this isomorphism proves the group axioms,
the identity statement, and ground-field rationality.  A morphism between
projective varieties is represented on a finite affine cover by regular
coordinate tuples; homogenizing those tuples gives the asserted finite
projective atlas.  Without a chosen rational point there is no distinguished
identification of a genus-one curve with its Jacobian, which proves the
torsor qualification.

It remains to compare the two origins.  On any genus-one curve carrying the
law \(+_D\), the operation
\(R_1\star R_2=R_1+_DR_2-_DQ\) has identity \(Q\).  Taking
\(Q=O_J^Q\), Proposition~\ref{prop:QRT-full-point-group-law-summary}
identifies \(\star\) with the law transported from the Jacobi quartic.
This proves \eqref{eq:QRT-Jacobi-origin-alignment-law} and the final
assertion.
\end{proof}

\begin{remark}[Jacobi-core operation counts]
The costs in \eqref{eq:QRT-projective-cost-table} are exact upper bounds for
the displayed Jacobi core under the stated square-root hypothesis.  A full
pointed-state evaluation combines this core with the QRT--Jacobi interface,
the complementary exceptional charts, and, when the two origins differ,
the fixed translation in
\eqref{eq:QRT-Jacobi-origin-alignment-law}.
\end{remark}

\begin{proposition}[Available bounds for full-point addition]
\label{prop:QRT-full-addition-full-addition-bounds}
For every pointed adjacent state satisfying
Theorem~\ref{thm:QRT-adjacent-Kummer-state}, generic full-point addition is
defined over \(k\) and admits a finite complete projective atlas.  If
\(\charac(k)\ne2\), the state is in the exact form
\(\Q_{\alpha,\beta,\gamma}\), and
\(c_J^2=-4\alpha\gamma\) has a solution in the working field, then the
extended-Jacobi formulas of
\eqref{eq:QRT-projective-X}--\eqref{eq:QRT-projective-TZ} give the rigorous
core upper bounds in \eqref{eq:QRT-projective-cost-table}.

For the original pointed-state law, the complete cost consists of the core
plus the interfaces, complementary exceptional charts, and, when
\(O_J^Q\ne O_S\), the alignment translation
\eqref{eq:QRT-Jacobi-origin-alignment-law}.  When
\(-4\alpha\gamma\in k^{\times2}\), all of these components are defined over
\(k\).  In the nonsplit case, full addition is evaluated on the original
pointed Weierstrass model and transported by
\eqref{eq:QRT-pointed-state-ground-field-law}.  This gives a ground-field
complete atlas in both square classes and an explicit Jacobi-core schedule
on the split locus.
\end{proposition}

\begin{proof}
The ground-field existence, completeness atlas, and torsor qualification
are Proposition~\ref{prop:QRT-pointed-state-ground-field-law}.  Under the
square-root hypothesis, the complete QRT--Jacobi interface
Proposition~\ref{prop:QRT-Jacobi-original-atlas} and the projective
calculation in Section~\ref{subsec:QRT-extended-Jacobi} prove the displayed
core formulas and operation counts.  Equation
\eqref{eq:QRT-Jacobi-origin-alignment-law} computes the precise difference
between the two group laws, so omitting that translation is valid exactly
when their identities agree.  If the square root is absent, the definition
of the Jacobi ordinate divides by \(c_J\), and hence does not define a
\(k\)-rational interface; the transported law
\eqref{eq:QRT-pointed-state-ground-field-law} remains available because it
uses the original \(k\)-isomorphism \(\mathcal S_D\).  Together these
calculations prove the stated ground-field atlas, the split Jacobi schedule,
and the origin-alignment rule.
\end{proof}

\begin{theorem}[State doubling and the QRT binary ladder]
\label{thm:QRT-state-doubling-ladder}
Under the hypotheses of
Theorem~\ref{thm:QRT-adjacent-Kummer-state}, define
\begin{equation}
  \Delta_D
  =\mathcal S_D\circ[2]\circ\mathcal S_D^{-1}.
  \label{eq:QRT-state-doubling-definition}
\end{equation}
Then
\begin{equation}
  \Delta_D\circ T_D=T_D^2\circ\Delta_D.
  \label{eq:QRT-state-semiconjugacy}
\end{equation}
For
\begin{equation}
  S_n=\mathcal S_D([n]D)
  =\bigl(\kappa([n]D),\kappa([n+1]D)\bigr),
  \label{eq:QRT-adjacent-multiple-state}
\end{equation}
one has
\begin{align}
  T_D(S_n)&=S_{n+1},
  \label{eq:QRT-state-index-plus-one}\\
  \Delta_D(S_n)&=S_{2n},
  \label{eq:QRT-state-index-double}\\
  T_D^{\varepsilon}\Delta_D(S_n)&=S_{2n+\varepsilon},
  \qquad \varepsilon\in\{-1,0,1\}.
  \label{eq:QRT-state-index-signed}
\end{align}
Consequently, if $n=(b_{r-1}\cdots b_0)_2$, the recurrence
\begin{equation}
  R\leftarrow T_D^{b_i}\Delta_D(R),
  \qquad i=r-1,r-2,\ldots,0,
  \label{eq:QRT-state-binary-Horner}
\end{equation}
started at $R=S_0$ ends at $S_n$.  It computes the adjacent Kummer
state of $[n]D$ in $O(\log n)$ state steps.
\end{theorem}

\begin{proof}
Let $\tau_D:E\to E$ denote translation by $D$,
$\tau_D(P)=P+D$.  Theorem~\ref{thm:QRT-adjacent-Kummer-state} states that
$\mathcal S_D:E\to\mathcal B_D$ is an isomorphism of smooth projective
curves and that
\[
  T_D=\mathcal S_D\circ\tau_D\circ\mathcal S_D^{-1}.
\]
Because $\mathcal S_D$ is an isomorphism of smooth projective curves,
$\Delta_D=\mathcal S_D\circ[2]\circ\mathcal S_D^{-1}$ is a globally defined morphism of the state curve.

For $P\in E$, compute both sides of the claimed semiconjugacy after
applying them to $\mathcal S_D(P)$.  On the one hand,
\begin{align*}
  (\Delta_D\circ T_D)(\mathcal S_D(P))
  &=\Delta_D(\mathcal S_D(P+D))\\
  &=\mathcal S_D([2](P+D))\\
  &=\mathcal S_D(2P+2D).
\end{align*}
On the other hand,
\begin{align*}
  (T_D^2\circ\Delta_D)(\mathcal S_D(P))
  &=T_D^2(\mathcal S_D(2P))\\
  &=\mathcal S_D(2P+2D).
\end{align*}
The two morphisms agree on every state because $\mathcal S_D$ is
surjective.  This proves
$\Delta_D\circ T_D=T_D^2\circ\Delta_D$.

Now let $n\in\mathbb Z$.  By the definition of $S_n$ and the fact that
$T_D$ is translation by $D$ on the state curve,
\begin{align*}
  T_D(S_n)
  &=T_D(\mathcal S_D([n]D))\\
  &=\mathcal S_D([n]D+D)\\
  &=\mathcal S_D([n+1]D)=S_{n+1}.
\end{align*}
Likewise,
\begin{align*}
  \Delta_D(S_n)
  &=\mathcal S_D([2]([n]D))\\
  &=\mathcal S_D([2n]D)=S_{2n}.
\end{align*}
For $\varepsilon=1,0,-1$, respectively, application of $T_D$, the
identity, or $T_D^{-1}$ to $S_{2n}$ gives
$S_{2n+1}$, $S_{2n}$, or $S_{2n-1}$.  This proves
\eqref{eq:QRT-state-index-signed}, including the negative branch.

It remains to verify the binary algorithm without suppressing the
prefix induction.  For $j=0,1,\ldots,r$, let $m_j$ be the integer
represented by the first $j$ processed bits, with $m_0=0$.  Thus
\[
  m_{j+1}=2m_j+b_{r-1-j}.
\]
We prove by induction on $j$ that the stored state after $j$ updates is
$S_{m_j}$.  The assertion is true for $j=0$ because the algorithm starts
at $S_0$.  If it is true after $j$ updates, the next update gives
\begin{align*}
  T_D^{b_{r-1-j}}\Delta_D(S_{m_j})
  &=T_D^{b_{r-1-j}}(S_{2m_j})\\
  &=S_{2m_j+b_{r-1-j}}\\
  &=S_{m_{j+1}}.
\end{align*}
The induction therefore reaches $S_{m_r}=S_n$.  The number of processed
bits is $r=O(\log n)$ for $n>0$; the case $n=0$ is the initial state.
The output is an adjacent Kummer state.  Recovering a signed full point,
when required, is a separate final operation and is not used in the
correctness proof of the ladder.
\end{proof}

\begin{figure}[H]
\centering
\resizebox{0.95\textwidth}{!}{%
\begin{tikzpicture}[node distance=10mm and 18mm,>=Latex,every node/.style={font=\small}]
  \node at (0,1.8) {A comparison of two adjacent-state viewpoints};
  \node[draw,rounded corners,inner sep=5pt] (mont0) at (-4.0,0.35) {$([n]D,[n+1]D)$};
  \node[draw,rounded corners,inner sep=5pt] (mont1) at (-4.0,-1.35) {$([2n+b]D,[2n+b+1]D)$};
  \draw[->,thick] (mont0) -- node[left]{Montgomery ladder step} (mont1);
  \node[align=center] at (-4.0,-2.4) {full points maintained\\on the elliptic curve};

  \node[draw,rounded corners,inner sep=5pt] (s0) at (3.8,0.6) {$S_n=(\kappa([n]D),\kappa([n+1]D))$};
  \node[draw,rounded corners,inner sep=5pt] (s1) at (1.6,-1.2) {$S_{2n}=\Delta_D(S_n)$};
  \node[draw,rounded corners,inner sep=5pt] (s2) at (6.0,-1.2) {$S_{2n+1}=T_D\Delta_D(S_n)$};
  \draw[->,thick] (s0)--(s1);
  \draw[->,thick] (s0)--(s2);
  \node[align=center] at (3.8,-2.5) {adjacent Kummer state\\on a symmetric biquadratic curve};
\end{tikzpicture}%
}
\caption{Montgomery ladders and QRT state ladders organize the same binary recursion in two different ambient spaces.  The former keeps a neighboring pair of full points, whereas the latter packages the neighboring Kummer values themselves into a single symmetric biquadratic state.}
\label{fig:QRT-Montgomery-state-correspondence}
\end{figure}
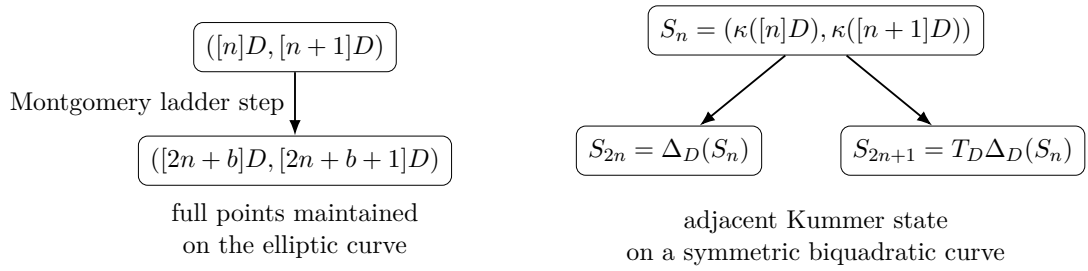

\subsection{The affine-index semigroup generated by \texorpdfstring{$T_D$ and $\Delta_D$}{TD and DeltaD}}
\label{subsec:QRT-affine-index-semigroup}

The three maps emphasized by the adjacent-state interpretation are
\begin{equation}
  T_D:\mathcal S_D(P)\longmapsto\mathcal S_D(P+D),
  \qquad
  \Delta_D:\mathcal S_D(P)\longmapsto\mathcal S_D(2P),
  \label{eq:QRT-three-state-maps-first-two}
\end{equation}
and
\begin{equation}
  T_D\Delta_D:\mathcal S_D(P)
  \longmapsto\mathcal S_D(2P+D).
  \label{eq:QRT-three-state-maps-third}
\end{equation}
They are not three unrelated formulas.  They are the first three elements of
an exact action of the semigroup of integral affine index maps.

\begin{theorem}[Affine-index semigroup on the state curve]
\label{thm:QRT-affine-index-semigroup}
Under the hypotheses of
Theorem~\ref{thm:QRT-adjacent-Kummer-state}, let $m\geq1$ and $r\in\mathbb Z$,
and define
\begin{equation}
  \Delta_{D,m}=\mathcal S_D\circ[m]\circ\mathcal S_D^{-1},
  \qquad
  \mathfrak F_{m,r}=T_D^r\circ\Delta_{D,m}.
  \label{eq:QRT-general-state-multiplication}
\end{equation}
Then, for every $P\in E$,
\begin{equation}
  \mathfrak F_{m,r}\bigl(\mathcal S_D(P)\bigr)
  =\mathcal S_D([m]P+[r]D).
  \label{eq:QRT-affine-index-action}
\end{equation}
The composition law is
\begin{equation}
  \mathfrak F_{m,r}\circ\mathfrak F_{n,s}
  =\mathfrak F_{mn,\,ms+r},
  \label{eq:QRT-affine-index-composition}
\end{equation}
which is the composition law of the affine maps
$t\mapsto mt+r$.  In particular,
\begin{align}
  \Delta_{D,m}\Delta_{D,n}&=\Delta_{D,mn},
  \label{eq:QRT-state-multiplication-composition}\\
  \Delta_{D,m}T_D^s&=T_D^{ms}\Delta_{D,m},
  \label{eq:QRT-state-general-semiconjugacy}\\
  \deg(\mathfrak F_{m,r})&=m^2.
  \label{eq:QRT-affine-index-degree}
\end{align}
If $D$ has infinite order, the maps $\mathfrak F_{m,r}$ are pairwise
distinct.  If $D$ has exact order $N$, then
\begin{equation}
  \mathfrak F_{m,r}=\mathfrak F_{n,s}
  \quad\Longleftrightarrow\quad
  m=n\ \text{ and }\ r\equiv s\pmod N.
  \label{eq:QRT-affine-index-faithfulness}
\end{equation}
\end{theorem}

\begin{proof}
Because $\mathcal S_D:E\to\mathcal B_D$ is an isomorphism, conjugating a
morphism of $E$ by $\mathcal S_D$ produces a globally defined morphism
of the state curve.  In particular, both $\Delta_{D,m}$ and
$\mathfrak F_{m,r}$ in
\eqref{eq:QRT-general-state-multiplication} are global morphisms.

For $P\in E$, first apply $\Delta_{D,m}$ and then the translation
$T_D^r$:
\begin{align*}
  \mathfrak F_{m,r}(\mathcal S_D(P))
  &=T_D^r\bigl(\mathcal S_D([m]P)\bigr)\\
  &=\mathcal S_D([m]P+[r]D).
\end{align*}
This is \eqref{eq:QRT-affine-index-action}.

To compute a composition, start with a state $\mathcal S_D(P)$.  The
inner map gives
\[
  \mathfrak F_{n,s}(\mathcal S_D(P))
  =\mathcal S_D([n]P+[s]D).
\]
Applying the outer map and using that $[m]$ is a group homomorphism,
\begin{align*}
  \mathfrak F_{m,r}
     \bigl(\mathfrak F_{n,s}(\mathcal S_D(P))\bigr)
  &=\mathcal S_D\bigl([m]([n]P+[s]D)+[r]D\bigr)\\
  &=\mathcal S_D\bigl([mn]P+[ms+r]D\bigr)\\
  &=\mathfrak F_{mn,ms+r}(\mathcal S_D(P)).
\end{align*}
Surjectivity of $\mathcal S_D$ proves equality of the two morphisms and
hence \eqref{eq:QRT-affine-index-composition}.  Taking $r=s=0$ gives
$\Delta_{D,m}\Delta_{D,n}=\Delta_{D,mn}$.  Taking the inner map to be
$T_D^s=\mathfrak F_{1,s}$ and the outer map to be
$\Delta_{D,m}=\mathfrak F_{m,0}$ gives
\[
  \Delta_{D,m}T_D^s
  =\mathfrak F_{m,0}\mathfrak F_{1,s}
  =\mathfrak F_{m,ms}
  =T_D^{ms}\Delta_{D,m}.
\]

Translations and isomorphisms have degree one.  The multiplication map
$[m]:E\to E$ has total degree $m^2$ for every $m\geq1$.  This statement
remains valid when $\operatorname{char}(k)$ divides $m$: in that case
$[m]$ may have a nontrivial inseparable degree, but the product of its
separable and inseparable degrees is still $m^2$.  Degree is invariant
under conjugation by an isomorphism and under composition with a
translation.  Therefore
\[
  \deg(\mathfrak F_{m,r})=\deg([m])=m^2.
\]

We finally prove the equality criterion in both directions.  Suppose
$\mathfrak F_{m,r}=\mathfrak F_{n,s}$.  Equality of degrees gives
$m^2=n^2$, and the assumptions $m,n\geq1$ imply $m=n$.  Conjugating the
remaining equality by $\mathcal S_D^{-1}$ yields
\[
  \tau_{[r]D}\circ[m]=\tau_{[s]D}\circ[m],
\]
where $\tau_Q(P)=P+Q$.  The finite morphism $[m]$ is surjective after
base change to an algebraic closure, including in the inseparable case.
Thus, for every geometric point $R\in E(\bar k)$, choose
$P\in E(\bar k)$ with $[m]P=R$ and obtain
\[
  R+[r]D=R+[s]D.
\]
Cancelling $R$ gives $[r-s]D=O$.  If $D$ has infinite order, this forces
$r=s$.  If $D$ has exact order $N$, it is equivalent to
$r\equiv s\pmod N$.

Conversely, if $m=n$ and $[r-s]D=O$, then
$[r]D=[s]D$.  Formula \eqref{eq:QRT-affine-index-action} shows that the
two maps have identical values on every state, hence they are equal.
This proves the stated faithfulness assertions and completes the
proof.
\end{proof}

\begin{corollary}[Affine reachability of a fixed-displacement state]
\label{cor:QRT-fixed-state-affine-reachability}
Any word generated by \(P\mapsto2P\) and translations
\(P\mapsto P+[r]D\) that contains exactly \(s\) doublings has the form
\begin{equation}
 P\longmapsto[2^s]P+[r']D
 \label{eq:QRT-single-state-affine-word}
\end{equation}
for some integer \(r'\).  Multi-scalar multiplication extends this
fixed-displacement engine by adjoining compatible state curves,
differential additions, or a higher-dimensional state for the selected
independent bases.
\end{corollary}

\begin{proof}
Before any operation the assertion holds with \(s=0\) and \(r'=0\).
A translation changes only \(r'\).  A doubling sends
\([2^s]P+[r']D\) to
\([2^{s+1}]P+[2r']D\).  Induction on the length of the word proves
\eqref{eq:QRT-single-state-affine-word}.  For two independent bases, the
corresponding construction augments the reachable family
\([2^s]P+[r']D\) by a second marked displacement or by a state carrying the
required differential data.  This is exactly the extension described in
the statement.
\end{proof}

For $m=2$, put
\begin{equation}
  \mathscr L_0=\Delta_D,
  \qquad
  \mathscr L_1=T_D\Delta_D.
  \label{eq:QRT-abstract-binary-branches}
\end{equation}
Then
\begin{equation}
  \mathscr L_b(S_n)=S_{2n+b},
  \qquad b\in\{0,1\}.
  \label{eq:QRT-abstract-binary-branch-index}
\end{equation}
More generally, if $b_0,b_1,\ldots,b_{\ell-1}$ are processed in that
order, then
\begin{equation}
  \mathscr L_{b_{\ell-1}}\circ\cdots\circ\mathscr L_{b_0}(S_n)
  =S_{2^\ell n+\sum_{j=0}^{\ell-1}b_j2^{\ell-1-j}}.
  \label{eq:QRT-binary-word-action}
\end{equation}
Thus the binary ladder realizes exactly the affine index monoid generated
by $n\mapsto2n$ and $n\mapsto2n+1$.  Replacing $2$ by any
$m\geq2$ gives a radix-$m$ state system based on
$T_D^r\Delta_{D,m}$, $0\leq r<m$; explicit low-cost formulas for
$\Delta_{D,m}$ are a separate arithmetic problem.

The state reflection supplies a second useful set of identities.
Let $\sigma(x,y)=(y,x)$ on the adjacent Kummer state curve and put
\begin{equation}
  \nu_D=T_D\sigma.
  \label{eq:QRT-state-negation-nu}
\end{equation}
By Theorem~\ref{thm:QRT-adjacent-Kummer-state},
$\sigma\mathcal S_D(P)=\mathcal S_D(-P-D)$, so
\begin{equation}
  \nu_D\mathcal S_D(P)=\mathcal S_D(-P).
  \label{eq:QRT-state-negation-nu-action}
\end{equation}
Consequently,
\begin{align}
  \sigma T_D\sigma&=T_D^{-1},
  \label{eq:QRT-state-reversibility}\\
  \nu_D\Delta_{D,m}&=\Delta_{D,m}\nu_D,
  \label{eq:QRT-state-negation-commutes-multiplication}\\
  \sigma\Delta_{D,m}\sigma&=T_D^{m-1}\Delta_{D,m}.
  \label{eq:QRT-state-reflection-multiplication}
\end{align}
For $m=2$, the last equality becomes the branch-conjugacy identity
\begin{equation}
  \boxed{\mathscr L_1=\sigma\mathscr L_0\sigma.}
  \label{eq:QRT-abstract-branch-conjugacy}
\end{equation}
It follows that a constant-pattern implementation needs only one branch
core: conditionally swap the two input Kummer coordinates, evaluate
$\mathscr L_0$, and conditionally swap the outputs.  Conditional swaps
are not field multiplications, so this organization does not change the
$(\M,\Sqr,\Dpar)$ count.

\begin{proposition}[Fixed states and exact periods]
\label{prop:QRT-state-fixed-periodic}
The map $\Delta_D$ has the unique fixed state $\mathcal S_D(O)$, and
$T_D\Delta_D$ has the unique fixed state $\mathcal S_D(-D)$.  If $D$
has exact order $N$, then $T_D$ has exact order $N$ and the state
sequence $(S_n)$ has exact period $N$:
\begin{equation}
  S_{n+N}=S_n\quad\text{for all }n,
  \qquad
  S_{n+t}=S_n\ \text{for all }n
  \Longrightarrow N\mid t.
  \label{eq:QRT-state-exact-period}
\end{equation}
\end{proposition}

\begin{proof}
Every state has a unique representation $\mathcal S_D(P)$ because
$\mathcal S_D$ is an isomorphism.

First consider $\Delta_D$.  The state $\mathcal S_D(P)$ is fixed by
$\Delta_D$ if and only if
\[
  \mathcal S_D(2P)=\mathcal S_D(P).
\]
Injectivity of $\mathcal S_D$ gives $2P=P$.  Subtracting $P$ in the
elliptic-curve group gives $P=O$.  Conversely,
$\Delta_D(\mathcal S_D(O))=\mathcal S_D(2O)=\mathcal S_D(O)$, so this
fixed state exists and is unique.

Next consider $T_D\Delta_D$.  One has
\[
  (T_D\Delta_D)(\mathcal S_D(P))=\mathcal S_D(2P+D).
\]
It is equal to $\mathcal S_D(P)$ if and only if $2P+D=P$, again by
injectivity.  Cancelling $P$ gives $P=-D$.  Conversely,
$2(-D)+D=-D$, so $\mathcal S_D(-D)$ is fixed.  This proves uniqueness of
the second fixed state.

For every nonnegative integer $t$, induction on $t$ gives
\[
  T_D^t(\mathcal S_D(P))=\mathcal S_D(P+[t]D).
\]
For $t<0$, the inverse $T_D^{-1}$ is translation by $-D$, and the same
induction applied to $-t$ gives the identical formula.  Hence $T_D^t$
is the identity morphism if and only if
\[
  \mathcal S_D(P+[t]D)=\mathcal S_D(P)
  \quad\text{for every }P.
\]
By injectivity this is equivalent to $P+[t]D=P$ for every $P$, and
therefore to $[t]D=O$.  If $D$ has exact order $N$, the least positive
integer satisfying this condition is $N$; thus $T_D$ has exact order
$N$.

The state sequence satisfies
\begin{align*}
  S_{n+N}
  &=\mathcal S_D([n+N]D)\\
  &=\mathcal S_D([n]D+[N]D)\\
  &=\mathcal S_D([n]D)=S_n.
\end{align*}
Conversely, suppose $S_{n+t}=S_n$ for at least one integer $n$.  Then
injectivity of $\mathcal S_D$ gives
$[n+t]D=[n]D$, and cancellation gives $[t]D=O$.  Therefore $N$ divides
$t$.  In particular, the implication in
\eqref{eq:QRT-state-exact-period}, which assumes equality for every
$n$, follows.  This argument also explains why the full adjacent state
must be used: equality of only one Kummer coordinate may identify a
point with its negative and can therefore suggest a spurious shorter
coordinate period, whereas equality of full states is governed by the
injective map $\mathcal S_D$.
\end{proof}

For implementation, the expression $T_D^{b_i}\Delta_D$ need not be
evaluated as two successive rational maps.  Write
$\operatorname{xDBL}(r)=\kappa(2R)$ when $r=\kappa(R)$, and let
$\operatorname{xADD}_D(r,s)$ denote the differential sum coordinate
when $s=\kappa(R+D)$.  The two bit branches are
\begin{align}
  \mathfrak L_0(r,s)
  &=\bigl(\operatorname{xDBL}(r),
           \operatorname{xADD}_D(r,s)\bigr),
  \label{eq:QRT-state-ladder-branch-zero}\\
  \mathfrak L_1(r,s)
  &=\bigl(\operatorname{xADD}_D(r,s),
           \operatorname{xDBL}(s)\bigr).
  \label{eq:QRT-state-ladder-branch-one}
\end{align}
For the state $S_n$, these outputs are $S_{2n}$ and $S_{2n+1}$,
respectively.  Thus both branches require one differential addition and
one doubling; a constant-pattern implementation obtains them by the
usual conditional-swap organization of a Montgomery ladder.  The
conceptual identity $T_D\Delta_D=S_{2n+1}$ therefore does not impose a
separate McMillan evaluation in each bit.

The coordinate swap gives a particularly small negation operation.
Indeed, by \eqref{eq:QRT-adjacent-state-swap},
\begin{equation}
  (T_D\circ\sigma)\mathcal S_D(P)=\mathcal S_D(-P).
  \label{eq:QRT-state-negation-abstract}
\end{equation}
In the canonical QRT coordinates this is the Vieta involution
\begin{equation}
  [-1]_O(x,y)
  =\mathcal M_Q(y,x)
  =\left(x,-y-\frac{\beta x}{x^2+\alpha}\right),
  \label{eq:QRT-state-negation-canonical}
\end{equation}
where the minus sign refers to the origin used to define the state
model.  Formula \eqref{eq:QRT-state-negation-canonical} must be replaced
by its projective Vieta representative on the exceptional fibres.

It is equally important that $T_D$ is not an elliptic-curve
endomorphism unless $D=O$.  In the fixed group law,
\[
  T_D(P+R)=P+R+D,
  \qquad
  T_D(P)+T_D(R)=P+R+2D.
\]
Thus a cheap QRT translation does not provide a GLV- or GLS-type scalar
decomposition.  The logarithmic algorithm comes from the separate map
$\Delta_D$, not from treating $T_D$ as a homomorphism.

\subsection{An explicit odd-characteristic pointed model.}
The general state theorem becomes completely explicit on
\eqref{eq:EAB}.

\begin{theorem}[Symmetric biquadratic attached to a marked point]
\label{thm:QRT-pointed-EAB-biquadratic}
Assume $\charac(k)\ne2$.  Let
\[
  E_{A,B}:v^2=u^3+Au^2+Bu,
  \qquad B(A^2-4B)\ne0,
\]
and let $D=(d,e)\in E_{A,B}(k)$ satisfy $2D\ne O$.  For a variable
point $P$, put
\begin{equation}
  x=u(P),\qquad y=u(P+D).
  \label{eq:QRT-pointed-EAB-state}
\end{equation}
Then the state curve is
\begin{equation}
\begin{split}
  \mathcal B_{A,B;d}:\quad
  0={}&x^2y^2-2dxy(x+y)+d^2(x^2+y^2)\\
      &{}-2Bd(x+y)-2(B+d^2+2Ad)xy+B^2.
  \label{eq:QRT-pointed-EAB-relation}
\end{split}
\end{equation}
It is a smooth symmetric biquadratic and depends on the marked point
only through $d=u(D)=u(-D)$.  Its two McMillan orientations distinguish
$D$ from $-D$.

Suppose additionally that $c\in k$ satisfies $c^2=B$.  The common
M\"obius change
\begin{equation}
  X=\frac{x-c}{x+c},
  \qquad
  Y=\frac{y-c}{y+c}
  \label{eq:QRT-pointed-EAB-evenization}
\end{equation}
identifies \eqref{eq:QRT-pointed-EAB-relation} with
$\Q_{\alpha_D,\beta_D,\gamma_D}$, where
\begin{align}
  \alpha_D&=\frac{Ad+B+d^2}{d(2c-A)},
  \label{eq:QRT-pointed-alpha}\\
  \beta_D&=\frac{2(B-d^2)}{d(2c-A)},
  \label{eq:QRT-pointed-beta}\\
  \gamma_D&=-\frac{A+2c}{2c-A}.
  \label{eq:QRT-pointed-gamma}
\end{align}
Under this identification, the genuine translation
$P\mapsto P+D$ becomes the canonical McMillan shift
\begin{equation}
  (X,Y)\longmapsto
  \left(Y,-X-\frac{\beta_DY}{Y^2+\alpha_D}\right).
  \label{eq:QRT-pointed-McMillan-shift}
\end{equation}
If $B$ is not a square in $k$, the same normal form exists over
$k(\sqrt B)$; over $k$ itself,
\eqref{eq:QRT-pointed-EAB-relation} is the corresponding twisted
symmetric state model.
\end{theorem}

\begin{proof}
Let $z=u(P-D)$.  Applying
Theorem~\ref{thm:EAB-Kummer-identities} to $P$ and $D$ gives
\begin{align}
  yz&=\frac{(xd-B)^2}{(x-d)^2},
  \label{eq:QRT-pointed-EAB-product}\\
  y+z&=\frac{2\bigl((x+d)(xd+B)+2Axd\bigr)}{(x-d)^2}.
  \label{eq:QRT-pointed-EAB-sum}
\end{align}
The two values $y,z$ are the roots of
$T^2-(y+z)T+yz$.  Multiplication by $(x-d)^2$ and substitution of
\eqref{eq:QRT-pointed-EAB-product}--
\eqref{eq:QRT-pointed-EAB-sum} gives
\[
  (x-d)^2y^2
  -2\bigl((x+d)(xd+B)+2Axd\bigr)y
  +(xd-B)^2=0.
\]
Expanding and collecting the symmetric monomials yields exactly
\eqref{eq:QRT-pointed-EAB-relation}.  Smoothness follows from
Theorem~\ref{thm:QRT-adjacent-Kummer-state}, because $2D\ne O$.
The ordinate $e$ does not occur in the Kummer identities, so the
relation depends only on the unordered pair $\{D,-D\}$.

Since $c^2=B$, the inverse of
\eqref{eq:QRT-pointed-EAB-evenization} is
\[
  x=c\frac{1+X}{1-X},
  \qquad
  y=c\frac{1+Y}{1-Y}.
\]
The four linear factors that occur in the substitution are
\[
  x-c=\frac{2cX}{1-X},
  \quad x+c=\frac{2c}{1-X},
  \quad y-c=\frac{2cY}{1-Y},
  \quad y+c=\frac{2c}{1-Y}.
\]
After multiplying by $(1-X)^2(1-Y)^2$, expanding only in the
symmetric basis
$1$, $X^2+Y^2$, $XY$, and $X^2Y^2$, and using $B=c^2$, the respective
coefficients are
\[
  -4c^2d(A+2c),
  \quad 4c^2(Ad+B+d^2),
  \quad 8c^2(B-d^2),
  \quad 4c^2d(2c-A).
\]
Consequently,
\begin{align*}
 &(1-X)^2(1-Y)^2
 \mathcal B_{A,B;d}
 \left(c\frac{1+X}{1-X},c\frac{1+Y}{1-Y}\right)\\
 &\qquad=4c^2d(2c-A)
 \bigl(X^2Y^2+\alpha_D(X^2+Y^2)+\beta_DXY+\gamma_D\bigr),
\end{align*}
where the three ratios of coefficients are precisely
\eqref{eq:QRT-pointed-alpha}--\eqref{eq:QRT-pointed-gamma}.
Here $d\ne0$ because the only point of \eqref{eq:EAB} with
$u=0$ is the two-torsion point $(0,0)$, and $2c-A\ne0$ because
$A^2-4B\ne0$.  Division by the nonzero scalar proves the parameter
formulas.  Finally, the state shift
$(u(P),u(P+D))\mapsto(u(P+D),u(P+2D))$ is the coordinate swap followed
by the second Vieta root; conjugating by the common M\"obius map gives
\eqref{eq:QRT-pointed-McMillan-shift}.
\end{proof}

The evenization has a simple group-theoretic explanation.  Translation
by the marked two-torsion point $(0,0)$ on $E_{A,B}$ sends
\begin{equation}
  u\longmapsto\frac{B}{u}.
  \label{eq:QRT-EAB-two-torsion-Kummer}
\end{equation}
Indeed, the line through $(u,v)$ and $(0,0)$ has slope $v/u$, and the
addition formula gives $(v/u)^2-A-u=B/u$.  When $B=c^2$, the M\"obius
coordinate $(u-c)/(u+c)$ conjugates $u\mapsto B/u$ to negation.  This
is why the general symmetric state relation becomes the even normal
form \eqref{eq:symmetric-QRT} after
\eqref{eq:QRT-pointed-EAB-evenization}.

\begin{proposition}[Odd-characteristic state doubling]
\label{prop:QRT-pointed-EAB-state-doubling}
In the setting of
Theorem~\ref{thm:QRT-pointed-EAB-biquadratic}, the affine state-doubling
map in the Kummer-adapted coordinates
$(x,y)=(u(P),u(P+D))$ is
\begin{equation}
  \Delta_D(x,y)=
  \left(
    \frac{(x^2-B)^2}{4x(x^2+Ax+B)},
    \frac{(xy-B)^2}{d(x-y)^2}
  \right)
  \label{eq:QRT-pointed-EAB-state-double-affine}
\end{equation}
where the displayed quotients are interpreted projectively.  For
projective input
$((X_0:Z_0),(X_1:Z_1))$, a homogeneous representative is
\begin{align}
  X_0'&=(X_0^2-BZ_0^2)^2,
  \label{eq:QRT-pointed-EAB-state-double-X0}\\
  Z_0'&=4X_0Z_0
  (X_0^2+AX_0Z_0+BZ_0^2),
  \label{eq:QRT-pointed-EAB-state-double-Z0}\\
  X_1'&=(X_0X_1-BZ_0Z_1)^2,
  \label{eq:QRT-pointed-EAB-state-double-X1}\\
  Z_1'&=d(X_0Z_1-Z_0X_1)^2.
  \label{eq:QRT-pointed-EAB-state-double-Z1}
\end{align}
The direct schedule costs
\begin{equation}
  6\M+5\Sqr+3\Dpar.
  \label{eq:QRT-pointed-EAB-state-double-cost}
\end{equation}
If $d$ is a compile-time curve constant rather than a run-time known
difference, the final multiplication in
\eqref{eq:QRT-pointed-EAB-state-double-Z1} is a fixed multiplication,
and the count becomes
\begin{equation}
  5\M+5\Sqr+4\Dpar.
  \label{eq:QRT-pointed-EAB-state-double-compiled-cost}
\end{equation}
\end{proposition}

\begin{proof}
The first coordinate of \eqref{eq:QRT-pointed-EAB-state-double-affine}
is \eqref{eq:EAB-Kummer-double-affine}.  For the second, apply
\eqref{eq:EAB-Kummer-product} to the pair $P$ and $P+D$.  Their sum is
$2P+D$, their difference is $-D$, and its Kummer coordinate is $d$.
Thus
\[
  u(2P+D)d=\frac{(xy-B)^2}{(x-y)^2},
\]
which proves the affine formula.  Homogenization gives
\eqref{eq:QRT-pointed-EAB-state-double-X0}--
\eqref{eq:QRT-pointed-EAB-state-double-Z1}.

The first output is the Kummer doubling schedule
\eqref{eq:EAB-xDBL}, of cost $2\M+3\Sqr+2\Dpar$.  The second is the
affine-known-difference product law
\eqref{eq:EAB-xADD-product}, of cost
$4\M+2\Sqr+1\Dpar$.  This gives
\eqref{eq:QRT-pointed-EAB-state-double-cost}.  Reclassifying the
multiplication by fixed $d$ changes one general multiplication into one
fixed multiplication and proves
\eqref{eq:QRT-pointed-EAB-state-double-compiled-cost}.  Because
$2D\ne O$, the known difference is neither the identity nor a
two-torsion state, so the fast product branch is the regular branch of
the complete differential atlas.

For completeness, the alternating cross-product used in this count can be
formed with three, rather than four, general products.  Put
\[
 P=X_0X_1,\qquad Q=Z_0Z_1,\qquad
 C=(X_0+Z_0)(X_1-Z_1).
\]
Then
\[
 P-Q-C=X_0Z_1-Z_0X_1.
\]
Thus \(P,Q,C\) supply both expressions in the second output pair.  This is
the alternating Karatsuba schedule underlying the displayed run-time and
compiled-constant counts.
\end{proof}

Transporting these formulas through
\eqref{eq:QRT-pointed-EAB-evenization} gives explicit state doubling on
the canonical $\Q_{\alpha_D,\beta_D,\gamma_D}$ model.  The operation
counts just stated refer to the Kummer-adapted representative
\eqref{eq:QRT-pointed-EAB-relation}; conjugating by a M\"obius change
may add fixed multiplications and should be counted separately in an
implementation.

On the exactly evenizable locus, the common M\"obius map can be represented
in the selected homogeneous coordinates by additions and subtractions only.
The standard Montgomery core therefore gives a sharper rigorous bound.

\begin{theorem}[Montgomery-conjugate complete state doubling]
\label{thm:QRT-Montgomery-conjugate-state-doubling}
Assume the hypotheses of
Theorem~\ref{thm:QRT-pointed-EAB-biquadratic}, and suppose
\(B=c^2\) with \(c\in k^\times\).  Put
\begin{equation}
  a=\frac Ac,\qquad
  \delta_D=\frac dc,\qquad
  A_{24}=\frac{a+2}{4},\qquad
  U=\frac uc,
  \label{eq:QRT-state-Montgomery-normalization}
\end{equation}
and use the exact-even coordinate
\begin{equation}
  z=\frac{U-1}{U+1}.
  \label{eq:QRT-state-even-coordinate-z}
\end{equation}
For a state represented by
\(z_i=(R_i:S_i)\), \(i=0,1\), make the input conversion
\begin{equation}
  (X_i:Z_i)=(R_i+S_i:S_i-R_i).
  \label{eq:QRT-state-z-to-Montgomery}
\end{equation}
Define
\begin{align*}
 A_0&=X_0+Z_0, & AA&=A_0^2,\\
 B_0&=X_0-Z_0, & BB&=B_0^2,\\
 E_0&=AA-BB,\\
 C_0&=X_1+Z_1, & D_0&=X_1-Z_1,\\
 DA&=D_0A_0, & CB&=C_0B_0.
\end{align*}
Then state doubling is represented by
\begin{align}
 X_{2P}&=AA\,BB,
 &Z_{2P}&=E_0(BB+A_{24}E_0),
 \label{eq:QRT-state-Montgomery-xDBL}\\
 X_{2P+D}&=(DA+CB)^2,
 &Z_{2P+D}&=\delta_D(DA-CB)^2.
 \label{eq:QRT-state-Montgomery-xADD}
\end{align}
The output conversion is
\begin{equation}
  (R_i':S_i')=(X_i'-Z_i':X_i'+Z_i').
  \label{eq:QRT-state-Montgomery-to-z}
\end{equation}
Equations \eqref{eq:QRT-state-z-to-Montgomery}--
\eqref{eq:QRT-state-Montgomery-to-z} define the state map at every
geometric point of the smooth adjacent-state curve.  If \(A_{24}\) and
\(\delta_D\) are compiled constants, their loop cost is
\begin{equation}
  \boxed{4\M+4\Sqr+2\Dpar}.
  \label{eq:QRT-state-Montgomery-compiled-cost}
\end{equation}
\end{theorem}

\begin{proof}
Substituting \(u=cU\) into the first coordinate of
\eqref{eq:QRT-pointed-EAB-state-double-affine} and using \(B=c^2\)
gives
\[
 U(2P)=\frac{(U^2-1)^2}{4U(U^2+aU+1)}.
\]
For the second coordinate, put \(U_0=x/c\) and \(U_1=y/c\).  Dividing
\[
 u(2P+D)d=\frac{(xy-c^2)^2}{(x-y)^2}
\]
by \(c^2\) gives the exact identity
\[
 U(2P+D)
 =\frac{(U_0U_1-1)^2}{\delta_D(U_0-U_1)^2}.
\]
Thus the two coordinates are precisely the homogeneous Montgomery doubling
and fixed-difference differential-addition identities
\eqref{eq:QRT-state-Montgomery-xDBL} and
\eqref{eq:QRT-state-Montgomery-xADD}
\cite{Montgomery1987,CostelloSmith2018}.  This normalization concerns only
the Kummer coordinate; no rescaling of the Weierstrass ordinate is used.

Solving \eqref{eq:QRT-state-even-coordinate-z} gives
\(U=(1+z)/(1-z)\).  Therefore
\eqref{eq:QRT-state-z-to-Montgomery} and
\eqref{eq:QRT-state-Montgomery-to-z} are inverse elements of
\(\operatorname{PGL}_2(k)\), evaluated with additions and subtractions
only.  Conjugating the Montgomery formulas by them gives the asserted state
map without adding a multiplication or a squaring.

We now prove completeness rather than infer it from affine formulas.  Up to
a nonzero common scalar, the doubling output pair is
\begin{equation}
 \bigl((X^2-Z^2)^2:
 4XZ(X^2+aXZ+Z^2)\bigr).
 \label{eq:QRT-state-Montgomery-doubling-basepair}
\end{equation}
If both entries vanished, the first would give \(X/Z=1\) or \(X/Z=-1\);
the second could then vanish only for \(a=-2\) or \(a=2\), respectively.
For \(a=-2\) the quadratic factor is \((U-1)^2\), and for \(a=2\) it
is \((U+1)^2\); in either case the cubic
\(U(U^2+aU+1)\) has a repeated root and the normalized Montgomery curve
is singular.  If \(Z=0\), the first entry is \(X^4\ne0\), and if \(X=0\),
it is \(Z^4\ne0\).  Hence the doubling pair has no geometric base point
on a smooth member.

For the differential output, removal of the harmless factor \(4\) gives
the homogeneous pair
\begin{equation}
 \bigl((X_0X_1-Z_0Z_1)^2:
 \delta_D(X_0Z_1-Z_0X_1)^2\bigr).
 \label{eq:QRT-state-Montgomery-addition-basepair}
\end{equation}
Here \(\delta_D\ne0\), since \(d=0\) would make
\(D=(0,0)\) the nonzero two-torsion point, contrary to \(2D\ne O\).
If one of \(X_0,Z_0,X_1,Z_1\) vanishes, the two bilinear forms in
\eqref{eq:QRT-state-Montgomery-addition-basepair} cannot vanish together.
Indeed, \(Z_0=0\) forces \(X_1=Z_1=0\), while \(X_0=0\) forces
\(Z_1=X_1=0\); interchanging the two inputs gives the contradictions for
\(Z_1=0\) and \(X_1=0\).  Each conclusion is impossible for a projective
point.  Thus a common zero must be affine.  Writing \(U_i=X_i/Z_i\), it
would satisfy \(U_0U_1=1\) and \(U_0=U_1\), hence
\(U_0=U_1=\pm1\).

Equality of the adjacent Kummer coordinates implies either \(D=O\) or
\(2P=-D\).  The first alternative is excluded.  Moreover \(U=\pm1\) is a
fixed point of the Kummer involution \(U\mapsto1/U\) induced by the marked
two-torsion point \(T_2=(0,0)\).  Since \(T_2\ne O\), the fixed-point
condition gives \(P+T_2=-P\), hence \(2P=-T_2\).  Comparing the two
equalities yields \(D=T_2\), again contradicting \(2D\ne O\).  Thus the
differential pair also has no geometric base point.  The input and output
PGL$_2$ maps are everywhere defined, so the complete four-form tuple
remains base-point-free on the state curve.

Finally, \(DA,CB,AA\,BB\), and
\(E_0(BB+A_{24}E_0)\) account for four general multiplications.  The
squares \(AA,BB,(DA+CB)^2,(DA-CB)^2\) account for four squarings, and
the products by \(A_{24}\) and \(\delta_D\) account for two fixed
multiplications.  This proves \eqref{eq:QRT-state-Montgomery-compiled-cost}.
\end{proof}

\begin{proposition}[State doubling is fixed-difference \texttt{xDBLADD}]
\label{prop:QRT-state-fixed-difference-xDBLADD}
Under the hypotheses of
Theorem~\ref{thm:QRT-adjacent-Kummer-state}, state doubling is intrinsically
the joint Kummer operation
\begin{equation}
 \bigl(\kappa(P),\kappa(P+D)\bigr)
 \longmapsto
 \bigl(\kappa(2P),\kappa(2P+D)\bigr).
 \label{eq:QRT-state-fixed-difference-xDBLADD}
\end{equation}
The first output is Kummer doubling.  The second output is differential
addition applied to \(P\) and \(P+D\), whose known difference is the fixed
point \(D\).  Hence \(\Delta_D\) is the fixed-difference
\texttt{xDBLADD} operation itself.

On the split exact-even locus of
Theorem~\ref{thm:QRT-Montgomery-conjugate-state-doubling}, the common
Kummer coordinate conjugates \eqref{eq:QRT-state-fixed-difference-xDBLADD}
to the standard Montgomery \texttt{xDBLADD} core.  This conjugation does
not change the numbers of general multiplications or squarings.  Moreover,
the other binary branch satisfies
\begin{equation}
 T_D\Delta_D=\sigma\Delta_D\sigma,
 \qquad \sigma(x,y)=(y,x),
 \label{eq:QRT-state-xDBLADD-branch-swap}
\end{equation}
so it uses the same core with input and output swaps.
\end{proposition}

\begin{proof}
By definition,
\[
 \Delta_D
 =\mathcal S_D\circ[2]\circ\mathcal S_D^{-1},
 \qquad
 \mathcal S_D(P)=
 \bigl(\kappa(P),\kappa(P+D)\bigr).
\]
Applying \(\Delta_D\) to \(\mathcal S_D(P)\) gives
\[
 \Delta_D(\mathcal S_D(P))
 =\mathcal S_D(2P)
 =\bigl(\kappa(2P),\kappa(2P+D)\bigr),
\]
which proves \eqref{eq:QRT-state-fixed-difference-xDBLADD}.  In its second
coordinate the two Kummer inputs represent \(P\) and \(P+D\); their
difference is \((P+D)-P=D\), independent of \(P\).  Thus the second
coordinate is a differential addition with compiled known difference.
Equations \eqref{eq:QRT-state-Montgomery-xDBL} and
\eqref{eq:QRT-state-Montgomery-xADD}, together with the addition-only
changes \eqref{eq:QRT-state-z-to-Montgomery} and
\eqref{eq:QRT-state-Montgomery-to-z}, prove the Montgomery-core assertion.

For the branch identity, use
\eqref{eq:QRT-adjacent-state-swap} successively.  On the elliptic parameter
\(P\), the three maps on the right-hand side act as
\[
 P\longmapsto -P-D
 \longmapsto -2P-2D
 \longmapsto 2P+D.
\]
The left-hand side sends \(P\) first to \(2P\) and then to \(2P+D\).
Since \(\mathcal S_D\) is an isomorphism, the two state morphisms are equal,
which proves \eqref{eq:QRT-state-xDBLADD-branch-swap}.
\end{proof}

\begin{corollary}[Minimum number of state-update formulas]
\label{cor:QRT-one-formula-complete-state-update}
Under the hypotheses of
Theorem~\ref{thm:QRT-Montgomery-conjugate-state-doubling}, a complete
projective atlas for \(\Delta_D\) requires exactly one formula.  The
bit-one branch \(T_D\Delta_D\) also requires one formula and uses the same
core with input and output swaps.
\end{corollary}

\begin{proof}
The preceding theorem gives one globally defined, base-point-free tuple, and
an atlas cannot contain fewer than one tuple.  Equation
\eqref{eq:QRT-abstract-branch-conjugacy} gives
\(T_D\Delta_D=\sigma\Delta_D\sigma\), so coordinate swaps transfer the
same completeness statement to the other branch without changing the field
operation count.
\end{proof}

\begin{remark}[Circuit class and sharp bound]
The cost \eqref{eq:QRT-state-Montgomery-compiled-cost} is the proved circuit
bound for the displayed adjacent-state update on the split exact-even locus.
The Lucas product factorization developed later represents a distinct
evaluation schedule; its completeness follows from its exact base-locus
calculation.  These two results give, respectively, the compiled state
circuit and the complete product-factorized circuit class.
\end{remark}

\subsection{An exact binary symmetric state model.}
Characteristic two gives an especially economical specialization.
It is not obtained by formally setting \(2=0\) in the odd-characteristic
descent theorem: the involution \(z\mapsto-z\) is the identity in
characteristic two and cannot encode a nontrivial split two-torsion action.
The result below is instead an independent Artin--Schreier construction
using the separable Kummer quotient of an ordinary binary elliptic curve.

\begin{theorem}[Binary pointed state equation and Artin--Schreier twist]
\label{thm:QRT-binary-pointed-state}
Assume $\charac(k)=2$, and let
\begin{equation}
  E_{a,b}^{(2)}:
  \qquad v^2+uv=u^3+au^2+b^2,
  \qquad b\ne0.
  \label{eq:QRT-binary-pointed-curve}
\end{equation}
Let $D=(d,e)\in E_{a,b}^{(2)}(k)$ satisfy $2D\ne O$, so $d\ne0$, and
put
\[
  x=u(P),\qquad y=u(P+D).
\]
Then the adjacent state curve is exactly
\begin{equation}
  x^2y^2+d^2(x^2+y^2)+dxy+b^2=0,
  \label{eq:QRT-binary-pointed-relation}
\end{equation}
that is,
\begin{equation}
  (\alpha,\beta,\gamma)=(d^2,d,b^2).
  \label{eq:QRT-binary-pointed-parameters}
\end{equation}
It is smooth because $bd\ne0$.  The coefficient $a$ does not occur in
the state relation.  More precisely, put
\begin{equation}
  a_D=d^2+\frac{b^2}{d^2},
  \qquad
  \lambda_D=d+\frac{e}{d}.
  \label{eq:QRT-binary-pointed-twist-parameters}
\end{equation}
Then
\begin{equation}
  a_D=a+\lambda_D^2+\lambda_D,
  \label{eq:QRT-binary-pointed-twist-relation}
\end{equation}
and the ordinate change
\begin{equation}
  (u,v)\longmapsto(u,v+\lambda_Du)
  \label{eq:QRT-binary-pointed-twist-isomorphism}
\end{equation}
identifies $E_{a,b}^{(2)}$ with
$E_{a_D,b}^{(2)}$, which is the shifted Weierstrass model attached to
$\Q_{d^2,d,b^2}$.
\end{theorem}

\begin{proof}
Let $z=u(P-D)$.  The binary Kummer identities
\eqref{eq:binary-Kummer-product}--
\eqref{eq:binary-Kummer-sum} give
\[
  yz=\frac{(xd+b)^2}{(x+d)^2},
  \qquad
  y+z=\frac{xd}{(x+d)^2}.
\]
In characteristic two, $y$ is a root of
$T^2+(y+z)T+yz$.  Multiplication by $(x+d)^2$ gives
\[
  (x+d)^2y^2+xdy+(xd+b)^2=0.
\]
Expanding the two squares yields
\eqref{eq:QRT-binary-pointed-relation}.

We next verify smoothness without appealing to the later general
characteristic-two criterion.  Homogenize the state equation as
\[
 X^2Y^2+d^2(X^2W^2+Z^2Y^2)+dXZYW+b^2Z^2W^2=0.
\]
Its four partial derivatives are
\[
 dZYW,\qquad dXYW,\qquad dXZW,\qquad dXZY.
\]
In the affine chart $Z=W=1$, a singular point would have $x=y=0$,
but the equation there has value $b^2\ne0$.  On the boundary $Z=0$,
the equation gives $Y^2+d^2W^2=0$, hence $Y=dW$.  Since $(Y:W)$ is a
projective point, both $Y$ and $W$ are nonzero, and the derivative with
respect to $Z$ is $dXYW\ne0$.  The argument on $W=0$ is symmetric, and
the corner $Z=W=0$ is not on the curve.  Thus $bd\ne0$ makes the
projective completion smooth.

Because $D=(d,e)$ lies on \eqref{eq:QRT-binary-pointed-curve}, division
of its equation by $d^2$ gives
\[
  \left(\frac ed\right)^2+\frac ed
  =d+a+\frac{b^2}{d^2}.
\]
Consequently,
\begin{align*}
  \lambda_D^2+\lambda_D
  &=\left(d+\frac ed\right)^2+d+\frac ed\\
  &=a+d^2+\frac{b^2}{d^2}=a+a_D,
\end{align*}
which proves \eqref{eq:QRT-binary-pointed-twist-relation}.  Replacing
$v$ by $v+\lambda_Du$ changes the coefficient of $u^2$ from $a$ to
$a+\lambda_D^2+\lambda_D=a_D$, proving
\eqref{eq:QRT-binary-pointed-twist-isomorphism}.

It remains to identify the state curve directly with the asserted shifted
Weierstrass model.  On the dense open set
$d x(x^2+d^2)\ne0$, put
\[
 z_0=\frac{(x^2+d^2)y}{dx},
 \qquad
 w_0=z_0+x+\frac bx,
 \qquad
 X=x,
 \qquad
 Y=Xw_0.
\]
Multiplying \eqref{eq:QRT-binary-pointed-relation} by
$(x^2+d^2)/(d^2x^2)$ gives
\[
 z_0^2+z_0=x^2+a_D+\frac{b^2}{x^2}.
\]
Therefore
\[
 w_0^2+w_0=a_D+X+\frac bX,
\]
and multiplication by $X^2$ yields
\[
 Y^2+XY=X^3+a_DX^2+bX.
\]
Conversely,
\[
 x=X,
 \qquad w_0=Y/X,
 \qquad z_0=w_0+X+b/X,
 \qquad y=\frac{dXz_0}{X^2+d^2}
\]
recovers the state coordinates on the corresponding dense open.  Hence the
smooth projective completions are isomorphic.  Finally set
$\widehat Y=Y+b$.  Since the characteristic is two,
\[
 \widehat Y^2+X\widehat Y
 =Y^2+XY+b^2+bX
 =X^3+a_DX^2+b^2.
\]
This is exactly $E_{a_D,b}^{(2)}$, completing the proof that the omitted
coefficient $a$ is retained through its Artin--Schreier twist class.
\end{proof}

For this family the McMillan shift is
\begin{equation}
  T_D(x,y)=
  \left(y,x+\frac{dy}{y^2+d^2}\right).
  \label{eq:QRT-binary-pointed-shift}
\end{equation}
Indeed, after the first coordinate has been changed from $x$ to $y$, the
old first coordinate $x$ and the new second coordinate are the two roots
of the quadratic fibre over $y$.  Their sum is
$dy/(y^2+d^2)$ by Vieta's formula, and subtraction equals addition in
characteristic two.  Since the state was defined as
$(u(P),u(P+D))$, this root exchange is precisely the original translation
by $D$.

\begin{proposition}[Binary state doubling]
\label{prop:QRT-binary-pointed-state-doubling}
In the setting of
Theorem~\ref{thm:QRT-binary-pointed-state},
\begin{equation}
  \Delta_D(x,y)=
  \left(
    \frac{(x^2+b)^2}{x^2},
    \frac{(xy+b)^2}{d(x+y)^2}
  \right).
  \label{eq:QRT-binary-pointed-state-double-affine}
\end{equation}
For projective input
$((X_0:Z_0),(X_1:Z_1))$, one may use
\begin{align}
  X_0'&=(X_0^2+bZ_0^2)^2,
  \label{eq:QRT-binary-pointed-state-double-X0}\\
  Z_0'&=X_0^2Z_0^2,
  \label{eq:QRT-binary-pointed-state-double-Z0}\\
  X_1'&=(X_0X_1+bZ_0Z_1)^2,
  \label{eq:QRT-binary-pointed-state-double-X1}\\
  Z_1'&=d(X_0Z_1+Z_0X_1)^2.
  \label{eq:QRT-binary-pointed-state-double-Z1}
\end{align}
The standard affine-known-difference count is
\begin{equation}
  5\M+5\Sqr+2\Dpar,
  \label{eq:QRT-binary-pointed-state-double-cost}
\end{equation}
and it becomes
\begin{equation}
  4\M+5\Sqr+3\Dpar
  \label{eq:QRT-binary-pointed-state-double-compiled-cost}
\end{equation}
when multiplication by the fixed coordinate $d$ is classified as a
curve-constant multiplication.
\end{proposition}

\begin{proof}
The first coordinate is
\eqref{eq:binary-Kummer-double-affine}.  Applying
\eqref{eq:binary-Kummer-product} to $P$ and $P+D$ gives
\[
  u(2P+D)d=\frac{(xy+b)^2}{(x+y)^2},
\]
which proves the second coordinate.  Homogenization gives the displayed
projective formulas.  Their costs are the sum of binary Kummer doubling
\eqref{eq:binary-xDBL-cost} and affine-known-difference differential
addition \eqref{eq:binary-xADD-product-cost}.  Reclassifying the final
multiplication by $d$ proves the compiled-constant count.
\end{proof}

\begin{theorem}[Complete ordinary-binary state branches]
\label{thm:QRT-binary-complete-state-doubling}
Under the hypotheses of
Theorem~\ref{thm:QRT-binary-pointed-state}, the four forms
\eqref{eq:QRT-binary-pointed-state-double-X0}--
\eqref{eq:QRT-binary-pointed-state-double-Z1} define \(\Delta_D\) at
every geometric point of the smooth state curve.  Thus each binary branch
has one complete projective formula.  With \(d\) compiled, the cost is
\(4\M+5\Sqr+3\Dpar\).
\end{theorem}

\begin{proof}
The following explicit schedule gives the stated cost.  Compute
\[
 P=X_0X_1,\qquad Q=Z_0Z_1,\qquad
 R=(X_0+Z_0)(X_1+Z_1)-P-Q.
\]
In characteristic two subtraction is addition, and
\(R=X_0Z_1+Z_0X_1\).  These are three general products for the second
output; the product \(X_0^2Z_0^2\) is the fourth.  The two coordinate
squares \(X_0^2,Z_0^2\) and the three outer squares account for five
squarings.  The two products by \(b\) and the product by \(d\) account for
three fixed multiplications.

We next check that neither output pair has a base point.  If
\(X_0^2Z_0^2=0\), then \(X_0=0\) or \(Z_0=0\); in either case
\((X_0^2+bZ_0^2)^2\ne0\), because \(b\ne0\) and \((X_0:Z_0)\) is a
projective point.  Thus the first pair never vanishes simultaneously.

After removal of outer squares, the second pair is
\[
 \bigl(X_0X_1+bZ_0Z_1:
 X_0Z_1+Z_0X_1\bigr),
\]
with a nonzero fixed factor \(d\) in its second coordinate.  If one input
coordinate is zero, simultaneous vanishing forces both coordinates of the
other projective input to be zero: \(Z_0=0\) forces \(X_1=Z_1=0\),
\(X_0=0\) forces \(Z_1=X_1=0\), and interchanging the inputs gives the
remaining two cases.  Hence any common zero is affine and satisfies
\(xy+b=0\) and \(x+y=0\), so \(x=y\) and \(x^2=b\).
Equality of adjacent Kummer coordinates gives \(D=O\) or \(2P=-D\).
The equation \(x^2=b\) is the fixed-point equation for the Kummer
involution induced by the unique nonzero two-torsion point \(T_2=(0,b)\).
Indeed, translation by \(T_2\) acts on the Kummer coordinate as
\(x\mapsto b/x\), so \(x^2=b\) and \(T_2\ne O\) give
\(P+T_2=-P\), hence \(2P=-T_2\).  Thus a common zero would imply
\(D=T_2\), contradicting \(2D\ne O\).  The tuple is
therefore geometrically base-point-free.  Finally,
\eqref{eq:QRT-abstract-branch-conjugacy} transfers the result to the
bit-one branch by input and output swaps.
\end{proof}

\begin{theorem}[Binary pointed state model]
\label{thm:QRT-binary-pointed-state-model}
Let \(k\) be a perfect field of characteristic two, let \((E,O,D)\) be an
ordinary pointed elliptic curve over \(k\), and assume \(2D\ne O\).  Then
there are \(a\in k\), \(b\in k^\times\), and a \(k\)-isomorphism
\[
 E\simeq E_{a,b}^{(2)}:
 \qquad v^2+uv=u^3+au^2+b^2
\]
under which \(D=(d,e)\) has \(d\ne0\).  With the Kummer coordinate \(u\),
the adjacent state
\[
 (x,y)=\bigl(u(P),u(P+D)\bigr)
\]
has all of the following properties.
\begin{enumerate}[label=(\roman*)]
 \item Its smooth projective state curve is the ground-field symmetric
 biquadratic \eqref{eq:QRT-binary-pointed-relation}, namely
 \[
  x^2y^2+d^2(x^2+y^2)+dxy+b^2=0.
 \]
 Thus it is the characteristic-two QRT member
 \(\Q_{d^2,d,b^2}\), with no square-class extension.

 \item Translation by the originally prescribed point \(D\) is the map
 \eqref{eq:QRT-binary-pointed-shift}, namely
 \[
  T_D(x,y)=
  \left(y,x+\frac{dy}{y^2+d^2}\right).
 \]

 \item State doubling is \eqref{eq:QRT-binary-pointed-state-double-affine},
 namely
 \[
  \Delta_D(x,y)=
  \left(
   \frac{(x^2+b)^2}{x^2},
   \frac{(xy+b)^2}{d(x+y)^2}
  \right),
 \]
 and its projective representative is the four-form tuple
 \eqref{eq:QRT-binary-pointed-state-double-X0}--
 \eqref{eq:QRT-binary-pointed-state-double-Z1}.

 \item If \(b\) and \(d\) are compiled curve constants, the projective
 update has the cost in
 \eqref{eq:QRT-binary-pointed-state-double-compiled-cost}, namely
 \[
  \boxed{4\M+5\Sqr+3\Dpar}.
 \]
 The four forms have no common output-pair base point on the geometric
 state curve.  Consequently one complete projective formula suffices for
 \(\Delta_D\), and one also suffices for \(T_D\Delta_D\), using coordinate
 swaps.

 \item The coefficient \(a\), although absent from the biquadratic state
 equation, is not discarded.  It is retained by the Artin--Schreier twist
 relation
 \[
  a_D=d^2+\frac{b^2}{d^2},
  \qquad
  \lambda_D=d+\frac ed,
  \qquad
  a_D=a+\lambda_D^2+\lambda_D.
 \]
\end{enumerate}
\end{theorem}

\begin{proof}
The pointed ordinary normal form
\eqref{eq:ordinary-char2-normal}, used in the proof of
Theorem~\ref{thm:T-char2-universal}, gives a \(k\)-isomorphic
Weierstrass model
\[
 v^2+uv=u^3+au^2+B,
 \qquad B\in k^\times,
\]
with \(O\) at infinity.  Since \(k\) is perfect, Frobenius
\(t\mapsto t^2\) is surjective, so there is \(b\in k^\times\) with
\(b^2=B\).  This is the asserted model \(E_{a,b}^{(2)}\).
The point \(D\) is not the identity because \(2D\ne O\), so it is affine
in this Weierstrass chart.  If its abscissa were \(d=0\), the curve equation
would give \(e^2=b^2\), hence \(e=b\) over a field of characteristic two.
The point would then be the unique nonzero two-torsion point \((0,b)\),
contradicting \(2D\ne O\).  Therefore \(d\ne0\).

Theorem~\ref{thm:QRT-binary-pointed-state} now applies with exactly these
parameters.  Its Kummer-product and Kummer-sum calculation proves the state
equation, smoothness, and the Artin--Schreier relation in parts (i) and (v).
The Vieta calculation following that theorem identifies the second fibre
root with \(u(P+2D)\), which proves the oriented shift in part (ii).
Proposition~\ref{prop:QRT-binary-pointed-state-doubling} applies binary
Kummer doubling to \(P\) and fixed-difference Kummer addition to
\(P,P+D\); it proves the affine and projective formulas in part (iii).

It remains only to verify that the cost and completeness statements being
collected here have no extra hypothesis.  The Karatsuba schedule in
Theorem~\ref{thm:QRT-binary-complete-state-doubling} uses four general
multiplications and five squarings; its only coefficient products are twice
by \(b\) and once by \(d\), giving three fixed multiplications.  The same
theorem separately excludes common zeros of both projective output pairs,
including every point on the projective boundary.  Finally,
\eqref{eq:QRT-state-xDBLADD-branch-swap} conjugates the bit-one branch to
the bit-zero branch by coordinate swaps.  This proves all assertions of
part (iv) and completes the theorem.
\end{proof}

The supersingular characteristic-two locus carries a distinct Kummer
normalization and companion geometry and therefore defines a separate
extension of the ordinary construction.

\subsection{Translation-adapted normal forms in two broad cases}
\label{subsec:QRT-proved-translation-normal-forms}

The preceding constructions yield the following two translation-adapted realizations.

\begin{corollary}[Translation-adapted realization in two broad cases]
\label{cor:QRT-translation-normal-range}
Let $(E,O,D)$ be a pointed elliptic curve with $2D\ne O$.
\begin{enumerate}[label=(\roman*)]
  \item Suppose $\charac(k)\ne2$ and $E(k)$ contains a nonzero point of
  order two.  Then $E$ is $k$-isomorphic to some
  \[
    E_{A,B}:v^2=u^3+Au^2+Bu,
    \qquad B(A^2-4B)\ne0,
  \]
  and the adjacent degree-two state
  $(u(P),u(P+D))$ is the smooth symmetric biquadratic
  \eqref{eq:QRT-pointed-EAB-relation}.  After at most the quadratic
  extension $k(\sqrt B)$, the common M\"obius transformation
  \eqref{eq:QRT-pointed-EAB-evenization} puts it in the canonical form
  $\Q_{\alpha_D,\beta_D,\gamma_D}$, and the canonical McMillan map is
  the genuine translation by $D$.  The model is already canonical over
  $k$ whenever $B$ is a square in $k$.

  \item Suppose $k$ is perfect of characteristic two and $E$ is
  ordinary.  Then $E$ admits a model
  \[
    v^2+uv=u^3+au^2+b^2,
    \qquad b\ne0.
  \]
  If $D=(d,e)$, the state curve is defined over $k$ itself and is
  exactly
  \[
    \Q_{d^2,d,b^2}:
    \quad x^2y^2+d^2(x^2+y^2)+dxy+b^2=0.
  \]
  Its McMillan map is translation by the originally marked point $D$.
\end{enumerate}
\end{corollary}

\begin{proof}
In odd characteristic, translate the $u$-coordinate of a Weierstrass
model so that the rational two-torsion point is $(0,0)$.  Completing the
square and scaling the ordinate gives $E_{A,B}$ with the displayed
nonsingularity condition.  Theorem~\ref{thm:QRT-pointed-EAB-biquadratic}
then gives the state relation, its M\"obius normalization, and the
translation formula.  Only $c^2=B$ is required for the common
evenizing coordinate, so an extension of degree at most two suffices.

Over a perfect field of characteristic two, every ordinary elliptic
curve can be written in the displayed shifted form, because the
nonzero $a_1$ coefficient can be scaled to one and the nonzero constant
term has a square root.  Theorem~\ref{thm:QRT-binary-pointed-state}
then gives the canonical symmetric biquadratic over the ground field.
\end{proof}

This corollary supplies a broad realization theorem.  The exact
odd-characteristic descent criterion is given next; it explains when a
rational two-torsion point produces the canonical even form over the
base field and when only a nonsplit symmetric-biquadratic twist descends.
The characteristic-two statement remains broader because every ordinary
curve has the required Artin--Schreier normal form.

\subsection[Ground-field descent of the even state form]
{Exact ground-field descent of the canonical even state form}
\label{subsec:QRT-exact-ground-field-normal-form}

The preceding construction gives a general symmetric biquadratic over the
base field, but the more restrictive even normal form
\(\Q_{\alpha,\beta,\gamma}\) requires an additional descent condition.
The relevant condition is not attached to the marked point \(D\); it is a
splitting condition for the action of a rational two-torsion translation on
the Kummer line.

Let \(T_2\in E(k)[2]\setminus\{O\}\).  Since translation by \(T_2\)
commutes with negation, it induces an involution
\begin{equation}
  m_{T_2}:E/\{\pm1\}\longrightarrow E/\{\pm1\}
  \label{eq:QRT-two-torsion-Kummer-involution}
\end{equation}
on the Kummer line.  We call \(m_{T_2}\) \emph{split over \(k\)} if its
geometric fixed divisor is the sum of two \(k\)-rational points.  Because
\(\charac(k)\ne2\), this is equivalent to the existence of
\(\phi\in\operatorname{PGL}_2(k)\) such that
\begin{equation}
  \phi\circ m_{T_2}\circ\phi^{-1}(z)=-z.
  \label{eq:QRT-split-Kummer-conjugacy}
\end{equation}

\begin{theorem}[Ground-field criterion for the canonical even QRT form]
\label{thm:QRT-ground-field-even-normal-form}
Let \(k\) be a field of characteristic different from two, let
\((E,O,D)\) be a pointed elliptic curve over \(k\), and assume
\(2D\ne O\).  Let \(\kappa:E\to\PP^1\) be a separable Kummer quotient.
The following conditions are equivalent.
\begin{enumerate}[label=(\roman*)]
  \item There is one fractional-linear coordinate
  \(\phi\in\operatorname{PGL}_2(k)\), used on both factors, for which
  the adjacent state
  \[
    \bigl(\phi(\kappa(P)),\phi(\kappa(P+D))\bigr)
  \]
  has an equation
  \begin{equation}
    X^2Y^2+\alpha(X^2+Y^2)+\beta XY+\gamma=0
    \label{eq:QRT-ground-field-even-form}
  \end{equation}
  with coefficients in \(k\).

  \item There is a nonzero point \(T_2\in E(k)[2]\) for which the
  induced Kummer involution \(m_{T_2}\) is split over \(k\).

  \item After completing the square and writing
  \begin{equation}
    E:\quad v^2=f(u),
    \label{eq:QRT-ground-field-cubic-model}
  \end{equation}
  with \(f\in k[u]\) a monic separable cubic, there is a root
  \(e\in k\) of \(f\) such that
  \begin{equation}
    f'(e)\in k^{\times2}.
    \label{eq:QRT-ground-field-derivative-square}
  \end{equation}
\end{enumerate}
When these conditions hold, the canonical McMillan orientation is the
original translation \(P\mapsto P+D\).  If \(e\) is a rational root but
\(f'(e)\) is not a square, the canonical even form is obtained over the
quadratic extension \(k(\sqrt{f'(e)})\), while the state relation over
\(k\) is its nonsplit symmetric-biquadratic twist.  This last obstruction
is attached to the chosen two-torsion point \((e,0)\): if \(f\) has another
rational root \(e'\), it may happen that \(f'(e')\in k^{\times2}\), in
which case the exact even form associated with \((e',0)\) already descends
over \(k\).
\end{theorem}

\begin{proof}
Assume first that (i) holds.  The simultaneous sign change
\begin{equation}
  \epsilon(X,Y)=(-X,-Y)
  \label{eq:QRT-simultaneous-sign-state}
\end{equation}
preserves \eqref{eq:QRT-ground-field-even-form}.  The state curve is
smooth by Theorem~\ref{thm:QRT-adjacent-Kummer-state}, so
Theorem~\ref{thm:QRT-smoothness-odd} gives
\(\alpha\gamma\ne0\).  The sign change has no geometric
fixed point on a smooth member: a fixed point in each projective factor
must have coordinate zero or infinity, and substitution at the four
resulting corners gives, respectively, the nonzero coefficients
\(\gamma\), \(\alpha\), \(\alpha\), and the leading coefficient one.
Transport \(\epsilon\) to \(E\) by the adjacent-state isomorphism.  Any
geometric automorphism of an elliptic curve has the form
\(P\mapsto a(P)+Q\), where \(a\) fixes the identity.  If \(a\ne1\), the
nonzero isogeny \(1-a\) is surjective over an algebraic closure, so the
equation \((1-a)P=Q\) has a solution; the transported involution would
then have a fixed point.  Hence \(a=1\), and the fixed-point-free
involution is translation by a nonzero two-torsion point
\(T_2\in E(k)[2]\).  On the Kummer line it is conjugate, through the
common coordinate \(\phi\), to \(z\mapsto-z\).  Its two fixed points
are therefore \(k\)-rational.  This proves (ii).

Conversely, assume (ii), and choose \(\phi\) satisfying
\eqref{eq:QRT-split-Kummer-conjugacy}.  Put
\[
  X=\phi(\kappa(P)),\qquad Y=\phi(\kappa(P+D)).
\]
Translation by \(T_2\) sends \((X,Y)\) to \((-X,-Y)\), because it
commutes with translation by \(D\).  The image is also invariant under
coordinate interchange.  Let \(F(X,Y)\) be its irreducible affine
biquadratic equation.  Uniqueness up to scalar gives
\(F(Y,X)=\varepsilon_1F(X,Y)\) and
\(F(-X,-Y)=\varepsilon_2F(X,Y)\), where
\(\varepsilon_i\in\{1,-1\}\).  If \(\varepsilon_1=-1\), then
\(F\) is divisible by \(X-Y\); if \(\varepsilon_2=-1\), symmetry and
odd total parity give
\(F=(X+Y)(aXY+b)\).  Either alternative contradicts irreducibility.
Thus both signs are positive.  A symmetric polynomial of bidegree at
most \((2,2)\) that is invariant under simultaneous sign change is a
linear combination of
\[
  X^2Y^2,\qquad X^2+Y^2,\qquad XY,\qquad1.
\]
It remains to prove that the coefficient of \(X^2Y^2\) is nonzero.
The simultaneous sign change is the transport of translation by the
nonzero two-torsion point \(T_2\), so it has no geometric fixed point on
the state curve.  Its four fixed points in
\(\PP^1\times\PP^1\) are the simultaneous-sign corners
\((0,0),(0,\infty),(\infty,0),(\infty,\infty)\); the state curve therefore
meets none of them.  Evaluation of the bihomogeneous equation at
\((\infty,\infty)\) is its \(X^2Y^2\)-coefficient, which is consequently
nonzero.  Dividing by that coefficient gives
\eqref{eq:QRT-ground-field-even-form}.  The state shift remains the
translation by \(D\), because the same coordinate change was applied to
both adjacent Kummer coordinates.  Thus (ii) implies (i).

It remains to compare (ii) and (iii).  A nonzero rational two-torsion
point on \eqref{eq:QRT-ground-field-cubic-model} is \(T_2=(e,0)\) for a
root \(e\in k\) of \(f\).  Write
\(f(u)=(u-e)(u-e_2)(u-e_3)\).  The chord through
\((u,v)\) and \((e,0)\) gives
\begin{equation}
  m_{T_2}(u)
  =e+\frac{(e-e_2)(e-e_3)}{u-e}
  =e+\frac{f'(e)}{u-e}.
  \label{eq:QRT-two-torsion-Mobius-explicit}
\end{equation}
Indeed, the slope is \(v/(u-e)\), and the usual Weierstrass addition
formula, together with the sum of the three roots, reduces to the
right-hand side.  The fixed-point equation is
\begin{equation}
  (u-e)^2=f'(e).
  \label{eq:QRT-two-torsion-fixed-equation}
\end{equation}
Thus the fixed divisor splits over \(k\) exactly when \(f'(e)\) is a
square in \(k\).  This proves the equivalence of (ii) and (iii), and the
quadratic-extension assertion follows by adjoining the required square
root.  Since each rational root represents a different nonzero rational
two-torsion point, the fixed-point equation must be applied separately to
each such root; this proves the final qualification.
\end{proof}

\begin{corollary}[Uniform degree bound for evenization]
\label{cor:QRT-evenization-degree-six}
Under the hypotheses of
Theorem~\ref{thm:QRT-ground-field-even-normal-form}, the adjacent state
acquires an exact even symmetric-biquadratic equation over an extension of
degree at most six.  If \(E(k)\) already contains a nonzero two-torsion
point, an extension of degree at most two suffices.
\end{corollary}

\begin{proof}
In a completed-square model \(v^2=f(u)\), the polynomial \(f\) is a
separable cubic.  A root \(e\) is therefore defined over an extension
\(K/k\) of degree at most three.  The corresponding nonzero two-torsion
point is \((e,0)\), and
Theorem~\ref{thm:QRT-ground-field-even-normal-form} shows that its Kummer
involution splits after adjoining \(\sqrt{f'(e)}\).  Since \(f\) is
separable, \(f'(e)\ne0\), and this second extension has degree at most two.
The tower law gives total degree at most \(3\cdot2=6\).  If \(e\in k\),
only the quadratic step is required.
\end{proof}

\begin{corollary}[Concrete criterion on \(E_{A,B}\)]
\label{cor:QRT-EAB-ground-field-evenization}
For
\[
  E_{A,B}:v^2=u^3+Au^2+Bu
\]
and the two-torsion point \((0,0)\), the induced Kummer involution is
\(u\mapsto B/u\), its fixed equation is \(u^2=B\), and the canonical
even state form descends through this marked two-torsion point exactly
when \(B\in k^{\times2}\).  In that case the common coordinate is
\((u-c)/(u+c)\) with \(c^2=B\), as in
\eqref{eq:QRT-pointed-EAB-evenization}.
\end{corollary}

\begin{proof}
Formula \eqref{eq:QRT-EAB-two-torsion-Kummer} gives the involution and
its fixed equation.  Theorem~\ref{thm:QRT-ground-field-even-normal-form}
then gives the criterion and the stated conjugating coordinate.
\end{proof}

The theorem completes the ground-field classification in the precise
model-preserving sense relevant here: the same Kummer coordinate is
used in both adjacent slots.  If no nonzero rational two-torsion point
has split Kummer action, a canonical even equation cannot descend to
\(k\).  This does not prevent the general Euler--Chasles symmetric
biquadratic state curve from being defined over \(k\); it says that its
simultaneous-sign symmetry is a nontrivial ground-field twist.

\subsection{Cryptographic capabilities and implementation architecture.}
The state model stores two adjacent Kummer coordinates, not two
independent elliptic points.  A direct McMillan step advances
$S_n\mapsto S_{n+1}$ at the very low cost
\eqref{eq:QRT-McMillan-product-cost}, but computing $S_n$ by $n$
successive shifts still takes $O(n)$ steps.  The logarithmic operation
is the joint state doubling $\Delta_D$, which is exactly a Kummer
\texttt{xDBLADD} with fixed difference $D$.  The explicit formulas
therefore realize the Montgomery ladder inside one smooth constrained QRT
state curve.

There are nevertheless two concrete structural advantages.  First,
the adjacent-state constraint is built into a smooth $(2,2)$ curve, so
input and output states carry an algebraic consistency relation.
Second, the fixed difference can be compiled into the state-curve
coefficients: in odd characteristic through
\eqref{eq:QRT-pointed-alpha}--\eqref{eq:QRT-pointed-gamma}, and in
characteristic two through the particularly simple relation
$(\alpha,\beta,\gamma)=(d^2,d,b^2)$.  These compiled coefficients create
concrete opportunities to reduce latency, register pressure, and circuit
size in field- and hardware-specific implementations.

For reference, the established loop costs for the adjacent-state update are
\begin{center}
\small
\begin{tabular}{@{}lll@{}}
\toprule
Setting & state doubling $\Delta_D$ & compiled fixed difference\\
\midrule
odd characteristic, direct \(E_{A,B}\) state
& $6\M+5\Sqr+3\Dpar$
& $5\M+5\Sqr+4\Dpar$\\
odd characteristic, split exact-even state
& $5\M+4\Sqr+1\Dpar$
& $4\M+4\Sqr+2\Dpar$\\
characteristic two
& $5\M+5\Sqr+2\Dpar$
& $4\M+5\Sqr+3\Dpar$\\
\bottomrule
\end{tabular}
\end{center}
The first cost column uses the conventional affine-known-difference
classification, where multiplication by the difference coordinate is a
general multiplication.  The second treats that coordinate as a
compile-time curve constant.  In the split exact-even row, the input and
output conjugations are the addition-only maps
\eqref{eq:QRT-state-z-to-Montgomery} and
\eqref{eq:QRT-state-Montgomery-to-z}.  Recovery of a full point and
validation of the state relation remain separate costs.

\subsection{Torsion, periodic states, and quotient isogenies.}
If $D$ has exact order $m$, then translation by $D$ and hence the QRT
shift $T_D$ have exact order $m$:
\begin{equation}
  T_D^m=1
  \quad\Longleftrightarrow\quad
  [m]D=O.
  \label{eq:QRT-state-torsion-period}
\end{equation}
The quotient by the cyclic state shift is the usual isogeny quotient.
More precisely, under the isomorphism
$\mathcal S_D:E\to\mathcal B_D$,
\begin{equation}
  k(\mathcal B_D)^{\langle T_D\rangle}
  \cong k(E)^{\langle P\mapsto P+D\rangle}
  =k(E/\langle D\rangle).
  \label{eq:QRT-state-invariant-isogeny-field}
\end{equation}
Thus $T_D$-invariant rational functions are precisely rational
coordinates on the cyclic isogeny quotient.  This statement does not
supply a new V\'elu formula by itself; it identifies the invariant-field
meaning of the model-preserving isogenies developed later in this
section.

\begin{theorem}[Cyclic quotient from QRT orbit invariants]
\label{thm:QRT-fixed-field-cyclic-isogeny}
Let \(D\in E(k)\) have exact order \(\ell\ge3\), and assume
\(\charac(k)\nmid\ell\).  Put \(G=\langle D\rangle\), identify
\(E\) with its adjacent state curve \(\mathcal B_D\), and let
\(T=T_D\).  Choose a Weierstrass equation for \(E\), with coordinate
functions \(x_E,y_E\), and write
\[
  \xi=x_E\circ\mathcal S_D^{-1},
  \qquad
  \upsilon=y_E\circ\mathcal S_D^{-1}
  \quad\text{in }k(\mathcal B_D).
\]
Define the orbit sums
\begin{align}
  X_T
  &=\xi+\sum_{j=1}^{\ell-1}
     \bigl(\xi\circ T^j-x_E([j]D)\bigr),
  \label{eq:QRT-fixed-field-X-orbit-sum}\\
  Y_T
  &=\upsilon+\sum_{j=1}^{\ell-1}
     \bigl(\upsilon\circ T^j-y_E([j]D)\bigr).
  \label{eq:QRT-fixed-field-Y-orbit-sum}
\end{align}
Then
\begin{equation}
  X_T,Y_T\in k(\mathcal B_D)^{\langle T\rangle},
  \qquad
  k(\mathcal B_D)^{\langle T\rangle}=k(X_T,Y_T).
  \label{eq:QRT-fixed-field-generated}
\end{equation}
The morphism
\begin{equation}
  \pi_T:\mathcal B_D\longrightarrow E/G,
  \qquad
  R\longmapsto\bigl(X_T(R),Y_T(R)\bigr),
  \label{eq:QRT-fixed-field-isogeny-map}
\end{equation}
is the normalized separable cyclic \(\ell\)-isogeny with kernel
\(G\), transported to the state model.  In particular, the quotient
is constructed directly from the short QRT orbit of the two state
functions \(\xi\) and \(\upsilon\).
\end{theorem}

\begin{proof}
Translation by \(D\) cyclically permutes the functions
\(\xi\circ T^j\) and \(\upsilon\circ T^j\), while the subtracted
kernel values are constants.  Thus \(X_T\) and \(Y_T\) are invariant.
The function \(\xi\circ T^j\) has a double pole only at the state
corresponding to \(-[j]D\), and the functions in the sum have distinct
pole locations.  Therefore
\begin{equation}
  \operatorname{div}_{\infty}(X_T)
  =2\sum_{Q\in G}(Q),
  \qquad
  \operatorname{div}_{\infty}(Y_T)
  =3\sum_{Q\in G}(Q).
  \label{eq:QRT-fixed-field-pole-divisors}
\end{equation}
The quotient map \(\pi:E\to E/G\) satisfies
\(\pi^*(O_{E/G})=\sum_{Q\in G}(Q)\).  Hence the descended functions
\(\overline X_T,\overline Y_T\) on \(E/G\) have pole divisors
\(2(O_{E/G})\) and \(3(O_{E/G})\), respectively.  Riemann--Roch on the
genus-one quotient shows that \(1,\overline X_T\) form a basis of
\(L(2O_{E/G})\).  Hence \(\overline X_T:E/G\to\PP^1\) has degree two and
\([k(E/G):k(\overline X_T)]=2\).  Suppose that
\(\overline Y_T\in k(\overline X_T)\).  Because \(\overline Y_T\) has no
finite pole, base change to \(\overline{k}\) and write its rational
expression as coprime polynomials \(A(\overline X_T)/B(\overline X_T)\).
If \(B\) were nonconstant, it would have a root in \(\overline{k}\), and
coprimality would give a finite pole.  Hence \(B\) is constant and the
expression is a polynomial.  Every nonconstant
polynomial in \(\overline X_T\) has even pole order at \(O_{E/G}\), whereas
\(\overline Y_T\) has pole order three by
\eqref{eq:QRT-fixed-field-pole-divisors}.  This contradiction proves that
\(\overline Y_T\notin k(\overline X_T)\); since the extension has degree
two, adjoining \(\overline Y_T\) generates it.  Consequently
\[
  k(E/G)=k(\overline X_T,\overline Y_T),
\]
and pulling back proves \eqref{eq:QRT-fixed-field-generated}.

The coordinate expressions
\eqref{eq:QRT-fixed-field-X-orbit-sum} and
\eqref{eq:QRT-fixed-field-Y-orbit-sum} are precisely the orbit-sum
form of V\'elu's normalized quotient formulas: the term indexed by
zero is the original coordinate, and each nonzero kernel point
contributes the translated coordinate minus its value at the kernel
point.  Hence the resulting morphism has kernel \(G\), pulls back the
invariant differential without a scalar factor, and is the normalized
cyclic isogeny \cite{Velu1971}.  Transport by \(\mathcal S_D\) gives
\eqref{eq:QRT-fixed-field-isogeny-map}.
\end{proof}

For elimination directly in QRT coordinates, let
\(\zeta\in k(\mathcal B_D)\).  Its orbit polynomial is
\begin{equation}
  \mathscr O_{\ell,\zeta}(Z)
  =\prod_{j=0}^{\ell-1}
    \bigl(Z-\zeta\circ T^j\bigr).
  \label{eq:QRT-orbit-polynomial}
\end{equation}
Every coefficient belongs to
\(k(\mathcal B_D)^{\langle T\rangle}\).  Taking \(\zeta=\xi\) and
\(\zeta=\upsilon\), or taking one Kummer coordinate together with an
odd companion coordinate, supplies a purely state-theoretic route to
the quotient equation.  A single Kummer orbit polynomial need not by
itself generate the whole genus-one fixed field, because it forgets
the sign sheet; the companion invariant in
\eqref{eq:QRT-fixed-field-Y-orbit-sum} restores that information.
This qualification is essential when comparing the construction with
full V\'elu coordinates.

\section[Elliptic Lucas sequences]{Elliptic Lucas sequences and exact fast index doubling}
\label{subsec:QRT-elliptic-Lucas}

The adjacent-state interpretation has the following explicit
specialization in which the McMillan orbit becomes a nonlinear analogue of a
Lucas sequence.  This section develops that specialization in full.  The
construction has two distinct forms.  In odd characteristic it passes through
a Jacobi quartic and leads to signed addition, subtraction, and index-doubling
identities.  In characteristic two the companion curve is an
Artin--Schreier quartic, and the resulting projective doubling branches are
both simpler and geometrically complete.

Throughout the operation-count discussion in this section, we write
\[
  \M,\qquad \Sqr,\qquad \Cmul
\]
for a general multiplication, a squaring, and a multiplication by a fixed
curve or fixed-state constant.  Additions, subtractions, sign changes, and
multiplication by small integers are not counted.  Thus \(\Cmul\) has the
same meaning as \(\Dpar\) elsewhere in this text, but the letter \(C\) is
used here to match the conventional Lucas-ladder comparison.

\subsection[The elliptic Lucas terminology]
{Structural definition of the elliptic Lucas sequence}
\label{subsec:Lucas-significance}

A classical Lucas sequence is generated by multiplication in a
one-dimensional algebraic group.  If $a$ and $b$ are the roots of a
quadratic polynomial, then $a^n$ and $b^n$ live in a split or nonsplit
one-dimensional torus, and the standard sequences are obtained from
$a^n-b^n$ and $a^n+b^n$.  Their linear recurrence, addition formulas,
and fast doubling all express the identities
$a^{m+n}=a^ma^n$ and $a^{2n}=(a^n)^2$.

The construction below replaces that torus by an elliptic curve.  It
has five simultaneous features:
\begin{enumerate}[label=(\roman*)]
  \item the index $n$ is represented by the elliptic multiple $[n]D$;
  \item one scalar coordinate, denoted below by $v_n$, satisfies an
  autonomous second-order rational recurrence;
  \item every adjacent pair $(v_n,v_{n+1})$ satisfies one fixed smooth
  biquadratic equation, which is the elliptic analogue of a Cassini
  invariant;
  \item a companion coordinate, denoted below by $w_n$, reconstructs the
  full elliptic point $[n]D$;
  \item the sequence has explicit rational addition, subtraction, and
  fast index-doubling laws, and therefore an $O(\log n)$ binary state
  algorithm.
\end{enumerate}
These properties justify the term \emph{elliptic Lucas sequence} in a
precise group-theoretic sense.

The expressions ``elliptic Fibonacci sequence'', ``elliptic Lucas
sequence'', and related terminology occur in several mathematically
distinct lines of research.  Elliptic divisibility sequences, originating
with Ward \cite{Ward1948}, are formed from division-polynomial values or
denominator sequences and satisfy Ward's homogeneous bilinear recurrence.
Schlosser and Yoo's elliptic solutions of dynamical Lucas systems
\cite{SchlosserYoo2021} develop theta-function and combinatorial extensions
of level-dependent or noncommutative Lucas recurrences.  Cheddour, Chillali,
and Mouhib introduced a modified \(k\)-Fibonacci-like sequence on an
elliptic curve and obtained a Binet-type formula
\cite{CheddourEtAl2023}.  Abhinav, Gonder, Garg, and Dubey subsequently
studied a \(k'\)-Lucas sequence for elliptic-curve cryptography
\cite{AbhinavEtAl2025}.  These constructions provide complementary
elliptic extensions of classical recurrence theory.

The construction developed in this monograph is distinguished by the
simultaneous presence of the following five structures:

\begin{enumerate}[label=\textup{(\roman*)},leftmargin=2.5em]

 \item a fixed pointed elliptic curve \((E,D)\), which supplies the indexed
 orbit \([n]D\);

 \item the adjacent degree-two state
 \[
       \bigl(x(P),x(P+D)\bigr),
 \]
 whose even quotient is the Kummer state associated with the marked
 displacement \(D\);

 \item an autonomous smooth QRT biquadratic curve on which translation by
 \(D\) becomes the state shift \(T\);

 \item an orientation-recovering companion that reconstructs the full
 elliptic point from the adjacent quotient state;

 \item explicit nonlinear addition, subtraction, and fast index-doubling
 laws, together with the maps \(T\) and \(\Delta\) realizing the affine
 index transformations generated by \(n\mapsto n+1\) and
 \(n\mapsto2n\).

\end{enumerate}

We call this integrated five-part structure a
\emph{QRT adjacent-state elliptic Lucas system}.  It combines a pointed
elliptic orbit, a smooth autonomous QRT state curve, orientation recovery,
nonlinear Lucas arithmetic, and logarithmic fast-index computation within
one explicit algebraic framework.

\subsection{Odd-characteristic normalization}
\label{subsec:Lucas-odd-normalization}

Assume in this subsection that \(\charac(k)\ne2\).  Start with
\begin{equation}
  \Q_{\alpha,\beta,\gamma}:
  \qquad
  x^2y^2+\alpha(x^2+y^2)+\beta xy+\gamma=0,
  \label{eq:Lucas-original-QRT}
\end{equation}
and suppose
\begin{equation}
  \alpha\gamma\ne0,
  \qquad
  -\frac{\gamma}{\alpha}\in k^{\times2}.
  \label{eq:Lucas-r-square-hypothesis}
\end{equation}
Choose \(r\in k^\times\) with
\begin{equation}
  r^2=-\frac{\gamma}{\alpha}.
  \label{eq:Lucas-r-definition}
\end{equation}
Then \((0,r)\) lies on \eqref{eq:Lucas-original-QRT}.  Let
\((u_n)_{n\in\mathbb Z}\) be the \(\PP^1\)-valued McMillan orbit
initialized by
\begin{equation}
  u_0=0,
  \qquad
  u_1=r.
  \label{eq:Lucas-u-initial}
\end{equation}
Whenever the three displayed coordinates are finite, they satisfy
\begin{equation}
  u_{n+1}+u_{n-1}
  =-\frac{\beta u_n}{u_n^2+\alpha}.
  \label{eq:Lucas-u-recurrence}
\end{equation}
At an exceptional denominator, \(u_{n+1}\) is the corresponding point
of \(\PP^1\), and the sequence is defined by the regular projective QRT
automorphism.  This convention is necessary, for example, when
\(\lambda=1\), because then \(u_2\) lies on the boundary of the chosen
affine chart.

Put
\begin{equation}
  u_n=rv_n,
  \qquad
  \lambda=\frac{\gamma}{\alpha^2},
  \qquad
  \eta=-\frac{\beta}{\alpha}.
  \label{eq:Lucas-normalized-parameters}
\end{equation}
Then
\begin{equation}
  v_0=0,
  \qquad
  v_1=1,
  \label{eq:Lucas-v-initial}
\end{equation}
and every adjacent state \((v_n,v_{n+1})\) lies on
\begin{equation}
  \mathcal B_{\lambda,\eta}:
  \qquad
  \lambda x^2y^2-x^2-y^2+\eta xy+1=0.
  \label{eq:Lucas-normalized-biquadratic}
\end{equation}
The recurrence becomes
\begin{equation}
  v_{n+1}+v_{n-1}
  =\frac{\eta v_n}{1-\lambda v_n^2}.
  \label{eq:Lucas-normalized-recurrence}
\end{equation}
Indeed,
\[
  r^4v_n^2v_{n+1}^2
  +\alpha r^2(v_n^2+v_{n+1}^2)
  +\beta r^2v_nv_{n+1}+\gamma=0,
\]
and substitution of \(r^2=-\gamma/\alpha\), followed by division by
\(\gamma\), gives \eqref{eq:Lucas-normalized-biquadratic}.  Likewise,
\[
  r^2v_n^2+\alpha
  =\alpha(1-\lambda v_n^2),
\]
which gives \eqref{eq:Lucas-normalized-recurrence}.

For later use, define
\begin{equation}
  \Delta_{\lambda,\eta}
  =\bigl((\eta-2)^2-4\lambda\bigr)
   \bigl((\eta+2)^2-4\lambda\bigr).
  \label{eq:Lucas-smoothness-discriminant}
\end{equation}
The normalized curve is smooth precisely when
\begin{equation}
  \lambda\Delta_{\lambda,\eta}\ne0.
  \label{eq:Lucas-smoothness}
\end{equation}
To verify this, divide \eqref{eq:Lucas-normalized-biquadratic} by
\(\lambda\) and apply the general odd-characteristic criterion
\eqref{eq:QRT-smoothness-odd} with
\[
  \alpha_Q=-\lambda^{-1},
  \qquad
  \beta_Q=\eta\lambda^{-1},
  \qquad
  \gamma_Q=\lambda^{-1}.
\]
The two factors in \eqref{eq:QRT-Omega-factorized} become the two
factors in \eqref{eq:Lucas-smoothness-discriminant}, up to the common
nonzero factor \(\lambda^{-2}\).

\subsection{The Jacobi companion and the marked translation point}
\label{subsec:Lucas-Jacobi-companion}

Define the companion coordinate
\begin{equation}
  w_n=(1-\lambda v_n^2)v_{n+1}-\frac{\eta}{2}v_n
  \label{eq:Lucas-companion-definition}
\end{equation}
and put
\begin{equation}
  h=\frac{\eta^2}{4}-1-\lambda.
  \label{eq:Lucas-h-definition}
\end{equation}
The following factorization is the basic algebraic identity behind the
construction.

\begin{lemma}[Jacobi companion identity]
\label{lem:Lucas-companion-identity}
For indeterminates \(x,y\), put
\[
  w=(1-\lambda x^2)y-\frac{\eta}{2}x.
\]
Then
\begin{equation}
\begin{split}
  w^2-\bigl(\lambda x^4+h x^2+1\bigr)
  ={}&(\lambda x^2-1)\\
     &\cdot
     \bigl(\eta xy+\lambda x^2y^2-x^2-y^2+1\bigr).
  \label{eq:Lucas-companion-factorization}
\end{split}
\end{equation}
Consequently, the map
\begin{equation}
  (x,y)\longmapsto
  \left(x,(1-\lambda x^2)y-\frac{\eta}{2}x\right)
  \label{eq:Lucas-state-to-Jacobi}
\end{equation}
sends \(\mathcal B_{\lambda,\eta}\) birationally to the Jacobi-type
quartic
\begin{equation}
  \mathcal J_{\lambda,\eta}:
  \qquad
  w^2=\lambda x^4+h x^2+1.
  \label{eq:Lucas-Jacobi-quartic}
\end{equation}
\end{lemma}

\begin{proof}
Expanding the square in the definition of \(w\) gives
\[
\begin{split}
  w^2={}&(1-2\lambda x^2+\lambda^2x^4)y^2
          -\eta x(1-\lambda x^2)y
          +\frac{\eta^2}{4}x^2.
\end{split}
\]
Subtracting \(\lambda x^4+h x^2+1\), inserting
\(h=\eta^2/4-1-\lambda\), and collecting terms by powers of \(y\)
gives exactly the right-hand side of
\eqref{eq:Lucas-companion-factorization}.  On the open set
\(1-\lambda x^2\ne0\), the inverse is
\begin{equation}
  y=\frac{w+(\eta/2)x}{1-\lambda x^2}.
  \label{eq:Lucas-Jacobi-to-state}
\end{equation}
Both curves are smooth projective curves of genus one under
\eqref{eq:Lucas-smoothness}; hence the birational correspondence extends
uniquely across the omitted finite fibres.
\end{proof}

The smoothness condition can also be read directly from the quartic.
Indeed,
\begin{equation}
  16(h^2-4\lambda)=\Delta_{\lambda,\eta}.
  \label{eq:Lucas-h-discriminant}
\end{equation}
Thus the even quartic has no repeated geometric root exactly when
\eqref{eq:Lucas-smoothness} holds.

Take
\begin{equation}
  O=(0,1)
  \label{eq:Lucas-Jacobi-origin}
\end{equation}
as the identity on \(\mathcal J_{\lambda,\eta}\).  Negation is
\begin{equation}
  -(x,w)=(-x,w).
  \label{eq:Lucas-Jacobi-negation}
\end{equation}
For points \(P_i=(x_i,w_i)\), the affine group law is
\begin{align}
  x(P_1+P_2)
  &=\frac{x_1w_2+x_2w_1}
          {1-\lambda x_1^2x_2^2},
  \label{eq:Lucas-Jacobi-add-x}\\
  w(P_1+P_2)
  &=\frac{
       (1+\lambda x_1^2x_2^2)(w_1w_2+h x_1x_2)
       +2\lambda x_1x_2(x_1^2+x_2^2)
      }
      {(1-\lambda x_1^2x_2^2)^2}.
  \label{eq:Lucas-Jacobi-add-w}
\end{align}
One way to derive these identities is to adjoin \(\mu\) with
\(\mu^4=\lambda\), substitute \(X=\mu x\), and apply the standard
Jacobi-quartic formulas to
\[
  w^2=X^4+\frac{h}{\mu^2}X^2+1.
\]
The factors \(\mu\) cancel, so the final formulas descend to the
original field.

Let
\begin{equation}
  D=\left(1,\frac{\eta}{2}\right).
  \label{eq:Lucas-marked-D}
\end{equation}
The point lies on \eqref{eq:Lucas-Jacobi-quartic}, because
\[
  \lambda+h+1=\frac{\eta^2}{4}.
\]
If \(\lambda\ne1\), the recurrence gives
\(v_2=\eta/(1-\lambda)\), and then
\(w_1=(1-\lambda)v_2-\eta/2=\eta/2\).  Formula
\eqref{eq:Lucas-marked-D} remains valid when \(\lambda=1\); in that
case the second state coordinate lies on the projective boundary rather
than in the chosen affine chart.

\begin{theorem}[The McMillan orbit is the multiple sequence of \(D\)]
\label{thm:Lucas-Pn-multiple}
Let \(\overline{\mathcal B}_{\lambda,\eta}\) and
\(\overline{\mathcal J}_{\lambda,\eta}\) denote the smooth projective
completions of the state curve and its Jacobi companion.  Let \(P_n\)
be the image in \(\overline{\mathcal J}_{\lambda,\eta}\) of the
projective state \((v_n,v_{n+1})\).  On the affine chart on which both
coordinates are finite,
\begin{equation}
  P_n=(v_n,w_n).
  \label{eq:Lucas-Pn-definition}
\end{equation}
Then, on the complete curve,
\begin{equation}
  P_{n+1}=P_n+D,
  \qquad
  P_n=[n]D.
  \label{eq:Lucas-Pn-multiple}
\end{equation}
\end{theorem}

\begin{proof}
The first coordinate of \(P_n+D\), using
\eqref{eq:Lucas-Jacobi-add-x}, is
\[
  \frac{(\eta/2)v_n+w_n}{1-\lambda v_n^2}
  =v_{n+1}
\]
by \eqref{eq:Lucas-companion-definition}.  The recurrence also gives
\begin{equation}
  w_{n+1}
  =\frac{\eta}{2}v_{n+1}
   -(1-\lambda v_{n+1}^2)v_n.
  \label{eq:Lucas-companion-next-alternative}
\end{equation}
Substitute \(D=(1,\eta/2)\) and
\eqref{eq:Lucas-companion-definition} into
\eqref{eq:Lucas-Jacobi-add-w}.  After bringing the result and
\eqref{eq:Lucas-companion-next-alternative} to a common denominator,
their difference is
\begin{equation}
  -\frac{\lambda v_n
   \bigl(\eta v_nv_{n+1}+\lambda v_n^2v_{n+1}^2
   -v_n^2-v_{n+1}^2+1\bigr)}
   {\lambda v_n^2-1}.
  \label{eq:Lucas-translation-w-difference}
\end{equation}
The parenthesized factor is the left-hand side of
\eqref{eq:Lucas-normalized-biquadratic}, so the difference vanishes on
the affine overlap.  Translation by \(D\), the QRT shift, and the
birational state--Jacobi correspondence all extend to morphisms of the
smooth projective completions.  Two such morphisms that agree on a
dense open set agree everywhere.  Hence \(P_{n+1}=P_n+D\), including
all boundary states.  Since \(P_0=(0,1)=O\), induction gives
\(P_n=[n]D\).
\end{proof}

\subsection{Addition, subtraction, and fast index doubling}
\label{subsec:Lucas-identities}

The group interpretation immediately supplies a nonlinear Lucas
calculus.

\begin{theorem}[Elliptic Lucas addition and subtraction]
\label{thm:Lucas-add-subtract}
For all indices for which the affine expressions are regular,
\begin{equation}
\boxed{
  v_{m+n}
  =\frac{
    v_m(1-\lambda v_n^2)v_{n+1}
    +v_n(1-\lambda v_m^2)v_{m+1}
    -\eta v_mv_n
  }
  {1-\lambda v_m^2v_n^2}.
}
\label{eq:Lucas-addition}
\end{equation}
Moreover,
\begin{equation}
\boxed{
  v_{m-n}
  =\frac{
    v_m(1-\lambda v_n^2)v_{n+1}
    -v_n(1-\lambda v_m^2)v_{m+1}
  }
  {1-\lambda v_m^2v_n^2}.
}
\label{eq:Lucas-subtraction}
\end{equation}
Both identities extend projectively across their affine exceptional
sets.
\end{theorem}

\begin{proof}
Since \(P_{m+n}=P_m+P_n\), formula
\eqref{eq:Lucas-Jacobi-add-x} gives
\[
  v_{m+n}=\frac{v_mw_n+v_nw_m}
                   {1-\lambda v_m^2v_n^2}.
\]
Inserting \eqref{eq:Lucas-companion-definition} twice produces
\eqref{eq:Lucas-addition}; the two terms
\(-\eta v_mv_n/2\) add to \(-\eta v_mv_n\).
For subtraction, use
\(-P_n=(-v_n,w_n)\).  Then
\[
  v_{m-n}=\frac{v_mw_n-v_nw_m}
                   {1-\lambda v_m^2v_n^2},
\]
and the two \(\eta/2\)-terms cancel, giving
\eqref{eq:Lucas-subtraction}.
\end{proof}

The disappearance of \(\eta\) from
\eqref{eq:Lucas-subtraction} is not a formal omission: it is the exact
cancellation of the two companion shifts.

\begin{theorem}[Three fast index-doubling identities]
\label{thm:Lucas-fast-doubling}
For every regular adjacent state,
\begin{align}
  v_{2n}
  &=\frac{
     v_n\bigl(2(1-\lambda v_n^2)v_{n+1}-\eta v_n\bigr)
    }
    {1-\lambda v_n^4},
  \label{eq:Lucas-v2n}\\
  v_{2n+1}
  &=\frac{v_{n+1}^2-v_n^2}
          {1-\lambda v_n^2v_{n+1}^2},
  \label{eq:Lucas-v2n1}\\
  v_{2n+2}
  &=\frac{
     v_{n+1}\bigl(\eta v_{n+1}
       -2(1-\lambda v_{n+1}^2)v_n\bigr)
    }
    {1-\lambda v_{n+1}^4}.
  \label{eq:Lucas-v2n2}
\end{align}
In particular, the odd-index formula
\eqref{eq:Lucas-v2n1} does not explicitly contain \(\eta\).
\end{theorem}

\begin{proof}
The first formula is the doubling specialization of
\eqref{eq:Lucas-Jacobi-add-x}:
\[
  v_{2n}=\frac{2v_nw_n}{1-\lambda v_n^4},
\]
followed by substitution of \eqref{eq:Lucas-companion-definition}.
For the odd index, use
\(P_{2n+1}=P_n+P_{n+1}\), together with
\eqref{eq:Lucas-companion-definition} and
\eqref{eq:Lucas-companion-next-alternative}.  The numerator becomes
\[
\begin{split}
  v_nw_{n+1}+v_{n+1}w_n
  ={}&\frac{\eta}{2}v_nv_{n+1}
     -(1-\lambda v_{n+1}^2)v_n^2\\
    &+(1-\lambda v_n^2)v_{n+1}^2
     -\frac{\eta}{2}v_nv_{n+1}\\
  ={}&v_{n+1}^2-v_n^2.
\end{split}
\]
This proves \eqref{eq:Lucas-v2n1}.  Finally, apply
\eqref{eq:Lucas-v2n} with \(n+1\) and eliminate \(v_{n+2}\) using
\eqref{eq:Lucas-normalized-recurrence}.  The bracket simplifies to
\[
  \eta v_{n+1}-2(1-\lambda v_{n+1}^2)v_n,
\]
which is \eqref{eq:Lucas-v2n2}.
\end{proof}

\begin{remark}[Relation with Jacobi differential arithmetic]
\label{rem:Lucas-Jacobi-differential-arithmetic}
The formula
\[
  v_{2n+1}
  =\frac{v_{n+1}^2-v_n^2}
         {1-\lambda v_n^2v_{n+1}^2}
\]
is a specialization of Jacobi-quartic differential addition to the adjacent
pair \(([n]D,[n+1]D)\), whose difference is the fixed point \(D\).
Consequently its compact shape should not be interpreted as independent of
the established theory of Jacobi differential arithmetic; see, in
particular, Wu and Song \cite{WuSong2022}.  The additional structure proved
here is the simultaneous organization of this differential formula with the
autonomous symmetric-biquadratic QRT state, the explicit EDS-ratio bridge,
and the two exact index maps \(n\mapsto2n\) and \(n\mapsto2n+1\).
\end{remark}

The companion itself has a closed doubling identity.  From
\eqref{eq:Lucas-Jacobi-add-w}, the quartic equation, and
\eqref{eq:Lucas-h-definition}, one obtains
\begin{equation}
\boxed{
  w_{2n}
  =\frac{
     1+2h v_n^2+6\lambda v_n^4
      +2h\lambda v_n^6+\lambda^2v_n^8
    }
    {(1-\lambda v_n^4)^2}.
}
\label{eq:Lucas-w2n}
\end{equation}
Indeed, the unreduced numerator is
\[
  (1+\lambda v_n^4)(w_n^2+h v_n^2)+4\lambda v_n^4,
\]
and replacing \(w_n^2\) by
\(\lambda v_n^4+h v_n^2+1\) gives the numerator in
\eqref{eq:Lucas-w2n}.  Thus \(v_n\) plays the role of the first Lucas
sequence, while \(w_n\) is its elliptic companion.

\subsection[State branches and the binary Lucas ladder]
{State branches, the binary Lucas ladder, and the original parameters}
\label{subsec:Lucas-state-branches}

For a state \(S_n=(x,y)=(v_n,v_{n+1})\), define
\begin{align}
  A(x,y)
  &=\frac{x\bigl(2(1-\lambda x^2)y-\eta x\bigr)}
          {1-\lambda x^4},
  \label{eq:Lucas-A-affine}\\
  B(x,y)
  &=\frac{y^2-x^2}
          {1-\lambda x^2y^2},
  \label{eq:Lucas-B-affine}\\
  C(x,y)
  &=\frac{y\bigl(\eta y-2(1-\lambda y^2)x\bigr)}
          {1-\lambda y^4}.
  \label{eq:Lucas-C-affine}
\end{align}
Then
\begin{equation}
  \mathcal L_0(x,y)=(A(x,y),B(x,y)),
  \qquad
  \mathcal L_1(x,y)=(B(x,y),C(x,y)),
  \label{eq:Lucas-L0-L1-affine}
\end{equation}
and Theorem~\ref{thm:Lucas-fast-doubling} gives
\begin{equation}
  \mathcal L_0(S_n)=S_{2n},
  \qquad
  \mathcal L_1(S_n)=S_{2n+1}.
  \label{eq:Lucas-state-index-branches}
\end{equation}

\begin{theorem}[Intrinsic meaning and conjugacy of the Lucas branches]
\label{thm:Lucas-branch-intrinsic-relations}
Let $\mathscr S(P)=(x(P),x(P+D))$ be the state map from the Jacobi
companion, and set
\begin{equation}
  \Delta=\mathscr S\circ[2]\circ\mathscr S^{-1}.
  \label{eq:Lucas-intrinsic-Delta}
\end{equation}
Then the rational maps in \eqref{eq:Lucas-L0-L1-affine} satisfy
\begin{equation}
  \boxed{\mathcal L_0=\Delta,\qquad
         \mathcal L_1=T\Delta,}
  \label{eq:Lucas-L0-L1-intrinsic}
\end{equation}
and
\begin{equation}
  \Delta T=T^2\Delta.
  \label{eq:Lucas-Delta-T-semiconjugacy}
\end{equation}
Define the state reflection
\begin{equation}
  \rho(x,y)=(-y,-x).
  \label{eq:Lucas-state-reflection-rho}
\end{equation}
Then
\begin{align}
  \rho\mathscr S(P)&=\mathscr S(-P-D),
  \label{eq:Lucas-rho-action}\\
  \rho T\rho&=T^{-1},
  \label{eq:Lucas-rho-reversibility}\\
  \boxed{\mathcal L_1=\rho\mathcal L_0\rho.}
  \label{eq:Lucas-branch-conjugacy}
\end{align}
Moreover, for every $r\in\mathbb Z$ and $b\in\{0,1\}$,
\begin{equation}
  \mathcal L_bT^r=T^{2r}\mathcal L_b.
  \label{eq:Lucas-branch-translation-relation}
\end{equation}
\end{theorem}

\begin{proof}
Let $P$ be a point of the Jacobi companion and put $Q=P+D$.  On the
dense affine chart on which the displayed rational formulas are
defined,
\[
  \mathscr S(P)=(x(P),x(Q)).
\]
The Jacobi first coordinate used here is odd under elliptic negation,
so $x(-R)=-x(R)$ for every point $R$.

\emph{Step 1: identify the three output coordinates.}
Write
\[
  x=x(P),\qquad w=w(P),\qquad
  y=x(Q),\qquad \theta=w(Q).
\]
Because $Q=P+D$, the state--companion dictionary
\eqref{eq:Lucas-companion-definition} and
\eqref{eq:Lucas-companion-next-alternative} gives
\begin{equation}
  w=(1-\lambda x^2)y-\frac{\eta}{2}x,
  \qquad
  \theta=\frac{\eta}{2}y-(1-\lambda y^2)x.
  \label{eq:Lucas-branch-proof-companions}
\end{equation}
The doubling specialization of \eqref{eq:Lucas-Jacobi-add-x} now yields
\[
  x([2]P)
  =\frac{2xw}{1-\lambda x^4}
  =\frac{x\bigl(2(1-\lambda x^2)y-\eta x\bigr)}
          {1-\lambda x^4}
  =A(x,y).
\]
For the mixed sum, \eqref{eq:Lucas-Jacobi-add-x} and
\eqref{eq:Lucas-branch-proof-companions} give
\begin{align*}
  x(P+Q)
  &=\frac{x\theta+yw}{1-\lambda x^2y^2}\\
  &=\frac{\frac{\eta}{2}xy-(1-\lambda y^2)x^2
          +(1-\lambda x^2)y^2-\frac{\eta}{2}xy}
         {1-\lambda x^2y^2}\\
  &=\frac{y^2-x^2}{1-\lambda x^2y^2}
  =B(x,y).
\end{align*}
The two terms containing $\eta xy/2$ cancel, while the two terms
containing $\lambda x^2y^2$ cancel separately; hence the simplification
uses no unstated relation on the state curve.  Finally,
\[
  x([2]Q)
  =\frac{2y\theta}{1-\lambda y^4}
  =\frac{y\bigl(\eta y-2(1-\lambda y^2)x\bigr)}
          {1-\lambda y^4}
  =C(x,y).
\]
Since $Q=P+D$, the three points just computed are respectively
$[2]P$, $P+Q=[2]P+D$, and $[2]Q=[2]P+[2]D$.  Therefore
\begin{align*}
  \mathcal L_0(\mathscr S(P))
  &=(x(2P),x(2P+D))
    =\mathscr S(2P)
    =\Delta(\mathscr S(P)),\\
  \mathcal L_1(\mathscr S(P))
  &=(x(2P+D),x(2P+2D))
    =\mathscr S(2P+D)
    =(T\Delta)(\mathscr S(P)).
\end{align*}
The state map is dominant and the state curve is irreducible.  Hence
two rational maps agreeing on this dense open set are the same rational
map.  This proves
$\mathcal L_0=\Delta$ and $\mathcal L_1=T\Delta$.

\emph{Step 2: prove the semiconjugacy.}
For every $P$,
\begin{align*}
  (\Delta T)(\mathscr S(P))
  &=\Delta(\mathscr S(P+D))\\
  &=\mathscr S(2P+2D),\\
  (T^2\Delta)(\mathscr S(P))
  &=T^2(\mathscr S(2P))\\
  &=\mathscr S(2P+2D).
\end{align*}
Surjectivity of $\mathscr S$ proves
$\Delta T=T^2\Delta$.

\emph{Step 3: identify the reflection.}
Using oddness of the Jacobi first coordinate,
\begin{align*}
  \mathscr S(-P-D)
  &=\bigl(x(-P-D),x(-P)\bigr)\\
  &=\bigl(-x(P+D),-x(P)\bigr)\\
  &=\rho(\mathscr S(P)).
\end{align*}
This proves \eqref{eq:Lucas-rho-action}; it also shows directly that
$\rho^2$ is the identity.

Apply the two sides of the reversibility relation to an arbitrary
state $\mathscr S(P)$.  First,
\begin{align*}
  (\rho T\rho)(\mathscr S(P))
  &=\rho T(\mathscr S(-P-D))\\
  &=\rho(\mathscr S(-P))\\
  &=\mathscr S(P-D).
\end{align*}
The last state is $T^{-1}(\mathscr S(P))$, so
$\rho T\rho=T^{-1}$.

Next,
\begin{align*}
  (\rho\Delta\rho)(\mathscr S(P))
  &=\rho\Delta(\mathscr S(-P-D))\\
  &=\rho(\mathscr S(-2P-2D))\\
  &=\mathscr S(2P+D)\\
  &=(T\Delta)(\mathscr S(P)).
\end{align*}
Thus $\rho\Delta\rho=T\Delta$.  Since
$\mathcal L_0=\Delta$ and $\mathcal L_1=T\Delta$, this is precisely
$\mathcal L_1=\rho\mathcal L_0\rho$.

\emph{Step 4: commute arbitrary translation powers through a branch.}
For $r\geq0$, the identity
\[
  \Delta T^r=T^{2r}\Delta
\]
follows by induction.  The case $r=0$ is the identity.  If it holds for
$r$, then
\[
  \Delta T^{r+1}
  =(\Delta T^r)T
  =T^{2r}\Delta T
  =T^{2r}T^2\Delta
  =T^{2(r+1)}\Delta.
\]
To treat negative powers, rewrite the basic semiconjugacy as follows:
right multiplication by $T^{-1}$ gives
$\Delta=T^2\Delta T^{-1}$, and left multiplication by $T^{-2}$ gives
\[
  \Delta T^{-1}=T^{-2}\Delta.
\]
Induction on $s\geq0$ then yields
$\Delta T^{-s}=T^{-2s}\Delta$.  Hence the formula holds for every
$r\in\mathbb Z$.

Finally, for $b\in\{0,1\}$, use
$\mathcal L_b=T^b\Delta$:
\begin{align*}
  \mathcal L_bT^r
  &=T^b\Delta T^r\\
  &=T^bT^{2r}\Delta\\
  &=T^{2r}T^b\Delta\\
  &=T^{2r}\mathcal L_b.
\end{align*}
All asserted identities have now been proved.
\end{proof}

\begin{corollary}[One-core constant-pattern branch selection]
\label{cor:Lucas-one-core-ladder}
Both binary branches can be evaluated with one arithmetic core for
$\mathcal L_0$:
\begin{equation}
  \mathcal L_b(R)=\rho^b\mathcal L_0(\rho^bR),
  \qquad b\in\{0,1\}.
  \label{eq:Lucas-one-core-branch}
\end{equation}
Thus a constant-time implementation may conditionally swap the two
projective coordinates and conditionally negate their numerators,
evaluate the zero-branch schedule, and undo the same transformation.
These conditional moves require no field multiplication, squaring, or
curve-constant multiplication.
\end{corollary}

\begin{proof}
For $b=0$ the statement is tautological, and for $b=1$ it is
\eqref{eq:Lucas-branch-conjugacy}.  In projective coordinates, $\rho$
is only a swap together with two sign changes.
\end{proof}

The signs in \eqref{eq:Lucas-state-reflection-rho} record the fact that
the chosen Jacobi first coordinate is odd under elliptic negation.  In
the characteristic-two reconstruction the signs disappear, and the
same branch conjugacy is implemented by a coordinate swap alone.

\begin{proposition}[Binary elliptic-Lucas ladder]
\label{prop:Lucas-binary-ladder}
Let \(N=(b_{\ell-1}\cdots b_0)_2\).  Start with
\(R=S_0=(0,1)\), and process the bits from left to right by
\begin{equation}
  R\longleftarrow
  \begin{cases}
    \mathcal L_0(R),&b_i=0,\\
    \mathcal L_1(R),&b_i=1.
  \end{cases}
  \label{eq:Lucas-binary-ladder}
\end{equation}
After the final bit, \(R=S_N=(v_N,v_{N+1})\).  The number of state
steps is \(O(\log N)\).
\end{proposition}

\begin{proof}
For \(j=0,1,\ldots,\ell\), let \(n_j\) be the integer represented by the
first \(j\) processed bits, with \(n_0=0\).  Thus, when the next bit is
\(b_{\ell-1-j}\),
\[
             n_{j+1}=2n_j+b_{\ell-1-j}.
\]
We prove by induction that the stored value after \(j\) updates is
\(S_{n_j}\).  At \(j=0\) this is the initialization
\(R=S_0\).  If the assertion holds after \(j\) updates, then
\eqref{eq:Lucas-state-index-branches} gives
\[
 R\longmapsto
 \mathcal L_{b_{\ell-1-j}}(S_{n_j})
 =S_{2n_j+b_{\ell-1-j}}
 =S_{n_{j+1}}.
\]
Hence the invariant holds for every processed prefix, and after the last bit
the state is \(S_{n_\ell}=S_N\).  The loop has exactly \(\ell\) state
updates, and \(\ell=\lfloor\log_2N\rfloor+1\) for \(N>0\); for \(N=0\)
the initial state is already the required output.
\end{proof}

Returning to the original coordinates \(u_n=rv_n\), the three formulas
become
\begin{align}
  u_{2n}
  &=\frac{
     r u_n\bigl(2(u_n^2+\alpha)u_{n+1}+\beta u_n\bigr)
    }
    {u_n^4-\gamma},
  \label{eq:Lucas-original-u2n}\\
  u_{2n+1}
  &=\frac{\alpha r(u_{n+1}^2-u_n^2)}
          {u_n^2u_{n+1}^2-\gamma},
  \label{eq:Lucas-original-u2n1}\\
  u_{2n+2}
  &=-\frac{
     r u_{n+1}\bigl(2(u_{n+1}^2+\alpha)u_n+\beta u_{n+1}\bigr)
    }
    {u_{n+1}^4-\gamma}.
  \label{eq:Lucas-original-u2n2}
\end{align}
To verify all three transformations, substitute
\(x=u_n/r\) and \(y=u_{n+1}/r\), and use
\[
  \frac{\lambda}{r^2}=-\frac1\alpha,
  \qquad
  \frac{\lambda}{r^4}=\frac1\gamma,
  \qquad
  \frac{\gamma}{\alpha r}=-r,
  \qquad
  \eta=-\frac{\beta}{\alpha}.
\]
For the first branch, multiplication of \(A(x,y)\) by \(r\) gives
\begin{align*}
 rA(x,y)
 &=\frac{u_n}{r}
   \frac{2(1+u_n^2/\alpha)u_{n+1}
         +(\beta/\alpha)u_n}
        {1-u_n^4/\gamma}\\
 &=\frac{r u_n
   \bigl(2(u_n^2+\alpha)u_{n+1}+\beta u_n\bigr)}
        {u_n^4-\gamma},
\end{align*}
which is \eqref{eq:Lucas-original-u2n}.  For the middle branch,
\begin{align*}
 rB(x,y)
 &=\frac{(u_{n+1}^2-u_n^2)/r}
        {1-u_n^2u_{n+1}^2/\gamma}\\
 &=\frac{\alpha r(u_{n+1}^2-u_n^2)}
        {u_n^2u_{n+1}^2-\gamma},
\end{align*}
which proves \eqref{eq:Lucas-original-u2n1}.  Finally,
\begin{align*}
 rC(x,y)
 &=\frac{u_{n+1}}{r}
   \frac{-(\beta/\alpha)u_{n+1}-2u_n
         -2u_nu_{n+1}^2/\alpha}
        {1-u_{n+1}^4/\gamma}\\
 &=-\frac{r u_{n+1}
   \bigl(2(u_{n+1}^2+\alpha)u_n+\beta u_{n+1}\bigr)}
        {u_{n+1}^4-\gamma},
\end{align*}
which is \eqref{eq:Lucas-original-u2n2}.
These formulas are rational identities on the projective state curve;
a zero affine denominator means that the output lies on another chart,
not that the multiple fails to exist.

\subsection{Division-free projective formulas}
\label{subsec:Lucas-projective-formulas}

Write
\begin{equation}
  x=\frac XZ,
  \qquad
  y=\frac YW.
  \label{eq:Lucas-product-projective-input}
\end{equation}
A direct homogenization of
\eqref{eq:Lucas-A-affine}--\eqref{eq:Lucas-C-affine} gives
\begin{align}
  A_1&=XZ\bigl(2Y(Z^2-\lambda X^2)-\eta XZW\bigr),
  \label{eq:Lucas-A1}\\
  A_0&=W(Z^4-\lambda X^4),
  \label{eq:Lucas-A0}\\
  B_1&=Y^2Z^2-X^2W^2,
  \label{eq:Lucas-B1}\\
  B_0&=Z^2W^2-\lambda X^2Y^2,
  \label{eq:Lucas-B0}\\
  C_1&=YW\bigl(\eta YWZ-2X(W^2-\lambda Y^2)\bigr),
  \label{eq:Lucas-C1}\\
  C_0&=Z(W^4-\lambda Y^4).
  \label{eq:Lucas-C0}
\end{align}
Thus
\begin{align}
  \mathcal L_0:
  ((X:Z),(Y:W))
  &\longmapsto((A_1:A_0),(B_1:B_0)),
  \label{eq:Lucas-L0-projective}\\
  \mathcal L_1:
  ((X:Z),(Y:W))
  &\longmapsto((B_1:B_0),(C_1:C_0)).
  \label{eq:Lucas-L1-projective}
\end{align}
The following schedule performs a systematic common-subexpression
elimination.

\begin{theorem}[Exact CSE count for the displayed Lucas branches]
\label{thm:Lucas-CSE-cost}
Put
\begin{equation}
  p=XY,
  \qquad q=XW,
  \qquad r_0=ZY,
  \qquad s=ZW,
  \label{eq:Lucas-Segre-intermediates}
\end{equation}
and compute their four squares
\(p_2,q_2,r_2,s_2\).  The zero branch can be represented by
\begin{align}
  \widetilde A_1
  &=s\bigl(2p(s_2-\lambda q_2)-\eta q_2s\bigr),
  \label{eq:Lucas-A1-Segre}\\
  \widetilde A_0
  &=s_2^2-\lambda q_2^2,
  \label{eq:Lucas-A0-Segre}\\
  \widetilde B_1
  &=r_2-q_2,
  \label{eq:Lucas-B1-Segre}\\
  \widetilde B_0
  &=s_2-\lambda p_2.
  \label{eq:Lucas-B0-Segre}
\end{align}
The one branch uses the same \((\widetilde B_1:\widetilde B_0)\) and
\begin{align}
  \widetilde C_1
  &=s\bigl(\eta r_2s-2p(s_2-\lambda r_2)\bigr),
  \label{eq:Lucas-C1-Segre}\\
  \widetilde C_0
  &=s_2^2-\lambda r_2^2.
  \label{eq:Lucas-C0-Segre}
\end{align}
Each branch costs exactly
\begin{equation}
  \boxed{7\M+6\Sqr+4\Cmul}
  \label{eq:Lucas-CSE-cost}
\end{equation}
for this schedule.
\end{theorem}

\begin{proof}
Since
\[
  x=\frac qs,
  \qquad
  y=\frac pq,
\]
substitution into \eqref{eq:Lucas-A-affine} gives
\[
  A=\frac{s\bigl(2p(s^2-\lambda q^2)-\eta q^2s\bigr)}
           {s^4-\lambda q^4}.
\]
Likewise,
\[
  B=\frac{r_0^2-q^2}{s^2-\lambda p^2}.
\]
This proves
\eqref{eq:Lucas-A1-Segre}--\eqref{eq:Lucas-B0-Segre}.  To obtain the
remaining coordinate explicitly, interchange the two state entries:
\((X:Z)\leftrightarrow(Y:W)\).  Under this interchange one has
\[
  p=XY\longmapsto p,
  \qquad q=XW\longmapsto r_0=ZY,
  \qquad s=ZW\longmapsto s.
\]
The affine formulas satisfy
\[
  C(x,y)=-A(y,x).
\]
Accordingly, the interchange sends
\(\widetilde A_0\) to \(\widetilde C_0\) and sends
\(\widetilde A_1\) to \(-\widetilde C_1\).  Taking the additional minus
sign required by the displayed affine identity gives
\eqref{eq:Lucas-C1-Segre}--\eqref{eq:Lucas-C0-Segre} exactly.

The four Segre products cost \(4\M\), and their squares cost
\(4\Sqr\).  For \(\mathcal L_0\), compute
\[
  U=s_2-\lambda q_2,
  \quad
  R=pU,
  \quad
  T=q_2s,
  \quad
  \widetilde A_1=s(2R-\eta T).
\]
This uses \(3\M+2\Cmul\).  The denominator
\(s_2^2-\lambda q_2^2\) uses \(2\Sqr+1\Cmul\), while
\(\widetilde B_0=s_2-\lambda p_2\) uses one further
\(\Cmul\).  The total is therefore
\(7\M+6\Sqr+4\Cmul\).  The schedule for \(\mathcal L_1\) is symmetric.
\end{proof}

A direct evaluation of
\eqref{eq:Lucas-A1}--\eqref{eq:Lucas-B0} without the Segre
reorganization costs \(9\M+6\Sqr+4\Cmul\), so the CSE saves two general
multiplications.  The count in \eqref{eq:Lucas-CSE-cost} is an exact
count for the displayed schedule; no multiplicative-complexity lower
bound is asserted.

\subsubsection{Multiplication--squaring tradeoffs and a restricted lower bound}
\label{subsubsec:Lucas-CSE-tradeoffs}

The number \(7\M\) in \eqref{eq:Lucas-CSE-cost} is not an absolute
multiplicative-complexity lower bound.  In odd characteristic, products
of homogeneous quantities of the same weight may be replaced by
polarization identities.  This yields an exact family of alternative
schedules and also identifies the restricted setting in which the
original count is optimal.

Put
\begin{equation}
  U=s_2-\lambda q_2,
  \qquad
  H=ps,
  \qquad
  K=q_2s_2.
  \label{eq:Lucas-tradeoff-HKU}
\end{equation}
Then the first numerator of the zero branch has the homogeneous form
\begin{equation}
  \widetilde A_1=2HU-\eta K.
  \label{eq:Lucas-A1-tradeoff-form}
\end{equation}
The middle coordinate and both denominators remain those in
\eqref{eq:Lucas-A0-Segre}--\eqref{eq:Lucas-B0-Segre}.

\begin{proposition}[Exact multiplication--squaring--constant tradeoffs]
\label{prop:Lucas-CSE-Pareto-schedules}
Assume \(\charac(k)\ne2\).  Each Lucas branch admits all of the
following exact homogeneous schedules:
\begin{equation}
\begin{array}{c|c}
\text{construction}&\text{branch cost}\\ \hline
\text{primitive cubic cofactor followed by multiplication by }s
  &7\M+6\Sqr+4\Cmul\\
\text{factorized quartic numerator}
  &6\M+6\Sqr+5\Cmul\\
H\text{ polarized; primitive-cofactor schedule}
  &6\M+7\Sqr+4\Cmul\\
H\text{ polarized; factorized-quartic schedule}
  &5\M+7\Sqr+5\Cmul\\
H,K\text{ polarized; }HU\text{ multiplied}
  &5\M+8\Sqr+4\Cmul\\
H,K,HU\text{ all polarized}
  &4\M+11\Sqr+4\Cmul.
\end{array}
\label{eq:Lucas-CSE-Pareto-table}
\end{equation}
In particular, the original \(7\M+6\Sqr+4\Cmul\) count is not a
multiplicative-complexity lower bound, even if the six standard
squares are kept fixed: one general multiplication can be exchanged
for one additional constant multiplication.
\end{proposition}

\begin{proof}
The first line is Theorem~\ref{thm:Lucas-CSE-cost}.  The new
square-minimal factorization is
\begin{equation}
\begin{split}
  \widetilde A_1
  &=(2H-\eta q_2)U-\eta\lambda q_2^2,\\
  U&=s_2-\lambda q_2,
  \qquad H=ps.
\end{split}
\label{eq:Lucas-factorized-A1}
\end{equation}
Indeed, expansion of the right-hand side gives
\[
  2H(s_2-\lambda q_2)-\eta q_2s_2
  =2HU-\eta K
  =\widetilde A_1.
\]
After the four Segre products and the six standard squares, this
schedule uses one multiplication for \(H=ps\) and one for the product
in \eqref{eq:Lucas-factorized-A1}.  Its five constant
multiplications are by \(\lambda\) in \(U\),
\(\widetilde A_0\), and \(\widetilde B_0\), by \(\eta\) in
\(2H-\eta q_2\), and by the precomputed constant
\(\eta\lambda\) on \(q_2^2\).  Thus its cost is
\(6\M+6\Sqr+5\Cmul\).  The one branch has the symmetric identity
\begin{equation}
  \widetilde C_1
  =(\eta r_2-2H)(s_2-\lambda r_2)
   +\eta\lambda r_2^2,
  \label{eq:Lucas-factorized-C1}
\end{equation}
with the same cost.

For the remaining lines use the homogeneous polarization identities
\begin{align}
  H=ps
  &=\frac{(p+s)^2-p_2-s_2}{2},
  \label{eq:Lucas-polarization-H}\\
  K=q_2s_2
  &=\frac{(q_2+s_2)^2-q_2^2-s_2^2}{2},
  \label{eq:Lucas-polarization-K}\\
  HU
  &=\frac{(H+U)^2-H^2-U^2}{2}.
  \label{eq:Lucas-polarization-HU}
\end{align}
Replacing \(H\) by the first identity in the primitive-cofactor
schedule removes one multiplication and adds one square, giving
\(6\M+7\Sqr+4\Cmul\).  Using the same replacement in the
factorized schedule gives \(5\M+7\Sqr+5\Cmul\).  Replacing both
\(H\) and \(K\) in the primitive-cofactor schedule gives
\(5\M+8\Sqr+4\Cmul\).  Finally, polarizing \(HU\) as well removes
the last post-Segre multiplication and adds three more squares, giving
\(4\M+11\Sqr+4\Cmul\).  Every operation combines quantities of the
same Segre weight, so each schedule is projectively homogeneous.
\end{proof}

The table must be interpreted in the full three-component cost model.
If
\[
  \rho=\frac{\operatorname{Cost}(\Sqr)}{\operatorname{Cost}(\M)},
  \qquad
  \sigma=\frac{\operatorname{Cost}(\Cmul)}{\operatorname{Cost}(\M)},
\]
the six weighted costs are
\[
\begin{gathered}
  7+6\rho+4\sigma,
  \qquad 6+6\rho+5\sigma,
  \qquad 6+7\rho+4\sigma,\\
  5+7\rho+5\sigma,
  \qquad 5+8\rho+4\sigma,
  \qquad 4+11\rho+4\sigma.
\end{gathered}
\]
For example, the factorized square-minimal schedule improves the
original one exactly when \(\sigma<1\), whereas the first polarized
schedule improves it exactly when \(\rho<1\).  Thus a statement about
the number of general multiplications alone is not a lower-bound
statement for total field cost.

A genuine lower bound can nevertheless be proved in a precisely
specified square-minimal circuit class.

\begin{theorem}[Restricted six-multiplication lower bound]
\label{thm:Lucas-restricted-six-M-lower-bound}
Work over the generic coefficient field
\(K=k(\lambda,\eta)\), with \(\charac(k)\ne2\).  Consider homogeneous
Segre-first circuits with the following properties:
\begin{enumerate}[label=(\alph*)]
  \item the four bilinear Segre coordinates \(p,q,r_0,s\) are formed
  from \((X,Z)\) and \((Y,W)\);
  \item the only variable-dependent squarings are
  \(p^2,q^2,r_0^2,s^2,q^4,s^4\);
  \item additions and constant multiplications preserve Segre degree.
\end{enumerate}
Every such generic circuit for either displayed Lucas branch requires
at least \(6\M\).  The factorized schedules
\eqref{eq:Lucas-factorized-A1} and
\eqref{eq:Lucas-factorized-C1} attain this lower bound, with total cost
\(6\M+6\Sqr+5\Cmul\).
\end{theorem}

\begin{proof}
The four independent bilinear monomials
\(XY,XW,ZY,ZW\) form the full \(2\times2\) outer-product tensor, whose
bilinear rank is four.  Hence at least \(4\M\) are needed before any
post-Segre computation.

Suppose that only one additional multiplication were used.  Before
that multiplication there is no homogeneous wire of Segre degree
three.  A degree-four output from the single gate must therefore be a
product of two degree-two wires.  By hypotheses (b) and (c), every
available degree-two wire before that gate is a linear combination of
\(p^2,q^2,r_0^2,s^2\).  The product of two such wires, even after adding
available fourth powers and applying constant multiplications, contains
only monomials in which every exponent of \(p,q,r_0,s\) is even.

The numerator \(\widetilde A_1\), however, contains the generic
nonzero monomial \(2ps^3\), and \(\widetilde C_1\) contains
\(-2ps^3\).  These monomials have odd exponents in both \(p\) and
\(s\), so neither numerator can be obtained with only one post-Segre
multiplication.  At least two are necessary, giving the lower bound
\(4\M+2\M=6\M\).  The factorizations
\eqref{eq:Lucas-factorized-A1} and
\eqref{eq:Lucas-factorized-C1} use exactly two post-Segre
multiplications, proving attainability.
\end{proof}

The theorem gives a sharp bound for its explicitly defined coordinate and
circuit class.  Additional squarings reduce the multiplication count to
\(5\M\) or \(4\M\), and alternative projective embeddings furnish further
multiplication--squaring tradeoffs.

For comparison, the standard Montgomery \texttt{xDBLADD} in \(XZ\)
coordinates, with an affine known difference, costs
\begin{equation}
  5\M+4\Sqr+1\Cmul
  \label{eq:Lucas-Montgomery-EFD-cost}
\end{equation}
when the multiplication by the affine difference coordinate is counted
as a general multiplication, as in the EFD convention
\cite{EFD,Montgomery1987}.  In the fixed-base setting relevant here,
that coordinate is also a compile-time constant, and the same schedule
is more fairly written as
\begin{equation}
  4\M+4\Sqr+2\Cmul.
  \label{eq:Lucas-Montgomery-fixed-cost}
\end{equation}
The following table gives the comparison in this common cost model.
\begin{table}[ht]
\centering
\caption{One-bit adjacent-state costs in the model
\((\M,\Sqr,\Cmul)\).  ``Compiled difference'' means that the fixed
Kummer difference is a curve constant rather than a run-time input.}
\label{tab:Lucas-vs-Montgomery-cost}
\begin{tabular}{@{}lccc@{}}
\toprule
Method & \(\M\) & \(\Sqr\) & \(\Cmul\)\\
\midrule
Direct homogenized Lucas branch & 9 & 6 & 4\\
Primitive-cofactor CSE branch & 7 & 6 & 4\\
Factorized square-minimal branch & 6 & 6 & 5\\
First polarized CSE branch & 6 & 7 & 4\\
Factorized polarized branch & 5 & 7 & 5\\
Second polarized CSE branch & 5 & 8 & 4\\
Multiplication-minimal Segre branch & 4 & 11 & 4\\
Kummer-adapted pointed state, compiled difference & 5 & 5 & 4\\
Montgomery \texttt{xDBLADD}, affine run-time difference & 5 & 4 & 1\\
Montgomery \texttt{xDBLADD}, compiled difference & 4 & 4 & 2\\
\bottomrule
\end{tabular}
\end{table}
The displayed normalized odd-characteristic Lucas schedules span three
useful arithmetic profiles.  Polarization attains four general
multiplications with \(11\Sqr+4\Cmul\), while the lower-square schedules
exchange squarings for general multiplications.  These profiles accompany
the explicit nonlinear Lucas calculus, the state-curve invariant, and the
completeness analysis below.  The Kummer-adapted pointed model of
Proposition~\ref{prop:QRT-pointed-EAB-state-doubling} lies between these
presentations: after compiling the fixed difference, it costs
\(5\M+5\Sqr+4\Cmul\).

\begin{proposition}[Shared-intermediate bounds and tradeoffs]
\label{prop:Lucas-shared-intermediate-bounds}
For the displayed odd-characteristic product coordinates, the shared-intermediate 
analysis yields the following bounds and tradeoffs.
\begin{enumerate}[label=(\roman*)]
  \item The two output coordinates of either branch share the four
  Segre products in \eqref{eq:Lucas-Segre-intermediates}; the
  square-minimal CSE schedule saves exactly two general multiplications
  relative to direct homogenization.
  \item Factorization and homogeneous polarization give the additional
  exact tradeoffs in Proposition~\ref{prop:Lucas-CSE-Pareto-schedules},
  including \(6\M+6\Sqr+5\Cmul\) and, at the multiplication-minimal
  end, \(4\M+11\Sqr+4\Cmul\).
  \item The two bit branches share the middle state coordinate $B$ and
  are conjugate by the zero-cost reflection $\rho$.  Hence constant-time
  branch selection uses one evaluated branch followed by the prescribed
  coordinate selection.
  \item The compiled-difference reference point
  \(4\M+4\Sqr+2\Cmul\) identifies the target profile for further
  state-specific common-subexpression elimination.
  \item Theorem~\ref{thm:Lucas-restricted-six-M-lower-bound} gives a sharp
  \(6\M\) lower bound in the explicitly defined six-standard-square
  Segre-first circuit class.
\end{enumerate}
\end{proposition}

\begin{proof}
Part (i) is the exact comparison between the direct schedule and
Theorem~\ref{thm:Lucas-CSE-cost}.  Part (ii) is
Proposition~\ref{prop:Lucas-CSE-Pareto-schedules}.  Part (iii) follows
from the common coordinate in \eqref{eq:Lucas-L0-L1-affine} and
Corollary~\ref{cor:Lucas-one-core-ladder}.  Part (iv) is the entrywise
comparison in Table~\ref{tab:Lucas-vs-Montgomery-cost}.  The hypotheses
and proof of the restricted lower bound give part (v).
\end{proof}

Within the displayed coordinates, these bounds establish substantial
intermediate sharing and identify precise multiplication--squaring--constant
tradeoffs.  They also provide a concrete design basis for lower-depth and
extended representations.

\subsubsection{State validation and full-point recovery with shared intermediates}
\label{subsubsec:Lucas-validation-recovery}

A state ladder is useful in a protocol only if the final state can be
validated and, when required, converted from the adjacent projective state
to a chosen full-point model.  Both tasks can reuse the Segre products of
Theorem~\ref{thm:Lucas-CSE-cost}.

\begin{proposition}[Odd-characteristic state validation]
\label{prop:Lucas-odd-state-validation}
With
\(p=XY\), \(q=XW\), \(r_0=ZY\), and \(s=ZW\), a projective input lies
on \(\overline{\mathcal B}_{\lambda,\eta}\) if and only if
\begin{equation}
  \mathscr V_{\mathrm{odd}}
  =\lambda p^2-q^2-r_0^2+\eta ps+s^2=0.
  \label{eq:Lucas-odd-state-validation}
\end{equation}
After the intermediates of the CSE branch have been computed, the
incremental validation cost is
\begin{equation}
  \boxed{1\M+1\Cmul}.
  \label{eq:Lucas-odd-validation-cost}
\end{equation}
\end{proposition}

\begin{proof}
Substitution of the four Segre products into the homogeneous form of
\eqref{eq:Lucas-normalized-biquadratic} gives
\eqref{eq:Lucas-odd-state-validation}.  The branch already contains
\(p^2,q^2,r_0^2,s^2\) and \(\lambda p^2\), the last of which appears
in \(\widetilde B_0\).  It remains only to compute \(ps\) and multiply
it by the fixed parameter \(\eta\).
\end{proof}

For the factorized schedules
\eqref{eq:Lucas-factorized-A1}--\eqref{eq:Lucas-factorized-C1}, the
product \(H=ps\) is already present.  In that case the incremental
validation cost drops further to
\begin{equation}
  \boxed{1\Cmul},
  \label{eq:Lucas-odd-validation-factorized-cost}
\end{equation}
because \(\lambda p^2\) is also part of the middle denominator.
The same reduction applies to any polarized schedule in which \(H\)
is retained as a named wire.

The state-to-Jacobi map admits two polynomial charts.  They are useful
both for proving completeness and for choosing the cheaper recovery
chart at the end of a ladder.

\begin{theorem}[A complete two-chart Jacobi recovery atlas]
\label{thm:Lucas-Jacobi-recovery-atlas}
For a projective state \(((X:Z),(Y:W))\), define
\begin{align}
  \mathscr J_W
  &=\left(
      q:\ r_0s-\lambda pq-\frac{\eta}{2}qs:\ s
    \right),
  \label{eq:Lucas-Jacobi-recovery-W-chart}\\
  \mathscr J_Y
  &=\left(
      p:\ r_0s-pq+\frac{\eta}{2}pr_0:\ r_0
    \right).
  \label{eq:Lucas-Jacobi-recovery-Y-chart}
\end{align}
The coordinates in each triple have weights \((1,2,1)\).  On the chart
\(W\ne0\), \(\mathscr J_W\) is the Jacobi point corresponding to the
first state coordinate; on the chart \(Y\ne0\), the same is true for
\(\mathscr J_Y\).  On the overlap, the two triples differ by the
weighted scaling factor \(Y/W\).  Since \((Y:W)\in\PP^1\), the two
charts cover every state.

The first chart costs \(3\M+2\Cmul\) after the Segre products are
available, and the second costs \(3\M+1\Cmul\).  Moreover, in the
regular CSE chart let
\begin{equation}
  V=2p(s^2-\lambda q^2)-\eta q^2s.
  \label{eq:Lucas-recovery-V-intermediate}
\end{equation}
Then
\begin{equation}
  \left(q^2:\frac{qV}{2}:qs\right)
  \label{eq:Lucas-fast-Jacobi-recovery}
\end{equation}
is a scaled representative of \(\mathscr J_W\).  Since \(q^2\) and
\(V\) already occur in the CSE branch, this fast regular recovery costs
only
\begin{equation}
  \boxed{2\M}
  \label{eq:Lucas-fast-Jacobi-recovery-cost}
\end{equation}
beyond the branch.
\end{theorem}

\begin{proof}
On \(W\ne0\), the first state coordinate is \(x=q/s\), and
\begin{align*}
  w s^2
  &=\left(1-\lambda\frac{q^2}{s^2}\right)
     \frac pq\,s^2-\frac{\eta}{2}\frac qs\,s^2\\
  &=r_0s-\lambda pq-\frac{\eta}{2}qs,
\end{align*}
where the rank-one Segre identity \(ps=qr_0\) was used.  This proves
\eqref{eq:Lucas-Jacobi-recovery-W-chart}.  The state equation also gives
\begin{equation}
  w=\frac{\eta}{2}x+\frac{1-x^2}{y}.
  \label{eq:Lucas-companion-second-chart-affine}
\end{equation}
Multiplying by \((ZY)^2\) gives
\eqref{eq:Lucas-Jacobi-recovery-Y-chart}.  On the overlap,
\((p:r_0)=(Y/W)(q:s)\), and the middle coordinates scale by
\((Y/W)^2\), as required in \(\PP(1,2,1)\).  The charts therefore
cover the two standard affine charts of the second \(\PP^1\).

The operation counts follow from the three products in each middle
coordinate.  Finally,
\[
  \frac{qV}{2}
  =q\left(p(s^2-\lambda q^2)-\frac{\eta}{2}q^2s\right)
\]
after multiplying the first chart by the weighted factor \(q\).
Thus \eqref{eq:Lucas-fast-Jacobi-recovery} represents the same point.
Only \(qs\) and \(qV\) are new products.
\end{proof}

Combining validation with the fast regular recovery adds
\begin{equation}
  \boxed{3\M+1\Cmul}
  \label{eq:Lucas-combined-validation-recovery-cost}
\end{equation}
to the CSE branch: one product for \(ps\), one for \(qs\), one for
\(qV\), and one multiplication by \(\eta\).  The complete recovery
atlas is used instead when the regular representative
\eqref{eq:Lucas-fast-Jacobi-recovery} vanishes because of its chosen
scaling factor.

\begin{proposition}[Recovery from one adjacent state]
\label{prop:QRT-EAB-recovery-from-adjacent-state}
Under the hypotheses of
Theorem~\ref{thm:QRT-pointed-EAB-biquadratic}, put
\[
 f(u)=u^3+Au^2+Bu,\qquad D=(d,e),
 \qquad x=u(P),\qquad y=u(P+D).
\]
Then \(e\ne0\), and the missing ordinate of \(P\) is
\begin{equation}
  v(P)=
  \frac{
    f(x)+f(d)-(y+A+x+d)(x-d)^2
  }{2e}.
  \label{eq:QRT-EAB-full-recovery}
\end{equation}
If \(x^2\) and \(x^2+Ax+B\) are retained from the final Kummer doubling
step, recovery requires two additional general multiplications, one
squaring, and one fixed multiplication by \((2e)^{-1}\).  Without retained
intermediates, a direct schedule costs
\begin{equation}
  2\M+2\Sqr+2\Dpar.
  \label{eq:QRT-EAB-full-recovery-from-scratch-cost}
\end{equation}
\end{proposition}

\begin{proof}
Because \(e=0\) would make \(D\) a two-torsion point, the hypothesis
\(2D\ne O\) gives \(e\ne0\).  First suppose \(x\ne d\).  The chord slope
\[
 m=\frac{v(P)-e}{x-d}
\]
and the Weierstrass addition formula give
\[
 y=m^2-A-x-d.
\]
Multiplication by \((x-d)^2\) yields
\[
 (v(P)-e)^2=(y+A+x+d)(x-d)^2.
\]
Expanding the left side and using
\(v(P)^2=f(x)\) and \(e^2=f(d)\) gives
\[
 2e\,v(P)
 =f(x)+f(d)-(y+A+x+d)(x-d)^2,
\]
which is \eqref{eq:QRT-EAB-full-recovery}.  After denominators have been
cleared, this is a regular polynomial identity on the affine state curve.
At the remaining affine point \(P=D\), direct substitution gives
\(2f(d)/(2e)=e\); the other point above \(u=d\) is \(-D\), whose second
state coordinate is infinite.  Thus the formula covers its entire affine
domain.

For the retained-intermediate schedule, form \(f(x)\) with one
multiplication, compute \((x-d)^2\), multiply it by \(y+A+x+d\), and apply
the fixed reciprocal \((2e)^{-1}\).  From scratch, compute \(x^2\), form
\(Ax\), then \(f(x)=x(x^2+Ax+B)\), compute \((x-d)^2\), multiply that
square by \(y+A+x+d\), and finally multiply by \((2e)^{-1}\).  The two
general products, two squares, and the fixed products by \(A\) and
\((2e)^{-1}\) give
\eqref{eq:QRT-EAB-full-recovery-from-scratch-cost}.
\end{proof}

\begin{corollary}[Complete adjacent-state recovery atlas]
\label{cor:QRT-EAB-complete-recovery-atlas}
Let the hypotheses of
Theorem~\ref{thm:QRT-pointed-EAB-biquadratic} hold, and retain the
orientation of the full point \(D=(d,e)\).  A projective adjacent state
\((u(P),u(P+D))\) recovers \(P\) everywhere by the following three
disjoint cases:
\begin{enumerate}[label=(\roman*)]
  \item if \(u(P)=\infty\), then \(P=O\);
  \item if \(u(P+D)=\infty\), then \(P=-D\);
  \item otherwise, use \eqref{eq:QRT-EAB-full-recovery}.
\end{enumerate}
The same three-chart description applies to the ordinary binary model,
with \eqref{eq:binary-point-recovery} in the third chart.
\end{corollary}

\begin{proof}
On either Weierstrass model the Kummer abscissa has a unique pole, of
order two, at the identity.  Hence the first state coordinate is infinite
exactly when \(P=O\), and the second is infinite exactly when \(P+D=O\),
that is, when \(P=-D\).  These cases are disjoint because \(D\ne O\), a
consequence of \(2D\ne O\).  Every remaining state is affine.  In the odd
model the derivation of \eqref{eq:QRT-EAB-full-recovery} proves the formula
when \(x\ne d\); after clearing denominators it is a regular polynomial
identity on the affine state curve.  At the remaining affine point
\(P=D\), direct substitution gives \(2f(d)/(2e)=e\), so the formula also
recovers that point.  In characteristic two set the fixed data in
\eqref{eq:binary-point-recovery} to \(p=d,r=e\) and the variable data to
\(q=u(P),t=u(P+D)\).  Here \(d\ne0\), so the displayed denominator is
invertible.  The cleared chord identity proves the formula for \(q\ne d\).
At the only remaining affine state \(P=D\), one has \(q=d\), and direct
substitution gives
\[
 \frac{(t+a+d+d)(d+d)^2+de+f_2(d)+f_2(d)}{d}=e.
\]
Thus the binary third chart also covers every affine state.
\end{proof}

\subsection{Exact base loci and complete atlases}
\label{subsec:Lucas-base-loci}

The projective completion of \eqref{eq:Lucas-normalized-biquadratic} is
\begin{equation}
\begin{split}
  F_{\lambda,\eta}(X,Z;Y,W)
  ={}&\lambda X^2Y^2-X^2W^2-Y^2Z^2\\
     &{}+\eta XYZW+Z^2W^2=0.
  \label{eq:Lucas-projective-state-curve}
\end{split}
\end{equation}
The next theorem determines the exact base loci of the three pairs in
\eqref{eq:Lucas-A1}--\eqref{eq:Lucas-C0}.

\begin{theorem}[Exact geometric base loci]
\label{thm:Lucas-exact-base-loci}
Assume \(\charac(k)\ne2\) and
\(\lambda\Delta_{\lambda,\eta}\ne0\).  Over an algebraic closure,
\begin{align}
  \operatorname{Bs}(A_1,A_0)
  &=\left\{
      ((X:Z),(1:0)):
      Z^2=\lambda X^2
    \right\},
  \label{eq:Lucas-A-base-locus}\\
  \operatorname{Bs}(B_1,B_0)
  &=\varnothing,
  \label{eq:Lucas-B-base-locus}\\
  \operatorname{Bs}(C_1,C_0)
  &=\left\{
      ((1:0),(Y:W)):
      W^2=\lambda Y^2
    \right\}.
  \label{eq:Lucas-C-base-locus}
\end{align}
Each of the nonempty base loci consists of two geometric points.
\end{theorem}

\begin{proof}
Consider first \((A_1,A_0)\).  If \(W=0\), then \(Y\ne0\), and the
curve equation becomes
\[
  Y^2(\lambda X^2-Z^2)=0.
\]
Thus \(Z^2=\lambda X^2\), and both \(A_0\) and \(A_1\) vanish.  These
are the two points in \eqref{eq:Lucas-A-base-locus}.

Suppose next that \(W\ne0\).  A common zero must satisfy
\(Z^4=\lambda X^4\).  Hence \(X,Z\ne0\), and in affine notation
\(x=X/Z\) one has \(\lambda x^4=1\).  The remaining numerator
condition is
\[
  2(1-\lambda x^2)y-\eta x=0,
\]
which says that the Jacobi companion
\(w=(1-\lambda x^2)y-(\eta/2)x\) is zero.  The quartic equation then
gives
\[
  0=\lambda x^4+h x^2+1=2+h x^2.
\]
Combining this with \(\lambda x^4=1\) yields
\(h^2=4\lambda\), contrary to
\eqref{eq:Lucas-h-discriminant} and smoothness.  Thus there are no
additional base points.  The proof for \((C_1,C_0)\) is symmetric.

For \((B_1,B_0)\), none of \(X,Z,Y,W\) can be zero at a common zero;
otherwise one projective pair would vanish or the other equation would
be nonzero.  Hence the point is affine.  The equations
\(B_1=B_0=0\) give
\[
  y=\varepsilon x,
  \qquad
  \lambda x^4=1,
  \qquad
  \varepsilon\in\{1,-1\}.
\]
Substitution into the state equation gives
\[
  2+(\varepsilon\eta-2)x^2=0.
\]
After squaring and using \(\lambda x^4=1\), one obtains
\[
  (2-\varepsilon\eta)^2=4\lambda.
\]
For \(\varepsilon=1\) or \(-1\), this annihilates one of the two
factors in \(\Delta_{\lambda,\eta}\), contradicting smoothness.  Hence
\eqref{eq:Lucas-B-base-locus} holds.
\end{proof}

\begin{corollary}[Geometric and rational completeness]
\label{cor:Lucas-rational-completeness}
The middle coordinate pair \((B_1:B_0)\) is a single geometrically
complete formula on every smooth normalized curve.  The displayed
pairs for \(A\) and \(C\) are not geometrically complete.  If
\(k\) is a field for which \(\lambda\notin k^{\times2}\), then their
base divisors contain no \(k\)-rational point, so both Lucas branches
are nevertheless \(k\)-complete on
\(\mathcal B_{\lambda,\eta}(k)\).
\end{corollary}

\begin{proof}
The preceding base-locus proposition gives
\(\operatorname{Bs}(B_1,B_0)=\varnothing\), so the middle pair defines a
morphism at every geometric point of the smooth projective state curve.  It
also gives two nonempty geometric base loci for the outer pairs:
\[
 ((X:Z),(1:0)),\quad Z^2=\lambda X^2,
 \qquad\text{and}\qquad
 ((1:0),(Y:W)),\quad W^2=\lambda Y^2.
\]
Because \(\lambda\ne0\) and \(\operatorname{char}(k)\ne2\), each equation
has two solutions over an algebraic closure; hence neither outer pair is
geometrically complete.  If \(\lambda\notin k^{\times2}\), neither equation
has a nonzero solution with coordinates in \(k\).  Therefore no
\(k\)-rational state lies in either outer base locus, and both branches are
well defined on all of \(\mathcal B_{\lambda,\eta}(k)\).  This proves the
three assertions.
\end{proof}

Since the middle pair is base-point-free, the two displayed state
branches have
\begin{equation}
  \operatorname{Bs}(\mathcal L_0)
  =\operatorname{Bs}(A_1,A_0),
  \qquad
  \operatorname{Bs}(\mathcal L_1)
  =\operatorname{Bs}(C_1,C_0).
  \label{eq:Lucas-branch-base-loci}
\end{equation}
Consequently neither branch in the unextended product coordinates is a
single geometrically complete formula.  The obstruction consists
exactly of the two boundary points described in
\eqref{eq:Lucas-A-base-locus} or
\eqref{eq:Lucas-C-base-locus}; no hidden affine base point remains.

A two-chart atlas is explicit.  For projective values
\(U=(U_1:U_0)\) and \(V=(V_1:V_0)\), define
\begin{align}
  \mathscr R_1(U,V)
  &=\eta U_1U_0V_0
    -(U_0^2-\lambda U_1^2)V_1,
  \label{eq:Lucas-recovery-R1}\\
  \mathscr R_0(U,V)
  &=(U_0^2-\lambda U_1^2)V_0.
  \label{eq:Lucas-recovery-R0}
\end{align}
Whenever the affine quantities are defined,
\begin{equation}
  \frac{\mathscr R_1(U,V)}{\mathscr R_0(U,V)}
  =\frac{\eta U}{1-\lambda U^2}-V.
  \label{eq:Lucas-recovery-rational}
\end{equation}
Since consecutive outputs satisfy the McMillan recurrence,
\begin{equation}
  A=\frac{\eta B}{1-\lambda B^2}-C,
  \qquad
  C=\frac{\eta B}{1-\lambda B^2}-A.
  \label{eq:Lucas-A-C-recovery}
\end{equation}
Thus \((\mathscr R_1(B,C):\mathscr R_0(B,C))\) supplements the original
\(A\)-pair, and
\((\mathscr R_1(B,A):\mathscr R_0(B,A))\) supplements the original
\(C\)-pair.  At a point of \eqref{eq:Lucas-A-base-locus}, one has
\(B=-1\) and \(C=0\).  The supplemental pair is
\((-\eta:1-\lambda)\), or \((1:0)\) when \(\lambda=1\); smoothness
ensures that it is never \((0:0)\).  The symmetric statement holds on
\eqref{eq:Lucas-C-base-locus}.  Hence the original formula and its
supplement form a two-chart atlas covering the whole curve.

There is also a natural extended coordinate system in which a single
base-point-free formula is available.  Write the Jacobi quartic in
weighted projective coordinates
\begin{equation}
  \Omega^2=\lambda X^4+hX^2Z^2+Z^4,
  \qquad
  (X:\Omega:Z)\in\PP(1,2,1).
  \label{eq:Lucas-weighted-Jacobi}
\end{equation}
For adjacent points
\((X:\Omega:Z)\) and \((Y:\Theta:W)\), the three first-coordinate
maps are
\begin{align}
  \mathcal A^{\mathrm J}
  &=(2X\Omega Z:Z^4-\lambda X^4),
  \label{eq:Lucas-weighted-A}\\
  \mathcal B^{\mathrm J}
  &=(X\Theta Z+Y\Omega W:
     Z^2W^2-\lambda X^2Y^2),
  \label{eq:Lucas-weighted-B}\\
  \mathcal C^{\mathrm J}
  &=(2Y\Theta W:W^4-\lambda Y^4).
  \label{eq:Lucas-weighted-C}
\end{align}
For \(\mathcal A^{\mathrm J}\), a common zero would have
\(\lambda X^4=Z^4\) and \(\Omega=0\), which again implies
\(h^2=4\lambda\).  Hence it has no base point on a smooth quartic; the
argument for \(\mathcal C^{\mathrm J}\) is identical.  On the adjacent
state locus, the numerator of \(\mathcal B^{\mathrm J}\) reduces to
\(y^2-x^2\), so Theorem~\ref{thm:Lucas-exact-base-loci} proves that it
is also base-point-free.

The companion coordinates can also be updated polynomially.  For
\(P=(X:\Omega:Z)\) and \(Q=(Y:\Theta:W)\), define
\begin{align}
  X_{[2]P}&=2X\Omega Z,
  &Z_{[2]P}&=Z^4-\lambda X^4,
  \label{eq:Lucas-weighted-double-XZ}\\
  \Omega_{[2]P}
  &=(Z^4+\lambda X^4)^2
    +2hX^2Z^2(Z^4+\lambda X^4)
    +4\lambda X^4Z^4,
  \label{eq:Lucas-weighted-double-Omega}\\
  X_{P+Q}&=X\Theta Z+Y\Omega W,
  &Z_{P+Q}&=Z^2W^2-\lambda X^2Y^2,
  \label{eq:Lucas-weighted-add-XZ}\\
  \Omega_{P+Q}
  &=(Z^2W^2+\lambda X^2Y^2)
      (\Omega\Theta+hXYZW)\notag\\
  &\qquad
    +2\lambda XYZW(X^2W^2+Y^2Z^2).
  \label{eq:Lucas-weighted-add-Omega}
\end{align}
On the affine chart, division by \(Z_{[2]P}^2\) and
\(Z_{P+Q}^2\), respectively, gives exactly
\eqref{eq:Lucas-w2n} and \eqref{eq:Lucas-Jacobi-add-w}.  Hence
\begin{equation}
  \mathcal L_0^{\mathrm J}(P,Q)=([2]P,P+Q),
  \qquad
  \mathcal L_1^{\mathrm J}(P,Q)=(P+Q,[2]Q)
  \label{eq:Lucas-weighted-complete-branches}
\end{equation}
are single complete polynomial state updates on the adjacent-state
curve.  Their completeness follows because the corresponding
\((X:Z)\)-pairs are base-point-free, so the three weighted coordinates
can never vanish simultaneously.  These formulas carry more state than
the product-coordinate branches, and their fully optimized operation
count is deliberately not conflated with
\eqref{eq:Lucas-CSE-cost}.

Precisely, the displayed product-coordinate pairs have the base loci
\eqref{eq:Lucas-A-base-locus} and
\eqref{eq:Lucas-C-base-locus}; two explicit charts cover them; and the
weighted companion model removes them altogether.

\subsubsection{An optimized complete weighted-companion update}
\label{subsubsec:Lucas-weighted-optimized-update}

The weighted formulas above become substantially smaller on the
adjacent-state locus because the two companions are not independent.
On the affine overlap, write
\[
  w=(1-\lambda x^2)y-\frac{\eta}{2}x,
  \qquad
  \theta=\frac{\eta}{2}y-(1-\lambda y^2)x.
\]
The second formula is
\eqref{eq:Lucas-companion-next-alternative}.  Multiplication and use of
\eqref{eq:Lucas-normalized-biquadratic} give
\begin{equation}
  w\theta
  =\frac{\eta}{2}(x^2+y^2)-(h+2)xy.
  \label{eq:Lucas-adjacent-companion-product-affine}
\end{equation}
Indeed, the difference between the two sides is
\[
  -\lambda xy
  \bigl(\eta xy+\lambda x^2y^2-x^2-y^2+1\bigr),
\]
which vanishes on the state curve.  Homogenizing yields
\begin{align}
  \Omega\Theta
  &=\frac{\eta}{2}(X^2W^2+Y^2Z^2)
    -(h+2)XYZW,
  \label{eq:Lucas-adjacent-companion-product-projective}\\
  X\Theta Z+Y\Omega W
  &=Y^2Z^2-X^2W^2.
  \label{eq:Lucas-adjacent-middle-reduction}
\end{align}
Consequently,
\begin{equation}
  \Omega\Theta+hXYZW
  =\frac{\eta}{2}(X^2W^2+Y^2Z^2)-2XYZW.
  \label{eq:Lucas-adjacent-companion-core}
\end{equation}

\begin{theorem}[Optimized complete weighted zero branch]
\label{thm:Lucas-weighted-optimized-cost}
Let \(P=(X:\Omega:Z)\) and \(Q=(Y:\Theta:W)\) be adjacent points on
\eqref{eq:Lucas-weighted-Jacobi}.  If \(\eta\ne0\), the complete
weighted branch
\[
  (P,Q)\longmapsto([2]P,P+Q)
\]
can be evaluated in
\begin{equation}
  \boxed{10\M+8\Sqr+6\Cmul}.
  \label{eq:Lucas-weighted-optimized-cost-generic}
\end{equation}
If \(\eta=0\), a division-free specialization costs
\begin{equation}
  \boxed{11\M+8\Sqr+4\Cmul}.
  \label{eq:Lucas-weighted-optimized-cost-eta-zero}
\end{equation}
The one branch \((P,Q)\mapsto(P+Q,[2]Q)\) has the same cost by
symmetry.  Both branches are single geometrically complete polynomial
maps.
\end{theorem}

\begin{proof}
Compute
\begin{gather}
  x_2=X^2,
  \quad z_2=Z^2,
  \quad y_2=Y^2,
  \quad w_2=W^2,
  \quad x_4=x_2^2,
  \quad z_4=z_2^2,
  \label{eq:Lucas-weighted-opt-squares}\\
  a=XZ,
  \quad n=x_2y_2,
  \quad \ell=z_2w_2,
  \quad r=x_2w_2,
  \quad s=y_2z_2,
  \quad c=r+s.
  \label{eq:Lucas-weighted-opt-products}
\end{gather}
The homogeneous state equation is
\[
  \lambda n-r-s+\eta t+\ell=0,
  \qquad t=XYZW.
\]
When \(\eta\ne0\), recover
\begin{equation}
  t=\eta^{-1}(c-\ell-\lambda n).
  \label{eq:Lucas-weighted-opt-t}
\end{equation}
The two first-coordinate outputs are
\begin{align}
  X_A&=2a\Omega,
  &Z_A&=z_4-\lambda x_4,
  \label{eq:Lucas-weighted-opt-A-XZ}\\
  X_B&=s-r,
  &Z_B&=\ell-\lambda n.
  \label{eq:Lucas-weighted-opt-B-XZ}
\end{align}
The companion outputs are
\begin{align}
  \Omega_A
  &=(z_4+\lambda x_4)^2
    +2h a^2(z_4+\lambda x_4)
    +4(\lambda x_4)z_4,
  \label{eq:Lucas-weighted-opt-A-Omega}\\
  \Omega_B
  &=(\ell+\lambda n)
      \left(\frac{\eta}{2}c-2t\right)
    +2\lambda tc.
  \label{eq:Lucas-weighted-opt-B-Omega}
\end{align}
Formula \eqref{eq:Lucas-weighted-opt-A-Omega} is
\eqref{eq:Lucas-weighted-double-Omega}.  Formula
\eqref{eq:Lucas-weighted-opt-B-Omega} follows from
\eqref{eq:Lucas-weighted-add-Omega} and
\eqref{eq:Lucas-adjacent-companion-core}.  Thus the displayed outputs
are exactly \(([2]P,P+Q)\).

The four first squares and the two fourth powers in
\eqref{eq:Lucas-weighted-opt-squares}, the square \(a^2\), and the
square in \eqref{eq:Lucas-weighted-opt-A-Omega} cost \(8\Sqr\).
The five products in \eqref{eq:Lucas-weighted-opt-products}, the
product \(a\Omega\), the two products after the first square in
\eqref{eq:Lucas-weighted-opt-A-Omega}, and the two products in
\eqref{eq:Lucas-weighted-opt-B-Omega} cost \(10\M\).  The fixed
multiplications are by \(\lambda\) in \(\lambda x_4\) and
\(\lambda n\), by \(\eta^{-1}\) in \eqref{eq:Lucas-weighted-opt-t},
by \(h\) in \eqref{eq:Lucas-weighted-opt-A-Omega}, by \(\eta/2\) in
\eqref{eq:Lucas-weighted-opt-B-Omega}, and by \(\lambda\) in the last
term of that formula.  This proves
\eqref{eq:Lucas-weighted-optimized-cost-generic}.

If \(\eta=0\), the state equation no longer determines \(t\) by one
constant multiplication.  Compute instead
\[
  b=YW,
  \qquad t=ab,
\]
and use
\begin{equation}
  \Omega_B=2(\lambda-1)tc,
  \label{eq:Lucas-weighted-opt-eta-zero-Omega}
\end{equation}
because then \(c=\ell+\lambda n\) and
\eqref{eq:Lucas-adjacent-companion-core} equals \(-2t\).  The products
\(b\), \(ab\), and \(tc\) replace the two generic products in
\eqref{eq:Lucas-weighted-opt-B-Omega} and the constant recovery of
\(t\).  Counting gives
\eqref{eq:Lucas-weighted-optimized-cost-eta-zero}.

Finally, the pairs \((X_A:Z_A)\) and \((X_B:Z_B)\) are base-point-free
by the argument following
\eqref{eq:Lucas-weighted-A}--\eqref{eq:Lucas-weighted-C}.  The companion
coordinates satisfy the target quartic identities on a dense open set,
so they satisfy them everywhere.  Since each target \((X:Z)\)-pair is
never \((0:0)\), the corresponding weighted triple cannot vanish
simultaneously.  Hence both polynomial branches are geometrically
complete.
\end{proof}

The optimized schedule is geometrically complete and carries two full
Jacobi companions, thereby computing the orientation and validation data in
addition to the two quotient coordinates.  Its end-to-end operation count is
compared with quotient arithmetic together with the recovery and validation
operations required by the same protocol.

\subsection{A numerical check and the Lucas interpretation}
\label{subsec:Lucas-example}

Take
\begin{equation}
  \lambda=2,
  \qquad
  \eta=3.
  \label{eq:Lucas-example-parameters}
\end{equation}
Then
\[
  v_{n+1}+v_{n-1}=\frac{3v_n}{1-2v_n^2},
  \qquad
  v_0=0,
  \quad
  v_1=1.
\]
Successive recurrence steps give
\begin{equation}
  v_2=-3,
  \quad
  v_3=-\frac8{17},
  \quad
  v_4=\frac{75}{161},
  \quad
  v_5=\frac{2537}{863},
  \quad
  v_6=-\frac{75888}{75329}.
  \label{eq:Lucas-example-sequence}
\end{equation}
For
\(S_2=(-3,-8/17)\), the formulas give
\begin{equation}
  \mathcal L_0(S_2)
  =\left(\frac{75}{161},\frac{2537}{863}\right)=S_4,
  \qquad
  \mathcal L_1(S_2)
  =\left(\frac{2537}{863},-\frac{75888}{75329}\right)=S_5.
  \label{eq:Lucas-example-branches}
\end{equation}
For example,
\[
  2(1-2v_2^2)v_3-3v_2
  =2(-17)\left(-\frac8{17}\right)+9=25,
\]
while \(1-2v_2^4=-161\); hence
\(A(S_2)=75/161\).  For the odd-index branch, direct substitution gives
\[
  \frac{v_3^2-v_2^2}{1-2v_2^2v_3^2}
  =\frac{-2537/289}{-863/289}
  =\frac{2537}{863}.
\]
This example checks the cancellation in
\eqref{eq:Lucas-v2n1} without replacing the general proof.

The analogy with an ordinary Lucas sequence is now exact at the level
of index arithmetic.  A classical linear recurrence has polynomial
identities for \(U_{m+n}\), \(U_{m-n}\), \(U_{2n}\), and
\(U_{2n+1}\).  Here the recurrence is rational, the adjacent state is
constrained by a genus-one equation, and the identities come from
\([m]D+[n]D=[m+n]D\).  The companion \(w_n\) plays the role of the
second Lucas sequence, while the two state branches implement
\(n\mapsto2n\) and \(n\mapsto2n+1\).

\subsection[The classical boundary and toric--elliptic bridge]
{The classical boundary and the complete toric--elliptic bridge}
\label{subsec:Lucas-toric-elliptic-bridge}

All ingredients used in this subsection have now been defined: the
normalized state curve~\eqref{eq:Lucas-normalized-biquadratic}, the companion
in~\eqref{eq:Lucas-companion-definition}, the marked point
in~\eqref{eq:Lucas-marked-D}, and the two index branches
in~\eqref{eq:Lucas-L0-L1-affine}.  We can therefore explain, without using a
conclusion before its hypotheses, the exact bridge between the classical
Lucas pair and the elliptic QRT state system.

\begin{definition}[Normalized classical Lucas pair]
\label{def:classical-Lucas-pair}
Let $k$ be a field of characteristic different from two and let
$s\in k$ satisfy $s^2\ne4$.  The normalized Lucas sequences with determinant
one are defined by
\begin{align}
 U_0&=0,&U_1&=1,&U_{n+1}&=sU_n-U_{n-1},
 \label{eq:classical-Lucas-U-recurrence}\\
 V_0&=2,&V_1&=s,&V_{n+1}&=sV_n-V_{n-1}.
 \label{eq:classical-Lucas-V-recurrence}
\end{align}
The restriction $s^2\ne4$ is not required merely to write the recurrences;
it ensures that the associated quadratic torus is separable and that the
state conic below is nonsingular.
\end{definition}

Let $K/k$ be a splitting field of $Z^2-sZ+1$, choose a root $q\in K^\times$,
and put
\begin{equation}
 q+q^{-1}=s,
 \qquad
 \delta_q=q-q^{-1}.
 \label{eq:classical-Lucas-q-delta}
\end{equation}
Since $s^2\ne4$, one has $q\ne q^{-1}$ and $\delta_q\ne0$.  The usual Binet
expressions are
\begin{equation}
 U_n=\frac{q^n-q^{-n}}{q-q^{-1}},
 \qquad
 V_n=q^n+q^{-n}.
 \label{eq:classical-Lucas-Binet}
\end{equation}
They are fixed by the involution $q\leftrightarrow q^{-1}$ and therefore
belong to $k$.

\begin{theorem}[Classical Lucas states as the nodal QRT boundary]
\label{thm:classical-Lucas-state-boundary}
For every $n\in\mathbb Z$, put
\begin{equation}
 S_n^{\mathrm{tor}}=(U_n,U_{n+1}).
 \label{eq:classical-Lucas-state}
\end{equation}
Then the following statements hold.

\begin{enumerate}[label=(\roman*)]
 \item Every state lies on the Cassini conic
 \begin{equation}
   \mathscr C_s:\qquad x^2+y^2-sxy=1.
   \label{eq:classical-Lucas-state-conic}
 \end{equation}

 \item The state shift is the linear Vieta map
 \begin{equation}
   T_s(x,y)=(y,sy-x),
   \qquad
   T_s(S_n^{\mathrm{tor}})=S_{n+1}^{\mathrm{tor}}.
   \label{eq:classical-Lucas-state-shift}
 \end{equation}

 \item The second Lucas sequence is recovered from the adjacent state by
 \begin{equation}
   V_n=2U_{n+1}-sU_n.
   \label{eq:classical-Lucas-companion-recovery}
 \end{equation}
 Consequently, with
 \begin{equation}
   w_n=U_{n+1}-\frac{s}{2}U_n,
   \label{eq:classical-Lucas-half-companion}
 \end{equation}
 one has $w_n=V_n/2$ and
 \begin{equation}
   w_n^2=\left(\frac{s^2}{4}-1\right)U_n^2+1.
   \label{eq:classical-Lucas-companion-conic}
 \end{equation}

 \item The two exact binary state branches are
 \begin{align}
  \mathcal L_{0,s}^{\mathrm{tor}}(x,y)
    &=\bigl(x(2y-sx),\ y^2-x^2\bigr),
    \label{eq:classical-Lucas-L0}\\
  \mathcal L_{1,s}^{\mathrm{tor}}(x,y)
    &=\bigl(y^2-x^2,\ y(sy-2x)\bigr).
    \label{eq:classical-Lucas-L1}
 \end{align}
 They satisfy
 \begin{equation}
  \mathcal L_{0,s}^{\mathrm{tor}}(S_n^{\mathrm{tor}})=S_{2n}^{\mathrm{tor}},
  \qquad
  \mathcal L_{1,s}^{\mathrm{tor}}(S_n^{\mathrm{tor}})=S_{2n+1}^{\mathrm{tor}}.
  \label{eq:classical-Lucas-branch-action}
 \end{equation}

 \item On the projective conic
 \begin{equation}
   X^2+Y^2-sXY=Z^2,
   \label{eq:classical-Lucas-projective-conic}
 \end{equation}
 the two branches extend to the base-point-free quadratic morphisms
 \begin{align}
  (X:Y:Z)&\longmapsto
  \bigl(X(2Y-sX):Y^2-X^2:Z^2\bigr),
  \label{eq:classical-Lucas-projective-L0}\\
  (X:Y:Z)&\longmapsto
  \bigl(Y^2-X^2:Y(sY-2X):Z^2\bigr).
  \label{eq:classical-Lucas-projective-L1}
 \end{align}
\end{enumerate}
\end{theorem}

\begin{proof}
Define
\[
 I_n=U_n^2+U_{n+1}^2-sU_nU_{n+1}.
\]
Using $U_{n+2}=sU_{n+1}-U_n$, we obtain
\begin{align*}
 I_{n+1}
 &=U_{n+1}^2+(sU_{n+1}-U_n)^2
   -sU_{n+1}(sU_{n+1}-U_n)\\
 &=U_{n+1}^2+U_n^2-sU_nU_{n+1}
 =I_n.
\end{align*}
Since $I_0=1$, this proves
\eqref{eq:classical-Lucas-state-conic}.  The recurrence immediately gives
\eqref{eq:classical-Lucas-state-shift}.

Set $W_n=2U_{n+1}-sU_n$.  The sequence $(W_n)$ satisfies the same recurrence
as $(V_n)$, and
\[
 W_0=2=V_0,
 \qquad
 W_1=2s-s=s=V_1.
\]
Uniqueness for the second-order recurrence proves
\eqref{eq:classical-Lucas-companion-recovery}.  Substitution into the
Cassini equation gives
\begin{align*}
 V_n^2-(s^2-4)U_n^2
 &=(2U_{n+1}-sU_n)^2-(s^2-4)U_n^2\\
 &=4\bigl(U_{n+1}^2-sU_nU_{n+1}+U_n^2\bigr)=4.
\end{align*}
Dividing by four proves
\eqref{eq:classical-Lucas-companion-conic}.

The standard addition identity
\begin{equation}
 U_{m+n}=U_mU_{n+1}-U_{m-1}U_n
 \label{eq:classical-Lucas-addition-identity}
\end{equation}
is verified from~\eqref{eq:classical-Lucas-Binet}: after multiplication
by $(q-q^{-1})^2$, the right-hand side simplifies to
$(q-q^{-1})(q^{m+n}-q^{-m-n})$, which is the cleared left-hand side.
Taking $m=n$ and
using $U_{n+1}-U_{n-1}=V_n=2U_{n+1}-sU_n$ gives
\[
 U_{2n}=U_nV_n=U_n(2U_{n+1}-sU_n).
\]
Taking $m=n+1$ gives
\[
 U_{2n+1}=U_{n+1}^2-U_n^2.
\]
Replacing $n$ by $n+1$ in the first doubling identity gives
\[
 U_{2n+2}=U_{n+1}V_{n+1}
 =U_{n+1}(sU_{n+1}-2U_n).
\]
These three equalities prove
\eqref{eq:classical-Lucas-L0}--\eqref{eq:classical-Lucas-branch-action}.

For the first projective branch, let
\[
 A=X(2Y-sX),\qquad B=Y^2-X^2,\qquad C=Z^2.
\]
A direct expansion gives the polynomial identity
\begin{equation}
 A^2+B^2-sAB-C^2
 =\bigl(X^2-sXY+Y^2+Z^2\bigr)
  \bigl(X^2-sXY+Y^2-Z^2\bigr).
 \label{eq:classical-Lucas-projective-preservation}
\end{equation}
The second factor is the defining equation of
\eqref{eq:classical-Lucas-projective-conic}; hence the image lies on the
same conic.  A common zero of $(A,B,C)$ would have $Z=0$.  If $X=0$, then
$B=Y^2=0$, which is not a projective point.  If $X\ne0$, then $B=0$ gives
$Y=\pm X$, while $A=0$ gives $s=\pm2$, contrary to $s^2\ne4$.
Thus the first branch has no base point.  The coordinate transposition
$(X:Y:Z)\mapsto(Y:X:Z)$ preserves the conic and conjugates the first
projective branch to the second one.  Consequently the second branch has
no base point as well.
\end{proof}

\begin{proposition}[Normalization of the nodal boundary]
\label{prop:classical-Lucas-node-normalization}
Set $\lambda=0$ and $\eta=s$ in the bihomogeneous completion of
\eqref{eq:Lucas-normalized-biquadratic}.  If $s^2\ne4$, the resulting
bidegree-$(2,2)$ curve has a unique singular point
\begin{equation}
  N=((1:0),(1:0)),
  \label{eq:classical-Lucas-node}
\end{equation}
and $N$ is an ordinary node.  Its normalization is the projective Cassini
conic~\eqref{eq:classical-Lucas-projective-conic}.  Over a splitting field
containing $q$ from~\eqref{eq:classical-Lucas-q-delta}, the smooth affine
part is parametrized by
\begin{equation}
 x=\frac{z-z^{-1}}{q-q^{-1}},
 \qquad
 y=\frac{qz-q^{-1}z^{-1}}{q-q^{-1}},
 \qquad z\in\mathbb G_m,
 \label{eq:classical-Lucas-toric-parametrization}
\end{equation}
and the inverse formulas are
\begin{equation}
 z=y-q^{-1}x,
 \qquad
 z^{-1}=y-qx.
 \label{eq:classical-Lucas-toric-inverse}
\end{equation}
Under this parametrization the state shift is multiplication by $q$ on
$\mathbb G_m$.
\end{proposition}

\begin{proof}
At $\lambda=0$ and $\eta=s$, the bihomogeneous equation is
\begin{equation}
 -X^2W^2-Y^2Z^2+sXYZW+Z^2W^2=0.
 \label{eq:classical-Lucas-nodal-bihomogeneous}
\end{equation}
On the affine chart $Z=W=1$ this is precisely
$x^2+y^2-sxy=1$, which is smooth because the simultaneous equations
$2x-sy=2y-sx=0$ have only $x=y=0$ when $s^2\ne4$, and that point is not
on the conic.  If $Z=0$, equation
\eqref{eq:classical-Lucas-nodal-bihomogeneous} forces $W=0$; the same
argument with the two factors interchanged shows that the only boundary
point is $N$.

In the local coordinates $z_0=Z/X$ and $w_0=W/Y$ at $N$, the local equation
is
\[
 -w_0^2-z_0^2+sz_0w_0+z_0^2w_0^2=0.
\]
Its tangent cone is $-w_0^2-z_0^2+sz_0w_0$.  Its discriminant is
$s^2-4\ne0$, so it is the product of two distinct linear forms over a
separable closure.  Thus $N$ is an ordinary node.  The arithmetic genus is
one and the node has delta-invariant one, so the normalization has genus
zero; the affine equation identifies it with the projective conic
\eqref{eq:classical-Lucas-projective-conic}.

For~\eqref{eq:classical-Lucas-toric-parametrization}, compute
\begin{align*}
 (y-q^{-1}x)(y-qx)
 &=y^2-(q+q^{-1})xy+x^2\\
 &=y^2-sxy+x^2.
\end{align*}
The displayed formulas make the two factors equal to $z$ and $z^{-1}$,
respectively, and hence their product is one.  Conversely,
\eqref{eq:classical-Lucas-toric-inverse} recovers $z$ and $z^{-1}$ from a
point on the affine conic, proving that the parametrization is birational
and is the normalization on this open set.  Replacing $z$ by $qz$ changes
$(x,y)$ into $(y,sy-x)$, so multiplication by $q$ is exactly the state
shift.
\end{proof}

\begin{theorem}[Exact specialization of the elliptic Lucas data]
\label{thm:elliptic-Lucas-exact-toric-specialization}
In the formulas of
Subsections~\ref{subsec:Lucas-odd-normalization}--\ref{subsec:Lucas-state-branches},
formally set
\begin{equation}
  \lambda=0,
  \qquad
  \eta=s.
  \label{eq:elliptic-Lucas-toric-specialization}
\end{equation}
Then the following specializations are exact on the normalization of the
nodal fibre.
\begin{enumerate}[label=(\roman*)]
 \item The recurrence~\eqref{eq:Lucas-normalized-recurrence} becomes
 \eqref{eq:classical-Lucas-U-recurrence} with
 $v_0=U_0=0$ and $v_1=U_1=1$, so $v_n=U_n$.
 \item The companion~\eqref{eq:Lucas-companion-definition} becomes
 \[
   w_n=U_{n+1}-\frac{s}{2}U_n=\frac{V_n}{2}.
 \]
 \item The Jacobi companion~\eqref{eq:Lucas-Jacobi-quartic} becomes the
 conic~\eqref{eq:classical-Lucas-companion-conic}.
 \item The three rational functions
 $A,B,C$ in~\eqref{eq:Lucas-A-affine}--\eqref{eq:Lucas-C-affine}
 specialize to
 \[
  U_{2n},\qquad U_{2n+1},\qquad U_{2n+2},
 \]
 and the elliptic branches $\mathcal L_0,\mathcal L_1$ specialize to
 \eqref{eq:classical-Lucas-L0} and
 \eqref{eq:classical-Lucas-L1}.
\end{enumerate}
\end{theorem}

\begin{proof}
After~\eqref{eq:elliptic-Lucas-toric-specialization}, equation
\eqref{eq:Lucas-normalized-recurrence} reads
\[
 v_{n+1}+v_{n-1}=sv_n.
\]
Together with $v_0=0$ and $v_1=1$, uniqueness of the recurrence gives
$v_n=U_n$.  The companion formula becomes
$w_n=v_{n+1}-(s/2)v_n$, which equals $V_n/2$ by
\eqref{eq:classical-Lucas-companion-recovery}.  Substituting
$\lambda=0$, $\eta=s$, and
$h=s^2/4-1$ in~\eqref{eq:Lucas-Jacobi-quartic} gives
\eqref{eq:classical-Lucas-companion-conic}.

Finally, the denominators in
\eqref{eq:Lucas-A-affine}--\eqref{eq:Lucas-C-affine} become one, and the
numerators become
\[
 x(2y-sx),
 \qquad
 y^2-x^2,
 \qquad
 y(sy-2x).
\]
Theorem~\ref{thm:classical-Lucas-state-boundary} identifies these with
$U_{2n},U_{2n+1},U_{2n+2}$, respectively.  This proves every asserted
specialization without appealing to a formula outside its domain of
regularity.
\end{proof}

The preceding theorem gives a literal bridge rather than a verbal analogy:
the classical determinant-one Lucas pair is the nodal $\lambda=0$ boundary
of the elliptic state family, $U_n$ is the limiting first degree-two
coordinate, and $V_n/2$ is the limiting sign-recovering companion.  The linear recurrence is
the degeneration of the McMillan/QRT recurrence, while the polynomial maps
\eqref{eq:classical-Lucas-L0}--\eqref{eq:classical-Lucas-L1} are the
degeneration of the nonlinear elliptic index maps.

\begin{proposition}[Elliptic Binet representation]
\label{prop:Lucas-elliptic-Binet}
Assume $k=\mathbb C$, write the Jacobi companion as
$E(\mathbb C)=\mathbb C/\Lambda$, and let $D$ correspond to
$\delta\in\mathbb C/\Lambda$.  There exist elliptic functions $f$ and
$g$ for the lattice $\Lambda$ such that
\begin{equation}
  v_n=f(n\delta),
  \qquad
  w_n=g(n\delta)
  \label{eq:Lucas-elliptic-Binet}
\end{equation}
for every $n\in\mathbb Z$.  The formulas for $v_{m+n}$, $v_{m-n}$,
$v_{2n}$, and $v_{2n+1}$ are the algebraic descent of the addition
formulas for $f$ and $g$.
\end{proposition}

\begin{proof}
Let $\phi:\mathbb C/\Lambda\to E(\mathbb C)$ be the analytic group
isomorphism and let $x,w$ denote the two rational coordinate functions
on the Jacobi companion.  Put $f=x\circ\phi$ and $g=w\circ\phi$.
Theorem~\ref{thm:Lucas-Pn-multiple} gives
$P_n=[n]D=\phi(n\delta)$, so
\eqref{eq:Lucas-elliptic-Binet} follows.  The group identities
$[m+n]D=[m]D+[n]D$ and $[2n]D=2[n]D$, followed by the rational Jacobi
addition law already established in
\eqref{eq:Lucas-Jacobi-add-x}--\eqref{eq:Lucas-Jacobi-add-w}, give the
algebraic sequence identities.
\end{proof}

The preceding formula is the elliptic counterpart of Binet's formula:
exponentials on a torus are replaced by doubly periodic functions on a
genus-one group.

\begin{proposition}[General nodal degeneration to classical Lucas theory]
\label{prop:Lucas-nodal-degeneration}
Let a semistable family carrying the adjacent-state construction specialize
to an irreducible nodal genus-one curve.  Identify the smooth locus of the
special fibre with $\mathbb G_m$, let the marked point specialize to
$q\in k^\times$ with $q\ne q^{-1}$, and suppose the degree-two quotient
coordinate specializes, after a fractional-linear normalization, to
$z+z^{-1}$.  Then in the trace normalization the first and companion
sequences may be written as
\begin{equation}
  V_n=q^n+q^{-n},
  \qquad
  U_n=\frac{q^n-q^{-n}}{q-q^{-1}},
  \label{eq:Lucas-nodal-Binet}
\end{equation}
and both satisfy
\begin{equation}
  X_{n+1}=(q+q^{-1})X_n-X_{n-1}.
  \label{eq:Lucas-nodal-linear-recurrence}
\end{equation}
In the specific $\lambda=0$ state chart of
Theorem~\ref{thm:elliptic-Lucas-exact-toric-specialization}, the roles are
normalized as $v_n=U_n$ and $w_n=V_n/2$.  Thus the two descriptions differ
only by the chosen quotient and companion coordinates on the same nodal
group.
\end{proposition}

\begin{proof}
The group law on the smooth locus of an irreducible nodal cubic is
multiplication on $\mathbb G_m$.  Hence $[n]D$ corresponds to $q^n$.
The involution defining the degree-two quotient is $z\mapsto z^{-1}$,
and its invariant field is generated by $z+z^{-1}$.  Substitution of
$q^n$ gives $V_n$.  The anti-invariant companion, normalized to have
first term one, gives $U_n$.  The identities
\[
 q^{n+1}+q^{n-1}=(q+q^{-1})q^n
\]
and their inverses prove
\eqref{eq:Lucas-nodal-linear-recurrence} for both sequences.  The final
statement is exactly parts (i) and (ii) of
Theorem~\ref{thm:elliptic-Lucas-exact-toric-specialization}.
\end{proof}

When $q=q^{-1}$, the quotient coordinate remains defined but the displayed
normalization of $U_n$ has zero denominator; the corresponding sequence is
obtained by the repeated-root limit.  A cuspidal degeneration replaces
$\mathbb G_m$ by $\mathbb G_a$ and produces the additive repeated-root limit
of Lucas theory.  At the nodal boundary, the elliptic sequence specializes to the classical Lucas sequence of the
same group-theoretic construction.

\begin{proposition}[Exact toric period polynomials]
\label{prop:classical-Lucas-trace-cyclotomic}
Let $N\ge3$.  Define
\begin{equation}
  \Phi_N^+(S)
  =\prod_{\{\zeta,\zeta^{-1}\}}
    \bigl(S-(\zeta+\zeta^{-1})\bigr),
  \label{eq:classical-Lucas-trace-cyclotomic}
\end{equation}
where the product runs over inversion orbits of primitive $N$th roots of
unity.  Then $\Phi_N^+(S)\in\mathbb Z[S]$ is monic of degree
$\varphi(N)/2$.  If $k$ is a field with $\charac(k)\nmid2N$ and
$s\in\overline{k}$, the toric shift $T_s$ has exact order $N$ on the
normalized state curve if and only if
\begin{equation}
  \Phi_N^+(s)=0.
  \label{eq:classical-Lucas-exact-period-polynomial}
\end{equation}
For $N=3,4,5,6$ one has
\begin{equation}
 \Phi_3^+(S)=S+1,
 \quad
 \Phi_4^+(S)=S,
 \quad
 \Phi_5^+(S)=S^2+S-1,
 \quad
 \Phi_6^+(S)=S-1.
 \label{eq:classical-Lucas-low-trace-cyclotomic}
\end{equation}
\end{proposition}

\begin{proof}
The set of numbers $\zeta+\zeta^{-1}$ in
\eqref{eq:classical-Lucas-trace-cyclotomic} is stable under the absolute
Galois group of $\mathbb Q$; each number is an algebraic integer.  Therefore
the monic product has rational algebraic-integer coefficients and hence lies
in $\mathbb Z[S]$.  Inversion has no fixed primitive $N$th root for $N\ge3$,
so the number of factors is $\varphi(N)/2$.

Let $q$ be a root of $Z^2-sZ+1$.  The parametrization
\eqref{eq:classical-Lucas-toric-parametrization} identifies $T_s$ with
$z\mapsto qz$.  Consequently $T_s$ has exact order $N$ exactly when $q$ is
a primitive $N$th root.  This is equivalent to
$s=q+q^{-1}$ being a root of $\Phi_N^+$, proving
\eqref{eq:classical-Lucas-exact-period-polynomial}.  The four displayed
low-degree polynomials follow by pairing the primitive roots under inversion.

\end{proof}

The proposition gives the exact period theory on the toric boundary.  Its
comparison with the smooth-fibre state-division factors is made only after
those factors have been constructed and proved below.

\phantomsection
\label{subsec:Lucas-complete-bridge-diagram}
The exact specialization just proved fits into a single commutative bridge.
It identifies four compatible realizations of the same index calculus:
integer indices, points of a one-dimensional algebraic group, oriented
adjacent degree-two states, and full-point companion coordinates.  The next
proposition records that bridge while distinguishing genuine morphisms on a
fixed fibre from specialization in the total family.

Put
\begin{equation}
  \mathscr H_s:\qquad
  w^2=\left(\frac{s^2}{4}-1\right)x^2+1
  \label{eq:Lucas-classical-companion-conic-H}
\end{equation}
and define
\begin{align}
  \mathcal S_q(z)
  &=\left(
      \frac{z-z^{-1}}{q-q^{-1}},
      \frac{qz-q^{-1}z^{-1}}{q-q^{-1}}
    \right),
  \label{eq:Lucas-toric-state-map-Sq}\\
  \mathcal C_s(x,y)
  &=\left(x,y-\frac{s}{2}x\right).
  \label{eq:Lucas-toric-companion-map-Cs}
\end{align}
Thus \(\mathcal S_q(q^n)=S_n^{\mathrm{tor}}\) and
\(\mathcal C_s(S_n^{\mathrm{tor}})=(U_n,V_n/2)\).

Write $E_{\lambda,\eta}$ for the pointed elliptic curve underlying
$\mathcal B_{\lambda,\eta}$ before choosing either the adjacent-state or
the Jacobi coordinate presentation.

\begin{proposition}[Complete toric--elliptic bridge]
\label{prop:Lucas-complete-bridge-diagram}
Let the smooth normalized state fibre
\(\mathcal B_{\lambda,\eta}\), its Jacobi companion
\(\mathcal J_{\lambda,\eta}\), and the marked point \(D\) be as in
Subsections~\ref{subsec:Lucas-odd-normalization} and
\ref{subsec:Lucas-Jacobi-companion}.  In the explicit family obtained by
letting \((\lambda,\eta)\) specialize to \((0,s)\), the following diagram
commutes after replacing the singular special fibre by its normalization:
\begin{figure}[htbp]
\centering
\resizebox{\textwidth}{!}{%
\begin{tikzcd}[ampersand replacement=\&,column sep=large,row sep=large]
 \mathbb Z
   \arrow[r,"{n\mapsto[n]D}"]
   \arrow[d,equal]
 \& E_{\lambda,\eta}
   \arrow[r,"\mathcal S_D"]
   \arrow[d,dashed,"\operatorname{sp}"']
 \& \mathcal B_{\lambda,\eta}
   \arrow[r,"\mathcal C_{\lambda,\eta}"]
   \arrow[d,dashed,"\operatorname{sp}"']
 \& \mathcal J_{\lambda,\eta}
   \arrow[d,dashed,"\operatorname{sp}"]
 \\
 \mathbb Z
   \arrow[r,"n\mapsto q^n"']
 \& \mathbb G_m
   \arrow[r,"\mathcal S_q"']
 \& \mathscr C_s
   \arrow[r,"\mathcal C_s"']
 \& \mathscr H_s
\end{tikzcd}}
\caption{The complete toric--elliptic bridge.  Dashed arrows denote
specialization in the family followed, on the nodal fibre, by normalization
of its smooth locus; they are not morphisms from one fixed geometric fibre to
another.}
\label{fig:Lucas-complete-bridge}
\end{figure}
Here
\begin{equation}
  \mathcal C_{\lambda,\eta}(x,y)
  =\left(
      x,
      (1-\lambda x^2)y-\frac{\eta}{2}x
    \right)
  \label{eq:Lucas-smooth-companion-map}
\end{equation}
is the state-to-companion map.  On the indexed orbit, the four entries in
each column are respectively
\begin{equation}
 n,
 \qquad [n]D,
 \qquad (v_n,v_{n+1}),
 \qquad (v_n,w_n)
 \label{eq:Lucas-complete-bridge-smooth-column}
\end{equation}
in the smooth fibre and
\begin{equation}
 n,
 \qquad q^n,
 \qquad (U_n,U_{n+1}),
 \qquad (U_n,V_n/2)
 \label{eq:Lucas-complete-bridge-toric-column}
\end{equation}
in the normalized nodal fibre.
\end{proposition}

\begin{proof}
The top row is defined in the following order.  The first map sends the
integer \(n\) to the elliptic multiple \([n]D\).  The second is the
oriented adjacent-state isomorphism
\[
  \mathscr S(P)=\bigl(v(P),v(P+D)\bigr).
\]
Its even quotient is the corresponding Kummer state, but the displayed map
retains the orientation needed to recover the full Jacobi point.  The third
is the companion map
\eqref{eq:Lucas-smooth-companion-map}.  By
Theorem~\ref{thm:Lucas-Pn-multiple}, these maps send \(n\) to the four
objects in \eqref{eq:Lucas-complete-bridge-smooth-column}.

For the bottom row, Proposition~\ref{prop:classical-Lucas-node-normalization}
identifies the smooth locus of the normalized nodal fibre with
\(\mathbb G_m\), sends the marked point to multiplication by \(q\), and
gives the state map \(\mathcal S_q\) in
\eqref{eq:Lucas-toric-state-map-Sq}.  Therefore
\[
  \mathcal S_q(q^n)
  =\left(
      \frac{q^n-q^{-n}}{q-q^{-1}},
      \frac{q^{n+1}-q^{-n-1}}{q-q^{-1}}
    \right)
  =(U_n,U_{n+1}).
\]
Equation~\eqref{eq:classical-Lucas-companion-recovery} then gives
\[
  \mathcal C_s(U_n,U_{n+1})
  =\left(U_n,U_{n+1}-\frac{s}{2}U_n\right)
  =\left(U_n,\frac{V_n}{2}\right).
\]
This proves commutativity of the bottom row.

It remains to justify the vertical arrows.  Theorem
\ref{thm:elliptic-Lucas-exact-toric-specialization} proves, formula by
formula, that the normalized state, companion, and binary-branch data
specialize to the corresponding Lucas data at \((\lambda,\eta)=(0,s)\).
The group law on the smooth locus of the nodal fibre is multiplication on
\(\mathbb G_m\), so the section \([n]D\) specializes to \(q^n\).  These
facts prove every square of Figure~\ref{fig:Lucas-complete-bridge}.
The qualification in the caption is necessary because specialization is a
map in the total family, not an isomorphism between a chosen smooth generic
fibre and the singular special fibre.
\end{proof}

\subsection[Affine index chains and the EDS bridge]
{Affine index chains and the algebraic EDS bridge}
\label{subsec:Lucas-chain-diagram}

The complete bridge becomes especially transparent when the integer affine
semigroup and the state semigroup are placed in a single chain.  For
\(m\in\mathbb Z\) and \(r\in\mathbb Z\), define
\begin{equation}
  a_{m,r}(n)=mn+r,
  \qquad
  \mu_{m,r}(P)=[m]P+[r]D.
  \label{eq:Lucas-affine-index-and-point-maps}
\end{equation}
Recall that
\begin{equation}
  \mathfrak F_{m,r}=T_D^r\Delta_{D,m},
  \qquad
  \mathfrak F_{m,r}(S_n)=S_{mn+r}.
  \label{eq:Lucas-state-affine-map-recall}
\end{equation}
On the companion model define
\begin{equation}
  \widetilde\mu_{m,r}
  =\mathcal C_{\lambda,\eta}\circ\mathfrak F_{m,r}
   \circ\mathcal C_{\lambda,\eta}^{-1}.
  \label{eq:Lucas-companion-affine-map}
\end{equation}
This notation denotes a morphism of the smooth complete companion curve;
the displayed affine expression for $\mathcal C_{\lambda,\eta}^{-1}$ is
used only on its regular chart.

\begin{theorem}[Index--point--state--companion chain]
\label{thm:Lucas-chain-diagram}
On every open set on which the displayed affine coordinates are defined,
the diagram
\begin{figure}[htbp]
\centering
\resizebox{\textwidth}{!}{%
\begin{tikzcd}[ampersand replacement=\&,column sep=large,row sep=large]
 n
   \arrow[r,mapsto,"\iota_D"]
   \arrow[d,mapsto,"a_{m,r}"']
 \& {[n]D}
   \arrow[r,mapsto,"\mathcal S_D"]
   \arrow[d,mapsto,"\mu_{m,r}"']
 \& S_n=(v_n,v_{n+1})
   \arrow[r,mapsto,"\mathcal C_{\lambda,\eta}"]
   \arrow[d,mapsto,"\mathfrak F_{m,r}"']
 \& (v_n,w_n)
   \arrow[d,mapsto,"\widetilde\mu_{m,r}"]
 \\
 mn+r
   \arrow[r,mapsto,"\iota_D"']
 \& {[mn+r]D}
   \arrow[r,mapsto,"\mathcal S_D"']
 \& S_{mn+r}
   \arrow[r,mapsto,"\mathcal C_{\lambda,\eta}"']
 \& (v_{mn+r},w_{mn+r})
\end{tikzcd}}
\caption{The chain diagram for the affine index map \(n\mapsto mn+r\).
The rightmost vertical arrow is the same group map
\(P\mapsto[m]P+[r]D\), written in Jacobi companion coordinates.}
\label{fig:Lucas-chain-diagram}
\end{figure}
commutes.  In particular,
\begin{align}
 (m,r)=(1,1)&:\quad T_D:S_n\longmapsto S_{n+1},
 \label{eq:Lucas-chain-shift-case}\\
 (m,r)=(2,0)&:\quad \mathcal L_0=\Delta_D:S_n\longmapsto S_{2n},
 \label{eq:Lucas-chain-even-case}\\
 (m,r)=(2,1)&:\quad \mathcal L_1=T_D\Delta_D:S_n\longmapsto S_{2n+1}.
 \label{eq:Lucas-chain-odd-case}
\end{align}
\end{theorem}

\begin{proof}
The left square commutes because
\[
  \mu_{m,r}([n]D)=[m]([n]D)+[r]D=[mn+r]D
  =\iota_D(a_{m,r}(n)).
\]
For the middle square, the affine-index semigroup theorem
\eqref{eq:QRT-affine-index-action} gives
\[
  \mathfrak F_{m,r}(\mathcal S_D(P))
  =\mathcal S_D([m]P+[r]D).
\]
Taking \(P=[n]D\) proves commutativity.  The companion map
\(\mathcal C_{\lambda,\eta}\) is the birational identification between
the state curve and the Jacobi companion established in
Subsection~\ref{subsec:Lucas-Jacobi-companion}.  Consequently applying
\(\mathcal C_{\lambda,\eta}\) before or after the group map
\(P\mapsto[m]P+[r]D\) gives the Jacobi coordinates of the same point.
This proves the right square.  The three displayed special cases follow by
substituting \((m,r)=(1,1),(2,0),(2,1)\).
\end{proof}

\phantomsection
\label{subsec:Lucas-EDS-bridge}
The affine-index chain becomes computationally explicit after a
Jacobi--Weierstrass dictionary and a division-polynomial normalization are
chosen.  The Lucas state sequence and an elliptic divisibility sequence are
not the same sequence: they are different systems of functions evaluated on
the same indexed points \([n]D\).  For the Jacobi quartic
\(\mathcal J_{\lambda,\eta}\), their relationship can nevertheless be
made completely explicit.  We first construct the bridge algebraically over
the ground field, keeping Kummer-only recovery, oriented full-point recovery,
and Ward-window arithmetic distinct; the analytic sigma normalization is
then used to explain the same formulas without being needed for their
finite-field validity.

\paragraph*{Jacobi--Weierstrass coordinates and the marked point.}
\phantomsection\label{subsubsec:Lucas-specific-Weierstrass-dictionary}

Put
\begin{equation}
  d_E=\eta+2.
  \label{eq:Lucas-specific-dE}
\end{equation}
The identity
\begin{equation}
  4h=d_E^2-4d_E-4\lambda
  \label{eq:Lucas-specific-h-dE}
\end{equation}
follows by substituting \(\eta=d_E-2\) into
\eqref{eq:Lucas-h-definition}.

\begin{proposition}[Explicit Jacobi--Weierstrass isomorphism]
\label{prop:Lucas-specific-Weierstrass-dictionary}
Assume \(\charac(k)\ne2\) and
\(\lambda\Delta_{\lambda,\eta}\ne0\).  On the affine open set
\(v\ne0\), define
\begin{equation}
  X=\frac{2(w+1)}{v^2},
  \qquad
  Y=\frac{2(X+h)}{v}.
  \label{eq:Lucas-specific-J-to-E}
\end{equation}
Then \((X,Y)\) lies on
\begin{equation}
  E_{\lambda,\eta}:
  \qquad
  Y^2=(X+h)(X^2-4\lambda)
  =X^3+hX^2-4\lambda X-4h\lambda.
  \label{eq:Lucas-specific-Weierstrass}
\end{equation}
Conversely, on the regular affine overlap,
\begin{equation}
  v=\frac{2(X+h)}{Y},
  \qquad
  w=\frac{X^2+2hX+4\lambda}{X^2-4\lambda}.
  \label{eq:Lucas-specific-E-to-J}
\end{equation}
The two rational maps extend uniquely to an isomorphism of the smooth
projective completions, sending the Jacobi identity \(O=(0,1)\) to the
Weierstrass identity.
\end{proposition}

\begin{proof}
From the first equation in \eqref{eq:Lucas-specific-J-to-E},
\[
  w=\frac{Xv^2}{2}-1.
\]
Substitution into
\(w^2=\lambda v^4+hv^2+1\) gives
\[
  \frac{X^2v^4}{4}-Xv^2+1
  =\lambda v^4+hv^2+1.
\]
After cancelling \(1\) and dividing by \(v^2\), which is permitted on
the stated chart, one obtains
\[
  \left(\frac{X^2}{4}-\lambda\right)v^2=X+h.
\]
Hence
\begin{equation}
  v^2=\frac{4(X+h)}{X^2-4\lambda}.
  \label{eq:Lucas-specific-v-square}
\end{equation}
Using the second equation in \eqref{eq:Lucas-specific-J-to-E},
\[
  Y^2
  =\frac{4(X+h)^2}{v^2}
  =(X+h)(X^2-4\lambda),
\]
which proves \eqref{eq:Lucas-specific-Weierstrass}.

Conversely, solving \(Y=2(X+h)/v\) gives the first formula in
\eqref{eq:Lucas-specific-E-to-J}.  Furthermore,
\[
\begin{aligned}
  w
  &=\frac{X}{2}\frac{4(X+h)}{X^2-4\lambda}-1\\
  &=\frac{X^2+2hX+4\lambda}{X^2-4\lambda},
\end{aligned}
\]
which is the second formula in
\eqref{eq:Lucas-specific-E-to-J}.  To check that the inverse recovers
\(X\), note that
\[
  w+1=\frac{2X(X+h)}{X^2-4\lambda}
\]
and combine this equality with
\eqref{eq:Lucas-specific-v-square}; then
\[
  \frac{2(w+1)}{v^2}=X.
\]
Substituting the recovered \(X\) into \(2(X+h)/v\) gives the displayed
inverse-coordinate value of \(Y\), by the defining relation
\(Y=2(X+h)/v\) in \eqref{eq:Lucas-specific-J-to-E}.  Thus the maps are
inverse on a dense open set.  Both
models are smooth projective genus-one curves under the stated smoothness
hypothesis.  A birational map between smooth projective curves extends
uniquely to an isomorphism.  Finally, as \(v\to0\) and \(w\to1\), the
function \(X=2(w+1)/v^2\) has a pole, so \(O=(0,1)\) maps to the point at
infinity of \(E_{\lambda,\eta}\).
\end{proof}

\begin{corollary}[The marked point in Weierstrass coordinates]
\label{cor:Lucas-specific-marked-point}
The point \(D=(1,\eta/2)\) corresponds to
\begin{equation}
  D_W=(X_D,Y_D)
  =\left(
      d_E,
      \frac{d_E^2-4\lambda}{2}
    \right).
  \label{eq:Lucas-specific-DW}
\end{equation}
Moreover,
\begin{equation}
  4(d_E+h)=d_E^2-4\lambda.
  \label{eq:Lucas-specific-W2-preidentity}
\end{equation}
\end{corollary}

\begin{proof}
At \(D\), the first formula in
\eqref{eq:Lucas-specific-J-to-E} gives
\[
  X_D=2\left(\frac{\eta}{2}+1\right)=\eta+2=d_E.
\]
Since \(v(D)=1\), the second formula gives
\[
  Y_D=2(d_E+h).
\]
Using \eqref{eq:Lucas-specific-h-dE},
\[
  4(d_E+h)=d_E^2-4\lambda,
\]
so \(Y_D=(d_E^2-4\lambda)/2\), as claimed.
\end{proof}

\paragraph*{Division-polynomial normalization and the two coordinate bridges.}
\phantomsection\label{subsubsec:Lucas-explicit-algebraic-EDS-bridge}

Let \(\psi_n\) denote the standard division polynomials for
\eqref{eq:Lucas-specific-Weierstrass}, with
\(\psi_0=0\), \(\psi_1=1\), and \(\psi_2=2Y\).  Define the algebraic
elliptic divisibility sequence attached to \((E_{\lambda,\eta},D_W)\) by
\begin{equation}
  W_n=\psi_n(D_W).
  \label{eq:Lucas-algebraic-EDS-definition}
\end{equation}
Thus \(W_0=0\), \(W_1=1\), and
\begin{equation}
  W_2=2Y_D=d_E^2-4\lambda
  =(\eta+2)^2-4\lambda=4(d_E+h).
  \label{eq:Lucas-explicit-W2}
\end{equation}

Unless a statement explicitly says otherwise, every assertion below that
identifies the zero set of \(W_n\), treats the \(n\)-division divisor as
reduced, or recovers a projective point from a finite EDS window is made
under the separability hypothesis
\[
  \operatorname{char}(k)=0
  \quad\text{or}\quad
  \operatorname{char}(k)\nmid n.
\]
The rational division-polynomial identities themselves remain valid as
identities of rational functions in arbitrary characteristic whenever the
underlying Weierstrass model is nonsingular; the hypothesis is needed only
when multiplicities or zero-set equivalences are invoked.

The algebraic sequence \((W_n)\) satisfies Ward's bilinear recurrence
\cite{Ward1948}:
\begin{equation}
\begin{aligned}
 W_{m+n}W_{m-n}
 ={}&W_{m+1}W_{m-1}W_n^2\\
   &-W_{n+1}W_{n-1}W_m^2.
\end{aligned}
\label{eq:Lucas-Ward-recurrence}
\end{equation}
It is the standard division-polynomial recurrence and therefore holds over
every field on which the displayed nonsingular Weierstrass model and the
marked point are defined.

\begin{corollary}[Ward fast-doubling identities]
\label{cor:Lucas-Ward-fast-doubling}
For all admissible integers \(n\),
\begin{align}
 W_{2n+1}
 &=W_{n+2}W_n^3-W_{n-1}W_{n+1}^3,
 \label{eq:Lucas-Ward-odd-doubling}\\
 W_{2n}W_2
 &=W_n\left(
      W_{n+2}W_{n-1}^2-W_{n-2}W_{n+1}^2
    \right).
 \label{eq:Lucas-Ward-even-doubling}
\end{align}
\end{corollary}

\begin{proof}
In \eqref{eq:Lucas-Ward-recurrence}, put \(m=n+1\) while keeping the
second index equal to \(n\).  Since \(W_1=1\), the result is
\eqref{eq:Lucas-Ward-odd-doubling}.  Next put \(m=n+1\) and replace the
second index by \(n-1\).  The left-hand side is \(W_{2n}W_2\), and the
right-hand side factors as
\[
 W_n\left(
      W_{n+2}W_{n-1}^2-W_{n-2}W_{n+1}^2
    \right),
\]
proving \eqref{eq:Lucas-Ward-even-doubling}.
\end{proof}

On the open set \(W_n\ne0\), define
\begin{equation}
  R_n=\frac{W_{n-1}W_{n+1}}{W_n^2},
  \qquad
  \Upsilon_n=\frac{W_{2n}}{W_n^4}.
  \label{eq:Lucas-algebraic-EDS-ratios}
\end{equation}
The letter \(\Upsilon_n\) is used here rather than \(S_n\), because
\(S_n=(v_n,v_{n+1})\) already denotes the adjacent Lucas state.

\begin{theorem}[The first and second EDS bridges]
\label{thm:Lucas-explicit-two-EDS-bridges}
For every \(n\) for which the displayed affine ratios are defined,
\begin{align}
  X([n]D_W)&=d_E-R_n,
  \label{eq:Lucas-explicit-X-bridge}\\
  2Y([n]D_W)&=\Upsilon_n.
  \label{eq:Lucas-explicit-Y-bridge}
\end{align}
Equivalently, the oriented ratio pair \((R_n,\Upsilon_n)\) recovers the
full Weierstrass point \([n]D_W\).  If \(W_n=0\), then
\([n]D_W=O\), and the affine ratios are replaced by the corresponding
projective boundary point.
\end{theorem}

\begin{proof}
The standard division-polynomial identity
\begin{equation}
  \phi_n=X\psi_n^2-\psi_{n-1}\psi_{n+1}
  \label{eq:Lucas-division-phi-identity}
\end{equation}
and the multiplication formula
\[
  X([n]P)=\frac{\phi_n(P)}{\psi_n(P)^2}
\]
give, after evaluation at \(P=D_W\),
\[
\begin{aligned}
  X([n]D_W)
  &=\frac{d_EW_n^2-W_{n-1}W_{n+1}}{W_n^2}\\
  &=d_E-R_n.
\end{aligned}
\]
This proves \eqref{eq:Lucas-explicit-X-bridge}.

For the odd coordinate, the standard formula for a Weierstrass equation
with \(a_1=a_3=0\) is
\begin{equation}
  [n]P
  =\left(
     \frac{\phi_n(P)}{\psi_n(P)^2},
     \frac{\omega_n(P)}{\psi_n(P)^3}
    \right),
  \qquad
  \omega_n=
  \frac{\psi_{n+2}\psi_{n-1}^2-
        \psi_{n-2}\psi_{n+1}^2}{4Y}.
  \label{eq:Lucas-division-omega-identity}
\end{equation}
Therefore
\begin{equation}
  Y([n]D_W)
  =\frac{W_{n+2}W_{n-1}^2-W_{n-2}W_{n+1}^2}
         {4Y_DW_n^3}.
  \label{eq:Lucas-Y-before-even-identity}
\end{equation}
The even-index division-polynomial identity is
\begin{equation}
  \psi_{2n}
  =\frac{\psi_n}{\psi_2}
   \left(
    \psi_{n+2}\psi_{n-1}^2-
    \psi_{n-2}\psi_{n+1}^2
   \right).
  \label{eq:Lucas-even-division-identity}
\end{equation}
Evaluating at \(D_W\), using \(W_2=2Y_D\), and solving for the
parenthesized expression gives
\[
  W_{n+2}W_{n-1}^2-W_{n-2}W_{n+1}^2
  =\frac{W_{2n}W_2}{W_n}.
\]
Substitution into \eqref{eq:Lucas-Y-before-even-identity} yields
\[
  Y([n]D_W)
  =\frac{W_{2n}(2Y_D)}{4Y_DW_n^4}
  =\frac{W_{2n}}{2W_n^4}.
\]
This is \eqref{eq:Lucas-explicit-Y-bridge}.  The final statement follows
from the defining property of division polynomials:
\(W_n=0\) if and only if \([n]D_W=O\), provided the characteristic does
not divide \(n\); in arbitrary characteristic the displayed rational
coordinate formulas still describe the regular affine open on which
\(W_n\ne0\).
\end{proof}

\begin{corollary}[Direct recovery of the elliptic Lucas coordinates]
\label{cor:Lucas-direct-EDS-coordinate}
On the regular oriented chart,
\begin{align}
  v_n
  &=\frac{W_2-4R_n}{\Upsilon_n},
  \label{eq:Lucas-direct-v-from-ratios}\\
  &=\frac{
     W_n^2\bigl(W_2W_n^2-4W_{n-1}W_{n+1}\bigr)
    }{W_{2n}},
  \label{eq:Lucas-direct-v-from-EDS}\\
  w_n
  &=\frac{(d_E-R_n)^2+2h(d_E-R_n)+4\lambda}
          {(d_E-R_n)^2-4\lambda}.
  \label{eq:Lucas-direct-w-from-ratio}
\end{align}
The ratio variables satisfy the genus-one equation
\begin{equation}
  \Upsilon_n^2
  =(W_2-4R_n)
   \bigl((d_E-R_n)^2-4\lambda\bigr).
  \label{eq:Lucas-ratio-state-curve}
\end{equation}
Furthermore,
\begin{equation}
  v_n^2
  =\frac{W_2-4R_n}{(d_E-R_n)^2-4\lambda}
  =\frac{4(d_E+h-R_n)}{(d_E-R_n)^2-4\lambda}.
  \label{eq:Lucas-v-square-from-R}
\end{equation}
Thus \(R_n\) contains only Kummer information, whereas
\(\Upsilon_n\) supplies the orientation needed to distinguish
\([n]D_W\) from \(-[n]D_W\).
\end{corollary}

\begin{proof}
Apply the inverse map
\eqref{eq:Lucas-specific-E-to-J} to
\((X,Y)=(d_E-R_n,\Upsilon_n/2)\).  Its first formula gives
\[
  v_n
  =\frac{2(d_E+h-R_n)}{\Upsilon_n/2}
  =\frac{4(d_E+h-R_n)}{\Upsilon_n}.
\]
By \eqref{eq:Lucas-explicit-W2},
\(4(d_E+h)=W_2\), proving
\eqref{eq:Lucas-direct-v-from-ratios}.  Substituting
\eqref{eq:Lucas-algebraic-EDS-ratios} and clearing denominators gives
\eqref{eq:Lucas-direct-v-from-EDS}.  The second inverse formula in
\eqref{eq:Lucas-specific-E-to-J} gives
\eqref{eq:Lucas-direct-w-from-ratio}.

The Weierstrass equation at
\(X=d_E-R_n\), \(2Y=\Upsilon_n\), is
\[
  \frac{\Upsilon_n^2}{4}
  =(d_E+h-R_n)
    \bigl((d_E-R_n)^2-4\lambda\bigr).
\]
Multiplication by \(4\) and use of
\(4(d_E+h-R_n)=W_2-4R_n\) prove
\eqref{eq:Lucas-ratio-state-curve}.  Dividing the square of
\eqref{eq:Lucas-direct-v-from-ratios} by
\eqref{eq:Lucas-ratio-state-curve} yields
\eqref{eq:Lucas-v-square-from-R}.

Finally, \(R_n=d_E-X([n]D_W)\) is unchanged by elliptic negation,
whereas \(\Upsilon_n=2Y([n]D_W)\) changes sign.  The coordinate \(v\)
is odd under negation by \eqref{eq:Lucas-Jacobi-negation}.  Consequently
\(R_n\) alone determines only the Kummer class, while the pair
\((R_n,\Upsilon_n)\) determines the oriented point.
\end{proof}

\paragraph*{Low-index EDS values and ratio-state dynamics.}
\phantomsection\label{subsubsec:Lucas-W2-W3-ratio-QRT}

The first nontrivial EDS values have the following explicit expressions in the parameters.  The formula for \(W_2\) has already been obtained in
\eqref{eq:Lucas-explicit-W2}.  For \(W_3\), the Weierstrass coefficients
are
\[
 a_1=a_3=0,
 \qquad a_2=h,
 \qquad a_4=-4\lambda,
 \qquad a_6=-4h\lambda.
\]
Hence
\[
 b_2=4h,
 \quad b_4=-8\lambda,
 \quad b_6=-16h\lambda,
 \quad b_8=-16h^2\lambda-16\lambda^2,
\]
and therefore
\begin{equation}
  \psi_3(X)
  =3X^4+4hX^3-24\lambda X^2-48h\lambda X
   -16h^2\lambda-16\lambda^2.
  \label{eq:Lucas-explicit-psi3}
\end{equation}

\begin{proposition}[Explicit low-index EDS values]
\label{prop:Lucas-explicit-W2-W3}
For the marked point \(D_W\),
\begin{align}
  W_2&=d_E^2-4\lambda,
  \label{eq:Lucas-explicit-W2-restated}\\
  W_3&=(d_E-\lambda-1)(d_E^2-4\lambda)^2
      =(\eta-\lambda+1)W_2^2.
  \label{eq:Lucas-explicit-W3}
\end{align}
Moreover,
\begin{align}
  W_4&=-\eta(\lambda-1)W_2^4,
  \label{eq:Lucas-explicit-W4}\\
  W_5&=-W_2^6\left[
      \eta(\lambda-1)W_2+(\eta-\lambda+1)^3
    \right].
  \label{eq:Lucas-explicit-W5}
\end{align}
Assume, in addition, that the relevant index is prime to
\(\operatorname{char}(k)\).  Then the vanishing of the primitive factor in
\eqref{eq:Lucas-explicit-W3}, \eqref{eq:Lucas-explicit-W4}, or
\eqref{eq:Lucas-explicit-W5} is equivalent to \(D_W\) having exact order
\(3\), \(4\), or \(5\), respectively.
\end{proposition}

\begin{proof}
The first equality is \eqref{eq:Lucas-explicit-W2}.  Substitute
\(X=d_E\) and
\(4h=d_E^2-4d_E-4\lambda\) into
\eqref{eq:Lucas-explicit-psi3}.  Term-by-term expansion gives
\[
\begin{aligned}
 W_3={}&d_E^5-(\lambda+1)d_E^4-8\lambda d_E^3
       +8\lambda(\lambda+1)d_E^2\\
     &+16\lambda^2d_E-16\lambda^2(\lambda+1).
\end{aligned}
\]
On the other hand,
\[
\begin{aligned}
 &(d_E-\lambda-1)(d_E^2-4\lambda)^2\\
 &\quad=(d_E-\lambda-1)
       (d_E^4-8\lambda d_E^2+16\lambda^2),
\end{aligned}
\]
and expanding the right-hand side gives the same polynomial.  This proves
\eqref{eq:Lucas-explicit-W3}.

For \(W_4\), the tangent slope at \(D_W=(d_E,W_2/2)\) is
\[
 m_D=\frac{3d_E^2+2hd_E-4\lambda}{W_2}.
\]
Using \eqref{eq:Lucas-specific-h-dE}, the numerator factors as
\[
 3d_E^2+2hd_E-4\lambda
 =\frac{\eta+4}{2}W_2,
\]
so \(m_D=(\eta+4)/2\).  The duplication formulas on
\eqref{eq:Lucas-specific-Weierstrass} give
\[
 X(2D_W)=m_D^2-h-2d_E=\lambda+1
\]
and
\[
\begin{aligned}
 Y(2D_W)
 &=-\frac{W_2}{2}
   +m_D\bigl(d_E-(\lambda+1)\bigr)\\
 &=-\frac{\eta(\lambda-1)}{2}.
\end{aligned}
\]
By \eqref{eq:Lucas-explicit-Y-bridge} with \(n=2\),
\(W_4/W_2^4=2Y(2D_W)\), proving
\eqref{eq:Lucas-explicit-W4}.

Finally, Ward's odd-index identity
\eqref{eq:Lucas-Ward-odd-doubling}, applied with \(n=2\), gives
\[
 W_5=W_4W_2^3-W_3^3.
\]
Substitution of \eqref{eq:Lucas-explicit-W3} and
\eqref{eq:Lucas-explicit-W4} yields
\eqref{eq:Lucas-explicit-W5}.

It remains to justify the exact-order interpretation without using any
later classification.  On the smooth family, \(W_2\ne0\), so \(D_W\) is
neither the identity nor a point of order \(2\).  If
\(\operatorname{char}(k)\nmid3\), the division-polynomial criterion gives
\(W_3=0\) if and only if \([3]D_W=O\); the preceding exclusion forces exact
order \(3\).  If \(\operatorname{char}(k)\nmid5\), then
\(W_5=0\) is equivalent to \([5]D_W=O\).  Since \(5\) is prime and
\(D_W\ne O\), this is equivalent to \(D_W\) having exact order
\(5\).  If
\(\operatorname{char}(k)\nmid4\), then \(W_4=0\) means that the order of
\(D_W\) divides \(4\); because orders \(1\) and \(2\) have already been
excluded, the order is exactly \(4\).  Finally, the nonzero powers of
\(W_2\) in \eqref{eq:Lucas-explicit-W3}--\eqref{eq:Lucas-explicit-W5}
do not change these zero loci, which proves the last assertion.
\end{proof}

\begin{proposition}[Multiplicative QRT recurrence for the EDS ratio]
\label{prop:Lucas-R-multiplicative-QRT}
On every regular orbit segment,
\begin{equation}
  R_{n+1}R_{n-1}
  =\frac{W_2^2R_n-W_3}{R_n^2}
  =W_2^2
    \frac{R_n-(\eta-\lambda+1)}{R_n^2}.
  \label{eq:Lucas-R-multiplicative-QRT}
\end{equation}
\end{proposition}

\begin{proof}
Take \(m=n\) and the second Ward index equal to \(2\) in
\eqref{eq:Lucas-Ward-recurrence}.  This gives
\[
 W_{n+2}W_{n-2}
 =W_{n+1}W_{n-1}W_2^2-W_3W_n^2.
\]
By definition,
\[
 R_{n+1}R_{n-1}
 =\frac{W_{n+2}W_n}{W_{n+1}^2}
  \frac{W_{n-2}W_n}{W_{n-1}^2}.
\]
Substituting the preceding Ward identity into the numerator and using
\[
 R_n=\frac{W_{n-1}W_{n+1}}{W_n^2}
\]
gives
\[
\begin{aligned}
 R_{n+1}R_{n-1}
 &=\frac{W_n^2
    \left(W_2^2W_{n-1}W_{n+1}-W_3W_n^2\right)}
    {W_{n+1}^2W_{n-1}^2}\\
 &=\frac{W_2^2R_n-W_3}{R_n^2}.
\end{aligned}
\]
The second equality in \eqref{eq:Lucas-R-multiplicative-QRT} follows
from \eqref{eq:Lucas-explicit-W3}.
\end{proof}

\phantomsection\label{subsubsec:Lucas-ratio-state-arithmetic}

The multiplicative ratio dynamics has an oriented lift.  The pair
\((R,\Upsilon)\) on \eqref{eq:Lucas-ratio-state-curve} is the affine
Weierstrass point
\[
  X=d_E-R,
  \qquad
  2Y=\Upsilon,
\]
one can write state doubling and the fixed translation by \(D_W\)
directly in these variables.

\begin{proposition}[Direct ratio-state doubling and translation]
\label{prop:Lucas-ratio-state-arithmetic}
For a regular point \((R,\Upsilon)\) of
\eqref{eq:Lucas-ratio-state-curve}, put
\begin{equation}
  A_R=3(d_E-R)^2+2h(d_E-R)-4\lambda.
  \label{eq:Lucas-ratio-tangent-numerator}
\end{equation}
Then doubling is
\begin{align}
  R^{(2)}
  &=3d_E+h-2R-\frac{A_R^2}{\Upsilon^2},
  \label{eq:Lucas-ratio-doubling-R}\\
  \Upsilon^{(2)}
  &=-\Upsilon+
    \frac{2A_R}{\Upsilon}\bigl(R^{(2)}-R\bigr).
  \label{eq:Lucas-ratio-doubling-Upsilon}
\end{align}
For fixed translation by \(D_W\), put
\begin{equation}
  M_R=\frac{W_2-\Upsilon}{2R}.
  \label{eq:Lucas-ratio-translation-slope}
\end{equation}
Then
\begin{align}
  R^+&=3d_E+h-R-M_R^2,
  \label{eq:Lucas-ratio-translation-R}\\
  \Upsilon^+&=-\Upsilon+2M_R(R^+-R).
  \label{eq:Lucas-ratio-translation-Upsilon}
\end{align}
These formulas are affine charts of morphisms of the complete genus-one
curve; the cases \(R=0\) or \(\Upsilon=0\) are handled by the corresponding
projective group-law chart.
\end{proposition}

\begin{proof}
Write
\(x=d_E-R\) and \(y=\Upsilon/2\).  On the Weierstrass model
\eqref{eq:Lucas-specific-Weierstrass}, the tangent slope at \((x,y)\) is
\[
  m=\frac{3x^2+2hx-4\lambda}{2y}
   =\frac{A_R}{\Upsilon}.
\]
The doubled first coordinate is
\(x^{(2)}=m^2-h-2x\).  Since
\(R^{(2)}=d_E-x^{(2)}\),
\[
\begin{aligned}
 R^{(2)}
 &=d_E-m^2+h+2(d_E-R)\\
 &=3d_E+h-2R-\frac{A_R^2}{\Upsilon^2},
\end{aligned}
\]
which proves \eqref{eq:Lucas-ratio-doubling-R}.  The doubled odd
coordinate is
\[
 y^{(2)}=-y+m(x-x^{(2)}).
\]
Because \(x-x^{(2)}=R^{(2)}-R\), multiplication by \(2\) gives
\eqref{eq:Lucas-ratio-doubling-Upsilon}.

For translation by
\(D_W=(d_E,W_2/2)\), the chord slope is
\[
 m_D=\frac{\Upsilon/2-W_2/2}{(d_E-R)-d_E}
     =\frac{W_2-\Upsilon}{2R}=M_R.
\]
The sum has first coordinate
\(x^+=M_R^2-h-x-d_E\).  Hence
\[
 R^+=d_E-x^+=3d_E+h-R-M_R^2,
\]
which proves \eqref{eq:Lucas-ratio-translation-R}.  The formula for the
odd coordinate is
\[
 2y^+=-\Upsilon+2M_R(x-x^+)
       =-\Upsilon+2M_R(R^+-R),
\]
proving \eqref{eq:Lucas-ratio-translation-Upsilon}.  The final
qualification follows because the Weierstrass group law is regular on the
smooth projective curve even where a chosen affine slope formula is not.
\end{proof}

\phantomsection
\label{subsubsec:Lucas-squared-QRT-semiconjugacy}
Forgetting the odd orientation coordinate produces a second, degree-two
quotient of the same indexed elliptic orbit.  Set
\begin{equation}
  z_n=v_n^2.
  \label{eq:Lucas-z-definition}
\end{equation}
Squaring the state relation
\eqref{eq:Lucas-normalized-biquadratic} gives
\begin{equation}
  \bigl(\lambda z_nz_{n+1}-z_n-z_{n+1}+1\bigr)^2
  =\eta^2z_nz_{n+1}.
  \label{eq:Lucas-z-biquadratic}
\end{equation}
For fixed \(z_n\), its two roots are \(z_{n-1}\) and \(z_{n+1}\).
Vieta's formulas therefore give
\begin{align}
  z_{n+1}z_{n-1}
  &=\left(\frac{1-z_n}{1-\lambda z_n}\right)^2,
  \label{eq:Lucas-z-multiplicative-QRT}\\
  z_{n+1}+z_{n-1}
  &=\frac{2(1-\lambda z_n)(1-z_n)+\eta^2z_n}
          {(1-\lambda z_n)^2}.
  \label{eq:Lucas-z-additive-QRT}
\end{align}

\begin{proposition}[The Kummer ratio quotient]
\label{prop:Lucas-squared-QRT-semiconjugacy}
Define
\begin{equation}
  \Phi_R(R)
  =\frac{W_2-4R}{(d_E-R)^2-4\lambda}.
  \label{eq:Lucas-R-to-z-map}
\end{equation}
Then \(z_n=\Phi_R(R_n)\).  The rational map \(\Phi_R:\PP^1\to\PP^1\)
has degree two for generic smooth parameters.  On every common regular
orbit, it sends the multiplicative ratio-QRT sequence
\eqref{eq:Lucas-R-multiplicative-QRT} to the squared QRT sequence
\eqref{eq:Lucas-z-multiplicative-QRT}.  After adjoining the orientation
coordinates \(\Upsilon_n\) and \(w_n\), the corresponding genus-one
covers are birational rather than merely semiconjugate.
\end{proposition}

\begin{proof}
The identity \(z_n=\Phi_R(R_n)\) is exactly
\eqref{eq:Lucas-v-square-from-R}.  The numerator of \(\Phi_R\) has
degree one and its denominator has degree two.  They have no common factor
for generic smooth parameters, so the induced map of projective lines has
degree two.

Both \(R_n\) and \(z_n\) are obtained from the same elliptic point
\([n]D\) by rational functions, and the equality
\(z_n=\Phi_R(R_n)\) holds for every integer \(n\) in the regular orbit.
Applying it at \(n-1,n,n+1\), and using the already proved recurrences
\eqref{eq:Lucas-R-multiplicative-QRT} and
\eqref{eq:Lucas-z-multiplicative-QRT}, proves the claimed semiconjugacy on
the invariant elliptic orbit.  Its domain is the orbit closure carrying the
two recurrences, and the equality intertwines their three consecutive
indexed values there.

Finally, Theorem~\ref{thm:Lucas-explicit-two-EDS-bridges} identifies
\((R,\Upsilon)\) with the full Weierstrass point, while
Proposition~\ref{prop:Lucas-specific-Weierstrass-dictionary} identifies
that point with \((v,w)\).  Their composition is a birational map of the
oriented genus-one covers.  The loss of degree two occurs only after the
odd orientation coordinates are forgotten.
\end{proof}

\phantomsection\label{subsubsec:Lucas-analytic-EDS-normalization}
The preceding formulas are algebraic.  Over \(\mathbb C\), they also arise
from the sigma function.  Write the same curve in analytic Weierstrass form
\begin{equation}
  E(\mathbb C):
  \qquad
  Y_a^2=4X_a^3-g_2X_a-g_3
  \simeq\mathbb C/\Lambda,
  \label{eq:Lucas-EDS-analytic-Weierstrass}
\end{equation}
with \(X_a=\wp(z)\) and \(Y_a=\wp'(z)\).  Let \(z_D\) represent the
marked point.  Define the analytic sigma sequence
\begin{equation}
  \widehat W_n(z_D)
  =\frac{\sigma(nz_D)}{\sigma(z_D)^{n^2}}.
  \label{eq:Lucas-EDS-sigma-definition}
\end{equation}
For the monic Weierstrass normalization used above, the standard algebraic
division-polynomial sequence differs from this analytic normalization by
the parity sign
\begin{equation}
  W_n=(-1)^{n-1}\widehat W_n.
  \label{eq:Lucas-algebraic-analytic-EDS-sign}
\end{equation}
The sign is immaterial in the even ratio
\(W_{n-1}W_{n+1}/W_n^2\), but it is essential for the odd coordinate
\(W_{2n}/W_n^4\).

\begin{theorem}[Sigma coordinate recovery]
\label{thm:Lucas-EDS-coordinate-recovery}
For every integer \(n\) for which the denominators are nonzero,
\begin{align}
  \wp(nz_D)
  &=\wp(z_D)-
    \frac{\widehat W_{n-1}\widehat W_{n+1}}
         {\widehat W_n^2},
  \label{eq:Lucas-EDS-x-recovery}\\
  \wp'(nz_D)
  &=-\frac{\widehat W_{2n}}{\widehat W_n^4}.
  \label{eq:Lucas-EDS-y-recovery}
\end{align}
Moreover,
\begin{equation}
  \widehat W_n=0
  \quad\Longleftrightarrow\quad
  nz_D\in\Lambda
  \quad\Longleftrightarrow\quad
  [n]D=O.
  \label{eq:Lucas-EDS-zero-torsion}
\end{equation}
After the affine change from \((\wp,\wp')\) to the monic model
\eqref{eq:Lucas-specific-Weierstrass}, equations
\eqref{eq:Lucas-EDS-x-recovery}--\eqref{eq:Lucas-EDS-y-recovery}, together
with \eqref{eq:Lucas-algebraic-analytic-EDS-sign}, become precisely
\eqref{eq:Lucas-explicit-X-bridge}--\eqref{eq:Lucas-explicit-Y-bridge}.
If an auxiliary symmetric state uses the even fractional-linear Kummer
coordinate \(\kappa=\mu\circ X\), then its coordinate value
\(\kappa_n=\kappa([n]D)\) is
\begin{equation}
  \kappa_n
  =\mu\left(
      X(D)-\frac{W_{n-1}W_{n+1}}{W_n^2}
    \right)
  \label{eq:Lucas-EDS-to-Kummer-coordinate}
\end{equation}
on the corresponding chart.  To retain the notation of the earlier
Kummer-only bridge without identifying it with the oriented Jacobi coordinate,
we record the same relation once more as
\begin{equation}
  \kappa_n
  =\mu\left(
      X(D)-\frac{W_{n-1}W_{n+1}}{W_n^2}
    \right).
  \label{eq:Lucas-EDS-to-vn}
\end{equation}
In the earlier notation the symbol attached to this formula denoted an even
Kummer coordinate.  It must not be confused with the present oriented
Jacobi coordinate \(v_n\), which also requires the second, odd EDS bridge
\eqref{eq:Lucas-explicit-Y-bridge}.
\end{theorem}

\begin{proof}
The Weierstrass sigma addition identity is
\begin{equation}
  \wp(v)-\wp(u)
  =\frac{\sigma(u+v)\sigma(u-v)}
         {\sigma(u)^2\sigma(v)^2}.
  \label{eq:Lucas-sigma-addition-identity}
\end{equation}
Substitute \(u=nz_D\) and \(v=z_D\).  The powers of
\(\sigma(z_D)\) cancel because
\[
  (n+1)^2+(n-1)^2-2n^2-2=0,
\]
which proves \eqref{eq:Lucas-EDS-x-recovery}.

For the odd coordinate, put \(v=u+\varepsilon\) in
\eqref{eq:Lucas-sigma-addition-identity} and let
\(\varepsilon\to0\).  Since
\[
  \sigma(-\varepsilon)=-\varepsilon+O(\varepsilon^3),
\]
one obtains
\[
  \wp'(u)=-\frac{\sigma(2u)}{\sigma(u)^4}.
\]
Taking \(u=nz_D\) and cancelling the powers of \(\sigma(z_D)\) gives
\eqref{eq:Lucas-EDS-y-recovery}.  The sigma function has simple zeros
exactly at lattice points, proving
\eqref{eq:Lucas-EDS-zero-torsion}.

It remains to justify the sign in
\eqref{eq:Lucas-algebraic-analytic-EDS-sign}.  The functions
\(\psi_n(D)\) and \(\widehat W_n(z_D)\) have the same divisor as
functions of the marked point.  Their quotient is therefore constant.
Near the identity, where \(z_D\to0\), one has
\(\widehat W_n\sim n z_D^{1-n^2}\).  In the monic coordinates
\(X=\wp-h/3\), \(Y=\wp'/2\), the leading term of the standard division
polynomial is the same for odd \(n\) and its negative for even \(n\).
Consequently the constant quotient is \((-1)^{n-1}\).  Substituting this
sign into the analytic formulas yields the algebraic first and second
bridges.  Applying the fractional-linear map \(\mu\) proves
\eqref{eq:Lucas-EDS-to-Kummer-coordinate}.  The final qualification follows
because \(X(P)=X(-P)\), whereas \(v(-P)=-v(P)\).
\end{proof}

\begin{proposition}[The algebraic EDS--Lucas closure diagram]
\label{prop:Lucas-algebraic-EDS-closure-diagram}
On the common regular open set, all arrows in
\begin{figure}[htbp]
\centering
\resizebox{\textwidth}{!}{%
\begin{tikzcd}[ampersand replacement=\&,column sep=large,row sep=large]
 \{W_j\}_{j\sim n}
  \arrow[r,mapsto,"{(R_n,\Upsilon_n)}"]
 \& (R_n,\Upsilon_n)
  \arrow[r,mapsto,"{(d_E-R_n,\Upsilon_n/2)}"]
  \arrow[d,mapsto,"{\Phi_R}"']
 \& {[n]D_W\in E_{\lambda,\eta}}
  \arrow[r,leftrightarrow,"{\text{Jacobi--Weierstrass}}"]
 \& (v_n,w_n)\in\mathcal J_{\lambda,\eta}
  \arrow[r,leftrightarrow,"{\text{companion}}"]
 \& (v_n,v_{n+1})\in\mathcal B_{\lambda,\eta}
 \\
 \& z_n=v_n^2
\end{tikzcd}}
\caption{The explicit algebraic closure from a Ward EDS window to the
oriented elliptic point, Jacobi companion, and adjacent QRT state.  The
downward arrow forgets orientation and has generic degree two.}
\label{fig:Lucas-algebraic-EDS-closure}
\end{figure}
are given by the formulas proved above.  The horizontal arrows are
birational after the odd orientation coordinate is retained.
\end{proposition}

\begin{proof}
The first arrow is definition
\eqref{eq:Lucas-algebraic-EDS-ratios}.  The second is
Theorem~\ref{thm:Lucas-explicit-two-EDS-bridges}.  The third is the
isomorphism of
Proposition~\ref{prop:Lucas-specific-Weierstrass-dictionary}.  The fourth
is the companion isomorphism
\eqref{eq:Lucas-state-to-Jacobi}--\eqref{eq:Lucas-Jacobi-to-state}.  The
downward arrow is \eqref{eq:Lucas-R-to-z-map}, whose generic degree is two
by Proposition~\ref{prop:Lucas-squared-QRT-semiconjugacy}.  Every square
therefore commutes by substitution into the displayed formulas, and the
horizontal compositions recover the same point \([n]D\).
\end{proof}

\begin{corollary}[Compressed bridge and model chains]
\label{cor:Lucas-compressed-bridge-model-chains}
The algebraic closure diagram admits the following two compressed forms.
First, the model chain from the original three-parameter curve to the EDS is
\begin{equation}
\boxed{
\begin{gathered}
 \Q_{\alpha,\beta,\gamma}
 \xrightarrow{\ (x,y)=(r\xi,r\zeta)\ }
 \mathcal B_{\lambda,\eta},\qquad
 r^2=-\frac{\gamma}{\alpha},\quad
 \lambda=\frac{\gamma}{\alpha^2},\quad
 \eta=-\frac{\beta}{\alpha},\\[1mm]
 \mathcal B_{\lambda,\eta}
 \xleftrightarrow{\text{companion}}
 \mathcal J_{\lambda,\eta}
 \xleftrightarrow{\text{Jacobi--Weierstrass}}
 E_{\lambda,\eta}
 \xrightarrow{\;D_W\mapsto\{\psi_n(D_W)\}\;}
 \{W_n\}.
\end{gathered}}
\label{eq:Lucas-compressed-model-chain}
\end{equation}
Second, the indexed arithmetic chain is
\begin{equation}
\boxed{
\begin{aligned}
 \{W_j\}_{j\sim n}
 &\longrightarrow (R_n,\Upsilon_n)
 \longrightarrow [n]D_W,\\
 [n]D_W
 &\longleftrightarrow (v_n,w_n)
 \longleftrightarrow (v_n,v_{n+1}),\\
 (v_n,v_{n+1})
 &\longrightarrow
 \bigl[\lambda v_n^2v_{n+1}^2-v_n^2-v_{n+1}^2
       +\eta v_nv_{n+1}+1=0\bigr],\\
 &\longrightarrow
 \left[v_{n+1}+v_{n-1}
       =\frac{\eta v_n}{1-\lambda v_n^2}\right]
 \longrightarrow (\mathcal L_0,\mathcal L_1).
\end{aligned}}
\label{eq:Lucas-compressed-arithmetic-chain}
\end{equation}
The even EDS branch is simultaneously summarized by
\begin{equation}
 \{W_j\}_{j\sim n}
 \longrightarrow R_n
 \longrightarrow
 R_{n+1}R_{n-1}
 =W_2^2\frac{R_n-(\eta-\lambda+1)}{R_n^2}.
\label{eq:Lucas-compressed-ratio-QRT-chain}
\end{equation}
Thus the Ward EDS, the multiplicative QRT on \(R_n\), the oriented
Jacobi lift, the symmetric McMillan recurrence, and Lucas-type fast index
arithmetic are different presentations of one indexed elliptic orbit.
\end{corollary}

\begin{proof}
For \eqref{eq:Lucas-compressed-model-chain}, substitution of
\((x,y)=(r\xi,r\zeta)\) into \eqref{eq:Lucas-original-QRT}, followed by
division by \(\gamma\), gives
\eqref{eq:Lucas-normalized-biquadratic} with the parameters in
\eqref{eq:Lucas-normalized-parameters}.  The next arrow is the companion
isomorphism
\eqref{eq:Lucas-state-to-Jacobi}--\eqref{eq:Lucas-Jacobi-to-state}; the
following arrow is Proposition~\ref{prop:Lucas-specific-Weierstrass-dictionary};
and the last arrow is the definition
\(W_n=\psi_n(D_W)\) in
\eqref{eq:Lucas-algebraic-EDS-definition}.  Hence every arrow in
\eqref{eq:Lucas-compressed-model-chain} has already been defined over the
stated ground field or over the explicitly indicated square-class
extension.

For \eqref{eq:Lucas-compressed-arithmetic-chain}, the first two arrows are
Theorem~\ref{thm:Lucas-explicit-two-EDS-bridges}; the Jacobi and companion
arrows are the two isomorphisms just cited.  The symmetric biquadratic
identity is \eqref{eq:Lucas-normalized-biquadratic}; applying its Vieta root
exchange at the middle coordinate gives
\eqref{eq:Lucas-normalized-recurrence}; and the formulas
\(\mathcal L_0(S_n)=S_{2n}\) and
\(\mathcal L_1(S_n)=S_{2n+1}\) are
\eqref{eq:Lucas-state-index-branches}.  Finally,
\eqref{eq:Lucas-compressed-ratio-QRT-chain} is
Proposition~\ref{prop:Lucas-R-multiplicative-QRT}.  This verifies every
link rather than treating the displayed chains as merely mnemonic.
\end{proof}

\subsection[EDS doubling and the commutative calculus]
{EDS repeated doubling and the total commutative calculus}
\label{subsec:Lucas-EDS-repeated-doubling-diagram}

The EDS identities and the Lucas state identities compute different data.
Nevertheless, after coordinate recovery they must agree because both compute
the same elliptic multiples.  We now make this assertion precise and check
that a finite EDS window is sufficient.  Repeated-doubling evaluation of EDS
windows was developed in the computational work of Shipsey and is discussed
in the subsequent EDS literature; see
\cite{Shipsey2000,EverestWard2001}.  Here we establish the exact commutative
comparison with the adjacent-state Lucas branches.

Assume \(2D\ne O\), so \(W_2\ne0\).  Define the radius-four EDS window
\begin{equation}
  \mathbf W_n
  =\bigl(W_{n-4},W_{n-3},\ldots,W_{n+3},W_{n+4}\bigr).
  \label{eq:Lucas-EDS-radius-four-window}
\end{equation}
For an integer \(j\), put
\begin{align}
  \mathscr O_j
  &=W_{j+2}W_j^3-W_{j-1}W_{j+1}^3,
  \label{eq:Lucas-EDS-odd-template}\\
  \mathscr E_j
  &=\frac{W_j}{W_2}
    \left(W_{j+2}W_{j-1}^2-W_{j-2}W_{j+1}^2\right).
  \label{eq:Lucas-EDS-even-template}
\end{align}
By Corollary~\ref{cor:Lucas-Ward-fast-doubling},
\begin{equation}
  \mathscr O_j=W_{2j+1},
  \qquad
  \mathscr E_j=W_{2j}.
  \label{eq:Lucas-EDS-template-values}
\end{equation}
Consequently the templates with
\(j\in\{n-2,n-1,n,n+1,n+2\}\) construct the two radius-four
windows \(\mathbf W_{2n}\) and \(\mathbf W_{2n+1}\).  Denote the
resulting rational window maps by
\begin{equation}
  \mathscr D_{W,0}(\mathbf W_n)=\mathbf W_{2n},
  \qquad
  \mathscr D_{W,1}(\mathbf W_n)=\mathbf W_{2n+1}.
  \label{eq:Lucas-EDS-window-branch-maps}
\end{equation}

For every integer \(j\), define
\begin{align}
  \Phi_j
  &=d_EW_j^2-W_{j-1}W_{j+1},
  \label{eq:Lucas-EDS-projective-Phi}\\
  \Xi_j
  &=W_{j+2}W_{j-1}^2-W_{j-2}W_{j+1}^2.
  \label{eq:Lucas-EDS-projective-Xi}
\end{align}
The division-polynomial multiplication formula gives the homogeneous
Weierstrass point
\begin{equation}
  \mathscr P_j
  =\bigl(
      2W_2\Phi_jW_j:
      \Xi_j:
      2W_2W_j^3
    \bigr)
  \in E_{\lambda,\eta}\subset\PP^2.
  \label{eq:Lucas-EDS-projective-point-recovery}
\end{equation}
Indeed, if \(W_j\ne0\), then
\[
  \frac{2W_2\Phi_jW_j}{2W_2W_j^3}
  =\frac{\Phi_j}{W_j^2}
  =X([j]D_W),
\]
and, by \eqref{eq:Lucas-division-omega-identity},
\[
  \frac{\Xi_j}{2W_2W_j^3}
  =\frac{\omega_j(D_W)}{W_j^3}
  =Y([j]D_W).
\]
When \(W_j=0\), the same homogeneous triple specializes to the point at
infinity, provided the multiplication-by-\(j\) map is separable at the
marked point; the general inseparable case is handled by the complete
morphism \([j]\) rather than by an affine division-polynomial chart.

The full point \(\mathscr P_j\) refines projective Kummer recovery by
retaining the odd coordinate.  Forgetting that coordinate gives
\begin{equation}
  \mathscr X_j
  =\bigl(d_EW_j^2-W_{j-1}W_{j+1}:W_j^2\bigr)\in\PP^1.
  \label{eq:Lucas-EDS-projective-Kummer-recovery}
\end{equation}
This pair is never \((0:0)\) when \(D_W\ne O\).  If \(W_j\ne0\), its
affine value is \(X([j]D_W)\).  If \(W_j=0\), then \([j]D_W=O\); if also
\(W_{j-1}=0\) or \(W_{j+1}=0\), subtracting the corresponding torsion
relations would give \(D_W=O\), a contradiction.  Thus the first entry is
nonzero at such a torsion index.  Applying a fractional-linear Kummer
normalization to \(\mathscr X_j\) recovers the even adjacent state, whereas
\(\mathscr P_j\) and the companion inverse retain the orientation of the full
Jacobi state.

Let
\[
  \overline\iota:
  \overline{\mathcal J}_{\lambda,\eta}
  \xrightarrow{\sim}
  E_{\lambda,\eta}
\]
be the projective isomorphism of
Proposition~\ref{prop:Lucas-specific-Weierstrass-dictionary}, and let
\(\overline{\mathcal C}_{\lambda,\eta}^{-1}\) be the complete inverse of
the companion map.  Define the projective state recovery map by
\begin{equation}
  \mathscr R(\mathbf W_n)
  =\overline{\mathcal C}_{\lambda,\eta}^{-1}
   \overline\iota^{-1}(\mathscr P_n)
  \in\overline{\mathcal B}_{\lambda,\eta}.
  \label{eq:Lucas-EDS-window-state-recovery-projective}
\end{equation}
The inverse companion morphism returns the complete adjacent state attached
to the recovered point \([n]D\), rather than only its Kummer coordinate.  On
the regular oriented chart this is exactly
\begin{equation}
\begin{aligned}
 \mathscr R(\mathbf W_n)
 =\Bigg(&
  \frac{W_n^2
    (W_2W_n^2-4W_{n-1}W_{n+1})}{W_{2n}},\\
 &\frac{W_{n+1}^2
    (W_2W_{n+1}^2-4W_nW_{n+2})}{W_{2n+2}}
 \Bigg)
 =(v_n,v_{n+1}).
\end{aligned}
\label{eq:Lucas-EDS-window-state-recovery}
\end{equation}
Thus the recovery uses both EDS bridges.  The even ratio recovers the
Kummer coordinate, while the odd division-polynomial numerator in
\(\mathscr P_j\) selects the correct orientation sheet.

\begin{theorem}[EDS/Lucas doubling square]
\label{thm:Lucas-EDS-Lucas-doubling-square}
For \(b\in\{0,1\}\), the diagram
\begin{figure}[htbp]
\centering
\begin{tikzcd}[column sep=huge,row sep=large]
 \mathbf W_n
   \arrow[r,"\mathscr D_{W,b}"]
   \arrow[d,"\mathscr R"']
 & \mathbf W_{2n+b}
   \arrow[d,"\mathscr R"]
 \\
 S_n
   \arrow[r,"\mathcal L_b"']
 & S_{2n+b}
\end{tikzcd}
\caption{The commutative square comparing Ward/EDS fast doubling with
Lucas-type adjacent-state fast doubling.}
\label{fig:Lucas-EDS-doubling-square}
\end{figure}
commutes for every index for which the Ward window formulas are defined,
when \(\mathscr R\) is interpreted projectively by
\eqref{eq:Lucas-EDS-window-state-recovery-projective}.  In particular,
for every \(r\ge0\),
\begin{equation}
  \mathscr R\bigl(\mathscr D_{W,0}^{\,r}(\mathbf W_n)\bigr)
  =\mathcal L_0^{\,r}(S_n)
  =S_{2^r n}.
  \label{eq:Lucas-EDS-repeated-doubling-commutation}
\end{equation}
\end{theorem}

\begin{proof}
First we verify that the radius-four input is sufficient.  If an output
index is odd, write it as \(2j+1\).  Formula
\eqref{eq:Lucas-EDS-odd-template} uses only the four terms from
\(W_{j-1}\) through \(W_{j+2}\).  If an output index is even, write it
as \(2j\).  Formula \eqref{eq:Lucas-EDS-even-template} uses only the
five terms from \(W_{j-2}\) through \(W_{j+2}\), together with the fixed
constant \(W_2\).  For the indices from \(2n-4\) through \(2n+5\),
which contain both radius-four windows centred at \(2n\) and \(2n+1\),
all required input indices lie between \(n-4\) and \(n+4\).  Therefore
\eqref{eq:Lucas-EDS-window-branch-maps} is a well-defined rational
construction from \(\mathbf W_n\).

By the projective recovery construction
\eqref{eq:Lucas-EDS-projective-point-recovery}--
\eqref{eq:Lucas-EDS-window-state-recovery-projective}, together with
Theorem~\ref{thm:Lucas-explicit-two-EDS-bridges},
\(\mathscr R(\mathbf W_n)=S_n\).  Applying the same construction to the
output window gives
\[
  \mathscr R(\mathscr D_{W,b}(\mathbf W_n))=S_{2n+b}.
\]
On the other hand, the state-branch theorem gives
\[
  \mathcal L_b(\mathscr R(\mathbf W_n))
  =\mathcal L_b(S_n)=S_{2n+b}.
\]
Thus the square commutes on the affine oriented chart.  At a boundary or
torsion state, the homogeneous point
\eqref{eq:Lucas-EDS-projective-point-recovery} and the complete projective
isomorphisms defining \(\mathscr R\) recover the full elliptic point, not
only its Kummer class.  Ward's identities compute the same multiples
\([2n+b]D\).  Both recovered paths are morphisms to the smooth projective
state curve and agree on the dense regular open, so they agree at every
state covered by the projective Ward construction.  In this
odd-characteristic section, \([2]\) is separable because
\(\charac(k)\ne2\), so no inseparable doubling case remains.
Finally,
\eqref{eq:Lucas-EDS-repeated-doubling-commutation} follows by induction on
\(r\): the case \(r=0\) is the definition of \(\mathscr R\), and the
inductive step is one application of the commutative square with \(b=0\).
\end{proof}

\begin{remark}[Data retained by the doubling bridge]
\label{rem:Lucas-EDS-doubling-data-cost}
Figure~\ref{fig:Lucas-EDS-doubling-square} identifies the elliptic point
recovered by both paths.  An EDS window retains division-function values for
several neighbouring multiples, whereas a Lucas branch retains the oriented
adjacent state.  The two arithmetic costs therefore correspond to different
retained data, and the diagram records their exact correctness
correspondence.
\end{remark}

\phantomsection
\label{subsec:Lucas-total-commutative-diagram}
The doubling square belongs to a larger calculus in which the integer index,
full elliptic point, adjacent state, companion, and EDS window are displayed
simultaneously.  For \(b\in\{0,1\}\), let
\begin{equation}
  \nu_b(P)=[2]P+[b]D.
  \label{eq:Lucas-total-diagram-point-map}
\end{equation}
Let \(\mathcal C_{\lambda,\eta}\) denote the companion identification in
\eqref{eq:Lucas-smooth-companion-map}.  When the state curve itself is
regarded as the chosen model of \(E\), the inverse arrow
\(\mathcal C_{\lambda,\eta}^{-1}\) is understood projectively through the
complete atlas, not only through the affine quotient formula.

\begin{theorem}[Total binary bridge]
\label{thm:Lucas-total-commutative-diagram}
The following indexed diagram commutes.  The two middle rows are related by
morphisms of the smooth complete elliptic/state curve, the bottom row records
the integer index update, and the upper row is the Ward recurrence on finite
EDS windows:
\begin{figure}[htbp]
\centering
\resizebox{0.88\textwidth}{!}{%
\begin{tikzcd}[ampersand replacement=\&,column sep=huge,row sep=large]
 \mathbf W_n
   \arrow[r,"\mathscr D_{W,b}"]
   \arrow[d,"\mathscr R"']
 \& \mathbf W_{2n+b}
   \arrow[d,"\mathscr R"]
 \\
 S_n
   \arrow[r,"\mathcal L_b=T_D^b\Delta_D"']
   \arrow[d,"\mathcal C_{\lambda,\eta}"']
 \& S_{2n+b}
   \arrow[d,"\mathcal C_{\lambda,\eta}"]
 \\
 {[n]D}
   \arrow[r,"\nu_b"']
 \& {[2n+b]D}
 \\
 n
   \arrow[u,mapsto,"\iota_D"]
   \arrow[r,mapsto,"n\mapsto2n+b"']
 \& 2n+b
   \arrow[u,mapsto,"\iota_D"']
\end{tikzcd}}
\caption{The total commutative diagram linking Ward/EDS windows,
adjacent QRT states, full elliptic points, and binary index arithmetic.}
\label{fig:Lucas-total-commutative-diagram}
\end{figure}
The diagram uses the fixed algebraic normalization
\(W_n=\psi_n(D_W)\) from
\eqref{eq:Lucas-algebraic-EDS-definition}.  The analytic sigma sequence is
related to it by the parity sign
\eqref{eq:Lucas-algebraic-analytic-EDS-sign}.  No independent EDS rescaling
is performed inside the diagram.
\end{theorem}

\begin{proof}
The top square is Theorem~\ref{thm:Lucas-EDS-Lucas-doubling-square}.
For the middle square, apply the state-to-point identification to the branch
identity
\[
  \mathcal L_b(\mathcal S_D(P))
  =\mathcal S_D([2]P+[b]D).
\]
Taking \(P=[n]D\) shows that both paths give the full point
\([2n+b]D\).  The bottom square commutes because
\[
  \nu_b([n]D)
  =[2]([n]D)+[b]D
  =[2n+b]D.
\]
These three checks prove commutativity of every square.  The fixed
algebraic normalization removes the usual EDS scaling ambiguity.  The
state-to-companion and Jacobi--Weierstrass arrows are complete morphisms, so
the lower part of the diagram includes all boundary states covered by the
separable multiplication maps.
\end{proof}

The diagram uses the fixed algebraic normalization \(W_n=\psi_n(D_W)\).
Under the rank-one EDS scale equivalence
\begin{equation}
  W_j\longmapsto c^{j^2-1}W_j,
  \qquad c\in k^\times,
  \label{eq:Lucas-EDS-scale-equivalence}
\end{equation}
one has
\begin{equation}
  R_j\longmapsto c^2R_j,
  \qquad
  \Upsilon_j\longmapsto c^3\Upsilon_j.
  \label{eq:Lucas-EDS-ratio-scale-equivariance}
\end{equation}
The exponents are
\[
 ((j-1)^2-1)+((j+1)^2-1)-2(j^2-1)=2
\]
and
\[
 ((2j)^2-1)-4(j^2-1)=3.
\]
They match the Weierstrass weights \(X\mapsto c^2X\) and
\(Y\mapsto c^3Y\).  If \(d_E=X(D_W)\) is changed simultaneously to
\(c^2d_E\), then
\[
  \Phi_j\longmapsto c^{2j^2}\Phi_j,
  \qquad
  \Xi_j\longmapsto c^{3j^2+3}\Xi_j.
\]
Consequently \(\mathscr P_j\) changes by the Weierstrass coordinate weights,
and \(\mathscr X_j\) changes by the weight of the Kummer \(X\)-coordinate.
Recovery is therefore equivariant under EDS scaling together with the
corresponding coordinate change; it is not invariant under rescaling the
sequence while holding a fixed Weierstrass equation unchanged.  The chosen
division-polynomial normalization removes this ambiguity from the diagrams.

The diagram separates three logically different facts.  The lower square is
pure group arithmetic.  The middle square is the oriented adjacent-state
realization of that arithmetic.  The upper square is the division-function
realization.  Conflating these levels would incorrectly identify an EDS term
with a Kummer coordinate.  The following comparison makes their distinction
exact before the common sigma-function source is introduced.

\subsection{Three sequence layers, sigma functions, and elliptic nets}
\label{subsec:Lucas-three-sequence-distinction}

The words ``Lucas sequence'', ``elliptic divisibility sequence'', and
``elliptic Lucas sequence'' refer to related but different constructions.
The following table records the precise distinction used in this text.

Before comparing the elliptic constructions with the classical toric Lucas
pair, it is useful to distinguish three sequences that live on the same
smooth elliptic orbit:
\begin{equation}
  W_n=\psi_n(D_W),
  \qquad
  R_n=\frac{W_{n-1}W_{n+1}}{W_n^2},
  \qquad
  v_n=v([n]D).
  \label{eq:Lucas-three-elliptic-sequences}
\end{equation}
They are linked, but they are not interchangeable.

\begin{proposition}[EDS, Kummer ratio, and oriented Lucas coordinate]
\label{prop:Lucas-three-elliptic-sequences}
The three sequences in \eqref{eq:Lucas-three-elliptic-sequences} have the
following precise meanings.
\begin{enumerate}[label=(\roman*)]
  \item \(W_n\) is a division-section value.  Its vanishing detects
  \([n]D=O\), and its natural recurrence is Ward's homogeneous bilinear
  recurrence.
  \item \(R_n=d_E-X([n]D_W)\) is an even Kummer coordinate.  It is
  unchanged by \([n]D_W\mapsto-[n]D_W\) and satisfies the multiplicative
  QRT recurrence \eqref{eq:Lucas-R-multiplicative-QRT}.
  \item \(v_n\) is an oriented degree-two elliptic coordinate.  It changes
  sign under elliptic negation, satisfies the symmetric McMillan recurrence
  \eqref{eq:Lucas-normalized-recurrence}, and is recovered from the EDS by
  \eqref{eq:Lucas-direct-v-from-EDS}.
\end{enumerate}
Thus neither \(R_n\) nor \(v_n\) is an EDS term, and \(R_n\) alone cannot
recover the orientation retained by \(v_n\).
\end{proposition}

\begin{proof}
Part (i) follows from the definition
\eqref{eq:Lucas-algebraic-EDS-definition}, the zero criterion
\eqref{eq:Lucas-EDS-zero-torsion}, and Ward's identity
\eqref{eq:Lucas-Ward-recurrence}.  Part (ii) is
\eqref{eq:Lucas-explicit-X-bridge}; since the Weierstrass \(X\)-coordinate
is invariant under negation, so is \(R_n\).  Its recurrence is
Proposition~\ref{prop:Lucas-R-multiplicative-QRT}.  For part (iii), the
Jacobi negation formula \eqref{eq:Lucas-Jacobi-negation} gives
\(v(-P)=-v(P)\), while the recurrence and EDS recovery are
\eqref{eq:Lucas-normalized-recurrence} and
\eqref{eq:Lucas-direct-v-from-EDS}.  These different divisor and parity
properties rule out an equality of the three sequences in general.
\end{proof}

\begin{definition}[Elliptic net]
\label{def:Lucas-elliptic-net}
Let \(A\) be a finite-rank free abelian group and \(R\) an integral
domain.  An elliptic net is a function \(W:A\to R\) with \(W(0)=0\)
satisfying, for all \(p,q,r,s\in A\),
\begin{equation}
\begin{aligned}
 &W(p+q+s)W(p-q)W(r+s)W(r)\\
 &\quad+W(q+r+s)W(q-r)W(p+s)W(p)\\
 &\quad+W(r+p+s)W(r-p)W(q+s)W(q)=0.
\end{aligned}
\label{eq:Lucas-elliptic-net-recurrence}
\end{equation}
\end{definition}

In rank one, the recurrence in
\eqref{eq:Lucas-elliptic-net-recurrence} is Ward's EDS recurrence after
normalization.  Higher rank retains several independent elliptic-point
indices.  After the three sequence theories have been distinguished, the
sigma-function construction below proves that these nets exist and explains
how adjacent-state coordinates arise from them.

\begingroup
\small
\setlength{\tabcolsep}{3pt}
\renewcommand{\arraystretch}{1.12}
\begin{longtable}{L{0.13\textwidth}L{0.27\textwidth}L{0.27\textwidth}L{0.27\textwidth}}
\caption{The exact distinction among the classical Lucas pair, an elliptic
 divisibility sequence, and the QRT adjacent-state elliptic Lucas system.}
\label{tab:Lucas-three-sequence-distinction}\\
\toprule
Feature
& Classical Lucas pair
& Elliptic divisibility sequence / rank-one elliptic net
& QRT adjacent-state elliptic Lucas system\\
\midrule
\endfirsthead
\multicolumn{4}{c}{\small\itshape Table~\thetable\ continued}\\
\toprule
Feature
& Classical Lucas pair
& Elliptic divisibility sequence / rank-one elliptic net
& QRT adjacent-state elliptic Lucas system\\
\midrule
\endhead
\midrule
\multicolumn{4}{r}{\small\itshape Continued on the next page}\\
\endfoot
\bottomrule
\endlastfoot
Ambient group
& A one-dimensional split or nonsplit torus; equivalently the smooth locus
  of a nodal genus-one curve
& A smooth elliptic curve with a marked point \(D\)
& The same kind of pointed elliptic curve, together with a chosen degree-two
  state coordinate and the adjacent-state embedding; its square may descend
  to a Kummer quotient\\

Term
& \(U_n=(q^n-q^{-n})/(q-q^{-1})\) and
  \(V_n=q^n+q^{-n}\)
& A division-function value
  \(\widehat W_n=\sigma(nz)/\sigma(z)^{n^2}\), or the algebraically
  normalized value \(W_n=(-1)^{n-1}\widehat W_n=\psi_n(D)\)
& The coordinate pair \(S_n=(v_n,v_{n+1})\), with full-point companion
  \(w_n\)\\

Defining recurrence
& A second-order linear recurrence
& Ward's homogeneous bilinear recurrence
& An autonomous second-order rational QRT recurrence, supplemented by a
  companion identity\\

Invariant geometry
& The Cassini conic \(x^2+y^2-sxy=1\)
& The zero divisor of the division section and the elliptic curve supporting
  it; no adjacent-state conic is part of the definition
& A fixed smooth symmetric biquadratic state curve and its Jacobi companion\\

Vanishing
& Governed by the order of \(q\) in the torus, with the usual caveats for the
  chosen trace or anti-trace normalization
& \(W_n=0\) exactly when \([n]D=O\)
& \(v_n=0\) means that \([n]D\) lies in the zero divisor of the chosen
  state coordinate; by itself it is not a torsion-to-identity criterion\\

Fast index arithmetic
& Polynomial addition and doubling identities
& Ward/EDS window formulas such as
  \eqref{eq:Lucas-Ward-odd-doubling}--\eqref{eq:Lucas-Ward-even-doubling}
& Rational branches \(\mathcal L_0,\mathcal L_1\), complete projective
  atlases, and the QRT binary ladder\\

Recovery of the point
& The pair \((U_n,V_n)\) determines the torus element up to the declared
  normalization
& Ratios of neighbouring \(W_j\)'s recover \(X([n]D)\), and \(W_{2n}/W_n^4\)
  recovers the odd coordinate
& The pair \((v_n,w_n)\) recovers the full elliptic point, while
  \(v_n^2\) is the associated even quotient coordinate\\

Higher-rank extension
& Multivariable toric recurrences may be formed, but are not the object
  defined here
& Elliptic nets indexed by \(\mathbb Z^r\)
& Multidirectional adjacent states obtained by projecting the same
  \(\mathbb Z^r\)-lattice of elliptic points\\
\end{longtable}
\endgroup

\begin{proposition}[The bridge is by ratios, not by equality of terms]
\label{prop:Lucas-EDS-not-identical-sequences}
For the present Jacobi normalization,
\begin{equation}
  v_n
  =\frac{
     W_n^2\bigl(W_2W_n^2-4W_{n-1}W_{n+1}\bigr)
    }{W_{2n}}
  \label{eq:Lucas-three-sequence-ratio-bridge}
\end{equation}
on the regular oriented chart, but in general \(v_n\ne W_n\).  The
torsion test is \(W_n=0\), or equivalently the full-state test
\(S_n=S_0\); the single condition \(v_n=0\) is not a torsion test.
\end{proposition}

\begin{proof}
Equation~\eqref{eq:Lucas-three-sequence-ratio-bridge} is
\eqref{eq:Lucas-direct-v-from-EDS}.  The two sides have different
divisor-theoretic meanings.  The term \(W_n\) is the value at \(D_W\)
of the division section \(\psi_n\), whose zero divisor detects
\([n]D=O\).  The term \(v_n\) is the value at \([n]D\) of the fixed rational
function \(v=2(X+h)/Y\).  Directly from
\(Y^2=(X+h)(X^2-4\lambda)\), the function \(v\) vanishes at the identity
and at \(T_h=(-h,0)\), and has simple poles at the two points
\(T_\pm=(\pm2\sqrt\lambda,0)\).  Hence its divisor is
\[
  (O)+(T_h)-(T_+)-(T_-).
\]
This fixed coordinate divisor is fundamentally different from the
\(n\)-division divisor of \(\psi_n\), so equality of the two sequences is
neither part of the construction nor preserved by changing the state
coordinate.

The equivalence \(W_n=0\Longleftrightarrow[n]D=O\), under the separability
hypothesis on \([n]\), is \eqref{eq:Lucas-EDS-zero-torsion}; the exact
state-period theorem gives
\([n]D=O\Longleftrightarrow S_n=S_0\).  By contrast,
\(v_n=0\) merely says that \([n]D\) lies in the zero divisor
\((O)+(T_h)\) of \(v\).  It may therefore also represent the nonzero
two-torsion point \(T_h\), and is not by itself a torsion-to-identity test.
\end{proof}

Classical Lucas theory meets both elliptic theories on a singular boundary.
The sequence \(U_n\) is simultaneously the limiting QRT coordinate sequence
and, after the standard singular normalization, a toric divisibility
sequence.  This common boundary explains the similarities displayed in the
preceding comparison without erasing the different definitions of the
three systems on a smooth elliptic fibre.

\phantomsection
\label{subsec:Lucas-sigma-elliptic-net}
The preceding distinctions do not separate the three theories completely.
The sigma function provides a common analytic source for division sequences,
coordinate addition laws, and elliptic nets, and the resulting identities
produce the maps appearing in the total diagram.

\begin{proposition}[Divisor of the Jacobi coordinate]\label{prop:Lucas-divisor-v}
The sigma description is especially concrete for the Jacobi coordinate
\(v=2(X+h)/Y\) arising in the explicit bridge above.  Over an algebraic
closure, let
\begin{equation}
  T_h=(-h,0),
  \qquad
  T_+=(2\sqrt\lambda,0),
  \qquad
  T_-=(-2\sqrt\lambda,0)
  \label{eq:Lucas-specific-two-torsion}
\end{equation}
be the three nonzero two-torsion points on
\eqref{eq:Lucas-specific-Weierstrass}.  The coordinate functions have
divisors
\begin{align}
  \operatorname{div}(X+h)&=2(T_h)-2(O),
  \label{eq:Lucas-divisor-X-plus-h}\\
  \operatorname{div}(Y)&=(T_h)+(T_+)+(T_-)-3(O).
  \label{eq:Lucas-divisor-Y}
\end{align}
Consequently,
\begin{equation}
  \operatorname{div}(v)
  =(O)+(T_h)-(T_+)-(T_-).
  \label{eq:Lucas-divisor-v}
\end{equation}
\end{proposition}

\begin{proof}
The root \(X=-h\) of the cubic in
\eqref{eq:Lucas-specific-Weierstrass} is simple because the curve is
smooth.  Near \(T_h\), the function \(Y\) is a local parameter, and
\[
  X+h=\frac{Y^2}{X^2-4\lambda}
\]
has a zero of order two because \(h^2-4\lambda\ne0\).  The function
\(X\) has a pole of order two at \(O\), proving
\eqref{eq:Lucas-divisor-X-plus-h}.  The function \(Y\) vanishes simply
at each of the three simple roots of the cubic and has a pole of order
three at \(O\), proving \eqref{eq:Lucas-divisor-Y}.  Since
\(v=2(X+h)/Y\), subtraction of divisors gives
\eqref{eq:Lucas-divisor-v}.
\end{proof}

\begin{proposition}[Sigma quotient for the elliptic Lucas coordinate]
\label{prop:Lucas-v-sigma-quotient}
Assume \(k=\mathbb C\) and identify
\(E(\mathbb C)\simeq\mathbb C/\Lambda\).  Choose half-period
representatives \(\omega_h,\omega_+,\omega_-\) for
\(T_h,T_+,T_-\).  They may be adjusted by lattice elements so that
\begin{equation}
  \omega_h=\omega_++\omega_-.
  \label{eq:Lucas-half-period-relation}
\end{equation}
If \(z_D\) represents the marked point and \(v(z_D)=1\), then
\begin{equation}
\boxed{
  v_n=
  \frac{\sigma(nz_D)\sigma(nz_D-\omega_h)}
       {\sigma(nz_D-\omega_+)\sigma(nz_D-\omega_-)}
  \frac{\sigma(z_D-\omega_+)\sigma(z_D-\omega_-)}
       {\sigma(z_D)\sigma(z_D-\omega_h)}.
}
\label{eq:Lucas-vn-specific-sigma-quotient}
\end{equation}
\end{proposition}

\begin{proof}
Define
\begin{equation}
  F_v(z)=
  \frac{\sigma(z)\sigma(z-\omega_h)}
       {\sigma(z-\omega_+)\sigma(z-\omega_-)}.
  \label{eq:Lucas-v-sigma-quotient}
\end{equation}
The sigma function has a simple zero at every lattice point and no other
zeros.  Therefore the divisor of \(F_v\) on \(\mathbb C/\Lambda\) is
\[
  (O)+(T_h)-(T_+)-(T_-),
\]
which agrees with \eqref{eq:Lucas-divisor-v}.

It remains to prove that \(F_v\) is elliptic rather than merely
quasi-periodic.  For every period \(\Omega\in\Lambda\), the sigma
quasi-periodicity law has the form
\[
  \sigma(z+\Omega)
  =\varepsilon(\Omega)
   \exp\!\left(\eta_\Omega\left(z+\frac{\Omega}{2}\right)\right)
   \sigma(z),
  \qquad \varepsilon(\Omega)\in\{\pm1\}.
\]
Applying this identity to the four factors in
\eqref{eq:Lucas-v-sigma-quotient}, the two multipliers
\(\varepsilon(\Omega)\) in the numerator and the two in the denominator
cancel.  The terms involving \(z\) and \(\Omega/2\) also cancel because
the numerator and denominator each contain two factors.  The remaining
multiplier is
\[
  \exp\!\left(
    \eta_\Omega(\omega_++\omega_--\omega_h)
  \right)=1
\]
by \eqref{eq:Lucas-half-period-relation}.  Hence
\(F_v(z+\Omega)=F_v(z)\).

The quotient \(v/F_v\) is now an elliptic function with neither zeros nor
poles.  It is therefore constant on the compact Riemann surface
\(\mathbb C/\Lambda\).  Evaluating at \(z_D\), where \(v(z_D)=1\), gives
\[
  v(z)=\frac{F_v(z)}{F_v(z_D)}.
\]
Substitution of \(z=nz_D\) yields
\eqref{eq:Lucas-vn-specific-sigma-quotient}.
\end{proof}

Formula~\eqref{eq:Lucas-vn-specific-sigma-quotient} is the analytic
counterpart of the algebraic EDS expression
\eqref{eq:Lucas-direct-v-from-EDS}.  It also explains why the Lucas
coordinate is not itself an EDS term: it is a ratio of four translated sigma
sections whose divisor records the oriented degree-two coordinate.

Let \(z_1,\ldots,z_r\in\mathbb C\) represent points
\(P_1,\ldots,P_r\in E(\mathbb C)=\mathbb C/\Lambda\).  For
\(\mathbf v=(v_1,\ldots,v_r)\in\mathbb Z^r\), define Stange's analytic
net function \cite{StangeEllipticNets2011}
\begin{equation}
 \Omega_{\mathbf v}(\mathbf z;\Lambda)
 =\frac{
   \sigma(v_1z_1+\cdots+v_rz_r;\Lambda)
 }{
   \displaystyle
   \prod_{i=1}^r
     \sigma(z_i;\Lambda)^{
       2v_i^2-\sum_{j=1}^r v_iv_j
     }
   \prod_{1\le i<j\le r}
     \sigma(z_i+z_j;\Lambda)^{v_iv_j}
 }.
 \label{eq:Lucas-Stange-Omega-definition}
\end{equation}
For \(\mathbf v=0\), set \(\Omega_{\mathbf 0}=0\).  In rank one this
reduces exactly to the analytic normalization
\begin{equation}
  \Omega_n(z;\Lambda)
  =\frac{\sigma(nz;\Lambda)}{\sigma(z;\Lambda)^{n^2}}
  =\widehat W_n(z).
  \label{eq:Lucas-Stange-rank-one-EDS}
\end{equation}
The algebraic division-polynomial normalization used in the explicit bridge
is related by \(W_n=(-1)^{n-1}\widehat W_n\), as in
\eqref{eq:Lucas-algebraic-analytic-EDS-sign}.

\begin{theorem}[Sigma functions produce the net underlying the state lattice]
\label{thm:Lucas-sigma-net-state-lattice}
The function
\begin{equation}
  \mathbf v\longmapsto\Omega_{\mathbf v}(\mathbf z;\Lambda)
  \label{eq:Lucas-Omega-is-net-map}
\end{equation}
is an elliptic net.  For \(p,q\in\mathbb Z^r\) on the regular open set,
\begin{equation}
  \frac{
    \Omega_{p+q}\Omega_{p-q}
  }{
    \Omega_p^2\Omega_q^2
  }
  =\wp(q\cdot\mathbf z)-\wp(p\cdot\mathbf z),
  \label{eq:Lucas-net-coordinate-difference}
\end{equation}
where \(p\cdot\mathbf z=\sum_i p_iz_i\).  In particular, for the
rank-two tuple \((P,D)\),
\begin{equation}
  \mathcal W_{P,D}(m,r)
  :=\Omega_{(m,r)}(z_P,z_D;\Lambda)
  \label{eq:Lucas-rank-two-net-PD}
\end{equation}
encodes the entire lattice of points
\begin{equation}
  [m]P+[r]D.
  \label{eq:Lucas-rank-two-point-lattice}
\end{equation}
The adjacent state
\begin{equation}
 \mathcal S_D([m]P+[r]D)
 =\bigl(
    \kappa([m]P+[r]D),
    \kappa([m]P+[r+1]D)
  \bigr)
 \label{eq:Lucas-rank-two-net-adjacent-state}
\end{equation}
is therefore a Kummer-coordinate projection of two neighbouring columns of
this rank-two elliptic net.
\end{theorem}

\begin{proof}
Theorem~3.7 of Stange~\cite{StangeEllipticNets2011} applies to the
fixed lattice \(\Lambda\) and the fixed tuple
\((z_1,\ldots,z_r)\) and proves that
\eqref{eq:Lucas-Omega-is-net-map} satisfies
\eqref{eq:Lucas-elliptic-net-recurrence}.  Its construction is meromorphic
on the tuple space; throughout this theorem we work on the nondegenerate
open locus on which the denominators in
\eqref{eq:Lucas-Stange-Omega-definition} are nonzero, and the algebraic net
polynomials give the corresponding continuation outside that chart.  We
verify next the coordinate identity needed for the present bridge.  Apply
the sigma addition identity
\eqref{eq:Lucas-sigma-addition-identity} with
\[
  u=p\cdot\mathbf z,
  \qquad
  v=q\cdot\mathbf z.
\]
The numerator is
\(\sigma((p+q)\cdot\mathbf z)\sigma((p-q)\cdot\mathbf z)\).
When each sigma factor is replaced by the definition
\eqref{eq:Lucas-Stange-Omega-definition}, all factors involving the fixed
\(z_i\)'s and \(z_i+z_j\)'s cancel.  The cancellation follows from the
quadratic identity
\[
  f(p+q)+f(p-q)-2f(p)-2f(q)=0
\]
for every quadratic form \(f\).  The remaining quotient is exactly the
left-hand side of \eqref{eq:Lucas-net-coordinate-difference}; hence the
sigma addition formula gives the asserted right-hand side.

For the tuple \((P,D)\), the linear combination represented by
\((m,r)\) is \([m]P+[r]D\).  Let \(e_2=(0,1)\).  Definition
\eqref{eq:Lucas-Stange-Omega-definition} gives
\(\Omega_{e_2}=1\).  Taking \(p=(m,r)\) and \(q=e_2\) in
\eqref{eq:Lucas-net-coordinate-difference} yields the explicit recovery
formula
\begin{equation}
 \wp(mz_P+rz_D)
 =\wp(z_D)
  -\frac{
     \Omega_{(m,r+1)}\Omega_{(m,r-1)}
   }{
     \Omega_{(m,r)}^2
   }.
 \label{eq:Lucas-rank-two-net-coordinate-recovery}
\end{equation}
Thus the net values \(\mathcal W_{P,D}(m,r)\) recover the Kummer
coordinate of every regular point in
\eqref{eq:Lucas-rank-two-point-lattice}.  Replacing \(r\) by \(r+1\)
recovers the adjacent point differing by \(D\), which proves
\eqref{eq:Lucas-rank-two-net-adjacent-state}.  Identity
\eqref{eq:Lucas-rank-two-net-coordinate-recovery} has a homogeneous
projective interpretation at the poles, exactly as in
\eqref{eq:Lucas-EDS-window-state-recovery-projective}.
\end{proof}

\begin{remark}[How the affine-index semigroup appears in the net]
\label{rem:Lucas-affine-semigroup-in-net}
For fixed \((P,D)\), the map \((m,r)\mapsto[m]P+[r]D\) is the geometric
meaning of the rank-two net index.  The state map
\(\mathfrak F_{m,r}\) studied in
Subsection~\ref{subsec:QRT-affine-index-semigroup} applies the same affine
combination to a variable point.  Thus the composition law
\[
  (m,r)\circ(n,s)=(mn,ms+r)
\]
is visible both in the affine-index semigroup of state maps and in the
corresponding reindexing of the rank-two net.  When \(P=[a]D\), the
point lattice collapses along
\([m]P+[r]D=[ma+r]D\).  After passing to the corresponding rank-one
specialization of the algebraic net polynomials, the independent index is
therefore \(ma+r\), and the construction reduces to the rank-one EDS of
\(D\).  This formulation avoids substituting a linearly dependent tuple
into an analytic denominator that may vanish.
\end{remark}

The sigma formula finally clarifies the exact relationship among the three
sequence theories.  The EDS/net records a normalized division section
\(\Omega_{\mathbf v}\).  The elliptic Lucas system records rational
coordinates of the point represented by \(\mathbf v\), arranged in an
adjacent state.  Classical Lucas theory is the toric degeneration in which
the sigma-function geometry degenerates to multiplicative exponentials.
Thus the three sequences are different functorial shadows of the same group
indexing, and the bridge diagrams above state precisely which rational maps
connect them.

\subsection[Characteristic-two companion and complete branches]
{Characteristic two: an Artin--Schreier companion and complete branches}
\label{subsec:Lucas-characteristic-two}

We now assume that \(\charac(k)=2\), that \(k\) is perfect, and that
\begin{equation}
  \alpha\beta\gamma\ne0.
  \label{eq:Lucas-char2-smoothness}
\end{equation}
The nonzero mixed coefficient is essential.  To verify this without
appealing to a later classification theorem, homogenize the equation as
\[
 F=X^2Y^2+\alpha(X^2W^2+Y^2Z^2)
   +\beta XYZW+\gamma Z^2W^2.
\]
On the affine chart \(Z=W=1\), characteristic two gives
\[
 F_x=\beta y,\qquad F_y=\beta x.
\]
If \(\beta\ne0\), an affine singular point would therefore have
\(x=y=0\), which is impossible because \(F(0,0)=\gamma\ne0\).
At a boundary point with \(Z=0\), the equation forces
\(Y^2=\alpha W^2\).  Since \(\alpha\ne0\), the coordinates
\(X,Y,W\) are all nonzero, and
\[
  F_Z=\beta XYW\ne0.
\]
The chart \(W=0\) is symmetric.  Hence
\(\alpha\beta\gamma\ne0\) implies smoothness.  Conversely, if
\(\beta=0\), perfectness of \(k\) supplies
\(r^2=\alpha\) and \(s^2=\gamma\), and
\[
 F=(XY+rXW+rYZ+sZW)^2.
\]
The projective curve is then nonreduced.  Thus the mixed term is exactly
what keeps the characteristic-two double cover separable and the
biquadratic smooth.

Before normalizing, define the characteristic-two companion
\begin{equation}
  \mathcal W=(x^2+\alpha)y.
  \label{eq:Lucas-char2-original-companion}
\end{equation}
Multiplying the original QRT equation by \(x^2+\alpha\) gives the
binary quartic
\begin{equation}
  \mathcal W^2+\beta x\mathcal W
  =\alpha x^4+(\alpha^2+\gamma)x^2+\alpha\gamma.
  \label{eq:Lucas-char2-original-quartic}
\end{equation}
This identity shows exactly why the \(\beta xy\)-term cannot be removed:
it becomes the Artin--Schreier cross term \(\beta x\mathcal W\), whose
nonvanishing keeps the double cover separable.

Choose
\begin{equation}
  t^2=\frac{\gamma}{\alpha},
  \qquad
  \lambda=\frac{\gamma}{\alpha^2},
  \qquad
  \eta=\frac{\beta}{\alpha},
  \qquad
  \mu^2=\lambda.
  \label{eq:Lucas-char2-normalized-parameters}
\end{equation}
All square roots are unique in a perfect field of characteristic two.
Writing \(u_n=tv_n\), the initial state \((0,t)\) becomes
\((v_0,v_1)=(0,1)\), and the normalized state curve is
\begin{equation}
  \mathcal B_{\lambda,\eta}^{(2)}:
  \qquad
  \lambda x^2y^2+x^2+y^2+\eta xy+1=0.
  \label{eq:Lucas-char2-normalized-curve}
\end{equation}
The recurrence is
\begin{equation}
  v_{n+1}+v_{n-1}
  =\frac{\eta v_n}{1+\lambda v_n^2}.
  \label{eq:Lucas-char2-recurrence}
\end{equation}
Define
\begin{equation}
  w_n=(1+\lambda v_n^2)v_{n+1}.
  \label{eq:Lucas-char2-companion}
\end{equation}
Then
\begin{equation}
  w_n^2+\eta v_nw_n
  =\lambda v_n^4+(1+\lambda)v_n^2+1.
  \label{eq:Lucas-char2-companion-quartic}
\end{equation}
Indeed, multiply \eqref{eq:Lucas-char2-normalized-curve} by
\(1+\lambda x^2\) and use \(w=(1+\lambda x^2)y\).  The companion curve
contains
\begin{equation}
  O=(0,1),
  \qquad
  D=(1,\eta).
  \label{eq:Lucas-char2-O-D}
\end{equation}
Transporting the group law from the state curve with origin
\((0,1)\) makes the McMillan shift translation by \(D\).  Unlike the
odd-characteristic companion, the binary companion is an
Artin--Schreier quartic rather than a completed square.

\begin{theorem}[Characteristic-two elliptic-Lucas doubling]
\label{thm:Lucas-char2-fast-doubling}
Let \((x,y)=(v_n,v_{n+1})\) on
\(\mathcal B_{\lambda,\eta}^{(2)}\).  Then
\begin{align}
  v_{2n}
  &=\frac{\eta x^2}{(1+\mu x^2)^2},
  \label{eq:Lucas-char2-v2n}\\
  v_{2n+1}
  &=\frac{(x+y)^2}{(1+\mu xy)^2},
  \label{eq:Lucas-char2-v2n1}\\
  v_{2n+2}
  &=\frac{\eta y^2}{(1+\mu y^2)^2}.
  \label{eq:Lucas-char2-v2n2}
\end{align}
The equalities are projective identities, so the displayed affine
quotients extend through zero denominators by the homogeneous formulas
below.
\end{theorem}

\begin{proof}
Set
\begin{equation}
  p=\frac{\mu}{\eta}x,
  \qquad
  q=\frac{\mu}{\eta}y,
  \qquad
  d=\frac1\eta,
  \qquad
  b=\frac{\mu}{\eta^2}.
  \label{eq:Lucas-char2-exact-state-scaling}
\end{equation}
Substituting the four fractions and placing every term over
\(\eta^4\) gives
\begin{align*}
 p^2q^2+d^2(p^2+q^2)+dpq+b^2
 &=
 \frac{\mu^4x^2y^2+\mu^2(x^2+y^2)
       +\mu^2\eta xy+\mu^2}{\eta^4}\\
 &=\frac{\mu^2}{\eta^4}
   \bigl(\mu^2x^2y^2+x^2+y^2+\eta xy+1\bigr).
\end{align*}
Since \(\mu^2=\lambda\), this is
\begin{equation}
\begin{split}
  p^2q^2+d^2(p^2+q^2)+dpq+b^2
  =\frac{\lambda}{\eta^4}
   \bigl(\lambda x^2y^2+x^2+y^2+\eta xy+1\bigr).
  \label{eq:Lucas-char2-exact-state-relation}
\end{split}
\end{equation}
Thus \((p,q)\) is the exact binary adjacent-state model
\eqref{eq:QRT-binary-pointed-relation} with known-difference coordinate
\(d\) and curve constant \(b\).

On the shifted binary Weierstrass model, the point
\(T_2=(0,b)\) is the nonzero two-torsion point.  The state
\((p,q)=(0,\mu/\eta)\) is obtained from the elliptic point \(T_2\),
because translation by \(T_2\) sends a Kummer coordinate \(z\) to
\(b/z\), and \(b/d=\mu/\eta\).  Hence the Lucas origin is \(T_2\), not
the Weierstrass point at infinity.  With this translated origin,
doubling is \(P\mapsto2P+T_2\).

The binary Kummer doubling identity
\eqref{eq:binary-Kummer-double-affine} gives
\[
  x(2P)=\frac{(p^2+b)^2}{p^2}.
\]
Translation by \(T_2\) then gives
\[
  x(2P+T_2)=\frac{b}{x(2P)}
  =\frac{bp^2}{(p^2+b)^2}.
\]
Substitution of \eqref{eq:Lucas-char2-exact-state-scaling}, followed by
division by the scaling factor \(\mu/\eta\), yields
\eqref{eq:Lucas-char2-v2n}.

For the middle coordinate, the Lucas sum is the translated-group sum of
points with Kummer coordinates \(p\) and \(q\).  In the original
Weierstrass group this is \(P+Q+T_2\), while
\(x(P-Q)=d\).  The binary Kummer product identity gives
\[
  x(P+Q)=\frac{(pq+b)^2}{d(p+q)^2}.
\]
After translation by \(T_2\),
\[
  x(P+Q+T_2)
  =\frac{bd(p+q)^2}{(pq+b)^2}.
\]
Substitution of \eqref{eq:Lucas-char2-exact-state-scaling} reduces this
expression to \eqref{eq:Lucas-char2-v2n1}.  The formula for
\(v_{2n+2}\) follows from \eqref{eq:Lucas-char2-v2n} with \(y\) in
place of \(x\).
\end{proof}

Let \(x=X/Z\) and \(y=Y/W\).  The characteristic-two projective
branches are
\begin{align}
  A_1^{(2)}&=\eta X^2Z^2,
  &A_0^{(2)}&=(Z^2+\mu X^2)^2,
  \label{eq:Lucas-char2-A-pair}\\
  B_1^{(2)}&=(XW+YZ)^2,
  &B_0^{(2)}&=(ZW+\mu XY)^2,
  \label{eq:Lucas-char2-B-pair}\\
  C_1^{(2)}&=\eta Y^2W^2,
  &C_0^{(2)}&=(W^2+\mu Y^2)^2.
  \label{eq:Lucas-char2-C-pair}
\end{align}
They define
\begin{equation}
  \mathcal L_0^{(2)}=((A_1^{(2)}:A_0^{(2)}),
                      (B_1^{(2)}:B_0^{(2)})),
  \qquad
  \mathcal L_1^{(2)}=((B_1^{(2)}:B_0^{(2)}),
                      (C_1^{(2)}:C_0^{(2)})).
  \label{eq:Lucas-char2-projective-branches}
\end{equation}

\begin{theorem}[Complete binary branches and exact cost]
\label{thm:Lucas-char2-complete-cost}
On every smooth curve \eqref{eq:Lucas-char2-normalized-curve}, all three
pairs in
\eqref{eq:Lucas-char2-A-pair}--\eqref{eq:Lucas-char2-C-pair} are
geometrically base-point-free.  Each branch in
\eqref{eq:Lucas-char2-projective-branches} costs
\begin{equation}
  \boxed{4\M+5\Sqr+3\Cmul}.
  \label{eq:Lucas-char2-branch-cost}
\end{equation}
\end{theorem}

\begin{proof}
For the zero branch, compute
\[
  X_2=X^2,
  \qquad
  Z_2=Z^2,
  \qquad
  T=X_2Z_2,
\]
which costs \(1\M+2\Sqr\).  Then
\[
  A_1^{(2)}=\eta T,
  \qquad
  A_0^{(2)}=(Z_2+\mu X_2)^2
\]
uses \(1\Sqr+2\Cmul\).  For the middle coordinate, use the Karatsuba
schedule
\[
  P=XY,
  \qquad
  Q=ZW,
  \qquad
  R=(X+Z)(Y+W).
\]
Then
\[
  XW+YZ=R+P+Q,
  \qquad
  ZW+\mu XY=Q+\mu P.
\]
This uses \(3\M+2\Sqr+1\Cmul\).  The total is
\eqref{eq:Lucas-char2-branch-cost}.  The one branch is symmetric.

For base-point freeness of \((A_1^{(2)},A_0^{(2)})\), a common zero
would require \(XZ=0\) and \(Z^2+\mu X^2=0\).  If \(X=0\), then
\(Z\ne0\) and the second expression is nonzero; if \(Z=0\), then
\(X\ne0\) and it equals \(\mu X^2\ne0\).  Thus there is no common
zero.  The proof for \(C\) is identical.

For the middle pair, suppose
\[
  XW+YZ=0,
  \qquad
  ZW+\mu XY=0.
\]
No projective coordinate can vanish, so the point is affine.  The first
equality gives \(x=y\), and the second gives \(1+\mu x^2=0\).  Since
\(\lambda=\mu^2\), substitution into
\eqref{eq:Lucas-char2-normalized-curve} gives
\[
  \lambda x^4+x^2+x^2+\eta x^2+1
  =1+\eta x^2+1
  =\eta x^2,
\]
which is nonzero because \(\eta\ne0\).  This contradiction proves that
\((B_1^{(2)},B_0^{(2)})\) is base-point-free.
\end{proof}

In characteristic two, the state reflection is simply
\begin{equation}
  \rho_2(x,y)=(y,x),
  \label{eq:Lucas-char2-state-reflection}
\end{equation}
because signs disappear and the first coordinate is Kummer-even.  The
state reflection sends the parameter \(P\) to \(-P-D\).  Therefore
\(\rho_2\mathcal L_0^{(2)}\rho_2\) sends \(P\) successively to
\(-P-D\), then to \(-2P-2D\), and finally to \(2P+D\), which is the
index action of \(\mathcal L_1^{(2)}\).  Hence
\begin{equation}
  \mathcal L_1^{(2)}
  =\rho_2\mathcal L_0^{(2)}\rho_2.
  \label{eq:Lucas-char2-branch-conjugacy}
\end{equation}
Hence the complete binary formulas also admit a single arithmetic core
with conditional input and output swaps.

The count \eqref{eq:Lucas-char2-branch-cost} is exactly the compiled
fixed-difference count already obtained in
Proposition~\ref{prop:QRT-binary-pointed-state-doubling}.  The binary Lucas
reconstruction is therefore a particularly symmetric, single-formula
realization of binary Kummer \texttt{xDBLADD}.  Its projective branch is
geometrically complete and uses one arithmetic core with coordinate-swap
branch selection.

For completeness, write \(r_2^2=\alpha\) and \(s_2^2=\gamma\), so that
\(t=s_2/r_2\) and \(\mu=s_2/\alpha\).  Returning to the original
coordinates gives
\begin{align}
  u_{2n}
  &=\frac{\beta s_2u_n^2}
          {r_2(u_n^2+s_2)^2},
  \label{eq:Lucas-char2-original-u2n}\\
  u_{2n+1}
  &=\frac{r_2s_2(u_n+u_{n+1})^2}
          {(u_nu_{n+1}+s_2)^2},
  \label{eq:Lucas-char2-original-u2n1}\\
  u_{2n+2}
  &=\frac{\beta s_2u_{n+1}^2}
          {r_2(u_{n+1}^2+s_2)^2}.
  \label{eq:Lucas-char2-original-u2n2}
\end{align}
These identities follow by substituting \(u_n=tv_n\) into
\eqref{eq:Lucas-char2-v2n}--\eqref{eq:Lucas-char2-v2n2} and using
\[
  \frac{\mu}{t^2}=\frac1{s_2},
  \qquad
  \frac{s_2^2}{t}=r_2s_2,
  \qquad
  \frac{t\eta}{t^2}\,s_2^2=\frac{\beta s_2}{r_2}.
\]

\subsubsection{Exact small-constant orbits in characteristic two}
\label{subsubsec:Lucas-char2-small-constants}

The phrase ``normalize \(\mu\) or \(\eta\) to a small constant'' must
be tied to an explicit coordinate-equivalence group.  Here we use the
monomial equivalence generated by diagonal scalings, interchange of the
two \(\mathbf P^1\)-factors, simultaneous inversion, and inversion of
one coordinate followed by renormalization.  The following theorem shows
directly that a normalized parameter pair has an orbit of size at most
two.

\begin{theorem}[Monomial orbit of the binary Lucas parameters]
\label{thm:Lucas-char2-mu-eta-orbit}
Let
\[
  \mathcal B^{(2)}_{\mu^2,\eta}:
  \quad
  \mu^2x^2y^2+x^2+y^2+\eta xy+1=0,
  \qquad \mu\eta\ne0.
\]
Under diagonal scalings, factor interchange, simultaneous inversion,
and one-coordinate inversion followed by renormalization, the only
possibly distinct normalized parameter pair is
\begin{equation}
  (\mu,\eta)
  \longmapsto
  \iota(\mu,\eta)
  =\left(\mu^{-1},\frac{\eta}{\mu}\right).
  \label{eq:Lucas-char2-mu-eta-involution}
\end{equation}
The transformation is an involution.  It preserves the binary
Weierstrass invariants
\begin{equation}
  a_Q=\frac{1+\mu^2}{\eta^2},
  \qquad
  b_Q=\frac{\mu}{\eta^2},
  \label{eq:Lucas-char2-mu-eta-Weierstrass-invariants}
\end{equation}
and therefore preserves both the ordinary elliptic curve and its
Artin--Schreier twist class.

Let \(S_\mu,S_\eta\subset k^\times\) be any prescribed sets of
implementation-friendly constants.  An equivalent normalized model
with \((\mu',\eta')\in S_\mu\times S_\eta\) exists if and only if
\begin{equation}
  (\mu,\eta)\in S_\mu\times S_\eta
  \quad\text{or}\quad
  \left(\mu^{-1},\frac{\eta}{\mu}\right)
  \in S_\mu\times S_\eta.
  \label{eq:Lucas-char2-cheap-set-criterion}
\end{equation}
\end{theorem}

\begin{proof}
It is useful first to derive the nontrivial generator directly.  Write the
normalized binary model in the auxiliary form
\[
  \mathcal N_{b,g}:\qquad
  gx^2y^2+x^2+y^2+bxy+1=0,
  \qquad bg\ne0.
\]
Invert the first coordinate by setting \(x=X^{-1}\), multiply by
\(X^2\), and then put \(Y=hy\), where \(h^2=g\).  The transformed
equation is
\[
  g^{-1}X^2Y^2+X^2+Y^2+\frac b hXY+1=0.
\]
Thus one-coordinate inversion followed by renormalization sends
\[
  (b,g)\longmapsto\left(\frac b h,\frac1g\right).
\]
For the Lucas model, \(b=\eta\), \(g=\mu^2\), and we may take
\(h=\mu\).  Hence the induced transformation is precisely
\[
  (\mu,\eta)
  \longmapsto
  \left(\mu^{-1},\frac{\eta}{\mu}\right).
\]
Interchanging the two factors changes nothing because the equation is
symmetric.  Diagonal scalings that preserve the displayed normalization
and simultaneous inversion either fix the pair or reproduce the same
representative.  Applying the displayed transformation twice gives
\[
  \left(\mu^{-1},\frac{\eta}{\mu}\right)
  \longmapsto
  \left(\mu,\eta\right),
\]
so the orbit contains at most the two stated pairs.

Substitution gives
\[
  \frac{1+\mu^{-2}}{(\eta/\mu)^2}
  =\frac{1+\mu^2}{\eta^2},
  \qquad
  \frac{\mu^{-1}}{(\eta/\mu)^2}
  =\frac{\mu}{\eta^2},
\]
so \eqref{eq:Lucas-char2-mu-eta-Weierstrass-invariants} is unchanged.
The orbit has precisely the two displayed representatives, possibly
coincident; the set criterion follows.
\end{proof}

For the exact pointed state model of
Theorem~\ref{thm:QRT-binary-pointed-state},
\begin{equation}
  \mu=\frac b{d^2},
  \qquad
  \eta=\frac1d.
  \label{eq:Lucas-char2-mu-eta-pointed}
\end{equation}
The theorem gives the following concrete consequences.

\begin{corollary}[Concrete constant-normalization criteria]
\label{cor:Lucas-char2-concrete-small-constants}
Let \(c\in k^\times\).  Under the monomial equivalence above:
\begin{enumerate}[label=(\roman*)]
  \item \(\mu\) can be normalized to \(1\) if and only if
  \(b=d^2\).
  \item \(\eta\) can be normalized to \(1\) if and only if
  \(d=1\) or \(b=d\).
  \item both \(\mu\) and \(\eta\) can be normalized to \(1\) if and
  only if \(b=d=1\).
  \item \(\mu\) can be normalized to \(c\) if and only if
  \(b=cd^2\) or \(b=c^{-1}d^2\).
  \item \(\eta\) can be normalized to \(c\) if and only if
  \(d=c^{-1}\) or \(b=d/c\).
\end{enumerate}
If \(k_0\subseteq k\) is a subfield, an equivalent representative has
both parameters in \(k_0\) if and only if
\(\mu,\eta\in k_0\) already.
\end{corollary}

\begin{proof}
Insert \eqref{eq:Lucas-char2-mu-eta-pointed} into the two orbit
representatives
\[
  \left(\frac b{d^2},\frac1d\right),
  \qquad
  \left(\frac{d^2}{b},\frac d b\right).
\]
Each assertion is obtained by setting the relevant coordinate equal to
the prescribed constant.  For the subfield statement, if the second orbit
representative belongs to \(k_0\times k_0\), then \(\mu^{-1}\in k_0\), so \(\mu\in k_0\), and
\(\eta/\mu\in k_0\), so \(\eta=\mu(\eta/\mu)\in k_0\).  If instead the
first representative belongs to \(k_0\times k_0\), the conclusion
\(\mu,\eta\in k_0\) is part of the assumption.  Conversely, when
\(\mu,\eta\in k_0\), both
\((\mu,\eta)\) and
\((\mu^{-1},\eta/\mu)\) have coordinates in \(k_0\), because \(k_0\) is
a field and \(\mu\ne0\).
\end{proof}

Thus small constants occur on explicit lower-dimensional parameter
loci; they are not available for a generic pointed binary curve.  This
classification preserves the marked-state isomorphism class and does
not silently change the Artin--Schreier twist in order to obtain a
cheaper constant.

\subsubsection{Binary validation and companion recovery}
\label{subsubsec:Lucas-char2-validation-recovery}

The characteristic-two branch also admits low-incremental-cost state
checks.  Put
\begin{equation}
  P=XY,
  \qquad Q=ZW,
  \qquad K=XW+YZ.
  \label{eq:Lucas-char2-validation-intermediates}
\end{equation}
Then the homogeneous state equation is equivalent to
\begin{equation}
  \mathscr V_2=(Q+\mu P)^2+K^2+\eta PQ=0.
  \label{eq:Lucas-char2-state-validation}
\end{equation}
The first two squares are exactly the denominator and numerator cores
of the middle branch in
\eqref{eq:Lucas-char2-B-pair}.  Hence validation adds only
\begin{equation}
  \boxed{1\M+1\Cmul}.
  \label{eq:Lucas-char2-validation-cost}
\end{equation}

The Artin--Schreier companion
\(w=(1+\lambda x^2)y\) has the two projective recovery charts
\begin{align}
  \mathscr J_{2,W}
  &=\bigl(XW:\ ZY\,ZW+\lambda XY\,XW:\ ZW\bigr),
  \label{eq:Lucas-char2-recovery-W}\\
  \mathscr J_{2,Y}
  &=\bigl(XY:\ XY\,XW+ZY\,ZW+\eta XY\,ZY:\ ZY\bigr).
  \label{eq:Lucas-char2-recovery-Y}
\end{align}
They are weighted triples for the binary quartic
\eqref{eq:Lucas-char2-companion-quartic}.  The first is valid on
\(W\ne0\), the second on \(Y\ne0\), and together they cover the whole
state curve.  Indeed, the first follows from
\[
  w=(1+\lambda x^2)y,
\]
while the state equation gives the complementary identity
\begin{equation}
  w=\eta x+\frac{x^2+1}{y}.
  \label{eq:Lucas-char2-companion-second-chart}
\end{equation}
On an affine final state for which \(\lambda x^2\) is already
available, recovering \(w\) costs one additional multiplication.  In
projective coordinates, the displayed atlas avoids inversions and makes
the exceptional \(y=0\) and \(y=\infty\) states explicit.

\subsection[Torsion, periods, and state-division polynomials]
{Torsion tests, exact state periods, and state-division polynomials}
\label{subsec:Lucas-torsion-periods}

The affine-index interpretation turns torsion into an exact periodicity
statement for the two-coordinate Lucas state.

\begin{proposition}[Exact Lucas-state period and torsion test]
\label{prop:Lucas-exact-period-torsion}
Let $D$ be the marked point defining the elliptic Lucas sequence and
let $S_n=(v_n,v_{n+1})$.  If $D$ has exact order $N$, then $(S_n)$ has
exact period $N$.  In the odd-characteristic normalization,
\begin{equation}
  [N]D=O
  \quad\Longleftrightarrow\quad
  S_N=S_0=(0,1).
  \label{eq:Lucas-torsion-state-test}
\end{equation}
The order is exactly $N$ if the equality fails for every proper
positive divisor of $N$.  A binary word in $\mathcal L_0,\mathcal L_1$
computes $S_N$ in $O(\log N)$ state steps.
\end{proposition}

\begin{proof}
The state map is an isomorphism and $S_n=\mathscr S([n]D)$.  Hence
$S_{n+t}=S_n$ for every $n$ if and only if translation by $[t]D$ is the
identity, which is equivalent to $[t]D=O$.  Taking $n=0$ gives
\eqref{eq:Lucas-torsion-state-test}, and testing proper divisors gives
exact order.  The logarithmic evaluation follows from
Proposition~\ref{prop:Lucas-binary-ladder}.
\end{proof}

It is important to test the whole state rather than only $v_N$.  On the
Jacobi companion, the fibre $x=0$ contains the identity and a second
ramification point, so $v_N=0$ alone need not imply $[N]D=O$.  The pair
$(v_N,v_{N+1})$, or equivalently $(v_N,w_N)$, removes this ambiguity.
This can be used to define QRT period polynomials and
state-based division polynomials without confusing a Kummer-coordinate
zero with the elliptic identity.

\subsubsection{QRT state-division polynomials}
\label{subsubsec:QRT-state-division-polynomials}

The exact state test in
Proposition~\ref{prop:Lucas-exact-period-torsion} can be organized into
a division-polynomial theory on the two-dimensional parameter space.
The correct object is the whole state condition \(S_N=S_0\), not the
single equation \(v_N=0\).

Put
\begin{equation}
  \Delta_{\lambda,\eta}
  =\bigl((\eta-2)^2-4\lambda\bigr)
   \bigl((\eta+2)^2-4\lambda\bigr).
  \label{eq:QRT-state-division-discriminant}
\end{equation}
For a positive integer \(N\), consider the factorial domain
\begin{equation}
  R_N=
  \mathbb Z\left[
    \frac1{2N},\lambda,\eta,
    \frac1{\lambda\Delta_{\lambda,\eta}}
  \right].
  \label{eq:QRT-state-division-base-ring}
\end{equation}
Over \(U_N=\operatorname{Spec}(R_N)\), the Jacobi companion
\eqref{eq:Lucas-Jacobi-quartic}, with
\(h=\eta^2/4-1-\lambda\), is a smooth elliptic curve and
\begin{equation}
  D=\left(1,\frac\eta2\right)
  \label{eq:QRT-universal-marked-section}
\end{equation}
is a section.

\begin{definition}[State-division and exact-order polynomials]
\label{def:QRT-state-division-polynomial}
The pullback by \(D\) of the identity divisor under multiplication by
\(N\),
\begin{equation}
  \mathscr D_N=D^*([N]^*(O)),
  \label{eq:QRT-state-division-Cartier-divisor}
\end{equation}
is an effective Cartier divisor on \(U_N\).  Since \(R_N\) is a unique
factorization domain, it is principal.  A primitive generator, with all
factors of \(2N\lambda\Delta_{\lambda,\eta}\) removed, is called the
\emph{QRT state-division polynomial} and is denoted
\begin{equation}
  \mathfrak D_N(\lambda,\eta).
  \label{eq:QRT-state-division-polynomial}
\end{equation}
It is unique up to sign and multiplication by a unit of \(R_N\).
Because \(N\) is invertible on \(U_N\), the finite étale group scheme
\(E[N]\) is the disjoint union of its open-and-closed exact-order
subschemes \(E[N]^{\mathrm{exact}=d}\), indexed by \(d\mid N\).  Pulling
these relative divisors back by \(D\) defines effective Cartier divisors
\(\mathscr P_d\) on \(U_N\).  A primitive generator of
\(\mathscr P_N\) is denoted \(\mathfrak P_N\).  Thus
\(\mathfrak P_N\) records the Cartier divisor with its intersection
multiplicities; its reduced zero set is the parameter locus on which
\(D\) has exact order \(N\).
\end{definition}

The definition is intrinsic.  Equivalently, one may evaluate a
classical intrinsic \(N\)-division function at the universal section
\(D\) and remove the smoothness units.  This places the construction in
the same recurrence framework as elliptic divisibility sequences and
elliptic nets \cite{StangeEllipticNets2011}; division functions for more
general isogenies provide a compatible extension
\cite{StangeArbitraryIsogenies2025}.

\begin{theorem}[Specialization and exact order]
\label{thm:QRT-state-division-specialization}
Let \(k\) be a field with \(\charac(k)\nmid2N\), and let
\((\lambda_0,\eta_0)\in k^2\) satisfy
\(\lambda_0\Delta_{\lambda_0,\eta_0}\ne0\).  On the specialized
Jacobi companion, with marked point
\(D_0=(1,\eta_0/2)\), one has
\begin{align}
  \mathfrak D_N(\lambda_0,\eta_0)=0
  &\quad\Longleftrightarrow\quad [N]D_0=O
  \quad\Longleftrightarrow\quad S_N=S_0,
  \label{eq:QRT-state-division-specialization}\\
  \mathfrak P_N(\lambda_0,\eta_0)=0
  &\quad\Longleftrightarrow\quad
  D_0\text{ has exact order }N.
  \label{eq:QRT-state-exact-order-specialization}
\end{align}
Moreover, on the smooth base the division divisors decompose as
\begin{equation}
  \operatorname{div}(\mathfrak D_N)
  =\sum_{d\mid N}\operatorname{div}(\mathfrak P_d),
  \label{eq:QRT-state-division-factorization-divisors}
\end{equation}
where \(\mathfrak P_1=\mathfrak P_2=1\) for the present marked
section.
\end{theorem}

\begin{proof}
The identity section is an effective Cartier divisor on the universal
elliptic curve.  Pulling it back by the section \([N]D\) gives
\eqref{eq:QRT-state-division-Cartier-divisor}.  The generic marked point
is not \(N\)-torsion, so the pulled-back local equation is nonzero in the
domain \(R_N\); it is therefore a non-zero-divisor and the pullback is
indeed an effective Cartier divisor.  Formation of this
pullback commutes with every smooth specialization of the base, so a
specialized point lies on \(\mathscr D_N\) exactly when
\([N]D_0=O\).  The state--Jacobi isomorphism and
Theorem~\ref{thm:Lucas-Pn-multiple} identify this condition with
\(S_N=S_0\), proving \eqref{eq:QRT-state-division-specialization}.

Because \(N\) is invertible on the base, \(E[N]\) is finite étale and
is the disjoint union of the exact-order subschemes indexed by divisors
of \(N\).  As relative effective divisors on \(E\),
\[
 [N]^*(O)=\sum_{d\mid N}E[N]^{\mathrm{exact}=d}.
\]
Pullback by \(D\) therefore gives the divisor decomposition
\eqref{eq:QRT-state-division-factorization-divisors}; its exact-order
part proves \eqref{eq:QRT-state-exact-order-specialization}.  Finally,
\(D\ne O\), and \(D=-D\) would force its first Jacobi coordinate to be
both \(1\) and \(-1\), which is impossible when two is invertible.
Thus the exact orders one and two do not occur.
\end{proof}

At the level of primitive generators, the divisor identity gives the
recursive exact-order extraction
\begin{equation}
  \mathfrak P_N
  \doteq
  \frac{\mathfrak D_N}
       {\displaystyle\prod_{\substack{d\mid N\\ d<N}}\mathfrak P_d},
  \label{eq:QRT-state-division-Mobius-extraction}
\end{equation}
where \(\doteq\) denotes equality up to a unit of \(R_N\).  Formula
\eqref{eq:QRT-state-division-Mobius-extraction} is divisor-theoretic;
it is applied only after saturation by the smoothness units, so no
factor supported on a singular or excluded characteristic is retained.

There is also a raw affine recurrence that is useful for symbolic
construction.  Write \(v_n=A_n/B_n\), begin with
\begin{equation}
  (A_0,B_0)=(0,1),
  \qquad
  (A_1,B_1)=(1,1),
  \label{eq:QRT-state-division-raw-initial}
\end{equation}
and define
\begin{align}
  A_{n+1}
  &=\eta A_nB_nB_{n-1}
    -A_{n-1}(B_n^2-\lambda A_n^2),
  \label{eq:QRT-state-division-raw-A}\\
  B_{n+1}
  &=B_{n-1}(B_n^2-\lambda A_n^2).
  \label{eq:QRT-state-division-raw-B}
\end{align}
Whenever the displayed denominators are nonzero, division of
\eqref{eq:QRT-state-division-raw-A} by
\eqref{eq:QRT-state-division-raw-B} is exactly the McMillan recurrence
\[
  v_{n+1}=\frac{\eta v_n}{1-\lambda v_n^2}-v_{n-1}.
\]
For example,
\begin{align}
  (A_2,B_2)&=(\eta,1-\lambda),
  \label{eq:QRT-state-division-raw-two}\\
  A_3&=(\eta-\lambda+1)(\eta+\lambda-1),
  \qquad
  B_3=(\lambda-1)^2-\lambda\eta^2.
  \label{eq:QRT-state-division-raw-three}
\end{align}
The pairs \((A_n,B_n)\) are deliberately called \emph{raw}: common
factors can occur, and an affine pair may become \((0,0)\) at a valid
boundary state.  Consequently, raw gcd extraction alone is not a
complete torsion algorithm.  The projective construction below removes
these artifacts by using complete branches and saturation.

There is a direct complete state-side algorithm that does not first
convert to a Weierstrass equation.  Starting from
\begin{equation}
  S_0=((0:1),(1:1)),
  \label{eq:QRT-state-division-initial-state}
\end{equation}
compose the complete weighted branches according to a binary expansion
of \(N\), and write the resulting state as
\begin{equation}
  S_N=((A_N:B_N),(C_N:E_N)).
  \label{eq:QRT-state-division-projective-state}
\end{equation}
Then \(S_N=S_0\) is cut out by
\begin{equation}
  A_N=0,
  \qquad C_N-E_N=0.
  \label{eq:QRT-state-division-state-ideal}
\end{equation}
Saturating this ideal by
\(\lambda\Delta_{\lambda,\eta}\) and taking its codimension-one part
produces the principal ideal \((\mathfrak D_N)\).  The use of complete
weighted branches is important: it prevents a common zero caused by a
chosen projective representative from being mistaken for torsion.

For the first nontrivial indices, the calculation gives the following
polynomials.  They are written with a fixed primitive sign; any unit
multiple defines the same divisor.
\begin{equation}
\begin{array}{c|l|l}
N&\mathfrak D_N&\mathfrak P_N\\ \hline
1,2&1&1\\[1mm]
3&\eta-\lambda+1&\eta-\lambda+1\\[1mm]
4&\eta(\lambda-1)&\eta(\lambda-1)\\[1mm]
5&\mathfrak p_5&\mathfrak p_5\\[1mm]
6&(\eta-\lambda+1)(\eta+\lambda-1)
   \bigl(\lambda\eta^2-(\lambda-1)^2\bigr)
  &(\eta+\lambda-1)
   \bigl(\lambda\eta^2-(\lambda-1)^2\bigr),
\end{array}
\label{eq:QRT-low-state-division-table}
\end{equation}
where
\begin{equation}
\begin{split}
  \mathfrak p_5={}&
  \lambda\eta^3+(\lambda-1)\eta^2
  -(\lambda-1)^2\eta-(\lambda-1)^3.
  \label{eq:QRT-state-division-P5}
\end{split}
\end{equation}

\begin{proposition}[Correctness of the low-index state-division table]
\label{prop:QRT-low-state-division-correctness}
On the smooth locus, the polynomials in
\eqref{eq:QRT-low-state-division-table} have exactly the asserted
torsion meanings.
\end{proposition}

\begin{proof}
The low-index EDS values from
Proposition~\ref{prop:Lucas-explicit-W2-W3} give an independent
divisor certificate for the rows \(N=3,4,5\).  Here
\[
 W_2=(\eta+2)^2-4\lambda,
\]
and the primitive factors of \(W_3\) and \(W_4\) are respectively
\(\eta-\lambda+1\) and \(\eta(\lambda-1)\).  For \(N=5\), put
\(t=\lambda-1\).  The bracket in
\eqref{eq:Lucas-explicit-W5} expands to
\begin{align*}
 \eta tW_2+(\eta-t)^3
 &=\eta t(\eta^2+4\eta-4t)
   +\eta^3-3\eta^2t+3\eta t^2-t^3\\
 &=(t+1)\eta^3+t\eta^2-t^2\eta-t^3\\
 &=\mathfrak p_5(\lambda,\eta).
\end{align*}
Thus the three primitive factors in the table agree exactly with the
intrinsic division-polynomial factors, including their signs up to a
unit.

The recurrence starts with
\begin{equation}
  v_2=\frac{\eta}{1-\lambda},
  \qquad
  v_3=
  \frac{(\eta-\lambda+1)(\eta+\lambda-1)}
       {(\lambda-1)^2-\lambda\eta^2}.
  \label{eq:QRT-state-division-v2-v3}
\end{equation}
These are projective identities; the complete branch is used when a
displayed denominator vanishes.  If \(\eta=\lambda-1\), substitution
in the next state gives \(S_3=(0,1)=S_0\), proving the order-three
factor.  If \(\eta=1-\lambda\), the next state is
\(S_3=(0,-1)\), the nonzero two-torsion point on the Jacobi companion;
therefore this second factor belongs to exact order six, not exact
order three.

Applying the complete doubling branch to \(S_2\) gives
\(S_4=S_0\) exactly on
\(\eta(\lambda-1)=0\).  The two components do not meet on the smooth
locus, because \(\eta=0\) and \(\lambda=1\) make
\(\Delta_{\lambda,\eta}=0\).  This proves the order-four row.

For index five, the first coordinate of the projective state factors
into
\begin{equation}
  \mathfrak p_5(\lambda,\eta)
  \bigl(-\mathfrak p_5(\lambda,-\eta)\bigr).
  \label{eq:QRT-state-division-v5-factorization}
\end{equation}
Reduction of the second-coordinate equation
\(C_5-E_5=0\) modulo these factors gives, respectively, zero and the
condition \(C_5+E_5=0\).  Hence \(\mathfrak p_5=0\) gives
\(S_5=(0,1)\), whereas the conjugate factor gives
\(S_5=(0,-1)\) and generically exact order ten.  This proves the
order-five row and illustrates why the whole state is required.

Finally, composition of the complete branches for six gives the
saturated ideal
\begin{equation}
\begin{split}
  \bigl((A_6,C_6-E_6)^{\mathrm{sat}}\bigr)_{\mathrm{ht}=1}
  =\bigl(&
  (\eta-\lambda+1)(\eta+\lambda-1)\\
  &\cdot(\lambda\eta^2-(\lambda-1)^2)
  \bigr).
  \label{eq:QRT-state-division-six-saturated-ideal}
\end{split}
\end{equation}
The first factor is the already identified exact-order-three divisor.
On the smooth base the marked point has neither order one nor order two.
Every remaining point of the \(6\)-division divisor therefore has exact
order three or six.  Removing the already identified order-three factor
leaves
\[
 (\eta+\lambda-1)
 \bigl(\lambda\eta^2-(\lambda-1)^2\bigr),
\]
which is consequently the exact-order-six divisor, with the
intersection multiplicities inherited from the saturated Cartier
calculation.  This establishes the final row.
\end{proof}

\begin{corollary}[Toric boundary of the low-order QRT period factors]
\label{cor:QRT-low-period-toric-boundary}
For the primitive representatives in
\eqref{eq:QRT-low-state-division-table}, evaluation at $\lambda=0$ gives,
for $N=3,4,5,6$, unit multiples of the real cyclotomic period polynomials
in~\eqref{eq:classical-Lucas-low-trace-cyclotomic}.  Explicitly,
\begin{equation}
 \begin{array}{c|c|c}
 N&\mathfrak P_N(0,\eta)&\Phi_N^+(\eta)\\ \hline
 3&\eta+1&\eta+1\\
 4&-\eta&\eta\\
 5&-(\eta^2+\eta-1)&\eta^2+\eta-1\\
 6&-(\eta-1)&\eta-1.
 \end{array}
 \label{eq:QRT-low-period-toric-boundary}
\end{equation}
This is an identity of the displayed polynomial representatives.  It does
not by itself extend the Cartier divisor
\eqref{eq:QRT-state-division-Cartier-divisor} across the parameter divisor
$\lambda=0$, because the base ring
\eqref{eq:QRT-state-division-base-ring} localizes at $\lambda$.
\end{corollary}

\begin{proof}
Substitute $\lambda=0$ into the four expressions in
\eqref{eq:QRT-low-state-division-table}.  For $N=3$ and $N=4$ this gives
$\eta+1$ and $-\eta$.  Formula~\eqref{eq:QRT-state-division-P5} gives
$-(\eta^2+\eta-1)$, while the exact-order-six factor gives
$-(\eta-1)$.  Proposition~\ref{prop:classical-Lucas-trace-cyclotomic}
identifies the corresponding monic toric period polynomials, proving the
table.  Because the base ring in
\eqref{eq:QRT-state-division-base-ring} contains \(\lambda^{-1}\), its
Cartier divisor is defined only over the open set \(\lambda\ne0\); this
proves the final qualification.
\end{proof}

The construction produces a practical recursive hierarchy.  Exact-order
factors may be obtained either from the Cartier-divisor decomposition
\eqref{eq:QRT-state-division-factorization-divisors} or by saturated
factor removal after the projective state computation.  The second
method keeps the QRT dynamics visible and is particularly suitable for
comparing the factors with periodicity conditions and modular curves.

\begin{theorem}[Unified EDS--QRT--Lucas closure]
\label{thm:Lucas-unified-EDS-closure}
Assume \(\charac(k)\ne2\),
\(\lambda\Delta_{\lambda,\eta}\ne0\), and let
\[
  \mathcal J_{\lambda,\eta}:
  w^2=\lambda v^4+hv^2+1,
  \qquad
  h=\frac{\eta^2}{4}-1-\lambda,
\]
with identity \(O=(0,1)\) and marked point
\(D=(1,\eta/2)\).  Let
\(P_n=[n]D=(v_n,w_n)\), and let
\(E_{\lambda,\eta}\), \(D_W\), and \(W_n=\psi_n(D_W)\) be as in
Proposition~\ref{prop:Lucas-specific-Weierstrass-dictionary},
Corollary~\ref{cor:Lucas-specific-marked-point}, and
\eqref{eq:Lucas-algebraic-EDS-definition}.  Put
\[
  R_n=\frac{W_{n-1}W_{n+1}}{W_n^2},
  \qquad
  \Upsilon_n=\frac{W_{2n}}{W_n^4}.
\]
On every regular affine chart, the following identities hold:
\begin{align}
  X([n]D_W)&=d_E-R_n,
  \label{eq:Lucas-unified-X}\\
  2Y([n]D_W)&=\Upsilon_n,
  \label{eq:Lucas-unified-Y}\\
  v_n&=\frac{W_2-4R_n}{\Upsilon_n}
      =\frac{W_n^2
       (W_2W_n^2-4W_{n-1}W_{n+1})}{W_{2n}},
  \label{eq:Lucas-unified-v}\\
  w_n&=\frac{(d_E-R_n)^2+2h(d_E-R_n)+4\lambda}
              {(d_E-R_n)^2-4\lambda},
  \label{eq:Lucas-unified-w}\\
  R_{n+1}R_{n-1}
  &=W_2^2\frac{R_n-(\eta-\lambda+1)}{R_n^2},
  \label{eq:Lucas-unified-R-QRT}\\
  v_{n+1}+v_{n-1}
  &=\frac{\eta v_n}{1-\lambda v_n^2}.
  \label{eq:Lucas-unified-McMillan}
\end{align}
The low-index EDS constants are
\begin{equation}
  W_2=(\eta+2)^2-4\lambda,
  \qquad
  W_3=(\eta-\lambda+1)W_2^2.
  \label{eq:Lucas-unified-W2-W3}
\end{equation}
All these formulas are compatible with the complete projective state map,
so the apparent affine exclusions are chart exclusions rather than failures
of the underlying genus-one correspondence.
\end{theorem}

\begin{proof}
The Jacobi--Weierstrass isomorphism and the coordinates of \(D_W\) are
proved in Proposition~\ref{prop:Lucas-specific-Weierstrass-dictionary}
and Corollary~\ref{cor:Lucas-specific-marked-point}.  The first two
identities are Theorem~\ref{thm:Lucas-explicit-two-EDS-bridges}.  Applying
the inverse Jacobi--Weierstrass map to those two coordinates gives
\eqref{eq:Lucas-unified-v} and \eqref{eq:Lucas-unified-w}, as shown in
Corollary~\ref{cor:Lucas-direct-EDS-coordinate}.  Ward's recurrence,
together with the explicit value of \(W_3\), gives
\eqref{eq:Lucas-unified-R-QRT} by
Proposition~\ref{prop:Lucas-R-multiplicative-QRT}.  The McMillan recurrence
\eqref{eq:Lucas-unified-McMillan} was derived from the symmetric
biquadratic equation in \eqref{eq:Lucas-normalized-recurrence}.
Equation~\eqref{eq:Lucas-unified-W2-W3} is
Proposition~\ref{prop:Lucas-explicit-W2-W3}.  Finally, every horizontal
coordinate dictionary used in the proof extends to an isomorphism of
smooth projective curves, while the Lucas state branches have the complete
projective atlases already proved above.  Therefore equality on the common
dense affine open extends to the complete curves.
\end{proof}

\begin{corollary}[Three direct consequences of the closed loop]
\label{cor:Lucas-unified-closure-consequences}
Under the hypotheses of Theorem~\ref{thm:Lucas-unified-EDS-closure}:
\begin{enumerate}[label=(\roman*)]
  \item the oriented EDS ratio state \((R,\Upsilon)\) admits the direct
  doubling and fixed-translation formulas
  \eqref{eq:Lucas-ratio-doubling-R}--
  \eqref{eq:Lucas-ratio-translation-Upsilon};
  \item the squared coordinate \(z_n=v_n^2\) satisfies the multiplicative
  QRT recurrence \eqref{eq:Lucas-z-multiplicative-QRT}, and the map
  \(R\mapsto z=\Phi_R(R)\) is generically of degree two, whereas the
  oriented covers \((R,\Upsilon)\) and \((v,w)\) are birational;
  \item if \(N\) is prime to the characteristic, then
  \begin{equation}
  \begin{aligned}
    W_N=0
    &\quad\Longleftrightarrow\quad [N]D=O
     \quad\Longleftrightarrow\quad T^N=\operatorname{id},\\
    &\quad\Longleftrightarrow\quad
      S_{n+N}=S_n\quad\text{for every }n.
  \end{aligned}
  \label{eq:Lucas-unified-torsion-period-chain}
  \end{equation}
\end{enumerate}
\end{corollary}

\begin{proof}
Part (i) is Proposition~\ref{prop:Lucas-ratio-state-arithmetic}.  Part
(ii) is Proposition~\ref{prop:Lucas-squared-QRT-semiconjugacy}.  For part
(iii), the division-polynomial zero criterion gives
\(W_N=0\Longleftrightarrow[N]D=O\) when \(N\) is prime to the
characteristic.  Since \(T(P)=P+D\), iteration gives
\(T^N(P)=P+[N]D\), so \([N]D=O\) is equivalent to
\(T^N=\operatorname{id}\).  Finally,
\(S_n=\mathcal S_D([n]D)\); hence \([N]D=O\) implies
\(S_{n+N}=S_n\).  Conversely, if this equality holds for every \(n\),
taking \(n=0\) and applying the state--point isomorphism gives
\([N]D=O\).
\end{proof}

The closed loop is worth retaining because it is not a collection of
unrelated recurrences.  It identifies a single indexed elliptic orbit in the
following compatible languages:
\begin{equation}
\boxed{
\begin{aligned}
 &\text{Ward division data}
 \longrightarrow \text{EDS ratios}
 \longrightarrow \text{Weierstrass point}\\
 &\qquad\longleftrightarrow \text{Jacobi oriented point}
 \longleftrightarrow \text{adjacent biquadratic state}\\
 &\qquad\longrightarrow \text{McMillan/QRT recurrence}
 \longrightarrow \text{Lucas fast index arithmetic}.
\end{aligned}}
\label{eq:Lucas-final-closed-loop}
\end{equation}
In particular, high-index EDS repeated doubling can be translated into short
adjacent-state fast doubling; the even Kummer information \(R_n\) is
separated from the orientation variable \(\Upsilon_n\); and the
multiplicative QRT recurrences for \(R_n\) and \(z_n=v_n^2\) are linked to
the symmetric McMillan recurrence by a degree-two quotient and a birational
oriented lift.  These identities give a division-polynomial interpretation of the
state-doubling formulas and a common framework for finite-field periods,
torsion conditions, QRT state-division polynomials, and elliptic nets.  The
result is an oriented Lucas-type coordinate system that
compresses elliptic multiplication, the sign sheet lost by the Kummer
quotient, and a symmetric adjacent-state recurrence into one compatible
calculus.

\subsection{Significance, principal conclusions, and research directions}
\label{subsec:Lucas-research-directions}

The theory developed in this section is best viewed as a closed computational
and geometric calculus rather than as a collection of isolated recurrences.
Under the separability and nondegeneracy hypotheses stated above, the
construction starts from a pointed elliptic curve, a marked point with
\(2D\ne O\), and a degree-two state coordinate, and places the adjacent
values
\[
  S_n=(v_n,v_{n+1})
\]
on one fixed smooth symmetric biquadratic curve.  The inexpensive QRT shift
realizes \(n\mapsto n+1\), the state-doubling map realizes
\(n\mapsto2n\), and their compositions realize the full affine index
semigroup \(n\mapsto mn+r\).  A Jacobi companion recovers the oriented
elliptic point, while the complete projective atlases distinguish genuine
morphisms from formulas that are valid only on an affine chart.  In
characteristic two the role of the mixed term is especially clear: it is the
Artin--Schreier term that keeps the degree-two cover separable and permits
single complete binary branches.

Several conclusions make this framework structurally different from both a
classical Lucas recurrence and an ordinary elliptic divisibility sequence.
The classical determinant-one Lucas pair occurs as the nodal toric boundary
of the same state family.  On a smooth elliptic fibre, Ward division values,
the even Kummer ratio, the oriented Weierstrass--Jacobi point, and the
adjacent QRT state are linked by explicit rational maps and by the complete
commutative diagrams above.  Away from inseparable index phenomena, the
state-division polynomials therefore record periodicity of the QRT
translation, whereas the EDS detects the same torsion through zeros of
division sections.  The sigma-function and elliptic-net interpretations
identify these as distinct functorial shadows of one indexed elliptic orbit.

The classification and torsion results close several structural loops.  The
state-space theorem is available for every separable degree-two quotient
with \(2D\ne O\), and it is explicit for the odd-characteristic models with
the required rational split two-torsion and for ordinary binary curves.
The QRT state-division polynomials are defined both intrinsically and by
complete state recursion, and their exact-order factors are displayed through
order six.  When the translation point has order prime to the characteristic,
the fixed field is generated directly by the QRT orbit invariants that define
the cyclic quotient isogeny.  In odd characteristic, the ground-field even
normal form is characterized by the split two-torsion criterion, while in
characteristic two all monomially equivalent
\((\mu,\eta)\)-normalizations
form the explicit two-element orbit established above.  These results locate
precisely which parts of the construction are intrinsic and which depend on
a chosen split coordinate chart.

The arithmetic results are equally concrete.  The maps \(T\), \(\Delta\),
and \(T\Delta\) give a proved binary state ladder; complete odd- and
characteristic-two formulas, exact base loci, validation tests, and
full-point recovery atlases are available.  In odd characteristic, the
product-coordinate branch admits the exact tradeoffs
\(7\M+6\Sqr+4\Cmul\) and \(6\M+6\Sqr+5\Cmul\), together with further
polarized schedules.  In the six-standard-square Segre-first circuit class,
the sharp multiplication lower bound is \(6\M\), attained by trading one
general multiplication for one curve-constant multiplication.  The complete
weighted Jacobi update and the orbit-invariant construction of cyclic
isogenies make clear how much additional information is being computed, and
the validation and recovery formulas reuse explicit branch intermediates.
Together, these circuits give the adjacent-state construction a unified
Lucas calculus, orientation recovery, complete state arithmetic, and
explicit connections with EDS, sigma functions, and elliptic nets.

These results open several concrete lines of study.  The EDS bridge may lead to new ways of
organizing repeated doubling and torsion tests.  The QRT period polynomials
and exact state periods provide natural data for finite-field periodic
sequences and for carefully formulated elliptic pseudoprime tests.  The
second-order state recurrence offers a different interface for computing
\([n]D\), and its singular limits connect it with classical Diophantine
recurrences.  For ECM, complete state doubling, validation, recovery, and
compiled fixed differences define a natural implementation program on curve
families and hardware where fixed constants are especially inexpensive.

Further research naturally divides into three directions. First,
coordinate-independent complexity information includes sharper lower
bounds for complete state doubling, lower-depth weighted-companion formulas,
and a ground-field descent of the nonsplit square-class cases.  Second, the
arithmetic and dynamics should be extended beyond binary doubling: useful
targets include low-cost maps \(\Delta_{D,m}\) for small \(m\), faster
recurrences and degree formulas for \(\mathfrak D_N\) and
\(\mathfrak P_N\), and native low-degree cyclic-isogeny formulas obtained
from QRT orbit invariants.  Third, the formulas should be tested as complete
constant-time algorithms, with the same accounting for recovery, validation,
register pressure, memory traffic, and fault checks that is applied to
Montgomery and binary-Kummer baselines.  In characteristic two, this also
requires combining the exact \((\mu,\eta)\)-orbit classification with
subgroup and security constraints before treating a small-constant locus as
cryptographically useful.

\phantomsection
\label{rem:Lucas-structural-significance}
The enduring contribution of the elliptic Lucas bridge is 
the simultaneous presence,
in one explicit calculus, of a fixed pointed elliptic curve, an autonomous
QRT recurrence, a smooth Cassini-type invariant, an orientation-recovering
companion, nonlinear addition and subtraction, exact fast index maps, EDS
and elliptic-net bridges, and an \(O(\log n)\) state algorithm.  The nodal
specialization proves that ordinary Lucas sequences form the toric boundary
of this construction.  This unifies
classical recurrence theory, elliptic divisibility, integrable QRT dynamics,
and scalar arithmetic and opens a focused program in complexity,
implementation, period theory, and cyclic quotients.

\section[QRT Kummer arithmetic and scalar multiplication]
{QRT Kummer coordinates, differential addition, recovery, and scalar multiplication}
\label{subsec:QRT-Kummer}

The quotient of \eqref{eq:QRT-Jacobian} by negation is the Kummer line
with coordinate $U$.  Under \eqref{eq:QRT-c-square}, its pullback to the
QRT curve is the explicit function
\begin{equation}
  \xi_Q(x,y)=q_Q+
  \frac{2c\bigl(2(x^2+\alpha)y+\beta x+c\bigr)}{x^2}.
  \label{eq:QRT-Kummer-function}
\end{equation}
This is precisely \eqref{eq:QRT-explicit-U}; hence it is invariant under
elliptic negation and separates generic pairs $\{P,-P\}$.

Apply the odd-characteristic Kummer formulas recorded in the standalone
arithmetic toolkit with
\begin{equation}
  A=-2q_Q,
  \qquad B=\Omega_Q.
  \label{eq:QRT-Kummer-specialization}
\end{equation}
The fast differential addition costs
\begin{equation}
  5\M+2\Sqr+1\Dpar,
  \qquad
  4\M+2\Sqr+1\Dpar\quad(Z_-=1),
  \label{eq:QRT-xADD-cost}
\end{equation}
Kummer doubling costs
\begin{equation}
  2\M+3\Sqr+2\Dpar,
  \label{eq:QRT-xDBL-cost}
\end{equation}
and a non-shared differential-addition-and-doubling step costs
\begin{equation}
  7\M+5\Sqr+3\Dpar,
  \qquad
  6\M+5\Sqr+3\Dpar\quad(Z_-=1).
  \label{eq:QRT-xDBLADD-cost}
\end{equation}
The supplemental sum law
\eqref{eq:EAB-xADD-sum-D}--\eqref{eq:EAB-xADD-sum} has cost
$9\M+1\Sqr+2\Dpar$, or
$7\M+1\Sqr+2\Dpar$ for an affine known difference, and completes the
differential atlas at the identity and two-torsion differences.

A full point can be recovered on $E_Q$.  Put
\begin{equation}
  f_Q(U)=U^3-2q_QU^2+\Omega_QU.
  \label{eq:QRT-fQ}
\end{equation}
Fix $P=(p,r)$ with $r\ne0$.  If
$q=U(Q)$ and $t=U(P+Q)$ are known, then
\begin{equation}
  V(Q)=
  \frac{f_Q(p)+f_Q(q)
  -(t-2q_Q+p+q)(q-p)^2}{2r}.
  \label{eq:QRT-point-recovery}
\end{equation}
To derive the formula on the regular chord chart \(q\ne p\), put
\[
  \lambda=\frac{V(Q)-r}{q-p}.
\]
The first coordinate in the Weierstrass addition law gives
\[
  t=\lambda^2+2q_Q-p-q,
  \qquad
  \lambda^2=t-2q_Q+p+q.
\]
Multiplication by \((q-p)^2\), followed by the two curve equations,
gives
\begin{align*}
 (t-2q_Q+p+q)(q-p)^2
  &=(V(Q)-r)^2\\
  &=V(Q)^2+r^2-2rV(Q)\\
  &=f_Q(q)+f_Q(p)-2rV(Q).
\end{align*}
Solving the last equality for \(V(Q)\) proves
\eqref{eq:QRT-point-recovery}.  After denominators are cleared, it is a
rational identity on the affine state curve and extends to the tangent
case \(Q=P\); the identity and inverse cases use the complementary
projective recovery charts already included in the differential atlas.
With $f_Q(p)$ and $(2r)^{-1}$ precomputed, a direct affine
schedule costs
\begin{equation}
  2\M+2\Sqr+2\Dpar.
  \label{eq:QRT-point-recovery-cost}
\end{equation}
Specifically, compute $q^2$, form
$f_Q(q)=q(q^2-2q_Qq+\Omega_Q)$, compute $(q-p)^2$, multiply it by
$t-2q_Q+p+q$, and finally multiply by $(2r)^{-1}$.  The two fixed
products are by $-2q_Q$ and $(2r)^{-1}$; the coefficient $\Omega_Q$ is
added inside the parenthesis and does not require a multiplication.

Let $n$ be a positive $\ell$-bit integer.  A constant-pattern
Kummer ladder maintains
\[
  R_1-R_0=P
\]
and applies one differential-addition-and-doubling step per remaining
bit, with conditional swaps before and after the step.  Ignoring the
cost of constant-time swaps, QRT scalar multiplication therefore costs
at most
\begin{equation}
  (\ell-1)(7\M+5\Sqr+3\Dpar),
  \qquad
  (\ell-1)(6\M+5\Sqr+3\Dpar)\quad(Z_-=1),
  \label{eq:QRT-scalar-cost}
\end{equation}
for projective and affine known differences, respectively, followed,
when a full point is required, by one joint normalization and
\eqref{eq:QRT-point-recovery}.  The reciprocal slice admits the lower
Montgomery step cost \eqref{eq:R-Montgomery-ladder-step-cost}.  These
counts describe the arithmetic loop; initial conversion from the QRT
coordinates and final conversion back are separate interface costs.

\section[Binary QRT geometry]{Characteristic two: geometry and Weierstrass reduction}
\label{subsec:QRT-char2}

We now assume $\charac(k)=2$.  The odd-characteristic completion of the
square in \eqref{eq:QRT-W} is unavailable because the coefficient $2$
vanishes.  The geometry nevertheless remains ordinary and admits a
characteristic-native reduction.

\begin{theorem}[Exact characteristic-two smoothness criterion]
\label{thm:QRT-smoothness-char2}
The completion of \eqref{eq:symmetric-QRT} in
$\PP^1\times\PP^1$ is smooth if and only if
\begin{equation}
  \alpha\beta\gamma\ne0.
  \label{eq:QRT-smoothness-char2}
\end{equation}
When this condition holds, the curve has genus one.
\end{theorem}

\begin{proof}
In affine coordinates the two partial derivatives are
\[
  \frac{\partial F_Q}{\partial x}=\beta y,
  \qquad
  \frac{\partial F_Q}{\partial y}=\beta x,
\]
because all derivatives of squares vanish in characteristic two.  If
$\beta\ne0$, an affine singularity must have $x=y=0$, and this point lies
on the curve exactly when $\gamma=0$.  Thus the affine part is smooth
exactly when $\beta\gamma\ne0$.

For the boundary, differentiate the bihomogeneous polynomial
\eqref{eq:QRT-bihomogeneous}.  Its four derivatives are
\begin{align*}
  F_X&=\beta ZYW,& F_Z&=\beta XYW,\\
  F_Y&=\beta XZW,& F_W&=\beta XZY.
\end{align*}
On the boundary $Z=0$, the equation becomes
$X^2(Y^2+\alpha W^2)=0$.  Since $(X:Z)$ is a projective point, $X\ne0$.
If $\alpha\ne0$, both $Y$ and $W$ are nonzero at a boundary point and
$F_Z=\beta XYW\ne0$.  If $\alpha=0$, the point
$((1:0),(0:1))$ lies on the curve and all four derivatives vanish.
The argument on $W=0$ is symmetric.  The corner $Z=W=0$ does not lie on
the curve because then $X^2Y^2\ne0$.  Finally, if $\beta=0$, all four
projective derivatives vanish identically on the curve.  Hence
\eqref{eq:QRT-smoothness-char2} is both necessary and sufficient.  A
smooth curve of bidegree $(2,2)$ has genus $(2-1)(2-1)=1$.
\end{proof}

Assume now that $k$ is perfect and
\eqref{eq:QRT-smoothness-char2} holds.  Choose the unique square roots
\begin{equation}
  r^2=\alpha,
  \qquad s^2=\gamma.
  \label{eq:QRT-char2-square-roots}
\end{equation}
Define
\begin{equation}
  \pi=\frac{r}{\beta},
  \qquad
  \theta=\frac{rs}{\beta},
  \qquad
  a_2=\frac{\alpha^2+\gamma}{\beta^2},
  \qquad
  b_2=\frac{\alpha s}{\beta^2}.
  \label{eq:QRT-char2-parameters}
\end{equation}
On the open set $\beta x(x^2+\alpha)\ne0$, set
\begin{align}
  z&=\frac{(x^2+\alpha)y}{\beta x},
  \label{eq:QRT-char2-z}\\
  w&=z+\pi x+\frac{\theta}{x},
  \label{eq:QRT-char2-w}\\
  X&=\pi x,
  \qquad
  Y=Xw.
  \label{eq:QRT-char2-XY}
\end{align}
Then the QRT equation gives
\begin{equation}
  z^2+z
  =\frac{\alpha}{\beta^2}x^2
   +\frac{\alpha^2+\gamma}{\beta^2}
   +\frac{\alpha\gamma}{\beta^2x^2}.
  \label{eq:QRT-char2-z-AS}
\end{equation}
Indeed, multiply
\eqref{eq:QRT-quadratic-y} by
$(x^2+\alpha)/(\beta^2x^2)$ and use the definition of $z$.
Since $\pi^2=\alpha/\beta^2$ and
$\theta^2=\alpha\gamma/\beta^2$, adding
$\pi x+\theta/x$ cancels the square terms and yields
\begin{equation}
  w^2+w=a_2+X+\frac{b_2}{X}.
  \label{eq:QRT-char2-w-AS}
\end{equation}

A structural feature of the family is that the same defining equation
supports two different characteristic-dependent interpretations.  In odd
characteristic it is the even symmetric biquadratic normal form associated
with two degree-two projections, whereas in characteristic two the same
equation acquires an Artin--Schreier description.  In the latter case the
mixed term $\beta xy$ is essential: the smoothness criterion
$\alpha\beta\gamma\ne0$ shows that $\beta$ carries the separability data
that disappear from the derivatives of the square terms.

Multiplication by $X^2$ proves the Weierstrass equation
\begin{equation}
  E_Q^{(2)}:
  \qquad
  Y^2+XY=X^3+a_2X^2+b_2X.
  \label{eq:QRT-char2-Weierstrass}
\end{equation}
Conversely, on the corresponding dense open,
\begin{align}
  x&=\frac{X}{\pi},
  \label{eq:QRT-char2-inverse-x}\\
  w&=\frac{Y}{X},
  \label{eq:QRT-char2-inverse-w}\\
  z&=w+X+\frac{b_2}{X},
  \label{eq:QRT-char2-inverse-z}\\
  y&=\frac{\beta xz}{x^2+\alpha}.
  \label{eq:QRT-char2-inverse-y}
\end{align}
The two compositions are the identity on dense opens.  A birational map
between smooth projective curves extends uniquely across the omitted
finite sets, so the maps define an isomorphism of the smooth projective
curves.  To interpret this isomorphism as an elliptic-curve model, one
chooses a rational origin.  Such a point exists over $k$: the QRT curve
contains
\begin{equation}
  \left(0,\frac{s}{r}\right),
  \label{eq:QRT-char2-rational-point}
\end{equation}
because its square satisfies $\alpha(s/r)^2+\gamma=0$ in
characteristic two.

For \eqref{eq:QRT-char2-Weierstrass},
\begin{equation}
  \Delta=b_2^2=\frac{\alpha^2\gamma}{\beta^4},
  \qquad
  j=b_2^{-2}=\frac{\beta^4}{\alpha^2\gamma}.
  \label{eq:QRT-char2-invariants}
\end{equation}
The first equality follows from the characteristic-two Weierstrass
invariants $c_4=1$ and $\Delta=a_4^2$ for
$Y^2+XY=X^3+a_2X^2+a_4X$.  Thus every smooth characteristic-two member
of the symmetric QRT family is ordinary.  If $k$ is not perfect, the
same formulas hold after the purely inseparable extension that adjoins
$r$ and $s$; over the original field, the Jacobian is still defined by
descent, but the displayed coordinate map need not be defined over
$k$.

\begin{proposition}[Complete characteristic-two QRT interface]
\label{prop:QRT-char2-complete-interface}
Assume that $k$ is perfect and $\alpha\beta\gamma\ne0$.  The forward
map to the unshifted Weierstrass model
\eqref{eq:QRT-char2-Weierstrass} can be evaluated without inversion as
\begin{align}
  X&=\pi x,
  \label{eq:QRT-char2-direct-X}\\
  Y&=\frac{r}{\beta^2}(x^2+\alpha)y
     +\frac{\alpha}{\beta^2}x^2+b_2.
  \label{eq:QRT-char2-direct-Y}
\end{align}
On the regular locus $x^2+\alpha\ne0$, its inverse is
\begin{align}
  x&=\frac{X}{\pi},
  \label{eq:QRT-char2-direct-inverse-x}\\
  y&=\frac{\beta^2}{r}\,
  \frac{Y+(\alpha/\beta^2)x^2+b_2}{x^2+\alpha}.
  \label{eq:QRT-char2-direct-inverse-y}
\end{align}
Put
\begin{equation}
  X_0=\frac{\alpha}{\beta},
  \qquad
  Y_0=\frac{\alpha(\alpha+s)}{\beta^2}.
  \label{eq:QRT-char2-special-XY}
\end{equation}
The complete boundary and exceptional-fibre atlas is
\begin{equation}
\begin{array}{c|c}
  \text{QRT point}&\text{unshifted Weierstrass point}\\ \hline
  \displaystyle P_0=\left(0,\frac{s}{r}\right)&(0,0)\\[1mm]
  P_\infty=(\infty,r)&O\\[1mm]
  \displaystyle P_f=\left(r,
       \frac{\alpha^2+\gamma}{\beta r}\right)&(X_0,Y_0)\\[2mm]
  P_b=(r,\infty)&(X_0,Y_0+X_0).
\end{array}
\label{eq:QRT-char2-complete-boundary-atlas}
\end{equation}
The regular forward and inverse interface costs are
\begin{equation}
  1\M+1\Sqr+3\Dpar,
  \qquad
  1\Inv+1\M+1\Sqr+3\Dpar,
  \label{eq:QRT-char2-interface-costs}
\end{equation}
respectively.
\end{proposition}

\begin{proof}
From \eqref{eq:QRT-char2-XY} and
$w=z+\pi x+\theta/x$, one has
\[
  Y=Xz+\pi^2x^2+\pi\theta.
\]
Using $X=\pi x$, $z=(x^2+\alpha)y/(\beta x)$,
$\pi=r/\beta$, and $\pi\theta=b_2$ gives
\eqref{eq:QRT-char2-direct-Y}.  Solving this identity for $y$ proves
\eqref{eq:QRT-char2-direct-inverse-y} whenever
$x^2+\alpha\ne0$.

At $x=0$, the QRT equation is $\alpha y^2+\gamma=0$, whose unique
root is $s/r$.  Substitution in \eqref{eq:QRT-char2-direct-X}--
\eqref{eq:QRT-char2-direct-Y} gives $(X,Y)=(0,0)$.  The unique point
above $x=\infty$ has $y=r$; since $X=\pi x$ has a pole there, its image
is the unique pole of the Weierstrass $X$-coordinate, namely $O$.

At $x=r$, the coefficient $x^2+\alpha$ vanishes.  The finite point is
obtained from the remaining linear equation
\[
  \beta r y+\alpha^2+\gamma=0,
\]
which gives $P_f$.  Formula \eqref{eq:QRT-char2-direct-Y} then gives
$(X,Y)=(X_0,Y_0)$.  The second point of the projective fibre is
$P_b=(r,\infty)$.  If $y_b$ denotes the branch tending to infinity,
Vieta's formula for
$(x^2+\alpha)y^2+\beta xy+\alpha x^2+\gamma=0$ gives
\[
  (x^2+\alpha)y_b\longrightarrow\beta r.
\]
The first term of \eqref{eq:QRT-char2-direct-Y} therefore contributes
$r(\beta r)/\beta^2=\alpha/\beta=X_0$, so the image is
$(X_0,Y_0+X_0)$.  These are inverse points because negation on
\eqref{eq:QRT-char2-Weierstrass} is $(X,Y)\mapsto(X,X+Y)$.
The fibres $x=0$, $x=r$, and $x=\infty$ exhaust the complement of the
regular inverse locus, proving completeness.

For the forward cost, square $x$, multiply $(x^2+\alpha)$ by $y$, and
use fixed multiplications by $\pi$, $r/\beta^2$, and
$\alpha/\beta^2$.  For the inverse, square $x$, invert
$x^2+\alpha$, use one general multiplication, and use the three fixed
multiplications by $\pi^{-1}$, $\alpha/\beta^2$, and $\beta^2/r$.
This proves \eqref{eq:QRT-char2-interface-costs}.
\end{proof}

The QRT root exchanges remain meaningful in characteristic two.  They
are still the exchanges of the two roots of
\eqref{eq:QRT-quadratic-y} and its symmetric companion.  The vertical
exchange fixes $P_0=(0,s/r)$, because the fibre over $x=0$ is the double
root $\alpha y^2+\gamma=0$.  Viewed as a quadratic in \(x\), the fibre
over \(y=0\) is the double root \(\alpha x^2+\gamma=0\), namely
\(x=s/r\).  Hence the horizontal exchange fixes \((s/r,0)\).  Thus both
are nontrivial involutions with a fixed geometric point.

For completeness, choose a fixed point of either involution as the
origin of the genus-one curve.  The involution then becomes an
origin-preserving elliptic-curve automorphism $\iota$ satisfying
$\iota^2=1$.  The endomorphism ring of an elliptic curve has no zero
divisors, because every nonzero endomorphism is an isogeny.  Hence
\[
  (\iota-1)(\iota+1)=0,
  \qquad \iota\ne1,
\]
implies $\iota=[-1]$.  Translations act trivially on $\Pic^0$, so the
same conclusion is independent of the chosen origin.  It follows that
each root exchange induces $[-1]$ on $\Pic^0$, their composition
induces the identity, and any automorphism of a genus-one curve acting
trivially on $\Pic^0$ is a translation.  The exchanges are nevertheless
wild rather than semisimple at their fixed points, as is unavoidable
for order-two automorphisms in characteristic two.

\section{Characteristic-two QRT arithmetic}
\label{subsec:QRT-char2-arithmetic}

Insert the parameters $a=a_2$ and $b=b_2$ from
\eqref{eq:QRT-char2-parameters} into the characteristic-two arithmetic
recorded in the standalone toolkit.  That arithmetic is written
on the shifted model
\begin{equation}
  \widehat Y^2+X\widehat Y
  =X^3+a_2X^2+b_2^2,
  \qquad
  \widehat Y=Y+b_2.
  \label{eq:QRT-char2-shifted-model}
\end{equation}
Consequently, an input from the QRT model is first converted to the
unshifted pair $(X,Y)$ by
Proposition~\ref{prop:QRT-char2-complete-interface} and then changed to
$(X,\widehat Y)=(X,Y+b_2)$.  After full-point arithmetic, one must undo
the shift by $Y=\widehat Y+b_2$ before applying the inverse QRT
interface.  The Kummer coordinate $X$ is unaffected by this ordinate
shift.

The resulting arithmetic is as follows:
\begin{itemize}
  \item affine full addition and doubling are given by
  Proposition~\ref{prop:binary-affine-law}, with costs
  $3\M+1\Sqr+1\Inv$ and $2\M+2\Sqr+1\Inv$;
  \item general, mixed, and affine L\'opez--Dahab addition cost
  $13\M+4\Sqr$, $8\M+5\Sqr+1\Dpar$, and
  $5\M+3\Sqr+1\Dpar$, respectively;
  \item L\'opez--Dahab doubling costs
  $3\M+5\Sqr+2\Dpar$, or
  $1\M+3\Sqr+2\Dpar$ for an affine input;
  \item Kummer differential addition and doubling are
  \eqref{eq:binary-xADD-product} and \eqref{eq:binary-xDBL}.  A
  differential-addition-and-doubling step costs
  $6\M+5\Sqr+2\Dpar$, or
  $5\M+5\Sqr+2\Dpar$ when the known difference is affine;
  \item the supplemental law \eqref{eq:binary-xADD-sum} completes the
  differential atlas, and \eqref{eq:binary-point-recovery} recovers a
  full point.
\end{itemize}
On the regular affine locus, the QRT interface contributes
$1\M+1\Sqr+3\Dpar$ on input and
$1\Inv+1\M+1\Sqr+3\Dpar$ on output.  Formula
\eqref{eq:QRT-char2-complete-boundary-atlas} handles every exceptional
point without evaluating a zero denominator.  During scalar
multiplication one stores $(X:Y:Z)$ in L\'opez--Dahab coordinates for
full-point arithmetic, or $(X:Z)$ on the Kummer line for a ladder.  An
$\ell$-bit binary Kummer ladder therefore has the upper bound
\begin{equation}
  (\ell-1)(6\M+5\Sqr+2\Dpar),
  \label{eq:QRT-char2-scalar-cost}
\end{equation}
or
$(\ell-1)(5\M+5\Sqr+2\Dpar)$ when the fixed difference is affine,
followed by one point recovery if a full output is required.

\section{The binary McMillan map and its translation point}
\label{subsec:QRT-char2-McMillan}

In characteristic two the same root exchange produces the recurrence
\begin{equation}
  \mathcal M_Q(x,y)
  =\left(y,x+\frac{\beta y}{y^2+\alpha}\right),
  \label{eq:QRT-char2-McMillan-map}
\end{equation}
and therefore
\begin{equation}
  x_{n+2}+x_n
  =\frac{\beta x_{n+1}}{x_{n+1}^2+\alpha}.
  \label{eq:QRT-char2-McMillan-recurrence}
\end{equation}
The invariant is still
\eqref{eq:QRT-McMillan-invariant}; subtraction and addition coincide.
The product-root projective formula
\eqref{eq:QRT-McMillan-projective-product} remains valid without any
change, while the numerator in the supplemental sum formula becomes
$N=XH+\beta ZC$.

Assume $k$ is perfect and use the notation
$r^2=\alpha$, $s^2=\gamma$ from
\eqref{eq:QRT-char2-square-roots}.  Select the boundary point
$P_\infty=(\infty,r)$ as origin.  The McMillan map sends it to
\begin{equation}
  P_f=\left(r,\frac{\alpha^2+\gamma}{\beta r}\right).
  \label{eq:QRT-char2-McMillan-point-Q}
\end{equation}
Under the complete interface
\eqref{eq:QRT-char2-complete-boundary-atlas}, this point becomes
\begin{equation}
  P_M^{\mathrm{unshifted}}
  =\left(
      \frac{\alpha}{\beta},
      \frac{\alpha(\alpha+s)}{\beta^2}
    \right).
  \label{eq:QRT-char2-McMillan-unshifted}
\end{equation}
After the ordinate shift
$\widehat Y=Y+b_2$, it is
\begin{equation}
  P_M^{(2)}
  =\left(X_0,X_0^2\right),
  \qquad
  X_0=\frac{\alpha}{\beta}.
  \label{eq:QRT-char2-McMillan-shifted}
\end{equation}

\begin{proposition}[Doubling the binary McMillan increment]
\label{prop:QRT-char2-McMillan-double}
On the shifted curve
\eqref{eq:QRT-char2-shifted-model},
\begin{equation}
  2P_M^{(2)}
  =\left(a_2,X_0^2+a_2\right).
  \label{eq:QRT-char2-McMillan-double}
\end{equation}
The point $P_M^{(2)}$ has exact order four if and only if
\begin{equation}
  a_2=0
  \quad\Longleftrightarrow\quad
  \gamma=\alpha^2.
  \label{eq:QRT-char2-McMillan-order-four}
\end{equation}
\end{proposition}

\begin{proof}
For $P=(x,y)=(X_0,X_0^2)$, the binary tangent slope in
\eqref{eq:binary-affine-double-lambda} is
\[
  \lambda=X_0+\frac{X_0^2}{X_0}=0.
\]
The doubling formulas therefore give
$x(2P)=a_2$ and $y(2P)=a_2+X_0^2$, which is
\eqref{eq:QRT-char2-McMillan-double}.  The unique nonzero two-torsion
point on the shifted model is $(0,b_2)$.  Hence $P_M$ has exact order
four exactly when $a_2=0$.  Since
$a_2=(\alpha^2+\gamma)/\beta^2$, this is equivalent to
$\gamma=\alpha^2$.  In that case $s=\alpha$ because the square root is
unique, and $b_2=\alpha^2/\beta^2=X_0^2$, so the double is indeed the
nonzero two-torsion point.
\end{proof}

Thus the smooth binary family has a genuinely nontrivial McMillan
increment because $\beta\ne0$ is forced by smoothness.  The special
order-four locus is $\gamma=\alpha^2$; the odd-characteristic branch
$\beta=0$ disappears because it is singular in characteristic two.
For scalar multiplication, the Kummer difference is the affine value
\begin{equation}
  X(P_M)=\frac{\alpha}{\beta}.
  \label{eq:QRT-char2-McMillan-known-difference}
\end{equation}
Consequently the binary QRT ladder uses the affine-difference cost
\begin{equation}
  5\M+5\Sqr+2\Dpar
  \label{eq:QRT-char2-McMillan-ladder-cost}
\end{equation}
per differential-addition-and-doubling step, before recovery.  Direct
iteration by the projective McMillan formula still costs
$2\M+2\Sqr+3\Dpar$ per step, but requires $O(n)$ steps to reach the
$n$th translate.

\section[Odd finite-field equivalence]{Finite-field model equivalence in odd characteristic}
\label{subsec:QRT-finite-field-model-equivalence-odd}

Let $k=\F_q$ with $q$ odd, and let $\chi:k^\times\to\{\pm1\}$ be the
quadratic character.  It is important to distinguish three notions:

\begin{enumerate}[label=(\roman*)]
  \item equality or equivalence of coefficient triples under a
  prescribed coordinate group;
  \item isomorphism of the completed curves as abstract genus-one
  curves;
  \item isomorphism of elliptic curves after an origin and, possibly, a
  marked two-torsion point have been selected.
\end{enumerate}

The counts below state explicitly which notion is being used.

For a smooth odd-characteristic member, introduce the normalized
parameters
\begin{equation}
  b=\frac{\beta}{\alpha},
  \qquad
  g=\frac{\gamma}{\alpha^2}.
  \label{eq:QRT-normalized-bg-odd}
\end{equation}
Then the factorized smoothness condition is
\begin{equation}
  g\ne0,
  \qquad
  \bigl((b+2)^2-4g\bigr)
  \bigl((b-2)^2-4g\bigr)\ne0.
  \label{eq:QRT-normalized-smooth-odd}
\end{equation}

\begin{proposition}[Number of smooth coefficient triples]
\label{prop:QRT-count-smooth-triples-odd}
The number of triples $(\alpha,\beta,\gamma)\in\F_q^3$ for which
$\Q_{\alpha,\beta,\gamma}$ is smooth is
\begin{equation}
  (q-1)(q^2-3q+3).
  \label{eq:QRT-count-smooth-triples-odd}
\end{equation}
\end{proposition}

\begin{proof}
Choose $\alpha\in\F_q^\times$ and use
\eqref{eq:QRT-normalized-bg-odd}.  There are $q(q-1)$ initial pairs
$(b,g)$ with $g\ne0$.  The first factor in
\eqref{eq:QRT-normalized-smooth-odd} vanishes for
\[
  g=(b+2)^2/4.
\]
Among pairs with $g\ne0$, this gives $q-1$ pairs, because $b=-2$ gives
$g=0$.  The second factor gives another $q-1$ pairs.  The two bad sets
intersect only at $(b,g)=(0,1)$: equality of the two displayed values
of $g$ gives $8b=0$, hence $b=0$, and then $g=1$.  Therefore the number
of bad pairs is $2q-3$, and the number of good pairs is
$q(q-1)-(2q-3)=q^2-3q+3$.  Multiplication by the $q-1$ choices of
$\alpha$ proves the formula.
\end{proof}

\begin{definition}[Diagonal symmetric model equivalence]
\label{def:QRT-diagonal-model-equivalence}
Two smooth odd-characteristic models are \emph{diagonally symmetric
model-equivalent} if one is obtained from the other by a transformation
\begin{equation}
  x=sX,\qquad y=\varepsilon sY,
  \qquad s\in\F_q^\times,
  \quad \varepsilon\in\{\pm1\},
  \label{eq:QRT-diagonal-model-change}
\end{equation}
possibly followed by interchange of $X$ and $Y$, and multiplication of
the equation by a nonzero scalar.
\end{definition}

\begin{theorem}[Exact diagonal classification and count]
\label{thm:QRT-diagonal-classification-odd}
Let $q$ be odd.  Two smooth models
$\Q_{\alpha,\beta,\gamma}$ and
$\Q_{\alpha',\beta',\gamma'}$ are diagonally symmetric
model-equivalent if and only if
\begin{equation}
  \chi(\alpha)=\chi(\alpha'),
  \qquad
  \frac{\gamma}{\alpha^2}
    =\frac{\gamma'}{\alpha'^2},
  \qquad
  \frac{\beta'}{\alpha'}
    =\pm\frac{\beta}{\alpha}.
  \label{eq:QRT-diagonal-classification-odd}
\end{equation}
The number of equivalence classes is
\begin{equation}
  (q-1)^2.
  \label{eq:QRT-diagonal-class-count-odd}
\end{equation}
\end{theorem}

\begin{proof}
Under \eqref{eq:QRT-diagonal-model-change}, division by $s^4$ gives
\[
  \alpha'=\frac{\alpha}{s^2},
  \qquad
  \beta'=\varepsilon\frac{\beta}{s^2},
  \qquad
  \gamma'=\frac{\gamma}{s^4}.
\]
Hence the three quantities in
\eqref{eq:QRT-diagonal-classification-odd} are invariant.  Conversely,
if those conditions hold, $\alpha/\alpha'$ is a square; choose $s$ with
$s^2=\alpha/\alpha'$ and choose $\varepsilon$ to match the signs of the
normalized $b$-parameters.  The equality of the normalized
$g$-parameters then gives the required equality of $\gamma'$.

It remains to count.  There are two square classes for $\alpha$.  The
involution $b\mapsto-b$ acts on the $q^2-3q+3$ smooth pairs $(b,g)$.
Its fixed pairs have $b=0$ and $g\in\F_q^\times\setminus\{1\}$, so
there are $q-2$ of them.  Burnside's lemma gives
\[
  \frac{(q^2-3q+3)+(q-2)}2=\frac{(q-1)^2}{2}
\]
orbits in $(b,g)$.  Multiplying by the two square classes of $\alpha$
gives \eqref{eq:QRT-diagonal-class-count-odd}.
\end{proof}

A larger but still geometrically natural coordinate group is obtained
by normalizing the involution $z\mapsto-z$ on each projective line.
Its normalizer in $\operatorname{PGL}_2(k)$ consists of the maps
$z\mapsto sz$ and $z\mapsto s/z$.

\begin{definition}[Normalizer equivalence]
\label{def:QRT-projection-normalizer-equivalence}
Two smooth odd-characteristic QRT models are
\emph{projection-normalizer equivalent} if they are related by
independent maps of the forms
\[
  x\mapsto sx\quad\hbox{or}\quad x\mapsto s/x,
  \qquad
  y\mapsto ty\quad\hbox{or}\quad y\mapsto t/y,
\]
followed by interchange of the two factors and renormalization of the
coefficient of $x^2y^2$ to one, provided the transformed coefficients
of $x^2$ and $y^2$ agree.
\end{definition}

The four types of transformation are explicit.  Besides
\eqref{eq:QRT-diagonal-model-change}, simultaneous inversion with
$t=\varepsilon s$ gives
\begin{equation}
  (\alpha,\beta,\gamma)
  \longmapsto
  \left(
    \frac{\alpha s^2}{\gamma},
    \varepsilon\frac{\beta s^2}{\gamma},
    \frac{s^4}{\gamma}
  \right).
  \label{eq:QRT-simultaneous-inversion-parameters}
\end{equation}
Inverting only the first coordinate is possible exactly when one can
choose $s,t\in k^\times$ with $s^2t^2=\gamma$; thus it is possible over
$k$ exactly when $\gamma$ is a square.  In that case
\begin{equation}
  (\alpha,\beta,\gamma)
  \longmapsto
  \left(
    \frac{s^2}{\alpha},
    \frac{\beta s}{\alpha t},
    \frac{s^2}{t^2}
  \right),
  \qquad s^2t^2=\gamma.
  \label{eq:QRT-one-inversion-parameters}
\end{equation}
The fourth type is obtained by interchanging the two coordinates.

In normalized data $(\delta,b,g)$, where
$\delta=\chi(\alpha)$, these operations are generated by
\begin{align}
  (\delta,b,g)&\longmapsto(\delta,-b,g),
  \label{eq:QRT-normalizer-generator-sign}\\
  (\delta,b,g)&\longmapsto(\delta\chi(g),b,g),
  \label{eq:QRT-normalizer-generator-both-inv}\\
  (\delta,b,g)&\longmapsto
  \left(\delta,\frac{b}{h},\frac1g\right),
  \qquad h^2=g,
  \label{eq:QRT-normalizer-generator-one-inv}
\end{align}
where the last generator is available only when $g$ is a square; the
choice $h\mapsto-h$ is absorbed by
\eqref{eq:QRT-normalizer-generator-sign}.

\begin{theorem}[Exact projection-normalizer count]
\label{thm:QRT-normalizer-count-odd}
Let $q$ be odd and put
\begin{equation}
  \varepsilon_4(q)=
  \begin{cases}
    1,&-1\text{ is a square in }\F_q,\\
    0,&-1\text{ is a nonsquare in }\F_q.
  \end{cases}
  \label{eq:epsilon4-definition}
\end{equation}
The number of projection-normalizer equivalence classes of smooth
models is
\begin{equation}
  N_{\mathrm{norm}}^{\mathrm{odd}}(q)
  =\frac{(q+1)(q-2)}2+\varepsilon_4(q).
  \label{eq:QRT-normalizer-count-odd}
\end{equation}
Moreover, two models are equivalent if and only if their normalized
triples lie in the same orbit under
\eqref{eq:QRT-normalizer-generator-sign}--
\eqref{eq:QRT-normalizer-generator-one-inv}.
\end{theorem}

\begin{proof}
The orbit criterion follows from the four explicit normalizer types
above, so only the count remains.

Suppose first that $g$ is a nonsquare.  Neither factor in
\eqref{eq:QRT-normalized-smooth-odd} can vanish, since it would express
$g$ as a square.  Hence every $b\in\F_q$ is allowed.  Simultaneous
inversion identifies the two values of $\delta$, while
$b\mapsto-b$ has $(q+1)/2$ orbits.  Since there are $(q-1)/2$
nonsquare values of $g$, this contributes
\begin{equation}
  \frac{q^2-1}{4}
  \label{eq:QRT-normalizer-nonsquare-contribution}
\end{equation}
to the orbit count.

Now suppose that $g$ is a square.  The total number of smooth raw pairs
$(b,g)$ with square $g$ is obtained by subtracting the
$q(q-1)/2$ nonsquare-$g$ pairs from $q^2-3q+3$; it equals
$(q-2)(q-3)/2$.  Among these, the fixed pairs for $b\mapsto-b$ have
$b=0$ and square $g\ne1$, so there are $(q-3)/2$.  The number of
sign-orbits is therefore
\begin{equation}
  \frac{(q-2)(q-3)/2+(q-3)/2}{2}
  =\frac{(q-1)(q-3)}4.
  \label{eq:QRT-square-g-sign-orbits}
\end{equation}
On this set, one-coordinate inversion induces an involution.  Its fixed
sign-orbits have $g=g^{-1}$.  For $g=1$, all smooth sign-orbits are
fixed and their number is $(q-3)/2$.  If $-1$ is a square, the case
$g=-1$ contributes one additional fixed sign-orbit, namely $b=0$; no
nonzero $b$ is fixed modulo sign because $b/h=\pm b$ would force
$h=\pm1$, contrary to $h^2=-1$.  Burnside's lemma, followed by the two
unchanged values of $\delta$, gives the square-$g$ contribution
\begin{equation}
  \frac{(q-3)(q+1)}4+\varepsilon_4(q).
  \label{eq:QRT-normalizer-square-contribution}
\end{equation}
Adding \eqref{eq:QRT-normalizer-nonsquare-contribution} and
\eqref{eq:QRT-normalizer-square-contribution} yields
\eqref{eq:QRT-normalizer-count-odd}.
\end{proof}

\section[Odd finite-field isomorphism]{Abstract finite-field isomorphism in odd characteristic}
\label{subsec:QRT-abstract-isomorphism-odd}

Model equivalence is finer than abstract isomorphism.  For
$\charac(k)\ne2,3$, the invariants of the Jacobian
\eqref{eq:QRT-Jacobian} are
\begin{align}
  c_4(Q)&=16\bigl(q_Q^2+192\alpha^2\gamma\bigr),
  \label{eq:QRT-c4-explicit}\\
  c_6(Q)&=64q_Q\bigl(576\alpha^2\gamma-q_Q^2\bigr),
  \label{eq:QRT-c6-explicit}\\
  \Delta(Q)&=4096\alpha^2\gamma\Omega_Q^2.
  \label{eq:QRT-Delta-explicit-again}
\end{align}
The first and third formulas agree with
\eqref{eq:QRT-j} and \eqref{eq:QRT-Delta-odd}; the middle formula
follows from $c_6=-b_2^3+36b_2b_4$ for
$V^2=U^3-2q_QU^2+\Omega_QU$.

\begin{theorem}[Abstract isomorphism criterion over $\F_q$]
\label{thm:QRT-abstract-isomorphism-odd}
Let $q$ have characteristic greater than three, and let $Q,Q'$ be
smooth members of the symmetric QRT family over $\F_q$.  Then their
smooth projective completions are $\F_q$-isomorphic if and only if
there exists $u\in\F_q^\times$ such that
\begin{equation}
  c_4(Q')=u^{-4}c_4(Q),
  \qquad
  c_6(Q')=u^{-6}c_6(Q).
  \label{eq:QRT-abstract-isomorphism-c4c6}
\end{equation}
\end{theorem}

\begin{proof}
Every smooth genus-one curve over a finite field has a rational point,
so both QRT curves are elliptic curves after choosing origins and are
isomorphic to their Jacobians.  In characteristic greater than three,
two Weierstrass curves are isomorphic over the ground field exactly
when their invariants are related by the indicated $u$-scaling.  Apply
this standard criterion to the explicit Jacobians
\eqref{eq:QRT-Jacobian}.  Conversely, such a scaling gives a
Weierstrass isomorphism, and composing with the two genus-one
identifications gives an isomorphism of the QRT completions.
\end{proof}

\subsection{Characteristic-three abstract classification}
The preceding $c_4,c_6$ scaling criterion uses the usual
characteristic-greater-than-three normal form.  In characteristic three
one can nevertheless give an equally explicit classification by using
the coefficient $q_Q$ itself.

\begin{theorem}[Abstract classification in characteristic three]
\label{thm:QRT-abstract-isomorphism-char3}
Let $q=3^m$, and let $Q,Q'$ be smooth symmetric QRT curves over
$\F_q$.  Write $q_Q,\Omega_Q$ and $q_{Q'},\Omega_{Q'}$ as in
\eqref{eq:QRT-q-Omega}.

\begin{enumerate}[label=(\roman*)]
  \item If $q_Q=q_{Q'}=0$, then $Q$ and $Q'$ are supersingular and
  \begin{equation}
    Q\cong_{\F_q}Q'
    \quad\Longleftrightarrow\quad
    \frac{\Omega_{Q'}}{\Omega_Q}\in(\F_q^\times)^4.
    \label{eq:QRT-char3-supersingular-criterion}
  \end{equation}

  \item If $q_Qq_{Q'}\ne0$, then $Q$ and $Q'$ are ordinary and
  \begin{equation}
    Q\cong_{\F_q}Q'
    \quad\Longleftrightarrow\quad
    j(Q)=j(Q')
    \quad\hbox{and}\quad
    \chi(q_Q)=\chi(q_{Q'}),
    \label{eq:QRT-char3-ordinary-criterion}
  \end{equation}
  where
  \begin{equation}
    j(Q)=
    \frac{q_Q^6}{\Omega_Q^2(q_Q^2-\Omega_Q)}.
    \label{eq:QRT-char3-j}
  \end{equation}

  \item A supersingular member is never isomorphic to an ordinary
  member.
\end{enumerate}
\end{theorem}

\begin{proof}
In characteristic three, \eqref{eq:QRT-Jacobian} becomes
\begin{equation}
  V^2=U^3+q_QU^2+\Omega_QU.
  \label{eq:QRT-char3-Jacobian}
\end{equation}
If $q_Q=0$, this is
$V^2=U^3+\Omega_QU$.  Its discriminant is nonzero because
$\Omega_Q\ne0$, and its $j$-invariant is zero.  A change
$U=u^2X$, $V=u^3Y$ replaces $\Omega_Q$ by
$\Omega_Q/u^4$.  Conversely, an admissible isomorphism between two
models of this shape has this effect on the nonzero linear
coefficient; an $X$-translation can preserve the zero constant term
but cannot alter the fourth-power class of the linear coefficient.
This proves \eqref{eq:QRT-char3-supersingular-criterion}.

Suppose $q_Q\ne0$ and put
\begin{equation}
  r_Q=\frac{\Omega_Q}{q_Q},
  \qquad
  c_Q=\frac{\Omega_Q^2(\Omega_Q-q_Q^2)}{q_Q^3}.
  \label{eq:QRT-char3-r-c}
\end{equation}
The substitution $U=X+r_Q$ removes the linear term and gives
\begin{equation}
  V^2=X^3+q_QX^2+c_Q.
  \label{eq:QRT-char3-ordinary-normal-form}
\end{equation}
Both $q_Q$ and $c_Q$ are nonzero by smoothness.  For a curve
$Y^2=X^3+aX^2+c$ with $ac\ne0$, an admissible change between two
normal forms of the same shape must have zero $X$-translation: after
$X=u^2X'+r$, the coefficient of $X'$ is $-ar/u^4$, and hence $r=0$.
The remaining scaling sends
\[
  (a,c)\longmapsto(a/u^2,c/u^6).
\]
Therefore two such curves are isomorphic exactly when their
$j=-a^3/c$ values agree and $a,a'$ have the same square class.  Applied
to \eqref{eq:QRT-char3-ordinary-normal-form}, this gives
\eqref{eq:QRT-char3-ordinary-criterion}; substituting
\eqref{eq:QRT-char3-r-c} gives \eqref{eq:QRT-char3-j}.
Finally, $q_Q=0$ is equivalent to $c_4=0$, whereas $q_Q\ne0$ gives
$c_4=q_Q^2\ne0$, so the two cases cannot be isomorphic.
\end{proof}

\begin{theorem}[Which elliptic curves admit a symmetric QRT model]
\label{thm:QRT-representation-by-point-odd}
Let $k$ have characteristic different from two.  Let
\begin{equation}
  E_{q_0,\Omega_0}:
  V^2=U(U^2-2q_0U+\Omega_0)
  \label{eq:QRT-representation-source-E}
\end{equation}
be smooth, so $\Omega_0(q_0^2-\Omega_0)\ne0$.  Then $E_{q_0,\Omega_0}$
is the Jacobian of a smooth model $\Q_{\alpha,\beta,\gamma}$ with
$q_Q=q_0$ and $\Omega_Q=\Omega_0$ if and only if
$E_{q_0,\Omega_0}(k)$ contains a point
$P=(D,V)$ with $DV\ne0$.

Given such a point, one may take
\begin{equation}
  \alpha=\frac{V}{4D},
  \qquad
  t=\frac{q_0-D}{4\alpha},
  \qquad
  \beta=t+2\alpha,
  \qquad
  \gamma=\frac{t^2-D}{4}.
  \label{eq:QRT-reconstruction-from-point-odd}
\end{equation}
Conversely, every smooth QRT model supplies the point
\begin{equation}
  P_M=(D_-,4\alpha D_-)
  \label{eq:QRT-representation-converse-point}
\end{equation}
on its Jacobian, and both coordinates are nonzero.
\end{theorem}

\begin{proof}
Assume first that $P=(D,V)$ with $DV\ne0$ is given and define the
parameters by \eqref{eq:QRT-reconstruction-from-point-odd}.  Then
\[
  (\beta-2\alpha)^2-4\gamma=t^2-(t^2-D)=D.
\]
Moreover,
\[
  q_Q=\beta^2-4\alpha^2-4\gamma
  =4\alpha t+D=q_0.
\]
Using the curve equation at $P$,
\[
  \frac{V^2}{D}=D^2-2q_0D+\Omega_0,
\]
and $16\alpha^2D=V^2/D$, one obtains
\begin{align*}
  64\alpha^2\gamma
  &=16\alpha^2(t^2-D)\\
  &=(q_0-D)^2-\frac{V^2}{D}
   =q_0^2-\Omega_0.
\end{align*}
Consequently $\Omega_Q=q_Q^2-64\alpha^2\gamma=\Omega_0$.  Smoothness
of $E_{q_0,\Omega_0}$ implies $\alpha\gamma\Omega_Q\ne0$, so the QRT
model is smooth.

Conversely, Theorem~\ref{thm:QRT-McMillan-point} provides
\eqref{eq:QRT-representation-converse-point}.  Smoothness gives
$\alpha D_-\ne0$, so its ordinate is nonzero.  This proves both
necessity and sufficiency.
\end{proof}

\begin{theorem}[Exact characteristic-three class count]
\label{thm:QRT-abstract-count-char3}
Let $q=3^m$.  The number of $\F_q$-isomorphism classes of elliptic
curves having a nonzero rational two-torsion point is
\begin{equation}
  N_2^{(3)}(q)
  =\frac{4q-6}{3}+\gcd(4,q-1).
  \label{eq:QRT-char3-N2}
\end{equation}
The number represented by smooth symmetric QRT models is
\begin{equation}
  N_Q^{\mathrm{abs}}(q)=
  \begin{cases}
    2,&q=3,\\
    13,&q=9,\\
    \displaystyle\frac{4q-6}{3}+\gcd(4,q-1),&q\ge27.
  \end{cases}
  \label{eq:QRT-char3-Q-count}
\end{equation}
\end{theorem}

\begin{proof}
We first count the ordinary classes.  Every ordinary curve in
characteristic three has a normal form
\begin{equation}
  E_{a,c}:\qquad y^2=x^3+ax^2+c,
  \qquad ac\ne0,
  \label{eq:QRT-char3-count-normal-form}
\end{equation}
and two such forms are isomorphic exactly under
$(a,c)\mapsto(a/u^2,c/u^6)$.

Count coefficient pairs in \eqref{eq:QRT-char3-count-normal-form} for
which the cubic has a rational root.  A marked root $r$ must be
nonzero, and then
\[
  c=-r^3-ar^2=-r^2(r+a).
\]
For each $r\in\F_q^\times$, all $a\in\F_q^\times\setminus\{-r\}$
are allowed, giving $(q-1)(q-2)$ marked-root pairs.  A separable cubic
of this form has either one or three rational roots.  To count the
completely split cases, choose an ordered pair of distinct nonzero
roots $r,s$ with $s\ne-r$.  The vanishing of the $x$-coefficient forces
\[
  t=-\frac{rs}{r+s}.
\]
The three roots are then distinct and their sum is
$(r-s)^2/(r+s)\ne0$.  Conversely every completely split polynomial is
obtained from its six ordered pairs of roots.  Hence their number is
\[
  \frac{(q-1)(q-3)}6.
\]
If $N_1$ and $N_3$ denote the numbers with one and three rational
roots, then
\[
  N_1+3N_3=(q-1)(q-2),
  \qquad
  N_3=\frac{(q-1)(q-3)}6.
\]
Thus
\begin{equation}
  N_1+N_3=\frac{(q-1)(2q-3)}3.
  \label{eq:QRT-char3-root-pair-count}
\end{equation}
In the scaling action, a pair with $ac\ne0$ is fixed only by
$u=\pm1$.  Burnside's lemma applied to
\eqref{eq:QRT-char3-root-pair-count} gives
\[
  \frac{2}{q-1}\cdot\frac{(q-1)(2q-3)}3
  =\frac{4q-6}{3}
\]
ordinary isomorphism classes with rational two-torsion.

The supersingular classes with rational two-torsion have the form
$y^2=x^3+Bx$, $B\ne0$.  Their isomorphism classes are the fourth-power
classes of $B$, so there are $\gcd(4,q-1)$ of them.  This proves
\eqref{eq:QRT-char3-N2}.

For $q\ge27$, Hasse's bound gives
$\#E(\F_q)>4$.  Hence every class counted by
\eqref{eq:QRT-char3-N2} contains a rational point outside its rational
$2$-torsion subgroup, and
Theorem~\ref{thm:QRT-representation-by-point-odd} constructs a smooth
QRT model.

It remains to justify the two finite exceptions.  Over $\F_3$, direct
use of \eqref{eq:QRT-smoothness-odd} gives six smooth coefficient
triples and exactly two abstract classes, represented by
\[
  y^2=x^3+x,
  \qquad
  y^2=x^3+x^2+1.
\]
The two classes with rational two-torsion that are not represented are
\[
  y^2=x^3+2x,
  \qquad
  y^2=x^3+2x^2+2;
\]
their rational groups contain no point outside $E[2]$.

For $\F_9=\F_3(\omega)$ with $\omega^2=-1$, the only omitted class is
\begin{equation}
  y^2=x^3+\omega x.
  \label{eq:QRT-char3-F9-exception}
\end{equation}
Indeed, if $\rho=1+\omega$, then $\rho^2=-\omega$ and
\[
  E(\F_9)=\{O,(0,0),(\rho,0),(-\rho,0)\}=E[2](\F_9).
\]
A curve over $\F_9$ with four rational points and full rational
$2$-torsion has trace $6$ and is supersingular.  The four
supersingular fourth-power twists can be represented by
$B=1,\rho,\omega,\rho^3$.  Evaluating the quadratic-character sum
$1+\sum_{x\in\F_9}(1+\chi(x^3+Bx))$ gives the exact table
\[
\begin{array}{c|cccc}
  B&1&\rho&\omega&\rho^3\\ \hline
  \#\{y^2=x^3+Bx\}(\F_9)&16&10&4&10.
\end{array}
\]
The count $4$ therefore occurs only for the fourth-power class of
$B=\omega$.  Thus \eqref{eq:QRT-char3-F9-exception} is the unique
unrepresented class, and $N_Q^{\mathrm{abs}}(9)=14-1=13$.
This proves \eqref{eq:QRT-char3-Q-count}.
\end{proof}

The modular parameter becomes especially compact if
\begin{equation}
  t_Q=\frac{q_Q^2}{64\alpha^2\gamma}.
  \label{eq:QRT-modular-t}
\end{equation}
Then smoothness is $t_Q\ne1$, and
\begin{equation}
  j(Q)=64\frac{(t_Q+3)^3}{(t_Q-1)^2}.
  \label{eq:QRT-j-modular-t}
\end{equation}
This is a rational parameter on the moduli problem with a marked
nonzero two-torsion point, as is already visible from the point
$(0,0)\in E_Q[2]$.

Thus the symmetric QRT normal form represents an elliptic curve with a
marked rational two-torsion point together with enough rational group
structure to choose a point outside the two-torsion subgroup.  This
extra point is precisely the McMillan increment.

\begin{theorem}[Exact number of abstract QRT isomorphism classes]
\label{thm:QRT-abstract-count-odd}
Let $q$ have characteristic greater than three.  Put
\begin{equation}
  \varepsilon_3(q)=
  \begin{cases}
    1,&3\mid(q-1),\\
    0,&3\nmid(q-1).
  \end{cases}
  \label{eq:epsilon3-definition}
\end{equation}
Equivalently, $\varepsilon_3(q)=1$ exactly when $-3$ is a square in
$\F_q$.  The number of $\F_q$-isomorphism classes of elliptic curves
having a nonzero rational two-torsion point is
\begin{equation}
  N_2(q)
  =\frac{4(q-2+\varepsilon_3(q))}{3}
   +\gcd(4,q-1).
  \label{eq:number-elliptic-classes-rational-2torsion}
\end{equation}
For $q\ge11$, every one of these classes has a smooth symmetric QRT
model, so the number of abstract isomorphism classes represented by
$\Q_{\alpha,\beta,\gamma}$ is exactly $N_2(q)$.  For the two remaining
fields of characteristic greater than three,
\begin{equation}
  N_Q^{\mathrm{abs}}(5)=6,
  \qquad
  N_Q^{\mathrm{abs}}(7)=9.
  \label{eq:QRT-small-odd-abstract-counts}
\end{equation}
\end{theorem}

\begin{proof}
We first count all elliptic-curve classes with rational two-torsion.
In characteristic greater than three, use short equations
\[
  E_{A,B}:y^2=x^3+Ax+B,
  \qquad 4A^3+27B^2\ne0.
\]
The group $\F_q^\times$ acts by
$(A,B)\mapsto(u^{-4}A,u^{-6}B)$, and its orbits are the
$\F_q$-isomorphism classes.

Let $\mathcal S$ be the set of nonsingular coefficient pairs for which
the cubic has a rational root.  Count first pairs with a marked root
$r$.  Writing $B=-r^3-Ar$, the cubic is singular if either the marked
root is repeated, which gives $A=-3r^2$, or the other two roots
coincide, which gives $A=-3r^2/4$.  For $r=0$ these exclusions coincide;
for $r\ne0$ they are distinct.  Hence the number of smooth pairs with a
marked root is
\begin{equation}
  (q-1)+(q-1)(q-2)=(q-1)^2.
  \label{eq:marked-root-pair-count}
\end{equation}
The remaining quadratic factor splits when its discriminant
$-4A-3r^2$ is a nonzero square.  For $r=0$ there are $(q-1)/2$ choices;
for $r\ne0$ there are $(q-3)/2$, because the square $9r^2$ must be
excluded.  Thus the number of marked roots belonging to completely
split cubics is $(q-1)(q-2)/2$.  Dividing by three gives
$(q-1)(q-2)/6$ completely split coefficient pairs.  It follows from
\eqref{eq:marked-root-pair-count} that
\begin{equation}
  |\mathcal S|=\frac{(q-1)(2q-1)}{3}.
  \label{eq:S-rational-root-size}
\end{equation}

Apply Burnside's lemma to the weighted scaling action.  The elements
$u=\pm1$ fix every pair in $\mathcal S$.  An element satisfying
$u^4=1$ but $u^6\ne1$ fixes the $q-1$ nonsingular pairs with $B=0$.
There are $\gcd(4,q-1)-2$ such elements.  An element satisfying
$u^6=1$ but $u^4\ne1$ fixes the pairs with $A=0$ for which $-B$ is a
cube.  This case occurs exactly when $3\mid q-1$; there are four such
group elements and $(q-1)/3$ fixed pairs for each.  No other element
fixes a nonsingular pair.  Substitution of
\eqref{eq:S-rational-root-size} into Burnside's formula gives exactly
\eqref{eq:number-elliptic-classes-rational-2torsion}.

Now suppose $q\ge11$.  Hasse's bound gives
\[
  \#E(\F_q)\ge q+1-2\sqrt q>4.
\]
Since $E[2](\F_q)$ has at most four points, every elliptic curve counted
by $N_2(q)$ has a rational point outside its two-torsion subgroup.
Theorem~\ref{thm:QRT-representation-by-point-odd} therefore constructs
a QRT model.  Conversely every smooth QRT model has rational
two-torsion, so the counts agree.

For $q=5$, the two excluded classes can be represented by
$y^2=x^3+x$ and $y^2=x^3+2x$.  Direct evaluation of the five abscissas
gives
\[
\begin{array}{c|c|c}
 E&\#E(\F_5)&E(\F_5)\setminus\{O\}\\ \hline
 y^2=x^3+x&4&(0,0),(2,0),(3,0)\\
 y^2=x^3+2x&2&(0,0).
\end{array}
\]
Thus every rational point is two-torsion.  For $q=7$, the unique
excluded class is $y^2=x^3-1$, and direct substitution gives
\[
  E(\F_7)=\{O,(1,0),(2,0),(4,0)\}=E[2](\F_7).
\]
Subtracting these exceptions from
\eqref{eq:number-elliptic-classes-rational-2torsion} gives
\eqref{eq:QRT-small-odd-abstract-counts}.
\end{proof}

\section{Finite-field classification in characteristic two}
\label{subsec:QRT-finite-field-classification-char2}

Let $k=\F_q$ with $q=2^m$.  Smoothness is exactly
$\alpha\beta\gamma\ne0$.  Because the square map is an automorphism of
$k$, diagonal scaling normalizes $\alpha$ uniquely.  With
\begin{equation}
  b=\frac{\beta}{\alpha},
  \qquad
  g=\frac{\gamma}{\alpha^2},
  \label{eq:QRT-normalized-bg-char2}
\end{equation}
we have $b,g\in k^\times$.

\begin{theorem}[Characteristic-two model counts]
\label{thm:QRT-char2-model-counts}
The number of smooth coefficient triples is
\begin{equation}
  (q-1)^3.
  \label{eq:QRT-char2-raw-count}
\end{equation}
The number of diagonal model-equivalence classes is
\begin{equation}
  (q-1)^2.
  \label{eq:QRT-char2-diagonal-count}
\end{equation}
Under the characteristic-two monomial projection equivalence generated
by diagonal scalings, factor interchange, simultaneous inversion, and
one-coordinate inversion, the number is
\begin{equation}
  N_{\mathrm{mon}}^{(2)}(q)=\frac{q(q-1)}2.
  \label{eq:QRT-char2-monomial-count}
\end{equation}
\end{theorem}

\begin{proof}
The raw count follows from the three independent nonzero parameters.
Under $x=sX$, $y=sY$, the normalized pair $(b,g)$ is unchanged, and the
unique square root of $\alpha/\alpha'$ supplies the converse scaling.
This proves \eqref{eq:QRT-char2-diagonal-count}.

Simultaneous inversion has no additional effect on $(b,g)$ in
characteristic two.  A one-coordinate inversion is always defined,
because every $g$ has a unique square root $h$, and it induces
\begin{equation}
  (b,g)\longmapsto\left(\frac{b}{h},\frac1g\right),
  \qquad h^2=g.
  \label{eq:QRT-char2-normalizer-involution}
\end{equation}
This is an involution on the $(q-1)^2$ normalized pairs.  Its fixed
points have $g=1$, hence $h=1$, while $b$ may be arbitrary nonzero;
there are $q-1$ fixed points.  Burnside's lemma gives
\[
  \frac{(q-1)^2+(q-1)}2=\frac{q(q-1)}2.
\]
\end{proof}

\begin{remark}[Meaning of monomial equivalence in characteristic two]
Here ``monomial'' is deliberate.  In characteristic two the map
$z\mapsto-z$ is the identity, whose full normalizer is all of
$\operatorname{PGL}_2$; therefore the odd-characteristic phrase
``full projection normalizer'' would be misleading.  The stated group
is exactly the characteristic-two specialization of the explicit
monomial transformations used in
Definition~\ref{def:QRT-projection-normalizer-equivalence}.
\end{remark}

For abstract isomorphism, use the shifted binary model
\eqref{eq:QRT-char2-shifted-model}.  Set
\begin{equation}
  a_Q=\frac{\alpha^2+\gamma}{\beta^2},
  \qquad
  b_Q=\frac{\alpha\sqrt\gamma}{\beta^2}.
  \label{eq:QRT-char2-abstract-invariants}
\end{equation}
Then
\[
  \widehat Y^2+X\widehat Y
  =X^3+a_QX^2+b_Q^2,
  \qquad j=b_Q^{-2}.
\]

\begin{theorem}[Binary abstract classification]
\label{thm:QRT-char2-abstract-classification}
Two smooth characteristic-two QRT curves over $\F_q$ are
$\F_q$-isomorphic if and only if
\begin{equation}
  b_Q=b_Q',
  \qquad
  \Tr_{\F_q/\F_2}(a_Q)
  =\Tr_{\F_q/\F_2}(a_Q').
  \label{eq:QRT-char2-abstract-criterion}
\end{equation}
For $q\ge8$, every ordinary elliptic-curve class over $\F_q$ has a
smooth symmetric QRT model.  Consequently,
\begin{equation}
  N_Q^{\mathrm{abs}}(q)=2(q-1)
  \qquad(q\ge8).
  \label{eq:QRT-char2-abstract-count-large}
\end{equation}
For the two smaller fields,
\begin{equation}
  N_Q^{\mathrm{abs}}(2)=1,
  \qquad
  N_Q^{\mathrm{abs}}(4)=5.
  \label{eq:QRT-char2-abstract-count-small}
\end{equation}
\end{theorem}

\begin{proof}
We first justify the isomorphism criterion without restricting in
advance the shape of an admissible change of variables.  Write
\[
  E_{a,b}:\quad Y^2+XY=X^3+aX^2+b^2,
  \qquad b\ne0.
\]
A general Weierstrass change from $E_{a,b}$ to another equation is
\[
  X=u^2X'+r,
  \qquad
  Y=u^3Y'+su^2X'+t,
  \qquad u\in\F_q^\times.
\]
After substitution and division by $u^6$, the coefficients of
$X'Y'$, $Y'$, and $X'$ are, respectively,
\[
  u^{-1},\qquad ru^{-3},\qquad
  u^{-4}(r^2+t+rs).
\]
For the target to have again $a_1=1$, $a_3=a_4=0$, these three
coefficients force successively
\[
  u=1,\qquad r=0,\qquad t=0.
\]
The remaining transformation is therefore
$Y=Y'+sX'$, and comparison of the $X'^2$ and constant coefficients
gives
\[
  a'=a+s^2+s,
  \qquad
  b'^2=b^2.
\]
The square map on $\F_q$ is bijective, so $b'=b$.  Conversely, every
$s\in\F_q$ gives precisely this isomorphism.  Hence two normal forms
are isomorphic exactly when their $b$-parameters agree and their
$a$-parameters differ by an Artin--Schreier element.  The image of
$s\mapsto s^2+s$ is the trace-zero subspace, which proves
\eqref{eq:QRT-char2-abstract-criterion}.  Equivalently, $b_Q$ is also
determined by $j=b_Q^{-2}$, and for each
$b_Q\in\F_q^\times$ there are exactly two trace classes for $a_Q$.
Thus the number of ordinary elliptic-curve classes is $2(q-1)$.

It remains to show representability for $q\ge8$.  Every ordinary
elliptic curve has the unique nonzero rational two-torsion point
$(0,b_Q)$.  Hasse's bound gives $\#E(\F_q)>2$ for $q\ge8$, and the group
order is even, so there is a point $P=(X_0,Y_0)$ outside
$\{O,(0,b_Q)\}$.  In particular $X_0\ne0$.  Replace
$Y$ by $Y+rX$ with
$r=(Y_0+X_0^2)/X_0$; then $P$ becomes $(X_0,X_0^2)$ and $a_Q$ changes
within the same isomorphism class.  The point equation gives
\begin{equation}
  X_0^4=a_QX_0^2+b_Q^2.
  \label{eq:QRT-char2-point-reconstruction-identity}
\end{equation}
Choose
\begin{equation}
  \beta=1,
  \qquad
  \alpha=X_0,
  \qquad
  \sqrt\gamma=\frac{b_Q}{X_0},
  \qquad
  \gamma=\frac{b_Q^2}{X_0^2}.
  \label{eq:QRT-char2-reconstruction}
\end{equation}
Then the quantities in \eqref{eq:QRT-char2-abstract-invariants} are
$b_Q$ and
$X_0^2+b_Q^2/X_0^2$, which equals $a_Q$ by
\eqref{eq:QRT-char2-point-reconstruction-identity}.  Thus the QRT model
is smooth and represents the original class.

For $q=2$, the only normalized pair is $(b,g)=(1,1)$, giving one
class.  For $q=4$, write
$\F_4=\F_2(\omega)$ with $\omega^2+\omega+1=0$, and use $a=0$ and
$a=\omega$ as representatives of the trace-zero and trace-one
classes.  The six ordinary classes and one point outside $E[2]$, when
such a point exists, are displayed in the following table:
\[
\begin{array}{c|c|c}
 b_Q&\Tr(a_Q)&\text{a point with }X\ne0\\ \hline
 1&0&(1,0)\\
 1&1&\text{none}\\
 \omega&0&(\omega^2,1)\\
 \omega&1&(1,0)\\
 \omega^2&0&(\omega,1)\\
 \omega^2&1&(1,\omega).
\end{array}
\]
Each entry follows by substitution in
$Y^2+XY=X^3+a_QX^2+b_Q^2$.  For the omitted class
$b_Q=1$, $\Tr(a_Q)=1$, a conceptual verification is also available:
for $X\ne0$, division by $X^2$ shows that an ordinate exists only if
\[
  \Tr\left(X+a_Q+X^{-2}\right)=0.
\]
Since $X^3=1$ on $\F_4^\times$, one has $X^{-2}=X$, so the displayed
trace equals $\Tr(a_Q)=1$.  Thus this curve has only the identity and
its unique nonzero two-torsion point.  The remaining five rows supply
points outside $E[2]$ and are represented by the construction above.
\end{proof}

The formulas
\eqref{eq:QRT-char2-abstract-invariants} also give an exact
coefficient-level realization test.  After normalizing $\alpha=1$, put
$t=\sqrt g$.  Then
\[
  a_Q=\frac{1+t^2}{b^2},
  \qquad
  b_Q=\frac{t}{b^2}.
\]
For a prescribed target
$Y^2+XY=X^3+aX^2+b_0^2$, the equality $b_Q=b_0$ first forces
$b^2=t/b_0$.  Substitution in $a_Q=a$ and multiplication by $t$ then
give the following necessary equation; reversing these two steps shows
that it is also sufficient.  Thus a normalized QRT model exists with
those specific coefficients if and only if there is
$t\in\F_q^\times$ such that
\begin{equation}
  t^2+\frac{a}{b_0}t+1=0,
  \qquad
  b^2=\frac{t}{b_0}.
  \label{eq:QRT-char2-realization-quadratic}
\end{equation}
If $a/b_0\ne0$, divide the quadratic by $(a/b_0)^2$ and put
$z=t/(a/b_0)$.  It becomes
\[
  z^2+z=\left(\frac{b_0}{a}\right)^2.
\]
An Artin--Schreier equation $z^2+z=c$ over $\F_q$ is soluble exactly
when $\Tr(c)=0$, and $\Tr(c^2)=\Tr(c)$.  Hence the quadratic is soluble
exactly when
\begin{equation}
  \Tr_{\F_q/\F_2}\left(\frac{b_0}{a}\right)=0.
  \label{eq:QRT-char2-realization-trace}
\end{equation}
If $a=0$, it is $(t+1)^2=0$ and has the unique root $t=1$.  Once $t$
is known, the equation $b^2=t/b_0$ has a unique solution because the
square map is an automorphism of $\F_q$.  This coefficient test is
finer than abstract isomorphism, because one may first change $a$ by an
Artin--Schreier element without changing the elliptic curve.

\begin{corollary}[Global finite-field isomorphism count]
\label{cor:QRT-global-finite-field-count}
Let $N_Q^{\mathrm{abs}}(q)$ denote the number of $\F_q$-isomorphism
classes of smooth projective curves represented by
$\Q_{\alpha,\beta,\gamma}$.  Then
\begin{equation}
N_Q^{\mathrm{abs}}(q)=
\begin{cases}
  1,&q=2,\\
  5,&q=4,\\
  2(q-1),&q=2^m\ge8,\\
  2,&q=3,\\
  13,&q=9,\\
  \dfrac{4q-6}{3}+\gcd(4,q-1),&q=3^m\ge27,\\
  6,&q=5,\\
  9,&q=7,\\
  \dfrac{4(q-2+\varepsilon_3(q))}{3}
       +\gcd(4,q-1),
       &\charac(\F_q)>3,\ q\ge11.
\end{cases}
\label{eq:QRT-global-finite-field-count}
\end{equation}
For two given smooth triples, an exact isomorphism decision is supplied
by Theorem~\ref{thm:QRT-abstract-isomorphism-odd} in characteristic
larger than three, Theorem~\ref{thm:QRT-abstract-isomorphism-char3}
in characteristic three, and
Theorem~\ref{thm:QRT-char2-abstract-classification} in characteristic
two.
\end{corollary}

\begin{proof}
The formula is the disjoint combination of
Theorems~\ref{thm:QRT-abstract-count-odd},
\ref{thm:QRT-abstract-count-char3}, and
\ref{thm:QRT-char2-abstract-classification}.  Those theorems also
prove the cited isomorphism criteria and account separately for every
small finite field at which the Hasse-bound representability argument
does not apply.
\end{proof}

\section{Model-preserving isogenies}
\label{subsec:QRT-model-preserving-isogenies}

An isogeny is intrinsic to the elliptic curve, whereas the equation
\eqref{eq:symmetric-QRT} contains additional coordinate data.  The
phrase \emph{model-preserving isogeny} will therefore have the following
precise meaning.

Before fixing coefficient normalizations, the adjacent-state construction
has an intrinsic functoriality statement.

\begin{theorem}[Coordinatewise transport of adjacent states]
\label{thm:QRT-coordinatewise-state-isogeny}
Let \(\varphi:E\to E'\) be an isogeny, let \(D\in E(k)\), put
\(D'=\varphi(D)\), and let
\(\kappa:E\to\PP^1\) and \(\kappa':E'\to\PP^1\) be separable degree-two
Kummer quotients.
There is a unique morphism
\(\overline\varphi:\PP^1\to\PP^1\) satisfying
\begin{equation}
 \overline\varphi\circ\kappa=\kappa'\circ\varphi.
 \label{eq:QRT-state-Kummer-isogeny-descent}
\end{equation}
It has degree \(\deg\varphi\).

If \(2D\ne O\) and \(2D'\ne O'\), the induced map between the two smooth
state curves is
\begin{equation}
 \Phi_\varphi
 =\mathcal S'_{D'}\circ\varphi\circ\mathcal S_D^{-1},
 \qquad
 \boxed{\Phi_\varphi(x,y)=
 (\overline\varphi(x),\overline\varphi(y))}.
 \label{eq:QRT-state-coordinatewise-isogeny}
\end{equation}
It has degree \(\deg\varphi\) and satisfies
\begin{align}
 \Phi_\varphi T_D&=T_{D'}\Phi_\varphi,
 \label{eq:QRT-state-isogeny-shift-commute}\\
 \Phi_\varphi\Delta_{D,m}&=
 \Delta_{D',m}\Phi_\varphi\qquad(m\in\mathbb Z),
 \label{eq:QRT-state-isogeny-multiplication-commute}
\end{align}
where
\(\Delta_{D,m}=\mathcal S_D\circ[m]\circ\mathcal S_D^{-1}\).
\end{theorem}

\begin{proof}
An isogeny is a group homomorphism and therefore commutes with negation.
The composite \(\kappa'\circ\varphi\) is consequently constant on the
generic fibres \(\{P,-P\}\) of \(\kappa\), so it descends uniquely to a
rational map \(\overline\varphi\) on the Kummer quotient.  A rational map
between smooth projective curves extends across its finitely many apparent
indeterminacies, hence \(\overline\varphi\) is a morphism.  Comparing
degrees in \eqref{eq:QRT-state-Kummer-isogeny-descent} gives
\[
 2\deg\varphi
 =\deg(\kappa'\circ\varphi)
 =2\deg\overline\varphi,
\]
including inseparable degrees, and proves
\(\deg\overline\varphi=\deg\varphi\).

For \(P\in E\), use \(\varphi(P+D)=\varphi(P)+D'\) to compute
\begin{align*}
 \mathcal S'_{D'}(\varphi(P))
 &=\bigl(\kappa'(\varphi(P)),
          \kappa'(\varphi(P)+D')\bigr)\\
 &=\bigl(\overline\varphi(\kappa(P)),
          \overline\varphi(\kappa(P+D))\bigr).
\end{align*}
This is \eqref{eq:QRT-state-coordinatewise-isogeny}.  Both state maps are
isomorphisms under the two nondegeneracy hypotheses, so conjugation proves
the degree statement.  The identities
\(\varphi(P+D)=\varphi(P)+D'\) and
\(\varphi([m]P)=[m]\varphi(P)\) prove
\eqref{eq:QRT-state-isogeny-shift-commute} and
\eqref{eq:QRT-state-isogeny-multiplication-commute}, respectively.
\end{proof}

If \(D'=O'\), the coordinatewise formula still maps to the diagonal
Kummer image.  If \(D'\ne O'\) but \(2D'=O'\), it maps to the corresponding
two-torsion state image, which also factors through the Kummer line.  In
either case \eqref{eq:QRT-state-coordinatewise-isogeny} remains a valid
coordinate identity, but the target state map is not an elliptic-state
isomorphism.  When exact even coordinates exist on both nondegenerate
state curves, their two PGL$_2$ normalizations merely conjugate
\(\overline\varphi\); the same one-variable rational function is still
applied to both state slots.

\begin{definition}[QRT model-preserving isogeny]
\label{def:QRT-model-preserving-isogeny}
Let $Q$ and $Q'$ be smooth QRT models with selected origins, and let
$\psi_Q:Q\to E_Q$ and $\psi_{Q'}:Q'\to E_{Q'}$ be explicit elliptic
identifications.  A separable isogeny
$\varphi:E_Q\to E_{Q'}$ is \emph{QRT model-preserving} relative to these
choices if the rational map
\begin{equation}
  \varphi_Q
  =\psi_{Q'}^{-1}\circ\varphi\circ\psi_Q
  \label{eq:QRT-model-preserving-composition}
\end{equation}
is defined over the ground field.  If only the Jacobians, but not the
torsors, are identified over the ground field, the term refers to the
Jacobian isogeny and not to a map between the original torsors.
\end{definition}

This definition records the ground-field data needed by the declared QRT
interface.  The codomain Weierstrass curve is realized in the selected QRT
coefficient normalization, and a genus-one torsor is supplied with the
rational origin and chart used by its birational interface.  Over finite
fields, the rational origin exists and the chosen chart determines the
explicit formula.

\subsection{Exact re-embedding criterion in odd characteristic.}
Let
\begin{equation}
  E':Y^2=X(X^2+A'X+B')
  \label{eq:QRT-isogeny-target-marked-2torsion}
\end{equation}
be smooth.  Fix $u\in k^\times$ and put
\begin{equation}
  q'=-\frac{A'}{2u^2},
  \qquad
  \Omega'=\frac{B'}{u^4},
  \qquad
  h'=\frac{q'^2-\Omega'}{64}.
  \label{eq:QRT-reembedding-qOmega}
\end{equation}
Then $E'$ is isomorphic, by $X=u^2U$, $Y=u^3V$, to the Jacobian of a
QRT model $\Q_{\alpha',\beta',\gamma'}$ if and only if there exists
$\alpha'\in k^\times$ for which
\begin{equation}
  \gamma'=\frac{h'}{\alpha'^2}
  \label{eq:QRT-reembedding-gamma}
\end{equation}
and
\begin{equation}
  \beta'^2=q'+4\alpha'^2+4\gamma'
  \label{eq:QRT-reembedding-beta-square}
\end{equation}
has a solution $\beta'\in k$.  These equations are necessary because
$q'=\beta'^2-4\alpha'^2-4\gamma'$ and
$h'=\alpha'^2\gamma'$.  Conversely, substituting the displayed
$\gamma'$ and $\beta'$ gives
\[
  q_{Q'}=\beta'^2-4\alpha'^2-4\gamma'=q'
\]
and
\[
  \Omega_{Q'}=q_{Q'}^2-64\alpha'^2\gamma'
  =q'^2-64h'=\Omega'.
\]
Thus the equations are also sufficient.  Over an algebraic closure
they are always soluble; over $k$ they are a genuine
field-of-definition test.

\subsection{The rational two-isogeny.}
For
\[
  E_Q:V^2=U(U^2-2q_QU+\Omega_Q),
\]
the point $(0,0)$ generates a rational subgroup of order two.  The
normalized quotient isogeny is
\begin{equation}
  \phi_2(U,V)
  =\left(
     U-2q_Q+\frac{\Omega_Q}{U},
     V\left(1-\frac{\Omega_Q}{U^2}\right)
   \right).
  \label{eq:QRT-explicit-2-isogeny}
\end{equation}
Equivalently, its first coordinate is $V^2/U^2$.  The codomain is
\begin{equation}
  E_Q^{(2)}:
  Y^2=X^3+4q_QX^2+256\alpha^2\gamma X.
  \label{eq:QRT-explicit-2-isogeny-target}
\end{equation}

\begin{proof}
The standard quotient of
$y^2=x^3+Ax^2+Bx$ by the subgroup generated by $(0,0)$ is
\[
  x'=x+A+B/x,
  \qquad
  y'=y(1-B/x^2),
\]
with codomain
$y'^2=x'^3-2Ax'^2+(A^2-4B)x'$.  Substituting
$A=-2q_Q$ and $B=\Omega_Q$ gives
\eqref{eq:QRT-explicit-2-isogeny}.  Moreover,
\[
  A^2-4B=4(q_Q^2-\Omega_Q)=256\alpha^2\gamma,
\]
which proves \eqref{eq:QRT-explicit-2-isogeny-target}.  If
$(X,Y)=\phi_2(U,V)$, then the target equation satisfies the exact
identity
\[
  Y^2-\bigl(X^3-2AX^2+(A^2-4B)X\bigr)
  =\frac{(B-U^2)^2}{U^4}
   \bigl(V^2-U^3-AU^2-BU\bigr).
\]
Hence every nonkernel point of the source maps to the codomain.  The
only finite pole is $U=0$, namely $(0,0)$, and both it and $O$ map to
the quotient identity; the map has degree two, so its kernel is
exactly $\{O,(0,0)\}$.
\end{proof}

For the target, the QRT invariants may be chosen as
\begin{equation}
  q_2=-2q_Q,
  \qquad
  \Omega_2=256\alpha^2\gamma,
  \qquad
  \alpha_2^2\gamma_2=\frac{\Omega_Q}{16}.
  \label{eq:QRT-2-isogeny-target-invariants}
\end{equation}
Thus, after choosing $\alpha_2\in k^\times$, one sets
\begin{equation}
  \gamma_2=\frac{\Omega_Q}{16\alpha_2^2},
  \qquad
  \beta_2^2
  =-2q_Q+4\alpha_2^2+\frac{\Omega_Q}{4\alpha_2^2}.
  \label{eq:QRT-2-isogeny-target-parameters}
\end{equation}
Whenever the last quantity is a square in $k$, these formulas produce
an explicit QRT codomain.  On the square-constant charts of
Section~\ref{subsec:QRT-odd-geometry}, the model-preserving map is
obtained by inserting \eqref{eq:QRT-explicit-2-isogeny} between the
explicit QRT--Weierstrass interfaces in
\eqref{eq:QRT-model-preserving-composition}.

\subsection{An explicit three-isogeny.}
Assume now that $\charac(k)\ne2,3$.  Write
\begin{equation}
  A=-2q_Q,
  \qquad B=\Omega_Q,
  \label{eq:QRT-3isog-AB}
\end{equation}
and make the shift
\begin{equation}
  U=x-\frac A3.
  \label{eq:QRT-shift-to-short}
\end{equation}
The resulting short equation is
\begin{equation}
  E_s:y^2=x^3+a_sx+b_s,
  \label{eq:QRT-short-for-isogeny}
\end{equation}
where
\begin{equation}
  a_s=B-\frac{A^2}{3},
  \qquad
  b_s=\frac{2A^3}{27}-\frac{AB}{3}.
  \label{eq:QRT-short-coefficients-isogeny}
\end{equation}
Let $K=(x_K,y_K)\in E_s(\overline{k})$ have exact order three and
assume that the subgroup $\{O,K,-K\}$ is Galois stable.  Equivalently,
$y_K\ne0$ and
\begin{equation}
  \psi_3(x_K)
  =3x_K^4+6a_sx_K^2+12b_sx_K-a_s^2=0.
  \label{eq:QRT-3division-polynomial}
\end{equation}
The Galois-stability assumption implies
$x_K,y_K^2\in k$, although $y_K$ itself need not belong to $k$.  Put
\begin{equation}
  t_K=3x_K^2+a_s,
  \qquad
  u_K=2y_K^2,
  \qquad
  w_K=u_K+t_Kx_K.
  \label{eq:QRT-3isog-tuw}
\end{equation}
Then the degree-three isogeny with kernel
$\{O,K,-K\}$ is
\begin{align}
  R_3(x)
  &=x+2\left(
       \frac{t_K}{x-x_K}
       +\frac{u_K}{(x-x_K)^2}
     \right),
  \label{eq:QRT-3isog-R}\\
  \phi_3(x,y)&=(R_3(x),R_3'(x)y),
  \label{eq:QRT-explicit-3isogeny}
\end{align}
with codomain
\begin{equation}
  E_3:y^2=x^3+a_3x+b_3,
  \qquad
  a_3=a_s-10t_K,
  \qquad
  b_3=b_s-14w_K.
  \label{eq:QRT-3isog-target-short}
\end{equation}
These are the degree-three specialization of V\'elu's formulas
\cite{Velu1971}.

The original rational two-torsion point has short $x$-coordinate
\begin{equation}
  e_0=\frac A3.
  \label{eq:QRT-original-2torsion-short-x}
\end{equation}
Since the kernel has odd order, its image remains a nonzero rational
two-torsion point.  Put
\begin{equation}
  e_3=R_3(e_0).
  \label{eq:QRT-image-2torsion-three-isog}
\end{equation}
The shift $x=X+e_3$ changes \eqref{eq:QRT-3isog-target-short} into
\begin{equation}
  Y^2=X(X^2+A_3X+B_3),
  \qquad
  A_3=3e_3,
  \qquad
  B_3=3e_3^2+a_3.
  \label{eq:QRT-3isog-marked-target}
\end{equation}
Equations \eqref{eq:QRT-reembedding-qOmega}--
\eqref{eq:QRT-reembedding-beta-square} then give the exact condition
and explicit parameters for a QRT codomain over $k$.

\subsection{Odd prime degree.}
Let $\ell\ne\charac(k)$ be an odd prime and let
$G\subset E_s(\overline{k})$ be a cyclic Galois-stable subgroup of
order $\ell$.  For $Q=(x_Q,y_Q)\in G\setminus\{O\}$ define
\begin{equation}
  t_Q=3x_Q^2+a_s,
  \qquad
  u_Q=2y_Q^2,
  \qquad
  w_Q=u_Q+t_Qx_Q.
  \label{eq:QRT-l-isog-local-data}
\end{equation}
Set
\begin{equation}
  t_G=\sum_{Q\in G\setminus\{O\}}t_Q,
  \qquad
  w_G=\sum_{Q\in G\setminus\{O\}}w_Q,
  \label{eq:QRT-l-isog-global-data}
\end{equation}
and
\begin{equation}
  R_G(x)=x+
  \sum_{Q\in G\setminus\{O\}}
  \left(
    \frac{t_Q}{x-x_Q}
    +\frac{u_Q}{(x-x_Q)^2}
  \right).
  \label{eq:QRT-l-isog-R}
\end{equation}
Then
\begin{equation}
  \phi_G(x,y)=(R_G(x),R_G'(x)y)
  \label{eq:QRT-explicit-l-isogeny}
\end{equation}
has kernel $G$ and codomain
\begin{equation}
  E_G:
  y^2=x^3+(a_s-5t_G)x+(b_s-7w_G).
  \label{eq:QRT-l-isog-target}
\end{equation}
The sums are Galois invariant, so the map and codomain are defined over
$k$.  Map the rational two-torsion point $e_0$ through $R_G$, shift its
image to $X=0$, and apply
\eqref{eq:QRT-reembedding-qOmega}--
\eqref{eq:QRT-reembedding-beta-square}.  This gives an explicit
model-preserving $\ell$-isogeny whenever the target realization
equations are soluble over $k$.  There is no universal three-parameter
coefficient formula independent of the kernel: the kernel polynomial
is essential input.  For large $\ell$, fast variants of V\'elu's
method may be used, but that changes the isogeny-evaluation algorithm,
not the model-theoretic re-embedding criterion.

\section{Isogenies in characteristic two}
\label{subsec:QRT-isogenies-char2}

Start from the unshifted binary model
\begin{equation}
  E_{a,b}^{\mathrm{lin}}:
  Y^2+XY=X^3+aX^2+bX,
  \qquad b\ne0.
  \label{eq:QRT-char2-linear-model-isogeny}
\end{equation}
The unique nonzero two-torsion point is $(0,0)$.  The separable
quotient by this point is
\begin{align}
  X_2&=X+\frac bX,
  \label{eq:QRT-char2-2isog-X}\\
  Y_2&=Y\left(1+\frac b{X^2}\right)+\frac bX,
  \label{eq:QRT-char2-2isog-Y}
\end{align}
with codomain
\begin{equation}
  Y_2^2+X_2Y_2=X_2^3+aX_2^2+b.
  \label{eq:QRT-char2-2isog-target-constant}
\end{equation}
Indeed, in characteristic two the difference between the two sides of
\eqref{eq:QRT-char2-2isog-target-constant} is
\[
 \frac{X^4+b^2}{X^4}
 \bigl(Y^2+XY-X^3-aX^2-bX\bigr).
\]
It therefore vanishes on the source curve wherever $X\ne0$.  The
finite pole $X=0$ is the nonzero two-torsion point $(0,0)$; together
with $O$ it is sent to the quotient identity, and the resulting
separable map has kernel $\{O,(0,0)\}$.

Let $r_b^2=b$.  The ordinate shift
$\widetilde Y_2=Y_2+r_b$ changes the codomain to
\begin{equation}
  \widetilde Y_2^2+X_2\widetilde Y_2
  =X_2^3+aX_2^2+r_bX_2,
  \label{eq:QRT-char2-2isog-target-linear}
\end{equation}
which is again of the linear form
\eqref{eq:QRT-char2-Weierstrass}.  For a QRT input,
$b=b_2$ in the notation of
\eqref{eq:QRT-char2-parameters}; hence the target has
\begin{equation}
  a_2'=a_2,
  \qquad
  b_2'=\sqrt{b_2}.
  \label{eq:QRT-char2-2isog-target-parameters}
\end{equation}
The reconstruction in
\eqref{eq:QRT-char2-reconstruction}, or equivalently the coefficient
test \eqref{eq:QRT-char2-realization-quadratic}, decides and constructs
a QRT codomain over the ground field.

\subsection{An explicit binary three-isogeny.}
Let $K=(r,s)\in E_{a,b}^{\mathrm{lin}}(k)$ have exact order three.
Since
\begin{equation}
  -K=(r,r+s),
  \label{eq:QRT-char2-neg-K}
\end{equation}
the kernel is $\{O,K,-K\}$.  The assumption $K\in E(k)$ ensures that
every quantity in the explicit formulas below belongs to the ground
field.  A non-split Galois-stable kernel of order three is handled by
the characteristic-uniform V\'elu sums following the proposition,
or equivalently by eliminating the quadratic kernel ordinate from
those sums.  For a point $P=(X,Y)$ with $X\ne r$, put
\begin{align}
  D_K&=X+r,
  &d_K&=\frac r{D_K},
  &\lambda_K&=\frac{Y+s}{D_K},
  \label{eq:QRT-char2-3isog-precomp}\\
  h_K&=d_K^2+d_K,
  &X_+&=\lambda_K^2+\lambda_K+a+X+r.
  \label{eq:QRT-char2-3isog-precomp-two}
\end{align}
Then the normalized degree-three V\'elu isogeny is
\begin{align}
  X_3&=X+h_K,
  \label{eq:QRT-char2-explicit-3isog-X}\\
  Y_3&=Y+r+d_KX_+
       +(\lambda_K+d_K+1)h_K+Xd_K.
  \label{eq:QRT-char2-explicit-3isog-Y}
\end{align}
It has kernel $\{O,K,-K\}$ and codomain
\begin{equation}
  E_3^{(2)}:
  \qquad
  Y_3^2+X_3Y_3
  =X_3^3+aX_3^2+(b+r)X_3+r.
  \label{eq:QRT-char2-explicit-3isog-target}
\end{equation}
The formulas extend across $X=r$ by sending the two nonzero kernel
points to the identity.

\begin{proof}
For distinct abscissas, the binary addition law gives
\[
  X(P+K)=\lambda_K^2+\lambda_K+a+X+r=X_+.
\]
For $-K=(r,r+s)$, the corresponding slope is
$\lambda_K+d_K$, so
\[
  X(P-K)=X_++d_K^2+d_K=X_++h_K.
\]
The V\'elu abscissa is
\[
  X_3=X+\bigl(X(P+K)-r\bigr)
        +\bigl(X(P-K)-r\bigr).
\]
In characteristic two the two occurrences of $r$ and the two
occurrences of $X_+$ cancel, leaving $X_3=X+h_K$.

For the ordinates, write
\[
  \nu_K=\frac{Xs+rY}{D_K}.
\]
Then
\[
  Y(P+K)=(\lambda_K+1)X_++\nu_K.
\]
For $P-K$, the line intercept is $\nu_K+Xd_K$ and the abscissa is
$X_++h_K$; hence
\[
  Y(P-K)
  =(\lambda_K+d_K+1)(X_++h_K)+\nu_K+Xd_K.
\]
Since $Y(K)+Y(-K)=s+(r+s)=r$, the V\'elu ordinate
\[
  Y_3=Y+\bigl(Y(P+K)-Y(K)\bigr)
       +\bigl(Y(P-K)-Y(-K)\bigr)
\]
reduces exactly to
\eqref{eq:QRT-char2-explicit-3isog-Y}.

It remains to calculate the codomain.  For the generalized Weierstrass
coefficients
$a_1=1$, $a_2=a$, $a_3=0$, $a_4=b$, $a_6=0$, V\'elu's local quantities
at $K$ are
\[
  g_K^x=3r^2+2ar+b-s=r^2+b+s,
  \qquad
  g_K^y=-2s-r=r,
\]
so
\[
  t_K=2g_K^x-a_1g_K^y=r,
  \qquad
  u_K=(g_K^y)^2=r^2,
  \qquad
  w_K=u_K+rt_K=0.
\]
The quotient coefficients are therefore
\[
  a_4'=a_4-5t_K=b+r,
  \qquad
  a_6'=a_6-(a_1^2+4a_2)t_K-7w_K=r,
\]
while $a_1,a_2,a_3$ are unchanged.  This proves
\eqref{eq:QRT-char2-explicit-3isog-target} and completes the proof.
\end{proof}

Over a perfect field, let $t_3^2=r$ and set
\begin{equation}
  \widetilde Y_3=Y_3+t_3.
  \label{eq:QRT-char2-3isog-target-shift}
\end{equation}
Then the codomain becomes
\begin{equation}
  \widetilde Y_3^2+X_3\widetilde Y_3
  =X_3^3+aX_3^2+b_3^{\mathrm{lin}}X_3,
  \qquad
  b_3^{\mathrm{lin}}=b+r+t_3.
  \label{eq:QRT-char2-3isog-target-linear}
\end{equation}
The coefficient $b_3^{\mathrm{lin}}$ is nonzero because the quotient
curve is smooth.  Thus Theorem~\ref{thm:QRT-char2-abstract-classification}
and the realization test
\eqref{eq:QRT-char2-realization-quadratic} apply directly to the
three-isogenous codomain.  This ordinate shift is essential: the
constant term $r$ in
\eqref{eq:QRT-char2-explicit-3isog-target} must not be silently treated
as though it were already in the QRT linear normal form.

For a separable odd-order subgroup $G$ in characteristic two, the
characteristic-uniform V\'elu map can be written directly on a general
Weierstrass model as
\begin{align}
  X_G(P)&=X(P)+
  \sum_{Q\in G\setminus\{O\}}
  \bigl(X(P+Q)-X(Q)\bigr),
  \label{eq:QRT-char2-Velu-X}\\
  Y_G(P)&=Y(P)+
  \sum_{Q\in G\setminus\{O\}}
  \bigl(Y(P+Q)-Y(Q)\bigr).
  \label{eq:QRT-char2-Velu-Y}
\end{align}
The binary addition law makes these rational functions explicit.  For
$\ell=3$, take $G=\{O,K,-K\}$; for general odd $\ell$, use the full
kernel or its kernel polynomial.  After computing the generalized
Weierstrass coefficients of the quotient, normalize $a_1=1$ and apply
Theorem~\ref{thm:QRT-char2-abstract-classification} and
\eqref{eq:QRT-char2-realization-quadratic}.  This is the correct
characteristic-two analogue of the odd-characteristic model-preserving
pipeline; the short-Weierstrass formulas
\eqref{eq:QRT-explicit-3isogeny}--
\eqref{eq:QRT-explicit-l-isogeny} must not be specialized by setting
$2=0$ or $3=0$.

\section[QRT dynamics on product families]{Transporting QRT dynamics to \texorpdfstring{$C_d$ and $\C_{a,b,d}$}{Cd and Cabd}, and reciprocal curves}
\label{subsec:QRT-transport-other-families}

The QRT mechanism is not tied to the affine equation
\eqref{eq:symmetric-QRT}.  Theorem~\ref{thm:QRT-adjacent-Kummer-state}
shows that it is the state-space form of a degree-two Kummer coordinate
and a fixed elliptic translation.

\begin{corollary}[Translation states on the $C$-curve families]
\label{cor:QRT-from-translation}
Let $E/k$ be the smooth completion of a member of $C_d$,
$\T_{a,d}$, $\C_{a,b,d}$, or the reciprocal $C$-curve family.  Fix a
$k$-rational origin, a separable even degree-two function
$\kappa:E\to\PP^1$, and a point $D\in E(k)$ with $2D\ne O$.  Then
\begin{equation}
  P\longmapsto
  \bigl(\kappa(P),\kappa(P+D)\bigr)
  \label{eq:QRT-C-family-adjacent-state}
\end{equation}
identifies $E$ with a smooth symmetric $(2,2)$ state curve, and
translation by $D$ is its Vieta--McMillan shift.  In odd characteristic,
whenever the relevant common involution on the two projective coordinates
has split fixed divisor over $k$, a common M\"obius transformation places
this state curve in the normal form $\Q_{\alpha,\beta,\gamma}$; the
nonsplit fixed divisor gives its twisted symmetric form.  Over a perfect
field of characteristic two, every ordinary pointed member has the explicit
normal form $\Q_{d^2,d,b^2}$ of
Theorem~\ref{thm:QRT-binary-pointed-state-model}.
\end{corollary}

\begin{proof}
Apply Theorem~\ref{thm:QRT-adjacent-Kummer-state} to the chosen elliptic
curve, Kummer function, and translation point.  In odd characteristic,
Theorem~\ref{thm:QRT-ground-field-even-normal-form} gives the stated
fixed-divisor criterion and the common-involution normalization.  In
characteristic two, Theorem~\ref{thm:QRT-binary-pointed-state-model} gives
the displayed ordinary normal form over the ground field.
\end{proof}

The corollary distinguishes the native coordinates of a $C$-curve from
a translation-adapted state coordinate.  A native projection can be
used when its associated translation is the desired one.  For a
specified point of larger order, one instead selects an even
Kummer function and forms \eqref{eq:QRT-C-family-adjacent-state}.  Over
an infinite field the point $D$ may have infinite order; over a finite
field it has finite order, and the order of the resulting QRT shift is
exactly the order of $D$ by
\eqref{eq:QRT-state-torsion-period}.  State doubling and logarithmic
scalar multiplication are then given by
Theorem~\ref{thm:QRT-state-doubling-ladder}, with the model-specific
Kummer formulas supplying the actual operation count.

The native coordinate involutions of the distinguished slices behave
more rigidly.  For centered $C_d$ in odd characteristic,
$\beta=0$, and
\begin{equation}
  \mathcal M_Q(r,s)=(s,-r).
  \label{eq:Cd-native-McMillan}
\end{equation}
It has order four, in agreement with
\eqref{eq:QRT-McMillan-order-four}.  For the normalized reciprocal
slice, $\gamma=\alpha^2$, and the curve equation reduces the product
form to
\begin{equation}
  \mathcal M_Q(x,y)=(y,-1/x),
  \label{eq:R-native-McMillan}
\end{equation}
which again has order four.  Thus the native coordinate QRT maps on these two highly symmetric
slices are torsion translations.  To obtain
non-torsion dynamics on the same abstract curve, one must choose a
non-torsion increment and generally a different degree-two coordinate,
as in Corollary~\ref{cor:QRT-from-translation}.

For $\C_{a,b,d}$, the two product involutions arise explicitly as
the deck transformations of the two degree-two projections.  Their composition is a translation whose
increment is the difference of the two reflection centres in
$\Pic^0$.  On special split-boundary subfamilies this divisor class is
torsion.  Across the full three-parameter family, its order is computed from
the corresponding class in $\Pic^0$.  The general translation construction
therefore separates the special boundary symmetries from the broader QRT
principle: every elliptic translation has a symmetric-biquadratic state
encoding, while the original product coordinates display the translation
attached to their two native projections.

From an algorithmic viewpoint, this gives two strategies.

\begin{enumerate}[label=(\roman*)]
  \item Use the native product or reciprocal coordinates when their
  exceptionally simple fixed-step maps are the desired operation.

  \item For large scalar multiplication, construct the Kummer line of
  the underlying Weierstrass/Jacobi model and use a differential
  ladder.  A nonnative QRT coordinate may be selected when it lowers
  the fixed-step cost or supports parallel ECM arithmetic.
\end{enumerate}

The QRT representation supplies complementary arithmetic layers: direct
fixed-step recurrence, full addition, Kummer ladders, point recovery, and
model conversion.  Each layer has its own exact operation count, and together
they form the complete state interface.

\section{The distinguished slices and their interpretation}
\label{subsec:QRT-slices}

In odd characteristic, centered coordinates
$r=2u+1$, $s=2v+1$ place $C_d$ in the QRT family:
\begin{equation}
  (r^2-1)(s^2-1)=16d
  \quad\Longleftrightarrow\quad
  r^2s^2-r^2-s^2+1-16d=0.
  \label{eq:Cd-QRT-slice}
\end{equation}
Thus
\begin{equation}
  C_d:
  \qquad
  (\alpha,\beta,\gamma)=(-1,0,1-16d).
  \label{eq:Cd-QRT-parameters}
\end{equation}
The root exchanges reduce to the two split sign changes
$r\mapsto-r$ and $s\mapsto-s$.  This embedding is intrinsically an
odd-characteristic centered description; it must not be used in
characteristic two, where the coordinate $2u+1$ degenerates.

The normalized reciprocal model is
\begin{equation}
  (x^2-1)(y^2-1)=\kappa xy
  \quad\Longleftrightarrow\quad
  x^2y^2-x^2-y^2-\kappa xy+1=0,
  \label{eq:R-QRT-slice}
\end{equation}
so
\begin{equation}
  \R_\kappa:
  \qquad
  (\alpha,\beta,\gamma)=(-1,-\kappa,1).
  \label{eq:R-QRT-parameters}
\end{equation}
The QRT root exchanges become
\begin{equation}
  (x,y)\longmapsto(-1/x,y),
  \qquad
  (x,y)\longmapsto(x,-1/y),
  \label{eq:R-QRT-involutions}
\end{equation}
exactly the reciprocal involutions.  In characteristic two the minus
signs disappear and these maps become the wild reciprocal involutions
described by the characteristic-two reciprocal arithmetic in the
standalone toolkit.

For the general reciprocal curve, expansion gives
\begin{equation}
  x^2y^2-\sigma x^2-\tau y^2-\kappa xy+\tau\sigma=0.
  \label{eq:R-general-expanded}
\end{equation}
After adjoining square roots of $\tau$ and $\sigma$ and scaling the two
coordinates separately, the coefficients of $x^2$ and $y^2$ become
equal, producing the symmetric reciprocal slice.  Over the ground
field, the failure of this scaling is precisely the nonsplit reciprocal
twist.  Consequently, the reciprocal family simultaneously connects:
\begin{enumerate}[label=(\roman*)]
  \item the Jacobi quartic obtained from its even-quartic reduction;
  \item the symmetric QRT root-exchange dynamics;
  \item the split or nonsplit four-torsion normal form determined by the
  order-four automorphism;
  \item the reciprocal twist of $C_d$ supplied by the Cayley
  transformation.
\end{enumerate}
These are four descriptions of the same genus-one geometry, not four
independent claims of isomorphism over every ground field.  The field
extensions and rational-point hypotheses have been stated explicitly in
the corresponding constructions.

\begin{center}
\small
\setlength{\tabcolsep}{3pt}
\renewcommand{\arraystretch}{1.24}
\begin{tabular}{@{}>{\raggedright\arraybackslash}p{0.16\textwidth}
                >{\raggedright\arraybackslash}p{0.21\textwidth}
                >{\raggedright\arraybackslash}p{0.26\textwidth}
                >{\raggedright\arraybackslash}p{0.27\textwidth}@{}}
\toprule
Model & QRT parameters & Root exchanges & Arithmetic interpretation\\
\midrule
Centered $C_d$ (odd characteristic)
& $(-1,0,1-16d)$
& $(r,s)\mapsto(-r,s)$ and $(r,s)\mapsto(r,-s)$
& Split reflections; one-sided twists give the twisted-Edwards and
Montgomery descriptions developed earlier.\\
\addlinespace
Normalized reciprocal model
& $(-1,-\kappa,1)$
& $(x,y)\mapsto(-1/x,y)$ and $(x,y)\mapsto(x,-1/y)$
& Jacobi full-point formulas, a rational four-torsion structure, and a
Montgomery Kummer ladder.\\
\addlinespace
General smooth QRT member
& $\alpha\gamma\Omega_Q\ne0$ in odd characteristic;
  $\alpha\beta\gamma\ne0$ in characteristic two
& Rational root exchanges
\eqref{eq:QRT-vertical-involution}--
\eqref{eq:QRT-horizontal-involution}
& Weierstrass/Jacobi arithmetic in odd characteristic and
L\'opez--Dahab/Kummer arithmetic in characteristic two.\\
\bottomrule
\end{tabular}
\end{center}

\begin{figure}[H]
\centering
\begin{tikzpicture}[x=1cm,y=1cm,>=Latex, box/.style={draw,rounded corners,align=center,inner sep=4pt,font=\small,text width=3.2cm}]
  \node[box] (P) at (0,2.8) {$[n]D$\\on the elliptic curve};
  \node[box] (EDS) at (-4.2,0.7) {$W_n$\\EDS / division values};
  \node[box] (Lucas) at (4.2,0.7) {$(v_n,w_n)$\\elliptic Lucas data};
  \node[box] (State) at (0,-1.2) {$(v_n,v_{n+1})$\\adjacent QRT state};
  \node[box] (Idx) at (0,-3.5) {$n\mapsto n+1,\ 2n,\ 2n+1$};
  \draw[->,thick] (P) -- node[above left,font=\scriptsize,fill=white,inner sep=1pt]{division data} (EDS);
  \draw[->,thick] (P) -- node[above right,font=\scriptsize,fill=white,inner sep=1pt]{Kummer / Lucas data} (Lucas);
  \draw[->,thick] (EDS) -- node[below left,font=\scriptsize,fill=white,inner sep=1pt]{EDS bridge} (State);
  \draw[->,thick] (Lucas) -- node[below right,font=\scriptsize,fill=white,inner sep=1pt]{state extraction} (State);
  \draw[->,thick] (State) -- node[right,font=\scriptsize,fill=white,inner sep=1pt]{state ladder} (Idx);
  \draw[->,thick,bend left=20] (Lucas) to node[above,font=\scriptsize,fill=white,inner sep=1pt]{same point sequence} (EDS);
\end{tikzpicture}
\caption{The Chapter~\ref{ch:symmetric-QRT-envelope} closed loop among point multiples, elliptic divisibility sequences, elliptic Lucas data, and adjacent QRT states.  The same sequence of multiples $[n]D$ can therefore be read through three different but explicitly bridged families of recurrence data.}
\label{fig:QRT-Lucas-EDS-closed-loop}
\end{figure}
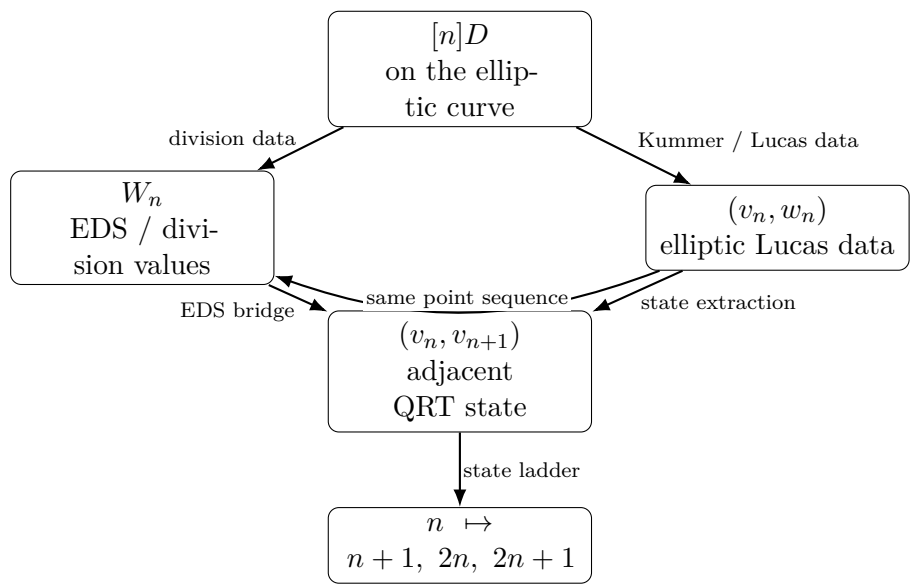

Figure~\ref{fig:QRT-Lucas-EDS-closed-loop} collects the bridges developed in
this chapter into a single closed-loop picture.

\section[Cryptographic prospects for QRT states]{Cryptographic prospects for symmetric biquadratic state models}
\label{subsec:QRT-cryptographic-prospects}

From a cryptographic viewpoint, the significance of
\(\Q_{\alpha,\beta,\gamma}\) is not that it supplies one more equation
birational to a Weierstrass curve.  Its distinctive feature is that a fixed
elliptic displacement is incorporated into a smooth \((2,2)\) state curve:
the two coordinates represent
\[
  \bigl(\kappa(P),\kappa(P+D)\bigr),
\]
and therefore carry both a Kummer value and the fixed difference required by
differential addition.  This turns scalar multiplication into arithmetic on
one constrained state space, makes the two ladder branches conjugate by
coordinate exchange, and permits the difference coordinate to be compiled
into the curve parameters.  These structural facts identify several lines of
research whose cryptographic value is more substantial than the mere
presentation of a new curve equation.

\paragraph{Native constant-time scalar multiplication.}
The central implementation problem is to determine the best complete state-doubling circuit in
coordinates intrinsic to \(\Q_{\alpha,\beta,\gamma}\).  On the split
odd-characteristic locus, the present formulas transport the Montgomery
\texttt{xDBLADD} core with the same numbers of general multiplications and
squarings.  This establishes a rigorous reference circuit from which native
formulas can reduce fixed multiplications, critical-path depth, live
registers, memory traffic, or conversion overhead.  End-to-end evaluation
includes complete formulas, conditional swaps, normalization, validation,
and final recovery in both software and hardware.

\paragraph{Full-point arithmetic and multi-scalar multiplication.}
An adjacent state is a native fixed-displacement arithmetic engine centered
at the marked point \(D\).  Short complete formulas for the transported full
group law, mixed addition, and variable-difference addition provide the
foundation for extending this engine to multi-base computations while
retaining the symmetry of the state model.  Generic multi-scalar
multiplication can be organized through two or more compatible state curves
whose fixed displacements correspond to the selected base points and whose
arithmetic circuits share intermediate products.  This framework creates
natural QRT-state realizations of Shamir, Straus, Pippenger, and
endomorphism-assisted multi-scalar multiplication.

\paragraph{Complete encodings, validation, and side-channel resistance.}
The state equation supplies an intrinsic algebraic membership test for every
encoded adjacent state, and the complete projective update gives a uniform
ladder circuit across all valid projective inputs.  Protocol deployment
requires a canonical state encoding, precise decoding and subgroup rules,
curve-and-twist classification, exceptional-state validation, constant-time
full-point recovery, and a proof that every accepted byte string has a unique
protocol interpretation.  These components combine naturally into an
end-to-end constant-pattern interface.  Fault injection, exceptional-state
propagation, timing behavior, cache access, and power consumption can then be
analyzed within the same validated state-transition system.

\paragraph{Ordinary binary curves.}
Over every perfect field of characteristic two, an ordinary pointed elliptic
curve admits the explicit adjacent-state member
\[
  \Q_{d^2,d,b^2}:
  \qquad
  x^2y^2+d^2(x^2+y^2)+dxy+b^2=0.
\]
The associated companion data retain the full Artin--Schreier twist
information and recover the original pointed elliptic curve from its state
model.  This universal ordinary-binary construction is particularly well
suited to normal-basis arithmetic, bit-sliced software, and compact hardware
in which squaring is inexpensive.  Its implementation program comprises a
uniform supersingular companion theory, basis-specific operation schedules,
complete encoding and recovery, and security evaluation alongside binary
Kummer and L\'opez--Dahab implementations.

\paragraph{Recovery and protocol interfaces.}
The adjacent coordinate \(\kappa(P+D)\), together with an oriented lift of
\(D\), restores the sign information omitted by the single Kummer coordinate
\(\kappa(P)\).  This structure supports ladders that terminate in a full
point, protocols requiring a canonical sign, and implementations in which
state validation and recovery reuse final-loop intermediates.  End-to-end
benchmarks can measure the complete cost of the additional state coordinate,
projective normalization, recovery, and serialization.  Quotient-only
Diffie--Hellman interfaces and full-point signature, proof, and pairing
interfaces thereby receive separate optimized schedules within the same
state framework.

\paragraph{Isogenies and model transport.}
An isogeny acts on an adjacent state by applying its induced Kummer function
to both state coordinates.  Low-degree model-preserving isogenies therefore
support efficient curve changes, cyclic quotients, cofactor management, and
families of compatible state parameters.  Protocol deployment specifies the
nondegeneracy of the image displacement, the ground-field descent of the
target normalization, the induced subgroup map, and the behavior of all
exceptional fibres.  Coordinatewise transport supplies the arithmetic layer
for these constructions.  Post-quantum security, when sought, is supplied by
the protocol-level isogeny assumption together with its parameter,
subgroup, encoding, and implementation analysis.

\paragraph{ECM, recurrence arithmetic, and batch computation.}
The QRT shift, elliptic Lucas data, elliptic divisibility sequences, and
state-division polynomials furnish applications beyond conventional
discrete-logarithm protocols.  They support parallel and batched stages for
elliptic-curve factorization, exact period and torsion tests, and recurrence
computations in which a fixed translation is reused many times.  In these
applications, small constants, inexpensive curve generation, smooth group
orders, product-tree compatibility, batch inversion, and parallel depth form
the principal selection criteria.  Dedicated parameterizations and
smoothness analyses can therefore exploit the factorized state equation and
its compiled translation data directly.

\paragraph{Security inheritance and parameter specification.}
The security of a
\(\Q_{\alpha,\beta,\gamma}\) instance is inherited from its underlying
elliptic-curve group and completed by its encoding, validation rules, subgroup
policy, and constant-time implementation.  Birational re-embedding preserves
the underlying elliptic discrete-logarithm problem, while the state geometry
adds an intrinsic membership equation, complete state transitions,
fixed-displacement compilation, and explicit recovery data.  A cryptographic
parameter set therefore specifies the group order and cofactor, twist order,
exceptional loci, available endomorphisms, subgroup checks, canonical
encoding, recovery procedure, and constant-pattern arithmetic interface.
These data make the geometric and computational advantages of each
\(\Q_{\alpha,\beta,\gamma}\) instance explicit and independently
verifiable.

\begin{theorem}[Arithmetic-state synthesis]
\label{thm:QRT-arithmetic-state-synthesis}
Let \(k\) be a field, let \((E,O,D)\) be a pointed elliptic curve over
\(k\), let
\(\kappa:E\to\PP^1\) be a separable degree-two Kummer quotient with
generic fibres \(\{P,-P\}\), and assume \(2D\ne O\).  Define
\[
  \mathcal S_D(P)=\bigl(\kappa(P),\kappa(P+D)\bigr).
\]
Then the following assertions hold.

\begin{enumerate}[label=(\roman*)]
  \item The map \(\mathcal S_D\) is a \(k\)-isomorphism from \(E\) onto a
  smooth symmetric curve \(\mathcal B_D\) of bidegree \((2,2)\) in
  \(\PP^1\times\PP^1\).  Under this isomorphism, translation by \(D\)
  becomes the adjacent-state shift \(T_D\), and multiplication by two
  becomes
  \[
    \Delta_D:\mathcal S_D(P)\longmapsto
    \bigl(\kappa(2P),\kappa(2P+D)\bigr).
  \]

  \item Suppose \(\charac(k)\ne2\).  One common
  \(\operatorname{PGL}_2(k)\) change in the two Kummer coordinates puts
  \(\mathcal B_D\) into the exact even form
  \[
    X^2Y^2+\alpha(X^2+Y^2)+\beta XY+\gamma=0
  \]
  if and only if \(E\) has a nonzero \(k\)-rational two-torsion point
  whose induced involution on the Kummer line has split fixed divisor.
  Equivalently, in a completed-square model \(v^2=f(u)\), there is a
  rational root \(r\) of \(f\) for which \(f'(r)\in k^{\times2}\).

  \item On this split odd-characteristic locus, move \((r,0)\) to
  \((0,0)\), write
  \[
    E_{A,B}:v^2=u^3+Au^2+Bu,\qquad B=c^2,\quad c\in k^\times,
  \]
  and use the common coordinate
  \(z=(u-c)/(u+c)\).  Each of the two scalar-ladder branches
  \(\Delta_D\) and \(T_D\Delta_D\) is represented by one complete
  projective formula.  It is an addition-only
  \(\operatorname{PGL}_2(k)\) conjugate of fixed-difference Montgomery
  \texttt{xDBLADD}; with its two coefficient multiplications compiled,
  its loop cost is
  \[
    \boxed{4\M+4\Sqr+2\Dpar}.
  \]

  \item Suppose that \(k\) is perfect of characteristic two and that
  \(E\) is ordinary.  Then \(E\) has a \(k\)-model
  \[
    E_{a,b}^{(2)}:\quad v^2+uv=u^3+au^2+b^2,\qquad b\ne0,
  \]
  in which \(D=(d,e)\) has \(d\ne0\).  After a common
  \(\operatorname{PGL}_2(k)\) change of the given Kummer coordinate,
  take \(\kappa=u\); then the adjacent state is
  \[
    \Q_{d^2,d,b^2}:\quad
    x^2y^2+d^2(x^2+y^2)+dxy+b^2=0.
  \]
  Translation by \(D\) is
  \[
    T_D(x,y)=\left(y,x+\frac{dy}{y^2+d^2}\right),
  \]
  and each scalar-ladder branch has one complete projective formula of
  compiled cost
  \[
    \boxed{4\M+5\Sqr+3\Dpar}.
  \]
  The omitted coefficient \(a\) is retained by
  \[
    d^2+\frac{b^2}{d^2}
    =a+\left(d+\frac ed\right)^2+\left(d+\frac ed\right).
  \]

  \item In the explicit odd and ordinary-binary models of parts (iii)
  and (iv), the orientation of \(D\) gives a complete three-chart
  recovery atlas: the two boundary charts return \(O\) and \(-D\), and
  the affine chart recovers the missing ordinate.  Thus the adjacent
  state retains enough information for full-point recovery once an
  oriented lift of \(D\) has been fixed.

  \item Let \(\varphi:E\to E'\) be an isogeny, let
  \(\kappa':E'\to\PP^1\) be a separable degree-two Kummer quotient, put
  \(D'=\varphi(D)\), and assume \(2D'\ne O'\).  Let
  \(\overline\varphi:\PP^1\to\PP^1\) be the unique induced Kummer map
  satisfying
  \(\overline\varphi\circ\kappa=\kappa'\circ\varphi\).  If
  \(\mathcal S'_{D'}(P')=(\kappa'(P'),\kappa'(P'+D'))\), then the state
  isogeny
  \(\Phi_\varphi=\mathcal S'_{D'}\circ\varphi\circ\mathcal S_D^{-1}\)
  is coordinatewise:
  \[
    \Phi_\varphi(x,y)
    =\bigl(\overline\varphi(x),\overline\varphi(y)\bigr).
  \]
  It has degree \(\deg\varphi\) and commutes with both the state shift and
  every transported multiplication map \(\Delta_{D,m}\).
\end{enumerate}
\end{theorem}

\begin{proof}
\smallskip
\noindent\emph{Step 1: the state curve.}
We first verify the assertions after base change to an algebraic closure;
because every map used below is defined over \(k\), the resulting identities
and morphisms descend to \(k\).  Suppose that two generic points
\(P,R\in E\) have the same adjacent
state.  Equality of their first Kummer coordinates gives \(R=P\) or
\(R=-P\).  In the second case, equality of the second coordinates gives
\[
  \kappa(D-P)=\kappa(D+P).
\]
The generic fibre description of \(\kappa\) then gives either
\(D-P=D+P\), which is possible only on the finite set \(E[2]\), or
\(D-P=-D-P\), which would imply \(2D=O\).  The latter is excluded.
Hence \(\mathcal S_D\) is generically injective and is birational onto
its integral image.  Its two projections are \(\kappa\) and
\(\kappa\circ\tau_D\), both of degree two, so the image has bidegree
\((2,2)\).  Such a curve has arithmetic genus
\((2-1)(2-1)=1\).  The morphism
\(E\to\mathcal B_D\) is finite because \(E\) is proper and its fibres
are finite.  Since it is also birational and \(E\) is normal, it is the
normalization map.  Therefore the genus formula
\[
 p_a(\mathcal B_D)-g(E)
 =\sum_{Q\in\mathcal B_D(\overline k)}\delta_Q
\]
has left-hand side zero.  Each \(\delta_Q\) is a nonnegative integer, so
every one vanishes.  The image is therefore smooth,
and the birational morphism \(\mathcal S_D\) is an isomorphism.

The evenness of the Kummer quotient gives
\[
  \mathcal S_D(-P-D)
  =\bigl(\kappa(P+D),\kappa(P)\bigr),
\]
which proves symmetry under coordinate exchange.  Direct substitution gives
\[
  \mathcal S_D(P+D)
   =\bigl(\kappa(P+D),\kappa(P+2D)\bigr)
\]
and
\[
  \mathcal S_D(2P)
   =\bigl(\kappa(2P),\kappa(2P+D)\bigr).
\]
These identities prove part (i), including the interpretations of \(T_D\)
and \(\Delta_D\).

\smallskip
\noindent\emph{Step 2: odd-characteristic descent.}
Assume first that one common Kummer coordinate gives an exact even equation.
The simultaneous sign change \((X,Y)\mapsto(-X,-Y)\) preserves that
equation.  On a smooth member it has no fixed point: its only fixed points
in \(\PP^1\times\PP^1\) are the four sign corners.  In homogeneous
coordinates the equation is
\[
 X_1^2Y_1^2
 +\alpha(X_1^2Y_0^2+X_0^2Y_1^2)
 +\beta X_0X_1Y_0Y_1
 +\gamma X_0^2Y_0^2=0.
\]
Its values at \((\infty,\infty)\), \((\infty,0)\), \((0,\infty)\), and
\((0,0)\) are respectively \(1,\alpha,\alpha,\gamma\).  Smoothness of an
even member forces \(\alpha\gamma\ne0\), so none of the four corners lies
on the curve.  Pull this involution back
to \(E\).  An elliptic-curve automorphism has the form
\(P\mapsto a(P)+Q\), with \(a(O)=O\).  If \(a\ne1\), the nonzero isogeny
\(1-a\) is surjective over \(\overline k\), so
\((1-a)P=Q\) has a solution and the automorphism has a fixed point.  Hence
\(a=1\), and the pulled-back fixed-point-free involution is translation by
a nonzero point \(T_2\in E(k)[2]\).  Its action on the Kummer line is
conjugate over \(k\) to \(z\mapsto-z\), whose two fixed points are rational;
the induced Kummer involution is therefore split.

Conversely, suppose a nonzero \(T_2\in E(k)[2]\) has split Kummer action.
Choose \(\phi\in\operatorname{PGL}_2(k)\) that conjugates that action to
\(z\mapsto-z\), and use \(\phi\circ\kappa\) in both state coordinates.
Translation by \(T_2\) then acts as \((X,Y)\mapsto(-X,-Y)\), while
\(P\mapsto-P-D\) exchanges \(X\) and \(Y\).  The irreducible bidegree
\((2,2)\) equation is consequently an eigenvector for both operations.
If it were anti-invariant under exchange, it would vanish on \(X=Y\) and
would contain \(X-Y\).  If it were symmetric under exchange but
anti-invariant under simultaneous sign, its affine part would have the
form
\[
 a(X^2Y+XY^2)+b(X+Y)=(X+Y)(aXY+b).
\]
Either case contradicts irreducibility.  The equation is therefore
invariant under both operations.  The invariant symmetric monomials are
exactly
\[
  X^2Y^2,\qquad X^2+Y^2,\qquad XY,\qquad 1.
\]
The coefficient of \(X^2Y^2\) cannot vanish.  Indeed, if the state curve
passed through \((\infty,\infty)\), then for some \(P\) one would have
\(\kappa(P)=\kappa(P+D)=\infty\).  Since \(\infty\) is a fixed value of
the \(T_2\)-action, this would give simultaneously \(2P=-D\) and
\(2P=-T_2\), hence \(D=T_2\), contrary to \(2D\ne O\).  Dividing by the
nonzero leading coefficient gives the asserted exact even equation.

To obtain the cubic criterion, write \(T_2=(r,0)\) on
\(v^2=f(u)\).  If
\(f(u)=(u-r)(u-r_2)(u-r_3)\), the chord formula gives the induced
Kummer involution
\[
  u\longmapsto r+\frac{(r-r_2)(r-r_3)}{u-r}
   =r+\frac{f'(r)}{u-r}.
\]
Its fixed equation is \((u-r)^2=f'(r)\).  Thus its fixed divisor is split
over \(k\) exactly when \(f'(r)\in k^{\times2}\), proving part (ii).

\smallskip
\noindent\emph{Step 3: the split odd-characteristic update.}
Translate \(r\) to zero.  The cubic becomes
\[
  v^2=u^3+Au^2+Bu,\qquad B=f'(r)=c^2.
\]
Writing the transformed marked point as \(D=(d,\varepsilon)\), one has
\(d\ne0\); otherwise \(D=(0,0)=T_2\), contrary to \(2D\ne O\).  Put
\(U=u/c\), \(a=A/c\), \(\delta_D=d/c\), and
\(A_{24}=(a+2)/4\).  Kummer doubling and fixed-difference addition give
\[
 U(2P)=\frac{(U^2-1)^2}{4U(U^2+aU+1)},
 \qquad
 U(2P+D)=
 \frac{(U_0U_1-1)^2}{\delta_D(U_0-U_1)^2},
\]
where \(U_0=U(P)\) and \(U_1=U(P+D)\).  These are precisely the
Montgomery \texttt{xDBL} and fixed-difference \texttt{xADD} identities.
The change
\[
 z=\frac{U-1}{U+1},\qquad
 (X:Z)=(R+S:S-R)\quad\text{for }z=(R:S),
\]
and its inverse use only additions and subtractions.  Therefore conjugating
the Montgomery four-form tuple does not add a general multiplication or a
squaring.

For completeness, consider first its doubling output.  Up to a nonzero
scalar the output pair is
\[
 \bigl((X^2-Z^2)^2:
 4XZ(X^2+aXZ+Z^2)\bigr).
\]
A common zero with \(XZ\ne0\) would have \(X/Z=\pm1\) and then
\(a=\mp2\); in either case \(U(U^2+aU+1)\) has a repeated root, so the
curve is singular.  If \(X=0\) or \(Z=0\), the first coordinate is
nonzero.  Thus the doubling pair is base-point-free on every smooth member.
The differential output is, again up to nonzero scalars,
\[
 \bigl((X_0X_1-Z_0Z_1)^2:
 \delta_D(X_0Z_1-Z_0X_1)^2\bigr).
\]
A common zero cannot contain a zero projective input coordinate.  If
\(Z_0=0\), then \(X_0\ne0\), and the two bilinear forms reduce to
\(X_0X_1\) and \(X_0Z_1\), forcing \(X_1=Z_1=0\).  If \(X_0=0\), then
\(Z_0\ne0\), and they reduce to \(-Z_0Z_1\) and \(-Z_0X_1\), again
forcing \(X_1=Z_1=0\).  If \(Z_1=0\), then \(X_1\ne0\), and their
vanishing forces \(X_0=Z_0=0\); if \(X_1=0\), then \(Z_1\ne0\), and
their vanishing has the same conclusion.  Each conclusion contradicts
the projectivity of one input.  In the affine chart, simultaneous
vanishing would imply
\(U_0U_1=1\) and \(U_0=U_1\), hence
\(U_0=U_1=\pm1\).  Equality of adjacent Kummer coordinates gives
\(2P=-D\), whereas \(U=\pm1\) is fixed by the Kummer involution induced by
\(T_2\) and gives \(2P=-T_2\).  Hence \(D=T_2\), again contradicting
\(2D\ne O\).  The full tuple is therefore base-point-free.

Using the intermediate names
\[
 \begin{gathered}
 AA=(X_0+Z_0)^2,\qquad BB=(X_0-Z_0)^2,\qquad E_0=AA-BB,\\
 DA=(X_1-Z_1)(X_0+Z_0),\qquad
 CB=(X_1+Z_1)(X_0-Z_0),
 \end{gathered}
\]
the products \(DA\), \(CB\), \(AA\,BB\), and
\(E_0(BB+A_{24}E_0)\) in
\eqref{eq:QRT-state-Montgomery-xDBL}--
\eqref{eq:QRT-state-Montgomery-xADD} account for \(4\M\); the two input
squares and the two differential-output squares account for \(4\Sqr\);
and multiplication by \(A_{24}\) and \(\delta_D\) accounts for
\(2\Dpar\).  Finally,
\(T_D\Delta_D=\sigma\Delta_D\sigma\), where
\(\sigma(X,Y)=(Y,X)\): on the elliptic parameter, the right-hand side sends
\(P\) successively to \(-P-D\), \(-2P-2D\), and \(2P+D\), which is also
the image under the left-hand side.  Thus the same complete tuple gives
both ladder branches and proves part (iii).

\smallskip
\noindent\emph{Step 4: ordinary characteristic two.}
An ordinary elliptic curve in characteristic two has a \(k\)-model
\[
  v^2+uv=u^3+au^2+B,\qquad B\ne0.
\]
Because \(k\) is perfect, Frobenius is surjective, so \(B=b^2\) for some
\(b\in k^\times\).  If \(D=(d,e)\) had \(d=0\), then \(e^2=b^2\), hence
\(D=(0,b)\), the unique nonzero two-torsion point; this contradicts
\(2D\ne O\).  Thus \(d\ne0\).

Let \(x=u(P)\), \(y=u(P+D)\), and \(z=u(P-D)\).  The binary Kummer
identities give
\[
 yz=\frac{(xd+b)^2}{(x+d)^2},
 \qquad
 y+z=\frac{xd}{(x+d)^2}.
\]
Since \(y\) is a root of \(T^2+(y+z)T+yz\), clearing the denominator and
expanding squares gives
\[
  x^2y^2+d^2(x^2+y^2)+dxy+b^2=0.
\]
Vieta exchange in the next fibre gives
\(T_D(x,y)=(y,x+dy/(y^2+d^2))\).  Binary Kummer doubling and
fixed-difference addition give
\[
  \Delta_D(x,y)=
  \left(\frac{(x^2+b)^2}{x^2},
        \frac{(xy+b)^2}{d(x+y)^2}\right),
\]
whose homogeneous output is
\[
 \bigl((X_0^2+bZ_0^2)^2:X_0^2Z_0^2\bigr),
 \qquad
 \bigl((X_0X_1+bZ_0Z_1)^2:
       d(X_0Z_1+Z_0X_1)^2\bigr).
\]
The first pair cannot vanish simultaneously because \(b\ne0\).  Write the
two unsquared forms in the second pair as
\[
 L=X_0X_1+bZ_0Z_1,\qquad R=X_0Z_1+Z_0X_1.
\]
If \(Z_0=0\), then \(X_0\ne0\), and \(L=R=0\) forces
\(X_1=Z_1=0\).  If \(X_0=0\), then \(Z_0\ne0\); since \(b\ne0\),
\(L=R=0\) again forces \(Z_1=X_1=0\).  If \(Z_1=0\), then
\(X_1\ne0\), and the two equations force \(X_0=Z_0=0\).  Finally, if
\(X_1=0\), then \(Z_1\ne0\), and they again force \(Z_0=X_0=0\).
Thus a common zero is affine and would satisfy
\(xy+b=0\) and \(x+y=0\), hence \(x=y\) and \(x^2=b\).  The first equality
\(x=y\) gives \(2P=-D\); the equation \(x^2=b\) is the fixed equation for
translation by the unique two-torsion point \(T_2=(0,b)\), and gives
\(2P=-T_2\).  Hence \(D=T_2\), a contradiction.  The tuple is complete.
For the cross term, compute
\[
 P_0=X_0X_1,\qquad Q_0=Z_0Z_1,\qquad
 R_0=(X_0+Z_0)(X_1+Z_1)-P_0-Q_0.
\]
Then \(R_0=X_0Z_1+Z_0X_1\).  These are three general products, and
\(X_0^2Z_0^2\) is the fourth.  The two input-coordinate squares and the
three outer squares give \(5\Sqr\); the two products by \(b\) and the one
by \(d\) give \(3\Dpar\).  Coordinate exchange gives the other branch.

Finally, division of the equation of \(D\) by \(d^2\) yields
\[
  \left(\frac ed\right)^2+\frac ed
  =d+a+\frac{b^2}{d^2}.
\]
Adding \(d^2+d\) to both sides gives
\[
  \left(d+\frac ed\right)^2+\left(d+\frac ed\right)
  =a+d^2+\frac{b^2}{d^2},
\]
which is the asserted Artin--Schreier twist relation.  This proves part
(iv).

\smallskip
\noindent\emph{Step 5: recovery.}
In the odd model write \(D=(d,\varepsilon)\) and
\(f(u)=u^3+Au^2+Bu\).  Since \(2D\ne O\), one has
\(\varepsilon\ne0\).  For \(x=u(P)\), \(y=u(P+D)\), and \(x\ne d\),
the chord slope \(m=(v(P)-\varepsilon)/(x-d)\) satisfies
\(y=m^2-A-x-d\).  Multiplying by \((x-d)^2\), using
\(v(P)^2=f(x)\) and \(\varepsilon^2=f(d)\), and solving for \(v(P)\)
gives
\[
  v(P)=
  \frac{f(x)+f(d)-(y+A+x+d)(x-d)^2}{2\varepsilon}.
\]
After denominators are cleared this is a regular identity on the affine
state curve; at \(P=D\) direct substitution returns \(\varepsilon\).
For the binary model put
\(f_2(t)=t^3+at^2+b^2\), write \(P=(x,v)\), and let
\(y=u(P+D)\).  When \(x\ne d\), the binary chord slope
\(\lambda=(e+v)/(d+x)\) satisfies
\[
 y+a+d+x=\lambda^2+\lambda.
\]
After multiplication by \((d+x)^2\), substitution of
\(e^2+de=f_2(d)\) and \(v^2+xv=f_2(x)\), and collection of the terms
linear in \(v\), this becomes
\[
 d\,v=(y+a+d+x)(d+x)^2+xe+f_2(d)+f_2(x).
\]
Since \(d\ne0\), this recovers \(v\); after clearing the chord denominator
the identity extends regularly to the remaining affine state \(P=D\).
Indeed, substituting \(x=d\) and \(v=e\) makes the squared-factor term
zero and the two values \(f_2(d)\) cancel, leaving \(de\) on both sides.
The Kummer abscissa has its unique pole at \(O\), so an infinite first
coordinate identifies \(P=O\), while an infinite second coordinate
identifies \(P=-D\).  These two boundary charts and the affine formula are
disjoint and exhaustive, proving part (v).

\smallskip
\noindent\emph{Step 6: isogeny functoriality.}
An isogeny commutes with negation.  Therefore
\(\kappa'\circ\varphi\) is constant on every generic fibre
\(\{P,-P\}\) of \(\kappa\), and it factors uniquely through a morphism
\(\overline\varphi:\PP^1\to\PP^1\) satisfying
\(\overline\varphi\circ\kappa=\kappa'\circ\varphi\).  Comparing degrees
gives
\[
  2\deg\varphi
  =\deg(\kappa'\circ\varphi)
  =2\deg\overline\varphi,
\]
so \(\deg\overline\varphi=\deg\varphi\).  Since
\(\varphi(P+D)=\varphi(P)+D'\),
\[
 \mathcal S'_{D'}(\varphi(P))
 =\bigl(\overline\varphi(\kappa(P)),
        \overline\varphi(\kappa(P+D))\bigr).
\]
Both state maps are isomorphisms when \(2D\ne O\) and \(2D'\ne O'\),
which proves the coordinatewise formula and its degree.  The identities
\(\varphi(P+D)=\varphi(P)+D'\) and
\(\varphi([m]P)=[m]\varphi(P)\) give commutation with \(T_D\) and
\(\Delta_{D,m}\), respectively.  This proves part (vi) and completes the
proof.
\end{proof}

\section{Synthesis and research directions}
\label{subsec:arithmetic-scope-synthesis}

The symmetric biquadratic family
\[
  \Q_{\alpha,\beta,\gamma}:
  \qquad
  x^2y^2+\alpha(x^2+y^2)+\beta xy+\gamma=0
\]
provides the common geometric object around which this chapter is organized.
A smooth member is an anticanonical \((2,2)\)-curve in
\(\PP^1\times\PP^1\), hence a genus-one curve, and its two degree-two
projections produce the Vieta involutions whose composition is the QRT map.
This description is intrinsic before an origin is chosen: the natural datum
is a genus-one curve together with a distinguished automorphism \((C,T)\).
After a rational origin is selected, the same automorphism becomes
translation by the point \(D=T(O)\).  Conversely, under the separability hypothesis on the degree-two Kummer
map and the nondegeneracy condition \(2D\ne O\), the adjacent-state
construction
\[
  P\longmapsto\bigl(\kappa(P),\kappa(P+D)\bigr)
\]
shows how a genuinely preassigned marked point can be built into a symmetric
biquadratic model.  The chapter therefore connects the classical
Euler--Chasles and QRT descriptions with the arithmetic of pointed elliptic
curves, while keeping the distinction between rechoosing an origin and
encoding a fixed point into the model.

In odd characteristic, the even-quartic reduction, the extended Jacobi
quartic, and the Weierstrass model give complementary arithmetic interfaces.
When \(c^2=-4\alpha\gamma\) is defined over the ground field, the
QRT--Jacobi isomorphism identifies the pointed QRT curve with identity
\(O_Q=(0,c/(2\alpha))\); the transported Jacobi law is consequently the
full-point group law on that pointed curve.  The affine formulas, extended
projective formulas, boundary charts, and Weierstrass atlas together form a
complete arithmetic system.  A single Jacobi tuple is \(k\)-complete over an
odd finite field under the stated nonsquare condition; outside that regime,
completeness means an explicit finite atlas rather than an unproved extension
of a generic denominator formula.  The Kummer coordinate, differential
addition, recovery, and scalar-multiplication formulas are valid without
requiring a strict Montgomery normalization, whereas the latter additionally
requires the relevant square class to split over the base field.

Characteristic two is not obtained by setting \(2=0\) in the odd formulas.
Here the mixed term \(\beta xy\) becomes the Artin--Schreier term that keeps
the degree-two cover separable, and smoothness is governed by
\(\alpha\beta\gamma\ne0\).  The binary Weierstrass, L\'opez--Dahab, Kummer,
and companion formulas therefore constitute a separate arithmetic branch.
This branch supplies complete state-doubling formulas and makes precise why
the general three-parameter model, rather than the \(\beta=0\) slice, is the
natural binary symmetric biquadratic model.  The same symmetric biquadratic
equation therefore accommodates the even-quartic geometry of odd characteristic and the Artin--Schreier
geometry of characteristic two without changing its defining form.

The finite-field and isogeny results place these formulas in a moduli-theoretic
framework.  The chapter distinguishes coefficient triples, diagonal or
projective-normalizer equivalence, and abstract elliptic-curve isomorphism;
the corresponding counting theorems apply only to the equivalence relation
specified in their statements.  Likewise, the explicit two-, three-, and
odd-prime-degree isogenies are transported through Weierstrass models, while
a quotient is called model-preserving only after the codomain realization
equations have been solved over the ground field.  The product family
\(\C_{a,b,d}\), the original \(C_d\) family, and reciprocal \(C\)-curves
appear as distinguished slices or twists of the same symmetric
biquadratic geometry, with the necessary square-root and descent hypotheses
kept explicit.

A second contribution of the chapter is the adjacent-Kummer-state calculus.
The shift \(T_D\), state multiplication \(\Delta_{D,m}\), and the maps
\(T_D^r\Delta_{D,m}\) realize the affine index transformations
\(n\mapsto mn+r\) on one fixed genus-one state curve.  For \(m=2\), the two
binary branches compute the states indexed by \(2n\) and \(2n+1\), giving a
proved logarithmic ladder rather than the linear-time iteration of a cheap
fixed translation.  The Jacobi companion restores the orientation discarded
by the Kummer quotient, and the resulting elliptic Lucas theory links the
QRT recurrence with classical Lucas sequences at the nodal boundary, Ward
elliptic divisibility sequences, division polynomials, sigma functions, and
elliptic nets.  State-division polynomials and QRT orbit invariants then
connect exact periods with torsion and, for translation order prime to the
characteristic, with the cyclic quotient isogeny. 

The operation counts in this chapter are exact upper bounds for the displayed
coordinate schedules.  The complete end-to-end count combines each headline
circuit with its listed encoding, decoding, subgroup validation, affine
normalization, and one-time transport components.  The direct McMillan update performs one
fixed translation and reaches its \(n\)-th iterate in \(O(n)\) steps; the
\(O(\log n)\) statements use state doubling or a Kummer ladder.  The
odd-characteristic Lucas schedules exhibit substantial common-subexpression
sharing and a sharp \(6\M\) lower bound within the explicitly defined
Segre-first circuit class.  On exactly the split evenizable locus, the input and
output changes are projectively linear, so
Theorem~\ref{thm:QRT-Montgomery-conjugate-state-doubling} transports the
Montgomery core with unchanged general multiplication and squaring counts,
giving \(4\M+4\Sqr+2\Dpar\).  

Natural extensions of the present theory include direct canonical-QRT circuits
that build on the Montgomery-conjugate upper bound while reducing
fixed-constant products, parallel depth, register
pressure, or completeness overhead.  Low-degree maps
\(\Delta_{D,m}\), faster recurrences for the state-division polynomials, and
short native formulas for cyclic quotients would extend the theory in
arithmetically meaningful directions.  Supersingular
characteristic-two states, binary small-constant orbits, and nonsplit
odd-characteristic states form the corresponding characteristic-specific
normal-form, descent, and security branches.  Applications to ECM, finite-field period detection,
elliptic pseudoprimes, and Diophantine recurrences lead naturally to complete
constant-time implementations that account for validation, recovery, memory
traffic, and exceptional branches.  These directions arise directly from the
proved geometry and arithmetic of the chapter.

\part[Supersingular-isogeny cryptography]
{Supersingular-Isogeny Cryptography on the \texorpdfstring{\(\mathcal C_d\)}{Cd} Platform and Its Extensions}
\label{part:supersingular-isogeny-cryptography}
\partoverview{Parts I--VI have established the native \(\mathcal C_d\) model, the full product and one-sided families, reciprocal \(C\)-curves, and the symmetric QRT envelope.  This part brings those structures together as one isogeny-cryptographic platform.  The chapter restates the native Kummer, target-parameter, full-point, degree-three, and odd- and even-characteristic degree-two formulas needed for an independent cryptographic reading; it then compares every declared output with the corresponding Edwards interface at the same arithmetic level.}

\chapter[Supersingular-Isogeny Cryptography]
{Supersingular-Isogeny Cryptography on the \texorpdfstring{\(\mathcal C_d\)}{Cd} Platform and Its Extensions}
\label{ch:Cd-isogeny-cryptography}

\section{Purpose, scope, and comparison discipline}
\label{sec:Cd-iso-purpose}

This chapter synthesizes the isogeny arithmetic of the native model and the
extension families established in Parts~I--VI.  Its principal object is
\[
  \mathcal C_d:\qquad (u^2+u)(v^2+v)=d,
  \qquad \rho=1-16d,
\]
and the one-sided family \(\mathcal T_{a,d}\), the full product family
\(\mathcal C_{a,b,d}\), reciprocal \(C\)-curves, and the QRT-Jacobian
envelope are used whenever the quotient leaves the one-parameter locus.
Every output is declared as a point, a Kummer class, or a target parameter
on one of these models.  Edwards, Montgomery, generalized Montgomery, and
binary Weierstrass coordinates serve as exact dictionaries or optimized
internal states; they do not silently replace the declared endpoint.

In characteristic two, the smooth models
\(\mathcal C_{a,b,d}\) and \(\mathcal T_{a,d}\) are ordinary, exactly as
binary Edwards is an ordinary binary model.  They therefore provide the
proper platform for the separable binary two-isogeny without asserting a
characteristic-two supersingular member.  For a self-contained reading, the
chapter restates the native \(u\)-Kummer, target-parameter, full-point,
degree-three, and both characteristic branches of the degree-two formulas.
The earlier chapters supply their wider geometric setting; the central
identities used here are reproved or verified by direct substitution.

Three logically distinct levels are kept separate throughout.
\begin{enumerate}[label=\textup{(\arabic*)}]
 \item \emph{Protocol security}: the public data and hardness assumption
       belong to the protocol and are unchanged by a coordinate dictionary.
 \item \emph{Algebraic availability}: the quotient, image parameter, and
       full or Kummer point maps require explicit separability and
       field-of-definition hypotheses.
 \item \emph{Implementation performance}: every comparison fixes the input
       and output representations, preprocessing policy, field backend, and
       platform.  Symbolic operation counts and cycle or nanosecond
       measurements describe complementary layers.
\end{enumerate}
In particular, an \(x\)-only or \(w\)-only evaluation is compared with
another Kummer evaluation, and a full-point map is compared with another
full-point map.  Kernel enumeration, kernel normalization, image-curve
recovery, and protocol-level validation are charged separately.

\section{The current supersingular-isogeny landscape}
\label{sec:Cd-iso-landscape}

The 2022--2023 attacks on SIDH recover the secret isogeny by exploiting the
auxiliary torsion images carried by the SIDH public key
\cite{CastryckDecru2023}.  They invalidate that public-key design, not the
algebraic identities for evaluating a separable isogeny.  The surviving and
developing branches use different assumptions and different arithmetic
layers:

\begin{longtable}{L{2.4cm}L{3.7cm}L{4.0cm}L{3.2cm}}
\caption{Protocol landscape and the role of the \(C_d\) platform}
\label{tab:Cd-current-isogeny-landscape}\\
\toprule
branch & current status and principal arithmetic & direct \(C_d\) or
extension-model interface & role in the \(C_d\) platform\\
\midrule
\endfirsthead
\toprule
branch & current status and principal arithmetic & direct \(C_d\) or
extension-model interface & role in the \(C_d\) platform\\
\midrule
\endhead
SIDH/SIKE
& the auxiliary-torsion-image public-key construction is broken
  \cite{CastryckDecru2023}
& the odd-degree product identities remain valid algebraic tools, but are
  not a security repair
& retain formulas; do not infer protocol security from a model change\\
CSIDH, CTIDH, and SQALE
& commutative class-group actions over \(\F_p\), with constant-time and
  square-root-V\'elu variants; deterministic batching, dummy-free
  hardening, and projective masking now refine the implementation layer
  \cite{CSIDH2018,CTIDH2021,SQALE2022,CamposEtAl2024,dCTIDH2025,
  HardenedCTIDH2025,SideChannelCSIDH2026}
& native \(u\)-Kummer for low and medium odd degrees; \(C_d\)-\(w\) for
  high degrees; \(\mathcal T_{a,d}\) for twisted-Edwards coverage; the
  projective pair \((u+1:u)\) supports quotient masking
& a direct target for the optimized kernels proved below\\
SQIsign and SQIsign2D
& SQIsign is in the NIST additional-signature Round~2 process; current
  variants combine quaternion algorithms with one- or two-dimensional
  isogeny representations
  \cite{NISTRound2Signatures2026,SQIsign2DWest2024,SQIsign2DPush2025,
  SQIsignFixedPrecision2026}
& \(C_d\) supplies compact elliptic input, output, and one-dimensional
  edge charts; a two-dimensional theta or product-isogeny layer retains
  its own polarization data
& use \(C_d\) at elliptic boundaries and factors, with an explicit
  product/polarization interface\\
two- and four-dimensional group actions
& KLaPoTi and qt-PEGASIS give effective higher-dimensional actions,
  while CORAL specializes to a fast restricted action built from
  two-dimensional two-isogenies
  \cite{KLaPoTi2025,qtPEGASIS2025,CORAL2026}
& \(C_d\), \(\mathcal T_{a,d}\), and reciprocal charts represent the
  elliptic factors and boundary curves; the QRT interface transports
  adjacent Kummer states coordinatewise
& retain the native elliptic endpoints while carrying the rank-two kernel
  and polarization in the natural product/theta layer\\
PRISM
& the 2026 revision uses large-prime-degree isogenies and gives a
  standard-model signature construction \cite{PRISM2026}
& the \(O(\ell)\) native \(u\) backend and the
  \(\widetilde O(\sqrt\ell)\) \(C_d\)-\(w\) backend directly address its
  large-degree evaluation step
& especially relevant to the dual-backend design below\\
oriented large-discriminant actions
& deterministic and compressed implementations now use twisted-Edwards
  isogeny arithmetic and wide-vector parallelism \cite{CSIDHLDO2026}
& the exact dictionary
  \(\mathcal T_{a,d}\leftrightarrow\) twisted Edwards gives a native
  \(C\)-curve realization with the same arithmetic state
& a principal application of the one-sided extension model\\
radical and structured actions
& current work includes radical \(3\)-isogenies on
  \((2,\varepsilon)\)-structures \cite{RadicalStructures2026}
& reciprocal and QRT extensions supply marked two-isogeny and
  re-embedding interfaces; a structure-specific descent proof is still
  required before identifying parameters
& a concrete extension program rather than an unsupported identification\\
\bottomrule
\end{longtable}

The Algebraic Isogeny Model gives a 2026 formal setting for arbitrary
isogenies over \(\F_{p^2}\), proves an ID-soundness result for the SQIsign
identification protocol in that model, and relates discrete-logarithm and
Diffie--Hellman problems for several SIDH-derived exchanges
\cite{AIM2026}.  This reinforces the separation just made: the security
model describes the adversary and public interface, whereas the \(C_d\)
results below optimize a mathematically equivalent arithmetic
representation.

\section[Algebraic and algorithmic criteria]
{Algebraic and algorithmic criteria for \texorpdfstring{\(C_d\)}{Cd} and its extension families}
\label{sec:Cd-four-tests}

Let \(K\) be a separable cyclic subgroup of odd order
\(\ell=2m+1\), and choose one representative
\(Q_i=(u_i,v_i)\) from each pair \(\{\pm Q_i\}\),
\(1\le i\le m\).  Four core properties organize the odd-degree analysis:
\begin{enumerate}[label=\textbf{C\arabic*.}]
 \item Does the quotient return to a member \(\mathcal C_{d'}\) over the
       declared field?
 \item Can the half-kernel be enumerated in the native Kummer coordinate
       \(\kappa_d=(u+1:u)\)?
 \item Can the isogeny be evaluated directly on that Kummer line, without
       constructing a full auxiliary-model point?
 \item Can \(d'\) be recovered by a simple product of kernel coordinates,
       preferably without an online inversion?
\end{enumerate}

They are necessary but not by themselves a complete cryptographic
assessment.  The complete comparison also charges kernel enumeration,
full-point evaluation, target re-embedding, the degree-three specialization,
and separable degree-two quotients in both odd and even characteristic.
It treats \(C_d\), the one-sided family \(\mathcal T_{a,d}\), the full
product family \(\mathcal C_{a,b,d}\), and the reciprocal/QRT closure
envelope separately.  This prevents an odd-degree closure theorem from
being misread as an even-degree closure theorem.

\begin{theorem}[Odd-degree closure, native kernel, evaluation, and parameter]
\label{thm:Cd-four-tests}
Assume \(\charac k\ne2\),
\(d\rho\ne0\), and
\(\gcd(\ell,\charac k)=1\).  Then \(C_d\) satisfies all four core
odd-degree properties.
More precisely:
\begin{enumerate}[label=\textup{(\roman*)}]
 \item the odd quotient has the canonical target
       \(\mathcal C_{d'}\);
 \item the points \(\kappa_d([i]Q)\), \(1\le i\le m\), are obtained by
       the native \(x\)DBL and \(x\)ADD laws of
       Theorems~\ref{thm:xdbl-thesis} and~\ref{thm:xadd-thesis};
 \item if
       \[
        A(u)=\prod_{i=1}^m(u+u_i+1),\qquad
        B(u)=\prod_{i=1}^m(u-u_i),
       \]
       then the native Kummer image is
       \[
         u'=\frac{uB(u)^2}
          {(u+1)A(u)^2-uB(u)^2};
       \]
 \item with \(R_K=\prod_i(2u_i+1)\),
       \[
         \rho'=\frac{\rho^\ell}{R_K^8},
         \qquad d'=\frac{1-\rho'}{16}.
       \]
\end{enumerate}
\end{theorem}

\begin{proof}
Part~(i) is Theorem~\ref{thm:odd-velu-thesis}; its hypotheses are exactly
the present smoothness and separability hypotheses.  For part~(ii), start
with \(\kappa_d(Q)\), compute
\(\kappa_d(2Q)\) by \(x\)DBL, and for \(i\ge2\) compute
\(\kappa_d((i+1)Q)\) from
\(\kappa_d(iQ),\kappa_d(Q),\kappa_d((i-1)Q)\) by \(x\)ADD.  The known
difference is \((iQ)-Q=(i-1)Q\), so every invocation satisfies the
oriented differential-addition hypothesis.

For part~(iii), use the ordinary Kummer coordinate
\(x=(u+1)/u\) and write \(x_i=(u_i+1)/u_i\).  The normalized odd-kernel
product on this Kummer line is
\[
 x'=x\prod_{i=1}^m
       \left(\frac{xx_i-1}{x-x_i}\right)^2.
\]
The numerator and denominator of one factor are
\[
 xx_i-1=\frac{u+u_i+1}{uu_i},\qquad
 x-x_i=\frac{u_i-u}{uu_i},
\]
and therefore the factor is
\(- (u+u_i+1)/(u-u_i)\).  Squaring removes the sign and gives
\[
 x'=\frac{u+1}{u}\frac{A(u)^2}{B(u)^2}.
\]
Since the target native coordinate satisfies \(u'=1/(x'-1)\), a single
common-denominator calculation yields
\[
 u'=\frac{uB(u)^2}{(u+1)A(u)^2-uB(u)^2},
\]
as asserted.  Finally, Theorem~\ref{thm:odd-velu-thesis} gives
\(\rho'=\rho^\ell/R_K^8\), and the fixed parameter dictionary
\(\rho'=1-16d'\) gives part~(iv).  Thus each property has an explicit
construction rather than only an existence assertion.
\end{proof}

\begin{longtable}{L{2.6cm}L{3.6cm}L{3.8cm}L{3.4cm}}
\caption{Scope of the complete \(C\)-curve isogeny analysis}
\label{tab:Cd-complete-isogeny-audit}\\
\toprule
interface & closure statement & native quotient data & low-degree role\\
\midrule
\endfirsthead
\toprule
interface & closure statement & native quotient data & low-degree role\\
\midrule
\endhead
\(\mathcal C_d\), odd degree
& canonical \(d'=(1-\rho')/16\) for every separable odd kernel
& \((u+1:u)\), two monic products, and one kernel product for \(\rho'\)
& optimized degree three is obtained by setting \(m=1\) and collapsing the
  derivative product\\
\(\mathcal C_d\), odd-characteristic degree two
& the quotient always has a root-free generalized-Montgomery/QRT endpoint;
  return to \(\mathcal T_{a,d}\) requires \(\sqrt\rho\), and return to
  \(\mathcal C_d\) has the additional stated square-class condition
& the marked kernel is fixed and the native source factor is
  \(t=u(u+1)\)
& one-square Kummer core; no half-kernel enumeration\\
\(\mathcal T_{a,d}\)
& exact twisted-Edwards dictionary in odd characteristic and twist-stable
  degree-two closure in characteristic two
& the split factor supplies a native Montgomery Kummer coordinate
& root-selected odd degree two and one-operation binary degree two\\
\(\mathcal C_{a,b,d}\)
& diagonal-Edwards transport in odd characteristic; Jacobian closure in
  characteristic two
& two boundary square classes and a Jacobian Kummer interface
& descent envelope when a pure \(C_d\) endpoint is unavailable\\
reciprocal and QRT extensions
& exact re-embedding criterion over the ground field
& general marked-two-torsion Kummer coordinate
& root-free odd-characteristic two-isogeny closure and coordinatewise state
  transport\\
Edwards comparator
& ordinary/twisted Edwards closure for odd kernels; even-degree return may
  require a square root
& Edwards full coordinates, \(y\)-Kummer, and \(w\)-Kummer
& strongest published degree-three and transported degree-two baselines\\
\bottomrule
\end{longtable}

\section{Exact dictionaries needed by the optimized formulas}
\label{sec:Cd-crypto-dictionaries}

\begin{proposition}[The \(C_d\)--Montgomery dictionary]
\label{prop:Cd-crypto-Montgomery}
Put
\[
 x=\frac{u+1}{u},\qquad y=x(2v+1),\qquad
 B_d=\frac1{4d},\qquad A_d=B_d-2.
\]
Then the smooth completion of \(\mathcal C_d\) is isomorphic to the
smooth completion of
\[
       M_d:\qquad B_dy^2=x^3+A_dx^2+x.
\]
On the common affine chart the inverse is
\[
       u=\frac1{x-1},\qquad
       v=\frac12\left(\frac yx-1\right).
\]
\end{proposition}

\begin{proof}
From \(x=(u+1)/u\) one obtains
\[
 u=\frac1{x-1},\qquad u^2+u=\frac{x}{(x-1)^2}.
\]
Writing \(V=2v+1=y/x\) gives
\[
 v^2+v=\frac{V^2-1}{4}
       =\frac{y^2-x^2}{4x^2}.
\]
Substitution in the \(C_d\) equation and multiplication by
\(4x(x-1)^2\) give
\[
 y^2=x^2+4dx(x-1)^2
     =4dx^3+(1-8d)x^2+4dx.
\]
Division by \(4d\) yields the displayed Montgomery equation because
\((1-8d)/(4d)=1/(4d)-2=A_d\).  The two rational substitutions are
inverse wherever their denominators are nonzero.  A birational map
between smooth projective curves extends uniquely to an isomorphism,
which supplies the omitted boundary points.
\end{proof}

\begin{proposition}[The \(C_d\)--Edwards and \(w\) dictionaries]
\label{prop:Cd-crypto-Edwards-w}
Let
\[
 U=2u+1,\qquad V=2v+1,\qquad
 \xi=V^{-1},\qquad \eta=U^{-1}.
\]
Then
\[
 \xi^2+\eta^2=1+\rho\xi^2\eta^2.
\]
The Edwards \(w\)-function and its shifted form pull back to
\begin{align}
 w_C
 &=\frac{\rho}{U^2V^2}
  =\frac{\rho(U^2-1)}{U^2(U^2-\rho)},
 \label{eq:Cd-crypto-w}\\
 \eta_C
 &=\frac{w_C+1}{2}
  =\frac{U^4-\rho}{2U^2(U^2-\rho)}.
 \label{eq:Cd-crypto-eta}
\end{align}
If \(U=P/Q\), an inversion-free representative of \(w_C\) is
\begin{equation}
 (W:Z)=
 \bigl(\rho Q^2(P^2-Q^2):
       P^2(P^2-\rho Q^2)\bigr).
 \label{eq:Cd-crypto-w-projective}
\end{equation}
\end{proposition}

\begin{proof}
The centered equation is
\[
       (U^2-1)(V^2-1)=16d=1-\rho.
\]
Expanding and dividing by \(U^2V^2\) gives the Edwards equation.
Solving the centered equation for \(V^2\) gives
\[
 V^2-1=\frac{1-\rho}{U^2-1},\qquad
 V^2=\frac{U^2-\rho}{U^2-1}.
\]
Substitution in \(\rho/(U^2V^2)\) proves
\eqref{eq:Cd-crypto-w}.  Moreover
\[
 \frac{w_C+1}{2}
 =\frac{\rho(U^2-1)+U^2(U^2-\rho)}
        {2U^2(U^2-\rho)}
 =\frac{U^4-\rho}{2U^2(U^2-\rho)},
\]
which proves \eqref{eq:Cd-crypto-eta}.  Finally insert \(U=P/Q\) in
\eqref{eq:Cd-crypto-w}, multiply numerator and denominator by \(Q^4\),
and read off \eqref{eq:Cd-crypto-w-projective}.
\end{proof}

\section[Edwards--Cd product correspondence]
{Term-by-term correspondence with the Edwards products}
\label{sec:Cd-Edwards-termwise-products}

The comparison with Edwards arithmetic is exact at the level of rational
functions, not only at the level of \(j\)-invariants.  Put
\[
 r=2u+1,\qquad s=2v+1,\qquad
 \xi=s^{-1},\qquad \eta=r^{-1}.
\]
Then \(C_d\) is the ordinary Edwards curve
\[
             \xi^2+\eta^2=1+\rho\xi^2\eta^2.
\]
For \(Q_i=(u_i,v_i)\), write
\[
 r_i=2u_i+1,\quad s_i=2v_i+1,\quad
 \alpha_i=s_i^{-1},\quad\beta_i=r_i^{-1},\quad
 B_K=\prod_i\beta_i.
\]
The Moody--Shumow products are
\begin{align}
 \xi'&=\frac{\xi}{B_K^2}
 \prod_i
 \frac{\beta_i^2\xi^2-\alpha_i^2\eta^2}
      {1-\rho^2\alpha_i^2\beta_i^2\xi^2\eta^2},
 \label{eq:Cd-crypto-Edwards-product-xi}\\
 \eta'&=\frac{\eta}{B_K^2}
 \prod_i
 \frac{\beta_i^2\eta^2-\alpha_i^2\xi^2}
      {1-\rho^2\alpha_i^2\beta_i^2\xi^2\eta^2}.
 \label{eq:Cd-crypto-Edwards-product-eta}
\end{align}
Substitution of the four reciprocal definitions gives, factor by factor,
\begin{align*}
 1-\rho^2\alpha_i^2\beta_i^2\xi^2\eta^2
 &=\frac{(rr_iss_i)^2-\rho^2}
         {(rr_iss_i)^2},\\
 \beta_i^2\eta^2-\alpha_i^2\xi^2
 &=\frac{s_i^2s^2-r_i^2r^2}{(rr_iss_i)^2},\\
 \beta_i^2\xi^2-\alpha_i^2\eta^2
 &=\frac{s_i^2r^2-r_i^2s^2}{(rr_iss_i)^2}.
\end{align*}
Thus the three Edwards factors are exactly the
\(\Delta_i,A_i,B_i\) factors of
Theorem~\ref{thm:odd-velu-thesis}; no factor and no normalization constant
is lost in the \(C_d\) pullback.

The smaller native \(u\)-Kummer formula is obtained by taking the same
quotient one level farther down.  Since
\[
 x=\frac{1+\eta}{1-\eta}=\frac{u+1}{u},\qquad
 x_i=\frac{1+\beta_i}{1-\beta_i}=\frac{u_i+1}{u_i},
\]
one Edwards/Montgomery kernel factor becomes
\begin{equation}
 \frac{xx_i-1}{x-x_i}
 =-\frac{u+u_i+1}{u-u_i}.
\label{eq:Cd-crypto-termwise-linear-factor}
\end{equation}
The sign disappears after squaring.  Moreover
\[
 B_K=\prod_i(2u_i+1)^{-1}=R_K^{-1},
\]
so the Edwards image parameter
\(\rho^\ell B_K^8\) is literally
\(\rho^\ell/R_K^8\).

\begin{longtable}{L{2.9cm}L{3.9cm}L{3.9cm}L{2.4cm}}
\caption{Term-by-term Edwards--\(C_d\) correspondence}
\label{tab:Cd-Edwards-termwise-products}\\
\toprule
object & Edwards expression & \(C_d\) expression & consequence\\
\midrule
\endfirsthead
\toprule
object & Edwards expression & \(C_d\) expression & consequence\\
\midrule
\endhead
full coordinates
& \((\xi,\eta)=(s^{-1},r^{-1})\)
& \((u,v)=((r-1)/2,(s-1)/2)\)
& exact birational endpoints\\
sign Kummer
& \((1+\eta:1-\eta)\)
& \((u+1:u)\)
& the projective pairs differ only by the common scale \(2/r\)\\
one kernel factor
& \((xx_i-1)/(x-x_i)\)
& \(-(u+u_i+1)/(u-u_i)\)
& two monic linear factors in the native coordinate\\
kernel normalization
& \(B_K=\prod_i\beta_i\)
& \(R_K=\prod_i(2u_i+1)\)
& \(B_K=R_K^{-1}\)\\
image parameter
& \(\rho'=\rho^\ell B_K^8\)
& \(\rho'=\rho^\ell/R_K^8\)
& identical invariant, direct \(d'=(1-\rho')/16\)\\
\(w\)-coordinate
& \(w=\rho\xi^2\eta^2\)
& \(w_C=\rho/(r^2s^2)\)
& identical function and identical square-root-V\'elu graph\\
full product
& equations
  \eqref{eq:Cd-crypto-Edwards-product-xi}--%
  \eqref{eq:Cd-crypto-Edwards-product-eta}
& \(\Delta_i,A_i,B_i\) followed by reciprocation
& global model closure and boundary completion agree\\
\bottomrule
\end{longtable}

This table separates two kinds of gain.  The \(w\)-backend is a
zero-loss transport and therefore has equal field-operation cost.
The \(u\)-backend is a further native specialization: after an affine
half-kernel has been made reusable, its two monic products cost
\(2m\M+2\Sqr\), while the full \(C_d\) product--derivative map and the
published Edwards full-point product are compared separately in
Section~\ref{sec:Cd-Edwards-full-comparison}.  No Kummer-only count is
used there as if it were a full-point count.

\section[Native odd-degree Kummer evaluation]
{Native half-kernel enumeration and odd-degree Kummer evaluation}
\label{sec:Cd-native-kernel-evaluation}

\subsection{Enumeration on the native Kummer line}

Let \(Q\) generate \(K\), set \(K_i=\kappa_d([i]Q)\), and retain
projective pairs throughout enumeration.  For \(m\ge2\), the direct
sequential schedule is
\[
 K_1=\kappa_d(Q),\qquad
 K_2=x\mathrm{DBL}(K_1),\qquad
 K_{i+1}=x\mathrm{ADD}(K_i,K_1,K_{i-1})
 \quad(2\le i<m).
\]
The third input is valid because
\([i]Q-Q=[i-1]Q\).  With the costs of
Theorems~\ref{thm:xdbl-thesis} and~\ref{thm:xadd-thesis}, this simple
schedule costs
\begin{equation}
 (4m-6)\M+(2m-2)\Sqr+\Dpar
 \qquad(m\ge2),
\label{eq:Cd-sequential-kernel-cost}
\end{equation}
excluding the given \(K_1\).  Indeed, one double contributes
\(2\M+2\Sqr+\Dpar\), and the \(m-2\) additions contribute
\((m-2)(4\M+2\Sqr)\).  A differential addition chain can replace this
sequential schedule without changing any later formula.

If the kernel is to be shared by several evaluations, convert
\[
 (X_i:Z_i)\longmapsto
 u_i=\frac{Z_i}{X_i-Z_i}
\]
with one batch inversion.  Inverting the \(m\) denominators costs
\(3(m-1)\M+\Inv\), and multiplying the \(m\) inverse denominators by
the \(Z_i\) costs another \(m\M\).  Thus the exact conversion cost is
\begin{equation}
                  (4m-3)\M+\Inv.
\label{eq:Cd-kernel-affinization-cost}
\end{equation}

\subsection[The specialized u-Kummer product]
{The specialized \(u\)-Kummer product}

\begin{theorem}[Native odd-degree \(C_d\) Kummer evaluation]
\label{thm:Cd-native-odd-xeval}
Under the hypotheses of Theorem~\ref{thm:Cd-four-tests}, define
\[
 A(u)=\prod_{i=1}^m(u+u_i+1),\qquad
 B(u)=\prod_{i=1}^m(u-u_i),
\]
\[
 N=(u+1)A(u)^2,\qquad D=uB(u)^2.
\]
Then the quotient isogeny, normalized to the target
\(\mathcal C_{d'}\), induces
\begin{equation}
   \boxed{\quad u'=\frac{D}{N-D},\qquad
          (u'_0:u'_1)=(D:N-D).\quad}
\label{eq:Cd-native-odd-xeval}
\end{equation}
The formula depends only on the unsigned kernel classes
\(\{Q_i,-Q_i\}\).
\end{theorem}

\begin{proof}
In the Montgomery coordinate \(x=(u+1)/u\), put
\[
       x_i=\frac{u_i+1}{u_i}.
\]
The normalized odd-degree Kummer quotient is
\[
 f(x)=x\prod_{i=1}^m
       \left(\frac{xx_i-1}{x-x_i}\right)^2
 \cite{CostelloHisil2017}.
\]
Every factor can be simplified inside the \(C_d\) function field:
\[
 xx_i-1
 =\frac{(u+1)(u_i+1)-uu_i}{uu_i}
 =\frac{u+u_i+1}{uu_i},
\]
and
\[
 x-x_i
 =\frac{(u+1)u_i-(u_i+1)u}{uu_i}
 =-\frac{u-u_i}{uu_i}.
\]
After squaring, the signs and the common denominators cancel.  Hence
\[
 f(x)=\frac{u+1}{u}
       \left(\frac{A(u)}{B(u)}\right)^2
      =\frac ND.
\]
The target native coordinate satisfies
\[
 x'=\frac{u'+1}{u'},\qquad
 u'=\frac1{x'-1}.
\]
Substitution of \(x'=N/D\) gives
\[
 u'=\frac1{N/D-1}=\frac{D}{N-D},
\]
and clearing its single denominator gives the projective pair in
\eqref{eq:Cd-native-odd-xeval}.  Since negation fixes the \(u\)-Kummer
coordinate, replacing \(Q_i\) by \(-Q_i\) leaves \(u_i,A,B\) unchanged.
\end{proof}

\begin{proposition}[Exact Kummer evaluation costs]
\label{prop:Cd-native-xeval-costs}
Ignore additions, subtractions, and multiplication by \(2\).
\begin{enumerate}[label=\textup{(\roman*)}]
 \item Affine \(u\), affine \(u_i\), and projective \(u'\) require
       \[
                     \boxed{2m\M+2\Sqr}.
       \]
 \item Affine normalization of the output adds
       \(\M+\Inv\).
 \item General projective input and kernel pairs require
       \[
                     \boxed{4m\M+2\Sqr}.
       \]
\end{enumerate}
\end{proposition}

\begin{proof}
For part~(i), the products \(A\) and \(B\) each use \(m-1\)
multiplications.  The two squares use \(2\Sqr\), and multiplication by
\(u+1\) and \(u\) uses \(2\M\).  The total is
\[
 2(m-1)\M+2\M+2\Sqr=2m\M+2\Sqr.
\]
Normalizing \(D:N-D\) uses one inverse and one final multiplication,
which proves part~(ii).

For part~(iii), write the source and kernel Kummer pairs as
\((X:Z)\) and \((X_i:Z_i)\), and set
\[
 P=\prod_i(XX_i-ZZ_i),\qquad
 Q=\prod_i(XZ_i-ZX_i).
\]
The output is
\[
                  (X':Z')=(XP^2:ZQ^2).
\]
Both bilinear factors for a fixed \(i\) are obtained with two
multiplications by computing
\[
 a_i=(X-Z)(X_i+Z_i),\qquad
 b_i=(X+Z)(X_i-Z_i);
\]
then \(a_i+b_i=2(XX_i-ZZ_i)\) and
\(a_i-b_i=2(XZ_i-ZX_i)\).  The common powers of \(2\) are projectively
irrelevant.  The \(m\) factor pairs cost \(2m\M\), their two
accumulators cost \(2(m-1)\M\), and final assembly costs
\(2\M+2\Sqr\).  The total is \(4m\M+2\Sqr\).
\end{proof}

\subsection{When kernel affinization pays}

Suppose one half-kernel is used to evaluate \(n\) Kummer points.
The general projective route costs
\[
            n(4m\M+2\Sqr).
\]
The affine-kernel route, including
\eqref{eq:Cd-kernel-affinization-cost}, costs
\[
 (4m-3)\M+\Inv+n(2m\M+2\Sqr).
\]
Consequently the affine route is cheaper precisely when
\begin{equation}
 \operatorname{cost}(\Inv)
  <(2nm-4m+3)\operatorname{cost}(\M).
\label{eq:Cd-affinization-threshold}
\end{equation}
The threshold depends only on the public degree, batch size, and field
backend, so a constant-time implementation can choose it offline.

\section{Image parameter and model closure}
\label{sec:Cd-image-parameter-crypto}

\begin{theorem}[Kernel-product formula for the target parameter]
\label{thm:Cd-image-parameter-product}
Let the hypotheses of Theorem~\ref{thm:Cd-four-tests} hold.
For affine kernel coordinates,
\begin{equation}
 \boxed{\qquad
 \rho'=\frac{\rho^\ell}
 {\left(\prod_{i=1}^m(2u_i+1)\right)^8},
 \qquad
 d'=\frac{1-\rho'}{16}.
 \qquad}
\label{eq:Cd-image-parameter-affine}
\end{equation}
For projective Kummer points \((X_i:Z_i)\),
\begin{equation}
 \rho'=\rho^\ell
  \left(\prod_{i=1}^m
       \frac{X_i-Z_i}{X_i+Z_i}\right)^8.
\label{eq:Cd-image-parameter-projective}
\end{equation}
\end{theorem}

\begin{proof}
The already proved full-coordinate quotient in
Theorem~\ref{thm:odd-velu-thesis} gives
\[
       \rho'=\frac{\rho^\ell}{R_K^8},
       \qquad R_K=\prod_i(2u_i+1).
\]
This is \eqref{eq:Cd-image-parameter-affine}.  To obtain the projective
form, use the native relation
\[
 (X_i:Z_i)=(u_i+1:u_i).
\]
It gives
\[
 \frac{X_i-Z_i}{X_i+Z_i}
 =\frac{(u_i+1)-u_i}{(u_i+1)+u_i}
 =\frac1{2u_i+1}.
\]
Multiplying these identities and raising to the eighth power proves
\eqref{eq:Cd-image-parameter-projective}.  Finally,
\(\rho'=1-16d'\) is equivalent to
\(d'=(1-\rho')/16\).
\end{proof}

Let \(R=\prod_i(2u_i+1)\).  No inversion is needed if the target
parameter is kept projectively:
\begin{equation}
       (\rho'_0:\rho'_1)=(\rho^\ell:R^8),
       \qquad
       (d'_0:d'_1)=(R^8-\rho^\ell:16R^8).
\label{eq:Cd-projective-image-parameter}
\end{equation}
The variable part of the affine-kernel computation is
\[
 (m-1)\M+3\Sqr+\operatorname{Pow}(\rho,\ell),
\]
where the public exponent \(\ell\) has a fixed addition chain.
For general projective kernel pairs, put
\[
 R_-=\prod_i(X_i-Z_i),\qquad
 R_+=\prod_i(X_i+Z_i).
\]
Then
\[
 (\rho'_0:\rho'_1)=(\rho^\ell R_-^8:R_+^8),
\]
at variable cost
\[
 (2m-1)\M+6\Sqr+\operatorname{Pow}(\rho,\ell).
\]
Indeed the two length-\(m\) products cost \(2m-2\) multiplications, the
eighth powers cost six squarings, and multiplication of
\(\rho^\ell\) by \(R_-^8\) costs the remaining multiplication.  Keeping
the parameter as the pair \((\rho^\ell:R^8)\) in the affine-kernel case
avoids this last multiplication because the numerator needs no kernel
factor.

\begin{corollary}[Odd-degree model closure]
\label{cor:Cd-odd-model-closure-crypto}
The quotient in Theorem~\ref{thm:Cd-image-parameter-product} is a smooth
member \(\mathcal C_{d'}\).  No square root or field extension is needed
to write its canonical parameter.
\end{corollary}

\begin{proof}
The quotient of a smooth elliptic curve by a finite separable subgroup is
smooth.  Theorem~\ref{thm:odd-velu-thesis} identifies its canonical
ordinary-Edwards parameter with the \(\rho'\) in
\eqref{eq:Cd-image-parameter-affine}.  A smooth ordinary Edwards curve
has \(\rho'\ne0,1\); hence
\[
        d'=\frac{1-\rho'}{16}\ne0,\qquad
        1-16d'=\rho'\ne0.
\]
These are exactly the \(C_d\) smoothness conditions.  Formula
\eqref{eq:Cd-projective-image-parameter} uses only field addition,
multiplication, squaring, and a public fixed exponent, so no root
extraction is present.
\end{proof}

\section[Cd-specific full-point product--derivative formula]
{A \(C_d\)-specific full-point product--derivative formula}
\label{sec:Cd-full-isogeny-crypto}

Here \emph{full point} means that both native coordinates are returned.
It does not mean that one affine rational pair covers every boundary
point.  The fast main chart below is accompanied by a projective chart and
the already proved global product map of
Theorem~\ref{thm:odd-velu-thesis}; together they represent the complete
isogeny morphism.

\begin{theorem}[\(C_d\) full-point product--derivative formula]
\label{thm:Cd-full-product-derivative}
Retain the hypotheses and notation of
Theorem~\ref{thm:Cd-native-odd-xeval}.  Put
\[
 V=2v+1,\qquad
 \dot A=\frac{dA}{du},\qquad
 \dot B=\frac{dB}{du},
\]
\begin{equation}
 J=AB-2u(u+1)(\dot A B-A\dot B),
 \qquad
 c_K=1+2\sum_{i=1}^m(2u_i+1).
\label{eq:Cd-full-J-cK}
\end{equation}
Then \(c_K\ne0\), and the marked, normalized isogeny to
\(\mathcal C_{d'}\) has the main-chart formula
\begin{equation}
 \boxed{\qquad
 u'=\frac{D}{N-D},\qquad
 V'=2v'+1=V\frac{J}{c_KAB}.
 \qquad}
\label{eq:Cd-full-product-derivative}
\end{equation}
Two independent inversion-free output pairs are
\begin{equation}
       (u'_0:u'_1)=(D:N-D),\qquad
       (V'_0:V'_1)=(VJ:c_KAB).
\label{eq:Cd-full-productive-output}
\end{equation}
\end{theorem}

\begin{proof}
The proof has five explicit steps.

\smallskip
\noindent\emph{Step 1: determine the possible second coordinate.}
On \(M_d\), negation fixes \(x\) and sends \(y\) to \(-y\).  An isogeny
commutes with negation, so once \(x'=f(x)\) is fixed, its second coordinate
has the form
\[
                    y'=\lambda\,y f'(x)
\]
for a nonzero constant \(\lambda\).  To verify this standard shape without
assuming it, write the invariant differentials as
\[
 \omega=\frac{dx}{2B_dy},\qquad
 \omega'=\frac{dx'}{2B_{d'}y'}.
\]
For a separable isogeny, \(\varphi^*\omega'=\alpha\omega\) for one
\(\alpha\in k^\times\).  Since \(dx'=f'(x)dx\),
\[
 \frac{f'(x)dx}{2B_{d'}y'}
 =\alpha\frac{dx}{2B_dy},
\]
and therefore
\[
 y'=\frac{B_d}{\alpha B_{d'}}\,y f'(x).
\]
The prefactor is a nonzero constant, proving the asserted shape.

\smallskip
\noindent\emph{Step 2: fix the marked normalization.}
The point \(x=1,y=1\) is the selected point of order four in the
Montgomery chart: indeed
\[
 B_d\cdot1^2=1+A_d+1=B_d.
\]
An odd-degree quotient preserves its exact order.  Choose the target
orientation so that its image is the target marked point
\((x',y')=(1,1)\).  The Kummer formula has \(f(1)=1\), because each
factor at \(x=1\) is
\[
 \left(\frac{x_i-1}{1-x_i}\right)^2=1.
\]
Consequently \(1=\lambda f'(1)\), so
\[
                  \lambda=\frac1{f'(1)}.
\]
This is the only normalization used below.

\smallskip
\noindent\emph{Step 3: express the ordinate ratio through a logarithmic
derivative.}
Since \(y=xV\), one has
\[
 V'=\frac{y'}{x'}
 =V\,\frac{x f'(x)}{f(x)f'(1)}.
\]
Let
\[
 F(u)=f\!\left(\frac{u+1}{u}\right)
     =\frac{u+1}{u}\left(\frac AB\right)^2.
\]
Because \(dx/du=-u^{-2}\),
\begin{align*}
 H(u):=\frac{x f'(x)}{f(x)}
 &= -u(u+1)\frac{F'(u)}{F(u)},\\
 \frac{F'}F
 &=\frac1{u+1}-\frac1u
    +2\left(\frac{\dot A}{A}-\frac{\dot B}{B}\right).
\end{align*}
Multiplication by \(-u(u+1)\) gives
\[
 H
 =1-2u(u+1)
       \left(\frac{\dot A}{A}-\frac{\dot B}{B}\right)
 =\frac{AB-2u(u+1)(\dot A B-A\dot B)}{AB}
 =\frac J{AB}.
\]

\smallskip
\noindent\emph{Step 4: compute \(f'(1)\) without a limit ambiguity.}
As \(u\) tends to infinity,
\begin{align*}
 A(u)
 &=u^m\left(1+\frac{\sum_i(u_i+1)}u+O(u^{-2})\right),\\
 B(u)
 &=u^m\left(1-\frac{\sum_i u_i}u+O(u^{-2})\right).
\end{align*}
Thus
\[
 \frac AB
 =1+\frac{\sum_i(2u_i+1)}u+O(u^{-2})
\]
and
\[
 F(u)
 =1+\frac{1+2\sum_i(2u_i+1)}u+O(u^{-2})
 =1+\frac{c_K}{u}+O(u^{-2}).
\]
Also \(x=1+u^{-1}\).  Therefore
\[
 f'(1)=\lim_{u\to\infty}\frac{F(u)-1}{x(u)-1}=c_K.
\]
The point of order four is not a branch point of the degree-two Kummer
map.  A separable isogeny is \'etale, and its image again has order four;
hence the induced local Kummer derivative \(f'(1)\) is nonzero.  This
proves \(c_K\ne0\).

\smallskip
\noindent\emph{Step 5: assemble the native output.}
Substitution of \(H=J/(AB)\) and \(f'(1)=c_K\) gives the second equation
in \eqref{eq:Cd-full-product-derivative}.  The first is
Theorem~\ref{thm:Cd-native-odd-xeval}.  Clearing the two denominators
separately gives \eqref{eq:Cd-full-productive-output}.
\end{proof}

\subsection{Simultaneous products and derivatives}

The products and their derivatives are evaluated without expanding a
kernel polynomial.  For \(a_i=u+u_i+1\), update
\begin{equation}
       (A,\dot A)\longleftarrow
       (Aa_i,\dot A\,a_i+A);
\label{eq:Cd-product-derivative-recurrence-A}
\end{equation}
for \(b_i=u-u_i\), use the identical recurrence for
\((B,\dot B)\).  The right-hand \(A\) or \(B\) in each derivative update
is the value before the update.  The product rule verifies the recurrence:
\[
       \frac d{du}(Aa_i)=\dot A\,a_i+A\frac{da_i}{du}
                        =\dot A\,a_i+A.
\]
Starting with the first factor and derivative \(1\), the four quantities
\(A,\dot A,B,\dot B\) cost exactly \(4(m-1)\M\).

\subsection{Homogeneous main chart and projective completion}

Write \(u=(U_0:U_1)\) and define
\[
 A_h=\prod_i\bigl(U_0+(u_i+1)U_1\bigr),\qquad
 B_h=\prod_i(U_0-u_iU_1),
\]
\[
 A_0=\frac{\partial A_h}{\partial U_0},\qquad
 B_0=\frac{\partial B_h}{\partial U_0}.
\]
The homogeneous polynomial
\[
 J_h=U_1A_hB_h
 -2U_0(U_0+U_1)(A_0B_h-A_hB_0)
\]
specializes to \(J\) at \(U_1=1\).  Put
\[
 N_h=(U_0+U_1)A_h^2,\qquad D_h=U_0B_h^2.
\]
For \(V=(V_0:V_1)\), the homogeneous main chart is
\begin{align}
 (u'_0:u'_1)&=(D_h:N_h-D_h),
 \label{eq:Cd-full-homogeneous-u}\\
 (V'_0:V'_1)&=
 (V_0J_h:V_1c_KU_1A_hB_h).
 \label{eq:Cd-full-homogeneous-V}
\end{align}
Every pair is bihomogeneous in its own projective input.  At a base point
of the displayed pair, use the full centered-coordinate product map
\eqref{eq:velu-r-thesis}--\eqref{eq:velu-s-thesis}.  Both constructions
agree with \eqref{eq:Cd-full-product-derivative} on a dense open set.
Two morphisms from a smooth projective curve to a separated projective
curve that agree on a dense open set agree everywhere.  Hence this
fallback is not a different isogeny: it is an additional chart of the
same global morphism.  In a fixed-time implementation, all required
charts can be evaluated and selected by zero masks.

\begin{theorem}[Exact full-point main-chart cost]
\label{thm:Cd-full-point-cost}
For affine \(u,V,u_i\), precomputed \(c_K\), and the two projective
outputs in \eqref{eq:Cd-full-productive-output}, the cost is
\begin{equation}
             \boxed{(4m+4)\M+2\Sqr+\mathbf C_K,}
\label{eq:Cd-full-point-cost-C}
\end{equation}
where \(\mathbf C_K\) is one multiplication by the precomputed
kernel-dependent constant \(c_K\).  If \(\mathbf C_K\) is charged as a
general multiplication, the cost is
\begin{equation}
             \boxed{(4m+5)\M+2\Sqr.}
\label{eq:Cd-full-point-cost-M}
\end{equation}
A projective input \(V_0:V_1\) adds one \(\M\).  Affine normalization of
both outputs adds \(5\M+\Inv\).
\end{theorem}

\begin{proof}
The simultaneous recurrences cost \(4(m-1)\M\).
The values \(N,D\) cost \(2\M+2\Sqr\).  To form \(J\), compute
\[
 AB,\quad \dot A B,\quad A\dot B,\quad u(u+1),\quad
 u(u+1)(\dot A B-A\dot B),
\]
using \(5\M\).  The numerator \(VJ\) costs \(1\M\), and the denominator
\(c_KAB\) costs \(\mathbf C_K\).  Therefore the number of general
multiplications is
\[
 4(m-1)+2+5+1=4m+4.
\]
This proves \eqref{eq:Cd-full-point-cost-C}; charging
\(\mathbf C_K=\M\) proves \eqref{eq:Cd-full-point-cost-M}.
For two affine outputs, one simultaneous inversion uses one denominator
product, two multiplications to split its inverse, and two numerator
multiplications, for \(5\M+\Inv\).
\end{proof}

\section{The optimized separable three-isogeny}
\label{sec:Cd-optimized-three-isogeny}

The degree-three case should not be charged by blindly substituting
\(m=1\) in a generic derivative recurrence.  Let
\[
 Q=(q,v_Q),\qquad [3]Q=O,\qquad Q\ne O,
\]
and put
\[
 a=u+q+1,\qquad b=u-q,\qquad
 c=2q+1,\qquad c_3=1+2c=4q+3,\qquad t=u(u+1).
\]

\begin{theorem}[Optimized native \(C_d\) three-isogeny]
\label{thm:Cd-optimized-three-isogeny}
Assume \(\charac k\ne2,3\), \(d\rho\ne0\), and
\(\langle Q\rangle\) is Galois stable.  Define
\[
 N=(u+1)a^2,\qquad D=ub^2,\qquad
 J=ab+2ct,\qquad V=2v+1.
\]
Then the normalized degree-three quotient has
\begin{equation}
 \boxed{\qquad
 (u'_0:u'_1)=(D:N-D),\qquad
 (V'_0:V'_1)=(VJ:c_3ab).
 \qquad}
\label{eq:Cd-optimized-three-isogeny-full}
\end{equation}
Its target parameter is represented without inversion by
\begin{equation}
 \boxed{\qquad
 (\rho'_0:\rho'_1)=(\rho^3:c^8),\qquad
 (d'_0:d'_1)=(c^8-\rho^3:16c^8).
 \qquad}
\label{eq:Cd-optimized-three-isogeny-parameter}
\end{equation}
\end{theorem}

\begin{proof}
For \(m=1\), the products in
Theorem~\ref{thm:Cd-full-product-derivative} are \(A=a\), \(B=b\), and
\(\dot A=\dot B=1\).  Hence
\[
 \dot A B-A\dot B=b-a=-(2q+1)=-c.
\]
Substitution in \eqref{eq:Cd-full-J-cK} gives
\[
 J=ab-2u(u+1)(-c)=ab+2ct,
\]
and the normalization constant is
\[
 c_K=1+2(2q+1)=c_3.
\]
This proves \eqref{eq:Cd-optimized-three-isogeny-full}.  The general
kernel-product formula has \(R_K=2q+1=c\) and \(\ell=3\); therefore
\(\rho'=\rho^3/c^8\).  Clearing the denominator and then applying
\(d'=(1-\rho')/16\) proves
\eqref{eq:Cd-optimized-three-isogeny-parameter}.
\end{proof}

\begin{proposition}[Exact degree-three costs]
\label{prop:Cd-optimized-three-isogeny-cost}
With affine \(u,V,q\) and two projective output pairs:
\begin{enumerate}[label=\textup{(\roman*)}]
 \item Kummer evaluation alone costs
       \[
                         \boxed{2\M+2\Sqr}.
       \]
 \item Full-point evaluation has the two equivalent schedules
       \[
       \boxed{5\M+2\Sqr+2\Cmul}
       \quad\hbox{or}\quad
       \boxed{4\M+3\Sqr+2\Cmul},
       \]
       according as \(t\) is formed by one multiplication or by
       \(t=u^2+u\).
 \item The projective target parameter in
       \eqref{eq:Cd-optimized-three-isogeny-parameter} costs
       \[
                         \boxed{\M+4\Sqr}.
       \]
 \item Consequently, Kummer image plus target parameter costs
       \[
                         \boxed{3\M+6\Sqr}.
       \]
\end{enumerate}
\end{proposition}

\begin{proof}
The Kummer row consists of the two squares \(a^2,b^2\) and the two final
multiplications by \(u+1,u\).  For the full row, compute \(ab\) once,
form \(t\) by either \(\M\) or \(\Sqr\), multiply \(t\) by the
kernel constant \(c\), form \(N,D\), multiply \(VJ\), and multiply
\(ab\) by \(c_3\).  Excluding the chosen realization of \(t\), this is
\(4\M+2\Sqr+2\Cmul\), proving both schedules.  Finally
\(\rho^3\) costs one square and one multiplication, while
\(c^8\) costs three successive squarings.  The two projective
coordinates require no final multiplication or inversion.
\end{proof}

The strongest published Edwards comparators must be separated by output.
Moody--Shumow record a degree-three projective full-point specialization
of
\[
                       5\M+4\Sqr+3\Cmul
\]
\cite{MoodyShumow2016}.  Against the first \(C_d\) full-point schedule,
the exact Edwards-minus-\(C_d\) difference is
\[
                         \boxed{2\Sqr+\Cmul};
\]
against the square realization of \(t\), it is
\(\M+\Sqr+\Cmul\).  Both differences are positive for every
positive operation pricing.

Kim--Yoon--Kwon--Hong give an Edwards \(YZ\)-Kummer image together with
projective image-curve coefficients in \(6\M+5\Sqr\)
\cite{KimYoonKwonHong2018}.  The equal-output \(C_d\) row is
\(3\M+6\Sqr\); it is cheaper exactly when
\begin{equation}
                         \operatorname{cost}(\Sqr)
                         <3\operatorname{cost}(\M).
\label{eq:Cd-Edwards-three-isogeny-threshold}
\end{equation}
This condition includes every conventional prime-field backend in which a
square is no more expensive than a general multiplication.  The comparison
does not use the Edwards \(YZ\)-only result as a full-point result.

\begin{table}[H]
\centering
\small
\caption{Equal-output degree-three comparison}
\label{tab:Cd-Edwards-three-isogeny-cost}
\begin{tabular}{L{4.1cm}L{3.4cm}L{4.9cm}}
\toprule
interface & exact cost & conclusion\\
\midrule
\(C_d\) Kummer point only
& \(2\M+2\Sqr\)
& native monic-linear specialization\\
\(C_d\) Kummer plus projective \(d'\)
& \(3\M+6\Sqr\)
& no inversion; one kernel coordinate \(q\)\\
Edwards \(YZ\)-Kummer plus curve coefficients
& \(6\M+5\Sqr\)
& \(C_d\) wins under
  \eqref{eq:Cd-Edwards-three-isogeny-threshold}\\
\(C_d\) full point, multiplication schedule
& \(5\M+2\Sqr+2\Cmul\)
& saves \(2\Sqr+\Cmul\) against the next row\\
published Edwards full-point specialization
& \(5\M+4\Sqr+3\Cmul\)
& strongest like-for-like full-point baseline used here\\
\bottomrule
\end{tabular}
\end{table}

\section[Full-point comparison with Edwards]
{Full-point comparison with the published Edwards formula}
\label{sec:Cd-Edwards-full-comparison}

Moody and Shumow give the projective Edwards full-point upper bound
\begin{equation}
           (3m+3)\M+4\Sqr+3m\Cmul
\label{eq:Moody-Shumow-published-cost}
\end{equation}
for an odd kernel of order \(2m+1\)
\cite{MoodyShumow2016}.  It was presented as a general formula rather
than an implementation lower bound.  Let
\[
 \sigma=\frac{\operatorname{cost}(\Sqr)}
              {\operatorname{cost}(\M)},\qquad
 \chi=\frac{\operatorname{cost}(\Cmul)}
             {\operatorname{cost}(\M)}.
\]
Comparing \eqref{eq:Cd-full-point-cost-C}, with its one
kernel-constant multiplication also costed by \(\chi\), against
\eqref{eq:Moody-Shumow-published-cost} gives
\begin{equation}
 \boxed{\quad
 C_d\text{ is no more expensive precisely when }
 (3m-1)\chi+2\sigma\ge m+1.
 \quad}
\label{eq:Cd-Edwards-cost-threshold}
\end{equation}
Indeed, subtracting the normalized \(C_d\) cost
\(4m+4+2\sigma+\chi\) from the Edwards cost
\(3m+3+4\sigma+3m\chi\) gives
\[
       (3m-1)\chi+2\sigma-(m+1).
\]

When a kernel-dependent field constant is an arbitrary field element,
take \(\chi=1\).  The two costs become
\[
 C_d:\ (4m+5)\M+2\Sqr,\qquad
 \text{Edwards baseline}:\ (6m+3)\M+4\Sqr,
\]
and the Edwards-minus-\(C_d\) difference is
\begin{equation}
                    (2m-2)\M+2\Sqr.
\label{eq:Cd-Edwards-symbolic-difference}
\end{equation}
Thus the \(C_d\) main chart has the smaller symbolic count for every
\(m\ge1\) under this cost assignment.  When a platform gives the
\(3m\) Edwards constants a special low-cost representation, its measured
\(\chi\) must instead be inserted in
\eqref{eq:Cd-Edwards-cost-threshold}.

\begin{table}[H]
\centering
\small
\caption{Degree-by-degree full-point counts when \(\Cmul=\M\)}
\label{tab:Cd-full-cost-by-degree}
\begin{tabular}{rrrrr}
\toprule
\(m\)&\(\ell\)&\(C_d\) \(u\)-Kummer&
\(C_d\) full point&published Edwards baseline\\
\midrule
1&3&\(2\M+2\Sqr\)&\(9\M+2\Sqr\)&\(9\M+4\Sqr\)\\
2&5&\(4\M+2\Sqr\)&\(13\M+2\Sqr\)&\(15\M+4\Sqr\)\\
3&7&\(6\M+2\Sqr\)&\(17\M+2\Sqr\)&\(21\M+4\Sqr\)\\
4&9&\(8\M+2\Sqr\)&\(21\M+2\Sqr\)&\(27\M+4\Sqr\)\\
8&17&\(16\M+2\Sqr\)&\(37\M+2\Sqr\)&\(51\M+4\Sqr\)\\
16&33&\(32\M+2\Sqr\)&\(69\M+2\Sqr\)&\(99\M+4\Sqr\)\\
32&65&\(64\M+2\Sqr\)&\(133\M+2\Sqr\)&\(195\M+4\Sqr\)\\
64&129&\(128\M+2\Sqr\)&\(261\M+2\Sqr\)&\(387\M+4\Sqr\)\\
\bottomrule
\end{tabular}
\end{table}

\section[Cd-w and square-root V\'elu]
{\(C_d\)-\(w\), generalized Montgomery coordinates, and
square-root V\'elu}
\label{sec:Cd-w-sqrt-Velu}

\subsection[Zero-loss transfer of the Edwards w backend]
{Zero-loss transfer of the Edwards \(w\) backend}

\begin{theorem}[\(C_d\)-\(w\) transfer]
\label{thm:Cd-w-zero-loss-transfer}
Consider an odd-degree isogeny algorithm whose field operations use only
the Edwards parameter \(\rho\), projective \(w\)-coordinates, homogeneous
differential polynomials, product trees, remainder trees, and resultants.
After the substitution \(w=w_C\) from
\eqref{eq:Cd-crypto-w}, the algorithm is a \(C_d\) algorithm with the
same sequence of field operations.  If its output parameter is
\(\rho'\), the native target parameter is obtained by
\[
                         d'=\frac{1-\rho'}{16}.
\]
\end{theorem}

\begin{proof}
Proposition~\ref{prop:Cd-crypto-Edwards-w} gives an isomorphism of the
smooth source curves and identifies \(w_C\) with the pullback of the
Edwards \(w\)-function.  Pullback by an isomorphism preserves products,
homogeneous identities, divisors, resultants, and the Kummer
differential relation.  The projective pair
\eqref{eq:Cd-crypto-w-projective} represents the same function value
used by the Edwards algorithm; hence every multiplication, squaring,
addition, product-tree node, and remainder-tree node has the same
operands after relabeling.  The only output conversion is the linear
parameter identity \(\rho'=1-16d'\), proving the claim.
\end{proof}

In the normalization of Kim--Yoon--Park--Hong
\cite{KimEtAl2019}, let
\(\eta_i=(w_C(Q_i)+1)/2\).  Their image-parameter formula transfers as
\begin{equation}
              \rho'=\rho^\ell\prod_{i=1}^m\eta_i^8.
\label{eq:Cd-w-image-parameter}
\end{equation}
Comparison with \eqref{eq:Cd-image-parameter-affine} gives the useful
kernel identity
\begin{equation}
 \prod_{i=1}^m\eta_i^8
 =\left(\prod_{i=1}^m(2u_i+1)\right)^{-8}.
\label{eq:Cd-w-u-kernel-identity}
\end{equation}
This is an equality of the complete kernel products; it does not assert a
termwise equality between \(\eta_i\) and \((2u_i+1)^{-1}\).

The Edwards-\(w\) literature has already integrated this coordinate with
CSIDH, CTIDH, an Elligator-type sampling step, and square-root V\'elu;
the reported Edwards implementation is comparable to, or slightly faster
than, its selected Montgomery comparator
\cite{MoriyaOnukiTakagi2023}.  Theorem
\ref{thm:Cd-w-zero-loss-transfer} imports that entire \(w\)-arithmetic
layer into \(C_d\), while the \(u\)-Kummer and full-point paths remain
available as additional native choices.

The generalized Montgomery-coordinate framework unifies Montgomery
\(x\), Montgomery-minus \(x\), Edwards \(w\), Huff \(w\), and twisted
Jacobi \(\omega\), together with isogeny evaluation and image-coefficient
recovery \cite{GeneralizedMontgomery2022}.  The \(C_d\) model carries two
of these optimized states at once:
\[
 \kappa_d=(u+1:u)
 \quad\hbox{is Montgomery \(x\),}\qquad
 w_C=\frac{\rho}{(2u+1)^2(2v+1)^2}
 \quad\hbox{is Edwards \(w\).}
\]
Both are native functions on the same \(C_d\) function field.  This is
why a \(C_d\) implementation can select a low-degree \(u\) product or a
high-degree \(w\) product without changing its curve object.

\subsection{Exact square-root-V\'elu parameter formulas}

Let \(m=(\ell-1)/2\), and choose an index system
\((I,J,K)\) for the unsigned kernel in the sense of the square-root
V\'elu construction \cite{BernsteinDeFeoLerouxSmith2020}.  Write the
source parameter as \(\rho=D/C\), and let
\[
 h_S(T_0,T_1)=\prod_{s\in S}(T_0Z_s-T_1W_s)
\]
for the projective \(w_C\)-points \((W_s:Z_s)\).  Let
\(E_J(-1,1,T)\) be the standard giant-step polynomial obtained from the
homogeneous \(w\)-differential addition polynomial.  With the
normalizations of Takahashi--Onuki--Takagi
\cite{TakahashiOnukiTakagi2022}, the exact resultants are
\begin{align}
 D'&=D^{\,2\#K+1}
 \left[
 h_K(-1,1)
 \operatorname{Res}_T
 \bigl(h_I(T,1),E_J(-1,1,T)\bigr)
 \right]^8,
 \label{eq:Cd-sqrt-Velu-D}\\
 C'&=C^\ell D^{\,4\#I\#J}
 \left[
 2^m h_K(1,0)
 \operatorname{Res}_T
 \bigl(h_I(T,1),h_J(1,T)\bigr)^2
 \right]^8.
 \label{eq:Cd-sqrt-Velu-C}
\end{align}
The target values are
\begin{equation}
       \rho'=\frac{D'}{C'},\qquad
       (d'_0:d'_1)=(C'-D':16C').
\label{eq:Cd-sqrt-Velu-target}
\end{equation}

The two resultants in
\eqref{eq:Cd-sqrt-Velu-D}--\eqref{eq:Cd-sqrt-Velu-C} have different
second arguments and generally acquire different scaling factors in a
pseudo-remainder tree.  This can be handled without an inversion.
Suppose an implementation returns
\[
 R_1=\gamma_1\operatorname{Res}_1,\qquad
 R_2=\gamma_2\operatorname{Res}_2,
\]
and tracks the nonzero scale factors \(\gamma_1,\gamma_2\).  Define
\[
 \widehat D=
 D^{\,2\#K+1}[h_K(-1,1)R_1]^8,\qquad
 \widehat C=
 C^\ell D^{\,4\#I\#J}
 [2^m h_K(1,0)R_2^2]^8.
\]
Since the true pair is
\[
 \left(\frac{\widehat D}{\gamma_1^8}:
       \frac{\widehat C}{\gamma_2^{16}}\right),
\]
an exact division-free representative is
\begin{equation}
             (D':C')=
       (\widehat D\,\gamma_2^{16}:
        \widehat C\,\gamma_1^8).
\label{eq:Cd-sqrt-Velu-scaled-resultants}
\end{equation}
Equation~\eqref{eq:Cd-sqrt-Velu-scaled-resultants} follows by multiplying
both projective coordinates by
\(\gamma_1^8\gamma_2^{16}\); it shows explicitly why the two scale
factors must not be silently identified.

With fast polynomial multiplication, the index-system construction gives
\[
 \widetilde O(\sqrt\ell)
\]
field operations, in contrast to the \(O(\ell)\) linear products.
The redundant-integer index system of Otsuki--Onuki--Takagi improves the
square-root-V\'elu schedule for \(44\%\) of the 367 tested primes between
97 and 2689 and reports about \(6.6\%\) improvement at 1279 and 2687
\cite{OtsukiOnukiTakagi2023}.  By
Theorem~\ref{thm:Cd-w-zero-loss-transfer}, the same index selection and
operation sequence apply to the \(C_d\)-\(w\) backend.

\section{Separable two-isogenies in odd characteristic}
\label{sec:Cd-odd-characteristic-two-isogeny}

Even degree is structurally different from the odd product formulas.  The
marked point of order four on \(C_d\) doubles to the rational two-torsion
point \(T_0=(0,0)\) on \(M_d\).  Quotienting by \(T_0\) destroys that
chosen order-four marking, so unconditional closure in the one-parameter
\(C_d\) family must not be asserted.

\begin{theorem}[Root-free marked two-isogeny]
\label{thm:Cd-odd-marked-two-isogeny}
Assume \(\charac k\ne2\) and \(d\rho\ne0\).  On
\[
 M_d:\qquad B_dy^2=x^3+A_dx^2+x,\qquad
 B_d=\frac1{4d},\quad A_d=B_d-2,
\]
put
\[
 C_d^{(2)}=A_d^2-4=\frac{\rho}{16d^2}.
\]
Then
\begin{equation}
 \boxed{\qquad
 X=x+A_d+\frac1x,\qquad
 Y=y\left(1-\frac1{x^2}\right)
 \qquad}
\label{eq:Cd-odd-two-isogeny-M}
\end{equation}
defines a separable two-isogeny with kernel
\(\{O,T_0\}\) and root-free codomain
\begin{equation}
 M_d^{(2)}:\qquad
 B_dY^2=X^3-2A_dX^2+C_d^{(2)}X.
\label{eq:Cd-odd-two-isogeny-target}
\end{equation}
For a projective source Kummer pair \(x=K_0/K_1\), its induced map is
\begin{equation}
 \boxed{\quad
 (K'_0:K'_1)=
 \bigl((K_0+K_1)^2+(A_d-2)K_0K_1:K_0K_1\bigr).
 \quad}
\label{eq:Cd-odd-two-isogeny-projective}
\end{equation}

For the raw native pair \((K_0:K_1)=(u+1:u)\), put
\[
             t=u(u+1),\qquad U=2u+1,\qquad V=2v+1.
\]
Then
\begin{align}
 (K'_0:K'_1)&=(t+4d:4dt),
 \label{eq:Cd-odd-two-isogeny-native-Kummer}\\
 (Y_0:Y_1)&=(UV:t).
 \label{eq:Cd-odd-two-isogeny-native-Y}
\end{align}
Equations
\eqref{eq:Cd-odd-two-isogeny-native-Kummer} and
\eqref{eq:Cd-odd-two-isogeny-native-Y} are two independent projective
ratios representing the full point on
\eqref{eq:Cd-odd-two-isogeny-target}.
\end{theorem}

\begin{proof}
The standard quotient of
\(By^2=x^3+Ax^2+bx\) by \((0,0)\) is
\[
 X=x+A+\frac b x,\qquad
 Y=y\left(1-\frac b{x^2}\right),
\]
with codomain
\[
 BY^2=X^3-2AX^2+(A^2-4b)X.
\]
The identity can also be checked without invoking a template:
\begin{align*}
 &B\,y^2\left(1-\frac1{x^2}\right)^2\\
 &\quad=(x^3+A_dx^2+x)\frac{(x^2-1)^2}{x^4}\\
 &\quad=\left(x+A_d+\frac1x\right)^3
       -2A_d\left(x+A_d+\frac1x\right)^2
       +(A_d^2-4)\left(x+A_d+\frac1x\right).
\end{align*}
The only finite pole is at \(x=0\), so \(O,T_0\) map to the target
identity.  The rational function \(X\) has degree two and its derivative
\[
                       \frac{dX}{dx}=1-\frac1{x^2}
\]
is not zero because \(\charac k\ne2\).  The map is therefore separable
of degree two and has exactly the stated kernel.

Homogenizing
\[
 X=\frac{x^2+A_dx+1}{x}
\]
gives
\[
 (K_0^2+A_dK_0K_1+K_1^2:K_0K_1).
\]
Since
\[
 K_0^2+A_dK_0K_1+K_1^2
 =(K_0+K_1)^2+(A_d-2)K_0K_1,
\]
this is \eqref{eq:Cd-odd-two-isogeny-projective}.

Finally \(x=(u+1)/u\) gives
\begin{align*}
 x+A_d+x^{-1}
 &=\frac{2t+1+A_dt}{t}
   =\frac{t+4d}{4dt},\\
 y(1-x^{-2})
 &=V(x-x^{-1})
   =\frac{UV}{t}.
\end{align*}
Clearing the two denominators separately proves
\eqref{eq:Cd-odd-two-isogeny-native-Kummer}--%
\eqref{eq:Cd-odd-two-isogeny-native-Y}.
\end{proof}

The root-free codomain is already a marked-two-torsion
generalized-Montgomery endpoint.  It also has a root-free QRT-Jacobian
normalization.  Set
\[
 \widetilde X=B_dX,\qquad \widetilde Y=B_d^2Y.
\]
Then \eqref{eq:Cd-odd-two-isogeny-target} becomes
\begin{equation}
 \widetilde Y^2
 =\widetilde X
  \left(\widetilde X^2-2q_2\widetilde X+\Omega_2\right),
 \qquad
 q_2=A_dB_d,\quad
 \Omega_2=C_d^{(2)}B_d^2.
\label{eq:Cd-odd-two-isogeny-QRT-endpoint}
\end{equation}
Thus the QRT/reciprocal Jacobian layer records the quotient without a
square root.  Returning the raw genus-one QRT equation over \(k\) is then
governed by the exact re-embedding criterion
\eqref{eq:QRT-reembedding-qOmega}--%
\eqref{eq:QRT-reembedding-beta-square}.

\begin{theorem}[Return to the one-sided and pure \(C_d\) families]
\label{thm:Cd-odd-two-isogeny-reembedding}
In addition to the hypotheses of
Theorem~\ref{thm:Cd-odd-marked-two-isogeny}, suppose
\(s\in k^\times\) satisfies \(s^2=\rho\).  Define
\begin{equation}
 \alpha_2=-\frac{16s^2}{(1+s)^2},\qquad
 \Delta_2=-\frac{16s^2}{(1-s)^2},
 \label{eq:Cd-odd-two-isogeny-alpha-delta}
\end{equation}
\[
 a_2=\frac{1-\alpha_2}{4},\qquad
 d_2=\frac{\alpha_2-\Delta_2}{16}
     =\frac{4s^3}{(1-s^2)^2}.
\]
Then \eqref{eq:Cd-odd-two-isogeny-target} is isomorphic over \(k\) to
\[
                \mathcal T_{a_2,d_2}.
\]
For a source \(C_d\) point, an inversion-free native target is
\begin{align}
 (u_{2,0}:u_{2,1})
 &=\bigl(sUV-(t+4d):2(t+4d)\bigr),
 \label{eq:Cd-odd-two-isogeny-T-u}\\
 (v_{2,0}:v_{2,1})
 &=\bigl(st:4d+(1-s)t\bigr).
 \label{eq:Cd-odd-two-isogeny-T-v}
\end{align}

If, moreover, \(\imath^2=-1\) in \(k\), then
\(\alpha_2=(4\imath s/(1+s))^2\), and a final constant scaling returns
to the pure family \(C_{d_2^{C}}\), where
\begin{equation}
 \rho_2=\frac{\Delta_2}{\alpha_2}
       =\left(\frac{1+s}{1-s}\right)^2,\qquad
 d_2^{C}=\frac{1-\rho_2}{16}
       =-\frac{s}{4(1-s)^2}.
\label{eq:Cd-odd-two-isogeny-pure-Cd-target}
\end{equation}
\end{theorem}

\begin{proof}
Put \(r=s/(4d)\).  Since \(s^2=\rho\),
\[
 r^2=\frac{\rho}{16d^2}=C_d^{(2)}.
\]
Set \(x_2=X/r\) and \(y_2=Y\).  Division of
\eqref{eq:Cd-odd-two-isogeny-target} by \(r^3\) gives
\[
 \frac{B_d}{r^3}y_2^2
 =x_2^3-\frac{2A_d}{r}x_2^2+x_2.
\]
The Montgomery dictionary of \(\mathcal T_{a_2,d_2}\) requires
\[
 B_2=\frac1{4d_2},\qquad
 A_2=\frac{\alpha_2}{4d_2}-2.
\]
Direct substitution of \(d=(1-s^2)/16\) verifies
\[
 \frac1{4d_2}=\frac{B_d}{r^3},\qquad
 \frac{\alpha_2}{4d_2}-2=-\frac{2A_d}{r},
\]
and
\[
 \alpha_2-16d_2=\Delta_2.
\]
This proves the target equation.

For the point map, the inverse
\(\mathcal T_{a_2,d_2}\)-to-Montgomery dictionary in centered native
coordinates is
\[
 U_2=\frac{y_2}{x_2}=\frac{rY}{X},\qquad
 V_2=\frac{x_2+1}{x_2-1}=\frac{X+r}{X-r}.
\]
Using
\[
 X=\frac{t+4d}{4dt},\qquad Y=\frac{UV}{t},
\]
one obtains the projective pairs
\[
 (U_{2,0}:U_{2,1})=(sUV:t+4d),
\]
\[
 (V_{2,0}:V_{2,1})
 =\bigl(4d+(1+s)t:4d+(1-s)t\bigr).
\]
Applying \(u_2=(U_2-1)/2\) and \(v_2=(V_2-1)/2\) gives
\eqref{eq:Cd-odd-two-isogeny-T-u}--%
\eqref{eq:Cd-odd-two-isogeny-T-v}.

Finally, if \(\imath^2=-1\), then the displayed square root of
\(\alpha_2\) scales the twisted Edwards coefficient \(\alpha_2\) to one.
The remaining coefficient is \(\Delta_2/\alpha_2\), and the two
equalities in \eqref{eq:Cd-odd-two-isogeny-pure-Cd-target} follow by
expanding the two squares.  This proves every re-embedding assertion.
\end{proof}

Once \(s\) has been selected, no inversion is required merely to retain
the target coefficients.  One may store
\[
 (\alpha_{2,0}:\alpha_{2,1})
   =(-16\rho:(1+s)^2),\qquad
 (\Delta_{2,0}:\Delta_{2,1})
   =(-16\rho:(1-s)^2),
\]
and
\[
 (d_{2,0}:d_{2,1})
   =(4s\rho:(1-\rho)^2).
\]
Because \(s^2=\rho\), the two squared denominators are formed by additions
from \(1,\rho,s\); the final \(d_2\) numerator is one multiplication by
the public constant \(s\).  Affine coefficient normalization may be
batched with other target normalizations.  The square-root selection, not
an unreported inversion, is therefore the genuine re-embedding cost.

\subsection{Exact costs and the Edwards endpoint comparison}

Ignoring additions and small integer multiplications,
\eqref{eq:Cd-odd-two-isogeny-projective} costs
\[
                         \M+\Sqr+\Dpar.
\]
Indeed, form \(K_0K_1\), square \(K_0+K_1\), and multiply the product by
\(A_d-2\).  On a raw native \(u\), form
\(t=u^2+u\) and \(4dt\), giving
\begin{equation}
              \boxed{\min(\M,\Sqr)+\Dpar}
\label{eq:Cd-odd-two-isogeny-native-cost}
\end{equation}
for the Kummer output.  The root-free full output adds \(UV\), for
\[
                         \boxed{\M+\Sqr+\Dpar}.
\]
The native \(\mathcal T_{a_2,d_2}\) output in
\eqref{eq:Cd-odd-two-isogeny-T-u}--%
\eqref{eq:Cd-odd-two-isogeny-T-v} uses \(t\), \(UV\), \(st\), and
\(sUV\), and therefore costs
\[
                         \boxed{\M+\Sqr+2\Dpar}.
\]

For an Edwards source point \((\xi,\eta)\), put
\(p=1-\eta^2\).  The same root-free quotient has
\begin{align}
 (X_0:X_1)&=\bigl(4+(A_d-2)p:p\bigr),
 \label{eq:Edwards-odd-two-isogeny-X}\\
 (Y_0:Y_1)&=(4\eta:\xi p).
 \label{eq:Edwards-odd-two-isogeny-Y}
\end{align}
Thus its Kummer and full-point costs are respectively
\(\Sqr+\Dpar\) and
\(\M+\Sqr+\Dpar\): the Kummer core ties the native \(C_d\) square
schedule, and the root-free full core is the same transported operation
graph.

When \(s^2=\rho\) and a standard twisted-Edwards target is required, one
may use
\begin{align}
 (x_{2,0}:x_{2,1})
 &=\bigl(\xi(1-s^2\eta^2):4s\eta\bigr),\\
 (y_{2,0}:y_{2,1})
 &=\bigl(1-s^2\eta^2-s(1-\eta^2):
          1-s^2\eta^2+s(1-\eta^2)\bigr).
\label{eq:Edwards-odd-two-isogeny-rooted-output}
\end{align}
This schedule costs
\(\M+\Sqr+3\Dpar\).  The \(C_d\)-native one-sided endpoint saves one
curve-constant multiplication, while both models face the same
\(\sqrt\rho\) field-of-definition condition.  This is the precise
square-root obstruction observed in the Edwards two-isogeny literature
\cite{Morain2009,KimYoonKwonHong2018}; it is avoided, rather than hidden,
by retaining the root-free QRT-Jacobian endpoint
\eqref{eq:Cd-odd-two-isogeny-QRT-endpoint}.

\begin{longtable}{L{3.3cm}L{3.4cm}L{3.1cm}L{3.3cm}}
\caption{Odd-characteristic two-isogeny costs at equal endpoints}
\label{tab:Cd-odd-two-isogeny-costs}\\
\toprule
interface & declared input and output & exact cost & conclusion\\
\midrule
\endfirsthead
\toprule
interface & declared input and output & exact cost & conclusion\\
\midrule
\endhead
general projective Kummer
& \(M_d\) pair to root-free target pair
& \(\M+\Sqr+\Dpar\)
& common transported core\\
\(C_d\), native Kummer
& affine \(u\) to root-free target pair
& \(\min(\M,\Sqr)+\Dpar\)
& fixed marked kernel and one native quadratic \(t\)\\
Edwards, affine Kummer
& affine \(\eta\) to the same root-free pair
& \(\Sqr+\Dpar\)
& ties the \(C_d\) square schedule\\
\(C_d\), native full
& affine \(u,v\) to two root-free target ratios
& \(\M+\Sqr+\Dpar\)
& direct \(t,UV\) circuit\\
Edwards, native full
& affine \(\xi,\eta\) to the same root-free target
& \(\M+\Sqr+\Dpar\)
& exact transport equality\\
\(C_d\to\mathcal T_{a_2,d_2}\)
& native full source and native full target
& \(\M+\Sqr+2\Dpar\)
& requires \(s^2=\rho\)\\
Edwards \(\to E_{\alpha_2,\Delta_2}\)
& standard affine source to projective full target
& \(\M+\Sqr+3\Dpar\)
& same root; one additional constant multiplication in the displayed
  native schedule\\
root-free target coefficients
& generalized Montgomery or QRT-Jacobian
& additions and public-constant products
& no square root and no kernel enumeration\\
pure \(C_d\)/ordinary Edwards target
& parameter and endpoint normalization
& \(s^2=\rho\), \(\imath^2=-1\)
& exactly the conditions in
  \eqref{eq:Cd-odd-two-isogeny-pure-Cd-target}\\
\bottomrule
\end{longtable}

If \(A_d^2-4\) is a square, the source has two additional rational
two-torsion points
\[
                  x_\pm=\frac{-A_d\pm\sqrt{A_d^2-4}}2.
\]
For either one, shift \(z=x-x_\pm\), put
\[
 A_\pm=3x_\pm+A_d,\qquad
 B_\pm=3x_\pm^2+2A_dx_\pm+1=x_\pm^2-1,
\]
and apply
\[
 X=z+A_\pm+\frac{B_\pm}{z},\qquad
 Y=y\left(1-\frac{B_\pm}{z^2}\right).
\]
A general projective \(z\)-Kummer evaluation costs
\(\M+2\Sqr+2\Dpar\).  These are not the marked-kernel specialization
\eqref{eq:Cd-odd-two-isogeny-native-Kummer}; the same two extra kernels
and the same shifted costs occur on the Edwards side under the exact
dictionary.

\section[Characteristic branches]
{Characteristic-unified statements with exact branch conditions}
\label{sec:Cd-isogeny-characteristics}

The defining \(C_d\) equation is characteristic-uniform, but a correct
isogeny theorem must retain the natural moduli parameter in each
characteristic.  In particular, \(\rho=1-16d\) becomes \(1\) in
characteristic two and therefore cannot encode the binary moduli.

\begin{theorem}[Prime-to-characteristic \(C_d\) Kummer products]
\label{thm:Cd-characteristic-branch-Kummer}
Let \(\ell=2m+1\) with \(m\ge1\), let the kernel isogeny be
separable, and choose \(Q_i=(u_i,v_i)\), \(1\le i\le m\), from the
nonzero pairs \(\{\pm Q_i\}\).
\begin{enumerate}[label=\textup{(\roman*)}]
 \item If \(\charac k\ne2\) and
       \(\gcd(\ell,\charac k)=1\), then
       \[
       F(u)=\frac{u+1}{u}
       \prod_{i=1}^m
       \left(\frac{u+u_i+1}{u-u_i}\right)^2,
       \qquad u'=\frac1{F(u)-1},
       \]
       and
       \[
       \rho'=\rho^\ell
       \prod_i(2u_i+1)^{-8}.
       \]
 \item If \(\charac k=2\), then
       \[
       F(u)=\frac{u+1}{u}
       \prod_{i=1}^m
       \left(\frac{u+u_i+1}{u+u_i}\right)^2,
       \qquad u'=\frac1{F(u)+1},
       \]
       while the moduli parameter is \(d\) and
       \[
       d'=d\prod_i\left(\frac{u_i}{u_i+1}\right)^2.
       \]
\end{enumerate}
Both branches have affine-kernel, projective-output cost
\(2m\M+2\Sqr\).
\end{theorem}

\begin{proof}
Part~(i) is Theorems~\ref{thm:Cd-native-odd-xeval} and
\ref{thm:Cd-image-parameter-product}.

For part~(ii), put \(z=(u+1)/u\) and
\(z_i=(u_i+1)/u_i\).  The characteristic-two curve equation gives
\[
 u^2+u=\frac{z}{(z+1)^2},\qquad
 v^2+v=d(z+z^{-1}).
\]
Translation by the unique nonzero two-torsion point acts on the
\(z\)-line by \(z\mapsto z^{-1}\).  The normalized quotient therefore
has a double pole at every \(z_i\), a double zero at every \(z_i^{-1}\),
and the simple zero and pole represented by the leading factor \(z\).
Thus it has the form
\[
 f(z)=c z\prod_{i=1}^m
       \left(\frac{zz_i+1}{z+z_i}\right)^2.
\]
At the marked value \(z=1\), each fraction is \(1\); the normalization
\(f(1)=1\) therefore gives \(c=1\).  Substituting
\(z=(u+1)/u\) and \(z_i=(u_i+1)/u_i\), then using
\(u'=1/(f(z)+1)\), gives the coordinate formula in part~(ii).

It remains to determine the binary parameter.  Let
\(r=\prod_i z_i\).  At \(z=0\),
\[
 f(z)=r^{-2}z+O(z^2).
\]
After pulling back the target Artin--Schreier equation, its simple-pole
term is \(d'f(z)^{-1}=d'r^2z^{-1}+O(1)\), while the source term is
\(dz^{-1}+O(1)\).  An ordinate correction changes the right side by
\(g^2+g\).  If \(g\) has a pole of order \(n>0\), the leading pole of
\(g^2+g\) has even order \(2n\); hence it cannot change a nonzero
simple-pole coefficient.  Consequently \(d'r^2=d\), and therefore
\[
 d'=d\prod_i z_i^{-2}
    =d\prod_i\left(\frac{u_i}{u_i+1}\right)^2.
\]
This proves the binary parameter formula without reducing the
odd-characteristic parameter \(\rho\).

In either characteristic, evaluate two monic products of length \(m\),
square both results, and form two final products.  This costs
\(2(m-1)\M+2\Sqr+2\M=2m\M+2\Sqr\).  For \(m=0\) the kernel is
trivial, the products are empty, the map is the identity, and the charged
cost is zero after operations on the unit constant are suppressed.
\end{proof}

\subsection{Direct characteristic-two dictionary}

Assume \(\charac k=2\) and \(d\ne0\).  Put
\[
 z=\frac{u+1}{u},\qquad Y=zv,\qquad x=dz,\qquad y=dY.
\]
Then
\begin{equation}
            y^2+xy=x^3+d^2x.
\label{eq:Cd-crypto-binary-Weierstrass}
\end{equation}
Indeed, \(u=1/(z+1)\), so
\[
 u^2+u=\frac{z}{(z+1)^2},\qquad
 v^2+v=\frac{Y^2+zY}{z^2}.
\]
The \(C_d\) equation gives
\[
 Y^2+zY=dz(z+1)^2=d(z^3+z).
\]
Multiplication by \(d^2\) and substitution of \(x,y\) prove
\eqref{eq:Cd-crypto-binary-Weierstrass}.  The inverse on the main chart is
\[
                     u=\frac d{x+d},\qquad v=\frac yx.
\]

\begin{theorem}[Binary \(C_d\) Kummer and target products]
\label{thm:Cd-binary-product-parameter}
Let \(\charac k=2\), let \(K\) be a separable odd-order kernel, and put
\[
              z_i=\frac{u_i+1}{u_i}.
\]
The normalized Kummer quotient is
\begin{equation}
 f(z)=z\prod_{i=1}^m
 \left(\frac{zz_i+1}{z+z_i}\right)^2,
\label{eq:Cd-binary-Kummer-product}
\end{equation}
and its target parameter is
\begin{equation}
 \boxed{\qquad
 d'=d\left(\prod_{i=1}^m z_i\right)^{-2}
    =d\prod_{i=1}^m
       \left(\frac{u_i}{u_i+1}\right)^2.
 \qquad}
\label{eq:Cd-binary-parameter-product}
\end{equation}
\end{theorem}

\begin{proof}
Translation by the unique nonzero two-torsion point acts on the
\(z\)-Kummer line as \(z\mapsto z^{-1}\).  The quotient map has a double
pole at each unsigned kernel value \(z_i\) and a double zero at its
two-torsion translate \(z_i^{-1}\).  It also has the simple zero and pole
carried by the leading factor \(z\).  Therefore its divisor determines it
up to a constant as
\[
 c\,z\prod_i\left(\frac{zz_i+1}{z+z_i}\right)^2.
\]
At \(z=1\), every squared fraction is \(1\), so the marked normalization
\(f(1)=1\) gives \(c=1\), proving
\eqref{eq:Cd-binary-Kummer-product}.

The source function field is the Artin--Schreier cover
\[
                 v^2+v=d(z+z^{-1}).
\]
Let \(r=\prod_i z_i\).  At \(z=0\),
\[
 f(z)=\frac{z}{r^2}+O(z^2).
\]
The target cover pulled back by \(f\) has simple-pole term
\[
 d'f(z)^{-1}=\frac{d'r^2}{z}+O(1),
\]
whereas the source has simple-pole term \(d/z\).  Their difference is
\(g^2+g\) for the ordinate correction \(v'=v+g(z)\).  If a rational
function \(g\) has a pole of order \(n>0\), then \(g^2+g\) has leading
pole order \(2n\); hence an Artin--Schreier coboundary cannot change a
nonzero simple-pole coefficient.  Therefore
\[
                       d'r^2=d,
\]
which proves the first equality in
\eqref{eq:Cd-binary-parameter-product}.  The second follows by inserting
\(z_i=(u_i+1)/u_i\).
\end{proof}

The full binary map is unambiguously specified by
\begin{equation}
 z'=f(z),\qquad v'=v+g(z),
\label{eq:Cd-binary-full-map}
\end{equation}
where
\begin{equation}
 g(z)^2+g(z)
 =d'\bigl(f(z)+f(z)^{-1}\bigr)
  +d\bigl(z+z^{-1}\bigr).
\label{eq:Cd-binary-g-equation}
\end{equation}
The two solutions differ by \(1\); the marked choice \(g(1)=0\) fixes
one.  Substitution of \eqref{eq:Cd-binary-g-equation} into
\[
 (v+g)^2+(v+g)
 =v^2+v+g^2+g
\]
gives \(v'^2+v'=d'(z'+z'^{-1})\), which verifies the target equation
directly.

\begin{corollary}[Compatibility with the binary V\'elu normalization]
\label{cor:Cd-binary-parameter-compatibility}
Use the notation of Theorem~\ref{thm:binary-velu-thesis}:
\[
 a_i=dz_i,\qquad
 \Sigma_K=\sum_i a_i,\qquad
 \tau_K^4=\Sigma_K.
\]
Then
\begin{equation}
 d+\tau_K^2+\tau_K
 =d\left(\prod_i z_i\right)^{-2}.
\label{eq:Cd-binary-two-parameter-formulas}
\end{equation}
\end{corollary}

\begin{proof}
Theorem~\ref{thm:binary-velu-thesis} constructs the quotient by \(K\)
and sends the image of the source marked order-four pair to the canonical
order-four pair on
\(W^+_{d+\tau_K^2+\tau_K}\).  Theorem
\ref{thm:Cd-binary-product-parameter} constructs the same quotient and
normalizes \(f(1)=1\), so the same marked Kummer point \(z=1\) maps to
the target marked point \(z'=1\).  In
\[
 W_e^+:\quad y^2+xy=x^3+e^2x,
\]
the \(x\)-coordinate of either point above the marked Kummer value
\(z=1\) is \(e\), because \(x=ez\).  Therefore equality of the two
normalized quotient maps forces equality of those target
\(x\)-coordinates.  The two displayed target parameters are thus equal,
which is exactly \eqref{eq:Cd-binary-two-parameter-formulas}.
\end{proof}

\begin{proposition}[Binary degree-three specialization]
\label{prop:Cd-binary-three-isogeny-cost-crypto}
Assume \(\charac k=2\), let \(Q=(q,v_Q)\) have exact order three, and
put
\[
 A_3(u)=u+q+1,\qquad B_3(u)=u+q.
\]
Then the native Kummer quotient is
\begin{equation}
 \boxed{\quad
 (u'_0:u'_1)=
 \bigl(uB_3(u)^2:
       (u+1)A_3(u)^2+uB_3(u)^2\bigr),
 \quad}
\label{eq:Cd-binary-three-isogeny-Kummer-crypto}
\end{equation}
and the target parameter has the inversion-free pair
\begin{equation}
 \boxed{\qquad
 (d'_0:d'_1)=\bigl(dq^2:(q+1)^2\bigr)
                      =\bigl(dq^2:q^2+1\bigr).
 \qquad}
\label{eq:Cd-binary-three-isogeny-parameter-crypto}
\end{equation}
Kummer evaluation costs \(2\M+2\Sqr\), parameter recovery costs
\(\Sqr+\Dpar\), and the combined output costs
\[
                         \boxed{2\M+3\Sqr+\Dpar}.
\]
\end{proposition}

\begin{proof}
Set \(m=1\) in
Theorem~\ref{thm:Cd-characteristic-branch-Kummer}.  In characteristic two
the two linear factors are \(u+q+1\) and \(u+q\), while subtraction of
the final numerator and denominator is addition.  This proves
\eqref{eq:Cd-binary-three-isogeny-Kummer-crypto} and its two-square,
two-multiplication count.  Formula
\eqref{eq:Cd-binary-parameter-product} gives
\[
                         d'=d\left(\frac q{q+1}\right)^2.
\]
Clearing the denominator proves the first parameter pair.  The Frobenius
identity \((q+1)^2=q^2+1\) shows that the two denominator squares share
the single value \(q^2\); multiplication of that value by \(d\) is one
curve-constant multiplication.  This proves the remaining costs.
\end{proof}

The proposition is the complete Kummer-plus-parameter degree-three core.
When a full binary point is required, the ordinate companion
\eqref{eq:Cd-binary-g-equation} or the explicit full wrapper
\eqref{eq:binary-native-3iso-wrapper} is additionally charged.  A binary
Edwards implementation that first enters the same \(W\)-Kummer state has
the identical intrinsic quotient; its standard endpoint conversion is a
separate interface cost, just as in the binary degree-two comparison
below.

\subsection{Characteristic two and characteristic three in the
supersingular program}

On \(W_d^+\), the point \((d,0)\) has exact order four and
\([2](d,0)=(0,0)\).  Thus a smooth \(C_d\) curve in characteristic two
has \(2\)-rank one and belongs to the ordinary locus.  This gives the
binary branch a precise role: it is a characteristic-uniform companion
for ordinary binary arithmetic, isogenies, and transport, while the
supersingular \(C_d\) cryptographic backend is naturally placed in odd
characteristic.

In characteristic three, \(2\) and \(16\) remain invertible.  All
odd-characteristic formulas in this chapter therefore apply when
\(3\nmid\ell\).  An isogeny containing a \(3\)-power inseparable part is
factored into separable and Frobenius/Verschiebung components as in
Chapter~\ref{ch:char3}.  The same statement holds in any odd
characteristic \(p\): the product formulas describe the separable
prime-to-\(p\) component, while the \(p\)-power component is represented
by Frobenius and Verschiebung.

\begin{table}[H]
\centering
\small
\setlength{\tabcolsep}{4pt}
\caption{Characteristic branches of the \(C_d\) isogeny interface}
\label{tab:Cd-isogeny-characteristic-branches}
\begin{tabular}{L{2.3cm}L{3.2cm}L{4.1cm}L{4.0cm}}
\toprule
characteristic&moduli parameter&separable isogeny formula&
cryptographic placement\\
\midrule
\(0\) or \(p>3\)&\(\rho=1-16d\)&
Theorems~\ref{thm:Cd-native-odd-xeval},
\ref{thm:Cd-full-product-derivative}, and
\ref{thm:Cd-image-parameter-product}, \(p\nmid\ell\);
Theorem~\ref{thm:Cd-odd-marked-two-isogeny} for degree two&
full supersingular odd-characteristic backend when the chosen
isogeny class has the marked four-torsion model; root-free QRT endpoint
when a two-isogenous target loses that marking\\
\(3\)&\(\rho=1-d\)&same odd formulas for \(3\nmid\ell\), including the
separable degree-two theorem&
prime-to-three steps; Frobenius/Verschiebung for the supersingular
three-primary part\\
\(2\)&\(d\)&
\eqref{eq:Cd-binary-Kummer-product}--%
\eqref{eq:Cd-binary-g-equation} for odd degree;
\eqref{eq:Cd-crypto-binary-two-isogeny-W} for degree two&
ordinary binary backend on \(C_d,\mathcal T_{a,d}\), and
\(\operatorname{Jac}(\mathcal C_{a,b,d})\)\\
\bottomrule
\end{tabular}
\end{table}

\section[Binary extension models and two-isogenies]
{Characteristic-two extension models and the separable two-isogeny}
\label{sec:Cd-crypto-binary-two-isogeny}

\subsection{Correct moduli placement}

In characteristic two the comparison model is not the odd-characteristic
Edwards equation reduced modulo two.  The relevant product models are
\[
 \mathcal C_{a,b,d}:
 (u^2+u+a)(v^2+v+b)=d,
 \qquad
 \mathcal T_{a,d}:
 (u^2+u+a)(v^2+v)=d.
\]
Their Jacobian and pointed one-sided equations are, respectively,
\begin{equation}
 W_{a+b,d}:\ y^2+xy=x^3+(a+b)x^2+d^2x,
 \qquad
 W_{a,d}:\ y^2+xy=x^3+ax^2+d^2x.
\label{eq:Cd-crypto-binary-extension-W}
\end{equation}
Both have \(c_4=1\), \(\Delta=d^4\), and \(j=d^{-4}\ne0\).  Hence every
smooth member is ordinary.  This is also the exact setting of binary
Edwards curves: their defining paper introduces them as models of ordinary
binary elliptic curves \cite{BinaryEdwards2008}, and the current complete
differential formulas retain that ordinary scope
\cite{FarashahiHosseini2023}.  Therefore no smooth member of these product
families, and no binary Edwards comparator, is being presented as a
characteristic-two supersingular equation.  Their cryptographic role here is
an ordinary binary isogeny and implementation branch; the supersingular
\(C\)-curve branch remains in odd characteristic.

The 2026 Binary Kummer Lines work gives two-way SIMD implementations for
fixed- and variable-base scalar multiplication on BKL251 and BKL257 and
connects those Kummer lines with binary Edwards endpoints
\cite{KaratiHajraSen2026}.  Here \emph{two-way} means two-lane vectorization,
not a degree-two isogeny.  It is nevertheless directly relevant to the
implementation architecture below: two independent instances of the
one-line quotient circuit can occupy the same vector lanes.

\subsection{Self-contained two-isogeny theorem}

\begin{theorem}[Binary twist-stable two-isogeny]
\label{thm:Cd-crypto-binary-two-isogeny}
Let \(k\) be perfect of characteristic two, \(d\ne0\), and \(e^2=d\).
For every \(A\in k\), the formulas
\begin{equation}
 \boxed{\qquad
 X=x+\frac{d^2}{x},\qquad
 Y=y+\frac{d^2(y+x)}{x^2}+d
 \qquad}
\label{eq:Cd-crypto-binary-two-isogeny-W}
\end{equation}
define a separable degree-two isogeny
\[
 \Phi_{A,d}:W_{A,d}\longrightarrow W_{A,e},
 \qquad
 \ker\Phi_{A,d}=\{O,(0,0)\}.
\]
For the Kummer coordinate \(z=x/d=K_0/K_1\), the induced map is
\begin{equation}
 \boxed{\qquad
 (K'_0:K'_1)=
 \bigl(e(K_0+K_1)^2:K_0K_1\bigr).
 \qquad}
\label{eq:Cd-crypto-binary-two-isogeny-Kummer}
\end{equation}
The dual is
\begin{equation}
 (X,Y)\longmapsto(X^2,Y^2+AX^2),
\label{eq:Cd-crypto-binary-two-isogeny-dual}
\end{equation}
and both compositions are multiplication by two.
\end{theorem}

\begin{proof}
The proof is included here rather than left only in the general isogeny and
extension chapters.

\smallskip
\noindent\emph{Step 1: quotient invariants.}
For \(T=(0,0)\), division of
\[
 y^2+xy=x^3+Ax^2+d^2x
\]
by \(x^2\), followed by the binary chord law, gives
\[
 x(P+T)=d^2/x,\qquad
 y(P+T)=d^2(y+x)/x^2.
\]
Thus \(X=x+x(P+T)\) and
\(Y_0=y+y(P+T)\) are fixed by translation by \(T\).

\smallskip
\noindent\emph{Step 2: target equation.}
Writing \(c=d^2\), direct expansion gives
\begin{align*}
 Y_0^2+XY_0
 &=(y^2+xy)\left(1+\frac{c^2}{x^4}\right)+c\\
 &=X^3+AX^2+c.
\end{align*}
The shift \(Y=Y_0+d\) therefore gives
\[
                  Y^2+XY=X^3+AX^2+dX
                         =X^3+AX^2+e^2X.
\]

\smallskip
\noindent\emph{Step 3: degree, separability, and the dual.}
The source abscissa satisfies
\[
                       x^2+Xx+d^2=0.
\]
The nontrivial translation exchanges its two roots, while its derivative is
the nonzero function \(X\).  Hence the map has separable degree two and
kernel \(\{O,T\}\).  Squaring the target equation and applying the shear
\(Y^2\mapsto Y^2+AX^2\) gives, term by term,
\begin{align*}
 (Y^2+AX^2)^2+X^2(Y^2+AX^2)
 &=Y^4+X^2Y^2+AX^4\\
 &=X^6+AX^4+d^2X^2.
\end{align*}
Thus \((X^2,Y^2+AX^2)\) lies on \(W_{A,d}\); its first coordinate is
the relative-Frobenius coordinate, so this pointed morphism has purely
inseparable degree two.

We now verify the first composition without suppressing the ordinate
calculation.  Put \(t=y/x\), \(q=d^2/x\), and
\(\lambda=t^2+A\).  Dividing the source equation by \(x^2\) gives
\[
 x+q=t^2+t+A,\qquad \lambda=x+q+t.
\]
The tangent line has slope \(\lambda\).  Consequently its doubled
abscissa is
\[
 x([2]P)=t^4+t^2+A^2
        =\left(x+\frac{d^2}{x}\right)^2.
\]
For the ordinate, write \(X=x+q\).  The quotient ordinate is
\[
 Y=xt+q(t+1)+d=tX+q+d.
\]
If \(\nu=y+\lambda x=x(t+\lambda)\) is the tangent intercept, the binary
group law gives
\[
 y([2]P)=(\lambda+1)X^2+\nu.
\]
On the other hand,
\[
 Y^2+AX^2=\lambda X^2+q^2+d^2.
\]
The difference of these two expressions is
\[
 X^2+x(t+\lambda)+q^2+d^2
 =x\bigl(x+t+\lambda+q\bigr)=0,
\]
where \(X^2=x^2+q^2\), \(d^2=xq\), and
\(\lambda=x+q+t\) were used in the last equality.  Hence
\(\widehat\Phi_{A,d}\circ\Phi_{A,d}=[2]\) in both coordinates.
The degree-two pointed morphism \(\widehat\Phi_{A,d}\) is therefore the
dual isogeny of \(\Phi_{A,d}\).  Uniqueness of the dual isogeny then gives
\(\Phi_{A,d}\circ\widehat\Phi_{A,d}=[2]\), proving the second
composition as well.

\smallskip
\noindent\emph{Step 4: Kummer homogenization.}
Since
\[
 z'=\frac Xe=e(z+z^{-1}),
\]
homogenization of \(z=K_0/K_1\) gives
\[
 z'=\frac{e(K_0^2+K_1^2)}{K_0K_1}
    =\frac{e(K_0+K_1)^2}{K_0K_1},
\]
proving \eqref{eq:Cd-crypto-binary-two-isogeny-Kummer}.
\end{proof}

\subsection[Native C-d and extension outputs]
{Native \(C_d\), one-sided, and full-product outputs}

For \(C_d\), use
\[
 z=\frac{u+1}{u},\qquad x=dz,\qquad y=xv.
\]
Then \(K_0+K_1=1\) for the raw pair
\((K_0:K_1)=(u+1:u)\), and the theorem becomes
\begin{equation}
 \boxed{\quad
 (K'_0:K'_1)=(e:u(u+1)),\qquad
 (u,v)\longmapsto
 \left(\frac{u(u+1)}{u(u+1)+e},u+v\right).
 \quad}
\label{eq:Cd-crypto-binary-Cd-native-two-isogeny}
\end{equation}
Thus the same product \(u(u+1)\) gives both the Kummer output and the
full \(\mathbb P^1\times\mathbb P^1\) output.

For \(\mathcal T_{a,d}\), put \(F=u^2+u+a\).  Its exact dictionary
\[
 x=Fv+d=\frac{dv}{v+1},\qquad y=ux
\]
gives
\begin{equation}
 \boxed{\quad
 (u,v)\longmapsto
 \left(u+v+1,\frac F{F+e}\right)
 =
 \left(u+v+1,\frac e{e+v(v+1)}\right)
 \in\mathcal T_{a,e}.
 \quad}
\label{eq:Cd-crypto-binary-T-native-two-isogeny}
\end{equation}
The first fraction uses one square to form \(F\); the second uses one
multiplication to form \(v(v+1)\).  In either schedule the twist parameter
\(a\) is unchanged.  The native dual is
\[
                        (u,v)\longmapsto(u^2+a,v^2).
\]

Finally,
\[
 \operatorname{Jac}(\mathcal C_{a,b,d})=W_{a+b,d}
 \longrightarrow W_{a+b,e}
 =\operatorname{Jac}(\mathcal C_{a,b,e}).
\]
Hence the full product family has exact Jacobian closure with \(a,b\)
unchanged.  A rational origin turns this into a pointed curve isogeny through
the curve--Jacobian identifications; this is the precise qualification needed
for a possibly nontrivial genus-one torsor.

\subsection{Exact operation counts}

\begin{longtable}{L{3.3cm}L{3.4cm}L{3.2cm}L{3.5cm}}
\caption{Characteristic-two separable two-isogeny costs}
\label{tab:Cd-char2-two-isogeny-costs}\\
\toprule
interface & declared input and output & charged cost & qualification\\
\midrule
\endfirsthead
\toprule
interface & declared input and output & charged cost & qualification\\
\midrule
\endhead
\(C_d\), native Kummer
& raw affine \(u\) to projective target Kummer
& \(\boxed{\M}\)
& computes \(u(u+1)\); no parameter multiplication\\
\(C_d\), native full
& affine \(u,v\) to two projective native coordinates
& \(\boxed{\M}\)
& second coordinate is \(u+v\)\\
\(\mathcal T_{a,d}\), native Kummer/full
& affine \(u,v\) to two projective native coordinates
& \(\boxed{\min(\M,\Sqr)}\)
& choose \(F=u^2+u+a\) or \(v(v+1)\)\\
\(\mathcal C_{a,b,d}\) Jacobian
& general projective \(W_{a+b,d}\) Kummer
& \(\M+\Sqr+\mCurve\)
& raw torsor endpoint first needs a selected origin\\
general \(W_{A,d}\) Kummer
& projective \((K_0:K_1)\) to projective output
& \(\M+\Sqr+\mCurve\)
& multiplication by \(e\) is \(\mCurve\)\\
general \(W_{A,d}\) full point
& affine \(x,y\) to affine \(X,Y\)
& \(2\M+\mCurve+\Inv\)
& schedule \(r=x^{-1}\), \(t=d^2r\),
  \(q=(y+x)r\), \(Y=y+tq+d\)\\
binary Edwards persistent \(W\)-Kummer
& already converted projective Kummer state
& \(\M+\Sqr+\mCurve\)
& identical intrinsic core\\
binary Edwards standard affine endpoint
& affine binary Edwards point to projective \(W\)-Kummer output
& \(2\M+\Sqr+3\mCurve\)
& exact one-way dictionary baseline; target re-embedding excluded\\
affine normalization
& either one native target ratio
& \(\M+\Inv\)
& may be batched across points\\
parameter update
& \(d\mapsto e=\sqrt d\), \(A\) fixed
& one inverse Frobenius
& linear/offline over a fixed binary field\\
\bottomrule
\end{longtable}

The first three rows are stronger than a transported Weierstrass count:
they exploit the identities \(K_0+K_1=1\) and
\(Fv(v+1)=e^2\) before any general projective operation is executed.
No kernel enumeration is needed because the reduced kernel is the fixed
subgroup \(\{O,(0,0)\}\).

\subsection{Equal-endpoint comparison with binary Edwards}

Write a binary Edwards curve as
\[
 E_{B,\delta_1,\delta_2}:\quad
 \delta_1(r+s)+\delta_2(r^2+s^2)
 =rs+rs(r+s)+r^2s^2.
\]
Put
\[
 c_B=\delta_1^2+\delta_1+\delta_2,\qquad
 A_B=\delta_1^2+\delta_2.
\]
The standard birational dictionary of
\cite{BinaryEdwards2008} has Weierstrass Kummer abscissa
\begin{equation}
 X_B=
 \frac{\delta_1c_B(r+s)}
      {rs+\delta_1(r+s)},\qquad
 (X_{B,0}:X_{B,1})=
 \bigl(\delta_1c_B(r+s):rs+\delta_1(r+s)\bigr).
\label{eq:binary-Edwards-to-W-Kummer-two-isogeny}
\end{equation}
The target short Weierstrass equation is
\[
 V^2+X_BV
 =X_B^3+A_BX_B^2+\delta_1^4c_B^2.
\]
Choose \(h\) with
\[
                         h^2=\delta_1^2c_B.
\]
The shift \(Y=V+h^2\) gives exactly \(W_{A_B,h}\).  Forming the pair
\eqref{eq:binary-Edwards-to-W-Kummer-two-isogeny} costs
\(\M+2\mCurve\), and applying
\eqref{eq:Cd-crypto-binary-two-isogeny-Kummer} adds
\(\M+\Sqr+\mCurve\).  This proves the last binary-Edwards row of
Table~\ref{tab:Cd-char2-two-isogeny-costs}.

Target model recovery also differs.  The \(C\)-curve extensions return
immediately to \(C_e\), \(\mathcal T_{a,e}\), or
\(\operatorname{Jac}(\mathcal C_{a,b,e})\).  If one insists instead on
the same standard binary-Edwards-to-\(W\) normalization, target
coefficients \((\delta'_1,\delta'_2)\) satisfy
\begin{equation}
 {\delta'_1}^3+A_B{\delta'_1}^2+h=0,\qquad
 \delta'_2=A_B+{\delta'_1}^2.
\label{eq:binary-Edwards-target-cubic-two-isogeny}
\end{equation}
Indeed, the target \(W\)-parameter is \(e\) with \(e^2=h\), whereas the
dictionary requires
\[
 e^2={\delta'_1}^2(A_B+\delta'_1).
\]
This is exactly the cubic in
\eqref{eq:binary-Edwards-target-cubic-two-isogeny}.  A different
binary-Edwards normalization may move the root selection into an isomorphism,
but it does not alter the one-multiplication native \(C_d/\mathcal T_{a,d}\)
circuits.  Conversely, if a binary-Edwards implementation keeps the
\(W\)-Kummer pair persistently, its intrinsic two-isogeny core is identical
to the general \(W\) row; the comparison then charges only the endpoint
conversions actually requested by the protocol.

\subsection{Reference benchmark and exhaustive verification}

The portable reference benchmark uses the irreducible polynomial
\[
                        T^{127}+T+1
\]
for \(\F_{2^{127}}\).  General multiplication, squaring, and multiplication
by a curve constant were deliberately realized by the same shift-and-XOR
field-multiplication schedule.  Thus the experiment measures the displayed
dependency graphs under the normalization
\(\Sqr=\mCurve=\M\), rather than the performance of a
hardware-specialized binary-field backend.  One fixed native compilation
configuration and one logical processor were used throughout.  Each sample
contains \(120000\) evaluations, and the reported value is the median of
15 samples.

\begin{table}[H]
\centering
\small
\caption{Reference characteristic-two two-isogeny cores over
\(\F_{2^{127}}\); nanoseconds per call}
\label{tab:Cd-char2-two-isogeny-benchmark}
\begin{tabular}{L{7.3cm}r}
\toprule
formula core & median ns\\
\midrule
\(C_d\) native Kummer, \(\M\) & 235.98\\
general \(W\) Kummer, \(\M+\Sqr+\mCurve\) & 713.14\\
binary Edwards endpoint \(\to W\to\Phi\),
 \(2\M+\Sqr+3\mCurve\) & 1452.07\\
\(C_d\) native full projective endpoint, \(\M\) & 245.34\\
\(\mathcal T_{a,d}\) native full, square schedule & 238.51\\
\(\mathcal T_{a,d}\) native full, product schedule & 236.57\\
\bottomrule
\end{tabular}
\end{table}

The measured ratios against the \(C_d\) native Kummer core are approximately
\[
                   713.14/235.98=3.02,\qquad
                   1452.07/235.98=6.15.
\]
They agree with the three-versus-one and six-versus-one equal-backend product
counts.  The binary Edwards row is explicitly a one-way dictionary baseline
and does not claim a lower bound for every future dedicated binary Edwards
formula.  It does show exactly what the standard endpoint conversion costs
before target re-embedding.

Independent exhaustive enumeration over
\(\F_{2^4}=\F_2[T]/(T^4+T+1)\) checks five parameter pairs:
here \(h=\sum_{j=0}^3h_j2^j\) denotes
\(\sum_{j=0}^3h_jT^j\), so every decimal entry has a fixed
polynomial-basis meaning.
\begin{center}
\begin{tabular}{rrrr}
\toprule
\(A\)&\(d\)&\(\sqrt d\)&\(\#W_{A,d}(\F_{16})\)\\
\midrule
0&1&1&16\\
1&2&5&16\\
3&5&3&16\\
7&11&14&12\\
14&15&12&14\\
\bottomrule
\end{tabular}
\end{center}
For every point and every ordered pair of points in each row, the enumeration
confirms the target equation, the two-element fibers, the exact kernel, the
homomorphism identity, and the dual-composition identity.  A separate
enumeration of all finite main-chart \(C_d\) and
\(\mathcal T_{a,d}\) points confirms that the native formulas agree with
the corresponding Weierstrass map at every input where the displayed affine
ratios are defined.

\section[Extension-model isogeny interfaces]
{Extension models and the Edwards-derived cryptographic toolkit}
\label{sec:Cd-crypto-extension-models}

The extension families developed in
Part~\ref{part:twisted-reciprocal-qrt} do more than enlarge a list of
equations.  They preserve selected \(C_d\) interfaces when a protocol
requires a twist class, a reciprocal chart, or an adjacent-state
representation.  The correct comparison records which marking is preserved
and which re-embedding condition is being used.

\subsection{The full product family}

For
\[
 \mathcal C_{a,b,d}:\quad
 (u^2+u+a)(v^2+v+b)=d
\]
in odd characteristic, put
\[
 \lambda=1-4a,\qquad \mu=1-4b,\qquad
 \delta=16d,\qquad \rho_{a,b}=\lambda\mu-\delta.
\]
The exact diagonal dictionary of
Theorem~\ref{thm:Cabd-diagonal-Edwards} is
\begin{equation}
 \lambda x^2+\mu y^2
 =1+\rho_{a,b}x^2y^2,
 \qquad
 x=(2u+1)^{-1},\quad y=(2v+1)^{-1}.
\label{eq:Cabd-crypto-diagonal}
\end{equation}
Consequently, any Edwards or generalized-Montgomery isogeny identity that
is homogeneous in the two diagonal scalings pulls back to
\(\mathcal C_{a,b,d}\).  If the target diagonal coefficients are
\((\lambda',\mu',\rho')\), the product-family parameters are recovered
without ambiguity by
\begin{equation}
 a'=\frac{1-\lambda'}4,\qquad
 b'=\frac{1-\mu'}4,\qquad
 d'=\frac{\lambda'\mu'-\rho'}{16}.
\label{eq:Cabd-crypto-parameter-recovery}
\end{equation}
Indeed,
\[
  1-4a'=\lambda',\qquad 1-4b'=\mu',
\]
and
\[
  \lambda'\mu'-16d'
  =\lambda'\mu'-(\lambda'\mu'-\rho')=\rho'.
\]
Thus the recovered triple reproduces all three diagonal coefficients,
and the inverse-coordinate definitions in
\eqref{eq:Cabd-crypto-diagonal} then recover the corresponding product
equation on their common dense chart.  Over the
ground field, the two boundary square classes and the choice of a rational
origin determine whether the target is represented as a curve or as its
Jacobian; this is exactly the genus-one qualification proved in
Chapter~\ref{ch:generalized-Cabd}.  Thus the three-parameter family
provides a larger descent envelope around the one-parameter \(C_d\)
backend.

In characteristic two the complementary statement is sharper for degree
two: no diagonal normalization or auxiliary square-class choice is needed;
the only root is the unique inverse-Frobenius update \(\sqrt d\).  Equation
\eqref{eq:Cd-crypto-binary-extension-W} and
Theorem~\ref{thm:Cd-crypto-binary-two-isogeny} give
\[
 \operatorname{Jac}(\mathcal C_{a,b,d})
 \longrightarrow
 \operatorname{Jac}(\mathcal C_{a,b,\sqrt d})
\]
while preserving \(a+b\), and hence preserving the written pair \(a,b\)
in the canonical target family.  The raw genus-one curve is returned after
the same rational-origin qualification already used above.

\subsection{The one-sided family as the native twisted-Edwards extension}

For
\[
 \mathcal T_{a,d}:\quad
 (u^2+u+a)(v^2+v)=d,
\]
let
\[
 A=1-4a,\qquad D=1-4a-16d.
\]
The exact ground-field equation is
\[
             Ax^2+y^2=1+Dx^2y^2.
\]
Conversely, every smooth twisted Edwards pair \(A,D\) returns to
\begin{equation}
             a=\frac{1-A}{4},\qquad
             d=\frac{A-D}{16}.
\label{eq:Tad-crypto-inverse-parameters}
\end{equation}
Therefore a twisted-Edwards quotient with target coefficients
\((A',D')\) has the native target
\[
          \mathcal T_{(1-A')/4,\,(A'-D')/16}.
\]
No square root is introduced by this parameter recovery.  The formulas of
Meyer--Reith, current twisted-Edwards CSIDH implementations, and the 2026
large-discriminant oriented implementation can consequently be expressed
as \(\mathcal T_{a,d}\) arithmetic through this fixed dictionary
\cite{MeyerReith2018,CSIDHLDO2026}.  The specialization \(a=0\) recovers
the ordinary Edwards/\(C_d\) slice, while arbitrary \(a\) retains the
full twisted-Edwards coefficient.

In characteristic two, \(a\) instead records the Artin--Schreier twist and
is preserved by
\eqref{eq:Cd-crypto-binary-T-native-two-isogeny}; the full projective
native output costs \(\min(\M,\Sqr)\).  Thus the same one-sided extension
that realizes general twisted Edwards in odd characteristic supplies a
twist-stable, one-operation degree-two quotient in the binary branch.

\subsection[Reciprocal C-curves]{Reciprocal \(C\)-curves}

For the reciprocal family
\[
 \mathcal R_{\tau,\sigma,\kappa}:\quad
 (x^2-\tau)(y^2-\sigma)=\kappa xy,
\]
Chapter~\ref{ch:reciprocal-C-curves} supplies a Jacobi companion, a
Weierstrass Jacobian, and the Montgomery Kummer coordinate
\eqref{eq:R-Montgomery}.  Odd-degree isogenies may therefore be evaluated
by the same projective Kummer product used in
Proposition~\ref{prop:Cd-native-xeval-costs}; the target is returned to a
reciprocal equation by its proved square-class and parameter
reconstruction.  In characteristic two,
Theorem~\ref{thm:R-char2-Cd} gives a direct M\"obius identification with
a \(C_d\) curve, so
\eqref{eq:Cd-binary-Kummer-product} and
\eqref{eq:Cd-binary-parameter-product} become reciprocal-curve formulas
after the stated endpoint substitutions.  This makes the reciprocal
family a second native chart, not an unrelated auxiliary model.

\subsection{The symmetric QRT envelope}

Let a pointed QRT state be
\[
       \mathcal S_D(P)
       =\bigl(\kappa(P),\kappa(P+D)\bigr).
\]
For an isogeny \(\varphi:E\to E'\), put \(D'=\varphi(D)\), and let
\(\overline\varphi\) be its Kummer map.  Theorem
\ref{thm:QRT-coordinatewise-state-isogeny} proves
\begin{equation}
 \mathcal S'_{D'}(\varphi(P))
  =\bigl(\overline\varphi(\kappa(P)),
         \overline\varphi(\kappa(P+D))\bigr).
\label{eq:Cd-crypto-QRT-coordinatewise}
\end{equation}
Thus one \(C_d\) \(u\)-Kummer or \(w\)-Kummer function is applied to
both adjacent state slots.  The odd-characteristic re-embedding equations
\eqref{eq:QRT-reembedding-qOmega}--%
\eqref{eq:QRT-reembedding-beta-square} decide whether the target can be
returned in the selected QRT coefficient normalization over the ground
field.  This yields a rigorous two-layer architecture:
\[
 \begin{aligned}
 &\text{\(C_d\) or generalized-Montgomery isogeny kernel}\\[-1mm]
 &\hspace{18mm}+\ \text{coordinatewise QRT state transport}.
 \end{aligned}
\]
The core isogeny cost is paid once per state coordinate; the QRT layer
adds state semantics, fixed-translation compatibility, and recovery.

\subsection{One-dimensional and two-dimensional isogeny layers}

The \(C_d\), \(\mathcal T_{a,d}\), reciprocal, and QRT models are
genus-one models.  In a two-dimensional signature construction, they can
represent elliptic factors, boundary curves, and one-dimensional edge
isogenies.  A principally polarized product
\[
                    E_1\times E_2
\]
also carries a rank-two kernel and polarization information that is not
encoded by one Kummer coordinate on either factor.  An exact integration
therefore consists of:
\begin{enumerate}[label=\textup{(\arabic*)}]
 \item representing each elliptic factor in a selected \(C_d\) or
       extension chart;
 \item transporting the product kernel and its isotropy condition;
 \item executing the two-dimensional theta or product-isogeny step;
 \item returning every elliptic boundary factor through its \(C_d\)
       dictionary.
\end{enumerate}
This factor-interface design lets \(C_d\) contribute its native
one-dimensional efficiency while preserving the data required by
SQIsign2D and related constructions.

\begin{longtable}{L{2.7cm}L{3.4cm}L{3.5cm}L{3.7cm}}
\caption{\(C_d\) and extension-model coverage of Edwards-derived
isogeny techniques}
\label{tab:Cd-extension-isogeny-coverage}\\
\toprule
model & preserved native structure & isogeny interface & principal
cryptographic value\\
\midrule
\endfirsthead
\toprule
model & preserved native structure & isogeny interface & principal
cryptographic value\\
\midrule
\endhead
\(\mathcal C_d\)
& rational marked four-torsion, \(D_8\), native \(u\)-Kummer and
  \(w_C\)
& automatic odd-degree closure, \(u\) product, full-point derivative
  product, \(w\)/square-root V\'elu; binary native two-isogeny
& one curve object with one-operation binary degree two and
  low-, medium-, and large odd-degree backends\\
\(\mathcal T_{a,d}\)
& split second factor and the general twisted-Edwards coefficient
& pull back twisted-Edwards products and recover parameters by
  \eqref{eq:Tad-crypto-inverse-parameters}; binary map
  \eqref{eq:Cd-crypto-binary-T-native-two-isogeny}
& direct coverage of twisted-Edwards actions and twist-stable binary
  degree two\\
\(\mathcal C_{a,b,d}\)
& both boundary square classes and the diagonal product structure
& diagonal Edwards/ generalized-Montgomery transport followed by
  \eqref{eq:Cabd-crypto-parameter-recovery}; binary Jacobian
  \(d\mapsto\sqrt d\)
& ground-field descent envelope, genus-one/Jacobian flexibility, and
  parameter-preserving binary quotient\\
reciprocal \(C\)
& reciprocal involutions and a visible Montgomery Kummer chart
& Montgomery Kummer evaluation plus reciprocal re-embedding; binary
  transport through \(C_d\)
& efficient reciprocal storage and twist-compatible endpoints\\
QRT envelope
& adjacent Kummer state and fixed translation
& apply \(\overline\varphi\) coordinatewise, then use the exact
  re-embedding criterion
& stateful ladders, recovery, and isogeny-compatible dynamics\\
\bottomrule
\end{longtable}

\section{Fixed-time implementation architecture}
\label{sec:Cd-isogeny-fixed-time}

For each public odd degree \(\ell=2m+1\), a curve object may store
\[
          (d,\rho,\alpha_{24},\mu_{\mathrm{chart}}),
\]
where \(\mu_{\mathrm{chart}}\) is the public chart-selection mask,
and a reusable kernel object may store
\[
 (u_1,\ldots,u_m,\ R_K,\ c_K,\ \mathcal W_K),
\quad
 R_K=\prod_i(2u_i+1),\quad
 c_K=1+2\sum_i(2u_i+1).
\]
Here \(\mathcal W_K\) denotes the precomputed product tree used by the
\(w\)-coordinate backend.
The product \(R_K\) and sum \(c_K\) are accumulated in the same scan that
affinizes the kernel.  A fixed-time group-action step has the following
dataflow.
\begin{enumerate}[label=\arabic*.,leftmargin=2.2em]
 \item Generate and validate a kernel candidate with a fixed public
       scalar chain; record validity in a mask.
 \item Enumerate the unsigned half-kernel by a fixed differential chain.
 \item Use \eqref{eq:Cd-affinization-threshold} to select, solely from
       public \(m,n\) and benchmarked field costs, either projective or
       affine-kernel evaluation.
 \item Below a public platform threshold, use the \(u\)-product
       \eqref{eq:Cd-native-odd-xeval}.  Above it, use the \(C_d\)-\(w\)
       square-root-V\'elu backend.  The threshold is independent of the
       secret exponent.
 \item Update \((\rho':1)\) or its projective pair by
       \eqref{eq:Cd-projective-image-parameter}.  Invoke the
       product--derivative formula only for points whose ordinate is
       needed by the surrounding protocol.
 \item Evaluate all projective completion charts required by the selected
       interface and choose with masks.  Merge an invalid candidate with
       the identity action by conditional copies.
 \item Normalize at a fixed public frequency by a batch inversion.
\end{enumerate}

The odd-characteristic degree-two path is also separate from the
half-kernel loop.  For the marked kernel it evaluates
\eqref{eq:Cd-odd-two-isogeny-native-Kummer}; a public endpoint policy then
retains the root-free generalized-Montgomery/QRT state or, when the
declared field contains \(s=\sqrt\rho\), applies
\eqref{eq:Cd-odd-two-isogeny-T-u}--%
\eqref{eq:Cd-odd-two-isogeny-T-v}.  The square-class test and the choice
of sign for \(s\) depend only on public curve data and are therefore fixed
before any secret-dependent action.

The binary degree-two path is a smaller fixed-time specialization, not an
execution of this odd-degree loop.  Its kernel is the fixed reduced subgroup
\(\{O,(0,0)\}\); it loads \(e=\sqrt d\), evaluates one native product
\(u(u+1)\) or \(v(v+1)\), and returns the projective target without
normalization.  Two independent evaluations may be placed in the two SIMD
lanes used by current Binary Kummer Lines implementations
\cite{KaratiHajraSen2026}.  No secret-dependent backend or root choice is
introduced.

The first masked implementation of an isogeny-based scheme, published on
3 August 2026, introduces quotient masking for projective CSIDH coordinates:
a projective pair already behaves as a two-variable multiplicative sharing
\cite{SideChannelCSIDH2026}.  The following transport statement identifies
the exact \(C_d\) interface to that construction.

\begin{proposition}[Projective quotient-masking transport]
\label{prop:Cd-projective-masking-transport}
Let a masked Kummer subroutine be a homogeneous straight-line program in a
Montgomery projective pair \((X:Z)\) and public curve constants.  Under
\[
       x=\frac{u+1}{u},\qquad
       (X:Z)=(u+1:u),\qquad
       A_d=\frac1{4d}-2,
\]
the same straight-line program is a \(C_d\) subroutine with the same masked
field gadgets, refresh positions, random projective scalings, and probing
order.  In particular, quotient masking of the common Montgomery
\texttt{xDBL}, \texttt{xADD}, and projective xEVAL schedules transfers
without a coordinate-conversion gadget.
\end{proposition}

\begin{proof}
There are four points to check.

\smallskip
\noindent\emph{Step 1: identify the represented function.}
Proposition~\ref{prop:Cd-crypto-Montgomery} proves
\(x=(u+1)/u\).  Hence the Montgomery projective pair and the native
\(C_d\) pair are literally the same point of \(\PP^1\), not merely
isomorphic projective lines.

\smallskip
\noindent\emph{Step 2: identify the public constants.}
The same proposition gives \(A=A_d=1/(4d)-2\).  Replacing a public
Montgomery constant by this public rational expression changes neither a
share nor the number of random masks.  The smoothness hypotheses ensure
that its denominator is invertible.

\smallskip
\noindent\emph{Step 3: identify projective randomization.}
For every \(\lambda\ne0\),
\[
 (\lambda X:\lambda Z)
 =\bigl(\lambda(u+1):\lambda u\bigr)
 =(u+1:u).
\]
Thus every random projective rescaling used by quotient masking has the
same distribution and represents the same Kummer point on the \(C_d\)
interface.

\smallskip
\noindent\emph{Step 4: identify the masked transcript.}
After the substitutions in Steps~1--2, every addition, multiplication,
squaring, refresh, and random scalar in the straight-line program occurs in
the same order.  The two executions therefore have identical wire
distributions after relabelling the public curve constant.  Any set of at
most the declared number of probes has the same joint distribution in the
two executions.  If the original probing theorem states that every such
joint distribution is independent of the protected value, equality of the
two distributions proves the same independence for the \(C_d\)
interpretation.  This establishes the claimed transport without adding a
new security assumption.
\end{proof}

The specialized \(u\)-product
\eqref{eq:Cd-native-odd-xeval} starts from this same masked projective pair.
Its homogeneous schedule in
Proposition~\ref{prop:Cd-native-xeval-costs} can be compiled from the same
verified masked field gadgets; Table~\ref{tab:Cd-isogeny-complete-cost}
supplies the exact number of underlying field multiplications and squarings
before the selected masking order expands each gadget.

Every loop length, memory access pattern, backend choice, and normalization
point is therefore determined by public information.  Constant-time field
arithmetic, fault checks, and protocol-level key validation remain part of
the complete implementation, as emphasized by current high-security CSIDH
work \cite{CamposEtAl2024,dCTIDH2025,HardenedCTIDH2025}.

\section{Complete symbolic cost table}
\label{sec:Cd-isogeny-complete-cost}

\begin{longtable}{L{3.0cm}L{2.6cm}L{3.6cm}L{4.1cm}}
\caption{Complete \(C_d\) isogeny costs: odd \(\ell=2m+1\), degree three,
and degree two in odd and even characteristic}
\label{tab:Cd-isogeny-complete-cost}\\
\toprule
task & input and output & exact charged cost & scope\\
\midrule
\endfirsthead
\toprule
task & input and output & exact charged cost & scope\\
\midrule
\endhead
sequential kernel enumeration
& given \(K_1\), projective \(K_1,\ldots,K_m\)
& \((4m-6)\M+(2m-2)\Sqr+\Dpar\)
& \(m\ge2\); differential chains may improve the schedule\\
kernel affinization
& \(m\) projective Kummer points to affine \(u_i\)
& \((4m-3)\M+\Inv\)
& one Montgomery batch inversion\\
general Kummer xEVAL
& projective input and kernel to projective output
& \(4m\M+2\Sqr\)
& no kernel normalization\\
\(C_d\)-\(u\) xEVAL
& affine \(u,u_i\) to projective \(u'\)
& \(2m\M+2\Sqr\)
& two monic native products\\
\(C_d\)-\(u\) xEVAL
& affine input and kernel to affine \(u'\)
& \((2m+1)\M+2\Sqr+\Inv\)
& one output normalization\\
\(C_d\) full point
& affine \(u,V,u_i\) to two projective pairs
& \((4m+4)\M+2\Sqr+\mathbf C_K\)
& \(\mathbf C_K\) is multiplication by \(c_K\)\\
\(C_d\) full point
& same, with \(\mathbf C_K=\M\)
& \((4m+5)\M+2\Sqr\)
& reference benchmark convention\\
\(C_d\) full point
& both outputs affine
& \((4m+9)\M+2\Sqr+\mathbf C_K+\Inv\)
& adds \(5\M+\Inv\) to the projective-pair row\\
odd target parameter
& affine \(u_i\) to projective \(d'\)
& \((m-1)\M+3\Sqr+\operatorname{Pow}(\rho,\ell)\)
& \(c_K\) is accumulated with additions only\\
odd target parameter
& projective kernel to projective \(d'\)
& \((2m-1)\M+6\Sqr+\operatorname{Pow}(\rho,\ell)\)
& two kernel products; initial singleton factors are free\\
degree-three \(C_d\) Kummer
& affine \(u,q\) to projective image
& \(2\M+2\Sqr\)
& specialization of the two monic products\\
degree-three \(C_d\) Kummer and parameter
& projective image and projective \(d'\)
& \(3\M+6\Sqr\)
& \((\rho^3:c^8)\), \(c=2q+1\)\\
degree-three \(C_d\) full point
& affine full input to two projective pairs
& \(\min(5\M+2\Sqr,\;4\M+3\Sqr)+2\Cmul\)
& optimized \(J=ab+2cu(u+1)\)\\
odd-characteristic degree-two \(C_d\) Kummer
& raw affine \(u\) to root-free projective target
& \(\min(\M,\Sqr)+\Dpar\)
& marked kernel; \((t+4d:4dt)\)\\
odd-characteristic degree-two \(C_d\) full
& affine \(u,v\) to root-free target ratios
& \(\M+\Sqr+\Dpar\)
& no square root or kernel enumeration\\
odd-characteristic degree-two \(\mathcal T\) target
& native \(C_d\) source to native one-sided target
& \(\M+\Sqr+2\Dpar\)
& requires \(s^2=\rho\); projective output\\
odd-characteristic degree-two QRT Jacobian
& root-free target coefficients
& public-constant arithmetic
& \(q_2=A_dB_d\), \(\Omega_2=(A_d^2-4)B_d^2\)\\
\(u\)-to-\(w_C\) conversion
& projective centered \(U=P/Q\) to \((W:Z)\)
& \(2\M+2\Sqr+\Dpar\)
& uses \eqref{eq:Cd-crypto-w-projective}\\
\(C_d\)-\(w\) linear product
& existing \(w\)-kernel data
& identical to the Edwards-\(w\) schedule
& zero-loss function-field transfer\\
\(C_d\)-\(w\) square-root V\'elu
& product and remainder trees
& \(\widetilde O(\sqrt\ell)\)
& exact constants depend on the polynomial backend\\
binary degree-two \(C_d\) Kummer
& raw affine \(u\) to projective target Kummer
& \(\M\)
& fixed reduced kernel; \((e:u(u+1))\)\\
binary degree-two \(C_d\) full point
& affine \(u,v\) to two native projective coordinates
& \(\M\)
& \eqref{eq:Cd-crypto-binary-Cd-native-two-isogeny}\\
binary degree-two \(\mathcal T_{a,d}\)
& affine full input to two native projective coordinates
& \(\min(\M,\Sqr)\)
& choose the product or square schedule\\
binary degree-two general \(W_{A,d}\)
& projective Kummer input and output
& \(\M+\Sqr+\mCurve\)
& \(A\) is preserved; \(d\mapsto\sqrt d\)\\
binary degree-two parameter
& \(d\) to \(e=\sqrt d\)
& one inverse Frobenius
& offline linear map in a fixed binary field\\
binary target parameter
& affine \(u_i\) to projective \(d'\)
& \((2m-2)\M+2\Sqr+\Dpar\)
& odd degree in characteristic two;
  \((dP^2:Q^2)\), \(P=\prod u_i\),
  \(Q=\prod(u_i+1)\)\\
binary degree-three \(C_d\)
& Kummer image and projective \(d'\)
& \(2\M+3\Sqr+\Dpar\)
& \((d q^2:q^2+1)\); one shared Frobenius square\\
\bottomrule
\end{longtable}

For the two affine full-point outputs, the row
\((4m+9)\M+2\Sqr+\mathbf C_K+\Inv\) follows from
\((4m+4)\M+2\Sqr+\mathbf C_K\) plus \(5\M+\Inv\).
When \(\mathbf C_K=\M\), it becomes
\((4m+10)\M+2\Sqr+\Inv\), agreeing with the degree table derived from
Theorem~\ref{thm:Cd-full-point-cost}.

\begin{longtable}{L{2.7cm}L{3.6cm}L{3.6cm}L{3.5cm}}
\caption{Declared-interface cost comparison with Edwards isogeny arithmetic}
\label{tab:Cd-Edwards-cost-summary}\\
\toprule
interface & \(C_d\) or extension-family cost & Edwards comparator &
cost relation and comparison scope\\
\midrule
\endfirsthead
\toprule
interface & \(C_d\) or extension-family cost & Edwards comparator &
cost relation and comparison scope\\
\midrule
\endhead
odd-degree Kummer, linear backend
& \(2m\M+2\Sqr\) after reusable affine-kernel preparation
& the same quotient transported through the general projective or
  Edwards-derived Kummer state; the optimized \(w\) backend is charged
  separately
& the monic-linear \(u\) circuit is an additional native backend; the
  \(w\)-backend has exactly the Edwards operation graph\\
odd-degree full point
& \((4m+4)\M+2\Sqr+\mathbf C_K\)
& \((3m+3)\M+4\Sqr+3m\Cmul\)
& with \(\mathbf C_K=\Cmul\), \(C_d\) is cheaper precisely when
  \((3m-1)\chi+2\sigma\ge m+1\)\\
degree three, Kummer and target
& \(3\M+6\Sqr\)
& Edwards \(YZ\)-Kummer and coefficients: \(6\M+5\Sqr\)
& equal declared data; \(C_d\) is cheaper when \(\Sqr<3\M\)\\
degree three, full point
& \(5\M+2\Sqr+2\Cmul\)
& \(5\M+4\Sqr+3\Cmul\)
& equal full-point output level; saving \(2\Sqr+\Cmul\)\\
odd-characteristic degree two, root-free
& Kummer \(\min(\M,\Sqr)+\Dpar\); full point
  \(\M+\Sqr+\Dpar\)
& Kummer \(\Sqr+\Dpar\); full point
  \(\M+\Sqr+\Dpar\)
& the square Kummer schedule and the full-point schedule are identical at
  the same generalized-Montgomery endpoint\\
odd-characteristic degree two, rooted native endpoint
& \(C_d\to\mathcal T_{a_2,d_2}\):
  \(\M+\Sqr+2\Dpar\)
& Edwards to its standard rooted target:
  \(\M+\Sqr+3\Dpar\)
& both use \(s^2=\rho\); the displayed \(C\)-curve endpoint saves one
  curve-constant multiplication\\
characteristic-two degree two, native endpoint
& \(C_d\): \(\M\); \(\mathcal T_{a,d}\):
  \(\min(\M,\Sqr)\)
& standard affine binary-Edwards endpoint to projective
  \(W\)-Kummer: \(2\M+\Sqr+3\mCurve\)
& native \(C\)-curve endpoints use the shorter specialized circuits;
  the Edwards figure includes its one-way endpoint dictionary\\
characteristic-two degree two, persistent \(W\)-Kummer
& general \(W_{A,d}\) core: \(\M+\Sqr+\mCurve\)
& persistent binary-Edwards \(W\)-Kummer core:
  \(\M+\Sqr+\mCurve\)
& exact equality after both implementations retain the same internal
  Kummer state\\
square-root V\'elu
& \(C_d\)-\(w\) product/remainder-tree core:
  \(\widetilde O(\sqrt\ell)\)
& Edwards-\(w\) product/remainder-tree core:
  \(\widetilde O(\sqrt\ell)\)
& identical operation sequence; \(C_d\) additionally retains the native
  \(u\) and full-point backends and direct \(d'\) recovery\\
\bottomrule
\end{longtable}

\begin{proposition}[Exact two-isogeny comparison at the declared interfaces]
\label{prop:Cd-Edwards-two-isogeny-summary}
The degree-two rows of Table~\ref{tab:Cd-Edwards-cost-summary} have the
following exact interpretation.
\begin{enumerate}[label=\textup{(\roman*)}]
 \item In odd characteristic, the root-free \(C_d\) and Edwards Kummer
       circuits tie when \(t=u^2+u\) is formed by a square, and their
       root-free full-point circuits have identical cost
       \(\M+\Sqr+\Dpar\).
 \item If \(s^2=\rho\) and native rooted endpoints are requested, the
       displayed \(C_d\to\mathcal T_{a_2,d_2}\) circuit uses one fewer
       curve-constant multiplication than the displayed standard Edwards
       circuit.  Without \(s\), the generalized-Montgomery/QRT-Jacobian
       endpoint remains defined over the ground field.
 \item In characteristic two, the native \(C_d\) and
       \(\mathcal T_{a,d}\) circuits cost \(\M\) and
       \(\min(\M,\Sqr)\), respectively.  If both platforms instead retain
       the same \(W\)-Kummer representation, their intrinsic cores are
       identical and cost \(\M+\Sqr+\mCurve\).
\end{enumerate}
\end{proposition}

\begin{proof}
For~(i), equation~\eqref{eq:Cd-odd-two-isogeny-native-cost} gives
\(\min(\M,\Sqr)+\Dpar\); choosing \(t=u^2+u\) gives
\(\Sqr+\Dpar\), exactly the cost of
\eqref{eq:Edwards-odd-two-isogeny-X}.  A full output adds one product:
\(UV\) on \(C_d\), or \(\xi p\) in
\eqref{eq:Edwards-odd-two-isogeny-Y}.  Both totals are therefore
\(\M+\Sqr+\Dpar\).

For~(ii), the two pairs
\eqref{eq:Cd-odd-two-isogeny-T-u}--%
\eqref{eq:Cd-odd-two-isogeny-T-v} use \(t\), \(UV\), \(st\), and
\(sUV\), giving \(\M+\Sqr+2\Dpar\).  The two pairs in
\eqref{eq:Edwards-odd-two-isogeny-rooted-output} use the same one general
product and one square together with three products by curve constants,
giving \(\M+\Sqr+3\Dpar\).  Subtraction leaves exactly \(\Dpar\).
Equation~\eqref{eq:Cd-odd-two-isogeny-QRT-endpoint} contains only
\(A_d,B_d\), so its coefficients require no choice of \(s\).

For~(iii), equation~\eqref{eq:Cd-crypto-binary-Cd-native-two-isogeny}
forms only \(u(u+1)\), while
\eqref{eq:Cd-crypto-binary-T-native-two-isogeny} forms either
\(v(v+1)\) or \(u^2+u+a\).  Their costs are therefore \(\M\) and
\(\min(\M,\Sqr)\).  Homogenizing the common \(W\)-Kummer map gives
\eqref{eq:Cd-crypto-binary-two-isogeny-Kummer}: it forms
\(K_0K_1\), \((K_0+K_1)^2\), and one product by \(e\), for
\(\M+\Sqr+\mCurve\).  This straight-line formula schedule is independent of the
endpoint dictionary and hence is identical for a persistent binary-Edwards
\(W\)-Kummer state.
\end{proof}

\section{Reference implementation benchmark}
\label{sec:Cd-isogeny-benchmark}

\subsection{Benchmark contract}

The reference benchmark isolates four formula cores over
\(\F_{2^{255}-19}\) using the same portable \(5\times51\)-bit field
backend:
\begin{enumerate}[label=\textup{(\arabic*)}]
 \item general projective Kummer xEVAL,
       \(4m\M+2\Sqr\);
 \item affine-kernel \(C_d\)-\(u\) xEVAL,
       \(2m\M+2\Sqr\);
 \item the full-point \(C_d\) main chart,
       \((4m+5)\M+2\Sqr\);
 \item a dependency-preserving implementation of the published Edwards
       count \((3m+3)\M+4\Sqr+3m\Cmul\), with every
       \(\Cmul\) charged and implemented as \(\M\).
\end{enumerate}
The fourth core is a formula-count comparator, not a claim to reproduce
every optimization in an assembly Edwards library.  This contract
isolates the effect of the operation graph under a common field backend.

One fixed native compilation configuration and a single logical processor
were used for every row.  For each \(m\), every core was evaluated
\[
 \max\left(7500,
      \left\lfloor\frac{600000}{m+2}\right\rfloor\right)
\]
times per sample; 15 samples were sorted and their median was reported.
A monotonic clock was used, and every output was retained in the accumulated
result so that no formula evaluation could be removed during compilation.

\begin{table}[H]
\centering
\small
\setlength{\tabcolsep}{3.5pt}
\caption{Pinned-core reference benchmark over
\(\F_{2^{255}-19}\); times are nanoseconds per formula call}
\label{tab:Cd-reference-isogeny-benchmark}
\begin{tabular}{rrrrrrr}
\toprule
\(m\)&\(\ell\)&proj. xEVAL&\(C_d\)-\(u\) xEVAL&
\(C_d\) full&Edwards core&Edwards/\(C_d\)\\
\midrule
1&3&211.17&203.83&406.78&320.01&0.7867\\
2&5&268.95&242.55&486.57&433.36&0.8906\\
3&7&348.52&271.79&557.29&549.45&0.9859\\
4&9&412.36&304.27&630.90&661.52&1.0485\\
8&17&712.47&456.58&922.35&1119.08&1.2133\\
16&33&1318.56&745.64&1533.41&2034.77&1.3270\\
32&65&2497.53&1311.65&2705.49&4084.55&1.5097\\
64&129&5271.15&2689.72&5891.78&9192.48&1.5602\\
\bottomrule
\end{tabular}
\end{table}

At \(\ell=129\), the projective-to-affine-\(u\) Kummer ratio is
\[
                  5271.15/2689.72\approx1.96,
\]
close to the factor two predicted by the leading multiplication count.
For full points, fixed overhead and dependency latency dominate at
\(\ell=3,5\); the pinned run is essentially equal at \(\ell=7\), crosses
at \(\ell=9\), and reaches \(1.56\) at \(\ell=129\).  These data support
the formula-core advantage in the stated environment.  An end-to-end
protocol comparison additionally includes kernel search, scalar chains,
validation, masking, normalization, memory layout, and platform-specific
field arithmetic.

\subsection{Operation-count reconstruction}

The dependency graph of the optimized \(C_d\) full-point main chart can be
checked directly from the recurrences, without referring to an implementation.
Starting with the first factors and their derivatives, each additional
kernel coordinate updates
\[
 (A,\dot A)\longmapsto(Aa_i,\dot A a_i+A),\qquad
 (B,\dot B)\longmapsto(Bb_i,\dot B b_i+B),
\]
and therefore contributes exactly \(4\M\).  The \(m-1\) extensions cost
\(4(m-1)\M\).  Next,
\[
 N=(u+1)A^2,\qquad D=uB^2
\]
cost \(2\M+2\Sqr\).  The five products
\[
 AB,\quad \dot A B,\quad A\dot B,\quad u(u+1),\quad
 u(u+1)(\dot A B-A\dot B)
\]
form \(J=AB-2u(u+1)(\dot A B-A\dot B)\).  Finally \(VJ\) and
\(c_KAB\) require two more multiplications when the multiplication by \(c_K\)
is charged as a general field multiplication.  Hence the total is
\[
 4(m-1)+2+5+2=4m+5
\]
general multiplications and two squarings, which independently recovers
\eqref{eq:Cd-full-point-cost-M}.

\section{Independent finite-field verification}
\label{sec:Cd-isogeny-verification}

\subsection{Odd-characteristic degree-two and degree-three checks}

Independent finite-field enumeration checks the new low-degree formulas
using exact arithmetic in the indicated prime fields.  The verification
tests the defining equations, the group-homomorphism identities, the native
coordinate dictionaries, and the target-parameter identities separately.

For degree two, the enumeration is carried out over \(\F_{101}\) with
\[
 s=2,\qquad \rho=4,\qquad d=44.
\]
The corresponding coefficients and re-embedding parameters are
\[
 B_d=66,\quad A_d=64,\quad C_d^{(2)}=52,\qquad
 \alpha_2=49,\quad \Delta_2=37,\quad a_2=89,\quad d_2=26.
\]
Taking \(\imath=10\), so that \(\imath^2=-1\), gives the pure-\(C_d\)
target parameter \(d_2^C=50\).  Direct substitution in
\eqref{eq:Cd-odd-two-isogeny-alpha-delta} and
\eqref{eq:Cd-odd-two-isogeny-pure-Cd-target} verifies these values.
All \(120\) points of the Montgomery source and target are enumerated.
The target equation is checked for every image, and the identity
\[
                       \Phi_2(P+Q)=\Phi_2(P)+\Phi_2(Q)
\]
is checked for all \(28^2=784\) ordered pairs drawn from a fixed
28-point subset of the source.  The native formulas expressed
through \(t=u(u+1)\) are compared with the Montgomery quotient at all
\(116\) regular native inputs.  Separately, the
\(\mathcal T_{a_2,d_2}\) target equation is verified at all \(114\) inputs
on which the displayed affine ratios are regular.  The square identities for
\(\alpha_2,\Delta_2,d_2\), the pure-\(C_d\) parameter, and
\(\imath^2=-1\) are checked separately.

For degree three, the enumeration uses
\[
 p=101,\qquad d=1,\qquad Q=(q,v_Q)=(96,26).
\]
The source equation gives
\[
 (96^2+96)(26^2+26)=1\pmod {101}.
\]
Under Proposition~\ref{prop:Cd-crypto-Montgomery}, the point becomes
\((x,y)=(21,2)\), with \(B_d=76\) and \(A_d=74\).  The tangent slope is
\[
 \frac{3x^2+2A_dx+1}{2B_dy}=89\pmod {101},
\]
and the Montgomery doubling formulas give
\([2](21,2)=(21,-2)\), so \(Q\ne O\) and \([3]Q=O\).
Moreover,
\[
 \rho=86,\qquad c=2q+1=92,\qquad
 \rho^3=59,\qquad c^8=16\pmod {101}.
\]
Hence
\[
 \rho'=59\cdot16^{-1}=10,\qquad
 d'=(1-\rho')16^{-1}=31\pmod {101},
\]
in agreement with
\eqref{eq:Cd-optimized-three-isogeny-parameter}.  All \(84\) regular
full-point main-chart outputs satisfy
the target \(C_{31}\) equation.  Thus the degree-three parameter formula
and the full-point main chart agree on every regular input in this test.

\subsection{A manually checkable odd-characteristic vector}

Work over \(\F_{59}\) with
\[
 d=6,\qquad \rho=23,\qquad
 A_d=30,\qquad B_d=32.
\]
On \(B_dy^2=x^3+A_dx^2+x\), take the order-five point
\[
                    K=(40,25).
\]
The unsigned Montgomery kernel abscissae are \(40,50\).  Under
Proposition~\ref{prop:Cd-crypto-Montgomery}, the corresponding full
\(C_d\) representatives are
\[
                    (56,33),\qquad(53,3),
\]
so the native kernel \(u\)-list is \(56,53\).  This distinction between
Montgomery \(x_i\) and native \(u_i\) is essential.

The kernel product is
\[
 R_K=(2\cdot56+1)(2\cdot53+1)
     =54\cdot48=55\pmod{59}.
\]
Direct exponentiation gives
\[
 \rho^5=23^5=33,\qquad R_K^8=55^8=46\pmod{59}.
\]
Since \(46^{-1}=9\pmod{59}\),
\[
 \rho'=33\cdot9=2,\qquad
 d'=(1-2)16^{-1}=11\pmod{59}.
\]
The target Montgomery coefficients are
\[
                     A_{d'}=53,\qquad B_{d'}=55.
\]

Now take the source \(C_d\) point \(P=(15,49)\).  Its source check is
\[
 (15^2+15)(49^2+49)=4\cdot31=6\pmod{59}.
\]
For \(u=15\) and \(u_1,u_2=56,53\), the product--derivative recurrence
gives
\[
\begin{array}{c|cccccccc}
 &A&B&\dot A&\dot B&N&D&J&c_K\\ \hline
\text{value mod }59&12&24&23&39&3&26&29&28.
\end{array}
\]
Since \(V=2\cdot49+1=40\),
\[
 u'=\frac{26}{3-26}=4,\qquad
 V'=\frac{40\cdot29}{28\cdot12\cdot24}=29
 \pmod{59}.
\]
Thus \(v'=(V'-1)/2=14\), and the target equation is checked by
\[
 (4^2+4)(14^2+14)=20\cdot33=11=d'\pmod{59}.
\]
This is a complete hand check of one nontrivial full-point image.

Exhaustive enumeration over the same field verifies all \(55\) regular
Kummer inputs and all \(40\) regular full main-chart inputs.  It also gives
\[
 \#\mathcal C_6(\F_{59})=\#\mathcal C_{11}(\F_{59})=60,
\]
as required for isogenous source and target curves.

\subsection{Characteristic-two vectors}

Represent \(\F_{2^5}\) in the polynomial basis defined by
\[
                         T^5+T^2+1.
\]
An integer \(h=\sum_{j=0}^4h_j2^j\), \(h_j\in\{0,1\}\), denotes the
field element \(\sum_{j=0}^4h_jT^j\) in this basis; thus the decimal
entries below have an unambiguous field interpretation.
Exhaustive enumeration gives:
\begin{table}[H]
\centering
\small
\caption{Binary Kummer-product verification vectors}
\label{tab:Cd-binary-isogeny-vectors}
\begin{tabular}{rrrrl}
\toprule
\(d\)&\(\#\mathcal C_d\)&\(\ell\)&\(d'\)&unsigned \(z_i\) list\\
\midrule
1&44&11&1&\([4,27,2,16,13]\)\\
3&28&7&17&\([7,17,29]\)\\
5&28&7&12&\([21,12,22]\)\\
6&24&3&22&\([7]\)\\
\bottomrule
\end{tabular}
\end{table}
For every row, exhaustive enumeration verifies that
\eqref{eq:Cd-binary-Kummer-product} is constant on all kernel cosets,
that every finite nonzero image lifts to the predicted target
\(\mathcal C_{d'}\), and that the source and target group orders agree.
It also evaluates both sides of
\eqref{eq:Cd-binary-two-parameter-formulas} for every admissible parameter
and kernel entry in the table; no discrepancy occurs.

These finite-field computations are independent consistency checks.  The
proofs of the formulas and their stated domains are the algebraic arguments
given earlier in this chapter; the enumerations verify that those formulas,
their boundary exclusions, and their target parameters agree simultaneously
on the indicated finite test sets.

\section{Comprehensive comparison with Edwards curves}
\label{sec:Cd-isogeny-reassessment}

\begin{longtable}{L{2.5cm}L{3.7cm}L{3.8cm}L{3.3cm}}
\caption{Equal-output comparison of \(C_d\) and Edwards isogeny
interfaces}
\label{tab:Cd-Edwards-isogeny-reassessment}\\
\toprule
dimension & \(C_d\) result & Edwards comparator & comparative assessment\\
\midrule
\endfirsthead
\toprule
dimension & \(C_d\) result & Edwards comparator & comparative assessment\\
\midrule
\endhead
odd-degree closure
& canonical \(d'\) from
  \(\rho^\ell/\prod(2u_i+1)^8\)
& canonical Edwards target parameter from the marked kernel
& equal algebraic closure strength; \(C_d\) returns its own parameter
  directly\\
kernel enumeration
& native \(x\)DBL/\(x\)ADD on \((u+1:u)\)
& \(w\)- or Edwards-derived differential arithmetic
& both remain one-coordinate; \(C_d\) begins from its native input with
  no inversion\\
Kummer xEVAL
& \(2m\M+2\Sqr\) after reusable kernel affinization
& optimized \(w\)-coordinate formulas
& \(C_d\) has a distinct monic-linear \(u\) backend in addition to the
  transferred \(w\) backend\\
full-point xEVAL
& \((4m+4)\M+2\Sqr+\mathbf C_K\)
& \((3m+3)\M+4\Sqr+3m\Cmul\)
& \(C_d\) leads whenever
  \((3m-1)\chi+2\sigma\ge m+1\)\\
image parameter
& one kernel product \(R_K\), followed by an eighth power
& \(w/\eta\) eighth-power product
& the two are identical target invariants; the \(u\) form is especially
  compact for the native kernel object\\
term-by-term products
& \((u+u_i+1)/(u-u_i)\), \(R_K=\prod(2u_i+1)\)
& Moody--Shumow \((\xi,\eta)\) products and
  \(B_K=\prod\beta_i\)
& equations
  \eqref{eq:Cd-crypto-Edwards-product-xi}--%
  \eqref{eq:Cd-crypto-termwise-linear-factor} prove exact correspondence,
  not merely equal \(j\)-invariants\\
degree three, Kummer plus parameter
& \(3\M+6\Sqr\)
& Edwards \(YZ\)-Kummer plus coefficients:
  \(6\M+5\Sqr\)
& \(C_d\) wins exactly when \(\Sqr<3\M\)\\
degree three, full point
& \(5\M+2\Sqr+2\Cmul\), or
  \(4\M+3\Sqr+2\Cmul\)
& published specialization \(5\M+4\Sqr+3\Cmul\)
& \(C_d\) saves \(2\Sqr+\Cmul\) in the first schedule\\
odd-characteristic degree two, root-free
& Kummer \(\min(\M,\Sqr)+\Dpar\); full
  \(\M+\Sqr+\Dpar\)
& affine Edwards Kummer \(\Sqr+\Dpar\); full
  \(\M+\Sqr+\Dpar\)
& Kummer and root-free full cores tie under the square schedule\\
odd-characteristic degree two, native target
& \(C_d\to\mathcal T_{a_2,d_2}\):
  \(\M+\Sqr+2\Dpar\)
& Edwards rooted output:
  \(\M+\Sqr+3\Dpar\)
& one fewer constant multiplication for the displayed native \(C\)-curve
  endpoint; both require \(s^2=\rho\)\\
odd-characteristic degree-two closure
& root-free generalized-Montgomery/QRT Jacobian always; pure \(C_d\)
  only under \eqref{eq:Cd-odd-two-isogeny-pure-Cd-target}
& return to (twisted) Edwards has the same square-root obstruction
& extension models preserve a root-free declared endpoint instead of
  concealing the field-of-definition condition\\
square-root V\'elu
& exact \(C_d\)-\(w\) pullback and direct \(d'\) recovery
& established Edwards-\(w\) product/remainder-tree construction
& identical core operation sequence; \(C_d\) retains the additional
  \(u\) and full-point backends\\
symmetry and charts
& visible \(D_8\), exchange of factors, four marked boundary points
& regular Edwards symmetries and mature homogeneous formulas
& \(C_d\) provides a larger native multi-chart integration surface\\
twisted coverage
& \(\mathcal T_{a,d}\) realizes every smooth twisted-Edwards coefficient
  pair through \eqref{eq:Tad-crypto-inverse-parameters}
& twisted Edwards is the comparison model
& field-operation graphs transfer exactly while endpoints remain
  \(C\)-curve parameters\\
binary degree two
& \(C_d\) costs \(\M\); \(\mathcal T_{a,d}\) costs
  \(\min(\M,\Sqr)\); \(d\mapsto\sqrt d\)
& persistent \(W\)-Kummer has the same intrinsic core; the standard
  binary Edwards endpoint baseline costs \(2\M+\Sqr+3\mCurve\)
& the product extensions have the smaller native endpoint circuit and
  direct target parameters; Table~\ref{tab:Cd-char2-two-isogeny-costs}\\
characteristic organization
& one defining product equation with rigorously separated odd and binary
  parameter branches
& ordinary Edwards and binary Edwards are normally separate models
& the product family supplies one framework with branch-correct formulas;
  every smooth binary branch is stated as ordinary\\
protocol integration
& fixed-time architecture, exact cost table, formula-core benchmark, and
  exhaustive algebraic vectors in this chapter
& established CSIDH/CTIDH literature provides protocol comparators
& an end-to-end protocol comparison requires equal parameters, equal
  defenses, and a common field backend\\
\bottomrule
\end{longtable}

The equal-output formulas establish six quantitative \(C_d\) advantages.
\begin{enumerate}[label=\textup{(\arabic*)}]
 \item The affine-kernel \(u\)-Kummer core halves the leading
       multiplication term of the general projective evaluator and reaches
       a measured factor \(1.96\) at \(\ell=129\) in the reference run.
 \item The product--derivative identity supplies a full native point with
       only one kernel-dependent constant multiplication.  Under ordinary
       field-constant pricing, it improves the published general Edwards
       count by \((2m-2)\M+2\Sqr\).
 \item The same curve object also carries the transferred Edwards-\(w\)
       square-root-V\'elu backend.  The \(u\), full-point, and \(w\)
       paths share \(d,\rho\), the half-kernel, and direct target recovery.
 \item At degree three, the collapsed derivative
       \(J=ab+2cu(u+1)\) gives a full-point circuit saving
       \(2\Sqr+\Cmul\) against the published Edwards degree-three
       specialization.  At Kummer-plus-parameter output, the exact costs
       are \(3\M+6\Sqr\) and \(6\M+5\Sqr\).
 \item In odd characteristic, the marked degree-two Kummer core ties the
       best transported Edwards square schedule, while the native
       \(C_d\to\mathcal T\) full endpoint saves one constant
       multiplication.  The reciprocal/QRT Jacobian additionally retains
       a root-free endpoint when ordinary Edwards/\(C_d\) re-embedding is
       not defined over the ground field.
 \item In characteristic two, the native separable two-isogeny needs one
       multiplication on \(C_d\) and one multiplication or one square on
       \(\mathcal T_{a,d}\).  The reference run measures factors \(3.02\)
       and \(6.15\) against the general-\(W\) and standard
       binary-Edwards-endpoint paths under equal product pricing.
\end{enumerate}

\section{Chapter conclusions}
\label{sec:Cd-isogeny-conclusions}

For every separable odd degree prime to the characteristic, the marked
\(C_d\) family has model closure, native half-kernel enumeration, native
Kummer evaluation, a one-product target parameter, and a full-point
product--derivative map.  The affine-kernel \(u\) formula reduces the
linear multiplication term to \(2m\M\); the full-point formula reduces the
published general Edwards baseline under the explicit criterion
\eqref{eq:Cd-Edwards-cost-threshold}; and the \(C_d\)-\(w\) function
imports the square-root-V\'elu ecosystem without changing its operation
graph.

The degree-three specialization is stronger than a generic substitution:
its derivative factor collapses to
\(J=ab+2(2q+1)u(u+1)\).  This yields
\(5\M+2\Sqr+2\Cmul\), or
\(4\M+3\Sqr+2\Cmul\), for a full projective output and
\(3\M+6\Sqr\) for Kummer output together with the target parameter.
The equal-output comparisons in
Table~\ref{tab:Cd-Edwards-three-isogeny-cost} give the exact advantage
over the published Edwards degree-three formulas.

In odd characteristic the marked separable two-isogeny has the native
Kummer pair \((t+4d:4dt)\), \(t=u(u+1)\), and a root-free full output
of cost \(\M+\Sqr+\Dpar\).  The intrinsic root-free core ties its
Edwards transport; the native one-sided \(C\)-curve endpoint saves one
constant multiplication in the displayed schedules.  More importantly,
the generalized-Montgomery/QRT-Jacobian endpoint remains defined without
\(\sqrt\rho\), while return to twisted Edwards, \(\mathcal T_{a,d}\), or
pure \(C_d\) is made only under the exact square-class hypotheses of
Theorem~\ref{thm:Cd-odd-two-isogeny-reembedding}.

In characteristic two, every smooth \(C_d\), \(\mathcal T_{a,d}\), and
\(\mathcal C_{a,b,d}\) Jacobian considered here is ordinary, just as the
binary Edwards comparator is ordinary.  On this correctly delimited branch,
the reduced two-torsion quotient preserves \(a\) or \(a+b\), updates only
\(d\mapsto\sqrt d\), and costs one multiplication on native \(C_d\) or
one multiplication/one square on \(\mathcal T_{a,d}\).  The full proof,
equal-endpoint binary Edwards dictionary, exact cost table, reference
benchmark, and exhaustive dual-composition checks are all contained in this
chapter.

The extension families make this platform broader.  The one-sided model
is the native \(C\)-curve form of twisted Edwards, the full product family
retains both boundary square classes, reciprocal \(C\)-curves expose a
Montgomery Kummer chart, and the QRT envelope transports an isogeny
coordinatewise on adjacent states.  Together these structures allow
current CSIDH/CTIDH, large-prime-degree, oriented-action, and elliptic-factor
workflows to choose the most efficient internal coordinate while keeping
an explicit \(C_d\) or extended-\(C\) endpoint.

Consequently, the \(C_d\) platform is more than an Edwards transcription.
It combines a native monic-linear \(u\)-Kummer backend for low and medium
odd degrees, a zero-loss \(w\)-backend for square-root V\'elu, an optimized
degree-three full-point circuit, a root-free QRT-Jacobian degree-two
closure in odd characteristic, and one-operation native degree-two
circuits in characteristic two.  At equal output level it preserves the
best transported core whenever the common quotient forces equality, while
the displayed schedules improve the Edwards baselines for degree-three
full points, rooted odd-characteristic degree-two endpoints, and native
binary endpoints.  Direct target-parameter recovery and the availability
of \(\mathcal T_{a,d}\), \(\mathcal C_{a,b,d}\), reciprocal, and QRT
endpoints make these gains usable without suppressing field-of-definition
conditions.  This combination of exact dictionaries and genuinely native
specializations is the central arithmetic advantage of the extended
\(C_d\) family for supersingular-isogeny implementations in odd
characteristic and for its ordinary binary two-isogeny branch.

\chapter{Conclusion}
\label{ch:conclusion}

\section{Synthesis}

The equation
\[
        (u^2+u)(v^2+v)=d
\]
supports a coherent theory from projective geometry to explicit
cryptographic operations.  Characteristic uniformity means that the equation, marked
identity, boundary, inverse, native quotient, and declared input/output
interface survive unchanged, while each characteristic uses its optimized
arithmetic realization.  The invariant skeleton is
\[
\begin{gathered}
\text{smooth \((2,2)\)-model}
 \longrightarrow
\text{intrinsic \(D_8\) square symmetry}\\[-1mm]
\Downarrow\\[-1mm]
\text{rational four-torsion}
 \longrightarrow
\text{inverse fixing \(u\)}
 \longrightarrow
\kappa=(u+1:u).
\end{gathered}
\]
The square symmetry is exact rather than decorative.  In characteristic
different from two, translating the affine origin to
\((-1/2,-1/2)\) turns the eight natural transformations into all signed
permutations of two coordinates, namely the four rotations and four
reflections of a square.  In every characteristic, including characteristic
two where no common affine center exists, the full action has fixed field
\[
 k(\mathcal C_d)^{D_8}
   =k\bigl((u^2+u)+(v^2+v)\bigr).
\]
This full square quotient is distinct from the Kummer quotient by inversion;
the former has generic degree eight and the latter degree two.

In odd characteristic, centered variables expose Edwards symmetry and give a
Kummer circuit with the Montgomery dependency graph.  In characteristic two, the same factors
become Artin--Schreier operators and lead to \(Z/4\mathbb Z\)-normal and
binary Kummer arithmetic.  Characteristic three belongs to the odd branch
for its group law and differential addition, but its tripling map factors as
relative Frobenius followed by Verschiebung.  This factorization is kept
separate from a separable three-isogeny with a chosen reduced kernel; in
particular, inexpensive cubing is not identified with scalar tripling or
with a V\'elu quotient.  These classical dictionaries explain and optimize
the formulas, while the marked \(\mathcal C_d\) completion defines the model.
Each branch begins with addition and doubling on \(\mathcal C_d\), and each optimized operation returns either
a full \(\mathcal C_d\)-point or its declared native quotient.

The reciprocal equation
\[
                  (u+1)(v+1)=du^2v^2
\]
is a complementary internal chart of the same marked model.  Its
completion is carried to \(\overline{\mathcal C}_d\) by an ambient
involution of \(\PP^1\times\PP^1\), and its affine function \(v+1\)
represents exactly the native Kummer point \((u_0+1:u_0)\).  Consequently a
point already stored in reciprocal coordinates supplies an affine known
difference without an inversion, whereas the original native input retains
the cheaper odd-characteristic first double.  This endpoint-sensitive
comparison exhibits two complementary interfaces of the same Cd model, each
with its own stored-input and normalization profile.

The completeness results make the distinction between a chart and a model
precise.  No formula with an everywhere-finite affine output can handle
\(P+(-P)=O=(0,\infty)\).  On the smooth \((2,2)\)-completion, however, the
native centered Segre coordinates
\[
       (X:Y:T:Z)=(RV_0:SU_0:U_0V_0:RS)
\]
support
\[
                    P_1+P_2=(EF:GH:EH:FG).
\]
This formula costs \(9\M+\Dpar\), its mixed specialization costs
\(8\M+\Dpar\), and its dedicated double costs \(4\M+4\Sqr\).
Over an odd finite field this single tuple is \(\F_q\)-complete exactly when
\(\rho=1-16d\) is a nonsquare.  For the remaining parameters, and in
characteristic two, finitely many native laws have empty common base locus.
The ordinary differential formulas and their exceptional-difference charts
likewise form complete native differential-addition atlases.  Thus
\(\mathcal C_d\) has a complete arithmetic system of its own.  Here the word
``complete'' refers to the smooth model or to an atlas; a finite affine chart
necessarily omits the boundary point \(O\).

Therefore \(\mathcal C_d\) is an explicit biquadratic arithmetic model with
its own smooth completion, marked boundary, function-field presentation,
addition-law spaces, and Kummer normalization.  The isomorphisms with
Edwards, Montgomery, and Weierstrass equations provide coordinate
dictionaries for proofs and optimization, while the native input and output
objects remain
\[
P_n\in\overline{\mathcal C}_d
\qquad\text{or}\qquad
\kappa_d([n]P)=(u_n+1:u_n).
\]
When \(P_n\) lies in the finite affine chart it is represented by
\((u_n,v_n)\); boundary outputs remain on the native Segre completion.
The operation counts stated in this monograph therefore refer to explicitly
specified native or quotient outputs.

For the specialization
\[
d_{\rm C}=(16\cdot121666)^{-1}
\qquad\text{over}\qquad \F_{2^{255}-19},
\]
the native quotient coordinate is exactly the X25519 \(U\)-coordinate.
The distinguished square root of \(-1\) gives an invertible diagonal
recoding of the native Segre completion, with complete addition costing
\(8\M+\Dpar\), mixed addition \(7\M+\Dpar\), and the native-input first
doubling \(\M+\Sqr+\Dpar\).  At the quotient level, the X25519 input space
covers the main Curve25519 curve and its quadratic twist; a full Cd25519 point
refines this quotient input by supplying the data of a main-curve lift.  The
RFC decoding, clamping, cofactor, and all-zero-output conditions combine with
the algebraic coordinate dictionary to give the complete protocol interface.

For odd-degree separable isogenies, the same native quotient admits a closed
product evaluator rather than merely an imported Montgomery routine.  If
\(m=(\ell-1)/2\), the projective \(u\)-Kummer output costs
\(2m\M+2\Sqr\) for affine kernel coordinates, and the full affine-chart
point map is obtained from two products and their derivatives.  The target
parameter is recovered as
\[
       \rho'=\frac{\rho^\ell}
       {\prod_{i=1}^{m}(2u_i+1)^8},
       \qquad d'=\frac{1-\rho'}{16}.
\]
Thus odd-degree model closure, kernel enumeration, Kummer evaluation, and
image-parameter recovery all remain on the \(\mathcal C_d\) interface.
The \(w\)-coordinate and square-root V\'elu constructions transfer without
loss through explicit rational dictionaries, while
\(\mathcal T_{a,d}\), reciprocal \(C\)-curves, and adjacent QRT states
extend the same backend to the corresponding twisted or stateful interfaces.
The characteristic-two branch has its own Artin--Schreier product formula;
the odd-characteristic supersingular backend and the binary ordinary backend
are thereby stated with their correct, distinct hypotheses.

The three-parameter family
\[
\mathcal C_{a,b,d}:\qquad
(u^2+u+a)(v^2+v+b)=d
\]
separates the properties of general factorized biquadratic genus-one curves
from those specific to \(\mathcal C_d\).  Every smooth member is a
\((2,2)\)-curve of genus one and retains the two factor involutions, although
a rational point is required before the curve itself becomes an elliptic
curve over the ground field.  In odd characteristic the centered reciprocal
form is diagonal Edwards, whereas in characteristic two the parameters
\(a\) and \(b\) contribute Artin--Schreier twist data.  The loci
\[
b=0,\qquad a=b,\qquad d=ab
\]
give, respectively, a rational Edwards boundary, coordinate-exchange
symmetry, and four explicit rational affine points.  In particular, the
specialization \(a=b=0\) combines characteristic-uniform geometry, rational
marked four-torsion, factor-exchange symmetry, and a native arithmetic
interface on the same \((2,2)\)-model.

The extension part separates three further structures.  The one-sided family
\(\mathcal T_{a,d}\) separates the moduli parameter from the quadratic or
Artin--Schreier twist.  Reciprocal \(\mathcal C\)-curves replace the split affine
reflections by reciprocal Möbius involutions while retaining explicit Jacobi,
Weierstrass, and Kummer interfaces.  The symmetric QRT family adds an
intrinsic displacement class to the genus-one curve.  After an origin is
chosen, this class is represented by the translation point
\[
D=T(O),
\]
while for a genuinely pointed elliptic curve \((E,O,D)\) the adjacent-state
map
\[
P\longmapsto
\bigl(\kappa(P),\kappa(P+D)\bigr)
\]
realizes the same translation on a symmetric biquadratic state curve.
These three levels---an unpointed genus-one automorphism, its
origin-dependent translation point, and a pointed adjacent-Kummer
state---are kept distinct in the arithmetic developed in the extension
chapters.

\section{Principal conclusions}

The main results can be summarized as nineteen conclusions.
\begin{enumerate}[label=\textup{(\arabic*)}]
 \item The smooth completion and its marked boundary determine the identity,
 inverse, rational four-torsion, exceptional divisors, and Kummer
 normalization.  These data constitute the geometric core of the arithmetic
 model and govern every native formula.

 \item The elementary symmetries
 \[
 (u,v)\mapsto(-u-1,v),\qquad
 (u,v)\mapsto(u,-v-1),\qquad
 (u,v)\mapsto(v,u)
 \]
 extend across the boundary and generate a faithful intrinsic
 \(D_8\)-subgroup.  The resulting eight coordinate maps are exactly the
 full rigid-motion group of a square after centring at
 \((-1/2,-1/2)\) in odd characteristic; the center itself belongs only to
 the singular member \(d=1/16\).  The quarter-turn cycles the marked
 boundary as \(O\mapsto -R\mapsto T\mapsto R\mapsto O\), and the invariant
 \(J=(u^2+u)+(v^2+v)\) generates the full fixed field in every
 characteristic.  In odd characteristic the maps are also
 translation--inversion combinations involving \(T\) and \(R\).  This
 symmetry is present for every smooth parameter, independently of the
 larger automorphism groups at CM parameters.

 \item Full-point arithmetic is intrinsic.  Direct affine formulas provide
 the basic group law; the centered Segre law provides a low-cost complete
 formula on the nonsquare-\(\rho\) subfamily; and a finite native atlas is
 complete in every characteristic.  Mixed addition and dedicated doubling
 are stated with operation counts separate from that of the general
 complete law.

 \item Quotient arithmetic is also intrinsic.  Since inversion fixes \(u\),
 the map \(\kappa_d=(u+1:u)\) supports direct doubling, differential
 addition, recovery, and complete differential atlases.  Its native input
 relation \(X-Z=1\) gives the first-double cost
 \(\M+\Sqr+\Dpar\); the book explicitly distinguishes this one-time saving
 from the recurring ladder-step cost.

 \item The reciprocal presentation
 \[
                  \mathcal R_d:(u+1)(v+1)=du^2v^2
 \]
 is an ambient chart of the same marked \((2,2)\)-completion.  It conjugates
 the four marked points and the full \(D_8\)-action explicitly, and its
 affine coordinate \((v+1:1)\) equals \(\kappa_d\) projectively.  The chart
 exposes the odd Montgomery and reduced binary Weierstrass abscissae,
 supplies inversion-free affine known differences when used as the stored
 input interface, and admits general, mixed, fully affine, and complete-atlas
 differential formulas.  Exact endpoint accounting shows at the same time
 that converting a raw native point incurs the usual normalization cost and
 that the native odd first double remains strictly cheaper.

 \item Characteristic uniformity is geometric rather than circuit-level.
 Centered reciprocal formulas are natural in odd characteristic,
 Frobenius--Verschiebung tripling is natural in characteristic three, and
 Artin--Schreier formulas are natural in characteristic two.  In the
 characteristic-three branch, \([3]\) has degree nine, relative Frobenius
 has degree three and is purely inseparable, and a point-sum V\'elu formula
 applies only to a separately specified reduced cyclic kernel.  All three
 branches retain the same equation and quotient interface.

 \item The native coordinate links local arithmetic with global division
 theory.  Besides the pulled-back Weierstrass and Edwards division
 polynomials, the paired recursions
 \(\mathbf K_n=(\mathsf X_n:\mathsf Z_n)\) construct
 \(\kappa_d([n]P)\) directly from native \texttt{xDBL} and
 \texttt{xADD}.  After
 primitive cancellation the pair has degree \(n^2\); its identity coordinate
 gives the native division polynomial.  Thus one addition-chain circuit
 serves scalar multiplication, torsion testing, and division-polynomial
 evaluation.

 \item The finite-field statistics are exact in both characteristic
 branches.  In addition to the geometric and rational isomorphism-class
 counts and the first trace moment, the native character sums give
 \[
 \sum_{d\in\mathcal D_q}a_d^2
 =\begin{cases}
 q^2-2q-3,&q\text{ odd},\\
 q^2-q-1,&q\text{ even}.
 \end{cases}
 \]
 The same Frobenius traces determine the zeta function and all
 extension-field point counts by a second-order recurrence.

 \item Division fibers satisfy exact trace identities.  For every
 \(n\) prime to the characteristic,
 \[
  \frac1{n^2}\sum_{[n]P=Q}\frac{u(P)+1}{u(P)}
     =\frac{u(Q)+1}{u(Q)}.
 \]
 In odd characteristic the reciprocal centered coordinates have means
 \(\xi(Q)/n\) and
 \((-1)^{(n-1)/2}\eta(Q)/n\) for odd \(n\), and zero for even \(n\).
 These statements include characteristic three when \(3\nmid n\) and are
 interpreted as function-field trace identities at pole fibers.

 \item The same native functions organize model-preserving isogenies and
 pairings.  Odd and binary V\'elu-type formulas return to a member of the
 \(\mathcal C_d\) family; kernel-polynomial evaluation, Galois descent,
 relative Frobenius--Verschiebung, composition, and dual isogenies remain on
 that family.  The ordinary and supersingular characteristic-three cases
 are separated by their kernel schemes and inseparable degrees.  Native
 Miller line numerators, Kummer-dependent vertical denominators, an exact
 denominator-elimination criterion, product pairings, and Tate and Weil
 wrappers keep the pairing computation on the \((u,v)\)-interface.

 \item The supersingular-isogeny specialization is constructive at every
 declared interface.  For odd \(\ell=2m+1\), the quotient remains in
 \(\mathcal C_{d'}\), kernel pairs are enumerated natively as
 \((u_i+1:u_i)\), affine-kernel Kummer evaluation costs
 \(2m\M+2\Sqr\), and
 \[
       \rho'=\frac{\rho^\ell}{\prod_i(2u_i+1)^8},
       \qquad d'=\frac{1-\rho'}{16}.
 \]
 The dedicated degree-three circuits improve the displayed Edwards
 full-point baseline by \(2\Sqr+\Cmul\).  In odd characteristic, the
 root-free degree-two Kummer and full-point cores match the corresponding
 Edwards cores, while the native
 \(\mathcal T_{a_2,d_2}\) endpoint saves one multiplication by a curve
 constant relative to the displayed rooted Edwards endpoint.  In
 characteristic two, the native \(C_d\) and \(\mathcal T_{a,d}\)
 degree-two circuits cost \(\M\) and \(\min(\M,\Sqr)\), respectively;
 persistent common \(W\)-Kummer states retain the identical shared core.
 The \(w\)-coordinate and square-root V\'elu backends transfer without
 loss at higher degree.  The complete symbolic cost table, controlled
 reference benchmarks, and exhaustive finite-field checks agree with these
 formula-level conclusions under their explicitly stated output and
 field-of-definition conditions.

 \item Parameter selection and encoding retain the native boundary.  The
 odd \(c=\alpha_{24}\) parameter has an exact smoothness and completeness
 test; affine decoding uses a square-and-sign criterion in odd
 characteristic and a trace-and-linear-functional criterion in
 characteristic two.  Four reserved tags encode the four boundary points,
 and the canonical decoder invariant separates algebraic correctness,
 subgroup policy, and constant-time failure handling.

 \item Special parameters add efficient non-scalar endomorphisms without
 changing that interface.  The \(j=1728\) point \(d=1/8\) has
 \(\phi^2=[-1]\), while the \(j=0\) locus \(A^2=3\) has
 \(\psi^2+\psi+[1]=[0]\).  Their sign-change and affine-projective Kummer
 actions yield the eigenvalue equations required for GLV decomposition.

 \item Cd25519 realizes the Curve25519/Edwards25519 isomorphism class
 within the native family.  Its quotient is the X25519 \(U\)-line; its
 \(\imath\)-twisted native coordinates support complete
 \(8\M+\Dpar\) addition, \(7\M+\Dpar\) mixed addition, and
 \(4\M+4\Sqr\) doubling; and its \(U=9\) fixed-difference ladder step
 costs \(4\M+4\Sqr+\Dpar+\Dnine\).  The leading-bit optimization and the
 RFC vector are verified, while curve/twist and protocol boundaries are
 kept explicit.
 \item The symmetric generalization
 \[
       \mathcal C_{a,b,d}:(u^2+u+a)(v^2+v+b)=d
 \]
 has an exact characteristic-dependent classification.  In odd
 characteristic its smoothness is
 \((1-4a)(1-4b)(16d)
 ((1-4a)(1-4b)-16d)\ne0\), its reciprocal equation is diagonal Edwards,
 and its base-field Jacobian is explicitly Weierstrass.  In characteristic
 two smoothness is simply \(d\ne0\), and the Jacobian is
 \(Y^2+XY=X^3+(a+b)X^2+d^2X\).  The cases \(b=0\), \(a=b\), and
 \(d=ab\) isolate, respectively, a ground-field twisted Edwards chart,
 extra \(D_8\)-symmetry, and a rationally pointed Huff--Weierstrass chart.

 \item The one-sided family \(\mathcal T_{a,d}\) represents the Jacobian
 classes of the full product family and, over finite fields, the curves
 themselves.  In odd characteristic it is a twisted-Edwards/Montgomery
 normal form with a distinguished native \(A=-1\) slice at \(a=1/2\); in
 characteristic two it gives every ordinary elliptic curve over a perfect
 field and separates \(j=d^{-4}\) from the Artin--Schreier class of \(a\).

 \item Reciprocal \(C\)-curves have an exact projective smoothness theory,
 reciprocal quotient involutions, a visible order-four action, Jacobi and
 Weierstrass models, and complete full-point and Kummer arithmetic atlases.
 Over a perfect field of characteristic two,
 \(\mathcal R_{\tau,\sigma,\kappa}\) is M\"obius-equivalent to
 \(\mathcal C_{\sqrt{\tau\sigma}/\kappa}\), so the apparent reciprocal
 involutions become Artin--Schreier translations.

 \item The symmetric QRT envelope has exact odd- and binary-characteristic
 smoothness conditions, explicit \(j\)-invariants and Weierstrass models,
 finite-field model and abstract-isomorphism counts, and model-preserving
 two-, three-, and odd-prime-degree isogeny criteria.  Isogenies act
 coordinatewise on nondegenerate adjacent states through their induced
 Kummer maps.  Generic rational formulas are kept distinct from complete
 atlases, degenerate image states, and square-root normalizations that
 require field-of-definition hypotheses.

 \item The McMillan map is the composition of Vieta root exchange and
 coordinate swap and acts as a fixed translation on every smooth QRT fiber.
 Before an origin is chosen, its invariant is a class in the Jacobian; after
 choosing \(O\), the translation point is \(D=T(O)\).  For a genuinely
 pointed curve \((E,D)\), the state
 \((\kappa(P),\kappa(P+D))\) supports explicit state doubling,
 differential addition, recovery, and binary or signed-digit ladder
 recurrences satisfying \(\Delta_D T_D=T_D^2\Delta_D\).  Exact evenization
 is equivalent to a split rational two-torsion Kummer involution and is
 available after an extension of degree at most six.  On the split locus,
 each state branch is a one-formula complete PGL$_2$ conjugate of compiled
 Montgomery \texttt{xDBLADD}, at \(4\M+4\Sqr+2\Dpar\); the ordinary-binary
 pointed state model is universal over perfect fields, retains its
 Artin--Schreier twist parameter, and costs
 \(4\M+5\Sqr+3\Dpar\).  The three-chart recovery atlas and the
 affine-index theorem distinguish fixed-base use from generic MSM.  The
 resulting elliptic Lucas calculus includes nonlinear addition and
 subtraction, state-division polynomials, a nodal specialization to the
 classical Lucas pair, and exact bridges to Ward EDS, sigma functions, and
 elliptic nets.
\end{enumerate}

Within the coordinate systems and cost models developed in this monograph,
the arithmetic of \(\mathcal C_d\) attains the optimized Edwards
field-operation counts on the corresponding linearly equivalent loci.
The \(C_d\) model simultaneously retains a
characteristic-uniform defining equation, an intrinsic \(D_8\)-action, four
rational boundary points, a native Kummer quotient, the same defining
equation in characteristic \(2\), and explicit extensions to twisted,
reciprocal, and symmetric-biquadratic QRT families.

\section{Final perspective}

The results above separate three levels that are often conflated in explicit
elliptic-curve arithmetic: the abstract elliptic curve, the chosen arithmetic
model, and its quotient representations.  For \(\mathcal C_d\), the smooth
marked \((2,2)\)-completion fixes the identity, rational four-torsion,
inverse, and native Kummer quotient before any auxiliary coordinate system is
introduced.  Edwards, Montgomery, Weierstrass, and related normal forms then
serve as explicit dictionaries for proofs and optimization.  This separation
makes the domain, completeness, and output type of each arithmetic circuit
part of the statement of the circuit itself.

The same marked geometry persists across characteristics, while the optimal
realization of its arithmetic changes with the field.  In odd characteristic
it yields Edwards-compatible full-point arithmetic and the native Montgomery
Kummer line; in characteristic two the quadratic factors become
Artin--Schreier operators and the same equation supports binary full-point and
Kummer arithmetic; in characteristic three the multiplication-by-three map
exhibits its Frobenius--Verschiebung factorization.  The reciprocal, twisted,
and symmetric-QRT extensions show how these structures change when the marked
branch data or the distinguished translation is varied.
For the QRT extension, the cryptographic conclusions are exact at the state
level: descent is controlled by split two-torsion, the pointed group law is
distinguished from the Jacobi choice of identity, the split odd and
ordinary-binary ladder branches have one complete formula, and the latter
arise for every ordinary pointed elliptic curve over a perfect binary field.
Recovery has two explicit boundary charts, and isogenies act coordinatewise
until the marked state degenerates.  The exact interfaces and cost ledgers
established here provide the basis for platform-matched performance
comparisons of canonical-QRT circuits.

Cd25519 provides a concrete specialization in which the native Kummer
coordinate is exactly the X25519 coordinate while complete full-point
arithmetic remains available on the native completion.  More generally,
within the coordinate systems and cost models developed in this monograph,
\(\mathcal C_d\) arithmetic attains the optimized Edwards dependency graphs
on the corresponding loci while retaining a characteristic-uniform defining equation, an
intrinsic \(D_8\)-action, four rational boundary points, native Kummer
coordinates, and direct extensions to twisted, reciprocal, and symmetric
biquadratic QRT models.  These properties establish \(\mathcal C_d\) as a
characteristic-uniform arithmetic model in its own right and provide a common
setting for the geometric, arithmetic, and dynamical constructions developed
throughout the monograph.

\appendix

\chapter{Algorithmic Templates}
\label{app:algorithms}
\section[Odd-characteristic native Cd Kummer ladder]
{Odd-characteristic native \(\mathcal C_d\) Kummer ladder}

\begin{enumerate}[label=\arabic*.,leftmargin=2.2em]
 \item Input \(P=(u_P,v_P)\in\mathcal C_d\), a scalar
 \(m=(m_{\ell-1}\cdots m_0)_2\), and
 \(\alpha_{24}=1/(16d)\).
 \item Form the native quotient without inversion,
 \[
       \Delta=\kappa_d(P)=(u_P+1:u_P).
 \]
 Set \(R_0=(1:0)\), \(R_1=\Delta\).  Normalize \(\Delta\) only when an
 affine known-difference implementation is desired.
 \item For \(i=\ell-1,\ldots,0\), conditionally swap \(R_0,R_1\);
 evaluate one \(x\)ADD with difference \(\Delta\) and one \(x\)DBL;
 conditionally swap back.
 \item If \(R_0=(X_m:Z_m)\), return the native quotient
 \(\kappa_d([m]P)=R_0\), with
 \[
       u([m]P)=Z_m/(X_m-Z_m).
 \]
 Retain \(R_1\) and recover \(v([m]P)\) when a full
 \(\mathcal C_d\)-point is required.
\end{enumerate}

\section{Cd25519 leading-bit Kummer ladder}

For a clamped X25519 scalar
\(m=2^{254}+\sum_{i=3}^{253}m_i2^i\), the following template removes the
otherwise redundant processing of the fixed leading bit while retaining a
uniform loop for all secret bits.
\begin{enumerate}[label=\arabic*.,leftmargin=2.2em]
 \item Decode the quotient input as \(\Delta=(X_P:Z_P)\).  For a genuine
 Cd25519 affine input \((u_P,v_P)\), use the inversion-free lift
 \(\Delta=(u_P+1:u_P)\).  For the standard base point take
 \(\Delta=(9:1)\).
 \item Set \(R_0=\Delta\).  Compute \(R_1=2\Delta\) once with
 \eqref{eq:Cd25519-xdbl}; when the source is an affine Cd25519 point, use
 the native first-double circuit \eqref{eq:Cd25519-first-double}, costing
 \(\M+\Sqr+\Dpar\).
 \item For \(i=253,\ldots,0\), perform the usual constant-time conditional
 swap, one \(x\)ADD with known difference \(\Delta\), one \(x\)DBL, and
 the closing swap.  The invariant is
 \(R_1-R_0=\Delta\), and after bit \(i\) the pair represents the two
 adjacent multiples selected by the processed prefix.
 \item Return \(R_0\).  If a canonical Cd25519 affine quotient is required,
 perform one final inversion and recover
 \[
       u([m]P)=\frac{Z_m}{X_m-Z_m}.
 \]
 A full point additionally requires the sign information discarded by the
 quotient; it cannot be reconstructed from an arbitrary X25519 byte string
 alone.
\end{enumerate}
For the fixed difference \((9:1)\), each of the (254) loop iterations
costs \(4\M+4\Sqr+\Dpar+\Dnine\).  This count excludes decoding, the final
inversion, and constant-time swap or memory costs.

\section{CM-assisted joint scalar multiplication}

\begin{enumerate}[label=\arabic*.,leftmargin=2.2em]
 \item Select one of the two parameter loci of
 Chapter~\ref{ch:native-endomorphisms}: either
 \(d=1/8\) with \(i^2=-1\), or \(A^2=3\) with
 \(\zeta^2+\zeta+1=0\).  Verify that the chosen prime-order subgroup
 \(G=\langle P\rangle\) is stable under \(\vartheta=\phi\) or \(\psi\).
 \item Determine the subgroup eigenvalue
 \[
 \vartheta(P)=[\lambda]P,
 \]
 with \(\lambda^2+1\equiv0\pmod n\) in the \(j=1728\) case or
 \(\lambda^2+\lambda+1\equiv0\pmod n\) in the \(j=0\) case.
 \item Reduce the lattice
 \[
 \{(a,b)\in\mathbb Z^2:a+b\lambda\equiv0\pmod n\}
 \]
 once for the fixed subgroup.  Decompose each input scalar as
 \(m\equiv m_1+m_2\lambda\pmod n\) with
 \(|m_1|,|m_2|=O(\sqrt n)\).
 \item For a full-point output, compute the full point \(\vartheta(P)\)
 by \eqref{eq:j1728-native-endomorphism} or
 \eqref{eq:j0-native-endomorphism}.  For quotient-only arithmetic, the
 sparse Kummer maps \eqref{eq:j1728-Kummer-action} or
 \eqref{eq:j0-Kummer-action} may be used only as part of a differential
 joint-multiplication procedure supplied with the required oriented
 differences and recovery data; two independent Kummer images do not
 determine the sign of their sum.
 \item In the full-point path, evaluate
 \[
              [m_1]P+[m_2]\vartheta(P)
 \]
 with a constant-time joint signed recoding and either a complete native
 full-point atlas or a single complete law whose parameter hypothesis has
 been verified.  In particular, do not infer completeness of the
 nonsquare-\(\rho\) Segre tuple merely from the existence of the CM map.
 Return the result in \(u,v\), not in an auxiliary CM model.  In the
 quotient-only path, return only the declared Kummer class unless a
 separately specified full-coordinate recovery step is available.
\end{enumerate}

\section{Characteristic-two native ladder}

\begin{enumerate}[label=\arabic*.,leftmargin=2.2em]
 \item Input \(P=(u_P,v_P)\in\mathcal C_d\), form
 \(K_P=(u_P+1:u_P)\), and set \(d_{\rm K}=d^{-1}\).
 \item Normalize to \(K_P=(1:t)\) only if the chosen circuit requires it;
 this costs \(\Inv\) and additions from a raw native input, using
 \(t=1+(u_P+1)^{-1}\).  Otherwise retain the
 inversion-free projective pair.  Initialize \(R_0=(1:0)\), \(R_1=K_P\).
 \item At every scalar bit, use constant-time swaps, binary \(x\)DBL
 from \eqref{eq:binary-xdbl-thesis}, and either the general projective
 law \eqref{eq:binary-xadd-thesis} for an unnormalized difference or
 \eqref{eq:binary-xadd-a-thesis}--\eqref{eq:binary-xadd-b-thesis} after
 affine normalization.
 \item Return \(R_0=(X_0:X_1)\), interpreted natively as
 \[
       u([m]P)=X_1/(X_0+X_1).
 \]
 Retain the adjacent state for
 Proposition~\ref{prop:binary-recovery-thesis} and convert the recovered
 Weierstrass point back to \(u,v\).
\end{enumerate}

\section{Characteristic-two native separable two-isogeny}

\begin{enumerate}[label=\arabic*.,leftmargin=2.2em]
 \item Input a perfect binary field, \(d\ne0\), the precomputed inverse
       Frobenius value \(e=\sqrt d\), and either a \(C_d\) point or a
       \(\mathcal T_{a,d}\) point.
 \item For \(C_d\), compute
       \[
                    t=u(u+1).
       \]
       Return the target Kummer pair \((e:t)\).  For a full native output,
       return
       \[
             (u'_0:u'_1)=(t:t+e),\qquad v'=u+v.
       \]
       This step costs one multiplication.
 \item For \(\mathcal T_{a,d}\), choose one of the two public
       field-backend schedules:
       \[
       F=u^2+u+a,\quad
       (v'_0:v'_1)=(F:F+e),
       \]
       costing one square, or
       \[
       h=v(v+1),\quad
       (v'_0:v'_1)=(e:e+h),
       \]
       costing one multiplication.  In both cases return
       \(u'=u+v+1\).
 \item Keep the outputs projective through a chain and normalize by batch
       inversion only at the declared endpoint.  Boundary pairs are valid
       points of the smooth \(\mathbb P^1\times\mathbb P^1\) completion
       and must not be sent to an affine division.
 \item If the dual is required, apply
       \[
       (u',v')\mapsto({u'}^2,{v'}^2)
       \quad\text{for }C_d,
       \qquad
       (u',v')\mapsto({u'}^2+a,{v'}^2)
       \quad\text{for }\mathcal T_{a,d}.
       \]
\end{enumerate}

\section{Ternary Miller update}

For a signed ternary digit, maintain
\(R=[m]P=(u_R,v_R)\in\mathcal C_d\), its double, and the
accumulator \(F_m(Q)\).  The tripling update is
\[
F\leftarrow F^3
g_{R,R}^{-}(Q)
g_{2R,R}^{-}(Q),\qquad
R\leftarrow3R.
\]
The point update may use the full \(\mathcal C_d\)-formula
\eqref{eq:char3-full-tripling-native}; the two factors are evaluated by
\eqref{eq:odd-native-Miller-thesis}, specialized to characteristic three.
A Kummer \(x\)TPL alone is insufficient because the two line
functions require ordinates and signs.

\chapter{Notation and Formula Index}
\label{app:notation-index}
\begin{longtable}{L{3.2cm}L{10.5cm}}
\toprule
symbol & meaning\\
\midrule
\endfirsthead
\toprule
symbol & meaning\\
\midrule
\endhead
\(\Cd\) & \((u^2+u)(v^2+v)=d\), with smooth completion understood\\
\(\iota_u,\iota_v,\sigma\) & the three intrinsic involutions
\((-u-1,v)\), \((u,-v-1)\), and \((v,u)\), generating \(D_8\)\\
\(\rho\) & \(1-16d\) for \(\mathcal C_d\); in
the odd-characteristic part of Chapter~\ref{ch:generalized-Cabd},
\(\rho=\lambda\mu-\delta\)\\
\(r,s\) & \(2u+1,2v+1\), the odd-characteristic polynomial shifts of
the native coordinates\\
\(\beta\) & \((16d)^{-1}\), also \(\alpha_{24}\), in odd
characteristic\\
\(A\) & \((4d)^{-1}-2\), the odd-characteristic Montgomery coefficient\\
\(U\) & \((u+1)/u\), the affine native Kummer coordinate\\
\(\kappa\) & \((u+1:u)\), the projective native Kummer map\\
\(\omega\) & \(\rho/((2u+1)^2(2v+1)^2)\), the degree-eight
post-isogeny Kummer coordinate in odd characteristic\\
\(e_\omega\) & \(4/\rho\), the auxiliary constant for the
post-isogeny Kummer formulas\\
\(T_4\) & the Edwards point \((1,0)\), corresponding to the marked
native boundary point \(R\), used in the order-four quotient\\
\(i,\imath,\zeta\) & \(i^2=\imath^2=-1\) in their declared local
fields, and \(\zeta^2+\zeta+1=0\); \(i,\zeta\) are the CM constants,
whereas \(\imath\) is also used for the Cd25519 diagonal recoding and
temporary scalar extensions\\
\(d_{\rm K}\) & \(d^{-1}\), the characteristic-two Kummer constant\\
\(N_d,a_d\) & \(\#\mathcal C_d(\F_q)\) and its Frobenius trace
\(a_d=q+1-N_d\)\\
\(k_{\rm emb}\) & extension or embedding degree used in the pairing chapter\\
\(h_{\rm red}\) & \((q^{k_{\rm emb}}-1)/n\), the reduced-Tate final exponent\\
\(\M,\Sqr,\Dpar,\Cmul,\Dnine,\Inv\) & multiplication, squaring,
fixed parameter multiplication in the general monograph cost model,
curve-constant multiplication in the elliptic Lucas schedules,
multiplication by \(9\), and inversion\\
\(\gamma_{\rm S},\gamma_{\rm D}\) & normalized costs
\(\Sqr/\M\) and \(\Dpar/\M\)\\
\(\gamma_{2,{\rm S}},\gamma_{2,{\rm C}}\) & binary normalized costs
\(\Sqr/\M\) and \(\mCurve/\M\) when \(\mBase=\M\)\\
\(\lambda,\mu,\delta\) & as the odd-characteristic chapter-level
parameters in
Chapter~\ref{ch:generalized-Cabd},
\(1-4a,1-4b,16d\), respectively\\
\(\mathcal T_{a,d}\) & the one-sided twisted family
\((u^2+u+a)(v^2+v)=d\)\\
\(\mathcal R_{\tau,\sigma,\kappa}\) & the reciprocal \(C\)-curve
\((x^2-\tau)(y^2-\sigma)=\kappa xy\), distinct from the internal
one-parameter reciprocal chart \(\mathcal R_d\)\\
\(\mathcal Q_{\alpha,\beta,\gamma}\) & the symmetric QRT biquadratic
\(x^2y^2+\alpha(x^2+y^2)+\beta xy+\gamma=0\)\\
\(q_Q,\Omega_Q\) & \(\beta^2-4\alpha^2-4\gamma\) and
\(q_Q^2-64\alpha^2\gamma\), the odd-characteristic QRT invariants\\
\(T_Q\) & the McMillan/QRT translation
\((x,y)\mapsto(y,-x-\beta y/(y^2+\alpha))\)\\
\(\delta_T\) & the intrinsic displacement class
\([T(P)-P]\in\operatorname{Pic}^0(C)\), independent of \(P\)\\
\(\mathcal S_D,\Delta_D\) & adjacent Kummer state
\(P\mapsto(\kappa(P),\kappa(P+D))\) and its induced state-doubling map\\
\(v_n,w_n\) & oriented elliptic Lucas coordinate and its Jacobi companion,
representing the marked point \([n]D\) in Chapter~\ref{ch:symmetric-QRT-envelope}\\
\(W_n,R_n,\Upsilon_n\) & Ward EDS term and the even/odd EDS ratios used to
recover the Weierstrass and oriented Lucas coordinates\\
\(\Sigma_K,\tau_K\) & characteristic-two half-kernel sum
\(\sum_i a_i\) and its unique fourth root\\
\(W_{A,d}\) & characteristic-two sparse ordinary equation
\(y^2+xy=x^3+Ax^2+d^2x\); its reduced two-torsion quotient preserves
\(A\) and sends \(d\) to \(\sqrt d\)\\
\(A,B,J,c_K\) & in Chapter~\ref{ch:Cd-isogeny-cryptography}, the two
native odd-isogeny kernel products, their product--derivative ordinate
factor, and the marked normalization constant; these local \(A,B\) are
distinct from the Montgomery coefficient above\\
\(w_C,\eta_C\) & the \(\mathcal C_d\)-native pullbacks of the Edwards
\(w\)-coordinate and its shifted square-root-V\'elu coordinate\\
\(\mBase,\mCurve\) & binary base-point and curve-constant multiplications\\
\(\Cube\) & characteristic-three Frobenius cubing\\
\(\Full,\Ffin,\Ktwo,\Kfour\) & completed, finite-chart, and
quotient-output types\\
\eqref{eq:native-affine-add-thesis} &
basic odd-characteristic addition on \(\mathcal C_d\)\\
\eqref{eq:native-affine-add-uv} &
the same addition written only in \(u_i,v_i,d\)\\
\eqref{eq:native-affine-double-uv} &
odd-characteristic doubling written only in \(u,v,d\)\\
\eqref{eq:ordinary-native-projective-add} &
ordinary native \((U:V:Z)\) projective addition\\
\eqref{eq:ordinary-native-projective-double} &
ordinary native \((U:V:Z)\) dedicated doubling\\
\eqref{eq:xdbl-thesis} &
odd-characteristic native Kummer doubling; first input costs
\(\M+\Sqr+\Dpar\)\\
\eqref{eq:binary-native-affine-double} &
basic characteristic-two doubling on \(\mathcal C_d\)\\
\eqref{eq:native-full-tripling-uv} &
full native tripling on \(\mathcal C_d\)\\
\eqref{eq:odd-full-projective-tripling-thesis} &
dedicated native full-point projective tripling\\
\eqref{eq:native-closed-2PQ} &
closed native \(2P+Q\) on \(\mathcal C_d\)\\
\eqref{eq:j1728-native-endomorphism} &
fourth-order \(j=1728\) endomorphism on \(\mathcal C_{1/8}\)\\
\eqref{eq:j1728-Kummer-action} &
\(j=1728\) Kummer action \((X:Z)\mapsto(-X:Z)\)\\
\eqref{eq:j0-native-endomorphism} &
third-order \(j=0\) endomorphism on the locus \(A^2=3\)\\
\eqref{eq:j0-Kummer-action} &
\(j=0\) Kummer action \((X:Z)\mapsto(\zeta X+cZ:Z)\)\\
\eqref{eq:odd-native-Miller-thesis} &
odd-characteristic native Miller factor in \(u,v\)\\
\eqref{eq:binary-native-Miller-thesis} &
characteristic-two native Miller factor in \(u,v\)\\
\eqref{eq:native-Miller-wrapper} &
raw and locally normalized Miller evaluations
\(\mathsf M_d^{\rm raw}(n;P,Q)\), \(\mathsf M_{d,t}(n;P,Q)\)\\
\eqref{eq:Tate-shifted-Cd} &
normalization-free reduced Tate pairing wrapper on \(\mathcal C_d\)\\
\eqref{eq:Weil-shifted-Cd} &
normalization-free Weil pairing wrapper on \(\mathcal C_d\)\\
\eqref{eq:T-char2-Weierstrass} &
characteristic-two Weierstrass form of \(\mathcal T_{a,d}\)\\
\eqref{eq:Cd-char2-two-isogeny-native-general} &
native characteristic-two separable two-isogeny on \(C_d\)\\
\eqref{eq:T-char2-two-isogeny-native} &
twist-stable native separable two-isogeny on \(\mathcal T_{a,d}\)\\
\eqref{eq:reciprocal-general} &
general reciprocal \(C\)-curve\\
\eqref{eq:R-Jacobi-quartic} &
normalized reciprocal Jacobi quartic\\
\eqref{eq:symmetric-QRT} &
symmetric QRT biquadratic\\
\eqref{eq:QRT-Jacobian} &
odd-characteristic Weierstrass model of the QRT envelope\\
\eqref{eq:QRT-McMillan-sum-form} &
Vieta--McMillan fixed-translation map\\
\eqref{eq:QRT-state-doubling-definition} &
adjacent-Kummer-state doubling \(\Delta_D\)\\
\eqref{eq:Cd-native-odd-xeval} &
odd-degree native \(u\)-Kummer isogeny product evaluator\\
\eqref{eq:Cd-image-parameter-affine} &
odd-degree image parameter
\(\rho'=\rho^\ell/\prod_i(2u_i+1)^8\)\\
\eqref{eq:Cd-full-product-derivative} &
native full-point product--derivative isogeny formula on its main chart\\
\eqref{eq:Cd-binary-Kummer-product} &
characteristic-two native Artin--Schreier Kummer isogeny product\\
\eqref{eq:Cd-crypto-binary-two-isogeny-Kummer} &
general characteristic-two degree-two Kummer quotient
\((e(K_0+K_1)^2:K_0K_1)\)\\
\eqref{eq:Cd-crypto-w}--\eqref{eq:Cd-crypto-eta} &
exact \(\mathcal C_d\) pullbacks of the Edwards \(w\)- and
\(\eta\)-coordinates\\
\bottomrule
\end{longtable}

\section*{Acknowledgments}
The author acknowledges the use of ChatGPT (OpenAI) as an interactive
assistant during the preparation of this monograph.  It was used to explore
alternative presentations, check symbolic computations, improve the
exposition, and edit the English and \LaTeX{} source.

The mathematical contributions, algorithms, proofs, and conclusions were
conceived, verified, and finalized by the author, who assumes full
responsibility for the accuracy of the work.

The author hopes that this monograph will serve as a useful reference for
researchers, educators, and students interested in elliptic curve theory and
its applications.

\backmatter

\phantomsection

\end{document}